\documentclass[10pt,reqno]{amsart}
\usepackage{kpfonts}
\usepackage{marvosym}
\usepackage{amsfonts,latexsym,amsmath,amssymb,mathrsfs,verbatim,cancel}
\usepackage[alphabetic]{amsrefs}
\usepackage{multirow}
\usepackage{url,color}
\usepackage{makecell}
\usepackage{pifont}
\usepackage{upgreek}
\usepackage{fancyhdr}
\usepackage{hyperref}
\usepackage[percent]{overpic}
\usepackage{caption}
\usepackage{subcaption}
\usepackage{pict2e}
\usepackage[hypcap=false]{caption}
\usepackage{xcolor}
\usepackage{soul}
\usepackage{wasysym}
\usepackage{scalerel}
\usepackage{accents}
\usepackage{slashed}
\usepackage{enumitem}
\usepackage{tabularx}
\usepackage{float}
\setlist[enumerate]{label*=\arabic*.}
\usepackage{array}
\newcolumntype{C}[1]{>{\centering\arraybackslash}m{#1}} 
\usepackage{multirow}
\allowdisplaybreaks

\topmargin - .4 in

\newtheorem{theorem}{Theorem}[section]
\newtheorem{proposition}[theorem]{Proposition}
\newtheorem{lemma}[theorem]{Lemma}
\newtheorem{corollary}[theorem]{Corollary}

\theoremstyle{definition}
\newtheorem{definition}[theorem]{Definition}
\newtheorem{remark}[theorem]{Remark}
\newtheorem{assumption}[theorem]{Assumption}
\newtheorem{notation}{Notation}[section]

\newcommand*\eqdef{\overset{\mbox{\tiny{def}}}{=}}
\newcommand{\R}{\mathbb{R}}
\newcommand{\T}{\mathbb{T}}
\newcommand{\p}{\partial}
\newcommand{\mr}{\mathring} 
\newcommand{\rmD}{\mathrm{d}}
\newcommand{\rmd}{\mathrm{d}}

\DeclareMathOperator{\Span}{span}
\newcommand{\mytr}{{\mbox{\upshape tr}}} 
\newcommand{\bm}{\boldsymbol}
\newcommand{\otimesarray}{\overset{\diamond}{\otimes}}

\newcommand{\tander}{\mathscr{P}}
\newcommand{\topordertancom}{\mathscr{O}}

\newcommand{\belowtopordercomforpartialmodifiedchi}{\widetilde{\mathscr{B}}}

\newcommand{\tandersmall}{\mathscr{P}_*}
\newcommand{\comder}{\mathscr{Z}}
\newcommand{\comdersmall}{\mathscr{Z}_*}
\newcommand{\comderdoublesmall}{\mathscr{Z}_{**}}
\newcommand{\tanderY}{\mathscr{Y}}
\newcommand{\tandergeneral}{P}

\DeclareMathOperator{\diag}{diag}

\DeclareMathOperator*{\esssup}{ess\,sup}
\DeclareMathOperator*{\essinf}{ess\,inf}

\DeclareMathOperator{\mydet}{det}

\newcommand{\RRiemannPS}{\mathcal{R}_{(+)}^{PS}}
\newcommand{\LRiemannPS}{\mathcal{R}_{(-)}^{PS}}

\newcommand{\muxmulevelsetvalue}{\mathfrak{n}}

\newcommand{\initialsmall}{\mathring{\upepsilon}}
\newcommand{\fundbootsmall}{\varepsilonup}
\newcommand{\auxbootsmall}{\varepsilonup^{1/2}}
\newcommand{\initialsmalldoublenull}{\underline{\initialsmall}}

\newcommand{\ubar}{\underline{u}}
\newcommand{\ingoingmu}{\underline{\upmu}}
\newcommand{\gnulltori}{\breve{\gtorus}}
\newcommand{\nulltorusproject}{\breve{\torusproject}}
\newcommand{\nullgeop}[1]{\breve{\partial}_{#1}}

\newcommand{\outgoingcharacteristicsurfacetwoarg}[2]{\mathcal{P}_{#1}^{#2}}
\newcommand{\ingoingcharacteristicsurfacearg}[1]{\underline{\mathcal{P}}_{#1}}
\newcommand{\ingoingcharacteristicsurfacetwoarg}[2]{\underline{\mathcal{P}}_{#1}^{#2}}
\newcommand{\doublenulltoritwoarg}[2]{\underline{\ell}_{#1,#2}}
\newcommand{\characteristicdiamondtwoarg}[2]{\mathcal{M}_{#1,#2}}

\newcommand{\blowupratetoporderacoustic}{K_*}
\newcommand{\blowupratetoporderwave}{M_*}
\newcommand{\blowuprateofwaveWRTnumberofcommutations}{\mathfrak{p}}
\newcommand{\blowuprateofacousticgeoWRTnumberofcommutations}{\mathfrak{r}}
\newcommand{\blowuprateoftransportWRTnumberofcommutations}{\mathfrak{q}}
\newcommand{\blowuprateoftopordertransportWRTnumberofcommutations}{\mathfrak{s}}

\newcommand{\Lunit}{L}
\newcommand{\uLunit}{\underline{L}}
\newcommand{\Lgeo}{L_{(Geo)}}
\newcommand{\uLgeo}{\underline{L}_{(Geo)}}
\newcommand{\newL}{\breve{L}}
\newcommand{\newuL}{\breve{\underline{L}}}
\newcommand{\olduLunit}{\underline{\mathscr{L}}}

\newcommand{\ReciprocalLunitAppliedtoTimeFunction}{\iota}
\newcommand{\ReciprocaluLunitAppliedtoTimeFunction}{\underline{\iota}}
\newcommand{\MagnitueofinnerproductofnewLandnewuL}{\uplambda}
\newcommand{\MagnitueofinnerproductofLunitanduLunit}{\Lambda}

\newcommand{\ToriTangentVectorfieldAssociatedToDoubleNullFolliations}{U}
\newcommand{\gtorusdoublenullCOV}{\underline{\mathbf{M}}}
\newcommand{\angularcoeffinmuweightedspatialcartesian}{\mathbf{Z}}
\newcommand{\nullangDiv}{\breve{\slashed{\Flatdiv}}}
\newcommand{\nullangD}{\breve{\angD}}
\newcommand{\nullangpartial}{\breve{\slashed{\partial}}}
\newcommand{\nullangcurl}{\breve{\Flatcurl}}
\newcommand{\nullangrmD}{\breve{\rmD} \mkern-8mu / \,}
\newcommand{\nullangChristoff}{\breve{\Upgamma}}
\newcommand{\Riemnulltori}{\breve{\mathfrak{Riem}}}
\newcommand{\Ricnulltori}{\breve{\mathfrak{Ric}}}
\newcommand{\Scalarnulltori}{\breve{\mathfrak{R}}}
\newcommand{\Gaussnulltori}{\breve{\mathfrak{K}}}
\newcommand{\UnitTimeNormalizedUbarDataSurfaceTangentVectorfield}{\mathring{\mathscr{V}}}
\newcommand{\coefficientsInVectorfieldsTangetToUbarDataSurface}{\varkappa}
\newcommand{\datasurfacegeometricp}[1]{\mathring{\partial}_{#1}} 
\newcommand{\datasurfacep}[1]{\mathring{\mathfrak{d}}_{#1}}

\newcommand{\volingoingnullhypersurface}{\mathrm d \underline{\varpi}}
\newcommand{\voloutgoingnullhypersurface}{\mathrm d \overline{\varpi}}
\newcommand{\voldoublenulltori}{\mathrm d \varpi_{\gnulltori}}
\newcommand{\voldiamond}{\mathrm d \boldsymbol{\varpi}}
\newcommand{\volcanonical}{\mathrm{d} \text{vol}}

\newcommand{\ingoingfluxwave}[1]{\underline{\mathbf{F}}_{#1}}
\newcommand{\outgoingfluxwave}[1]{\mathbf{F}_{#1}}
\newcommand{\ingoingfluxtransport}[1]{\underline{\mathbf{F}}_{(Transport; #1)}}
\newcommand{\outgoingfluxtransport}[1]{\mathbf{F}_{(Transport; #1)}}

\newcommand{\spacetimecoercive}[1]{\mathbf{K}_{#1}}
\newcommand{\newspacetimecoercive}[1]{\underline{\mathbf{K}}_{#1}}

\newcommand{\Lsmall}{L_{(\textnormal{Small})}}
\newcommand{\Xsmall}{X_{(\textnormal{Small})}}
 
\newcommand{\Yvf}[1]{Y_{(#1)}}
\newcommand{\Yvfsmall}[1]{Y_{(#1;Small)}}
\newcommand{\Yvfsmallcoeff}[1]{y_{(#1)}}

\newcommand{\Lrough}{\widetilde{L}}

\newcommand{\controlvars}{\upgamma}
\newcommand{\badcontrolvars}{\underline{\upgamma}}

\newcommand{\wavearray}{\vec{\Psi}}

\newcommand{\wavearraypartial}{{\vec{\Psi}}_{\textnormal{(Partial)}}}
\newcommand{\velocityarray}{\vec{V}}
\newcommand{\velocityarraypartial}{\vec{V}_{\textnormal{(Small)}}}

\newcommand{\Wtransarg}[1]{{ ^{(#1)} \mkern-2mu \breve{W}}}

\newcommand{\gtorusroughfirstfund}{\widetilde{\gtorus}}
\newcommand{\gtorusroughinversefirstfund}{\widetilde{\gtorus} \,^{-1}}

\newcommand{\timefunction}{\uptau}
\newcommand{\timefunctionarg}[1]{{^{(#1)} \mkern-1mu \uptau}}

\newcommand{\timefunctionboot}{\timefunction_{\textnormal{Boot}}}

\newcommand{\mupositive}{\mulevelsetvalue_0}
\newcommand{\boringregionmupositive}{\mulevelsetvalue_1}
\newcommand{\awayfromsmallneightborhoodnmupositive}{\upmu_2}
\newcommand{\Ntop}{N_{\textnormal{top}}}

\newcommand{\HarmlessWave}[1]{\textnormal{Harmless}_{(\textnormal{Wave})}^{#1}}

\newcommand{\zetatan}{\upzeta^{(\text{Tan--}\wavearray)}} 
\newcommand{\zetatrans}{\upzeta^{(\text{Trans--}\wavearray)}} 
\newcommand{\angktan}{\angk^{(\text{Tan--}\wavearray)}} 
\newcommand{\angktrans}{\angk^{(\text{Trans--}\wavearray)}}

\newcommand{\totalfluxcontrolvelocity}{\mathbf{V}} 

\newcommand{\strongangularcontrolvelocity}{\mathbf{K}\slashed{\mathbf{V}}}

\newcommand{\totalcontrolwave}{\mathbf{W}}

\newcommand{\totalfluxcontrolprecisefullymodifiedchi}{\mathbf{X}^{(\textnormal{Precise})}}
\newcommand{\totalfluxcontrolimprecisefullymodifiedchi}{\mathbf{X}^{(\textnormal{Imprecise})}}
\newcommand{\totalfluxcontrolpartialmodifiedchi}{\widetilde{\mathbf{X}}}
\newcommand{\totalbulkcontrolpartialmodifiedchi}{\mathbf{K}\widetilde{\mathbf{X}}}
\newcommand{\totalbulkcontrolprecisefullymodifiedchi}{\mathbf{K}\mathbf{X}^{(\textnormal{Precise})}}
\newcommand{\totalbulkcontrolimprecisefullymodifiedchi}{\mathbf{K}\mathbf{X}^{(\textnormal{Imprecise})}}
\newcommand{\weakspacetimetransvelocitycontrol}{\underline{\mathbf{K}\breve{\mathbf{V}}}}

\newcommand{\weakspacetimeangularvelocitycontrol}{\underline{\mathbf{K}\slashed{\mathbf{V}}}}

\newcommand{\fluxcontrolVort}{\mathbf{O}}
\newcommand{\bulkcontrolVort}{\mathbf{KO}}
\newcommand{\fluxcontrolGradEnt}{\mathbf{S}}
\newcommand{\bulkcontrolGradEnt}{\mathbf{KS}}
\newcommand{\fluxcontrolVortVort}{\mathbf{C}}
\newcommand{\bulkcontrolVortVort}{\mathbf{KC}}
\newcommand{\fluxcontrolDivGradEnt}{\mathbf{D}}
\newcommand{\bulkcontrolDivGradEnt}{\mathbf{KD}}
\newcommand{\toricontrolVort}{\,^{(\textnormal{Null Tori})}\mathbf{O}}
\newcommand{\toricontrolGradEnt}{\,^{(\textnormal{Null Tori})}\mathbf{S}}
\newcommand{\toricontrolVortVort}{\,^{(\textnormal{Null Tori})}\mathbf{C}}
\newcommand{\toricontrolDivGradEnt}{\,^{(\textnormal{Null Tori})}\mathbf{D}}
\newcommand{\mastercontroltop}{\mathbf{M}^{(\textnormal{Top})}}

\newcommand{\fullymodquant}[1]{ ^{(#1)}\mathscr{X}}
\newcommand{\partialmodquant}[1]{ ^{(#1)}\widetilde{\mathscr{X}}}

\newcommand{\partialmodquantinhom}[1]{ ^{(#1)}\widetilde{\mathfrak{X}}}

\newcommand{\hypthreearg}[3]{{^{(#3)}\widetilde{\Sigma}_{#1}^{#2}}}
\newcommand{\hyptwoarg}[2]{{^{(#2)}}\widetilde{\Sigma}_{#1}}

\newcommand{\twoargMrough}[2]{{^{(#2)}\mathcal{M}_{#1}}}

\newcommand{\roughtori}[1]{\widetilde{\ell}_{#1}}
\newcommand{\twoargroughtori}[2]{{^{(#2)}\widetilde{\ell}_{#1}}}

\newcommand{\singhyp}{\mathcal B}
\newcommand{\crease}{\partial_- \mathcal B}

\newcommand{\Cauchyhor}{\underline{\mathcal{C}}}

\newcommand{\COVframe}{\Lambda}
\newcommand{\COVL}{\lambda}

\newcommand{\geop}[1]{\frac{\partial}{\partial #1}}
\newcommand{\twogeop}[2]{\frac{\partial^2}{\partial{#1}\partial{#2}}}
\newcommand{\roughgeop}[1]{\frac{\widetilde{\partial}}{\widetilde{\partial} #1}}

\newcommand{\Flatdiv}{\mbox{\upshape div}\mkern 1mu}
\newcommand{\Flatcurl}{\mbox{\upshape curl}\mkern 1mu}

\newcommand{\RRiemann}{\mathcal{R}_{(+)}}
\newcommand{\LRiemann}{\mathcal{R}_{(-)}}
\newcommand{\almostRiemann}{\mathcal{R}}

\newcommand{\Ent}{s}
\newcommand{\GradEnt}{S}

\newcommand{\Density}{\varrho}

\newcommand{\LogDensity}{\uprho}

\newcommand{\vort}{\omega}
\newcommand{\vortrenormalized}{\Omega}

\newcommand{\Speed}{c}

\newcommand{\VortVort}{\mathcal{C}}
\newcommand{\DivGradEnt}{\mathcal{D}}

\DeclareMathOperator{\ErrorIBP}{\mbox{\upshape Error}}

\newcommand{\Errortop}{\mbox{\upshape Error}_{N}^{(\textnormal{Top})}} 

\newcommand{\Errortoparg}[1]{\mbox{\upshape Error}_{#1}^{(\mbox{\upshape \footnotesize Top})}} 
\newcommand{\Errorsubcriticalarg}[1]{\mbox{\upshape Error}^{(\mbox{\upshape \footnotesize Sub-critical})}_{#1}}

\newcommand{\enmomem}{\mathbf{Q}}

\newcommand{\Ricfour}{\mathbf{Ric}}
\newcommand{\Riemfour}{\mathbf{Riem}}

\newcommand{\Transport}{\mathbf{B}}

\newcommand{\muX}{\breve X}

\newcommand{\Jen}[1]{^{(#1)} \mathbf{J}}
\newcommand{\Jenarg}[2]{{^{(#1)}  \mathbf{J}^{#2}}}

\newcommand{\deform}[1]{{^{(#1)} \mkern-1mu \pmb{\pi}}}
\newcommand{\deformarg}[3]{{^{(#1)} \pmb{\pi}_{#2 #3}}}

\newcommand{\angdeform}[1]{{^{(#1)} \mkern-1mu \pi \mkern-8mu /}}
\newcommand{\angdeformarg}[2]{\angdeform{#1}_{#2}}

\newcommand{\gfour}{\mathbf{g}}
\newcommand{\hfour}{\mathbf{h}}

\newcommand{\gtorus}{g \mkern-8mu / }

\newcommand{\gtorusdoublearg}[2]{\gtorus_{#1#2}}

\newcommand{\Chfour}{\pmb{\Gamma}}

\newcommand{\torusproject}{{\Pi \mkern-12mu / } \,\, }

\newcommand{\Sigmatproject}{\Pi}

\newcommand{\smoothtorusproject}{{\Pi \mkern-12mu / } \,\, }

\newcommand{\Dfour}{\mathbf{D}}

\newcommand{\angD}{ {\nabla \mkern-10mu / \, \,} }

\newcommand{\angDarg}[1]{{\angD_{\mkern-3mu #1}}}

\newcommand{\angDsquared}{ {\angD_{\mkern-3mu}}^{\, 2} }

\newcommand{\angdiv}{\mbox{\upshape{div} $\mkern-17mu /$\,\,}}

\newcommand{\angLap}{ {\Delta \mkern-9mu / \, } }
\newcommand{\angrmD}{\rmD \mkern-8mu / \,}
\newcommand{\argangrmd}[1]{\angrmD_{#1}}
\newcommand{\angrmd}{{\rmd \mkern-8mu / \, }}

\newcommand{\newangD}{ {\nabla \mkern-10mu / \,} }
\newcommand{\newangDarg}[1]{ {\nabla \mkern-10mu / \,}_{#1} }

\newcommand{\angk}{ { {k \mkern-7mu /} \, } }

\newcommand{\angG}{{\vec{G} \mkern -10mu /} \,\,}

\newcommand{\weight}{\mathscr{W}}

\newcommand{\Lie}{\mathcal{L}}
\newcommand{\SigmatLie}{\underline{\mathcal{L}}}
\newcommand{\angLie}{ { \mathcal{L} \mkern-8mu / } }

\newcommand{\smoothfunction}{\mathrm{f}}

\newcommand{\coercivequadraticformintoriboundaryterms}{\mathfrak{P}}

\newcommand{\Currentboundaryerrorhavetocontrolprincipal}{\mathfrak{E}_{(\textnormal{Principal})}}
\newcommand{\Currentboundaryerrorhavetocontrollowerorder}{\mathfrak{E}_{(\textnormal{Lower-order})}}
\newcommand{\Currentboundaryerrorhavetocontrolprincipalweighted}[1]{\mathfrak{E}_{(\textnormal{Principal};#1)}}
\newcommand{\Currentboundaryerrorhavetocontrollowerorderweighted}[1]{\mathfrak{E}_{(\textnormal{Lower-order};#1)}}

\newcommand{\EllipticHyperbolicCurrentIntegralIdentityTotalSpacetimeErrorTerm}{\mathfrak{M}}

\newcommand{\Nullhypersurfaceproject}{\overline{\Pi}}

\newcommand{\SigmatTan}{W}

\newcommand{\ehcurrent}{\mathscr{J}}
\newcommand{\flatcurrent}{\mathscr{J}_{(\textnormal{Flat})}}

\newcommand{\ellipticCoerciveQuadratic}{\mathscr{Q}}

\newcommand{\Singletan}{P}

\newcommand{\MInteresting}{\mathcal{M}_{\textnormal{Interesting}}}

\newcommand{\smallneighborhoodofcreasearg}[1]{\mathcal{N}_{#1}}

\newcommand{\mulevelsetvalue}{\mathfrak{m}}

\newcommand{\sumofbrevexvisquared}{\breve{\mathfrak{V}}}

\newcommand{\almostRiemannfunction}{F}

\newcommand{\twoargmumuxtorus}[2]{\breve{\mathbf{T}}_{#1,#2}}
\newcommand{\argubarnewultorus}[1]{\underline{\breve{\mathbf{T}}}_{#1}}
\newcommand{\Cartesiantisafunctiononubarnewultori}{\underline{\mathfrak{T}}}
\newcommand{\Eikonalisafunctiononubarnewultori}{\underline{\mathfrak{U}}}

\newcommand{\Cartesiantisafunctiononlevelsetsofroughtimefunctionarg}[2]{\mathfrak{t}_{#1,#2}}
\newcommand{\Cartesiantisafunctiononlevelsetsofingoingcharacteristicsarg}[1]{\underline{\mathfrak{t}}_{#1}}

\newcommand{\Cartesiantisafunctiononmumxtoriarg}[2]{\mathfrak{T}_{#1,#2}}
\newcommand{\Eikonalisafunctiononmumuxtoriarg}[2]{\mathfrak{U}_{#1,#2}}

\newcommand{\embeddatahypersurfacearg}[1]{{^{(#1)}E}}

\newcommand{\embeddednewuLhypersurface}{\underline{E}}

\newcommand{\datahypfortimefunctionarg}[1]{\breve{\mathbb{X}}_{#1}}
\newcommand{\datahypfortimefunctiontwoarg}[2]{\breve{\mathbb{X}}_{#1}^{#2}}

\newcommand{\newuLmulevelset}{\underline{\breve{\mathbb{L}}}}

\newcommand{\mulevelsetarg}[1]{\breve{\mathbb{M}}_{#1}}
\newcommand{\mulevelsettwoarg}[2]{{\breve{\mathbb{M}}_{#1}^{#2}}}

\newcommand{\domainforembeddingdatahypfortimefunctiontwoarg}[2]{{^{(#1)} \mkern-2mu \mathscr{H}_{#2}}}

\newcommand{\embeddingdatahypfortimefunctionarg}[1]{{^{(#1)} \mkern-2mu H}}

\newcommand{\scalarembeddingdatahypfortimefunctionarg}[1]{{^{(#1)} \mkern-2mu h}}

\newcommand{\flowmapWtransargtwoarg}[2]{{^{(#1)} \mkern-2mu \iota_{#2}}}

\newcommand{\nullhyp}{\mathcal{P}}
\newcommand{\nullhyparg}[1]{\mathcal{P}_{#1}}

\newcommand{\nullhypthreearg}[3]{{^{(#1)} \mkern-2mu \mathcal{P}_{#2}^{#3}}}

\newcommand{\multipliervectorfield}{\breve{T}}

\newcommand{\Fullset}{\mathscr{Z}}
\newcommand{\Tanset}{\mathscr{P}}
\newcommand{\Angularset}{\mathscr{Y}}

\newcommand{\CHOVgeotorough}[1]{{^{(#1)} \mkern-2mu \mathscr{T}}}
\newcommand{\CHOVgeotodoublenull}{\underline{\mathscr{U}}}
\newcommand{\CHOVJacobiangeotodoublenull}{{^{(\CHOVgeotodoublenull)} \mkern-2mu \mathbf{J}}}
\newcommand{\CHOVdoublenulltoubarnewulmu}{\underline{\Phi}}

\newcommand{\CHOVroughtomumuxmu}[1]{{^{(#1)} \mkern-1mu \Phi}}
\newcommand{\InverseCHOVroughtomumuxmu}[1]{{^{(#1)} \mkern-1mu {\Phi^{-1}}}}
\newcommand{\CHOVJacobianroughtomumuxmu}[1]{{^{(\CHOVroughtomumuxmu{#1})} \mkern-2mu \mathbf{J}}}

\newcommand{\CHOVgeotomumuxmu}{\breve{\mathscr{M}}}

\newcommand{\CHOVJacobiangeotomumuxmu}{{^{(\CHOVgeotomumuxmu)} \mkern-2mu \mathbf{J}}}

\newcommand{\nullform}{\mathfrak{Q}}

\newcommand{\CH}{\underline{\mathcal{C}}}

\newcommand{\mainnullform}{\mathfrak{M}}

\newcommand{\FlowmapnewL}{\breve{\Lambda}}

\newcommand{\InverseFlowmapnewL}{\breve{\Lambda}^{-1}}
\newcommand{\FlowmapnewLarg}[1]{\breve{\Lambda}^{#1}}

\newcommand{\argLrough}[1]{{^{(#1)} \mkern-2mu \widetilde{L}}}

\newcommand{\mydiam}{{\mkern-1mu \scaleobj{.75}{\blacklozenge}}}

\newcommand{\leftubar}{\underline{U}_0}
\newcommand{\ubarboot}{\ubar_{\textnormal{Boot}}}
\newcommand{\ubarlocal}{\ubar_{\textnormal{Local}}}
\newcommand{\leftu}{U_2}
\newcommand{\moreinterestingu}{U^{\mbox{\small \Lightning}}}
\newcommand{\interestingu}{U_{\mbox{\tiny \Radioactivity}}}
\newcommand{\rightu}{U_1}

\newcommand{\blowupdelta}{\mr \updelta_*}
\newcommand{\blowupdeltadoublenull}{\underline{\blowupdelta}}

\newcommand{\secondtransversalderivativemulowerbound}{M_2}

\newcommand{\flowmapofnewularg}[1]{\underline{\Lambda}_{#1}}

\numberwithin{equation}{section}

\begin{document}
	\title{The emergence of the Cauchy horizon from the crease for $3D$ compressible Euler flow}
	\author[LA,JS]{Leonardo Abbrescia$^{* \dagger}$ and Jared Speck$^{** \dagger\dagger}$}

\thanks{$^{*}$Georgia Institute of Technology, Atlanta, GA, USA.
\texttt{abbrescia@math.gatech.edu}}

\thanks{$^{**}$Vanderbilt University, Nashville, TN, USA.
\texttt{jared.speck@vanderbilt.edu}}

\thanks{$^\dagger$ L. Abbrescia gratefully acknowledges support from a Travel Support for Mathematicians grant from the Simons Foundation.}

\thanks{$^{\dagger\dagger}$ J. Speck gratefully acknowledges support from NSF grants DMS-2349575 and DMS-2054184.
}

	\begin{abstract}
We study open sets of nearly plane symmetric initial data for the $3D$ compressible Euler equations with dynamic entropy and non-trivial vorticity. In our prior work \cite{abbrescia2022emergence}, we proved that the solutions develop a gradient singularity along a \emph{singular boundary}, 
which is a hypersurface that emanates from a co-dimension 2, acoustically spacelike submanifold called the \emph{crease}. In the present paper, we prove that a null hypersurface, called a \emph{Cauchy horizon}, also emanates from the crease, and propagates in a direction transversal to the singular boundary. The solution remains smooth up to the Cauchy horizon, even though it ``feels'' the influence of the gradient singularity at the crease. The union of the Cauchy horizon, the crease, and the singular boundary, make up a connected portion of the boundary of a \emph{maximal globally hyperbolic development} (MGHD) of the data. One can roughly view an MGHD as a ``largest'' way that a smooth data set can evolve into a classical solution. Known examples show that the only way one can guarantee uniqueness of an MGHD is by proving that certain structural properties hold along its \emph{entire} boundary. We prove that the needed properties are satisfied along the portion we construct in this paper and \cite{abbrescia2022emergence}.  

Our work provides the first description of the formation, structure, and stability of an $O(1)$-size portion of the Cauchy horizon in $3D$ without any simplifying assumptions. Our results hold for an arbitrary equation of state with positive speed of sound, aside from that of a Chaplygin gas, for which shocks are not expected to form. The presence of dynamic entropy and vorticity stretching, both of which are absent in simpler settings such as the $2D$ isentropic case, introduces significant difficulties that we resolve with novel techniques. Our approach relies on new foliations of spacetime, dynamically adapted to the shape of the crease, and a PDE framework based on double-null foliation and the remarkable energy identities we derived in \cite{abbrescia2025remarkable}. Collectively, our techniques allow us to handle a challenging new solution regime in which the crease lacks strict convexity.

		\bigskip
		
		\bigskip

	\end{abstract}
	
	\maketitle

	\tableofcontents

	\newpage

\section{Introduction} \label{S:INTRO}
We study the Cauchy problem for the $3D$ compressible Euler equations:
\begin{subequations} \label{E:INTRO:COMPRESSIBLEEULER}
		\begin{align}
		\Transport v^i 
		& = 
		- \frac{\partial_i p}{\varrho},
		&& (i=1,2,3),
		\label{E:INTROTRANSPORTVI}
			\\
		 \label{E:INTROTRANSPORTDENSITY}
		\Transport \varrho
		& = - \varrho \partial_a v^a,
		&&
			\\
		\Transport \Ent
		& = 0,
		&&
		\label{E:INTROBS}
		\end{align}
	\end{subequations}
where $v : \R \times  \Sigma \to \R^3$ is the velocity, $\varrho : \R\times\Sigma \to [0,\infty)$ is the density, $s : \R\times\Sigma \to \R$ denotes the entropy, and $\Transport$ is the \emph{material vectorfield}:
	\begin{align}
		\Transport 
		& \eqdef \partial_t + v^a \partial_a. \label{E:MATERIALDERIVATIVEVECOTRFIELD}
	\end{align}
We close the under-determined system \eqref{E:INTROTRANSPORTVI}--\eqref{E:INTROBS} by assuming that the pressure $p = p(\varrho,s)$ is given as a function of $\varrho$ and $s$, called the equation of state. Our results hold \emph{for any equation of state except for that of the Chaplygin gas} -- for which shocks are not expected to form -- with positive sound speed $\Speed$ defined by: 
\begin{align}
	\Speed \eqdef \sqrt{p_{;\Density}}, \label{E:SOUNDSPEED}
\end{align}
where $p_{;\Density}$ is the partial derivative of the pressure with respect to the density at fixed entropy. In this paper, $\Sigma \eqdef \R\times \T^2$ denotes the ``spatial manifold.'' 
Throughout, we denote the standard Cartesian coordinates on $\R\times \Sigma$ by
$(t,x^1,x^2,x^3)$, where $t = x^0$ denotes the Cartesian time function, $x^1$ is the Cartesian coordinate on $\R$, and $(x^2,x^3)$ denote Cartesian coordinates on $\T^2 \eqdef [-\pi,\pi]^2$ with the endpoints identified. Our main results, which are described in detail below, can be roughly summarized as follows:

\begin{quote}
We provide a rigorous description of the dynamic formation, geometry, and stability of an $\mathcal{O}(1)$-portion of 
a \emph{Cauchy horizon} $\Cauchyhor$; see Fig.\,\ref{F:INTROPICTURESOFMAINRESULTS}. 
We prove that such a Cauchy horizon develops in classical solutions launched by open sets 
of small, smooth, nearly plane-symmetric initial data (exact symmetry is not assumed) that lead to finite-time gradient-blowup (i.e., shock formation).
The data are allowed to have non-trivial entropy and vorticity. 
Roughly, $\Cauchyhor$ is a manifold-with-boundary that lies in the domain of influence 
of the ``initial points'' where the blowup occurs. Even though $\Cauchyhor$ is a boundary of the classical solution region, 
the solution remains smooth up to it.
In \cite{abbrescia2022emergence}, for the same initial data,
we constructed a related manifold-with-boundary, the singular boundary $\singhyp$, and proved that
the fluid variables' gradients blow up on $\singhyp$.
The union $\Cauchyhor \cup \singhyp$ delineates a portion of the boundary of a maximal region in which one can uniquely solve for  
a classical solution to \eqref{E:INTROTRANSPORTVI}--\eqref{E:INTROBS}. 

Our initial data are assumed to satisfy a ``transversal convexity'' property. Roughly speaking, 
it is a stable form of convexity \underline{only} in the approximately plane-symmetric direction; even in plane symmetry, 
transversal convexity is needed to prevent a continuum of characteristics from collapsing into a single point. 
In particular, we do not assume any convexity in the ``symmetry breaking directions.''
Other approaches to studying similar problems relied on a full convexity assumption on the initial data, which does not hold not for data with symmetry or their perturbations. In particular, \emph{this is the first stability result that derives the totality of the structures of $\Cauchyhor\cup\singhyp$ in an open neighborhood of plane-symmetry}. 
We develop many techniques to overcome the substantial difficulties introduced by the lack of full convexity, the
vorticity stretching term, and the presence dynamic entropy, which have not previously been handled in the context of Cauchy horizon formation.
Our methods are not constrained to perturbations of plane-symmetry or even the compressible Euler equations (e.g., they are applicable to 
other solution regimes and to quasilinear wave equations that fail to satisfy the null condition). 
In addition, our work has implications for classical determinism and classical uniqueness for
the fluid equations \eqref{E:INTROTRANSPORTVI}--\eqref{E:INTROBS}. Despite being a revered set of equations,
even the most basic question of \underline{global} uniqueness of classical solutions 
arising from smooth data is very far from being understood in full generality.
\end{quote}

\subsection{Outline of the remainder of the introduction} \label{SS:OUTLINEFORREMAINDEROFINTRO}
The rest of this introduction is devoted to fleshing out the main results and methods stated above, as well as placing them in the broader context of the existing literature for the compressible Euler equations. In Sect.\,\ref{SS:BROADPERSPECTIVEOFOURWORK}, we give a rough overview of our work, focusing on its connection to singularity formation and \emph{maximal globally hyperbolic developments} (MGHDs). 
In Sect.\,\ref{SS:CLASSICALDETERMINISM}, we briefly explain how our work makes progress on understanding the (still wide open) problem of global uniqueness of classical solutions (classical determinism). In Sect.\,\ref{SS:ALESSROUGHOVERVIEWOFOURWORK}, we give a less rough overview of our work, and highlight some of the main new ideas in the paper paper. In Sect.\,\ref{SS:CAUCHYPROBLEM}, 
we set up the Cauchy problem for the compressible Euler equations. 
In Sect.\,\ref{SS:GEOMETRICWAVETRANSPORTINTRO}, we reformulate the equations as a geometric wave-transport system, one that is well-suited to our PDE analysis. In Sect.\,\ref{SS:MGHDS}, we further discuss MGHDs in the context of the compressible Euler equations. In Sect.\,\ref{SS:RESULTSFROMSINGULARBOUNDARYPAPER}, we briefly discuss some of the main ideas and techniques used in \cite{abbrescia2022emergence} that are necessary for the current paper. In Sect.\,\ref{SS:NEWIDEASANDMETHODS}, we discuss in great detail the new ideas and methods developed in this paper. 
In Sect.\,\ref{SS:ABBREVIATEDMAINRESULTS}, we provide an abbreviated version of the main results.

\subsection{A rough overview of our work}  \label{SS:BROADPERSPECTIVEOFOURWORK}
In our previous paper \cite{abbrescia2022emergence}, we studied open sets of smooth initial data for the $3D$ compressible Euler equations with dynamic entropy and vorticity and gave a detailed description of the stable formation of shock singularities in a compact subset of spacetime. Within that compact subset, we gave a complete description of the co-dimension-one manifold-with-boundary, called the \emph{singular boundary}, where the fluid's gradient blows up. 
In the present paper, we study those classical solutions in a companion region, that is, up to a bounded portion of the \emph{Cauchy horizon}. The Cauchy horizon is a characteristic hypersurface that emerges from the gradient singularity and delineates a portion of the boundary of the region of classical existence. We prove that the fluid remains smooth up to the Cauchy horizon, even though it ``feels the influence'' of the gradient-blowup in its causal past. Taken together, the two papers provide a local-in-spacetime description of a portion of a \emph{maximal globally hyperbolic development of the data} (MGHD), including an $\mathcal{O}(1)$-size neighborhood of its boundary. Below, we will discuss MGHDs in detail. 
In the last several decades, there 
has been a lot of work on the study of gradient-blowup for multi-dimensional compressible fluids and related PDE systems, starting with foundational work by Alinhac and Christodoulou. What distinguishes our work here is the generality of our results, the generality of the data that our methods can treat, the largeness of the size of the region that our methods can handle, and the potential of the tools that we have developed for understanding the full structure of the MGHD.
In particular, the case of the $3D$ compressible Euler equations with vorticity and entropy is much more difficult to treat than the $2D$ isentropic case,
and our methods allow us to study singular sets that do not enjoy strict convexity. Such general methods are needed, for example, to study perturbations of symmetric solutions, where strict convexity does not hold.

Roughly, an MGHD is a largest possible classical solution + region that is uniquely determined by the initial data. The MGHD existence region can have a complicated shape, one that is not easily described by standard coordinates, such as the Cartesian ones. In \emph{any} spatial dimension, uniqueness of MGHDs for shock-forming data \emph{remains an open problem.} As we will explain below, the problem is fundamentally global, and cannot be answered by studying only compact portions of spacetime. Here and throughout, the word ``global'' does not mean ``for all time,'' but rather refers to the idea of the 
``totality of the classical solution.'' 
Although here we have treated only $\mathcal{O}(1)-$size compact regions of spacetime, 
the tools that we developed mesh well with tools that have proven to be exceptionally useful for studying the global structure of solutions, such as the geometric vectorfield method.
Before we describe our work in detail, we first place it in a general context.

\subsection{The breakdown of classical determinism -- a fundamental open problem}  \label{SS:CLASSICALDETERMINISM}
Classical determinism -- whether smooth initial conditions of a PDE determine \emph{unique} classical solutions -- is a fundamental notion of immense mathematical and physical significance. The breakdown of classical determinism, whether it be from singularity formation or a lack of uniqueness, has continuously driven mathematical progress for centuries. For the $3D$ compressible Euler equations governing gas dynamics, and more generally, quasilinear wave equations, determining whether any breakdown in classical determinism occurs is extraordinarily subtle. 
For a general quasilinear hyperbolic PDE, it is even difficult to rigorously \emph{formulate} the problem. A key difficulty is that proper formulations certainly require restrictions on the structure of the domains where uniqueness of classical solutions might hold, e.g., 
one cannot expect classical uniqueness on domains that are not in sufficient causal contact with the initial data. 
Another key difficulty, which we elaborate on 
later, is that the problem is global in nature,
i.e., there is a fundamental difference between ``locally unique'' classical solutions and ``globally unique'' ones. 
In Sect.\,\ref{SS:MGHDS}, for the $3D$ compressible Euler equations, we provide a precise formulation of the problem of classical determinism and discuss many of the subtleties. In this paper and its predecessor \cite{abbrescia2022emergence}, we solve many ``local'' aspects of the 
problem for open sets of initial data that lead to shock formation. We again highlight that it \emph{remains a fundamental outstanding open problem of the field to prove existence and uniqueness of a largest possible classical solution}.

The first mathematician to prove singularity formation in compressible fluid dynamics was Riemann in his historic 1860 work \cite{bR1860}. There, he proved that the \emph{classical evolution} of solutions launched by smooth initial data for the $1D$ compressible Euler equations can \emph{terminate} due to the formation of shocks. Extending this result to a constructive proof in multiple spatial dimensions without symmetry remained an open problem for 139 years until the celebrated results 
of Alinhac \cite{sA1999a,sA1999b}.\footnote{Prior to Alinhac, Sideris \cite{tS1985} had proved that singularities form for open sets of $3D$ compressible Euler solutions, but the proof was through a contradiction argument and non-constructive.} 
More precisely, Alinhac studied multi-dimensional quasilinear wave equations that essentially model irrotational and isentropic fluids.
Alinhac's data was non-degenerate in the sense that at the time of first blowup,
a gradient singularity occurred at an \emph{isolated} point. Moreover, his methods only allowed him to follow the solution up to the Cartesian time 
of first blowup; see Fig.\,\ref{F:MULTIDSHOCKFORMATION}. However, the fundamental principle of finite speed of propagation implies that the formation of a singularity at a given spacetime location does not preclude one from continuing the \emph{classical} evolution elsewhere to a larger region. The \emph{largest} region of spacetime on which a smooth solution's classical evolution is determined by the initial conditions is called the \emph{maximal globally hyperbolic development} (MGHD) of the data; see Sect.\,\ref{SS:MGHDS}.\footnote{MGHDs are sometimes loosely referred to as ``the  classical development'' of the data.} By definition, spacetime points on which the solution is singular are part of the \emph{boundary} of the MGHD. Therefore, describing the solution up to the boundary of the MGHD for shock-forming solutions, and proving the \emph{uniqueness} of such an MGHD, is a rigorous way to give a definitive answer to what causes the breakdown of the classical evolution. In 2007, Christodoulou published a breakthrough monograph \cite{dC2007} that developed methods to study shock formation for a much larger set of data than Alinhac's. In addition, the methods did not have the limitation of halting at the Cartesian time of first shock formation. In particular, in \cite{dC2007}*{Chapter 15}, Christodoulou derived a large, though implicit, portion of the MGHD's boundary and correctly predicted the structure and behavior of the solution on other portions of the boundary. To varying degrees of generality (described in more detail in Sect.\,\ref{SS:RESULTSFROMSINGULARBOUNDARYPAPER}), the latest progress \cite{abbrescia2022emergence,shkoller2024geometry} has yielded a rigorous and constructive proof of the behavior and state of the solution on \emph{portions} of the boundary components initially described in \cite{dC2007}*{Chapter 15}. However, the MGHD portions derived in \cite{abbrescia2022emergence,shkoller2024geometry} and in this paper are \textbf{\emph{extendible}}, and hence, as we explain in Sect.\,\ref{SS:MGHDS}, the uniqueness of the classical solution in the regions thus-far constructed is an open problem;
known proofs of uniqueness require one to construct the entire MGHD.  

\subsection{A less rough overview of our work and methods} \label{SS:ALESSROUGHOVERVIEWOFOURWORK}
\begin{figure}[ht]
\centering
\begin{subfigure}{.4\textwidth}
 \centering
\begin{overpic}[scale=.36,grid=false]{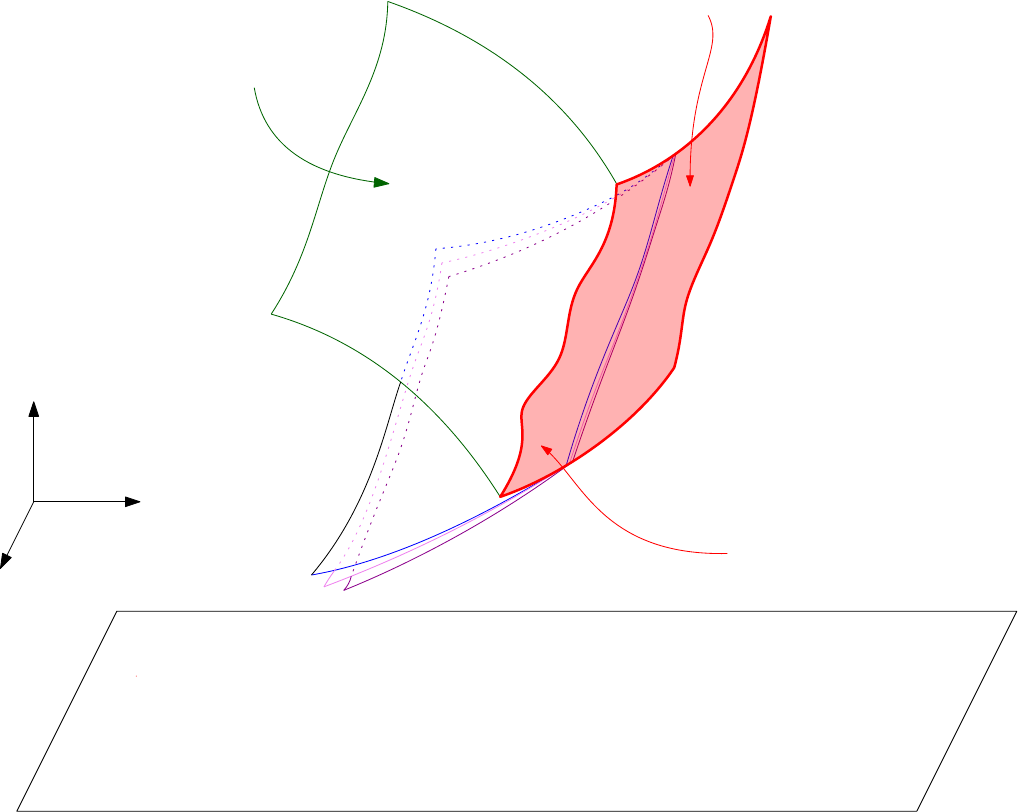} 
			\put (50,8) {$\Sigma_0$}
			\put (1.2,22) {$(x^2,x^3) \in \mathbb{T}^2$}
			\put (0,36) {$t$}
			\put (41,32) {$\nullhyparg{u}$}
			\put (14.5,29) {$x^1 \in \mathbb{R}$}
			\put (68,80) {$\mathcal{B}$}
			\put (72,24) {$\crease$}
			\put (23,74) {$\underline{\mathcal{C}}$}
\end{overpic}
  \caption{Infinite density of the characteristics $\nullhyparg{u}$ on $\mathcal{B}$, shaded in red, and constructed in \cite{abbrescia2022emergence}.}
  \label{F:INFININTEDENSITYOFCHARACTERISTICSONSINGULARBOUNDARY}
\end{subfigure}\, \, \, \,
\begin{subfigure}{.4\textwidth}
  \centering
  	\begin{overpic}[scale=.36, grid = false, tics=5, trim=-.5cm 0cm -1cm -.5cm, clip]{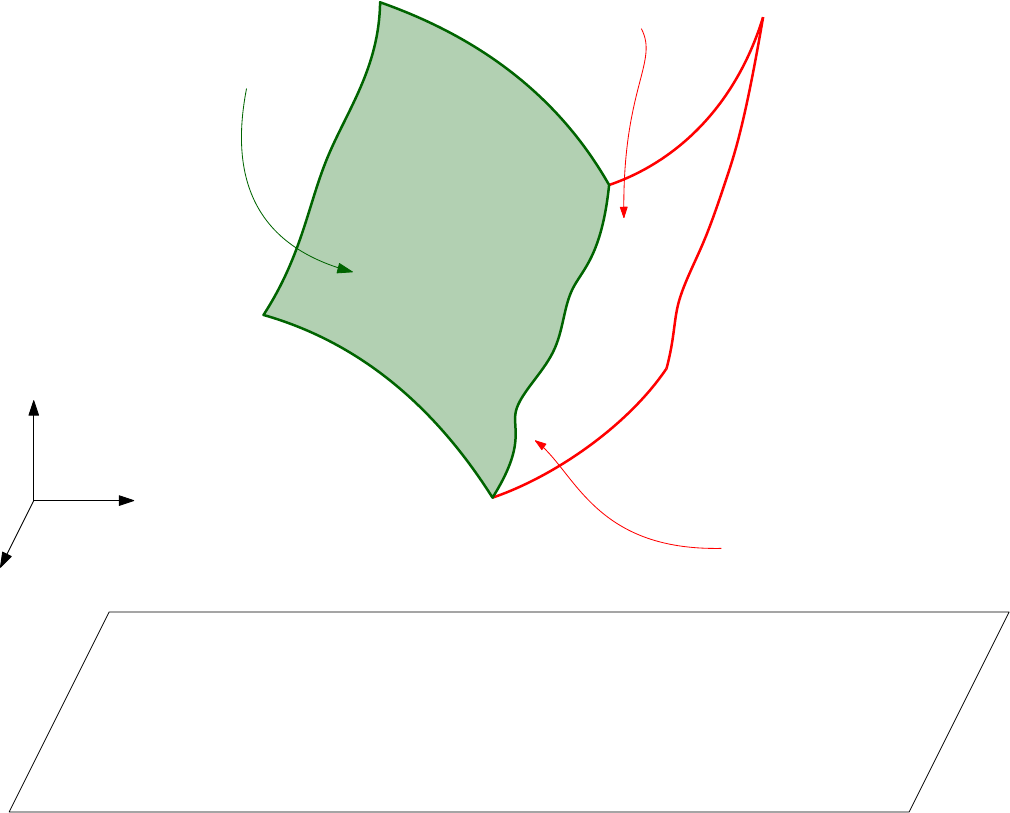}
			\put (50,7) {$\Sigma_0$}
			\put (4,20) {$(x^2,x^3) \in \mathbb{T}^2$}
			\put (2,34) {$t$}
			\put (16,27) {$x^1 \in \mathbb{R}$}
			\put (59,72.5) {$\mathcal{B}$}
			\put (70,22) {$\crease$}			
			\put (23,68.5) {$\underline{\mathcal{C}}$}
		\end{overpic}
		\caption{The Cauchy horizon, shaded in green, emanating from the crease and the subject of this paper.}
		\label{F:MAXDEVELOPMENTINCARTESIAN}
\end{subfigure}%
\caption{Cartesian coordinate space illustrations of the main results}
\label{F:INTROPICTURESOFMAINRESULTS}
\end{figure}
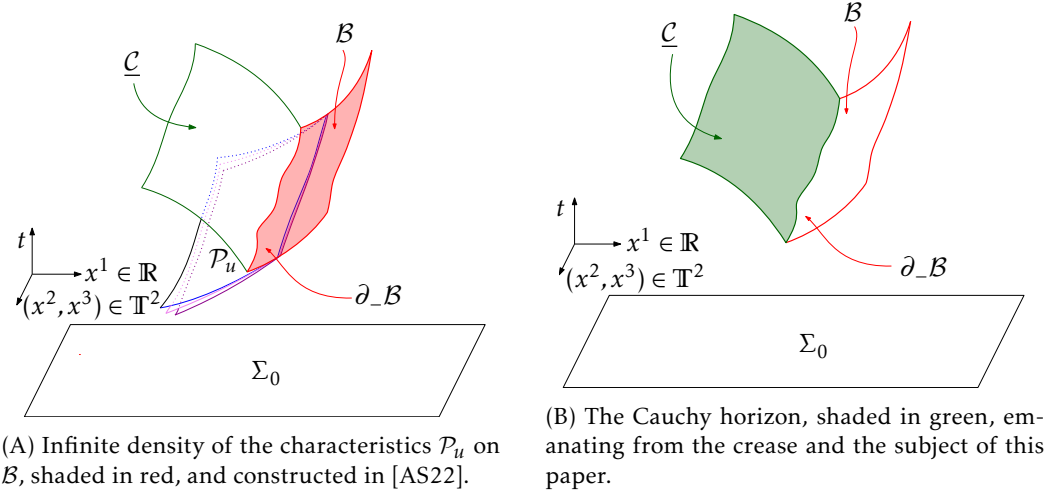

We now give a less rough overview of our work and method, which we remind the reader hold for solutions to \eqref{E:INTROTRANSPORTVI}--\eqref{E:INTROBS} with non-trivial vorticity, dynamic entropy, and without any symmetry reductions or strict convexity assumptions. The results mark the conclusion of the local-in-spacetime aspects of our research program initiated in \cite{abbrescia2022emergence} where, for \emph{general} small perturbations of smooth, non-degenerate, simple, plane-symmetric initial data, we constructed a large, localized portion of the maximal globally hyperbolic development (MGHD) of the data up to the boundary. In \cite{abbrescia2022emergence}, for these open sets, 
we proved that certain gradient components of the density $\Density$ and velocity $v$ become infinite along a distinguished co-dimension 2 submanifold of spacetime, denoted $\crease$, which can be thought of as the true ``initial singularity.'' We refer to 
$\crease$ as the \emph{crease}. These gradients then \emph{continue} to blow-up along a co-dimension 1 submanifold-with-boundary of spacetime, denoted $\singhyp$, which emanates from $\crease$.
We refer to $\singhyp$ as the \emph{singular boundary}. The gradient-blowup mechanism is an infinitely dense collapse of the characteristics along 
$\singhyp$, much as in the case of the $1D$ Burgers' equation; see Fig.\,\ref{F:INFININTEDENSITYOFCHARACTERISTICSONSINGULARBOUNDARY}. Also as in the case of Burgers' equation, the components of $\pmb{\partial}(\Density,v)$ that blow-up are \emph{transverse} to the collapsing characteristics, while the derivatives of the solution in directions tangent to the characteristics remain bounded.

\renewcommand{\thesubfigure}{\Alph{subfigure}} 

The purpose of this paper is to construct the Cauchy horizon $\Cauchyhor$, which is another co-dimension 1 submanifold-with-boundary of spacetime that, 
like $\singhyp$, also emanates from the crease; see Fig.\,\ref{F:MAXDEVELOPMENTINCARTESIAN}. We prove that the solution can be extended smoothly to $\Cauchyhor$, even though it still ``feels the influence'' of the gradient-blowup from $\crease$. By this, we mean that any $p \in \crease$ can be reached by a past-directed acoustically null geodesic in $\Cauchyhor$, i.e. $\Cauchyhor$ lies in the causal future domain of influence of $\crease$. These and all other notions of causality in the paper are measured with respect to the \emph{acoustical metric} $\gfour$ (see \eqref{E:ACOUSTICALMETRICANDMATERIALDER}), the solution-dependent Lorentzian metric that governs the geometry of sound waves. The reader can jump to Theorem~\ref{T:ABBREVIATEDMAINRESULTS} for an abbreviated version of the main results of this paper, 
and to Theorem~\ref{T:EXISTENCEUPTOCAUCHYHORIZONBYCONTINUATIONCRITERIA} for precise, extended statements.  

 Theorem~\ref{T:EXISTENCEUPTOCAUCHYHORIZONBYCONTINUATIONCRITERIA} and \cite[Theorem 34.1]{abbrescia2022emergence} are the first results to construct and fully justify 
Fig.\,\ref{F:INTROPICTURESOFMAINRESULTS} for open sets of initial data without symmetry, irrotationality, isentropicity, or 
\emph{strict} non-degeneracy. From a technical perspective, our main new contributions can be summarized as follows:

\begin{itemize}
	\item
	Compared to previous works, e.g., papers on the $2D$ isentropic or irrotational cases or papers that assumed strict convexity on the shape of the 	
	singular set, in the general $3D$ case with vorticity and entropy that we treat, 
	it is \emph{much harder} to follow the solution all the way to $\singhyp \cup \Cauchyhor$.
	A key difficulty is that the energy estimates are extremely intricate, in part because the foliations on which we must derive estimates are not pre-determined and have \underline{less regularity} than the solution, and 
	in part because of the presence of a vorticity stretching term (which is absent in $2D$) and the coupling to entropy.
	Consequently, to close the proof and avoid derivative-loss, we complement hyperbolic energy estimates with elaborate 
	\emph{elliptic-hyperbolic estimates} for the acoustic geometry and fluid, described below. 
	These estimates allow us to gain one derivative on the acoustic geometry, the vorticity, and the entropy compared to standard estimates.
	These elliptic-hyperbolic estimates can be avoided in the $2D$ isentropic or irrotational cases, but in the general $3D$ case, 
	it is not clear how to close the problem without them.
	\item We employ a double-null foliation of spacetime by characteristic surfaces to follow the solution up to the Cauchy horizon $\Cauchyhor$. A new aspect of the analysis that \emph{was not at all present in \cite{abbrescia2022emergence}} is that the $L^2$-hyperbolic estimates on characteristic surfaces are only  positive \emph{semi}-definite, compared to the positive definite $L^2$-hyperbolic estimates on the $\gfour$-spacelike hypersurfaces of \cite{abbrescia2022emergence}. This loss in coercivity implied that some dangerous top-order singular terms present in the commuted equations can no longer be handled as in \cite{abbrescia2022emergence}; see Sect.\,\ref{PAR:ALMOSTRIEMANNINVARIANTSNOTUSABLEFORCH}. In this paper, we overcome this new difficulty by uncovering a new hidden divergence structure that, upon integrating by parts, allows the singular terms to be handled perturbatively; see Sect.\,\ref{PAR:VELOCITYASAFUNDAMENTALUNKNOWN}. 
	This hidden structure was only uncovered due to the precise structure of the equations derived in \cite{jS2019c}. These issues were not present in \cite{shkoller2024geometry} because they did not utilize a double-null foliation to access the Cauchy horizon. Instead, they relied on ``approximate'' ingoing characteristic surfaces that are $\gfour$-spacelike and hence they did not encounter the above coercivity issues. We stress that the energy methods we developed specifically for double-null foliations are robust and have applications to other PDE problems, see \cite{speckyucharacteristicIVP}.
\end{itemize}

The portion of $\singhyp\cup\Cauchyhor$ accessed by our methods is of size $\mathcal{O}(1)$, and it contains an \emph{entire connected component of the crease} 
$\crease$. As we described in \cite[Sect.\,1.4]{abbrescia2022emergence} in detail, these theorems resolve several open problems and allow one to fully set up the \emph{shock-development problem} in $3D$, which is the problem of (locally) describing the transition of the flow from classical, to developing the crease as the ``initial singularity'', to developing a 
\emph{shock}, which, under a suitable weak formulation of the flow,
is a co-dimension $1$ hypersurface $\mathcal{K}$ emanating from $\crease$ across which the fluid variables experience discontinuous jumps. In other words, the initial data for the shock development problem \emph{is} the state of the fluid on $\crease\cup \Cauchyhor$. While the shock development problem is still open in $3D$ in full generality, there has been exciting progress, which we described in \cite[Sect. 1.9.12]{abbrescia2022emergence}. We refer the reader to the introduction of \cite{abbrescia2022emergence} for a focused, but far from exhaustive, literature review on recent multi-dimensional progress on singularities in compressible fluids. In this paper, we focus the discussion more on shock-formation results and how they relate to \cite[Theorem 34.1]{abbrescia2022emergence} and Theorem\,\ref{T:EXISTENCEUPTOCAUCHYHORIZONBYCONTINUATIONCRITERIA}.

We stress that our results \emph{apply to situations where $\crease$ is \textbf{not} ``convex in the $(x^2,x^3)$-directions''}. Full convexity along the crease characterizes what we meant by the \emph{strictly} non-degenerate case mentioned earlier in the context of Alinhac's work; 
see Fig.\,\ref{F:SHKOLLERMGHD} below. We also often refer to this setting as \textbf{\emph{strict convexity}}. In contrast, in the \emph{non-strict} case that we treat here, the crease $\crease$ does not need to have a single global lowest point. More generally, for symmetric shock-forming solutions and their perturbations, there will not be a single global lowest point in the singular set, 
and therefore, such solution regimes cannot be treated with techniques that rely on the strict convexity of $\crease$.
On the other hand, our assumptions on the data \emph{do} result in the ``upwards bending'' of $\singhyp$ in directions
transversal to the characteristics, e.g., in the direction of $x^1$ in Fig.\,\ref{F:INFININTEDENSITYOFCHARACTERISTICSONSINGULARBOUNDARY}. 
We call this weakened form of convexity \textbf{\emph{transversal convexity}}. Without transversal convexity, a continuum of characteristics could collapse into the same co-dimension $2$ set, and if that happens, it is not clear how to construct $\singhyp$ or whether the shock development problem is well-posed. 

We mention here the recent results of \cite{shkoller2024geometry} which, given strictly convex (i.e. strictly non-degenerate) data of large gradient size $\frac{1}{\epsilon}$ for the \emph{isentropic} $2D$ compressible Euler equations, followed the solution up until a small disconnected $\mathcal{O}(\epsilon)$ portion of $\crease\cup\singhyp\cup\Cauchyhor$, as in Fig.\,\ref{F:SHKOLLERMGHD}.  We discuss the similarities and differences between the methods and results of \cite{shkoller2024geometry} and our 2022 paper \cite{abbrescia2022emergence} (as well as the this one) in more detail in Appendix\,\,\ref{S:DICTIONARY}. 
In particular, there we prove (see Prop.\,\ref{P:ALETOGEOCOORD}) 
that the ALE coordinate systems used in \cite{shkoller2024geometry} are identical to the geometric coordinates afforded by Christodoulou's framework. We will also highlight some of the inspiration that \cite{abbrescia2022emergence, shkoller2024geometry} and this current work took from the breakthrough monograph of Christodoulou \cite{dC2007}, where a vision for how to study such problems was developed.

We close this section of the introduction by highlighting again that, in spite of the progress by \cite{dC2007,abbrescia2022emergence,shkoller2024geometry} and this current work, the problem of describing the entire MGHD of shock-forming solutions up to the boundary, as well as proving its extraordinarily subtle uniqueness, \emph{remains open in \underline{any} spatial dimension} for the compressible Euler equations. The only result constructing a provably unique MGHD up to its boundary for a system of hyperbolic PDEs with shock-forming data is our recent work \cite{abbresciaBlueSierskiSpeck2025quasilinear}, joint with P.\ Blue and J.\ Sbierski, for a model problem tailored to emulate the compressible Euler equations in $1D$. We elaborate on notions of MGHDs in Sect.\,\ref{SS:MGHDS}. We again stress that the question of MGHD uniqueness is global and generally cannot be answered if one has only constructed a portion of an MGHD.

\begin{center}
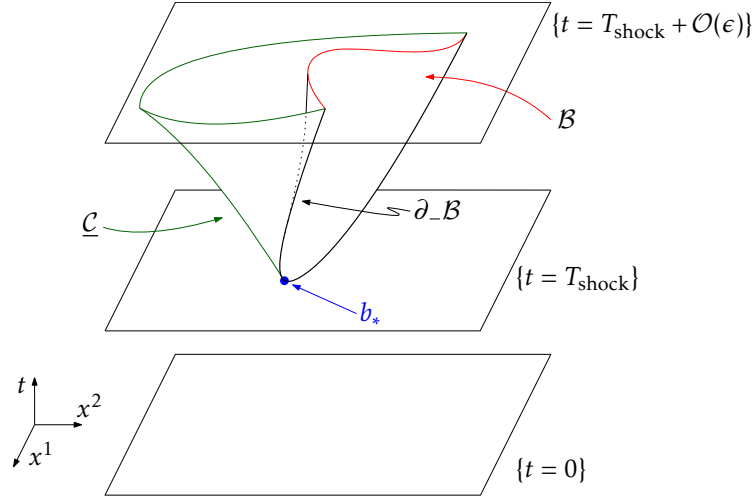
\begin{figure}  
\begin{overpic}[scale=.55, grid = false, tics=3, trim=-.5cm -1cm -1cm -.5cm, clip]{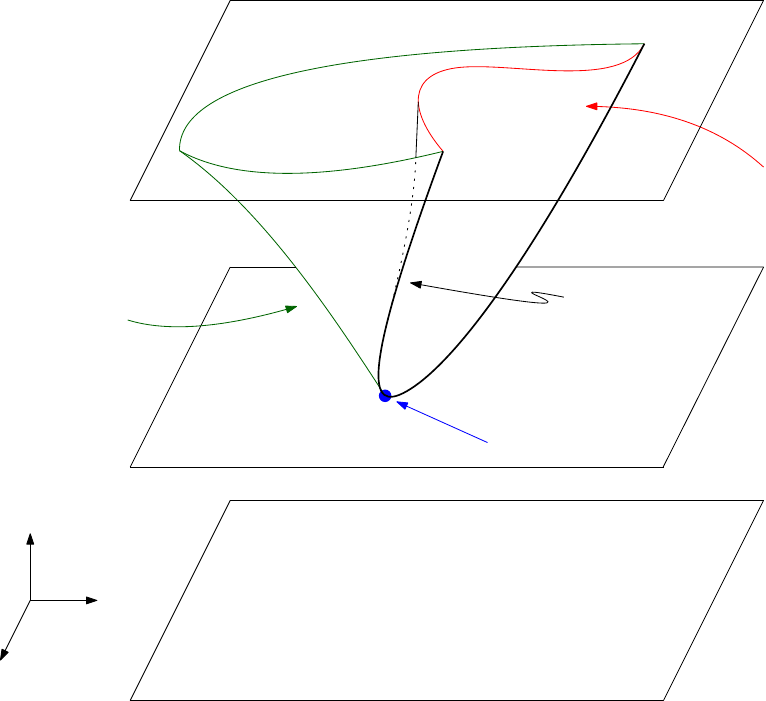}
	\put (6,12) {$x^1$}
	\put (4,24) {$t$}
	\put (14,20) {$x^2$}
	\put (87,10) {$\{t=0\}$}
	\put (87,42) {$\{t=T_{\textnormal{shock}}\}$}
	\put (93,84) {$\{t=T_{\textnormal{shock}} + \mathcal{O}(\epsilon)\}$}	
	\put (61,36) {{\color{blue}$b_*$}}
	\put (94,68) {$\mathcal{B}$}
	\put (70,53) {$\crease$}
	\put (15,51) {$\underline{\mathcal{C}}$}
\end{overpic} 
\caption{The crease, singular boundary, and Cauchy horizon derived in \cite{shkoller2024geometry} for \emph{strictly} non-degenerate shock forming data.}
\label{F:SHKOLLERMGHD}                                                                              
\end{figure}
\end{center}

\subsection{The Cauchy initial value problem for the Euler equations} \label{SS:CAUCHYPROBLEM}

 In \cite{abbrescia2022emergence}, we prescribed Cauchy data for \eqref{E:INTROTRANSPORTVI}--\eqref{E:INTROBS} on $\Sigma_0 \eqdef \{t=0\}= \{0\}\times\R\times\T^2$ which was a smooth (asymmetric) perturbation of simple plane-symmetric isentropic shock-forming data to \eqref{E:INTROTRANSPORTVI}--\eqref{E:INTROBS}. Of particular importance is that the background data led to a ``background solution'' which featured \emph{transversal convexity} along a large portion of the singular boundary. Transversal convexity along the singular boundary was preserved by the perturbations by choosing them to be sufficiently small in a sufficiently high-order Sobolev space.\footnote{In \cite{abbrescia2022emergence}, different solution variables had distinct, directionally dependent amounts of regularity.}
 
 To construct $\crease \cup \singhyp$, we controlled the fluid on a large region of classical existence $\MInteresting$ whose boundary contained  $\crease \cup \singhyp$.  It is within this region that we set up the data for the Cauchy horizon.  Although it is described in more detail below in Sect.\,\ref{SSS:CONSTRUCTIONOFEIKONALFUNCTION}, we mention here already that the data for the Cauchy horizon will be posed on a pair of transversally intersecting $\gfour$-null hypersurfaces. We construct these $\gfour$-null hypersurfaces from within $\MInteresting$ \emph{using the control of the solution from the already-proved estimates for the solution from \cite{abbrescia2022emergence}}. The Cauchy horizon is then obtained from a breakdown criterion of local solutions to the \emph{characteristic initial value problem}, with the data being the state of the fluid restricted to the aforementioned $\gfour$-null hypersurfaces. We provide our precise data assumptions in Sect.\,\ref{S:ASSUMPTIONSONTHEDATA}, 
and in Appendix~\ref{A:DATASSUMPTIONS}, we prove that these data assumptions are induced by the already-known solution in $\MInteresting$. 
 
We highlight the following subtle fact: \emph{it was imperative that we already controlled the solution up to and including $\crease$ in \cite{abbrescia2022emergence}}. Since the crease intersects the causal past of any point of the Cauchy horizon, the fluid on $\Cauchyhor$ ``feels the influence'' of the singularity on $\crease$, even though it is regular on it. For this reason, the pair of initial data $\gfour$-null hypersurfaces mentioned in the previous paragraph are constructed directly from the singularity at the $\crease$. We describe this construction in detail in Sect.\,\ref{S:CONSTRUCTIONANDESTIMATESFORINGOINGEIKONALFUNCTION}. 
It is not clear whether the fluid can be followed up to Cauchy horizon from data on $\Sigma_0$ without \emph{first} constructing $\crease$ as we did in \cite{abbrescia2022emergence}.

 \subsection{The geometric wave--transport formulation} \label{SS:GEOMETRICWAVETRANSPORTINTRO}
 Although we prescribe initial data for the Cauchy initial value problem of \eqref{E:INTROTRANSPORTVI}--\eqref{E:INTROBS}, we do not directly use the first order formulation for most of the paper.\footnote{More precisely, we do not directly commute the first order formulation to derive energy estimates. However, we \emph{do} rely on the first order formulation to express derivatives of the density in terms of the other fluid variables; see Sect.\,\ref{S:DENSITYDERIVATIVESINTERMSOFOTHERS}. This is convenient for the $L^2$-type energy method we employ.}

 It was proven by the second author \cite{jS2019c} that for sufficiently regular\footnote{The solutions in this paper are $C^{\infty}$ away from $\singhyp$.} solutions to \eqref{E:INTROTRANSPORTVI}--\eqref{E:INTROBS}, 
the fluid variables \emph{also} satisfy a system of covariant wave equations with respect to $\gfour$, coupled to transport-div-curl equations for the specific vorticity $\vortrenormalized^i \eqdef \frac{1}{\varrho} (\Flatcurl \, v)^i$ and gradient entropy $S^i\eqdef\p_i s$. In addition to being amenable to Klainerman's vectorfield method \cite{sK1985}, this formulation also makes the nonlinear interactions between sound waves and the advection of vorticity and entropy readily apparent. For $\Psi \in \vec{\Psi} \eqdef \{s,\uprho,v^1,v^2,v^3\}$ (where $\LogDensity = \ln(\varrho/\overline{\varrho})$ is the logarithmic density) and $\boldsymbol{\p}= (\p_t, \p_{x^1},\p_{x^2},\p_{x^3})$, the equations are, schematically:\footnote{In \eqref{E:ACOUSTICALMETRICANDMATERIALDER}, $\mathcal{Q}(\boldsymbol{\p} \wavearray,\boldsymbol{\p}\wavearray)$ are $\gfour$-null forms, which are quadratic nonlinearities with crucial nonlinear cancellations of self-resonant terms. The absence of such resonant terms on the RHS of the equations is fundamental for closing the estimates.}
\begin{align} \label{E:ACOUSTICALMETRICANDMATERIALDER}
	\gfour & \eqdef - \mathrm{d} t \otimes \mathrm{d} t + \Speed^{-2} \sum_{a=1}^3 (\mathrm{d} x^a - v^a \mathrm{d}t)\otimes  (\mathrm{d} x^a - v^a \mathrm{d}t), & & \Transport = \p_t + v^a \p_a,
\end{align}
\begin{subequations}
	\begin{align} 
		\Box_{\gfour(\vec{\Psi})} v^i  & \simeq (\Flatcurl \, \vortrenormalized)^i + \mathcal{Q}(\boldsymbol{\p}\vec{\Psi},\boldsymbol{\p}\vec{\Psi}), \label{E:COVWAVEEQv} \\
		\Box_{\gfour(\vec{\Psi})} \uprho & \simeq \Flatdiv S  + \mathcal{Q}(\boldsymbol{\p}\vec{\Psi},\boldsymbol{\p}\vec{\Psi}), \label{E:COVWAVEEQRHO} \\
		\Box_{\gfour(\vec{\Psi})} s & \simeq \Flatdiv S + \mathcal{Q}(\boldsymbol{\p}\vec{\Psi},\boldsymbol{\p}\vec{\Psi}), \label{E:COVWAVEEQS} \\
		 \Transport \vortrenormalized^i \simeq \boldsymbol{\p} \wavearray,  &  \qquad \qquad \qquad  \Transport S^i \simeq  \boldsymbol{\p} \wavearray,  \label{E:TRANSPORTVORTANDS}
	\end{align}
	\begin{align}
		\begin{split}  \label{E:MODIFIEDFLUIDVARIABLES} 
			\Transport(\Flatcurl\,\vortrenormalized)^i \simeq \mathcal{Q}(\boldsymbol{\p} \vortrenormalized , \boldsymbol{\p}\wavearray) + \mathcal{Q}(\boldsymbol{\p} \GradEnt , \boldsymbol{\p}\wavearray)  + \mathcal{Q}(\boldsymbol{\p}\vec{\Psi},\boldsymbol{\p}\vec{\Psi}),& \qquad \qquad \qquad  \Transport(\Flatdiv\,\GradEnt)  \simeq \mathcal{Q}(\boldsymbol{\p} \GradEnt , \boldsymbol{\p}\wavearray) + \mathcal{Q}(\boldsymbol{\p}\vec{\Psi},\boldsymbol{\p}\vec{\Psi}).  \\
			 \Flatdiv\, \vortrenormalized  \simeq \partial \wavearray, & \qquad \qquad \qquad   \Flatcurl\, S  = 0, 
		\end{split}
	\end{align}
\end{subequations}
Due to \eqref{E:COVWAVEEQv}--\eqref{E:MODIFIEDFLUIDVARIABLES}, $\wavearray$ are typically called the \emph{wave variables} and $(\vortrenormalized, \GradEnt, \Flatcurl\, \vortrenormalized, \Flatdiv S)$ are called the \emph{transport variables}.

\begin{remark}[Modified fluid variables]  \label{R:MODIFIEDFLUIDVARS} The higher order transport variables used in our work are not actually  $\Flatcurl\, \vortrenormalized, \Flatdiv S$, but rather a nonlinear perturbations of these up to lower-order-terms. Specifically, we use the following variables, which satisfy equations with surprisingly good structures:
	\begin{subequations}
		\begin{align} \label{E:MODIFIEDCURLOFVORTICITYINTRO}
			\VortVort^i
			& \eqdef
			\exp(-\LogDensity) (\Flatcurl \vortrenormalized)^i
			+
			\exp(-3\LogDensity) \Speed^{-2} \frac{p_{;\Ent}}{\overline{\varrho}} \GradEnt^a \partial_a v^i
			-
			\exp(-3\LogDensity) \Speed^{-2} \frac{p_{;\Ent}}{\overline{\varrho}} (\Flatdiv v) \GradEnt^i,
				\\
			\DivGradEnt 
			& \eqdef 
			\exp(-2 \LogDensity) \Flatdiv \GradEnt 
			-
			\exp(-2 \LogDensity) \GradEnt^a \partial_a \LogDensity.
			\label{E:MODIFIEDDIVERGENCEOFENTROPYGRADIENTINTRO}
		\end{align}
	\end{subequations}
In \eqref{E:MODIFIEDCURLOFVORTICITYINTRO}--\eqref{E:MODIFIEDDIVERGENCEOFENTROPYGRADIENTINTRO}, $p_{;\Ent}$ denotes differentiation of $p(\LogDensity,\Ent)$ with respect to $\Ent$ at fixed $\LogDensity$.
\end{remark}

\subsection{MGHDs} \label{SS:MGHDS}
We now elaborate on the concept of MGHDs in the context of the compressible Euler equations. We begin with the following definition, which introduces some causal notions of curves.

\begin{definition}[Causality and extendibility of curves] \label{D:CAUSALITYANDEXTENDIBLITYOFCURVES}
	Let $\mathcal{M} \subset \R^{1+3}$ be an open subset of spacetime on which the fluid variables
	$\vec{\Psi} := \{s,\uprho,v^1,v^2,v^3\}$ are defined and $C^1$,
	and let $\gfour(\vec{\Psi})$ be the acoustical metric defined in \eqref{E:ACOUSTICALMETRICANDMATERIALDER}. Let $\upgamma : I \to \mathcal{M}$ be a $C^1$ curve, where $I \subseteq \R$ is an interval of the form $[a,b], (a,b), [a,b)$, or $(a,b]$ where $-\infty \le a < b \le \infty$.
		\begin{itemize}
			\item We say that $\upgamma$ is $\gfour$-timelike ($\gfour$-null) if $\gfour(\dot \upgamma(s), \dot \upgamma(s)) < 0$ $(\gfour(\dot \upgamma(s), \dot \upgamma(s)) =0)$ for all $s \in I$. We say $\upgamma$ is $\gfour$-causal if it is either $\gfour$-timelike or $\gfour$-null.
			\item We say that $\upgamma$ is \emph{extendible} in $\mathcal{M}$ if at least one of the limits $\lim_{s \to a} \gamma(s)$ or $\lim_{s \to b} \gamma(s)$ exist \emph{in $\mathcal{M}$}. We note that this notion of extendibility is with respect to the \emph{domain} $\mathcal{M}$. 
			\item If $\upgamma$ is not extendible, it is called inextendible.
		\end{itemize}
\end{definition}

\begin{definition}[Cauchy hypersurfaces and Globally hyperbolic regions] \label{D:CAUCHYHYPERSURFACESANDGLOBALLYHYPERBOLICREGIONS}
	Let $\mathcal{M} \subset \R^{1+3}$ be an open subset of spacetime on which the fluid variables
	$\vec{\Psi} := \{s,\uprho,v^1,v^2,v^3\}$ are defined and $C^1$,
	and let $\gfour(\vec{\Psi})$ be the acoustical metric defined in \eqref{E:ACOUSTICALMETRICANDMATERIALDER}.
	\begin{itemize}
		\item  We say that a surface $\Sigma \subset \mathcal{M}$ is a \emph{Cauchy hypersurface} for $\mathcal{M}$ (with respect to $\gfour$) if every $\gfour$-timelike inextendible curve in $\mathcal{M}$ intersects $\Sigma$ exactly once.
		\item We say that $(\mathcal{M},\gfour(\vec{\Psi}))$ is \emph{globally hyperbolic} whenever $\mathcal{M}$ contains a Cauchy hypersurface.
	\end{itemize}
\end{definition}

\begin{definition}[Globally hyperbolic development] \label{D:GHD}
		A \emph{globally hyperbolic development} (GHD) of initial data  
		$\mathring{\vec{\Psi}} \eqdef \{\mathring{s},\mathring{\uprho},\mathring{v}^1,\mathring{v}^2,\mathring{v}^3\}$ 
		on a surface $\Sigma$ is a pair $(\mathcal{M},\vec{\Psi})$, where $\mathcal{M} \subset \R^{1+3}$ 
		is a spacetime region containing $\Sigma$ and $\vec{\Psi} \eqdef \{s,\uprho,v^1,v^2,v^3\}$ is a classical solution
		to \eqref{E:INTROTRANSPORTVI}--\eqref{E:INTROBS} on $\mathcal{M}$
		such that $\vec{\Psi} \restriction_{\Sigma} = \mathring{\vec{\Psi}}$ and such that
		$(\mathcal{M},\gfour(\vec{\Psi}))$ is globally hyperbolic and $\Sigma$ is a Cauchy hypersurface.
\end{definition}

	\begin{definition}[Maximal globally hyperbolic developments] \label{D:MGHD} \hfill
		\begin{itemize} 
			\item We say that a GHD $(\mathcal{M},\vec{\Psi})$ of the initial data $\mathring{\vec{\Psi}}$ on $\Sigma$ is extendible as a GHD
		 if there is another GHD $(\widetilde{\mathcal{M}},\widetilde{\vec{\Psi}})$ of $\mathring{\vec{\Psi}}$
		 such that $\mathcal{M} \subsetneq \widetilde{\mathcal{M}}$
		 and such that $\widetilde{\vec{\Psi}} \restriction_{\mathcal{M}} = \vec{\Psi} $.
		\item We say that $(\mathcal{M},\vec{\Psi})$ is a \emph{maximal globally hyperbolic development} (MGHD)
			of the initial data $\mathring{\vec{\Psi}}$ on $\Sigma$ if it is not extendible.
		\end{itemize}
	\end{definition}
	
The key subtle point from the definition of extendibility is that $\widetilde{\mathcal{M}}$ is \emph{also} required to be globally hyperbolic in its own right. In general, if $\widetilde{\mathcal{M}}$ is not, then there are many (often infinitely many) ways to extend a solution from $\mathcal{M}$ to $\widetilde{\mathcal{M}}$. 

For use below, we note that Definitions \ref{D:CAUCHYHYPERSURFACESANDGLOBALLYHYPERBOLICREGIONS}--\ref{D:MGHD}
extend in the natural fashion to the following Cauchy initial value problems on $\R^{1+d}$:
			\begin{subequations}
				\begin{align}
					(\hfour^{-1}(\Phi,\pmb{\p}\Phi))^{\alpha\beta} \p_\alpha\p_\beta \Phi & = N(\Phi,\pmb{\p}\Phi), \label{E:JANQNLW} \\
					(\Phi,\p_t \Phi)|_{t=0} & = (f,g), \label{E:JANQNLWDATA}
				\end{align}
			\end{subequations}
		whenever $\hfour(\Phi,\pmb{\p}\Phi)$ is a solution-dependent Lorentzian metric. 

For the compressible Euler equations or equations of type \eqref{E:JANQNLW}, 
the \emph{construction} of an MGHD of initial data on $\lbrace t= 0 \rbrace$ that terminates
due to shock-formation on its boundary, as well as the \emph{uniqueness} of any such MGHD, remains a fundamental open problem. More precisely, 
the only such result is our joint work with P.\ Blue and J.\ Sbierski \cite{abbresciaBlueSierskiSpeck2025quasilinear} for a model problem tailored to emulate the compressible Euler equations in $1D$. For hyperbolic equations posed on a fixed background such as the Euler equations, \emph{uniqueness} of an MGHD is an extremely subtle issue. 

Definitions\,\ref{D:CAUCHYHYPERSURFACESANDGLOBALLYHYPERBOLICREGIONS}--\ref{D:MGHD} fundamentally rely 
on having a \emph{Lorentzian} metric to define what it means to be ``causal.'' 
For hyperbolic systems \emph{without} a Lorentzian structure, there is no natural, compelling way to delineate spacetime regions
on which the uniqueness of classical solutions holds. The phenomenon of classical non-uniqueness is well-known for
model equations, such as Burgers' equation.  

	\begin{proposition}[Non-unique classical development for Burgers' equation, \cite{abbrescia2022emergence}*{Appendix A}, \cite{abbresciaspecknotices}*{Sect.\,2.1.1}] \label{P:NONUNIQUEBURGERS}
		Consider Burgers' equation $\partial_t \psi + (1+\psi)\partial_x \psi = 0$. Given initial data on $\{t=0\}$, we say $(\mathcal{M},\psi)$  is a \emph{classical development} of the data if $\psi : \mathcal{M}\subset\R^{1+1} \to \R$ is a $C^1$ solution to Burgers' equation. In particular, for every $p \in M$, $\psi(p)$ is determined by integrating back along a unique characteristic of $\partial_t + (1+\psi) \partial_x$ to $\{t =0\}$. 
		
	Then for the fixed initial condition $\psi(0,x) := -x+\tfrac{1}{3}x^3$,
	there are two distinct classical developments $(\mathcal{M}_1,\psi_1)$ and $(\mathcal{M}_2,\psi_2)$ such that 
	$\mathcal{M}_1\cap\mathcal{M}_2 \neq \emptyset$ and such that there exists a $p \in \mathcal{M}_1\cap\mathcal{M}_2 $ with $\psi_1(p) \neq \psi_2(p)$. 
	Both solutions $\psi_i$ develop gradient catastrophes on a portion of their boundary $\p \mathcal{M}_i$.  
	\end{proposition}
	
Remarkably, the following theorem shows that global hyperbolicity alone is \emph{not} sufficient to ensure uniqueness of an MGHD!  

	\begin{theorem}[\cite{fEhRjS2019}*{Theorem 3.45 and Remark 3.55}]  \label{T:NONUNIQUEMGHD} Let $d= 1$, $N = 0$, and 
	consider the quasilinear wave equation \eqref{E:JANQNLW} with the following metric:
		\begin{align}
			\hfour & := \begin{pmatrix}
				-(1+(\p_{x^1}\Phi)^2) & \p_t\Phi \p_{x^1} \Phi \\ 
				\p_t\Phi \p_{x^1} \Phi & 1 - (\p_t \Phi)^2
				\end{pmatrix}. \label{E:BORNINFELDSCALARMETRIC}
		\end{align}
	There exists smooth initial data on $\{t=0\}$ and two MGHDs $(\mathcal{M}_1,\Phi_1)$ and $(\mathcal{M}_2,\Phi_2)$ of the same data such that there exists a $p \in \mathcal{M}_1\cap \mathcal{M}_2$ with $\Phi_1(p) \neq \Phi_2(p)$.
	\end{theorem}
The results of Theorem~\ref{T:NONUNIQUEMGHD} are rather striking. This is because solutions of the initial value problem to Einstein's equations of general relativity -- where the notion of MGHDs was first developed -- that develop in a globally hyperbolic manner are unique! This fundamental theorem of mathematical relativity was proved in Choquet-Bruhat--Geroch’s celebrated and historic work \cite{cBgR1969}.

The equation defined by \eqref{E:BORNINFELDSCALARMETRIC} is the Born-Infeld scalar field equation. This is a well-studied equation in mathematical physics that enjoys a double-null condition with respect to $\hfour$. In particular, the non-uniqueness mechanism of Theorem~\ref{T:NONUNIQUEMGHD} is not due to the formation of shocks, which is contrast to Prop.\,\ref{P:NONUNIQUEBURGERS}. In the same paper \cite{fEhRjS2019}, the authors proved a sufficient condition for MGHD uniqueness, which we recall here.

	\begin{theorem}[\cite{fEhRjS2019}*{Theorem 4.82}] \label{T:UNIQUECONDITION}
		Let $(\mathcal{M},\Phi)$ be an MGHD for the quasilinear wave initial value problem \eqref{E:JANQNLW}--\eqref{E:JANQNLWDATA}. If $\mathcal{M}$ lies on one side of its boundary,\footnote{More precisely, for every $p \in \p \mathcal{M}$, there exists a neighborhood $U \subset \R^{1+d}$ of $p$ and a coordinate chart $f : U \to (-\epsilon,\epsilon)^{d+1}$, $\epsilon > 0$, and a continuous function $g : (-\epsilon,\epsilon)^d \to (-\epsilon,\epsilon)$ such that $\p \mathcal{M} \cap U = f^{-1}(\textnormal{graph $g$})$, all points \emph{below} the graph of $g$ are mapped into $\mathcal{M}$, and all points \emph{above} the graph of $g$ are mapped to $\R^{1+d} \setminus \mathcal{M}$.} then $(M,\Phi)$ is the unique MGHD of the data.
	\end{theorem}

We note that in order to apply Theorem~\ref{T:UNIQUECONDITION}, \emph{the MGHD must be known in its entirety}, and then the global one-sidedness property must be proved \emph{a posteriori} \underline{globally}. To date, no MGHD of shock-forming data to the compressible Euler equations has ever been constructed in \emph{any spatial dimension}. We will address this issue in future work.

In particular, since there is a large portion of the data whose domain of dependence lies outside the portion of spacetime between $\Sigma_0$ and $\Cauchyhor\cup\crease\cup\singhyp$, neither \cite{dC2007, abbrescia2022emergence, shkoller2024geometry} nor this paper construct an MGHD. Although \cite{shkoller2024geometry} do not construct an MGHD, they define a ``MGHD in a box''. If the first time of gradient blow-up is $T_{\textnormal{shock}}$, which exists due to their strict convexity restrictions, \cite{shkoller2024geometry}*{Sect.\,1.3} defines ``MGHD in a box'' as:
	\begin{quote}
		``The intersection of the MGHD of [strictly convex] data posed on $\{t=T_{\textnormal{shock}} - C\epsilon\}$  with $[T_{\textnormal{shock}} - C\epsilon, T_{\textnormal{shock}} + C\epsilon]\times \T^2$.''
	\end{quote}
	The ``MGHDs in a box'' from \cite{shkoller2024geometry}
	are $\mathcal{O}(\epsilon)$-thick GHDs that contain an $\mathcal{O}(\epsilon)$-size portion of a convex singular boundary 
	and an $\mathcal{O}(\epsilon)$-size portion of a corresponding Cauchy horizon.

 \subsection{Portion of the MGHD up to the singular boundary} \label{SS:RESULTSFROMSINGULARBOUNDARYPAPER}
 In this section, we give a terse overview of the methods and results from \cite{abbrescia2022emergence} that we need to construct the Cauchy horizon. We give a precise statement of the needed results in Sect.\,\ref{S:CLASSICALEXISTENCEUPTOTHECREASE}, specifically in Theorem~\ref{T:MAINRESULTSFROMSINGULARBOUNDARYPAPER}.
 
 \subsubsection{Gradient blow-up by nonlinear geometric optics} \label{SSS:NONLINEARGEOMETRICOPTICS}
 The gradient singularities along $\crease$ are caused by an infinitely dense collapse of a family of ``outgoing'' characteristics; see Fig.\,\ref{F:INFININTEDENSITYOFCHARACTERISTICSONSINGULARBOUNDARY}. The analysis we use to prove gradient blow-up, called \emph{nonlinear geometric optics}, is carried out in a \emph{geometric coordinate system} $(t,u,x^2,x^3)$, where $u$ is the solution to the following \emph{eikonal equation} initial value problem: 
 	\begin{subequations}
		\begin{align}
			(\gfour^{-1})^{\alpha\beta} \partial_\alpha u \partial_\beta u & = 0, \label{E:EIKONALEQUATIONINTRO} \\
			\partial_t u & > 0, \label{E:EIKONALBRANCHINTRO} \\
			u|_{t=0} & = -x^1. \label{E:EIKONALDATAINTRO}
		\end{align}
	\end{subequations}
The level sets of $u$ in Cartesian space $\R\times\R\times\T^2$, which we denote by $\nullhyparg{u}$, are called \emph{characteristics} or (outgoing) \emph{null hypersurfaces}.
Our analysis shows that the solution to \eqref{E:EIKONALDATAINTRO} is such that the characteristics become infinitely dense along $\singhyp$ 
precisely when the fluid gradient blows up. The geometric coordinate system $(t,u,x^2,x^3)$ can be thought of as a gauge that
``unfolds'' the collapsing characteristics. 
For this reason, the solution remains rather smooth relative to the differential structure of the geometric coordinates. 
The price one pays is that the geometric coordinate system degenerates relative to the Cartesian one as the singularity forms.
  
 Although the characteristics do \emph{not} become infinitely dense in geometric coordinates, one can still identify their collapse by the vanishing of a non-negative function $\upmu \approx |\partial_x u|^{-1}$ called the \emph{inverse foliation density}. $\upmu$ is positive everywhere in the MGHD and vanishes on $\crease\cup\singhyp$ both in Cartesian and geometric coordinates. Roughly, $\frac{1}{\upmu}$ measures the density of $\nullhyparg{u}$ relative to constant-time-hypersurfaces and thus: 
\begin{quote}
	The study of shock formation essentially entails following the solution all the way up to the level set $\lbrace \upmu = 0 \rbrace$,
	which contains the singular boundary $\singhyp$.
\end{quote}
The relationship between the fluid variables and $\upmu$ in the geometric and Cartesian differential structures is given by the following schematic estimate:
	\begin{align}
  		|\partial_{x^i} (v,\LogDensity)| \approx \frac{1}{\upmu} \left|\geop{u}(v,\LogDensity)\right|. \label{E:SCHEMATICGRADIENTBLOWUP}
	\end{align}
In \eqref{E:SCHEMATICGRADIENTBLOWUP}, $\geop{u}$ denotes partial differentiation with respect to $u$ at fixed $(t,x^2,x^3)$. The data assumptions and ``smoothness'' in geometric coordinates imply that $ 0 < C \le  \left|\geop{u}(v,\varrho)\right| < \infty$, where $C$ depends on the data. From this and \eqref{E:SCHEMATICGRADIENTBLOWUP}, one recovers the gradient singularity in the Cartesian coordinate differential structure as $\upmu \downarrow 0$. Here, we write ``smoothness'' in quotes because the $L^2$-analysis leaves the potential for the mid-to-top-order energy of the fluid to blow up as $\upmu \downarrow 0$, \emph{even} in geometric coordinates.

We highlight the following under-appreciated fact regarding the methods described in the previous two paragraphs, pioneered by Christodoulou in his 2007 monograph \cite{dC2007}:

	\begin{quote}
		Nonlinear geometric optics can be used to prove shock formation up to the first Cartesian time of gradient blow-up \emph{for any perturbation of a plane-symmetric shock-forming background solution, \textbf{even in the very degenerate case where transversal convexity is \underline{not present at all.}}} This was accomplished for the full $3D$ compressible Euler equations in \cite{LukSpeck2024stability} (see also \cite{luk2018shock} for $2D$ isentropic Euler). When transversal convexity \textbf{is} present, 
		nonlinear geometric optics can be supplemented by the methods we introduced in \cite{abbrescia2022emergence} to go \emph{past} the first Cartesian time of 
		gradient blow-up, allowing one to access the entirety of $\crease$ and an $\mathcal{O}(1)$ portion of $\singhyp$.  
	\end{quote}

\noindent Techniques using modulated self-similar analysis were recently developed \cite{tBsSvV2022,tBsSvV2019a,tBsSvV2020} to also follow the solutions \emph{with \textbf{strict} convexity} up to the first Cartesian time of gradient blow-up, say $T_{\textnormal{shock}}$, but \text{not} further. For the $3D$ compressible Euler equations, these papers derived an analog of Alinhac’s shock formation results \cite{sA1999a,sA1999b} for quasilinear wave equations failing the classical null condition. Alinhac's solution also featured \emph{\textbf{strict} convexity}. We emphasize that \cite{sA1999a,sA1999b,tBsSvV2019a,tBsSvV2020,tBsSvV2022}'s requirement for strict convexity implied that the results do not hold for perturbations of symmetric shock-forming solutions. In contrast, \cite{LukSpeck2024stability} proved shock formation for the full compressible Euler equations in $3D$ for a much larger class of data \emph{without} even \textbf{transversal} convexity and can incorporate perturbations of symmetric shock-forming solutions. A consequence of this is that gradient-blow up in  \cite{sA1999a,sA1999b,tBsSvV2022,tBsSvV2019a,tBsSvV2020} occurs at an isolated point, i.e.  $\Sigma_{T_{\textnormal{shock}}} \cap \{\upmu=0\} = b_*$, while \cite{LukSpeck2024stability} allowed for much more general singular sets; see Fig.\,\ref{F:MULTIDSHOCKFORMATION}.

	\begin{remark}[$C^{0,1/3}$ cusps]
	When the solution \emph{does} have transversal or strict convexity, a sharper picture of the singular behavior at points where $\upmu = 0$ can be obtained. In \cite{sA1999a,sA1999b,tBsSvV2019a,tBsSvV2020,tBsSvV2022}, \textbf{strict} convexity was used to prove that the fluid variables in fact developed a $C^{0,1/3}$ cusp at an isolated point. \cite[Cor. 4.5]{LukSpeck2024stability} proved that the weaker assumption of \textbf{transversal} convexity is sufficient for  $C^{0,1/3}$ cusps to develop on $\Sigma_{T_{\textnormal{shock}}} \cap \{\upmu=0\}$.
	\end{remark}

\begin{center}
\begin{figure}  
\begin{overpic}[scale=.4, grid = false, tics=3, trim=-.5cm -1cm -1cm -.5cm, clip]{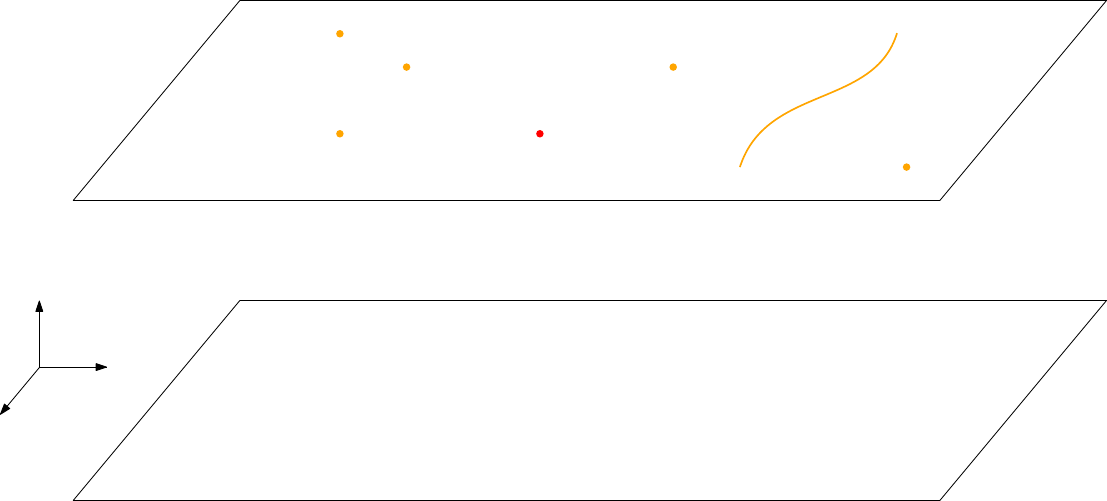} 
	\put (12,16) {$x^1$}
	\put (7,19) {$t$}
	\put (88,10) {$\Sigma_0$}
	\put (87,30) {$\Sigma_{T_{\textnormal{shock}}}$}	
	\put (48,33) {{\color{red}$b_*$}}
\end{overpic} 
\caption{The set of points where shock-formation occurs in \cite{sA1999a,sA1999b,tBsSvV2022,tBsSvV2019a,tBsSvV2020,LukSpeck2024stability}. The isolated point $\Sigma_{T_{\textnormal{shock}}} \cap \{\upmu=0\} = b_*$ in \cite{sA1999a,sA1999b,tBsSvV2022,tBsSvV2019a,tBsSvV2020} is shown as red, whereas the geometry of $\Sigma_{T_{\textnormal{shock}}} \cap \{\upmu=0\}$ from \cite{LukSpeck2024stability}  is allowed to be much more general (i.e. it can even feature connected curves) and is shown as orange.}
\label{F:MULTIDSHOCKFORMATION}
\end{figure}
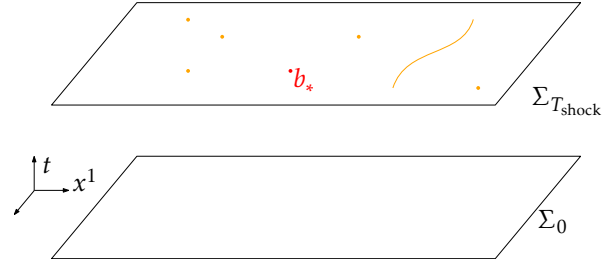
\end{center}

The following table summarizes the latest breakthroughs in multi-dimensional settings without symmetry up until the time of first shock formation. 

\begin{table}[H]
\centering
\begin{tabular}{|c|c|c|c|}
\hline
 & \makecell{$2D$ isentropic \\ Euler  up to $T_{\textnormal{shock}}$} &  \makecell{$3D$ quasilinear scalar \\  equations up to $T_{\textnormal{shock}}$} & \makecell{$3D$ full Euler \\ system up to $T_{\textnormal{shock}}$} \\ \hline
 Strict convexity & \cite{tBsSvV2019a,luk2018shock} & \cite{dC2007,jS2016b} & \cite{tBsSvV2022,LukSpeck2024stability} \\ \hline
Transversal convexity  & \cite{luk2018shock} & \cite{dC2007,jS2016b} & \cite{LukSpeck2024stability} \\ \hline
No convexity  present & \cite{luk2018shock}&  \cite{dC2007,jS2016b} & \cite{LukSpeck2024stability} \\ \hline
\end{tabular}
\end{table}

The following table summarizes the latest breakthroughs in multi-dimensional settings describing portions of the boundary of the MGHD, beyond the time for first blowup.

\begin{table}[H]
\centering
\begin{tabular}{|c|c|c|c|c|c|c|}
\hline
 &   \makecell{$2D$ isentropic \\ Euler up to $\singhyp$} &   \makecell{$2D$ isentropic \\ Euler up to $\Cauchyhor$} &  \makecell{$3D$ quasilinear scalar \\  equations up to $\singhyp$} & \makecell{$3D$ full Euler \\ system up to $\singhyp$ } & \makecell{$3D$ full Euler\\ system up to $\Cauchyhor$} \\ \hline
 Strict convexity & \cite{shkoller2024geometry} & \cite{shkoller2024geometry} & \cite{dC2007} & &\\ \hline
Transversal convexity  &    &  &  & \cite{abbrescia2022emergence}  & This paper \\ \hline
\end{tabular}
\end{table}

\subsubsection{Rough foliations and time functions up to $\crease\cup\singhyp$} \label{SSS:ROUGHFOLIATIONSANDTIMEFUNCTION} 
The points on $\Sigma_{T_{\textnormal{shock}}} \cap \{\upmu=0\}$ described in Sect.\,\ref{SSS:NONLINEARGEOMETRICOPTICS} constitute some points of $\crease$. Alinhac used nonlinear geometric optics in \cite{sA1999a,sA1999b}, but unlike Christodoulou, he used Nash-Moser estimates to recover a loss of regularity generated when derivatives in the commuted equations fell on error terms from the acoustic geometry. It was this technical implementation of Nash-Moser that prevented Alinhac from going past the first Cartesian time of gradient-blow-up $T_{\textnormal{shock}}$. In \cite{abbrescia2022emergence}, we instead recover this loss of derivatives through an exhaustive study of the acoustic geometry near the shock using a variant of a method introduced in \cite{dC2007}. Roughly, Christodoulou found a nonlinear combination of acoustic geometry + fluid variables which solved a transport equation called \emph{Raychaudhuri's equation} \cite{aR1955} that could be explicitly integrated to recover the loss.  Christodoulou's motivation to use Raychaudhuri's equation came from his earlier groundbreaking work with Klainerman on the stability of Minkowski space as solutions to the Einstein vacuum equations \cite{christodoulou1994global}.

In order to control the fluid up to $\crease$, we constructed in \cite{abbrescia2022emergence} a new quasilinear \emph{rough foliation}   $\cup_{\timefunction \in [\timefunction_0,0]} \hyptwoarg{\timefunction}{0}$ of spacetime which was able to accurately ``predict'' the geometry of the crease. By this, we mean that the crease was completely contained in the final leaf of the foliation, i.e., $\crease \subset \hyptwoarg{0}{0}$. By ``quasilinear,'' we mean that each of these rough hypersurfaces were dynamically constructed from the state of the fluid near the expected time of shock formation. For each $\timefunction \in [\timefunction_0,0]$ where $\timefunction_0 < 0$ is a fixed constant that depends only on the background solution, they are level sets of a \emph{rough time function} $\timefunctionarg{0} = \timefunctionarg{0}(t,u,x^2,x^3)$, i.e. $\hyptwoarg{\timefunction}{0} \eqdef \{ (t,u,x^2,x^3) \, : \, \timefunctionarg{0}(t,u,x^2,x^3) = \timefunction\}$. We mention here that in \cite{abbrescia2022emergence}, we actually constructed a continuum of rough foliations (or, equivalently, a continuum of rough time functions). For any $\muxmulevelsetvalue \in [0,\muxmulevelsetvalue_0]$ where $\muxmulevelsetvalue_0 > 0$ is a fixed constant that depends only on the background solution, we constructed a rough foliation  $\cup_{\timefunction \in [\timefunction_0,0]} \hyptwoarg{\timefunction}{\muxmulevelsetvalue}$ of spacetime which accurately ``predicts'' a distinguished two-dimensional submanifold $\twoargmumuxtorus{0}{-\muxmulevelsetvalue}$
of $\singhyp$. As a set,\footnote{In defining  $\twoargmumuxtorus{0}{-\muxmulevelsetvalue}$, $\muX = \upmu X$ where $X$ is the $\Sigma_t$-tangent $\gfour$-spacelike unit vectorfield (i.e., $\gfour(X,X) = 1$) which is $\gfour$-orthogonal to the tori $\Sigma_t\cap \nullhyparg{u}$. Schematically, $\muX \simeq \geop{u}$.} $\twoargmumuxtorus{0}{-\muxmulevelsetvalue} \eqdef \{\upmu = 0 \} \cap \{\muX \upmu = -\muxmulevelsetvalue\}$ was always a subset of the last leaf, i.e. $\twoargmumuxtorus{0}{-\muxmulevelsetvalue} \subset  \hyptwoarg{0}{\muxmulevelsetvalue}$. In particular, $\crease \eqdef \{\upmu = 0 \} \cap \{\muX \upmu = 0\} = \twoargmumuxtorus{0}{0}$ and the portion of the singular boundary constructed in \cite{abbrescia2022emergence} is precisely $\singhyp = \cup_{\muxmulevelsetvalue \in [0,\muxmulevelsetvalue_0]} \twoargmumuxtorus{0}{-\muxmulevelsetvalue}$.

\begin{figure}[ht]
\centering
\begin{subfigure}{.5\textwidth}
  \centering
  	\begin{overpic}[scale=.36, grid = false,trim=-.5cm -1cm -1cm -1.9cm, clip=true]{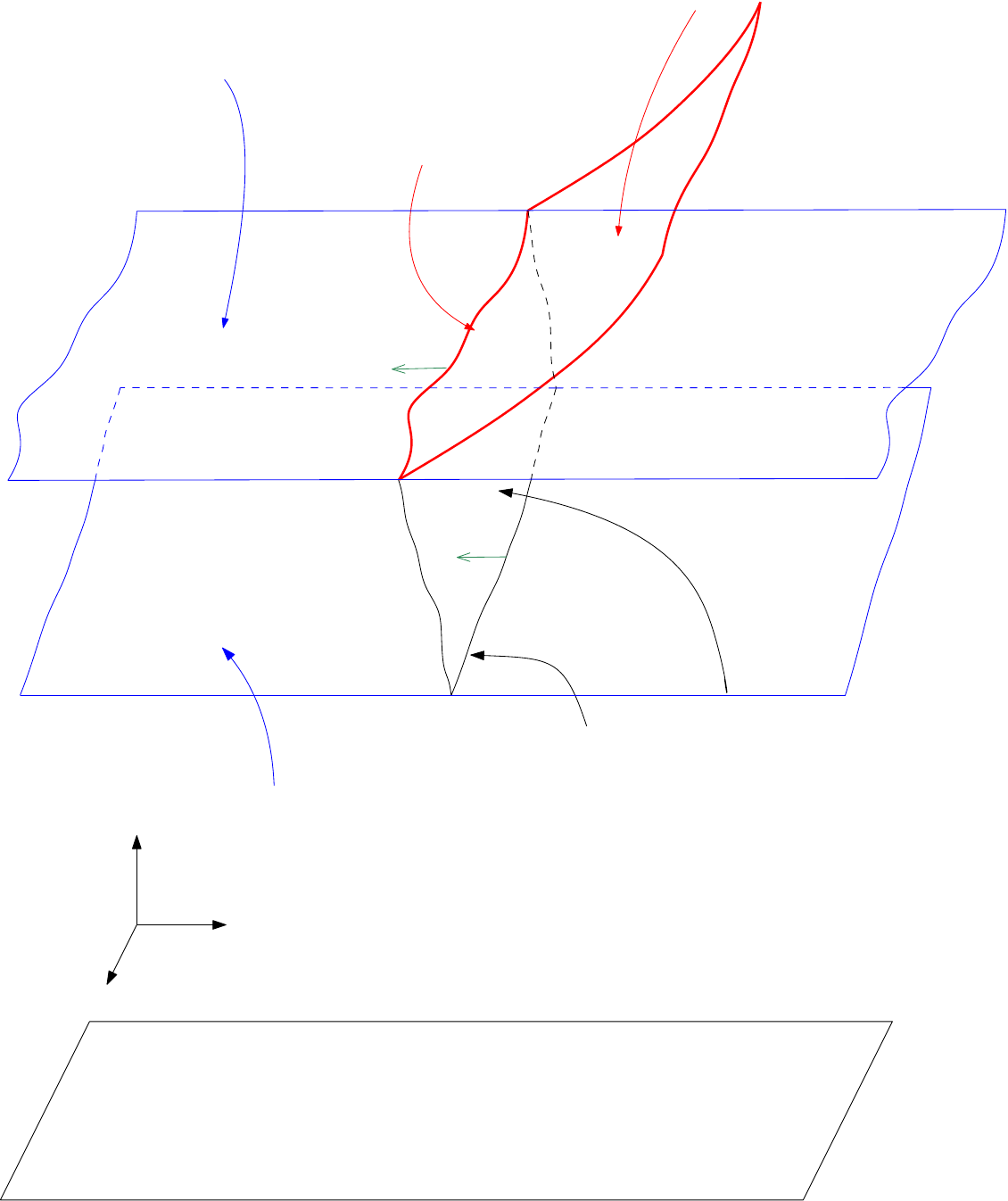}
			\put (36,11) {$\Sigma_0$}
			\put (11,19) {$(x^2,x^3) \in \mathbb{T}^2$}
			\put (10,29) {$t$}
			\put (20,23.5) {$x^1 \in \mathbb{R}$}
			\put (31,90) {$\twoargmumuxtorus{0}{0}$}
			\put (34,85) {\rotatebox{90}{$=$}}
			\put (31,82) {$\partial_- \mathcal{B}$}
			\put (54,93) {$\mathcal{B}$}
			\put (12,89) {$\hyptwoarg{0}{0}$}
			\put (42.5,36.5) {$\twoargmumuxtorus{-\timefunction}{0}$}
			\put (55,37) {$\datahypfortimefunctionarg{0}$}
			\put (18,30) {$\hyptwoarg{\timefunction}{0}$}
			\put (24.5,65) {$\Wtransarg{0}$}
			\put (31,52) {$\Wtransarg{0}$}
		\end{overpic}
		\caption{Rough foliations adapted to the crease}
		\label{F:CARTESIANROUGHFOLIATIONCREASE}
\end{subfigure}%
\begin{subfigure}{.5\textwidth}
 \centering
\begin{overpic}[scale=.36,grid=false,trim=-.5cm -1cm -1cm -.5cm, clip=true]{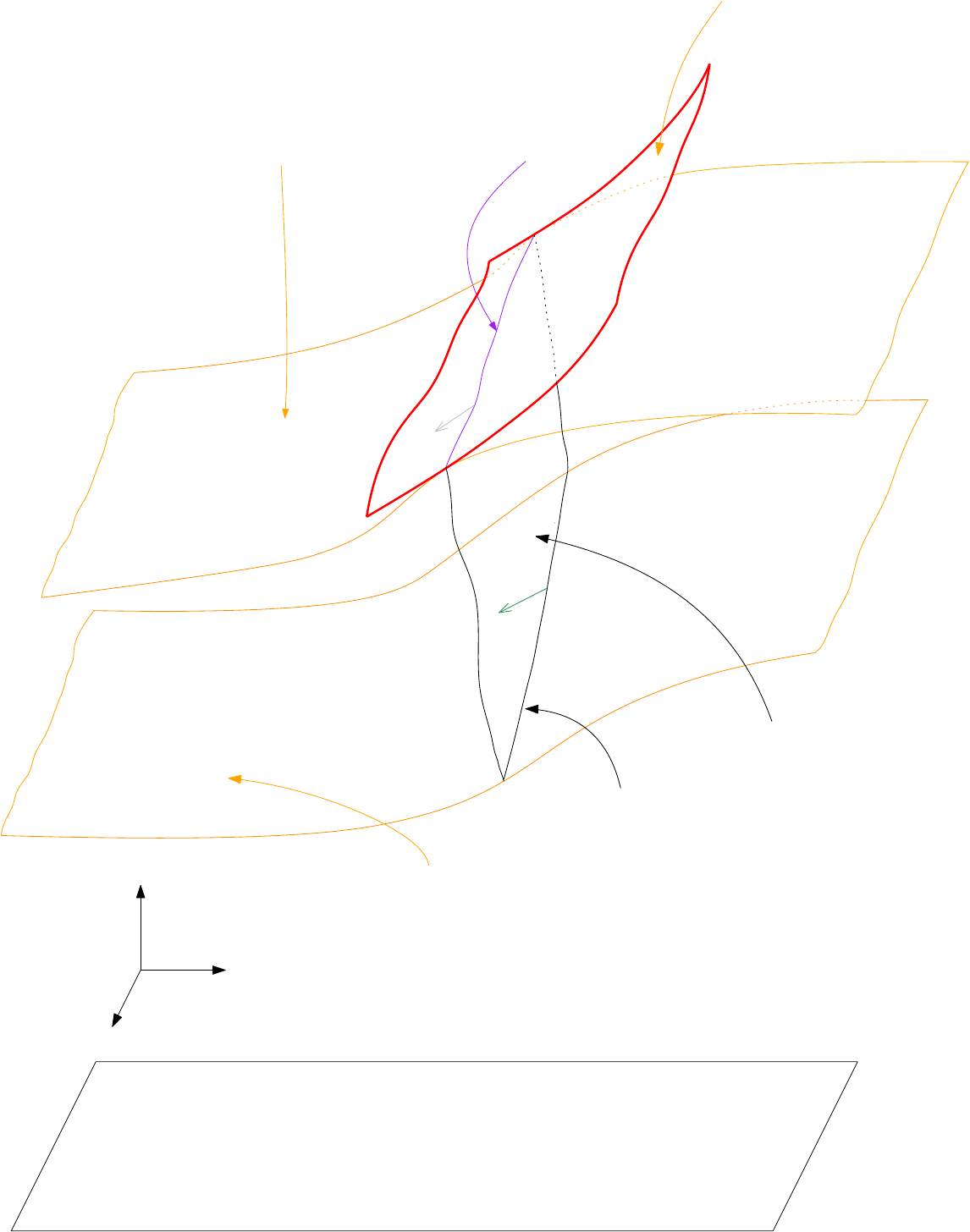} 
			\put (36,11) {$\Sigma_0$}
			\put (11.5,19) {$(x^2,x^3) \in \mathbb{T}^2$}
			\put (10,29) {$t$}
			\put (20,23.5) {$x^1 \in \mathbb{R}$}
			\put (39,88) {$\twoargmumuxtorus{0}{-\muxmulevelsetvalue}$}
			\put (58,98) {$\mathcal{B}$}
			\put (46,35) {$\twoargmumuxtorus{-\timefunction}{-\muxmulevelsetvalue}$}
			\put (15,90) {$\hyptwoarg{0}{\muxmulevelsetvalue}$}
			\put (59,39) {$\datahypfortimefunctionarg{-\muxmulevelsetvalue}$}
			\put (30,27.5) {$\hyptwoarg{\timefunction}{\muxmulevelsetvalue}$}
			\put (29,66) {$\Wtransarg{\muxmulevelsetvalue}$}
			\put (34,52.5) {$\Wtransarg{\muxmulevelsetvalue}$}
\end{overpic}
  \caption{Rough foliations adapted to a non-crease torus $\twoargmumuxtorus{0}{-\muxmulevelsetvalue} \subset \mathcal{B}$}
	\label{F:CARTESIANROUGHFOLIATIONFUTUREOFCREASE}
\end{subfigure}
\caption{Rough foliations adapted to the singular boundary, depicted in Cartesian coordinate space}
\label{F:INTROCARTESIANROUGHFOLIATIONS}
\end{figure}

Of fundamental importance in our work \cite{abbrescia2022emergence}, and particularly fundamental in future applications to the shock development problem, was that $\{ \upmu = - \timefunction\}$ and $\{\muX \upmu = -\muxmulevelsetvalue\}$ intersected \emph{transversally}\footnote{This transversality is a consequence of transversal convexity.} in a two-dimensional manifold $\twoargmumuxtorus{-\timefunction}{-\muxmulevelsetvalue}\eqdef\{ \upmu = - \timefunction\}\cap\{\muX \upmu = -\muxmulevelsetvalue\}$ that was diffeomorphic to $\T^2$ for $(\timefunction,\muxmulevelsetvalue)\in [\timefunction_0,0]\times[0,\muxmulevelsetvalue_0]$. To see why, fix $\muxmulevelsetvalue \in [0,\muxmulevelsetvalue_0]$ and define $\Wtransarg{\muxmulevelsetvalue} \eqdef \muX + \frac{\muxmulevelsetvalue}{\Lunit \upmu} \Lunit$, 
where the quantities on the RHS are rigorously constructed in Sect.\,\ref{SSS:IMPORTANTACOUSTICVECTORFIELDS}.
We constructed the rough time functions using the state of the fluid on $\twoargmumuxtorus{-\timefunction}{-\muxmulevelsetvalue}$ by solving the following Cauchy initial value problem:
\begin{subequations} \label{E:CAUCHYIVPFORROUGHTIMEFUNCTION}
	\begin{align}
		\Wtransarg{\muxmulevelsetvalue} \timefunctionarg{\muxmulevelsetvalue} & = 0, \label{E:EQUATIONFORROUGHTIMEFUNCTION}   \\
		\timefunctionarg{\muxmulevelsetvalue}\big|_{\twoargmumuxtorus{-\timefunction}{-\muxmulevelsetvalue}} & = -\upmu\big|_{\twoargmumuxtorus{-\timefunction}{-\muxmulevelsetvalue}} = \timefunction. \label{E:DATAFORROUGHTIMEFUNCTION}
	\end{align}
\end{subequations}
The Cauchy IVP \eqref{E:CAUCHYIVPFORROUGHTIMEFUNCTION} is well-posed because of transversal convexity. Indeed, $\Wtransarg{\muxmulevelsetvalue}$ is tangent to level sets of $\upmu$ but \emph{transverse} to level sets of $\muX \upmu$ due to transversal convexity, meaning that $\{\muX \upmu = -\muxmulevelsetvalue_0\}$ is ``not characteristic'' to $\Wtransarg{\muxmulevelsetvalue}$ in the sense of \cite{hormander1997lectures}.

\subsubsection{Energy descent scheme in rough geometric coordinates} \label{SSS:ENERGYDESCENTSCHEMEINROUGHCOORDINATES}
We now elaborate on what we referred to as ``remains rather smooth'' in Sect.\,\ref{SSS:NONLINEARGEOMETRICOPTICS} when using nonlinear geometric optics.
Our analysis is based on a descent scheme whose broad strategy originated in
the works of Alinhac \cite{sA1999a,sA1999b} and Christodoulou \cite{dC2007}.
We rely on many technical innovations tied to the special structure of equations \eqref{E:COVWAVEEQv}--\eqref{E:MODIFIEDFLUIDVARIABLES}
in order to implement the scheme.
For convenience and brevity, we focus on the case $\muxmulevelsetvalue = 0$. The precise tensorial components which blow-up in the schematic expression $|\partial (\LogDensity,v^i)|\sim \frac{1}{\upmu}$ are derivatives which are \emph{transverse} to the outgoing characteristics. Hence, the set of commutators $\mathscr{P}$ applied to equations \eqref{E:COVWAVEEQv}--\eqref{E:COVWAVEEQS} consists of derivative operators which are \emph{tangent} to the characteristics. These are $\{\Lunit, \Yvf{2},\Yvf{3}\}$, where $\Lunit^\alpha = - \upmu (\gfour^{-1})^{\alpha\beta}\p_\beta u$ is a $\gfour$-null vectorfield tangent to $\nullhyparg{u}$ normalized by $\Lunit t = 1$, and $\Yvf{A}$ is the $\Sigma_t$-tangent $\gfour$-orthogonal projection of $\p_{x^A}$ onto $\nullhyparg{u}$ (which removes the transverse components) for $A = 2,3$. The energy estimates for the commuted wave equations \eqref{E:COVWAVEEQv}--\eqref{E:COVWAVEEQS} are coupled to the transport equations \eqref{E:TRANSPORTVORTANDS}--\eqref{E:MODIFIEDFLUIDVARIABLES} (which we also commute with with $\mathscr{P}$), but we do not elaborate on that aspect here; see Sects.\,\ref{S:BELOWTOPORDERHYPERBOLICL2ESTIMATESFORSPECIFICVORTICITYANDENTROPYGRADIENT}--\ref{S:TOPORDERELIPTICHYPERBOLICL2ESTIMATESFORSPECIFICVORTICITYANDENTROPYGRADIENT}. We also note that to close the proof, energy estimates for 
the wave equations must be supplemented with elliptic-hyperbolic estimates, which we elaborate on in Sect.\ref{PAR:THENEEDFORUSINGCANDD}.

In \cite{abbrescia2022emergence}, our $L^2$ analysis was localized in the rough geometric coordinate region $[\timefunction_0,0]\times [-\rightu,\leftu]$, where $\leftu,\rightu > 0$ are constants that depended on the background solution. The $L^2$ estimates were based on coercive energies on truncated rough hypersurfaces $\hypthreearg{\timefunction}{[-\rightu,u']}{0} \eqdef \hyptwoarg{\timefunction}{0} \cap \{u \in [-\rightu,u']\}$, fluxes on truncated characteristics $\nullhypthreearg{0}{u}{[\timefunction_0,\timefunction']} \eqdef \nullhyparg{u}\cap\{\timefunctionarg{0} \in [\timefunction_0,\timefunction']\}$, and spacetime bulk integrals on $\twoargMrough{[\timefunction_0,\timefunction'],[-\rightu,u']}{0} \eqdef \bigcup_{\timefunction \in [\timefunction_0,\timefunction']} \hypthreearg{\timefunction}{[-\rightu,u']}{0}$. Denoting $\tander^N$ an $N$-th order string of $\nullhyparg{u}$-tangential differential operators and $\mathfrak{P}^{(N)}$ the set of all such operators, these can be schematically depicted as follows, where $\muX$ is an $\mathcal{O}(\upmu)$-size vectorfield that is $\Sigma_t$-tangent, transversal to the characteristics 
$\nullhypthreearg{0}{u}{[\timefunction_0,\timefunction]}$, and satisfies $\muX u = 1$
(see Def.\,\ref{D:COMVECTORFIELDS}):
	\begin{subequations} \label{E:ROUGHL2CONTROLINGQUANTITIES}
		\begin{align}
			\widetilde{\mathbf{E}}_N(\timefunction,u) & \simeq \sum_{\tander^N \in \mathfrak{P}^{(N)}} \int_{\hypthreearg{\timefunction}{[-\rightu,u]}{0}} \upmu \left|\Lunit \tander^N \wavearray\right|^2 + \sum_{A=2}^3  \upmu \left| \Yvf{A} \tander^N \wavearray \right|^2 + \left| \muX \tander^N \wavearray \right|^2, \label{E:INTRO:ROUGHHYPERSURFACEENERGIES} \\ 
			 \widetilde{\mathbf{F}}_N(\timefunction,u) & \simeq \sum_{\tander^N \in \mathfrak{P}^{(N)}} \int_{\nullhypthreearg{0}{u}{[\timefunction_0,\timefunction]}} \left|\Lunit \tander^N \wavearray\right|^2 +  
\sum_{A=2}^3  \upmu \left| \Yvf{A} \tander^N \wavearray \right|^2, \label{E:INTRO:CHARACTERISTICHYPERSURFACEFLUXES} \\
			\widetilde{\mathbf{K}}_N(\timefunction,u) & \simeq \sum_{\tander^N \in \mathfrak{P}^{(N)}} \int_{\twoargMrough{[\timefunction_0,\timefunction],[-\rightu,u]}{0} } \sum_{A=2}^3 \left|\Yvf{A}\tander^N \wavearray \right|^2, \label{E:INTRO:ROUGHBULKINTEGRAL} 
		\end{align}
	\end{subequations}
where we have suppressed the volume forms in \eqref{E:INTRO:ROUGHHYPERSURFACEENERGIES}--\eqref{E:INTRO:ROUGHBULKINTEGRAL}. We highlight the following two important facts:
	\begin{enumerate}
		\item The rough hypersurface energies $\widetilde{\mathbf{E}}_N(\timefunction,u)$ control $\Lunit \tander^N \wavearray$ in $L^2$, but with a weighted $\upmu$-factor that degenerates as the shock forms, i.e., as $\upmu \downarrow 0$. 
		They also control the angular derivatives $\Yvf{A} \tander^N \wavearray$ in $L^2$ with the same factor of $\upmu$.
		\item The null-fluxes $\widetilde{\mathbf{F}}_N(\timefunction,u)$ also control $\Lunit \tander^N \wavearray$ in $L^2$, but this coercivity is much \emph{stronger} than the one from energies on the rough hypersurfaces because there is no factor of $\upmu$. Similarly, the bulk spacetime integrals $\widetilde{\mathbf{K}}_N(\timefunction,u)$ provide a stronger control over the angular derivatives $\Yvf{A}\tander^N \wavearray$ than either the null-fluxes or the hypersurface energies.
	\end{enumerate}

The top order energy estimate derived in \cite{abbrescia2022emergence} is schematically of the following form:
\begin{align} \label{E:INTRO:SCHEMATICTOPORDERESTIMTEFORROUGHFOLIATIONS}
		\widetilde{\mathbf{E}}_{\Ntop}(\timefunction,u) +  \widetilde{\mathbf{F}}_{\Ntop}(\timefunction,u) + \widetilde{\mathbf{K}}_{\Ntop}(\timefunction,u)  \le C \mathring{\upepsilon}^2 + \boxed{A} \int_{\timefunction' = \timefunction_0}^{\timefunction} \frac{1}{|\timefunction'|} \widetilde{\mathbf{E}}_{\Ntop}(\timefunction',u) \, \mathrm{d} \timefunction' + \cdots,
	\end{align}
where $\mathring{\upepsilon} > 0$ is the small size of the perturbation of the initial data away from the background solution, and $\boxed{A} >0$ is a \emph{universal constant} that is independent of the equation state, and $\cdots$ denotes easier\footnote{In reality, there are other difficult error terms present in $\cdots$ that require a further integration by parts and are just as singular as $ \int_{\timefunction' = \timefunction_0}^{\timefunction} \frac{1}{|\timefunction'|} \widetilde{\mathbf{E}}_N(\timefunction',u) \, \mathrm{d} \timefunction'$. We have chosen to stick with the RHS\,\eqref{E:INTRO:SCHEMATICTOPORDERESTIMTEFORROUGHFOLIATIONS} for the sake of simplicity in the introduction.} or less singular terms. Applying Gr\"onwall's inequality to \eqref{E:INTRO:SCHEMATICTOPORDERESTIMTEFORROUGHFOLIATIONS} yields:
\begin{align} \label{E:INTRO:SCHEMATICTOPORDERESTIMTEFORROUGHFOLIATIONSAFTERGRONWALL}
		\widetilde{\mathbf{E}}_{\Ntop}(\timefunction,u) +  \widetilde{\mathbf{F}}_{\Ntop}(\timefunction,u) + \widetilde{\mathbf{K}}_{\Ntop}(\timefunction,u)  \lesssim \mathring{\upepsilon}^2 |\timefunction|^{-A}.
	\end{align}
We note that \eqref{E:INTRO:SCHEMATICTOPORDERESTIMTEFORROUGHFOLIATIONSAFTERGRONWALL} is \emph{singular} as $\timefunction \uparrow 0$, which reflects that $\upmu = 0$ on the crease, which is contained in $\hypthreearg{0}{[-\rightu,\leftu]}{0}$. \

The singular nature of the estimate \eqref{E:INTRO:SCHEMATICTOPORDERESTIMTEFORROUGHFOLIATIONSAFTERGRONWALL} begs the question of whether ``hiding'' the singularity with nonlinear geometric optics will lead to fruitful results. The answer to this is affirmative, and it lies in the differing methods used to derive the top-order estimate \eqref{E:INTRO:SCHEMATICTOPORDERESTIMTEFORROUGHFOLIATIONS} and the below-top-order estimate:
	\begin{align}
		\begin{split} \label{E:INTRO:SCHEMATICBELOWTOPORDERESTIMTEFORROUGHFOLIATIONS} 
			\widetilde{\mathbf{E}}_{\Ntop-1}(\timefunction,u) & +  \widetilde{\mathbf{F}}_{\Ntop-1}(\timefunction,u) + \widetilde{\mathbf{K}}_{\Ntop-1}(\timefunction,u) \\
			&  \le C \mathring{\upepsilon}^2 +  C\int_{\timefunction' = \timefunction_0}^{\timefunction} \frac{1}{|\timefunction'|^{1/2}} \widetilde{\mathbf{E}}_{\Ntop-1}(\timefunction',u) \int_{\timefunction'' = \timefunction_0}^{\timefunction'} \frac{1}{|\timefunction''|^{1/2}} \widetilde{\mathbf{E}}_{\Ntop}(\timefunction'',u) \, \rmd \timefunction'' \rmd \timefunction' + \cdots.
		\end{split}
	\end{align}
Although we do not elaborate on the details of the derivation, we note that the singular nature of \eqref{E:INTRO:SCHEMATICTOPORDERESTIMTEFORROUGHFOLIATIONS} is a consequence of using Raychaudhuri's equation (mentioned in Sect.\,\ref{SSS:ROUGHFOLIATIONSANDTIMEFUNCTION}) to prevent derivative-loss. In contrast, in deriving the below-top-order estimate \eqref{E:INTRO:SCHEMATICBELOWTOPORDERESTIMTEFORROUGHFOLIATIONS}, we instead choose to accept the loss of one derivative, which is reflected in the coupling of RHS~\eqref{E:INTRO:SCHEMATICBELOWTOPORDERESTIMTEFORROUGHFOLIATIONS}  to top-order energy $\widetilde{\mathbf{E}}_{\Ntop}$. 
Inserting \eqref{E:INTRO:SCHEMATICTOPORDERESTIMTEFORROUGHFOLIATIONSAFTERGRONWALL} into \eqref{E:INTRO:SCHEMATICBELOWTOPORDERESTIMTEFORROUGHFOLIATIONS} and applying 
Gr\"onwall's inequality, we obtain the following \emph{less singular estimate}: 
\begin{align} \label{E:INTRO:SCHEMATICBELOWTOPORDERESTIMTEFORROUGHFOLIATIONSAFTERGRONWALL}
		\widetilde{\mathbf{E}}_{\Ntop-1}(\timefunction,u) +  \widetilde{\mathbf{F}}_{\Ntop-1}(\timefunction,u) + \widetilde{\mathbf{K}}_{\Ntop-1}(\timefunction,u)  \lesssim \mathring{\upepsilon}^2 |\timefunction|^{-A+2}.
	\end{align}
This ``descent'' can be continued downward, with each estimate becoming less singular, e.g.,
	\begin{align} 
		\widetilde{\mathbf{E}}_{\Ntop-2}(\timefunction,u) +  \widetilde{\mathbf{F}}_{\Ntop-2}(\timefunction,u) + \widetilde{\mathbf{K}}_{\Ntop-2}(\timefunction,u)  & \lesssim \mathring{\upepsilon}^2 |\timefunction|^{-A+4}, \label{E:INTRO:SCHEMATICTWOBELOWTOPORDERESTIMTEFORROUGHFOLIATIONSAFTERGRONWALL} \\
		\widetilde{\mathbf{E}}_{\Ntop-3}(\timefunction,u) +  \widetilde{\mathbf{F}}_{\Ntop-3}(\timefunction,u) + \widetilde{\mathbf{K}}_{\Ntop-3}(\timefunction,u)  & \lesssim \mathring{\upepsilon}^2 |\timefunction|^{-A+6}.  \label{E:INTRO:SCHEMATICTHREEBELOWTOPORDERESTIMTEFORROUGHFOLIATIONSAFTERGRONWALL}
	\end{align} 
One continues the descent as in \eqref{E:INTRO:SCHEMATICBELOWTOPORDERESTIMTEFORROUGHFOLIATIONSAFTERGRONWALL}--\eqref{E:INTRO:SCHEMATICTHREEBELOWTOPORDERESTIMTEFORROUGHFOLIATIONSAFTERGRONWALL} until finally arriving at the \emph{non-singular estimate}:
\begin{align}
		\widetilde{\mathbf{E}}_{\Ntop-\frac{A}{2}}(\timefunction,u) +  \widetilde{\mathbf{F}}_{\Ntop-\frac{A}{2}}(\timefunction,u) + \widetilde{\mathbf{K}}_{\Ntop-\frac{A}{2}}(\timefunction,u)  & \lesssim \mathring{\upepsilon}^2. \label{E:INTRO:SCHEMATICNONSINGULARENERGYESTIMATEFORROUGHFOLIATIONSAFTERGRONWALL}
	\end{align}
	 
From \eqref{E:INTRO:SCHEMATICNONSINGULARENERGYESTIMATEFORROUGHFOLIATIONSAFTERGRONWALL} and Sobolev embedding, we finally conclude that the solution is bounded with respect to the rough geometric coordinates $(\timefunction,u,x^2,x^3)$. This is what we meant by ``rather'' regular in geometric or rough geometric coordinates: although the top order energies could potentially blow-up, at the low-derivative levels, the solution's $C^{\le \Ntop - \frac{A}{2} - 2}$-norms are finite when measuring regularity in the $\nullhyparg{u}$-directions. 
We finally recover the desired statements in the original background Cartesian coordinate space by proving quantitative estimates for the change of variables map $(\timefunction,u,x^2,x^3) \mapsto (t,u,x^2,x^3) \mapsto (t,x^1,x^2,x^3)$. In particular, we show that these coordinate transformations are diffeomorphisms everywhere on $\twoargMrough{[\timefunction_0,\timefunction],[-\rightu,\leftu]}{0}$, except at the crease $\crease$, where they are homeomorphisms.

\begin{remark}[Using $(t,u,x^3,x^3)$-adapted commutators in $(\timefunction,u,x^2,x^3)$-coordinates] \label{R:ROUGHCOORDINATESCOULDNOTBEUSEDASCOMMUTATORS}
	Since our estimates are done in the ($\timefunction,u,x^2,x^3$) coordinate system, one might wonder why we used commutators $\mathscr{Z} = \{\Lunit,\Yvf{2},\Yvf{3}\}$ adapted to the geometric coordinates $(t,u,x^2,x^3)$ instead of, perhaps, coordinate differentiations $\{\roughgeop{\timefunction},\roughgeop{x^2},\roughgeop{x^3}\}$ in the rough geometric coordinates $(\timefunction,u,x^2,x^3)$ that are also $\nullhyparg{u}$-tangent. The reason is that the regularity of the time function is constrained by the regularity of $\upmu$ by \eqref{E:DATAFORROUGHTIMEFUNCTION}, which depends on \emph{transversal} derivatives of the solution, i.e. $\muX \wavearray$. From Sobolev embedding and the coercivity of the $L^2$-controlling quantities, \eqref{E:ROUGHL2CONTROLINGQUANTITIES} implies that $\wavearray$ is smooth up to order $\Ntop - \frac{A}{2} -2$ \emph{in $\nullhyparg{u}$-tangential directions}. The regularity of $\wavearray$ in \emph{transverse} directions is low and limited to $C^{3,1}$. This in turn implies that the best estimate possible  for the rough time function is $\timefunctionarg{0} \in C^{2,1}$. This brings us back to the reason as to why we did not derive estimates relative to $\{\roughgeop{\timefunction},\roughgeop{x^2},\roughgeop{x^3}\}$: when expressed in $(t,u,x^2,x^3)$, it is easy to see that  $\roughgeop{\timefunction} = \smoothfunction(\tander \timefunctionarg{0}) \Lunit +\smoothfunction^A(\tander \timefunctionarg{0}) \Yvf{A}$  where $\smoothfunction$ is a smooth function (and similarly for $\roughgeop{x^A}$). Hence, Sobolev $L^2$ estimates derived from using $\{\roughgeop{\timefunction},\roughgeop{x^2},\roughgeop{x^3}\}$  as commutators in rough geometric coordinates would require controlling  $\Ntop - \frac{A}{2} -1$ derivatives of $\timefunctionarg{0}$. This is inconsistent with the $C^{2,1}$ regularity afforded in $(t,u,x^2,x^3)$, which is restricted fundamentally from the Cauchy data's regularity posed at $\{t=0\}$ in Cartesian coordinate space.
	\end{remark}
	
\begin{remark}[Gauss curvature estimates in rough geometric coordinates] \label{R:GAUSSCURVATUREOFROUGHTORI}
	One of the technical steps we omitted in deriving \eqref{E:INTRO:SCHEMATICTOPORDERESTIMTEFORROUGHFOLIATIONS} was an $L^\infty$ bound on the Gauss curvature of the tori $\roughtori{\timefunction,u} \eqdef \hypthreearg{\timefunction}{[-\rightu,\leftu]}{0} \cap \nullhypthreearg{0}{u}{[\timefunction_0,0]}$ induced by the acoustical metric $\gfour$. Since the $\timefunctionarg{0}\in C^{2,1}$, this bound was attainable. We highlight this because Gauss curvature estimates analogous to those from \cite{abbrescia2022emergence} are needed in this paper, but require new methods which we discuss in Sect.\,\ref{SS:NEWIDEASANDMETHODS}.
\end{remark}

\subsection{New ideas and methods used in the main proof} \label{SS:NEWIDEASANDMETHODS}
We now summarize the main new ideas and methods we introduce in this paper to derive an $\mathcal{O}(1)$-portion of a Cauchy horizon.

\subsubsection{Construction and regularity subtleties of the ingoing eikonal function} \label{SSS:CONSTRUCTIONOFEIKONALFUNCTION} 
The main way we control the fluid up to the Cauchy horizon is by introducing a new \emph{ingoing} eikonal function $\ubar$ and a corresponding \emph{double-null coordinate system} $(\ubar,u,x^2,x^3)$. The gauge for the ingoing eikonal function is precisely chosen such that the Cauchy horizon and the $\crease$ are subsets of the 0-level set of $\ubar$, i.e. $(\Cauchyhor \cup \crease )\subset \{ (t,u,x^2,x^3) \, : \, \ubar(t,u,x^2,x^3) = 0\}$. By setting $ \datahypfortimefunctiontwoarg{0}{[\timefunction_0,0]}\eqdef \{ \muX\upmu = 0\} \cap\{ \timefunctionarg{0} \in [\timefunction_0,0]\}$, we consider the following fully non-linear Cauchy initial value problem:
\begin{subequations} \label{E:INTRO:INGOINGEIKONALEQUATIONINITIALVALUEPROBLEM}
		\begin{align}
			(\gfour^{-1})^{\alpha\beta} \p_\alpha  \ubar \p_\beta \ubar &  = 0, \label{E:INTROINGOINGEIKONALEQUATION} \\
			\Lunit \ubar & >  0, \label{E:INTROTRANSVERSALITYCONDITIONFORINGOINGEIKONALWITHRESPECTOU} \\
			\ubar|_{\datahypfortimefunctiontwoarg{0}{[\timefunction_0,0]}} & = -\upmu|_{\datahypfortimefunctiontwoarg{0}{[\timefunction_0,0]}}. \label{E:INTRODATAFORINGOINGEIKONALINITIALVALUEPROBLEM}
		\end{align}
	\end{subequations}
To solve for $\ubar$, we derive a second order quasilinear system of \emph{transport} equations for the derivatives of the eikonal function.  Let  $\newuL$ be the $\gfour$-null vector field proportional to the $\gfour$-gradient vectorfield $\uLgeo^\alpha = (\gfour^{-1})^{\alpha\beta} \p_\beta \ubar$ of $\ubar$ normalized by $\newuL u = 1$. Then, for any $Z \in \{\Lunit,\muX,\Yvf{2},\Yvf{3}\}$, 
upon Lie differentiating \eqref{E:INTROINGOINGEIKONALEQUATION},
it is not difficult to derive the following:
\begin{align}
	\begin{split}
		 \newuL ( Z \ubar) & = \frac{1}{2 \gfour^{-1}(\mathrm{d} u,\mathrm{d} \ubar)} (\Lie_Z\gfour^{-1})(\mathrm{d} \ubar,\mathrm{d}\ubar). \label{E:INTRODIFFERENTIATEDEIKONALEQUATIONGENERAL}
	\end{split}
\end{align}Now, to solve \eqref{E:INTRODIFFERENTIATEDEIKONALEQUATIONGENERAL}, one must  derive the initial data for $Z \ubar|_{\datahypfortimefunctiontwoarg{0}{[\timefunction_0,0]}}$ and ensure that $\newuL$ is transverse to the initial data surface $\datahypfortimefunctiontwoarg{0}{[\timefunction_0,0]}$. Since $\ubar = -\upmu$ on $\datahypfortimefunctiontwoarg{0}{[\timefunction_0,0]}$, all derivatives \emph{tangent} to $\datahypfortimefunctiontwoarg{0}{[\timefunction_0,0]}$ are equal to the negative of those of $\upmu$, while \emph{transversal} derivatives along $\datahypfortimefunctiontwoarg{0}{[\timefunction_0,0]}$ are solved for using the eikonal equation \eqref{E:INTROINGOINGEIKONALEQUATION}. Solving for these derivatives, one finds that the initial data for $Z\ubar$ in terms of $\upmu$ is schematically as follows:
	\begin{subequations} \label{E:INTROSCHEMADICDATAFORZUBAR}
		\begin{align}
			\Lunit \ubar \big|_{\datahypfortimefunctiontwoarg{0}{[\timefunction_0,0]}} & = - \left(1 + \upmu\smoothfunction(Z^2\upmu)\right)\Lunit \upmu +  \mathcal{O}(\mathring{\upepsilon})  \smoothfunction(Z^{\le 2}  \upmu)  \big|_{\datahypfortimefunctiontwoarg{0}{[\timefunction_0,0]}}, \label{E:INTRODATAFORINGOINGLUBAR} \\
			\muX \ubar \big|_{\datahypfortimefunctiontwoarg{0}{[\timefunction_0,0]}} & = \upmu\left(1 + \upmu \smoothfunction(Z^2\upmu)\right)\Lunit \upmu  +\upmu \mathcal{O}(\mathring{\upepsilon})  \smoothfunction(Z^{\le 2}  \upmu) \big|_{\datahypfortimefunctiontwoarg{0}{[\timefunction_0,0]}}, \label{E:INTRODATAFORINGOINGMUXUBAR} \\
			\Yvf{A} \ubar \big|_{\datahypfortimefunctiontwoarg{0}{[\timefunction_0,0]}}  & = - \Yvf{A} \upmu +  \mathcal{O}(\mathring{\upepsilon})  \smoothfunction(Z^{\le 2}\upmu) \big|_{\datahypfortimefunctiontwoarg{0}{[\timefunction_0,0]}};\label{E:INTRODATAFORINGOINGYUBAR}
		\end{align}
	\end{subequations}
see Lemmas\,\ref{L:IDENTITYFORLUBARANDMUXUBARINTERMSOFDATASURFACETANGETDERIVATIVESSOVLEDBYEIKONALEQUAITON} and\,\ref{L:DATAREGULARITYFORINGOINGEIKONALFUNCTION}, Remark\,\ref{R:LEADINGORDERFORUBAR}, and equations \eqref{E:IDENTITYFORLAPXUBARINTERMSOFDATASURFACEDERIVATIVES}--\eqref{E:IDENTITYFORMUXAPXUBARINTERMSOFDATASURFACEDERIVATIVES}. The reason one should expect two derivatives of $\upmu$ on the RHS\,\eqref{E:INTROSCHEMADICDATAFORZUBAR} is that the initial data set $\datahypfortimefunctiontwoarg{0}{[\timefunction_0,0]}$ already depends on \emph{one} derivative of $\upmu$. For example, to solve for $\Lunit \ubar|_{\datahypfortimefunctiontwoarg{0}{[\timefunction_0,0]}}$, we differentiate the data $\ubar|_{\datahypfortimefunctiontwoarg{0}{[\timefunction_0,0]}}   = -\upmu|_{\datahypfortimefunctiontwoarg{0}{[\timefunction_0,0]}}$ with respect to the $\datahypfortimefunctiontwoarg{0}{[\timefunction_0,0]}$-tangent vector field $\Lunit - \frac{\Lunit \muX \upmu}{\muX \muX \upmu} \muX$ to find $\Lunit \ubar - \frac{\Lunit \muX \upmu}{\muX \muX \upmu} \muX \ubar |_{\datahypfortimefunctiontwoarg{0}{[\timefunction_0,0]}} = \Lunit \upmu - \frac{\Lunit \muX \upmu}{\muX \muX \upmu} \muX \upmu |_{\datahypfortimefunctiontwoarg{0}{[\timefunction_0,0]}} $. We then use the eikonal equation \eqref{E:INTROINGOINGEIKONALEQUATION} to solve for $\muX \ubar$ in terms of $\Lunit\ubar$ and $\datahypfortimefunctiontwoarg{0}{[\timefunction_0,0]}$-tangential derivatives of $\ubar$ (the latter of which can be expressed purely in terms of $\upmu$). The denominator $\muX\muX\upmu$ is always bounded away from zero and infinity near the crease due to transversal convexity. Similarly, one can prove that $\newuL \muX \upmu |_{\crease} =\muX\muX\upmu |_{\crease}$, which is a necessary condition to prove that the system \eqref{E:INTROSCHEMADICDATAFORZUBAR} is well-posed because it shows that $\newuL$ is transverse to the initial data surface, even at the crease.

We now highlight several crucial consequences of \eqref{E:INTROSCHEMADICDATAFORZUBAR} 
that are not immediately apparent. First, from \eqref{E:INTRODATAFORINGOINGMUXUBAR}, it follows that\footnote{The following bound holds everywhere in the spacetime regions we study: $\Lunit \upmu \approx -1$.} $X \ubar|_{\crease} \approx -1$ and so $\ubar$ morally behaves like $x^1$ at the crease. From \eqref{E:INTRODATAFORINGOINGLUBAR}--\eqref{E:INTRODATAFORINGOINGYUBAR} and the eikonal equation, one can also compute that the inverse foliation density $\ingoingmu \approx |\partial_x \ubar|^{-1}$ for the \emph{ingoing} characteristics $\ingoingcharacteristicsurfacearg{\ubar'} \eqdef \{ (t,u,x^2,x^3) \, : \, \ubar(t,u,x^2,x^3) = \ubar'\}$ satisfies  $\ingoingmu \approx 1$. These crucial estimates imply that the level sets of $\ubar$ do not collapse themselves to form a shock in the ingoing characteristic directions.

The data estimates \eqref{E:INTROSCHEMADICDATAFORZUBAR} highlight a striking difference between using $\ubar$ as a double-null coordinate function and using $\timefunctionarg{0}$ as a coordinate time function. As we mentioned in Remark\,\ref{R:ROUGHCOORDINATESCOULDNOTBEUSEDASCOMMUTATORS}, $\timefunctionarg{0} \in C^{2,1}$ because it was just as regular as $\upmu \in C^{2,1}$. The reason for this is that \eqref{E:DATAFORROUGHTIMEFUNCTION} implies that $\timefunctionarg{0}$ is as smooth as $\upmu$ in directions \emph{tangential} to $\datahypfortimefunctiontwoarg{0}{[\timefunction_0,0]}$, while \eqref{E:EQUATIONFORROUGHTIMEFUNCTION} (which reads $\muX \timefunctionarg{0} = 0$) implies that it is much smoother in directions transverse to $\datahypfortimefunctiontwoarg{0}{[\timefunction_0,0]}$. Contrastingly, even though $\timefunctionarg{0}|_{\datahypfortimefunctiontwoarg{0}{[\timefunction_0,0]}} = \ubar|_{\datahypfortimefunctiontwoarg{0}{[\timefunction_0,0]}} = -\upmu |_{\datahypfortimefunctiontwoarg{0}{[\timefunction_0,0]}}$ on the initial data hypersurface, the eikonal equation \eqref{E:INTROINGOINGEIKONALEQUATION} forces there to be a \emph{loss} of derivatives in the data estimates \eqref{E:INTROSCHEMADICDATAFORZUBAR} and $\ubar \in C^{1,1}$ is the most regularity allowable from the initial conditions. Thus, for an even stronger reason than in Remark\,\ref{R:ROUGHCOORDINATESCOULDNOTBEUSEDASCOMMUTATORS}, the double-null coordinate vector fields $\{\nullgeop{\ubar},\nullgeop{u},\nullgeop{x^2},\nullgeop{x^3}\}$ cannot be used as commutators when deriving energy estimates. For this reason, in this paper, we \emph{continue} to use the same commutators as in \cite{abbrescia2022emergence}: $\{\Lunit,\Yvf{2},\Yvf{3}\}$.

\subsubsection{Raychaudhuri's equation and cancellations in Gauss curvature estimates to avoid derivative-loss in the acoustic geometry} \label{SSS:GAUSSCURVATUREESTIMATESNOTREQUIREDIN2D} 
We now explain why using $\{\Lunit,\Yvf{2},\Yvf{3}\}$ as commutators poses novel difficulties in the current paper compared to all prior $3D$ shock formation works \cite{dC2007,jS2016b,LukSpeck2024stability}, as well as why these difficulties are not present in the simpler $2D$ setting of \cite{luk2018shock,shkoller2024geometry}. For the same reasons as in Remark\,\ref{R:GAUSSCURVATUREOFROUGHTORI}, which we expand on below, using  $\{\Lunit,\Yvf{2},\Yvf{3}\}$ as commutators requires an estimate on the Gauss curvature of the \emph{double-null tori} $\doublenulltoritwoarg{\ubar}{u}\eqdef \ingoingcharacteristicsurfacearg{\ubar}\cap \nullhyparg{u}$, where $\ingoingcharacteristicsurfacearg{\ubar'} := \{ (t,u,x^2,x^3) : \ubar(t,u,x^2,x^3) = \ubar'\}$ denotes the \emph{ingoing} characteristics. A priori, this estimate requires controlling \emph{three} derivatives of $\ubar$. The sharp regularity $\ubar \in C^{1,1}$ discussed in Sect.\,\ref{SSS:CONSTRUCTIONOFEIKONALFUNCTION} thus brings into question whether these double-null coordinates are compatible with using $\{\Lunit,\Yvf{2},\Yvf{3}\}$ as commutators. This issue is completely resolved in Sect.\,\ref{SSS:GAUSSCURVATUREOFROUGHTORI}, where we prove that there is a \emph{complete} cancellation of \emph{all} thrice differentiated factors of $\ubar$ when expanding the Gauss curvature in the intrinsic coordinates on $\doublenulltoritwoarg{\ubar}{u}$, leaving at most twice differentiated factors of $\ubar$. These terms can be controlled since $\ubar \in C^{1,1}$ by Rademacher's theorem.

As we mentioned above, the Gauss curvature estimates from the previous paragraph, as well as the ones stated in Remark\,\ref{R:GAUSSCURVATUREOFROUGHTORI} that were proved in \cite{abbrescia2022emergence}, \emph{would not be needed if our results were restricted to two spatial dimensions}.  In particular, the $2D$ works \cite{luk2018shock,shkoller2024geometry} did not require them. The fundamental reason is that symmetric $\binom{0}{2}$-tensors in two dimensions have a pure-trace part and a trace-free part, while they are \emph{always} pure-trace in one dimension. To see why this is the reason, we now elaborate on how Raychaudhuri's equation is used to avoid derivative-loss in the commutator method of this paper as well as \cite{abbrescia2022emergence}.

Denoting \eqref{E:COVWAVEEQv}--\eqref{E:COVWAVEEQS} schematically as\footnote{We always commute the covariant wave equations \eqref{E:COVWAVEEQv}--\eqref{E:COVWAVEEQS} \emph{after} multiplying them by $\upmu$. The reason for this is that for any $P  \in\{\Lunit,\Yvf{2},\Yvf{3}\}$, \eqref{E:SCHEMATICGRADIENTBLOWUP} implies $[P,\partial_{x^\alpha}]\sim \frac{1}{\upmu}\partial_{x^\alpha}$ and so commuting $P$ with $\Box_{\gfour} \sim \gfour^{\alpha\beta} \partial_{\alpha}\partial_{\beta}$ would lead to uncontrollable powers of $\frac{1}{\upmu}$. By multiplying $\Box_{\gfour}$ by $\upmu$, we can instead express $\upmu \Box_{\gfour} \sim Z^2$ where $Z \in \{\Lunit,\muX,\Yvf{2},\Yvf{3}\}$. We prove in Lemma\,\ref{L:SIMPLECOMMUTATORIDENTITY} the much better schematic relation $[P,Z]\sim \Yvf{A}$.}  $\upmu \Box_{\gfour} \wavearray = \vec{F}$ (where $\vec{F}$ denotes $\upmu$ times the RHS\,\eqref{E:COVWAVEEQv}--\eqref{E:COVWAVEEQS}), we have:
\begin{align} \label{E:SCHEMATICCOMMUTEDWAVEEQUATION}
		\upmu \Box_{\gfour}\tander^{\Ntop} \wavearray =  \tander^{\Ntop} \mytr_{\gtorus}\upchi 
		+ \tander^{\Ntop} \left(\vec{F}\right) + \cdots.
	\end{align}
In \eqref{E:SCHEMATICCOMMUTEDWAVEEQUATION},  $\gtorus$ is the Riemannian metric on $\ell_{t,u} \eqdef \Sigma_t \cap \nullhyparg{u}$ induced by the acoustical metric $\gfour$, $\upchi = \frac{1}{2} \angLie_\Lunit \gtorus$ is the null second fundamental form of the outgoing null hypersurfaces, and $ \mytr_{\gtorus}\upchi \eqdef (\gtorus^{-1})^{\alpha\beta}\upchi_{\alpha\beta}$ is the trace of $\upchi$ with respect to $\gtorus$.
It turns out that the term $ \tander^{\Ntop} \mytr_{\gtorus}\upchi $, 
is, at first glance, uncontrollable from the point of view of regularity. Since $\upchi \sim \pmb{\p}^2 u \sim \pmb{\p}^2 \wavearray$, it follows that $\tander^{\Ntop} \mytr_{\gtorus}\upchi$ 
 will depend on $\Ntop+2$ derivatives of $\wavearray$ and therefore cannot be treated as a source term in the energy estimates. Indeed, for solutions to \eqref{E:SCHEMATICCOMMUTEDWAVEEQUATION}, Raychaudhuri's equation (which is the governing evolution equation for $\mytr_{\gtorus}\upchi$) roughly implies: 
\begin{align} 
		\Lunit \left(\upmu \tander^{\Ntop} \mytr_{\gtorus}\upchi\right) = \Lunit\left( \muX  \tander^{\Ntop} \wavearray + \tander^{\Ntop+1}\wavearray\right) + \tander^N\left( |\upchi |_{\gtorus}^2\right) + \tander^N\left(\vec{F}\right) \cdots , \label{E:SCHEMATICRAYCHAUDHURIINTRO}
	\end{align}
where $|\upxi|^2_{\gtorus} \eqdef \gtorus^{\alpha\beta}\gtorus^{\delta\gamma}\upxi_{\alpha\delta}\upxi_{\beta\gamma}$ is the square norm of a symmetric type $\binom{0}{2}$ $\ell_{t,u}$-tangent tensor. In particular, \eqref{E:SCHEMATICRAYCHAUDHURIINTRO} shows that the evolution of $\tander^{\Ntop} \mytr_{\gtorus}\upchi$ does depend on $\Ntop+2$ derivatives of $\wavearray$. However, there is an extra structure\footnote{The reason for these structures is that, even without knowing that $\wavearray$ solves $\upmu \Box_{\gfour} \wavearray = \vec{F}$, Raychaudhuri's equation schematically reads  $\Lunit (\upmu \mytr_{\gtorus}\upchi) \simeq \upmu \Ricfour(\Lunit,\Lunit) +\upmu |\upchi|_{\gtorus}^2$, where $\Ricfour(\Lunit,\Lunit)$ is the double contraction of the Ricci curvature of the acoustical metric $\gfour$ and $\Lunit$. Straightforward calculations show that $\upmu \Ricfour(\Lunit,\Lunit) = \Lunit \tander\wavearray + \upmu \angLap \wavearray + \cdots$, where $\angLap$ is the Laplace--Beltrami operator with respect to $\gtorus$ and $\cdots$ are lower order terms. This is where solving $\upmu \Box_{\gfour} \wavearray = \vec{F}$ comes in: one can decompose $\upmu \Box_{\gfour}$ with respect to the frame $\{\Lunit,\muX,\Yvf{2},\Yvf{3}\}$ as $\upmu \Box_{\gfour}\wavearray = \Lunit (\muX \wavearray  + \Lunit \wavearray) + \upmu \angLap \wavearray + \cdots$ and solve for $\upmu\angLap \wavearray$ in terms of factors with $\Lunit$ as its outer-most derivative operator and $\vec{F}$.} in that \emph{$\Lunit$ is the outer-most derivative operator acting on $\Ntop+1$ derivatives of $\wavearray$!} We can then define a new \emph{modified} top-order acoustic variable as $\fullymodquant{\tander^N}\eqdef \upmu \tander^{\Ntop} \mytr_{\gtorus}\upchi - \muX \tander^{\Ntop}\wavearray - \tander^{\Ntop+1}\wavearray$, which solves: 
\begin{align} 
		\Lunit \fullymodquant{\tander^{\Ntop}} 
		& = \tander^{\Ntop}\left( |\upchi |_{\gtorus}^2\right) + \tander^{\Ntop}\left(\vec{F}\right) \cdots. 
		\label{E:SCHEMATICRAYCHAUDHURIFORMODIFIEDCHIINTRO}
	\end{align}
One can then re-write \eqref{E:SCHEMATICCOMMUTEDWAVEEQUATION} as:
\begin{align} \label{E:SCHEMATICCOMMUTEDWAVEEQUATIONWITHMODIFIEDCHI}
		\upmu \Box_{\gfour}\tander^{\Ntop} \wavearray = \frac{1}{\upmu}\, \fullymodquant{\tander^N} + \frac{1}{\upmu} \left( \muX \tander^{\Ntop}\wavearray + \tander^{\Ntop+1}\wavearray\right) + \tander^{\Ntop} \left( \vec{F}\right) + \cdots
	\end{align}Now, \emph{if} one had sufficient $L^2$ control over RHS\,\eqref{E:SCHEMATICRAYCHAUDHURIFORMODIFIEDCHIINTRO}, 
then there would no longer be a loss of derivatives 
coming from the term $\frac{1}{\upmu}\, \fullymodquant{\tander^N}$ on~RHS\,\eqref{E:SCHEMATICCOMMUTEDWAVEEQUATIONWITHMODIFIEDCHI}, as $\fullymodquant{\tander^N}$ would be a controllable error term from the point of view of regularity. 
We can finally state the key difficulty in controlling RHS\,\eqref{E:SCHEMATICRAYCHAUDHURIFORMODIFIEDCHIINTRO}: the Leibniz rule implies that $\tander^{\Ntop} |\upchi|_{\gtorus}^2$ features top order derivatives of \emph{the entire tensor $\upchi$}, not just the pure-trace part $\tander^{\Ntop}\mytr_{\gtorus}\upchi$. Unfortunately, Raychaudhuri's equation only governs the evolution of the pure-trace part $\tander^{\Ntop}\mytr_{\gtorus}\upchi$. The standard way to overcome these issues is to derive elliptic estimates. Given a closed $2$-dimensional manifold $(M,h)$ and a symmetric type $\binom{0}{2}$-tensorfield $\upxi$, it is well known that one can derive elliptic $L^2$-estimates which bound the full tensor $\nabla \upxi$ (where $\nabla$ is the Levi-Civita connection of $h$) in terms of up to one derivative of the trace $\mytr_h\upxi$ \emph{if one has control of the Gauss curvature of $M$}. In the setting of this paper, the correct $2$-dimensional manifolds to consider are the double-null tori  $\doublenulltoritwoarg{\ubar}{u}$ because understanding the solution on them will give a precise understanding of the solution on the ingoing-null hypersurface $\ingoingcharacteristicsurfacearg{\ubar} = \cup_u \doublenulltoritwoarg{\ubar}{u}$. Recall that the ultimate goal of this paper is to understand the solution on the Cauchy horizon, which is a strict subset of $\ingoingcharacteristicsurfacearg{0}$. All in all, \emph{if} we are able to control the Gauss curvature of $\doublenulltoritwoarg{\ubar}{u}$, we would be able to derive elliptic estimates for $\doublenulltoritwoarg{\ubar}{u}$-tangent $\binom{0}{2}$-type tensorfields. Although we are able to control the Gauss curvature using the novel cancellations discussed above, $\upchi$ is \emph{not} $\doublenulltoritwoarg{\ubar}{u}$-tangent! To close argument, we prove that the difference between derivatives of $\upchi$ and derivatives of its projection onto $\doublenulltoritwoarg{\ubar}{u}$ are purely $\Lunit$-derivatives of $\upchi$,  which are better. This must be done while being mindful about crucial factors of $\upmu$, whose vanishing signifies the gradient singularity at the crease. This delicate analysis is located in Sect.\,\ref{S:ELLIPTICESTIMATESFORTHEACOUSTICGEOMETRYONTHEDOUBLENULLTORI}, where we construct two distinct frames on $\nullhyparg{u}$, one adapted to the acoustic tori $\ell_{t,u}$, and one adapted to the double-null tori $\doublenulltoritwoarg{\ubar}{u}$. We control the relationship between the two frames as a key step in obtaining the elliptic estimates.

There is still one remaining source of potential derivative-loss hiding in the terms $\tander^{\Ntop}(\vec{F})$ 
in equations \eqref{E:SCHEMATICRAYCHAUDHURIFORMODIFIEDCHIINTRO} and \eqref{E:SCHEMATICCOMMUTEDWAVEEQUATIONWITHMODIFIEDCHI},
tied to the difficulty in estimating the top-order derivatives of the specific vorticity $\vortrenormalized$ 
and entropy gradient $\GradEnt$. In Sect.\,\ref{SSS:INTROENERGYESTIMATES}, we describe how we overcome these difficulties.

We conclude this section by re-emphasizing that in two spatial dimensions as in, $\ell_{t,u}$ is one dimensional and hence $\upchi = (\mytr_{\gtorus} \upchi) \gtorus$. This completely removes the need for deriving elliptic estimates to control the RHS\,\eqref{E:SCHEMATICRAYCHAUDHURIFORMODIFIEDCHIINTRO}.

\subsubsection{Overview of the $L^2$-analysis and its novel difficulties} \label{SSS:INTROENERGYESTIMATES}

There are many new aspects of the $L^2$ analysis present in this paper compared to all other prior results on shock formation. Since the gauge choice of $\ubar$ is made so that the Cauchy horizon is a subset of the ingoing characteristic surface $\ingoingcharacteristicsurfacearg{0}$, our energy estimates are done on a characteristic diamond foliated by truncated null hypersurfaces. The relevant spacetime sets are:
	\begin{align} \label{E:INTROSPACETIMEREGIONS}
		\outgoingcharacteristicsurfacetwoarg{u}{[\leftubar,0]} \eqdef \bigcup_{\ubar \in [\leftubar,0]} \doublenulltoritwoarg{\ubar}{u}, \qquad \ingoingcharacteristicsurfacetwoarg{\ubar}{[\moreinterestingu_1,\moreinterestingu_2]}  \eqdef \bigcup_{u \in [-\rightu,\leftu]}  \doublenulltoritwoarg{\ubar}{u}, \qquad 		\characteristicdiamondtwoarg{[\leftubar,0]}{[\moreinterestingu_1,\moreinterestingu_2]} = \bigcup_{\ubar \in [\leftubar,0]} \ingoingcharacteristicsurfacetwoarg{\ubar}{[\moreinterestingu_1,\moreinterestingu_2]} ,
	\end{align}
where $\leftubar \in [\tfrac{3}{4} \timefunction_0,\frac{1}{4}\timefunction_0]$, and $-\rightu < \moreinterestingu_1 < 0 < \moreinterestingu_2<U_2$ are  chosen so that all of the singular behavior near the Cauchy horizon is contained in $u \in [\moreinterestingu_1, \moreinterestingu_2]$, e.g. $\crease \subset \{ u \in [\moreinterestingu_1,\moreinterestingu_2]\}$. We point out that our actual analysis is done on a bootstrap sub-region of \eqref{E:INTROSPACETIMEREGIONS} where $\ubar \in [\leftubar,\ubarboot)$ with $\ubarboot < 0$, and extended up to and including $\ubar = 0$ through a continuation criterion in Sect.\,\ref{S:EXISTENCEUPTOCAUCHYHORIZONBYCONTINUATIONCRITERIA}.

\begin{center}
	\begin{figure}  
		\begin{overpic}[scale=.7, grid = false, tics=5]{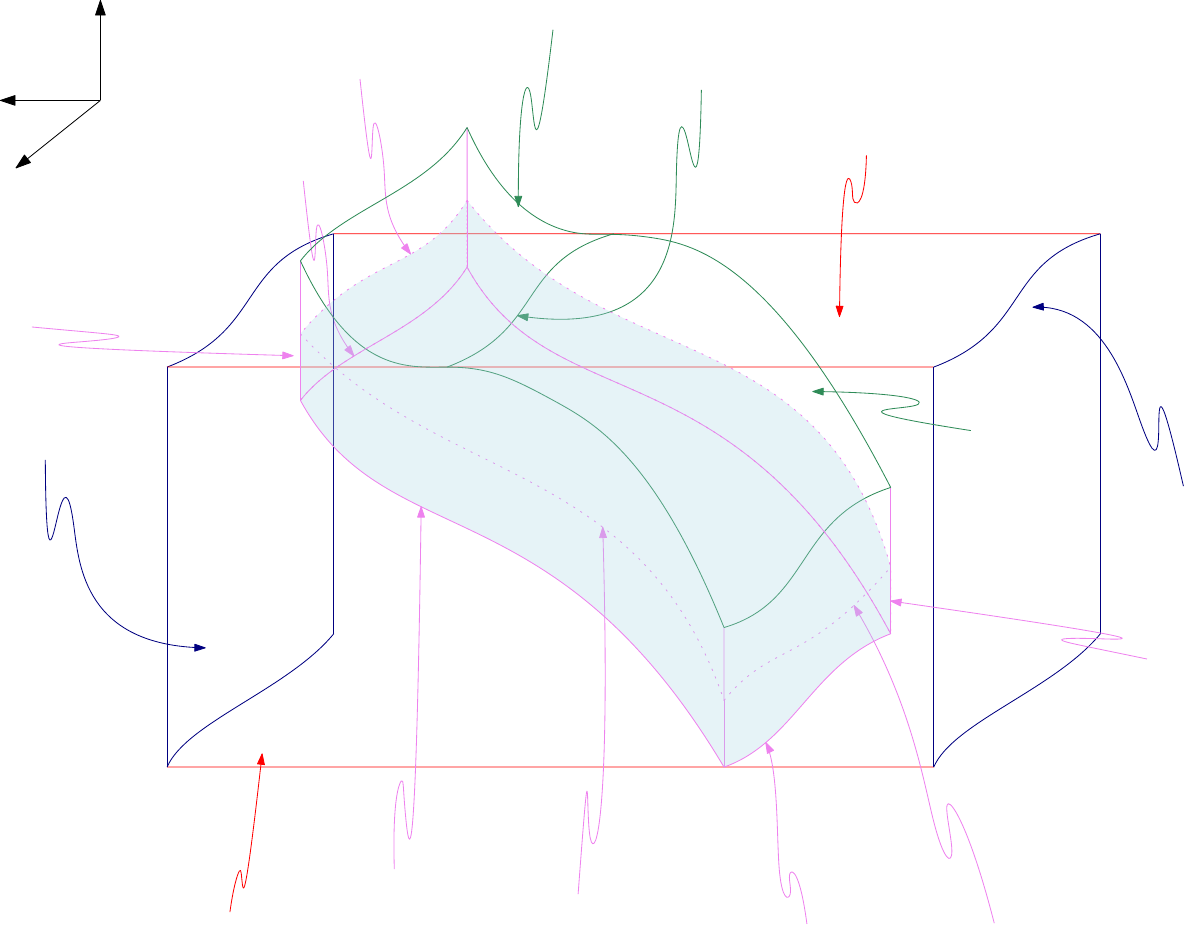}
			\put (55,72) {$\crease = \twoargmumuxtorus{0}{0}$}
			\put (68,66) {$\hypthreearg{0}{[-\interestingu,\interestingu]}{0}$}
			\put (47,76) {$\Cauchyhor$}
			\put (83,40) {$\ingoingcharacteristicsurfacetwoarg{0}{[\moreinterestingu_1,\moreinterestingu_2]}$}
			\put (95.5,34) {$\nullhypthreearg{0}{-\interestingu}{[\timefunction_0,0]}$}
			\put (-2.5,51.5) {$\outgoingcharacteristicsurfacetwoarg{\moreinterestingu_2}{[\leftubar,0]}$}
			\put (98,21) {$\outgoingcharacteristicsurfacetwoarg{\moreinterestingu_1}{[\leftubar,0]}$}
			\put (82,-1) {$\doublenulltoritwoarg{\ubarboot}{\moreinterestingu_1}$}
			\put (66,-1) {$\doublenulltoritwoarg{\leftubar}{\moreinterestingu_1}$}
			\put (46,-1) {$\ingoingcharacteristicsurfacetwoarg{\ubarboot}{[\moreinterestingu_1,\moreinterestingu_2]}$}
			\put (30,1.5) {$\ingoingcharacteristicsurfacetwoarg{\leftubar}{[\moreinterestingu_1,\moreinterestingu_2]}$}
			\put (22,65.5) {$\doublenulltoritwoarg{\leftubar}{\moreinterestingu_2}$}
			\put (27,74.5) {$\doublenulltoritwoarg{\ubarboot}{\moreinterestingu_2}$}
			\put (-1.5,42) {$\nullhypthreearg{0}{\interestingu}{[\timefunction_0,0]} $}
			\put (13,-2) {$\hypthreearg{\timefunction_0}{[-\interestingu,\interestingu]}{0}$}
			\put (4,64.5) {$(x^2,x^3) \in \mathbb{T}^2$}
			\put (9.5,78) {$t$}
			\put (-1,72) {$u \in \mathbb{R}$}
		\end{overpic}
		\vspace{0.5cm}
		\caption{The shaded characteristic diamond is the bootstrap region $\characteristicdiamondtwoarg{[\leftubar,\ubarboot)}{[\moreinterestingu_1,\moreinterestingu_2]}$. This region is extended to the larger characteristic diamond $\characteristicdiamondtwoarg{[\leftubar,0]}{[\moreinterestingu_1,\moreinterestingu_2]}$, which contains $\ingoingcharacteristicsurfacetwoarg{0}{[\moreinterestingu_1,\moreinterestingu_2]}$, through a continuation criterion. The portion of the Cauchy horizon $\Cauchyhor$ we construct is the region above the crease $\crease$ within $\ingoingcharacteristicsurfacetwoarg{0}{[\moreinterestingu_1,\moreinterestingu_2]}$.}
		\label{F:BOOTSTRAPDOMAINFORINGOINGEIKONALFUNCTIONANDTHECAUCHYHORIZON}
	\end{figure}
\end{center}

The top-order $L^2$-controlling quantities for the modified fluid variables are schematically as follows:
	\begin{subequations} \label{E:INTROVELOCITYL2CONTROLLINGQUANTITIES}
		\begin{align} 
			\ingoingfluxwave{\Ntop}[V](\ubar,u) & = \sum_{i=1}^3 \sum_{\tander^{\Ntop} \in \mathfrak{P}^{(\Ntop)}}\int_{\ingoingcharacteristicsurfacetwoarg{\ubar}{[\moreinterestingu_1,u]}} \left\{ \left| \newuL \tander^{\Ntop} v^i\right|^2 +\sum_{B=2}^3 \upmu \left| \nullgeop{x^B} \tander^{\Ntop}  v^i \right|^2 \right\} |\ubar|^A, \label{E:INTROINGOINGVELOCITYFLUX} \\
			\outgoingfluxwave{\Ntop}[V](\ubar,u) & = \sum_{i=1}^3 \sum_{\tander^{\Ntop} \in \mathfrak{P}^{(\Ntop)}}\int_{\outgoingcharacteristicsurfacetwoarg{u}{[\leftubar,\ubar]}} \left\{ \left| \Lunit  \tander^{\Ntop} v^i\right|^2 +\sum_{B=2}^3 \upmu \left| \Yvf{B} \tander^{\Ntop}  v^i \right|^2 \right\} |\ubar|^A, \label{E:INTROOUTGOINGVELOCITYFLUX} \\
			\newspacetimecoercive{\Ntop}[V](\ubar,u) & = \sum_{i=1}^3 \sum_{\tander^{\Ntop} \in \mathfrak{P}^{(\Ntop)}}\int_{\characteristicdiamondtwoarg{[\leftubar,\ubar]}{[\moreinterestingu_1,u]}} \left\{ \left| \newuL \tander^{\Ntop} v^i\right|^2 +\sum_{B=2}^3 \upmu \left| \nullgeop{x^B} \tander^{\Ntop}  v^i \right|^2 \right\} |\ubar|^{A-1},  \label{E:INTROWEAKCOERCIVEVELOCITYBULK} \\
			\spacetimecoercive{\Ntop}[V](\ubar,u) & = \sum_{i=1}^3 \sum_{\tander^{\Ntop} \in \mathfrak{P}^{(\Ntop)}}\int_{\characteristicdiamondtwoarg{[\leftubar,\ubar]}{[\moreinterestingu_1,u]}} \sum_{B=1}^3 \left| \Yvf{B} \tander^{\Ntop} v^i \right|^2 |\ubar|^{A},  \label{E:INTROSTRONGCOERCIVEVELOCITYBULK}
		\end{align}
	\end{subequations}
	\begin{subequations}  \label{E:INTROVORTANDGRADENTL2CONTROLLINGQUANTITIES}
		\begin{align} 
			\ingoingfluxwave{\Ntop}[(\vortrenormalized,\GradEnt)](\ubar,u) & = \sum_{\tander^{\Ntop} \in \mathfrak{P}^{(\Ntop)}}\int_{\ingoingcharacteristicsurfacetwoarg{\ubar}{[\moreinterestingu_1,u]}} \upmu \left|  \tander^{\Ntop} (\vortrenormalized,\GradEnt) \right|^2 |\ubar|^{A-1}, \label{E:INTROINGOINGVORTANDGRADENTFLUX} \\		
			\outgoingfluxwave{\Ntop}[(\vortrenormalized,\GradEnt)](\ubar,u) & = \sum_{\tander^{\Ntop} \in \mathfrak{P}^{(\Ntop)}}\int_{\outgoingcharacteristicsurfacetwoarg{u}{[\leftubar,\ubar]}}  \left|  \tander^{\Ntop} (\vortrenormalized,\GradEnt) \right|^2 |\ubar|^{A-1}, \label{E:INTROOUTGOINGVORTANDGRADENTFLUX} \\
			\newspacetimecoercive{\Ntop}[(\vortrenormalized,\GradEnt)](\ubar,u) & = \sum_{\tander^{\Ntop} \in \mathfrak{P}^{(\Ntop)}}\int_{\characteristicdiamondtwoarg{[\leftubar,\ubar]}{[\moreinterestingu_1,u]}} \upmu \left| \tander^{\Ntop} (\vortrenormalized,\GradEnt)\right|^2  |\ubar|^{A-2}, 	\label{E:INTROVORTANDGRADENTBULK}
		\end{align}
	\end{subequations}
	\begin{subequations} \label{E:INTROVORTVORTANDDIVGRADENTL2CONTROLLINGQUANTITIES}
		\begin{align}
			\ingoingfluxwave{\Ntop}[(\VortVort,\DivGradEnt)](\ubar,u) & = \sum_{\tander^{\Ntop} \in \mathfrak{P}^{(\Ntop)}}\int_{\ingoingcharacteristicsurfacetwoarg{\ubar}{[\moreinterestingu_1,u]}} \upmu \left|  \tander^{\Ntop} (\VortVort,\DivGradEnt) \right|^2 |\ubar|^{A + \frac{3}{2}}, \label{E:INTROINGOINGVORTVORTDIVGRADENTFLUX} \\		
			\outgoingfluxwave{\Ntop}[(\VortVort,\DivGradEnt)](\ubar,u) & = \sum_{\tander^{\Ntop} \in \mathfrak{P}^{(\Ntop)}}\int_{\outgoingcharacteristicsurfacetwoarg{u}{[\leftubar,\ubar]}}  \left|  \tander^{\Ntop} (\VortVort,\DivGradEnt) \right|^2 |\ubar|^{A+\frac{3}{2}},	\label{E:INTROOUTGOINGVORTVORTDIVGRADENTFLUX} \\
			\newspacetimecoercive{\Ntop}[(\VortVort,\DivGradEnt)](\ubar,u) & = \sum_{\tander^{\Ntop} \in \mathfrak{P}^{(\Ntop)}} \int_{\characteristicdiamondtwoarg{[\leftubar,\ubar]}{[\moreinterestingu_1,u]}} \upmu   \left|  \tander^{\Ntop} (\VortVort,\DivGradEnt) \right|^2 |\ubar|^{A + \frac{1}{2}}, 	\label{E:INTROOUTGOINGVORTVORTDIVGRADENTBULK} 			
		\end{align}
	\end{subequations}
where in \eqref{E:INTROINGOINGVELOCITYFLUX} and \eqref{E:INTROWEAKCOERCIVEVELOCITYBULK}, $\newuL^\alpha$ is a $\gfour$-null vectorfield proportional to $(\gfour^{-1})^{\alpha\beta} \p_\beta \ubar$ normalized so that\footnote{In particular, the reader could roughly think that $\newuL$ is the $\gfour$-null analog of the $\gfour$-spacelike vectorfield  $\muX$, which satisfies $\muX u = 1$} $\newuL u = 1$, and $\nullgeop{x^B}$ is the unique vectorfield satisfying $\nullgeop{x^B}(x^C) = \updelta_B^C$ for $B,C\in \{2,3\}$ and $\nullgeop{x^B}\ubar = \nullgeop{x^B} u = 0$. In \eqref{E:INTROVORTVORTANDDIVGRADENTL2CONTROLLINGQUANTITIES}, $\VortVort \simeq \Flatcurl \, \vortrenormalized + \cdots $ and $\DivGradEnt \simeq \Delta s + \cdots$ are the top-order modified fluid variables of Remark\,\ref{R:MODIFIEDFLUIDVARS}.

Even before mentioning anything about the proof of the $L^2$estimates, we first make a series of remarks concerning the $L^2$-controlling quantities \eqref{E:INTROVELOCITYL2CONTROLLINGQUANTITIES}--\eqref{E:INTROVORTVORTANDDIVGRADENTL2CONTROLLINGQUANTITIES}.

	\begin{enumerate}
		\item For reasons that we will describe below, it is important that the $L^2$-controlling quantities of the velocity is a sum of each component $v^i$. This is different compared to the $L^2$-controlling quantity of, say, $\vortrenormalized$ or $\VortVort$ in \eqref{E:INTROVORTANDGRADENTL2CONTROLLINGQUANTITIES}--\eqref{E:INTROVORTVORTANDDIVGRADENTL2CONTROLLINGQUANTITIES}, which have only been expressed schematically in this introduction.\footnote{In the precise analysis, the $L^2$-controlling quantities take the maximum over $i \in \{1,2,3\}$, see Def.\,\ref{D:FUNDAMENTALTRANSPORTL2CONTROLLINGQUANTITIES}.} 
		\item For reasons that we will describe below, \emph{we do not derive $L^2$-estimates for the logarithmic density $\LogDensity$} (and therefore do not need to define an $L^2$-controlling quantity for it). To control this fundamental quantity, we instead derive a new geometric decomposition of all first derivatives $Z \LogDensity$ for $Z \in \{ \Lunit,\muX,\Yvf{2},\Yvf{3}\}$ expressed in terms of first derivatives of the velocity, vorticity, and entropy. The precise details are in Sect.\,\ref{S:DENSITYDERIVATIVESINTERMSOFOTHERS}. 
		
		For related reasons, we do not treat the entropy $\Ent$ as a wave variable, even though it solves the covariant wave equation \eqref{E:MODIFIEDFLUIDVARIABLES}. Instead, using that $\Ent$ is transported (see \eqref{E:INTROBS}), we can derive sufficient estimates for the entropy by instead controlling its gradient
		$\GradEnt$, which we can treat as a transport variable; see \eqref{E:MODIFIEDFLUIDVARIABLES}.
		\item The constant $A$ in \eqref{E:INTROVELOCITYL2CONTROLLINGQUANTITIES}--\eqref{E:INTROVORTVORTANDDIVGRADENTL2CONTROLLINGQUANTITIES} can be interpreted as the same uniform constant $A$ in \eqref{E:INTRO:SCHEMATICTOPORDERESTIMTEFORROUGHFOLIATIONS}.
		\item The ingoing flux for the velocity \eqref{E:INTROINGOINGVELOCITYFLUX} is only coercive in derivative directions that are \emph{tangent} to the ingoing characteristics $\ingoingcharacteristicsurfacetwoarg{\ubar}{[\moreinterestingu_1,u]}$. Thus, while at a formal level the null flux $\ingoingfluxwave{\Ntop}[V](\ubar,u)$ plays the role of the $L^2$-energies of $\widetilde{\mathbf{E}}_N(\timefunction,u)$ in \eqref{E:INTRO:ROUGHHYPERSURFACEENERGIES} from our prior work \cite{abbrescia2022emergence}, the null-fluxes here are only semi-definite while the energies there were positive-definite. This leads to numerous new novel difficulties in the analysis discussed below.
	\end{enumerate}

One of the main differences between the $L^2$-controlling quantities from \eqref{E:INTROVELOCITYL2CONTROLLINGQUANTITIES}--\eqref{E:INTROVORTVORTANDDIVGRADENTL2CONTROLLINGQUANTITIES} and those in \eqref{E:ROUGHL2CONTROLINGQUANTITIES} are the weights in powers of $|\ubar| = -\ubar$. Consequently, the fundamental energy inequality for the velocity takes the following form in this paper:
\begin{align}
		\begin{split} \label{E:INTROSCHEMATICENERGYINEQUALITYFORVELOCITY}
			& \ingoingfluxwave{\Ntop}[V](\ubar,u) + \outgoingfluxwave{\Ntop}[V](\ubar,u)  +  \boxed{A}  \newspacetimecoercive{\Ntop}[V](\ubar,u) + \spacetimecoercive{\Ntop}[V](\ubar,u) \\
			& \ \ \le \text{data terms} +  \sum_{\tander^{\Ntop} \in \mathfrak{P}^{(\Ntop)}} \int_{\characteristicdiamondtwoarg{[\leftubar,\ubar]}{[\moreinterestingu_1,u]}} \sum_{i=1}^3 \left\{ \tander^{\Ntop} \mytr_{\gtorus}\upchi + \tander^{\Ntop} \left(\upmu \VortVort^i \right) \right\}  (\multipliervectorfield \tander^{\Ntop}v^i  )|\ubar|^A + \cdots, 
		\end{split}
	\end{align}
where $ \tander^{\Ntop} \mytr_{\gtorus}\upchi + \tander^{\Ntop} \left(\upmu \VortVort^i \right) $ denote the worst terms from the RHS $\upmu \Box_{\gfour} \tander^{\Ntop}v^i$. Compared to the energy method from Sect.\,\ref{SSS:ENERGYDESCENTSCHEMEINROUGHCOORDINATES}, where the mid- to top-order $L^2$-controlling quantities have the potential to blow-up via the descent scheme, \emph{every single controlling quantity in this paper remains bounded and small all the way up to the Cauchy horizon}, even at top-order! The descent scheme in this paper is no longer indexed with singular powers, but with \emph{weighted} powers of $|\ubar|$, e.g.: 
	\begin{subequations} \label{E:INTRODESCENTSCHEMEFORCH}
		\begin{align}
			\ingoingfluxwave{\Ntop-1}[V](\ubar,u) & = \sum_{i=1}^3 \sum_{\tander^{\Ntop} \in \mathfrak{P}^{(\Ntop)}}\int_{\ingoingcharacteristicsurfacetwoarg{\ubar}{[\moreinterestingu_1,u]}} \left\{ \left| \newuL \tander^{\Ntop} v^i\right|^2 +\sum_{B=2}^3 \upmu \left| \nullgeop{x^B} \tander^{\Ntop}  v^i \right|^2 \right\} |\ubar|^{A-2},  \label{E:INTRODESCENTJUSTBELOWTOP} \\
			\ingoingfluxwave{\Ntop-2}[V](\ubar,u) & = \sum_{i=1}^3 \sum_{\tander^{\Ntop} \in \mathfrak{P}^{(\Ntop)}}\int_{\ingoingcharacteristicsurfacetwoarg{\ubar}{[\moreinterestingu_1,u]}} \left\{ \left| \newuL \tander^{\Ntop} v^i\right|^2 +\sum_{B=2}^3 \upmu \left| \nullgeop{x^B} \tander^{\Ntop}  v^i \right|^2 \right\} |\ubar|^{A-4},  \label{E:INTRODESCENTTWOJUSTBELOWTOP}
		\end{align}
	\end{subequations}
and similarly for the remaining quantities from \eqref{E:INTROVELOCITYL2CONTROLLINGQUANTITIES}--\eqref{E:INTROVORTVORTANDDIVGRADENTL2CONTROLLINGQUANTITIES}. This is of particular interest for future applications to the shock development problem, where the initial data for shock development is the state of the solution on the Cauchy horizon, the crease, and the singular boundary. The upshot is that using a weighted descent scheme would provide \emph{finite} initial data to the shock development problem, instead of potentially infinite initial data. These ideas could have been implemented in \cite{abbrescia2022emergence}, and with very minor additional difficulty, one could adapt the analysis of our prior work \cite{abbrescia2022emergence} to prove that weighted $L^2$-controlling quantities tailored to the shape of the singular boundary, e.g.  
	\begin{align}
		\int_{\hypthreearg{\timefunction}{[-\rightu,u]}{0}} \left\{ \upmu \left|\Lunit \tander^N \wavearray\right|^2 + \sum_{A=2}^3  \upmu \left| \Yvf{A} \tander^N \wavearray \right|^2 + \left| \muX \tander^N \wavearray \right|^2\right\}|\timefunction|^A  \label{E:ALTERNATIVEWEIGHTEDENERGYIESFORSINGULARBOUNDARYPAPER}
	\end{align}
also remain bounded. We also mention that the weighted powers of the transport variables in \eqref{E:INTROVORTANDGRADENTL2CONTROLLINGQUANTITIES}--\eqref{E:INTROVORTVORTANDDIVGRADENTL2CONTROLLINGQUANTITIES} reflect the same singular powers of potential blow-up for the $L^2$ controlling quantities in \cite{abbrescia2022emergence}. For example, in that paper, we had:
	\begin{align}
		  \int_{\hypthreearg{\timefunction}{[-\rightu,u]}{0}}  \upmu \left|\tander^{\Ntop} \vortrenormalized \right|^2 \lesssim \mr{\upepsilon}^2 |\timefunction|^{-A+1},\qquad  \int_{\hypthreearg{\timefunction}{[-\rightu,u]}{0}} \upmu \left|\tander^{\Ntop} \VortVort\right|^2 \lesssim \mr{\upepsilon}^2 |\timefunction|^{-A-\frac{3}{2}}. \label{E:INTROSAMESINGULARPOWERSINSINGULARBOUNDARYASCHWEIGHTS}
	\end{align}

\begin{remark}[Weights tailored to the geometry] 
To study the solution near the Cauchy horizon, we use weighted $L^2$-energies whose weights are tailored to the geometry of the foliations; see, e.g., 
\eqref{E:INTROVELOCITYL2CONTROLLINGQUANTITIES}--\eqref{E:INTROVORTVORTANDDIVGRADENTL2CONTROLLINGQUANTITIES}. 
The use of such weights outside of symmetry appears to have originated in Christodoulou's earlier work \cite{dC2019}*{Sect.\,9.7} on the restricted shock development problem. Analogous weights were used in \cite{shkoller2024geometry}. 
\end{remark}

\paragraph{\emph{Almost Riemann invariants are not suitable as solution variables for the Cauchy horizon}} \label{PAR:ALMOSTRIEMANNINVARIANTSNOTUSABLEFORCH} \hfill

Another main difference between this work and \cite{abbrescia2022emergence} is that, in our prior work, the main solution wave variables were in fact $\{ \RRiemann,\LRiemann,v^2,v^3,s\}$, where $\mathcal{R}_{\pm}:= v \pm F(\LogDensity,\Ent)$ are ``almost'' Riemann invariants, and $F$ is determined by the equation of state (see \eqref{E:ALMOSTRIEMANNINVARIANTS}). These are natural choices for the fundamental unknowns because the main result of \cite{abbrescia2022emergence} constructs the singular boundary for an open set of data centered at a \emph{simple, plane-symmetric, isentropic} background solution $(\RRiemannPS, \LRiemannPS =0, v^2 = 0, v^3 = 0, \Ent = \Ent_0)$, where $\RRiemannPS$ was explicitly constructed to feature quantitative transversal convexity. We now explain why the fundamental solution wave variables in this paper are $\{v^1,v^2,v^3\}$ (with $\Ent$ being treated purely as a transport variable through its gradient $\GradEnt$, and $\LogDensity$ being controlled in terms of the velocity, vorticity, and entropy), even though the class of data for this paper is the same as in \cite{abbrescia2022emergence}, i.e. perturbations of simple, plane-symmetric, isentropic transversally convex data. If we \emph{were} to use the almost Riemann invariants, the top-order weighted null flux-energy for inequality would schematically read:
\begin{align}
		\begin{split} \label{E:INTROSCHEMATICENERGYINEQUALITYFORRRIEMANNONDIAMOND}
			& \ingoingfluxwave{\Ntop}[\RRiemann](\ubar,u) + \outgoingfluxwave{\Ntop}[\RRiemann](\ubar,u)  +  \boxed{A}  \newspacetimecoercive{\Ntop}[\RRiemann](\ubar,u) + \spacetimecoercive{\Ntop}[\RRiemann](\ubar,u) \\
			& \ \ \le \text{data terms} +  \sum_{\tander^{\Ntop} \in \mathfrak{P}^{(\Ntop)}} \int_{\characteristicdiamondtwoarg{[\leftubar,\ubar]}{[\moreinterestingu_1,u]}} \left\{ \tander^{\Ntop} \mytr_{\gtorus}\upchi + \tander^{\Ntop} \left(\upmu \VortVort +  \upmu \DivGradEnt \right) \right\}  (\multipliervectorfield \tander^{\Ntop}\RRiemann  )|\ubar|^A + \cdots, 
		\end{split} \\
		\begin{split} \label{E:INTROSCHEMATICENERGYINEQUALITYFORVORTVORT}
			& \ingoingfluxwave{\Ntop}[\VortVort](\ubar,u) + \outgoingfluxwave{\Ntop}[\VortVort](\ubar,u)  + \boxed{A+\tfrac{3}{2}} \newspacetimecoercive{\Ntop}[\VortVort](\ubar,u) \\
			& \ \ \le \text{data terms} +  \sum_{\tander^{\Ntop} \in \mathfrak{P}^{(\Ntop)}} \int_{\characteristicdiamondtwoarg{[\leftubar,\ubar]}{[\moreinterestingu_1,u]}} \left\{ \pmb{\p} \tander^{\Ntop}(\vortrenormalized,\GradEnt)  \right\}( \tander^{\Ntop} \VortVort )|\ubar|^{A+\frac{3}{2}} + \cdots,
		\end{split}\\
		\begin{split} \label{E:INTROSCHEMATICENERGYINEQUALITYFORDIVGRADENT}
			& \ingoingfluxwave{\Ntop}[\DivGradEnt](\ubar,u) + \outgoingfluxwave{\Ntop}[\DivGradEnt](\ubar,u)  + \boxed{A+\tfrac{3}{2}} \newspacetimecoercive{\Ntop}[\DivGradEnt](\ubar,u) \\
			& \ \ \le \text{data terms} +  \sum_{\tander^{\Ntop} \in \mathfrak{P}^{(\Ntop)}} \int_{\characteristicdiamondtwoarg{[\leftubar,\ubar]}{[\moreinterestingu_1,u]}} \left\{ \pmb{\p} \tander^{\Ntop}(\vortrenormalized,\GradEnt)  \right\}( \tander^{\Ntop} \DivGradEnt )|\ubar|^{A+\frac{3}{2}} + \cdots.
		\end{split}		
	\end{align}
	In \eqref{E:INTROSCHEMATICENERGYINEQUALITYFORRRIEMANNONDIAMOND}, the $L^2$ controlling quantities for the almost Riemann invariant controls the same derivatives as in the case for the velocity in \eqref{E:INTROINGOINGVELOCITYFLUX}, e.g. 
	\begin{align}
		\outgoingfluxwave{\Ntop}[\RRiemann](\ubar,u)  = \sum_{\tander^{\Ntop} \in \mathfrak{P}^{(\Ntop)}}\int_{\outgoingcharacteristicsurfacetwoarg{u}{[\leftubar,\ubar]}} \left\{ \left| \Lunit \tander^{\Ntop} \RRiemann\right|^2 +\sum_{B=2}^3 \upmu \left| \Yvf{A} \tander^{\Ntop}  \RRiemann \right|^2 \right\} |\ubar|^A, \label{E:INTROSCHEMATICOUTGOINGFLUXFORRRIEMANN}
	\end{align}
	and similarly for $\ingoingfluxwave{\Ntop}[\RRiemann](\ubar,u), \newspacetimecoercive{\Ntop}[\RRiemann](\ubar,u)$, and $\spacetimecoercive{\Ntop}[\RRiemann](\ubar,u)$. The multiplier vectorfield $\breve{T}:= (1+2\upmu)\Lunit + 2\muX$ is the same one used in \cite{dC2007,jS2016b,LukSpeck2024stability}. The term $\tander^{\Ntop} \mytr_{\gtorus}\upchi$ appears on the RHS\,\eqref{E:INTROSCHEMATICENERGYINEQUALITYFORRRIEMANNONDIAMOND} for the same reason it did in \eqref{E:SCHEMATICCOMMUTEDWAVEEQUATION}, and that reason $\tander^{\Ntop} (\upmu \VortVort) \simeq\tander^{\Ntop}( \upmu \Flatcurl \vortrenormalized + \cdots)$ and $\tander^{\Ntop} (\upmu \DivGradEnt) \simeq \tander^{\Ntop}(\upmu \Flatdiv \GradEnt + \cdots)$ both appear on the on RHS\,\eqref{E:INTROSCHEMATICENERGYINEQUALITYFORRRIEMANNONDIAMOND} is because \eqref{E:COVWAVEEQv}--\eqref{E:COVWAVEEQRHO} imply $ \Box_{\gfour}\RRiemann  \simeq  \VortVort +  \DivGradEnt + \cdots$. The top-ordered error integrals generated from  $\tander^{\Ntop}\mytr_{\gtorus} \upchi$ are handled in a similar way as in Sect.\,\ref{SSS:GAUSSCURVATUREESTIMATESNOTREQUIREDIN2D}. We do not elaborate on it further and only note that, just like in Sect.\,\ref{SSS:GAUSSCURVATUREESTIMATESNOTREQUIREDIN2D}, the use of Raychaudhuri allows us to avoid a loss of derivative at top order at the expense of the weighted descent scheme of \eqref{E:INTRODESCENTSCHEMEFORCH}. The term $\pmb{\p}\tander^{\Ntop}(\vortrenormalized,\GradEnt)$ in braces featuring in \eqref{E:INTROSCHEMATICENERGYINEQUALITYFORVORTVORT} comes from commuting the transport equation satisfied by $\VortVort$ at top order: $\upmu \Transport \tander^{\Ntop}\VortVort \simeq  \pmb{\p}\tander^{\Ntop}(\vortrenormalized,\GradEnt) + \cdots$, where $\cdots$ are easier error terms.

We first focus on the following error term, generated from the un-weighted $\Lunit$ factor of the multiplier vectorfield and applying the Leibniz rule: 
\begin{align} \label{E:IMPOSSIBLEERRORTERMFROMRIEMANNINVARIANTS}
		 \int_{\characteristicdiamondtwoarg{[\leftubar,\ubar]}{[\moreinterestingu_1,u]}} \left\{\upmu \tander^{\Ntop}\VortVort +  \upmu \tander^{\Ntop} \DivGradEnt  \right\}  (\Lunit \tander^{\Ntop}\RRiemann  )|\ubar|^A.
	\end{align}
As both factors are top-order, we simply bound these terms as follows using H\"older's and Young's inequality:
	\begin{align} 
		\begin{split} \label{E:TRYINGTOBOUNDIMPOSSIBLEERRORTERMFROMRIEMANNINVARIANTS}
			 \int_{\characteristicdiamondtwoarg{[\leftubar,\ubar]}{[\moreinterestingu_1,u]}} \left\{\upmu \tander^{\Ntop}\VortVort +  \upmu \tander^{\Ntop} \DivGradEnt  \right\}  (\Lunit \tander^{\Ntop}\RRiemann  )|\ubar|^A & \lesssim  \int_{\characteristicdiamondtwoarg{[\leftubar,\ubar]}{[\moreinterestingu_1,u]}} \left\{ \upmu^2 \left| \tander^{\Ntop} \VortVort\right|^2 +  \upmu^2 \left| \tander^{\Ntop} \DivGradEnt \right|^2 \right\} |\ubar|^A \\
			 & \ \ + \int_{\characteristicdiamondtwoarg{[\leftubar,\ubar]}{[\moreinterestingu_1,u]}}  \left| \Lunit \tander^{\Ntop} \RRiemann \right|^2 |\ubar|^A \\
		\end{split}
	\end{align}
The last integral on the RHS\,\eqref{E:TRYINGTOBOUNDIMPOSSIBLEERRORTERMFROMRIEMANNINVARIANTS} featuring $\tander^{\Ntop}\RRiemann$ is $= \int_{u' = \moreinterestingu_1}^u \outgoingfluxwave{\Ntop}[\RRiemann](\ubar,u') \, \rmd  u'$, which can be controlled using Gr\"onwall's inequality in \eqref{E:INTROSCHEMATICENERGYINEQUALITYFORRRIEMANNONDIAMOND}. The other error integral featuring $\tander^{\Ntop}(\VortVort,\DivGradEnt)$ on RHS\,\eqref{E:TRYINGTOBOUNDIMPOSSIBLEERRORTERMFROMRIEMANNINVARIANTS} \emph{does not have sufficient $|\ubar|$-weights to be Gr\"onwalled} using \eqref{E:INTROSCHEMATICENERGYINEQUALITYFORVORTVORT}--\eqref{E:INTROSCHEMATICENERGYINEQUALITYFORDIVGRADENT}. The bound $\upmu \lesssim |\ubar|$ is provably false on the characteristic diamond and hence there is no hope to try to absorb the error integrals featuring $\tander^{\Ntop}(\VortVort,\DivGradEnt)$ by the bulk terms $\boxed{A+\tfrac{3}{2}} \newspacetimecoercive{\Ntop}[\VortVort](\ubar,u)$ and $\boxed{A+\tfrac{3}{2}} \newspacetimecoercive{\Ntop}[\DivGradEnt](\ubar,u)$ on LHS\,\eqref{E:INTROSCHEMATICENERGYINEQUALITYFORVORTVORT}--\eqref{E:INTROSCHEMATICENERGYINEQUALITYFORDIVGRADENT}. Contrastingly, if one tries to ``borrow'' half a factor of $|\ubar|$ from the $\Lunit \tander^{\Ntop}\RRiemann$ error integral, one could bound: 
	\begin{align} 
		\begin{split} \label{E:SECONDTRYTOBOUNDIMPOSSIBLEERRORTERMFROMRIEMANNINVARIANTS}
			 \int_{\characteristicdiamondtwoarg{[\leftubar,\ubar]}{[\moreinterestingu_1,u]}} \left\{\upmu \tander^{\Ntop}\VortVort +  \upmu \tander^{\Ntop} \DivGradEnt  \right\}  (\Lunit \tander^{\Ntop}\RRiemann  )|\ubar|^A & \lesssim  \varsigma \int_{\characteristicdiamondtwoarg{[\leftubar,\ubar]}{[\moreinterestingu_1,u]}} \left\{ \upmu^2 \left| \tander^{\Ntop} \VortVort\right|^2 +  \upmu^2 \left| \tander^{\Ntop} \DivGradEnt \right|^2 \right\} |\ubar|^{A+\frac{1}{2}}\\
			 & \ \ +  \varsigma^{-1}\int_{\characteristicdiamondtwoarg{[\leftubar,\ubar]}{[\moreinterestingu_1,u]}}  \left| \Lunit \tander^{\Ntop} \RRiemann \right|^2 |\ubar|^{A-\tfrac{1}{2}},
		\end{split}
	\end{align}
for some small constant $\varsigma$. Since $\upmu \lesssim 1$, the spacetime error integrals on RHS\,\eqref{E:SECONDTRYTOBOUNDIMPOSSIBLEERRORTERMFROMRIEMANNINVARIANTS} featuring  $\tander^{\Ntop}(\VortVort,\DivGradEnt)$ could be absorbed by the bulk terms $\boxed{A+\tfrac{3}{2}} \newspacetimecoercive{\Ntop}[\VortVort](\ubar,u)$ and $\boxed{A+\tfrac{3}{2}} \newspacetimecoercive{\Ntop}[\DivGradEnt](\ubar,u)$ on LHS\,\eqref{E:INTROSCHEMATICENERGYINEQUALITYFORVORTVORT}--\eqref{E:INTROSCHEMATICENERGYINEQUALITYFORDIVGRADENT} if $\varsigma$ is sufficiently small. However, the second error integral $\Lunit \tander^{\Ntop} \RRiemann$ can no longer be Gr\"onwalled with the outgoing flux $\outgoingfluxwave{\Ntop}[\RRiemann](\ubar,u)$ on LHS\,\eqref{E:INTROSCHEMATICENERGYINEQUALITYFORRRIEMANNONDIAMOND} because there is no uniform upper bound for $|\ubar|^{-1}$ on the outgoing null hypersurfaces $\outgoingcharacteristicsurfacetwoarg{u}{[\leftubar,\ubar]}$.

This issue was not present in our work \cite{abbrescia2022emergence}\footnote{As we mentioned above, we did not use $|\timefunction|$-weighted energies and fluxes in \cite{abbrescia2022emergence}. Regardless, an analogous estimate as in \eqref{E:VALIDTOPORDERESTIMATEINVOLVINGVORTVORTANDDIVGRADENTANDRIEMANNVAIRABLEINSINGULARBOUNDARY} holds for the potentially singular energies introduced in Sect.\,\ref{SSS:ENERGYDESCENTSCHEMEINROUGHCOORDINATES}.} because the corresponding error integral could be controlled instead by:
\begin{align}
		\begin{split} \label{E:VALIDTOPORDERESTIMATEINVOLVINGVORTVORTANDDIVGRADENTANDRIEMANNVAIRABLEINSINGULARBOUNDARY}
			\int_{\twoargMrough{[\timefunction_0,\timefunction],[-\rightu,u]}{0}}  \left\{\upmu \tander^{\Ntop}\VortVort +  \upmu \tander^{\Ntop} \DivGradEnt  \right\}  (\Lunit \tander^{\Ntop}\RRiemann  )|\timefunction|^A & \lesssim  \int_{\twoargMrough{[\timefunction_0,\timefunction],[-\rightu,u]}{0}} \left\{ \upmu \left|\tander^N \VortVort\right|^2 + \upmu \left| \tander^N \DivGradEnt \right|^2 \right\} |\timefunction|^{A + \frac{3}{4}}  \\
			& \ \ \ \ + \int_{\twoargMrough{[\timefunction_0,\timefunction],[-\rightu,u]}{0}}  \upmu \left|\Lunit \tander^N \RRiemann \right|^2 |\timefunction|^{A-\frac{3}{4}} \, \\
			& \lesssim \int_{\timefunction'=\timefunction_0}^\timefunction \frac{1}{|\timefunction|^{3/4}} \left\{\widetilde{\mathbf{E}}_{\Ntop}[\VortVort](\timefunction',u) +  \widetilde{\mathbf{E}}_{\Ntop}[\DivGradEnt](\timefunction',u)\right\} d\timefunction' \\
			& \ \ \ \ +  \int_{\timefunction'=\timefunction_0}^\timefunction \frac{1}{|\timefunction|^{3/4}} \widetilde{\mathbf{E}}_{\Ntop}[\RRiemann](\timefunction',u)  d\timefunction',
		\end{split}
	\end{align}
where the energies $\widetilde{\mathbf{E}}_{\Ntop}[\RRiemann](\timefunction,u), \widetilde{\mathbf{E}}_{\Ntop}[\DivGradEnt](\timefunction,u)$, and $\widetilde{\mathbf{E}}_{\Ntop}[\VortVort](\timefunction,u)$ for the top-order Riemann invariant and  modified fluid variables on the rough hypersurfaces   are given by:
	\begin{subequations}	
		\begin{align}
			\widetilde{\mathbf{E}}_{\Ntop}[\RRiemann](\timefunction,u) & =  \int_{\hypthreearg{\timefunction}{[-\rightu,u]}{0}} \left\{\upmu \left|\Lunit \tander^{\Ntop} \RRiemann\right|^2 + \sum_{A=2}^3  \upmu \left| \Yvf{A} \tander^{\Ntop} \RRiemann \right|^2 + \left| \muX \tander^{\Ntop} \RRiemann \right|^2\right\}|\timefunction|^A, \label{E:SPACELIKEENERGYRIEMANNONROUGHHYPERSURFACES} \\		
			\widetilde{\mathbf{E}}_{\Ntop}[\VortVort](\timefunction,u) & =  \int_{\hypthreearg{\timefunction}{[-\rightu,u]}{0}} \upmu |\tander^{\Ntop}\VortVort|^2|\timefunction|^{A+\tfrac{3}{2}},  \label{E:SPACELIKEENERGYVORTVORTONROUGHHYPERSURFACES} \\
			\widetilde{\mathbf{E}}_{\Ntop}[\DivGradEnt](\timefunction,u) & =  \int_{\hypthreearg{\timefunction}{[-\rightu,u]}{0}} \upmu |\tander^{\Ntop}\DivGradEnt|^2|\timefunction|^{A+\tfrac{3}{2}}. \label{E:SPACELIKEENERGYGRADENTONROUGHHYPERSURFACES} 
		\end{align}
	\end{subequations}			
The upshot is that $|\timefunction|^{-3/4}$ is integrable on $[\timefunction_0,0]$, and so the RHS\,\eqref{E:VALIDTOPORDERESTIMATEINVOLVINGVORTVORTANDDIVGRADENTANDRIEMANNVAIRABLEINSINGULARBOUNDARY} could be controlled using Gr\"onwall. 

The fundamental difference between \eqref{E:VALIDTOPORDERESTIMATEINVOLVINGVORTVORTANDDIVGRADENTANDRIEMANNVAIRABLEINSINGULARBOUNDARY}  and the failed attempts in \eqref{E:TRYINGTOBOUNDIMPOSSIBLEERRORTERMFROMRIEMANNINVARIANTS}--\eqref{E:SECONDTRYTOBOUNDIMPOSSIBLEERRORTERMFROMRIEMANNINVARIANTS}  is that level sets of $\timefunction$ provide a $\gfour$-\emph{spacelike} foliation of $\twoargMrough{[\timefunction_0,\timefunction],[-\rightu,u]}{0}$, while level sets of $\ubar$ provide a $\gfour$-\emph{null} foliation of $\characteristicdiamondtwoarg{[\leftubar,\ubar]}{[\moreinterestingu_1,u]}$. Indeed, the hypersurface energies given by \eqref{E:SPACELIKEENERGYRIEMANNONROUGHHYPERSURFACES} are positive definite in all derivatives of $ \tander^{\Ntop} \RRiemann$, \emph{including $\Lunit$-derivatives}, while the hypersurface fluxes on $\ingoingcharacteristicsurfacetwoarg{\ubar}{[\moreinterestingu_1,\moreinterestingu_2]} $ are only positive semi-definite and do not control $|\Lunit \tander^{\Ntop}\RRiemann|^2$. In particular, we could not multiply and divide  \eqref{E:IMPOSSIBLEERRORTERMFROMRIEMANNINVARIANTS} by powers of $|\ubar|$ as in \eqref{E:VALIDTOPORDERESTIMATEINVOLVINGVORTVORTANDDIVGRADENTANDRIEMANNVAIRABLEINSINGULARBOUNDARY} and apply Gr\"onwall for the same reason that Gr\"onwall was not applicable to RHS\,\eqref{E:SECONDTRYTOBOUNDIMPOSSIBLEERRORTERMFROMRIEMANNINVARIANTS}.

\paragraph{\emph{Using the velocity as a fundamental unknown}} \label{PAR:VELOCITYASAFUNDAMENTALUNKNOWN} \hfill

We now explain why using $\{v^1,v^2,v^3\}$ as the fundamental wave variables overcomes the following glaring issue:  the fluxes \eqref{E:INTROINGOINGVELOCITYFLUX}--\eqref{E:INTROOUTGOINGVELOCITYFLUX} for the velocity are also semi-definite and feature the same $\upmu$ and $|\ubar|$ weights as those of the almost Riemann invariants, leading to seemingly uncontrollable top-order error integrals as in Sect.\,\ref{PAR:ALMOSTRIEMANNINVARIANTSNOTUSABLEFORCH}. Since the $L^2$-controlling quantity for the velocity features the \emph{sum} of each component, the dangerous integral is, schematically:
\begin{align}
		\int_{\characteristicdiamondtwoarg{[\leftubar,\ubar]}{[\moreinterestingu_1,u]}} \sum_{i=1}^3  \tander^{\Ntop}\left(\upmu \VortVort^i \right) \left( \Lunit \tander^{\Ntop} v^i\right) |\ubar|^A. \label{E:INTRODANGEROUSTOPORDERINTEGRALBUTUSINGVELOCITY}
	\end{align}We note here already some differences between \eqref{E:INTRODANGEROUSTOPORDERINTEGRALBUTUSINGVELOCITY} and \eqref{E:IMPOSSIBLEERRORTERMFROMRIEMANNINVARIANTS}. Most noticeably, there is no $\DivGradEnt$ term present in \eqref{E:INTRODANGEROUSTOPORDERINTEGRALBUTUSINGVELOCITY} because the geometric wave equation satisfied by $v^i$ do not feature them (see \eqref{E:COVWAVEEQv} and Rmk.\,\ref{R:MODIFIEDFLUIDVARS}). The second main difference is that we have not applied the Leibniz rule to $\upmu \VortVort^i$ in \eqref{E:INTRODANGEROUSTOPORDERINTEGRALBUTUSINGVELOCITY}. Instead, using $\upmu \VortVort^i \simeq \upmu \Flatcurl(\vortrenormalized)^i + \cdots$, we commute\footnote{The commutation had to be done in this order. Writing $\tander^{\Ntop} (\upmu\VortVort) \simeq \upmu \tander^{\Ntop} \VortVort + \cdots \simeq \upmu \Flatcurl \tander^{\Ntop} \vortrenormalized + [\tander^{\Ntop},\Flatcurl] \vortrenormalized + \cdots$ would lead to additional singular powers of $\upmu^{-1}$. Commuting with a $\upmu$-weighted $\Flatcurl$ renormalizes the singularity. \label{F:COMMUTINGFIRSTORDEROPERATORSWITHMU}} $\tander^{\Ntop}(\upmu \VortVort)^i = \upmu \Flatcurl\left(\tander^{\Ntop}\vortrenormalized\right)^i + [\tander^{\Ntop},\upmu\Flatcurl]\vortrenormalized + \cdots$. By $\tander^{\Ntop}\vortrenormalized$, we mean the $\Sigma_t$-tangent vectorfield $(\tander^{\Ntop}\vortrenormalized^1,\tander^{\Ntop}\vortrenormalized^2,\tander^{\Ntop}\vortrenormalized^3)$. Ignoring many easier error terms, \eqref{E:INTRODANGEROUSTOPORDERINTEGRALBUTUSINGVELOCITY} reduces to:
\begin{align}
		\int_{\characteristicdiamondtwoarg{[\leftubar,\ubar]}{[\moreinterestingu_1,u]}} \sum_{i=1}^3 \upmu \upepsilon_{ijk} \partial_j  ( \tander^{\Ntop} \vortrenormalized)^k \left( \Lunit \tander^{\Ntop} v^i\right) |\ubar|^A, \label{E:INTRODANGEROUSTOPORDERINTEGRALBUTUSINGVELOCITYAFTEREXPANDINGVORTVORT}
	\end{align}
 where $\upepsilon_{ijk}$ is the fully anti-symmetric symbol normalized to $\upepsilon_{123} = 1$. We note that \eqref{E:INTRODANGEROUSTOPORDERINTEGRALBUTUSINGVELOCITYAFTEREXPANDINGVORTVORT} also features a summation over $j,k\in \{1,2,3\}$, but we choose to only highlight the summation over $i\in\{1,2,3\}$ for convenience and exposition. We then perform \emph{two} integration by parts: to swap the curl differential operator $ \upmu \upepsilon_{ijk} \partial_j $ and the $\Lunit$, leading to:
	\begin{align}
		\begin{split}	\label{E:INTRODANGEROUSTOPORDERINTEGRALBUTAFTERSWITCHINGMUCURLANDL}
			& \int_{\characteristicdiamondtwoarg{[\leftubar,\ubar]}{[\moreinterestingu_1,u]}} \sum_{i=1}^3 \upmu \upepsilon_{ijk} \partial_j  ( \tander^{\Ntop} \vortrenormalized)^k \left( \Lunit \tander^{\Ntop} v^i\right) |\ubar|^A =  \int_{\characteristicdiamondtwoarg{[\leftubar,\ubar]}{[\moreinterestingu_1,u]}} \sum_{i=1}^3 \Lunit ( \tander^{\Ntop} \vortrenormalized)^k \upmu \upepsilon_{ijk} \partial_j   \left(  \tander^{\Ntop} v^i\right) |\ubar|^A \\
			& \ \ +  \int_{\ingoingcharacteristicsurfacetwoarg{\ubar}{[\moreinterestingu_1,u]}}  \underline{\text{Dangerous boundary terms}}  + \int_{\outgoingcharacteristicsurfacetwoarg{u}{[\leftubar,\ubar]}} \text{Dangerous boundary terms} + \cdots,
		\end{split}
	\end{align}
where $\cdots$ are either boundary terms controlled by the initial data or harmless bulk terms. The boundary terms are dangerous because the above procedure yields boundary terms on null hypersurfaces of the form $Z \tander^{\Ntop} v^i$ for $Z \in \{\Lunit,\muX,\Yvf{2},\Yvf{3}\}$, while the null fluxes are semi-definite and only control \emph{tangential} derivatives of $\tander^{\Ntop} v^i$, see \eqref{E:INTROINGOINGVELOCITYFLUX}--\eqref{E:INTROOUTGOINGVELOCITYFLUX}. Before explaining how to control the dangerous boundary terms, we first explain how to control the bulk spacetime integral on RHS\,\eqref{E:INTRODANGEROUSTOPORDERINTEGRALBUTAFTERSWITCHINGMUCURLANDL}. Upon commuting and applying the Leibniz rule, we have $\upmu\upepsilon_{ijk} \partial_j\tander^{\Ntop} v =\upmu \tander^{\Ntop} (\upepsilon_{ijk}\p_j v^i) + [\upmu\epsilon_{ijk}\partial_j,\tander^{\Ntop}] v^k + \cdots$. Next, we use  the anti-symmetry of $\upepsilon_{ijk}$ to write $\upepsilon_{ijk}\p_j v^i = - \exp(\LogDensity) \vortrenormalized^k$. This proves that the bulk integral features \emph{a perfect derivative} up to error terms:
	\begin{align}
		  \int_{\characteristicdiamondtwoarg{[\leftubar,\ubar]}{[\moreinterestingu_1,u]}} \sum_{i=1}^3 \Lunit ( \tander^{\Ntop} \vortrenormalized)^k \upmu \upepsilon_{ijk} \partial_j   \left(  \tander^{\Ntop} v^i\right) |\ubar|^A = -  \frac{1}{2}   \int_{\characteristicdiamondtwoarg{[\leftubar,\ubar]}{[\moreinterestingu_1,u]}} \Lunit \left( \upmu \exp(\LogDensity)  \left|\tander^{\Ntop} \vortrenormalized\right|^2 |\ubar|^A \right) + \cdots, \label{E:PERFECTDERIVATIVEINTHETOPORDERTERM}
	\end{align}
where $\cdots$ in \eqref{E:PERFECTDERIVATIVEINTHETOPORDERTERM} features easier error terms. The term with the perfect derivative even has a good sign\footnote{The good sign featured in RHS\,\eqref{E:THEGOODSIGN} is not strictly needed for this paper, though we expect it could have potential use in future projects.}  (note that \eqref{E:PERFECTDERIVATIVEINTHETOPORDERTERM} is present in the RHS of the energy identity):
\begin{align}
		\begin{split} \label{E:THEGOODSIGN}
			&  - \frac{1}{2}   \int_{\characteristicdiamondtwoarg{[\leftubar,\ubar]}{[\moreinterestingu_1,u]}} \Lunit \left( \upmu \exp(\LogDensity)  \left|\tander^{\Ntop} \vortrenormalized\right|^2 |\ubar|^A \right)  \\
			& \ \ = -  \frac{1}{2} \int_{\ingoingcharacteristicsurfacetwoarg{\ubar}{[\moreinterestingu_1,u]}}  \upmu \exp(\LogDensity)  \left|\tander^{\Ntop} \vortrenormalized\right|^2 |\ubar|^A +  \frac{1}{2} \int_{\ingoingcharacteristicsurfacetwoarg{\leftubar}{[-\rightu,u]}}  \upmu \exp(\LogDensity)  \left|\tander^{\Ntop} \vortrenormalized\right|^2 |\ubar|^A,
		\end{split}
	\end{align}
and the first integral on RHS\,\eqref{E:THEGOODSIGN} can be moved to the LHS of the fundamental energy identity for the velocity. This perfect derivative, which depended on the sum over the indices $i,j,k \in \{1,2,3\}$ is the exact reason as to why we needed to have the $L^2$-controlling quantities for the velocity be the \emph{sum} of each component; see \eqref{E:INTROVELOCITYL2CONTROLLINGQUANTITIES}. 

We now return to the dangerous boundary terms on RHS\,\eqref{E:INTRODANGEROUSTOPORDERINTEGRALBUTAFTERSWITCHINGMUCURLANDL}, specifically those on $\ingoingcharacteristicsurfacetwoarg{\ubar}{[\moreinterestingu_1,u]}$. Part of the boundary terms were generated from integrating by parts off the $\upmu\epsilon_{ijk}\partial_j$ from the vorticity term onto the velocity term in deriving \eqref{E:INTRODANGEROUSTOPORDERINTEGRALBUTAFTERSWITCHINGMUCURLANDL}\footnote{This step must come first as integrating by parts off the $\Lunit$ first would yield hypersurface boundary integrals with $\pmb{\p}\tander^{\Ntop}\vortrenormalized$ terms that we do not have control over.}.  Since $\newuL$ is the $\gfour$-normal vector to $\ingoingcharacteristicsurfacetwoarg{\ubar}{[\moreinterestingu_1,u]}$, this would yield $-  (\tander^{\Ntop} \vortrenormalized)^k ( \Lunit \tander^{\Ntop} v^i) \gfour(\upmu \upepsilon_{ijk}\partial_j, \newuL) |\ubar|^A$ as the integrand on the ingoing characteristic hypersurface. On the other hand, the resulting spacetime bulk term from this first integration by parts features the error term $- (\tander^{\Ntop}\vortrenormalized)^k \upmu \upepsilon_{ijk}\p_j \Lunit \tander^{\Ntop}v^i |\ubar|^A$. By commuting $- \upmu \upepsilon_{ijk}\p_j \Lunit \tander^{\Ntop}v^i  = - \Lunit (\upmu \upepsilon_{ijk}\p_j  \tander^{\Ntop}v^i) +[\Lunit, \upmu\upepsilon_{ijk}\p_j]\tander^{\Ntop}v^i$, we now integrate by parts off the $\Lunit$, leading to both the bulk term on RHS\,\eqref{E:INTRODANGEROUSTOPORDERINTEGRALBUTAFTERSWITCHINGMUCURLANDL} as well as another boundary integral $\ingoingcharacteristicsurfacetwoarg{\ubar}{[\moreinterestingu_1,u]}$ featuring $ (\tander^{\Ntop} \vortrenormalized)^k( \upmu \upepsilon_{ijk}\p_j \tander^{\Ntop} v^i) \gfour(\Lunit, \newuL) |\ubar|^A$. The upshot is that, by combining both of these boundary terms, there is a complete cancellation of all derivatives of $\tander^{\Ntop} v^i$ \emph{transverse} to $\ingoingcharacteristicsurfacetwoarg{\ubar}{[\moreinterestingu_1,u]}$, i.e.,
	\begin{align}
		\begin{split}
			& (\tander^{\Ntop} \vortrenormalized)^k( \upmu \upepsilon_{ijk}\p_j \tander^{\Ntop} v^i) \gfour(\Lunit, \newuL) -  (\tander^{\Ntop} \vortrenormalized)^k ( \Lunit \tander^{\Ntop} v^i) \gfour(\upmu \upepsilon_{ijk}\partial_j, \newuL) \\
			& \qquad \qquad = (\tander^{\Ntop} \vortrenormalized)^k \left\{ \alpha \newuL \tander^{\Ntop} v^i + \beta^A \nullgeop{x^A} \tander^{\Ntop} v^i \right\},
		\end{split}
	\end{align}
for some explicit constants $\alpha,\beta^A$. This null structure, which also holds for the dangerous boundary integrals on the outgoing characteristics $\outgoingcharacteristicsurfacetwoarg{u}{[\leftubar,\ubar]}$ in RHS\,\eqref{E:INTRODANGEROUSTOPORDERINTEGRALBUTAFTERSWITCHINGMUCURLANDL}, has little to do with compressible fluids and is instead a statement about swapping first order differential operators via integrating by parts on a region covered by a double-null foliation. We provide the full details of this newly observed null structure in Lemma\,\ref{L:IBPTWODIFFERENTIALOPERATORS}. The application of this method 
to the study of the Cauchy horizon is located in Sects.\,\ref{S:IBPIDENTITIESNEEDEDFORTOPORDERVORTICITY} and\,\ref{SS:L2ESTIMATESFORTOPORDERVORTICITYINRHSOFWAVEEQUATIONS}.

\paragraph{\emph{Analysis of the logarithmic density}} \label{PAR:ANALYSISOFTHELOGARITHMICDENSITY} \hfill
The previous section showed that using the sum of the velocity vectorfield components as the fundamental $L^2$-controlling quantities as in \eqref{E:INTROVELOCITYL2CONTROLLINGQUANTITIES} led to a perfect derivative of otherwise uncontrollable derivative losing error terms. Formally, this was possible because RHS\,\eqref{E:COVWAVEEQv} could be viewed schematically as $\Flatcurl \Flatcurl v$, and hence one of the $\Flatcurl$'s could be integrated by parts to the multiplier factor $\multipliervectorfield v$ in the multiplier energy method. Note, however, that \emph{this would never be possible for the logarithmic density}. This is because the derivative losing term on the RHS\,\eqref{E:COVWAVEEQRHO} is $\Flatdiv \GradEnt$, and hence one could never hope to find a perfect derivative by integrating by parts the $\Flatdiv$ to the multiplier factor $\multipliervectorfield \LogDensity$. 

In this paper, this is overcome by \emph{never deriving energy estimates for the logarithmic density}. In Sect.\,\ref{S:DENSITYDERIVATIVESINTERMSOFOTHERS}, we use the first-order formulation of the compressible Euler equations \eqref{E:INTRO:COMPRESSIBLEEULER} to derive algebraic identities of the form: 
	\begin{subequations} \label{E:DERIVATIVESOFLOGDENSITY}
		\begin{align}
			\tandergeneral \LogDensity & = \smoothfunction(\wavearray) \tandergeneral v+ \smoothfunction(\wavearray) (\vortrenormalized,\GradEnt) \label{E:INTROTANGENTIALDERIVATIVESOFLOGDENSITY} \\
			\muX \LogDensity & = \smoothfunction(\wavearray) (\muX v + \tandergeneral v) +  \smoothfunction(\wavearray) \GradEnt,\label{E:INTROINTROTRANSVERSEDERIVATIVESOFLOGDENSITY}
		\end{align}
	\end{subequations}
where $P \in \{\Lunit,\Yvf{2},\Yvf{3}\}$. Once we derive pointwise estimates for the velocity, vorticity, and gradient entropy, 
\eqref{E:DERIVATIVESOFLOGDENSITY} automatically implies corresponding estimates for the logarithmic density.

It is important to note that inserting the algebraic identities \eqref{E:DERIVATIVESOFLOGDENSITY} into the $\gfour$-null forms $\mathcal{Q}(\boldsymbol{\p}\vec{\Psi},\boldsymbol{\p}\vec{\Psi})$ on the RHS\,\eqref{E:COVWAVEEQv}--\eqref{E:MODIFIEDFLUIDVARIABLES}, and multiplying by $\upmu$, we obtain the following:
\begin{align}
		\upmu \mathcal{Q}(\boldsymbol{\p}\vec{\Psi},\boldsymbol{\p}\vec{\Psi}) \simeq \upmu \mathcal{Q}(\boldsymbol{\p}v,\boldsymbol{\p}v) + \smoothfunction(Z\wavearray) (\vortrenormalized,\GradEnt) +  \upmu \smoothfunction(\wavearray) (\vortrenormalized,\GradEnt)^2 \label{E:NULLFORMSREMAINNULLFORMS}
	\end{align}
A fruitful way to interpret \eqref{E:NULLFORMSREMAINNULLFORMS} is that null-forms of the full solution variables can be re-expressed as null-forms of the velocity up to harmless error terms. This is crucial because if this structure were not preserved, it would imply that the quadratic error terms on RHS\,\eqref{E:COVWAVEEQv}--\eqref{E:MODIFIEDFLUIDVARIABLES} are \emph{not} perturbative near the formation of the shock.

\paragraph{\emph{The need for the modified fluid variables $\VortVort$, $\DivGradEnt$, and elliptic-hyperbolic integral identities}} \label{PAR:THENEEDFORUSINGCANDD} 
As we mentioned at the start of Sect.\,\ref{SSS:INTROENERGYESTIMATES}, we do not derive energy estimates for the entropy $\Ent$ from the wave equation \eqref{E:COVWAVEEQS}. Although we \emph{could} derive a similar perfect derivative as in \eqref{E:THEGOODSIGN} by integrating the $\Flatdiv$ from the RHS\,\eqref{E:COVWAVEEQS} (as opposed to the logarithmic density equation \eqref{E:COVWAVEEQRHO}), for the sake of brevity, we instead derive energy estimates solely for its gradient $\GradEnt^i := \p_i \Ent$ by treating it as a transport variable (see \eqref{E:TRANSPORTVORTANDS}). 

The reader might then wonder why we even introduce the modified fluid variables $\VortVort$ and $\DivGradEnt$ in the first place. After all, \textbf{I)} $(\VortVort,\DivGradEnt)$ play no role in the transport equations \eqref{E:TRANSPORTVORTANDS} for $\vortrenormalized$ and $\GradEnt$; \textbf{II)}  for the velocity wave equations, $\VortVort$ is handled via the perfect derivative method of Sect.\,\ref{PAR:VELOCITYASAFUNDAMENTALUNKNOWN} and never estimated directly; \textbf{III)} $\Ent$ is not estimated via the wave equation and hence $\DivGradEnt$ seemingly does not enter into the equations. The reason we still have to control \emph{both} $\VortVort$ and $\DivGradEnt$ is that \emph{we still employ Raychaudhuri to control the top-order $\tander^{\Ntop}\mytr_{\gtorus}\upchi$ terms in \eqref{E:INTROSCHEMATICENERGYINEQUALITYFORVELOCITY}}. As we explained in Sect.\,\ref{SSS:GAUSSCURVATUREESTIMATESNOTREQUIREDIN2D}, the RHS of the wave equations \eqref{E:COVWAVEEQv}--\eqref{E:COVWAVEEQS} for \emph{all of the wave variables $(v,\LogDensity,\Ent)$} enter as source terms for the modified acoustic variable $\fullymodquant{\tander^{\Ntop}}$, see the $\tander^{\Ntop}(\vec{F})$ term on RHS \eqref{E:SCHEMATICRAYCHAUDHURIINTRO}--\eqref{E:SCHEMATICRAYCHAUDHURIFORMODIFIEDCHIINTRO}. As there is no way to find a perfect divergence between bilinear terms of the form $(\, \fullymodquant{\tander^{\Ntop}})( \tander^{\Ntop}(\VortVort,\DivGradEnt))$, we must still independently control $(\VortVort,\DivGradEnt)$ as fundamental unknowns. 

The most difficult error integral generated by commuting \,\eqref{E:MODIFIEDFLUIDVARIABLES} is the following:
\begin{align} 
		\begin{split} \label{E:INTROMOSTDIFFICULTINTEGRALFORCOMMUTEDCANDDEQUATIONS}
			& \int_{\characteristicdiamondtwoarg{[\leftubar,\ubar]}{[\moreinterestingu_1,u]}} \left\{ \pmb{\p}\tander^{\Ntop}(\vortrenormalized,\GradEnt)\right\} \left\{ \tander^{\Ntop}(\VortVort,\DivGradEnt)\right\} |\ubar|^{A + 1.5} \\
			& \lesssim \varsigma \int_{\characteristicdiamondtwoarg{[\leftubar,\ubar]}{[\moreinterestingu_1,u]}} \left| \pmb{\p}\tander^{\Ntop}(\vortrenormalized,\GradEnt)\right|^2 |\ubar|^{A + 1.5} + \varsigma^{-1}  \int_{\characteristicdiamondtwoarg{[\leftubar,\ubar]}{[\moreinterestingu_1,u]}} \left| \tander^{\Ntop}(\VortVort,\DivGradEnt)\right|^2  |\ubar|^{A + 1.5},
		\end{split}
	\end{align}
where $\varsigma$ is a small constant. Using Fubini, the second integral on RHS\,\eqref{E:INTROMOSTDIFFICULTINTEGRALFORCOMMUTEDCANDDEQUATIONS} can be expressed as the sum $\varsigma^{-1}\int_{u' = \moreinterestingu_1}^u \left\{ \outgoingfluxwave{\Ntop}[\VortVort](\ubar,u') + \outgoingfluxwave{\Ntop}[\DivGradEnt](\ubar,u') \right\} \rmd u'$, which can handled using Gr\"onwall's inequality. 

In order to control the first spacetime integral on RHS\,\eqref{E:INTROMOSTDIFFICULTINTEGRALFORCOMMUTEDCANDDEQUATIONS}, we derive a new localized elliptic-hyperbolic integral identity adapted to the characteristic diamond. This identity builds on the framework we developed in \cite{abbrescia2025remarkable}, valid for \emph{any} compact globally hyperbolic regions, as well as the one in \cite[Sect.\,21]{abbrescia2022emergence} which introduced the concept of \emph{characteristic currents}. We note already that the identities derived in \cite{abbrescia2025remarkable} used a different \emph{non}-characteristic current which, when translated to the current setting of this paper, would generate dangerous boundary integrals that are uncontrollable at first glance.\footnote{We believe that there are cancellations of the uncontrollable terms present in these boundary integrals, but do not investigate this carefully in this paper.} The characteristic current of \cite{abbrescia2022emergence} is tailored to the $\gfour$-null frame $\{\Lunit, \olduLunit\}$, where $\olduLunit = \Lunit + 2X$ and is transverse to the ingoing characteristics $\ingoingcharacteristicsurfacetwoarg{\ubar}{[\moreinterestingu_1,u]}$ except at the crease (where it is tangent), whereas the one in this paper is tailored to the $\gfour$-null frame $\{\Lunit, \uLunit\}$, where $\uLunit$ is everywhere tangent to the ingoing characteristics $\ingoingcharacteristicsurfacetwoarg{\ubar}{[\moreinterestingu_1,u]}$.

For convenience, we denote $\SigmatTan = \tander^{\Ntop}\vortrenormalized$ or $\tander^{\Ntop} \GradEnt$.  In Sect.\,\ref{S:ELLIPTICHYPERBOLICIDENTITIES}, we construct a $\outgoingcharacteristicsurfacetwoarg{u}{[\leftubar,\ubar]}$-tangent vectorfield $\ehcurrent^\alpha$ (hence the name \emph{characteristic} current) depending on $\SigmatTan$ and $\pmb{\p}\SigmatTan$ that satisfies the following crucial identity, expressed schematically:
\begin{align} \label{E:INTROSCHEMATICEHCURRENTDIVERGENCE}
		\Dfour_\alpha \ehcurrent^\alpha = \ellipticCoerciveQuadratic(\pmb{\p}\SigmatTan,\pmb{\p}\SigmatTan) + (\Flatcurl \SigmatTan)^2 + (\Flatdiv \SigmatTan)^2 + (\Transport \SigmatTan)^2 + \cdots,
	\end{align}
where $\Dfour$ is the Levi-Civita connection of $\gfour$, $\ellipticCoerciveQuadratic(\pmb{\p}\SigmatTan,\pmb{\p}\SigmatTan)$ is a coercive quadratic form satisfying $\ellipticCoerciveQuadratic(\pmb{\p}\SigmatTan,\pmb{\p}\SigmatTan) \approx |\pmb{\p} \SigmatTan|^2$, and $\cdots$ are easier terms which are at most linear in $\pmb{\p} \SigmatTan$. Part of the detailed structure of $\ehcurrent$ implies that  all second derivatives $\pmb{\p}^2 \SigmatTan$ cancel in $\Dfour_\alpha \ehcurrent^\alpha$ when deriving  \eqref{E:INTROSCHEMATICEHCURRENTDIVERGENCE}. Upon integrating \eqref{E:INTROSCHEMATICEHCURRENTDIVERGENCE} over the characteristic diamond, the coercivity of $\ellipticCoerciveQuadratic(\pmb{\p}\SigmatTan,\pmb{\p}\SigmatTan)$ implies that we are able to bound the first integral on RHS\,\eqref{E:INTROMOSTDIFFICULTINTEGRALFORCOMMUTEDCANDDEQUATIONS} so long as we can control the spacetime integral of $(\Flatcurl \SigmatTan)^2 + (\Flatdiv \SigmatTan)^2 + (\Transport \SigmatTan)^2 + \cdots$ and the boundary terms generated from the divergence theorem applied to the LHS\,\eqref{E:INTROSCHEMATICEHCURRENTDIVERGENCE}. 

The spacetime integrals of $(\Flatcurl \SigmatTan)^2 + (\Flatdiv \SigmatTan)^2 + (\Transport \SigmatTan)^2 + \cdots$ are not difficult to bound in a way that is consistent with the $|\ubar|$-weighted hierarchy by using directly the div-curl transport system of \eqref{E:TRANSPORTVORTANDS}--\eqref{E:MODIFIEDFLUIDVARIABLES}, see also \eqref{E:ANTISYMMETRICVORTICITYELLIPTICHYPERBOLICIDENTITYERORTERMPOINTWISE}--\eqref{E:NULLGEOMETRYENTROPYGRADIENTELLIPTICHYPERBOLICIDENTITYERORTERMPOINTWISE} for the details. 

The most delicate part of the analysis, and arguably of the entire paper, is controlling the boundary terms upon applying the divergence theorem to LHS\,\eqref{E:INTROSCHEMATICEHCURRENTDIVERGENCE}. Note firstly that as $\ehcurrent$ is $\outgoingcharacteristicsurfacetwoarg{u}{[\leftubar,\ubar]}$-tangent, there are no boundary terms on the outgoing characteristics because of the vanishing contraction $\ehcurrent^\alpha \Lunit_\alpha = 0$. Hence, we need only focus on the $\ingoingcharacteristicsurfacetwoarg{\ubar}{[\moreinterestingu_1,u]}$ boundary terms, which feature the contraction$\ehcurrent^\alpha \newuL_\alpha$.  As we noted above, $J \simeq V \cdot \pmb{\p} V$ and so these \emph{hypersurface} integrals could potentially be fatal, as they would require us to control   $\pmb{\p}\tander^{\Ntop}(\vortrenormalized,\GradEnt)$ terms on a hypersurface. In Sect.\,\ref{SS:CONTRACTIONOFCHARACTERISTICCURRENTSANDNEWUL}, using the definition of the modified fluid variables \eqref{E:MODIFIEDCURLOFVORTICITYINTRO}--\eqref{E:MODIFIEDDIVERGENCEOFENTROPYGRADIENTINTRO}, the first order formulation of the equations \eqref{E:INTRO:COMPRESSIBLEEULER}, and careful geometric decompositions of various tensors, we derive the following miraculous null structure (expressed schematically):
	\begin{align} \label{E:MIRACULOUSNULLSTRUCTUREINTHEBOUNDARYTERMS}
		\ehcurrent^\alpha \newuL_\alpha & = - \newuL \left( \coercivequadraticformintoriboundaryterms(\SigmatTan,\SigmatTan) \right) + \nullangDiv\left(( \SigmatTan\cdot\SigmatTan) \breve{\slashed{\Transport}} \right) +  \SigmatTan\cdot \tander^{\Ntop}(\VortVort,\DivGradEnt) + \SigmatTan \cdot \newuL \tander^{\Ntop} v + \SigmatTan \cdot \breve{\slashed{\Transport}}\tander^{\Ntop} v,
	\end{align}\noindent where $\coercivequadraticformintoriboundaryterms(\SigmatTan,\SigmatTan)  \approx |\SigmatTan|^2_e$ (where $e$ is the three dimensional Euclidean metric) is a coercive quadratic form, $ \breve{\slashed{\Transport}}$ denotes the $\gfour$-orthogonal projection of the material vectorfield $\Transport$ onto $\doublenulltoritwoarg{\ubar}{u}$, and $\nullangDiv$ denotes the divergence with respect to the Riemannian metric $\gnulltori$  induced on $\doublenulltoritwoarg{\ubar}{u}$ by $\gfour$. In total, \emph{all} derivatives of $\tander^{\Ntop} v$ on RHS \eqref{E:MIRACULOUSNULLSTRUCTUREINTHEBOUNDARYTERMS} are $\ingoingcharacteristicsurfacetwoarg{\ubar}{[\moreinterestingu_1,u]}$-tangent, which is necessary to control them by the velocity fluxes.

The derivatives of $\SigmatTan$, which are the principal terms from the point of view of regularity are either: \textbf{I)} a perfect $\nullangDiv$-divergence, which vanishes upon integrating; \textbf{II)} the modified fluid variables $\tander^{\Ntop}(\VortVort, \DivGradEnt)$, which are easy to estimate; \textbf{III)} feature a $\newuL$ derivative, and are the most delicate boundary terms. The crucial observation is that $\newuL \left( \coercivequadraticformintoriboundaryterms(\SigmatTan,\SigmatTan) \right)$ is a \emph{perfect $\newuL$ derivative} which can be integrated by parts (one can think of $\newuL \approx \nullgeop{u}$), and, miraculously, \emph{features a good sign} so that the corresponding co-dimension two double-null tori integrals are: 
\begin{align} \label{E:THEMIRACULOUSSIGN}
		 - \int_{\ingoingcharacteristicsurfacetwoarg{\ubar}{[\moreinterestingu_1,u]}} \newuL \left(\coercivequadraticformintoriboundaryterms(\SigmatTan,\SigmatTan) \right) = - \int_{\doublenulltoritwoarg{\ubar}{u}} \coercivequadraticformintoriboundaryterms(\SigmatTan,\SigmatTan)  + \underbrace{ \int_{\doublenulltoritwoarg{\ubar}{U_0}} \coercivequadraticformintoriboundaryterms(\SigmatTan,\SigmatTan) }_{\textnormal{= data}}.
	\end{align}
The upshot is that the first error integral on RHS\,\eqref{E:THEMIRACULOUSSIGN} can be brought to the other side of the spacetime integral of the divergence identity \eqref{E:INTROSCHEMATICEHCURRENTDIVERGENCE}, meaning the elliptic-hyperbolic integral identities provide control over: 
\begin{align} \label{E:WHATTHEELLIPTICHYPERBOLICIDENTITIESCONTROL}
		\int_{\characteristicdiamondtwoarg{[\leftubar,\ubar]}{[\moreinterestingu_1,u]}} \ellipticCoerciveQuadratic(\pmb{\p}\SigmatTan,\pmb{\p}\SigmatTan) +  \int_{\doublenulltoritwoarg{\ubar}{u}} \coercivequadraticformintoriboundaryterms(\SigmatTan,\SigmatTan).
	\end{align}
 
\begin{remark}[Final remarks on the elliptic-hyperbolic integral identities] \label{R:FINALREMARKONINTEGRALIDENTITIES}
	A direct analog of the remarkable null structure from \eqref{E:MIRACULOUSNULLSTRUCTUREINTHEBOUNDARYTERMS}, including the miraculous good sign, was present in our paper \cite{abbrescia2022emergence}. However, in that paper, the boundary terms were on the $\gfour$-\emph{spacelike} rough hypersurfaces $\hypthreearg{\timefunction}{[-\rightu,u]}{0}$ unlike the $\gfour$-null ingoing characteristics $\ingoingcharacteristicsurfacetwoarg{\ubar}{[\moreinterestingu_1,u]}$. Since the null-fluxes are only semi-definite and control only tangential derivatives of the velocity, while the energies on the rough hypersurfaces are positive definite in all derivatives of the velocity, we needed to uncover new null structures from highly non-trivial tensorial decompositions to get the structure featured in $ \SigmatTan \cdot \newuL \tander^{\Ntop} v + \SigmatTan \cdot \breve{\slashed{\Transport}}\tander^{\Ntop} v$.  This highly technical aspect of the analysis was proved in Sect.\,\ref{SS:GEOMETRICSTRUCTUREOFEXTERIORDERIVATIVES}--\ref{SS:CONTRACTIONOFCHARACTERISTICCURRENTSANDNEWUL}.

In this introductory section, we suppressed the role of the $\upmu$-weights and $|\ubar|$-weights. This aspect requires careful attention as it is in fact the weights of some of the boundary terms described above which dictate the $|\ubar|^A$, and $|\ubar|^{A+3/2}$ hierarchy between the velocity $v$ and modified fluid variables $(\VortVort,\DivGradEnt)$, see \eqref{E:INTROVELOCITYL2CONTROLLINGQUANTITIES} and \eqref{E:INTROVORTVORTANDDIVGRADENTL2CONTROLLINGQUANTITIES}. To be clear, the constant $A$ is completely determined by the velocity estimates. Here, we mean that the careful analysis of the weights suppressed in this introduction prove that the difference between the $|\ubar|$-weight of the top-order $(\VortVort,\DivGradEnt)$ and the velocity $v$ $L^2$-controlling quantities is always $\tfrac{3}{2}$, regardless of $A$.  The details can be found int Sects.\,\ref{S:POINTWISEESTIMATESFORCONTROLLINGSPECIFICVORTICITYANDENTROPYGRADIENT} and\,\ref{S:TOPORDERELIPTICHYPERBOLICL2ESTIMATESFORSPECIFICVORTICITYANDENTROPYGRADIENT}. 
 \end{remark}

\subsection{Abbreviated version of the main results} \label{SS:ABBREVIATEDMAINRESULTS}
In this section, we provide an abbreviated statement of the main results; see Theorem\,\ref{T:EXISTENCEUPTOCAUCHYHORIZONBYCONTINUATIONCRITERIA} 
for a more detailed version. We begin by providing the precise definition of the Cauchy horizon.

\begin{definition}[The Cauchy horizon] \label{D:THECAUCHYHORIZON} Let $\Cartesiantisafunctiononmumxtoriarg{0}{0},\Eikonalisafunctiononmumuxtoriarg{0}{0} : \T^2 \to \R$ be the defining scalar functions for describing the crease as a graph, namely  $\crease = \twoargmumuxtorus{0}{0} = \{ \left(\Cartesiantisafunctiononmumxtoriarg{0}{0}(x^2,x^3),\Eikonalisafunctiononmumuxtoriarg{0}{0}(x^2,x^3),x^2,x^3 \right) \, : \, (x^2,x^3) \in \T^2\}$. Then we define the  \textbf{Cauchy horizon} to be the following
subset of spacetime, where the null hypersurface portion $\ingoingcharacteristicsurfacetwoarg{0}{[\moreinterestingu_1,\moreinterestingu_2]}$ 
is defined in  \eqref{E:INTROSPACETIMEREGIONS}:
	\begin{align}
		\Cauchyhor \eqdef \left\{ (t,u,x^2,x^3)  \in  \ingoingcharacteristicsurfacetwoarg{0}{[\moreinterestingu_1,\moreinterestingu_2]} \, \Big| \, u > \Eikonalisafunctiononmumuxtoriarg{0}{0}(x^2,x^3) \, \textnormal{for all $(x^2,x^3) \in \T^2$}\right\}. \label{E:INTROCAUCHYHORIZON}
	\end{align}
\end{definition}

	\begin{theorem}[Abbreviated version of the main results] \label{T:ABBREVIATEDMAINRESULTS}
		\hfill
		
		\begin{itemize} 
			\item \underline{\textbf{Assumptions.}} Let $(\LogDensity^1|_{\Sigma_0},v^1|_{\Sigma_0}, v^2|_{\Sigma_0},v^3|_{\Sigma_0},\Ent|_{\Sigma_0}) \in C^\infty(\Sigma_0)^{\times 5}$ be  smooth functions satisfying the same assumptions (e.g. quantitative transversal convexity, Sobolev regularity, etc.) as in \cite{abbrescia2022emergence}. In particular, there exists an acoustical eikonal function $u$ solving the eikonal equation $(\gfour^{-1})^{\alpha\beta}\p_\alpha u \p_\beta u = 0$, and a rough time function 
				$\timefunctionarg{0} = \timefunctionarg{0}(t,u,x^2,x^3)$ solving $\muX \timefunctionarg{0} =0$
				with range $[\timefunction_0,0] = [-\mupositive,0]$ such that $(\timefunctionarg{0},u,x^2,x^3)$ form a global coordinate system on 
				$\twoargMrough{[\timefunction_0,0],[- \rightu,\leftu]}{0}
				= \bigcup_{\timefunction \in [\timefunction_0,0]} \hypthreearg{\timefunction}{[- \rightu,\leftu]}{0}
				$. This is a compact portion of the maximal classical development 
				of the data up to its boundary with respect to the differential structure of the geometric coordinates $(t,u,x^2,x^3)$, where $\rightu,\leftu > 0$ depend only on the background solution. With respect to this differential structure of the geometric coordinates, the $C^{3,1}$-norm of the solution variables $\wavearray = (\RRiemann,\LRiemann,v^2,v^3,\Ent)$ remain bounded on $\twoargMrough{[\timefunction_0,0],[- \rightu,\leftu]}{0}$, and $\upmu > 0$ everywhere on  $\twoargMrough{[\timefunction_0,0],[- \rightu,\leftu]}{0}$ except for the crease $\twoargmumuxtorus{0}{0} = \crease$, where it vanishes: $\upmu|_{\twoargmumuxtorus{0}{0}} = 0$. The change of variables map $\Upsilon(t,u,x^2,x^3) = (t,x^1,x^2,x^3)$ is a diffeomorphism everywhere on $\twoargMrough{[\timefunction_0,0],[- \rightu,\leftu]}{0}$ onto its image except for the crease $ \crease$, where it is merely a homeomorphism. There is gradient blow-up along $\crease$ relative to Cartesian differential structure of $(t,x^1,x^2,x^3)$ in that the following estimate holds on $\Upsilon\left(\twoargMrough{[\timefunction_0,0],[- \interestingu,\interestingu]}{0} \right)$, where $0 < \interestingu <  \min\{1,\leftu,\rightu,\frac{1}{22}\blowupdelta^{-1}\}$:
				\begin{align}
					|\partial_{x^\alpha} \RRiemann|, |\partial_{x^\alpha} v^1|,|\partial_{x^\alpha} \LogDensity| \gtrsim \frac{1}{\upmu}. \label{E:ABBREVIATEDTHEOREMSHOCKFORMATIONONCREASE}
				\end{align}
				The solution is smooth on 	$\Upsilon\left(\twoargMrough{[\timefunction_0,0],[- \interestingu,\interestingu]}{0} \setminus \crease \right)$, while for $\alpha = 0,1,2,3$, $i=1,2,3$, and $A=2,3$,
		the fluid quantities 
		$\LogDensity$, 
		$v^i$, 
		$\Ent$, 
		$\Flatcurl v^i$, 
		$\partial_i \Ent$, 
		and 
		$\gfour_{ab} \Yvf{A}^a \partial_{\alpha} v^b$
		remain bounded on all of $\Upsilon(\twoargMrough{[\timefunction_0,0],[- \interestingu,\interestingu]}{0})$. 
			\item \underline{\textbf{Ingoing eikonal function, double-null coordinates, and the characteristic diamond.}}
				There exists an \emph{ingoing} acoustical eikonal function $\ubar$ solving $(\gfour^{-1})^{\alpha\beta}\p_\alpha \ubar \p_\beta \ubar = 0$ with initial data $\ubar|_{\datahypfortimefunctiontwoarg{0}{[\timefunction_0,0]}} = - \upmu|_{\datahypfortimefunctiontwoarg{0}{[\timefunction_0,0]}}$. With respect to the differential structure of the geometric coordinates, $\ubar \in C^{1,1}$. The  range of $\ubar$ is $[\leftubar,0]$, where $\leftubar \in [\frac{3}{4}\timefunction_0,\frac{1}{4}\timefunction_0]$. The level sets $\ingoingcharacteristicsurfacetwoarg{\ubar'}{[-\rightu,\leftu]} := \{(t,u,x^2,x^3) \, : \, \ubar(t,u,x^2,x^3) = \ubar', \, u_1\le u\le u_2\}$ of $\ubar$ are $\gfour$-null hypersurfaces. There exists $\moreinterestingu_1 \in (-\interestingu,-\tfrac{3}{4} \interestingu)$ and $ \moreinterestingu_2 \in (\frac{3}{4} \interestingu,\interestingu)$ such that $(\ubar,u,x^2,x^3)$ form a global coordinate system on the characteristic diamond $\characteristicdiamondtwoarg{[\leftubar,0]}{[\moreinterestingu_1,\moreinterestingu_2]}$. 
			\item \underline{\textbf{The singular behavior of the fluid on the characteristic diamond.}} The crease is contained in the final ingoing characteristic, that is, $\crease \subset \ingoingcharacteristicsurfacetwoarg{0}{[\moreinterestingu_1,\moreinterestingu_2]}$. In particular, $\upmu > 0$ everywhere on $\characteristicdiamondtwoarg{[\leftubar,0]}{[\moreinterestingu_1,\moreinterestingu_2]}$ except for the crease, where it vanishes.
			\item \underline{\textbf{Cauchy horizon emanates from the crease.}}
			The following estimate holds:
				\begin{align}
					\moreinterestingu_2 - \sup_{\twoargmumuxtorus{0}{0}} u \ge \frac{1}{4}\interestingu. \label{E:INTRONONTRIVIALSETOFUVALUESINLEVELSETOFU}
				\end{align}
			Namely, \eqref{E:INTROCAUCHYHORIZON} and \eqref{E:INTRONONTRIVIALSETOFUVALUESINLEVELSETOFU} imply that the portion of the Cauchy horizon contained in $ \ingoingcharacteristicsurfacetwoarg{0}{[\moreinterestingu_1,\moreinterestingu_2]}$ is non-trivial. In addition, the following identities hold:
			\begin{subequations}
				\begin{align}
					\crease & =  \left\{ p   \in  \ingoingcharacteristicsurfacetwoarg{0}{[\moreinterestingu_1,\moreinterestingu_2]} \, \Big| \upmu(p) = \newuL \upmu = 0\right\}, \label{E:INTROCREASEWITHRESPECTTOULUNIT}  \\
					\Cauchyhor & =  \left\{ p   \in  \ingoingcharacteristicsurfacetwoarg{0}{[\moreinterestingu_1,\moreinterestingu_2]} \, \Big| \newuL \upmu (p) > 0\right\}, \label{E:INTROCAUCHYHORWITHRESPECTTOULUNIT}
				\end{align}
			\end{subequations}
			where $\newuL^\alpha $ is the future-directed $\gfour$-null vectorfield proportional to $(\gfour^{-1})^{\alpha\beta}\p_\beta \ubar$ normalized by $\newuL u = 1$. Hence, \eqref{E:INTROCREASEWITHRESPECTTOULUNIT}--\eqref{E:INTROCAUCHYHORWITHRESPECTTOULUNIT} imply that $\Cauchyhor$ emanates from the crease $\crease$, which the Cauchy horizon's past boundary.
			\item \underline{\textbf{Behavior of the change of variables map.}} The change of variables map $\Upsilon(t,u,x^2,x^3)=(t,x^1,x^2,x^3)$ is a diffeomorphism everywhere onto the image of the characteristic diamond $\characteristicdiamondtwoarg{[\leftubar,0]}{[\moreinterestingu_1,\moreinterestingu_2]}$, except for the crease, where it is a homeomorphism. 
			\item \underline{\textbf{Smooth extension up to the Cauchy horizon in Cartesian coordinates.}} Relative to the Cartesian coordinate differential structure, the  fluid solution is smooth on $\Upsilon\left( \characteristicdiamondtwoarg{[\leftubar,0]}{[\moreinterestingu_1,\moreinterestingu_2]}\setminus\crease\right)$. In particular, it is smooth up to the Cauchy horizon.
						 
		\end{itemize}
	\end{theorem}

\section{Basic setup, compressible Euler flow, and its geometric reformulation} 
\label{S:COMPRESSIBLEEULERFLOWANDITSGEOMETRIC}

In this section, we first introduce some basic notational conventions and definitions. 
We then provide a standard first-order quasilinear hyperbolic formulation of compressible Euler flow.
Next, we define a series of additional fluid variables and geometric tensors associated to the flow.
Finally, we recall the new formulation of the flow derived in \cite{jS2019c}. 

\subsection{Basic notation and conventions}
\label{SS:NOTATIONANDCONVENTIONS}
The precise definitions of some of the concepts referred to
here are provided later in the article.

\begin{itemize}
	\item  (\textbf{Cartesian coordinates}) 
	Our analysis takes place on subsets of the spacetime manifolds 
	$\mathbb{R} \times \Sigma$,
	where $\Sigma \eqdef \mathbb{R} \times \mathbb{T}^2$ is the spatial manifold. 
	We fix a standard Cartesian coordinate system $\lbrace x^{\alpha} \rbrace_{\alpha = 0,1,2,3}$
on $\mathbb{R} \times \Sigma$, where $t \eqdef x^0 \in \mathbb{R}$
is the time coordinate and $(x^1,x^2,x^3) \in \Sigma$
are the spatial coordinates (where $(x^2,x^3)$ are standard coordinates\footnote{While the coordinates $x^2,x^3$ on $\mathbb{T}^2$ are only locally defined, the corresponding partial derivative vectorfields $\partial_2,\partial_3$ can be extended
so as to form a global smooth frame on $\mathbb{T}^2$. Similar remarks apply to the one-forms $\rmd x^2,\rmd x^3$
These simple observations are relevant for this paper because when we derive estimates, the coordinate functions $x^2,x^3$
themselves are never directly relevant; what matters are estimates for the components of various tensorfields
with respect to the frame $\lbrace \partial_t, \partial_1, \partial_2, \partial_3 \rbrace$
and the basis dual co-frame $\lbrace \rmd t, \rmd x^1, \rmd x^2, \rmd x^3 \rbrace$, which are everywhere smooth. 
\label{FN:COORDINATESARENOTGLOBAL}}  
on $\mathbb{T}^2$). We sometimes refer to $t$ as the ``Cartesian time function.'' 
$\Sigma_{t'} \eqdef \lbrace (t,x^1,x^2,x^3) \in \R \times \R \times \T^2 | \ t = t' \rbrace$
denotes the standard flat hypersurface of constant Cartesian time.
\item (\textbf{Cartesian coordinate partial derivatives}) 
We use the notation $\lbrace \partial_{\alpha} \rbrace_{\alpha = 0,1,2,3}$ 
(or $\partial_t \eqdef \partial_0$) to denote the Cartesian coordinate partial derivative
vectorfields. 
	\item (\textbf{Lowercase Greek index conventions}) Lowercase Greek spacetime indices 
	$\alpha$, $\beta$, etc.\ correspond to the Cartesian coordinate spacetime coordinates
	and vary over $0,1,2,3$. All lowercase Greek indices are lowered and raised with the acoustical metric
	$\gfour$ (see definition~\eqref{E:ACOUSTICALMETRIC}) 
	and its inverse $\gfour^{-1}$, and \emph{not with the Minkowski metric}.
	Throughout the article, if $\upxi$ is a type $\binom{m}{n}$ spacetime tensorfield,
	then unless we indicate otherwise, in our identities and estimates,
	\textbf{$\left\lbrace \upxi_{\beta_1 \cdots \beta_n}^{\alpha_1 \cdots \alpha_m} \right\rbrace_{
	\alpha_1, \cdots \alpha_m, \beta_1, \cdots, \beta_n = 0,1,2,3}$ denotes its 
	components with respect to the Cartesian coordinates}. 
	This is important because some of our identities
	and estimates hold only with respect to the Cartesian coordinates.
	\item (\textbf{Lowercase Latin index conventions})
	Lowercase Latin spatial indices $a$, $b$, etc.\ correspond to the Cartesian spatial coordinates and vary over $1,2,3$.
	Much like in the previous point, 
	if $\upxi$ is a type $\binom{m}{n}$ $\Sigma_t$-tangent tensorfield 
	(see Def.~\ref{D:PROJECTIONTENSORFIELDSANDTANGENCYTOHYPERSURFACES}),
	then unless we indicate otherwise, in our identities and estimates,
	\textbf{$\left\lbrace \upxi_{b_1 \cdots b_n}^{a_1 \cdots a_m} \right\rbrace_{a_1,\cdots,a_m,b_1,\cdots,b_n = 1,2,3}$ 
	denotes its components with respect to the Cartesian spatial coordinates}. 
	\item (\textbf{Uppercase Latin index conventions})
	Uppercase Latin spatial indices $A,B$, etc.\ correspond to the coordinates 
	$(x^2,x^3)$ on $\mathbb{T}^2$ and vary over $2,3$. In particular, if $V$ is a vectorfield,
	and $A \in \lbrace 2,3 \rbrace$, then $V^A = V^{\alpha} \partial_{\alpha} x^A$, 
	where $(x^2,x^3)$ are the standard Cartesian coordinates on $\mathbb{T}^2$.
\item (\textbf{tilded indices}) We use tilded indices such as $\widetilde{\alpha}$ in the same way as their
	non-tilded counterparts.
\item (\textbf{Einstein summation}) 
	We use Einstein's summation convention in that repeated indices are summed,
	e.g., $\Lunit^A X^A \eqdef \Lunit^2 X^2 + \Lunit^3 X^3$.
\item (\textbf{Use of ``$\cdot$''})  
		We sometimes use ``$\cdot$'' to denote the natural contraction between two tensors. 
		For example, if $\upxi$ is a spacetime one-form and $V$ is a 
		spacetime vectorfield,
		then $\upxi \cdot V \eqdef \upxi_{\alpha} V^{\alpha}$.
		At other times, we use ``$\cdot$'' to schematically denote products,
		e.g., $A_1 \cdot A_2 \cdot A_3$ is a trilinear form in $A_1, A_2, A_3$.
\item (\textbf{Tensor contractions})  
	If $V$ and $W$ are vectorfields, then $V_W \eqdef V^{\alpha} W_{\alpha} = \gfour_{\alpha \beta} V^{\alpha} W^{\beta}$.
	If $\upxi$ is a one-form and $V$ is a vectorfield, then $\upxi_V \eqdef \upxi_{\alpha} V^{\alpha}$.
	We use similar notation when contracting higher-order tensorfields against vectorfields.
	For example, if $\upxi$ is a type $\binom{0}{2}$ tensorfield and
	$V$ and $W$ are vectorfields, then $\upxi_{VW} \eqdef \upxi_{\alpha \beta} V^{\alpha} W^{\beta}$.
\item (\textbf{Commutator of operators})  If $Q_1$ and $Q_2$ are two operators, then $[Q_1,Q_2] \eqdef Q_1 Q_2 - Q_2 Q_1$ denotes their commutator. 
\item (\textbf{Constants}) We establish conventions for constants (such as ``$C$'') in Sect.~\ref{SS:CONVENTIONSFORCONSTANTS}.

\end{itemize}

\subsection{Basic differential operators}
\label{SS:BASICDIFFERENTIALOPERATORS}
In our analysis, we will encounter many kinds of differential operators.
Here, we define some basic operators.

\begin{definition}[Gradient one-form of a scalar function]
	\label{D:GRADIENTONEFORMOFSCALARFUNCTION}
	If $f$ is a scalar function, then $\rmd f$
	denotes the gradient one-form associated to $f$,
	e.g., $(\rmd f)_{\alpha} \eqdef \rmd f \cdot \partial_{\alpha} = \partial_{\alpha} f$.
\end{definition}	

\begin{definition}[Vectorfield derivative of scalar functions]
	\label{D:VECTORFIELDDERIVATIVEOFSCALARFUNCTION}
	If $V$ is a vectorfield and $f$ is a scalar function, then 
	$V f \eqdef V^{\alpha} \partial_{\alpha} f = V \cdot \rmd f$
	denotes the derivative of $f$ in the direction $V$.
\end{definition}

\begin{definition}[Euclidean divergence and curl]
\label{D:EUCLIDEANDIVERGENCEANDCURL}
$\Flatdiv$ and $\Flatcurl$ respectively denote the Euclidean spatial divergence and curl operators. 
That is, given a $\Sigma_t$-tangent 
vectorfield $V = V^a \partial_a$, we define, 
relative to the Cartesian spatial coordinates, 
$\Flatdiv V$ and $\Flatcurl V$ to be the following
scalar function and $\Sigma_t$-tangent vectorfield:
\begin{align} \label{E:FLATDIVANDCURL}
		\Flatdiv V
		& \eqdef \partial_a V^a,
		&
		(\Flatcurl V)^i
		& \eqdef \upepsilon_{iab} \partial_a V^b,
	\end{align}
	where $\upepsilon_{iab}$ is the fully antisymmetric symbol normalized by $\upepsilon_{123} = 1$.
\end{definition}

\subsection{A first-order formulation involving the logarithmic density}
\label{SS:FIRSTORDERFORMULATIONWITHLOGDENSITY}

\subsubsection{The logarithmic density, assumptions on the equation of state, and normalizations}
\label{SSS:LOGDENSITYASSUMPTIONSONEOSANDNORMALIZATIONS}
We find it convenient to work with the logarithmic density featured in the next
definition, rather than the density.
In the rest of the paper,
	\begin{align} \label{E:BACKGROUNDDENSITY}
		\overline{\varrho} & > 0
	\end{align}
	denotes a fixed constant ``background density.''
	
\begin{definition}[Logarithmic density] 
\label{D:LOGDENS}
 We define the \emph{logarithmic density} $\LogDensity$ as follows:
\begin{align} \label{E:LOGDENS}
	\LogDensity 
	& \eqdef \ln\left(\varrho / \overline{\varrho}\right).
	\end{align}
\end{definition}

In the rest of the paper,
we view the speed of sound $\Speed$ (which is defined in \eqref{E:SOUNDSPEED}) 
to be a function of $(\LogDensity,\Ent)$.
Note that by \eqref{E:SOUNDSPEED} and the chain rule, 
we have $\Speed(\LogDensity,\Ent)
= \sqrt{(\overline{\varrho})^{-1} \exp(-\LogDensity) p_{;\LogDensity}}$,
where $p_{;\LogDensity} \eqdef \tfrac{\p p}{\p \LogDensity}$ denotes the partial derivative of the equation of state 
with respect to the logarithmic density at fixed $\Ent$. 

\begin{notation}[Partial differentiation with respect to state-space variables]
\label{N:PARTIALDIFFERENTIATIONWITHRESPECTTOSTATESPACEVARIABLES}
In accordance with the above notation, 
for any scalar function $f = f(\LogDensity,\Ent)$, we use the notation 
$f_{;\LogDensity} \eqdef \frac{\partial f}{\partial \LogDensity}$ to denote
the partial derivative of $f$ 
with respect to the logarithmic density at fixed $\Ent$.
Similarly, we denote the partial derivative of $f$ with respect to $\Ent$ at fixed $\LogDensity$
by $f_{;\Ent} \eqdef \frac{\partial f}{\partial \Ent}$.
We also write $ f_{;\LogDensity;\Ent} \eqdef \frac{\partial^2 f}{\partial \Ent \partial \LogDensity}$, 
and we use similar notation for other higher partial derivatives of $f$ with respect to $\LogDensity, \Ent$.
\end{notation}

To ensure that shocks occur for solutions near static isentropic fluid states with constant density $\overline{\varrho} > 0$,
we assume the following non-degeneracy condition:
\begin{align}\label{E:NONDEGENCONDITION}
\overline{\Speed}^{-1} \overline{\Speed_{;\uprho}} + 1 
& \neq 0,
\end{align}
where LHS \eqref{E:NONDEGENCONDITION} is defined to be the constant
obtained by evaluating $\Speed^{-1} \Speed_{;\LogDensity} + 1$ at $\LogDensity = \Ent \equiv 0$.
Equation \eqref{E:NONDEGENCONDITION} ensures that the null condition fails to hold for perturbations
of the background solution $\LogDensity = \Ent \equiv 0$; see Sect.\,\ref{SS:THEFACTORDRIVINGTHESHOCKFORMATION}.
Our main results hold for all equations of state except for that of a Chaplygin gas,
namely $p = C_0 - C_1 \exp(- \LogDensity)$, where $C_0 \geq 0$ and $C_1 > 0$ are constants.
This equation of state is degenerate in the following sense: $\Speed^{-1} \Speed_{;\LogDensity} + 1 \equiv 0$.

By rescaling Cartesian time if necessary, we can assume the following convenient normalization condition:
\begin{align} \label{E:BACKGROUNDSOUNDSPEEDISUNITY}
	\Speed(\LogDensity = 0, \Ent = 0) 
	& = 1.
\end{align}
\subsubsection{A first-order formulation involving the logarithmic density}
\label{SSS:FIRSTORDERFORMULATIONINVOLVINGLOGDENSITY}
From definition~\eqref{E:LOGDENS}, equations \eqref{E:INTROTRANSPORTVI}--\eqref{E:INTROBS}, and the chain rule,
it follows that
relative to the standard Cartesian coordinates on $\mathbb{R} \times \Sigma$,
the compressible Euler equations can be expressed as the following system
in $\LogDensity$, $v$, and $\Ent$:
 \begin{subequations} \label{E:FIRSTORDERFORMULATION}
\begin{align}
	\Transport v^i 
	& = 
	- \Speed^2 \updelta^{ia} \partial_a \LogDensity
	- \exp(-\LogDensity) \frac{p_{;\Ent}}{\overline{\varrho}} \updelta^{ia} \partial_a \Ent,
	\label{E:BVIEVOLUTION}
		\\
		 \label{E:BLOGDENSITYEVOLUTION}
	\Transport \LogDensity
	& = - \Flatdiv v,
		\\
	\Transport \Ent
	& = 0.
	\label{E:BENTROPYEVOLUTION}
\end{align}
\end{subequations}\subsection{The almost Riemann invariants}
\label{SS:ALMOSTRIEMANNINVARIANTSANDLFLUIDVARIABLEARRAYS}
To study solutions close to simple isentropic plane-symmetric solutions, 
we find it convenient to replace
$\LogDensity$ and $v^1$ with a pair of ``almost Riemann invariants,''
denoted by $\RRiemann$ and $\LRiemann$. 


\begin{definition}[The almost Riemann invariants] \label{D:ALMOSTRIEMANNINVARIANTS}
We define the \emph{almost Riemann invariants away from symmetry} $\almostRiemann_{(\pm)}$ as follows:
\begin{align} \label{E:ALMOSTRIEMANNINVARIANTS}
\almostRiemann_{(\pm)} 
& \eqdef 
v^1 
\pm 
\almostRiemannfunction(\LogDensity,\Ent),
& 
\mbox{where } 
\almostRiemannfunction(\LogDensity,\Ent)
&
\eqdef 
\int_0^{\LogDensity} 
	\Speed(\LogDensity',\Ent) 
\, \rmd \LogDensity'.
\end{align}\end{definition}

\begin{remark}[\textbf{Clarification on our approach to estimating $\LogDensity$ and $v^1$}]
	\label{R:HOWWEESTIMATEDENSITYANDV1}
	We have introduced $\almostRiemann_{(\pm)}$ because they are convenient for
	studying perturbations of simple isentropic plane-waves (for which only $\RRiemann$ is non-vanishing);
	$\almostRiemann_{(\pm)}$ allow us to capture various kinds of smallness of the
	perturbations. 
	It is well-known that for isentropic plane-symmetric solutions, 
	one can use $\lbrace \RRiemann, \LRiemann \rbrace$ as
	the unknowns in place of $\lbrace \LogDensity, v^1 \rbrace$. Indeed, it is precisely this approach that we took in our prior work  \cite{abbrescia2022emergence}*{Appendix A} to construct isentropic plane-symmetric solutions with transversal convexity which serve as our background solution.
	\end{remark}


\subsection{The higher order fluid variables}
\label{SS:HIGHERORDERFLUIDVARIABLES}
The ``higher order'' fluid variables in the next definition appear in Theorem\,\ref{T:GEOMETRICWAVETRANSPORTSYSTEM},
which provides the formulation of compressible Euler flow that we use throughout our analysis.

\begin{definition}[The higher order fluid variables] \label{D:HIGHERORDERANDMODFLUIDVARSDEF} \hfill
\label{D:HIGHERORDERFLUIDVARIABLES}
\begin{enumerate}
\item We define the \emph{specific vorticity} to be the $\Sigma_t$-tangent vectorfield whose Cartesian spatial components are:
\begin{align} \label{E:SPECIFICVORTICITYDEF}
	\vortrenormalized^i 
	& \eqdef 
	\frac{(\Flatcurl v)^i}{\exp(\LogDensity)} 
	= 
	\frac{\upepsilon_{ijk} \updelta^{jl}\p_lv^k}{\exp(\LogDensity)},
\end{align}
where $\updelta^{jl}$ is the Kronecker delta.
\item We define the \emph{entropy gradient} to be the $\Sigma_t$-tangent vectorfield whose Cartesian spatial components are: 
	\begin{align} \label{E:GRADENTDEF}
		\GradEnt^i 
		& \eqdef \updelta^{ia}\p_a \Ent = \p_i \Ent.
	\end{align}
\item We define the \emph{modified fluid variables} to be the $\Sigma_t$-tangent vectorfield $\VortVort$ and the scalar function 
$\DivGradEnt$ whose Cartesian spatial components are:
\begin{subequations}
	\begin{align} \label{E:MODIFIEDCURLOFVORTICITY}
		\VortVort^i
		& \eqdef
			\exp(-\LogDensity) (\Flatcurl \vortrenormalized)^i
			+
			\exp(-3\LogDensity) \Speed^{-2} \frac{p_{;\Ent}}{\overline{\varrho}} \GradEnt^a \partial_a v^i
			-
			\exp(-3\LogDensity) \Speed^{-2} \frac{p_{;\Ent}}{\overline{\varrho}} (\Flatdiv v) \GradEnt^i,
				\\
		\DivGradEnt 
		& \eqdef 
			\exp(-2 \LogDensity) \Flatdiv \GradEnt 
			-
			\exp(-2 \LogDensity) \GradEnt^a \partial_a \LogDensity.
			\label{E:MODIFIEDDIVERGENCEOFENTROPYGRADIENT}
\end{align}
\end{subequations}
\end{enumerate}
\end{definition}

\subsection{Arrays of fluid variables and array norm notation}
\label{SS:FLUIDVARIABLEARRAYSANDNORMNOTATION}
We provide the next definition for notational convenience.

\begin{definition}[The fluid variable array $\wavearray$ and the partial array $\wavearraypartial$]
\label{D:ARRAYSOFWAVEVARIABLES} 
We define the \emph{array of wave\footnote{These ``wave-variables'' solve wave equations; 
see Theorem\,\ref{T:GEOMETRICWAVETRANSPORTSYSTEM}.} variables} as follows:
\begin{subequations}
\begin{align} \label{E:ARRAYOFWAVEVARIABLES} 
\wavearray 
& \eqdef (\Psi_0,\Psi_1,\Psi_2,\Psi_3,\Psi_4) 
\eqdef (\LogDensity,v^1,v^2,v^3,\Ent).
\end{align}We define the \emph{partial array of wave-variables} by:
\begin{align} \label{E:PARTIALWAVEARRAY}
	\wavearraypartial 
	& \eqdef (\Psi_1,\Psi_2,\Psi_3,\Psi_4)
	 = (\LRiemann,v^2,v^3,\Ent).
	\end{align}  
	\end{subequations}
We define the \emph{array of velocity variables} by:
	\begin{align} \label{E:ARRAYOFVELOCITYVARIABLES}
		\velocityarray \eqdef (v^1,v^2,v^3).
	\end{align}
	
We define the \emph{array of small velocity variables} by:
	\begin{align} \label{E:PARTIALVELOCITYARRAY}
		\velocityarraypartial \eqdef (v^2,v^3)
	\end{align}

\end{definition} 
We view $\wavearray$ to be an array of scalar functions $\Psi_{\iota}$, where 
$\iota = 0,\cdots,4$. We will not attribute any tensorial structure to the labeling index $\iota$ besides simple contractions, denoted by $\diamond$, corresponding to the chain rule; see Def.~\ref{D:DERIVATIVESOFARRAYS}. In view of Remark\,\ref{R:HOWWEESTIMATEDENSITYANDV1}, we often abuse notation and use $\wavearray = \{\RRiemann,\LRiemann,v^2,v^3,\Ent\}$. In most aspects of the analysis, this abuse of notation makes no difference. When the choice of solution variables is importantly, we explicitly state whether $\LogDensity,v^1$ are used over $\RRiemann,\LRiemann$ and vice-versa. 

In the next definition, we introduce notation for norms of arrays.
\begin{definition}[Conventions with variable arrays] \label{D:CONVENTIONESFORDIFFERENTIATION} \hfill
\begin{itemize}
\item Given the fluid variable array $\wavearray$ from Def.~\ref{D:ARRAYSOFWAVEVARIABLES}, we define:
\begin{align} \label{E:ABSOLUTVALUEOFPSIARRAY}
	|\wavearray| 
	& \eqdef \max_{\iota \in \{0,\cdots,4\}} | \Psi_\iota|.
\end{align}
For any norm $\|\cdot\|$ on scalar functions that appears in the paper, we set:
\begin{align} \label{E:NORMOFPSIARRAY}
	\| \wavearray\| 
	& \eqdef \max_{\iota \in \{0,\cdots,4\}} \| \Psi_\iota \|.
\end{align}
We use a similar convention for $\vortrenormalized$:
$|\vortrenormalized| = \max_{a = 1,2,3} |\vortrenormalized^a|$, 
and similarly for $\wavearraypartial$, $\GradEnt$, $\VortVort$,
etc. 
\item We use the following convention when taking norms of more than one variable at a time:
\begin{align} \label{E:NORMSOFMORETHANONEVARIABLE}
 \| (\vortrenormalized,\GradEnt)\| 
	& \eqdef \max\{  \| \vortrenormalized\|,  \| \GradEnt \|\}.
\end{align}
\end{itemize}
\end{definition}

\subsection{The acoustical metric and related geometric objects}
In the following definition, we introduce the acoustical metric and its inverse. 
This Lorentzian\footnote{By ``Lorentzian,'' we mean that viewed as a quadratic form, the symmetric $4 \times 4$ matrix
$(\gfour_{\alpha \beta})_{\alpha,\beta=0,1,2,3}$ has signature $(-,+,+,+)$.} 
metric drives the propagation of sound waves and is necessary to reveal the full geometry of the singular boundary.

\begin{definition}[The acoustical metric] \label{D:ACOUSTICALMETRICDEF}
Relative to the Cartesian coordinates $(t,x^1,x^2,x^3)$, 
we define the acoustical metric $\gfour$ and the inverse acoustical metric $\gfour^{-1}$ as 
follows, where the material derivative vectorfield $\Transport$ is defined in \eqref{E:MATERIALDERIVATIVEVECOTRFIELD}
and the speed of sound $\Speed$ is defined in \eqref{E:SOUNDSPEED}:
\begin{subequations}
\begin{align}
\gfour & = - \rmd t \otimes \rmd t + \Speed^{-2} \sum_{a=1}^3 (\rmd x^a - v^a \rmd t)\otimes  (\rmd x^a - v^a \rmd t), 
\label{E:ACOUSTICALMETRIC} \\
\gfour^{-1} & = - \Transport \otimes \Transport + \Speed^2 \sum_{a=1}^3 \p_a\otimes \p_a.
\label{E:INVERSEACOUSTICALMETRIC}
\end{align}
\end{subequations}
\end{definition}

Straightforward calculations yield that $\gfour_{\alpha \gamma} (\gfour^{-1})^{\gamma \beta} = \updelta_{\alpha}^{\beta}$,
where $\updelta_{\alpha}^{\beta}$ is the Kronecker delta, i.e., $\gfour^{-1}$ is indeed the inverse of $\gfour$.
In the remainder of the article, we silently lower and raise lowercase Greek indices with $\gfour$ and $\gfour^{-1}$, 
e.g., $V^{\alpha} = (\gfour^{-1})^{\alpha \beta} V_{\beta}$.


In our forthcoming analysis, the undifferentiated quantities
$v^i$ and $\Speed - 1$ will be small, where 
we quantify their smallness via the parameters $\mathring{\upalpha}$ and $\initialsmall$,
which we introduce in Sect.~\ref{S:PARAMETERSANDSIZEASSUMPTIONSANDCONVENTIONSFORCONSTANTS}.
Hence, in view of \eqref{E:ACOUSTICALMETRIC}, we find it convenient to introduce the following decomposition:
\begin{align} \label{E:SPLITMETRICINTOMINKOWSKIANDREMAINDERPART}
	\gfour_{\alpha \beta}(\wavearray) 
	& = m_{\alpha \beta} + \gfour_{\alpha \beta}^{(\textnormal{Small})}(\wavearray),
\end{align}
where $m_{\alpha \beta} = \diag(-1,1,1,1)$ is the Minkowski metric and $\gfour_{\alpha \beta}^{(\textnormal{Small})}(\wavearray)$ is a smooth function of $\wavearray$ satisfying: 
\begin{align} \label{E:METRICPERTURBATIONVANISHESATTRIVIALPSISOLUTION}
	\gfour_{\alpha \beta}^{(\textnormal{Small})}(\wavearray = 0) 
	& = 0.
\end{align}

The scalar functions $G^{\iota}_{\alpha \beta}$ in the following definition
will appear as coefficients in many of the equations that we study.

\begin{definition}[$\wavearray$-derivatives of $\gfour$] 
\label{D:DERIVATIVESOFMETRICWRTFLUID}
Viewing the Cartesian component functions
$\gfour_{\alpha \beta} = \gfour_{\alpha \beta}(\wavearray)$ 
as functions of the wave-variables, for $\alpha,\beta = 0,1,2,3$ and $\iota = 0,1,2,3,4$, 
we define: 
\begin{subequations}
\begin{align} \label{E:DERIVATIVESOFMETRICWRTFLUID}
G^{\iota}_{\alpha \beta}(\wavearray) 
& \eqdef \frac{\p}{\p \Psi_{\iota}} \gfour_{\alpha \beta}(\wavearray), 
	\\
\vec{G}_{\alpha \beta} 
& = \vec{G}_{\alpha \beta}(\wavearray)  
\eqdef 
\left(G_{\alpha \beta}^0 (\wavearray), 
	G_{\alpha \beta}^1 (\wavearray), 
	G_{\alpha \beta}^2 (\wavearray), 
	G_{\alpha \beta}^3 (\wavearray), 
	G_{\alpha \beta}^4 (\wavearray)
\right).
 \label{E:ARRAYVERSIONDERIVATIVESOFMETRICWRTFLUID}
\end{align}
\end{subequations}
\end{definition}

For each fixed $\iota \in \{0,\cdots,4\}$, 
we view $\{G^{\iota}_{\alpha \beta}\}_{\alpha,\beta = 0,\cdots,3}$ to be the Cartesian components of the spacetime tensorfield
``$G^{\iota}$.''
Similarly, we view 
$\{ \vec{G}_{\alpha \beta}\}_{\alpha,\beta = 0,\cdots,3}$ to be the Cartesian components of 
the array-valued spacetime tensorfield $\vec{G}$. 

\begin{remark}[No ambiguity in the vector $\vec{G}_{\alpha\beta}$]
	As we mentioned just below \eqref{E:PARTIALVELOCITYARRAY}, we often abuse notation when writing $\wavearray =\{\LogDensity,v^1,v^2,v^3,\Ent\}$ or $\wavearray = \{\RRiemann,\LRiemann,v^2,v^3,\Ent\}$. We sometimes require precise structure in the components of $\vec{G}^\iota_{\alpha\beta}$, so to prevent confusion in this paper, we always declare \eqref{E:DERIVATIVESOFMETRICWRTFLUID} \emph{is always defined with respect to \eqref{E:ARRAYOFWAVEVARIABLES}}. For example, ${G}^0_{\alpha\beta} = \frac{\p}{\p \LogDensity} \gfour_{\alpha\beta}$ and ${G}^1_{\alpha\beta} = \frac{\p}{\p v^1} \gfour_{\alpha\beta}$. Of course, expressions for $\frac{\p}{\p \RRiemann} \gfour_{\alpha\beta}$ could easily be obtained from $G^0_{\alpha\beta}$ and $G^1_{\alpha\beta}$ by the chain rule. In practice, the precise choice will be clear from context.
\end{remark}

\begin{definition}[Operators involving $\wavearray$] \label{D:DERIVATIVESOFARRAYS} Let $V_1,V_2$ be vectorfields, 
and let $D$ be a differential operator. We define:
\begin{align} \label{E:DERIVATIVEOFARRAY}
D \wavearray 
& \eqdef (D \Psi_0, D \Psi_1, D \Psi_2, D \Psi_3, D \Psi_4),
&
\vec{G}_{V_1 V_2} \diamond D \wavearray 
& \eqdef \sum_{\iota = 0}^4 G_{\alpha \beta}^{\iota} V_1^{\alpha}V_2^\beta D \Psi_\iota.
\end{align}
\end{definition}

\subsection{Covariant wave operator and $\gfour$-null forms}
\label{SS:COVARIANTWAVEOPERATORANDGNULLFORMS}
In this section, we provide some definitions that we need
to state Theorem\,\ref{T:GEOMETRICWAVETRANSPORTSYSTEM}, which provides
the geometric formulation of compressible Euler flow that we use throughout our analysis.

We start by recalling the standard definition of the covariant wave operator $\square_{\gfour}$.

\begin{definition}[Covariant wave operator of the acoustical metric]
\label{D:COVWAVEOP}
The covariant wave operator $\square_{\gfour}$ 
of the acoustical metric $\gfour = \gfour(\wavearray)$ 
acts on scalar-valued functions $\varphi$ as follows:\footnote{The formula \eqref{E:WAVEOPERATORARBITRARYCOORDINATES}
holds relative to arbitrary coordinates.}
\begin{align} \label{E:WAVEOPERATORARBITRARYCOORDINATES}
\square_{\gfour} \varphi
	& \eqdef 
	\frac{1}{\sqrt{|\mydet \gfour|}}
	\partial_{\alpha}
	\left\lbrace
			\sqrt{|\mydet \gfour|} (\gfour^{-1})^{\alpha \beta}
			\partial_{\beta} \varphi
	\right\rbrace.
\end{align}
\end{definition}

We now recall the definition of a standard null form with respect to the acoustical metric
(``$\gfour$-null form'' for short).

\begin{definition}[Standard $\gfour$-Null forms]
\label{D:STANDARDNULLFORMS}
Let $\varphi$ and $\widetilde{\varphi}$ be scalar functions.
We define
 $\nullform_{(\gfour)}(\partial\varphi,\partial\widetilde{\varphi})$ to be the following 
derivative-quadratic term:
\begin{subequations}
\begin{align}\label{E:Q0DEF}
\nullform_{(\gfour)}(\partial \varphi,\partial \widetilde{\varphi})
& \eqdef 
(\gfour^{-1})^{\alpha \beta} 
\partial_{\alpha} \varphi \partial_{\beta} \widetilde{\varphi}.
\end{align}For $0 \leq \alpha < \beta \leq 3$, 
we define $\nullform_{(\alpha \beta)}(\partial \varphi,\partial\widetilde{\varphi})$
to be the following derivative-quadratic term:
\begin{align}\label{E:QIJDEF}
\nullform_{(\alpha \beta)}(\partial\varphi,\partial\widetilde{\varphi})
& \eqdef 	\partial_{\alpha} \varphi \partial_{\beta} \widetilde{\varphi} 
			- 
			\partial_{\beta} \varphi \partial_{\alpha} \widetilde{\varphi}.
\end{align}
\end{subequations}
\end{definition}

In the rest of the paper, we use the terminology
\emph{null form relative to $\gfour$} 
or
\emph{$\gfour$-null form} 
to denote any linear combination
of the standard null forms
\eqref{E:Q0DEF}--\eqref{E:QIJDEF}
with (possibly solution dependent) coefficients that are controllable under the scope of our approach.


\subsection{The geometric wave-transport-divergence-curl formulation of the compressible Euler equations} 
\label{SS:GEOMETRICFORMULATIONOFFLOW}
Our main results fundamentally rely on the following formulation of the compressible Euler equations, first derived in \cite{jS2019c}.

\begin{theorem} [The geometric wave-transport-divergence-curl formulation of the compressible Euler equations]
	\label{T:GEOMETRICWAVETRANSPORTSYSTEM}
	Let $\overline{\varrho} > 0$ be any constant background density,\footnote{Recall
	that $\LogDensity$ depends on $\overline{\varrho}$; see Def.~\ref{D:LOGDENS}.}
	and assume that $(\LogDensity,v^1,v^2,v^3,\Ent)$ is
	a solution
	to the compressible Euler equations 
	\eqref{E:BVIEVOLUTION}--\eqref{E:BENTROPYEVOLUTION}
	in three spatial dimensions under an arbitrary equation of state $p = p(\varrho,\Ent)$ with positive sound speed 
	$\Speed$ (see \eqref{E:SOUNDSPEED}).
	Let $\Transport$ be the material derivative vectorfield defined in \eqref{E:MATERIALDERIVATIVEVECOTRFIELD},
	let $\gfour$ be the acoustical metric from Def.~\ref{D:ACOUSTICALMETRICDEF},
	let $\square_{\gfour}$ be the corresponding covariant wave operator from Def.~\ref{D:COVWAVEOP},
	let $\almostRiemannfunction  = \almostRiemannfunction(\LogDensity,\Ent)$ be the function from
	Def.~\ref{D:ALMOSTRIEMANNINVARIANTS},
	and let 
	$\vortrenormalized$,
	$\GradEnt$,
	$\VortVort$,
	and
	$\DivGradEnt$
	be the higher order variables from Def.~\ref{D:HIGHERORDERFLUIDVARIABLES}.
	Then the scalar-valued functions
	$\LogDensity$,
	$v^i$, 
	$\almostRiemann_{(\pm)}$, 
	$\Ent$,
$\vortrenormalized^i$,
$\GradEnt^i$,
$\Flatdiv \vortrenormalized$,
$\VortVort^i$,
$\DivGradEnt$,
and
$(\Flatcurl \GradEnt)^i$,
	($i=1,2,3$),
	also solve the following equations,
	where $\upepsilon_{ijk}$ is the fully antisymmetric symbol normalized by $\upepsilon_{123}=1$, 
	and \textbf{the Cartesian component functions $v^i$ are 
	treated as scalar-valued functions
	under covariant differentiation on LHS~\eqref{E:VELOCITYWAVEEQUATION}}:

\medskip

\noindent \underline{\textbf{\upshape Covariant wave equations}}.
	\begin{subequations} \label{E:COVARIANTWAVEEQUATIONSWAVEVARIABLES}
	\begin{align}
		\square_{\gfour(\wavearray)} v^i
		& = 
			- 
			\Speed^2 \exp(2 \LogDensity) \VortVort^i
			+ 
			\nullform_{(v)}^i
			+ 
			\mathfrak{L}_{(v)}^i,
			\label{E:VELOCITYWAVEEQUATION}	\\
	\square_{\gfour(\wavearray)} \LogDensity
	& = 
		-
		\exp(\LogDensity) \frac{p_{;\Ent}}{\overline{\varrho}} \DivGradEnt
		+
		\nullform_{(\LogDensity)}
		+
		\mathfrak{L}_{(\LogDensity)},
			\label{E:RENORMALIZEDDENSITYWAVEEQUATION} 
				\\
	\square_{\gfour(\wavearray)} \Ent
	& = 
		\Speed^2 \exp(2 \LogDensity)  \DivGradEnt
		+
		\mathfrak{L}_{(\Ent)}.
	\label{E:ENTROPYWAVEEQUATION}
\end{align}
\end{subequations}\medskip

\noindent \underline{\textbf{\upshape Transport equations}}.
\begin{subequations}
\begin{align}	\Transport \vortrenormalized^i
	& = \mathfrak{L}_{(\vortrenormalized)}^i,
		\label{E:RENORMALIZEDVORTICTITYTRANSPORTEQUATION}
		\\
	\Transport \Ent
	& = 0,	
		\label{E:ENTROPYTRANSPORTMAINSYSTEM}
			\\
	\Transport \GradEnt^i
	& = \mathfrak{L}_{(\GradEnt)}^i.
		\label{E:GRADENTROPYTRANSPORT}
\end{align}
\end{subequations}\medskip	
	
\noindent \underline{\textbf{\upshape Transport-divergence-curl system for the specific vorticity}}.
\begin{subequations}
\begin{align} \label{E:FLATDIVOFRENORMALIZEDVORTICITY}
	\Flatdiv \vortrenormalized
	& = 
		\mathfrak{L}_{(\Flatdiv \vortrenormalized)},
		\\
\Transport \VortVort^i
& = 
	\mainnullform_{(\VortVort)}^i 
	+ 
	\nullform_{(\VortVort)}^i
	+	
	\mathfrak{L}_{(\VortVort)}^i. 
	\label{E:EVOLUTIONEQUATIONFLATCURLRENORMALIZEDVORTICITY} 
\end{align}	
\end{subequations}\medskip

\noindent \underline{\textbf{\upshape Transport-divergence-curl system for the entropy gradient}}.
\begin{subequations}
\begin{align} 	
\Transport \DivGradEnt
	& =  \mainnullform_{(\DivGradEnt)} +
			\nullform_{(\DivGradEnt)},
 \label{E:TRANSPORTFLATDIVGRADENT}
			\\
	(\Flatcurl \GradEnt)^i & = 0.
	\label{E:CURLGRADENTVANISHES}
\end{align}
\end{subequations}
	Above, the main terms $\mainnullform_{(\VortVort)}^i$ and $\mainnullform_{(\DivGradEnt)}$ 
	in the transport equations for the modified fluid variables are 
	the \textbf{null forms relative to} $\gfour$ (see Def.~\ref{D:STANDARDNULLFORMS})
	defined by:\footnote{Actually, the last the last term on
 RHS\,\eqref{E:TRANSPORTDIVGRADENTMAINTERMS} is not a null form, but rather 
	a simpler harmless error term.}
	\begin{subequations}
	\begin{align}
	\begin{split} \label{E:TRANSPORTVORTVORTMAINTERMS}
		\mainnullform_{(\VortVort)}^i & \eqdef 
		- 
		2 \updelta_{jk} \upepsilon_{iab} \exp(-\LogDensity) (\partial_a v^j) \partial_b \vortrenormalized^k
		+
		\upepsilon_{ajk}
		\exp(-\LogDensity)
		(\partial_a v^i) 
		\partial_j \vortrenormalized^k
			\\
& \ \
		+ 
		\exp(-3 \LogDensity) \Speed^{-2} \frac{p_{;\Ent}}{\overline{\varrho}} 
		\left\lbrace
			(\Transport \GradEnt^a) \partial_a v^i
			-
			(\Transport v^i) \partial_a \GradEnt^a
		\right\rbrace
			\\
	& \ \
		+
		\exp(-3 \LogDensity) \Speed^{-2} \frac{p_{;\Ent}}{\overline{\varrho}}  
		\left\lbrace
			(\Transport v^a) \partial_a \GradEnt^i
			- 
			( \Transport \GradEnt^i) \partial_a v^a
		\right\rbrace,
	\end{split}
			\\
	\mainnullform_{(\DivGradEnt)} & = 	2 \exp(-2 \LogDensity) 
			\left\lbrace
				(\partial_a v^a) \partial_b \GradEnt^b
				-
				(\partial_a \GradEnt^b) \partial_b v^a
			\right\rbrace
			+
			\exp(-\LogDensity) \updelta_{ab} (\Flatcurl \vortrenormalized)^a \GradEnt^b. 
			\label{E:TRANSPORTDIVGRADENTMAINTERMS}  
	\end{align}
	\end{subequations} 
	Moreover, 
	$\nullform_{(v)}^i$,
	$\nullform_{(\pm)}$,
	$\nullform_{(\LogDensity)}$, 
	$\nullform_{(\VortVort)}^i$,
	and
	$\nullform_{(\DivGradEnt)}$
	are\footnote{The term 
	$\mainnullform_{(\VortVort)}^i$
	on RHS\,\eqref{E:EVOLUTIONEQUATIONFLATCURLRENORMALIZEDVORTICITY} 
	and the term $\mainnullform_{(\DivGradEnt)}$
	on RHS\,\eqref{E:TRANSPORTFLATDIVGRADENT} are also null forms relative to $\gfour$.
	We have isolated these two null forms with different notation because
	they 
	are more difficult to treat than 
	$\nullform_{(v)}^i$, $\nullform_{(\pm)}$
	$\nullform_{(\LogDensity)}$, 
	$\nullform_{(\VortVort)}^i$,
	and
	$\nullform_{(\DivGradEnt)}$;
	to bound the top-order derivatives of the ``$\mainnullform$'' terms,
	we rely on the delicate ``elliptic-hyperbolic'' identities
	that we derive in Sect.~\ref{S:ELLIPTICHYPERBOLICIDENTITIES}.
	\label{FN:MORENULLFORMS}} 
	the null forms relative to $\gfour$ defined by: 
	\begin{subequations}
		\begin{align}
		\nullform_{(v)}^i	
		& \eqdef 	-
					\left\lbrace
						1
						+
						\Speed^{-1} \Speed_{;\LogDensity}
					\right\rbrace
					(\gfour^{-1})^{\alpha \beta} (\partial_{\alpha} \LogDensity) \partial_{\beta} v^i,
			\label{E:VELOCITYNULLFORM} 
			\\
		\nullform_{(\pm)} 
			& 
			\eqdef \nullform_{(v)}^1 \mp 2 \Speed_{;\LogDensity} 
			(\gfour^{-1})^{\alpha \beta} \p_\alpha \LogDensity \p_\beta \LogDensity
			\pm 
			\Speed \left\lbrace (\p_a v^a) (\p_b v^b) - (\p_a v^b) \p_b v^a 
			\right\rbrace, 
				\label{E:RIEMANNINVARIANTWAVEEQUATIONNULLFORM}
		\\
		\nullform_{(\LogDensity)}
		& \eqdef 
		- 
		3 \Speed^{-1} \Speed_{;\LogDensity} 
		(\gfour^{-1})^{\alpha \beta} (\partial_{\alpha} \LogDensity) \partial_{\beta} \LogDensity
		+ 
		\left\lbrace
			(\partial_a v^a) \partial_b v^b
			-
			(\partial_a v^b) \partial_b v^a
		\right\rbrace,
			\label{E:DENSITYNULLFORM}
				\\
	\begin{split} \label{E:RENORMALIZEDVORTICITYCURLNULLFORM}
	\nullform_{(\VortVort)}^i
	& \eqdef
		\exp(-3 \LogDensity) \Speed^{-2} \frac{p_{;\Ent}}{\overline{\varrho}}  \GradEnt^i
		\left\lbrace
			(\partial_a v^b) \partial_b v^a
			-
			(\partial_a v^a) \partial_b v^b
		\right\rbrace
			 \\
& \ \
		+ 
		\exp(-3 \LogDensity) \Speed^{-2} \frac{p_{;\Ent}}{\overline{\varrho}}
		\GradEnt^b 
		\left\lbrace
			(\partial_a v^a) \partial_b v^i 
			- 
			(\partial_a v^i) \partial_b v^a 
		\right\rbrace
		\\
& \ \
		+ 
		2 \exp(-3 \LogDensity) \Speed^{-2} \frac{p_{;\Ent}}{\overline{\varrho}}
		\GradEnt^a
		\left\lbrace
			(\partial_a \LogDensity)  \Transport v^i
			 - 
		  (\partial_a v^i) \Transport \LogDensity
		\right\rbrace
			\\
	&  \ \
			+ 
			2 \exp(-3 \LogDensity) \Speed^{-3} \Speed_{;\LogDensity} \frac{p_{;\Ent}}{\overline{\varrho}}
			\GradEnt^a 
			\left\lbrace
				(\partial_a \LogDensity)  \Transport v^i
				- 
				(\partial_a v^i) \Transport \LogDensity
			\right\rbrace
				\\
	& \ \
		+ 
		\exp(-3 \LogDensity) \Speed^{-2} \frac{p_{;\Ent;\LogDensity}}{\overline{\varrho}}
		\GradEnt^a 
		\left\lbrace
			(\partial_a v^i)
			\Transport \LogDensity 
			- 
			(\partial_a \LogDensity) \Transport v^i
		\right\rbrace
			\\
	& \ \
		+
		\exp(-3 \LogDensity) \Speed^{-2} \frac{p_{;\Ent;\LogDensity}}{\overline{\varrho}} \GradEnt^i
		\left\lbrace
			(\Transport v^a) \partial_a \LogDensity
			-
			 (\Transport \LogDensity) \partial_a v^a
		\right\rbrace
				\\
	& \ \
		+
		2 \exp(-3 \LogDensity) \Speed^{-2} \frac{p_{;\Ent}}{\overline{\varrho}} \GradEnt^i
		\left\lbrace
			(\Transport \LogDensity) \partial_a v^a 
			- 
			(\Transport v^a) \partial_a \LogDensity
		\right\rbrace
				\\
	& \ \
		 	+
			2 \exp(-3 \LogDensity) \Speed^{-3} \Speed_{;\LogDensity} \frac{p_{;\Ent}}{\overline{\varrho}} \GradEnt^i
			\left\lbrace
				(\Transport \LogDensity) \partial_a v^a
		 		- 
		 		(\Transport v^a) \partial_a \LogDensity
		 	\right\rbrace,
\end{split} 
		 		\\
\label{E:DIVENTROPYGRADIENTNULLFORM}
\nullform_{(\DivGradEnt)} 
	& \eqdef
		2 \exp(-2 \LogDensity) 
		\GradEnt^a 
		\left\lbrace
			(\partial_a v^b) \partial_b \LogDensity
			-
			(\partial_a \LogDensity)
			\partial_b v^b 
		\right\rbrace.
\end{align}
\end{subequations}
	In addition, the terms
	$\mathfrak{L}_{(v)}^i$, $\mathfrak{L}_{(\pm)}$, 
	$\mathfrak{L}_{(\LogDensity)}$,
	$\mathfrak{L}_{(\Ent)}$,
	$\mathfrak{L}_{(\vortrenormalized)}^i$,
	$\mathfrak{L}_{(\GradEnt)}^i$,
		$\mathfrak{L}_{(\Flatdiv \vortrenormalized)}$,
		and
	$\mathfrak{L}_{(\VortVort)}^i$,
	which are at most linear in the derivatives of the unknowns, are defined as follows:
	\begin{subequations}
	\begin{align} 
	\begin{split} \label{E:VELOCITYILINEARORBETTER} 
		\mathfrak{L}_{(v)}^i
		& \eqdef 
		2 \exp(\LogDensity) \upepsilon_{iab} (\Transport v^a) \vortrenormalized^b
		-
		\frac{p_{;\Ent}}{\overline{\varrho}} \upepsilon_{iab} \vortrenormalized^a \GradEnt^b
			\\
	& \ \
		- 
		\frac{1}{2} \exp(-\LogDensity) \frac{p_{;\LogDensity;\Ent}}{\overline{\varrho}} \GradEnt^a \partial_a v^i
				\\
		& \ \
		- 
		2 \exp(-\LogDensity) \Speed^{-1} \Speed_{;\LogDensity} \frac{p_{;\Ent}}{\overline{\varrho}} 
		(\Transport \LogDensity) \GradEnt^i
		+
		\exp(-\LogDensity) \frac{p_{;\Ent;\LogDensity}}{\overline{\varrho}} (\Transport \LogDensity) \GradEnt^i,
		\end{split} 
			\\
		\mathfrak{L}_{(\LogDensity)}  \label{E:DENSITYLINEARORBETTER}
		& \eqdef
		-
		\frac{5}{2} \exp(-\LogDensity) \frac{p_{;\Ent;\LogDensity}}{\overline{\varrho}} \GradEnt^a \partial_a \LogDensity 
		-
		\exp(-\LogDensity) \frac{p_{;\Ent;\Ent}}{\overline{\varrho}} \updelta_{ab} \GradEnt^a \GradEnt^b,
		  \\
		\mathfrak{L}_{(\Ent)}
		& \eqdef 
			\Speed^2  \GradEnt^a \partial_a \LogDensity
			- 
			\Speed \Speed_{;\LogDensity} \GradEnt^a \partial_a \LogDensity
			- 
			\Speed \Speed_{;\Ent} \updelta_{ab} \GradEnt^a \GradEnt^b,
			\label{E:ENTROPYLINEARORBETTER} 
				\\
		\mathfrak{L}_{(\vortrenormalized)}^i
		& \eqdef 
		\vortrenormalized^a \partial_a v^i
		-
		\exp(-2 \LogDensity) \Speed^{-2} \frac{p_{;\Ent}}{\overline{\varrho}} \upepsilon_{iab} (\Transport v^a) \GradEnt^b,
		\label{E:SPECIFICVORTICITYLINEARORBETTER}
			\\
		\mathfrak{L}_{(\GradEnt)}^i
		& \eqdef
			- 
			\GradEnt^a \partial_a v^i
			+ 
			\upepsilon_{iab} \exp(\LogDensity) \vortrenormalized^a \GradEnt^b,
			\label{E:ENTROPYGRADIENTLINEARORBETTER}
				\\
	\mathfrak{L}_{(\Flatdiv \vortrenormalized)}
		& \eqdef - \vortrenormalized^a \partial_a \LogDensity,
		\label{E:RENORMALIZEDVORTICITYDIVLINEARORBETTER} 
		\\
	\begin{split} \label{E:RENORMALIZEDVORTICITYCURLLINEARORBETTER} 
	\mathfrak{L}_{(\VortVort)}^i	 		
	& \eqdef
		 		2 \exp(-3 \LogDensity) \Speed^{-3} \Speed_{;\Ent} \frac{p_{;\Ent}}{\overline{\varrho}} 
				(\Transport v^i) \updelta_{ab} \GradEnt^a \GradEnt^b
				\\
		& \ \
			-
				2 \exp(-3 \LogDensity) \Speed^{-3} \Speed_{;\Ent} \frac{p_{;\Ent}}{\overline{\varrho}} 
				\updelta_{ab} \GradEnt^a (\Transport v^b) \GradEnt^i
					\\
		& \ \
			+ 
			\exp(-3 \LogDensity) \Speed^{-2} \frac{p_{;\Ent;\Ent}}{\overline{\varrho}} \updelta_{ab} (\Transport v^a) \GradEnt^b \GradEnt^i
				\\
		& \ \ 
			- 
			\exp(-3 \LogDensity) \Speed^{-2} \frac{p_{;\Ent;\Ent}}{\overline{\varrho}} (\Transport v^i) \updelta_{ab} \GradEnt^a \GradEnt^b.
		\end{split}
		\end{align}
	\end{subequations}
\end{theorem}

The following corollary is immediate based on the definition of the modified curl of vorticity $\VortVort$ \eqref{E:MODIFIEDCURLOFVORTICITY} and the covariant wave equation satisfied by the velocity \eqref{E:VELOCITYWAVEEQUATION}. 

\begin{corollary}[The covariant wave equations for the velocity in terms of $\Flatcurl(\Flatcurl(v))$]
The fluid velocity components $v^i$ also satisfy the following covariant wave equation: 
	\begin{align}
		\begin{split}
			\square_{\gfour(\wavearray)} v^i
			& = 
			- 
			\Speed^2 \exp(\LogDensity)(\Flatcurl \, \vortrenormalized)^i
			 - \exp(-\LogDensity) \frac{p_{;\Ent}}{\overline{\varrho}} \GradEnt^a \p_a v^i + \exp(-\LogDensity) \frac{p_{;\Ent}}{\overline{\varrho}} (\Flatdiv v) \GradEnt^i +
			\nullform_{(v)}^i
			+ 
			\mathfrak{L}_{(v)}^i,
			\label{E:VELOCITYWAVEEQUATIONWITHCURLCURL} 
		\end{split}
	\end{align}
where $\nullform_{(v)}^i$ and $\mathfrak{L}_{(v)}^i$ are as in \eqref{E:VELOCITYNULLFORM} and \eqref{E:VELOCITYILINEARORBETTER}, respectively.

\end{corollary}


\section{The outgoing acoustic geometry and the arrays \texorpdfstring{$\controlvars$ and $\badcontrolvars$}{}}
\label{S:ACOUSTICGEOMETRYANDCOMMUTATORVECTORFIELDS}
In this section, we carefully construct what we call the \emph{outgoing} acoustic geometry, reveal its basic properties,
and provide the evolution equations satisfied by various geometric tensors. The fundamental object behind all the constructions is an acoustic eikonal function 
$u$ initially set to $u|_{\Sigma_0} = -x^1$. As in prior work on shock formation in perturbations of plane symmetry \cite{jS2016b,gHsKjSwW2016, luk2018shock, LukSpeck2024stability}, from $u$ we construct \emph{commutation vectorfields} that have just enough regularity to derive the a priori $L^2$ estimates. For this reason, we present in this section the construction of $u$ and the same relevant tensors as in \cite{jS2016b,gHsKjSwW2016, luk2018shock, LukSpeck2024stability}, which we refer to as \emph{the outgoing} acoustic geometry. We also introduce the solution variable arrays $\controlvars$ and $\badcontrolvars$,
which contain the wave-variables and various components of the outgoing acoustic geometry. We use these arrays throughout the paper to simplify the notation and to allow for convenient,
schematic expressions. 

\subsection{The outgoing eikonal function and inverse foliation density}
\label{SS:EIKONAFUNCTIONANDINVERSEFOLIATIONDENSITY}
In the following definition, we introduce the \emph{outgoing eikonal function} $u$ and the outgoing inverse foliation density $\upmu$. The level-sets of $u$ are the characteristics for the wave operator $\square_{\gfour}$, 
while the inverse (i.e., reciprocal) 
of the inverse foliation density measures the density of these characteristics.

\begin{definition}[Outgoing eikonal function and inverse foliation density]
\label{D:REGULAREIKONALFUNCTIONANDMU}
The \emph{outgoing eikonal function} $u$ is the solution of the following fully nonlinear hyperbolic initial value problem,
where $\gfour$ is the acoustical metric defined in \eqref{E:ACOUSTICALMETRIC}
and the PDE is known as the \emph{acoustic eikonal equation}:
\begin{align} \label{E:REGULAREIKONALEQUATION}
\begin{cases}
(\gfour^{-1})^{\alpha \beta} \partial_{\alpha} u \partial_{\beta} u = 0, 
	\\
\p_t u > 0, 
	\\
u|_{\Sigma_0} = - x^1.
\end{cases}
\end{align}
We define the \emph{outgoing inverse foliation density} $\upmu$ by:
\begin{align} \label{E:MUDEF}
	\upmu 
	& 
	\eqdef 
	- 
	\frac{1}{(\gfour^{-1})^{\alpha \beta} \partial_{\alpha} t \partial_{\beta} u} > 0.
\end{align}
\end{definition}

\begin{remark}[The term ``outgoing''] \label{R:WEDROPTERMOUTGOINGWHENCLEAR}
When it is clear from the context, we will drop the term ``outgoing'', e.g. when none of the ``ingoing'' acoustic geometry is present as in this section.
\end{remark}

\subsection{Acoustical subsets of spacetime}
\label{SS:ACOUSTICSUBSETSOFSPACETIME}

\begin{definition}[Acoustical subsets of spacetime] 
	\label{D:ACOUSTICSUBSETSOFSPACETIME}
We define the following ``acoustical subsets'' of spacetime:
\begin{subequations}
\begin{align}
\Sigma_{t'} 
& \eqdef 
\{(t,x^1,x^2,x^3) \in \R\times\R\times\T^2 \ | \ t = t'\}, 
	\label{E:SIGMAT}
	\\
\nullhyparg{u'} 
& \eqdef \{(t,x^1,x^2,x^3) \in \R\times\R\times\T^2 \ | \ u(t,x^1,x^2,x^3) = u'\}, 
	\label{E:NULLHYPERSURFACES}
	\\
\ell_{t',u'} 
& 
\eqdef 
\Sigma_{t'} \cap \nullhyparg{u'} 
= 
\{(t,x^1,x^2,x^3) \in \R\times\R\times\T^2 \ | \ t= t',\ u(t,x^1,x^2,x^3) = u'\}. 
\label{E:ELLTUSMOOTHTORI} 
\end{align}
\end{subequations}Given real numbers $u_1 \leq u_2$ and $t_1 \leq t_2$, we define the following ``truncated'' subsets of spacetime:
\begin{subequations}
\begin{align}
\Sigma_{t'}^{[-\rightu,\leftu]} 
& \eqdef \Sigma_{t'} \cap \{(t,x^1,x^2,x^3) \in \R\times\R\times\T^2 \ | \ u_1 \leq u(t,x^1,x^2,x^3) \leq u_2 \}, 
	\label{E:TRUNCATEDSIGMAT}
	\\
\nullhyparg{u'}^{[t_1,t_2]} 
& \eqdef \nullhyparg{u'} \cap \{(t,x^1,x^2,x^3) \in \R\times\R\times\T^2 \ | \ t_1 \leq t \leq t_2 \}.
	\label{E:TRUNCATEDNULLHYPERSURFACES}
\end{align}
\end{subequations}\end{definition}

We refer to the $\Sigma_t$ as ``constant Cartesian-time hypersurfaces,'' 
the $\nullhyparg{u}$ as ``null hypersurfaces,''
``acoustic characteristics,''  
or ``characteristics,'' and the $\ell_{t,u}$ as ``acoustic tori.'' 
We note that many of our estimates will take place on the rough hypersurfaces 
and the rough tori of Def.\,\ref{D:TRUNCATEDROUGHSUBSETS}.


\subsection{Geometric coordinates and important acoustic vectorfields}

\label{SS:REGULARGEOMETRICCOORDIANTESANDPARTIALDERIVATIVEVECTORFIELDS}
\subsubsection{Geometric coordinates}
\label{SSS:REGULARGEOMETRICCOORDINATES}

\begin{definition}[The  geometric coordinates and their corresponding partial derivative vectorfields] 
\label{D:GEOMETRICCOORDIANTESANDPARTIALDERIVATIVEVECTORFIELDS}
We refer to $(t,u,x^2,x^3)$ as the \emph{geometric coordinates}. We denote the
geometric coordinate partial derivative vectorfields by
$\left\lbrace\geop{t},\geop{u},\geop{x^2},\geop{x^3} \right\rbrace$.  
\end{definition}

	\begin{remark}[Coordinate systems on $\ell_{t,u}$ and $\nullhyparg{u}$]
	Note that $(x^2,x^3)$ form a coordinate system on the acoustic tori $\ell_{t,u}$ 
	and that $\left\lbrace \geop{x^2},\geop{x^3} \right\rbrace$ span the tangent space of $\ell_{t,u}$.
	Similarly, $(t,x^2,x^3)$ form a coordinate system on the null hypersurfaces $\nullhyparg{u}$ 
	and $\left\lbrace \geop{t}, \geop{x^2},\geop{x^3} \right\rbrace$ span the tangent space of 
	$\nullhyparg{u}$.
\end{remark}


\begin{notation}[Conventions used with $(x^2,x^3)$ and $\lbrace \geop{x^2}, \geop{x^3} \rbrace$]
	\label{N:NOTATIONINVOLVINGSMOOTHTORUSCOORDINATES}
	\hfill
\begin{enumerate}
\item If $V$ is a vectorfield, then for $A = 2,3$, $V^A \eqdef V x^A = V^{\alpha} \partial_{\alpha} x^A$.
	In particular, if $V$ is $\ell_{t,u}$-tangent, then $V = V^A \geop{x^A}$, and  
	$V^A$ are the components of $V$ with respect to geometric coordinates $(x^2,x^3)$ on $\ell_{t,u}$.
\item If $\upxi$ is a one-form, then we use the abbreviated notation 
	$\upxi_A \eqdef \upxi\left( \geop{x^A} \right) = \upxi_{\alpha} (\geop{x^A})^{\alpha}$ for $A = 2,3$. 
\item We use a similar convention for higher order tensorfields, 
	e.g., $\upxi_{AB} = \upxi(\geop{x^A},\geop{x^B})$ for $A,B = 2,3$. 
\item We sum repeated uppercase Latin indices over $A = 2,3$, e.g., 
$\upxi_{AA} \eqdef \upxi_{22} + \upxi_{33}$.
\end{enumerate}
\end{notation}

\subsubsection{The important acoustic vectorfields}
\label{SSS:IMPORTANTACOUSTICVECTORFIELDS}
The vectorfields in the next definition are fundamental for the rest of the paper.

\begin{definition}[The important acoustic vectorfields] \label{D:COMVECTORFIELDS} 
\hfill
\begin{enumerate}
\item We define the \emph{geodesic null vectorfield} by:
\begin{align} \label{E:LGEO}
\Lgeo^{\alpha}
& \eqdef - (\gfour^{-1})^{\alpha \beta}\p_\beta u,
\end{align}
and the \emph{rescaled null vectorfield} 
as follows, where $\upmu$ is the inverse foliation density defined in \eqref{E:MUDEF}:
\begin{align} \label{E:LUNIT}
\Lunit 
& \eqdef \upmu \Lgeo.
\end{align}
For $i = 1,2,3$, we define the scalar functions $\Lsmall^i$ as follows, 
where throughout the paper, $\Lunit^i$ denotes the Cartesian component $\Lunit x^i$:
\begin{align} \label{E:LSMALLDEF}
	\Lsmall^1 
	& \eqdef L^1 - 1, \qquad \Lsmall^2 \eqdef L^2,\qquad \Lsmall^3 \eqdef L^3.
\end{align}\item We define $X$ to be the unique vectorfield that is $\Sigma_t$-tangent, and $\gfour$-orthogonal to $\ell_{t,u}$, and normalized by: 
\begin{align} \label{E:X}
\gfour(\Lunit,X) 
& = -1,
\end{align}
and we define the rescaled vectorfield $\muX$ by: 
\begin{align} \label{E:MUX}
\muX 
& \eqdef \upmu X.
\end{align}
For $i = 1,2,3$, we define the scalar functions $\Xsmall^i$ as follows, 
where throughout the paper, $X^i$ denotes the Cartesian component $X x^i$:\footnote{For the solutions covered by our main results, the functions $\Lsmall^i$ and $\Xsmall^i$ will have magnitudes that are $\ll 1$.}  
\begin{align} \label{E:XSMALL}
\Xsmall^1 
& \eqdef X^1 + 1, \qquad \Xsmall^2 \eqdef X^2, \qquad \Xsmall^3 \eqdef X^3.
\end{align}\item We define $\olduLunit$ to be the ingoing null vectorfield given by:
	\begin{align} \label{E:OLDULUNIT}
		\olduLunit \eqdef \Lunit + 2X.
	\end{align}
\item We define $\Yvf{2}$, $\Yvf{3}$, to be the following $\ell_{t,u}$-tangent vectorfields:
\begin{align} \label{E:YCOMMUTATOR}
\Yvf{2} & \eqdef \p_2 - \gfour(\p_2,X) X,  & \Yvf{3} & \eqdef \p_3 - \gfour(\p_3,X) X.  
\end{align}
We also define $\Yvfsmall{2}$, $\Yvfsmall{3}$, to be the following vectorfields 
(which are not generally $\ell_{t,u}$-tangent):
\begin{align} \label{E:YSMALL}
\Yvfsmall{2} & \eqdef \Yvf{2} - \partial_2, & \Yvfsmall{3} & \eqdef \Yvf{3} - \partial_2. 
\end{align}
	We similarly define the Cartesian component functions 
	$\Yvf{2}^i$, $\Yvf{3}^i$, $\Yvfsmall{2}^i$, $\Yvfsmall{3}^i$ 
	analogously to \eqref{E:XSMALL}.
\item We define the \emph{commutation vectorfields} $\Fullset$, the 
$\nullhyparg{u}$-tangential subset $\Tanset$, and the $\ell_{t,u}$-tangential subset
$\Angularset$:
\begin{align} \label{E:COMMUTATIONVECTORFIELDS}
	\Fullset
	& \eqdef \lbrace \Lunit, \muX, \Yvf{2}, \Yvf{3} \rbrace,
	& 
	\Tanset
	& \eqdef \lbrace \Lunit, \Yvf{2}, \Yvf{3} \rbrace,
	&
	\Angularset
	& 
	\eqdef \lbrace \Yvf{2}, \Yvf{3} \rbrace.
\end{align}\end{enumerate}
\end{definition}

In order to obtain wave equation energy estimates that are sufficient to allow us to track the solution up to the Cauchy horizon, we use the following multiplier vectorfield:

\begin{definition}[The $\gfour$-timelike multiplier vectorfield $\multipliervectorfield $] 
	We define $\multipliervectorfield$ to be the following vectorfield:
	\begin{align} \label{E:MULTIPLIERVECTORFIELD}
		\multipliervectorfield 
		& \eqdef \Lunit + 2 \upmu \Transport =  (1 + 2 \upmu) \Lunit 
		+ 
		2  \muX.
	\end{align}
\end{definition}

Lemma\,\ref{L:COMMUTATORSTOCOORDINATES} shows that
$\Fullset$ spans the tangent spaces of spacetime equipped with
the differential structure corresponding to the geometric coordinates $(t,u,x^2,x^3)$.
We sometimes refer to $\Fullset$ as the \emph{rescaled frame}
because the vectorfield $\muX = \upmu X$ degenerates with respect to the Cartesian differential
structure as $\upmu \downarrow 0$, i.e., $\muX^i = \upmu X^i$ tends to $0$.
Similarly, the lemma shows
that $\Tanset$ spans the tangent spaces of the characteristics $\nullhyparg{u}$
and that $\Angularset$ spans the tangent spaces of the $\ell_{t,u}$.
To derive $L^{\infty}$ and H\"{o}lder estimates, we commute various PDEs with 
elements of $\Fullset$. To derive energy estimates,
we will commute various PDEs with the elements of $\Tanset$.
For a handful of key estimates, we will refer to the set $\Angularset$.
We also note that from definitions \eqref{E:PROJECTIONOFTENSORONTOFLATTORUS}
and \eqref{E:YCOMMUTATOR}, it is straightforward to check that
$\Yvf{A} = \smoothtorusproject \partial_A$, i.e., 
$\Yvf{A}$ is the $\gfour$-orthogonal projection of the Cartesian partial derivative vectorfield $\partial_A$
onto $\ell_{t,u}$; see also the first equality in \eqref{E:SMOOTHTORUSMETRICINTERMSOFSIGMATMETRICANDX}.


Throughout the paper, we will often silently use the identities
featured in the following lemma.

\begin{lemma}[Basic properties of the vectorfields] 
\label{L:BASICPROPERTIESOFVECTORFIELDS} 
The follow results hold.

\begin{enumerate}
\item The vectorfield $\Lgeo$ is geodesic and $\gfour$-null, i.e., with 
$\Dfour$ denoting the Levi-Civita connection of $\gfour$, we have:
 \begin{align} \label{E:LGEOISGEODESICANDNULL}
 \gfour(\Lgeo,\Lgeo) & = 0, 
	&  
\Dfour_{\Lgeo}\Lgeo & = 0.
\end{align}
The rescaled vectorfield $\Lunit$ is also $\gfour$-null:
\begin{align} \label{E:LUNITISNULL}
\gfour(\Lunit,\Lunit) 
& = 0,
\end{align}
and it satisfies the following identity,
where $\upmu$ is the inverse foliation density defined in \eqref{E:MUDEF}:
\begin{align} \label{E:COVARIANTLUNITDERIVATIVEOFLUNIT}
	\Dfour_{\Lunit} \Lunit
	& = \frac{(\Lunit \upmu)}{\upmu} \Lunit.
\end{align}\item $\Lunit$ is $\gfour$-orthogonal to the characteristics $\nullhyparg{u}$, 
that is, for any vectorfield $\Singletan$ tangent to $\nullhyparg{u}$,
	we have:
	\begin{align} \label{E:LUNITISNORMALTOCHARACTERISTICS}
		\gfour(\Lunit,\Singletan)
		& = 0.
	\end{align}\item The following identities hold:
\begin{align}
\Lunit u & = 0, & \Lunit t = \Lunit^0 & = 1, & \muX u & = 1, & \muX t = \muX^0 &= 0, \label{E:LULTMUXUMUXT} \\
\gfour(X,X) & = 1, & \gfour(\muX,\muX) & = \upmu^2, & \gfour(\Lunit,X) & = -1, & \gfour(\Lunit,\muX) &= -\upmu. \label{E:GOFLANDX}
\end{align}\item The vectorfield $\olduLunit$ is $\gfour$-null and transverse to $\nullhyparg{u}$. In particular, we have:
	\begin{align}
		\gfour(\olduLunit,\olduLunit) = 0, & & \gfour(\Lunit,\olduLunit) = -2. \label{E:LPLUS2XISNULLANDTRANSVERSETOPU}
	\end{align}

\item The material vectorfield $\Transport$ is future directed, $\gfour$-orthogonal to $\Sigma_t$ (and hence also to $\ell_{t,u}$), 
and it is of $\gfour$-unit size: 
\begin{align} \label{E:TRANSPORTISUNITLENGTH}
\gfour(\Transport,\Transport) 
& = -1.
\end{align}
Moreover, we have:
\begin{align} \label{E:BISLPLUSX}
\Transport & = \Lunit + X, 
\end{align}
and relative to the Cartesian coordinates, we have:
\begin{align}
\Transport_{\alpha}
& = 
-
\updelta_\alpha^0, \label{E:MATERIALDERIVATIVELOWEREDCARTESIANCOORDINATES}
\end{align}
where $\updelta_{\alpha}^{\beta}$ is the Kronecker delta. 
\item Finally, the following identities hold for $i=1,2,3$ and $A=2,3$:
\begin{subequations}
\begin{align}
	\Xsmall^i  & = - \Lsmall^i + v^i, 
	\label{E:XSMALLINTERMSOFLSMALLANDVELOCITY} 
		\\
	\Yvfsmall{A}^i
	& = 
			-\Speed^{-2} \Xsmall^A
			X^i
			=
			-\Speed^{-2} 
			(- \Lsmall^A + v^A)
			(-\Lunit^i + v^i).
			\label{E:YSMALLINTERMSOFLSMALLANDVELOCITY}
\end{align}
\end{subequations} 

\end{enumerate}
\end{lemma}


\begin{proof}
All aspects of the lemma except for 
\eqref{E:COVARIANTLUNITDERIVATIVEOFLUNIT}, \eqref{E:LPLUS2XISNULLANDTRANSVERSETOPU}
and
\eqref{E:YSMALLINTERMSOFLSMALLANDVELOCITY} follow from
minor modifications of the proofs of \cite{jSgHjLwW2016}*{(2.12), (2.13) and Lemma 2.1}.
The identity \eqref{E:COVARIANTLUNITDERIVATIVEOFLUNIT} follows from
definition \eqref{E:LUNIT},
\eqref{E:LGEOISGEODESICANDNULL},
and the Leibniz rule for the connection $\Dfour$. The identity \eqref{E:LPLUS2XISNULLANDTRANSVERSETOPU} follows from applying \eqref{E:LUNITISNULL} and \eqref{E:GOFLANDX} to \eqref{E:OLDULUNIT}.
The identity \eqref{E:YSMALLINTERMSOFLSMALLANDVELOCITY} follows from definitions
\eqref{E:XSMALL}--\eqref{E:YCOMMUTATOR}, the form \eqref{E:ACOUSTICALMETRIC} of $\gfour_{\alpha \beta}$,
and \eqref{E:XSMALLINTERMSOFLSMALLANDVELOCITY}.
\end{proof}


\subsection{First fundamental forms and metric duality}
\label{SS:FIRSTFUDNAMENTALFORMS}

\begin{definition}[First fundamental forms] \hfill
\label{D:FIRSTFUNDAMENTALFORMS}
\begin{enumerate}
\item We define $g$, the first fundamental form of $\Sigma_t$ relative to $\gfour$, 
	to be the symmetric type $\binom{0}{2}$ tensorfield $\Sigmatproject \gfour$. Note that
	$g(Y,Z) = \gfour(Y,Z)$ for
	all pairs $(Y,Z)$ of $\Sigma_t$-tangent vectorfields.
	We define the corresponding inverse first fundamental form $g^{-1}$
	to be the symmetric type $\binom{2}{0}$ tensorfield that is $\gfour$-dual to $g$
	i.e., 
	$(g^{-1})^{\alpha \beta} 
	\eqdef 
	(\gfour^{-1})^{\alpha \widetilde{\alpha}} 
	(\gfour^{-1})^{\beta \widetilde{\beta}}
	g_{\widetilde{\alpha} \widetilde{\beta}}$.
	Note that the restriction of $g$ to $\Sigma_t$-tangent tensorfields is the Riemannian\footnote{It is Riemannian
	is because the $\gfour$-normal to $\Sigma_t$ is the $\gfour$-timelike vectorfield $\Transport$; see
	Lemma\,\ref{L:BASICPROPERTIESOFVECTORFIELDS}.} 
	metric on $\Sigma_t$ induced by $\gfour$. In particular, 
	relative to the Cartesian spatial coordinates, we have: 
	\begin{align} \label{E:CARTESIANCOMPONENTSOFFIRSTFUNDAMENTLAFORMOFSIGMAT}
		g_{ij} 
		& = \gfour_{ij}
		= \Speed^{-2} \updelta_{ij},
	\end{align}
	where $\updelta_{ij}$ denotes the Kronecker delta,
	and to obtain the last equality in \eqref{E:CARTESIANCOMPONENTSOFFIRSTFUNDAMENTLAFORMOFSIGMAT}, 
	we have used \eqref{E:ACOUSTICALMETRIC}.
 \item We define $\gtorus$, the first fundamental form of the acoustic tori $\ell_{t,u}$ relative to $\gfour$, 
	to be the symmetric type $\binom{0}{2}$ tensorfield $\smoothtorusproject \gfour$. Note that
	$\gtorus(Y,Z) = \gfour(Y,Z)$ for
	all pairs $(Y,Z)$ of $\ell_{t,u}$-tangent vectorfields.
	Note that the restriction of $\gtorus$ to $\ell_{t,u}$-tangent tensorfields is the Riemannian\footnote{It is Riemannian
	is because $\ell_{t,u}$ is a sub-manifold of the spacelike hypersurface $\Sigma_t$.} 
	metric on $\ell_{t,u}$ induced by $\gfour$.
	We define the corresponding inverse first fundamental form $\gtorus^{-1}$
	to be the symmetric type $\binom{2}{0}$ tensorfield that is $\gfour$-dual to $\gtorus$
	i.e., 
	$(\gtorus^{-1})^{\alpha \beta} 
	\eqdef 
	(\gfour^{-1})^{\alpha \widetilde{\alpha}} 
	(\gfour^{-1})^{\beta \widetilde{\beta}}
	\gtorus_{\widetilde{\alpha} \widetilde{\beta}}$.
\end{enumerate}
\end{definition}

\subsubsection{Metric duality and musical notation}
\label{SSS:METRICDUALITYANDMUSICALNOTATION}

\begin{definition}[Metric duality and musical notation]
	\label{D:METRICDUALITY}
	If $V$ is a spacetime vectorfield, then $V_{\flat}$ denotes the one-form that is
	$\gfour$-dual to $V$, i.e., $(V_{\flat})_{\alpha} = \gfour_{\alpha \beta} V^{\beta}$.
	Consistent with the conventions of Sect.\,\ref{SS:NOTATIONANDCONVENTIONS},
	we typically write $V_{\alpha}$ instead of $(V_{\flat})_{\alpha}$.
	
	If $Y = Y^B \geop{x^B}$ is an $\ell_{t,u}$-tangent vectorfield, 
	then $Y_{\flat}$ denotes the $\ell_{t,u}$-tangent one-form that is $\gtorus$-dual to $Y$,
	i.e., for $A=2,3$, $(Y_{\flat})_A = \gtorusdoublearg{A}{B} Y^B$. 
	Since $Y$ is $\ell_{t,u}$-tangent, $Y_{\flat}$ can also be viewed as the dual of
	$Y$ with respect to $\gfour$, i.e., 
	$(Y_{\flat})_{\alpha} \eqdef Y_{\flat} \cdot \partial_{\alpha} = \gfour(Y,\partial_{\alpha})
	= \gfour_{\alpha \beta} Y^{\beta}
	$.
	Similarly, if $\upxi$ is an $\ell_{t,u}$-tangent one-form, 
	then $\upxi^{\sharp}$ denotes the $\ell_{t,u}$-tangent vectorfield that is $\gtorus$-dual to $\upxi$,
	i.e., for $A=2,3$, $(\upxi^{\#})^A = (\gtorus^{-1})^{AB} \upxi_B$, where 
	$\upxi_B = \upxi \cdot \geop{x^B}$.
	Similarly, if $\upxi$ is a symmetric type $\binom{0}{2}$ $\ell_{t,u}$-tangent tensorfield, 
	then $\upxi^{\sharp}$ denotes the type $\binom{1}{1}$ $\ell_{t,u}$-tangent tensorfield
	obtained by raising one index of $\upxi$ with $\gtorus^{-1}$,
	while
	$\upxi^{\sharp \sharp}$ denotes the type $\binom{2}{0}$ $\ell_{t,u}$-tangent tensorfield
	obtained by raising both indices of $\upxi$ with $\gtorus^{-1}$.
\end{definition}


\subsection{Projection tensorfields and hypersurface tangency}
\label{SS:PROJECTIONSANDRELATEDDIFFERENTIALOPERATORS}

\begin{definition}[Projection tensorfields and hypersurface tangency]
\label{D:PROJECTIONTENSORFIELDSANDTANGENCYTOHYPERSURFACES}
\hfill
\begin{enumerate}
\item We define the type $\binom{1}{1}$
$\Sigma_t$-projection tensorfield $\Sigmatproject$ and the type $\binom{1}{1}$ $\ell_{t,u}$-projection tensorfield 
$\smoothtorusproject$ as follows, where $\updelta_{\beta}^{\alpha}$ denotes the Kronecker delta: 
\begin{subequations}
\begin{align}
\Sigmatproject_{\beta}^{\ \alpha}
& \eqdef 
	\updelta_{\beta}^{\alpha}
	+ 
	\Transport^{\alpha} 
	\Transport_\beta,  
	\label{E:CARTESIANSIGMAPROJECT} 
	\\
\smoothtorusproject^{\ \alpha}_\beta 
& 
\eqdef 
\updelta_{\beta}^{\alpha}
+  
\Transport^{\alpha} \Transport_\beta 
- 
X^{\alpha} X_\beta 
= 
\updelta_{\beta}^{\alpha}
- 
\Lunit^{\alpha} \updelta^0_\beta 
+
X^{\alpha} \Lunit_{\beta}.
 \label{E:SMOOTHTORUSPROJECT}
\end{align} 
\end{subequations}\item Given any type $\binom{m}{n}$ spacetime tensorfield $\upxi$, 
we respectively define its $\gfour$-orthogonal projection onto $\Sigma_t$, denoted by $\Sigmatproject \upxi$, its $\gfour$-orthogonal projection onto $\ell_{t,u}$, denoted by $\smoothtorusproject \upxi$, as follows:
\begin{subequations}
\begin{align}
(\Sigmatproject \upxi)_{\beta_1 \cdots \beta_n}^{\alpha_1 \cdots \alpha_m}
& 
\eqdef 
\Sigmatproject_{\widetilde{\alpha}_1}^{\ \alpha_1}
\cdots 
\Sigmatproject_{\widetilde{\alpha}_m}^{\ \alpha_m}
\Sigmatproject_{\beta_1}^{\ \widetilde{\beta}_1}
\cdots 
\Sigmatproject_{ \beta_n}^{\ \widetilde{\beta}_n}
\upxi_{\widetilde{\beta}_1 \cdots \widetilde{\beta}_n}^{\widetilde{\alpha}_1 \cdots \widetilde{\alpha}_m}, 
	\label{E:PROJECTIONOFTENSORONTOCARTESIANSIGMAT} 
	\\
(\smoothtorusproject \upxi)_{\beta_1\cdots\beta_n}^{\alpha_1\cdots\alpha_m} 
& \eqdef 
\smoothtorusproject_{\widetilde{\alpha}_1}^{\ \alpha_1} 
\cdots 
\smoothtorusproject_{\widetilde{\alpha}_m}^{\ \alpha_m} 
\smoothtorusproject_{ \beta_1}^{\ \widetilde{\beta}_1}
\cdots 
\smoothtorusproject_{ \beta_n}^{\ \widetilde{\beta}_n}
\upxi_{\widetilde{\beta}_1 \cdots \widetilde{\beta}_n}^{\widetilde{\alpha}_1 \cdots \widetilde{\alpha}_m}.
	\label{E:PROJECTIONOFTENSORONTOFLATTORUS} 
\end{align}
\end{subequations}\item
We say that a spacetime tensorfield $\upxi$ is $\Sigma_t$-tangent if $\Sigmatproject \upxi = \upxi$. 
We say that a spacetime tensorfield $\upxi$ is $\ell_{t,u}$-tangent if $\smoothtorusproject \upxi = \upxi$.  
\item If $\upxi$ is a symmetric type $\binom{0}{2}$-spacetime tensor and $V$ is a vectorfield, 
then we define $\slashed{\upxi}_V \eqdef \smoothtorusproject(\upxi \cdot V)$, 
where $\upxi \cdot V$ is the one-form with components
$(\upxi \cdot V)_{\alpha} \eqdef \upxi_{\alpha \beta} V^{\beta}$.
\item If $\upxi$ is a spacetime tensor, then we define $\slashed{\upxi} = \smoothtorusproject \upxi$.
From 3 above, it follows that $\xi$ is $\ell_{t,u}$-tangent if and only if $\slashed{\upxi} = \xi$.
\end{enumerate}
\end{definition}

\subsection{Coordinate partial derivative transformations} \label{SS:COORDINATEPARTIALDERIVATIVETRANSFORMATIONS}


\begin{lemma}[\cite{abbrescia2022emergence}*{Lemma 6.4}; Geometric coordinate vectorfields in terms of the Cartesian ones] The following identities hold,
where 
$\left\lbrace 
	\geop{t}, \geop{u}, \geop{x^2}, \geop{x^3} 
\right\rbrace
$ 
are the geometric coordinate partial derivative vectorfields, $\lbrace \partial_t, \partial_1, \partial_2, \partial_3 \rbrace$ are the Cartesian coordinate partial derivative vectorfields:
\label{L:GEOMETRICVECTORFIELDSINTERMSOFCARTESIANONES}
\begin{align} \label{E:GEOMETRICVECTORFIELDSINTERMSOFCARTESIANONES}
\geop{t} 
& 
= \p_t 
+ 
\left\lbrace \frac{L^1 X^1 + L^2 X^2 + L^3 X^3}{X^1} \right\rbrace \p_1, 
& 
\geop{u} 
& 
=  \frac{\upmu  \Speed^2}{X^1} \p_1, 
& 
\geop{x^2} & = \p_2 - \frac{X^2}{X^1} \p_1, 
&
\geop{x^3} 
& = \p_3 - \frac{X^3}{X^1} \p_1.
\end{align}
\end{lemma}

\begin{lemma}[\cite{abbrescia2022emergence}*{Lemma 6.5}; Relationship between $\left\lbrace\geop{t},\geop{u},\geop{x^2},\geop{x^3}\right\rbrace$ and $\{\Lunit,X,\Yvf{2},\Yvf{3}\}$] 
\label{L:COMMUTATORSTOCOORDINATES}
The commutation vectorfields from Def.\,\ref{D:COMVECTORFIELDS}
can be expressed as follows:
\begin{subequations} 
\begin{align}
\Lunit 
& = \geop{t} + \Lunit^A \geop{x^A},  
	\label{E:LUNITINTERMSOFGEOMETRICCOORDINATEVECTORFIELDS} 
	\\
\muX 
& = \geop{u} + \upmu X^A \geop{x^A}, 
	\label{E:MUXINTERMSOFGEOMETRICCOORDINATEVECTORFIELDS} 
	\\
\Yvf{2} 
& = \left\lbrace 1 - \Speed^{-2}(X^2)^2 \right\rbrace 
		\geop{x^2}
		- \Speed^{-2} X^2 X^3 \geop{x^3}, 
	\label{E:Y2INTERMSOFGEOMETRICCOORDINATEVECTORFIELDS}
		\\
\Yvf{3} 
& = \left\lbrace 1 - \Speed^{-2} (X^3)^2 \right\rbrace 
		\geop{x^3}
		- 
		\Speed^{-2} X^2 X^3 \geop{x^2}, 
	\label{E:Y3INTERMSOFGEOMETRICCOORDINATEVECTORFIELDS}
\end{align}
\end{subequations}Moreover, the following identities hold:
\begin{subequations}
\begin{align}
	\geop{t} 
	& = \Lunit - \Lunit^A \Yvf{A} 
	- 
		\frac{L^2 X^2 + L^3 X^3}{(X^1)^2} X^A \Yvf{A}, 
		\label{E:GEOPTOCOMMUTATORS} 
		\\
	\geop{u} 
	& = \muX - \frac{1}{(X^1)^2} \upmu \Speed^2 X^A \Yvf{A}, 
		\label{E:GEOPUTOCOMMUTATORS} 
		\\
	\geop{x^2} 
	& = 
	\left\{ \frac{(X^1)^2 + (X^2)^2}{(X^1)^2} \right\} \Yvf{2} +  \frac{X^2 X^3}{(X^1)^2} \Yvf{3},
		\label{E:GEOP2TOCOMMUTATORS} 
		\\
	\geop{x^3} 
	& = 
	\left\{ \frac{(X^1)^2 + (X^3)^2}{(X^1)^2} \right\} \Yvf{3} +  \frac{X^2 X^3}{(X^1)^2} \Yvf{2}.		\label{E:GEOP3TOCOMMUTATORS}
\end{align}
\end{subequations}
\end{lemma}

\begin{lemma}[\cite{abbrescia2022emergence}*{Lemma 6.6}; Relationship between $\{\p_t,\p_1,\p_2,\p_3\}$ and $\{\Lunit,X, \Yvf{2},\Yvf{3}\}$] \label{L:RELATIONSHIPBETWEENCARTESIANPARTIALDERIVATIVESANDSMOOTHGEOMETRICCOMMUTATORS} The following identities hold:
\begin{subequations} 
\begin{align}
\p_t 
& = \Lunit  - \frac{L^1X^1 + L^2 X^2 + L^3 X^3}{\Speed^2} X + \frac{\Lunit^1}{X^1} X^C \Yvf{C} - \Lunit^C \Yvf{C}, 
	\label{E:CARTESIANPARTIALTTOCOMMUTATORS} 
		\\
\p_1 
& = \frac{X^1}{\Speed^2} X - \frac{1}{X^1} X^A \Yvf{A}.
	\label{E:CARTESIANPARTIAL1TOCOMMUTATORS}
	\\
\p_A 
& = \frac{X^A}{\Speed^2} X + \Yvf{A}.
 \label{E:CARTESIANPARTIAL23TOCOMMUTATORS}
\end{align}
\end{subequations}
\end{lemma}

\begin{corollary}[\cite{abbrescia2022emergence}*{Corollary 6.7}; Expressions for $\smoothtorusproject_{\beta}^{\ \alpha} \partial_{\alpha}$ in terms of $\lbrace \Yvf{2}, \Yvf{3} \rbrace$]
	\label{C:ELLTUPROJECTEDVERSIONOFCARTESIANPARTIALDERIVATIVES}
	The following identities hold in Cartesian coordinates:
	\begin{align} \label{E:ELLTUPROJECTEDVERSIONOFCARTESIANPARTIAL0}
		\smoothtorusproject_0^{\ \alpha} \partial_{\alpha}
		& 
		= 
		\frac{\Lunit^1}{X^1} X^A \Yvf{A} - \Lunit^A \Yvf{A},
			\\
			\smoothtorusproject_1^{\ \alpha} \partial_{\alpha}
		& 
		= 
		- \frac{1}{X^1} X^A \Yvf{A},
			\label{E:ELLTUPROJECTEDVERSIONOFCARTESIANPARTIAL1} 
				\\
		\smoothtorusproject_2^{\ \alpha} \partial_{\alpha}
		& = \Yvf{2},
		&
		\smoothtorusproject_3^{\ \alpha} \partial_{\alpha}
		& = \Yvf{3}.
		\label{E:ELLTUPROJECTEDVERSIONOFCARTESIANPARTIAL2AND3}
	\end{align}
\end{corollary}

\subsection{Differential operators associated with the projections and metrics}
\label{SS:DIFFERENTIALOPERATORSASSOCIATEDWITHPROJECTIONSANDMETRICS}

\begin{definition}[$\ell_{t,u}$-differential] 
\label{D:ANGULARDIFFERENTIAL}
If $\varphi$ is a scalar function, then we define $\angrmd \varphi$ 
to be the following $\ell_{t,u}$-tangent one-form:
\begin{align} \label{E:ANGULARDIFFERENTIAL}
\angrmd \varphi 
& \eqdef \smoothtorusproject \rmd \varphi.
\end{align}
\end{definition}
Note that $\argangrmd{A} \varphi = \angrmd \varphi \cdot \geop{x^A}  = \geop{x^A} \varphi$ for $A = 2,3$
and that $\argangrmd{\alpha} \varphi 
= 
\angrmd \varphi \cdot \partial_{\alpha} 
=  
\smoothtorusproject^{\ \beta}_{\alpha} \partial_{\beta} \varphi$ 
for $\alpha = 0,1,2,3$.

\begin{definition}[Levi-Civita connections and associated differential operators] \hfill
\label{D:CONNECTIONSANDDIFFERENTIALOPERATORS}
\begin{enumerate}
\item We denote the Levi-Civita connection of $\gfour$ by $\Dfour$. 
\item We denote the Levi-Civita connection of $\gtorus$ by $\newangD$. In particular,
	for $\ell_{t,u}$-tangent tensorfields $\upxi$, we have $\newangD \upxi = \smoothtorusproject \Dfour \upxi$.
\item If $\upxi$ is an $\ell_{t,u}$-tangent one-form, then we define its $\ell_{t,u}$-divergence 
	to be the scalar function
	$\angdiv \upxi \eqdef \mytr_{\gtorus} \angD \upxi =  (\gtorus^{-1})^{AB}\cdot \newangD_A \upxi_B$. Similarly, if $Y$ is an $\ell_{t,u}$-tangent vectorfield, 
	then we define its $\ell_{t,u}$-divergence to be the scalar function $\angdiv Y \eqdef \gtorus^{-1} \cdot \newangD Y_{\flat}$.
	where $Y_{\flat}$ is the $\ell_{t,u}$-tangent one-form that is $\gfour$-dual to $Y$.
\item If $\upxi$ is an $\ell_{t,u}$-tangent one-form, we define its $\ell_{t,u}$-curl to be the scalar function given by $\in^{AB} \angD_A \upxi_B$, where $\in_{AB}$ are the components of the area form element of the Riemannian manifold $(\ell_{t,u},\gtorus)$ relative to an arbitrary frame $e_A$, $A = 2,3$.
\item If $\upxi$ is a symmetric type $\binom{0}{2}$ $\ell_{t,u}$-tangent tensorfield, 
	then we define its $\ell_{t,u}$-divergence 
	to be the $\ell_{t,u}$-tangent one-form with the following $\ell_{t,u}$-components for $A=2,3$:
	$(\angdiv \upxi)_A \eqdef (\gtorus^{-1})^{BC} \cdot \newangDarg{B} \upxi_{CA}$. 
\item We denote the covariant wave operator of $\gfour$ by $\Box_{\gfour} \eqdef \gfour^{-1}\cdot\Dfour^2
	= (\gfour^{-1})^{\alpha \beta} \Dfour_{\alpha} \Dfour_{\beta}
$. 
\item We denote the angular Hessian and angular Laplacian on $\ell_{t,u}$ associated to $\gtorus$ by $\angD^2$ and 
$\angLap \eqdef \gtorus^{-1} \cdot \newangD^2
=
(\gtorus^{-1})^{AB} \angDarg{A} \angDarg{B}
$, respectively.
\end{enumerate}
\end{definition}


\begin{definition}[Projected Lie derivatives] 
\label{D:PROJECTEDLIEDERIVATIVES}
Given a spacetime tensorfield $\upxi$
and a vectorfield $Z$,
we define $\SigmatLie_Z \upxi$ and $\angLie_Z \upxi$ to respectively be the following $\Sigma_t$-tangent and $\ell_{t,u}$-tangent tensorfields:
\begin{align} \label{E:PROJECTEDLIEDERIVATIVES}
\SigmatLie_Z \upxi 
	&	\eqdef \Sigmatproject \Lie_Z \upxi, 
&
\angLie_Z \upxi 
& \eqdef \smoothtorusproject \Lie_Z \upxi. 
	\end{align}
\end{definition}

\begin{lemma} [\cite{jSgHjLwW2016}*{Lemma 2.10}; Angular differential $\angrmd$ commutes with $\angLie$] 
\label{L:ANGULARDIFFERENTIALCOMMUTESWITHANGLIE}
Let $f$ be a scalar function and let 
$Z \in \Fullset$ (see definition \eqref{E:COMMUTATIONVECTORFIELDS}). 
Then the following identity holds:
\begin{align} \label{E:ANGULARDIFFERENTIALCOMMUTESWITHANGLIE}
\angLie_Z \angrmd f 
& = \angrmd Z f.
\end{align}
\end{lemma}


\subsection{Controlling quantities $\controlvars$ and $\badcontrolvars$} 
\label{SS:CONTROLVARS}
We use the arrays 
$\controlvars$ and $\badcontrolvars$ from the next definition to present
schematic formulas in contexts where precise structure is not important.

\begin{definition}[The controlling quantities] \label{D:CONTROLVARS} 
We define $\controlvars$ and $\badcontrolvars$ to be the following arrays of scalar functions:
\begin{align} 	\label{E:CONTROLVARSGOODANDBAD} 
\controlvars 
& 
\eqdef (\velocityarray,\Ent, \Lsmall^1, \Lsmall^2, \Lsmall^3), 
&
\badcontrolvars 
& \eqdef  
	(\velocityarray,\Ent, \upmu - 1, \Lsmall^1, \Lsmall^2, \Lsmall^3).
\end{align}
Notice that $\LogDensity$ is \emph{not} part of $\controlvars$ or $\badcontrolvars$.
\end{definition} 

For the solutions that we study in our main results, along the data-hypersurface 
$\Sigma_0$, $\controlvars$ and $\badcontrolvars$ are small in $L^{\infty}$.

\subsection{Identities for the $\ell_{t,u}$-projection tensorfield and the first fundamental forms}
\label{SS:IDENTITIESFORFIRSTFUNDSANDSMOOTHTORUSPROJECTION}

\begin{lemma}[Useful identities for the first fundamental forms; \cite{abbrescia2022emergence}*{Lemma\,3.15}]
\label{L:USEFULIDENTITIESFORFIRSTFUNDAMENTALFORM} 
Let $\smoothtorusproject$, 
$\gtorus$,
$g$,
and $X$ be as in Defs.\,\ref{D:FIRSTFUNDAMENTALFORMS},
\ref{D:PROJECTIONTENSORFIELDSANDTANGENCYTOHYPERSURFACES},
and\,\ref{D:COMVECTORFIELDS}.
Then the following identities hold relative to the geometric coordinates ($A,B = 2,3$):
\begin{subequations}
\begin{align}
\gtorus_{AB} & = \Speed^{-2} \updelta_{AB} +\Speed^{-2} \frac{X^AX^B}{(X^1)^2}, 
	\label{E:SMOOTHTORIGABEXPRESSION} 
		\\ 
(\gtorus^{-1})^{AB}
& 
= \Speed^2 \updelta^{AB} - X^A X^B, 
	\label{E:SMOOTHGINVERSEABEXPRESSION} 
		\\
\mydet \gtorus & = \frac{1}{\Speed^2 (X^1)^2}. 
\label{E:DETERMSMOOTHGTORUSRELTOGEOMETRICCOORDS}
\end{align}
\end{subequations}Moreover, the following identities hold relative to arbitrary coordinates:
\begin{subequations}
	\begin{align} \label{E:FIRSTFUNDOFSIGMATINTERMSOFACOUSTICALMETRICANDTRANSPORT}
		g_{\alpha \beta}
		& =
		\gfour_{\alpha \beta} 
		+
		\Transport_{\alpha} \Transport_{\beta},
			\\
		(g^{-1})^{\alpha \beta}
		& 
		=
		(\gfour^{-1})^{\alpha \beta} 
		+
		\Transport^{\alpha} \Transport^{\beta}.
			\label{E:INVERSEFIRSTFUNDOFSIGMATINTERMSOFACOUSTICALMETRICANDTRANSPORT}
	\end{align}
\end{subequations}Furthermore, the following identities hold relative to the Cartesian coordinates:
\begin{subequations}
	\begin{align}
		g
		& = \Speed^{-2} \sum_{a=1}^3 (\rmd x^a - v^a \rmd t)\otimes  (\rmd x^a - v^a \rmd t),
			\label{E:FIRSTFUNDOFSIGMATIDENTITY} 
				\\
		g^{-1}
		& = \Speed^2 \sum_{a=1}^3 \p_a \otimes \p_a.
		\label{E:INVERSEFIRSTFUNDOFSIGMATIDENTITY}
	\end{align}
\end{subequations}In addition, the following identities hold relative to arbitrary coordinates:
\begin{subequations}
	\begin{align}
			\gtorus_{\alpha \beta}
			& = 
			g_{\alpha \beta} 
			-
			X_{\alpha} X_{\beta}
			=
			\gfour_{\alpha \beta} 
			+
			\Transport_{\alpha} \Transport_{\beta}
			-
			X_{\alpha} X_{\beta}
			= \gfour_{\alpha \beta} 
				+
				\Lunit_{\alpha} \Lunit_{\beta}
				+
				\Lunit_{\alpha} X_{\beta}
				+
				X_{\alpha} \Lunit_{\beta},
				\label{E:SMOOTHTORUSMETRICINTERMSOFSIGMATMETRICANDX} 
				\\
		(\gtorus^{-1})^{\alpha \beta}
		& = 
			(g^{-1})^{\alpha \beta} 
			-
			X^{\alpha} X^{\beta}
			=
			(\gfour^{-1})^{\alpha \beta} 
			+
			\Transport^{\alpha} \Transport^{\beta}
			-
			X^{\alpha} X^{\beta}
			=
				(\gfour^{-1})^{\alpha \beta}
				+
				\Lunit^{\alpha} \Lunit^{\beta}
				+
				\Lunit^{\alpha} X^{\beta}
				+
				X^{\alpha} \Lunit^{\beta},
			\label{E:SMOOTHTORUSINVERSEMETRICINTERMSOFINVERSESIGMATMETRICANDX}
				\\
		\smoothtorusproject^{\ \alpha}_\beta 
		& = 
		\gfour_{\beta \gamma}
		(\gtorus^{-1})^{\alpha \gamma}
		=
		\gtorus_{\beta \gamma}
		(\gtorus^{-1})^{\alpha \gamma}.
		\label{E:LTUSPROJECTIONISONEINDEXLOWERINGOFGTORUSINVERSE}
	\end{align}
\end{subequations}Finally, relative to the Cartesian coordinates, the following identities hold for $\alpha, \beta = 0,1,2,3$:
\begin{subequations}
\begin{align} \label{E:GTORUSINVERSE0COMPONENTSVANISH}
		(\gtorus^{-1})^{0 \alpha}
		& 
		=
		(\gtorus^{-1})^{\alpha 0}
		= 0,
			\\
		(g^{-1})^{0 \alpha}
		& 
		=
		(g^{-1})^{\alpha 0}
		= 0,
		\label{E:FIRSTFUNDSIGMAT0COMPONENTSVANISH}
			\\
	\smoothtorusproject^{\ 0}_{\beta}
	& = 0.
	\label{E:SMOOTHTORUSUPPER0COMPONENTSVANISH}
\end{align}
\end{subequations}\end{lemma}

\subsection{Identities for $\gtorus$-gradients} \label{SS:GTORUSGRADIENTS}

\begin{lemma}[Identities for $\gtorus$-gradients] \label{L:GTORUSGRADIENTS}
The following identity holds for scalar functions $\phi$ and $\psi$: 
	\begin{align}
		\gtorus^{-1}(\angrmD \phi, \angrmD \psi) & = \Speed^2 \sum_{A=2,3} (\Yvf{A} \phi) \Yvf{A} \psi + \frac{\Speed^2}{(X^1)^2} (X^A \Yvf{A} \phi) X^B \Yvf{B} \psi. \label{E:GTORUSINNERPRODUCTFORSMOOTHTORIGRADIENTS}
	\end{align}
	
\end{lemma}

\begin{proof}
Identities \eqref{E:GEOP2TOCOMMUTATORS}--\eqref{E:GEOP3TOCOMMUTATORS} and \eqref{E:SMOOTHTORIGABEXPRESSION} imply $\geop{x^A} = \Speed^2 \sum_{B=2,3} \gtorus_{AB} \Yvf{B}$. Applying this identity to the LHS \eqref{E:GTORUSINNERPRODUCTFORSMOOTHTORIGRADIENTS}, it follows that $\gtorus^{-1}(\angrmD \phi, \angrmD \psi) = (\gtorus^{-1})^{AB} \geop{x^A}\phi \geop{x^B} \psi = \sum_{A,B= 2,3} \Speed^4 \gtorus_{AB} \Yvf{A} \phi \Yvf{B} \psi$. The desired identity \eqref{E:GTORUSINNERPRODUCTFORSMOOTHTORIGRADIENTS} then follows from \eqref{E:SMOOTHTORIGABEXPRESSION}.
\end{proof}


\subsection{Traces of tensorfields} 
\label{SS:TRACEOFTENSORS}

\begin{definition}[Traces of tensorfields]
\label{D:TRACEOFTENSORS} 
If $\upxi$ is a type $\binom{0}{2}$ spacetime tensorfield, 
then we define its $\gfour$-trace and its $\gtorus$-trace as follows:
\begin{align} \label{E:TRACES}
	\mytr_{\gfour}\upxi 
	& \eqdef (\gfour^{-1})^{\alpha \beta} \upxi_{\alpha \beta},
	&
	\mytr_{\gtorus} \upxi 
& \eqdef 
(\gtorus^{-1})^{\alpha \beta} \upxi_{\alpha \beta}.
\end{align}
\end{definition} 

\subsection{Pointwise norms and semi-norms of tensorfields}
\label{SS:POINTWISESEMINORMSOFTENSORS}
In the next definition, we define various pointwise norms and semi-norms
that we will use to measure the size of tensorfields.

\begin{definition}[Pointwise norms] \label{D:POINTWISESEMINORMS} 
If $\upxi$ is a type $\binom{m}{n}$ spacetime tensorfield,
then we define\footnote{While RHS\,\eqref{E:SQUAREPOINTWISELORENTZIANNORMWITHRESPECTTOACOUSTICALMETRIC} can be
negative, we will only ever use the notation ``$|\upxi|_{\gfour}$'' when it is positive. \label{FN:NONEGATIVENORMS}}
$|\upxi|_{\gfour} \geq 0$,
$|\upxi|_{\gtorus} \geq 0$,
and
$|\upxi|_g \geq 0$ by:
\begin{subequations}
\begin{align} \label{E:SQUAREPOINTWISELORENTZIANNORMWITHRESPECTTOACOUSTICALMETRIC}
|\upxi|_{\gfour}^2 
& \eqdef \gfour_{\alpha_1 \widetilde{\alpha}_1} \cdots \gfour_{\alpha_m \widetilde{\alpha}_m} (\gfour^{-1})^{\beta_1 \widetilde{\beta}_1} \cdots (\gfour^{-1})^{\beta_n \widetilde{\beta}_n} \upxi_{\beta_1\cdots\beta_n}^{\alpha_1\cdots\alpha_n} \upxi_{\widetilde{\beta}_1 \cdots \widetilde{\beta}_n}^{\widetilde{\alpha}_1\cdots\widetilde{\alpha}_m},
	\\
|\upxi|_{\gtorus}^2 
& \eqdef \gtorus_{\alpha_1 \widetilde{\alpha}_1} \cdots \gtorus_{\alpha_m \widetilde{\alpha}_m} (\gtorus^{-1})^{\beta_1 \widetilde{\beta}_1} \cdots (\gtorus^{-1})^{\beta_n \widetilde{\beta}_n} \upxi_{\beta_1\cdots\beta_n}^{\alpha_1\cdots\alpha_n} \upxi_{\widetilde{\beta}_1 \cdots \widetilde{\beta}_n}^{\widetilde{\alpha}_1\cdots\widetilde{\alpha}_m},
	\label{E:SQUAREPOINTWISESEMINORMWITHRESPECTTOFIRSTFUNDOFSMOOTHTORI} \\
|\upxi|_{g}^2 
& \eqdef g_{\alpha_1 \widetilde{\alpha}_1} \cdots g_{\alpha_m \widetilde{\alpha}_m} (g^{-1})^{\beta_1 \widetilde{\beta}_1} \cdots (g^{-1})^{\beta_n \widetilde{\beta}_n} \upxi_{\beta_1\cdots\beta_n}^{\alpha_1\cdots\alpha_n} \upxi_{\widetilde{\beta}_1 \cdots \widetilde{\beta}_n}^{\widetilde{\alpha}_1\cdots\widetilde{\alpha}_m}.
	\label{E:SQUAREPOINTWISESEMINORMWITHRESPECTTOFIRSTFUNDOFSIGMAT}
\end{align}
\end{subequations}
\end{definition}

\begin{remark}[Norms vs.\ semi-norms]
\label{R:NORMSANDSEMINORMS}
$|\cdot|_{\gfour}$ is a pointwise norm on the space of $\gfour$-spacelike tensorfields.
$|\cdot|_g$ is a pointwise norm on the space of $\Sigma_t$-tangent tensorfields and a pointwise semi-norm on the space of all 
tensorfields.
$|\cdot|_{\gtorus}$ is a pointwise norm on the space of $\ell_{t,u}$-tangent tensorfields and a pointwise semi-norm on the space of all 
tensorfields.

\end{remark}

\begin{remark}[Omitting the $0$ component in Cartesian coordinates]
\label{R:OMITTINGZEROCOMPONENTINCARTESIANCOORDINATES}
In view of 
\eqref{E:MATERIALDERIVATIVELOWEREDCARTESIANCOORDINATES},
Def.\,\ref{D:PROJECTIONTENSORFIELDSANDTANGENCYTOHYPERSURFACES},
definitions 
\eqref{E:SQUAREPOINTWISESEMINORMWITHRESPECTTOFIRSTFUNDOFSMOOTHTORI}--\eqref{E:SQUAREPOINTWISESEMINORMWITHRESPECTTOFIRSTFUNDOFSIGMAT}, 
and
\eqref{E:GTORUSINVERSE0COMPONENTSVANISH}--\eqref{E:FIRSTFUNDSIGMAT0COMPONENTSVANISH},
we see that relative to the Cartesian coordinates, we have:
\begin{align} \label{E:NO0COMPONENTSQUAREPOINTWISESEMINORMWITHRESPECTTOFIRSTFUNDOFSMOOTHTORI}
|\upxi|_{\gtorus}^2 
& 
= 
\gtorus_{a_1 \widetilde{a}_1} 
\cdots 
\gtorus_{a_m \widetilde{a}_m} 
(\gtorus^{-1})^{b_1 \widetilde{b}_1} 
\cdots 
(\gtorus^{-1})^{b_n \widetilde{b}_n} 
\upxi_{b_1 \cdots b_n}^{a_1 \cdots a_n}
\upxi_{\widetilde{b}_1 \cdots \widetilde{b}_n}^{\widetilde{a}_1 \cdots \widetilde{a}_m},
	\\
|\upxi|_g^2 
& = 
g_{a_1 \widetilde{a}_1} 
\cdots 
g_{a_m \widetilde{a}_m} 
(g^{-1})^{b_1 \widetilde{b}_1} 
\cdots 
(g^{-1})^{b_n \widetilde{b}_n} 
\upxi_{b_1 \cdots b_n}^{a_1 \cdots a_n} 
\upxi_{\widetilde{b}_1 \cdots \widetilde{b}_n}^{\widetilde{a}_1 \cdots\widetilde{a}_m},
\label{E:NO0COMPONENTSQUAREPOINTWISESEMINORMWITHRESPECTTOFIRSTFUNDOFSIGMAT}
\end{align}
i.e., we can omit all ``$0$'' components on 
RHSs~\eqref{E:NO0COMPONENTSQUAREPOINTWISESEMINORMWITHRESPECTTOFIRSTFUNDOFSMOOTHTORI}--\eqref{E:NO0COMPONENTSQUAREPOINTWISESEMINORMWITHRESPECTTOFIRSTFUNDOFSIGMAT}.
Similarly, taking into account definition~\eqref{E:PROJECTIONOFTENSORONTOFLATTORUS} and \eqref{E:SMOOTHTORUSUPPER0COMPONENTSVANISH}, 
we see that if $\SigmatTan$ is a $\Sigma_t$-tangent vectorfield, 
then relative to the Cartesian coordinates, we have
$\smoothtorusproject_{\beta}^{\ \alpha} \partial_{\alpha} \SigmatTan^{\beta}
=
\smoothtorusproject_b^{\ a} \partial_a \SigmatTan^b
$.
We will silently use these basic facts throughout.
\end{remark}


\subsection{Second fundamental forms and the torsion}
\label{SS:SECONDFUNDAMENTALFORMSANDTORSION}

\begin{definition}[The second fundamental forms $k$ and $\upchi$, and the one-form $\upzeta$] 
\label{D:SECONDFUNDAMENTALFORMSANDZETAONEFORM}
\hfill
\begin{enumerate}
\item We define the \emph{second fundamental form} $k$ of $\Sigma_t$ as follows:
\begin{align} \label{E:SIGMATSECONDFUND}
	k 
	& 
	\eqdef \frac{1}{2} \SigmatLie_{\Transport} g.
\end{align}
\item We define the \emph{null second fundamental form} of $\ell_{t,u}$ as follows: 
\begin{align} \label{E:NULLSECONDFUND}
	\upchi  
	& 
	\eqdef \frac{1}{2} \angLie_\Lunit \gtorus.
\end{align}
\item We define $\upzeta$ to be the $\ell_{t,u}$-tangent one-form with the following components:
\begin{align} \label{E:TORISONTENSORFIELD}
	\upzeta_A 
	& 
	\eqdef 
	\gfour(\Dfour_A \Lunit, X).
\end{align}
\end{enumerate}
\end{definition}

\subsection{Transport equations for the eikonal function quantities}
\label{SS:TRANSPORTFORMUANDLUNITI}
To control the eikonal function quantities $\upmu$ and $\Lunit^i$, we will use the following transport equations.

\begin{lemma}[\cite{jSgHjLwW2016}*{Lemma 2.12}; Transport equations satisfied by $\upmu$ and $\Lunit^i$] \label{L:TRANSPORTMUANDLUNITI}
The scalar functions $\upmu$ and $\Lunit^i$ satisfy the following transport equations
(see Def.\,\ref{D:DERIVATIVESOFARRAYS} regarding the notation):
\begin{align} \label{E:MUTRANSPORT}
\Lunit \upmu 
& 
= 
\frac{1}{2} \vec{G}_{\Lunit \Lunit} \diamond \muX \wavearray 
- 
\frac{1}{2} \upmu \vec{G}_{\Lunit \Lunit} \diamond \Lunit \wavearray 
- 
\upmu \vec{G}_{\Lunit X} \diamond \Lunit \wavearray, 
	\\
\Lunit \Lsmall^i 
& 
= \frac{1}{2} (\vec{G}_{\Lunit \Lunit}\diamond \Lunit \wavearray) X^i 
	- 
	(\angG_{\Lunit}^{\#} \diamond \Lunit \wavearray)\cdot \angrmd x^i 
	+ 
	\frac{1}{2} (\vec{G}_{\Lunit \Lunit} \diamond \angrmd^{\#}\wavearray)\cdot \angrmd x^i. 
	\label{E:LUNITITRANSPORT}
\end{align}
\end{lemma}

\subsection{The factor driving the shock formation and formulas involving $G_{\Lunit \Lunit}$}
\label{SS:THEFACTORDRIVINGTHESHOCKFORMATION}
In the following lemma, we compute
an expression for the product $\frac{1}{2} \vec{G}_{\Lunit \Lunit} \diamond \muX \wavearray$
on the RHS of the evolution equation \eqref{E:MUTRANSPORT} for $\upmu$.
For every smooth equation of state besides that of a Chaplygin gas, 
there exist open sets of background densities $\overline{\varrho} > 0$ such that
the non-degeneracy condition \eqref{E:NONDEGENCONDITION} holds.
The identity \eqref{E:IDENTITYFORMAINTERMDRIVINGTHESHOCK} then shows
that for solutions that are close to the trivial solution $\wavearray \equiv 0$,
the expansion of $\frac{1}{2} \vec{G}_{\Lunit \Lunit} \diamond \muX \wavearray$
features a non-zero term proportional to $\muX \RRiemann$; 
the presence of this term is crucial for our main results, as it drives the formation of the shock, 
i.e., it drives $\upmu$ to $0$.
In contrast, for the equation of state $p = C_0 - C_1 \exp(- \LogDensity)$ of a Chaplygin gas, 
one can compute that $\Speed^{-1} \Speed_{;\LogDensity} + 1 \equiv 0$, and the non-degeneracy condition
\eqref{E:NONDEGENCONDITION} is therefore impossible. 
In this case, equation \eqref{E:IDENTITYFORMAINTERMDRIVINGTHESHOCK} shows that
the product $\frac{1}{2} \vec{G}_{\Lunit \Lunit} \diamond \muX \wavearray$ 
does not depend on the solution's $\muX$ derivative, and hence our main results do not apply. 
We note that one can show that for irrotational and isentropic solutions,
the equation $\Speed^{-1} \Speed_{;\LogDensity} + 1 = 0$ 
is equivalent to the statement that the quasilinear wave equation for a potential function
satisfies Klainerman's null condition \cite{sK1984}.

\begin{lemma}[Identity for the factor driving the shock formation] 
\label{L:IDENTIFYFORFACTORDRIVINGSHOCK}
For solutions to the compressible Euler equations
\eqref{E:BVIEVOLUTION}--\eqref{E:BENTROPYEVOLUTION},
the following identity holds,
where $\almostRiemannfunction(\LogDensity,\Ent)$ 
is the scalar function from \eqref{E:ALMOSTRIEMANNINVARIANTS}:
\begin{align} \label{E:IDENTITYFORMAINTERMDRIVINGTHESHOCK}
\begin{split}
\frac{1}{2} \vec{G}_{\Lunit \Lunit}  \muX \wavearray 
& = - \frac{1}{2} \Speed^{-1}(\Speed^{-1} \Speed_{;\LogDensity} + 1)  
			\left\lbrace 
				\muX \RRiemann - \muX \LRiemann 
			\right\rbrace 
	\\
& - \frac{1}{2} \upmu \Speed^{-1} X^1 
\left\lbrace \Lunit \RRiemann + \Lunit \LRiemann\} - \upmu \Speed^{-2} \lbrace X^2 \Lunit v^2 + X^3 \Lunit v^3
\right\rbrace 
	\\
& - \upmu \Speed^{-1} \Speed_{;\Ent} X^a \GradEnt^a + \upmu \Speed^{-1} (\Speed^{-1} \Speed_{;\LogDensity} + 1) 
\almostRiemannfunction_{;\Ent} X^a \GradEnt^a.
\end{split}
\end{align}
\end{lemma}

\begin{proof}
This is the same as \cite{jLjS2018}*{Lemmas 2.45, 2.46}, except for minor
modifications incorporating the third dimension and the entropy
(via the $\Speed_{;\Ent}$-dependent and $\almostRiemannfunction_{;\Ent}$-dependent products).
\end{proof}

In the next lemma, we derive expressions for $G_{\Lunit \Lunit}^2$ and $G_{\Lunit \Lunit}^3$.
When deriving estimates, we will use the expressions to track smallness.
\begin{lemma}[Formulas for $G_{\Lunit \Lunit}^A$; \cite{abbrescia2022emergence}*{Lemma 3.23}] 
\label{L:GLLAEXPRESSION}
The following identities hold for $A= 2,3$:
\begin{align}  \label{E:GLLAEXPRESSION}
G_{\Lunit \Lunit}^A 
& 
= 2 \Speed^{-2}(v^A - \Lunit^A) 
= 2\Speed^{-2} X^A 
= 2\Speed^{-2} \Xsmall^A.
\end{align}
\end{lemma}
%

\subsection{Useful geometric decompositions}
\label{SS:USEFULGEOMETRIC}
\begin{lemma}[\cite{jS2016b}*{Lemma 3.61}; Alternate expressions for $\upchi$ and $\angk$]
\label{L:ALTERNATEEXPRESSIONSFORSECONDFUNDAMENTALFORMS}
The second fundamental forms from Def.\,\ref{D:SECONDFUNDAMENTALFORMSANDZETAONEFORM}
satisfy the following identities:
\begin{align} \label{E:SECONDFUNDSALTERNATE} 
\upchi_{AB} 
& 
= 
\gfour\left(\Dfour_A \Lunit, \geop{x^B}\right), 
&
\angk_{AB} 
& = \gfour\left(\Dfour_A \Transport, \geop{x^B}\right). 
\end{align}
\end{lemma}

\begin{lemma}[Useful identities and decompositions for $\upchi$, $\angk$, and $\upzeta$] 
\label{L:USEFULIDENTITIESANDDECOMPOSITIONSFORSECONDFUNDAMENTALFORMSANDTORSION} 
Let $\upchi$, $k$, and $\upzeta$ be the tensorfields from Def.\,\ref{D:SECONDFUNDAMENTALFORMSANDZETAONEFORM}.
Then the following\footnote{Here, $\otimesarray$ is defined by 
$\angG_{\Lunit} \otimesarray \angrmd \wavearray \eqdef \sum_{\iota = 0}^4 
\angG_{\Lunit}^{\iota} \otimes \angrmd \Psi_\iota$, 
and similarly for $\angrmd \Psi \otimesarray\angG_{\Lunit}$, $\angG_X \otimesarray \angrmd \wavearray$, etc.} identities hold:
\begin{subequations}
\begin{align}
	\upchi & = \gfour_{ab} \angrmd \Lunit^a \otimes \angrmd x^b + \frac{1}{2} \angG \diamond \Lunit \wavearray + \frac{1}{2} \angrmd \wavearray \otimesarray \angG_{\Lunit} - \frac{1}{2} \angG_{\Lunit} \otimesarray \angrmd \wavearray, 
		\label{E:CHIEXPRESSSIONINTERMSOFDERIVATIVESOFLUNITI}
	\\
	\mytr_{\gtorus}\upchi & = \gfour_{ab} \gtorus^{-1} \cdot \left\lbrace \angrmd \Lunit^a \otimes \angrmd x^b\right\rbrace + \frac{1}{2} \gtorus^{-1} \cdot \angG \diamond \Lunit \wavearray. \label{E:TRCHIEXPRESSSIONINTERMSOFDERIVATIVESOFLUNITI}
\end{align}
\end{subequations}Moreover, we can decompose $\angk$ and $\upzeta$ into $\upmu^{-1}$-singular and $\upmu^{-1}$-regular pieces as follows:
\begin{subequations}
\begin{align} \label{E:CONNECTIONCOEFFICIENTDECOMPOSITIONS}
\upzeta 
& = \zetatan + \upmu^{-1} \zetatrans, & \angk & = \angktan + \upmu^{-1}\angktrans,
\end{align}
where:
\begin{align}
\angktan & \eqdef  
\frac{1}{2} \angG \diamond \Lunit \wavearray 
- 
\frac{1}{2} \angG_{\Lunit} \otimesarray \angrmd \wavearray 
- 
\frac{1}{2} \angrmd \wavearray \otimesarray \angG_{\Lunit} 
- 
\frac{1}{2} \angG_X \otimesarray \angrmd \wavearray 
- 
\frac{1}{2} \angrmd \wavearray \otimesarray \angG_X, \label{E:ANGKTAN}
	\\
\angktrans & \eqdef 
\frac{1}{2} \angG\diamond \muX \wavearray, 
\label{E:ANGKTRANS} 
\\
\zetatan 
& \eqdef 
\frac{1}{2} \angG_X \diamond \Lunit \wavearray 
-
\frac{1}{2} \vec{G}_{\Lunit X} \diamond \angrmd \wavearray 
- 
\frac{1}{2} \vec{G}_{XX} \diamond \angrmd \wavearray, 
	\label{E:ZETATAN} 
	\\
\zetatrans 
& 
\eqdef 
- \frac{1}{2}  \angG_{\Lunit} \diamond \muX \wavearray. 
\label{E:ZETATRANS}
\end{align}
\end{subequations}
\end{lemma}

\begin{proof}
The same proofs of \cite{jLjS2018}*{Lemmas 2.13, 2.15} holds with minor modifications accounting for the third spatial dimension.
\end{proof}

\begin{lemma}[\cite{jS2016b}*{Proposition~5.4}; Frame decomposition of $\upmu \square_{\gfour(\wavearray)} f$] 
\label{L:FRAMEDCOMPOSITOINOFMUBOXG}
For scalar functions $f$, we have:
\begin{subequations}
\begin{align} \label{E:BOXDECOMPLOUTSIDE} 
\upmu \square_{\gfour(\wavearray)} f 
& = - \Lunit (\upmu \Lunit f + 2 \muX f) 
	+ 
	\upmu \angLap f 
	- 
	(\mytr_{\gtorus} \upchi) \muX f
	- 
	\upmu \mytr_{\gtorus} \angk \Lunit f 
	- 
	2 \upmu \upzeta^{\#} \cdot \angrmd f,
	\\
\upmu \square_{\gfour(\wavearray)} f 
& = 
- (\upmu \Lunit + 2 \muX)(\Lunit f) 
+ 
\upmu \angLap f 
- 
(\mytr_{\gtorus} \upchi)\muX f 
- 
(\Lunit \upmu) \Lunit f 
\label{E:BOXDECOMPMULBAROUTSIDE} 
- 
\upmu \mytr_{\gtorus} \angk \Lunit f 
+ 
2 \upmu \upzeta^{\#} \cdot \angrmd f 
+ 
2 (\angrmd^{\#} \upmu) \cdot \angrmd f. 
\end{align}
\end{subequations}
\end{lemma}

\section{Classical existence in a region of spacetime containing the crease} \label{S:CLASSICALEXISTENCEUPTOTHECREASE}
In this section, we adapt the necessary results concerning the $3D$ compressible Euler equations from \cite{abbrescia2022emergence} into the current setting. 
Here we only provide the results from \cite{abbrescia2022emergence} that are crucially used in the current paper. We invite the reader to read \cite{abbrescia2022emergence} for more details. We note that in that setting, we used a continuous family of time functions $\timefunctionarg{\mathfrak{n}}$ for $\mathfrak{n} \in [0,\mathfrak{n}_0]$ in order to construct a large portion of the singular boundary, see \cite{abbrescia2022emergence}. In the current work, we are only interested in the crease, which is embedded in the $0$-level set of $\timefunctionarg{0}$. Hence, we use the simplified notation $\timefunction \eqdef \timefunctionarg{0}$. 

\subsection{Basic constructions} 
\label{SS:BASICCONSTRUCTIONSFORROUGHTIMEFUNCTION}
We now introduce some basic ingredients that were used to construct $\timefunction$.

\subsubsection{$\upmu$-adapted subsets}
\label{SSS:MUADAPATEDSUBSETS}
To follow the solution up to the singular boundary, we analyzed it
on the following subsets (and others as well), 
which are adapted to the shape of the singular boundary.

\begin{definition}[Level sets of $\upmu$ and $\muX \upmu$ and the $\upmu$-adapted tori $\twoargmumuxtorus{\mulevelsetvalue}{0}$] 
	\label{D:LEVELSETSOFMUANDXMUANDMUMUXTORI}
	Recall that $\muX$ is defined in \eqref{E:MUX}.
	Given a real number $\mulevelsetvalue\ge 0$, we define:
	\begin{subequations}
		\begin{align} \label{E:LEVELSETSOFMU}
			\mulevelsetarg{\mulevelsetvalue}
			& \eqdef
			\left\lbrace
			(t,u,x^2,x^3) 
			\in \R \times \R \times \T^2 
			\ | \
			\upmu(t,u,x^2,x^3) = \mulevelsetvalue
			\right\rbrace
			\cap
			\lbrace |u| \leq \interestingu \rbrace,
			\\
			\datahypfortimefunctionarg{0}
			& \eqdef
			\left\lbrace
			(t,u,x^2,x^3) 
			\in \R \times \R \times \T^2 
			\ | \
			\muX \upmu(t,u,x^2,x^3) = 0
			\right\rbrace
			\cap
			\lbrace |u| \leq \interestingu \rbrace,
			\label{E:LEVELSETSOFMUXMU} \\
			\twoargmumuxtorus{\mulevelsetvalue}{0}
			& \eqdef 
			\mulevelsetarg{\mulevelsetvalue} 
			\cap
			\datahypfortimefunctionarg{0}.
			\label{E:MUXMUTORI}
		\end{align}
	\end{subequations}
\end{definition}

\begin{remark}[Differential structure with respect to the geometric coordinates] 
	\label{R:DIFFERENTIALSTRUCTUREWITHRESPECTTOGEOMETRICCOORDINATES} 
	In the rest of the paper, unless we explicitly mention otherwise, 
	we implicitly consider $\mulevelsetarg{\mulevelsetvalue}$, $\datahypfortimefunctionarg{0}$,
	and $\twoargmumuxtorus{\mulevelsetvalue}{0}$
	to be subsets of spacetime with the differential structure 
	\emph{induced by the geometric coordinates $(t,u,x^2,x^3)$};
	this is already apparent from Def.\,\ref{D:LEVELSETSOFMUANDXMUANDMUMUXTORI}.
	Similar remarks apply to the 
	sets
	$\hypthreearg{\timefunction}{I}{0}$,
	$\twoargroughtori{\timefunction,u}{0}$,
	$\nullhypthreearg{0}{u}{I}$,
	$\twoargMrough{I,J}{0}$,
	$\mulevelsettwoarg{\mulevelsetvalue}{I}$,
	and $\datahypfortimefunctiontwoarg{0}{I}$
	defined below. 
\end{remark}

\subsection{Rough time functions, the parameter $\timefunction_0$, and rough adapted coordinates}
\label{SS:ROUGHTIMEFUNCTIONANDROUGHADAPATEDCOORDINATES}
We now provide the transport equation initial value problem whose solution is the rough time function.

\subsubsection{Rough time functions and the parameter $\timefunction_0$}
\label{SSS:ROUGHTIMEFUNCTIONSANDINITIALROUGHSLICE}

\begin{definition}[The rough time function $\timefunction$]  
	\label{D:ROUGHTIMEFUNCTION}
	Let $\mupositive = -\timefunction_0 > 0$ be the real number from \cite{abbrescia2022emergence}*{Theorem~2.8} and let $\mulevelsetvalue \in [0,\mupositive]$
	(see Sect.\,\ref{S:PARAMETERSANDSIZEASSUMPTIONSANDCONVENTIONSFORCONSTANTS} for further discussion of how 
	the specific choice of $\mupositive$ that we make in our main theorem is tied to the initial data).
	Let $\muX$ be the vectorfield defined in \eqref{E:MUX},
	and let $\twoargmumuxtorus{\mulevelsetvalue}{0}$ 
	be the $\upmu$-adapted torus defined by \eqref{E:MUXMUTORI}.
	For $\mulevelsetvalue \in [0,\mupositive]$, 
	we define the \textbf{rough time function} $\timefunction$ 
	to be the solution to the following transport equation 
	initial value problem:
	\begin{subequations} \label{E:IVPFORROUGHTTIMEFUNCTION}
		\begin{align} 	
			\muX \timefunction 
			& = 0,  
			\label{E:TRANSPORTEQUATIONFORROUGHTIMEFUNCTION}
			\\
			\timefunction|_{\twoargmumuxtorus{\mulevelsetvalue}{0}} 
			& =  
			- \mulevelsetvalue 
			= 
			- \upmu|_{\twoargmumuxtorus{\mulevelsetvalue}{0}}. 
			\label{E:INITIALCONDITIONFORROUGHTIMEFUNCTION}
		\end{align} 
	\end{subequations}
\end{definition}

We sometimes refer to $\twoargmumuxtorus{\mulevelsetvalue}{0}$ 
as a \emph{primal torus} for $\timefunction$ because 
$\timefunction$  is ``flowed out from'' it.

\begin{remark}[$\muX$ is tangent to the level sets of $\timefunction$]
	\label{R:WTRANSISTANGENTTOLEVELSETSOFROUGHTIMEFUNCTION}
	Note that equation \eqref{E:TRANSPORTEQUATIONFORROUGHTIMEFUNCTION} implies that
	$\muX$ is tangent to the level sets of
	$\timefunction$.
\end{remark}

\begin{definition}[The parameter $\timefunction_0$]
	\label{D:RELATIONBETWEENINITIALROUGHTIMESLICEANDSMALLMUVALUE}
	We define the parameter $\timefunction_0 < 0$ as follows:
	\begin{align} \label{E:RELATIONBETWEENINITIALROUGHTIMESLICEANDSMALLMUVALUE}
		\timefunction_0 
		& \eqdef 
		- 
		\mupositive.
	\end{align}
\end{definition}

In \cite{abbrescia2022emergence}, portions of the rough hypersurface
$\lbrace \timefunction = \timefunction_0 \rbrace$
played the role of an ``initial data'' hypersurface near 
the singularity.
Note that by construction, we have $\timefunction_0 \leq \timefunction \leq 0$.

\subsubsection{Rough adapted coordinates}
\label{SSS:ROUGHADAPTEDSUBSETS}
Having constructed the eikonal function $u$ and the
rough time function $\timefunction$, 
we now define a system 
of coordinates adapted to them.

\begin{definition}[The rough adapted coordinates and their partial derivative vectorfields] 
	\label{D:ROUGHCOORDINATESANDPARTIALDERIVATIVES}
	We call $(\timefunction,u,x^2,x^3)$ the \emph{rough adapted coordinates}. We 
	denote the corresponding rough adapted coordinate partial derivative vectorfields by 
	$\left\lbrace\roughgeop{\timefunction},\roughgeop{u},\roughgeop{x^2},\roughgeop{x^3} \right\rbrace$. 
\end{definition}

\subsection{Rough subsets}
\label{SS:ROUGHSUBSETS}
In this section, we define various subsets of spacetime that are tied to
$u$ and $\timefunction$. Most of the delicate PDE analysis in the bulk
of \cite{abbrescia2022emergence} took place on these subsets.

\subsubsection{Truncated $\timefunction$-adapted subsets}
\label{SSS:TRUNCATEDROUGHSUBSETS}

\begin{definition}[Truncated $\timefunction$-adapted subsets]  \label{D:TRUNCATEDROUGHSUBSETS}
	Given intervals $I, J \subset \mathbb{R}$ and real numbers 
	$\timefunction, u \in \mathbb{R}$, we define:
	\begin{subequations}
		\begin{align} \label{E:ROUGHHYPERSURFACES}
			\hypthreearg{\timefunction'}{J}{0} 
			& \eqdef \{ (t,u,x^2,x^3) \in \R\times\R\times\T^2 \ | \ \timefunction(t,u,x^2,x^3) 
			= 
			\timefunction', \ u \in J \},
			\\
			\twoargroughtori{\timefunction',u}{0}
			& 
			\eqdef 
			\{ (t,u,x^2,x^3) \ | \ 
			(t,x^2,x^3) \in \mathbb{R} \times \mathbb{T}^2,
			\,
			\timefunction(t,u,x^2,x^3) 
			= 
			\timefunction' 
			\},
			\label{E:ROUGHTORI}
			\\
			\nullhypthreearg{0}{u}{I} 
			& 
			\eqdef 
			\bigcup_{\timefunction' \in I} \twoargroughtori{\timefunction',u}{0},
			\label{E:NULLHYPERSURFACEROUGHTRUNCATED} 
			\\
			\twoargMrough{I,J}{0} 
			& 
			\eqdef \bigcup_{\timefunction' \in I} \hypthreearg{\timefunction'}{J}{0} 
			= 
			\bigcup_{u' \in J} \nullhypthreearg{0}{u'}{I}.
			\label{E:TRUNCATEDMROUGH}
		\end{align}
	\end{subequations}
\end{definition}

We refer to the $\hypthreearg{\timefunction}{J}{0}$ as
\emph{rough hypersurfaces}. 
We sometimes refer to $\hypthreearg{\timefunction_0}{J}{0}$ as the \emph{initial rough hypersurface},
where $\timefunction_0 < 0$ is the parameter from Sect.\,\ref{SSS:ROUGHTIMEFUNCTIONSANDINITIALROUGHSLICE}.
We refer to the $\twoargroughtori{\timefunction,u}{0}$ as
\emph{rough tori}.
We also note that $\nullhypthreearg{0}{u}{I}$ is a portion of the $\gfour$-null surface $\nullhyparg{u}$.

From Defs.\,\ref{D:ROUGHCOORDINATESANDPARTIALDERIVATIVES} and\,\ref{D:TRUNCATEDROUGHSUBSETS},
it follows that $\left\lbrace\roughgeop{u},\roughgeop{x^2},\roughgeop{x^3} \right\rbrace$
spans the tangent space of $\hypthreearg{\timefunction}{I}{0}$,
that $\left\lbrace\roughgeop{x^2},\roughgeop{x^3} \right\rbrace$
spans the tangent space of $\twoargroughtori{\timefunction,u}{0}$,
that 
$\left\lbrace\roughgeop{t},\roughgeop{x^2},\roughgeop{x^3} \right\rbrace$
spans the tangent space of $\nullhypthreearg{0}{u}{I}$,
and that $\left\lbrace\roughgeop{\timefunction}, \roughgeop{u},\roughgeop{x^2},\roughgeop{x^3} \right\rbrace$
spans the tangent space of $\twoargMrough{I,J}{0}$.

\subsubsection{Truncated $\upmu$-adapted subsets}
\label{SSS:TRUNCATEDMUADAPATEDSUBSETS}

\begin{definition}[Truncated level sets of $\upmu$ and $\muX \upmu$] 
	\label{D:TRUNCATEDLEVELSETSOFMUANDXMUANDMUMUXTORI}
	Let $I \subset [\timefunction_0,0]$ be an interval,
	let $\mulevelsetvalue \in [0,\upmu_0$], and let $\mulevelsetarg{\mulevelsetvalue}$
	and $\datahypfortimefunctionarg{0}$ be the sets from Def.\,\ref{D:LEVELSETSOFMUANDXMUANDMUMUXTORI}.
	We define:
	\begin{subequations}
		\begin{align} \label{E:TRUNCATEDLEVELSETSOFMU}
			\mulevelsettwoarg{\mulevelsetvalue}{I}
			& 
			\eqdef
			\mulevelsetarg{\mulevelsetvalue}
			\cap
			\left\lbrace (t,u,x^2,x^3) \in \R \times \R \times \T^2 
			\ | \ 
			\timefunction(t,u,x^2,x^3) \in I 
			\right\rbrace,
			\\
			\datahypfortimefunctiontwoarg{0}{I}
			& \eqdef
			\datahypfortimefunctionarg{0}
			\cap
			\left\lbrace (t,u,x^2,x^3) \in \R \times \R \times \T^2 
			\ | \ 
			\timefunction(t,u,x^2,x^3) \in I 
			\right\rbrace.
			\label{E:TRUNCATEDLEVELSETSOFMUXMU} 
		\end{align}
	\end{subequations}
\end{definition}

Just as in Remark\,\ref{R:DIFFERENTIALSTRUCTUREWITHRESPECTTOGEOMETRICCOORDINATES}, 
we view 
$\hypthreearg{\timefunction}{J}{0}$,
$\nullhypthreearg{0}{u}{I}$,
$\twoargMrough{I,J}{0}$, 
$\mulevelsettwoarg{\mulevelsetvalue}{I}$,
and
$\datahypfortimefunctiontwoarg{0}{I}$
as submanifolds of spacetime equipped with the differential structure induced by the geometric coordinates $(t,u,x^2,x^3)$.

\subsection{Change of variables maps}
\label{SS:ALLTHECHOVMAPSFORSINGULARBOUNDARY}
In this short section, we define various change of variables maps that we use to prove our main results in \cite{abbrescia2022emergence}.

\begin{definition}[Change of variables maps]
	\label{D:ALLTHECHOVMAPS}
	We define the change of variables map from geometric to Cartesian coordinates as follows:
	\begin{align} \label{E:CHOVGEOTOCARTESIAN}
		\Upsilon(t,u,x^2,x^3) 
		& 
		\eqdef (t,x^1,x^2,x^3).
	\end{align}
	
	We define the change of variables map from geometric coordinates 
	to rough adapted coordinates as follows:
	\begin{align} \label{E:CHOVGEOTOROUGH}
		\CHOVgeotorough{0}(t,u,x^2,x^3) 
		& 
		\eqdef (\timefunctionarg{0},u,x^2,x^3).
	\end{align}
	
	We define the map $\CHOVgeotomumuxmu$ 
	from geometric coordinates to ``$(\upmu,\muX \upmu,x^2,x^3)$-space''
	and its Jacobian 
	$\CHOVJacobiangeotomumuxmu$ as follows:
	\begin{subequations}
		\begin{align} \label{E:CHOVFROMGEOMETRICCOORDINATESTOMUWEGIGHTEDXMUCOORDINATES}
			\CHOVgeotomumuxmu(t,u,x^2,x^3)
			& \eqdef (\upmu,\muX \upmu,x^2,x^3),
			\\
			\CHOVJacobiangeotomumuxmu(t,u,x^2,x^3)
			& \eqdef 
			\frac{\partial \CHOVgeotomumuxmu(t,u,x^2,x^3)}{\partial(t,u,x^2,x^3)}
			=
			\frac{\partial (\upmu,\muX \upmu,x^2,x^3)}{\partial(t,u,x^2,x^3)}.
			\label{E:JACOBIANMATRIXFORCHOVFROMGEOMETRICCOORDINATESTOMUWEGIGHTEDXMUCOORDINATES}
		\end{align}
	\end{subequations}
	
	We define the map $\CHOVroughtomumuxmu{0}$ 
	from rough adapted coordinates to ``$(\upmu,\muX \upmu,x^2,x^3)$-space''
	and its Jacobian 
	$\CHOVJacobianroughtomumuxmu{0}$ as follows:
	\begin{subequations}
		\begin{align} \label{E:CHOVFROMROUGHCOORDINATESTOMUWEGIGHTEDXMUCOORDINATES}
			\CHOVroughtomumuxmu{0}(\timefunctionarg{0},u,x^2,x^3)
			& \eqdef (\upmu,\muX \upmu,x^2,x^3),
			\\
			\CHOVJacobianroughtomumuxmu{0}(\timefunctionarg{0},u,x^2,x^3)
			& \eqdef 
			\frac{\partial \CHOVroughtomumuxmu{0}(\timefunction,u,x^2,x^3)}{\partial(\timefunctionarg{0},u,x^2,x^3)}
			=
			\frac{\partial (\upmu,\muX \upmu,x^2,x^3)}{\partial(\timefunctionarg{0},u,x^2,x^3)}.
			\label{E:JACOBIANMATRIXFORCHOVFROMROUGHCOORDINATESTOMUWEGIGHTEDXMUCOORDINATES}
		\end{align}
	\end{subequations}
	
\end{definition}
We note the following identity: 
\begin{align} \label{E:SIMPLERELATIONSHIPBETWEENCHOVMAPS}
	\CHOVgeotomumuxmu 
	& = \CHOVroughtomumuxmu{0} \circ \CHOVgeotorough{0}.
\end{align}\subsection{Existence in a region containing the crease $\twoargmumuxtorus{0}{0}$} \label{SS:ESTIMATESFORROUGHTIMEFUNCTION}
We begin with the following theorem, which is a summarized and adapted version of some of the main results of \cite{abbrescia2022emergence}. Importantly, we record here some of the geo-analytic structures of the rough subsets of spacetime.

\begin{theorem}[Existence in a region containing the singular boundary torus $\twoargmumuxtorus{0}{0}$] \label{T:MAINRESULTSFROMSINGULARBOUNDARYPAPER}
	Assume that $(\RRiemann,\LRiemann,v^2,v^3,\Ent)|_{t=0} \in C^\infty(\Sigma_0)^{\times 5}$ is the same smooth initial datas as in \cite{abbrescia2022emergence}. In particular, the same transversal convexity and Sobolev inequality as \cite{abbrescia2022emergence} holds. Then, the following hold:
	
	\noindent \underline{\textbf{The rough time function and classical existence relative to the geometric coordinates}}.
	
	\begin{itemize}
		\item There exists a rough time function 
		$\timefunctionarg{0} = \timefunctionarg{0}(t,u,x^2,x^3)$ 
		with range $[\timefunction_0,0] = [-\mupositive,0]$ (recall that $\timefunction_0 = - \mupositive$),
		and we denote its level set portions,
		viewed as subsets of geometric coordinate space $\mathbb{R} \times \mathbb{R} \times \mathbb{T}^2$,
		as follows:
		$\hypthreearg{\timefunction}{[-\rightu,\leftu]}{0}
		= 
		\lbrace 
		(t,u,x^2,x^3) 
		\ | \
		\timefunctionarg{0}(t,u,x^2,x^3) = \timefunction,
		\,
		u_1 \leq u \leq u_2,
		\,
		(x^2,x^3) \in \mathbb{T}^2
		\rbrace
		$.
		More precisely,
		$\timefunctionarg{0}$ is
		defined on the portion 
		$\twoargMrough{[\timefunction_0,0],[- \rightu,\leftu]}{0}
		= \bigcup_{\timefunction \in [\timefunction_0,0]} \hypthreearg{\timefunction}{[- \rightu,\leftu]}{0}
		$
		of the maximal classical development 
		of the data with respect to the differential structure of the geometric coordinates $(t,u,x^2,x^3)$.
		\item The change of variables map 
		$\CHOVgeotorough{0}(t,u,x^2,x^3) = (\timefunctionarg{0},u,x^2,x^3)$
		is a diffeomorphism from $\twoargMrough{[\timefunction_0,0],[- \rightu,\leftu]}{0}$
		onto its image $[-\mupositive,0] \times [- \rightu,\leftu] \times \mathbb{T}^2$ 
		satisfying: 
		\begin{align} 
		\| \CHOVgeotorough{0} \|_{C_{\textnormal{geo}}^{2,1}(\twoargMrough{[\timefunction_0,0],[- \rightu,\leftu]}{0})} & \leq C, \label{E:C21ESTIMATESOFTHEROUGHTIMEFUNCTIONINGEOMETRICCOORDINATES} \\
		\left\| \geop{x^A} \timefunctionarg{0} \right\|_{C_{\textnormal{geo}}^{1,1}(\twoargMrough{[\timefunction_0,0],[- \rightu,\leftu]}{0})} & \lesssim \initialsmall,\label{E:GEOMETRICPARTIALANGULARDERIVATIVEOFTIMEFUNCTIONARESMALL} 
		\end{align} 
		Moreover, 
		\begin{align}
			\geop{t} \timefunctionarg{0} & \approx 1, &\text{on $\twoargMrough{[\timefunction_0,0],[- \rightu,\leftu]}{0}$}.\label{E:GEOMETRICPARTIALTDERIVATIVEOFTIMEFUNCTIONISUNITSIZED}
		\end{align}
		\item The wave variables
		$\wavearray$, the transport variables $\{\vortrenormalized,\GradEnt,\VortVort,\DivGradEnt\}$,
		the eikonal function $u$,
		$\upmu$, $\Lunit^i$,
		the Cartesian coordinate functions
		$(t,x^1,x^2,x^3)$ 
		and all of the auxiliary quantities constructed out of these quantities
		exist classically
		with respect to the geometric coordinates $(t,u,x^2,x^3)$
		on $\twoargMrough{[\timefunction_0,0],[- \rightu,\leftu]}{0}$.
		\item The following quantities extend as solutions to the compact set 
	$\twoargMrough{[\timefunction_0,0],[- \rightu,\leftu]}{0}$
	as elements of the following spacetime H\"{o}lder spaces
	with respect to the geometric coordinates,
	and their corresponding spacetime H\"{o}lder norms on $\twoargMrough{[\timefunction_0,0],[- \rightu,\leftu]}{0}$
	are bounded by $\leq C$:
	\begin{itemize}
		\item $\wavearray, \, \vortrenormalized^i, \, \GradEnt^i, \, \VortVort^i, \, \DivGradEnt 
 \in C_{\textnormal{geo}}^{3,1}(\twoargMrough{[\timefunction_0,0],[- \rightu,\leftu]}{0})$,
		\item $\Upsilon \in C_{\textnormal{geo}}^{3,1}(\twoargMrough{[\timefunction_0,0],[- \rightu,\leftu]}{0})$,
		\item $\Lunit^i, \, \upmu \in C_{\textnormal{geo}}^{2,1}(\twoargMrough{[\timefunction_0,0],[- \rightu,\leftu]}{0})$.
	\end{itemize} 
		\item The following quantities are small:
			\begin{align}
				\left\| \geop{x^A} \wavearray \right\|_{C_{\textnormal{geo}}^{2,1}(\twoargMrough{[\timefunction_0,0],[- \rightu,\leftu]}{0})} & \lesssim \initialsmall, \label{E:WAVEANDTRANSPORTVARIABLESARESMALLINCLASSICALDEVELOPMENTFROMROUGHCOORDINATES} \\
				\left\| \geop{x^A} \Lunit^i, \,  \geop{x^A}  \upmu ^i \right\|_{C_{\textnormal{geo}}^{1,1}(\twoargMrough{[\timefunction_0,0],[- \rightu,\leftu]}{0})} & \lesssim \initialsmall. \label{E:ACOUSTICVARIABLESARESMALLINCLASSICALDEVELOPMENTFROMROUGHCOORDINATES} 
			\end{align}
		\item For $\timefunction \in [\timefunction_0,0]$,
		there exists a function 
		$\Cartesiantisafunctiononlevelsetsofroughtimefunctionarg{\timefunction}{0}: 
		[- \rightu,\leftu] \times \mathbb{T}^2 \rightarrow \mathbb{R}$,
		depending on $\timefunction$ and $0$,
		such that:
		\begin{align} \label{E:C21BOUNDFORCONSTANTTIMEFUNCTIONGRAPH}
			\| \Cartesiantisafunctiononlevelsetsofroughtimefunctionarg{\timefunction}{0} \|_{C^{2,1}([- \rightu,\leftu] \times \mathbb{T}^2)}
			& \leq C
		\end{align}
		and such that relative to the geometric coordinates, we have:
		\begin{align} \label{E:LEVELSETSOFTIMEFUNCTIONAREAGRAPH}
			\hypthreearg{\timefunction}{[- \rightu,\leftu]}{0}
			& =
			\left\lbrace
			(t,u,x^2,x^3)
			\ | \
			t = \Cartesiantisafunctiononlevelsetsofroughtimefunctionarg{\timefunction}{0}(u,x^2,x^3),
			\,
			(u,x^2,x^3) \in [- \rightu,\leftu] \times \mathbb{T}^2
			\right\rbrace.
		\end{align}
		\item 
		The following estimates hold for $\timefunction \in [\timefunction_0,0]$:
	\begin{subequations}
\begin{align} \label{E:SIZEOFCARTESIANT}
		\frac{1}{3 \mathring{\updelta}_*}
		& 
		\leq
		\min_{\hypthreearg{\timefunction}{[- \rightu,\leftu]}{0}} t
		\leq
		\sup_{\hypthreearg{\timefunction}{[- \rightu,\leftu]}{0}} t
		\leq 
		\frac{3}{\mathring{\updelta}_*},
			\\
		- \leftu
		+
		\frac{1}{3 \mathring{\updelta}_*}
		& 
		\leq
		\min_{\hypthreearg{\timefunction}{[- \rightu,\leftu]}{0}} x^1
		\leq
		\sup_{\hypthreearg{\timefunction}{[- \rightu,\leftu]}{0}} x^1
		\leq 
		\rightu
		+
		\frac{3}{\mathring{\updelta}_*}.
		\label{E:SIZEOFCARTESIANX1}
	\end{align}
	\end{subequations}	\end{itemize}		
	
	\medskip
	
	\noindent \underline{\textbf{Structure of the rough adapted tori $\twoargmumuxtorus{-\timefunction}{0}$ and properties of the flow map of $\muX$}}.
	
	\begin{itemize}
		\item There exists scalar functions $\Cartesiantisafunctiononmumxtoriarg{\mulevelsetvalue}{0}, \, \Eikonalisafunctiononmumuxtoriarg{\mulevelsetvalue}{0} \colon \T^2 \to \R$ such that the map: 	
		\begin{equation} \label{E:CARTESIANTANDEIKONALFUNCTIONAREGRAPHSONMUMUXTORI}
		\embeddatahypersurfacearg{0}(\mulevelsetvalue,x^2,x^3)
		=
		\left(\Cartesiantisafunctiononmumxtoriarg{\mulevelsetvalue}{0}(x^2,x^3),\Eikonalisafunctiononmumuxtoriarg{\mulevelsetvalue}{0}(x^2,x^3),x^2,x^3 \right)
		\in
		\twoargmumuxtorus{\mulevelsetvalue}{0}
	\end{equation}
		 extends to a $C^{1,1}$ embedding
		from $[0,\mupositive] \times \mathbb{T}^2$ 
		onto its image, which is $\datahypfortimefunctiontwoarg{0}{[\timefunction_0,0]}$.
		In addition, there is a $C > 1$ such that the extended embedding satisfies:
		\begin{align} \label{E:EMBEDDINGOFDATAHYPERSURFACEC11EXTENDEDC11BOUND}
			\| \embeddatahypersurfacearg{0} \|_{C^{1,1}([0,\mupositive] \times \mathbb{T}^2)} 
			& \leq C, 
		\end{align}
		\begin{align}
			\left\| \left(\geop{x^2} \Cartesiantisafunctiononmumxtoriarg{\mulevelsetvalue}{0},
			\geop{x^3} \Cartesiantisafunctiononmumxtoriarg{\mulevelsetvalue}{0} \right) \right\|_{C^{0,1}(\mathbb{T}^2)},
			\,
		\left\| \left(\geop{x^2} \Eikonalisafunctiononmumuxtoriarg{\mulevelsetvalue}{0},
			\geop{x^3} \Eikonalisafunctiononmumuxtoriarg{\mulevelsetvalue}{0} \right) \right\|_{C^{0,1}(\mathbb{T}^2)}
		& \leq C \initialsmall,
		\label{E:C01SMALLNESSOFANGULARDERIVATIVESOFTANDUGRAPHS}
		\end{align}
		and: 
		\begin{align} \label{E:CLOSEDIMPROVEMENTMUEQUALSMINUSKAPPAEMBEDDINGCARTESIANTIMEFUNCTIONNEGATIVEMUDERIVATIVE}
			- 
			C
			& 
			<
			\min_{(\mulevelsetvalue,x^2,x^3) \in [0,\mupositive] \times \mathbb{T}^2} 
			\frac{\partial}{\partial \mulevelsetvalue} \Cartesiantisafunctiononmumxtoriarg{\mulevelsetvalue}{0}(x^2,x^3) 
			\leq
			\max_{(\mulevelsetvalue,x^2,x^3) \in [0,\mupositive] \times \mathbb{T}^2} 
			\frac{\partial}{\partial \mulevelsetvalue} \Cartesiantisafunctiononmumxtoriarg{\mulevelsetvalue}{0}(x^2,x^3) 
			< 
			- \frac{1}{C}.
		\end{align}
		
		Furthermore, for $\mulevelsetvalue \in [0,\mupositive]$,
		we have: 
		\begin{align} \label{E:GRAPHDESCRIPTIONOFMUXMUTORUS}
			\twoargmumuxtorus{\mulevelsetvalue}{0}
			& = 
			\left\lbrace
			\left(\Cartesiantisafunctiononmumxtoriarg{\mulevelsetvalue}{0}(x^2,x^3),\Eikonalisafunctiononmumuxtoriarg{\mulevelsetvalue}{0}(x^2,x^3),x^2,x^3 \right)
			\ | \
			(x^2,x^3) \in \mathbb{T}^2
			\right\rbrace
		\end{align}
		and:
		\begin{align}
			\datahypfortimefunctiontwoarg{0}{[\timefunction_0,0]}
			& 
			=
			\bigcup_{\mulevelsetvalue \in [0,\mupositive]}
			\twoargmumuxtorus{\mulevelsetvalue}{0}.
			\label{E:CLOSEDIMPROVEMENTLEVELSETSTRUCTUREOFMUXEQUALSMINUSKAPPA}
		\end{align}
		\item 	Let $\muX$ be the vectorfield defined in \eqref{E:MUX}, and recall that
		$\muX \timefunction = 0$
		and
		$\muX u = 1$.
		Let $(\Delta u,t,u,x^2,x^3) \rightarrow \flowmapWtransargtwoarg{0}{\Delta u}(t,u,x^2,x^3)$ 
		denote the flow map of $\muX$,
		i.e., at each fixed 
		$
		(t,u,x^2,x^3) \in \twoargMrough{[\timefunction_0,0],[- \rightu,\leftu]}{0}
		$,
		the components of $\flowmapWtransargtwoarg{0}{\Delta u}(t,u,x^2,x^3)$ 
		solve the following ODE system initial value problem
		on the flow interval $\Delta u \in [- \rightu - u,\leftu-u]$:
		\begin{align} \label{E:FLOWMAPFORGENERATOROFROUGHTIMEFUNCTION}
			\frac{\partial}{\partial \Delta u} \flowmapWtransargtwoarg{0}{\Delta u}(t,u,x^2,x^3) 
			& = \muX \circ \flowmapWtransargtwoarg{0}{\Delta u}(t,u,x^2,x^3),
			&
			\flowmapWtransargtwoarg{0}{0}(t,u,x^2,x^3)
			& = 
			(t,u,x^2,x^3).
		\end{align}
		
		Then for each fixed $\timefunction \in [\timefunction_0,0]$
		and each pair $u_1, u_2 \in [- \rightu, \leftu]$,
		$\flowmapWtransargtwoarg{0}{u_2 - u_1}$ is a diffeomorphism from 
		the rough torus $\twoargroughtori{\timefunction,u_1}{0}$ onto the rough torus $\twoargroughtori{\timefunction,u_2}{0}$.
		In particular, the integral curves of $\muX$ thread $\hypthreearg{\timefunction}{[- \rightu,\leftu]}{0}$.
		Moreover, for every fixed $\timefunction \in [\timefunction_0,0]$,
		each integral curve of $\muX$ passes through precisely one point on
		the $\upmu$-adapted torus
		$
		\twoargmumuxtorus{-\timefunction}{0}$
		defined in \eqref{E:MUXMUTORI}.
		
		Moreover, with $d_{\textnormal{geo}}$ denoting the differential with respect to the geometric coordinates,
	we have the following bounds, where the implicit constants in \eqref{E:BOUNDFORDETERMINANTWFLOWMAP}
	are independent of all $\Delta u$ such that $|\Delta u| \leq \rightu + \leftu$:
	\begin{align} \label{E:C21FORWFLOWMAP} 
		\sup_{|\Delta u| \leq \rightu + \leftu}
		\| d_{\textnormal{geo}} \flowmapWtransargtwoarg{0}{\Delta u} \|_{C^{1,1}
		(
		\twoargMrough{[\timefunction_0,0],[- \rightu,\leftu]}{0}
		\cap
		\twoargMrough{[\timefunction_0,0],[- \rightu - \Delta u,\leftu - \Delta u]}{0})
		}
		& \leq C,
	\end{align}
	\begin{align} \label{E:BOUNDFORDETERMINANTWFLOWMAP}
		\det \left(d_{\textnormal{geo}} \flowmapWtransargtwoarg{0}{\Delta u} \right)
		&
		\approx 1
		\mbox{ on }
		\twoargMrough{[\timefunction_0,0],[- \rightu,\leftu]}{0}
		\cap
		\twoargMrough{[\timefunction_0,0],[- \rightu - \Delta u,\leftu - \Delta u]}{0}.
	\end{align}	\end{itemize}
	\medskip
	
	\noindent \underline{\textbf{The behavior of $\upmu$ and properties of $\Upsilon$}}.
	
	\begin{itemize}
		\item For $\timefunction \in [-\mupositive,0]$,
		we have: 
		\begin{align} \label{E:MAINRESULTSMINVALUEOFMUONROUGHFOLIATION}
			\min_{\hypthreearg{\timefunction}{[- \rightu,\leftu]}{0}} \upmu
			& = - \timefunction.
		\end{align}
		Moreover, within $\hypthreearg{\timefunction}{[- \rightu,\leftu]}{0}$,
		the minimum value $- \timefunction$
		in \eqref{E:MAINRESULTSMINVALUEOFMUONFOLIATION} is achieved by $\upmu$ precisely on the
		set $\twoargmumuxtorus{-\timefunction}{0}$ from definition\,\ref{E:MUXMUTORI},
		which is a $C^{1,1}$-embedded torus.
		In particular, in $\twoargMrough{[\timefunction_0,0],[- \rightu,\leftu]}{0}$,
		$\upmu$ vanishes precisely along the torus $\twoargmumuxtorus{0}{0}$,
		which is a subset of $\hypthreearg{0}{[- \rightu,\leftu]}{0}$
		that is contained in the singular boundary.
		\item  On $\twoargMrough{[\timefunction_0,0],[- \rightu,\leftu]}{0}$,
		the change of variables map 
		$\Upsilon(t,u,x^2,x^3) = (t,x^1,x^2,x^3)$ is an injection onto its image in Cartesian coordinate space
		satisfying $\| \Upsilon \|_{C_{\textnormal{geo}}^{3,1}(\twoargMrough{[\timefunction_0,0],[- \rightu,\leftu]}{0})} \leq C$.
		In particular, $\Upsilon$ is a homeomorphism from the compact set
		$\twoargMrough{[\timefunction_0,0],[- \rightu,\leftu]}{0}$ onto its image.
		\item With $d_{\textnormal{geo}} \Upsilon$ denoting the Jacobian matrix of $\Upsilon$,
		we have:
		\begin{align} \label{E:FIXEDKAPPAMAINTHEOREMGEOTOCARTESIANJACOBIANDETERMINANTESTIMATE}
			\mbox{\upshape det} d_{\textnormal{geo}} \Upsilon
			& \approx
			- 
			\upmu.
		\end{align}
		Hence, on $\twoargMrough{[\timefunction_0,0],[- \rightu,\leftu]}{0} \backslash \twoargmumuxtorus{0}{0}$,
		$\Upsilon$ is a diffeomorphism.
		\item For $\mulevelsetvalue \in [0,\mupositive]$, 
		$\Upsilon\left(\twoargmumuxtorus{\mulevelsetvalue}{0} \right)$
		is an embedded two-dimensional $C^{1,1}$ torus in Cartesian coordinate space.
		In particular, the restriction of $\Upsilon$ to $\twoargmumuxtorus{\mulevelsetvalue}{0}$ 
		is a diffeomorphism from $\twoargmumuxtorus{\mulevelsetvalue}{0}$ onto its image
		$\Upsilon\left(\twoargmumuxtorus{\mulevelsetvalue}{0} \right)$.
	\end{itemize}
	
	\medskip
	
	\noindent \underline{\textbf{A description of the singular and regular behavior with respect to Cartesian coordinates}}.
	\begin{itemize}
		\item (Region without singularities).
		On the subset
		$\Upsilon\left(\twoargMrough{[\timefunction_0,0],[- \rightu,\leftu]}{0} \backslash \twoargmumuxtorus{0}{0} \right)$
		of Cartesian coordinate space,
		the solution exists classically with respect to the Cartesian coordinates.
		\item (The shock singularity).
		The following lower bound holds in 
		$\Upsilon\left(\twoargMrough{[\timefunction_0,0],[- \interestingu,\interestingu]}{0} \right)$: 
		\begin{align} \label{E:FIXEDKAPPAMAINTHEOREMBLOWUPLOWERBOUND}
			|X \RRiemann|
			&
			\geq
			\frac{\mathring{\updelta}_*}{\upmu |\bar{\Speed}_{;\LogDensity} + 1|\upmu},
		\end{align}
				where $\mathring{\updelta}_* > 0$ is the data-parameter:
\begin{align} \label{E:DELTASTARDEF}
	\mathring{\updelta}_* 
	& 
	\eqdef 
	\sup_{\hypthreearg{\timefunction_0}{[- \rightu,\leftu]}{\muxmulevelsetvalue}} 
	\frac{1}{2} 
	\left[ 
		(\Speed^{-1} \Speed;_{\LogDensity} + 1) \muX \RRiemann 
	\right]_+,
\end{align}
		$\bar{\Speed}_{;\LogDensity} \eqdef \Speed_{;\LogDensity}(\LogDensity = 0,\Ent=0)$
		is $\Speed_{;\LogDensity}$ evaluated at the trivial solution, 
		$\bar{\Speed}_{;\LogDensity} + 1$ is a \underline{non-zero} constant
, and the $\Sigma_t$-tangent vectorfield $X$ has Euclidean length satisfying
		$\sqrt{\sum_{a=1}^3 (X^a)^2} = 1 + \mathcal{O}(\mathring{\upalpha})$,
		where $\mathring{\upalpha}$ is the small parameter from \cite{abbrescia2022emergence}*{Sect.12.3}.
		In particular, if $q \in \Upsilon\left(\twoargmumuxtorus{0}{0} \right)$,
		then since 
		$\twoargmumuxtorus{0}{0} \in \twoargMrough{[\timefunction_0,0],[- \frac{\interestingu}{2},\frac{\interestingu}{2}]}{0}$ 
		by \eqref{E:IMPROVEDLEVELSETSTRUCTUREANDLOCATIONOFMIN}, 
		and since $\upmu = 0$ along 
		$\Upsilon\left(\twoargmumuxtorus{0}{0} \right)$,
		it follows that
		$|X \RRiemann|(q') \to \infty$ as $q' \rightarrow q$ in
		$\Upsilon\left(\twoargMrough{[\timefunction_0,0],[- \rightu,\leftu]}{0} 
		\backslash \twoargmumuxtorus{0}{0} \right)$.
		\item (Regular behavior along the characteristics).
		The derivatives of
		$\wavearray$,
		up to order $\Ntop - 6$
		with respect to the vectorfields in the 
		$\nullhyparg{u}$-tangent commutation set $\Tanset$ defined in \eqref{E:COMMUTATIONVECTORFIELDS}
		and the derivatives of
		with respect to the elements of $\Tanset$ are $L^{\infty}$-bounded 
		on $\Upsilon\left(\twoargMrough{[\timefunction_0,0],[- \rightu,\leftu]}{0} \right)$.
	\end{itemize}
\end{theorem}

In Sect.\,\ref{S:DOUBLENULLFOLLIATIONSETUP} we will define the ingoing eikonal function $\ubar$ as a solution to the eikonal equation $\gfour(\Dfour \ubar,\Dfour \ubar) = 0$ with data posed on a subset of $\datahypfortimefunctiontwoarg{0}{[\timefunction_0,0]}$, see Def.\,\ref{D:INGOINGEIKONALFUNCTION}. Hence, we record here a result from our work \cite{abbrescia2022emergence}*{Sect.\,15} that will be used in the analysis of $\ubar$ carried out in Sect.\,\ref{S:CONSTRUCTIONANDESTIMATESFORINGOINGEIKONALFUNCTION}.

\begin{proposition}[Properties of the surface $\datahypfortimefunctiontwoarg{0}{[\timefunction_0,0]}$] \label{P:PROPERTIESOFMUXMUZEROLEVELSET}

Working under the assumptions and conclusions of Theorem\,\ref{T:MAINRESULTSFROMSINGULARBOUNDARYPAPER}, let: 
	\begin{align} \label{E:DOMAINOFEMBEDDINGFORXMUEQUALSMINUSKAPPAYSURFACE}
	\domainforembeddingdatahypfortimefunctiontwoarg{0}{[0,\mupositive]}
	&
	\eqdef 
	\left\lbrace 
		(t,x^2,x^3) \in \mathbb{R} \times \mathbb{T}^2 
		\ | \
		\Cartesiantisafunctiononmumxtoriarg{\mupositive}{0}(x^2,x^3) \leq t 
		\le 
		\Cartesiantisafunctiononmumxtoriarg{0}{0}(x^2,x^3)
	\right\rbrace.
	\end{align}
	Then $\domainforembeddingdatahypfortimefunctiontwoarg{0}{[0,\mupositive]}$ is precompact,
	and there exists an embedding
	$\embeddingdatahypfortimefunctionarg{0}: \domainforembeddingdatahypfortimefunctiontwoarg{0}{[0,\mupositive]} 
	\rightarrow \twoargMrough{[\timefunction_0,0],[-\frac{3}{4} \interestingu,\frac{3}{4} \interestingu]}{0}$
	of the form
	$\embeddingdatahypfortimefunctionarg{0}(t,x^2,x^3) 
	=\left(t,\scalarembeddingdatahypfortimefunctionarg{0}(t,x^2,x^3),x^2,x^3 \right)$
	such that $\embeddingdatahypfortimefunctionarg{0}\in C^{1,1}(\domainforembeddingdatahypfortimefunctiontwoarg{0}{[0,\mupositive]})$ 
	and such that $\embeddingdatahypfortimefunctionarg{0}$
	is a diffeomorphism from $\domainforembeddingdatahypfortimefunctiontwoarg{0}{[0,\mupositive]}$ onto 
	$\datahypfortimefunctiontwoarg{0}{[\timefunction_0,0]}$. In particular, relative to the geometric coordinates $(t,u,x^2,x^3)$, we have: 
	\begin{align} \label{E:XMUISMINUSCAPPISAGRAPH} 
		\datahypfortimefunctiontwoarg{0}{[\timefunction_0,0]}
		& = 
		\left\lbrace 
			\left(t,\scalarembeddingdatahypfortimefunctionarg{0}(t,x^2,x^3),x^2,x^3 \right) 
			\in \twoargMrough{[\timefunction_0,0],[- \rightu,\leftu]}{0}
			\ | \
			(t,x^2,x^3) 
			\subset 
			\domainforembeddingdatahypfortimefunctiontwoarg{0}{[0,\mupositive]}
		\right\rbrace.
	\end{align}
	Moreover, 
	the following estimates hold:
	\begin{subequations}
	\begin{align} \label{E:BOUNDONEMBEDDINGOFXMUEQUALSMINUSKAPPA}
		\left\| 
			\embeddingdatahypfortimefunctionarg{0} 
		\right\|_{C^{1,1}\left(\domainforembeddingdatahypfortimefunctiontwoarg{0}{[0,\mupositive]}\right)}
		& \leq C,
			\\
		\left\| 	
			\left(\datasurfacep{x^2} \scalarembeddingdatahypfortimefunctionarg{0},
			\datasurfacep{x^3} \scalarembeddingdatahypfortimefunctionarg{0}\right) 
		\right\|_{C^{1,1}\left(\domainforembeddingdatahypfortimefunctiontwoarg{0}{[0,\mupositive]}\right)}
		& \leq C \initialsmall.
		\label{E:SMALLNESSBOUNDONEMBEDDINGOFXMUEQUALSMINUSKAPPA}
	\end{align}
	\begin{align} \label{E:PARTIALTEMBEDDINGFUNCTIONKEYESTIMATE}
		\datasurfacep{t} \scalarembeddingdatahypfortimefunctionarg{0}
		& = - \frac{\Lunit \muX \upmu \circ \embeddingdatahypfortimefunctionarg{0}}{\muX \muX \upmu 
			\circ 
			\embeddingdatahypfortimefunctionarg{0}} 
			+ 
			\mathcal{O}(\initialsmall).
	\end{align}
	\end{subequations}
	In \eqref{E:SMALLNESSBOUNDONEMBEDDINGOFXMUEQUALSMINUSKAPPA}--\eqref{E:PARTIALTEMBEDDINGFUNCTIONKEYESTIMATE}, $\{\datasurfacep{t},\datasurfacep{x^2},\datasurfacep{x^3}\}$ denote coordinate partial differentiation on $\domainforembeddingdatahypfortimefunctiontwoarg{0}{[0,\mupositive]}$.
\end{proposition}

\subsection{Sharp control of $\upmu$ and its derivatives}
\label{SS:SHARPCONTROLOFMUANDDERIVATIVESFORSINGULARBOUNDARY}

\begin{proposition}[Sharp control of $\upmu$ and its derivatives]
\label{P:SHARPCONTROLOFMUANDDERIVATIVESROUGHFOLIATIOSNREGION}
Suppose the statements and conclusions of Theorem\,\ref{T:MAINRESULTSFROMSINGULARBOUNDARYPAPER} hold. Then there exists a constant $0 < \interestingu < 1$ small enough such that the following estimates hold.

\medskip

\noindent \underline{\textbf{Minima of $\upmu$ occur precisely along $\twoargmumuxtorus{-\timefunction}{0}$}}.
For each fixed $\timefunction \in [\timefunction_0,0] = [-\mupositive,0]$,
we have:
		\begin{align} \label{E:MINVALUEOFMUONFOLIATION}
			\min_{\hypthreearg{\timefunction}{[- \rightu,\leftu]}{0}}
			\upmu
			& = 
			 - \timefunction,
		\end{align}
		and the minimum value of $- \timefunction$
		in \eqref{E:MINVALUEOFMUONFOLIATION} is achieved by $\upmu$ precisely on the $\upmu$-adapted torus 
		$
		\twoargmumuxtorus{-\timefunction}{0}
		$ defined in \eqref{E:MUXMUTORI}. In particular, the following holds along the crease:
		\begin{equation} \label{E:ALLROUGHSPATIALDERIVATIVESVANISHALONGTHECREASE}
			\roughgeop{u}  \upmu (p) = \roughgeop{x^2} \upmu (p) = \roughgeop{x^3} \upmu (p) = 0, \qquad \forall p \in \twoargmumuxtorus{0}{0}.
		\end{equation}

\medskip		
		
\noindent \underline{\textbf{$\upmu$ is large when $|u| \geq \interestingu$}}.
	The following lower bound holds, where $\boringregionmupositive > 0$ is the constant appearing in 
	\eqref{E:DATAMUISLARGEINBORINGREGION}:
	\begin{align} \label{E:MUISLARGEINBORINGREGIONFORSINGULARBOUNDARY}
			\min_{\twoargMrough{[\timefunction_0,0],[- \rightu,\leftu]}{0}
			\backslash 
			\twoargMrough{[\timefunction_0,0],[-\interestingu,\interestingu]}{0}} \upmu
			& \geq \frac{\boringregionmupositive}{2}.
\end{align}
	
\medskip

\noindent \underline{\textbf{Location of $\datahypfortimefunctiontwoarg{0}{[\timefunction_0,0]}$ 
and $\twoargmumuxtorus{-\timefunction}{0}$}}.
With
$\datahypfortimefunctiontwoarg{0}{[\timefunction_0,0]}$
and
$\twoargmumuxtorus{-\timefunction}{0}$
denoting the sets defined in \eqref{E:LEVELSETSOFMUXMU}--\eqref{E:MUXMUTORI},
we have:
\begin{subequations}
\begin{align} \label{E:MUXMUKAPPALEVELSETLOCATION}
			\datahypfortimefunctiontwoarg{0}{[\timefunction_0,0]}
			& 
			\subset
			\twoargMrough{[\timefunction_0,0],[-\frac{1}{2}\interestingu,\frac{1}{2}\interestingu]}{0},
				\\
		\text{for each $ \timefunction \in [\timefunction_0,0]$,} \qquad
		\twoargmumuxtorus{-\timefunction}{0}, 
		& 
		\subset
		\hypthreearg{\timefunction}{[-\frac{1}{2}\interestingu,\frac{1}{2}\interestingu]}{0}.
			\label{E:IMPROVEDLEVELSETSTRUCTUREANDLOCATIONOFMIN} 
\end{align}
\end{subequations}Moreover, with $\InverseCHOVroughtomumuxmu{0}$
denoting the inverse function of the function
$\CHOVroughtomumuxmu{0}$ defined in \eqref{E:CHOVFROMROUGHCOORDINATESTOMUWEGIGHTEDXMUCOORDINATES},
we have:
\begin{align} \label{E:INTERMSOFPHICHOVMAPIMPROVEDLEVELSETSTRUCTUREANDLOCATIONOFMINFORSINGULARBOUNDARY} 
\text{For each $\mulevelsetvalue \in [\mupositive,0]$,} \qquad
	\InverseCHOVroughtomumuxmu{0}
	\left(\lbrace \mulevelsetvalue \rbrace \times \lbrace - 0 \rbrace \times \mathbb{T}^2 \right)
	\subset
	\lbrace - \mulevelsetvalue \rbrace \times \left[-\frac{\interestingu}{2},\frac{\interestingu}{2}\right] \times \mathbb{T}^2.
\end{align}\medskip

\noindent \underline{\textbf{Transversal convexity of $\upmu$ and its consequences}}.
The following estimates hold:
\begin{align} 
	\begin{split} \label{E:MUTRANSVERSALCONVEXITYFORSINGULARBOUNDARY}
			\frac{\secondtransversalderivativemulowerbound}{2}
			& \leq 
			\min_{\twoargMrough{[\timefunction_0,0],[-\interestingu,\interestingu]}{0}}
				\left\lbrace			
				\muX \muX \upmu,
					\,
				\muX \muX \upmu - \frac{(\muX \upmu) \Lunit \muX \upmu}{\Lunit \upmu},
					\,
				\geop{u} \muX \upmu,
					\,
				\geop{u} \muX \upmu - \frac{(\geop{u} \upmu) \geop{t} \muX \upmu}{\geop{t} \upmu}, 
					\,
				\roughgeop{u} \muX \upmu
				\right\rbrace
					\\
			&
			\leq 
			\max_{\twoargMrough{[\timefunction_0,0],[-\interestingu,\interestingu]}{0}}
			\left\lbrace
				\muX \muX \upmu,
					\,
				\muX \muX \upmu - \frac{(\muX \upmu) \Lunit \muX \upmu}{\Lunit \upmu},
					\,
				\geop{u} \muX \upmu,
					\,
				\geop{u} \muX \upmu - \frac{(\geop{u} \upmu) \geop{t} \muX \upmu}{\geop{t} \upmu}, 
					\,
				\roughgeop{u} \muX \upmu
				\right\rbrace
			\leq 
			\frac{2}{\secondtransversalderivativemulowerbound},
	\end{split}
	\end{align}\begin{align} 
	\min_{\twoargMrough{[\timefunction_0,0],[-\interestingu,\interestingu]}{0}
			\backslash
			\mathcal{M}^{(0)}_{[\timefunction_0,0],[-\frac{\interestingu}{2},\frac{\interestingu}{2}]}} 
			|\muX \upmu | 
			& \geq \frac{\secondtransversalderivativemulowerbound \interestingu}{8}.
			\label{E:REGIONWHEREMUXMUKAPPALEVELSETISNOTLOCATED}
\end{align}Moreover, the following pointwise estimates hold on
$
\twoargMrough{[\timefunction_0,0],[-\interestingu,\interestingu]}{0}
$:
\begin{align} \label{E:WMANDRTRANSMUUBOUNDEDBYSQRTMU}
	|\muX \upmu|,
	& \leq 
	C \sqrt{\upmu}.
	\end{align}\medskip

\noindent \underline{\textbf{Rough control of null and almost null derivatives of $\upmu$ in the interesting region}}.
The following estimates hold,
where $\argLrough{0} = \frac{1}{\Lunit \timefunction} \Lunit$:
\begin{subequations}
\begin{align} \label{E:BOUNDSONLMUINTERESTINGREGIONFROMTAUFOLIATIONS} 
	- 
	\frac{9}{8}
	\mathring{\updelta}_*
	\leq
	\min_{\twoargMrough{[\timefunction_0,0],[-\interestingu,\interestingu]}{0}} \Lunit \upmu
	\leq 
	\max_{\twoargMrough{[\timefunction_0,0],[-\interestingu,\interestingu]}{0}} \Lunit \upmu
	\leq 
	- 
	\frac{7}{8}
	\mathring{\updelta}_*,
		\\
	- 
	\frac{9}{8}
	\mathring{\updelta}_*
	\leq
	\min_{\twoargMrough{[\timefunction_0,0],[-\interestingu,\interestingu]}{0}} \geop{t} \upmu
	\leq 
	\max_{\twoargMrough{[\timefunction_0,0],[-\interestingu,\interestingu]}{0}} \geop{t} \upmu
	\leq 
	- 
	\frac{7}{8}
	\mathring{\updelta}_*,
		\label{E:BOUNDSONGEOMETRICTDERIVATIVEMUINTERESTINGREGIONFORSINGULARBOUNDARY}
		\\
	- 
	\frac{3}{2}
	\leq
	\min_{\twoargMrough{[\timefunction_0,0],[-\interestingu,\interestingu]}{0}} \argLrough{0} \upmu
	\leq 
	\max_{\twoargMrough{[\timefunction_0,0],[-\interestingu,\interestingu]}{0}} \argLrough{0} \upmu
	\leq 
	- \frac{2}{3}.
	\label{E:SIZEBOUNDSONWIDETILDELMUINTERESTINGREGION}
	\end{align}
	\end{subequations}
	
	\noindent \underline{\textbf{Sharp control of $\Lunit \timefunctionarg{0}$}}.
	The following estimates hold:
	\begin{align} \label{E:ROUGHTIMEFUNCTIONLDERIVATIVEBOUNDS}
	\frac{7}{9} 
	\mathring{\updelta}_*
	\leq
	\min_{\twoargMrough{[\timefunction_0,0],[- \rightu,\leftu]}{0}} 
		\Lunit \timefunctionarg{0}
	\leq 
	\max_{\twoargMrough{[\timefunction_0,0],[- \rightu,\leftu]}{0}} 
		\Lunit \timefunctionarg{0}
	\leq 
	\frac{9}{7} \mathring{\updelta}_*.
	\end{align}\end{proposition}

\section{The double-null folliations} \label{S:DOUBLENULLFOLLIATIONSETUP}

In this section we begin the construction of the double-null folliations of spacetime. With $u$ serving as one of the double-null coordinates, we construct $\ubar$ in Sect.\,\ref{SS:CONSTRUCTIONOFUBAR} as the coordinate for the other transversal null direction, while still being adapted to the crease. That is, we will construct a function $\ubar$ on a subset of $\twoargMrough{[\timefunction_0,0],[-\rightu,\leftu]}{0}$ satisfying $\gfour(\Dfour \ubar,\Dfour\ubar) = 0$ and such that $\twoargmumuxtorus{0}{0} \subset \{ \ubar = 0\}$. Through a bootstrap argument, the ingoing eikonal function will be extended past $\twoargMrough{[\timefunction_0,0],[-\rightu,\leftu]}{0}$ so that a nontrivial part of the Cauchy horizon is also a subset of $\{\ubar = 0\}$.

\subsection{The ingoing eikonal function} \label{SS:CONSTRUCTIONOFUBAR}

In this section we merely provide the definition of $\ubar$. We will justify the validity of the Cauchy problem the eikonal equation defined here in Sect.\,\ref{S:CONSTRUCTIONANDESTIMATESFORINGOINGEIKONALFUNCTION}.

\begin{definition}[The ingoing eikonal function] \label{D:INGOINGEIKONALFUNCTION}
	Let $\datahypfortimefunctiontwoarg{0}{[\timefunction_0,0]}$ be the truncated level set of $\muX \upmu$ defined in \eqref{E:TRUNCATEDLEVELSETSOFMUXMU}. We define \textbf{the ingoing eikonal function} $\ubar$ to be the scalar function on the differential structure of the geometric coordinates $(t,u,x^2,x^3)$ solving the following eikonal equation initial value problem:
	\begin{subequations} \label{E:INGOINGEIKONALEQUATIONINITIALVALUEPROBLEM}
		\begin{align}
			(\gfour^{-1})^{\alpha\beta} \p_\alpha  \ubar \p_\beta \ubar &  = 0, \label{E:INGOINGEIKONALEQUATION} \\
			\Lunit \ubar & >  0, \label{E:TRANSVERSALITYCONDITIONFORINGOINGEIKONALWITHRESPECTOU} \\
			\ubar|_{\datahypfortimefunctiontwoarg{0}{[\timefunction_0,0]}} & = -\upmu|_{\datahypfortimefunctiontwoarg{0}{[\timefunction_0,0]}}. \label{E:DATAFORINGOINGEIKONALINITIALVALUEPROBLEM}
		\end{align}
	\end{subequations}
\end{definition}

\begin{remark}[The data for $\ubar$ is the same as that for $\timefunction$] \label{R:DATAFORINGOINGEIKONALISTHESAMEASROUGHTIMEFUNCTION} From \eqref{E:INITIALCONDITIONFORROUGHTIMEFUNCTION} and \eqref{E:DATAFORINGOINGEIKONALINITIALVALUEPROBLEM} we see that $\ubar$ and $\timefunction$ agree on $\datahypfortimefunctiontwoarg{0}{[\timefunction_0,0]}$. 
\end{remark}


\subsection{Double-null coordinates}
\label{SS:DOUBLENULLCOORDINATES}
Having constructed both the ingoing  and outgoing eikonal functions, we now define a coordinate system adapted to them.

\begin{definition}[The double-null coordinates and their partial derivative vector fields] \label{D:DOUBLENULLCOORDINATES}
	We call $(\ubar,u,x^2,x^3)$ the \textbf{double-null coordinates}. We denote the corresponding coordinate partial derivative vector fields by $\{\nullgeop{\ubar},\nullgeop{u},\nullgeop{x^2},\nullgeop{x^3}\}$. 
\end{definition}

Our analysis will show that the map 
$(t,x^1,x^2,x^3) 
\rightarrow 
\left(\ubar,u,x^2,x^3 \right)
$
is a homeomorphism -- all the way up to the Cauchy horizon -- and that it
is a diffeomorphism away from the crease.
Moreover, the map
$(t,u,x^2,x^3) 
\rightarrow 
\left(\ubar,u,x^2,x^3 \right)
$
is a diffeomorphism all the way up to the Cauchy horizon boundary;
see Theorem\,\ref{T:EXISTENCEUPTOCAUCHYHORIZONBYCONTINUATIONCRITERIA}.

\subsection{Corresponding subsets of spacetime} \label{SS:SUBSETSOFSPACETIME}

\begin{definition} \label{D:SUBSETSOFSPACETIME}
	Given intervals $I,J \subset \R$ and real numbers $\ubar,u \in \R$, we define:
	\begin{subequations}
		\begin{align}
			\doublenulltoritwoarg{\ubar}{u}  & \eqdef \{ (t,u,x^2,x^3) \in \R\times\R\times\T^2 \ | \ \ubar(t,u,x^2,x^3) = u\}, \label{E:DOUBLENULLTORI} \\
			\outgoingcharacteristicsurfacetwoarg{u}{I} & \eqdef \bigcup_{\ubar \in I} \doublenulltoritwoarg{\ubar}{u}, \label{E:OUTGOINGCHRACTERISTICSURFACE} \\
			\ingoingcharacteristicsurfacetwoarg{\ubar}{J} & \eqdef \bigcup_{u \in J}  \doublenulltoritwoarg{\ubar}{u},  \label{E:INGOINGCHARACTERISTICSURFACE} \\
			\characteristicdiamondtwoarg{I}{J} & \eqdef \bigcup_{\ubar \in I} \ingoingcharacteristicsurfacetwoarg{\ubar}{J} = \bigcup_{u \in J} \outgoingcharacteristicsurfacetwoarg{u}{I}. \label{E:CHARACTERISTICDIAMOND}
		\end{align}
	\end{subequations}
\end{definition}

We refer to the $\ingoingcharacteristicsurfacetwoarg{\ubar}{J}$ as \emph{ingoing characteristic surfaces} or simply \emph{ingoing characteristics}. We respectively refer to $\outgoingcharacteristicsurfacetwoarg{u}{I}$ as the \emph{outgoing characteristic surfaces} or simply \emph{outgoing characteristics}. We refer to the $\doublenulltoritwoarg{\ubar}{u}$ as \emph{double-null tori}. We also note that $\outgoingcharacteristicsurfacetwoarg{u}{I}$ are part of the $\gfour$-null surface $\nullhyparg{u}$. Moreover, we refer to the spacetime region $\characteristicdiamondtwoarg{I}{J}$ as \emph{characteristic diamond}. 

From Defs.\,\ref{D:DOUBLENULLCOORDINATES} and\,\ref{D:SUBSETSOFSPACETIME}, it follows that $\left\{\nullgeop{u},\nullgeop{x^2},\nullgeop{x^3}\right\}$ spans the tangent space of $\ingoingcharacteristicsurfacetwoarg{\ubar}{J}$, that $\left\{\nullgeop{x^2},\nullgeop{x^3}\right\}$ spans the tangent space of $\doublenulltoritwoarg{\ubar}{u}$, that $\left\{\nullgeop{\ubar},\nullgeop{x^2},\nullgeop{x^3}\right\}$ spans the tangent space of $\outgoingcharacteristicsurfacetwoarg{u}{I}$, and that $\left\{\nullgeop{\ubar},\nullgeop{u},\nullgeop{x^2},\nullgeop{x^3}\right\}$ spans the tangent space of $\characteristicdiamondtwoarg{I}{J}$.
  
\begin{definition}[Level sets of $\newuL \upmu$ and the $\newuL \upmu$-adapted tori] \label{D:LEVELSETOFNEWULMUANDNEWULADAPTEDTORI}
Given a non-negative number $\ubar' \le 0$, we define:
	\begin{subequations}
		\begin{align}
			\newuLmulevelset & \eqdef \{ (t,u,x^2,x^3) \in \R \times \R \times \T^2 \ | \ \newuL \upmu(t,u,x^2,x^3) = 0\} \cap \{ |u| \le \interestingu\}, \label{E:DEFOFLEVELSETOFNEWULMU} \\
			\argubarnewultorus{\ubar'} & \eqdef \{ (t,u,x^2,x^3) \in \R\times\R\times \T^2 \ | \ \ubar(t,u,x^2,x^3) = \ubar'\} \cap \newuLmulevelset. \label{E:DEFOFNEWULADAPTEDTORI}
		\end{align}
	\end{subequations}
We sometimes find it useful to truncate $\newuLmulevelset$ as: 
	\begin{align}
		\newuLmulevelset^{I} \eqdef \newuLmulevelset \cap \{ \ubar \in I\}. \label{E:DEFOFTRUNCATEDLEVELSETOFNEWULMU}
	\end{align}
\end{definition}

\section{Coordinate Transformations} \label{S:COORDINATETRANSFORMATIONS}	

To prove our main results, we will have to control the transformations between the Cartesian coordinates $(t,x^1,x^2,x^3)$,
the geometric coordinates $(t,u,x^2,x^3)$, and
the double-null coordinates $(\ubar,u,x^2,x^3)$.
In this section, we define the relevant change of variables maps and derive some basic
relationships between the partial derivative vectorfields in the different coordinate systems as well as the geometric vectorfields.

\subsection{Change of variables maps}
\label{SS:ALLTHECHOVMAPS}
In this short section, we define several of the change of variables maps that we use to prove our main results. Recall that the map from geometric to Cartesian coordinates was defined in \eqref{E:CHOVGEOTOCARTESIAN}. 

\begin{definition}[Change of variables maps] \label{D:CHOVINVOLVINGTHEDOUBLENULLCOORD} We define the change of variables map from geometric coordinates to double-null coordinates and its Jacobian as follows:
	\begin{subequations}
		\begin{align}
			\CHOVgeotodoublenull(t,u,x^2,x^3) & \eqdef (\ubar,u,x^2,x^3),  \label{E:CHOVFROMGEOTODOUBLENULLCOORDINATES} \\ 
			\CHOVJacobiangeotodoublenull(t,u,x^2,x^3) & \eqdef \frac{\p(\ubar,u,x^2,x^3)}{\p(t,u,x^2,x^3)}. \label{E:JACOBIANOFCHOVFROMGEOTODOUBLENULLCOORDINATES}
		\end{align}
	\end{subequations}
We define the change of variables map from double-null coordinates to ``$(\ubar,\newuL \upmu,x^2,x^3)$-space''
as follows:
		\begin{align}
			\CHOVdoublenulltoubarnewulmu(\ubar,u,x^2,x^3) & \eqdef (\ubar, \newuL\upmu(\ubar,u,x^2,x^3), x^2,x^3), \label{E:CHOVFROMDOUBLENULLTOUBARNEWULMUCOORDINATES} 
		\end{align}
			
\end{definition}

\subsection{Coordinate double-null partial derivative transformations} \label{SS:COORDINATEDOUBLENULLPARTIALDERIVATIVETRANSFORMATIONS}

\begin{lemma}[Relationship between $\{\p_1,\p_{x^2},\p_{x^3}\}$ and $\{\nullgeop{\ubar},\nullgeop{u},\nullgeop{x^2},\nullgeop{x^3}\}$] \label{L:RELATIONSHIPBETWEENSPATIALCARTESIANDERIVATIVESSANDDOUBLENULL}
The following identities hold:
\end{lemma}

\begin{lemma}[Relationship between $\{\geop{t},\geop{u},\geop{x^2},\geop{x^3}\}$ and $\{\nullgeop{\ubar},\nullgeop{u},\nullgeop{x^2},\nullgeop{x^3}\}$] \label{L:RELATIONSHIPBETWEENGEOMETRICPARTIALDERIVATIVESANDDOUBLENULL}
	The following identities hold:
	\begin{subequations}
		\begin{align}
			\nullgeop{\ubar} & = \frac{1}{\geop{t} \ubar} \geop{t}, \label{E:NULLCOORDINATEPARTIALINGOINGEIKONALINTERMSOFGEOMETRICVF}
			\\
			\nullgeop{u} & = \geop{u} - \frac{\geop{u} \ubar}{\geop{t} \ubar} \geop{t}, \label{E:NULLCOORDINATEPARTIALUINTERMSOFGEOMETRICVF} \\
			\nullgeop{x^A} & = \geop{x^A} - \frac{\geop{x^A} \ubar}{\geop{t} \ubar} \geop{t} =  \geop{x^A} +  \frac{\geop{x^A} \ubar}{\geop{t} \ubar}  L^B \geop{x^B} -  \frac{\geop{x^A} \ubar}{\geop{t} \ubar} \Lunit . \label{E:NULLCOORDINATEPARTIALXAINTERMSOFGEOMETRICVF} 
		\end{align}
	\end{subequations}
	
	Moreover, the following identity holds:
\begin{align} \label{E:SMOOTHANGULARVECTORFIELDSINTERMSOFNULLONES}
\geop{x^A} 
& = 
\nullgeop{x^A} 
- 
\frac{\geop{x^A} \ubar}{\Lunit \ubar} \Lunit^B \nullgeop{x^B} 
+ 
\frac{\geop{x^A} \ubar}{\Lunit \ubar} \Lunit.
\end{align}
\end{lemma}
\begin{proof}
	The identity \eqref{E:NULLCOORDINATEPARTIALINGOINGEIKONALINTERMSOFGEOMETRICVF} follows from the chain rule identity $\geop{t} = (\geop{t} \ubar) \nullgeop{\ubar} + (\geop{t} u) \nullgeop{u} + (\geop{t} x^A) \nullgeop{x^A}$ and the fact that $\geop{t} u = \geop{t} x^A = 0$. The identities \eqref{E:NULLCOORDINATEPARTIALUINTERMSOFGEOMETRICVF}--\eqref{E:NULLCOORDINATEPARTIALXAINTERMSOFGEOMETRICVF} follow from similar arguments and \eqref{E:NULLCOORDINATEPARTIALINGOINGEIKONALINTERMSOFGEOMETRICVF}.
	
	Finally, with the help of \eqref{E:LUNITANDULUNITAPPLIEDTOEIKONALANDCARTESIANTIME},
it is straightforward to confirm the identity \eqref{E:SMOOTHANGULARVECTORFIELDSINTERMSOFNULLONES} 
by checking that both sides evaluate to the same values when acting on the double-null coordinate functions 
$\ubar,u,x^2,x^3$. 
\end{proof}

\section{The double-null acoustic geometry and curvature tensors} \label{S:DOUBLENULLGEOMETRICFRAMEWORK}

In this section, we set up the acoustic geometry of the double-null foliations. That is, the geometry associated to the ingoing characteristics $\ingoingcharacteristicsurfacetwoarg{\ubar}{J}$, the outgoing characteristics $\outgoingcharacteristicsurfacetwoarg{u}{I}$, and the double-null tori $\doublenulltoritwoarg{\ubar}{u}$. In particular, we construct several vectorfields from the ingoing and outgoing eikonal functions $(\ubar,u)$. We also construct several null-geometry scalars which can be viewed as $(\ubar,u,x^2,x^3)$-analogs of the several $(t,u,x^2,x^3)$-acoustic scalars seen in prior works \cite{dC2007, gHsKjSwW2016, jLjS2018, LukSpeck2024stability, abbrescia2022emergence} such as $\upmu$. We also define the first fundamental forms induced by the spacetime metric $\gfour$ on $\doublenulltoritwoarg{\ubar}{u}$, exhibit various geometric decompositions, and introduce the Riemann and Ricci curvature of the spacetime metric $\gfour$ and various  curvature tensors of the first fundamental form of $\doublenulltoritwoarg{\ubar}{u}$.

\subsection{Geometric vectorfields and related scalar functions}
\label{SS:GEOMETRICVECTORFIELDSANDSCALARFUNCTIONS}

\subsubsection{$\gfour$-null vectorfields} \label{SSS:NULLVECTORFIELDS}

Associated to $u$ and $\underline{u}$, we define the geodesic vectorfields:\footnote{More precisely, 
	using the eikonal equations, is straightforward to show
	that $\Dfour_{\Lgeo} \Lgeo = \Dfour_{\uLgeo} \uLgeo = 0$, where $\Dfour$ is the Levi--Civita connection of $\gfour$.} 
\begin{align} \label{E:NULLGEODESIVECTORFIELDS}
	\Lgeo^{\alpha}
	& \eqdef - (\gfour^{-1})^{\alpha \beta} \partial_{\beta} u,
	&
	\uLgeo^{\alpha}
	& \eqdef - 	
	(\gfour^{-1})^{\alpha \beta} \partial_{\beta} \underline{u}.
\end{align}
From \eqref{E:NULLGEODESIVECTORFIELDS}, it follows that $\Lgeo$ is $\gfour$-orthogonal to $\outgoingcharacteristicsurfacetwoarg{u}{I} $,
while $\uLgeo$ is $\gfour$-orthogonal to $\ingoingcharacteristicsurfacetwoarg{\ubar}{J}$.
The equations in \eqref{E:REGULAREIKONALEQUATION} and \eqref{E:INGOINGEIKONALEQUATIONINITIALVALUEPROBLEM} imply that $\Lgeo$ and $\uLgeo$ are $\gfour$-null:
\begin{align} \label{E:NULLGEODESIVECTORFIELDSAREINFACTNULL}
	\gfour(\Lgeo,\Lgeo)
	&
	= \gfour(\uLgeo,\uLgeo)
	= 0.
\end{align}Next, we define the following scalar functions on $\characteristicdiamondtwoarg{I}{J}$:
\begin{subequations} \label{E:NULLACOUSTICSCALARS}
	\begin{align} \label{E:LFOLIATIONDENSITY}
		\upmu
		& \eqdef 
		\frac{- 1}{
			(\gfour^{-1})^{\alpha \beta}
			\partial_{\alpha} t 
			\partial_{\beta} u}
		=
		\frac{1}{\Lgeo^0} = \frac{1}{\Transport u}, 
		\\
		\underline{\upmu}
		& \eqdef 
		\frac{-1}{
			(\gfour^{-1})^{\alpha \beta}
			\partial_{\alpha} t 
			\partial_{\beta} \underline{u}}
		=
		\frac{1}{\uLgeo^0} = \frac{1}{\Transport \ubar},
		\label{E:ULFOLIATIONDENSITY} \\
		\MagnitueofinnerproductofnewLandnewuL
		& \eqdef 
		\frac{-1}{
			(\gfour^{-1})^{\alpha \beta}
			\partial_{\alpha} u 
			\partial_{\beta} \underline{u}}
		=
		\frac{-1}{\gfour(\Lgeo,\uLgeo)},
\label{E:RECIPROCALOFNEGATIVENULLGEODESICVECTORFIELDINNERPRODUCT}
		\\
		\ReciprocalLunitAppliedtoTimeFunction 
		& \eqdef \frac{\MagnitueofinnerproductofnewLandnewuL}{\upmu},
		&
		\ReciprocaluLunitAppliedtoTimeFunction
		& \eqdef \frac{\MagnitueofinnerproductofnewLandnewuL}{\underline{\upmu}}.
		\label{E:RATIOOFNULLGEOSICINNERPRODUCTANDFOLIATIONDENSITY}
	\end{align}
\end{subequations}
We also define: 
	\begin{align}
		\MagnitueofinnerproductofLunitanduLunit \eqdef \frac{\MagnitueofinnerproductofnewLandnewuL}{\ReciprocalLunitAppliedtoTimeFunction \ReciprocaluLunitAppliedtoTimeFunction}. \label{E:DEFOFINNERPRODUCTOFLUNITANDULUNIT}
	\end{align}

	The bootstrap assumptions of Sect.\,\ref{S:BOOTSTRAPEVERYTHINGEXCEPTENERGIES} will imply that:
\begin{align} \label{E:POSITIVITYOFFOLIATIONDENSITYANDNULLGEODESICINNERPRODUCTETC}
	\upmu 
	& > 0,
	&
	\underline{\upmu}
	& > 0,
	&
	\MagnitueofinnerproductofnewLandnewuL
	& > 0,
	&
	\ReciprocalLunitAppliedtoTimeFunction
	& > 0,
	&
	\ReciprocaluLunitAppliedtoTimeFunction
	& > 0, 
	& \MagnitueofinnerproductofLunitanduLunit > 0.
\end{align}Next, we define the following vectorfields, which are rescaled versions of $\Lgeo$ and $\uLgeo$:
\begin{subequations}
	\begin{align}
		\Lunit^{\alpha}
		& \eqdef - \upmu (\gfour^{-1})^{\alpha \beta} \partial_{\beta} u,
		&
		\uLunit^{\alpha}
		& \eqdef - 	
		\underline{\upmu} (\gfour^{-1})^{\alpha \beta} \partial_{\beta} \underline{u}
		\label{E:DOUBLENULLCARTESIANTIMENORMALIZEDNULLVECTORFIELDS} \\
		\newL^{\alpha}
		& \eqdef - \MagnitueofinnerproductofnewLandnewuL (\gfour^{-1})^{\alpha \beta} \partial_{\beta} u,
		&
		\newuL^{\alpha}
		& \eqdef - \MagnitueofinnerproductofnewLandnewuL (\gfour^{-1})^{\alpha \beta} \partial_{\beta} \underline{u}.
		\label{E:DOUBLENULLEIKONALFUNCTIONNORMALIZEDNULLVECTORFIELDS}
	\end{align}
\end{subequations}The following identities easily follow from the above definitions:
\begin{align} \label{E:DOUBLENULLALLTHENEWVECTORFIELDSARENULL}
	\gfour(\Lunit,\Lunit)
	& = \gfour(\uLunit,\uLunit)
	= \gfour(\newL,\newL)
	= \gfour(\newuL,\newuL)
	= 0,
\end{align}\begin{subequations}
	\begin{align} \label{E:INNERPRODUCTOFCARTESIANNORMALIZEDNULLVECTORFIELDS}
		\gfour(\Lunit,\uLunit)
		& = - \frac{ \upmu \underline{\upmu}}{\MagnitueofinnerproductofnewLandnewuL}
		= -\frac{ \MagnitueofinnerproductofnewLandnewuL}{\ReciprocalLunitAppliedtoTimeFunction \ReciprocaluLunitAppliedtoTimeFunction} =-\MagnitueofinnerproductofLunitanduLunit,
		\\
		\gfour(\newL,\newuL)
		& = - \MagnitueofinnerproductofnewLandnewuL,
		\label{E:INNERPRODUCTOFEIKONALFUNCTIONNORMALIZEDNULLVECTORFIELDS}
	\end{align}
\end{subequations}
\begin{subequations}
	\begin{align}
		\Lunit u 
		& = \uLunit \underline{u}
		= 0,
		&
		\Lunit \underline{u}
		& = \frac{1}{\ReciprocalLunitAppliedtoTimeFunction},
		\,
		\uLunit u
		= \frac{1}{\ReciprocaluLunitAppliedtoTimeFunction},
		&
		\Lunit t 
		& = \uLunit t 
		= 1,
		\label{E:LUNITANDULUNITAPPLIEDTOEIKONALANDCARTESIANTIME}
	\end{align}
	\begin{align}
		\newL u 
		& = \newuL \underline{u}
		= 0,
		&
		\newL \underline{u}
		& =
		\newuL u
		= 1,
		\label{E:GEOMETRICCOORDINATENULLVECTORFIELDSAPPLIEDTOGEOMETRICCOORDINATES}
	\end{align}
\end{subequations}\begin{align} \label{E:RELATIONBETWEENCARTESIANNORMALIZEDNULLVECTORFIELDSANDEIKONALFUNCTIONORMALIZEDNULLVECTORFIELDS}
	\newL
	& = \ReciprocalLunitAppliedtoTimeFunction \Lunit,
	&
	\newuL
	& = \ReciprocaluLunitAppliedtoTimeFunction \uLunit.
\end{align}Note that by \eqref{E:MATERIALDERIVATIVELOWEREDCARTESIANCOORDINATES}, the last two equalities in \eqref{E:LUNITANDULUNITAPPLIEDTOEIKONALANDCARTESIANTIME} are equivalent to:
\begin{align} \label{E:EQUIVALENTLUNITANDULUNITAPPLIEDTOEIKONALANDCARTESIANTIME}
	\gfour(\Lunit,\Transport)
	& 
	= \gfour(\uLunit,\Transport)
	= - 1.
\end{align}\begin{lemma}[Relationship between $\{\Lunit,\uLunit,\newL,\newuL\}$ and $\{\nullgeop{\ubar},\nullgeop{u},\nullgeop{x^2},\nullgeop{x^3}\}$] \label{L:RELATIONSHIPBETWEENNULLCOORDINATEPARTIALSANDDOUBLENULLVECTORFIELDS}
	Let $\Lunit^A = \Lunit x^A$ and $\uLunit^A \eqdef \uLunit x^A$. Then the following identities hold:
	\begin{subequations}
		\begin{align}
			\Lunit & = \frac{1}{\ReciprocalLunitAppliedtoTimeFunction} \nullgeop{\ubar} + L^A \nullgeop{x^A}, \label{E:LUNITINTERMSOFNULLCOORDINATEPARTIALS} \\
			\uLunit & = \frac{1}{\ReciprocaluLunitAppliedtoTimeFunction} \nullgeop{u} + \uLunit^A \nullgeop{x^A}, \label{E:ULUNITINTERMSOFNULLCOORDINATEPARTIALS} \\
			\newL & = \nullgeop{\ubar} + \ReciprocalLunitAppliedtoTimeFunction L^A \nullgeop{x^A}, \label{E:NEWLINTERMSOFNULLCOORDINATEPARTIALS} \\
			\newuL & =  \nullgeop{u} + \ReciprocaluLunitAppliedtoTimeFunction \uLunit^A \nullgeop{x^A}. \label{E:NEWULINTERMSOFNULLCOORDINATEPARTIALS}
		\end{align}
	\end{subequations}
\end{lemma}

\begin{proof}
	We express $\Lunit$ in terms the double-null coordinate differential structure as $\Lunit = \Lunit \ubar \nullgeop{\ubar} + \Lunit u \nullgeop{u} + L^a \nullgeop{x^A}$. Using \eqref{E:LUNITANDULUNITAPPLIEDTOEIKONALANDCARTESIANTIME}, \eqref{E:LUNITINTERMSOFNULLCOORDINATEPARTIALS} follows. Multiplying \eqref{E:LUNITINTERMSOFNULLCOORDINATEPARTIALS} by $\ReciprocalLunitAppliedtoTimeFunction$ and using \eqref{E:RELATIONBETWEENCARTESIANNORMALIZEDNULLVECTORFIELDSANDEIKONALFUNCTIONORMALIZEDNULLVECTORFIELDS} proves \eqref{E:NEWLINTERMSOFNULLCOORDINATEPARTIALS}. The analogous proof for $\uLunit$ and $\newuL$ holds from identical arguments; we omit the details. 
\end{proof}

%
%

\subsection{First fundamental forms and projections}
\label{SS:DOUBLENULLFIRSTFUNDANDPROJECTIONS}

We now introduce the first fundamental forms of the double-null tori $\doublenulltoritwoarg{\ubar}{u}$ induced by the spacetime metric $\gfour$ and the $\gfour$-orthogonal projection tensor onto $\doublenulltoritwoarg{\ubar}{u}$.  

\begin{definition}[The first fundamental forms of $\doublenulltoritwoarg{\ubar}{u}$] \label{D:FIRSTFUNDAMENTALFORMSOFDOUBLENULLTORI} \hfill
\begin{itemize}
	\item We define the \emph{first fundamental forms of} $\doublenulltoritwoarg{\ubar}{u}$ relative to $\gfour$ to be the symmetric type $\binom{0}{2}$ tensorfield $\gnulltori$ such that $\gnulltori(Y,Z) = \gfour(Y,Z)$ for all pairs of vectorfields tangent to $\doublenulltoritwoarg{\ubar}{u}$ and such that $\gnulltori(V,\cdot) = \gnulltori(\cdot,V) = 0$ if $V$ is $\gfour$-orthogonal to $\doublenulltoritwoarg{\ubar}{u}$. 
	\item We define the \emph{inverse first fundamental form} $\gnulltori^{-1}$ to be the dual of $\gnulltori$ relative to $\gfour$, i.e., relative to arbitrary coordinates, it is the symmetric type $\binom{2}{0}$ tensorfield wit the following components:
	\begin{equation} \label{E:INVERSEFIRSTFUNDFORMOFDOUBLENULLTORI} 
		(\gnulltori^{-1})^{\alpha\beta} = (\gfour^{-1})^{\alpha\gamma} (\gfour^{-1})^{\beta\delta} \gnulltori_{\gamma\delta}.
	\end{equation}
\end{itemize}
\end{definition}

\begin{definition}[$\gfour$-orthogonal projection onto the double-null tori $\doublenulltoritwoarg{\ubar}{u}$ and $\doublenulltoritwoarg{\ubar}{u}$-tangency] \label{D:DOUBLENULLORTHOGONALPROJECTIONANDTANGENCY} \hfill
	\begin{itemize}
		\item We define the $\doublenulltoritwoarg{\ubar}{u}$-projection tensorfield $\nulltorusproject$ as follows, where $\MagnitueofinnerproductofnewLandnewuL$ is the scalar function from \eqref{E:RECIPROCALOFNEGATIVENULLGEODESICVECTORFIELDINNERPRODUCT} and $\newL$, $\newuL$ are the vectorfields from \eqref{E:DOUBLENULLEIKONALFUNCTIONNORMALIZEDNULLVECTORFIELDS}:
			\begin{align}
				\nulltorusproject_\beta^\alpha & \eqdef \updelta_{\beta}^{\alpha}
				+
				\frac{1}{\MagnitueofinnerproductofnewLandnewuL}
				\newL^{\alpha} \newuL_{\beta}
				+
				\frac{1}{\MagnitueofinnerproductofnewLandnewuL}
				\newuL^{\alpha} \newL_{\beta}, 
				\label{E:DOUBLENULLTORIPROJECTIONDEFININGEQUATION} \\
				& = \updelta_\beta^\alpha + \frac{1}{\MagnitueofinnerproductofLunitanduLunit} \Lunit^\alpha \uLunit_\beta + \frac{1}{\MagnitueofinnerproductofLunitanduLunit} \uLunit^\alpha \Lunit_\beta. \label{E:DOUBLENULLTORIPROJECTIONDEFININGEQUATIONWRTUNITL}
			\end{align}
		\item Given any type $\binom{m}{n}$ spacetime tensorfield $\xi$, we define its $\gfour$-orthogonal projection onto $\doublenulltoritwoarg{\ubar}{u}$, denoted by $\nulltorusproject \xi$, as follows: 
			\begin{equation} \label{E:TENSORFIELDPROJECTIONTONULLTORIDEFININGEQUATION}
				(\nulltorusproject \xi)_{\beta_1 \cdots \beta_n}^{\alpha_1\cdots\alpha_m} \eqdef \nulltorusproject_{\widetilde{\alpha}_1}^{\alpha_1} \cdots \nulltorusproject_{\widetilde{\alpha}_m}^{\alpha_m} \nulltorusproject^{\widetilde{\beta}_1}_{\beta_1} \cdots \nulltorusproject^{\widetilde{\beta}_n}_{\beta_n} \xi_{\widetilde{\beta}_1 \cdots \widetilde{\beta}_n}^{\widetilde{\alpha}_1 \cdots \widetilde{\alpha}_m}.
			\end{equation}
		\item Let $f$ be a scalar function. We define $\nullangrmD f$ to be the $\doublenulltoritwoarg{\ubar}{u}$-tangent one-form given by:
			\begin{align} \label{E:DOUBLENULLTORIDIFFERENTIAL}
				\nullangrmD f \eqdef \nulltorusproject \mathrm{d} f.
			\end{align}
			Note that $\nullangrmD f \cdot \nullgeop{x^A} = \nullgeop{x^A} f$ for $A = 2,3$.
	\end{itemize}
\end{definition} 

\begin{lemma}[Decomposition of $\gfour$ and $\gfour^{-1}$] \label{L:DECOMPOSITIONOFACOUSTICALMETRICANDINVERSEACOUSTICALMETRICINTERMSOFNULLVECTORFIELDSANDNULLTORIMETRIC}
The spacetime metric $\gfour$ and the inverse spacetime metric $\gfour^{-1}$ can be expressed as follows relative to the vectorfields $\newL$ and $\newuL$ from \eqref{E:DOUBLENULLEIKONALFUNCTIONNORMALIZEDNULLVECTORFIELDS} and the first fundamental forms $\gnulltori$ and $\gnulltori^{-1}$ from Def.\,\ref{D:FIRSTFUNDAMENTALFORMSOFDOUBLENULLTORI}: 
	\begin{subequations}
		\begin{align}
			\gfour^{-1}
			& = 
			-
			\frac{1}{\MagnitueofinnerproductofnewLandnewuL}
			\newL \otimes \newuL
			-
			\frac{1}{\MagnitueofinnerproductofnewLandnewuL}
			\newuL \otimes \newL 
			+
			\gnulltori^{-1},		\label{E:ACOUSTICALINVERSEMETRICINTERMSOFNULLVECTORFIELDSANDDOUBLENULLSPHEREINVERSEFIRSTFUND} \\
			& = 
			-
			\frac{1}{\MagnitueofinnerproductofLunitanduLunit}
			\Lunit \otimes \uLunit
			-
			\frac{1}{\MagnitueofinnerproductofnewLandnewuL}
			\uLunit \otimes \Lunit 
			+
			\gnulltori^{-1},	\label{E:ACOUSTICALINVERSEMETRICINTERMSOFUNITNULLVECTORFIELDSANDDOUBLENULLSPHEREINVERSEFIRSTFUND} \\
			\gfour 
			& = 
			-
			\frac{1}{\MagnitueofinnerproductofnewLandnewuL}
			\newL_{\flat} \otimes \newuL_{\flat}
			-
			\frac{1}{\MagnitueofinnerproductofnewLandnewuL} 
			\newuL_{\flat} \otimes \newL_{\flat}
			+
			\gnulltori,
			\label{E:ACOUSTICALMETRICINTERMSOFNULLVECTORFIELDSANDDOUBLENULLSPHEREFIRSTFUND} \\
			& = 
			-
			\frac{1}{\MagnitueofinnerproductofLunitanduLunit}
			\Lunit_{\flat} \otimes \uLunit_{\flat}
			-
			\frac{1}{\MagnitueofinnerproductofLunitanduLunit} 
			\uLunit_{\flat} \otimes \Lunit_{\flat}
			+
			\gnulltori, \label{E:ACOUSTICALMETRICINTERMSOFUNITNULLVECTORFIELDSANDDOUBLENULLSPHEREFIRSTFUND}			
		\end{align}
	\end{subequations}
where $\newL_{\flat},\, \newuL_{\flat}, \Lunit_{\flat}, \uLunit_{\flat}$  in \eqref{E:ACOUSTICALMETRICINTERMSOFNULLVECTORFIELDSANDDOUBLENULLSPHEREFIRSTFUND} denote the one-form $\gfour$-dual to $\newL, \, \newuL, \, \Lunit,\, \uLunit$, respectively. Identity \eqref{E:ACOUSTICALINVERSEMETRICINTERMSOFUNITNULLVECTORFIELDSANDDOUBLENULLSPHEREINVERSEFIRSTFUND} then follows from In particular, the following holds:
\begin{equation} \label{E:NULLTORIFIRSTFUNDISPROJECTIONOFACOUSTICALMETRIC}
	\gnulltori = \nulltorusproject \gfour.
\end{equation}

\end{lemma}
\begin{proof}
	Recall that $\mbox{\upshape span} \lbrace \newuL, \newL \rbrace$ is equal to the $\gfour$-orthogonal complement of $\doublenulltoritwoarg{\ubar}{u}$.
With the help of \eqref{E:INNERPRODUCTOFEIKONALFUNCTIONNORMALIZEDNULLVECTORFIELDS}, \eqref{E:ACOUSTICALINVERSEMETRICINTERMSOFNULLVECTORFIELDSANDDOUBLENULLSPHEREINVERSEFIRSTFUND} follows from contracting both sides with combinations of $\{\newuL, \newL, \nullgeop{x^2},\nullgeop{x^3}\}$. Identity \eqref{E:ACOUSTICALINVERSEMETRICINTERMSOFUNITNULLVECTORFIELDSANDDOUBLENULLSPHEREINVERSEFIRSTFUND} follows from the now proved \eqref{E:ACOUSTICALINVERSEMETRICINTERMSOFNULLVECTORFIELDSANDDOUBLENULLSPHEREINVERSEFIRSTFUND}, \eqref{E:DEFOFINNERPRODUCTOFLUNITANDULUNIT} and \eqref{E:RELATIONBETWEENCARTESIANNORMALIZEDNULLVECTORFIELDSANDEIKONALFUNCTIONORMALIZEDNULLVECTORFIELDS}.

The metric identity
\eqref{E:ACOUSTICALMETRICINTERMSOFNULLVECTORFIELDSANDDOUBLENULLSPHEREFIRSTFUND} follows from expressing \eqref{E:ACOUSTICALINVERSEMETRICINTERMSOFNULLVECTORFIELDSANDDOUBLENULLSPHEREINVERSEFIRSTFUND} in terms of Cartesian coordinates and lowering the indeces with $\gfour$. Identity \eqref{E:NULLTORIFIRSTFUNDISPROJECTIONOFACOUSTICALMETRIC} readily follows from \eqref{E:ACOUSTICALMETRICINTERMSOFNULLVECTORFIELDSANDDOUBLENULLSPHEREFIRSTFUND}. Finally, \eqref{E:ACOUSTICALMETRICINTERMSOFUNITNULLVECTORFIELDSANDDOUBLENULLSPHEREFIRSTFUND} follows from the same arguments as in the proof of \eqref{E:ACOUSTICALINVERSEMETRICINTERMSOFUNITNULLVECTORFIELDSANDDOUBLENULLSPHEREINVERSEFIRSTFUND}.

\end{proof}

The following lemma reveals the relationship between geometry of the double-null tori $\doublenulltoritwoarg{\ubar}{u}$ and the flat tori $\ell_{t,u}$.

\begin{lemma}[Relationship between $\gtorus$ and $\gnulltori$] \label{L:RELATIONSHIPBETWEENFLATTORIANDDOUBLENULLTORIMETRICS}
	When restricted to the tangent space of $\doublenulltoritwoarg{\ubar}{u}$, we have the following identity for $\gnulltori$ relative to the double-null coordinates:
		\begin{equation} \label{E:DOUBLENULLTORIFIRSTFUNDRESTRICTEDTODOUBLENULLTORI}
			\gnulltori = \gnulltori\left(\nullgeop{x^A},\nullgeop{x^B}\right) \mathrm{d} x^A \otimes \mathrm{d}x^B,
		\end{equation}
	where with $\gtorus_{AB} = \gtorus(\geop{x^A},\geop{x^B})$, 
	we have:
		\begin{align} \label{E:GNULLTORICOMPONENTS}
			\gnulltori\left(\nullgeop{x^A},\nullgeop{x^B}\right) & = 
			\gtorus_{AB} 
			+   
			\frac{\geop{x^A}\ubar}{\geop{t} \ubar} \gtorus_{BC} \Lunit^C  
			+   
			\frac{\geop{x^B}\ubar}{\geop{t} \ubar} \gtorus_{AC} \Lunit^C 
			+ 
			\frac{(\geop{x^A} \ubar) \geop{x^B} \ubar}{(\geop{t} \ubar)^2} \gtorus_{CD} \Lunit^C \Lunit^D.
		\end{align}
	Moreover, relative to the double-null coordinates, the following identity holds:
		\begin{align} \label{E:GNULLTORIINVERSECOMPONENDS}
			\gnulltori^{-1} = \gnulltori^{-1}\left(\mathrm{d} x^A,\mathrm{d} x^B\right) \nullgeop{x^A}\otimes \nullgeop{x^B},
		\end{align}
	where: 
		\begin{align} \label{E:DOUBLENULLTORIMETRICCOMPONENTSANDTHEINVERSECOMPONENTSRELATION}
			\gnulltori^{-1}\left(\mathrm{d} x^A,\mathrm{d} x^C \right)\gnulltori\left(\nullgeop{x^C},\nullgeop{x^B}\right) = \updelta_B^A,
		\end{align}
	and $\updelta_B^A$ is the Kronecker delta.

	Next, let $\gtorusdoublenullCOV$ be the $2\times 2$ matrix with the following entries:
	\begin{equation}
		\gtorusdoublenullCOV^B_A \eqdef \updelta_A^B - \frac{\geop{x^A} \ubar}{\Lunit \ubar} \Lunit^B, \label{E:CHOVCOEFFICIENTSSMOOTHANGULARDERIVATIVESINTERMSOFDOUBLENULLONESANDL}
	\end{equation}
	where $\updelta_A^B$ is the Kronecker delta. Then the following identity holds:
	\begin{equation}
		\geop{x^A} = \gtorusdoublenullCOV^B_A \nullgeop{x^B} + \frac{\geop{x^A} \ubar}{\Lunit \ubar} \Lunit. \label{E:SMOOTHANGULARDERIVATIVESINTERMSOFDOUBLENULLONESANDL}
	\end{equation}
	Furthermore, 
	we have the following relationship between 
	$\gtorus_{AB}$	and 
	$
	\gnulltori\left(\nullgeop{x^A},\nullgeop{x^B} \right)
	$:
	\begin{subequations}
	\begin{align} \label{E:SMOOTHTORUSFIRSTFUNDCOMPONENTSINTERMSOFDOUBLENULLTORUSFIRSTFUNDCOMPONENTS}
		\gtorus_{AB}
		& 
		= \gtorusdoublenullCOV_A^C \gtorusdoublenullCOV_B^D \gnulltori\left(\nullgeop{x^C},\nullgeop{x^D} \right),
		\\
		(\gtorus^{-1})^{AB} 
		& 
		= (\gtorusdoublenullCOV^{-1})_C^A (\gtorusdoublenullCOV^{-1})_D^B \gnulltori^{-1}\left(\mathrm{d} x^C,\mathrm{d} x^D \right).
		\label{E:SMOOTHTORUSINVERSEFIRSTFUNDCOMPONENTSINTERMSOFDOUBLENULLTORUSINVERSEFIRSTFUNDCOMPONENTS}
	\end{align}
	\end{subequations}
	In addition, the inverse $(\gtorusdoublenullCOV^{-1})_A^B$ of $\gtorusdoublenullCOV_B^A$, defined by
	$(\gtorusdoublenullCOV^{-1})_B^C \gtorusdoublenullCOV_C^A = \updelta_B^A$, can be expressed as follows:
	\begin{align} \label{E:INVERSECHOVCOEFFICIENTSSMOOTHANGULARDERIVATIVESINTERMSOFDOUBLENULLONESANDL}
		(\gtorusdoublenullCOV^{-1})_A^B 
		& 
		= 
		\updelta_A^B 
		+ 
		\frac{\geop{x^A}\ubar}{\geop{t} \ubar} \Lunit^B.
	\end{align}
	Finally, we have the following identities:
	\begin{subequations}
		\begin{align}
			\gnulltori\left(\nullgeop{x^A},\nullgeop{x^B} \right) & = (\gtorusdoublenullCOV^{-1})_A^C (\gtorusdoublenullCOV^{-1})_B^D \gtorus_{CD}, \label{E:DOUBLENULLFIRSTFUNDCOMPONENTSINTERMSOFSMOOTHTORUSCOMPONENTS} \\
			 \gnulltori^{-1}\left(\mathrm{d} x^A,\mathrm{d} x^B \right) & = \gtorusdoublenullCOV_C^A \gtorusdoublenullCOV_D^B (\gtorus^{-1})^{CD}. \label{E:DOUBLENULLINVERSEFIRSTFUNDCOMPONENENTSINTERMSOFSMOOTHTORUSCOMPONENETS}
		\end{align}
	\end{subequations}
\end{lemma}
\begin{proof}
\eqref{E:DOUBLENULLTORIFIRSTFUNDRESTRICTEDTODOUBLENULLTORI} is a simple consequence of the fact that $(x^2,x^3)$ are coordinates on  $\doublenulltoritwoarg{\ubar}{u}$. The identity \eqref{E:GNULLTORICOMPONENTS} follows from $\gnulltori(\nullgeop{x^A},\nullgeop{x^B}) = \gfour(\nullgeop{x^A},\nullgeop{x^B})$, \eqref{E:NULLCOORDINATEPARTIALXAINTERMSOFGEOMETRICVF}, and the fact that $\gfour(\Lunit,\Lunit) = \gfour(\Lunit,\nullgeop{x^A}) = 0$. The identity \eqref{E:GNULLTORIINVERSECOMPONENDS} is a simple consequence of the fact that the coordinate vectorfields $\{\nullgeop{x^2},\nullgeop{x^3}\}$ span the tangent space of $\doublenulltoritwoarg{\ubar}{u}$. Next, we note that it is straightforward to check, using \eqref{E:INVERSEFIRSTFUNDFORMOFDOUBLENULLTORI}, that the type $\binom{1}{1}$ tensorfield with components $\gnulltori^{-1})^{\alpha\gamma} \gnulltori_{\gamma\beta}$ is precisely the $\gfour$-orthogonal projection tensorfield \eqref{E:DOUBLENULLTORIPROJECTIONDEFININGEQUATION}. From this fact and the fact that $(x^2,x^3)$ are coordinates on $\doublenulltoritwoarg{\ubar}{u}$, \eqref{E:DOUBLENULLTORIMETRICCOMPONENTSANDTHEINVERSECOMPONENTSRELATION} readily follows. 

The identity \eqref{E:SMOOTHANGULARDERIVATIVESINTERMSOFDOUBLENULLONESANDL} is a restatement of \eqref{E:SMOOTHANGULARVECTORFIELDSINTERMSOFNULLONES}. To derive \eqref{E:INVERSECHOVCOEFFICIENTSSMOOTHANGULARDERIVATIVESINTERMSOFDOUBLENULLONESANDL}, we note that \eqref{E:SMOOTHANGULARDERIVATIVESINTERMSOFDOUBLENULLONESANDL} implies $\nullgeop{x^A} = (\gtorusdoublenullCOV^{-1})_A^B \geop{x^B} + f L$ for some scalar function $f$. We then note that the coefficients of $(\gtorusdoublenullCOV^{-1})_A^B$ are given by \eqref{E:NULLCOORDINATEPARTIALXAINTERMSOFGEOMETRICVF}. Next, using \eqref{E:INVERSECHOVCOEFFICIENTSSMOOTHANGULARDERIVATIVESINTERMSOFDOUBLENULLONESANDL}, we see that \eqref{E:DOUBLENULLFIRSTFUNDCOMPONENTSINTERMSOFSMOOTHTORUSCOMPONENTS} follows from \eqref{E:GNULLTORICOMPONENTS}. \eqref{E:SMOOTHTORUSFIRSTFUNDCOMPONENTSINTERMSOFDOUBLENULLTORUSFIRSTFUNDCOMPONENTS} follows from applying two factors of $\gtorusdoublenullCOV$ to each side of \eqref{E:DOUBLENULLFIRSTFUNDCOMPONENTSINTERMSOFSMOOTHTORUSCOMPONENTS}. \eqref{E:SMOOTHTORUSINVERSEFIRSTFUNDCOMPONENTSINTERMSOFDOUBLENULLTORUSINVERSEFIRSTFUNDCOMPONENTS} follows from taking the inverse of \eqref{E:SMOOTHTORUSFIRSTFUNDCOMPONENTSINTERMSOFDOUBLENULLTORUSFIRSTFUNDCOMPONENTS} and using \eqref{E:DOUBLENULLTORIMETRICCOMPONENTSANDTHEINVERSECOMPONENTSRELATION}. Finally, \eqref{E:DOUBLENULLINVERSEFIRSTFUNDCOMPONENENTSINTERMSOFSMOOTHTORUSCOMPONENETS} follows from applying two factors of $\gtorusdoublenullCOV$ to \eqref{E:SMOOTHTORUSINVERSEFIRSTFUNDCOMPONENTSINTERMSOFDOUBLENULLTORUSINVERSEFIRSTFUNDCOMPONENTS}.
\end{proof}

\subsection{The pointwise semi-norms of tensors with respect to $\gnulltori$ and the $\gnulltori$-trace}
\label{SS:POINTWISESEMINORMSOFTENSORSONDOUBLENULLTORI}

\begin{definition}[Pointwise norms] \label{D:POINTWISESEMINORMWITHRESPECTTOFIRSTFUNDOFDOUBLENULLTORI} 
If $\upxi$ is a type $\binom{m}{n}$ tensorfield, then we define $|\upxi|_{\gnulltori} \geq 0$ by:
\begin{align} \label{E:SQUAREPOINTWISESEMINORMWITHRESPECTTOFIRSTFUNDOFDOUBLENULLTORI}
|\upxi|_{\gnulltori}^2 
& \eqdef 
\gnulltori_{\alpha_1 \widetilde{\alpha}_1} \cdots 
\gnulltori_{\alpha_m \widetilde{\alpha}_m} 
(\gnulltori^{-1})^{\beta_1 \widetilde{\beta}_1} 
\cdots 
(\gnulltori^{-1})^{\beta_n \widetilde{\beta}_n} 
\upxi_{\beta_1 \cdots \beta_n}^{\alpha_1 \cdots \alpha_m} 
\upxi_{\widetilde{\beta}_1 \cdots \widetilde{\beta}_n}^{\widetilde{\alpha}_1\cdots\widetilde{\alpha}_m}.
\end{align}
\end{definition}

\begin{definition} [$\gnulltori$-trace]
\label{D:TRACEOFDOUBLENULLTORITANGENT02TENSORS}
If $\upxi$ is a type $\binom{0}{2}$ tensorfield, 
then we define its $\gnulltori$-trace $\mytr_{\gnulltori} \upxi$
as follows:
\begin{align} \label{E:TRACEOFNULLTORITANGENT02TENSORS}
\mytr_{\gnulltori} \upxi
& 
\eqdef 
(\gnulltori^{-1})^{\alpha \beta}
\upxi_{\alpha \beta}.
\end{align}
\end{definition}

\subsection{Identities involving the double-null frame and the null-geometry scalars} \label{SS:IDENTITIESFORDOUBLENULLFRAMEANDACOUSTICGEOMETRY}

\begin{lemma}[Identities involving the double-null frame and the null-geometry scalars] \label{L:IDENTITIESFORDOUBLENULLFRAMEANDACOUSTICGEOMETRY} Let $\ReciprocalLunitAppliedtoTimeFunction, \ReciprocaluLunitAppliedtoTimeFunction, \ingoingmu, \MagnitueofinnerproductofnewLandnewuL$ be the double-null acoustic geometry scalar functions defined in Sect.\,\ref{SSS:NULLVECTORFIELDS}.  Then the following identities hold:
	\begin{subequations}
		\begin{align}
			\muX \ubar & = \frac{\upmu}{\ingoingmu} - \frac{\upmu}{\ReciprocalLunitAppliedtoTimeFunction}, \label{E:MUXACTINGONINGOINGEIKONAL} \\
			0  & = \frac{\ReciprocaluLunitAppliedtoTimeFunction}{\ReciprocalLunitAppliedtoTimeFunction} + \frac{\upmu}{\ingoingmu} - \frac{\upmu}{\ReciprocalLunitAppliedtoTimeFunction}  - \MagnitueofinnerproductofnewLandnewuL |\angD \ubar|^2_{\gtorus}, \label{E:IDENTITYFORNULLSCALARS} \\
			\ReciprocaluLunitAppliedtoTimeFunction & = \frac{1}{2}\upmu + \frac{1}{2}\ReciprocalLunitAppliedtoTimeFunction^2 \upmu |\angD \ubar|^2_{\gtorus}\label{E:IDENTITYFORRECIPROCALULUNITAPPLIEDTOTIMEFUNCTION}.
		\end{align}
	\end{subequations}
	
	Moreover, the vectorfields $\uLunit, \newuL$ can be decomposed in terms of $L,\muX$, and $\geop{x^A}$ as follows, where $\angD^A = \gtorus^{-1}(\mathrm{d} x^A,\mathrm{d} x^B) \geop{x^B}$:
	\begin{align}
		\uLunit = \Lunit + \frac{1}{\ReciprocaluLunitAppliedtoTimeFunction} \muX - \ingoingmu \angD^A \ubar \geop{x^A} & & \newuL = \ReciprocaluLunitAppliedtoTimeFunction \Lunit + \muX- \MagnitueofinnerproductofnewLandnewuL \angD^A \ubar \geop{x^A}.  \label{E:CONVENIENTIDENTITYFORULUNIT} 
	\end{align}	Next, define:  
	\begin{equation}
		\ToriTangentVectorfieldAssociatedToDoubleNullFolliations \eqdef \angD^A \ubar \gtorusdoublenullCOV^B_A \nullgeop{x^B}, \label{E:TORITANGENTVECTORFIELDASSOCIATEDTODOUBLENULLFRAME}
	\end{equation}
	where  $\gtorusdoublenullCOV^B_A$ is as in \eqref{E:CHOVCOEFFICIENTSSMOOTHANGULARDERIVATIVESINTERMSOFDOUBLENULLONESANDL}. Then the following identities hold:
	\begin{align}
		\newuL & = -\frac{\upmu}{\ingoingmu} \newL  + \upmu \Transport - \MagnitueofinnerproductofnewLandnewuL \ToriTangentVectorfieldAssociatedToDoubleNullFolliations = -\frac{\MagnitueofinnerproductofnewLandnewuL}{\ingoingmu} \Lunit  + \upmu \Transport - \MagnitueofinnerproductofnewLandnewuL \ToriTangentVectorfieldAssociatedToDoubleNullFolliations,  \label{E:CONVENIENTIDENTITYFORNEWUL} \\
		& = \left( -\frac{\MagnitueofinnerproductofnewLandnewuL}{\ingoingmu} + \upmu\right)\Lunit + \muX - \MagnitueofinnerproductofnewLandnewuL \ToriTangentVectorfieldAssociatedToDoubleNullFolliations. \label{E:NEWULWITHRESPECTTOLOLDULUNIT}
	\end{align}
In addition, we have the following for $\alpha = 0,1,2,3$:
	\begin{subequations}
		\begin{align} \label{E:DECOMPOOFPARTIALALPHAINTODOUBLENULLFRAME}
			\p_\alpha & = -\frac{1}{\MagnitueofinnerproductofnewLandnewuL} \newuL_\alpha \newL - \frac{1}{\MagnitueofinnerproductofnewLandnewuL} \newL_\alpha \newuL  + \nullangpartial_\alpha, \\
			& =  -\frac{1}{\MagnitueofinnerproductofLunitanduLunit} \uLunit_\alpha \Lunit - \frac{1}{\MagnitueofinnerproductofLunitanduLunit} \Lunit_\alpha \uLunit  + \nullangpartial_\alpha, \label{E:DECOMPOOFPARTIALALPHAINTOUNITDOUBLENULLFRAME}
		\end{align}
	\end{subequations}
where $\nullangpartial_\alpha \eqdef \nulltorusproject \cdot \partial_\alpha$.\footnote{Note that $\nullangpartial_{A} = \nullgeop{x^A}$ for $A = 2,3$.}
\end{lemma}

\begin{proof}
	Since $\muX = \upmu \Transport - \upmu L$, \eqref{E:MUXACTINGONINGOINGEIKONAL} follows from \eqref{E:ULFOLIATIONDENSITY} and \eqref{E:DOUBLENULLCARTESIANTIMENORMALIZEDNULLVECTORFIELDS}.
	
	Next we prove \eqref{E:CONVENIENTIDENTITYFORULUNIT}. Since $\{L,\muX,\geop{x^2},\geop{x^3}\}$ span the tangent space to any given point in spacetime, we may write $\uLunit = \alpha \Lunit + \beta \muX + c^A \geop{x^A}$. Acting both sides of this expression on $t$, we find that $\alpha = 1$ by \eqref{E:LUNITANDULUNITAPPLIEDTOEIKONALANDCARTESIANTIME} and Lemma\,\ref{L:BASICPROPERTIESOFVECTORFIELDS}. Acting both sides on $u$ proves that $\beta = \frac{1}{\ReciprocaluLunitAppliedtoTimeFunction}$, see again \eqref{E:LUNITANDULUNITAPPLIEDTOEIKONALANDCARTESIANTIME}. Finally, taking the inner product with $\geop{x^B}$, using the definition \eqref{E:DOUBLENULLCARTESIANTIMENORMALIZEDNULLVECTORFIELDS} of $\uLunit$, and using that $\geop{x^B}$ is $\gfour$-orthogonal to $\Lunit$ and $\muX$ proves $c^A = - \ingoingmu  \gtorus^{-1}(\mathrm{d} x^A,\mathrm{d} x^B) \geop{x^B} \ubar$ and so the first identity in \eqref{E:CONVENIENTIDENTITYFORULUNIT} is proved. The second identity in \eqref{E:CONVENIENTIDENTITYFORULUNIT} follows from multiplying the first identity by $\ReciprocaluLunitAppliedtoTimeFunction$. 
	
	Next, using $\newuL \ubar = 0$ and the now proved \eqref{E:CONVENIENTIDENTITYFORULUNIT},  \eqref{E:MUXACTINGONINGOINGEIKONAL} shows \eqref{E:IDENTITYFORNULLSCALARS}. The now proved \eqref{E:IDENTITYFORNULLSCALARS} implies \eqref{E:IDENTITYFORRECIPROCALULUNITAPPLIEDTOTIMEFUNCTION} from the identities $\frac{{\ReciprocaluLunitAppliedtoTimeFunction}}{{\ReciprocalLunitAppliedtoTimeFunction}} = \frac{\upmu}{{\ingoingmu}}$ and ${\MagnitueofinnerproductofnewLandnewuL} = \upmu {\ReciprocalLunitAppliedtoTimeFunction}$ and straightforward algebraic manipulation.

	To prove \eqref{E:CONVENIENTIDENTITYFORNEWUL}, we expand $\geop{x^A}$ in \eqref{E:CONVENIENTIDENTITYFORULUNIT} using  \eqref{E:SMOOTHANGULARDERIVATIVESINTERMSOFDOUBLENULLONESANDL}, add and subtract a copy of $\upmu \Lunit = \frac{\upmu}{\ReciprocalLunitAppliedtoTimeFunction}\newL$, and simplify with \eqref{E:LUNITANDULUNITAPPLIEDTOEIKONALANDCARTESIANTIME} and \eqref{E:IDENTITYFORNULLSCALARS}. Moreover, \eqref{E:NEWULWITHRESPECTTOLOLDULUNIT} follows from \eqref{E:CONVENIENTIDENTITYFORNEWUL} and $\upmu \Transport =\upmu \Lunit + \muX $.

	Finally, \eqref{E:DECOMPOOFPARTIALALPHAINTODOUBLENULLFRAME} follows immediately from the definition of \eqref{E:DOUBLENULLTORIPROJECTIONDEFININGEQUATION}. Identity  \eqref{E:DECOMPOOFPARTIALALPHAINTOUNITDOUBLENULLFRAME} follows from \eqref{E:DEFOFINNERPRODUCTOFLUNITANDULUNIT} and \eqref{E:RELATIONBETWEENCARTESIANNORMALIZEDNULLVECTORFIELDSANDEIKONALFUNCTIONORMALIZEDNULLVECTORFIELDS} applied to \eqref{E:DECOMPOOFPARTIALALPHAINTODOUBLENULLFRAME}.
	
\end{proof}

\begin{remark}[Equivalent ways of writing the inverse foliation density and the eikonal equation for $\ubar$]
We highlight that the identities \eqref{E:MUXACTINGONINGOINGEIKONAL}--\eqref{E:MUXACTINGONINGOINGEIKONAL}, which could appear unmotivated at first glance, are simply an equivalent formulation of the definition of $\ingoingmu$ and the eikonal equation for $\ubar$.  Indeed, the first identity is equivalent to $\ingoingmu = \frac{1}{\Transport \ubar}$ since $\Transport = L + X = L + \frac{1}{\upmu} \muX$. Similarly, using \eqref{E:SMOOTHTORUSMETRICINTERMSOFSIGMATMETRICANDX}, the eikonal equation is equivalent to $0 = - \upmu (\Lunit \ubar)^2 - 2 \muX \ubar \Lunit \ubar + \upmu |\angD \ubar|_{\gtorus}^2$ and hence: 
\begin{equation} \label{E:MUXINGOINGEIKONALDETERMINEDBYOTHERDERIVATIVES}
	\muX {\ubar} = -\frac{\upmu}{2} \Lunit {\ubar} + \frac{\upmu}{2 \Lunit {\ubar} }|\angD {\ubar}|^2_{\gtorus}.
\end{equation} 
This, in turn, implies: 
\begin{equation} \label{E:INGOINGMUDETERMINEDBYEVERYTHINGELSE}
	\frac{1}{{\ingoingmu}} = \frac{1}{2 {\ReciprocalLunitAppliedtoTimeFunction}}  + \frac{{\ReciprocalLunitAppliedtoTimeFunction}}{2} |\angD {\ubar}|^2_{\gtorus} =  \frac{1}{2} \Lunit {\ubar} + \frac{1}{2 \Lunit {\ubar} }|\angD {\ubar}|^2_{\gtorus}, 
\end{equation}
which is equivalent to \eqref{E:IDENTITYFORNULLSCALARS} using $\frac{{\ReciprocaluLunitAppliedtoTimeFunction}}{{\ReciprocalLunitAppliedtoTimeFunction}} = \frac{\upmu}{{\ingoingmu}}$ and ${\MagnitueofinnerproductofnewLandnewuL} = \upmu {\ReciprocalLunitAppliedtoTimeFunction}$.

\end{remark}

\subsection{Identities involving $\Transport$ and the double-null frame}

\begin{lemma} \label{L:BINTERMSOFDOUBLENULLFRAME}
		The material vectorfield $\Transport$ admits the following decompositions in terms of the double-null frame and the vectorfield from \eqref{E:TORITANGENTVECTORFIELDASSOCIATEDTODOUBLENULLFRAME}:
		\begin{subequations}
			\begin{align}
				\Transport &= \frac{1}{\upmu} \newuL + \frac{1}{\ingoingmu} \newL + \ReciprocalLunitAppliedtoTimeFunction \ToriTangentVectorfieldAssociatedToDoubleNullFolliations  \label{E:BINTERMSOFDOUBLENULLFRAME} \\
				& = \frac{1}{\MagnitueofinnerproductofLunitanduLunit} \uLunit + \frac{1}{\MagnitueofinnerproductofLunitanduLunit} \Lunit + \breve{\slashed{\Transport}}. \label{E:BINTERMSOFUNITDOUBLENULLFRAME}
			\end{align}
		\end{subequations}
In particular, we have: 
	\begin{align}
		\ReciprocalLunitAppliedtoTimeFunction \ToriTangentVectorfieldAssociatedToDoubleNullFolliations  = \breve{\slashed{\Transport}} \label{E:TORITANGENTVECTORFIELDASSOCIATEDTODOUBLENULLISPROJECTIONOFTRANSPORT}
	\end{align}In addition, there exists unique $\Sigma_t$-tangent vectorfields $\Lunit^\perp$ and $\uLunit^\perp$ such that: 
	\begin{subequations}
		\begin{align}
			\Lunit & = \Transport + \Lunit^\perp, \label{E:LUNITINTERMSOFBANDSIGMATTAN} \\
			\uLunit & = \Transport + \uLunit^\perp. \label{E:ULUNITINTERMSOFBANDSIGMATTAN}
		\end{align}
	\end{subequations}
Moreover, $\Lunit^\perp$ and $\uLunit^\perp$ satisfy:
	\begin{align}
		 \gfour(\Lunit^\perp,\Lunit^\perp) & = 1, & \gfour(\uLunit^\perp,\uLunit^\perp) & = 1, & \gfour(\Lunit^\perp,\uLunit^\perp) & = 1 - \MagnitueofinnerproductofLunitanduLunit. \label{E:SIZEOFLUBITPERPANDULUNITPERPANDTHEIRINNERPRODUCT}
	\end{align}\end{lemma}

	\begin{proof}
		The identity, \eqref{E:BINTERMSOFDOUBLENULLFRAME} follows from solving for $\Transport$ in \eqref{E:CONVENIENTIDENTITYFORNEWUL} and using \eqref{E:RATIOOFNULLGEOSICINNERPRODUCTANDFOLIATIONDENSITY}. The related identity \eqref{E:BINTERMSOFUNITDOUBLENULLFRAME} follows from \eqref{E:EQUIVALENTLUNITANDULUNITAPPLIEDTOEIKONALANDCARTESIANTIME} and \eqref{E:DOUBLENULLTORIPROJECTIONDEFININGEQUATIONWRTUNITL}. Next, \eqref{E:TORITANGENTVECTORFIELDASSOCIATEDTODOUBLENULLISPROJECTIONOFTRANSPORT} follows from $\frac{1}{\MagnitueofinnerproductofLunitanduLunit} \uLunit = \frac{1}{\upmu} \newuL$ and $\frac{1}{\MagnitueofinnerproductofLunitanduLunit}\Lunit = \frac{1}{\ingoingmu} \newL$.
		
		The existence and uniqueness of $\Lunit^\perp$ and $\uLunit^\perp$ follow from the fact that $\Transport t = \Lunit t  = \uLunit t = 1$ and that, for any $p \in \Sigma_t$,  $\Span\{\Lunit\} \oplus T_p \Sigma_t$ and $\Span\{\uLunit\} \oplus T_p\Sigma_t$ span the entire spacetime tangent space to $p$.
		
		To prove the first two identities in \eqref{E:SIZEOFLUBITPERPANDULUNITPERPANDTHEIRINNERPRODUCT}, we note $\gfour(\Lunit^\perp,\Lunit^\perp) = \gfour(\Lunit - \Transport, \Lunit-\Transport) = -2 \gfour(\Lunit,\Transport) + \gfour(\Transport,\Transport) = 1$. An identical argument holds for $\uLunit^\perp$.
		
		The final identity in \eqref{E:SIZEOFLUBITPERPANDULUNITPERPANDTHEIRINNERPRODUCT} follows from nearly identical arguments supplemented by \eqref{E:INNERPRODUCTOFCARTESIANNORMALIZEDNULLVECTORFIELDS}.
	\end{proof}

\subsection{A convenient identity for $\nullgeop{u}$} \label{SS:CONVENIENTIDENTITYFORDOUBLENULLUDERIVATIVE}

\begin{lemma}[An identity for $\nullgeop{u}$]
		\label{L:DOUBLENULLUDERIVATIVEINTERMSOFGEOMETRICVECTORFIELDSANDDERIVATIVESOFUBAR}
	The following identity holds: 
	\begin{align} \label{E:DOUBLENULLUDERIVATIVEINTERMSOFGEOMETRICVECTORFIELDSANDDERIVATIVESOFUBAR}
	\nullgeop{u} 
		& = \newuL 
				- \left\{ \ReciprocaluLunitAppliedtoTimeFunction L^A + \upmu X^A - \MagnitueofinnerproductofnewLandnewuL \angD^A \ubar\right\} \left\{ \geop{x^A} - \frac{\geop{x^A} \ubar}{\geop{t} \ubar} \geop{t} \right\}
		\end{align}
\end{lemma}

\begin{proof}
The identity follows from decomposing $\nullgeop{u} = \newuL - \newuL^A \nullgeop{x^A}$ using \eqref{E:NEWULINTERMSOFNULLCOORDINATEPARTIALS} as well as using identity \eqref{E:NULLCOORDINATEPARTIALXAINTERMSOFGEOMETRICVF}.
\end{proof}

\subsection{Some geometric differential operators in terms of the double-null frame} \label{SS:SOMEDIFFERENTIALOPERTORSINTHEDOUBLENULLFRAME}

We will encounter several differential operators in our $L^2$ analysis whose decomposition relative to the double-null frame is useful to have. 

\begin{lemma}[Decomposition of some differential operators relative to the double-null frame] \label{L:MOREUSEFULFRAMEDECOMPOSITIONS}
Let $\multipliervectorfield$ be the multiplier vectorfield from \eqref{E:MULTIPLIERVECTORFIELD}, and let $\ToriTangentVectorfieldAssociatedToDoubleNullFolliations$ be the $\doublenulltoritwoarg{\ubar}{u}$-tangent vectorfield from \eqref{E:TORITANGENTVECTORFIELDASSOCIATEDTODOUBLENULLFRAME}. Then the following identity holds: 
	\begin{align}
		\multipliervectorfield = 2  \newuL + \left( 1 + 2\frac{ \MagnitueofinnerproductofnewLandnewuL }{\ingoingmu}\right) \Lunit + 2  \MagnitueofinnerproductofnewLandnewuL  \ToriTangentVectorfieldAssociatedToDoubleNullFolliations. \label{E:MULTIPLIERINTERMSOFDOUBLENULLFRAME}
	\end{align}In addition, the following identity holds for $\Yvf{A}$:
	\begin{align}
		\Yvf{A} = \frac{1}{\Speed^2} (\gtorus^{-1})^{AB} \gtorusdoublenullCOV_B^C \nullgeop{x^C} + \frac{\Yvf{A} \ubar}{\Lunit \ubar} \Lunit.  \label{E:YCOMMUTATORINTERMSOFNULLCOORDINATEVECTORFIELDSANDL}
	\end{align}
	
Finally set $\angularcoeffinmuweightedspatialcartesian_i^A$ to be the following $3\times 2$ matrix:
	\begin{align} \label{E:ANGULARCOEFFINMUWEIGHTEDCARTESIAN}
		\angularcoeffinmuweightedspatialcartesian_i^A \eqdef  \frac{-1}{X^1}\upepsilon_{i23} X^A + \updelta_i^A, & & i = 1,2,3, \quad A = 2,3,
	\end{align}
where $\upepsilon_{ijk}$ is the fully antisymmetric tensor normalized by $\upepsilon_{123} = 1$. Then the following identities hold:
\begin{align}
		\upmu \p_i & = \frac{X^i}{\Speed^2}\newuL + \left\{ \left(\frac{\MagnitueofinnerproductofnewLandnewuL}{\ingoingmu} - \upmu\right)\frac{X^i}{\Speed^2} + \upmu \angularcoeffinmuweightedspatialcartesian_i^A \frac{\Yvf{A} \ubar}{\Lunit \ubar}\right\} \Lunit + \frac{X^i}{\Speed^2} \MagnitueofinnerproductofnewLandnewuL \ToriTangentVectorfieldAssociatedToDoubleNullFolliations + \frac{\upmu}{\Speed^2}\angularcoeffinmuweightedspatialcartesian_i^A (\gtorus^{-1})^{AB}  \gtorusdoublenullCOV_B^C \nullgeop{x^C}. \label{E:MUWEIGHTEDCARTESIANINTERMSOFDOUBLENULLFRAME}
	\end{align}
\end{lemma}

\begin{proof}
The expression for $\multipliervectorfield$ in \eqref{E:MULTIPLIERINTERMSOFDOUBLENULLFRAME} follows from \eqref{E:NEWULWITHRESPECTTOLOLDULUNIT}.

Using \eqref{E:SMOOTHGINVERSEABEXPRESSION} and \eqref{E:Y2INTERMSOFGEOMETRICCOORDINATEVECTORFIELDS}--\eqref{E:Y3INTERMSOFGEOMETRICCOORDINATEVECTORFIELDS}, it follows that $\Yvf{A} = \frac{1}{\Speed^2} (\gtorus^{-1})^{AB} \geop{x^B}$. The desired identity \eqref{E:YCOMMUTATORINTERMSOFNULLCOORDINATEVECTORFIELDSANDL} follows from this and \eqref{E:SMOOTHANGULARDERIVATIVESINTERMSOFDOUBLENULLONESANDL}. 

To prove \eqref{E:MUWEIGHTEDCARTESIANINTERMSOFDOUBLENULLFRAME}, we first note that Lemma\,\ref{L:RELATIONSHIPBETWEENCARTESIANPARTIALDERIVATIVESANDSMOOTHGEOMETRICCOMMUTATORS} implies the following expression for the $\upmu$-weighted Cartesian spatial derivatives:
	\begin{align} \label{E:EXPRESSIONFORMUWEIGHTEDSPATIALCARTESIANDERIVATIVESWRTCOMMUTATORS}
		\upmu \p_i = \frac{X^i}{\Speed^2} \muX + \upmu \angularcoeffinmuweightedspatialcartesian_i^A \Yvf{A}.
	\end{align}
The result follows from \eqref{E:NEWULWITHRESPECTTOLOLDULUNIT} and the now proved \eqref{E:YCOMMUTATORINTERMSOFNULLCOORDINATEVECTORFIELDSANDL}.
\end{proof}

\subsection{The Levi-Civita connection $\nullangD$ of $\gnulltori$ and related differential operators}
\label{SS:CONNECTIONSANDDIFFERENTIALOPERATORSONROUGHTORI}

\begin{definition}[The Levi-Civita connection $\nullangD$ of $\gnulltori$ and related differential operators] 
\hfill
\label{D:CONNECTIONSANDDIFFERENTIALOPERATORSONNULLTORI}
\begin{enumerate}
\item We denote the Levi-Civita connection of $\gnulltori$ by $\nullangD$.
	In particular,
	for $\doublenulltoritwoarg{\ubar}{u}$-tangent tensorfields $\upxi$, 
	we have $\nullangD \upxi = \nulltorusproject \Dfour \upxi$.
\item If $\Xi$ is an $\doublenulltoritwoarg{\ubar}{u}$-tangent one-form, 
	then we define its $\doublenulltoritwoarg{\ubar}{u}$-divergence 
	to be the scalar function
	$\nullangDiv \Xi \eqdef \gnulltori^{-1} \cdot \nullangD \Xi$. 
	Similarly, if $V$ is an $\doublenulltoritwoarg{\ubar}{u}$-tangent vectorfield, 
	then we define its $\doublenulltoritwoarg{\ubar}{u}$-divergence to be the scalar function 
	$\nullangDiv V \eqdef \gnulltori^{-1} \cdot \nullangD V_{\flat}$,
	where $V_{\flat}$ is the one-form that is $\gfour$-dual to $V$.
\item If $\Xi$ is an $\doublenulltoritwoarg{\ubar}{u}$-tangent one-form, we define its curl to be the scalar function given by $\nullangcurl \eqdef \breve{\in}^{AB} \nullangD_A \Xi_B$, where $\breve{\in}$ is the area element of the Riemannian manifold $(\doublenulltoritwoarg{\ubar}{u},\gnulltori)$ with respect to an arbitrary frame $f_A, A = 2,3$.
\item If $\upxi$ is a symmetric type $\binom{0}{2}$ $\doublenulltoritwoarg{\ubar}{u}$-tangent tensorfield, 
	then we define its $\doublenulltoritwoarg{\ubar}{u}$-divergence $\nullangDiv \upxi$ 
	to be the $\doublenulltoritwoarg{\ubar}{u}$-tangent one-form with the following 
	$\doublenulltoritwoarg{\ubar}{u}$-components for $A=2,3$:
	\begin{align} \label{E:NULLANGDIVOFTYPE02NULLTORITANGENTTENSORFIELD}
	[\nullangDiv \upxi]\left(\nullgeop{x^A} \right) 
	\eqdef 
	(\gnulltori^{-1})(\rmd x^B,\rmd x^C) 
	\left[\nullangD_{\nullgeop{x^B}} \upxi \right]\left(\nullgeop{x^C},\nullgeop{x^A} \right).
	\end{align}
\end{enumerate}
\end{definition}

\subsection{Curvature tensors of $\gfour$ and $\gnulltori$}
\label{SS:CURVATURETENSORS}
The Riemann curvature tensors of $\gfour$ and $\gnulltori$ 
play a central role in the geometric analysis of the acoustic geometry. 


\begin{definition}[Curvature tensors of $\gfour$ and $\gnulltori$] \label{D:CURVATURETENSORS} The \emph{Riemann curvature tensor} $\Riemfour$ of the spacetime metric $\gfour$ is the type $\binom{0}{4}$ spacetime tensorfield 
defined by:
\begin{align} \label{E:ACOUSTICALCURVATURETENSOR}
\Riemfour(\mathbf{X},\mathbf{Y},\mathbf{Z},\mathbf{W}) 
& \eqdef \gfour(- \Dfour^2_{\mathbf{X} \mathbf{Y}} \mathbf{Z} 
+ 
\Dfour^2_{\mathbf{Y} \mathbf{X}} \mathbf{Z}, \mathbf{W}),
\end{align}
where $\mathbf{X},\mathbf{Y},\mathbf{Z},\mathbf{W}$ are arbitrary spacetime vectorfields, and
$\Dfour^2_{\mathbf{X} \mathbf{Y}} \mathbf{Z} 
\eqdef 
\mathbf{X}^{\alpha} \mathbf{Y}^\beta \Dfour_\alpha \Dfour_\beta \mathbf{Z}$. 

The \emph{Ricci curvature tensor} $\Ricfour$ of the spacetime metric $\gfour$ is the type $\binom{0}{2}$ spacetime tensor defined 
relative to arbitrary coordinates as follows: 
\begin{align} \label{E:SPACETIMERICCITENSOR}
\Ricfour_{\alpha \beta} 
& \eqdef 
(\gfour^{-1})^{\kappa \lambda} \Riemfour_{\alpha \kappa \beta \lambda}.
\end{align}Similarly, the Riemann curvature tensor $\Riemnulltori$ of the Riemannian metric $\gnulltori$ on 
$\doublenulltoritwoarg{\ubar}{u}$ is the type $\binom{0}{4}$ $\doublenulltoritwoarg{\ubar}{u}$ tensorfield defined 
as follows: 
\begin{align} \label{E:RIEMANNCURVATURETENSOROFDOUBLENULLTORUS}
\Riemnulltori(X,Y,Z,W) 
& \eqdef \gnulltori(- \nullangD^2_{XY}Z + \nullangD^2_{YX}Z,W),
\end{align}
where $X,Y,Z,W$ are arbitrary $\doublenulltoritwoarg{\ubar}{u}$-tangent vectorfields, $\nullangD$ is the Levi-Civita connection of $\gnulltori$, and $\nullangD_{XY}^2 Z= X^\alpha Y^\beta \nullangD_\alpha \nullangD_\beta Z$.

The \emph{Ricci curvature tensor} $\Ricnulltori$ of $\gnulltori$ is the type $\binom{0}{2}$ spacetime tensor defined 
relative to arbitrary coordinates as follows: 
\begin{align} \label{E:ROUGHTORUSRICCITENSOR}
\Ricnulltori_{\alpha \beta} 
& \eqdef 
(\gnulltori^{-1})^{\kappa \lambda} \Riemnulltori_{\alpha \kappa \beta \lambda}.
\end{align}The \emph{scalar curvature} $\Scalarnulltori$ of $\gnulltori$ 
is the scalar function defined 
relative to arbitrary coordinates as follows: 
\begin{align} \label{E:ROUGHTORUSSCALARCURVATURE}
\Scalarnulltori
& \eqdef 
(\gnulltori^{-1})^{\alpha \beta} \Ricnulltori_{\alpha \beta}.
\end{align}
\end{definition}

It is well-known that because $\doublenulltoritwoarg{\ubar}{u}$ is two-dimensional,
the \emph{Gauss curvature} $\Gaussnulltori$ of $\gnulltori$ can be expressed as follows
in terms of it scalar curvature:
\begin{align} \label{E:2DGAUSSCURVATUREISTWICESCALARCURVATURE}
	\Gaussnulltori
	& = \frac{1}{2} \Scalarnulltori.
\end{align}\section{Deformation tensors of the acoustic and double-null geometry} \label{S:DEFORMATIONTENSORSOFACOUSTICANDDOUBLENULLGEOMETRY}

In our analysis, we encounter the deformation tensors of various vectorfields that we now introduce.

\subsection{Spacetime and $\ell_{t,u}$-tangent deformation tensors}
\label{SS:DEFORMATIONTENSORS}

\begin{definition}[Deformation tensors] 
Let $Z$ be a spacetime vectorfield.
We define the \emph{deformation tensor} $\deform{Z}$  of $Z$ 
(with respect to $\gfour$)
to be the following symmetric type $\binom{0}{2}$ tensorfield:
\begin{align} \label{E:DEFORMATIONTENSORDEF}
\deformarg{Z}{\alpha}{\beta} 
& 
\eqdef \Lie_Z \gfour_{\alpha \beta} 
= 
\Dfour_\alpha Z_\beta + \Dfour_\alpha Z_\beta,
\end{align}
where the final equality in \eqref{E:DEFORMATIONTENSORDEF} follows from the torsion-free property of $\Dfour$. We define $\angdeform{Z} \eqdef \torusproject \deform{Z}$ as in Def.\,\ref{D:PROJECTIONTENSORFIELDSANDTANGENCYTOHYPERSURFACES}. 
\end{definition}

\begin{lemma}[Identities for $\ell_{t,u}$-projection of $\deform{\comder}$]\label{L:IDENTITIESFORFLATOTURSPROJECTIONOFDEFORM}
The following identities hold: 
\begin{align} \label{E:IDENTITIESFORFLATOTURSPROJECTIONOFDEFORM}
	\angLie_Z \gtorus = \angdeform{Z}, & & \angLie_Z \gtorus^{-1} = - \angdeform{Z}^{\#\#}.
\end{align}\end{lemma}

\begin{proof}
Applying $\Lie_Z$ to \eqref{E:SMOOTHTORUSMETRICINTERMSOFSIGMATMETRICANDX} yields $\deform{Z} = -( \Lie_Z \Transport)_\flat \otimes \Transport_\flat - \Transport_\flat\otimes( \Lie_Z \Transport)_\flat + (\Lie_Z X)_\flat \otimes X_\flat + X_\flat \otimes( \Lie_Z X)_\flat + \Lie_Z \gtorus$. Since $\torusproject \Transport = \torusproject X = 0$, the first identity in \eqref{E:IDENTITIESFORFLATOTURSPROJECTIONOFDEFORM} follows from \eqref{E:PROJECTEDLIEDERIVATIVES}. The second identity is a straightforward consequence of the identity for $\angLie_Z \gtorus$ and applying the Leibniz rule to $\gtorus_{AB} (\gtorus^{-1})^{BC} = \updelta_A^C$. 
\end{proof}

\begin{lemma}[The frame components of $\deform{Z}$] \label{L:ANGULARDEFORMATIONTENSORFRAMECOMPONENTS}
Let $\zetatan$, $\zetatrans$, $\angktan$, and $\angktrans$ be the $\ell_{t,u}$-tensorfields defined in \eqref{E:CONNECTIONCOEFFICIENTDECOMPOSITIONS}. Then the following identities for the frame components of the angular deformation tensor $\deform{Z}$ of $Z \in \comder = \{\Lunit,\muX,\Yvf{A},\Yvf{B}\}$:
\begin{subequations} \label{E:ANGDEFORMOFL}
	\begin{align}
		\deformarg{L}{L}{L} & = 0, & & \deformarg{L}{\muX}{X} = 2 \Lunit \upmu, & & \deformarg{L}{L}{\muX} = - \Lunit\upmu \\
		\angdeformarg{L}{L} & = 0, & & \angdeformarg{L}{\muX}  = \angrmD \upmu + 2 \zetatrans + 2 \upmu \zetatan,  & & \\
		\angdeform{L} & = 2 \upchi, 
	\end{align}
\end{subequations}\begin{subequations} \label{E:ANGDEFORMOFMUX} 
	\begin{align}
		\deformarg{\muX}{L}{L} & = 0, & & \deformarg{\muX}{\muX}{X}  = 2 \muX \upmu, & & \deformarg{\muX}{L}{\muX} = -  \muX \upmu, \\
		 \angdeformarg{\muX}{\Lunit} & = - \angrmD \upmu - 2 \zetatrans - 2 \upmu \zetatan, & & \angdeformarg{\muX}{\muX}  = 0 & & \\
		\angdeform{\muX} & = - 2 \upmu \upchi + 2 \angktrans + 2 \upmu \angktan,
	\end{align}
\end{subequations}
\begin{subequations}
	\begin{align}
		\begin{split} \label{E:ANGDEFORMOFY}
			\angdeform{\Yvf{A}} & = 2 \Speed^{-2} X^A \upchi + \angG \diamond \Yvf{A} \wavearray -  \Speed^{-2} X^A \angG \diamond \Lunit \wavearray \\
			& \ \ +  \Speed^{-2} X^A \angG_L \otimesarray \angrmD \wavearray + \Speed^{-2} X^A \angrmD \wavearray \otimesarray \angG_L  +  \Speed^{-2} X^A \angG_X \otimesarray \angrmD \wavearray + \Speed^{-2} X^A \angrmD \wavearray \otimesarray \angG_X.
		\end{split}
	\end{align}
\end{subequations}
\end{lemma}

\begin{proof}
A detailed proof in $2D$ was provided in \cite{jSgHjLwW2016}*{Lemma 2.18} in the case of a quasilinear wave equation. The minor modifications in $3D$ are minor and we omit the detials. 
\end{proof}

\subsection{Identities involving the double-null toroidal components of deformation tensors}
\label{S:ROUGHTORUSCOMPONENTSOFDEFORMATIONTENSORS}
In the next lemma, we derive simple identities relating the
$\doublenulltoritwoarg{\ubar}{u}$-components 
of a deformation tensor $\deform{Z}$ 
and the $\doublenulltoritwoarg{\ubar}{u}$-components of $\Lie_Z \gnulltori$.

\begin{lemma}[Relating the $\doublenulltoritwoarg{\ubar}{u}$ components of $\Lie_Z \gnulltori$ and $\deform{Z}$]
\label{L:LIEGNULLTORUS}
Let $Z$ be a spacetime vectorfield, and let
$\gnulltori$ be the first fundamental form of 
$\twoargroughtori{\timefunction}{u}$, as in \eqref{E:ACOUSTICALMETRICINTERMSOFNULLVECTORFIELDSANDDOUBLENULLSPHEREFIRSTFUND}. 
Then the following identities hold for $A,B=2,3$:
\begin{align} \label{E:LIEGNULLTORUS}
\deform{Z}\left(\nullgeop{x^A},\nullgeop{x^B}\right) 
& \eqdef
[\Lie_Z \gfour]\left(\nullgeop{x^A},\nullgeop{x^B}\right) 
=
[\Lie_Z \gnulltori]\left(\nullgeop{x^A},\nullgeop{x^B}\right).
\end{align}
\end{lemma}

\begin{proof}
First, using \eqref{E:ACOUSTICALMETRICINTERMSOFNULLVECTORFIELDSANDDOUBLENULLSPHEREFIRSTFUND} and the Leibniz rule for Lie differentiation,
we deduce:
\begin{align} 
\label{E:PROOFSTEPLIEGNULLTORUS}
\deform{Z}
\eqdef 
\Lie_{Z} \gfour
& = 
- 
 \frac{1}{\MagnitueofinnerproductofnewLandnewuL}  (\Lie_{Z} \newL_{\flat}) \otimes \newuL 
- 
 \frac{1}{\MagnitueofinnerproductofnewLandnewuL}  \newL_{\flat} \otimes \Lie_{Z}  \newuL_{\flat} - \frac{1}{\MagnitueofinnerproductofnewLandnewuL^2} \newuL(\MagnitueofinnerproductofnewLandnewuL) (\newL_{\flat} \otimes \newuL_{\flat} + \newuL_{\flat} \otimes \newL_{\flat}) + \Lie_{Z} \gnulltori,
\end{align}
where $\newL_{\flat}$ and $\newuL_{\flat}$ denote the one-form $\gfour$-dual to $\newL$ and $\newuL$, respectively. Contracting \eqref{E:PROOFSTEPLIEGNULLTORUS} against $\nullgeop{x^A} \otimes \nullgeop{x^B}$
and using that $\newL$ and $\newuL$ are $\gfour$-orthogonal to 
$\left\lbrace \nullgeop{x^2}, \nullgeop{x^3} \right\rbrace$, 
we conclude \eqref{E:LIEGNULLTORUS}.
\end{proof}

\section{Strings of differential operator and associated commutation formulas and identities}
\label{S:COMMUTATIONFORMULASIDENTITIESANDSTRINGSOFCOMMUTATORS}

\subsection{The basic commutator identity and strings of vectorfield differentiation}

The basic starting point is the following lemma.

\begin{lemma}\cite[Simple commutator identities; Lemma 8.5]{LukSpeck2024stability}
	\label{L:SIMPLECOMMUTATORIDENTITY}
	For any $\nullhyparg{u}$-tangent vectorfields 
	$\Singletan, \Singletan_1, \Singletan_2 \in \{\Lunit,\Yvf{2},\Yvf{3}\}$, 
	the commutators $[\Singletan_1,\Singletan_2]$ and $[\muX,\Singletan]$ are $\ell_{t,u}$-tangent.
	Moreover, we have the following identities, where $\smoothfunction$
	schematically denotes a smooth function of its arguments:
	\begin{subequations}
		\begin{align} \label{E:TANGENTIALSIMPLECOMMUTATORIDENTITY}
			[\Singletan_1, \Singletan_2] 
			& = 
			\smoothfunction(\tander^{\le 1}\controlvars) \Yvf{2} 
			+
			\smoothfunction(\tander^{\le 1}\controlvars) \Yvf{3}.
		\end{align}
		
		For each $\Singletan \in \{\Lunit,\Yvf{2},\Yvf{3}\}$,
		there exist smooth functions, all schematically denoted by ``$\smoothfunction$,''
		such that the following identity holds:
		\begin{align} \label{E:TRANSVERSALTANGENTIALSIMPLECOMMUTATORIDENTITY}
			[\Singletan, \muX] 
			& = 
			\smoothfunction(\tander^{\le 1} \badcontrolvars,\muX \wavearray) \Yvf{2}
			+
			\smoothfunction(\tander^{\le 1} \badcontrolvars,\muX \wavearray) \Yvf{3}.
		\end{align}
	\end{subequations}
\end{lemma}

\begin{definition}[Strings of commutation vectorfields and vectorfield seminorms] 
	\label{D:STRINGSOFCOMMUTATIONVECTORFIELDS}
	Recall that 
	$\Fullset = \{\Lunit, \muX, \Yvf{2},\Yvf{3}\}$, 
	$\Tanset = \{\Lunit,\Yvf{2},\Yvf{3}\}$,
	and 
	$\Angularset = \{\Yvf{2},\Yvf{3}\}$ are the sets of commutation vectorfields 
	defined in \eqref{E:COMMUTATIONVECTORFIELDS}.
	We adopt the following definitions for differential operators constructed out of strings of commutation vectorfields.
	In this section, $f$ denotes a scalar function.
	
	\begin{itemize}
		\item $\comder^{N;M} f$ denotes an arbitrary string of $N$ commutation vectorfields in $\comder$ applied to $f$, where the string contains \emph{precisely} $M$ factors of $\muX$. We set $\comder^{0,0}f = f$. We often write
		$\comder f$ instead of $\comder^{1;1} f$ and $\comder^{N;M} f$ when the string contains \emph{precisely} factors of $\muX$.
		\item $\tander^N f$ denotes an arbitrary string of $N$ commutation vectorfields in $\tander$ applied to $f$. We set $\tander^0 f = f$. We often write $\tander f$ instead of $\tander^1 f$.
		\item $\tanderY^N f$ denotes an arbitrary string of $N$ commutation vectorfields in $\mathscr{Y}$ applied to $f$. We set $\tanderY^0 f = f$.
		\item $\comdersmall^{N;M} f$ denotes an arbitrary string of $N$ commutation vectorfields in $\comder$ applied to $f$, where the string contains \emph{at least} one factor of $\tander$ and \emph{precisely} $M$ factors of $\muX$.
		\item $\tandersmall^N$ denotes an arbitrary string of $N$ commutation vectorfields in $\tander$ applied to $f$, where the string contains \emph{at least two factors} of $\Lunit$ or \emph{at least one factor of $\Yvf{2}, \Yvf{3}$}.
		\item $\comderdoublesmall^{N;M} f$ denotes an arbitrary string of $N$ commutation vectorfields in $\comder$ applied to $f$, 
		where the string contains at least two factors of $\Lunit$ or at least one factor of $\Yvf{A}$ and \emph{precisely} $M$ factors of $\muX$.
		\item $\mathfrak{P}^{(N)}$ denotes the set of all differential operators of the form $\tander^N$.
		\item $\mathfrak{Y}^{(N)}$ denotes the set of all differential operators of the form $\tanderY^N$.
		\item $\mathfrak{Z}_*^{(N;1)}$ and $\mathfrak{Z}_*^{(N;\le 1)}$ denote the set of all differential operators of the form $\comdersmall^{N;1}$ and $\comdersmall^{N;\le 1}$.
		\item  We define order $N$ strings of $\ell_{t,u}$-projected Lie derivatives
	such as $\angLie_{\tander}^N$ and $\angLie_{\comder}^{N;M}$ in an analogous fashion.  
	Such operators act on $\ell_{t,u}$-tangent tensorfields $\upxi$, e.g., $\angLie_{\tander}^N \upxi$.
		\item $\angLie_{\mathfrak{P}}^{(N)}$ denotes the set of all differential operators of the form $\angLie_{\tander}^N$.
		\item $\angLie_{\mathfrak{Y}}^{(N)}$ denotes the set of all differential operators of the form $\angLie_{\tanderY}^N$.
		\item  $\comder^{\leq N;M} f$ denotes the array of all terms of the form
		$\comder^{N';M} f$, where $0 \leq N' \leq N$. 
		\item If $N_1 < N_2$, then
		$\comder^{[N_1,N_2];M} f$ denotes the array of all terms of the form
		$\comder^{N';M} f$, where $N_1 \leq N' \leq N_2$.
		\item We define arrays such as $\tanderY^{\leq 2} f$ in an analogous fashion.
	\end{itemize}
	We also define pointwise seminorms:
	\begin{itemize}
		\item $|\comder^{N;M} f|$ denotes the magnitude of $\comder^{N;M}f$ as defined above (there is no summation).
		\item $|\comder^{\leq N;M} f|$ denotes the sum over all terms of the form $|\comder^{N';M} f|$ with $N' \leq N$. 
		\item $| \comder^{[N_1,N_2];M}f |$ is the sum over all terms of the form $|\comder^{N';M} f|$ with $N_1 \leq N' \leq N_2$. 
		\item Terms such as 
		$|\tandersmall^{[N_1,N_2]}f|$, 
		$|\angLie_{\comder}^{\leq N;M} \upxi|_{\gtorus}$, 
		$|\tanderY^{\leq N} f|$,
		etc., are defined analogously,
		e.g., $|\angLie_{\comder}^{\leq N;M} \upxi|_{\gtorus}$
		is the sum over all terms of the form
		$|\angLie_{\comder}^{N';M} \upxi|_{\gtorus}$
		with $N' \leq N$.
	\end{itemize}
	
	\begin{itemize}
		\item We will freely combine the above definitions with Def.\,\ref{D:CONVENTIONESFORDIFFERENTIATION}, e.g.,
		\begin{align} \label{E:NOTATIONFORSTRINGDERIVATIVESOFMULTIPLEARRAYS}
			|\tander^N \wavearray| 
			& \eqdef 
			\max_{k \in \{0,\cdots,3\}} | v^i|
		\end{align}
	\end{itemize}
\end{definition}

\section{Schematic structure capturing structure and schematic identities}
\label{S:SCHEMATICSTRUCTUREANDIDENTITIES}
In this section, we introduce schematic notation that will help us succinctly exhibit
the important qualitative features of various equations.
We then provide a collection of identities expressed in schematic form;
they will be helpful when we derive estimates.

\subsection{Some schematic notation}
\label{SS:SOMESCHEMATICANOTATION}
\begin{notation}[Schematic functional dependence] 
\label{N:SCHEMATICFUNCTIONALDEPENDENCE}
We often use the notation $\smoothfunction(\upxi_{(1)},\cdots,\upxi_{(m)})$ to schematically depict an expression 
(often tensorial and involving contractions) that depends smoothly on the $\ell_{t,u}$-tangent
tensorfields $\upxi_{(1)},\cdots,\upxi_{(m)}$. 
Note that in general, $\smoothfunction(0) \neq 0$. 
\end{notation}

\begin{notation}[Schematic use of the symbol $\tander$] 
\label{N:SCHEMATICUSEOFTANGENTIALDIFFERENTATION}
Throughout the rest of the paper, 
$\tander$ 
schematically denotes a differential operator that is tangent to 
the characteristics $\nullhyparg{u}$. 
For example, $\tander f$ might denote $\angrmd f, \Lunit f$, or $\Yvf{2}f$. 
We use such notation when the details of $\tander$ are unimportant. 
\end{notation}

We use the notation $\vec{x}$ to denote the array of spatial Cartesian coordinates,
i.e.,
\begin{align} \label{E:ARRAYOFOFCARTESIANSPATIALCOORDIANTES}
	\vec{x}
	& \eqdef (x^1, x^2, x^3).
\end{align}
We use the same conventions from Def.\,\ref{D:DERIVATIVESOFARRAYS} for
differential operators acting on $\vec{x}$, e.g.,
\begin{align} \label{E:ANDGULARDIFFERENTIALOFCARTESIANSPATIALCOORDIANTES}
	\angrmd \vec{x}
	& \eqdef (\angrmd x^1, \angrmd x^2, \angrmd x^3).
\end{align}
\subsection{Schematic structure of various tensorfields}
\label{SS:SCHEMATICSTRUCTUREOFVARIOUSTENSORFIELDS}

In this section we reveal schematic structure for many of the geometric objects that we will estimate in the paper. In Prop.\,\ref{P:SCHEMATICSTRUCTUREOFVARIOUSTENSORSINTERMSOFCONTROLVARS}, many of these structural dependencies feature both the array of controlling quantities $(\controlvars,\badcontrolvars)$ (recall \eqref{E:CONTROLVARSGOODANDBAD}) as well as $\LogDensity$. In Sect.\,\ref{S:DENSITYDERIVATIVESINTERMSOFOTHERS}, we will prove that all derivatives of $\LogDensity$ can be expressed in terms of $\velocityarray$, the specific vorticity $\vortrenormalized$, the gradient entropy $\GradEnt$, and lower order derivatives of $(\controlvars,\badcontrolvars)$. Hence, we only keep the dependence of $\LogDensity$ in Prop.\,\ref{P:SCHEMATICSTRUCTUREOFVARIOUSTENSORSINTERMSOFCONTROLVARS} for the sake of precison, but the reader should be aware that the bulk of the paper will never directly estimate $\LogDensity$.

\begin{proposition}[Schematic structure of various tensorfields] 
\label{P:SCHEMATICSTRUCTUREOFVARIOUSTENSORSINTERMSOFCONTROLVARS} 
Recall that $\controlvars$ and $\badcontrolvars$ are the arrays from Def.\,\ref{D:CONTROLVARS}.
The following schematic relations hold for scalar functions
($\alpha,\beta = 0,1,2,3$, $\kappa = 0,1,2,3,4$):
\begin{subequations} \label{E:SCHEMATICSTRUCTURALDEPENDENCE}
\begin{align}
\gfour_{\alpha \beta}, 
	\,
(\gfour^{-1})^{\alpha \beta}, 
	\,
\gtorus_{\alpha \beta}, 
	\,
(\gtorus^{-1})^{\alpha \beta}, 
	\,
G_{\alpha \beta}^{\kappa}, 
	\,
\smoothtorusproject_\beta^{\alpha}, 
	\,
\Lunit^{\alpha}, 
	\,
X^{\alpha}, 
	\,
\olduLunit^{\alpha}, 
	\,
\Yvf{2}^{\alpha}, 
	\,
\Yvf{3}^{\alpha}, 
& 
= \smoothfunction(\controlvars,\LogDensity),  
	\label{E:SCHEMATICSTRUCTUREOFMETRICETC} 
	\\
\left(\geop{t}\right)^{\alpha},
	\,
\left(\geop{x^2}\right)^{\alpha},
	\,
\left(\geop{x^3}\right)^{\alpha}
& 
= \smoothfunction(\controlvars,\LogDensity),
\label{E:SCHEMATICSTRUCTUREOFCARTESIANCOMPONENTSOFTANGENTIALGEOMETRICCOORDINATEVECTORFIELDS}
	\\
\left(\geop{u}\right)^{\alpha}
& 
= \smoothfunction(\badcontrolvars,\LogDensity),
	\label{E:SCHEMATICSTRUCTUREOFCARTESIANCOMPONENTSOFTRANSVERSALGEOMETRICCOORDINATEVECTORFIELD} 
		\\
G_{\Lunit \Lunit}^{\kappa}, 
	\,
G_{\Lunit X}^{\kappa}, 
	\,
G_{XX}^{\kappa} 
& = \smoothfunction(\controlvars,\LogDensity), 
	\label{E:SCHEMATICSTRUCTUREOFGFRAMESCALARS} \\
\Xsmall^{\alpha}, 
	\,
\Yvfsmall{2}^{\alpha},
	\,
\Yvfsmall{3}^{\alpha}, 
	\,
\Yvfsmallcoeff{A}
& 
= \smoothfunction(\controlvars,\LogDensity)(\controlvars,\LogDensity),
	\label{E:SCHEMATICSTRUCTUREOFXSMALL} \\
\muX^{\alpha}
& = \smoothfunction(\controlvars,\LogDensity)\upmu.
	\label{E:SCHEMATICSTRUCTUREOFMUX}
\end{align}
\end{subequations}Moreover, if $\varphi$ is a scalar function, then we have the following schematic relation for 
its $\angD$-Hessian (which is a symmetric type $\binom{0}{2}$-tangent tensorfield):
\begin{align} \label{E:ELLTUHESSIANOFSCALARFUNCTIONSCHEMATIC}
	\angDsquared \varphi
	& = \smoothfunction(\controlvars, \angrmd \vec{x} \, ) \tanderY^2 \varphi
			+
			\smoothfunction(\controlvars,\LogDensity, \angrmd \vec{x} \, ) \cdot \tanderY (\controlvars,\LogDensity) \cdot \tanderY \varphi.
\end{align}Finally, we have the following schematic relations for $\ell_{t,u}$-tangent tensorfields,
where $\angrmd \vec{x}$ is defined in \eqref{E:ANDGULARDIFFERENTIALOFCARTESIANSPATIALCOORDIANTES}:
\begin{subequations} \label{E:SCHEMATICSTRUCTURALDEPENDENCEOFTORITANGENTTENSORFIELDS}
\begin{align}
\gtorus, 
	\,
\angG_{\Lunit}, 
	\,
\angG_X,
	\,
\angG 
& = \smoothfunction(\controlvars,\LogDensity, \angrmd \vec{x} \, ), 
	\label{E:SCHEMATICSTRUCTUREOFGFRAMEELLTUTENSORS} 
		\\
\Yvf{2}, 
	\,
\Yvf{3}, 
& = \smoothfunction(\controlvars,\LogDensity, \gtorus^{-1}, \angrmd \vec{x} \, ), 
	\label{E:SCHEMTAICSTRUCTUREOFSMOOTHANGULARCOMMUTATORVECTORFIELDS} 
		\\
\upchi 
& = \smoothfunction(\controlvars,\LogDensity,  \angrmd \vec{x} \, ) \tander (\controlvars,\LogDensity), 
	\label{E:SCHEMATICSTRUCTUREOFNULLSECONDFUNDAMENTALFORM} 
		\\
\mytr_{\gtorus} \upchi 
& 
= 
\smoothfunction(\controlvars,\LogDensity, \gtorus^{-1}, \angrmd \vec{x} \, ) \tander  (\controlvars,\LogDensity), 
	\label{E:SCHEMATICSTRUCTUREOFTRACEOFNULLSECONDFUNDAMENTALFORM}  
		\\
\zetatan, 
	\,
\angktan 
& 
= \smoothfunction(\controlvars,\LogDensity, \angrmd \vec{x} \, ) \tander \wavearray, 
	\label{E:SCHEMATICSTRUCTUREOFGOODPARTSOFCONNECTIONCOEFFICIENTS} 
		\\
\zetatrans, 
	\,
\angktrans 
& = \smoothfunction(\controlvars,\LogDensity, \angrmd \vec{x} \, ) \muX \wavearray. \label{E:SCHEMATICSTRUCTUREOFTRANSVERSALDIFFERENTIATIONPARTSOFCONNECTIONCOEFFICIENTS}
\end{align}
\end{subequations}\end{proposition}

\begin{proof}
The proposition was proved in  \cite{abbrescia2022emergence}*{Prop.\,10.1}. However, \cite{abbrescia2022emergence} defined the array of solution variables to be $\controlvars = (\wavearray, \Lsmall^1, \Lsmall^2, \Lsmall^3)$ and $
\badcontrolvars = (\wavearray,\upmu - 1, \Lsmall^1, \Lsmall^2, \Lsmall^3)$, compared to \eqref{E:CONTROLVARSGOODANDBAD}. The only difference is that here we have redefined $\controlvars$ and $\badcontrolvars$ to not have any dependence on $\LogDensity$, which was extracted on the RHS\,\eqref{E:SCHEMATICSTRUCTURALDEPENDENCE}--\eqref{E:SCHEMATICSTRUCTURALDEPENDENCEOFTORITANGENTTENSORFIELDS}.
\end{proof}

\subsection{Structural properties of the $\gfour$-null forms}
\label{SS:TRANSVERSALDERIVATIVESINTERMSOFTANGENTIALONES}
In this section, we exhibit some crucial structural properties of the $\gfour$-null forms.

\begin{lemma}[Crucial structural properties of $\gfour$-null forms;\cite{LukSpeck2024stability}*{Lemma\,8.2}]
	\label{L:CRUCIALSTRUCTUREOFNULLFORMS}
	The product of
	$\upmu$ and the terms defined in \eqref{E:Q0DEF}--\eqref{E:QIJDEF} 
	satisfy the following pointwise estimates for any $\Psi_i, \Psi_j \in \wavearray = \{\RRiemann,\LRiemann,v^2,v^3,\Ent\}$:
\begin{align} 
\upmu \nullform_{(\gfour)}(\pmb\p \Psi_i,\pmb\p \Psi_j),
	\, 
\upmu \nullform_{(\alpha\beta)}(\pmb \p \Psi_i, \pmb \p \Psi_j) 
= \smoothfunction(\badcontrolvars, \comder \wavearray) 
	\cdot
	\tander \wavearray,
		\label{E:NULLFORMSTRUCTUREWAVEVARIABLES} 
	\end{align}
\end{lemma}

\subsection{Additional schematic identities involving differentiation}
\label{SS:ADDITIONALSCHEMATICIDENTITIESINVOLVINGDIFFERENTIATION}
For future use, in this section, we provide some additional schematic identities involving differentiation.

\begin{lemma}[Schematic identity for $\angLap \varphi$]
	\label{L:SCHEMATICEXPRESSIONFORANGULARLAPLACIAN}
	If $\varphi$ is a scalar function, then its angular Laplacian on $\ell_{t,u}$ can
	be expressed as follows, where the first two terms on RHS\,\eqref{E:SCHEMATICEXPRESSIONFORANGULARLAPLACIAN}
	is written precisely and the last two are written schematically:
	\begin{align} \label{E:SCHEMATICEXPRESSIONFORANGULARLAPLACIAN}
		\angLap \varphi
		& = 
			\Speed^{2} \sum_{A=2,3}
			\Yvf{A}(\Yvf{A} \varphi) + \Speed^2 \sum_{C,D = 2,3} \frac{X^C X^D}{(X^1)^2} \Yvf{C} \Yvf{D} \varphi 	+
			\smoothfunction(\controlvars,\LogDensity) 
			\cdot
			\tanderY \controlvars
			\cdot
			\tanderY \varphi.
	\end{align}
\end{lemma}
\begin{proof}
The angular Laplacian is given by $\angLap \varphi = (\gtorus^{-1})^{AB} \angD_A \geop{x^B} \varphi = (\gtorus^{-1})^{AB} \geop{x^A}\geop{x^B} \varphi - (\gtorus^{-1})^{CD} \gtorus^{FG} \slashed{\Gamma}_{CFD} \geop{x^G} \varphi$ in the intrinsic coordinates $(x^2,x^3)$ on $\ell_{t,u}$, where $\slashed{\Gamma}_{CFD} = \frac{1}{2}( \geop{x^C} \gtorus_{FD} + \geop{x^D} \gtorus_{CF} - \geop{x^F} \gtorus_{CD})$  is the fully lowered Christoffel symbol of $\gtorus$. By Prop.\,\ref{P:SCHEMATICSTRUCTUREOFVARIOUSTENSORSINTERMSOFCONTROLVARS} and \eqref{E:GEOP2TOCOMMUTATORS}--\eqref{E:GEOP3TOCOMMUTATORS}, the term $(\gtorus^{-1})^{CD} \gtorus^{FG} \slashed{\Gamma}_{CFD} \geop{x^G} \varphi$ is schematically represented by $\smoothfunction(\controlvars)\cdot \tanderY \controlvars \cdot \tanderY \varphi$. It is easy to see that the term $ (\gtorus^{-1})^{AB} \geop{x^A}\geop{x^B} \varphi$ is exactly the first two sums on RHS, \eqref{E:SCHEMATICEXPRESSIONFORANGULARLAPLACIAN} using \eqref{E:GEOP2TOCOMMUTATORS}--\eqref{E:GEOP3TOCOMMUTATORS} and \eqref{E:SMOOTHTORIGABEXPRESSION}--\eqref{E:SMOOTHGINVERSEABEXPRESSION}.
\end{proof}

\begin{lemma}[Identity satisfied by $\muX \Lunit^i$; \cite{abbrescia2022emergence}*{Lemma\,10.8}] 
\label{L:MUXLISCHEMATICIDENTITY}
There exist smooth functions,
all schematically denoted by ``$\smoothfunction$,''
such that the following identity holds:
\begin{align} \label{E:MUXLISCHEMATICIDENTITY}
	\muX \Lunit^i 
	& 
	= 
	\smoothfunction(\controlvars,\LogDensity) \cdot \muX \wavearray \cdot(-\updelta_1^i + \Xsmall^i)
	+
	\smoothfunction(\controlvars,\LogDensity) \cdot \muX \wavearraypartial
	+ 
	\upmu 
	\smoothfunction(\controlvars,\LogDensity) \tander \wavearray
	+ 
	\smoothfunction(\controlvars,\LogDensity) \tanderY \upmu.
\end{align}
\end{lemma}

\subsection{A schematic rewriting of the wave equations satisfied by $\wavearray$}
\label{SS:SCHEMATICREWRITINGOFGEOMETRICWAVEEQUATIONS}
The following lemma shows that $\muX \Psi$ obeys a transport equation with
source terms that are small but lose one derivative.
We will use it in Sect.\,\ref{S:MAINLINFINITYESTIMATES}, 
when we derive improvements of the auxiliary bootstrap assumptions.

\begin{lemma}[A schematic rewriting of the wave equations satisfied by $\wavearray$; \cite{abbrescia2022emergence}*{Lemma\,10.11}]
\label{L:SCHEMATICREWRITINGOFGEOMETRICWAVEEQUATIONS}
The covariant wave equations 
\eqref{E:COVARIANTWAVEEQUATIONSWAVEVARIABLES}
verified by $\Psi \in \{\LogDensity,v^1,v^2,v^3,\Ent\}$ 
can be expressed in the following schematic form:
\begin{align} \label{E:SCHEMATICREWRITINGOFWAVEEQUATIONSSATISFIEDBYWAVEVARIABLES}
\Lunit \muX \Psi 
& = 
\smoothfunction(\badcontrolvars,\LogDensity) \tander^2 \wavearray 
+ 
\smoothfunction(\badcontrolvars,\LogDensity,\comder \wavearray) \tander \controlvars  
+ 
\smoothfunction(\badcontrolvars, \GradEnt, \comder \wavearray)
\cdot
\tander^{\leq 1}(\vortrenormalized,\GradEnt).
\end{align}\end{lemma}

\section{Norms, area forms, and volume forms}
\label{S:NORMSANDVOLUMEFORMS}

In this section, we introduce the appropriate \emph{weighted} norms on regions of spacetime that are tied to the double-null-acoustic geometry. We also introduce the area and volume forms that we use in our $L^2$ analysis. 

\subsection{Weight function used in the analysis}

Our main a priori estimates will incorporate \emph{weighted} energies and null-fluxes. The weight $\weight$ must be properly chosen to capture the correct singular nature of the shock at the crease \emph{in the $\newuL$ directions}. Hence, we will call $\weight$ the \textbf{blow-up weight}, but emphasize that $\weight$ itself is regular. 

\begin{definition}[The blow-up weight] \label{D:BLOWUPWEIGHT} We define the \textbf{blow-up weight} $\weight$ to be the function defined in double-null coordinates as: 
\begin{equation} \label{E:BLOWUPWEIGHT}
\weight(\ubar,u,x^2,x^3) = - \ubar.
\end{equation}
We will only consider cases where $I \subset (-\infty,0]$ and so $\weight \ge 0$ on $\characteristicdiamondtwoarg{I}{J}$.
\end{definition}
 
\begin{remark}[Notation for the weight] \label{R:NOTATIONFORWEIGHT.} 
We will keep the notation $\weight$ whenever the scalar function $-\ubar$ is used as a weight, e.g. in our weighted $L^2$ energy-flux analysis. Although we could clearly write instead $-\ubar$, we use $\weight$ in those contexts with future applications in mind because we will derive numerous identities that hold true for arbitrary weight functions
\end{remark}

\subsection{Area forms, volume forms, and corresponding $L^2$ norms} 
\label{SS:FORMSANDL2NORMS}
We now define the area and volume forms on the rough subsets (see Def.\,\ref{D:SUBSETSOFSPACETIME})
$\doublenulltoritwoarg{\ubar}{u}, \, \outgoingcharacteristicsurfacetwoarg{u}{I}, \, \ingoingcharacteristicsurfacetwoarg{\ubar}{J}$, and $\characteristicdiamondtwoarg{I}{J}$ that we will use in our analysis.
We also define corresponding $L^2$-type norms.
Our definitions are in terms of the double-null coordinates $(\ubar,u,x^2,x^3)$
because those are the coordinates that we use
in our energy identities, where the forms arise.

\subsubsection{Geometric forms and related integrals}
\label{SSS:GEOMETRICFORMSANDINTEGRALS}

\begin{definition}[Geometric forms and related integrals]\label{D:VOLFORMS}
	\hfill
	\begin{itemize}
		\item Recall that $\gnulltori$ is the first fundamental form of the torus $\doublenulltoritwoarg{\ubar}{u}$.
		We define the canonical area form of $\doublenulltoritwoarg{\ubar}{u}$ 
		induced by $\gnulltori$
		in the double-null coordinates $(\ubar,u,x^2,x^3)$ by:
		\begin{align} \label{E:AREAFORMDOUBLENULLTORUS}
			\voldoublenulltori 
			& = \voldoublenulltori(\ubar',u',x^2,x^3) 
			\eqdef 
			\sqrt{\det \gnulltori (\ubar',u',x^2,x^3)} \, \mathrm{d} x^2 \mathrm{d} x^3,
		\end{align}
		where $\det \gnulltori(\ubar',u',x^2,x^3)$ is the determinant of the $2 \times 2$ matrix
		$\left(\gnulltori(\ubar',u',x^2,x^3) \left(\nullgeop{x^A},\nullgeop{x^B}\right) \right)_{A,B=2,3}$
		(see RHS\,\eqref{E:DOUBLENULLTORIFIRSTFUNDRESTRICTEDTODOUBLENULLTORI}).
		\item We define the (non-canonical) area form $\volingoingnullhypersurface$ of $\ingoingcharacteristicsurfacetwoarg{\ubar}{J}$ 
		in the double-null coordinates $(\ubar,u,x^2,x^3)$ by: 
		\begin{align} \label{E:VOLUMEFORMINGOINGNULLHYPERSURFACE}
			\volingoingnullhypersurface 
			& = \volingoingnullhypersurface(\ubar,u',x^2,x^3) 
			= \voldoublenulltori(\ubar,u',x^2,x^3) \, \mathrm{d} u'.
		\end{align}
		\item We define the (non-canonical) area form $\voloutgoingnullhypersurface$ on $\outgoingcharacteristicsurfacetwoarg{u}{I}$ 
		in the double-null coordinates $(\ubar,u,x^2,x^3)$ by:  
		\begin{align} \label{E:VOLUMEFORMOUTGOINGNULLHYPERSURFACE}
			\voloutgoingnullhypersurface 
			= 
			\voloutgoingnullhypersurface(\ubar',u,x^2,x^3) 
			& \eqdef \voldoublenulltori(\ubar',u,x^2,x^3) \, \mathrm{d} \ubar'.
		\end{align} 
		\item We define the (non-canonical) volume form $\voldiamond$ of $\characteristicdiamondtwoarg{I}{J}$ 
		in the double-null  coordinates $(\ubar,u,x^2,x^3)$ by:
		\begin{align} \label{E:VOLUMEFORMSPACTEIMROUGH}
			\voldiamond 
			= \voldiamond(\ubar',u',x^2,x^3) 
			& 
			\eqdef 
			\voldoublenulltori(\ubar',u',x^2,x^3) \, \mathrm{d} u' \, \mathrm{d} \ubar'.
		\end{align}
	\end{itemize}
	
	Unless we explicitly indicate otherwise, 
	all integrals along 
	$\doublenulltoritwoarg{\ubar}{u}, \, \outgoingcharacteristicsurfacetwoarg{u}{I}, \, \ingoingcharacteristicsurfacetwoarg{\ubar}{J}$, and $\characteristicdiamondtwoarg{I}{J}$
	are defined with respect to the above forms. 
	Moreover, we will often suppress the variables with respect to which we integrate, 
	e.g., we write:
	\begin{subequations} 
		\begin{align}
			\int_{\doublenulltoritwoarg{\ubar}{u}} 
			f 
			\, \voldoublenulltori 
			& 
			\eqdef 
			\int_{(x^2,x^3) \in \T^2}  
			f(\ubar,u,x^2,x^3) 
			\, \voldoublenulltori(\ubar,u,x^2,x^3),   
			\label{E:DOUBLENULLTORUSINTEGRAL} \\
			\int_{\outgoingcharacteristicsurfacetwoarg{u}{I}} 
			f \, 
			\voloutgoingnullhypersurface 
			& 
			\eqdef \int_{\ubar' \in I}  
			\int_{(x^2,x^3) \in \T^2} 
			f(\ubar',u,x^2,x^3) \, \voldoublenulltori(\ubar',u,x^2,x^3)
			\, \mathrm{d} \ubar', 
			\label{E:OUTGOINGNULLHYPINTEGRAL} \\
			\int_{\ingoingcharacteristicsurfacetwoarg{\ubar}{J}} 
			f 
			\, \volingoingnullhypersurface 
			& \eqdef 
			\int_{u' \in J} \int_{(x^2,x^3) \in \T^2} 
			f(\ubar,u',x^2,x^3) 
			\, \voldoublenulltori(\ubar,u',x^2,x^3) \, \mathrm{d} u', 
			\label{E:INGOINGNULLHYPSURFACEINTEGRAL} \\
			\int_{\characteristicdiamondtwoarg{I}{J}} 
			f 
			\, \voldiamond 
			& 
			\eqdef 
			\int_{\ubar' \in I} 
			\int_{u' \in J} 
			\int_{(x^2,x^3) \in \T^2} 
			f(\ubar',u',x^2,x^3) 
			\, \voldoublenulltori(\ubar',u',x^2,x^3)
			\, 
			\mathrm{d} u' \, \mathrm{d}\ubar'.
			\label{E:DIAMONDSPACETIMEINTEGRAL}
		\end{align}
	\end{subequations}
\end{definition}

In a few of our calculations, we will also refer to the canonical volume forms of $\gfour$ relative to the double null coordinates, which we provide
in the following definition.

\begin{definition}[Canonical volume forms relative to the double-null coordinates] 
	\label{D:CANONICALVOLFORMSINDOUBLECOORDINATES} 
	We define 
	$\volcanonical_{\gfour} \eqdef \sqrt{|\det \gfour|}\, \mathrm{d} x^2 \, \mathrm{d} x^3\, \mathrm{d} u' \, \mathrm{d} \ubar'$ 			to be the canonical area form on $\characteristicdiamondtwoarg{I}{J}$ induced by 
	the spacetime metric $\gfour$, 
	where $\mbox{\upshape det} \gfour$ is evaluated relative to the double-null
	$(\ubar',u',x^2,x^3)$ on $\characteristicdiamondtwoarg{I}{J}$.
	
\end{definition}

\subsubsection{Identities involving the forms}
\label{SSS:IDENTITIESINVOLVINGVOLUMEFORMSINROUGHADAPATEDCOORDINATES}

\begin{lemma}[Identities involving $\volcanonical_{\gfour}$ and $\voldoublenulltori,\, \voldiamond$]
	\label{L:IDENTITIESINVOLVINGVOLUMEFORMSINROUGHADAPATEDCOORDINATES}
	The following identities hold relative to the double-null coordinates $(\ubar,u,x^2,x^3)$,
	e.g., $\mbox{\upshape det} \gnulltori(\ubar',u',x^2,x^3)$ is the determinant of the $2 \times 2$ matrix
	$\left(\gnulltori(\ubar',u',x^2,x^3) \left(\nullgeop{x^A},\nullgeop{x^B}\right) \right)_{A,B=2,3}$:
	\begin{equation}
		\mbox{\upshape det} \gfour 
		= - \MagnitueofinnerproductofnewLandnewuL^2 \mbox{\upshape det} \gnulltori. 
		\label{E:DETACOUSTICALMETRICINDOUBLENULLCOORDS} 
	\end{equation}

	Moreover, with $
	\volcanonical_{\gfour}$
	and
	$\voldoublenulltori$ denoting the area and volume forms from
	Defs.\,\ref{D:VOLFORMS} and\,\ref{D:CANONICALVOLFORMSINDOUBLECOORDINATES}, 
	we have the following identities relative to the double-null coordinates:
	\begin{equation}
		\volcanonical_{\gfour}
		= \MagnitueofinnerproductofnewLandnewuL \, \voldiamond = \upmu \ReciprocalLunitAppliedtoTimeFunction \, \voldiamond
		= 
		\upmu \ReciprocalLunitAppliedtoTimeFunction \,
		\voldoublenulltori
		\, \mathrm{d} u' 
		\, \mathrm{d} \ubar'. 
		\label{E:VOLFORMACOUSTICALMETRICDOUBLENULLCOORDS}
	\end{equation}
	
\end{lemma}

\begin{proof}

\end{proof}

\subsubsection{Geometric $L^2$ and $L^{\infty}$ norms}
\label{SSS:GEOMETRICL2NORMS}

\begin{definition}[Geometric $L^2$ norms] 
	\label{D:GEOMETRICL2NORMS}
	We define the $L^2$ weighted norms with respect to the area and volume forms introduced 
	in Def.\,\ref{D:VOLFORMS}. 
	Recall that we measure the norm of $\ell_{t,u}$-tangent tensorfields with $\gtorus$, i.e.,
	if $\upxi$ is a type $\binom{0}{2}$ $\ell_{t,u}$-tangent tensorfield, then
	$|\upxi|_{\gtorus}^2 \eqdef (\gtorus^{-1})^{\alpha \beta}(\gtorus^{-1})^{\sigma \delta} \upxi_{\alpha \sigma} \upxi_{\beta \delta}$. 
	Then for scalar functions or $\ell_{t,u}$-tangent tensorfields $\upxi$ and positive numbers $q$, 
	we define: 
	\begin{subequations}
		\begin{align}
			\| \upxi \|_{L_{q}^2\left(\doublenulltoritwoarg{\ubar}{u}\right)}
			& 
			\eqdef 
			\left( 	\int_{\doublenulltoritwoarg{\ubar}{u}} 
			|\upxi|_{\gtorus}^2 \, \weight^{q}
			\, \voldoublenulltori 
			\right)^{1/2},
			& 
			\| \upxi \|_{L_{q}^2\left(\outgoingcharacteristicsurfacetwoarg{u}{I}\right)} 
			& \eqdef 
			\left( 
			\int_{\outgoingcharacteristicsurfacetwoarg{u}{I}} |\upxi|_{\gtorus}^2 \, \weight^{q} \, \voloutgoingnullhypersurface 
			\right)^{1/2}, 
			\label{E:GEOMETRICL2NORMSTORIANDOUTGOINGNULLHYPERSURFACES} 
			\\
			\| \upxi \|_{L_{q}^2\left(\ingoingcharacteristicsurfacetwoarg{\ubar}{J}\right)}
			& 
			\eqdef \left( 
			\int_{\ingoingcharacteristicsurfacetwoarg{\ubar}{J}} |\upxi|_{\gtorus}^2 \, \weight^{q}
			\, \volingoingnullhypersurface 
			\right)^{1/2}, 
			& 
			\| \upxi \|_{L_{q}^2\left(\characteristicdiamondtwoarg{I}{J}\right)} 
			& \eqdef 
			\left( 
			\int_{\characteristicdiamondtwoarg{I}{J}} |\upxi|_{\gtorus}^2 \, \weight^{q} \, \voldiamond 
			\right)^{1/2}.
			\label{E:GEOMETRICL2NORMSINGOINGNULLHYPERSURFACESANDSPACETIMEREGIONS}
		\end{align}
	\end{subequations}
\end{definition}

\begin{definition}[Geometric $L^{\infty}$ norms] 
	\label{D:GEOLINFINITYTORUS}
	For scalar functions or $\ell_{t,u}$-tangent tensorfields $\upxi$, we define 
	the following $L^{\infty}$ norm on the double-null tori $\doublenulltoritwoarg{\ubar}{u}$:
	\begin{align} \label{E:GEOLINFINITYDOUBLENULLTORUS}
		\| \upxi \|_{L^{\infty}\left(\doublenulltoritwoarg{\ubar}{u}\right)}
		& 
		\eqdef 
		\mbox{ess sup}_{(x^2,x^3) \in \mathbb{T}^2}
		|\upxi|_{\gtorus}(\ubar,u,x^2,x^3),
	\end{align}
	where on RHS\,\eqref{E:GEOLINFINITYDOUBLENULLTORUS}, we are viewing
	$\upxi$ as a function of the rough adapted coordinates
	$(\ubar,u,x^2,x^3)$.
\end{definition}

\begin{remark}[Carefully note the role of $\gtorus$] \label{R:ROLEOFFIRSTFUNDOFSMOOTHTORI}
	We stress that $\gtorus$ is the Riemannian metric on the smooth tori,
	even though the integrals defining
	$\| \cdot \|_{L_{q}^2\left(\doublenulltoritwoarg{\ubar}{u}\right)}$,
	$	\| \cdot \|_{L_{q}^2\left(\outgoingcharacteristicsurfacetwoarg{u}{I}\right)} $,
	$	\| \cdot \|_{L_{q}^2\left(\ingoingcharacteristicsurfacetwoarg{\ubar}{J}\right)}$,
	and
	$	\| \cdot \|_{L_{q}^2\left(\characteristicdiamondtwoarg{I}{J}\right)}$ 
	are over regions and with respect to forms
	tied to the rough adapted coordinates.
	
	Similarly, on RHS\,\eqref{E:GEOLINFINITYDOUBLENULLTORUS}, $|\upxi|_{\gtorus}$
	is the pointwise norm of $\upxi$ with respect to the 
	Riemannian metric $\gtorus$ on the smooth tori.
\end{remark}


\section{Modified quantities for controlling the acoustic geometry} \label{S:CONSTRUCTIONOFMODIFIEDQUANTITIES}
In this section we introduce and derive transport equations for ``modified'' versions of the eikonal function quantity $\mytr_{\gtorus} \upchi$. When combined with elliptic estimates on the tori $\doublenulltoritwoarg{\ubar}{u}$, the modified quantities will allow us to estimate top-order derivatives of $\mytr_{\gtorus} \upchi$ without losing derivatives.

\subsection{Definition of the modified quantities} 
\label{SS:DEFSOFMODQUANTITIES}
The definition of the modified quantities is motivated by the structure of the terms in Lemma\,\ref{L:RICLLDECOMPOSITIONS} below.


\begin{lemma}[The key identity verified by $\upmu \Ricfour_{\Lunit \Lunit}$] \label{L:RICLLDECOMPOSITIONS} Assume that the entries of $\wavearray = (v^1,v^2, v^2,v^3,\LogDensity,\Ent)$ solve the geometric wave equations \eqref{E:VELOCITYWAVEEQUATION}--\eqref{E:ENTROPYWAVEEQUATION}. Then the following identity holds:
	\begin{align}  
		\upmu \Ricfour_{\Lunit \Lunit} = & L \left\lbrace - \vec{G}_{\Lunit \Lunit} \diamond \muX \wavearray - \frac{1}{2} \upmu \mytr_{\gtorus} \angG \diamond \Lunit \wavearray - \frac{1}{2} \upmu \vec{G}_{\Lunit \Lunit} \diamond \Lunit \wavearray + \upmu \angG_{\Lunit}^{\#} \diamond \cdot \angrmD \wavearray \right\rbrace \notag \\
		& \ \ + \mathfrak{A},\label{E:RICCICONTRACTEDLANDL}
	\end{align}
	where $\mathfrak{A}$ has the following schematic structure:
	\begin{equation} \label{E:AINHOMRIC}
		\mathfrak{A}  = \smoothfunction(\badcontrolvars,\gtorus^{-1},\angrmD \vec{x}, \muX \wavearray, \tander \wavearray) \tander \wavearray + \smoothfunction(\controlvars) \upmu (\VortVort,\DivGradEnt) + \smoothfunction(\controlvars) (\upmu \tander \wavearray) (\Omega,\GradEnt) + \smoothfunction(\controlvars)\muX \wavearray (\Omega,\GradEnt) + \smoothfunction(\controlvars)(\Omega,\GradEnt)^2 
	\end{equation}
	
	Moreover, without assuming that the geometric wave equations \eqref{E:VELOCITYWAVEEQUATION}--\eqref{E:ENTROPYWAVEEQUATION} are satisfied, we have:
	\begin{align} 
		\Ricfour_{\Lunit \Lunit} = & \frac{\Lunit \upmu}{\upmu} \mytr_{\gtorus} \upchi + L \left\lbrace -\frac{1}{2} \mytr_{\gtorus} \angG \diamond \Lunit \wavearray + \angG_{\Lunit}^{\#} \diamond \cdot \angrmD \wavearray \right\rbrace \notag \\
		& - \frac{1}{2} \vec{G}_{\Lunit \Lunit} \diamond \angLap \wavearray + \mathfrak{B}, \label{E:RICCICONTRACTEDLANDLtype2}
	\end{align}
	where $\mathfrak{B}$ has the following schematic structure:
	\begin{equation} \label{E:BINHOMRIC}
		\mathfrak{B} = \text{$\smoothfunction(\controlvars, \gtorus^{-1}, \angrmD \vec{x} )$} (\tander \wavearray) \tander \controlvars.
	\end{equation}

\end{lemma}

\begin{proof}[Sketch of a proof] Identities \eqref{E:RICCICONTRACTEDLANDL} and \eqref{E:RICCICONTRACTEDLANDLtype2} were essentially proved in \cite[Lemma 6.1]{jLjS2018} in the case of two space dimensions. Nearly identical arguments found in \cite[Corollary 11.4]{jS2016b} hold true in the current case of three space dimensions. In fact, the proof of \eqref{E:RICCICONTRACTEDLANDLtype2} goes through without any significant changes. The only significant new feature in the present work is that RHS\,\eqref{E:RICCICONTRACTEDLANDL} features $\smoothfunction(\controlvars)(\upmu \VortVort + \upmu \DivGradEnt + \upmu \tander \wavearray (\Omega,S) + \muX \wavearray (\Omega,S) + (\Omega,S)^2)$. To see how these terms arise, one uses the identity \eqref{E:BOXDECOMPLOUTSIDE} and the wave equations \eqref{E:VELOCITYWAVEEQUATION}--\eqref{E:ENTROPYWAVEEQUATION} to express 
	
	\begin{multline}\label{E:DETAILEDINHOM}
		-\frac{1}{2} \upmu \vec{G}_{\Lunit \Lunit} \diamond \angLap \wavearray = - \frac{1}{2} \Lunit \left\lbrace \vec{G}_{\Lunit \Lunit} \diamond (\upmu \Lunit\wavearray + 2 \muX \wavearray) \right\rbrace - \frac{1}{2} \mytr_{\gtorus}\upchi \vec{G}_{\Lunit \Lunit} \diamond \muX \wavearray \\
		+ \smoothfunction(\controlvars) \cdot \text{Inhom} +  f(\badcontrolvars,\gtorus^{-1},\angrmD \vec{x}, \muX \wavearray, \tander \wavearray) \tander \wavearray. 
	\end{multline}
	In \eqref{E:DETAILEDINHOM}, Inhom denotes $\upmu\cdot($RHS\,\eqref{E:VELOCITYWAVEEQUATION}--\eqref{E:ENTROPYWAVEEQUATION}$)$.
	Decomposing the resulting $\upmu$-multiplied null forms $\upmu\mathfrak{Q}$ using Lemma\,\ref{L:CRUCIALSTRUCTUREOFNULLFORMS} and the linear terms $\upmu\mathfrak{L}$ schematically as $\smoothfunction(\controlvars)(\upmu \tander \wavearray (\Omega,S) + \muX \wavearray (\Omega,S) + (\Omega,S)^2)$ , we arrive to \eqref{E:RICCICONTRACTEDLANDL}--\eqref{E:AINHOMRIC}. In a detailed proof (see \cite[Lemma 6.1]{jSgHjLwW2016}), one finds that the term 
	\[ -\frac{1}{2} \mytr_{\gtorus}\upchi \vec{G}_{\Lunit \Lunit} \diamond \muX \wavearray\]
	on RHS\,\eqref{E:DETAILEDINHOM} is cancelled by another term and hence does not appear in \eqref{E:RICCICONTRACTEDLANDL}--\eqref{E:AINHOMRIC}.  
\end{proof}

We are now ready to define the modified quantities. Note that the RHS\,\eqref{E:MODIFIEDQUANTITYINHOM} is present in the perfect $\Lunit$ differentiation of \eqref{E:RICCICONTRACTEDLANDL}.

\begin{definition}[Modified versions of the derivatives of $\mytr_{\gtorus} \upchi$] \label{D:FULLYANDPARTIALLYMODIFIEDQUANTITIES} The \emph{fully modified quantity} $\fullymodquant{\tander^N}$ is defined to be: 
	\begin{subequations}
		\begin{align}
			\fullymodquant{\tander^N} & \eqdef \upmu \tander^N \mytr_{\gtorus} \upchi + \tander^N \mathfrak{X},
			\label{E:FULLYMODIFIEDQUANTITY} \\
			\mathfrak{X} & \eqdef - \vec{G}_{\Lunit \Lunit} \diamond \muX \wavearray - \frac{1}{2} \upmu \mytr_{\gtorus} \angG \diamond \Lunit \wavearray - \frac{1}{2} \upmu \vec{G}_{\Lunit \Lunit} \diamond \Lunit \wavearray + \upmu \angG_{\Lunit}^{\#} \diamond \cdot \angrmD \wavearray. \label{E:MODIFIEDQUANTITYINHOM}
		\end{align}
	\end{subequations}
	The \emph{partially modified quantity} $\partialmodquant{\tander^N}$ is defined to be:
	\begin{subequations}
		\begin{align}
			\partialmodquant{\tander^N} 
			& 
			\eqdef 
			\tander^N \mytr_{\gtorus} \upchi 
			+ 
			\partialmodquantinhom{\tander^N}, 
			\label{E:PARTIALMODIFIEDQUANTITY} 
			\\
			\partialmodquantinhom{\tander^N} 
			& \eqdef 
			- 
			\frac{1}{2} \mytr_{\gtorus} \angG \diamond \Lunit \tander^N \wavearray 
			+ 
			\angG_{\Lunit}^{\#} \diamond \cdot \angrmD \tander^N \wavearray. 
			\label{E:PARTIALMODIFIEDQUANTITYINHOMDEFINITION}
		\end{align}
	\end{subequations}
	The ``$0^{\text{th}}$-order'' version of \eqref{E:PARTIALMODIFIEDQUANTITYINHOM} is simply:
	\begin{align} \label{E:PARTIALMODIFIEDQUANTITYINHOMZEROORDER}
		\widetilde{\mathfrak{X}} 
		& \eqdef - \frac{1}{2} \mytr_{\gtorus} \angG \diamond \Lunit \wavearray + \angG_{\Lunit}^{\#} \diamond \cdot \angrmD \wavearray.
	\end{align} 
\end{definition}

\subsection{Transport equations for the modified quantities}

In this subsection we will derive transport equations for the fully modified quantities. We begin with the following lemma, which records the transport equation satisfied by $\mytr_{\gtorus}\upchi$ along the integral curves of $\Lunit$. This identity is the analog of well-known \emph{Raychaudhuri equation} in General Relativity \cite{aR1955}.


\begin{lemma}[Raychaudhuri-type transport equation for $\mytr_{\gtorus}\upchi$] The following identity holds:
	\begin{align}  \label{E:RAYCHAUDHURITRANSPORTCHI}
		\upmu \Lunit \mytr_{\gtorus}\upchi 
		& = \Lunit \upmu \mytr_{\gtorus}\upchi - \upmu \Ricfour_{\Lunit \Lunit} - \upmu |\upchi|^2. 
	\end{align}
\end{lemma}
\begin{proof}
	The same proof of \cite[(11.23)]{jS2016b} holds in the current setting.
\end{proof}


By differentiating \eqref{E:RAYCHAUDHURITRANSPORTCHI}, straightforward commutation identities can be used to prove the following transport equations for the fully modified quantity $\fullymodquant{\tander^N}$. Minor modifications to the proof of \cite[Prop. 6.2]{jSgHjLwW2016} imply the following proposition.


\begin{proposition}[Transport equation satisfied by $\fullymodquant{\tander^N}$] Assume that $\wavearray = (\RRiemann,\LRiemann,v^2,v^3,\Ent)$ verify the geometric wave equations \eqref{E:VELOCITYWAVEEQUATION}--\eqref{E:ENTROPYWAVEEQUATION}. Then the following identity holds:
	\begin{align}  \label{E:TRANSPORTEQUATIONFORFULLYMODIFIEDQUANTITY} \begin{split}
		L \fullymodquant{\tander^N}  &  = 2 \frac{\Lunit \upmu}{\upmu}  \fullymodquant{\tander^N} - 2 \frac{\Lunit \upmu}{\upmu} \tander^N \mathfrak{X} +  \upmu [\Lunit,\tander^N] \mytr_{\gtorus} \upchi    \\
		& \ \ + [\Lunit,\tander^N] \mathfrak{X} + [\upmu,\tander^N] \Lunit \mytr_{\gtorus}\upchi + [\tander^N,\Lunit \upmu] \mytr_{\gtorus}\upchi   \\ 
		&\ \  - \tander^N(\upmu |\upchi|^2)  - \tander^N \mathfrak{A},
		\end{split}
	\end{align}
	where $\mathfrak{A}$ is the term on the right hand side of \eqref{E:RICCICONTRACTEDLANDL} and $\tander^N$ is the same differential operator every time it appears in \eqref{E:TRANSPORTEQUATIONFORFULLYMODIFIEDQUANTITY}
\end{proposition}

\begin{proposition}[Transport equation satisfied by $\partialmodquant{\tanderY^{N-1}}$] The following identity holds:
	\begin{equation} \label{E:TRANSPORTEQUATIONFORPARTIALMODIFIEDQUANTITY}
		\Lunit \partialmodquant{\tanderY^{N-1}} =  \frac{1}{2} \vec{G}_{\Lunit \Lunit}\diamond \angLap \tanderY^{N-1} \wavearray +\, ^{(\tanderY^{N-1})}\mathfrak{B},
	\end{equation}
	where:
	\begin{align}
		\begin{split} \label{E:PARTIALMODIFIEDQUANTITYINHOM}
		^{(\tanderY^{N-1})}\mathfrak{B} & = - \tanderY^{N-1} \mathfrak{B} - \tanderY^{N-1} |\upchi|^2   \\
		& \ \ + \frac{1}{2} [\tanderY^{N-1},\vec{G}_{\Lunit \Lunit}] \diamond \angLap \wavearray + \frac{1}{2} \vec{G}_{\Lunit \Lunit} \diamond [ \tanderY^{N-1},\angLap]\wavearray  \\
		& \ \ + [\Lunit,\tanderY^{N-1}] \mytr_{\gtorus}\upchi + [\Lunit,\tanderY^{N-1}] \widetilde{\mathfrak{X}} + L \left\lbrace \partialmodquantinhom{\tanderY^{N-1}} - \tanderY^{N-1}\widetilde{\mathfrak{X}}\right\rbrace, 
		\end{split}
	\end{align}
	where $\mathfrak{B}$ is as in \eqref{E:BINHOMRIC} and $\tanderY^{N-1}$ is the same differential operator every time it appears in \eqref{E:TRANSPORTEQUATIONFORPARTIALMODIFIEDQUANTITY}.
	
\end{proposition}

\begin{proof}
	The proof begins with dividing the Raychaudhuri equation for $\mytr_{\gtorus}\upchi$ \eqref{E:RAYCHAUDHURITRANSPORTCHI} by $\upmu$ to find 
	\[\Lunit \mytr_{\gtorus}\upchi = \frac{\Lunit \upmu}{\upmu} \mytr_{\gtorus}\upchi - \Ricfour_{\Lunit \Lunit} - |\upchi|^2.\] 
	Inserting the Ricci identity \eqref{E:RICCICONTRACTEDLANDLtype2} and using the $0^{\text{th}}$ order $\widetilde{\mathfrak{\upchi}}$ defined in \eqref{E:PARTIALMODIFIEDQUANTITYINHOMZEROORDER}, we obtain: 
	\begin{equation} \label{E:TRANSPORTEQUATIONFORPARTIALMODIFIEDQUANTITYzerothorder}
		\Lunit \left(\mytr_{\gtorus} \upchi + \widetilde{\mathfrak{X}}\right) = \frac{1}{2} \vec{G}_{\Lunit \Lunit} \diamond \angLap \wavearray - |\upchi|^2 - \mathfrak{B}.
	\end{equation}
	The transport equation \eqref{E:TRANSPORTEQUATIONFORPARTIALMODIFIEDQUANTITY} then follows by differentiating \eqref{E:TRANSPORTEQUATIONFORPARTIALMODIFIEDQUANTITYzerothorder} by $\tanderY^{N-1}$ and straightforward commutations. 
	
\end{proof}


\section{Construction and estimates for the ingoing eikonal function}\label{S:CONSTRUCTIONANDESTIMATESFORINGOINGEIKONALFUNCTION}

The goal of this section is to provide the construction of and to derive pointwise estimates for the ingoing eikonal function $\ubar$. This section is organized as follows. In Sect.\,\ref{SS:SECONDORDERQUASILINEATRANSPORTSYSTEMFORINGOINGEIKONALFUNCTION}, we re-write the eikonal equation \eqref{E:INGOINGEIKONALEQUATION} as an equivalent \emph{second order} PDE system for the first derivatives of $\ubar$. The reason for this is that the eikonal equation is a \emph{fully nonlinear} hyperbolic equation, and so by differentiating, we derive a \emph{quasilinear} system of PDEs for $\Lunit  \ubar, \, \muX  \ubar, \, \geop{x^A} {\ubar}$ which is much simpler to analyze. As Lemma\,\ref{L:QUASILINEARTRANSPORTSYSTEMFORDERIVATIVESOFUBAR} reveals, the structure of the eikonal equation actually forces this quasilinear system to be a \emph{transport} system along integral curves of $\newuL$. In Sect.\,\ref{SS:GRADIENTDATAESTIMATESFORINGOINGEIKONALFUNCTIONONINITIALDATASURFACE}, we derive sharp estimates for all first derivatives ${\ubar}$ on the initial data surface $\datahypfortimefunctiontwoarg{0}{[\timefunction_0,0]}$ that are compatible with the eikonal equation. Then in Sect.\,\ref{SS:ANALYSISOFINGOINGEIKONALFUNCTION}, we actually construct a solution to \eqref{E:INGOINGEIKONALEQUATIONINITIALVALUEPROBLEM} through a fixed point argument.

\subsection{Second order quasilinear transport system for the ingoing eikonal function}

\label{SS:SECONDORDERQUASILINEATRANSPORTSYSTEMFORINGOINGEIKONALFUNCTION}

Since the eikonal equation is a fully nonlinear hyperbolic equation, in this section, we derive the equivalent quasilinear system of equations satisfied by derivatives of $\ubar$. Roughly, by taking a derivative of $(\gfour^{-1})^{\alpha\beta} \p_\alpha \ubar \p_\beta \ubar$ with respect to $\pmb{\partial}$, we obtain a PDE system for $\pmb{\partial}^2 \ubar$ where the coefficients for $\pmb{\partial}^2 \ubar$ are schematically $F(\pmb{\partial} \ubar)$. As Lemma\,\ref{L:QUASILINEARTRANSPORTSYSTEMFORDERIVATIVESOFUBAR} will reveal, this quasilinear system is actually a quasilinear \emph{transport} system along integral curves of $\newuL$ where the \emph{only} second order derivatives are of the form $\newuL \pmb{\partial} \ubar$. This fact is not difficult to see: for any vectorfield $Z$,
\begin{align}
	\begin{split}
		 0 & = \Lie_{Z} \left( \gfour^{-1}(\mathrm{d} \ubar, \mathrm{d} \ubar)\right) = 2 \gfour^{-1} (\mathrm{d}\Lie_Z \ubar, \mathrm{d} \ubar) + (\Lie_Z\gfour^{-1})(\mathrm{d} \ubar,\mathrm{d}\ubar)  \\
		 & = 2 (\gfour^{-1})^{\alpha\beta} \p_\alpha Z \ubar \p_\beta \ubar   + (\Lie_Z\gfour^{-1})(\mathrm{d} \ubar,\mathrm{d}\ubar)  \\
		 & = - 2 \uLgeo ( Z \ubar) + (\Lie_Z\gfour^{-1})(\mathrm{d} \ubar,\mathrm{d}\ubar). \label{E:DIFFERENTIATEDEIKONALEQUATIONGENERAL}
	\end{split}
\end{align}
The point of Lemma\,\ref{L:QUASILINEARTRANSPORTSYSTEMFORDERIVATIVESOFUBAR}  is then to obtain precise expressions for the lower order coefficients of the quasilinear transport system satisfied by $\newuL Z \ubar$ with respect to the null-geometry scalars. This is beneficial because it immediately reveals certain monotonicity properties of the coefficients which is harder to see from \eqref{E:DIFFERENTIATEDEIKONALEQUATIONGENERAL}. 

Because the bulk of the paper will use commutators tailored to the geometric coordinate system $(t,u,x^2,x^3)$, it will be convenient to work with $Z \in \{\Lunit, \muX, \geop{x^2},\geop{x^3}\}$. 

We begin with the following preliminary lemma whose proof is an immediate consequence of \eqref{E:SMOOTHTORUSMETRICINTERMSOFSIGMATMETRICANDX} and \eqref{E:INGOINGEIKONALEQUATIONINITIALVALUEPROBLEM}.

\begin{lemma}[The eikonal equation in terms of commutator vectorfields] \label{L:EIKONALEQUATIONWITHOUTTHEMUWEIGHT} 
	Let $\ubar$ be a solution to the eikonal equation \eqref{E:INGOINGEIKONALEQUATION}. Then ingoing eikonal function in rough adapted coordinates. Then the following identity holds:
	\begin{equation}
		0 = -(\Lunit \ubar)^2 - 2(\Lunit \ubar) X \ubar +  |\angD \ubar|^2_{\gtorus}. \label{E:EIKONALEQUATIONWITHOUTTHEMUWEIGHT}
	\end{equation}

\end{lemma} 

The following lemma is the main result of Sect.\,\ref{SS:SECONDORDERQUASILINEATRANSPORTSYSTEMFORINGOINGEIKONALFUNCTION},

\begin{lemma}[Quasilinear system satisfied by $Z\ubar$] \label{L:QUASILINEARTRANSPORTSYSTEMFORDERIVATIVESOFUBAR}
	Let $\newuL$ be the vectorfield defined in \eqref{E:DOUBLENULLEIKONALFUNCTIONNORMALIZEDNULLVECTORFIELDS}. If $\ubar$ is a solution to the eikonal equation $(\gfour^{-1})^{\alpha\beta} \p_\alpha \ubar \p_\beta \ubar = 0$ with $\Lunit \ubar > 0$, then the scalar functions $\Lunit \ubar, \, \muX \ubar,  \, \geop{x^A} \ubar$ solve the following quasilinear system of equations for $A =2,3$:
		\begin{align} 
			\newuL \Lunit \ubar & = - \frac{1}{2} (\Lunit \upmu) \Lunit \ubar + [\muX, \Lunit] \ubar -  \MagnitueofinnerproductofnewLandnewuL \angD^A \ubar \left[\geop{x^A},\Lunit\right] \ubar +  \frac{\Lunit \upmu}{2\Lunit \ubar} |\angD \ubar|_{\gtorus}^2 +\frac{ \upmu}{2\Lunit \ubar} \left( \Lunit (\gtorus^{-1})^{AB}\right) \geop{x^A} \ubar \geop{x^B} \ubar, \label{E:QUASILINEAREQUATIONSATISFIEDBYLUBAR}\\
			\newuL \muX \ubar & = - \frac{1}{2} (\muX \upmu) \Lunit \ubar + \ReciprocaluLunitAppliedtoTimeFunction[\Lunit,\muX] \ubar -  \MagnitueofinnerproductofnewLandnewuL \angD^A \ubar \left[\geop{x^A},\muX\right] \ubar +  \frac{\muX \upmu}{2\Lunit \ubar} |\angD \ubar|_{\gtorus}^2 +\frac{\upmu}{2\Lunit \ubar} \left( \muX (\gtorus^{-1})^{AB}\right) \geop{x^A} \ubar \geop{x^B} \ubar, \label{E:QUASILINEAREQUATIONSATISFIEDBYMUXUBAR} \\ 
			\begin{split} \label{E:QUASILINEAREQUATIONSATISFIEDBYSMOOTHANGULARDERIVATIVEOFUBAR}
				\newuL \geop{x^A} \ubar & = - \frac{1}{2} \left(\geop{x^A} \upmu\right) \Lunit \ubar +\ReciprocaluLunitAppliedtoTimeFunction \left[\Lunit,\geop{x^A}\right] \ubar + \left[\muX, \geop{x^A}\right] \ubar +  \frac{\geop{x^A} \upmu}{2\Lunit \ubar} |\angD \ubar|_{\gtorus}^2  +\frac{ \upmu}{2\Lunit \ubar} \left( \geop{x^A} (\gtorus^{-1})^{BC}\right) \geop{x^B} \ubar \geop{x^C} \ubar.  
			\end{split}
		\end{align} 
Conversely, if $\Lunit \ubar, \, \muX \ubar,  \, \geop{x^A} \ubar$ solve \eqref{E:QUASILINEAREQUATIONSATISFIEDBYLUBAR}--\eqref{E:QUASILINEAREQUATIONSATISFIEDBYSMOOTHANGULARDERIVATIVEOFUBAR}  and the following condition holds: 
\begin{align} \label{E:INITIALCONDITIONSFORSECONDORDERSYSTEMFORINGOINGEIKONALFUNCTION}
\muX \ubar \big|_{\twoargmumuxtorus{0}{0}} = 0,
\end{align}
then the eikonal equation $(\gfour^{-1})^{\alpha\beta} \p_\alpha \ubar \p_\beta \ubar = 0$ is satisfied. \end{lemma} 
\begin{proof}
	Suppose $\ubar$ is a solution to the eikonal equation, which we rewrite as: 
		\begin{equation} \label{E:MUWEIGHTEDEIKONALEQUATION}
			0 = -\upmu (\Lunit \ubar)^2 -2 \muX \ubar \Lunit \ubar + \upmu |\angD\ubar|_{\gtorus}^2
		\end{equation}
	using \eqref{E:EIKONALEQUATIONWITHOUTTHEMUWEIGHT}. For $Z \in \{\Lunit, \muX, \geop{x^2},\geop{x^3}\}$, we note the following identity satisfied by $Z \ubar$ (see \eqref{E:CONVENIENTIDENTITYFORULUNIT})
		\begin{equation} \label{E:NEWULDERIVATIVEOFZUBAR}
			\newuL Z \ubar = \ReciprocaluLunitAppliedtoTimeFunction \Lunit Z \ubar + \muX Z \ubar - \MagnitueofinnerproductofnewLandnewuL \angD^A \ubar \geop{x^A} Z \ubar.
		\end{equation}
	Differentiating \eqref{E:MUWEIGHTEDEIKONALEQUATION} with respect to $Z$, using the commutator identity $Z\muX = \muX Z + [Z,\muX]$, substituting RHS\,\eqref{E:MUXACTINGONINGOINGEIKONAL} for $\muX \ubar$, solving the resulting expression for $\muX Z \ubar$ and inserting it into \eqref{E:NEWULDERIVATIVEOFZUBAR}, and using the identity $\ReciprocaluLunitAppliedtoTimeFunction = \frac{\upmu \ReciprocalLunitAppliedtoTimeFunction}{\ingoingmu} = \frac{\upmu}{\ingoingmu \Lunit \ubar}$ (see \eqref{E:RATIOOFNULLGEOSICINNERPRODUCTANDFOLIATIONDENSITY} and \eqref{E:LUNITANDULUNITAPPLIEDTOEIKONALANDCARTESIANTIME}), we arrive to the following:
	\begin{equation}
		\label{E:NEWULDERIVATIVEOFZUBARWITHINSERTEDEIKONALEQUATION}
		\newuL Z \ubar  = - \frac{1}{2} Z \upmu \Lunit \ubar + \ReciprocaluLunitAppliedtoTimeFunction [\Lunit,Z] \ubar + [\muX, Z] \ubar - \MagnitueofinnerproductofnewLandnewuL \angD^A \ubar \left[\geop{x^A},Z\right] \ubar +  \frac{Z \upmu}{2\Lunit \ubar} |\angD \ubar|_{\gtorus}^2 +\frac{ \upmu}{2\Lunit \ubar} \left( Z (\gtorus^{-1})^{AB}\right) \geop{x^A} \ubar \geop{x^B} \ubar.
	\end{equation}
	The forward direction of the lemma follows. 
	
	Suppose now that \eqref{E:QUASILINEAREQUATIONSATISFIEDBYLUBAR}--\eqref{E:INITIALCONDITIONSFORSECONDORDERSYSTEMFORINGOINGEIKONALFUNCTION} hold. Then, since the LHS of \eqref{E:QUASILINEAREQUATIONSATISFIEDBYLUBAR}--\eqref{E:QUASILINEAREQUATIONSATISFIEDBYSMOOTHANGULARDERIVATIVEOFUBAR} are equivalent to \eqref{E:NEWULDERIVATIVEOFZUBARWITHINSERTEDEIKONALEQUATION} for the respective choice of $Z$, following the same arguments to prove the forward direction implies: 
\begin{equation} \label{E:SECONDORDERSYSTEMIMPLIESTHATGRADIENTOFUBARISCONSTANT}  - \upmu(\Lunit \ubar)^2 - 2 \muX \ubar \Lunit \ubar + \upmu |\angD\ubar|_{\gtorus}^2 = \text{constant}.
\end{equation}
The condition \eqref{E:INITIALCONDITIONSFORSECONDORDERSYSTEMFORINGOINGEIKONALFUNCTION} implies that the constant in \eqref{E:SECONDORDERSYSTEMIMPLIESTHATGRADIENTOFUBARISCONSTANT} is zero. 
\end{proof}

\begin{remark}[The need for a transport equation for $\muX \ubar$] \label{R:NEEDTODERIVETRANSPORTEQUATIONFORTHEFULL1JETOFINGOINGEIKONAL}

We highlight that a transport equation for $\muX \ubar$ in \eqref{E:QUASILINEAREQUATIONSATISFIEDBYMUXUBAR} is necessary \emph{even though} $\muX \ubar$ is \emph{completely determined} by $\Lunit\ubar, \, \geop{x^2} \ubar, \geop{x^3} \ubar$ for solutions to the eikonal equation, see \eqref{E:MUXINGOINGEIKONALDETERMINEDBYOTHERDERIVATIVES}. The reason is that, if we \emph{only} assume \eqref{E:QUASILINEAREQUATIONSATISFIEDBYLUBAR} and \eqref{E:QUASILINEAREQUATIONSATISFIEDBYSMOOTHANGULARDERIVATIVEOFUBAR}, then we could at best conclude that the quantity $-\upmu (\Lunit \ubar)^2 - 2 \muX \ubar \Lunit \ubar + \upmu |\angD \ubar|_{\gtorus}^2$ is only constant along the $\Lunit$ and $\geop{x^A}$ directions. That is to say, it was crucial for the backwards direction of the proof of Lemma\,\ref{L:QUASILINEARTRANSPORTSYSTEMFORDERIVATIVESOFUBAR} that $\{\Lunit,\muX,\geop{x^2},\geop{x^3}\}$ span the entire tangent space at every point on $\twoargMrough{[\timefunction_0,0],[- \rightu,\leftu]}{0}$.

\end{remark}

\subsection{Gradient data estimates of the ingoing eikonal function}

\label{SS:GRADIENTDATAESTIMATESFORINGOINGEIKONALFUNCTIONONINITIALDATASURFACE}

We now derive data estimates for first derivatives of $\ubar$ along the initial data surface $\datahypfortimefunctiontwoarg{0}{[\timefunction_0,0]}$ with respect to its intrinsic coordinates $(t,x^2,x^3)$ given by \eqref{E:DOMAINOFEMBEDDINGFORXMUEQUALSMINUSKAPPAYSURFACE}--\eqref{E:XMUISMINUSCAPPISAGRAPH}. By \eqref{E:DATAFORINGOINGEIKONALINITIALVALUEPROBLEM}, derivatives of $\ubar$ in directions \emph{tangent} to the co-dimension 1 surface $\datahypfortimefunctiontwoarg{0}{[\timefunction_0,0]}$ will inherit the regularity from $\upmu$. The remaining gradient component that is transverse to $\datahypfortimefunctiontwoarg{0}{[\timefunction_0,0]}$ can be solved for from the eikonal equation \eqref{E:INGOINGEIKONALEQUATION}.

\subsubsection{Intrinsic coordinates, vectorfields, and norms on $\datahypfortimefunctiontwoarg{0}{[\timefunction_0,0]}$} We now properly set up the notation and definitions needed to quantify the regularity of $\ubar$ on the initial data hypersurface. We begin with the following notation for the restriction of $\ubar$ on $\datahypfortimefunctiontwoarg{0}{[\timefunction_0,0]}$ . 

\begin{notation}[Restriction of $\ubar$ on the initial data hypersurface $\datahypfortimefunctiontwoarg{0}{[\timefunction_0,0]}$]\label{N:RESTRICTIONOFUBARONDATASURFACE}
	Let $\ubar$ be the solution to the eikonal function initial value problem \eqref{E:INGOINGEIKONALEQUATIONINITIALVALUEPROBLEM}. Then by $\mr \ubar$ we denote its restriction to the initial data surface:
	\begin{equation} \label{E:RESTRICTIONOFUBARONDATASURFACE}
		\mr \ubar \eqdef 
			\ubar|_{\datahypfortimefunctiontwoarg{0}{[\timefunction_0,0]}}  = -\upmu|_{\datahypfortimefunctiontwoarg{0}{[\timefunction_0,0]}}.
	\end{equation}
	Similarly, we adorn the null-geometry scalar functions defined in \eqref{E:NULLACOUSTICSCALARS} with a ring when considering their values on the initial data surface, e.g.  $\mr \ingoingmu \eqdef \ingoingmu |_{\datahypfortimefunctiontwoarg{0}{[\timefunction_0,0]}}$.
\end{notation}

Next we have the following definition, which introduces a useful set of vectorfields tangent to $\datahypfortimefunctiontwoarg{0}{[\timefunction_0,0]}$. 

\begin{definition}[Vectorfields tangent to the data hypersurface surface $\datahypfortimefunctiontwoarg{0}{[\timefunction_0,0]}$] \label{D:VECTORFIELDSTANGENTTOUBARDATASURFACE} \hfill
\begin{itemize}
	\item The intrinsic $(t,x^2,x^3)$ coordinate partial differentiation on $\datahypfortimefunctiontwoarg{0}{[\timefunction_0,0]}$ given by \eqref{E:DOMAINOFEMBEDDINGFORXMUEQUALSMINUSKAPPAYSURFACE}--\eqref{E:XMUISMINUSCAPPISAGRAPH} is: 
	\begin{subequations}
			\begin{align}
				\datasurfacep{t} & \eqdef \geop{t} - \frac{\geop{t} \muX \upmu}{\geop{u} \muX \upmu} \geop{u}, \label{E:TCOORDINATEPARTIALDERIVATIVEINTRINSICTODATASURFACE} \\
				\datasurfacep{x^2} & \eqdef \geop{x^2} - \frac{\geop{x^2} \muX \upmu}{\geop{u} \muX \upmu} \geop{u},   \label{E:X2COORDINATEPARTIALDERIVATIVEINTRINSICTODATASURFACE} \\
				\datasurfacep{x^3} & \eqdef \geop{x^3} - \frac{\geop{x^3} \muX \upmu}{\geop{u} \muX \upmu} \geop{u}.  \label{E:X3COORDINATEPARTIALDERIVATIVEINTRINSICTODATASURFACE}
			\end{align}
		\end{subequations}
		With $x^0 = t$ and a multiindex $\vec{\alpha} = (\alpha_0,\alpha_2,\alpha_3)$, we denote 
		\[ \mathring{\mathfrak{d}}^{\vec{\alpha}} \eqdef \datasurfacep{t}^{\alpha_0} \datasurfacep{x^2}^{\alpha_2} \datasurfacep{x^3}^{\alpha_3}.\]
	\item We define $\UnitTimeNormalizedUbarDataSurfaceTangentVectorfield$ to be the $\datahypfortimefunctiontwoarg{0}{[\timefunction_0,0]}$-tangent vectorfield given by: 
	\begin{equation} \label{E:FUTUREORIENTEDVECTORFIELDTANGENTTOUBARDATASURFACE}
		\UnitTimeNormalizedUbarDataSurfaceTangentVectorfield \eqdef \Lunit - \frac{\Lunit \muX \upmu}{\muX \muX \upmu}\muX.
	\end{equation} 
	For convenience, throughout, we denote $\coefficientsInVectorfieldsTangetToUbarDataSurface \eqdef \frac{\Lunit \muX \upmu}{\muX \muX \upmu}$ as the coefficient of $\muX$ in \eqref{E:FUTUREORIENTEDVECTORFIELDTANGENTTOUBARDATASURFACE}. That is, $\UnitTimeNormalizedUbarDataSurfaceTangentVectorfield = \Lunit - \coefficientsInVectorfieldsTangetToUbarDataSurface \muX$. 
	\item For $A = 2,3$, we define $\datasurfacegeometricp{x^A}$ to be the $\datahypfortimefunctiontwoarg{0}{[\timefunction_0,0]}$-tangent vectorfield given by: 
	\begin{equation} \label{E:ANGULARVECTORFIELDTANGETTOUBARDATASURFACE}
		\datasurfacegeometricp{x^A} \eqdef \geop{x^A} - \frac{\geop{x^A} \muX \upmu}{\muX \muX \upmu} \muX. 
	\end{equation}
	For convenience, throughout, we denote $\coefficientsInVectorfieldsTangetToUbarDataSurface_A \eqdef \frac{\geop{x^A} \muX \upmu}{\muX \muX \upmu}$ as the coefficient of $\muX$ in \eqref{E:ANGULARVECTORFIELDTANGETTOUBARDATASURFACE}. That is, $\datasurfacegeometricp{x^A} = \Lunit - \coefficientsInVectorfieldsTangetToUbarDataSurface_A \muX$.  
\end{itemize}

\end{definition}

\begin{remark}[The different uses of $\{\datasurfacep{t},\datasurfacep{x^2},\datasurfacep{x^3}\}$ and $\{\UnitTimeNormalizedUbarDataSurfaceTangentVectorfield, \datasurfacegeometricp{x^2},\datasurfacegeometricp{x^3}\}$] \label{E:TWODIFFERENTSPANNINGSETSOFDATAHYPERSURFACE}
We will use  $\{\datasurfacep{t},\datasurfacep{x^2},\datasurfacep{x^3}\}$ to define $C^{m,1}$ H\"older norms for the data on $\datahypfortimefunctiontwoarg{0}{[\timefunction_0,0]}$ since they are the intrinsic coordinate partial derivative vectorfields. The set $\{\UnitTimeNormalizedUbarDataSurfaceTangentVectorfield, \datasurfacegeometricp{x^2},\datasurfacegeometricp{x^3}\}$ is useful because it is easier to solve for $\Lunit \upmu$ and $\muX \upmu$ from the eikonal equation \eqref{E:INGOINGEIKONALEQUATION} with respect to it than it is with respect to $\{\datasurfacep{t},\datasurfacep{x^2},\datasurfacep{x^3}\}$. Nevertheless, we highlight that the following schematic identities hold:
\begin{subequations}
	\begin{align}
		\UnitTimeNormalizedUbarDataSurfaceTangentVectorfield & = \left(1 + \mathcal{O}(\initialsmall)\right) \datasurfacep{t} + \mathcal{O}(\initialsmall) (\datasurfacep{x^2} + \datasurfacep{x^3}), \label{E:UNITTIMENORMALIZEDDATASURFACEVECTORFIELDINTERMSOFDATASURFACECOORDINATEDERIVATIVES} \\
		\datasurfacegeometricp{x^2} & =\left(1 + \mathcal{O}(\initialsmall)\right) \datasurfacep{x^2} + \mathcal{O} (\initialsmall) (\datasurfacep{x^t} + \datasurfacep{x^3} ), \label{E:DATASURFACEANGULARX2GEOMETRICINTERMSOFDATASURFACECOORDINATEDERIVATIVES} \\
		\datasurfacegeometricp{x^3} & =\left(1 + \mathcal{O}(\initialsmall)\right) \datasurfacep{x^3} + \mathcal{O} (\initialsmall) (\datasurfacep{x^t} + \datasurfacep{x^2}).\label{E:DATASURFACEANGULARX3GEOMETRICINTERMSOFDATASURFACECOORDINATEDERIVATIVES}
	\end{align}
\end{subequations}
\end{remark}

\begin{definition}[H\"older norms on $\datahypfortimefunctiontwoarg{0}{[\timefunction_0,0]}$] \label{D:HOLDERNORMSONDATASURFACE} Let $f$ be a scalar function on $\datahypfortimefunctiontwoarg{0}{[\timefunction_0,0]}$, let $m \geq 0$ be an integer. Then we define the H\"older norm of $f$ relative to geometric coordinates (see the embedding \eqref{E:XMUISMINUSCAPPISAGRAPH}): 
	\begin{align} \label{E:HOLDERNORMSONDATASURFACE}
		\begin{split} 
		\| f \|_{C_{\textnormal{data}}^{m,1} \left(\datahypfortimefunctiontwoarg{0}{[\timefunction_0,0]}\right)} & \eqdef \sum_{|\vec{\alpha}| \le m} \max_{(t,x^2,x^3) \in \domainforembeddingdatahypfortimefunctiontwoarg{0}{[0,\mupositive]}} | \mathring{\mathfrak{d}}^{\vec{\alpha}} f(t,h(t,x^2,x^3),x^2,x^3) | \\
		& \ \  + \sum_{|\vec{\alpha}| = m} \sup_{\substack{(t,x^2,x^3),(t_0,x^2_0,x^3_0) \in \domainforembeddingdatahypfortimefunctiontwoarg{0}{[0,\mupositive]} \\ (t,x^2,x^3)\neq(t_0,x^2_0,x^3_0)} } \frac{\left|\mathring{\mathfrak{d}}^{\vec{\alpha}} f(t,h(t,x^2,x^3),x^2,x^3) - \mathring{\mathfrak{d}}^{\vec{\alpha}} f(t_0,h(t_0,x^2_0,x^3_0),x^2_0,x^3_0)\right|}{\textnormal{dist}_{\textnormal{euc}}\left(\big(t,h(t,x^2,x^3),x^2,x^3\big),\big(t_0,h(t_0,x^2_0,x^3_0),x^2_0,x^3_0\big)\right)}.
		\end{split}
	\end{align}\end{definition}

In the following lemma we show that $\Lunit \ubar$ and $\muX \ubar$ can be solved for in the eikonal equation \eqref{E:INGOINGEIKONALEQUATION} in terms of the $\datahypfortimefunctiontwoarg{0}{[\timefunction_0,0]}$-tangent derivatives $\UnitTimeNormalizedUbarDataSurfaceTangentVectorfield \ubar, \, \datasurfacegeometricp{x^2}\ubar, \,\datasurfacegeometricp{x^3} \ubar$.

\begin{lemma}[$\Lunit \ubar$ and $\muX \ubar$ are determined by $\datahypfortimefunctiontwoarg{0}{[\timefunction_0,0]}$-tangent derivatives] \label{L:IDENTITYFORLUBARANDMUXUBARINTERMSOFDATASURFACETANGETDERIVATIVESSOVLEDBYEIKONALEQUAITON}
Let $\{\UnitTimeNormalizedUbarDataSurfaceTangentVectorfield, \datasurfacegeometricp{x^2},\datasurfacegeometricp{x^3}\}$ be the vectorfields from Def.\,\ref{D:VECTORFIELDSTANGENTTOUBARDATASURFACE} and let $\ubar$ be a solution to \eqref{E:INGOINGEIKONALEQUATIONINITIALVALUEPROBLEM}. Then the following identities hold:
\begin{subequations} 
	\begin{align}
		\begin{split} \label{E:LUBARINTERMSOFDATASURFACEDERIVATIVES}
			\Lunit \ubar & =  \left\{ - 2 \coefficientsInVectorfieldsTangetToUbarDataSurface \UnitTimeNormalizedUbarDataSurfaceTangentVectorfield \ubar  + 2 \frac{\upmu}{(\muX \muX \upmu)^2} |\angD \muX \upmu|_{\gtorus}^2 - 2 \frac{\upmu \coefficientsInVectorfieldsTangetToUbarDataSurface}{\muX \muX \upmu} \angD^A \muX \upmu \datasurfacegeometricp{x^A} \ubar \right. \\
			& \ \ - \coefficientsInVectorfieldsTangetToUbarDataSurface
			\left\{\left(2 \UnitTimeNormalizedUbarDataSurfaceTangentVectorfield \ubar + 2\upmu \coefficientsInVectorfieldsTangetToUbarDataSurface \UnitTimeNormalizedUbarDataSurfaceTangentVectorfield \ubar - 2 \frac{\upmu}{\muX \muX \upmu} \angD^A \muX \upmu \datasurfacegeometricp{x^A} \ubar\right)^2  + 4\left( 2 \coefficientsInVectorfieldsTangetToUbarDataSurface  + \upmu \coefficientsInVectorfieldsTangetToUbarDataSurface^2  - \frac{\upmu}{(\muX \muX \upmu)^2} |\angD \muX \upmu|_{\gtorus}^2\right)\left(\upmu \gtorus^{AB} \datasurfacegeometricp{x^A} \ubar \datasurfacegeometricp{x^B} \ubar - \upmu (\UnitTimeNormalizedUbarDataSurfaceTangentVectorfield \ubar)^2\right) \right\}^{1/2}\\
			& \ \ \times \left\{ -4 \coefficientsInVectorfieldsTangetToUbarDataSurface  - 2\upmu \coefficientsInVectorfieldsTangetToUbarDataSurface^2 + \frac{2\upmu}{(\muX \muX \upmu)^2} |\angD \muX \upmu|_{\gtorus}^2\right\}^{-1}
		\end{split} \\ 
		\begin{split} \label{E:MUXUBARINTERMSOFDATASURFACEDERIVATIVES}
			\muX \ubar & = \left\{ 2 \UnitTimeNormalizedUbarDataSurfaceTangentVectorfield \ubar + 2\upmu \coefficientsInVectorfieldsTangetToUbarDataSurface \UnitTimeNormalizedUbarDataSurfaceTangentVectorfield \ubar - 2 \frac{\upmu}{\muX \muX \upmu} \angD^A \muX \upmu \datasurfacegeometricp{x^A} \ubar \right. \\
			& \ \ - 
			\left\{\left(2 \UnitTimeNormalizedUbarDataSurfaceTangentVectorfield \ubar + 2\upmu \coefficientsInVectorfieldsTangetToUbarDataSurface \UnitTimeNormalizedUbarDataSurfaceTangentVectorfield \ubar - 2 \frac{\upmu}{\muX \muX \upmu} \angD^A \muX \upmu \datasurfacegeometricp{x^A} \ubar\right)^2  + 4\left( 2 \coefficientsInVectorfieldsTangetToUbarDataSurface  + \upmu \coefficientsInVectorfieldsTangetToUbarDataSurface^2  - \frac{\upmu}{(\muX \muX \upmu)^2} |\angD \muX \upmu|_{\gtorus}^2\right)\left(\upmu \gtorus^{AB} \datasurfacegeometricp{x^A} \ubar \datasurfacegeometricp{x^B} \ubar - \upmu (\UnitTimeNormalizedUbarDataSurfaceTangentVectorfield \ubar)^2\right) \right\}^{1/2}\\
			& \ \ \times \left\{ -4 \coefficientsInVectorfieldsTangetToUbarDataSurface  - 2\upmu \coefficientsInVectorfieldsTangetToUbarDataSurface^2 + \frac{2\upmu}{(\muX \muX \upmu)^2} |\angD \muX \upmu|_{\gtorus}^2\right\}^{-1}
		\end{split}
	\end{align} 
\end{subequations}
\end{lemma}

\begin{proof} We begin with the following identity, which follows immediately from \eqref{E:FUTUREORIENTEDVECTORFIELDTANGENTTOUBARDATASURFACE}:
\begin{align} \label{E:DEFINITIONOFFUTUREORIENTEDDATASURFACEVECTORFIELDMULTIPLIEDBYLUBAR}
	(\Lunit \ubar)^2 - \coefficientsInVectorfieldsTangetToUbarDataSurface \Lunit \ubar \muX \ubar = (\Lunit \ubar) \UnitTimeNormalizedUbarDataSurfaceTangentVectorfield \ubar.
\end{align}
Solving for $(\Lunit \ubar)^2$ in the eikonal equation \eqref{E:EIKONALEQUATIONWITHOUTTHEMUWEIGHT} and inserting it into \eqref{E:DEFINITIONOFFUTUREORIENTEDDATASURFACEVECTORFIELDMULTIPLIEDBYLUBAR} yields:
\begin{align} \label{E:INSERTINGLUBARSQUAREDFROMEIKONALEQUATIONINTODEFINITIONOFDATASURFACEDERIVATIVE}
	-2 (\Lunit \ubar) X \ubar + |\angD \ubar|^2_{\gtorus} - \coefficientsInVectorfieldsTangetToUbarDataSurface \Lunit \ubar \muX \ubar = (\Lunit \ubar) \UnitTimeNormalizedUbarDataSurfaceTangentVectorfield \ubar.
\end{align}
Using \eqref{E:ANGULARVECTORFIELDTANGETTOUBARDATASURFACE}, we can express $|\angD \ubar|_{\gtorus}^2$ in terms of $\muX$ and $\datahypfortimefunctiontwoarg{0}{[\timefunction_0,0]}$-tangent derivatives as: 
	\begin{align}
		|\angD \ubar|_{\gtorus}^2 = \gtorus^{AB} \datasurfacegeometricp{x^A} \ubar \datasurfacegeometricp{x^B} \ubar  + \frac{2}{\muX \muX \upmu} \muX \ubar \angD^A \muX \upmu \datasurfacegeometricp{x^A} \ubar + \frac{1}{(\muX \muX \upmu)^2} |\angD \muX \upmu|_{\gtorus}^2 (\muX \ubar)^2. \label{E:ANGULARGRADIENTINTERMSOFBREVEXDATASURFACEDERIVATIVES}
	\end{align}
Inserting this expression for $|\angD \ubar|_{\gtorus}^2$ in \eqref{E:INSERTINGLUBARSQUAREDFROMEIKONALEQUATIONINTODEFINITIONOFDATASURFACEDERIVATIVE}, multiplying the resulting identity by $\upmu$, using $\muX = \upmu X$, and re-expressing $\Lunit \ubar =  \UnitTimeNormalizedUbarDataSurfaceTangentVectorfield \ubar + \coefficientsInVectorfieldsTangetToUbarDataSurface \muX \ubar$, we have the following quadratic equation for $\muX \ubar$:
\begin{align}
	\begin{split}
		 \bigg\{ - 2 \coefficientsInVectorfieldsTangetToUbarDataSurface  - \upmu \coefficientsInVectorfieldsTangetToUbarDataSurface^2 & \left.+ \frac{\upmu}{(\muX \muX \upmu)^2} |\angD \muX \upmu|_{\gtorus}^2\right\} (\muX \ubar)^2 + 2 \left\{ -  \UnitTimeNormalizedUbarDataSurfaceTangentVectorfield \ubar -  \upmu  \coefficientsInVectorfieldsTangetToUbarDataSurface \UnitTimeNormalizedUbarDataSurfaceTangentVectorfield \ubar   + \frac{\upmu}{\muX \muX \upmu} \angD^A \muX \upmu \datasurfacegeometricp{x^A} \ubar \right\} \muX \ubar \\
		& \ \ + \upmu \gtorus^{AB} \datasurfacegeometricp{x^A} \ubar \datasurfacegeometricp{x^B} \ubar - \upmu (\UnitTimeNormalizedUbarDataSurfaceTangentVectorfield \ubar)^2 = 0.
	\end{split}
\end{align}
A straightforward application of the quadratic formula yields \eqref{E:MUXUBARINTERMSOFDATASURFACEDERIVATIVES}. Using the now proven \eqref{E:MUXUBARINTERMSOFDATASURFACEDERIVATIVES} for $\muX \ubar$, identity \eqref{E:LUBARINTERMSOFDATASURFACEDERIVATIVES} follows from inserting the RHS\,\eqref{E:MUXUBARINTERMSOFDATASURFACEDERIVATIVES} into $\Lunit \ubar = \UnitTimeNormalizedUbarDataSurfaceTangentVectorfield \ubar  + \coefficientsInVectorfieldsTangetToUbarDataSurface \muX \ubar$ and straightforward algebraic manipulations. We highlight that this procedure yields a complete cancelation of $2 \upmu \coefficientsInVectorfieldsTangetToUbarDataSurface^2 \UnitTimeNormalizedUbarDataSurfaceTangentVectorfield \ubar$, which is a term resulting from multiplying the numerator of the RHS\,\eqref{E:MUXUBARINTERMSOFDATASURFACEDERIVATIVES} by $\coefficientsInVectorfieldsTangetToUbarDataSurface$. We also note that we took the negative root in the quadratic formula in order to be consistent with \eqref{E:TRANSVERSALITYCONDITIONFORINGOINGEIKONALWITHRESPECTOU}.

\end{proof}

\begin{remark}[Loss of regularity for $\ubar$] \label{R:LOSSOFREGULARITYINDATASURFACETRANSVERSALDIRECTIONS}
	From the RHS\,\eqref{E:LUBARINTERMSOFDATASURFACEDERIVATIVES}--\eqref{E:MUXUBARINTERMSOFDATASURFACEDERIVATIVES}, we see that $\Lunit \ubar $ loses a derivative with respect to $\upmu$, \emph{even though they agree on} $\datahypfortimefunctiontwoarg{0}{[\timefunction_0,0]}$ up to a sign by virtue of \eqref{E:DATAFORINGOINGEIKONALINITIALVALUEPROBLEM}. This is a fundamental difference between using $\ubar$ as an eikonal null-coordinate function and using $\timefunction$ from \cite{abbrescia2022emergence} as a time function, which does not lose a derivative \eqref{E:C21ESTIMATESOFTHEROUGHTIMEFUNCTIONINGEOMETRICCOORDINATES} with respect to $\upmu$ due to the compatability condition $0 = \muX \timefunction = \muX \upmu |_{\datahypfortimefunctiontwoarg{0}{[\timefunction_0,0]}}$, \emph{even though} $\mr \ubar$ and $\timefunction$ agree on the initial data defining surface. 
\end{remark}

\begin{remark}[The leading order behavior of $\Lunit \ubar$ and $\muX \ubar$] \label{R:LEADINGORDERFORUBAR} 
It might appear difficult to fully read off properties $\Lunit \ubar$ and $\muX \ubar$ from RHS\,\eqref{E:LUBARINTERMSOFDATASURFACEDERIVATIVES}--\eqref{E:MUXUBARINTERMSOFDATASURFACEDERIVATIVES}. To aid the reader, we refer them to identities \eqref{E:IDENTITYFORLAPXUBARINTERMSOFDATASURFACEDERIVATIVES}--\eqref{E:IDENTITYFORMUXAPXUBARINTERMSOFDATASURFACEDERIVATIVES}. Those identities correspond to an approximate eikonal function $ \, ^{(\textnormal{apx})} \ubar$ which solves $\left(\frac{1}{2} \upmu \Lunit + \muX \right) \, ^{(\textnormal{apx})} \ubar  = 0$ and agrees with $\mr \ubar$ on the initial data surface $\datahypfortimefunctiontwoarg{0}{[\timefunction_0,0]}$. The point is that, in plane symmetry, it is easy to see that $\newuL = \frac{1}{2} \upmu \Lunit + \muX$ and so the leading order behavior of $\Lunit \ubar$ and $\muX \ubar$ from RHS\,\eqref{E:LUBARINTERMSOFDATASURFACEDERIVATIVES}--\eqref{E:MUXUBARINTERMSOFDATASURFACEDERIVATIVES} is given by RHS\,\eqref{E:IDENTITYFORLAPXUBARINTERMSOFDATASURFACEDERIVATIVES}--\eqref{E:IDENTITYFORMUXAPXUBARINTERMSOFDATASURFACEDERIVATIVES} with $\ubar$ in place of $^{(\textnormal{apx})} \ubar$.
\end{remark}

\subsubsection{H\"older estimates for the ingoing eikonal function on the initial data surface} \label{SS:HOLDERESTIMATESFORUBARONDATASURFACE} In this section we derive H\"older estimates for first derivatives of $\ubar$ on $\datahypfortimefunctiontwoarg{0}{[\timefunction_0,0]}$.

\begin{lemma}[The ingoing eikonal function $\ubar$ on $\datahypfortimefunctiontwoarg{0}{[\timefunction_0,0]}$] \label{L:DATAREGULARITYFORINGOINGEIKONALFUNCTION} 

Let $\ubar$ be the ingoing eikonal function solving \eqref{E:INGOINGEIKONALEQUATIONINITIALVALUEPROBLEM} and let $\mr \ubar$ be its restriction to $\datahypfortimefunctiontwoarg{0}{[\timefunction_0,0]}$ as in \eqref{E:RESTRICTIONOFUBARONDATASURFACE}. Then the following identities hold for $A = 2,3$:
\begin{subequations}
	\begin{align} 
		\UnitTimeNormalizedUbarDataSurfaceTangentVectorfield \mr \ubar & = - \Lunit \upmu \big|_{\datahypfortimefunctiontwoarg{0}{[\timefunction_0,0]}},  \label{E:FUTUREDIRECTEDDATASURFACEDERIVATIVEOFOFUBAR} \\ 
		\datasurfacegeometricp{x^A} \mr \ubar & = - \geop{x^A} \upmu \big|_{\datahypfortimefunctiontwoarg{0}{[\timefunction_0,0]}}, \label{E:ANGULARDATASURFACETANGENTIALDERIVATIVESOFUBAR} 
	\end{align}
\end{subequations}
Moreover, the following sharp estimate for $\Lunit \mr \ubar$ and identity for $\muX \mr \ubar$ hold provided that $|\timefunction_0|$ is sufficiently small: 
\begin{subequations}
	\begin{align} 
		-1.01  \le \min_{\datahypfortimefunctiontwoarg{0}{[\timefunction_0,0]}} \frac{\Lunit \mr \ubar}{\Lunit \upmu} & \le \max_{\datahypfortimefunctiontwoarg{0}{[\timefunction_0,0]}} \frac{\Lunit \mr \ubar}{\Lunit \upmu}  \le -0.99, \label{E:LUBARISLIKEMINUSLMUALONGDATASURFACE} \\
		\muX \mr \ubar \big|_{\twoargmumuxtorus{0}{0}} &= 0, \label{E:MUXUBARVANISHESATTHECREASE}
	\end{align}
\end{subequations}
In addition, the following H\"older estimates hold:
\begin{subequations}
	\begin{align} 
		\left\| \Lunit \mr \ubar \right\|_{C_{\textnormal{data}}^{0,1} \left(\datahypfortimefunctiontwoarg{0}{[\timefunction_0,0]}\right)}, \, \left\| \newuL \Lunit \mr \ubar \right\|_{C_{\textnormal{data}}^{0,1} \left(\datahypfortimefunctiontwoarg{0}{[\timefunction_0,0]}\right)}, \, \left\| \newuL\newuL \Lunit \mr \ubar \right\|_{C_{\textnormal{data}}^{0,1} \left(\datahypfortimefunctiontwoarg{0}{[\timefunction_0,0]}\right)} & \le C,\label{E:C01INFINITYESTIMATESFORLUBARONINITIALDATASURFACE} \\
		 \left\|  \muX \mr \ubar \right\|_{C_{\textnormal{data}}^{0,1} \left(\datahypfortimefunctiontwoarg{0}{[\timefunction_0,0]}\right)},\,  \left\| \newuL \muX \mr \ubar \right\|_{C_{\textnormal{data}}^{0,1} \left(\datahypfortimefunctiontwoarg{0}{[\timefunction_0,0]}\right)}, \, \left\| \newuL\newuL \muX \mr \ubar \right\|_{C_{\textnormal{data}}^{0,1} \left(\datahypfortimefunctiontwoarg{0}{[\timefunction_0,0]}\right)} & \le C, \label{E:C01INFINITYESTIMATESFORMUXUBARONINITIALDATASURFACE} \\
		 \left\|  \geop{x^A}  \mr \ubar \right\|_{C_{\textnormal{data}}^{0,1} \left(\datahypfortimefunctiontwoarg{0}{[\timefunction_0,0]}\right)}, \,  \left\| \newuL \geop{x^A}  \mr \ubar \right\|_{C_{\textnormal{data}}^{0,1} \left(\datahypfortimefunctiontwoarg{0}{[\timefunction_0,0]}\right)}, \,  \left\| \newuL\newuL \geop{x^A}  \mr \ubar \right\|_{C_{\textnormal{data}}^{0,1} \left(\datahypfortimefunctiontwoarg{0}{[\timefunction_0,0]}\right)} & \le C \initialsmall.  \label{E:C01INFINITYESTIMATESFORANGULARUBARONINITIALDATASURFACE}
	\end{align}
\end{subequations}

\end{lemma}

\begin{proof}

To prove  \eqref{E:FUTUREDIRECTEDDATASURFACEDERIVATIVEOFOFUBAR}--\eqref{E:ANGULARDATASURFACETANGENTIALDERIVATIVESOFUBAR}, we recall that $\muX \upmu = 0$ on $\datahypfortimefunctiontwoarg{0}{[\timefunction_0,0]}$. This immediately implies $\UnitTimeNormalizedUbarDataSurfaceTangentVectorfield \upmu|_{\datahypfortimefunctiontwoarg{0}{[\timefunction_0,0]}} = \Lunit \upmu|_{\datahypfortimefunctiontwoarg{0}{[\timefunction_0,0]}}$ and $\datasurfacegeometricp{x^A} \upmu|_{\datahypfortimefunctiontwoarg{0}{[\timefunction_0,0]}} = \geop{x^A} \upmu |_{\datahypfortimefunctiontwoarg{0}{[\timefunction_0,0]}}$. The result then follows from \eqref{E:RESTRICTIONOFUBARONDATASURFACE} and straightforward calculations  using the definitions \eqref{E:FUTUREORIENTEDVECTORFIELDTANGENTTOUBARDATASURFACE}--\eqref{E:ANGULARVECTORFIELDTANGETTOUBARDATASURFACE}.

 Estimate \eqref{E:LUBARISLIKEMINUSLMUALONGDATASURFACE} follows from inserting \eqref{E:FUTUREDIRECTEDDATASURFACEDERIVATIVEOFOFUBAR}--\eqref{E:ANGULARDATASURFACETANGENTIALDERIVATIVESOFUBAR} into RHS\,\eqref{E:LUBARINTERMSOFDATASURFACEDERIVATIVES}, using estimate \eqref{E:ACOUSTICVARIABLESARESMALLINCLASSICALDEVELOPMENTFROMROUGHCOORDINATES}, the fact that: 
 \begin{equation} \label{E:MUISSMALLALONGUBARDATASURFACE}
 	\left\{- \upmu(p) \ | \ \datahypfortimefunctiontwoarg{0}{[\timefunction_0,0]}\right\} = \left\{ \timefunction(p) \ | \ \datahypfortimefunctiontwoarg{0}{[\timefunction_0,0]}\right\} \subset [0,\timefunction_0],
\end{equation}
 and taking $\initialsmall$ and $\timefunction_0$ sufficiently small. The proof of \eqref{E:MUXUBARVANISHESATTHECREASE} follows from similar arguments as well as the fact that $\upmu$ and $\muX\upmu$ vanish at the crease $\twoargmumuxtorus{0}{0}$.
 
 Estimates \eqref{E:C01INFINITYESTIMATESFORLUBARONINITIALDATASURFACE}--\eqref{E:C01INFINITYESTIMATESFORMUXUBARONINITIALDATASURFACE} for the undifferentiated $\Lunit \mr \ubar$ and $\muX \mr \ubar$ follow from inserting \eqref{E:FUTUREDIRECTEDDATASURFACEDERIVATIVEOFOFUBAR}--\eqref{E:ANGULARDATASURFACETANGENTIALDERIVATIVESOFUBAR} into RHS\,\eqref{E:LUBARINTERMSOFDATASURFACEDERIVATIVES}--\eqref{E:MUXUBARINTERMSOFDATASURFACEDERIVATIVES} as well as using the $C^{2,1}_{\textnormal{geo}}\left(\twoargMrough{[\timefunction_0,0],[- \rightu,\leftu]}{0}\right)$ bounds for $\upmu$ stated in Theorem\,\ref{T:MAINRESULTSFROMSINGULARBOUNDARYPAPER}. See also Rmk.\,\ref{R:LEADINGORDERFORUBAR}. Using the now proved bound for $\muX \mr \ubar$, the smallness from \eqref{E:ACOUSTICVARIABLESARESMALLINCLASSICALDEVELOPMENTFROMROUGHCOORDINATES}, and \eqref{E:ANGULARVECTORFIELDTANGETTOUBARDATASURFACE} prove \eqref{E:C01INFINITYESTIMATESFORANGULARUBARONINITIALDATASURFACE} for the undifferentiated $\geop{x^A} \mr \ubar$. Estimates \eqref{E:C01INFINITYESTIMATESFORLUBARONINITIALDATASURFACE}--\eqref{E:C01INFINITYESTIMATESFORANGULARUBARONINITIALDATASURFACE} for the $\newuL$-differentiated scalar functions $\Lunit \mr \ubar, \, \muX \mr \ubar, \, \geop{x^A} \mr \ubar$, follow from the undifferentiated estimates and the equations \eqref{E:QUASILINEAREQUATIONSATISFIEDBYLUBAR}--\eqref{E:QUASILINEAREQUATIONSATISFIEDBYSMOOTHANGULARDERIVATIVEOFUBAR}, which provide a precise expression for $\newuL$-derivatives. 
\end{proof}

\begin{remark}[$\ubar$ behaves like $x^1$ at the crease] \label{R:UNWEIGHTEDXUBARATTHECREASE} Using \eqref{E:SMOOTHTORUSMETRICINTERMSOFSIGMATMETRICANDX}, \eqref{E:INGOINGEIKONALEQUATION},  \eqref{E:BOUNDSONLMUINTERESTINGREGIONFROMTAUFOLIATIONS}, \eqref{E:LUBARISLIKEMINUSLMUALONGDATASURFACE}, and \eqref{E:C01INFINITYESTIMATESFORANGULARUBARONINITIALDATASURFACE}, we have the following estimate for $X \ubar$ along the crease:
	\[ X \ubar|_{\twoargmumuxtorus{0}{0}} \approx - \blowupdelta,\]
	which is reminiscent of how $x^1$ behaves near the crease with respect to derivatives of $X$ in the Cartesian differential structure.
\end{remark}

\begin{corollary}[Derivatives of the null-acoustic geometry scalars on $\datahypfortimefunctiontwoarg{0}{[\timefunction_0,0]}$] \label{C:DERIVATIVESOFNULLACOUSTICSCALARSONUBARDATASURFACE} 
Let $\ubar$ be the ingoing eikonal function that solves \eqref{E:INGOINGEIKONALEQUATIONINITIALVALUEPROBLEM} and let $\ingoingmu,\,  \ReciprocalLunitAppliedtoTimeFunction=  1/ \Lunit \ubar, \, \ReciprocaluLunitAppliedtoTimeFunction = 1/ \uLunit u, \, \MagnitueofinnerproductofnewLandnewuL$ be the null-acoustic geometry scalar functions defined in \eqref{E:NULLACOUSTICSCALARS}. 

Then, the following sharp estimates and identities hold, where we use Notation\,\ref{N:RESTRICTIONOFUBARONDATASURFACE}: 
\begin{subequations} 
\begin{align} \label{E:SHARPDATAESTIMATEFORINVERSELUBARONUBARDATASURFACE}
		0.87 \blowupdelta^{-1} \le \min_{\datahypfortimefunctiontwoarg{0}{[\timefunction_0,0]}} \mr \ReciprocalLunitAppliedtoTimeFunction & \le \max_{\datahypfortimefunctiontwoarg{0}{[\timefunction_0,0]}} \mr \ReciprocalLunitAppliedtoTimeFunction \le 1.16 \blowupdelta^{-1},  \\
		1.72 \blowupdelta \le \min_{\datahypfortimefunctiontwoarg{0}{[\timefunction_0,0]}} \mr  \ingoingmu & \le \max_{\datahypfortimefunctiontwoarg{0}{[\timefunction_0,0]}} \mr \ingoingmu \le  2.33 \blowupdelta^{-1}, \label{E:SHARPDATAESTIMATEFORINGOINGMUONUBARDATASURFACE}
\end{align} 
\end{subequations}
\begin{subequations}
\begin{align}
\mr \ReciprocaluLunitAppliedtoTimeFunction, \,  \mr \MagnitueofinnerproductofnewLandnewuL\big|_{\twoargmumuxtorus{0}{0}} & = 0 \label{E:NULLACOUSTICSCALARSTHATAREVANISHINGLIKEMUATTHECREASE} \\
\newuL \mr \ReciprocaluLunitAppliedtoTimeFunction, \, \newuL \mr\MagnitueofinnerproductofnewLandnewuL \big|_{\twoargmumuxtorus{0}{0}} & = 0  \label{E:NULLACOUSTICSCALARSTHATAREVANISHINGLIKEMUXMUATTHECREASE}
\end{align}
\end{subequations}Moreover, the following estimates hold for $M = 1,2$:
  \begin{align}
\left\| \mr \ingoingmu,\, \mr \ReciprocalLunitAppliedtoTimeFunction, \, \mr \ReciprocaluLunitAppliedtoTimeFunction, \, \mr \MagnitueofinnerproductofnewLandnewuL\right\|_{C_{\textnormal{data}}^{0,1} \left(\datahypfortimefunctiontwoarg{0}{[\timefunction_0,0]}\right)} & \le C, \label{E:HOLDERNORMSONUBARDATASURFACEOFTHENULLACOUSTICSCALARS} \\
 \left\| \newuL^M \mr \ingoingmu,\,  \newuL^M  \mr\ReciprocalLunitAppliedtoTimeFunction, \, \newuL^M \mr \ReciprocaluLunitAppliedtoTimeFunction, \, \newuL^M \mr \MagnitueofinnerproductofnewLandnewuL\right\|_{C_{\textnormal{data}}^{0,1} \left(\datahypfortimefunctiontwoarg{0}{[\timefunction_0,0]}\right)} & \le C. \label{E:HOLDERNORMSONUBARDATASURFACEOFNEWULDERIVATIVESOFTHENULLACOUSTICSCALARS}
\end{align}

\end{corollary}

\begin{proof}

Since $\ReciprocalLunitAppliedtoTimeFunction= 1/\Lunit \ubar$, estimate \eqref{E:SHARPDATAESTIMATEFORINVERSELUBARONUBARDATASURFACE} follows from \eqref{E:BOUNDSONLMUINTERESTINGREGIONFROMTAUFOLIATIONS} and \eqref{E:LUBARISLIKEMINUSLMUALONGDATASURFACE}. Similarly, \eqref{E:SHARPDATAESTIMATEFORINGOINGMUONUBARDATASURFACE} holds from \eqref{E:INGOINGMUDETERMINEDBYEVERYTHINGELSE}, the smallness of \eqref{E:C01INFINITYESTIMATESFORANGULARUBARONINITIALDATASURFACE}, and the now proved \eqref{E:SHARPDATAESTIMATEFORINVERSELUBARONUBARDATASURFACE}.

The identities in \eqref{E:NULLACOUSTICSCALARSTHATAREVANISHINGLIKEMUATTHECREASE} follow immediately from estimates \eqref{E:SHARPDATAESTIMATEFORINVERSELUBARONUBARDATASURFACE}--\eqref{E:SHARPDATAESTIMATEFORINGOINGMUONUBARDATASURFACE}, from the vanishing of $\upmu$ along the crease, and $\MagnitueofinnerproductofnewLandnewuL = \upmu \ReciprocalLunitAppliedtoTimeFunction, \, \ReciprocaluLunitAppliedtoTimeFunction = \MagnitueofinnerproductofnewLandnewuL/\ingoingmu$, see \eqref{E:RATIOOFNULLGEOSICINNERPRODUCTANDFOLIATIONDENSITY}. 

We now prove the identity in \eqref{E:NULLACOUSTICSCALARSTHATAREVANISHINGLIKEMUXMUATTHECREASE} for $\newuL \, \MagnitueofinnerproductofnewLandnewuL$. Upon differentiating $\MagnitueofinnerproductofnewLandnewuL = \upmu \ReciprocalLunitAppliedtoTimeFunction$ with respect to $\newuL$, it remains to show that $\newuL  \MagnitueofinnerproductofnewLandnewuL  = \newuL \upmu \ReciprocalLunitAppliedtoTimeFunction+ \upmu \newuL \ReciprocalLunitAppliedtoTimeFunction$ vanishes at the crease. Using \eqref{E:CONVENIENTIDENTITYFORULUNIT} and \eqref{E:NULLACOUSTICSCALARSTHATAREVANISHINGLIKEMUATTHECREASE}, one easily sees $\newuL \upmu|_{\twoargmumuxtorus{0}{0}}  = 0 $. The proof follows since $\newuL  \ReciprocalLunitAppliedtoTimeFunction= \newuL \left(\frac{1}{\Lunit \ubar}\right)$ is $\mathcal{O}(1)$ at the crease, which is an immediate consequence of \eqref{E:QUASILINEAREQUATIONSATISFIEDBYLUBAR} and the estimates of Lemma\,\ref{L:DATAREGULARITYFORINGOINGEIKONALFUNCTION}. 

Similarly, the identity in \eqref{E:NULLACOUSTICSCALARSTHATAREVANISHINGLIKEMUXMUATTHECREASE} for $\newuL \ReciprocaluLunitAppliedtoTimeFunction$ follows upon differentiating $\ReciprocaluLunitAppliedtoTimeFunction = \MagnitueofinnerproductofnewLandnewuL/\ingoingmu$ with respect to $\newuL$, using the already proved vanishing of $\newuL \upmu, \, \newuL \, \MagnitueofinnerproductofnewLandnewuL$ at the crease whenever the differentiation falls on $\upmu$ or $ \MagnitueofinnerproductofnewLandnewuL$, and using \eqref{E:INGOINGMUDETERMINEDBYEVERYTHINGELSE} and the identities \eqref{E:QUASILINEAREQUATIONSATISFIEDBYLUBAR}--\eqref{E:QUASILINEAREQUATIONSATISFIEDBYSMOOTHANGULARDERIVATIVEOFUBAR} whenever the derivative falls on $\ingoingmu$ to prove $\newuL \ingoingmu = \mathcal{O}(1)$ at the crease. 

Finally, the bounds in \eqref{E:HOLDERNORMSONUBARDATASURFACEOFTHENULLACOUSTICSCALARS}--\eqref{E:HOLDERNORMSONUBARDATASURFACEOFNEWULDERIVATIVESOFTHENULLACOUSTICSCALARS} follow from similar ideas to the above as well as the proven estimates \eqref{E:C01INFINITYESTIMATESFORLUBARONINITIALDATASURFACE}--\eqref{E:C01INFINITYESTIMATESFORMUXUBARONINITIALDATASURFACE}.
\end{proof}

In the following lemma we use the estimates derived for the ingoing eikonal function $\ubar$ and associated null-acoustic scalars to quantify the behavior of $\upmu$ at the crease near the singularity in terms of the null-acoustic geometry. In particular, we will show that $\newuL$ is \emph{transverse} to the data hypersurface $\datahypfortimefunctiontwoarg{0}{[\timefunction_0,0]}$, even at the crease.

\begin{lemma}[Behavior of $\upmu$ on $\twoargmumuxtorus{0}{0}$ and  $\datahypfortimefunctiontwoarg{0}{[\timefunction_0,0]}$ in terms of $\newuL$] \label{L:SHARPBEHAVIOROFMUANDROUGHTIMEFUNCTIONNEARTHECREASEINNEWULDIRECTIONS}
The following identities hold at the crease:
\begin{align}
	\newuL \upmu  \big|_{\twoargmumuxtorus{0}{0}} & = 0, & \newuL \newuL \upmu \big|_{\twoargmumuxtorus{0}{0}} = \muX \muX \upmu \big|_{\twoargmumuxtorus{0}{0}}. \label{E:BEHAVIOROFMUATTHECREASEINDIRECTIONSOFNEWUL}
\end{align}In addition, for $\mulevelsetvalue \in [0,\mupositive]$, the following estimate holds on $\twoargmumuxtorus{\mulevelsetvalue}{0}$:
\begin{align}
	\newuL \upmu\big|_{\twoargmumuxtorus{\mulevelsetvalue}{0}} =  \left\{\tfrac{1}{2}\Lunit \upmu + \mathcal{O}(\initialsmall)\right\} \mulevelsetvalue. \label{E:NEWULMUISNEGATIVEONDATAHYPERSURFACE}
\end{align}Moreover, the following estimates hold along the crease and data hypersurface for $\ubar$, where $\secondtransversalderivativemulowerbound > 0$ is the constant from \eqref{E:MUTRANSVERSALCONVEXITY}:
\begin{subequations}
	\begin{align}
		\frac{\secondtransversalderivativemulowerbound}{2} \le \min_{\twoargmumuxtorus{0}{0}} \newuL \newuL \upmu & \le \max_{\twoargmumuxtorus{0}{0}} \newuL \newuL \upmu \le  \frac{2}{\secondtransversalderivativemulowerbound},  \label{E:MUISSTRICTLYINCREASINGALONGTHECREASEALONGULUBAR} \\
		\frac{\secondtransversalderivativemulowerbound}{3} \le \min_{\datahypfortimefunctiontwoarg{0}{[\timefunction_0,0]}} \left\{ \newuL \muX \upmu, \newuL \newuL \upmu\right\} & \le \max_{\datahypfortimefunctiontwoarg{0}{[\timefunction_0,0]}} \left\{ \newuL \muX \upmu, \newuL \newuL \upmu\right\} \le  \frac{3}{\secondtransversalderivativemulowerbound}. \label{E:NEWULISTRANSVERSALTOUBARDATASURFACE}
	\end{align}
\end{subequations}In particular, \underline{\textbf{$\newuL$ is transverse to the data hypersurface $\datahypfortimefunctiontwoarg{0}{[\timefunction_0,0]}$, even at the crease.}}

\end{lemma}

\begin{proof}

Using \eqref{E:CONVENIENTIDENTITYFORULUNIT} and \eqref{E:NULLACOUSTICSCALARSTHATAREVANISHINGLIKEMUATTHECREASE}, one easily sees $\newuL \upmu|_{\twoargmumuxtorus{0}{0}}  = 0 $. Next we express $\newuL \newuL \upmu$ as $\newuL (\ReciprocaluLunitAppliedtoTimeFunction\Lunit \upmu + \muX \upmu - \MagnitueofinnerproductofnewLandnewuL \angD^A \ubar \geop{x^A}\upmu)$, again using \eqref{E:CONVENIENTIDENTITYFORULUNIT}. Next we apply the Leibniz with the following procedure: we express $\newuL$ as the RHS\,\eqref{E:CONVENIENTIDENTITYFORULUNIT} then a $\newuL$ differentiation falls on derivatives of $\upmu$, but make no substitutions when $\newuL$ falls on $\ReciprocaluLunitAppliedtoTimeFunction, \, \MagnitueofinnerproductofnewLandnewuL$, or $\newangD^A \ubar$. Then using \eqref{E:NULLACOUSTICSCALARSTHATAREVANISHINGLIKEMUATTHECREASE}--\eqref{E:NULLACOUSTICSCALARSTHATAREVANISHINGLIKEMUXMUATTHECREASE} and \eqref{E:QUASILINEAREQUATIONSATISFIEDBYSMOOTHANGULARDERIVATIVEOFUBAR}, the rest of \eqref{E:BEHAVIOROFMUATTHECREASEINDIRECTIONSOFNEWUL} follows. The desired estimate \eqref{E:MUISSTRICTLYINCREASINGALONGTHECREASEALONGULUBAR} immediately follows from \eqref{E:MUTRANSVERSALCONVEXITY}. 

\eqref{E:NEWULMUISNEGATIVEONDATAHYPERSURFACE} follows from \eqref{E:IDENTITYFORRECIPROCALULUNITAPPLIEDTOTIMEFUNCTION},  \eqref{E:CONVENIENTIDENTITYFORULUNIT}, the identity $\MagnitueofinnerproductofnewLandnewuL = \upmu \ReciprocalLunitAppliedtoTimeFunction$, the estimates \eqref{E:C01INFINITYESTIMATESFORLUBARONINITIALDATASURFACE}, \eqref{E:C01INFINITYESTIMATESFORANGULARUBARONINITIALDATASURFACE}, and the fact that $\upmu|_{\twoargmumuxtorus{\mulevelsetvalue}{0}} = \mulevelsetvalue$. 

We now prove \eqref{E:NEWULISTRANSVERSALTOUBARDATASURFACE}. We first express $\newuL \muX \upmu$ as $\ReciprocaluLunitAppliedtoTimeFunction\Lunit \muX \upmu + \muX\muX  \upmu - \MagnitueofinnerproductofnewLandnewuL \angD^A \ubar \geop{x^A}\muX \upmu$ using \eqref{E:CONVENIENTIDENTITYFORULUNIT}. The desired estimate for $\newuL \muX \upmu$ is then a consequence of  transversal convexity \eqref{E:MUTRANSVERSALCONVEXITY}, the identities $\MagnitueofinnerproductofnewLandnewuL = \upmu \ReciprocalLunitAppliedtoTimeFunction, \, \ReciprocaluLunitAppliedtoTimeFunction = \MagnitueofinnerproductofnewLandnewuL/\ingoingmu$ (see \eqref{E:RATIOOFNULLGEOSICINNERPRODUCTANDFOLIATIONDENSITY}), the bounds \eqref{E:SHARPDATAESTIMATEFORINVERSELUBARONUBARDATASURFACE}--\eqref{E:SHARPDATAESTIMATEFORINGOINGMUONUBARDATASURFACE}, and \eqref{E:MUISSMALLALONGUBARDATASURFACE} by making $|\timefunction_0|$ sufficiently small. To prove the estimate for $\newuL \newuL \upmu$, we again write $\newuL \newuL \upmu = \newuL (\ReciprocaluLunitAppliedtoTimeFunction\Lunit \upmu + \muX \upmu - \MagnitueofinnerproductofnewLandnewuL \angD^A \ubar \geop{x^A}\upmu)$ and focus on the first and last terms. The last term is of size $\mathcal{O}(\initialsmall^2)$ by \eqref{E:ACOUSTICVARIABLESARESMALLINCLASSICALDEVELOPMENTFROMROUGHCOORDINATES} and \eqref{E:C01INFINITYESTIMATESFORANGULARUBARONINITIALDATASURFACE}. To estimate the first term, we expand $\newuL (\ReciprocaluLunitAppliedtoTimeFunction \Lunit \upmu) = \ReciprocaluLunitAppliedtoTimeFunction \newuL \Lunit \upmu + (\newuL \ReciprocaluLunitAppliedtoTimeFunction) \Lunit \upmu$. The term $\ReciprocaluLunitAppliedtoTimeFunction \newuL \Lunit \upmu$ is of size $O(|\timefunction_0|)$ along $\datahypfortimefunctiontwoarg{0}{[\timefunction_0,0]}$. Using \eqref{E:RATIOOFNULLGEOSICINNERPRODUCTANDFOLIATIONDENSITY}, \eqref{E:INGOINGMUDETERMINEDBYEVERYTHINGELSE}, \eqref{E:MUISSMALLALONGUBARDATASURFACE}, and \eqref{E:NEWULMUISNEGATIVEONDATAHYPERSURFACE}, we see that $\newuL \ReciprocaluLunitAppliedtoTimeFunction = \frac{1}{2} \newuL \upmu + \mathcal{O}(\initialsmall) = \mathcal{O}(|\timefunction_0|) + \mathcal{O}(\initialsmall)$ along $\datahypfortimefunctiontwoarg{0}{[\timefunction_0,0]}$. The desired estimate for $\newuL\newuL\upmu$ then follows by making $|\timefunction_0|$ and $\initialsmall$ sufficiently small.

\end{proof}

\begin{lemma}[Properties and identities of $\UnitTimeNormalizedUbarDataSurfaceTangentVectorfield$ and $\datasurfacegeometricp{x^A}$] \label{L:PROPERTIESANDIDENTITIESOFVECTORFIELDSTANGENTTOUBARDATASURFACE} \hfill
\medskip

\noindent \underline{\textbf{Commutator identities:}} The following identities hold:
\begin{align}
	[\muX, \UnitTimeNormalizedUbarDataSurfaceTangentVectorfield] & = [\muX, L] - \muX(\coefficientsInVectorfieldsTangetToUbarDataSurface) \muX, & & [\muX, \datasurfacegeometricp{x^A}] = \left[\muX,\geop{x^A}\right] - \muX(\coefficientsInVectorfieldsTangetToUbarDataSurface_A) \muX. \label{E:COMMUTATORSOFMUXANDUBARDATASURFACEVECTORFIELDS}
\end{align}\medskip 

\noindent \noindent \underline{\textbf{Second order $\datahypfortimefunctiontwoarg{0}{[\timefunction_0,0]}$-transversal derivative operators:}} The following identities holds:
\begin{subequations}
	\begin{align} 
		\geop{x^A} \muX & = \datasurfacegeometricp{x^A} \muX + \coefficientsInVectorfieldsTangetToUbarDataSurface_A \muX \muX, \label{E:IDENTITYFORANGULARMUXINTERSOFUBARDATASURFACETANGTIALDERIVATIVES} \\
		\Lunit \muX & = \UnitTimeNormalizedUbarDataSurfaceTangentVectorfield \muX + \coefficientsInVectorfieldsTangetToUbarDataSurface \muX\muX, \label{E:IDENTITYFORLMUXINTERMSOFUBARDATASURFACETANGENTDERIVATIVES} \\
		\geop{x^A} \Lunit & = \datasurfacegeometricp{x^A} \Lunit + \coefficientsInVectorfieldsTangetToUbarDataSurface_A \UnitTimeNormalizedUbarDataSurfaceTangentVectorfield \muX +  \coefficientsInVectorfieldsTangetToUbarDataSurface \coefficientsInVectorfieldsTangetToUbarDataSurface_A \muX \muX + \coefficientsInVectorfieldsTangetToUbarDataSurface_A [\muX,\Lunit] \label{E:IDENTITYFORANGULARLINTERSOFUBARDATASURFACETANGENTDERIVATIVES}, \\ 
			\Lunit \Lunit &  = \UnitTimeNormalizedUbarDataSurfaceTangentVectorfield L +  \coefficientsInVectorfieldsTangetToUbarDataSurface \UnitTimeNormalizedUbarDataSurfaceTangentVectorfield \muX + \coefficientsInVectorfieldsTangetToUbarDataSurface [\muX,L]  +  2 \coefficientsInVectorfieldsTangetToUbarDataSurface \muX(\coefficientsInVectorfieldsTangetToUbarDataSurface) \muX, \\
			\geop{x^A}\geop{x^B} & = \datasurfacegeometricp{x^A}\geop{x^B} + \coefficientsInVectorfieldsTangetToUbarDataSurface_A \datasurfacegeometricp{x^B} \muX + \coefficientsInVectorfieldsTangetToUbarDataSurface_A \coefficientsInVectorfieldsTangetToUbarDataSurface_B \muX \muX +  \coefficientsInVectorfieldsTangetToUbarDataSurface_A \left[\muX,\geop{x^B}\right] \label{E:IDENTITYFORTWOANGULARINTERMSOFUBARDATASURFACETANGETDERIVATIVES} 
	\end{align}
\end{subequations}
 
\noindent \underline{\textbf{First and second order $\datahypfortimefunctiontwoarg{0}{[\timefunction_0,0]}$-tangential derivatives of $\mr \ubar$:}} The following identities hold, where we follow Notation\,\ref{N:RESTRICTIONOFUBARONDATASURFACE}:
\begin{subequations}
	\begin{align}
		\UnitTimeNormalizedUbarDataSurfaceTangentVectorfield \mr \ubar & = - \Lunit \upmu \big|_{\datahypfortimefunctiontwoarg{0}{[\timefunction_0,0]}},  \label{E:FUTUREDIRECTEDDATASURFACEDERIVATIVEOFOFUBARRESTATED} \\ 
		\datasurfacegeometricp{x^A} \mr \ubar & = - \geop{x^A} \upmu \big|_{\datahypfortimefunctiontwoarg{0}{[\timefunction_0,0]}}, \label{E:ANGULARDATASURFACETANGENTIALDERIVATIVESOFUBARRESTATED} \\
		\UnitTimeNormalizedUbarDataSurfaceTangentVectorfield \UnitTimeNormalizedUbarDataSurfaceTangentVectorfield \mr \ubar & = - \Lunit \Lunit \upmu + \coefficientsInVectorfieldsTangetToUbarDataSurface \muX \Lunit \upmu + \coefficientsInVectorfieldsTangetToUbarDataSurface  \Lunit \muX  \upmu -\coefficientsInVectorfieldsTangetToUbarDataSurface^2 \muX\muX \upmu \big|_{\datahypfortimefunctiontwoarg{0}{[\timefunction_0,0]}}, \label{E:TWOFUTUREDIRECTEDDATASURFACEDERIVATIVESOFUBAR} \\
		 \datasurfacegeometricp{x^A}  \UnitTimeNormalizedUbarDataSurfaceTangentVectorfield \mr \ubar & = - \geop{x^A} \Lunit \upmu + \coefficientsInVectorfieldsTangetToUbarDataSurface_A \muX \Lunit \upmu + \coefficientsInVectorfieldsTangetToUbarDataSurface \geop{x^A} \muX \upmu - \coefficientsInVectorfieldsTangetToUbarDataSurface \coefficientsInVectorfieldsTangetToUbarDataSurface_A \muX\muX \upmu  \big|_{\datahypfortimefunctiontwoarg{0}{[\timefunction_0,0]}}, \label{E:FUTUREDIRECTEDANDANGULARSECONDORDERDATASURFACETANGENTIALDERIVATIVESOFUBAR} \\
		\UnitTimeNormalizedUbarDataSurfaceTangentVectorfield \datasurfacegeometricp{x^A} \mr \ubar & = - \Lunit \geop{x^A} \upmu + \coefficientsInVectorfieldsTangetToUbarDataSurface \muX \geop{x^A} \upmu +  \coefficientsInVectorfieldsTangetToUbarDataSurface_A \Lunit \muX \upmu - \coefficientsInVectorfieldsTangetToUbarDataSurface \coefficientsInVectorfieldsTangetToUbarDataSurface_A \muX\muX \upmu\big|_{\datahypfortimefunctiontwoarg{0}{[\timefunction_0,0]}}, \label{E:ANGULARANDFUTUREDIRECTEDSECONDORDERDATASURFACETANGENTIALDERIVATIVESOFUBAR}  \\
	 	\datasurfacegeometricp{x^A}\datasurfacegeometricp{x^B} \mr \ubar & = - \geop{x^A}\geop{x^B} \upmu + \coefficientsInVectorfieldsTangetToUbarDataSurface_A \muX \geop{x^B} \upmu   + \coefficientsInVectorfieldsTangetToUbarDataSurface_B \geop{x^A}\muX  \upmu  - \coefficientsInVectorfieldsTangetToUbarDataSurface_A \coefficientsInVectorfieldsTangetToUbarDataSurface_B \muX \muX \upmu\big|_{\datahypfortimefunctiontwoarg{0}{[\timefunction_0,0]}}. \label{E:TWOANGULARDATASURFACEDERIVATIVESOFUBAR}
	\end{align}
\end{subequations}	

\end{lemma}
\begin{proof}
Identities \eqref{E:COMMUTATORSOFMUXANDUBARDATASURFACEVECTORFIELDS}--\eqref{E:IDENTITYFORTWOANGULARINTERMSOFUBARDATASURFACETANGETDERIVATIVES} follow from straightforward calculations using the definitions \eqref{E:FUTUREORIENTEDVECTORFIELDTANGENTTOUBARDATASURFACE}--\eqref{E:ANGULARVECTORFIELDTANGETTOUBARDATASURFACE}.

Identities \eqref{E:FUTUREDIRECTEDDATASURFACEDERIVATIVEOFOFUBARRESTATED}--\eqref{E:ANGULARDATASURFACETANGENTIALDERIVATIVESOFUBARRESTATED} are simply a restatement of \eqref{E:FUTUREDIRECTEDDATASURFACEDERIVATIVEOFOFUBAR}--\eqref{E:ANGULARDATASURFACETANGENTIALDERIVATIVESOFUBAR}. The other identities follow from similar calculations used to prove those, which we omit.
\end{proof}

\begin{remark}[Control of the full set of second order partial derivatives of $\ubar$ on $\datahypfortimefunctiontwoarg{0}{[\timefunction_0,0]}$] \label{R:CONTROLOFALLSECONDORDERPARTIALDERIVATIVESOFUBAR}
Identities \eqref{E:FUTUREDIRECTEDDATASURFACEDERIVATIVEOFOFUBAR}--\eqref{E:TWOANGULARDATASURFACEDERIVATIVESOFUBAR} demonstrate how all second order derivatives of $\ubar$ in directions \emph{tangent} to the initial data hypersurface $\datahypfortimefunctiontwoarg{0}{[\timefunction_0,0]}$ are controlled by derivatives of $\upmu$. Hence, to control the full set of second order partial derivatives of $\ubar$ on $\datahypfortimefunctiontwoarg{0}{[\timefunction_0,0]}$, identities \eqref{E:IDENTITYFORANGULARMUXINTERSOFUBARDATASURFACETANGTIALDERIVATIVES}--\eqref{E:TWOANGULARDATASURFACEDERIVATIVESOFUBAR} imply that it suffices to get control of $\UnitTimeNormalizedUbarDataSurfaceTangentVectorfield \muX \mr \ubar$, $\muX \muX \ubar$, and $\muX \mr \ubar$ in terms of $\datahypfortimefunctiontwoarg{0}{[\timefunction_0,0]}$-tangential derivatives. This can be done by differentiating equations \eqref{E:LUBARINTERMSOFDATASURFACEDERIVATIVES}--\eqref{E:MUXUBARINTERMSOFDATASURFACEDERIVATIVES}, which are equivalent to the eikonal equation \eqref{E:INGOINGEIKONALEQUATION}.
\end{remark}

\subsection{Construction of the ingoing eikonal function} \label{SS:ANALYSISOFINGOINGEIKONALFUNCTION}

The main goal of this section is to construct the ingoing eikonal function by a contraction mapping argument. Lemma\,\ref{L:THEFLOWOFREGULARINGOINGNEWUL} provides a  ``fixed neighborhood'' $^{(\textnormal{reg})}\mathscr{D}$ of $\datahypfortimefunctiontwoarg{0}{[\timefunction_0,0]}$ used to define the contraction space, see \eqref{E:DEFINITIONOFCONTRACTIONSPACEFORINGOINGEIKONAL}. The construction of  $^{(\textnormal{reg})}\mathscr{D}$ depends on the $\gfour$-null vectorfield $^{(\textnormal{reg})}\newuL$ defined below.

\begin{definition}[Regular ingoing $\gfour$-null vectorfield] \label{D:REGULARINGOINGNULLVECTORFIELD}
Let $\olduLunit$ be the vectorfield defined in \eqref{E:OLDULUNIT}. Then we define $^{(\textnormal{reg})}\newuL$ to be the vectorfield:
\begin{align} \label{E:REGULARINGOINGNULLVECTORFIELD}
	^{(\textnormal{reg})}\newuL \eqdef \frac{1}{2} \upmu  \olduLunit = \frac{1}{2}\upmu \Lunit + \muX.
\end{align}
\end{definition}

\begin{lemma}[Properties of $^{(\textnormal{reg})}\newuL$, its flow, and an approximate eikonal function] \label{L:THEFLOWOFREGULARINGOINGNEWUL} \hfill

\noindent \underline{\textbf{Properties of $^{(\textnormal{reg})}\newuL$:}}

Let $^{(\textnormal{reg})}\newuL$ be as in \eqref{E:REGULARINGOINGNULLVECTORFIELD}. Then $^{(\textnormal{reg})}\newuL$ is a $\gfour$-null vectorfield. Moreover, in the differential structure of the geometric coordinates, $^{(\textnormal{reg})}\newuL$ is a $\nullhyparg{u}$-transverse vectorfield of class $C^{1,1}_{\textnormal{geo}}\left(\twoargMrough{[\timefunction_0,0],[- \rightu,\leftu]}{0}\right)$. In addition, $^{(\textnormal{reg})}\newuL$ is also transverse to $\datahypfortimefunctiontwoarg{0}{[\timefunction_0,0]}$.

\medskip

\noindent \underline{\textbf{The flow of $^{(\textnormal{reg})}\newuL$:}}

Denote the flow map of $^{(\textnormal{reg})}\newuL$ by $^{(\textnormal{reg})}\Lambda_{\Delta s}$. That is, the $(t,u,x^2,x^3)$-coordinate components of $^{(\textnormal{reg})}\newuL$ solve: 
\begin{align}
	\frac{\p}{\p \Delta s} \, ^{(\textnormal{reg})}\Lambda_{\Delta s} (t,u,x^2,x^3) = \,^{(\textnormal{reg})}\newuL \circ \, ^{(\textnormal{reg})}\Lambda_{\Delta s}(t,u,x^2,x^3), & & ^{(\textnormal{reg})}\Lambda_{0}  (t,u,x^2,x^3) = (t,u,x^2,x^3). \label{E:FLOWMAPOFREGULARNEWUL}
\end{align}
Let $\domainforembeddingdatahypfortimefunctiontwoarg{0}{[0,\mupositive]}, \, \embeddingdatahypfortimefunctionarg{0}, \, \scalarembeddingdatahypfortimefunctionarg{0}$ be the set and maps defined in Prop.\,\ref{P:PROPERTIESOFMUXMUZEROLEVELSET}. Then there exists a unique maximal flowout domain $^{(\textnormal{reg})} \mathscr{F} \subset \R \times \domainforembeddingdatahypfortimefunctiontwoarg{0}{[0,\mupositive]}$ and flow $^{(\textnormal{reg})}F$ given by: 
\begin{subequations}
	\begin{align}
		\begin{split}\label{E:DEFINITIONOFFLOWOUTDOMAINFORREGULARNEWUL}
			^{(\textnormal{reg})}\mathscr{F} & \eqdef \left\{ (\Delta s, t,x^2,x^3) \ | \ (t,x^2,x^3) \in \domainforembeddingdatahypfortimefunctiontwoarg{0}{[0,\mupositive]}  \right. \\
			& \ \ \left. \text{ and } \Delta s \in [\, ^{(\textnormal{reg})}\Delta a(t,x^2,x^3) - \scalarembeddingdatahypfortimefunctionarg{0}(t,x^2,x^3), \, ^{(\textnormal{reg})}\Delta b(t,x^2,x^3) - \scalarembeddingdatahypfortimefunctionarg{0}(t,x^2,x^3)] \right\},
		\end{split} \\
		^{(\textnormal{reg})} F(\Delta s,t,x^2,x^3) & \eqdef \, ^{(\textnormal{reg})}\Lambda_{\Delta s} \circ  \embeddingdatahypfortimefunctionarg{0}(t, x^2,x^3), \label{E:DEFINITIONOFFLOWFORREGULARNEWUL} 
	\end{align} 
\end{subequations}
for some $C^{1,1}$ functions  $^{(\textnormal{reg})}\Delta b, \, ^{(\textnormal{reg})} \Delta a \colon \domainforembeddingdatahypfortimefunctiontwoarg{0}{[0,\mupositive]} \to \R$. Set: 
\begin{align}
	^{(\textnormal{reg})} \mathscr{D} \eqdef \, ^{(\textnormal{reg})} F \left(\, ^{(\textnormal{reg})}\mathscr{F} \right) \label{E:FLOWOUTOFINITIALDATAHYPERSURFACE}
\end{align}
as the flowout of $\datahypfortimefunctiontwoarg{0}{[\timefunction_0,0]}$ by $^{(\textnormal{reg})}\newuL$. Then $^{(\textnormal{reg})} F$ is a $C^{1,1}(^{(\textnormal{reg})}\mathscr{F})$-diffeomorphism with inverse  $^{(\textnormal{reg})} F^{-1} \in C^{1,1}_{\textnormal{geo}}(^{(\textnormal{reg})}\mathscr{D})$ satisfying:
\begin{subequations}
	\begin{align}
		\left\|\, ^{(\textnormal{reg})} F \right\|_{C^{1,1}(\, ^{(\textnormal{reg})}\mathscr{F})} \lesssim 1, \label{E:C11ESTIMATEFORFLOWOFREGULARNEWUL} \\ 
		\left\|\, ^{(\textnormal{reg})} F^{-1} \right\|_{C^{1,1}_{\textnormal{geo}}(\, ^{(\textnormal{reg})}\mathscr{D})} \lesssim 1. \label{E:C11ESTIMATEFORINVERSEFLOWOFREGULARNEWUL}
	\end{align}
\end{subequations}Moreover, if $|\timefunction_0|$ is sufficiently small, for any  $(t,x^2,x^3) \in \domainforembeddingdatahypfortimefunctiontwoarg{0}{[0,\mupositive]}$, we have: 
\begin{subequations}
	\begin{align}
		\timefunction \circ \, ^{(\textnormal{reg})}F(\, ^{(\textnormal{reg})}\Delta b(t,x^2,x^3),t,x^2,x^3)  & = 0,  \label{E:REQUIREMENTOFFLOWINTERVALFORREGULARNEWULTOMAKELEFTFLOWOUTENDINROUGHHYPZERO} \\
		\timefunction \circ \, ^{(\textnormal{reg})}F( \, ^{(\textnormal{reg})}\Delta a(t,x^2,x^3),t,x^2,x^3)  & = \timefunction_0, \label{E:REQUIREMENTOFFLOWINTERVALTOMAKERIGHTLOWOUTOFREGULARNEWULENDININITIALROUGHHYP}.
	\end{align}
\end{subequations}
That is, the flowouts to the left and right of $\datahypfortimefunctiontwoarg{0}{[\timefunction_0,0]}$ always terminate on $ \hypthreearg{0}{[- \rightu,\leftu]}{0}$ and $ \hypthreearg{\timefunction_0}{[- \rightu,\leftu]}{0}$, respectively.

Finally, if $|\timefunction_0|$ is smaller still, it follows that: 
\begin{align}
	 ^{(\textnormal{reg})}\mathscr{D} \subset \twoargMrough{[\timefunction_0,0],[- \interestingu, \interestingu]}{0} \label{E:FLOWOUTOFREGULARNEWULISCONTAINEDININTERESTINGREGION}
\end{align}\medskip 

\noindent \underline{\textbf{An approximate ingoing eikonal function:}}

Let $\mr \ubar \in {C_{\textnormal{data}}^{2,1} \left(\datahypfortimefunctiontwoarg{0}{[\timefunction_0,0]}\right)}$ be defined as in \eqref{E:RESTRICTIONOFUBARONDATASURFACE}. Let $^{(\textnormal{apx})}\ubar$ be the solution to the IVP: 
\begin{subequations} \label{E:APPROXIMATEEIKONALFUNCTION}
	\begin{align} 
		^{(\textnormal{reg})} \newuL ^{(\textnormal{apx})}\ubar & = 0, \label{E:APPROXIMATEEIKONALFUNCTIONEQUATION} \\
		^{(\textnormal{apx})} \ubar \big|_{\datahypfortimefunctiontwoarg{0}{[\timefunction_0,0]}} & =  \mr \ubar. \label{E:APPROXIMATEEIKONALFUNCTIONINITIALCONDITION}
	\end{align}
\end{subequations}
Then $^{(\textnormal{apx})} \ubar$ is an approximate ingoing eikonal function in the sense that the level sets of $^{(\textnormal{apx})}\ubar$ are $\gfour$-null and transverse to $\nullhyparg{u}$, but: 
\begin{align}
	\upmu \gfour^{-1}(\Dfour \, ^{(\textnormal{apx})}\ubar, \Dfour \, ^{(\textnormal{apx})}\ubar) = \upmu \mathcal{O}(\initialsmall) \not\equiv 0.  \label{E:APPROXEIKONALFUNCTIONHASSMALLGFOURGRADIENTNORM}
\end{align} Moreover, the following estimate holds:
\begin{align}
	\left\| ^{(\textnormal{apx})}\ubar \right\|_{ C_{\textnormal{geo}}^{1,1}\left(^{(\textnormal{reg})} \mathscr{D}\right)} \lesssim 1. \label{E:C11ESTIMATEFORAPPROXIMATEEIKONALFUNCTION}
\end{align}
Finally, the following identities hold:
\begin{subequations}
	\begin{align}
		\Lunit \, ^{(\textnormal{apx})} \ubar = \frac{2}{2 + \upmu \coefficientsInVectorfieldsTangetToUbarDataSurface} \UnitTimeNormalizedUbarDataSurfaceTangentVectorfield \, ^{(\textnormal{apx})}\ubar, \label{E:IDENTITYFORLAPXUBARINTERMSOFDATASURFACEDERIVATIVES} \\
		\muX \,  ^{(\textnormal{apx})} \ubar = - \frac{\upmu}{2 + \upmu \coefficientsInVectorfieldsTangetToUbarDataSurface} \UnitTimeNormalizedUbarDataSurfaceTangentVectorfield \, ^{(\textnormal{apx})}\ubar. \label{E:IDENTITYFORMUXAPXUBARINTERMSOFDATASURFACEDERIVATIVES}
	\end{align}
\end{subequations}
	

\end{lemma}

\begin{proof} \hfill

\noindent \underline{\textbf{Proof of the properties of $^{(\textnormal{reg})}\newuL$:}} The vectorfield $^{(\textnormal{reg})}\newuL$ is $\gfour$-null because $\olduLunit$ is. It is transverse to $\nullhyparg{u}$ because: 
\begin{align}
	^{(\textnormal{reg})}\newuL u = 1. \label{E:REGULARNEWULACTINGONU}
\end{align}
Finally, the $C^{1,1}_{\textnormal{geo}}\left(\twoargMrough{[\timefunction_0,0],[- \rightu,\leftu]}{0}\right)$ regularity of $^{(\textnormal{reg})}\newuL$ follows from the regularity of $\upmu$ proved in Thm\,\ref{T:MAINRESULTSFROMSINGULARBOUNDARYPAPER}. The transversality of $^{(\textnormal{reg})}\newuL$ to $\datahypfortimefunctiontwoarg{0}{[\timefunction_0,0]}$ follows from the same proof of \eqref{E:NEWULISTRANSVERSALTOUBARDATASURFACE} (whose main tool was the smallness of $|\timefunction_0|$) with $^{(\textnormal{reg})}\newuL$ in place of $\newuL$.

\medskip

\noindent \underline{\textbf{Proof of statements of the flow of $^{(\textnormal{reg})}\newuL$:}}

Since $^{(\textnormal{reg})}\newuL$ is transverse to  $\datahypfortimefunctiontwoarg{0}{[\timefunction_0,0]}$, by the $C^{1,1}_{\textnormal{geo}}$-regularity of $^{(\textnormal{reg})}\newuL$,
\eqref{E:DEFINITIONOFFLOWOUTDOMAINFORREGULARNEWUL}--\eqref{E:C11ESTIMATEFORINVERSEFLOWOFREGULARNEWUL}  follow from the Fundamental theorem of flows and flowouts, see e.g. \cite[Theorems 9.12,9.20]{lee2012smooth}.

We now prove \eqref{E:REQUIREMENTOFFLOWINTERVALFORREGULARNEWULTOMAKELEFTFLOWOUTENDINROUGHHYPZERO}--\eqref{E:REQUIREMENTOFFLOWINTERVALTOMAKERIGHTLOWOUTOFREGULARNEWULENDININITIALROUGHHYP}. Consider any point $ p_{\mulevelsetvalue} = (\Cartesiantisafunctiononmumxtoriarg{\mulevelsetvalue}{0}(x^2,x^3), \Eikonalisafunctiononmumuxtoriarg{\mulevelsetvalue}{0}(x^2,x^3),x^2,x^3) \in \datahypfortimefunctiontwoarg{0}{[\timefunction_0,0]}$ where $\mulevelsetvalue \in [0,\mupositive]$ and a flowout to the left, that is, consider $\Delta s > 0$. Then the $t$-value of the flowout from $p_{\mulevelsetvalue}$ is given by:
\begin{align}
	t \circ \,^{(\textnormal{reg})} \widehat{\underline{\Lambda}}_{\Delta s} (p_{\mulevelsetvalue}) = \Cartesiantisafunctiononmumxtoriarg{\mulevelsetvalue}{0}(x^2,x^3) + \int_{\Delta s' = 0}^{\Delta s} (\,^{(\textnormal{reg})} \hat{\uLunit} t)\circ \,^{(\textnormal{reg})} \widehat{\underline{\Lambda}}_{\Delta s'} (p_{\mulevelsetvalue}) \, \mathrm{d} \Delta s'. \label{E:CARTESIANTVALUEOFREGULARNEWULFLOWOUTS}
\end{align}
To estimate the RHS of \eqref{E:CARTESIANTVALUEOFREGULARNEWULFLOWOUTS}, we compute the following: 
\begin{align}
	^{(\textnormal{reg})}\newuL t \big|_{\datahypfortimefunctiontwoarg{0}{[\timefunction_0,0]}} & =  \frac{1}{2}\upmu\big|_{\datahypfortimefunctiontwoarg{0}{[\timefunction_0,0]}}, &  ^{(\textnormal{reg})}\newuL ^{(\textnormal{reg})}\newuL t \big|_{\datahypfortimefunctiontwoarg{0}{[\timefunction_0,0]}}  & = \frac{1}{2} \, ^{(\textnormal{reg})}\newuL \upmu \big|_{\datahypfortimefunctiontwoarg{0}{[\timefunction_0,0]}}  , &  ^{(\textnormal{reg})}\newuL \, ^{(\textnormal{reg})} \newuL \, ^{(\textnormal{reg})}\newuL t \big|_{\datahypfortimefunctiontwoarg{0}{[\timefunction_0,0]}}  & = \frac{1}{2} \muX \muX \upmu \big|_{\datahypfortimefunctiontwoarg{0}{[\timefunction_0,0]}}  + \mathcal{O}(|\timefunction_0|). \label{E:ESTIMATESFORCARTESIANTFORREGULARNEWUL}
\end{align}
From \eqref{E:ESTIMATESFORCARTESIANTFORREGULARNEWUL} and Taylor's theorem, there are constants $C_1,C_2> 0$ such that he following bound holds for any $\Delta s' \in [0,\Delta s]$:
\begin{align} 
	(\,^{(\textnormal{reg})}  \hat{\uLunit} t)\circ \,^{(\textnormal{reg})}  \widehat{\underline{\Lambda}}_{\Delta s'} (p_{\mulevelsetvalue}) & 
	\ge C_1 |\timefunction_0| - C_2|\timefunction_0|  \Delta s'+ \frac{\secondtransversalderivativemulowerbound}{6} (\Delta s')^2. \label{E:ESTIMATEFORREGULARNEWULDERIVATIVEOFCARTESIANTONFLOWOUTFORITERATES}
\end{align}
For $|\timefunction_0|$ sufficiently small, the quadratic polynomial in \eqref{E:ESTIMATEFORREGULARNEWULDERIVATIVEOFCARTESIANTONFLOWOUTFORITERATES} has no real roots and is uniformly bounded from below. That is, we have the crude bound: 
	\begin{align}
		t \circ \,^{(\textnormal{reg})}\widehat{\underline{\Lambda}}_{\Delta s} (p_{\mulevelsetvalue}) \ge \Cartesiantisafunctiononmumxtoriarg{\mulevelsetvalue}{0}(x^2,x^3)  + C \Delta s, \label{E:CRUDEBOUNDFORTFORCARTESIANTINLEFTFLOWOUTFORREGULARNEWUL}
	\end{align}
		from which \eqref{E:REQUIREMENTOFFLOWINTERVALFORREGULARNEWULTOMAKELEFTFLOWOUTENDINROUGHHYPZERO} follows from the maximality of the flowout $^{(\textnormal{reg})}\mathscr{D}$   since  $\hypthreearg{0}{[- \rightu,\leftu]}{0}$ lies on the boundary of $\twoargMrough{[\timefunction_0,0],[- \rightu,\leftu]}{0}$. A symmetric argument for $\Delta s \le 0$ proves \eqref{E:REQUIREMENTOFFLOWINTERVALTOMAKERIGHTLOWOUTOFREGULARNEWULENDININITIALROUGHHYP}. 
		
We now prove \eqref{E:FLOWOUTOFREGULARNEWULISCONTAINEDININTERESTINGREGION}. By \eqref{E:MUXMUKAPPALEVELSETLOCATION}, it suffices to show that: 
\begin{align}
	\sup_{\domainforembeddingdatahypfortimefunctiontwoarg{0}{[0,\mupositive]}} \, ^{(\textnormal{reg})} \Delta b \le \frac{1}{2} \interestingu & & \inf_{\domainforembeddingdatahypfortimefunctiontwoarg{0}{[0,\mupositive]}} \, ^{(\textnormal{reg})} \Delta a \ge - \frac{1}{2}\interestingu. \label{E:BOUNDSFORREGULARNEWULFLOWOUTDOMAININTERVALDEFININGFUNCTIONS}
\end{align}
By \eqref{E:REQUIREMENTOFFLOWINTERVALFORREGULARNEWULTOMAKELEFTFLOWOUTENDINROUGHHYPZERO}--\eqref{E:REQUIREMENTOFFLOWINTERVALTOMAKERIGHTLOWOUTOFREGULARNEWULENDININITIALROUGHHYP} as well as estimate \eqref{E:CRUDEBOUNDFORTFORCARTESIANTINLEFTFLOWOUTFORREGULARNEWUL} and the analogous one for $\Delta s < 0$, this will  follow once we show that the difference of Cartesian $t$-values between $\hypthreearg{0}{[- \rightu,\leftu]}{0}, \,\twoargmumuxtorus{\mupositive}{0}$ and between $\twoargmumuxtorus{0}{0}, \, \hypthreearg{\timefunction_0}{[- \rightu,\leftu]}{0}$ is of size $\mathcal{O}(|\timefunction_0|)$. First we note that $\muX \timefunction = 0$ implies:
\begin{align}
	\min_{\hypthreearg{0}{[- \rightu,\leftu]}{0}} t & = \min_{\twoargmumuxtorus{0}{0}} t,  & \max_{\hypthreearg{0}{[- \rightu,\leftu]}{0}} t & = \max_{\twoargmumuxtorus{0}{0}} t, \label{E:TIMEBOUNDSALONGFINALROUGHHYPISTHESAMEASTIMEBOUNDSALONGCREASE} \\
	\min_{\hypthreearg{\timefunction_0}{[- \rightu,\leftu]}{0}} t & = \min_{\twoargmumuxtorus{\mupositive}{0}} t, & \max_{\hypthreearg{\timefunction_0}{[- \rightu,\leftu]}{0}} t & = \max_{\twoargmumuxtorus{\mupositive}{0}} t. \label{E:TIMEBOUNDSALONGINITIALROUGHHYPISTHESAMEASTIMEBOUNDSALONGPRIMALTORUS}
\end{align}
The claim follows from \eqref{E:CLOSEDIMPROVEMENTMUEQUALSMINUSKAPPAEMBEDDINGCARTESIANTIMEFUNCTIONNEGATIVEMUDERIVATIVE}.
	
\medskip 

\noindent \underline{\textbf{Proof of the statements of the approximate ingoing eikonal function:}}

The evolution equation \eqref{E:APPROXIMATEEIKONALFUNCTIONEQUATION} for $^{(\textnormal{apx})}\ubar$ can be written equivalently as $\frac{\p}{\p \Delta s} \big(\, ^{(\textnormal{apx})}\ubar \circ \, ^{(\textnormal{reg})} F(\Delta s,t,x^2,x^3)\big) = 0$. The initial condition \eqref{E:APPROXIMATEEIKONALFUNCTIONINITIALCONDITION} can similarly be expressed as $^{(\textnormal{reg})} F(0,t,x^2,x^3) = \mr \ubar \circ \embeddingdatahypfortimefunctionarg{0}(t,x^2,x^3)$. Hence, the approximate eikonal function is given by $^{(\textnormal{apx})}\ubar \circ \, ^{(\textnormal{reg})} F(\Delta s,t,x^2,x^3) =  \mr \ubar \circ \embeddingdatahypfortimefunctionarg{0}(t,x^2,x^3)$. Composing with $^{(\textnormal{reg})} F^{-1}$ and using \eqref{E:C11ESTIMATEFORINVERSEFLOWOFREGULARNEWUL} proves \eqref{E:C11ESTIMATEFORAPPROXIMATEEIKONALFUNCTION}. Since $^{(\textnormal{reg})}\newuL$ is tangent to level sets of $^{(\textnormal{apx})}\ubar$, the statement about its level sets follow from the proved properties of $^{(\textnormal{reg})}\newuL$. Using \eqref{E:SMOOTHTORUSMETRICINTERMSOFSIGMATMETRICANDX} and \eqref{E:REGULARINGOINGNULLVECTORFIELD}, we express the product of $\upmu$ and the inverse spacetime metric as $\upmu \gfour^{-1} = - \Lunit \otimes\, ^{(\textnormal{reg})}\newuL - \, ^{(\textnormal{reg})}\newuL \otimes \Lunit + \upmu \gtorus^{-1}$. From this and \eqref{E:APPROXIMATEEIKONALFUNCTIONEQUATION}, the desired identity \eqref{E:APPROXEIKONALFUNCTIONHASSMALLGFOURGRADIENTNORM} holds. Finally, \eqref{E:IDENTITYFORLAPXUBARINTERMSOFDATASURFACEDERIVATIVES}--\eqref{E:IDENTITYFORMUXAPXUBARINTERMSOFDATASURFACEDERIVATIVES} hold from straightforward calculations using the definitions  \eqref{E:REGULARINGOINGNULLVECTORFIELD},\eqref{E:FUTUREORIENTEDVECTORFIELDTANGENTTOUBARDATASURFACE} as well as \eqref{E:APPROXIMATEEIKONALFUNCTIONEQUATION}.

\end{proof}

\begin{center}
	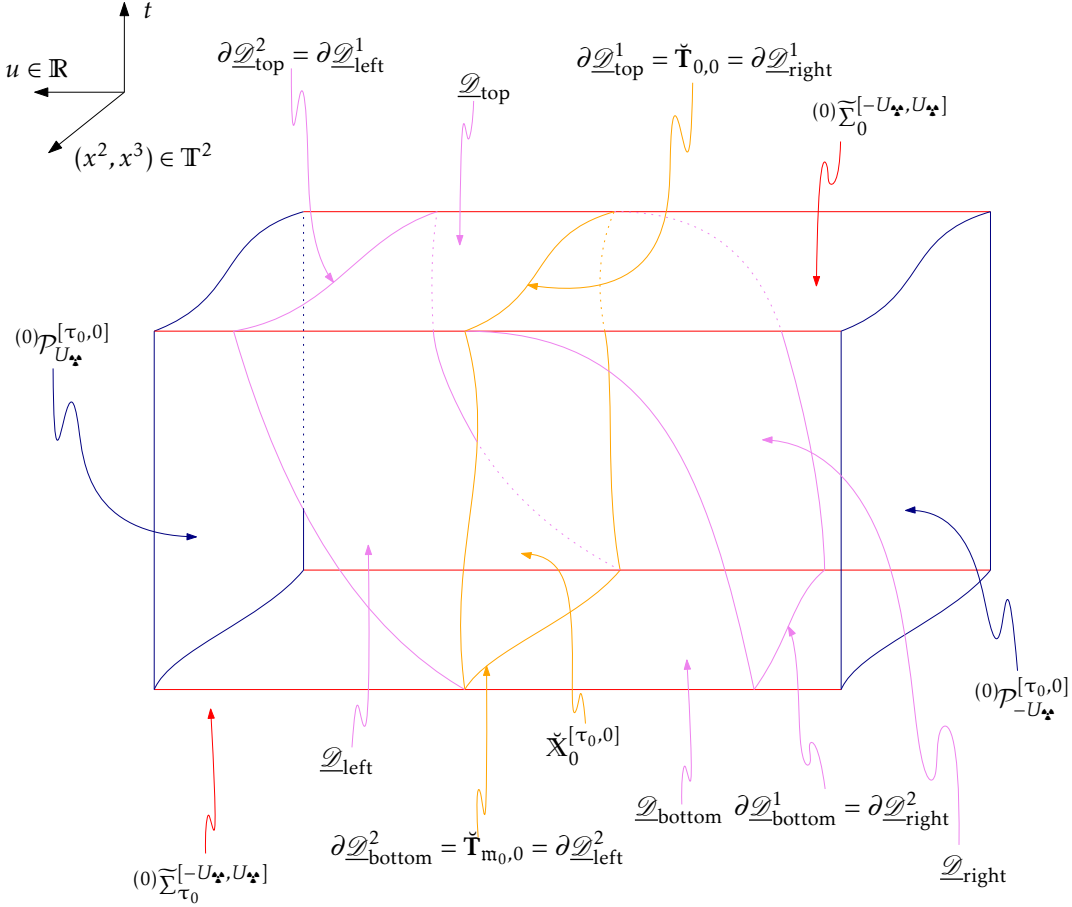
\begin{figure}  
		\begin{overpic}[scale=.7, grid = false, tics=5]{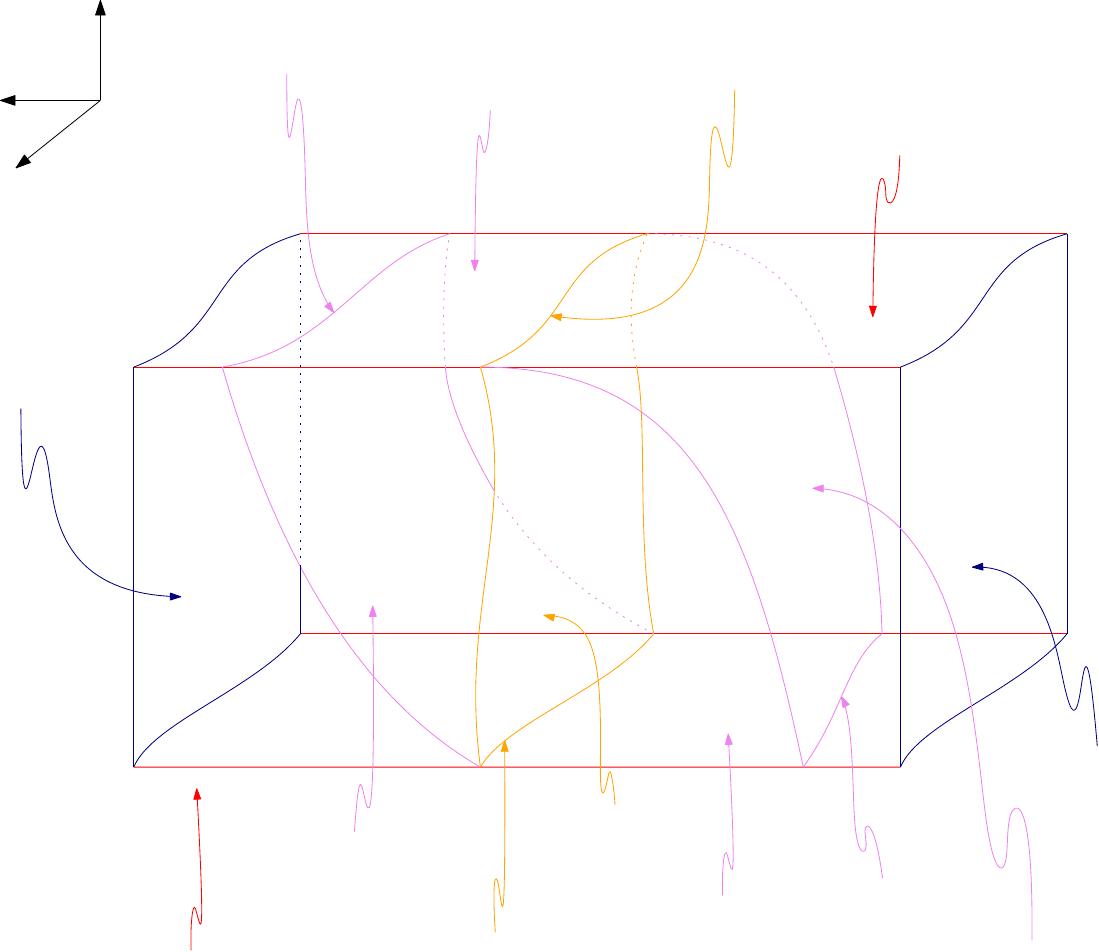}
			\put (43,77.5) {$\underline{\mathscr{D}}_{\textnormal{top}}$}
			\put (55,80) {$\p \underline{\mathscr{D}}_{\textnormal{top}}^1 = \twoargmumuxtorus{0}{0} = \p \underline{\mathscr{D}}_{\textnormal{right}}^1$}
			\put (18.5,80.5) {$\p \underline{\mathscr{D}}_{\textnormal{top}}^2 = \p \underline{\mathscr{D}}_{\textnormal{left}}^1$}
			\put(30,0) {$\p\underline{\mathscr{D}}_{\textnormal{bottom}}^2 =\twoargmumuxtorus{\mupositive}{0} = \p \underline{\mathscr{D}}_{\textnormal{left}}^2$}
			\put (29,9) {$\underline{\mathscr{D}}_{\textnormal{left}}$}
			\put (92,-2) {$\underline{\mathscr{D}}_{\textnormal{right}}$}
			\put (71,4) {$\p \underline{\mathscr{D}}_{\textnormal{bottom}}^1 = \p \underline{\mathscr{D}}_{\textnormal{right}}^2 $}
			\put (79,74) {$\hypthreearg{0}{[-\interestingu,\interestingu]}{0}$}
			\put (95.5,15) {$\nullhypthreearg{0}{-\interestingu}{[\timefunction_0,0]}$}
			\put (-2,51) {$\nullhypthreearg{0}{\interestingu}{[\timefunction_0,0]} $}
			\put (10,-4) {$\hypthreearg{\timefunction_0}{[-\interestingu,\interestingu]}{0}$}
			\put (61,4) {$\underline{\mathscr{D}}_{\textnormal{bottom}}$}
			\put (4,70) {$(x^2,x^3) \in \mathbb{T}^2$}
			\put (11,85) {$t$}
			\put (-3,79) {$u \in \mathbb{R}$}
			\put (52,10) {$\datahypfortimefunctiontwoarg{0}{[\timefunction_0,0]}$}
		\end{overpic}
		\vspace{0.5cm}
		\caption{The domain $\underline{\mathscr{D}}$ on which $\ubar$ is constructed in the differential structure of $(t,u,x^2,x^3)$.}
		\label{F:DOMAINOFCONSTRUCTIONFORINGOINGEIKONALFUNCTIONINGEOMETRICCOORDINATES}
	\end{figure}
\end{center}




\begin{theorem}[Construction and estimates of the ingoing eikonal function]\label{T:CONSTRUCTIONOFTHEINGOINGEIKONALFUNCTIONCONTRACTIONMAPPING} Assume the assumptions and conclusions of Theorem\,\ref{T:MAINRESULTSFROMSINGULARBOUNDARYPAPER} and Prop.\,\ref{P:PROPERTIESOFMUXMUZEROLEVELSET} hold. Let $\datahypfortimefunctiontwoarg{0}{[\timefunction_0,0]}$ be the truncated zero-level set $\muX \upmu$ as in \eqref{E:TRUNCATEDLEVELSETSOFMUXMU} and let $\domainforembeddingdatahypfortimefunctiontwoarg{0}{[0,\mupositive]}$ be the domain for the embedding $\embeddingdatahypfortimefunctionarg{0}$ as in \eqref{E:DOMAINOFEMBEDDINGFORXMUEQUALSMINUSKAPPAYSURFACE}--\eqref{E:XMUISMINUSCAPPISAGRAPH}. Let $\Cartesiantisafunctiononmumxtoriarg{\mulevelsetvalue}{0}, \, \Eikonalisafunctiononmumuxtoriarg{\mulevelsetvalue}{0}\colon \T^2 \to \R$ be the $(t,u)$-values of $\datahypfortimefunctiontwoarg{0}{[\timefunction_0,0]}$ for $\mulevelsetvalue \in [0, \mulevelsetvalue_0]$, see \eqref{E:CARTESIANTANDEIKONALFUNCTIONAREGRAPHSONMUMUXTORI} and  \eqref{E:CLOSEDIMPROVEMENTLEVELSETSTRUCTUREOFMUXEQUALSMINUSKAPPA}. Let $\flowmapofnewularg{\Delta s}$ be the flow map of $\newuL$, that is, the $(t,u,x^2,x^3)$-coordinates of $\flowmapofnewularg{\Delta s}$ solve:
	\begin{align} \label{E:FLOWMAPOFNEWUL}
		\frac{\p}{\p \Delta s} \flowmapofnewularg{\Delta s}(t,u,x^2,x^3) = \newuL \circ \flowmapofnewularg{\Delta s}(t,u,x^2,x^3), & & \flowmapofnewularg{0}(t,u,x^2,x^3) = (t,u,x^2,x^3).
	\end{align} 
Let $\interestingu$ denote the same constant as in \cite{abbrescia2022emergence}*{Equation\,(A.44)}, which satisfies the following smallness bound: 
\begin{equation}
	0 < \interestingu < \min\{1,\leftu,\rightu,\frac{1}{22}\blowupdelta^{-1}\}.  \label{E:BOUNDFORINTERESTINGU}
\end{equation}
Then if $|\timefunction_0|$ and $\initialsmall$ are sufficiently small, the following hold:

There exists a non-negative function $\Delta \underline b(t,x^2,x^3)$ and a non-positive function $\Delta \underline a(t,x^2,x^3)$ such that for any $(t,x^2,x^3) \in \domainforembeddingdatahypfortimefunctiontwoarg{0}{[0,\mupositive]}$, 
\begin{subequations} 
	\begin{align}
		\timefunction \left( \flowmapofnewularg{\Delta \underline{b}(t,x^2,x^3)}\circ\embeddingdatahypfortimefunctionarg{0}(t, x^2,x^3) \right) & = 0, \label{E:REQUIREMENTOFFLOWINTERVALFORNEWULTOMAKELEFTFLOWOUTENDINROUGHHYPZERO} \\
		\timefunction \left( \flowmapofnewularg{\Delta \underline{a}(t,x^2,x^3)}\circ\embeddingdatahypfortimefunctionarg{0}(t, x^2,x^3) \right) & = \timefunction_0 ,\label{E:REQUIREMENTOFFLOWINTERVALFORNEWULTOMAKERIGHTFLOWOUTENDININITIALROUGHHYP}		
	\end{align}
\end{subequations}
such that a unique solution $\ubar$ to \eqref{E:INGOINGEIKONALEQUATIONINITIALVALUEPROBLEM} exists on the flowout: 
\begin{align} \label{E:FLOWOUTOFDATASURFACE}
	\underline{\mathscr{D}} \eqdef \underline{F} \left(\underline{\mathscr{F}}\right),
\end{align} where $\underline F$ and $\underline{\mathscr{F}}$ are the map and set defined by:
\begin{subequations}
	\begin{align}
		\underline{F} (\Delta s,t,x^2,x^3) & \eqdef \flowmapofnewularg{\Delta s} \circ \embeddingdatahypfortimefunctionarg{0}(t, x^2,x^3), \label{E:DEFINITIONOFFLOWOUTMAP} \\
		\begin{split} \label{E:DEFINITIONOFFLOWOUTDOMAINFORUBAR}
			\underline{\mathscr{F}} & \eqdef \left\{ (\Delta s, t,x^2,x^3) \ | \ (t,x^2,x^3) \in \domainforembeddingdatahypfortimefunctiontwoarg{0}{[0,\mupositive]} 
			 \text{ and } \Delta s \in [\Delta \underline{a}(t,x^2,x^3) 
			 ,\Delta \underline{b}(t,x^2,x^3)
			 ] \right\}.
		\end{split}
	\end{align}
\end{subequations}
Moreover, the flowout domain of $\underline{\mathscr{D}}$ satisfies the following:

	\begin{itemize}
		\item The location of $\underline{\mathscr{D}}$ is in the interesting region:
			\begin{align}
				\underline{\mathscr{D}} \subset   \twoargMrough{[\timefunction_0,0],[- \interestingu, \interestingu]}{0} \label{E:LOCATIONOFUBARFLOWOUTDOMAIN}
			\end{align}
		\item The ``top" lateral boundary of $\underline{\mathscr{D}}$ is the hypersurface:
		\begin{align}
			\underline{\mathscr{D}}_{\textnormal{top}}  & \eqdef  \left\{ \flowmapofnewularg{\Delta \underline b(t,x^2,x^3)}\circ\embeddingdatahypfortimefunctionarg{0}(t, x^2,x^3) \ \bigg| \ (t,x^2,x^3) \in \domainforembeddingdatahypfortimefunctiontwoarg{0}{[0,\mupositive]} \right\}, \label{E:TOPLATERALBOUNDARYOFSPACETIMEREGIONFORINGOINGEIKONAL} \\
		\end{align}
		whose boundary components are: 
		\begin{subequations}
			\begin{align}
				\p \underline{\mathscr{D}}_{\textnormal{top}}^1 & \eqdef  \bigg\{ \flowmapofnewularg{\Delta \underline b(\Cartesiantisafunctiononmumxtoriarg{0}{0}(x^2,x^3),x^2,x^3)}\circ\embeddingdatahypfortimefunctionarg{0}(\Cartesiantisafunctiononmumxtoriarg{0}{0}(x^2,x^3),x^2,x^3)  \ \bigg| \ (x^2,x^3) \in \T^2 \bigg\}, \label{E:FIRSTBOUNDARYPIECEOFTOPBOUNDARYSURFACE} \\ 
				\p \underline{\mathscr{D}} _{\textnormal{top}}^2 & \eqdef   \bigg\{ \flowmapofnewularg{\Delta \underline b(\Cartesiantisafunctiononmumxtoriarg{\mupositive}{0}(x^2,x^3),x^2,x^3)}\circ\embeddingdatahypfortimefunctionarg{0}(\Cartesiantisafunctiononmumxtoriarg{\mupositive}{0}(x^2,x^3),x^2,x^3)  \ \bigg| \ (x^2,x^3) \in \T^2 \bigg\}.
			\label{E:SECONDBOUNDARYPIECEOFTOPBOUNDARYSURFACE}
			\end{align}
		\end{subequations} 
		 In particular, \eqref{E:REQUIREMENTOFFLOWINTERVALFORNEWULTOMAKELEFTFLOWOUTENDINROUGHHYPZERO} implies that $\p\underline{\mathscr{D}}_{\textnormal{top}}^1 = \twoargmumuxtorus{0}{0}$ and $\underline{\mathscr{D}}_{\textnormal{top}} \subset \hypthreearg{0}{[- \rightu,\leftu]}{0}$.
		\item The ``bottom" lateral boundary of $\mathscr{D}$ is the hypersurface: 		
		\begin{equation} 
			\underline{\mathscr{D}}_{\textnormal{bottom}} \eqdef   \left\{ \flowmapofnewularg{\Delta \underline a(t,x^2,x^3)}\circ\embeddingdatahypfortimefunctionarg{0}(t, x^2,x^3) \ \bigg| \ (t,x^2,x^3) \in \domainforembeddingdatahypfortimefunctiontwoarg{0}{[0,\mupositive]} \right\},\label{E:BOTTOMLATERALBOUNDARYOFSPACETIMEREGIONFORINGOINGEIKONAL}
		\end{equation}
			whose boundary components are: 
			\begin{subequations}
				\begin{align}
					\p\underline{\mathscr{D}}_{\textnormal{bottom}}^1 & \eqdef   \bigg\{ \flowmapofnewularg{\Delta \underline a(\Cartesiantisafunctiononmumxtoriarg{0}{0}(x^2,x^3),x^2,x^3)}\circ\embeddingdatahypfortimefunctionarg{0}(\Cartesiantisafunctiononmumxtoriarg{0}{0}(x^2,x^3),x^2,x^3)  \ \bigg| \ (x^2,x^3) \in \T^2 \bigg\}, \label{E:FIRSTBOUNDARYPIECEOFBOTTOMBOUNDARYSURFACE} \\
					\p \underline{\mathscr{D}}_{\textnormal{bottom}}^2 & \eqdef   \bigg\{ \flowmapofnewularg{\Delta \underline a(\Cartesiantisafunctiononmumxtoriarg{\mupositive}{0}(x^2,x^3),x^2,x^3)}\circ\embeddingdatahypfortimefunctionarg{0}(\Cartesiantisafunctiononmumxtoriarg{\mupositive}{0}(x^2,x^3),x^2,x^3)  \ \bigg| \ (x^2,x^3) \in \T^2 \bigg\}.\label{E:SECONDBOUNDARYPIECEOFBOTTOMBOUNDARYSURFACE}
				\end{align}
			\end{subequations} 			
			In particular, \eqref{E:REQUIREMENTOFFLOWINTERVALFORNEWULTOMAKERIGHTFLOWOUTENDININITIALROUGHHYP} implies that $\p \mathscr{D}_{\textnormal{bottom}}^2 =  \twoargmumuxtorus{\mupositive}{0}$ and $\mathscr{D}_{\textnormal{bottom}} \subset \hypthreearg{\timefunction_0}{[- \rightu,\leftu]}{0}$.
			
		\item The ``right" boundary is the hypersurface: 
		\begin{align} 
			\begin{split} \label{E:RIGHTBOUNDARYOFSPACETIMEREGIONFORINGOINGEIKONAL} 
				\underline{\mathscr{D}}_{\textnormal{right}} & \eqdef  \left\{ \flowmapofnewularg{\Delta s}\circ\embeddingdatahypfortimefunctionarg{0}(\Cartesiantisafunctiononmumxtoriarg{0}{0}(x^2,x^3),x^2,x^3) \ \bigg| \ (x^2,x^3) \in \T^2, 
				\text{ and } [\Delta \underline{a}(\Cartesiantisafunctiononmumxtoriarg{0}{0}(x^2,x^3),x^2,x^3)
				,\Delta \underline{b}(\Cartesiantisafunctiononmumxtoriarg{0}{0}(x^2,x^3),x^2,x^3)
				]\right\},
			\end{split}
		\end{align}
			whose boundary components are: 
			\begin{subequations}
				\begin{align}
					\p \underline{\mathscr{D}}_{\textnormal{right}}^1 & \eqdef  \bigg\{ \flowmapofnewularg{\Delta \underline b(\Cartesiantisafunctiononmumxtoriarg{0}{0}(x^2,x^3),x^2,x^3)}\circ\embeddingdatahypfortimefunctionarg{0}(\Cartesiantisafunctiononmumxtoriarg{0}{0}(x^2,x^3),x^2,x^3)  \ \bigg| \ (x^2,x^3) \in \T^2 \bigg\}, 
 \label{E:FIRSTBOUNDARYPIECEOFRIGHTPBOUNDARYSURFACE} \\
					\p \underline{\mathscr{D}}_{\textnormal{right}}^2 & \eqdef   \bigg\{ \flowmapofnewularg{\Delta \underline a(\Cartesiantisafunctiononmumxtoriarg{0}{0}(x^2,x^3),x^2,x^3)}\circ\embeddingdatahypfortimefunctionarg{0}(\Cartesiantisafunctiononmumxtoriarg{0}{0}(x^2,x^3),x^2,x^3)  \ \bigg| \ (x^2,x^3) \in \T^2 \bigg\}. \label{E:SECONDBOUNDARYPIECEOFRIGHTBOUNDARYSURFACE} 
				\end{align}
			\end{subequations}
			In particular, the first boundary component satisfies $\p \underline{\mathscr{D}}_{\textnormal{right}}^1  = \twoargmumuxtorus{0}{0}$. Moreover, $\ubar \equiv 0$ on $\mathscr{D}_{\textnormal{right}}$. 
		\item The ``left" boundary is the hypersurface: 
		\begin{align}
			\begin{split}\label{E:LEFTBOUNDARYOFSPACETIMEREGIONFORINGOINGEIKONAL}
				\underline{\mathscr{D}}_{\textnormal{left}} & \eqdef  \left\{ \flowmapofnewularg{\Delta s}\circ\embeddingdatahypfortimefunctionarg{0}(\Cartesiantisafunctiononmumxtoriarg{\mupositive}{0}(x^2,x^3),x^2,x^3) \ \bigg| \ (x^2,x^3) \in \T^2,
				\text{ and } [\Delta \underline{a}(\Cartesiantisafunctiononmumxtoriarg{\mupositive}{0}(x^2,x^3),x^2,x^3)
				,\Delta \underline{b}(\Cartesiantisafunctiononmumxtoriarg{\mupositive}{0}(x^2,x^3),x^2,x^3)
				]\right\},
			\end{split}
		\end{align}
		whose boundary components are: 
		\begin{subequations}
			\begin{align}
				\p \underline{\mathscr{D}}_{\textnormal{left}}^1 & \eqdef  \bigg\{ \flowmapofnewularg{\Delta \underline b(\Cartesiantisafunctiononmumxtoriarg{\mupositive}{0}(x^2,x^3),x^2,x^3)}\circ\embeddingdatahypfortimefunctionarg{0}(\Cartesiantisafunctiononmumxtoriarg{\mupositive}{0}(x^2,x^3),x^2,x^3)  \ \bigg| \ (x^2,x^3) \in \T^2 \bigg\}, \label{E:FIRSTBOUNDARYPIECEOLEFTBOUNDARYSURFACE} \\
				\p \underline{\mathscr{D}}_{\textnormal{left}}^2 & \eqdef    \bigg\{ \flowmapofnewularg{\Delta \underline a(\Cartesiantisafunctiononmumxtoriarg{\mupositive}{0}(x^2,x^3),x^2,x^3)}\circ\embeddingdatahypfortimefunctionarg{0}(\Cartesiantisafunctiononmumxtoriarg{\mupositive}{0}(x^2,x^3),x^2,x^3)  \ \bigg| \ (x^2,x^3) \in \T^2 \bigg\}. \label{E:SECONDBOUNDARYPIECEOLEFTBOUNDARYSURFACE}
			\end{align}
		\end{subequations}
		In particular, the second boundary piece satisfies $\p \mathscr{D}_{\textnormal{left}}^2  = \twoargmumuxtorus{\mupositive}{0}$ 
		Moreover, $\ubar \equiv \timefunction_0 = - \mupositive$ on $\underline{\mathscr{D}}_{\textnormal{left}}$. 
		\item The boundary components $\underline{\mathscr{D}}_{\textnormal{left}},  \underline{\mathscr{D}}_{\textnormal{right}}$ are co-dimension one $C^{1,1}_{geo}$-submanifolds of $\twoargMrough{[\timefunction_0,0],[- \rightu, \leftu]}{0}$. By \eqref{E:C21ESTIMATESOFTHEROUGHTIMEFUNCTIONINGEOMETRICCOORDINATES}, the boundary components $\underline{\mathscr{D}}_{\textnormal{top}}, \underline{\mathscr{D}}_{\textnormal{bottom}}$ are co-dimension one $C^{2,1}_{geo}$-submanifolds of $\twoargMrough{[\timefunction_0,0],[- \rightu, \leftu]}{0}$.
	\end{itemize}

In addition, the following estimates hold for $A = 2,3$:
	\begin{align}
		\left\|  \Lunit \ubar \right\|_{C_{\textnormal{geo}}^{0,1}(\underline{\mathscr{D}})} & \le C & \left\|  \muX \ubar \right\|_{C_{\textnormal{geo}}^{0,1}(\underline{\mathscr{D}})} & \le C, & \left\|   \geop{x^A}  \ubar \right\|_{C_{\textnormal{geo}}^{0,1}(\underline{\mathscr{D}})} & \le C \initialsmall. \label{E:ANISOTROPICPOINTWISEESTIMATESFORFIRSTDERIVATIVESOFINGOINGEIKONALONFLOWOUTDOMAIN}  \\
		\left\| \newuL \Lunit \ubar \right\|_{C_{\textnormal{geo}}^{0,1}(\underline{\mathscr{D}})} & \le C & \left\| \newuL \muX \ubar \right\|_{C_{\textnormal{geo}}^{0,1}(\underline{\mathscr{D}})} & \le C, & \left\| \newuL  \geop{x^A}  \ubar \right\|_{C_{\textnormal{geo}}^{0,1}(\underline{\mathscr{D}})} & \le C \initialsmall, \label{E:ANISOTROPICPOINTWISEESTIMATESFORHIGHERNEWULDERIVATIVESOFINGOINGEIKONALONFLOWOUTDOMAIN}
	\end{align}
	\begin{align}
		\frac{7}{9} \blowupdelta  \le \min_{\underline{\mathscr{D}}} \Lunit \ubar & \le \max_{\underline{\mathscr{D}}} \Lunit \ubar \le \frac{9}{7}\blowupdelta. \label{E:SHARPUNIFORMESTIMATEFORLUBARINFLOWOUTREGION}
	\end{align}
\end{theorem}
\begin{proof}
We will prove that \eqref{E:QUASILINEAREQUATIONSATISFIEDBYLUBAR}--\eqref{E:QUASILINEAREQUATIONSATISFIEDBYSMOOTHANGULARDERIVATIVEOFUBAR} have unique solutions $\Lunit \ubar, \, \muX \ubar, \, \geop{x^2} \ubar, \, \geop{x^3} \ubar$ through a contraction mapping argument. 

\medskip

\noindent \underline{\textbf{The admissibility class and contraction space:}}

Given a $4$-tuple of spacetime functions $(\hat{\mathcal{I}},\, \hat{J}, \, \hat{z}_2,\, \hat{z}_3)$, consider the following:
	\begin{subequations}
\begin{align}
			\hat \ReciprocalLunitAppliedtoTimeFunction & \eqdef \frac{1}{\hat{\mathcal{I}}}, & \hat \MagnitueofinnerproductofnewLandnewuL & \eqdef \upmu \hat \ReciprocalLunitAppliedtoTimeFunction,  \label{E:ITERATIONDEFINITIONOFLUBARANDINNERPRODUCTOFNEWLANDNEWUL} \\
			\hat \ReciprocaluLunitAppliedtoTimeFunction & \eqdef \frac{\upmu \hat \ReciprocalLunitAppliedtoTimeFunction}{ \hat \ingoingmu}, & \frac{1}{\hat{\ingoingmu} }& \eqdef \frac{1}{2} \hat{\mathcal{I}} + \frac{1}{2 \mathcal{\hat{I}}} (\gtorus^{-1})^{AB} \hat z_A \hat z_B, \label{E:ITERATIONDEFINITIONOFLBARUANDINGOINGMU}
		\end{align}
		\begin{equation} \label{E:ITERATIONDEFINITIONOFNEWUL}
			\hat{\uLunit} \eqdef \hat \ReciprocaluLunitAppliedtoTimeFunction \Lunit + \muX - \hat \MagnitueofinnerproductofnewLandnewuL \gtorus^{AB} \hat z_B \geop{x^A}.
		\end{equation}
	\end{subequations}
	We note that  \eqref{E:ITERATIONDEFINITIONOFLUBARANDINNERPRODUCTOFNEWLANDNEWUL}--\eqref{E:ITERATIONDEFINITIONOFLBARUANDINGOINGMU} are the respective definitions of the scalars introduced in Sect.\,\ref{SSS:NULLVECTORFIELDS} if $\hat{\mathcal{I}} = \Lunit \ubar$, where we used \eqref{E:INGOINGMUDETERMINEDBYEVERYTHINGELSE} to motivate the definition of $\hat{\ingoingmu}$. Likewise, $\hat{J}$ corresponds to $\muX \ubar$.  Also, \eqref{E:ITERATIONDEFINITIONOFNEWUL} is formally the vectorfield $\newuL$ defined by the scalar functions $\hat{\mathcal{I}},\, \hat{J}, \, \hat{z}_2,\, \hat{z}_3$, see \eqref{E:CONVENIENTIDENTITYFORULUNIT}. Denote further $\widehat{\underline{\Lambda}}_{\Delta s}$ as the flow map of $\hat{\uLunit}$, that is, the $(t,u,x^2,x^3)$-coordinate components of $\widehat{\underline{\Lambda}}_{\Delta s}$ solve: 
	\begin{align} \label{E:FLOWMAPOFNEWULITERATE}
		\frac{\p}{\p \Delta s} \widehat{\underline{\Lambda}}_{\Delta s}(t,u,x^2,x^3) = \hat{\uLunit} \circ \widehat{\underline{\Lambda}}_{\Delta s}(t,u,x^2,x^3), & & \widehat{\underline{\Lambda}}_{0}(t,u,x^2,x^3) = (t,u,x^2,x^3).
	\end{align} 
	
	Next we define an admissibility class $\mathcal{S}$ consisting of those $4$-tuples of functions $\hat{\mathbf{v}} \eqdef (\hat{\mathcal{I}},\hat{J},\hat{z}_2,\hat{z}_3)$ such that:
	\begin{align}
		\bullet & \, \, \text{There is a neighborhood $\mathscr{K}$ of $\datahypfortimefunctiontwoarg{0}{[\timefunction_0,0]}$ such that $\hat{\mathbf{v}}  \in C^{0,1}_{\textnormal{geo}}(\hat{\mathscr{K}})^{\times 4}$;} \tag{$\mathcal{S}_1$} \label{E:REQUIREDREGULARITYFORADMISSIBILITYCLASS} \\
		\bullet & \, \, \text{The vectorfield $\hat{\uLunit}$ is transverse to the data hypersurface $\datahypfortimefunctiontwoarg{0}{[\timefunction_0,0]}$.}  \tag{$\mathcal{S}_2$} \label{E:TRANSVERSALITYTODATAHYPERSURFACEOFADMISSIBILITYCLASS}
	\end{align}
The set $\mathcal{S}$ is non-empty since $^{(\textnormal{apx})} \mathbf{v} \eqdef \left(\Lunit \,^{(\textnormal{apx})}\ubar, \muX \, ^{(\textnormal{apx})}\ubar, \geop{x^2} \, ^{(\textnormal{apx})}\ubar, \geop{x^3} \, ^{(\textnormal{apx})}\ubar\right) \in \mathcal{S}$ where $^{(\textnormal{apx})}\ubar$ is the solution of \eqref{E:APPROXIMATEEIKONALFUNCTION}.

We now describe the spaces that we will run our contraction argument in. Let $d_{\mathcal{H}}$ be the Hausdorff distance on $\R^2\times\T^2$. For $\hat{\mathbf{v}} \in \mathcal{S}$, let $\hat{\mathscr{F}}$ and $\hat F$ be the unique maximal flowout domain of $\datahypfortimefunctiontwoarg{0}{[\timefunction_0,0]}$ and flow of $\hat{\uLunit}$ defined analogously to \eqref{E:DEFINITIONOFFLOWOUTDOMAINFORREGULARNEWUL}--\eqref{E:DEFINITIONOFFLOWFORREGULARNEWUL} with the appropriate $\Delta \hat a, \, \Delta \hat b, \, \widehat{\underline{\Lambda}}_{\Delta s}$ in place of $^{(\textnormal{reg})} \Delta a, \, ^{(\textnormal{reg})}\Delta b, \, ^{(\textnormal{reg})}\Lambda_{\Delta s}$. Set $\hat{\mathscr{D}} \eqdef \hat{F} (\hat{\mathscr{F}})$. Since \eqref{E:REQUIREDREGULARITYFORADMISSIBILITYCLASS} implies $\hat{\uLunit}$ is a $C^{0,1}$-vectorfield and \eqref{E:TRANSVERSALITYTODATAHYPERSURFACEOFADMISSIBILITYCLASS} implies that it is transverse to the initial data surface $\datahypfortimefunctiontwoarg{0}{[\timefunction_0,0]}$, it is well known (see e.g. the proof of \cite[Theorem 5.6]{rampazzo2007frobenius}) that $\hat F$ is a $C^{0,1}$-Lipeomorphism\footnote{A map $f \colon A \to B$ between Banach spaces is said to be a Lipeomorphism if it a homeomorphism and both $f,f^{-1}$ are Lipschitz.} with inverse $\hat{F}^{-1}$ and: 
	\begin{align}
		\left \| \hat{F} \right\|_{C^{0,1}( \hat{\mathscr{F})}} & \le C,\label{E:REQUIREDESTIMATEFORFLOWFORCONTRACTIONLASS} \\
		\left\| \hat{F}^{-1}\right\|_{C^{0,1}_{\textnormal{geo}}(\hat{\mathscr{D}})} & \le C. \label{E:REQUIREDESTIMATEFORINVERSEFLOWFORCONTRACTIONLASS}
	\end{align}  For some constants $1 < M_1,\, M_2 < \infty$ to be determined later, define the set:
\begin{subequations}
\begin{align}
	 \mathcal{S}^{M_1,M_2}_{\initialsmall} & \eqdef \Bigg\{  (\hat{\mathcal{I}},\, \hat{J}, \hat z_2, \hat z_3) \in  \mathcal{S}  \, \, \text{s.t.} \, :    \label{E:DEFINITIONOFCONTRACTIONSPACE}  \\
	& \bullet \, d_{\mathcal{H}}\left(\, ^{(\textnormal{reg})}(\mathscr{D}), \hat{\mathscr{D}}\right) \le \initialsmall^{1/2}, \label{E:FLOWOUTSOFINITIALDATASURFACEFORITERATESMUSTBECLOSETOAPPROXINGOINGEIKONALFLOWOUT} \\
	& \bullet \, \left\|  \hat{\mathcal{I}} \right\|_{C^{0,1}_{geo}(\hat{\mathscr{D}})}  \le M_1, \,  \, \left\|  \hat J \right\|_{C^{0,1}_{geo}(\hat{\mathscr{D}})} \le M_2; \label{E:C0BOUNDSFORLANDMUXDERIVATIVESOFINGOINGEIKONALCONTRACTIONS} \\ 
	& \bullet \, \left\|   \hat z_2 \right\|_{C^{0,1}_{geo}(\hat{\mathscr{D}})}, \,  \left\| \hat z_3 \right\|_{C^{0,1}_{geo}(\hat{\mathscr{D}})}   \le \initialsmall^{1/2}; \label{E:SMALLC0BOUNDSFORANGULARDERIVATIVESOFINGOINGEIKONALCONTRACTIONS} \\ 
	& \left. \bullet \, \frac{7}{9} \blowupdelta \le \min_{\hat{\mathscr{D}}} \hat{\mathcal{I}}  \le   \max_{\hat{\mathscr{D}}} \hat{\mathcal{I}} \le \frac{9}{7}\blowupdelta \right\} \label{E:LUNITOFITERATESISUNIFORMLYBOUNDED} 
\end{align}
\end{subequations}It is straight forward to check that $ \mathcal{S}^{M_1,M_2}_{\initialsmall}$ is non-empty as $^{(\textnormal{apx})}\mathbf{v} \in \mathcal{S}^{M_1,M_2}_{\initialsmall}$ for $M_1,\, M_2$ sufficiently large.

\medskip

\noindent \underline{\textbf{The mapping:}} 
Let $\mr \ubar$ be defined as in \eqref{E:RESTRICTIONOFUBARONDATASURFACE}, from which identities  \eqref{E:FUTUREDIRECTEDDATASURFACEDERIVATIVEOFOFUBAR}--\eqref{E:ANGULARDATASURFACETANGENTIALDERIVATIVESOFUBAR} hold for $\UnitTimeNormalizedUbarDataSurfaceTangentVectorfield \mr \ubar$ and $\datasurfacegeometricp{x^A}\mr \ubar$. Let $\Lunit \mr \ubar$ and $\muX \mr \ubar$ be the scalar functions defined on $\datahypfortimefunctiontwoarg{0}{[\timefunction_0,0]}$ given by inserting RHS\,\eqref{E:FUTUREDIRECTEDDATASURFACEDERIVATIVEOFOFUBAR}--\eqref{E:ANGULARDATASURFACETANGENTIALDERIVATIVESOFUBAR} into RHS\,\eqref{E:LUBARINTERMSOFDATASURFACEDERIVATIVES}--\eqref{E:MUXUBARINTERMSOFDATASURFACEDERIVATIVES}. Similarly, for $A = 2,3$, let $\geop{x^A} \mr \ubar$ be the scalar functions defined by $\geop{x^A} \mr \ubar \eqdef \datasurfacegeometricp{x^A} \mr \ubar +  \coefficientsInVectorfieldsTangetToUbarDataSurface_A \muX \mr \ubar$, see \eqref{E:ANGULARVECTORFIELDTANGETTOUBARDATASURFACE}. We note that $\Lunit \mr \ubar, \, \muX \mr \ubar$, and $\geop{x^A} \mr \ubar$ \emph{lose a derivative} from the data. That is, we have $\Lunit \mr \ubar, \, \muX \mr \ubar,\, \geop{x^A} \mr \ubar \in C_{\textnormal{data}}^{0,1} \left(\datahypfortimefunctiontwoarg{0}{[\timefunction_0,0]}\right)$ whereas $\UnitTimeNormalizedUbarDataSurfaceTangentVectorfield \mr \ubar, \, \datasurfacegeometricp{x^A} \mr \ubar \in C_{\textnormal{data}}^{1,1} \left(\datahypfortimefunctiontwoarg{0}{[\timefunction_0,0]}\right)$. 

We define $T$ as the mapping which takes $\hat{\mathbf{v}} = (\hat {\mathcal{I}},\, \hat J, \, \hat z_2, \, \hat z_3) \in \mathcal{S}$ as inputs and outputs a vector $T \hat{\mathbf{v}}$ whose components are solutions to the following integral equations for any $p \in \domainforembeddingdatahypfortimefunctiontwoarg{0}{[0,\mupositive]}$ and $(\Delta s,p) \in \hat{\mathscr{F}}$: 
\begin{subequations} \label{E:DEFINITIONOFCONTRACTIONSPACEFORINGOINGEIKONAL}
	\begin{align}
		\begin{split} \label{E:INTEGRALEQUATIONFORLUBARITERATES}
			T\hat{\mathcal{I}} \circ \hat{F}(\Delta s,p) & = \Lunit \mr \ubar \circ \embeddingdatahypfortimefunctionarg{0}(p) + \int_{\Delta s' = 0}^{\Delta s} \left\{ - \frac{1}{2} (\Lunit \upmu) \hat{\mathcal{I}} + [\muX, \Lunit]^A \hat{z}_A -  \hat{\MagnitueofinnerproductofnewLandnewuL} (\gtorus^{-1})^{AB} \hat z_B \left[\geop{x^A},\Lunit\right]^C \hat z_C \right. \\
			&  \ \ \left. + \frac{\Lunit \upmu}{2 \hat{\mathcal{I}}} |\hat z|_{\gtorus}^2 +\frac{ \upmu}{2\hat{\mathcal{I}}}  \left( \Lunit (\gtorus^{-1})^{AB}\right)\hat z_A \hat z_B \right\} \circ  \hat{F}(\Delta s',p) \, \mathrm{d} \Delta s',
		\end{split} \\
		\begin{split} \label{E:INTEGRALEQUATIONFORMUXUBARITERATES} 
			T\hat{J} \circ  \hat{F}(\Delta s,p)& =  \muX \mr \ubar \circ \embeddingdatahypfortimefunctionarg{0}(p) +  \int_{\Delta s' = 0}^{\Delta s}\left\{ - \frac{1}{2} (\muX \upmu) \hat{\mathcal{I}} + \hat{\ReciprocaluLunitAppliedtoTimeFunction}[\Lunit,\muX]^A \hat{z}_A -  \hat{\MagnitueofinnerproductofnewLandnewuL} (\gtorus^{-1})^{AB} \hat z_B \left[\geop{x^A},\muX\right]^C \hat z_C \right. \\
			&  \ \ \left. + \frac{\muX \upmu}{2 \hat{\mathcal{I}}} |\hat z|_{\gtorus}^2 +\frac{ \upmu}{2\hat{\mathcal{I}}} \left( \muX (\gtorus^{-1})^{AB}\right)\hat z_A \hat z_B \right\} \circ \hat{F}(\Delta s',p) \, \mathrm{d} \Delta s',
		\end{split} \\
		\begin{split} \label{E:INTEGRALEQUATIONFORANGULARDERIVATIVEOFUBARITERATES} 
			T\hat{z}_A  \circ \hat{F}(\Delta s,p) & = \geop{x^A} \mr \ubar \circ \embeddingdatahypfortimefunctionarg{0}(p) +  \int_{\Delta s' = 0}^{\Delta s} \left\{ - \frac{1}{2} \left(\geop{x^A} \upmu\right) \hat{\mathcal{I}} + \hat{\ReciprocaluLunitAppliedtoTimeFunction} \left[\Lunit,\geop{x^A}\right]^B \hat z_B  + \left[\muX, \geop{x^A}\right]^B \hat z_B  \right. \\
			& \ \  \left. +  \frac{\geop{x^A} \upmu}{2 \hat{\mathcal{I}}} |\hat z|_{\gtorus}^2  +\frac{ \upmu}{2 \hat{\mathcal{I}}} \left( \geop{x^A} (\gtorus^{-1})^{BC}\right)\hat z_B \hat z_C \right\} \circ  \hat{F}(\Delta s',p) \, \mathrm{d} \Delta s'.
		\end{split} 
	\end{align}
\end{subequations}
From the definition of $T$, we see that: 
\begin{align}
	\left(T \hat{\mathcal{I}}, T \hat{J}, T \hat{z}_2, T \hat{z}_3\right)\big|_{\datahypfortimefunctiontwoarg{0}{[\timefunction_0,0]}} = \left(\Lunit \mr \ubar, \muX \mr \ubar, \geop{x^2} \mr \ubar, \geop{x^3} \mr \ubar\right). \label{E:INITIALCONDITIONFORCONTRACTIONMAPOFITERATES}
\end{align}
We set: 
	\begin{subequations}
	\begin{align}
			T\hat \ReciprocalLunitAppliedtoTimeFunction & \eqdef \frac{1}{T\hat{\mathcal{I}}}, & T\hat \MagnitueofinnerproductofnewLandnewuL & \eqdef \upmu T \hat \ReciprocalLunitAppliedtoTimeFunction,  \label{E:IMAGEOFCONTRACTIONMAPDEFINITIONOFLUBARANDINNERPRODUCTOFNEWLANDNEWUL} \\
			T \hat \ReciprocaluLunitAppliedtoTimeFunction & \eqdef \frac{\upmu T \hat \ReciprocalLunitAppliedtoTimeFunction}{ T\hat \ingoingmu}, & \frac{1}{T\hat{\ingoingmu} }& \eqdef \frac{1}{2}T \hat{\mathcal{I}} + \frac{1}{2 T\mathcal{\hat{I}}} (\gtorus^{-1})^{AB} T \hat z_A T \hat z_B, \label{E:IMAGEOFCONTRACTIONMAPDEFINITIONOFLBARUANDINGOINGMU}
		\end{align}
		\begin{align}
			^{(T)}\hat{\uLunit} \eqdef (T \hat{\ReciprocaluLunitAppliedtoTimeFunction}) \Lunit + \muX - T \hat{\MagnitueofinnerproductofnewLandnewuL} (\gtorus^{-1})^{AB} T \hat{z}_B \geop{x^A}, \label{E:IMAGEOFCONTRACTIONMAPDEFINITIONOFNEWUL}
		\end{align}
	\end{subequations}
as the corresponding scalars and vectorfield from \eqref{E:ITERATIONDEFINITIONOFLUBARANDINNERPRODUCTOFNEWLANDNEWUL}--\eqref{E:ITERATIONDEFINITIONOFNEWUL} defined by the image components of the map $T$. Denote the flow map of $^{(T)}\hat{\uLunit}$ as $^{(T)}\widehat{\underline{\Lambda}}_{\Delta s}$. That is, the $(t,u,x^2,x^3)$-coordinate components of $^{(T)}\widehat{\underline{\Lambda}}_{\Delta s}$ solve:
	\begin{align}
		\frac{\p}{\p \Delta s} \,^{(T)}\widehat{\underline{\Lambda}}_{\Delta s}(t,u,x^2,x^3) = \, ^{(T)} \hat{\uLunit}\circ \, ^{(T)}\widehat{\underline{\Lambda}}_{\Delta s}(t,u,x^2,x^3), & & ^{(T)}\widehat{\underline{\Lambda}}_{0}(t,u,x^2,x^3) = (t,u,x^2,x^3). \label{E:IMAGEOFCONTRACTIONMAPDEFINITIONOFNEWULFLOWMAP}
	\end{align}
		
A fixed point of $T$ will then be a solution to the quasilinear transport system \eqref{E:QUASILINEAREQUATIONSATISFIEDBYLUBAR}--\eqref{E:QUASILINEAREQUATIONSATISFIEDBYSMOOTHANGULARDERIVATIVEOFUBAR} written in integral form. We will now show that for an appropriate choice of constants $M_1,\, M_2$, if $|\timefunction_0|, \, \initialsmall$ are sufficiently small, then: 
	\begin{align}
		T \colon  \mathcal{S}^{M_1,M_2}_{\initialsmall} \to  \mathcal{S}^{M_1,M_2}_{\initialsmall} \label{E:CONTRATIONMAPTAKESCONTRACTIONSETTOITSELF}
	\end{align}
and is a contraction relative to a suitable norm. Banach's fixed point theorem will then provide the desired unique fixed point.
		
\medskip

\noindent \underline{\textbf{$T$ is an endofunction:}} To prove \eqref{E:CONTRATIONMAPTAKESCONTRACTIONSETTOITSELF}, we first prove that $T \hat{\mathbf{v}} \in\mathcal{S}$ provided that $\hat{\mathbf{v}} \in \mathcal{S}^{M_1,M_2}_{\initialsmall}$. Under this assumptions on $\hat{\mathbf{v}}$, we have that the RHS\,\eqref{E:INTEGRALEQUATIONFORLUBARITERATES}--\eqref{E:INTEGRALEQUATIONFORANGULARDERIVATIVEOFUBARITERATES} are functions of class $C^{0,1}(\hat{\mathscr{F}})$. Composing $T \hat{\mathbf{v}} \circ \hat{F}$ with $\hat{F}^{-1}$ and using \eqref{E:REQUIREDESTIMATEFORINVERSEFLOWFORCONTRACTIONLASS} proves that \eqref{E:REQUIREDREGULARITYFORADMISSIBILITYCLASS} holds for $T \hat{\mathbf{v}}$ with $\mathscr{K} \eqdef \hat{\mathscr{D}}$. Next, since the initial conditions \eqref{E:INITIALCONDITIONFORCONTRACTIONMAPOFITERATES} hold, the same proof of \eqref{E:NEWULISTRANSVERSALTOUBARDATASURFACE} shows that $^{(T)} \hat{\uLunit}$ is transverse to $\datahypfortimefunctiontwoarg{0}{[\timefunction_0,0]}$ and so \eqref{E:TRANSVERSALITYTODATAHYPERSURFACEOFADMISSIBILITYCLASS} is also satisfied for $T\hat{\mathbf{v}}$. 

Under the assumptions that $\hat{\mathbf{v}} \in \mathcal{S}^{M_1,M_2}_{\initialsmall}$, estimates \eqref{E:C01INFINITYESTIMATESFORLUBARONINITIALDATASURFACE}--\eqref{E:C01INFINITYESTIMATESFORANGULARUBARONINITIALDATASURFACE},  \eqref{E:REQUIREDESTIMATEFORFLOWFORCONTRACTIONLASS}, and \eqref{E:C0BOUNDSFORLANDMUXDERIVATIVESOFINGOINGEIKONALCONTRACTIONS}--\eqref{E:LUNITOFITERATESISUNIFORMLYBOUNDED} imply the following bounds:
\begin{subequations}
	\begin{align} 
		\left\| T \hat{\mathcal{I}} \circ \hat{F} \right\|_{C^{0,1}(\hat{\mathscr{F}})}  \lesssim 1, & &
		\left\| T \hat{J} \circ \hat{F} \right\|_{C^{0,1}(\hat{\mathscr{F}})}  \lesssim 1,  & & \left\| T \hat{z}_A \circ \hat{F} \right\|_{C^{0,1}(\hat{\mathscr{F}})}   \lesssim \initialsmall. \label{E:FIRSTESTIMATEFORIMAGEOFUNDERCONTRACTIONMAP}
	\end{align}
Composing with $\hat{F}^{-1}$ and using \eqref{E:REQUIREDESTIMATEFORINVERSEFLOWFORCONTRACTIONLASS} proves 
	\begin{align} 
		\left\| T \hat{\mathcal{I}}\right\|_{C^{0,1}_{\textnormal{geo}}(\hat{\mathscr{D}})}  \lesssim 1, & &
		\left\| T \hat{J} \right\|_{C^{0,1}_{\textnormal{geo}}(\hat{\mathscr{D}})}  \lesssim 1,  & & \left\| T \hat{z}_A  \right\|_{C^{0,1}_{\textnormal{geo}}(\hat{\mathscr{D}})}   \lesssim \initialsmall. \label{E:SECONDESTIMATEFORIMAGEOFUNDERCONTRACTIONMAP}
	\end{align}
\end{subequations} 
The bounds in \eqref{E:SECONDESTIMATEFORIMAGEOFUNDERCONTRACTIONMAP} implies that $^{(T)} \hat{\uLunit}$ is a $C^{0,1}_{\textnormal{geo}}(\hat{\mathscr{D}})$-vectorfield and hence there exists a maximal flowout domain $^{(T)}\hat{\mathscr{F}} \subset \hat{\mathscr{F}}$, a flow $^{(T)}\hat{F}$, and a flowout $^{(T)}\hat{\mathscr{D}} \eqdef \, ^{(T)}\hat{F}\left( \, ^{(T)}\hat{\mathscr{F}}\right)$, such that $^{(T)}\hat{F}$ is a $C^{0,1}$-Lipeomorphism with inverse $^{(T)}\hat{F}^{-1}$ and:  
	\begin{align}
		\left \| \, ^{(T)}\hat{F} \right\|_{C^{0,1}(\, ^{(T)}\hat{\mathscr{F})}} & \le C,\label{E:ESTIMATEFORFLOWMAPOFCONTRACTIONMAPIMAGE} \\
		\left\| \, ^{(T)}\hat{F}^{-1}\right\|_{C^{0,1}_{\textnormal{geo}}(\, ^{(T)}\hat{\mathscr{D}})} & \le C. \label{E:ESTIMATEFORINVERSEFLOWMAPOFCONTRACTIONMAPIMAGE}
	\end{align}
Since $^{(T)}\hat{\mathscr{D}} \subset \hat{\mathscr{D}}$, \eqref{E:SECONDESTIMATEFORIMAGEOFUNDERCONTRACTIONMAP} implies that \eqref{E:C0BOUNDSFORLANDMUXDERIVATIVESOFINGOINGEIKONALCONTRACTIONS}--\eqref{E:SMALLC0BOUNDSFORANGULARDERIVATIVESOFINGOINGEIKONALCONTRACTIONS} hold for the components of $T \hat{\mathbf{v}}$, and in fact the $C^{0,1}_{\textnormal{geo}}$-norm of $T \hat{z}_A$ is even smaller than that of $\hat{z}_A$.

To prove that \eqref{E:LUNITOFITERATESISUNIFORMLYBOUNDED} holds for $T \hat{\mathcal{I}}$, we first note the following estimates, which follow from \eqref{E:BOUNDSFORREGULARNEWULFLOWOUTDOMAININTERVALDEFININGFUNCTIONS} and \eqref{E:FLOWOUTSOFINITIALDATASURFACEFORITERATESMUSTBECLOSETOAPPROXINGOINGEIKONALFLOWOUT}:
\begin{align}
	\sup_{\domainforembeddingdatahypfortimefunctiontwoarg{0}{[0,\mupositive]}} \, \Delta \hat b \le \frac{1}{2} \interestingu + \initialsmall^{1/2} & & \inf_{\domainforembeddingdatahypfortimefunctiontwoarg{0}{[0,\mupositive]}} \, \Delta \hat a \ge - \frac{1}{2}\interestingu - \initialsmall^{1/2}. \label{E:BOUNDSFORITERATESOFNEWULFLOWOUTDOMAININTERVALDEFININGFUNCTIONS}
\end{align}
We will actually prove \eqref{E:LUNITOFITERATESISUNIFORMLYBOUNDED} on $\hat{\mathscr{D}}$ for the $\Delta s \ge 0$. This will automatically prove the desired estimate on $^{(T)}\hat{\mathscr{D}}\subset \hat{\mathscr{D}}$. The proof for $\Delta s < 0$ holds from the similar arguments and we omit it. Since $\Lunit \upmu \approx - \blowupdelta < 0$ in $[-\interestingu,\interestingu]$, using the assumptions on the set $ \mathcal{S}^{M_1,M_2}_{\initialsmall}$, the integral equation \eqref{E:INTEGRALEQUATIONFORLUBARITERATES} implies
\[ T \hat{\mathcal{I}} \circ \hat{F}(\Delta s,p) \ge \Lunit \mr \ubar \circ \embeddingdatahypfortimefunctionarg{0}(p) + C \initialsmall^{1/2}.\]
The desired  bound $\min_{^{(T)}\hat{\mathscr{D}}} T \hat{\mathcal{I}} \ge \frac{7}{9} \blowupdelta$ then follows from \eqref{E:BOUNDSONLMUINTERESTINGREGIONFROMTAUFOLIATIONS}, \eqref{E:LUBARISLIKEMINUSLMUALONGDATASURFACE},  and choosing $\initialsmall$ sufficiently small. To prove that $ \max_{^{(T)}\hat{\mathscr{D}}} T \hat{\mathcal{I}}  \le \frac{9}{7} \blowupdelta$, we use the already known bound on $\hat{\mathcal{I}}$ (by virtue of being an element of the set  $ \mathcal{S}^{M_1,M_2}_{\initialsmall}$), the estimate \eqref{E:BOUNDSONLMUINTERESTINGREGIONFROMTAUFOLIATIONS}, and the integral equation \eqref{E:INTEGRALEQUATIONFORLUBARITERATES} to bound: 
\begin{align} \label{E:UPPERBOUNDFORIMAGEOFLUBARITERATESINPOSITIVEFLOWREGION}
	 T\hat{\mathcal{I}}\circ \hat{F}(\Delta s,p) \le  \Lunit \mr \ubar \circ \embeddingdatahypfortimefunctionarg{0}(p) +  \frac{81}{56} \blowupdelta^2 \Delta s + C \initialsmall^{1/2}.
\end{align}
The desired upper bound for $T \hat{\mathcal{I}}$ follows from \eqref{E:BOUNDSFORITERATESOFNEWULFLOWOUTDOMAININTERVALDEFININGFUNCTIONS} and the smallness of $\interestingu$ in \eqref{E:BOUNDFORINTERESTINGU}.

Next, we estimate using \eqref{E:IMAGEOFCONTRACTIONMAPDEFINITIONOFLUBARANDINNERPRODUCTOFNEWLANDNEWUL}--\eqref{E:IMAGEOFCONTRACTIONMAPDEFINITIONOFNEWUL} and the now improved bound \eqref{E:SECONDESTIMATEFORIMAGEOFUNDERCONTRACTIONMAP}:
\begin{align}
	\left(\, ^{(T)}\hat{\uLunit} - \, ^{(\textnormal{reg})}\newuL\right)t  = \mathcal{O}(\initialsmall), & & \left(\, ^{(T)}\hat{\uLunit} - \, ^{(\textnormal{reg})}\newuL\right)u = 0, & & \left(\, ^{(T)}\hat{\uLunit} - \, ^{(\textnormal{reg})}\newuL\right)x^A = \mathcal{O}(\initialsmall). \label{E:DIFFERENCEBETWEENIMAGEOFNEWULITERATESANDREGULARNEWUL}
\end{align}
Since the corresponding estimates in \eqref{E:DIFFERENCEBETWEENIMAGEOFNEWULITERATESANDREGULARNEWUL} for $\hat{\uLunit} -\, ^{(\textnormal{reg})} \newuL$ are at best $\mathcal{O}(\initialsmall^{1/2})$, it follows that $d_{\mathcal{H}}(\, ^{(\textnormal{reg})}\mathscr{D}, \, ^{(T)}\hat{\mathscr{D}}) \le d_{\mathcal{H}}(\, ^{(\textnormal{reg})}\mathscr{D},\hat{\mathscr{D}})$.

\medskip

\noindent \underline{\textbf{The contraction:}} 

We now prove that $T$ is a contraction on $ \mathcal{S}^{M_1,M_2}_{\initialsmall}$ using a suitably defined metric. For $i = 1,2$, let $ (\hat{\mathcal{I}}_i,\hat{J}_i,\hat{z}_{2,i}, \hat{z}_{3,i})$ be the corresponding components of an element $\hat{\mathbf{v}}_i \in \mathcal{S}^{M_1,M_2}_{\initialsmall}$. Let $\hat{\mathscr{D}}_i$ be their respective flowouts of $\datahypfortimefunctiontwoarg{0}{[\timefunction_0,0]}$. Consider the metric $d \colon  \mathcal{S}^{M_1,M_2}_{\initialsmall} \times  \mathcal{S}^{M_1,M_2}_{\initialsmall} \to \R_{\ge 0}$ defined by: 
\begin{align}
	\begin{split} \label{E:CONTRACTIONMAPPINGMETRIC}  
		 d( \hat{\mathbf{v}}_1, \hat{\mathbf{v}}_2)  &  \eqdef \max \left\{ d_{\mathcal{H}}\left(\hat{\mathscr{D}}_1,\hat{\mathscr{D}}_2\right), \, \left\|  \hat{\mathcal{I}}_1 -  \hat{\mathcal{I}}_2\right\|_{ C^{0,1}_{\textnormal{geo}}\left(\hat{\mathscr{D}}_1\cap \hat{\mathscr{D}}_2\right)}, \, \left\|   \hat{J}_1 -   \hat{J}_2 \right\|_{ C^{0,1}_{\textnormal{geo}}\left(\hat{\mathscr{D}}_1\cap \hat{\mathscr{D}}_2\right)}, \right. \\
		& \ \ \left. \left\| \hat{z}_{2,1} -  \hat{z}_{2,2} \right\|_{ C^{0,1}_{\textnormal{geo}}\left(\hat{\mathscr{D}}_1\cap \hat{\mathscr{D}}_2\right)}, \,  \left\|  \hat{z}_{3,1} -   \hat{z}_{3,2} \right\|_{ C^{0,1}_{\textnormal{geo}}\left(\hat{\mathscr{D}}_1\cap \hat{\mathscr{D}}_2\right)}\right\}.
	\end{split}
\end{align}
The contraction bound $d_{\mathcal{H}}(\, ^{(T)} \hat{\mathscr{D}}_1, \, ^{(T)}\hat{\mathscr{D}}_2) \le \frac{1}{2}  d_{\mathcal{H}}(\hat{\mathscr{D}}_1,\hat{\mathscr{D}}_2)$ follows from the estimates:
\begin{align}
	\left(\, ^{(T)}\hat{\uLunit}_1 - \, ^{(T)}\hat{\uLunit}_2\right)t  = \mathcal{O}(\initialsmall^4), & & \left(\, ^{(T)}\hat{\uLunit}_1 - \, ^{(T)} \hat{\uLunit}_2\right)u = 0, & & \left(\, ^{(T)}\hat{\uLunit}_1 - \, ^{(T)}\hat{\uLunit}_2\right)x^A = \mathcal{O}(\initialsmall). \label{E:DIFFERENCEBETWEENIMAGEOFNEWULITERATES} \\
	\left(\hat{\uLunit}_1 - \hat{\uLunit}_2\right)t  = \mathcal{O}(\initialsmall^2), & & \left(\hat{\uLunit}_1 - \hat{\uLunit}_2\right)u = 0, & & \left(\hat{\uLunit}_1 - \hat{\uLunit}_2\right)x^A = \mathcal{O}(\initialsmall^{1/2}). \label{E:DIFFERENCEBETWEENEWULITERATES}
\end{align}Next we note the following identity, which is a straightforward consequence of the definition \eqref{E:ITERATIONDEFINITIONOFNEWUL}:
\begin{align}
	\hat{\uLunit}_1 T \hat{\mathcal{I}}_2 & = \hat{\uLunit}_2 T \hat{\mathcal{I}}_2 + (\hat \ReciprocaluLunitAppliedtoTimeFunction_1 - \hat \ReciprocaluLunitAppliedtoTimeFunction_2) \Lunit T \hat{\mathcal{I}}_2 - \left( \hat \MagnitueofinnerproductofnewLandnewuL_1 (\gtorus^{-1})^{AB} \hat{z}_{A,1} -  \hat \MagnitueofinnerproductofnewLandnewuL_2 (\gtorus^{-1})^{AB} \hat{z}_{A,2}\right) \geop{x^B}  T \hat{\mathcal{I}}_2. \label{E:TRANSPORTEQUATIONFORI2WITHRESPECTTOL1ITERATE}
\end{align}
Using the definition of $T \hat{\mathcal{I}}_2$ as in \eqref{E:INTEGRALEQUATIONFORLUBARITERATES}, we have:
\begin{align}
	\begin{split} \label{E:TRANSPORTEQUATIONFORI2WITHRESPECTTOL1ITERATEWITHRHSINSERTED}
		\hat{\uLunit}_1 T \hat{\mathcal{I}}_2 & = - \frac{1}{2} (\Lunit \upmu) \hat{\mathcal{I}}_2 + [\muX, \Lunit]^A \hat{z}_{A,2} -  \hat{\MagnitueofinnerproductofnewLandnewuL}_2 (\gtorus^{-1})^{AB} \hat z_{B,2} \left[\geop{x^A},\Lunit\right]^C \hat z_{C,2} + \frac{\Lunit \upmu}{2 \hat{\mathcal{I}}_2} |\hat z_2|_{\gtorus}^2 +\frac{ \upmu}{2\hat{\mathcal{I}}_2}  \left( \Lunit (\gtorus^{-1})^{AB}\right)\hat z_{A,2} \hat z_{B,2} \\
		& \ \ + (\hat \ReciprocaluLunitAppliedtoTimeFunction_1 - \hat \ReciprocaluLunitAppliedtoTimeFunction_2) \Lunit T \hat{\mathcal{I}}_2 - \left( \hat \MagnitueofinnerproductofnewLandnewuL_1 (\gtorus^{-1})^{AB} \hat{z}_{A,1} -  \hat \MagnitueofinnerproductofnewLandnewuL_2 (\gtorus^{-1})^{AB} \hat{z}_{A,2}\right) \geop{x^B}  T \hat{\mathcal{I}}_2.
	\end{split}
\end{align}
Let $p \in \, ^{(T)}\hat{\mathscr{D}}_1 \cap \, ^{(T)} \hat{\mathscr{D}}_2$. In particular, it follows that $p \in \hat{\mathscr{D}}_1$ and $p \in \hat{\mathscr{D}}_2$.  Since $\hat{F}_1$ is a Lipeomorphism onto $\hat{\mathscr{D}}_1$, there is a $(\Delta s,t,x^2,x^3) \in \hat{\mathscr{F}}_1$ such that $p = \hat{F}_1(\Delta s,t,x^2,x^3)$. Note that by Rademacher's theorem and \eqref{E:SECONDESTIMATEFORIMAGEOFUNDERCONTRACTIONMAP}, $|\Lunit T \hat{\mathcal{I}}_2|, |\geop{x^A} T \hat{\mathcal{I}}_2| \lesssim 1$. Using the fact that the difference $T \hat{\mathcal{I}}_1 - T \hat{\mathcal{I}}_2$ vanishes along $\datahypfortimefunctiontwoarg{0}{[\timefunction_0,0]}$, we then have the following estimate along integral curves of $\hat{\uLunit}_1$ by subtracting the integral of RHS\,\eqref{E:TRANSPORTEQUATIONFORI2WITHRESPECTTOL1ITERATEWITHRHSINSERTED} from the RHS\,\eqref{E:INTEGRALEQUATIONFORLUBARITERATES} for $\hat{\mathcal{I}}_1$:
\begin{align} \label{E:BOUNDSONDIFFERENCEOFLMUINTERATES}
	\left| \left\{ T \hat{\mathcal{I}}_1 - T \hat{\mathcal{I}}_2\right\}(p) \right| \lesssim \left\| \Lunit \upmu \right\|_{L^\infty}\left\|  \hat{\mathcal{I}}_1 -  \hat{\mathcal{I}}_2 \right\|_{C^{0,1}_{\textnormal{geo}}\left(\hat{\mathscr{D}}_1 \cap \hat{\mathscr{D}}_2\right)}  |\Delta s|  + C|\Delta s|  \initialsmall^{1/2}  d(\hat{\mathbf{v}}_1,\hat{\mathbf{v}}_2).
\end{align}
By \eqref{E:BOUNDSFORITERATESOFNEWULFLOWOUTDOMAININTERVALDEFININGFUNCTIONS}, it follows that $\left\| T \hat{\mathcal{I}}_1 - T \hat{\mathcal{I}}_2\right\|_{ C^{0,1}_{\textnormal{geo}}\left( \, ^{(T)} \hat{\mathscr{D}}_1\cap \, ^{(T)}\hat{\mathscr{D}}_2\right)} \le \frac{1}{2} d(\hat{\mathbf{v}}_1,\hat{\mathbf{v}}_2)$ provided that $\initialsmall$ is sufficiently small due to the smallness of \eqref{E:BOUNDFORINTERESTINGU}. The rest of the differences appearing on the RHS\,\eqref{E:CONTRACTIONMAPPINGMETRIC} under the image of $T$ can similarly be shown to be bounded from above by $\frac{1}{2} d(\hat{\mathbf{v}}_1,\hat{\mathbf{v}}_2)$. We omit the details. 

\medskip

\noindent \underline{\textbf{Properties of the fixed point:}}

The preceding steps proved that the mapping $T$ given by \eqref{E:INTEGRALEQUATIONFORLUBARITERATES}--\eqref{E:INTEGRALEQUATIONFORANGULARDERIVATIVEOFUBARITERATES} has a unique fixed point. That is, we have a unique solution $(\Lunit \ubar, \muX \ubar, \geop{x^2} \ubar, \geop{x^3} \ubar)$ to the system of equations \eqref{E:QUASILINEAREQUATIONSATISFIEDBYLUBAR}--\eqref{E:QUASILINEAREQUATIONSATISFIEDBYSMOOTHANGULARDERIVATIVEOFUBAR} on a domain $\underline{\mathscr{D}}$ satisfying \eqref{E:TOPLATERALBOUNDARYOFSPACETIMEREGIONFORINGOINGEIKONAL}--\eqref{E:SECONDBOUNDARYPIECEOLEFTBOUNDARYSURFACE}, on which the estimates \eqref{E:ANISOTROPICPOINTWISEESTIMATESFORFIRSTDERIVATIVESOFINGOINGEIKONALONFLOWOUTDOMAIN} and \eqref{E:SHARPUNIFORMESTIMATEFORLUBARINFLOWOUTREGION} hold. The initial conditions \eqref{E:INITIALCONDITIONFORCONTRACTIONMAPOFITERATES} imply that \eqref{E:INITIALCONDITIONSFORSECONDORDERSYSTEMFORINGOINGEIKONALFUNCTION} holds and so the results of Lemma\,\ref{L:QUASILINEARTRANSPORTSYSTEMFORDERIVATIVESOFUBAR} prove that we have a solution $\ubar \in C^{1,1}_{\textnormal{geo}}(\underline{\mathscr{D}})$ of \eqref{E:INGOINGEIKONALEQUATIONINITIALVALUEPROBLEM}. The regularity of $\ubar$ implies the regularity of the boundary components $\underline{\mathscr{D}}_{\textnormal{left}}$ and $\underline{\mathscr{D}}_{\textnormal{right}}$. To prove \eqref{E:LOCATIONOFUBARFLOWOUTDOMAIN}, we note that $^{(\textnormal{reg})} \newuL t \le \newuL t$ and $^{(\textnormal{reg})} \newuL u =  \newuL u = 1$ imply $\underline{\mathscr{D}} \subset \, ^{(\textnormal{reg})}\mathscr{D}$. Hence, \eqref{E:LOCATIONOFUBARFLOWOUTDOMAIN} follows from \eqref{E:FLOWOUTOFREGULARNEWULISCONTAINEDININTERESTINGREGION}. Estimate \eqref{E:ANISOTROPICPOINTWISEESTIMATESFORHIGHERNEWULDERIVATIVESOFINGOINGEIKONALONFLOWOUTDOMAIN} holds by estimating directly the RHS of \eqref{E:QUASILINEAREQUATIONSATISFIEDBYLUBAR}--\eqref{E:QUASILINEAREQUATIONSATISFIEDBYSMOOTHANGULARDERIVATIVEOFUBAR}.

 This also proves \eqref{E:INITIALOUTGOINGCHARACTERISTICSURFACE}.  To see why \eqref{E:INITIALINGOINGCHARACTERISTICSURFACE} and \eqref{E:CONDITIONTHATGUARANTEESNONTRIVIALPARTOFCAUCHYHORIZON} hold, recall $\scalarembeddingdatahypfortimefunctionarg{0}\left(\domainforembeddingdatahypfortimefunctiontwoarg{0}{[0,\mupositive]}\right) \subset [ -\tfrac{1}{2} \interestingu, \tfrac{1}{2} \interestingu]$ by \eqref{E:MUXMUKAPPALEVELSETLOCATION}. Since $\scalarembeddingdatahypfortimefunctionarg{0}\left( \Cartesiantisafunctiononmumxtoriarg{\mathfrak{m}}{0}(x^2,x^3),x^2,x^3\right) = \Eikonalisafunctiononmumuxtoriarg{\mathfrak{m}}{0}(x^2,x^3)$ for any $(\mathfrak{m},x^2,x^3) \in [0,\mupositive]\times\T^2$, we have that $u \circ \ \flowmapofnewularg{\Delta s}(\Cartesiantisafunctiononmumxtoriarg{\frac{1}{2}\mupositive}{0}(x^2,x^3),x^2,x^3) = \Eikonalisafunctiononmumuxtoriarg{\frac{1}{2}\mupositive}{0}(x^2,x^3) + \Delta s$. The results follow from \eqref{E:BOUNDSFORITERATESOFNEWULFLOWOUTDOMAININTERVALDEFININGFUNCTIONS}. 
\end{proof}

\section{Parameters, their size assumptions, and conventions for constants}
\label{S:PARAMETERSANDSIZEASSUMPTIONSANDCONVENTIONSFORCONSTANTS}
In this section, we list and describe the ``size-parameters'' that appear throughout the paper.
These parameters will play a crucial role in Sect.\,\ref{S:ASSUMPTIONSONTHEDATA},
when we describe our assumptions on the data.
Then, in Sect.\,\ref{SS:CONVENTIONSFORCONSTANTS}, 
we state our conventions for how constants 
appearing in our analysis,
such as $C$ and $C_{\mydiam}$,
are allowed to depend on the parameters.

\subsection{Parameters} 
\label{SS:PARAMETERS}

\subsubsection{Parameters of the perturbed solutions}
\label{SSS:PARAMETERSOFPERTURBEDSOLUTIONS}

\subsubsection{Informal description of the parameters}
\label{SSS:INFORMALDESCRIPTIONOFTHEPARAMETERS}
To help guide the reader,
we now informally describe the role that various parameters
will play in our forthcoming analysis.
\begin{itemize}
		\item $\Ntop \geq 14$ represents the maximum number of times we commute the equations
			when we derive energy estimates. 
		\item The parameters $\timefunction_0$ and $\mupositive$ are related by $\timefunction_0 = - \mupositive$ 
			(see Def.\,\ref{D:RELATIONBETWEENINITIALROUGHTIMESLICEANDSMALLMUVALUE}).
			We view $\timefunction_0$ to be the ``initial rough time,'' i.e., the value of 
			$\timefunctionarg{0}$ corresponding to the initial state of the solution near the singularity.
		\item The parameter $\mupositive$ is the minimum value of $\upmu$ along the initial rough hypersurface portion
			$\hypthreearg{\timefunction_0}{[-\interestingu,\interestingu]}{0}$.
			For convenience, we assume that $\mupositive$ is small. While the smallness of
			 $\mupositive$ is not essential, it allows us to
			focus on studying the solution only near the singularity and allows us to give short proofs of various estimates.
		\item The parameter $\leftubar < 0$ corresponds to the initial value of the ingoing eikonal function $\ubar$ in the region where we will derive a priori estimates.
		\item The parameters $0 < \rightu$,
			$0 < \leftu$, 
			$0 < \interestingu < \min\{1,\leftu,\rightu,\frac{1}{22}\blowupdelta^{-1}\}$
			delineate ranges of values for the outgoing eikonal function $u$ used in \cite{abbrescia2022emergence}. The relevant parameters for the current paper are $- \interestingu \le \moreinterestingu_1 <  - \frac{3}{4} \interestingu$ and $\frac{3}{4} \interestingu < \moreinterestingu_2 \le \interestingu$. 
		\item The parameter $\mathring{\upalpha}$ measures the $L^{\infty}$-size of the amplitude of 
			$\RRiemann$ along the initial ingoing characteristic hypersurface $\ingoingcharacteristicsurfacetwoarg{\leftubar}{[\moreinterestingu_1,\moreinterestingu_2]}$.
		\item  The parameter $\boringregionmupositive$ quantifies the positivity of $\upmu$  
			when $u \in [-\interestingu,\interestingu]\setminus[\moreinterestingu_1,\moreinterestingu_2]$;
			see \eqref{E:DATAMUISLARGEINBORINGREGION}.
		\item  The parameter $\blowupdeltadoublenull$ measures the size of 
			the crucial factor that drives the blowup; see definition \eqref{E:DELTASTARDOUBLENULLDEF}. 
		\item The parameter $\mathring{\updelta}$ measures the $L^{\infty}$-size of the transversal derivatives of $\RRiemann$ along 
			the initial ingoing characteristic hypersurface $\ingoingcharacteristicsurfacetwoarg{\leftubar}{[\moreinterestingu_1,\moreinterestingu_2]}$. 
			It also controls the $L^{\infty}$ size of the transversal derivatives
			of various geometric quantities constructed out of the eikonal function.
			We make no smallness assumptions on $\mathring{\updelta}$. 
	\item The parameter $\initialsmall \geq 0$ measures the extent to which the solution
			``breaks the simple isentropic plane-symmetry'' relative to the rough foliations $ \hypthreearg{\timefunction}{[- \rightu,\leftu]}{0}$.
	\item The parameter $\initialsmalldoublenull \ge 0$ plays a similar role as $\initialsmall$, but measures the extent to which the solution
			``breaks the simple isentropic plane-symmetry'' relative to the double-null foliations. 
	\item $\fundbootsmall$ is a small ``bootstrap parameter'' first appearing in 
		Sect.\,\ref{SS:MAINQUANTITATIVEBOOTSTRAPASSUMPTIONS}.
	\item The parameter $\secondtransversalderivativemulowerbound$ 
		quantifies the transversal convexity of $\upmu$ 
		(namely, the positive size of various second-order $\nullhyparg{u}$-transversal derivatives of $\upmu$, 
		such as $\muX \muX \upmu$); 
		see, for example, \eqref{E:DATATASSUMPTIONMUTRANSVERSALCONVEXITY}.
\end{itemize}

\subsection{Parameter size assumptions}
\label{SS:PARAMETERSIZEASSUMPTIONS}
In this section, we state the size assumptions on the parameters that are sufficient
for our main results to hold, i.e., 
for Theorems~\ref{T:EXISTENCEUPTOCAUCHYHORIZONBYCONTINUATIONCRITERIA} to hold.

\begin{quote}
For the remainder of the article, 
when we say that ``$A$ is small relative to $B$,''
we mean that $A \geq 0$, that $B > 0$, and that
there exists a continuous increasing function\footnote{Although we do not specify their form, 
the functions $f$ could always be chosen to be 
polynomials with positive coefficients or exponentials of such polynomials.}
$f :(0,\infty) \rightarrow (0,\infty)$ 
such that 
$
A < f(B)
$.
The functions $f$ are allowed to depend on the equation of state.
\end{quote}

\begin{assumption}[Size assumptions on the parameters]
	\label{ASSUMPTIONS:RELATIVESIZEASSUMPTIONS} \hfill
\begin{itemize}
	\item To close our estimates, we assume that the regularity-parameter $\Ntop$ is an integer satisfying:
		\begin{align} \label{E:NTOPLARGENESSASSUMPTION}
			\Ntop \geq 15.
		\end{align}
	\item We assume that the following parameters are positive, 
		but they do not have to be small or large: 
		$\rightu$,
		$\interestingu$,
		$\leftu$, 
		$-\moreinterestingu_1$,
		$\moreinterestingu_2$
		$\mathring{\updelta}$, 
		$\mathring{\updelta}_*$, 
		$\blowupdeltadoublenull$, 
		and
		$\boringregionmupositive$.
	\item We assume that $0 < \interestingu < \min\{1,\leftu,\rightu,\frac{1}{22}\blowupdelta^{-1}\}$ and $- \interestingu \le \moreinterestingu_1 <  - \frac{3}{4} \interestingu$ and $\frac{3}{4} \interestingu < \moreinterestingu_2 \le \interestingu$. 
	\item We assume that $0 < \secondtransversalderivativemulowerbound < 1$. We make no other assumptions
		on $\secondtransversalderivativemulowerbound$.
	\item We assume that $\mathring{\upalpha} > 0$, and that $\mathring{\upalpha}$ is small relative to $1$.
	\item We assume that 
			$\mupositive$, $|\timefunction_0|$, $|\leftubar|$
			are small relative to 
			$1$,
			$\mathring{\updelta}^{-1}$,
			$\mathring{\updelta}_*$,
			and $\secondtransversalderivativemulowerbound$.
			In particular, we assume:
			\begin{align} \label{E:INITIALTIMEFUNCTIONVALUELESSTHANONE}
				|\leftubar|, |\timefunction_0|
				& < 1.
			\end{align}
	\item We assume that $\initialsmall$ is small relative to 
		$1$,
		$\mathring{\upalpha}$,
		$\rightu$,
		$\interestingu$,
		$\leftu$, 
		$\mathring{\updelta}^{-1}$, 
		$\mathring{\updelta}_*$,
		$\mupositive$, $-\leftubar$,
		and
		$\secondtransversalderivativemulowerbound$.
		\item 
			We assume that
			$\mupositive < \frac{\boringregionmupositive}{2}$,
			and
			$\timefunction_0 = - \mupositive$.
		\item We assume that $\initialsmalldoublenull = C \initialsmall$, where $\initialsmall$ is the constant from \cite{abbrescia2022emergence} and $C$ is the largest implicit constant appearing in the construction of $\ubar$ in Sect.\,\ref{S:CONSTRUCTIONANDESTIMATESFORINGOINGEIKONALFUNCTION}.
		\item Our main results will hold under the assumption that
			$\fundbootsmall = C \initialsmalldoublenull$ for some large constant $C$,
			where $\fundbootsmall$ is the bootstrap parameter first appearing in
			Sects.\,\ref{SS:MAINQUANTITATIVEBOOTSTRAPASSUMPTIONS}--\ref{SSS:AUXBOOTSTRAP}.
			This is consistent with the following parameter-size relations, which we assume in order to simplify our
			bootstrap argument:
			\begin{subequations}
			\begin{align} \label{E:DATAEPSILONISSMALLERTHANBOOTSTRAPEPSILONSMALLERTHANSQUAREOFDATAALPHA}
			\initialsmalldoublenull 
			& \leq 
			\fundbootsmall
			\leq 
			\mathring{\upalpha}^2,
				\\
			\fundbootsmall^{3/2}
			& 
			\leq
			\initialsmalldoublenull.
			\label{E:NONLINEARINEQUALITYRELATINGDATAEPSILONANDBOOTSTRAPEPSILON}
			\end{align}
			\end{subequations}
\end{itemize}
\end{assumption}

\subsection{Conventions for constants}
\label{SS:CONVENTIONSFORCONSTANTS}
In this section, we state our conventions for how the constants $C$, $C_*$, $\mathfrak{c}$, and $C_{\mydiam}$ 
appearing in our analysis are allowed to depend on the parameters introduced above.

\begin{itemize}
\item The constants $C$, $C_*$, and $\mathfrak{c}$ are free to vary from line to line,
	and we use them in a similar fashion; we mainly use ``$C$'' in our estimates, 
	introducing $C_*$ only in a few arguments in which multiple constants play a role.
These constants are allowed to depend on the nonlinearities (i.e., on the equation of state),
and they can continuously depend on the quantities
$\rightu$,
$\interestingu$, 
$\leftu$, $-\moreinterestingu_1$, $\moreinterestingu_2$,
$\mathring{\updelta}^{-1}$, 
$\mathring{\updelta}_*$,
$\boringregionmupositive^{-1}$, 
and
$\secondtransversalderivativemulowerbound^{-1}$
from Sect.\,\ref{SSS:INFORMALDESCRIPTIONOFTHEPARAMETERS}. 
In particular, $C$ and $\mathfrak{c}$ are allowed, in principle, 
to \emph{increase} with respect to 
$\rightu$,
$\interestingu$, 
$\leftu$, $-\moreinterestingu_1$, $\moreinterestingu_2$,
$\mathring{\updelta}^{-1}$, 
$\mathring{\updelta}_*$,
$\boringregionmupositive^{-1}$, 
and
$\secondtransversalderivativemulowerbound^{-1}$.
However, $C$, $C_*$, and can be chosen to be \textbf{independent of the parameters}
$\mathring{\upalpha}$, 
$\initialsmalldoublenull$, 
$\fundbootsmall$, 
$\leftubar$, and
$\mulevelsetvalue_0$.
under the smallness assumptions of Sect.\,\ref{SS:PARAMETERSIZEASSUMPTIONS}.
\item $A \lesssim B$ means that there exists a constant $C > 0$ (where $C$ has the properties described above) 
such that $A \leq C B$.
\item $A \approx B$ means that $A \lesssim B$ and $B \lesssim A$.
\item $A = \mathcal{O}(B)$ means that $|A| \lesssim |B|$. 
\item Constants $C_{\mydiam}$ are also allowed to vary from line to line.
\begin{quote}
However, unlike $C$ and $C_*$,
the $C_{\mydiam}$ are \textbf{universal} 
in the sense that under the smallness assumptions of Sect.\,\ref{SS:PARAMETERSIZEASSUMPTIONS}, 
they can be chosen to be \textbf{independent} of 
\begin{align}
\rightu, 
	\,
\interestingu, 
	\, 
\leftu, \, -\moreinterestingu_1, \, \moreinterestingu_2, 
	\,
\timefunction_0,
	\, 
\mupositive, \leftubar, 
	\, 
\boringregionmupositive,
	\, 
\mathring{\updelta}, 
	\, 
\mathring{\updelta}_*, \, \blowupdeltadoublenull, 
	\, 
\mathring{\upalpha},	
	\, 
\secondtransversalderivativemulowerbound, \initialsmalldoublenull \label{E:WHATDIAMONDCONSTANTSAREINDEPENDENTOF}
\end{align}
\end{quote}
\item $A = \mathcal{O}_{\mydiam}(B)$ means that there exists a constant $C_{\mydiam} > 0$
(where $C_{\mydiam}$ has the properties described above) such that 
$|A| \leq C_{\mydiam} |B|$.
\item As examples, we note that $\mathring{\updelta} \initialsmall \leq 1 \eqdef C_{\mydiam}$
and $\mathring{\updelta}^2 \initialsmall \leq 1 \eqdef C_{\mydiam}$
(because $\initialsmall$ is assumed to be small relative to $\mathring{\updelta}^{-1}$ and 
relative to increasing functions of $\mathring{\updelta}^{-1}$, such as $\mathring{\updelta}^{-2}$),
that $10 \mathring{\upalpha}^2 \leq C_{\mydiam} \mathring{\upalpha}$,
while we have only $\mathring{\updelta} \mathring{\upalpha} \leq C$
(i.e., our smallness assumptions on $\mathring{\upalpha}$ are not strong enough to ensure
that $\mathring{\updelta} \mathring{\upalpha}$ is small because $\mathring{\updelta}$ might be large).
\item As another example, in our estimates 
(e.g., the proof of \eqref{E:MUTRANSVERSALCONVEXITY}), 
by assuming that $|\leftubar|$ is small
and using that $C$ is independent of $|\leftubar|$
(in particular, $C$ does not implicitly contain any factors of $\frac{1}{|\leftubar|}$),
we can ensure that
$C |\leftubar| \leq \frac{1}{\secondtransversalderivativemulowerbound}$
and
$C |\leftubar| \leq \frac{\secondtransversalderivativemulowerbound}{20}$. 

\end{itemize}

\section{Assumptions on the data} 
\label{S:ASSUMPTIONSONTHEDATA}
In this section, we state our assumptions 
on the data in terms of the parameters listed in Sect.\,\ref{SS:PARAMETERS}.

We begin with the following definition, where we introduce the notation for the powers of $\weight$ used throughout the paper. 
\begin{definition}[Notation for powers of $\weight$] \label{D:NOTATIONFORPOWERSOFWEIGHT}
	For $N \le \Ntop$, we set:
	\begin{subequations}
		\begin{align}
			\blowuprateofwaveWRTnumberofcommutations(N) & \eqdef \max\{0,\blowupratetoporderwave - 2\Ntop + 2N\},\label{E:BLOWUPRATEOFWAVEWRTNUMBEROFCOMMUTATORS} \\
			\blowuprateoftransportWRTnumberofcommutations(N) & \eqdef \max\{0, \blowupratetoporderwave - 1 - 2 \Ntop + 2N\},\label{E:BLOWUPRATEOFTRANSPORTWRTNUMBEROFCOMMUTATORS} \\				
			\blowuprateofacousticgeoWRTnumberofcommutations(N) & \eqdef \begin{cases}
			\blowupratetoporderacoustic & N = \Ntop, \\
			\max\{0, \blowupratetoporderwave + 1 - 2 \Ntop + 2N\} & N \le \Ntop - 1. 
			\end{cases}, \label{E:BLOWUPRATEOFACOUSTICBELOWTOPWRTNUMBEROFCOMMUTATORS} 
		\end{align}
	\end{subequations}
	In this work, $\blowuprateofwaveWRTnumberofcommutations(\Ntop)$ will correspond to the necessary power of $\weight$ needed to close the top-order flux-energy estimates for the velocity $\tander^{\Ntop} \velocityarray$. The number $\blowuprateoftransportWRTnumberofcommutations(\Ntop)$ denotes the power of $\weight$ used in the $L^2$ flux-energy estimates for the vorticity and entropy. The number $\blowuprateofacousticgeoWRTnumberofcommutations(\Ntop)$ denotes the power of $\weight$ used in the $L^2$ flux-energy estimates for the acoustic geometry and modified fluid variables $\VortVort,\DivGradEnt$. 
	
We assume that: 
	\begin{subequations}
		\begin{align}
			\blowupratetoporderwave & = 6.13, \label{E:BLOWUPRATEDEFFORWAVE} \\
			\blowupratetoporderacoustic & = 7.63.\label{E:BLOWUPRATEDEFFORACOUSTICGEO}
		\end{align}
	\end{subequations}
\end{definition}

\subsection{The initial data assumptions}
\label{SS:ASSUMPTIONSONDATA}

\subsubsection{Localized assumptions of the initial ingoing characteristic surface} \label{SS:LOCALIZEDASSUMPTIONSINITIALINGOINGCHARACTERISTICSURFACE}

We now make assumptions on the $u$-width of the initial ingoing characteristic surface that are necessary in order derive estimates for the fluid all the way up to a nontrivial portion of the Cauchy horizon. In Appendix\,\ref{A:DATASSUMPTIONS}, we prove that the classical solution on $\twoargMrough{[\timefunction_0,0],[- \rightu,\leftu]}{0}$ induces data on the initial characteristic surfaces satisfying the assumptions of this section. 

We assume that there are constants: 
	\begin{subequations}
		\begin{align} 
			\leftubar & \in [\frac{3}{4}\timefunction_0,\tfrac{1}{4}\timefunction_0], \label{E:MOSTNEGATIVEVALUEOFUBAR} \\
			\moreinterestingu_1 & \in [- \interestingu, -\tfrac{3}{4}\interestingu), \label{E:MOSTNEGATIVEMOREINTERESTINGVALUEOFU}\\
			\moreinterestingu_2 & \in (\tfrac{3}{4} \interestingu, \interestingu], \label{E:MOSTPOSITIVEMOREINTERESTINGVALUEOFU}
		\end{align} 
	\end{subequations}
such that: 
	\begin{subequations}
		\begin{align}
			\ingoingcharacteristicsurfacetwoarg{\leftubar}{[\moreinterestingu_1,\moreinterestingu_2]} & \subset \underline{\mathscr{D}}\cap \twoargMrough{[\timefunction_0,\frac{1}{100}\timefunction_0],[- \rightu,u]}{0},  \label{E:INITIALINGOINGCHARACTERISTICSURFACE} \\
			\outgoingcharacteristicsurfacetwoarg{\moreinterestingu_1}{[\leftubar,0]} & \subset \underline{\mathscr{D}},  \label{E:INITIALOUTGOINGCHARACTERISTICSURFACE}
		\end{align}
	\end{subequations}
	and: 
	\begin{align} \label{E:CONDITIONTHATGUARANTEESNONTRIVIALPARTOFCAUCHYHORIZON}
		\moreinterestingu_2 - \moreinterestingu_1 > \sup_{(x^2,x^3) \in \T^2} \Eikonalisafunctiononmumuxtoriarg{0}{0}(x^2,x^3) - \moreinterestingu_1.
	\end{align}
\subsubsection{Localized assumptions on the data of $\upmu$ and its derivatives on $\ingoingcharacteristicsurfacetwoarg{\leftubar}{[\moreinterestingu_1,\moreinterestingu_2]}$}
\label{SSS:LOCALIZEDDATAASSUMPTIONSFORMUANDDERIVATIVES}

We now make localized quantitative and qualitative assumptions on the behavior of $\upmu$ along $\ingoingcharacteristicsurfacetwoarg{\leftubar}{[\moreinterestingu_1,\moreinterestingu_2]}$. We start by defining the data-parameter $\blowupdeltadoublenull$.

\begin{definition}[Key data-parameter tied to the Cartesian time of first blowup]
\label{D:DELTASTARDOUBLENULLDEF}
We define $\blowupdeltadoublenull$ by: 
\begin{align} \label{E:DELTASTARDOUBLENULLDEF}
	\blowupdeltadoublenull
	& 
	\eqdef 
	\sup_{\ingoingcharacteristicsurfacetwoarg{\leftubar}{[\moreinterestingu_1,\moreinterestingu_2]}} 
	 \frac{1}{2} \left[ G_{\Lunit \Lunit}^0 \muX \RRiemann\right]_-.
\end{align}
\end{definition}

We assume that:
\begin{align} \label{E:DELTASTARDOUBLENULLPOSITIVE}
	\blowupdeltadoublenull 
	& > 0.
\end{align}\begin{remark}[Connection between $\blowupdeltadoublenull$, $\mathring{\updelta}_*$, and the Cartesian time of first blowup]
\label{R:CONNECTOINBETWEENDELTASTARANDDELTASTARRBLOWUPTIMEANDSHOCKTIME}
For simple isentropic plane-symmetric solutions, the Cartesian time of first shock formation is precisely $\frac{1}{\mathring{\updelta}_*}$; 
see \cite{abbrescia2022emergence}*{Theorem~2.8}.
Our main results of \cite{abbrescia2022emergence} show that for the perturbed solutions under study, 
the Cartesian time of first blowup, which we denote here by $t_{\mbox{\tiny First shock}}$,
satisfies $t_{\mbox{\tiny First shock}} = \lbrace 1 + \mathcal{O}(\initialsmall) \rbrace \frac{1}{\mathring{\updelta}_*}$. The constant $\blowupdeltadoublenull$ is the analog of $\blowupdelta$ for the double-null foliations.
\end{remark}

\noindent  \underline{\textbf{Transversal convexity of $\upmu$:}} We assume that:
\begin{align} 
		\frac{\secondtransversalderivativemulowerbound}{4} \le  \min_{\ingoingcharacteristicsurfacetwoarg{\leftubar}{[\moreinterestingu_1,\moreinterestingu_2]}} \left\{  \newuL \newuL \upmu, \newuL \muX  \upmu, \, \muX \muX \upmu \right\}  & \le   \max_{\ingoingcharacteristicsurfacetwoarg{\leftubar}{[\moreinterestingu_1,\moreinterestingu_2]}} \left\{ \newuL \newuL \upmu,\, \newuL \muX \upmu, \,\muX \muX \upmu \right\} \le \frac{4}{\secondtransversalderivativemulowerbound}, \label{E:DATATASSUMPTIONMUTRANSVERSALCONVEXITY} \\
		\frac{\secondtransversalderivativemulowerbound}{4} \le  \essinf_{\ingoingcharacteristicsurfacetwoarg{\leftubar}{[\moreinterestingu_1,\moreinterestingu_2]}} \nullgeop{u}\newuL \upmu & \le   \esssup_{\ingoingcharacteristicsurfacetwoarg{\leftubar}{[\moreinterestingu_1,\moreinterestingu_2]}}  \nullgeop{u}\newuL \upmu \le \frac{4}{\secondtransversalderivativemulowerbound}. \label{E:DATATASSUMPTIONMUSLIGHTLYLESSREGULARTRANSVERSALCONVEXITY}
	\end{align}
\noindent \underline{\textbf{Quantitative negativity of $\Lunit \upmu$:}}
We assume that the following inequalities hold, 
where $\blowupdeltadoublenull$ is defined in \eqref{E:DELTASTARDOUBLENULLDEF}:
\begin{align} \label{E:DATATASSUMPTIONLMUQUANTITATIVENEGATIVITY}
	- 
	\frac{9}{8}
	\blowupdeltadoublenull
	\leq
	\min_{\ingoingcharacteristicsurfacetwoarg{\leftubar}{[\moreinterestingu_1,\moreinterestingu_2]}} \Lunit \upmu
	\leq 
	\max_{\ingoingcharacteristicsurfacetwoarg{\leftubar}{[\moreinterestingu_1,\moreinterestingu_2]}} \Lunit \upmu
	\leq 
	- 
	\frac{7}{8}
	\blowupdeltadoublenull.
\end{align}\medskip
	
\noindent  \underline{\textbf{$\argubarnewultorus{\leftubar}$
is located near $\lbrace u=0 \rbrace$.}}
With $\argubarnewultorus{\leftubar}$ as in \eqref{E:DEFOFNEWULADAPTEDTORI}, 
we assume that:
\begin{subequations}
\begin{align} \label{E:DATATASSUMPTIONUBARNEWULMULEVELSETLOCATION}
			\argubarnewultorus{\leftubar}
			& 
			\subset
			\ingoingcharacteristicsurfacetwoarg{\leftubar}{[-\frac{1}{2}\interestingu,\frac{1}{2}\moreinterestingu_2 + \frac{1}{8} \moreinterestingu]},
				\\
			\min_{\ingoingcharacteristicsurfacetwoarg{\leftubar}{[\moreinterestingu_1,\moreinterestingu_2]}
			\setminus
			\ingoingcharacteristicsurfacetwoarg{\leftubar}{[-\frac{1}{2}\interestingu,\frac{1}{2}\moreinterestingu_2 + \frac{1}{8} \moreinterestingu]}}
			|\newuL \upmu| 
			& \geq \frac{\secondtransversalderivativemulowerbound \interestingu}{4}.
			\label{E:DATATASSUMPTIONREGIONWHERENEWULMULEVELSETISNOTLOCATED}
\end{align}
\end{subequations}
\medskip

\noindent \underline{\textbf{Quantitative positivity of $\upmu$ away from $u = 0$:}}
We assume that there is a constant $\boringregionmupositive > 0$ such that:
\begin{align} \label{E:DATAMUISLARGEINBORINGREGION}
			\min_{\ingoingcharacteristicsurfacetwoarg{\leftubar}{[\moreinterestingu_1,\moreinterestingu_2]} \setminus \ingoingcharacteristicsurfacetwoarg{\leftubar}{[-\frac{3}{4}\interestingu,\frac{3}{4}\interestingu]}} \upmu
			& \geq \boringregionmupositive.
\end{align}
We also assume:
		\begin{align} \label{E:MU1BIGGERTHANU0}
			\frac{\boringregionmupositive}{2} > -\leftubar.
		\end{align}\subsubsection{Quantitative assumptions on the data of the fluid and eikonal function quantities along
$\ingoingcharacteristicsurfacetwoarg{\leftubar}{[\moreinterestingu_1,\moreinterestingu_2]}$ and
$\doublenulltoritwoarg{\ubar}{u}$}
	\label{SSS:QUANTITATIVEASSUMPTIONSONDATAAWAYFROMSYMMETRY} 
In this section, we state quantitative assumptions on the data of the wave variables and the eikonal
function quantities along
$\ingoingcharacteristicsurfacetwoarg{\leftubar}{[\moreinterestingu_1,\moreinterestingu_2]}$ and
$\doublenulltoritwoarg{\leftubar}{u}$.

We refer to Sect.\,\ref{S:COMMUTATIONFORMULASIDENTITIESANDSTRINGSOFCOMMUTATORS} for notation regarding
strings of commutation vectorfields and to Defs.\,\ref{D:GEOMETRICL2NORMS}
and\,\ref{D:GEOLINFINITYTORUS} for the definitions of our $L^2$ and $L^{\infty}$ norms. We also refer to Def.\,\ref{D:NOTATIONFORPOWERSOFWEIGHT} regarding the weights we will use in our defining $L^2$ norms.

\medskip

\noindent \underline{\textbf{$L^{\infty}$ assumptions on the wave variables and their pure transversal derivatives}}.
For $u \in [\moreinterestingu_1,\moreinterestingu_2]$ and $M = 1,2,3,4$, we assume 
(recall that $\wavearray$ and $\wavearraypartial$ are defined in Def.\,\ref{D:ARRAYSOFWAVEVARIABLES}):
\begin{subequations}
\begin{align}  
	\left \| \RRiemann \right \|_{L^{\infty}\left(\doublenulltoritwoarg{\leftubar}{u}\right)} 
	& \leq \mathring{\upalpha}, 
	\label{E:LINFINITYINITIALROUGHHYPERSURFACEBOUNDRRIEMANNAMPLITUDE}
		\\
	\left \| 
		\muX^M \RRiemann 
	\right\|_{L^{\infty}\left(\doublenulltoritwoarg{\leftubar}{u}\right)} 
	& 
	\leq \mathring{\updelta},  
		\label{E:LINFINITYINITIALROUGHHYPERSURFACEBOUNDRRIEMANNTRANSVERSALDERIVATIVES}
		\\
	\left\| 
		\wavearraypartial
	\right\|_{L^{\infty}\left(\doublenulltoritwoarg{\leftubar}{u}\right)},
		\,
	\left\| 
		\muX^M \wavearraypartial
	\right\|_{L^{\infty}\left(\doublenulltoritwoarg{\leftubar}{u}\right)}
	& 
	\leq \initialsmalldoublenull. 
	\label{E:LINFINITYINITIALROUGHHYPERSURFACEBOUNDSMALLPURETRANSVERSALDERIVATIVESSMALLWAVEVARIABLES} 
\end{align}
\end{subequations}\medskip

\noindent \underline{\textbf{$L^{\infty}$ assumptions involving tangential derivatives of the wave variables}}.
For $u \in [\moreinterestingu_1,\moreinterestingu_2]$, we assume:
	\begin{equation} 	\label{E:LINFINITYINITIALDOUBLENULLTORIBOUNDSWAVE} 
	\begin{split}
		&
		\left \| 
			\tander^{[1, \Ntop - 6]}\wavearray 
		\right \|_{L^{\infty}\left(\doublenulltoritwoarg{\leftubar}{u}\right)}, 
		\, 
		\left\| 
			\comdersmall^{[1,\Ntop - 7;1]} \wavearray 
		\right \|_{L^{\infty}\left(\doublenulltoritwoarg{\leftubar}{u}\right)}, 
		\, 
		\left\| 
			\comdersmall^{[1,6];2} \wavearray 
		\right \|_{L^{\infty}\left(\doublenulltoritwoarg{\leftubar}{u}\right)}, 
		\\
		&
		\left\| 
			\comdersmall^{[1,5];3} \wavearray 
		\right \|_{L^{\infty}\left(\doublenulltoritwoarg{\leftubar}{u}\right)},
			\,
		\left\| 
			\Lunit \muX \muX \muX \muX \wavearray 
		\right\|_{L^{\infty}\left(\doublenulltoritwoarg{\leftubar}{u}\right)} 
		\leq \initialsmalldoublenull. 
	\end{split}
	\end{equation}

\medskip

\noindent \underline{\textbf{$L^2$ assumptions along $\ingoingcharacteristicsurfacetwoarg{\leftubar}{[\moreinterestingu_1,\moreinterestingu_2]}$}}.
For $N \in \{1,\dots,\Ntop\}$, we assume:
\begin{subequations}
\begin{align}
\left\| 
	\newuL \tander^N\velocityarray 
\right\|_{L^2_{\blowuprateofwaveWRTnumberofcommutations(N)}\left(\ingoingcharacteristicsurfacetwoarg{\leftubar}{[\moreinterestingu_1,\moreinterestingu_2]} \right)} 
& \leq  \initialsmalldoublenull, \label{E:TRANSVERSALDERIVATIVEOFTANGENTIALL2NORMSOFWAVEVARIABLESSMALLALONGINITIALINGOINGNULLSURFACE}  
		\\
\left\| \sqrt{\upmu}
	\nullangD \tander^N\velocityarray 
\right\|_{L^2_{\blowuprateofwaveWRTnumberofcommutations(N)}\left(\ingoingcharacteristicsurfacetwoarg{\leftubar}{[\moreinterestingu_1,\moreinterestingu_2]} \right)} 
& \leq  \initialsmalldoublenull. \label{E:NULLANGULARDERIVATIVEOFTANGENTIALL2NORMSOFWAVEVARIABLESSMALLALONGINITIALROUGHHYPERSURFACE}
\end{align}
\end{subequations}In addition, for $N \in \{0,\dots,\Ntop\}$, we have: 
\begin{subequations}
	\begin{align}
		\left\| \sqrt{\ReciprocaluLunitAppliedtoTimeFunction} \, \tander^N(\vortrenormalized,\GradEnt)\right\|_{L^2_{\blowuprateoftransportWRTnumberofcommutations(N)}\left(\ingoingcharacteristicsurfacetwoarg{\leftubar}{[\moreinterestingu_1,\moreinterestingu_2]} \right)} 
		& \leq  \initialsmalldoublenull, 		\label{E:L2NORMSOFVORTANDGRADENTSMALLALONGINITIALINGOINGNULLSURFACE}  
		\\
		\left\| \sqrt{\ReciprocaluLunitAppliedtoTimeFunction} \, \tander^N(\VortVort,\DivGradEnt)\right\|_{L^2_{\blowuprateofacousticgeoWRTnumberofcommutations(N)}\left(\outgoingcharacteristicsurfacetwoarg{\leftubar}{[\moreinterestingu_1,\moreinterestingu_2]} \right)} 
		& \leq  \initialsmalldoublenull.\label{E:L2NORMSOFVORTVORTANDDIVGRADENTSMALLALONGINITIALINGOINGNULLSURFACE}  
\end{align}
\end{subequations}
\medskip

\noindent \underline{\textbf{$L^2$ assumptions along $\outgoingcharacteristicsurfacetwoarg{\moreinterestingu_1}{[\leftubar,0]}$}}.
For $N \in \{1,\dots,\Ntop\}$, we assume:
\begin{subequations}
\begin{align}
\left\| 
	L \tander^N\velocityarray 
\right\|_{L^2_{\blowuprateofwaveWRTnumberofcommutations(N)}\left(\outgoingcharacteristicsurfacetwoarg{\moreinterestingu_1}{[\leftubar,0]} \right)} 
& \leq  \initialsmalldoublenull, 
	\label{E:LDERIVATIVEOFTANGENTIALL2NORMSOFWAVEVARIABLESSMALLALONGINITIALOUTGOINGCHARACTERISTICHYPERSURFACE}  
		\\
\left\| \sqrt{\upmu}
	\angD \tander^N\velocityarray 
\right\|_{L^2_{\blowuprateofwaveWRTnumberofcommutations(N)}\left(\outgoingcharacteristicsurfacetwoarg{\moreinterestingu_1}{[\leftubar,0]} \right)} 
& \leq  \initialsmalldoublenull. \label{E:SMOOTHANGULARDERIVATIVEOFTANGENTIALL2NORMSOFWAVEVARIABLESSMALLALONGINITIALOUTGOINGCHARACTERISTICHYPERSURFACE}
\end{align}
\end{subequations}In addition, for $N \in \{0,\dots,\Ntop\}$, we have: 
\begin{subequations}
	\begin{align}
		\left\| \tander^N(\vortrenormalized,\GradEnt)\right\|_{L^2_{\blowuprateoftransportWRTnumberofcommutations(N)}\left(\outgoingcharacteristicsurfacetwoarg{\moreinterestingu_1}{[\leftubar,0]} \right)} 
		& \leq  \initialsmalldoublenull, 		\label{E:L2NORMSOFVORTANDGRADENTSMALLALONGINITIALOUTGOINGNULLSURFACE}  
		\\
		\left\|  \tander^N(\VortVort,\DivGradEnt)\right\|_{L^2_{\blowuprateofacousticgeoWRTnumberofcommutations(N)}\left(\outgoingcharacteristicsurfacetwoarg{\moreinterestingu_1}{[\leftubar,0]} \right)} 
		& \leq  \initialsmalldoublenull.\label{E:L2NORMSOFVORTVORTANDDIVGRADENTSMALLALONGINITIALOUTGOINGNULLSURFACE}  
\end{align}
\end{subequations}
\medskip

\noindent \underline{\textbf{$L^2$ assumptions along $\doublenulltoritwoarg{\leftubar}{u}$}}.
For $N \in \{1,\dots,\Ntop\}$ and $u \in [\moreinterestingu_1,\moreinterestingu_2]$, we assume:
	\begin{subequations}
		\begin{align}
			\left\| 
			\tander^{N} \velocityarray 
			\right\|_{L^2_{\blowuprateofwaveWRTnumberofcommutations(N)}(\doublenulltoritwoarg{\leftubar}{u})}
			& \leq \initialsmalldoublenull.  \label{E:SMALLDATAOFPSIONINITIALDOUBLENULLTORI} \\
			\left\| \tander^{ N-1} (\vortrenormalized,\GradEnt) \right\|_{L^2_{\blowuprateoftransportWRTnumberofcommutations(N)}(\doublenulltoritwoarg{\leftubar}{u})} & \leq  \initialsmalldoublenull.  \label{E:SMALLDATAOFVORTICITYANDGRADENTROPYONDOUBLENULLTORIONPU0}  \\
		\left\| \tander^{ N-1} (\VortVort,\DivGradEnt) \right\|_{L^2_{\blowuprateofacousticgeoWRTnumberofcommutations(N)}(\doublenulltoritwoarg{\leftubar}{u})} & \leq  \initialsmalldoublenull.  \label{E:SMALLDATAOFVORTVORTANDDIVGRADENTROPYONDOUBLENULLTORIONPU0}
		\end{align}  
	\end{subequations}

\medskip

\noindent \underline{\textbf{$L^2$ control of $\tander^{\Ntop}(\vortrenormalized,\GradEnt)$ for the elliptic hyperbolic identities.}}	
In addition, we assume the following for every $\ubar \in [\leftubar,0]$ and $u \in [\moreinterestingu_1,\moreinterestingu_2]$: 
	\begin{subequations}
		\begin{align}
			\left\| \frac{1}{\sqrt{\upmu}} \tander^{\Ntop} (\vortrenormalized,\GradEnt) \right\|_{L^2_{\blowupratetoporderacoustic}(\doublenulltoritwoarg{\leftubar}{u})} & \leq  \initialsmalldoublenull, \label{E:SMALLDATAOFTOPORDERVORTICITYANDGRADENTROPYONDOUBLENULLTORIONUBAR0FORELLIPTICHYPERBOLICIDENTITIES} \\
			\left\| \frac{1}{\sqrt{\upmu}} \tander^{\Ntop} (\vortrenormalized,\GradEnt) \right\|_{L^2_{\blowupratetoporderacoustic}(\doublenulltoritwoarg{\ubar}{\moreinterestingu_1})} & \leq \initialsmalldoublenull. \label{E:SMALLDATAOFTOPORDERVORTICITYANDGRADENTROPYONDOUBLENULLTORIONU1FORELLIPTICHYPERBOLICIDENTITIES}
		\end{align}
	\end{subequations}
	


\noindent \underline{\textbf{$L^{\infty}$ assumptions tied to transversal derivatives of the eikonal function quantities}}.
For $u \in [\moreinterestingu_1,\moreinterestingu_2]$ and $M = 0,1,2,3$, we assume:\footnote{Some of these ``assumptions'' 
can in fact be derived as a consequence of other assumptions,
up to constant factors that we absorb into the parameters
$\mathring{\updelta}$, 
$\mathring{\upalpha}$, 
and $\initialsmall$. 
For convenience, instead of deriving those ``assumptions,'' we just assume them.
The same is true for other assumptions on the eikonal function quantities stated below.
For example, even though the data of the $\muX$-derivatives of $\Lunit^i$
could be controlled via the identity \eqref{E:MUXLISCHEMATICIDENTITY},
we just assume \eqref{E:LINFINITYINITIALROUGHHYPERSURFACETRANSVERSALDERIVATIVESOFL1}.
We also refer to \eqref{E:MUTRANSPORT}
for intuition behind the form of 
RHS\,\eqref{E:LINFINITYINITIALROUGHHYPERSURFACELDERIVATIVEOFMUANDTRANSVERSALDERIVATIVES},
and to Remark\,\ref{R:CONNECTOINBETWEENDELTASTARANDDELTASTARRBLOWUPTIMEANDSHOCKTIME} for intuition
behind the factor $\blowupdeltadoublenull$
in the second term on RHS\,\eqref{E:LINFINITYINITIALROUGHHYPERSURFACEMUANDTRANSVERSALDERIVATIVES}.
\label{FN:INTUITIONONMUDATAASSUMPTIONS}}
	\begin{subequations} 
	\begin{align}  \label{E:LINFINITYINITIALROUGHHYPERSURFACELDERIVATIVEOFMUANDTRANSVERSALDERIVATIVES} 
	\left\| \Lunit \muX^M \upmu \right\|_{L^{\infty}(\doublenulltoritwoarg{\leftubar}{u})}
	& \leq 
			\frac{1}{2} 
			\left\|
				\muX^M
				\left\lbrace
					G_{\Lunit\Lunit}^0
					\muX \RRiemann 
				\right\rbrace
			\right\|_{L^{\infty}(\doublenulltoritwoarg{\leftubar}{u})}
		+
		\initialsmalldoublenull,
			\\
	\left\| \muX^M (\upmu -1)\right \|_{L^{\infty}(\doublenulltoritwoarg{\leftubar}{u})}
	& \leq
		\frac{1}{2 \blowupdeltadoublenull}
		\left\|
				\muX^M
				\left\lbrace
				G_{\Lunit\Lunit}^0
				\muX \RRiemann 
				\right\rbrace
			\right\|_{L^{\infty}(\doublenulltoritwoarg{\leftubar}{u})}
		+
		\initialsmalldoublenull.
		\label{E:LINFINITYINITIALROUGHHYPERSURFACEMUANDTRANSVERSALDERIVATIVES}  
\end{align}Moreover, for $u \in [- \rightu,\leftu]$ and $M = 1,2,3,4$, we assume:
\begin{align}
	 \left\| \Lunit \muX^M \Lsmall^i \right\|_{L^{\infty}(\doublenulltoritwoarg{\leftubar}{u})}
	& \leq 
			\initialsmalldoublenull,
			\label{E:LINFINITYINITIALROUGHHYPERSURFACELDERIVATIVEOFTRANSVERSALDERIVATIVESOFLI}  
				\\
	\left\| \muX^M \Lsmall^1 \right \|_{L^{\infty}(\doublenulltoritwoarg{\leftubar}{u})}
	& \leq
		\mathring{\updelta},
		\label{E:LINFINITYINITIALROUGHHYPERSURFACETRANSVERSALDERIVATIVESOFL1}
				\\
	\left\| \muX^M \Lsmall^A \right \|_{L^{\infty}(\doublenulltoritwoarg{\leftubar}{u})}
		& \leq
		\initialsmalldoublenull.
		\label{E:LINFINITYINITIALROUGHHYPERSURFACETRANSVERSALDERIVATIVESOFLA}
	\end{align}
	\end{subequations}\medskip

\noindent \underline{\textbf{$L^{\infty}$ assumptions involving tangential derivatives of the eikonal function quantities}}.
For $u \in [\moreinterestingu_1,\moreinterestingu_2]$, we assume:
	\begin{subequations}
	\begin{align}
	\left\| \tander_*^{[1,\Ntop-8]} \upmu \right\|_{L^{\infty}(\doublenulltoritwoarg{\leftubar}{u})},
			\,
	\left\| \comderdoublesmall^{[1,5];1} \upmu \right\|_{L^{\infty}(\doublenulltoritwoarg{\leftubar}{u})},
		\, 
	\left\| \comderdoublesmall^{[1,4];2} \upmu \right\|_{L^{\infty}(\doublenulltoritwoarg{\leftubar}{u})}
	& \leq \initialsmalldoublenull,
		\label{E:LINFINITYINITIALROUGHHYPERSURFACESMALLBOUNDSMUTRANSVERSALANDTANDERIVATIVE}  
			\\
	\left\|\Lsmall^1 \right\|_{L^{\infty}(\doublenulltoritwoarg{\leftubar}{u})} 
	& 
	\leq \mathring{\upalpha},	
		\label{E:LINFINITYINITIALROUGHHYPERSURFACEBOUNDL1SMALLAMPLITUDE}
			\\
	\begin{split}
	\left\| 
		\Lsmall^A 
	\right\|_{L^{\infty}(\doublenulltoritwoarg{\leftubar}{u})},
		\,
	\left\| 
		\tander^{[1, \Ntop-7]} \Lsmall^i 
	\right\|_{L^{\infty}(\doublenulltoritwoarg{\leftubar}{u})}, 
		\, 
	\left\| 
			\comdersmall^{[1,\Ntop-8];1} \Lsmall^i 
	\right\|_{L^{\infty}(\doublenulltoritwoarg{\leftubar}{u})}, 
		\label{E:LINFINITYINITIALROUGHHYPERSURFACESMALLBOUNDLISMALLMIXEDTRANSVERSALTANGENTIAL}
			\\
	\left\| 
		\comdersmall^{[1,5];2} \Lsmall^i 
	\right\|_{L^{\infty}(\doublenulltoritwoarg{\leftubar}{u})},
		\, 
	\left\| 
		\comdersmall^{[1,4];3} \Lsmall^i 
	\right\|_{L^{\infty}(\doublenulltoritwoarg{\leftubar}{u})}
	& 
	\leq \initialsmalldoublenull.
	\end{split}
	\end{align}
	\end{subequations}
\medskip

\noindent \underline{\textbf{$L^{2}$ assumptions involving tangential derivatives of the eikonal function quantities}}.

For $u \in [\moreinterestingu_1,\moreinterestingu_2]$ and $1 \le N \le \Ntop-1$, we assume:
\begin{subequations}
\begin{align}
	\left\|  
		\tandersmall^{[1,N+1]} \upmu 
	\right \|_{L^2_{\blowuprateofacousticgeoWRTnumberofcommutations(N)}\left( \ingoingcharacteristicsurfacetwoarg{\leftubar}{[\moreinterestingu_1,u]}\right)}, 
		\, 
	\left\| 
		\tander^{[1,N+1]} \Lsmall^i 
	\right \|_{L^2_{\blowuprateofacousticgeoWRTnumberofcommutations(N)}\left( \ingoingcharacteristicsurfacetwoarg{\leftubar}{[\moreinterestingu_1,u]}\right)}, 
		\, 
	\left\| 
		\comdersmall^{[1,N+1];1} \Lsmall^i 
	\right\|_{L^2_{\blowuprateofacousticgeoWRTnumberofcommutations(N)}\left( \ingoingcharacteristicsurfacetwoarg{\leftubar}{[\moreinterestingu_1,u]}\right)}
	& \leq \initialsmall,
		\label{E:L2DATASMALLNESSFOREIKONALFUNCTIONQUANTITIES} 
			\\
	\left\| 
		\angLie_{\comder}^{\leq N;1} \upchi 
	\right\|_{L^2_{\blowuprateofacousticgeoWRTnumberofcommutations(N)}\left( \ingoingcharacteristicsurfacetwoarg{\leftubar}{[\moreinterestingu_1,u]}\right)},
		\,
	\left\| 
		\angLie_{\tander}^{\leq N} \upchi 
	\right\|_{L^2_{\blowuprateofacousticgeoWRTnumberofcommutations(N)}\left( \ingoingcharacteristicsurfacetwoarg{\leftubar}{[\moreinterestingu_1,u]}\right)}
	& \leq \initialsmall.
		\label{E:L2DATASMALLNESSFORCHIUPTOBELOWTOPORDER}
\end{align} 
\end{subequations}
For $u \in [\moreinterestingu_1,\moreinterestingu_2]$ and $N = \Ntop$, we assume:
\begin{subequations}
	\begin{align}
		\left\| \fullymodquant{\tander^{N}} \right\|_{L^2_{\blowupratetoporderwave + 1}\left( \ingoingcharacteristicsurfacetwoarg{\leftubar}{[\moreinterestingu_1,u]}\right)}, \, \left\|\ReciprocalLunitAppliedtoTimeFunction \sqrt{(\muX v^1)^2 + (\muX v^2)^2 + (\muX v^3)^2} \, \fullymodquant{\tander^{N}}  \right\|_{L^2_{\blowupratetoporderwave}\left( \ingoingcharacteristicsurfacetwoarg{\leftubar}{[\moreinterestingu_1,u]}\right)} & \le \initialsmall, \label{E:L2DATASSUMPTIONFORFULLYMODIFIEDCHITOPORDER}\\
		 \left\| \sqrt{(\muX v^1)^2 + (\muX v^2)^2 + (\muX v^3)^2} \, \partialmodquant{\tanderY^{N-1}} \right\|_{L^2_{\blowupratetoporderwave - 1}\left( \ingoingcharacteristicsurfacetwoarg{\leftubar}{[\moreinterestingu_1,u]}\right)}  \le \initialsmall. \label{E:L2DATASSUMPTIONFORPARTIALLYMODIFIEDCHITOPORDER}
	\end{align} 
\end{subequations}
\section{The bootstrap assumptions, except those concerning the wave energies}
\label{S:BOOTSTRAPEVERYTHINGEXCEPTENERGIES}
In proving our main results, we rely on a continuity argument
based on deriving improvements of a set of bootstrap assumptions
for the solution on a region of the form $\characteristicdiamondtwoarg{[\leftubar,\ubarboot)}{[\moreinterestingu_1,\moreinterestingu_2]}$. In this section, we set up the bootstrap argument and 
state all the bootstrap assumptions, 
except for the ones concerning $L^2$-type energies, 
which we provide in Sect.\,\ref{SS:BOOTSTRAPASSUMPTIONSFORTHEWAVEENERGIES}.

\subsection{The start of the bootstrap argument: the bootstrap time interval $[\leftubar,\ubarboot)$}
\label{SS:BOOTSTRAPTIMEINTERVAL}
From now until Sect.\,\ref{S:EXISTENCEUPTOCAUCHYHORIZONBYCONTINUATIONCRITERIA}, 
we assume that there is a classical solution on an ``open-at-the-top'' region of the form
$\characteristicdiamondtwoarg{[\leftubar,\ubarboot)}{[\moreinterestingu_1,\moreinterestingu_2]}$,
and we will state our bootstrap assumptions on the same region. 
Here and throughout, 
\begin{align} \label{E:BOOTSTRAPTIME}
	\ubarboot & \in (\frac{3}{4}\leftubar,0)
\end{align}
is the ``bootstrap $\ubar$-time.''
Hence, at the start\footnote{In 
Lemma\,\ref{L:PROPERTIESANDDIFFEOMORPHICEXTENSIONOFDOUBLENULLCOORDINATES},
we will show that $\ubar$ can be suitably extended to have range $[\leftubar,\ubarboot]$,
and in Theorem\,\ref{T:EXISTENCEUPTOCAUCHYHORIZONBYCONTINUATIONCRITERIA},
we will show that it can be suitably extended to have range $[\leftubar,0]$.} 
of our bootstrap argument,
$\ubar$ has range $[\leftubar,\ubarboot)$,
which contains $[\leftubar, \frac{3}{4} \leftubar]$.

\begin{center}
	\begin{figure}  
		\begin{overpic}[scale=.7, grid = false, tics=5]{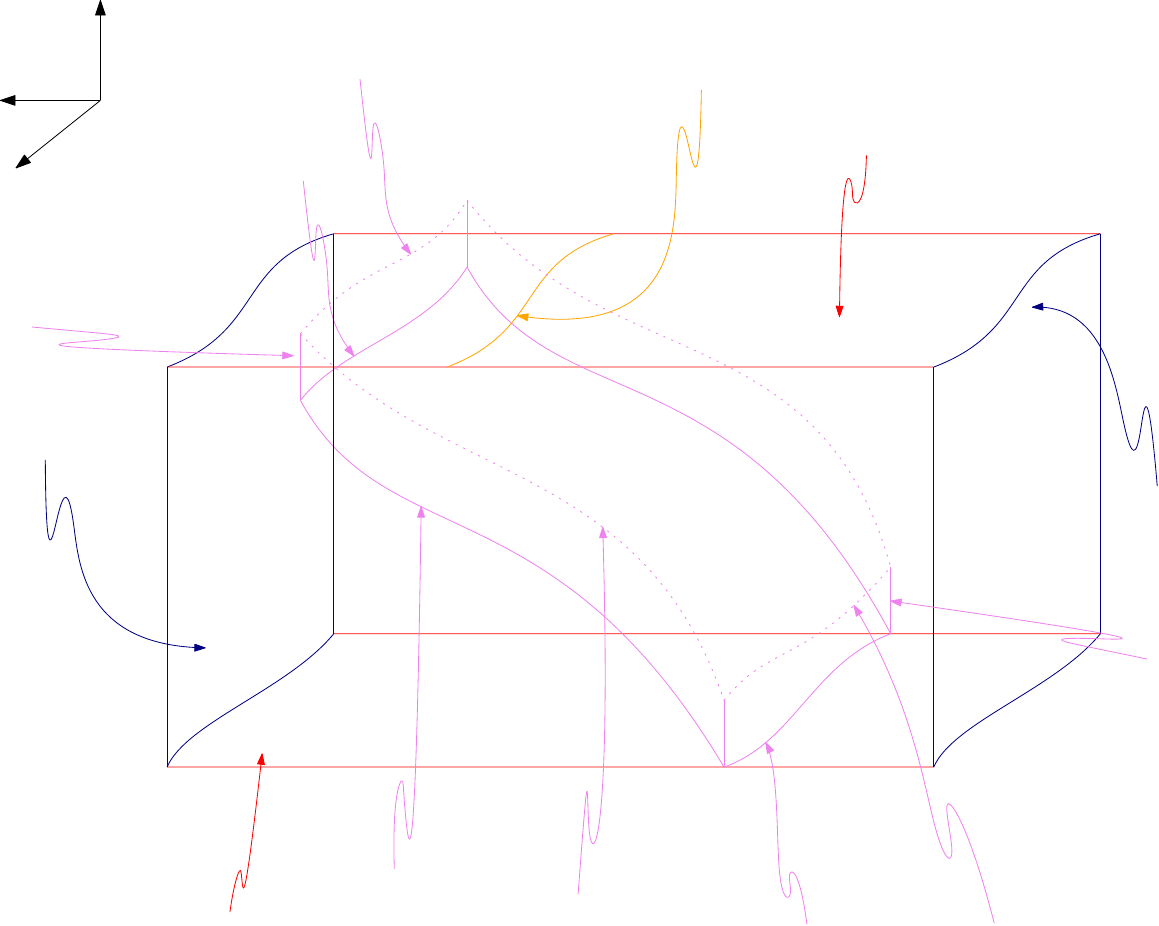}
			\put (55,73) {$\crease = \twoargmumuxtorus{0}{0}$}
			\put (68,66) {$\hypthreearg{0}{[-\interestingu,\interestingu]}{0}$}
			\put (95.5,34) {$\nullhypthreearg{0}{-\interestingu}{[\timefunction_0,0]}$}
			\put (-2.5,51.5) {$\outgoingcharacteristicsurfacetwoarg{\moreinterestingu_2}{[\leftubar,\ubarboot)}$}
			\put (99,22) {$\outgoingcharacteristicsurfacetwoarg{\moreinterestingu_1}{[\leftubar,\ubarboot)}$}
			\put (82,0) {$\doublenulltoritwoarg{\ubarboot}{\moreinterestingu_1}$}
			\put (66,0) {$\doublenulltoritwoarg{\leftubar}{\moreinterestingu_1}$}
			\put (46,-1) {$\ingoingcharacteristicsurfacetwoarg{\ubarboot}{[\moreinterestingu_1,\moreinterestingu_2]}$}
			\put (30,1.5) {$\ingoingcharacteristicsurfacetwoarg{\leftubar}{[\moreinterestingu_1,\moreinterestingu_2]}$}
			\put (22,66.5) {$\doublenulltoritwoarg{\leftubar}{\moreinterestingu_2}$}
			\put (27,75.5) {$\doublenulltoritwoarg{\ubarboot}{\moreinterestingu_2}$}
			\put (-1.5,42) {$\nullhypthreearg{0}{\interestingu}{[\timefunction_0,0]} $}
			\put (13,-2) {$\hypthreearg{\timefunction_0}{[-\interestingu,\interestingu]}{0}$}
			\put (4,65) {$(x^2,x^3) \in \mathbb{T}^2$}
			\put (9.5,78) {$t$}
			\put (-1,72) {$u \in \mathbb{R}$}
		\end{overpic}
		\vspace{0.5cm}
		\caption{The bootstrap region on which the solution is constructed. The region extends past  $\twoargMrough{[\timefunction_0,0],[- \rightu,\leftu]}{0}$ in order to access the Cauchy horizon.}
		\label{F:BOOTSTRAPDOMAINFORINGOINGEIKONALFUNCTION}
	\end{figure}
\end{center}

\subsection{Bootstrap assumptions tied to the fundamental scaffolding of the analysis}
\label{SS:BOOTSTRAPSCAFFOLDING}
The bootstrap assumptions in this section 
ensure that various fundamental aspects of our approach
(such as the change of variables maps from Sect.\,\ref{SS:ALLTHECHOVMAPS})
are well-defined and enjoy basic properties that
we use throughout the rest of the paper.

\subsubsection{Bootstrap assumptions for the inverse foliation density}
\label{SSS:BAFORINVERSEFOLIATIONDENSITY}

\begin{enumerate}
		\item  We assume that the following estimate holds on 
		$\characteristicdiamondtwoarg{[\leftubar,\ubarboot)}{[\moreinterestingu_1,\moreinterestingu_2]}$:
		\begin{align} \label{E:BAMUPOSITIVE}
			\upmu & > 0. \tag{\textbf{BA} $\upmu > 0$}
		\end{align}
	\item We assume that $\Lunit \upmu$ and $\geop{t} \upmu$ are quantitatively negative in 
	$\characteristicdiamondtwoarg{[\leftubar,\ubarboot)}{[\moreinterestingu_1,\moreinterestingu_2]}$,
	where $\blowupdeltadoublenull$ is defined in \eqref{E:DELTASTARDOUBLENULLDEF}:
	\begin{align} \label{E:BABOUNDSONLMUINTERESTINGREGION} \tag{\textbf{BA} $\Lunit \upmu$ neg}
	- 
	\frac{5}{4}
	\blowupdeltadoublenull
	& 
	\leq
	\min_{\characteristicdiamondtwoarg{[\leftubar,\ubarboot)}{[\moreinterestingu_1,\moreinterestingu_2]}} \Lunit \upmu
	\leq 
	\max_{\characteristicdiamondtwoarg{[\leftubar,\ubarboot)}{[\moreinterestingu_1,\moreinterestingu_2]}} \Lunit \upmu
	\leq 
	- 
	\frac{3}{4}
	\blowupdeltadoublenull,
		\\
	- 
	\frac{5}{4}
	\blowupdeltadoublenull
	& 
	\leq
	\min_{\characteristicdiamondtwoarg{[\leftubar,\ubarboot)}{[\moreinterestingu_1,\moreinterestingu_2]}}\geop{t} \upmu
	\leq 
	\max_{\characteristicdiamondtwoarg{[\leftubar,\ubarboot)}{[\moreinterestingu_1,\moreinterestingu_2]}} \geop{t} \upmu
	\leq 
	- 
	\frac{3}{4}
	\blowupdeltadoublenull.
		\label{E:BABOUNDSONGEOMETRICTDERIVATIVEMUINTERESTINGREGION}
		\tag{\textbf{BA} $\geop{t} \upmu$ neg}
\end{align}
		\item We assume that
		$\upmu$ is quantitatively convex in directions transversal to 
		$\datahypfortimefunctiontwoarg{0}{[\leftubar,0]}$. That is, we assume that
		the following estimates hold on $\characteristicdiamondtwoarg{[\leftubar,\ubarboot)}{[\moreinterestingu_1,\moreinterestingu_2]}$,
		where $0 < \secondtransversalderivativemulowerbound < 1$ is as in
		\eqref{E:DATATASSUMPTIONMUTRANSVERSALCONVEXITY}:
		\begin{align} \label{E:BAMUTRANSVERSALCONVEXITY} \tag{\textbf{BA} $\upmu$ cnvx}
			\begin{split}
			\frac{\secondtransversalderivativemulowerbound}{6}
			 \leq 
			\inf_{\characteristicdiamondtwoarg{[\leftubar,\ubarboot)}{[\moreinterestingu_1,\moreinterestingu_2]}}
				\Big\lbrace
				\newuL \newuL \upmu, \, \newuL \muX \upmu,  \,
				\muX \muX \upmu \Big\rbrace
			& \leq 
			\sup_{\characteristicdiamondtwoarg{[\leftubar,\ubarboot)}{[\moreinterestingu_1,\moreinterestingu_2]}}
			\Big\lbrace
				\newuL \newuL \upmu, \, \newuL \muX \upmu, \, \muX \muX \upmu
				\Big\rbrace
				\leq 
			\frac{6}{\secondtransversalderivativemulowerbound},
			\end{split} \\
			\frac{\secondtransversalderivativemulowerbound}{6} \le  \essinf_{\ingoingcharacteristicsurfacetwoarg{\leftubar}{[\moreinterestingu_1,\moreinterestingu_2]}} \nullgeop{u}\newuL \upmu & \le   \esssup_{\ingoingcharacteristicsurfacetwoarg{\leftubar}{[\moreinterestingu_1,\moreinterestingu_2]}}  \nullgeop{u}\newuL \upmu \le \frac{6}{\secondtransversalderivativemulowerbound}. \tag{\textbf{BA} $\upmu$ less reg cnvx}  \label{E:BAMULESSREGULARTRANSVERSALCONVEXITY}
		\end{align}	

\end{enumerate}

\subsubsection{Bootstrap assumptions for the ingoing eikonal function}
\label{SSS:BOOTSTRAPINGOINGEIKONALFUNCTION}
\begin{enumerate}
	\item We assume that
		$\geop t \ubar > 0$ in $\characteristicdiamondtwoarg{[\leftubar,\ubarboot)}{[\moreinterestingu_1,\moreinterestingu_2]}$.
	\item We assume that the following estimates hold, 
		where $\mathring{\updelta}_*$ is defined in \eqref{E:DELTASTARDOUBLENULLDEF}:
		\begin{align} \label{E:BALDERIVATIVEOFUBARAPPROXIMATELYUNITY} 
			\frac{3}{4} \blowupdeltadoublenull	
			&
			\leq
			\Lunit \ubar
			\leq 
			\frac{5}{4} \blowupdeltadoublenull,
			&&
			\mbox{on } \characteristicdiamondtwoarg{[\leftubar,\ubarboot)}{[\moreinterestingu_1,\moreinterestingu_2]}. \tag{\textbf{BA} $\Lunit\text{\underline{$u$}}$ pos}
		\end{align}
 \end{enumerate}
 
 \subsubsection{Bootstrap assumptions for the flow map of $\newuL$}
 \label{SSS:BOOTSTRAPNEWULFLOWMAP}
Let $\domainforembeddingdatahypfortimefunctiontwoarg{0}{[0,\mupositive]}$ be the domain for the embedding $\embeddingdatahypfortimefunctionarg{0}$ from Prop.\,\ref{P:PROPERTIESOFMUXMUZEROLEVELSET}, which is of the form $\embeddingdatahypfortimefunctionarg{0}(t,x^2,x^3) =	\left(t,\scalarembeddingdatahypfortimefunctionarg{0}(t,x^2,x^3),x^2,x^3 \right)$. Recall that $\flowmapofnewularg{\Delta s}$ is the flow map of $\newuL$, see \eqref{E:FLOWMAPOFNEWUL}. Let $\underline{F}$ and $\underline{\mathscr{F}}$ be the map and set defined by:
	\begin{subequations}
		\begin{align}
			\underline{F} (\Delta s,t,x^2,x^3) & \eqdef \flowmapofnewularg{\Delta s} \circ \embeddingdatahypfortimefunctionarg{0}(t, x^2,x^3), \label{E:BADEFINITIONOFFLOWOUTMAP} \\
			\begin{split} \label{E:BADEFINITIONOFFLOWOUTDOMAINFORUBAR}
				\underline{\mathscr{F}} & \eqdef \left\{ (\Delta s, t,x^2,x^3) \ | \ (t,x^2,x^3) \in \domainforembeddingdatahypfortimefunctiontwoarg{0}{(-\ubarboot,-\leftubar]}   \text{ and } \Delta s \in [\moreinterestingu_1 - \scalarembeddingdatahypfortimefunctionarg{0}(t,x^2,x^3), \moreinterestingu_2 - \scalarembeddingdatahypfortimefunctionarg{0}(t,x^2,x^3)] \right\}.
			\end{split} 
		\end{align}
	\end{subequations}
			
Then we assume that $\underline{F}$ is a $W^{1,\infty}$-bijection onto $\characteristicdiamondtwoarg{[\leftubar,\ubarboot)}{[\moreinterestingu_1,\moreinterestingu_2]}$ such that $\underline{F}$ and its inverse $\underline{F}^{-1}$ satisfy: 
	\begin{subequations}
		\begin{align}
			\left\| \underline{F} \right\|_{W^{1,\infty}(\underline{\mathscr{F}})} & < \infty, \label{E:BAFLOWMAPOFNEWULW1INFTYFINITE} \\
			\left\| \underline{F}^{-1} \right\|_{W^{1,\infty}_{\textnormal{geo}}\left( \characteristicdiamondtwoarg{[\leftubar,\ubarboot)}{[\moreinterestingu_1,\moreinterestingu_2]}\right)} & < \infty, \label{E:BAINVERSEFLOWMAPOFNEWULW1INFTYFINITE}
		\end{align}
	\end{subequations}
where we note the regularity measured by \eqref{E:BAFLOWMAPOFNEWULW1INFTYFINITE} is with respect to $(\Delta s,t,x^2,x^3)$-coordinate space.

\subsubsection{Bootstrap assumptions for change of variables maps}
\label{SSS:BACHOVMAPS}
 
We make the following assumptions for various change of variables maps.

\begin{enumerate}
	\item The change of variables map 
			$\Upsilon(t,u,x^2,x^3) = (t,x^1,x^2,x^3)$
			defined in \eqref{E:CHOVGEOTOCARTESIAN}
			satisfies
			$\left\| 
				\Upsilon
			\right\|_{C^{3,1}_{\textnormal{geo}}\left(\characteristicdiamondtwoarg{[\leftubar,\ubarboot)}{[\moreinterestingu_1,\moreinterestingu_2]} \right)} 
			< \infty
			$
			and is a diffeomorphism 
			from $\characteristicdiamondtwoarg{[\leftubar,\ubarboot)}{[\moreinterestingu_1,\moreinterestingu_2]}$
			onto its image.
	\item The change of variables map $\CHOVgeotodoublenull(t,u,x^2,x^3)  = (\ubar,u,x^2,x^3)$
			defined in \eqref{E:CHOVFROMGEOTODOUBLENULLCOORDINATES}
			satisfies
			$\left\| 
				\CHOVgeotodoublenull
			\right\|_{C^{1,1}_{\textnormal{geo}}\left(\characteristicdiamondtwoarg{[\leftubar,\ubarboot)}{[\moreinterestingu_1,\moreinterestingu_2]} \right)} 
			< \infty
			$
			and is a diffeomorphism from 
			$\characteristicdiamondtwoarg{[\leftubar,\ubarboot)}{[\moreinterestingu_1,\moreinterestingu_2]}$
			onto its image, which is $[\leftubar,\ubarboot) \times [\moreinterestingu_1,\moreinterestingu_2] \times \mathbb{T}^2$.
	\item The change of variables map $\CHOVdoublenulltoubarnewulmu$ defined in \eqref{E:CHOVFROMDOUBLENULLTOUBARNEWULMUCOORDINATES} satisfies $\left\| \CHOVdoublenulltoubarnewulmu \right\|_{C^{0,1}\left( [\leftubar,\ubarboot)\times[\moreinterestingu_1,\moreinterestingu_2]\times\T^2\right)} < \infty$ and is a Lipeomorphism from $[\leftubar,\ubarboot)\times[\moreinterestingu_1,\moreinterestingu_2]\times\T^2$ onto its image. 
	\end{enumerate}

\subsubsection{Bootstrap assumptions for the structure and location of $\argubarnewultorus{\ubar}$}
\label{SSS:BOOTSTRAPASSUMPTIONFORTORISTRUCTURE}
Recall that the $\newuL\upmu$-adapted tori $\argubarnewultorus{\ubar}$ are defined in \eqref{E:DEFOFNEWULADAPTEDTORI}. 
	\begin{enumerate}
		\item We assume that for each $\ubar \in [\leftubar,\ubarboot)$ there exist functions $\Cartesiantisafunctiononubarnewultori_{\ubar},\, \Eikonalisafunctiononubarnewultori_{\ubar} \in W^{1,\infty}(\T^2)$, depending on $\ubar$, such that relative to the geometric coordinates $(t,u,x^2,x^3)$, we have:
			\begin{align}
				\argubarnewultorus{\ubar} = \left\{ \left( \Cartesiantisafunctiononubarnewultori_{\ubar}(x^2,x^3), \Eikonalisafunctiononubarnewultori_{\ubar}(x^2,x^3),x^2,x^3\right) \ | \ (x^2,x^3) \in \T^2\right\} \label{E:BAUBARNEWULMUTORISTRUCTURE} \tag{\textbf{BA} $\text{\underline{$\newL$}}\upmu$ -- \textbf{TORI STRUCTURE}}
			\end{align}
			In particular, in geometric coordinates, $\argubarnewultorus{\ubar}$ is a $W^{1,\infty}$ graph over $\T^2$.
		\item We assume that for each $\ubar \in [\leftubar,\ubarboot)$, we have:
			\begin{align}
				\argubarnewultorus{\ubar} \subset \ingoingcharacteristicsurfacetwoarg{\ubar}{[-\frac{1}{2}\interestingu, \frac{3}{4}\moreinterestingu_2 + \frac{1}{2} \interestingu]}. \label{E:BAUBARNEWULMUTORILOCATION} \tag{\textbf{BA} $\text{\underline{$\newL$}}\upmu$ -- \textbf{TORI LOCATION}}
			\end{align}
		\item We assume that the map $\embeddednewuLhypersurface \colon [\leftubar,\ubarboot) \times \T^2 \to \characteristicdiamondtwoarg{[\leftubar,\ubarboot)}{[\moreinterestingu_1,\moreinterestingu_2]}$ defined by: 
			\begin{align}
				\embeddednewuLhypersurface(\ubar,x^2,x^3) =  \left( \Cartesiantisafunctiononubarnewultori_{\ubar}(x^2,x^3), \Eikonalisafunctiononubarnewultori_{\ubar}(x^2,x^3),x^2,x^3\right) \label{E:EMBEDDINGNEWULMUHYPERSURFACE}
			\end{align}
			is a Lipeomorphism from $[\leftubar,\ubarboot)\times\T^2$ onto $\newuLmulevelset^{[\leftubar,\ubarboot)}$ such that for every $\ubar' \in (\leftubar,\ubarboot)$, 
			we have:
	\begin{align} \label{E:EMBEDDEDEDNEWULMUHYPERSURFACEFINITEC01NORMONCOMPACTSUBSETS}
		\embeddednewuLhypersurface 
		& \in C^{0,1}([\leftubar,\ubar'] \times \mathbb{T}^2),
	\end{align}
	and such that, for any $(x^2,x^3) \in \T^2$, the map $\ubar \mapsto \Cartesiantisafunctiononubarnewultori_{\ubar}(x^2,x^3)$ is non-decreasing:
	\begin{align} \label{E:BANEWULMUHYPERSURFACEEMBEDDINGINCREASING}
		\Cartesiantisafunctiononubarnewultori_{\ubar_1}(x^2,x^3) \le \Cartesiantisafunctiononubarnewultori_{\ubar_2}(x^2,x^3), & &  \ubar_1 < \ubar_2. & & \ubar_1,\ubar_2 \in [\leftubar,\ubarboot).
	\end{align}	\end{enumerate}

\subsubsection{Soft bootstrap assumptions concerning regularity}
\label{SSS:SOFTBACONCERNINGREGULARITY}
We assume that for every $\ubar \in (\leftubar,\ubarboot)$, we have:
\begin{subequations}
\begin{align}
		\wavearray, 
			& \in C_{\textnormal{geo}}^{3,1}\left(\characteristicdiamondtwoarg{[\leftubar,\ubar]}{[\moreinterestingu_1,\moreinterestingu_2]}\right),
				\tag{\textbf{BA} \textbf{Fluid Regularity}} 
					\label{E:FLUIDC31ONCOMPACTSUBSETS} \\
		\Upsilon & \in C_{\textnormal{geo}}^{3,1}\left(\characteristicdiamondtwoarg{[\leftubar,\ubar]}{[\moreinterestingu_1,\moreinterestingu_2]}\right),	
			\tag{\textbf{BA} $\Upsilon$-\textbf{Regularity}} 
				\label{E:GEOTOCARTESIANCHOVC31ONCOMPACTSUBSETS} \\
		\Lunit^i, \, \upmu & \in C_{\textnormal{geo}}^{2,1}\left(\characteristicdiamondtwoarg{[\leftubar,\ubar]}{[\moreinterestingu_1,\moreinterestingu_2]}\right).
			 \tag{\textbf{BA} \textbf{Geometry Regularity}}  
			\label{E:GEOMETRYC31ONCOMPACTSUBSETS}
\end{align}
\end{subequations}\begin{remark}
	In \eqref{E:FLUIDC31ONCOMPACTSUBSETS}--\eqref{E:GEOMETRYC31ONCOMPACTSUBSETS},
	we are not making any quantitative assumptions on the size of the norms.
	That is, we are assuming only that the norms are finite,
	e.g., $\| \upmu \|_{C_{\textnormal{geo}}^{2,1}\left(\characteristicdiamondtwoarg{[\leftubar,\ubar]}{[\moreinterestingu_1,\moreinterestingu_2]}\right)} < \infty$.
	In Lemma\,\ref{L:CONTINUOUSEXTNESION}, we will show that all of these norms are $\leq C$.
\end{remark}

\subsection{The main quantitative bootstrap assumptions} 
\label{SS:MAINQUANTITATIVEBOOTSTRAPASSUMPTIONS}
We now state our main quantitative bootstrap assumptions.
In the rest of the paper, 
$\fundbootsmall \geq 0$ denotes a small ``bootstrap'' parameter whose smallness we described in Sect.\,\ref{SS:PARAMETERSIZEASSUMPTIONS}.
Later on, we will close our bootstrap argument by setting 
$\fundbootsmall = C \initialsmalldoublenull$ for some large constant $C$, 
where $\initialsmalldoublenull$ is the data-size parameter from Sect.\,\ref{SSS:QUANTITATIVEASSUMPTIONSONDATAAWAYFROMSYMMETRY};
see, in particular, Prop.\,\ref{P:IMPROVEMENTOFFUNDAMENTALQUANTITATIVEBOOTSTRAPASSUMPTIONS}.

\subsubsection{Fundamental quantitative bootstrap assumptions}
\label{SSS:FUNDAMENTALQUANTITATIVE}
Our fundamental quantitative bootstrap assumptions for 
$\wavearray$ are that the following inequalities hold for: 
$(\ubar,u) \in [\leftubar,\ubarboot) \times [\moreinterestingu_1,\moreinterestingu_2]$: 
\begin{align}
\left\| 
	\tander^{[1,\Ntop-6]} \wavearray
\right\|_{L^{\infty}\left(\doublenulltoritwoarg{\ubar}{u}\right)}, \,   \left\| 
	\tander^{[1,\Ntop-6]} (\vortrenormalized,\GradEnt)
\right\|_{L^{\infty}\left(\doublenulltoritwoarg{\ubar}{u}\right)},  \,   \left\| 
	\tander^{[1,\Ntop-7]} (\VortVort,\DivGradEnt)
\right\|_{L^{\infty}\left(\doublenulltoritwoarg{\ubar}{u}\right)}, \,  & \leq \fundbootsmall. \tag{\textbf{BA} $L^{\infty}$ \textbf{FUND}}
\label{E:FUNDAMENTALQUANTITATIVEBOOT} 
\end{align}

\subsubsection{Auxiliary bootstrap assumptions} \label{SSS:AUXBOOTSTRAP}
To derive pointwise and $L^{\infty}$ estimates, we find it convenient to make the following auxiliary bootstrap assumptions. 

\medskip

\noindent \underline{\textbf{Auxiliary bootstrap assumptions for small quantities}}.
 We assume that the following inequalities hold for 
$(\ubar,u) \in [\leftubar,\ubarboot) \times [\moreinterestingu_1,\moreinterestingu_2]$
(recall that $\wavearray$ and $\wavearraypartial$ are defined in Def.\,\ref{D:ARRAYSOFWAVEVARIABLES}
and that $\newL$ is defined in \eqref{E:DOUBLENULLEIKONALFUNCTIONNORMALIZEDNULLVECTORFIELDS}):
\begin{align} \label{BA:AUXR+} \tag{\textbf{AUX $\RRiemann$ SMALL}} 
	\left\|\RRiemann \right\|_{L^{\infty}\left(\doublenulltoritwoarg{\ubar}{u}\right)} 
	& \leq \mathring{\upalpha}^{1/2} + \auxbootsmall,  
		\\
	\left\| \wavearraypartial \right\|_{L^{\infty}\left(\doublenulltoritwoarg{\ubar}{u}\right)} 
	& \leq \auxbootsmall,
		 \label{BA:AUXWAVEARRAYPARTIALLINFINITY} 
			\tag{\textbf{AUX $\wavearraypartial$ SMALL}} 
			\\
	 \label{BA:AUXWAVEARRAY} \tag{\textbf{AUX $\wavearray$ SMALL}} 
	\begin{split}  
	\left\|\newL\comder^{\leq \Ntop-7;1} \wavearray \right\|_{L^{\infty}\left(\doublenulltoritwoarg{\ubar}{u}\right)},
		\,
	\left\|\comdersmall^{[1,\Ntop-7];1} \wavearray \right\|_{L^{\infty}\left(\doublenulltoritwoarg{\ubar}{u}\right)},
			  \\
	\left\| \newL \tander^{\leq 4} \muX \muX \wavearray \right\|_{L^{\infty}\left(\doublenulltoritwoarg{\ubar}{u}\right)},
		\,
	\left\| \comdersmall^{[1,6];2} \wavearray \right\|_{L^{\infty}\left(\doublenulltoritwoarg{\ubar}{u}\right)},
		&   \\
	\left\| \newL \tander^{\leq 2} \muX \muX \muX \wavearray \right\|_{L^{\infty}\left(\doublenulltoritwoarg{\ubar}{u}\right)},
		\,
	\left\| \comdersmall^{[1,5];3} \wavearray \right\|_{L^{\infty}\left(\doublenulltoritwoarg{\ubar}{u}\right)},
	&  \\
	\left\| \newL\muX \muX \muX \muX \wavearray \right\|_{L^{\infty}\left(\doublenulltoritwoarg{\ubar}{u}\right)}
	& \leq 
	\auxbootsmall,
	\end{split} 
		\\
	\label{BA:AUXMUSMALL} \tag{\textbf{AUX $\upmu$ SMALL}} 
	\begin{split}
	\left\| \newL\tander^{[1,\Ntop-8]} \upmu \right\|_{L^{\infty}\left(\doublenulltoritwoarg{\ubar}{u}\right)}, 
		\,
	\left\| \tander_*^{[1,\Ntop-8]} \upmu \right\|_{L^{\infty}\left(\doublenulltoritwoarg{\ubar}{u}\right)},
	& 
	 \\
	\left\| \newL\comdersmall^{[1,5];1}  \upmu \right\|_{L^{\infty}\left(\doublenulltoritwoarg{\ubar}{u}\right)}, 
		\,
	\left\| \comderdoublesmall^{[1,5];1} \upmu \right\|_{L^{\infty}\left(\doublenulltoritwoarg{\ubar}{u}\right)},
	&	
		\\
	\left\| \newL\comdersmall^{[1,4];2}  \upmu \right\|_{L^{\infty}\left(\doublenulltoritwoarg{\ubar}{u}\right)}, 
		\,
	\left\| \comderdoublesmall^{[1,4];2} \upmu \right\|_{L^{\infty}\left(\doublenulltoritwoarg{\ubar}{u}\right)}
	& \leq \auxbootsmall,
	\end{split}	
		\\
	\left\| \Lsmall^1 \right\|_{L^{\infty}\left(\doublenulltoritwoarg{\ubar}{u}\right)} 
	& 
	\leq \mathring{\upalpha}^{1/2},
		\label{BA:AUXL1SMALL} \tag{\textbf{AUX $\Lsmall^1$ SMALL}}
		\\
	\left\|\Lsmall^A \right\|_{L^{\infty}\left(\doublenulltoritwoarg{\ubar}{u}\right)}
	\label{BA:AUXLASMALL} \tag{\textbf{AUX $\Lsmall^A$ SMALL}} 
	& 
	\leq \auxbootsmall
				\\
	\label{BA:AUXTANGENTIALDERIVATIVESLISMALL} \tag{\textbf{AUX $\tander \Lsmall^i$ SMALL}} 
	\begin{split} 
	\left\| \newL\tander^{\leq \Ntop-7} \Lsmall^i \right\|_{L^{\infty}\left(\doublenulltoritwoarg{\ubar}{u}\right)}, 
			\,
	\left\| \tander^{[1,\Ntop-7]} \Lsmall^i \right\|_{L^{\infty}\left(\doublenulltoritwoarg{\ubar}{u}\right)}, 
			\\
	\left\| \newL\comder^{[1,\Ntop-8];1} \Lsmall^i \right\|_{L^{\infty}\left(\doublenulltoritwoarg{\ubar}{u}\right)}, 
		\,
	\left\| \comdersmall^{[1,\Ntop-8];1} \Lsmall^i \right\|_{L^{\infty}\left(\doublenulltoritwoarg{\ubar}{u}\right)}, 
			\\
	\left\| \newL\comder^{[1,5];2}\Lsmall^i \right\|_{L^{\infty}\left(\doublenulltoritwoarg{\ubar}{u}\right)}, 
		\,
	\left\| \comdersmall^{[1,5];2} \Lsmall^i \right\|_{L^{\infty}\left(\doublenulltoritwoarg{\ubar}{u}\right)},
		\\\
	\left\| \newL\comder^{[1,4];3}\Lsmall^i \right\|_{L^{\infty}\left(\doublenulltoritwoarg{\ubar}{u}\right)}, 
		\,
	\left\| \comdersmall^{[1,4];3} \Lsmall^i \right\|_{L^{\infty}\left(\doublenulltoritwoarg{\ubar}{u}\right)}
	& 
	\leq \auxbootsmall.
\end{split}
\end{align}\medskip

\noindent \underline{\textbf{Auxiliary bootstrap assumptions tied to pure transversal derivatives}}. 
We assume that the following inequalities hold for 
$(\ubar,u) \in [\leftubar,\ubarboot) \times [\moreinterestingu_1,\moreinterestingu_2]$:
\begin{align}
	\left 
		\| \muX^M \RRiemann
	\right\|_{L^{\infty}\left(\doublenulltoritwoarg{\ubar}{u}\right)} 
	& 
	\leq 
	\mathring{\updelta}
	+ 
	\auxbootsmall, 
	\label{BA:AUXTRANSVERSALPDERIVATIVESRRIEMANNLARGE} \tag{\textbf{AUX $\muX^M \RRiemann$ LARGE}} 
	&
	1 \leq M \leq 4,
		\\
	\left 
		\| \muX^M \wavearraypartial 
	\right\|_{L^{\infty}\left(\doublenulltoritwoarg{\ubar}{u}\right)} 
	& 
	\leq 
	\auxbootsmall, \label{BA:AUXTRANSVERSALPDERIVATIVESPARTIALWAVEARRAYSMALL}   
	\tag{\textbf{AUX $\muX^M \wavearraypartial$ SMALL}} 
	&
	1 \leq M \leq 4,
		\\
	\left\| 
		\muX^M \Lsmall^1 
	\right\|_{L^{\infty}\left(\doublenulltoritwoarg{\ubar}{u}\right)} 
	& 
	\leq	 
	\mathring{\updelta}
	+ 
	\auxbootsmall, 
		\label{BA:AUXTRANSVERSALPDERIVATIVESL1LARGE}
		\tag{\textbf{AUX $\muX^{[1,3]}\Lsmall^1$ LARGE}} 
	&
	1 \leq M \leq 3,
		\\
	\left\| 
		\muX^M \Lsmall^A 
	\right\|_{L^{\infty}\left(\doublenulltoritwoarg{\ubar}{u}\right)} 
	& \leq 
	\auxbootsmall, 
	\label{BA:AUXTRANSVERSALPDERIVATIVESLASMALL} \tag{\textbf{AUX $\muX^{[1,3]}\Lsmall^A$ SMALL}} 
	&
	1 \leq M \leq 3,
		\\
	\left\| 
		\muX^M \upmu 
	\right\|_{L^{\infty}\left(\doublenulltoritwoarg{\ubar}{u}\right)} 
	& 
	\leq 
	1 +
	\frac{3}{2 \blowupdeltadoublenull}
	\left\|
				\muX^M
				\left\lbrace
				G_{\Lunit\Lunit}^0
				\muX \RRiemann 
				\right\rbrace
	\right\|_{L^{\infty}\left(\doublenulltoritwoarg{\leftubar}{u}\right)}
	+ 
	\auxbootsmall, 
		\label{BA:AUXMULARGE} \tag{\textbf{AUX $\muX^{\le 3}\upmu$ LARGE}} 
	&
	0 \leq M \leq 3,
		\\
	\left\|
		\newL \muX^M \upmu 
	\right\|_{L^{\infty}\left(\doublenulltoritwoarg{\ubar}{u}\right)}
	& 
	\leq 
		\frac{2}{3 \blowupdeltadoublenull} 
			\left\|
				\muX^M
				\left\lbrace
				G_{\Lunit\Lunit}^0
				\muX \RRiemann 
				\right\rbrace
			\right\|_{L^{\infty}\left(\doublenulltoritwoarg{\leftubar}{u}\right)}
	+ 
	\auxbootsmall, \label{BA:AUXLMULARGE} 
	\tag{\textbf{AUX $\newL \muX^{\le 3}\upmu$ LARGE}}
	&
	0 \leq M \leq 3.
\end{align}
\subsection{Improvements of the bootstrap assumptions}

\label{SS:SUMMARYOFIMPROVEMENTOFBOOTSTRAPASSUMPTIONS}
In the subsequent sections, we will derive strict improvements of all the bootstrap assumptions
that we made throughout Sect.\,\ref{S:BOOTSTRAPEVERYTHINGEXCEPTENERGIES}.
For the reader's convenience, here we state all the forthcoming results that yield the desired
strict improvements. Here we clarify that by ``strict improvements,'' we mean one or more of the following
three things:
\begin{enumerate}
	\item (\textbf{Quantitative improvement}) By this, we mean 
		that some quantity $Q$ 
		was assumed to satisfy $A_1 \leq Q \leq A_2$ in the bootstrap assumptions
		(where $A_1,A_2$ are real numbers),
		and we derive the improved bound $B_1 \leq Q \leq B_2$, 
		where $A_1 < B_1 \leq B_2 < A_2$.
		\item (\textbf{From soft to quantitative}) By this, we mean that in the bootstrap assumptions,
		we assumed that some function $f$ belongs to some function space and has a finite norm in that space,
		and our improvement is a quantitative estimate for the norm of $f$.
	\item (\textbf{Extension to the closure})
		By this, we mean that our bootstrap assumptions involved an assumption on the ``open-at-the-top'' domain
			$\characteristicdiamondtwoarg{[\leftubar,\ubarboot)}{[\moreinterestingu_1,\moreinterestingu_2]}$,
			and we derive an improved result showing that the assumption holds
			on the closed domain $\characteristicdiamondtwoarg{[\leftubar,\ubarboot]}{[\moreinterestingu_1,\moreinterestingu_2]}$.
\end{enumerate}

Here are the precise spots in the article where we derive improvements of the bootstrap assumptions.
\begin{itemize}
	\item Regarding the bootstrap assumptions of Sect.\,\ref{SSS:BAFORINVERSEFOLIATIONDENSITY}:
		we derive improvements of \eqref{E:KEYESTIMATECONTROLLINGINVERSEMUBYINVERSEWEIGHT} in \eqref{E:MINVALUEOFMUONFOLIATION},
		of \eqref{E:BABOUNDSONLMUINTERESTINGREGION} in \eqref{E:BOUNDSONLMUINTERESTINGREGION},
		of \eqref{E:BABOUNDSONGEOMETRICTDERIVATIVEMUINTERESTINGREGION} in \eqref{E:BOUNDSONGEOMETRICTDERIVATIVEMUINTERESTINGREGION}, 
		and of
		\eqref{E:BAMUTRANSVERSALCONVEXITY} in \eqref{E:MUTRANSVERSALCONVEXITY}.
	\item We derive improvements of the bootstrap assumptions of Sect.\,\ref{SSS:BOOTSTRAPINGOINGEIKONALFUNCTION}
		in Lemma~\ref{L:PROPERTIESANDDIFFEOMORPHICEXTENSIONOFDOUBLENULLCOORDINATES}.
	\item We derive improvements of the bootstrap assumptions of Sect.\,\ref{SSS:BOOTSTRAPNEWULFLOWMAP} in Lemma\,\ref{L:PROPERTIESOFFLOWMAPOFNEWUL}.
	\item We derive improvements of the bootstrap assumptions of Sect.\,\ref{SSS:BACHOVMAPS}
				in Lemmas~\ref{L:PROPERTIESANDDIFFEOMORPHICEXTENSIONOFDOUBLENULLCOORDINATES}
				and
				\ref{L:CHOVFROMDOUBLENULLTONEWULMUCOORDINATES}
				and in Prop.\,\ref{P:HOMEOMORPHICANDDIFFEOMORPHICEXTENSIONOFCARTESIANCOORDINATES}.	
	\item Regarding the bootstrap assumptions of Sect.\,\ref{SSS:BOOTSTRAPASSUMPTIONFORTORISTRUCTURE}:
		we derive improvements of \eqref{E:BAUBARNEWULMUTORISTRUCTURE}
		in \eqref{E:C01BOUNDFORCARTESIANTCOORDINATEGRAPHTORI}--\eqref{E:C01BOUNDFORUDOUBLENULLCOORDINATEGRAPHTORI}, 
		of \eqref{E:BAUBARNEWULMUTORILOCATION} in \eqref{E:NEWULMULEVELSETIMPROVELOCATION}--\eqref{E:IMPROVEDLEVELSETSTRUCTUREANDLOCATIONOFMINONNEWULMUTORI}, and of Item 3 of in Lemma\,\ref{L:CHOVFROMDOUBLENULLTONEWULMUCOORDINATES}.
	\item We derive improvements of the bootstrap assumptions of Sect.\,\ref{SSS:SOFTBACONCERNINGREGULARITY}
			in Lemma~\ref{L:CONTINUOUSEXTNESION}.		
	\item We derive improvements of the fundamental quantitative bootstrap assumptions of Sect.\,\ref{SSS:FUNDAMENTALQUANTITATIVE}
		in Prop.\,\ref{P:IMPROVEMENTOFFUNDAMENTALQUANTITATIVEBOOTSTRAPASSUMPTIONS}.
	\item We derive improvements of the auxiliary bootstrap assumptions of Sect.\,\ref{SSS:AUXBOOTSTRAP}
		in Prop.\,\ref{P:IMPROVEMENTOFAUXILIARYBOOTSTRAP}.
	\item In Sect.\,\ref{SS:BOOTSTRAPASSUMPTIONSFORTHEWAVEENERGIES}, 
		we state bootstrap assumptions for the 
		$L^2$-type energies for the wave variables $\wavearray$. 
		We derive improvements of them in Prop.\,\ref{P:APRIORIL2ESTIMATESWAVEVARIABLES}.
\end{itemize}


\section{Estimates for the ingoing eikonal function, continuous extensions, and various diffeomorphisms} \label{S:ESTIMATESFORUBARCONTINUOUSEXTENSIONSANDDIFFEOS}

We continue to work under the assumptions of Sect.\,\ref{SS:SILENTFACTS}.
In this section, we derive estimates for the ingoing eikonal function $\ubar$. We then show that the map $\CHOVgeotodoublenull(t,u,x^2,x^3)  = (\ubar,u,x^2,x^3)$ extends to a diffeomorphism on the closure of $\characteristicdiamondtwoarg{[\leftubar,\ubarboot)}{[\moreinterestingu_1,\moreinterestingu_2]}$. Next, we show that various variables extend to the compact set $\characteristicdiamondtwoarg{[\leftubar,\ubarboot]}{[\moreinterestingu_1,\moreinterestingu_2]}$ as functions with substantial H\"older regularity, which we will exploit throughout the paper. 

\begin{lemma}[Estimates for $\ubar$, $\CHOVgeotodoublenull$, and $\CHOVgeotodoublenull^{-1}$]
		\label{L:LINFTYESTIMATESFORINGOINGEIKONALFUNCTIONANDDERIVATIVES}
		The following estimates hold.
		
		\medskip
		
	\noindent \underline{\textbf{Estimates for $\ubar$}}:
		\begin{subequations}
		\begin{align}
			\left\| 
				\ubar
			\right\|_{W_{\textnormal{geo}}^{2,\infty}\left(\characteristicdiamondtwoarg{(\leftubar,\ubarboot)}{(\moreinterestingu_1,\moreinterestingu_2)}\right)}
		& 
		\lesssim 1,
			\label{E:ALLDERIVATIVESLINFTYESTIMATESFORUBARANDDERIVATIVES}
			\\
			\left\| 
				\left(\geop{x^2} \ubar,\geop{x^3} \ubar \right)  
			\right\|_{W_{\textnormal{geo}}^{1,\infty}\left(\characteristicdiamondtwoarg{(\leftubar,\ubarboot)}{(\moreinterestingu_1,\moreinterestingu_2)}\right)}
			& 
			\lesssim \auxbootsmall. 
			\label{E:SMALLDERIVATIVESLINFTYESTIMATESFORUBARANDDERIVATIVES}
		\end{align}
		\end{subequations}
	
		Moreover, 
		\begin{subequations}
		\begin{align} \label{E:PARTIALTIMEDERIVATIVEOFUBARAPPROXIMATELYUNITY}
				\frac{1}{2} \blowupdeltadoublenull	
				&
				\leq
				\geop{t} \ubar
				\leq 
				\frac{3}{2} \blowupdeltadoublenull,
				&&
				\mbox{on } \characteristicdiamondtwoarg{[\leftubar,\ubarboot)}{[\moreinterestingu_1,\moreinterestingu_2]}.
		\end{align}
		\end{subequations}\medskip

\noindent \underline{\textbf{Estimates for $\CHOVgeotodoublenull$ and $\CHOVgeotodoublenull^{-1}$}}:
		The following estimate holds,
		where $\CHOVgeotodoublenull$ is the change of variables map from geometric coordinates
		to double-null coordinates defined in \eqref{E:CHOVFROMGEOTODOUBLENULLCOORDINATES}:
		\begin{align} \label{E:W2INFINITYBOUNDFORCHOVFROMGEOTODOUBLENULLOPEN}
			\left\| 
				\CHOVgeotodoublenull
			\right\|_{W_{\textnormal{geo}}^{2,\infty}\left(\characteristicdiamondtwoarg{(\leftubar,\ubarboot)}{(\moreinterestingu_1,\moreinterestingu_2)}\right)}
		& 
		\lesssim 1.
		\end{align}

		The following estimates hold, where
		$\CHOVJacobiangeotodoublenull$ 
		is the Jacobian matrix of $\CHOVgeotodoublenull$ defined in \eqref{E:JACOBIANOFCHOVFROMGEOTODOUBLENULLCOORDINATES}:
		\begin{align} \label{E:KEYJACOBIANDETERMINANTESTIMATECHOVGEOTODOUBLENULL}
				\mydet \left(\CHOVJacobiangeotodoublenull \right)
				& \approx 1,
				&&
				\mbox{on } \characteristicdiamondtwoarg{[\leftubar,\ubarboot)}{[\moreinterestingu_1,\moreinterestingu_2]},
					\\
			\left\|
				[\CHOVJacobiangeotodoublenull]^{-1}
			\right\|_{W^{1,\infty} \left(\characteristicdiamondtwoarg{(\leftubar,\ubarboot)}{[\moreinterestingu_1,\moreinterestingu_2]} \right)}
			& \leq C.
			\label{E:INVERSEOFJACOBIANMATRIXFROMGEOTODOUBLENULLW1INFINITYBOUND}
		\end{align}
		
		In addition, the following estimates hold, 
		where $\CHOVgeotodoublenull^{-1}$ is the change of variables map from double-null coordinates to geometric coordinates
		and
		$d_{\textnormal{null}}\CHOVgeotodoublenull^{-1}$ 
		is its Jacobian matrix:
			\begin{align} \label{E:KEYJACOBIANDETERMINANTESTIMATECHOVDOUBLENULLTOGEO}
				\nullgeop{\ubar} t 
				& =
				\mydet \left(d_{\textnormal{null}} \CHOVgeotodoublenull^{-1} \right)
				\approx 1,
				&&
				\mbox{on } [\leftubar,\ubarboot) \times [\moreinterestingu_1,\moreinterestingu_2] \times \mathbb{T}^2.
		\end{align}
		
		In addition, the following estimates hold:
		\begin{align} \label{E:W2INFINITYBOUNDFORINVERSECHOVGEOTODOUBLENULL}
			\left\| 
				\CHOVgeotodoublenull^{-1}
			\right\|_{W_{\textnormal{null}}^{2,\infty}\left((\leftubar,\ubarboot) \times (\moreinterestingu_1,\moreinterestingu_2) \times \mathbb{T}^2\right)}
			\leq C.
		\end{align}
		
		Finally, $\CHOVgeotodoublenull^{-1}$ extends to 
		a $C^{1,1}([\leftubar,\ubarboot] \times [\moreinterestingu_1,\moreinterestingu_2] \times \mathbb{T}^2)$
		map satisfying the estimate:
		\begin{align} \label{E:C11BOUNDFORINVERSECHOVGEOTODOUBLENULL}
			\left\| 
				\CHOVgeotodoublenull^{-1}
			\right\|_{C_{\textnormal{null}}^{1,1}\left([\leftubar,\ubarboot] \times [\moreinterestingu_1,\moreinterestingu_2] \times \mathbb{T}^2\right)}
			\leq C.
		\end{align}
		
\end{lemma}

\begin{proof}
Throughout the proof, we silently use the fact that functions $f \in W_{\textnormal{geo}}^{1,\infty}(\characteristicdiamondtwoarg{(\leftubar,\ubarboot)}{(\moreinterestingu_1,\moreinterestingu_2)})$,
	are locally Lipschitz and thus have (by Rademacher's theorem) 
	a.e.\ differentiable locally Lipschitz traces along the $C^{1,1}$ hypersurface
	$\datahypfortimefunctiontwoarg{0}{[\leftubar,\ubarboot)}$ (see Prop.\,\ref{P:PROPERTIESOFMUXMUZEROLEVELSET}) such that 	
	\[ 	\| f \|_{W^{1,\infty}\left(\datahypfortimefunctiontwoarg{0}{(\leftubar,\ubarboot)}\right)} = 
	\| f \|_{W^{1,\infty}\left(\mbox{\upshape int}(\datahypfortimefunctiontwoarg{0}{[\leftubar,\ubarboot)})\right)}
	\lesssim
	\| f \|_{W_{\textnormal{geo}}^{1,\infty}\left(\characteristicdiamondtwoarg{(\leftubar,\ubarboot)}{(\moreinterestingu_1,\moreinterestingu_2)}\right)}.\]

We being with the bounds $\| \upmu \|_{W_{\textnormal{geo}}^{3,\infty}(\characteristicdiamondtwoarg{(\leftubar,\ubarboot)}{(\moreinterestingu_1,\moreinterestingu_2})} 
	\lesssim 1, \, \left\|
		\left(\geop{x^2} \upmu,\geop{x^3} \upmu \right)  
	\right\|_{W_{\textnormal{geo}}^{2,\infty}(\characteristicdiamondtwoarg{(\leftubar,\ubarboot)}{(\moreinterestingu_1,\moreinterestingu_2})}
	\lesssim
	\auxbootsmall
	$, which follow from  Lemma\,\ref{L:COMMUTATORSTOCOORDINATES}, 
	Prop.\,\ref{P:SCHEMATICSTRUCTUREOFVARIOUSTENSORSINTERMSOFCONTROLVARS},
	and the bootstrap assumptions. From \eqref{E:C01INFINITYESTIMATESFORLUBARONINITIALDATASURFACE}--\eqref{E:C01INFINITYESTIMATESFORANGULARUBARONINITIALDATASURFACE} and $[\leftubar,\ubarboot]\subset [\timefunction_0,0]$, it follows that 
	$\| \Lunit \ubar, \muX \ubar \|_{W_{\textnormal{geo}}^{1,\infty}\left(\datahypfortimefunctiontwoarg{0}{(\leftubar,\ubarboot)}\right)} 
	\lesssim 1, \, \left\|
		\left(\geop{x^2} \ubar,\geop{x^3} \ubar \right)  
	\right\|_{W_{\textnormal{geo}}^{1,\infty}\left(\datahypfortimefunctiontwoarg{0}{(\leftubar,\ubarboot)}\right)}
	\lesssim
	\auxbootsmall$. Using the bootstrap assumption \eqref{E:BALDERIVATIVEOFUBARAPPROXIMATELYUNITY}, we then have that \eqref{E:QUASILINEAREQUATIONSATISFIEDBYLUBAR}--\eqref{E:QUASILINEAREQUATIONSATISFIEDBYSMOOTHANGULARDERIVATIVEOFUBAR} form a quasilinear system of transport equations with $W^{1,\infty}(\characteristicdiamondtwoarg{(\leftubar,\ubarboot)}{(\moreinterestingu_1,\moreinterestingu_2)})$-coefficients and with $W_{\textnormal{geo}}^{1,\infty}\left(\datahypfortimefunctiontwoarg{0}{(\leftubar,\ubarboot)}\right)$-data. Using the fact that we are assuming, by virtue of the bootstrap assumptions, that an ingoing eikonal function $\ubar$ exists on $\characteristicdiamondtwoarg{[\leftubar,\ubarboot)}{[\moreinterestingu_1,\moreinterestingu_2]}$, techniques similar to the ones used in Theorem\,\ref{T:CONSTRUCTIONOFTHEINGOINGEIKONALFUNCTIONCONTRACTIONMAPPING}, estimates \eqref{E:ALLDERIVATIVESLINFTYESTIMATESFORUBARANDDERIVATIVES}--\eqref{E:SMALLDERIVATIVESLINFTYESTIMATESFORUBARANDDERIVATIVES} follow. Moreover, using the identity $\geop{t}\ubar = \Lunit \ubar - L^A \geop{x^A} \ubar$ and the bound $|L^A \geop{x^A} \ubar| \lesssim \auxbootsmall$, \eqref{E:PARTIALTIMEDERIVATIVEOFUBARAPPROXIMATELYUNITY} follows from \eqref{E:BALDERIVATIVEOFUBARAPPROXIMATELYUNITY}. 
	
	\eqref{E:W2INFINITYBOUNDFORCHOVFROMGEOTODOUBLENULLOPEN} follows from \eqref{E:ALLDERIVATIVESLINFTYESTIMATESFORUBARANDDERIVATIVES} and the definition of $\CHOVgeotodoublenull$. 
	\eqref{E:KEYJACOBIANDETERMINANTESTIMATECHOVGEOTODOUBLENULL} now follows from \eqref{E:PARTIALTIMEDERIVATIVEOFUBARAPPROXIMATELYUNITY} and the simple identity $\geop{t} \ubar = \mydet \left(\CHOVJacobiangeotodoublenull \right)$, which follows easily from the definition of $\CHOVgeotodoublenull$. 
	\eqref{E:INVERSEOFJACOBIANMATRIXFROMGEOTODOUBLENULLW1INFINITYBOUND} now follows from \eqref{E:W2INFINITYBOUNDFORCHOVFROMGEOTODOUBLENULLOPEN} and \eqref{E:KEYJACOBIANDETERMINANTESTIMATECHOVGEOTODOUBLENULL}. 
	The ``$=$'' in \eqref{E:KEYJACOBIANDETERMINANTESTIMATECHOVDOUBLENULLTOGEO} follows easily from the definition of $\CHOVgeotodoublenull$. The ``$\approx$'' in \eqref{E:KEYJACOBIANDETERMINANTESTIMATECHOVDOUBLENULLTOGEO} follows from the identity
	\[ \mydet \left(\left[\mathrm{d}_{\textnormal{null}} \CHOVgeotodoublenull^{-1}\right] \circ \CHOVgeotodoublenull \right) = \left( \mydet \left[\CHOVJacobiangeotodoublenull\right] \right)^{-1}.\]
	and \eqref{E:INVERSEOFJACOBIANMATRIXFROMGEOTODOUBLENULLW1INFINITYBOUND}.
	\eqref{E:W2INFINITYBOUNDFORINVERSECHOVGEOTODOUBLENULL} follows from differentiating once, with the double-null coordinate partial derivative vectorfields, the identity $\mathrm{d}_{\textnormal{null}} \CHOVgeotodoublenull^{-1} = [\CHOVJacobiangeotodoublenull]^{-1} \circ \CHOVgeotodoublenull^{-1}$, \eqref{E:INVERSEOFJACOBIANMATRIXFROMGEOTODOUBLENULLW1INFINITYBOUND}, and the chain rule.
	\eqref{E:C11BOUNDFORINVERSECHOVGEOTODOUBLENULL} follows from \eqref{E:W2INFINITYBOUNDFORINVERSECHOVGEOTODOUBLENULL}, and the following Sobolev embedding result for scalar functions $f$ (see the proof of \cite{lE1998}*{Theorem\,5 in Section~5.6}),
which relies on the convexity of the domain $(\leftubar,\ubarboot) \times (\moreinterestingu_1,\moreinterestingu_2) \times \mathbb{T}^2$: $\| f\|_{C^{0,1}_{\textnormal{null}}\left(  [\leftubar,\ubarboot] \times [\moreinterestingu_1,\moreinterestingu_2] \times \mathbb{T}^2 \right)} \le C \| f\|_{W^{1,\infty}_{\textnormal{null}}\left((\leftubar,\ubarboot) \times (\moreinterestingu_1,\moreinterestingu_2) \times \mathbb{T}^2\right)}$.
	
 \end{proof}

\begin{lemma}[Properties of $\CHOVgeotodoublenull$ on the closure of $\characteristicdiamondtwoarg{[\leftubar,\ubarboot)}{[\moreinterestingu_1,\moreinterestingu_2]}$ and quasi-convexity] \label{L:PROPERTIESANDDIFFEOMORPHICEXTENSIONOFDOUBLENULLCOORDINATES}

		The following results hold.
		
		\begin{enumerate}
			\item $\characteristicdiamondtwoarg{[\leftubar,\ubarboot]}{[\moreinterestingu_1,\moreinterestingu_2]}$ 
				is precompact in the topology of the geometric coordinates $(t,u,x^2,x^3)$.
			\item The change of variables map $\CHOVgeotodoublenull(t,u,x^2,x^3) = (\ubar,u,x^2,x^3)$
				extends to a $C^{1,1}$ diffeomorphism on the closure of
				$\characteristicdiamondtwoarg{[\leftubar,\ubarboot)}{[\moreinterestingu_1,\moreinterestingu_2]}$,
				which we denote by $\mbox{\upshape cl} \left(\characteristicdiamondtwoarg{[\leftubar,\ubarboot)}{[\moreinterestingu_1,\moreinterestingu_2]} \right)$.
				Moreover,
				 $\mbox{\upshape cl} \left(\characteristicdiamondtwoarg{[\leftubar,\ubarboot)}{[\moreinterestingu_1,\moreinterestingu_2]} \right)
				= \characteristicdiamondtwoarg{[\leftubar,\ubarboot]}{[\moreinterestingu_1,\moreinterestingu_2]}$,
				and
				$\CHOVgeotodoublenull\left(\characteristicdiamondtwoarg{[\leftubar,\ubarboot]}{[\moreinterestingu_1,\moreinterestingu_2]}\right) 
			= [\leftubar,\ubarboot] \times [\moreinterestingu_1,\moreinterestingu_2] \times \mathbb{T}^2$.
		\item The following estimates hold for the extended maps:
		\begin{align}
				\frac{1}{2} \blowupdeltadoublenull	
				&
				\leq
				\geop{t} \ubar
				\leq 
				\frac{3}{2} \blowupdeltadoublenull,
				&&
				\mbox{on } \characteristicdiamondtwoarg{[\leftubar,\ubarboot]}{[\moreinterestingu_1,\moreinterestingu_2]}.  \label{E:CLOSEDVERSIONKEYJACOBIANDETERMINANTESTIMATECHOVGEOTODOUBLENULL} \\
				\frac{3}{4} \blowupdeltadoublenull	
				&
				\leq
				\Lunit \ubar
				\leq 
				\frac{5}{4} \blowupdeltadoublenull,
				&&
				\mbox{on } \characteristicdiamondtwoarg{[\leftubar,\ubarboot]}{[\moreinterestingu_1,\moreinterestingu_2]}.  \label{E:CLOSEDVERSIONLUBARESTIMATECHOVGEOTODOUBLENULL} 	
		\end{align}
		\item 	\begin{align} \label{E:CLOSEDKEYJACOBIANDETERMINANTESTIMATECHOVDOUBLENULLTOGEO}
				\nullgeop{\ubar} t 
				& =
				\mydet \left(d_{\textnormal{null}} \CHOVgeotodoublenull \right)
				\approx 1,
				&&
				\mbox{on } [\leftubar,\ubarboot] \times [\moreinterestingu_1,\moreinterestingu_2] \times \mathbb{T}^2,
		\end{align}
		\item \begin{align}
							\left\| 
								\CHOVgeotodoublenull
							\right\|_{C_{\textnormal{geo}}^{1,1}(\characteristicdiamondtwoarg{[\leftubar,\ubarboot]}{[\moreinterestingu_1,\moreinterestingu_2]})}
							& \leq C,
								\label{E:CLOSEDVERSIONC11BOUNDFORCHOVDOUBLENULLTOGEO} 
									\\
							\left\| 
								\CHOVgeotodoublenull^{-1}
							\right\|_{C_{\textnormal{null}}^{1,1}([\leftubar,\ubarboot] \times [\moreinterestingu_1,\moreinterestingu_2] \times \mathbb{T}^2)}
							& \leq C,
							\label{E:CLOSEDVERSIONC11BOUNDFORINVERSECHOVGEOTODOUBLENULL}
								\\
							\left\| 
								\left(\geop{x^2} \ubar,\geop{x^3}\ubar \right)
							\right\|_{C_{\textnormal{geo}}^{0,1}\characteristicdiamondtwoarg{[\leftubar,\ubarboot]}{[\moreinterestingu_1,\moreinterestingu_2]})}
							& \leq C \auxbootsmall. 
							\label{E:SMALLC01ESTIMATESFORUBAR}
						\end{align}
				\item (\textbf{Quasi-convexity of} 
				$\characteristicdiamondtwoarg{[\leftubar,\ubarboot]}{[\moreinterestingu_1,\moreinterestingu_2]}$). 
					For every pair of points 
					$q_1,q_2 \in [\leftubar,\ubarboot] \times [\moreinterestingu_1,\moreinterestingu_2] \times \mathbb{T}^2$,
					we have:
					\begin{align} \label{E:COMPARABLEDISTANCES}
						\textnormal{dist}_{\textnormal{flat}}
						\left(\CHOVgeotodoublenull^{-1}(q_1),\CHOVgeotodoublenull^{-1}(q_2) \right)
						& \approx
						\textnormal{dist}_{\textnormal{flat}}(q_1,q_2),
					\end{align}
					where $\textnormal{dist}_{\textnormal{flat}}(q_1,q_2)$ is the
					standard Euclidean distance between $q_1$ and $q_2$ in the flat space 
					$\mathbb{R}_{\ubar} \times \mathbb{R}_u \times \mathbb{T}^2$.
					
					Moreover, $\characteristicdiamondtwoarg{[\leftubar,\ubarboot]}{[\moreinterestingu_1,\moreinterestingu_2]}$ is quasi-convex.\footnote{Here, 
					we are not just interested in a qualitative version of quasi-convexity, but rather in obtaining control over the constants ``$C$.''
					Similar remarks apply for the Sobolev embedding result \eqref{E:SOBOELVEMBEDDINGRELYINGONQUASICONVEXITY} and for other 
					quasi-convexity and Sobolev embedding results derived throughout the paper. \label{FN:QUANTITATIVEQUASICONVEXITY}}
					That is,
					every pair of points
					$p_1,p_2 \in \characteristicdiamondtwoarg{[\leftubar,\ubarboot]}{[\moreinterestingu_1,\moreinterestingu_2]}$
					is connected by a $C_{\textnormal{geo}}^1$ curve in $\characteristicdiamondtwoarg{[\leftubar,\ubarboot]}{[\moreinterestingu_1,\moreinterestingu_2]}$
					whose length with respect to the standard flat Euclidean metric on geometric coordinate space
					$\mathbb{R}_t \times \mathbb{R}_u \times \mathbb{T}^2$
					is $\leq C \textnormal{dist}_{\textnormal{flat}}(p_1,p_2)$.
				\item (\textbf{Sobolev embedding}).
					There is a constant $C > 0$ that is \textbf{independent of $\timefunctionboot$}
					such that the following Sobolev embedding result holds for scalar functions $f$ on 
					$\characteristicdiamondtwoarg{(\leftubar,\ubarboot)}{(\moreinterestingu_1,\moreinterestingu_2)}$:
					\begin{align} \label{E:SOBOELVEMBEDDINGRELYINGONQUASICONVEXITY}
						\| f \|_{C_{\textnormal{geo}}^{0,1}(\characteristicdiamondtwoarg{[\leftubar,\ubarboot]}{[\moreinterestingu_1,\moreinterestingu_2]})}
						& \leq
						C
						\| f \|_{W_{\textnormal{geo}}^{1,\infty}(\characteristicdiamondtwoarg{(\leftubar,\ubarboot)}{(\moreinterestingu_1,\moreinterestingu_2)})}.
					\end{align}
				\item ($\ingoingcharacteristicsurfacetwoarg{\ubar}{[\moreinterestingu_1,\moreinterestingu_2]}$ \textbf{is a graph}).
					For $\ubar \in [\leftubar,\ubarboot]$,
					there exists a function 
					$\Cartesiantisafunctiononlevelsetsofingoingcharacteristicsarg{\ubar}: 
					[\moreinterestingu_1,\moreinterestingu_2] \times \mathbb{T}^2 \rightarrow \mathbb{R}$,
					depending on $\ubar$,
					such that:
					\begin{align} \label{E:C11BOUNDFORCONSTANTTIMEFUNCTIONGRAPH}
						\left\| \Cartesiantisafunctiononlevelsetsofingoingcharacteristicsarg{\ubar} \right\|_{C^{1,1}([\moreinterestingu_1,\moreinterestingu_2] \times \mathbb{T}^2)}  
						& \leq C
					\end{align}
					and such that relative to the geometric coordinates, we have:
					\begin{align} \label{E:LEVELSETSOFUBARAREAGRAPH}
						\ingoingcharacteristicsurfacetwoarg{\ubar}{[\moreinterestingu_1,\moreinterestingu_2]}
							& =
						\left\lbrace
							(t,u,x^2,x^3)
							\ | \
							t = \Cartesiantisafunctiononlevelsetsofingoingcharacteristicsarg{\ubar}(u,x^2,x^3),
								\,
							(u,x^2,x^3) \in [\moreinterestingu_1,\moreinterestingu_2] \times \mathbb{T}^2
						\right\rbrace.
					\end{align}
		\end{enumerate}

\end{lemma}

\begin{proof}
In \eqref{E:C11BOUNDFORINVERSECHOVGEOTODOUBLENULL} and just above it, we proved that $\CHOVgeotodoublenull^{-1}$ extends as a $C^{1,1}_{\textnormal{null}}$ function to the compact, convex domain $[\leftubar,\ubarboot]\times[\moreinterestingu_1,\moreinterestingu_2] \times \mathbb{T}^2$ such that \eqref{E:KEYJACOBIANDETERMINANTESTIMATECHOVDOUBLENULLTOGEO} holds on $[\leftubar,\ubarboot]\times[\moreinterestingu_1,\moreinterestingu_2] \times \mathbb{T}^2$. In particular, this yields \eqref{E:CLOSEDKEYJACOBIANDETERMINANTESTIMATECHOVDOUBLENULLTOGEO} and \eqref{E:CLOSEDVERSIONC11BOUNDFORINVERSECHOVGEOTODOUBLENULL}. From these facts, \eqref{E:PARTIALTIMEDERIVATIVEOFUBARAPPROXIMATELYUNITY}, and the fact that $(t,u,x^2,x^3) =  \CHOVgeotodoublenull^{-1}(\ubar,u,x^2,x^3)$, we conclude that the map $\CHOVgeotodoublenull^{-1}$ with domain $[\leftubar,\ubarboot]\times[\moreinterestingu_1,\moreinterestingu_2]$ has a global inverse, i.e., that $\CHOVgeotodoublenull$ extends to the compact domain $\characteristicdiamondtwoarg{[\leftubar,\ubarboot]}{[\moreinterestingu_1,\moreinterestingu_2]}$ as an inverstible map such that $\ubar$ satisfies \eqref{E:CLOSEDVERSIONKEYJACOBIANDETERMINANTESTIMATECHOVGEOTODOUBLENULL}. From this extension, \eqref{E:CLOSEDVERSIONLUBARESTIMATECHOVGEOTODOUBLENULL}, follows from \eqref{E:BALDERIVATIVEOFUBARAPPROXIMATELYUNITY}.

We now prove \eqref{E:COMPARABLEDISTANCES}.
	For $i=1,2$, we set 
	$q_i \eqdef (\ubar_i,u_i,x_i^2,x_i^3)$
	and
	$p_i \eqdef\CHOVgeotodoublenull^{-1}(q_i) \eqdef (t_i,u_i,x_i^2,x_i^3)$.
	We define 
	$\Delta \ubar \eqdef \ubar_2 - \ubar_1$,
	$\Delta u \eqdef u_2 - u_1$,
	and
	$
	|\Delta q|_{\textnormal{Taxi}}
	\eqdef
	|\Delta \ubar| + |\Delta u| + |\Delta x^2|_{\mathbb{T}} + |\Delta x^3|_{\mathbb{T}}
	$,
	where for $j=2,3$,
	$|\Delta x^j|_{\mathbb{T}}$ is the Euclidean distance between $x_2^j$ and $x_1^j$ in the torus.
	We similarly define
	$\Delta t \eqdef t_2 - t_1$ and
	$
	|\Delta p|_{\textnormal{Taxi}}
	\eqdef
	|\Delta t| + |\Delta u| + |\Delta x^2|_{\mathbb{T}} + |\Delta x^3|_{\mathbb{T}}
	$.
	Note that $|\Delta q|_{\textnormal{Taxi}} \approx \textnormal{dist}_{\textnormal{flat}}(q_1,q_2)$
	and
	$|\Delta p|_{\textnormal{Taxi}} \approx \textnormal{dist}_{\textnormal{flat}}(p_1,p_2)$.
	Without loss of generality, we assume 
	$\Delta \ubar \geq 0$.
	Then by 
	\eqref{E:CLOSEDKEYJACOBIANDETERMINANTESTIMATECHOVDOUBLENULLTOGEO}
	and
	\eqref{E:CLOSEDVERSIONC11BOUNDFORINVERSECHOVGEOTODOUBLENULL},
	there is a constant $C > 1$ such that
	$
	\frac{1}{C} \Delta \ubar 
	-
	C(|\Delta u| + |\Delta x^2|_{\mathbb{T}} + |\Delta x^3|_{\mathbb{T}})
	\leq
	\Delta t 
	\leq C \Delta \ubar + C(|\Delta u| + |\Delta x^2|_{\mathbb{T}} + |\Delta x^3|_{\mathbb{T}})
	$.
	Hence, we see that there exists a (possibly different) constant $C > 1$ such that:
	\begin{align} \label{E:QUANTITATIVESTEPINPROOFOFCOMPARABLEDISTANCES}
	\left| 
	\frac{1}{C} \Delta \ubar 
	-
	C\left(|\Delta u| + |\Delta x^2|_{\mathbb{T}} + |\Delta x^3|_{\mathbb{T}} \right)
	\right|
	+
	|\Delta u| + |\Delta x^2|_{\mathbb{T}} + |\Delta x^3|_{\mathbb{T}}
	& 
	\leq
	|\Delta p|_{\textnormal{Taxi}}
	\leq
	C\left(|\Delta \ubar| + |\Delta u| + |\Delta x^2|_{\mathbb{T}} + |\Delta x^3|_{\mathbb{T}} \right).
	\end{align}
	From \eqref{E:QUANTITATIVESTEPINPROOFOFCOMPARABLEDISTANCES}, 
	it follows that 
	$
	|\Delta p|_{\textnormal{Taxi}}
	\approx
	|\Delta q|_{\textnormal{Taxi}}
	$,
	which implies \eqref{E:COMPARABLEDISTANCES}.
	
	We now prove the quasi-convexity of $\characteristicdiamondtwoarg{[\leftubar,\ubarboot]}{[\moreinterestingu_1,\moreinterestingu_2]}$. Let $p_1,p_2 \in \characteristicdiamondtwoarg{[\leftubar,\ubarboot]}{[\moreinterestingu_1,\moreinterestingu_2]}$, and let $q_1,q_2 \in [\leftubar,\ubarboot] \times [\moreinterestingu_1,\moreinterestingu_2] \times \mathbb{T}^2$ be the unique points such that $p_i = \CHOVgeotodoublenull^{-1}(q_i)$ as above. Let $\ell$ be a straight line in $[\leftubar,\ubarboot] \times [\moreinterestingu_1,\moreinterestingu_2] \times \mathbb{T}^2$ whose flat length is equal to $\textnormal{dist}_{\textnormal{flat}}(q_1,q_2)$. From \eqref{E:CLOSEDVERSIONC11BOUNDFORINVERSECHOVGEOTODOUBLENULL}, it follows that the image curve $\CHOVgeotodoublenull^{-1}(\ell)$ has Euclidean length that is $\lesssim \textnormal{dist}_{\textnormal{flat}}(p_1,p_1)$, as desired.  
	
\eqref{E:SOBOELVEMBEDDINGRELYINGONQUASICONVEXITY} is a standard Sobolev embedding result (see, for example, \cite{pHpKhT2008}*{Theorem\,7}),
	which relies on the
	quantitative quasi-convexity of $\characteristicdiamondtwoarg{[\leftubar,\ubarboot]}{[\moreinterestingu_1,\moreinterestingu_2]}$. 

\eqref{E:CLOSEDVERSIONC11BOUNDFORCHOVDOUBLENULLTOGEO} follows from \eqref{E:W2INFINITYBOUNDFORCHOVFROMGEOTODOUBLENULLOPEN} and \eqref{E:SOBOELVEMBEDDINGRELYINGONQUASICONVEXITY}. 

The existence of the function $\Cartesiantisafunctiononlevelsetsofingoingcharacteristicsarg{\ubar}$ such that \eqref{E:LEVELSETSOFUBARAREAGRAPH} holds, as well as estimate \eqref{E:C11BOUNDFORCONSTANTTIMEFUNCTIONGRAPH}, follows from the fact that $(t,u,x^2,x^3) = \CHOVgeotodoublenull^{-1}(\ubar,u,x^2,x^3)$ and the estimate \eqref{E:CLOSEDVERSIONC11BOUNDFORINVERSECHOVGEOTODOUBLENULL}, i.e., $\Cartesiantisafunctiononlevelsetsofingoingcharacteristicsarg{\ubar}(u,x^2,x^3)$ is the first component of $\CHOVgeotodoublenull^{-1}(\ubar,u,x^2,x^3)$.
\end{proof}

\begin{corollary}[Preliminary estimates for the double-null geometry scalars] \label{C:PRELIMINARYESTIMATESFORTHEDOUBLENULLACOUSTICSCALARS}
Let $\ingoingmu, \, \MagnitueofinnerproductofnewLandnewuL, \, \ReciprocalLunitAppliedtoTimeFunction$, and $\ReciprocaluLunitAppliedtoTimeFunction$ be the double-null acoustic scalrs from \eqref{E:ULFOLIATIONDENSITY}--\eqref{E:RATIOOFNULLGEOSICINNERPRODUCTANDFOLIATIONDENSITY}. Then the following estimates hold on $\characteristicdiamondtwoarg{[\leftubar,\ubarboot]}{[\moreinterestingu_1,\moreinterestingu_2]}$:
\begin{subequations}
	\begin{align} 
		\ingoingmu & \approx 1, \label{E:INGOINGMUISAPPROX1} \\
		\MagnitueofinnerproductofnewLandnewuL & \approx \upmu, \label{E:MAGNITUDEOFINNERPRODUCTOFNEWLANDNEWULAPPROXIMATELYMU} \\
		\ReciprocaluLunitAppliedtoTimeFunction & \approx \upmu. \label{E:RECIPROCALULUNITAPPLIEDTOTIMEFUNCTIONAPPROXIMATELYMU}
	\end{align} 
\end{subequations}
Moreover, the following sharp estimate holds on $\characteristicdiamondtwoarg{[\leftubar,\ubarboot]}{[\moreinterestingu_1,\moreinterestingu_2]}$: 
	\begin{align} \label{E:SHARPESTIMATEFORRATIOOFRATIOOFNULLGEOSICINNERPRODUCTANDFOLIATIONDENSITYANDINGOINGMU}
		\ReciprocaluLunitAppliedtoTimeFunction = \frac{ \upmu \ReciprocalLunitAppliedtoTimeFunction}{\ingoingmu} =  \frac{1}{2} \upmu \left\{ 1 +  \mathcal{O}(\fundbootsmall)\right\}. 
	\end{align}
\end{corollary}
\begin{proof}
	Estimates \eqref{E:INGOINGMUISAPPROX1} and \eqref{E:SHARPESTIMATEFORRATIOOFRATIOOFNULLGEOSICINNERPRODUCTANDFOLIATIONDENSITYANDINGOINGMU} follow from \eqref{E:INGOINGMUDETERMINEDBYEVERYTHINGELSE}, the identity $\ReciprocalLunitAppliedtoTimeFunction = \frac{1}{\Lunit\ubar}$, and \eqref{E:CLOSEDVERSIONLUBARESTIMATECHOVGEOTODOUBLENULL}--\eqref{E:SMALLC01ESTIMATESFORUBAR}. Estimates \eqref{E:RECIPROCALULUNITAPPLIEDTOTIMEFUNCTIONAPPROXIMATELYMU}--\eqref{E:SHARPESTIMATEFORRATIOOFRATIOOFNULLGEOSICINNERPRODUCTANDFOLIATIONDENSITYANDINGOINGMU} follow from similar arguments and the definitions \eqref{E:RECIPROCALOFNEGATIVENULLGEODESICVECTORFIELDINNERPRODUCT}-\eqref{E:RATIOOFNULLGEOSICINNERPRODUCTANDFOLIATIONDENSITY}.
\end{proof}

\subsection{H\"{o}lder-space extensions to the compact set $\characteristicdiamondtwoarg{[\leftubar,\ubarboot]}{[\moreinterestingu_1,\moreinterestingu_2]}$}
\label{SS:HOLDERSPACEEXTENSIONS}
With the help of the bootstrap assumptions and
Lemma\,\ref{L:PROPERTIESANDDIFFEOMORPHICEXTENSIONOFDOUBLENULLCOORDINATES}, we now show that various wave variables extend to
the compact set $\characteristicdiamondtwoarg{[\leftubar,\ubarboot]}{[\moreinterestingu_1,\moreinterestingu_2]}$
as functions with substantial H\"{o}lder regularity relative to the geometric coordinates.

\begin{lemma}[H\"{o}lder-space extensions to the compact set 
$\characteristicdiamondtwoarg{[\leftubar,\ubarboot]}{[\moreinterestingu_1,\moreinterestingu_2]}$]
	\label{L:CONTINUOUSEXTNESION}
	The following quantities extend to the compact set 
	$\characteristicdiamondtwoarg{[\leftubar,\ubarboot]}{[\moreinterestingu_1,\moreinterestingu_2]}$
	as elements of the following H\"{o}lder spaces, and their  
	corresponding spacetime H\"{o}lder norms on $\characteristicdiamondtwoarg{[\leftubar,\ubarboot]}{[\moreinterestingu_1,\moreinterestingu_2]}$
	are bounded by $\leq C$:
	\begin{itemize}
		\item $\wavearray 
			\in C_{\textnormal{geo}}^{3,1}\left(\characteristicdiamondtwoarg{[\leftubar,\ubarboot]}{[\moreinterestingu_1,\moreinterestingu_2]}\right)$
		\item $\Upsilon \in C_{\textnormal{geo}}^{3,1}\left(\characteristicdiamondtwoarg{[\leftubar,\ubarboot]}{[\moreinterestingu_1,\moreinterestingu_2]}\right)$
		\item $\Lunit^i, \, \upmu \in C_{\textnormal{geo}}^{2,1}\left(\characteristicdiamondtwoarg{[\leftubar,\ubarboot]}{[\moreinterestingu_1,\moreinterestingu_2]}\right)$
		\item $\ubar \in C_{\textnormal{geo}}^{1,1}\left(\characteristicdiamondtwoarg{[\leftubar,\ubarboot]}{[\moreinterestingu_1,\moreinterestingu_2]}\right)$
	\end{itemize}
\end{lemma}

\begin{proof}
	The results for $\ubar$ were already proved as \eqref{E:CLOSEDVERSIONC11BOUNDFORCHOVDOUBLENULLTOGEO}.
	For the remaining results, we give the proof only for $\wavearray$ since the
	other wave variables can be handled using nearly identical arguments.
	To proceed, we note that
	Lemma\,\ref{L:COMMUTATORSTOCOORDINATES},
	Prop.\,\ref{P:SCHEMATICSTRUCTUREOFVARIOUSTENSORSINTERMSOFCONTROLVARS},
	and the bootstrap assumptions imply that
	$\| \wavearray \|_{W_{\textnormal{geo}}^{4,\infty}\left(\characteristicdiamondtwoarg{(\leftubar,\ubarboot)}{(\moreinterestingu_1,\moreinterestingu_2)}\right)}
	\leq C
	$.
	From this bound and \eqref{E:SOBOELVEMBEDDINGRELYINGONQUASICONVEXITY},
	we conclude that
	$
	\| \wavearray \|_{C_{\textnormal{geo}}^{3,1}\left(\characteristicdiamondtwoarg{[\leftubar,\ubarboot]}{[\moreinterestingu_1,\moreinterestingu_2]}\right)} 
	\leq C
	$
	as desired.
\end{proof}

\subsection{Properties of $\CHOVdoublenulltoubarnewulmu$}
\label{SS:CHOVFROMDOUBLENULLCOORDINATESTOUBARNEWULMUCOORDINATES}
In this section, 
we provide a detailed analysis of the change of variables map 
$\CHOVdoublenulltoubarnewulmu(\ubar,u,x^2,x^3) = (\ubar, \newuL \upmu,x^2,x^3)$.

\begin{lemma}[Properties of $\CHOVdoublenulltoubarnewulmu$]
	\label{L:CHOVFROMDOUBLENULLTONEWULMUCOORDINATES}
	\hfill
	
	\noindent \underline{\textbf{Estimates for and diffeomorphism properties of $\CHOVdoublenulltoubarnewulmu$}}.
	The map $\CHOVdoublenulltoubarnewulmu(\ubar,u,x^2,x^3) = (\ubar, \newuL \upmu,x^2,x^3)$
	from Def.\,\ref{D:CHOVINVOLVINGTHEDOUBLENULLCOORD}
	extends to a $C_{\textnormal{null}}^{0,1}$ map on $[\leftubar,\ubarboot]\times[\moreinterestingu_1,\moreinterestingu_2] \times \mathbb{T}^2$
	satisfying:
	\begin{align} \label{E:NEWULMUCHOVMAPC01BOUND}
		\left\| 	
			\CHOVdoublenulltoubarnewulmu 
		\right\|_{C_{\textnormal{null}}^{0,1}\left([\leftubar,\ubarboot]\times[\moreinterestingu_1,\moreinterestingu_2] \times \mathbb{T}^2 \right)}
		& \leq C.
	\end{align}
	Moreover, for each fixed $(\ubar,x^2,x^3) \in [\leftubar,\ubarboot]\times\T^2$, the map $(\ubar,u,x^2,x^3) \mapsto \newuL \upmu(\ubar,u,x^2,x^3)$ is strictly increasing on $u \in [\moreinterestingu_1,\moreinterestingu_2]$. In particular, $\CHOVdoublenulltoubarnewulmu$ 
	is a Lipeomorphism from the compact, convex set
	$[\leftubar,\ubarboot]\times[\moreinterestingu_1,\moreinterestingu_2] \times \mathbb{T}^2$
	onto its (compact) image 
	$\CHOVdoublenulltoubarnewulmu\left([\leftubar,\ubarboot]\times[\moreinterestingu_1,\moreinterestingu_2] \times \mathbb{T}^2 \right)$,
	where $\CHOVdoublenulltoubarnewulmu\left([\leftubar,\ubarboot]\times[\moreinterestingu_1,\moreinterestingu_2] \times \mathbb{T}^2 \right)$ 
	enjoys the following properties:
	\begin{enumerate}
	\item It contains
	$
	[\leftubar,\ubarboot]\times \{0\} \times \mathbb{T}^2
	$.
	\item It contains
	$
	(\leftubar,\ubarboot) \times \{0\} \times \mathbb{T}^2
	$
	in its interior.
	\item It
	is quasi-convex in the following sense:
	every pair of points
	$r_1, r_2 \in \CHOVdoublenulltoubarnewulmu \left([\leftubar,\ubarboot]\times[\moreinterestingu_1,\moreinterestingu_2] \times \mathbb{T}^2 \right)$
	is connected by a $C^1$ curve in 
	$\CHOVdoublenulltoubarnewulmu \left([\leftubar,\ubarboot]\times[\moreinterestingu_1,\moreinterestingu_2] \times \mathbb{T}^2 \right)$
	whose length with respect to the standard Euclidean metric on 
	$\mathbb{R} \times \mathbb{R} \times \mathbb{T}^2$
	is $\lesssim \textnormal{dist}_{\textnormal{flat}}(r_1,r_2)$,
	where $\textnormal{dist}_{\textnormal{flat}}(r_1,r_2)$ is the
	standard Euclidean distance between $r_1$ and $r_2$ in the flat space 
	$\mathbb{R} \times \mathbb{R} \times \mathbb{T}^2$.
	\end{enumerate}

	\medskip
	\noindent 
	\underline{\textbf{Properties of the $\newuL \upmu$-adapted tori as graphs over $\mathbb{T}^2$ in double-null coordinates}}.
	There exists a family of functions
	$\lbrace \Cartesiantisafunctiononubarnewultori_{\ubar}, \, \Eikonalisafunctiononubarnewultori_{\ubar}
	\rbrace_{\ubar \in [\leftubar,\ubarboot]}$
	on $\mathbb{T}^2$
	that, for $\ubar \in [\leftubar,\ubarboot)$, are equal to the functions 
	from Sect.\,\ref{SSS:BOOTSTRAPASSUMPTIONFORTORISTRUCTURE}, and
	such that
	for each $\ubar \in [\leftubar,\ubarboot]$,
	we have:
	\begin{align}
	\sup_{\ubar \in [\leftubar,\ubarboot]} 
	\|  \Cartesiantisafunctiononubarnewultori_{\ubar} \|_{C^{0,1}(\mathbb{T}^2)}
	& \leq C, \label{E:C01BOUNDFORCARTESIANTCOORDINATEGRAPHTORI} \\
	\sup_{\ubar \in [\leftubar,\ubarboot]} 
	\|  \Eikonalisafunctiononubarnewultori_{\ubar} \|_{C^{0,1}(\mathbb{T}^2)}
	& \leq C, \label{E:C01BOUNDFORUDOUBLENULLCOORDINATEGRAPHTORI}
	\end{align}
	\begin{align}
	\CHOVgeotodoublenull\left(\argubarnewultorus{\ubar}\right)
	& 
	=
	\left\lbrace
		\left(\ubar,\Eikonalisafunctiononubarnewultori_{\ubar}(x^2,x^3),x^2,x^3 \right) 
		\ | \ 
		(x^2,x^3) \in \mathbb{T}^2
	\right\rbrace,
	\label{E:LASTSLICETORIAREGRAPHSABOVEFLATTORIINGEOMETRICCOORDINATES}
	\end{align}
	where $\CHOVgeotodoublenull$ is the change of variables map defined in \eqref{E:CHOVFROMGEOTODOUBLENULLCOORDINATES}.
	In particular,
	\begin{align} \label{E:LASTSLICELEVELSETSTRUCTUREANDLOCATIONOFMIN} 
		\argubarnewultorus{\ubarboot}
		& 
		\subset
		 \ingoingcharacteristicsurfacetwoarg{\ubarboot}{[-\frac{1}{2}\interestingu, \frac{3}{4}\moreinterestingu_2 + \frac{1}{2} \interestingu]}.
	\end{align}

	\medskip
	\noindent 
	\underline{\textbf{Foliation of $\newuLmulevelset^{[\leftubar,\ubarboot]}$ by $\argubarnewultorus{\ubar}$:}} 
	For any $(x^2,x^3) \in \T^2$, the map $\ubar \mapsto \Cartesiantisafunctiononubarnewultori_{\ubar}(x^2,x^3)$ is strictly increasing:
	\begin{align} \label{E:NEWULMUHYPERSURFACEEMBEDDINGSTRICTLYINCREASING}
		\Cartesiantisafunctiononubarnewultori_{\ubar_1}(x^2,x^3) < \Cartesiantisafunctiononubarnewultori_{\ubar_2}(x^2,x^3), & &  \ubar_1 < \ubar_2. & & \ubar_1,\ubar_2 \in [\leftubar,\ubarboot].
	\end{align}
	In particular, we have: 
	\begin{align}
	\newuLmulevelset^{[\leftubar,\ubarboot]} = \bigcup_{\ubar \in [\leftubar,\ubarboot]} \argubarnewultorus{\ubar}. \label{E:LEVELSETOFNEWULMUFOLIATED}
	\end{align}

	\medskip
	\noindent \underline{\textbf{Properties of $\CHOVdoublenulltoubarnewulmu^{-1}$}}.
	With $\CHOVdoublenulltoubarnewulmu^{-1}$ denoting the inverse map
	and $\CHOVgeotodoublenull^{-1}$ denoting the inverse of 
	the change of variables map $\CHOVgeotodoublenull$ defined in \eqref{E:CHOVFROMGEOTODOUBLENULLCOORDINATES},
	the following holds for $\ubar \in [\leftubar,\ubarboot]$:
	\begin{align} \label{E:PHIINVERSEIMAGEOFTORUSISTORUSCONTAINEDINUBARFUNCTIONLEVELSET}
		\CHOVgeotodoublenull^{-1}
		\circ
		\CHOVdoublenulltoubarnewulmu^{-1}\left(\lbrace \ubar \rbrace \times \lbrace 0 \rbrace \times \mathbb{T}^2 \right)
		= 
		\argubarnewultorus{\ubar} 
		\subset
		\ingoingcharacteristicsurfacetwoarg{\ubar}{[\moreinterestingu_1,\moreinterestingu_2]} ,
	\end{align}	
	and:
	\begin{align} \label{E:PHIINVERSEIMAGEOFTORUSCROSSMUINTERVALISTORUSCROSSINTERVALCONTAINEDININTERESTINGREGION}
		\CHOVgeotodoublenull^{-1}
		\circ
		\CHOVdoublenulltoubarnewulmu^{-1}
		\left([\leftubar,\ubarboot] \times 
		\lbrace 
			0 
		\rbrace 
		\times \mathbb{T}^2 \right)
		= 
		\newuLmulevelset^{[\leftubar,\ubarboot]}
		\subset
		\characteristicdiamondtwoarg{[\leftubar,\ubarboot]}{[\moreinterestingu_1,\moreinterestingu_2]}.
	\end{align}	
	In addition, $\CHOVdoublenulltoubarnewulmu^{-1}$ satisfies the following estimate:
	\begin{align} \label{E:INVERSEOFNEWULMUCHOVMAPC01BOUND}
		\left\| 
			\CHOVdoublenulltoubarnewulmu^{-1} \right 
		\|_{C^{0,1}\left( 
		\CHOVdoublenulltoubarnewulmu\left([\leftubar,\ubarboot]\times[\moreinterestingu_1,\moreinterestingu_2] \times \mathbb{T}^2 \right)
		\right)}
		& \leq C.
	\end{align}\end{lemma}

\begin{proof} \hfill

\noindent \textbf{Proof of \eqref{E:NEWULMUCHOVMAPC01BOUND}:} Lemmas\,\ref{L:COMMUTATORSTOCOORDINATES},
\ref{P:SCHEMATICSTRUCTUREOFVARIOUSTENSORSINTERMSOFCONTROLVARS},
\ref{L:PROPERTIESANDDIFFEOMORPHICEXTENSIONOFDOUBLENULLCOORDINATES}, 
and\,\ref{L:CONTINUOUSEXTNESION} yield \eqref{E:NEWULMUCHOVMAPC01BOUND}.

\medskip
\noindent \textbf{Proof of the Lipeomorphism properties of $\CHOVdoublenulltoubarnewulmu$
and the quasi-convexity of its image:}
Since the Fundamental Theorem of Calculus holds for absolutely continuous functions, and hence, for fixed $(\ubar,x^2,x^3) \in [\leftubar,\ubarboot] \times \mathbb{T}^2$, the Lipschitz function $g(u) \eqdef \newuL \upmu(\ubar,u,x^2,x^3)$, it follows from \eqref{E:BAMULESSREGULARTRANSVERSALCONVEXITY} that $g(u)$
	is strictly increasing on $[\moreinterestingu_1,\moreinterestingu_2]$.
	Hence, the map $(\ubar,u,x^2,x^3) \rightarrow \left(\ubar,\newuL \upmu,x^2,x^3 \right)$
	is a $C^{0,1}$ Lipeomorphism from 
	$[\leftubar,\ubarboot] \times [\moreinterestingu_1,\moreinterestingu_2]\times \mathbb{T}^2$ onto its image. From this monotonicity, \eqref{E:NEWULMUCHOVMAPC01BOUND}, and the convexity of $[\leftubar,\ubarboot] \times [\moreinterestingu_1,\moreinterestingu_2]\times \mathbb{T}^2$, it follows that the image of $\CHOVdoublenulltoubarnewulmu$ is quasi-convex. 

Next, we use \eqref{E:BAUBARNEWULMUTORILOCATION}, the fact that $\CHOVdoublenulltoubarnewulmu$ is a Lipeomorphism on $[\leftubar,\ubarboot] \times [\moreinterestingu_1,\moreinterestingu_2]\times \mathbb{T}^2$, and the fact that $\newuL \upmu|_{\argubarnewultorus{\ubar}} = 0$ to deduce that $\CHOVdoublenulltoubarnewulmu\left([\leftubar,\ubarboot] \times [\moreinterestingu_1,\moreinterestingu_2]\times \mathbb{T}^2\right)$ contains $[\leftubar,\ubarboot] \times \{0\}\times \mathbb{T}^2$ and that $(\leftubar,\ubarboot) \times \{0\}\times \mathbb{T}^2$ is contained in the interior.

\medskip
\noindent \textbf{Proof of the properties of $\Cartesiantisafunctiononubarnewultori_{\ubar}, \, \Eikonalisafunctiononubarnewultori_{\ubar}$ 
and the estimates
\eqref{E:C01BOUNDFORUDOUBLENULLCOORDINATEGRAPHTORI}--\eqref{E:LASTSLICELEVELSETSTRUCTUREANDLOCATIONOFMIN}}:

All of the results but \eqref{E:C01BOUNDFORCARTESIANTCOORDINATEGRAPHTORI} follow from the form  
\eqref{E:CHOVFROMDOUBLENULLTOUBARNEWULMUCOORDINATES}
of $\CHOVdoublenulltoubarnewulmu$,
the fact that 
$\CHOVdoublenulltoubarnewulmu^{-1}$ 
is $C^{0,1}$ on
the compact, quasi-convex set
$
\CHOVdoublenulltoubarnewulmu\left([\leftubar,\ubarboot] \times [\moreinterestingu_1,\moreinterestingu_2]\times \mathbb{T}^2\right)
$,
and the bootstrap assumption \eqref{E:BAUBARNEWULMUTORILOCATION},
which implies that for $\ubar \in [\leftubar,\ubarboot]$
and $(x^2,x^3) \in \mathbb{T}^2$,
we have
$
\CHOVdoublenulltoubarnewulmu^{-1}(\ubar,0,x^2,x^3) \in 
\ingoingcharacteristicsurfacetwoarg{\ubar}{[-\frac{1}{2}\interestingu, \frac{3}{4}\moreinterestingu_2 + \frac{1}{2} \interestingu]}
$.

Estimate \eqref{E:C01BOUNDFORCARTESIANTCOORDINATEGRAPHTORI} follows from \eqref{E:C11BOUNDFORCONSTANTTIMEFUNCTIONGRAPH} and the fact that $\Cartesiantisafunctiononlevelsetsofingoingcharacteristicsarg{\ubar}(\Eikonalisafunctiononubarnewultori_{\ubar}(x^2,x^3),x^2,x^3) = \Cartesiantisafunctiononubarnewultori_{\ubar}(x^2,x^3)$. 

\medskip
\noindent \textbf{Proof of \eqref{E:NEWULMUHYPERSURFACEEMBEDDINGSTRICTLYINCREASING} and the foliation of $\newuLmulevelset^{[\leftubar,\ubarboot]}$ by $\argubarnewultorus{\ubar}$:} This is an immediate consequence of the Fundamental Theorem of Calculus and \eqref{E:CLOSEDKEYJACOBIANDETERMINANTESTIMATECHOVDOUBLENULLTOGEO}. 

\medskip
\noindent \textbf{Proof of \eqref{E:PHIINVERSEIMAGEOFTORUSISTORUSCONTAINEDINUBARFUNCTIONLEVELSET}--\eqref{E:INVERSEOFNEWULMUCHOVMAPC01BOUND}:}

Set identities and inclusions \eqref{E:PHIINVERSEIMAGEOFTORUSISTORUSCONTAINEDINUBARFUNCTIONLEVELSET}--\eqref{E:PHIINVERSEIMAGEOFTORUSCROSSMUINTERVALISTORUSCROSSINTERVALCONTAINEDININTERESTINGREGION} follow from the \eqref{E:BAUBARNEWULMUTORILOCATION}, \eqref{E:LASTSLICELEVELSETSTRUCTUREANDLOCATIONOFMIN}, \eqref{E:LEVELSETOFNEWULMUFOLIATED}, definitions \eqref{E:DEFOFNEWULADAPTEDTORI} and \eqref{E:DEFOFNEWULADAPTEDTORI}, definitions \eqref{E:CHOVFROMGEOTODOUBLENULLCOORDINATES} and \eqref{E:CHOVFROMDOUBLENULLTOUBARNEWULMUCOORDINATES}, the fact that $\CHOVgeotodoublenull$ is a Lipeomorphism on $\characteristicdiamondtwoarg{[\leftubar,\ubarboot]}{[\moreinterestingu_1,\moreinterestingu_2]}$ (see Lemma\,\ref{L:PROPERTIESANDDIFFEOMORPHICEXTENSIONOFDOUBLENULLCOORDINATES}),  and the fact that $\CHOVdoublenulltoubarnewulmu$ is a Lipeomorphism from $[\leftubar,\ubarboot] \times [\moreinterestingu_1,\moreinterestingu_2]\times \mathbb{T}^2$ onto a set containing $[\leftubar,\ubarboot] \times\{0\}\times \mathbb{T}^2$.

Estimate \eqref{E:INVERSEOFNEWULMUCHOVMAPC01BOUND} follows from \eqref{E:NEWULMUCHOVMAPC01BOUND} and the fact that $\CHOVdoublenulltoubarnewulmu$ is a Lipeomorphism from $[\leftubar,\ubarboot] \times [\moreinterestingu_1,\moreinterestingu_2]\times \mathbb{T}^2$ onto its image, and the quasi-convexity of $\CHOVdoublenulltoubarnewulmu\left([\leftubar,\ubarboot] \times [\moreinterestingu_1,\moreinterestingu_2]\times \mathbb{T}^2\right)$.

\end{proof}


\section{Derivatives of the logarithmic density in terms of velocity, vorticity, and entropy} \label{S:DENSITYDERIVATIVESINTERMSOFOTHERS}

In this section we express $Z\LogDensity$ in terms of the velocity, vorticity, and entropy for $Z \in \mathscr{Z} = \{ \Lunit,\muX,\Yvf{2},\Yvf{3}\}$. Importantly, whenever $Z \in \tander$, the derivatives of $v$ present are also tangent to $\outgoingcharacteristicsurfacetwoarg{u}{[\leftubar,\ubar)}$. This will guarantee that null forms in the commuted wave equation involving $\LogDensity$ are replaced by null forms involving only $v$ and below-top-order vorticity and entropy.

\begin{proposition}[Derivatives of logarithmic density] \label{P:ENTROPYDERIVATIVESINTERMSOFOTHERS}
Let $\angularcoeffinmuweightedspatialcartesian_i^A$ be the $3\times 2$ matrix from  \eqref{E:ANGULARCOEFFINMUWEIGHTEDCARTESIAN}. Then the following identities hold:
\begin{subequations}	
	\begin{align}
		\upmu \Transport \LogDensity & = - X_i \muX v^i - \upmu \angularcoeffinmuweightedspatialcartesian_i^A \Yvf{A} v^i, \label{E:MUBRHOINTERMSOFV} \\
		\Lunit \LogDensity & = X_i \Lunit v^i - \angularcoeffinmuweightedspatialcartesian_i^A \Yvf{A} v^i  + \exp(-\LogDensity)\frac{p_{;\Ent}}{\overline{\varrho}} \GradEnt^i X_i \label{E:LRHOINTERMSOFVANDENT}\\
		\muX \LogDensity & = - X_i \muX v^i - \upmu X_i \Lunit v^i - \upmu  \exp(-\LogDensity)\frac{p_{;\Ent}}{\overline{\varrho}} \GradEnt^i X_i \label{E:MUXRHOINTERMSOFVANDENT}
	\end{align}
Additionally, if $W$ is any $\ell_{t,u}$-tangent vectorfield, we have: 
\begin{align}
	W \LogDensity & = - X_i W v^i - W_i \Lunit v^i - \Speed^{-2}\exp(\LogDensity) X^i W^j \upepsilon_{ijk} \vortrenormalized^k - \exp(-\LogDensity)\frac{p_{;\Ent}}{\overline{\varrho}}W_i \GradEnt^i  \label{E:TANGENTIALDERIVATIVESOFRHOINTERMSOFVANDVORTANDENT}
\end{align} 
\end{subequations}Consequently, for any $1\le N \le \Ntop$, we have: 
\begin{subequations}
	\begin{align}
		\left| \tander^{[1,N+1]} \LogDensity \right| & \lesssim \left| \tander^{[1,N+1]} \velocityarray\right| +  \left|\tander^{\le N} (\vortrenormalized,\GradEnt)\right| + \fundbootsmall \left| \tander^{[1,N]} \controlvars\right|  \label{E:HIGHERORDERTANGENTIALDERIVATIVESOFDENSITY} \\
		\left| \comdersmall^{[1,N+1];1} \LogDensity \right| & \lesssim \left| \comdersmall^{[1,N+1];1} \velocityarray\right| + \upmu \left| \tander^{[1,N+1]}\velocityarray\right| +  \left|\tander^{\le N} (\vortrenormalized,\GradEnt)\right| + \fundbootsmall \left| \tandersmall^{[1,N]} \badcontrolvars\right|  \label{E:HIGHERORDEREXACTLYONEMUXINDERIVATIVESOFDENSITY} \\		
		\left| \comdersmall^{[1,N+1]; 2} \LogDensity \right| & \lesssim \left|\comdersmall^{[1,N+1]; 2}\velocityarray\right| + \upmu \left| \comdersmall^{[1,N+1];\le 1} \velocityarray\right| + \upmu  \left|\tander^{\le N} (\vortrenormalized,\GradEnt)\right| + \fundbootsmall \left|\comderdoublesmall^{[1,N];\le 1} \badcontrolvars\right|\label{E:HIGHERORDERMIXEDDERIVATIVESOFDENSITY}
	\end{align}
\end{subequations}
\end{proposition} 

\begin{proof}
Identity \eqref{E:MUBRHOINTERMSOFV} follows immediately from \eqref{E:BLOGDENSITYEVOLUTION} and \eqref{E:CARTESIANCOMPONENTSOFFIRSTFUNDAMENTLAFORMOFSIGMAT}. 

Next, we prove \eqref{E:MUXRHOINTERMSOFVANDENT}. By \eqref{E:CARTESIANCOMPONENTSOFFIRSTFUNDAMENTLAFORMOFSIGMAT}, upon contracting \eqref{E:BVIEVOLUTION} with $X_i$, we have:
\begin{align}\label{E:INTERMEDIATESTEPFORGRADIENTDECOMPOFRHO1}
	X_i \Lunit v^i + X_i X v^i = - X \LogDensity - \exp(-\LogDensity)\frac{p_{;\Ent}}{\overline{\varrho}} X_i S^i.
\end{align}
The result follows upon multiplying \eqref{E:INTERMEDIATESTEPFORGRADIENTDECOMPOFRHO1} by $\upmu$.

To prove \eqref{E:LRHOINTERMSOFVANDENT}, we use \eqref{E:BLOGDENSITYEVOLUTION} and \eqref{E:EXPRESSIONFORMUWEIGHTEDSPATIALCARTESIANDERIVATIVESWRTCOMMUTATORS} to see that $\Lunit \LogDensity + X \LogDensity = - X_i Xv^i - \angularcoeffinmuweightedspatialcartesian_i^A \Yvf{A} v^i$. Inserting \eqref{E:INTERMEDIATESTEPFORGRADIENTDECOMPOFRHO1} into this identity we observe a cancelation of $X_i X v^i$ and conclude the proof of \eqref{E:LRHOINTERMSOFVANDENT}. 

Finally, let $W$ be any $\ell_{t,u}$-tangent vectorfield. Contracting \eqref{E:BVIEVOLUTION} with $W_i$, we have:
\begin{align} \label{E:INTERMEDIATESTEPFORGRADIENTDECOMPOFRHO2}
	W \LogDensity = - W_i \Lunit v^i - W_i Xv^i - \exp(-\LogDensity) \frac{p_{;\Ent}}{\overline{\varrho}} W_i S^i 
\end{align}
The proof of \eqref{E:TANGENTIALDERIVATIVESOFRHOINTERMSOFVANDVORTANDENT} follows from the identities $W_i Xv^i = W_i X^j \p_j v^i = W_i X^j \p_i v^j + W_i X^j \upepsilon_{jik} \vort^k = W_i X^j \p_i v^j +\exp(\LogDensity) W_iX^j \upepsilon_{jik}\vortrenormalized^k$.

The proofs of \eqref{E:HIGHERORDERTANGENTIALDERIVATIVESOFDENSITY}--\eqref{E:HIGHERORDERMIXEDDERIVATIVESOFDENSITY} follow from \eqref{E:MUXRHOINTERMSOFVANDENT}--\eqref{E:TANGENTIALDERIVATIVESOFRHOINTERMSOFVANDVORTANDENT}, straightforward applications of the Leibniz rule, Prop.\,\ref{P:SCHEMATICSTRUCTUREOFVARIOUSTENSORSINTERMSOFCONTROLVARS}, and the bootstrap assumptions. We clarify that equations \eqref{E:RENORMALIZEDVORTICTITYTRANSPORTEQUATION}--\eqref{E:GRADENTROPYTRANSPORT} imply that $\muX(\vortrenormalized,\GradEnt) = \upmu \Lunit(\vortrenormalized,\GradEnt) + \smoothfunction(\comder^{\le 1} \wavearray) (\vortrenormalized,\GradEnt)$ and so even two transversal derivatives $\muX\muX(\vortrenormalized,\GradEnt)$ depend only on tangential derivatives $\tander^{[1,2]}(\vortrenormalized,\GradEnt)$.

\end{proof}

\section{Preliminary pointwise, commutator, and differential operator comparison estimates}
\label{S:PRELIMINARYPOINTWISECOMMUTATORANDOPERATORCOMPARISON}
In this section, we use the bootstrap assumptions to derive several preliminary estimates,
including pointwise estimates,
commutator estimates, 
and differential operator estimates comparing $\angLie$ and $\newangD$. The results of this section were proven in our companion work \cite{abbrescia2022emergence}. We simply state the results here for completeness, and refer the reader to \cite{abbrescia2022emergence} for the proofs.

\subsection{The norm of the $\ell_{t,u}$-tangent commutation vectorfields and simple comparison estimates}
\label{SS:NORMOFSMOOTHTORITANGENTCOMMUTATORSANDSIMPLECOMPARISON}

\begin{lemma}[The norm of the $\ell_{t,u}$-tangent commutation vectorfields and simple comparison estimates; \cite{abbrescia2022emergence}*{Lemma 14.1}]
	\label{L:NORMOFSMOOTHTORITANGENTCOMMUTATORSANDSIMPLECOMPARISON}
	The $\ell_{t,u}$-tangent commutation vectorfields $\Angularset = \lbrace \Yvf{2}, \Yvf{3} \rbrace$ 
	satisfy the following pointwise estimates on $\characteristicdiamondtwoarg{[\leftubar,\ubarboot)}{[\moreinterestingu_1,\moreinterestingu_2]}$:
	\begin{align} \label{E:POINTWISESEMINORMOFYVECTORFIELDS}
		|\Yvf{A}|_{\gtorus} 
		& 
		= 
		1 
		+ 
		\mathcal{O}_{\mydiam}(\mathring{\upalpha}^{1/2}).
	\end{align}	
	
	Moreover, for any type $\binom{m}{n}$ 
	$\ell_{t,u}$-tangent tensorfield $\upxi_{\beta_1 \cdots \beta_n}^{\alpha_1 \cdots \alpha_m}$,
	the following pointwise estimates hold on $\characteristicdiamondtwoarg{[\leftubar,\ubarboot)}{[\moreinterestingu_1,\moreinterestingu_2]}$,
	where $\Angularset$ is defined in \eqref{E:COMMUTATIONVECTORFIELDS}:
	\begin{align} \label{E:SMOOTHTORUSNORMCOMPARBLETOTANGENTIALCONTRACTIONS}
		\begin{split}	
			|\upxi|_{\gtorus}^2
			& 
			=
			\left\lbrace
			1 + \mathcal{O}_{\mydiam}(\mathring{\upalpha}^{1/2})
			\right\rbrace
			\sum_{\substack{ {^{(1)}U}, \cdots, {^{(m)}U} \in \Angularset \\ {^{(1)}V}, \cdots, {^{(n)}V} \in \Angularset}} 
			\left|
			{^{(1)}U_{\alpha_1}} \cdots {^{(m)}U_{\alpha_m}} 
			{^{(1)}V^{\beta_1}} \cdots {^{(n)}V^{\beta_n}}  
			\upxi_{\beta_1 \cdots \beta_n}^{\alpha_1 \cdots \alpha_n}
			\right|^2
			\\
			&
			=
			\left\lbrace
			1 + \mathcal{O}_{\mydiam}(\mathring{\upalpha}^{1/2})
			\right\rbrace
			\sum_{\substack{A_1,\cdots,A_m =2,3 \\ B_1,\cdots,B_n =2,3}}
			\left|
			\upxi_{B_1 \cdots B_n}^{A_1 \cdots A_n}
			\right|^2.
		\end{split}
	\end{align}
\end{lemma}

\subsection{Basic facts that we use silently when deriving estimates}
\label{SS:SILENTFACTS}
In the rest of the paper, we silently use the following basic facts.
\begin{enumerate}
	\item All of the estimates we derive hold on the bootstrap region 
	$\characteristicdiamondtwoarg{[\leftubar,\ubarboot)}{[\moreinterestingu_1,\moreinterestingu_2]}$. 
	Moreover, in deriving estimates, 
	we often rely on the parameter- and data-size assumptions
	of Sects.\,\ref{SS:PARAMETERSIZEASSUMPTIONS}
	and\,\ref{SS:ASSUMPTIONSONDATA}
	and the bootstrap assumptions of 
	Sects.\,\ref{SS:BOOTSTRAPSCAFFOLDING} and\,\ref{SS:MAINQUANTITATIVEBOOTSTRAPASSUMPTIONS}.
	\item All quantities that we estimate can be controlled in terms of the variables 
	$\badcontrolvars = \{\wavearray,\upmu-1,\Lsmall^1,\Lsmall^2,\Lsmall^3\}$.
	\item We use the Leibniz rule for the operators $\angLie_Z$ and $\newangD$ when deriving pointwise estimates for the $\angLie_Z$ and $\newangD$ derivatives of tensor products of the schematic form $\prod_{i=1}^m \upxi_{(i)}$, where the $\upxi_{(i)}$ are scalar functions or $\ell_{t,u}$-tangent tensors. Our derivative counts are such that all the $\upxi_{(i)}$ except at most one are uniformly bounded in $L^{\infty}$ on $\characteristicdiamondtwoarg{[\leftubar,\ubarboot)}{[\moreinterestingu_1,\moreinterestingu_2]}$. Thus, our pointwise estimates often explicitly feature (on the right-hand sides) only one factor with many derivatives on it, multiplied by a constant that uniformly bounds the other factors. In some estimates, the right-hand sides also gain smallness factor such as $\fundbootsmall$, generated by the remaining 
	$\upxi_{(i)}'s$. 
	\item We use the conventions for constants $C, \mathfrak{c}$, and $C_{\mydiam}$
	stated in Sect.\,\ref{SS:CONVENTIONSFORCONSTANTS}.
	\item We use the conventions for strings of commutation vectorfields stated in Sect.\,\ref{S:COMMUTATIONFORMULASIDENTITIESANDSTRINGSOFCOMMUTATORS}.
	\item We use the comparison estimates of Lemma\,\ref{L:NORMOFSMOOTHTORITANGENTCOMMUTATORSANDSIMPLECOMPARISON}.
\end{enumerate}

\subsection{Pointwise estimates for Cartesian components of geometric vectorfields}
\label{SS:POINTWISEFORCARTESIANCOMPONENTSOFGEOMETRICVECTORFIELDS}
In this section, we provide simple pointwise estimates for the Cartesian components 
of the vectorfields $\{\Lunit,\muX, \Yvf{2},\Yvf{3}\}$ and their derivatives.

\begin{lemma}[Pointwise estimates for $x^i$ and the Cartesian components of the vectorfields 
	$\{\Lunit,\muX, \Yvf{2},\Yvf{3}\}$; \cite{abbrescia2022emergence}*{Lemma 14.2}] \label{L:POINTWISEESTIMATESFORCARTESIANCOMPONENTSOFVECTORFIELDS} 
	Let $\Singletan \in \{\Lunit,\Yvf{2},\Yvf{3}\}$. For $i = 1,2,3$,
	the following pointwise estimates hold on $\characteristicdiamondtwoarg{[\leftubar,\ubarboot)}{[\moreinterestingu_1,\moreinterestingu_2]}$:
	\begin{subequations} 
		\begin{align} 
			|\Singletan^i| 
			& \lesssim 1 + |\controlvars|, 
			\label{E:TANGENTIALCOMMUTATORCARTESIANCOMPONENTPOINTWISEESTIMATE} 
			\\
			|\tander^{[1,N]} \Singletan^i| 
			& \lesssim  
			| \tander^{[1,N]} \controlvars|, 
			\label{E:DERIVATIVESOFCARTESIANCOMPONENTSOFCOMMUTATORS} 
			\\
			|\comdersmall^{[1,N];1} \Singletan^i| 
			& \lesssim  
			|\comdersmall^{[1,N];1} \controlvars|,
			\label{E:COMDERSMALLPCARTESIANCOMPEST} 
			\\ 
			|\comder^{[1,N];1} \Singletan^i| 
			& \lesssim  
			|\comder^{[1,N];1} \controlvars|, 
			\label{E:COMDERPCARTESIANCOMPEST} 
			\\ 
			|\muX^i| 
			& \lesssim 
			1 
			+ 
			|\badcontrolvars|, 
			\label{E:MUXCARTESIANCOMPEST} 
			\\
			|\tander^{[1,N]} \muX^i| 
			& \lesssim 
			|\tander^{[1,N]} \badcontrolvars|, 
			\label{E:DERMUXCARTESIANCOMPEST} 
			\\
			|\comdersmall^{[1,N];1} \muX^i| 
			& \lesssim 
			|\comdersmall^{[1,N];1} \badcontrolvars|, 
			\label{E:COMDERSMALLMUXCARTESIANCOMPEST}
			\\
			|\comder^{[1,N];1} \muX^i| 
			& \lesssim 
			|\comder^{[1,N];1} \badcontrolvars|, 
			\label{E:COMDERMUXCARTESIANCOMPEST}
			\\
			|\angrmD x^i|_{\gtorus} 
			& \lesssim 
			1 
			+ 
			|\controlvars|, 
			\label{E:ANGDCARTESIANCOORDINATEIPOINTWISESTIMATE} 
			\\
			|\angrmD \tander^{[1,N]} x^i |_{\gtorus} 
			& \lesssim 
			|\tander^{[1,N]} \controlvars|, 
			\label{E:ANGDTANDERXIEST} 
			\\
			|\angrmD \comdersmall^{[1,N];1} x^i|_{\gtorus} 
			& \lesssim  
			|\comdersmall^{[1,N];1} \controlvars| 
			+ 
			|\tandersmall^{[1,N]} \badcontrolvars|, 
			\label{E:ANGDCOMDERSMALLXIEST} 
			\\
			|\angrmD \comdersmall^{[1,N];1} x^i|_{\gtorus} 
			& \lesssim  
			| \tandersmall^{[1,N]} \controlvars  
			+ 
			|\tandersmall^{[1,N]} \badcontrolvars |. 
			\label{E:ANGDCOMDERxiEST}
		\end{align}
	\end{subequations}
\end{lemma}


\subsection{Pointwise estimates for various $\ell_{t,u}$-tangent tensorfields} 
In this section, we record several pointwise estimates of various $\ell_{t,u}$-tangent tensorfields. 

\begin{lemma}[Crude pointwise estimates for the Lie derivatives of $\gtorus$ and $\upchi$] 
	\label{L:CRUDEPOINTWISEESTIMATESFORTENSORFIELDS} 
	The following pointwise estimates hold on $\characteristicdiamondtwoarg{[\leftubar,\ubarboot)}{[\moreinterestingu_1,\moreinterestingu_2]}$:
	\begin{subequations} 
		\begin{align}  
			|\angLie_{\tander}^{N+1} \gtorus|_{\gtorus}, 
			\, 
			|\angLie_{\tander}^{N+1} \gtorus^{-1}|_{\gtorus},
			\, 
			|\angLie_{\tander}^N \upchi|_{\gtorus}, 
			\, 
			|\tander^N \mytr_{\gtorus}\upchi| 
			& \lesssim  
			|\tander^{[1,N+1]}\controlvars|,  
			\label{E:TANDERGANDCHIESTIMATE} 
			\\
			|\angLie_{\comdersmall}^{N+1;1} \gtorus|_{\gtorus},
			\,  
			|\angLie_{\comdersmall}^{N+1;1} \gtorus^{-1}|_{\gtorus},
			\,  
			|\angLie_{\comder}^{N;1} \upchi|_{\gtorus}, 
			\, 
			|\comder^{N;1} \mytr_{\gtorus}\upchi| 
			& \lesssim  
			|\comdersmall^{[1,N+1];1}\controlvars| 
			+ 
			|\tandersmall^{[1,N+1]} \badcontrolvars|, 
			\label{E:COMDERSMALLGANDCHIESTIMATE} 
			\\
			|\angLie_{\comder}^{N+1;1} \gtorus|_{\gtorus},
			\, 
			|\angLie_{\comder}^{N+1;1} \gtorus^{-1}|_{\gtorus} 
			& \lesssim  
			|\comder^{[1,N+1];1}\controlvars| 
			+ 
			|\tandersmall^{[1,N+1]} \badcontrolvars|. 
			\label{E:COMDERGESTIMATE}
		\end{align}
	\end{subequations}
\end{lemma}
\begin{proof}
	The same proof of \cite[Lemma 8.5]{jLjS2018} holds with minor modifications to account for the third spatial dimension.
\end{proof}

\subsection{Differential operator estimates comparing $\angLie$ and $\newangD$}
In this section, we provide several pointwise estimates comparing different differential operators. 

\begin{lemma}[Differential operator pointwise comparison estimates; \cite{abbrescia2022emergence}*{Lemma 14.6}]
	Let $\varphi$ be a scalar function on $\characteristicdiamondtwoarg{[\leftubar,\ubarboot)}{[\moreinterestingu_1,\moreinterestingu_2]}$. 
	Then the following pointwise estimates hold:
	\begin{subequations}
		\begin{align}
			|\angrmD \varphi|_{\gtorus}^2 
			& = \left\lbrace
			1 + \mathcal{O}_{\mydiam}(\mathring{\upalpha}^{1/2}) 
			\right\rbrace
			\sum_{A = 2}^3 |\Yvf{A} \varphi|^2, 
			\label{E:ANGDFPOINTWISEBOUNDEDBYCOMMUTATORVECTORFIELDS} 
			\\
			\left|\angD^2 \varphi\right|_{\gtorus}^2 
			& \leq 
			\left\{1 + C_{\mydiam} \mathring{\upalpha}^{1/2}\right\} 
			\sum_{A = 2}^3 |\angrmD \Yvf{A} \varphi|_{\gtorus}^2 
			+ 
			C \auxbootsmall |\angrmD \varphi|_{\gtorus}^2, \label{E:SMOOTHANGULARHESSIANOFFPOINTWISEBOUNDEDBYCOMMUTATORVECTORFIELDS} \\
			|\angLap \varphi|^2  
			& \leq 
			2(1 + C_{\mydiam} \mathring{\upalpha}^{1/2}) 
			\sum_{A = 2}^3 |\angrmD \Yvf{A} \varphi|_{\gtorus}^2 
			+ 
			C \auxbootsmall |\angrmD \varphi|_{\gtorus}^2. 
			\label{E:SMOOTHANGULARLAPLACIANPOINTWISEBOUNDEDBYCOMMUTATORVECTORFIELDS} 
		\end{align}
	\end{subequations}
\end{lemma}


\subsection{Commutator estimates}
\label{SS:COMMUTATORESTIMATES}

The following proposition provides various vectorfield commutator 
estimates that we use in our analysis.

\begin{proposition}[Pointwise commutator estimates involving $\angLie$ ; \cite{abbrescia2022emergence}*{Prop.\,14.5}] 
	\label{P:COMMUTATORESTIMATES}
	Let $1 \leq N \leq \Ntop$ be an integer, and let $\varphi$ be a scalar function. 
	For any $\Singletan \in \Tanset$ (see definition \eqref{E:COMMUTATIONVECTORFIELDS}),
	iterated commutators can be bounded pointwise as follows
	on $\characteristicdiamondtwoarg{[\leftubar,\ubarboot)}{[\moreinterestingu_1,\moreinterestingu_2]}$:
	\begin{subequations}
		\begin{align} \label{E:COMMUTATOROFTANGENTIALANDTANGENTIALCOMMUTATORS}
			\left|
			[\Singletan, \tander^N] \varphi
			\right|
			& \lesssim 
			\auxbootsmall \left|\tander^{[1,N]} \varphi \right| 
			+ 
			\underbrace{\sum_{\substack {N_1+N_2 \leq N+1 \\ N_1,\,N_2\leq N}} 
				\left|\tander^{[2,N_1]}\controlvars \right| \left|\tander^{[1,N_2]} \varphi \right|}_{\mbox{Absent if $N=1$}},
			\\
			\begin{split}
				\left|[\muX, \tander^N]\varphi \right|,
				\,
				\left|
				[\Singletan, \comder^{N;1}] \varphi
				\right|
				& \lesssim
				\left|\tander^{[1,N]} \varphi \right| 
				+ 
				\underbrace{\sum_{\substack {N_1+N_2 \leq N+1 \\ N_1,\,N_2\leq N}} 
					\left|\tander^{[2,N_1]} \badcontrolvars \right| \left|\tander^{[1,N_2]} \varphi \right|}_{\mbox{Absent if $N=1$}}
				\label{E:COMMUTATOROFMUXANDTANGENTIALCOMMUTATORS} 
				\\
				& \ \
				+ 
				\underbrace{\sum_{\substack {N_1+N_2 \leq N \\ N_1 \leq N-1}} 
					\left|\tander^{[1,N_1]} \muX \wavearray \right| \left|\tander^{[1,N_2]} \varphi \right|}_{\mbox{Absent if $N=1$}}.
			\end{split}
		\end{align}
	\end{subequations}
	
	In particular, 
	for any $\Singletan \in \{\Lunit,\Yvf{2},\Yvf{3}\}$,
	we have the following pointwise estimates
	on $\characteristicdiamondtwoarg{[\leftubar,\ubarboot)}{[\moreinterestingu_1,\moreinterestingu_2]}$:
	\begin{subequations}
		\begin{align} \label{E:POINTWISEBOUNDCOMMUTATORSTANGENTIALANDTANGENTIALDERIVATIVESONSCALARFUNCTION}
			\left|[\Singletan, \tander^N] \varphi \right| 
			& 
			\lesssim 
			\auxbootsmall \left| \tander^{[1,N]}\varphi \right|,
			&
			\mbox{if } 1 \leq N \leq \Ntop - 11,
			\\
			\left|[\muX, \tander^N]\varphi \right|,
			\,
			\left|
			[\Singletan, \comder^{N;1}] \varphi
			\right|
			& 
			\lesssim \left|\tander^{[1,N]} \varphi \right|,
			& 
			\mbox{if } 1 \leq N \leq \Ntop - 12.
			\label{E:POINTWISEBOUNDCOMMUTATORSMUXANDTANGENTIALDERIVATIVESONSCALARFUNCTION}
		\end{align}
	\end{subequations}
	
	Moreover, the following pointwise estimates hold:
	\begin{subequations}
		\begin{align} 
			\left|
			[\Singletan, \comder^{N;2}] \varphi
			\right|,
			\,
			\left|
			[\muX, \comder^{N;1}] \varphi 
			\right|
			& 
			\lesssim 
			\left| \comdersmall^{[1,N];1} \varphi \right|,
			&
			\mbox{if } 1 \leq N \leq 5,
			\label{E:POINTWISEBOUNDCOMMUTATORSUPTOTWOMUXANDTANGENTIALDERIVATIVESONSCALARFUNCTION}
			\\
			\left|
			[\Singletan, \comder^{N;3}] \varphi
			\right|,
			\,
			\left|
			[\muX, \comder^{N;2}] \varphi 
			\right|
			& 
			\lesssim 
			\left| \comdersmall^{[1,N];2} \varphi \right|,
			&
			\mbox{if } 1 \leq N \leq 4.
			\label{E:POINTWISEBOUNDCOMMUTATORSUPTOTHREEMUXANDTANGENTIALDERIVATIVESONSCALARFUNCTION}
		\end{align}
	\end{subequations}
	
	Finally, if $\upxi$ is an $\ell_{t,u}$-tangent one-form or an $\ell_{t,u}$-tangent type $\binom{0}{2}$-tensorfield 
	and $\Singletan \in \{\Lunit,\Yvf{2},\Yvf{3}\}$,
	then the following estimates hold:
	\begin{subequations}
\begin{align} \label{E:LANDANGLIECOMMUTATOR}
\left|
	[\angLie_{\Singletan}, \angLie_{\tander}^N] \upxi 
\right|_{\gtorus} 
&
\lesssim 
\auxbootsmall \left| \angLie_{\tander}^{\leq N} \upxi \right|_{\gtorus}  
+ 
\underbrace{\sum_{\substack{N_1+N_2 \leq N+1 \\ N_2\leq N}} 
\left|\tandersmall^{[2,N_1]} \controlvars \right| 
\left|\angLie_{\tander}^{[1,N_2]}  \upxi \right|_{\gtorus}}_{\mbox{Absent if $N=1$}},
	\\
\begin{split} \label{E:LANDANGLIEWITHUPTOONEMUXCOMMUTATOR}
\left|
	[\angLie_{\Singletan}, \angLie_{\comder}^{N;1}] \upxi 
\right|_{\gtorus} 
&
\lesssim 
\left| \angLie_{\tander}^{\leq N} \upxi \right|_{\gtorus}  
+ 
\underbrace{\sum_{\substack{N_1+N_2 \leq N+1 \\ N_2\leq N}} 
\left|\tandersmall^{[2,N_1]} \badcontrolvars \right| 
\left|\angLie_{\tander}^{[1,N_2]}  \upxi \right|_{\gtorus}}_{\mbox{Absent if $N=1$}}
	\\
& \ \
	+
	\underbrace{\sum_{\substack {N_1+N_2 \leq N \\ N_1 \leq N-1}} 
	\left|\tander^{[1,N_1]} \muX \wavearray \right| 
	\left|\angLie_{\tander}^{[1,N_2]}  \upxi \right|_{\gtorus}}_{\mbox{Absent if $N=1$}}.
\end{split}
\end{align}
\end{subequations}
		
\end{proposition}

\subsection{Differential operator estimates involving the double-null acoustic geometry} \label{SS:POINTWISEESTIMATESTIEDTOTHEDOUBLENULLACOUSTICGEOMETRY}

In this section, we provide estimates which relate derivatives of $\wavearray$ relative to the double-null frame $\{\newuL,\newL,\nullgeop{x^2},\nullgeop{x^3}\}$ to those of the smooth frame $\{\muX, \Lunit,\geop{x^2},\geop{x^3}\}$. These estimates are necessary because the former are used to derive coercive energy--null-flux estimates, while the latter are present manifestly present as revealed in the last section.

\begin{lemma}[Double-null differential operator pointwise estimates] \label{L:POINTWISEESTIMATESTIEDTOTHEDOUBLENULLACOUSTICGEOMETRY}
The following estimates hold on $\characteristicdiamondtwoarg{[\leftubar,\ubarboot)}{[\moreinterestingu_1,\moreinterestingu_2]}$. 
	For $\muX$ and any scalar function $\varphi$, we have:
		\begin{align}
			|\muX \varphi| & \le |\newuL \varphi| + C \upmu |L \varphi| + C \auxbootsmall \upmu |\angrmd \varphi|_{\gtorus}. \label{E:POINTWISESTIMATEFORBREVEXINTERMSOFNEWUL}
		\end{align}
	For any scalar function, $\varphi$, let $\angrmd \varphi$ be the $\ell_{t,u}$-tangent one-form from \eqref{E:ANGULARDIFFERENTIAL} and let $\nullangrmD \varphi$ be the $\doublenulltoritwoarg{\ubar}{u}$-tangent one-form given by Def.\,\ref{D:DOUBLENULLORTHOGONALPROJECTIONANDTANGENCY}. Then we have: 
		\begin{subequations}
			\begin{align}
				|\angrmd \varphi |_{\gtorus} & \le (1 + \mathcal{O}(\auxbootsmall)) |\nullangrmD \varphi |_{\gnulltori} + \mathcal{O}(\auxbootsmall)|L\varphi|, \label{E:POINTWISESMOOTHTORIDIFFERENTIALNORMINTERMSOFDOULBENULLTORIDIFFERENTIALINORMANDL}  \\
				|\nullangrmD \varphi |_{\gnulltori} & \le (1 + \mathcal{O}(\auxbootsmall))|\angrmd \varphi |_{\gtorus} + \mathcal{O}(\auxbootsmall)|L\varphi| \label{E:POINTWISEDOULBENULLTORIDIFFERENTIALINORMINTERMSOFSMOOTHTORIDIFFERENTIALNORMANDL}
			\end{align} 
		\end{subequations}
		
	In addition, we have the following for $A = 2,3$:
		\begin{align} \label{E:NULLCOORDINATEPARTIALINTERMSOFMETRICNORM}
			\left|\nullgeop{x^A} \varphi\right| \le \left(1 + \mathcal{O}(\mr \upalpha^{1/2}) + \mathcal{O}(\auxbootsmall)\right) \left|\nullangrmD \varphi \right|_{\gnulltori}.
		\end{align}
		
	Moreover, for an arbitrary string of commutator vectorfields, we have the following two estimates for $\comdersmall^{[1,N+1];1} \varphi$:
		\begin{subequations}
			\begin{align}
				\begin{split}\label{E:ARBITRARYZCOMMUTATORCOMMUTEDTOHAVEBREVEXLAST}
					\left|\comdersmall^{[1,N+1];1} \varphi\right|  & \lesssim \left|\muX \tander^{[1,N]} \varphi\right| + \left|\tander^{[1,N]} \varphi\right| + \auxbootsmall \left|\comdersmall^{[1,N];1} \controlvars \right| + \auxbootsmall \left| \tandersmall^{[1,N]} \badcontrolvars \right|.
				\end{split} \\
				\begin{split}  \label{E:ARBITRARYCOMMUTATORSTRINGESTIMATE}
					\left|\comdersmall^{[1,N+1];1} \varphi \right| &  \lesssim \left|\newuL \tander^{[1,N]} \varphi \right| + \upmu \left| \nullangrmD \tander^{[1,N]} \varphi \right|_{\gnulltori}  + \upmu \left| \tander^{[1,N+1]} \varphi \right|  \\
					& \ \ + \left|\tander^{[1,N]}\varphi \right| + \auxbootsmall \left|\comdersmall^{[1,N];1} \controlvars \right| + \auxbootsmall \left| \tandersmall^{[1,N]} \badcontrolvars \right|,
				\end{split} 
		\end{align}
	\end{subequations}\end{lemma}
\begin{proof}
	Estimate \eqref{E:POINTWISESTIMATEFORBREVEXINTERMSOFNEWUL} follows from \eqref{E:GEOP2TOCOMMUTATORS}--\eqref{E:GEOP3TOCOMMUTATORS}, \eqref{E:CONVENIENTIDENTITYFORULUNIT}, Prop\,\ref{P:SCHEMATICSTRUCTUREOFVARIOUSTENSORSINTERMSOFCONTROLVARS}, \eqref{E:MAGNITUDEOFINNERPRODUCTOFNEWLANDNEWULAPPROXIMATELYMU}--\eqref{E:RECIPROCALULUNITAPPLIEDTOTIMEFUNCTIONAPPROXIMATELYMU}, and \eqref{E:ANGDFPOINTWISEBOUNDEDBYCOMMUTATORVECTORFIELDS}.
	
	To prove \eqref{E:POINTWISESMOOTHTORIDIFFERENTIALNORMINTERMSOFDOULBENULLTORIDIFFERENTIALINORMANDL}, we use the identities $|\angrmd \varphi|_{\gtorus}^2 = (\gtorus^{-1})^{AB} \geop{x^A} \varphi \geop{x^B} \varphi$ for scalar functions, \eqref{E:SMOOTHANGULARDERIVATIVESINTERMSOFDOUBLENULLONESANDL}, \eqref{E:SMOOTHTORUSINVERSEFIRSTFUNDCOMPONENTSINTERMSOFDOUBLENULLTORUSINVERSEFIRSTFUNDCOMPONENTS}, estimates \eqref{E:CLOSEDVERSIONKEYJACOBIANDETERMINANTESTIMATECHOVGEOTODOUBLENULL} \eqref{E:SMALLC01ESTIMATESFORUBAR}, the bootstrap assumptions, the Cauchy-Schwarz inequality, and Young's inequality. Estimate \eqref{E:POINTWISEDOULBENULLTORIDIFFERENTIALINORMINTERMSOFSMOOTHTORIDIFFERENTIALNORMANDL} follows from similar arguments which we omit.
	
	To prove \eqref{E:NULLCOORDINATEPARTIALINTERMSOFMETRICNORM}, we use the definition of the $\gnulltori$-gradient and the $\gnulltori$-Cauchy--Schwarz inequality to see
	\[ \left| \nullgeop{x^A}\varphi\right| = \left| \nullangrmD \varphi \cdot \nullgeop{x^A}\right| = \left|\gnulltori\left( \nullangD \varphi, \nullgeop{x^A}\right)  \right| \le \left| \nullangD \varphi\right|_{\gnulltori} \cdot |\nullgeop{x^A}|_{\gnulltori}.\]
	Since $\left| \nullangD \varphi\right|_{\gnulltori} = \left| \nullangrmD \varphi\right|_{\gnulltori}$ and $|\nullgeop{x^A}|_{\gnulltori} = 1 + \mathcal{O}(\mr \upalpha^{1/2}) + \mathcal{O}(\auxbootsmall)$ (which is a consequence of \eqref{E:GNULLTORICOMPONENTS}), \eqref{E:NULLCOORDINATEPARTIALINTERMSOFMETRICNORM} follows.
	
	To prove \eqref{E:ARBITRARYZCOMMUTATORCOMMUTEDTOHAVEBREVEXLAST}, we repeatedly use the commutator estimates \eqref{E:COMMUTATOROFTANGENTIALANDTANGENTIALCOMMUTATORS}--\eqref{E:COMMUTATOROFMUXANDTANGENTIALCOMMUTATORS} and the bootstrap assumptions.

Using \eqref{E:POINTWISESTIMATEFORBREVEXINTERMSOFNEWUL}--\eqref{E:POINTWISESMOOTHTORIDIFFERENTIALNORMINTERMSOFDOULBENULLTORIDIFFERENTIALINORMANDL} and the now proven \eqref{E:ARBITRARYZCOMMUTATORCOMMUTEDTOHAVEBREVEXLAST}, estimate \eqref{E:ARBITRARYCOMMUTATORSTRINGESTIMATE}  follows.

\end{proof}


\subsection{Transport inequalities for the eikonal function quantities} \label{SS:TRANSPORTINEQUALITIES}
In this section, we provide transport inequalities satisfied by the eikonal function quantities 
$\upmu, \Lsmall^i, \upchi$, and $\mytr_{\gtorus}\upchi$.
We also provide pointwise estimates for the differentiated quantity $\angLie_{\Lunit} \angLie_{\tander}^{N-1} \upchi$. 
These estimates \emph{involve a loss of one order of differentiability} relative to $\wavearray$ 
in the sense that the right-hand sides of the transport equations that we use to derive the inequalities
depend on the first-order derivatives of $\wavearray$. 
In Sect.\,\ref{S:PRELIMINARYL2ESTIMATESFORBELOWTOPORDERDERIVATIVESOFACOUSTICGEOMETRYANDDERIVATIVELOSING}, 
we use the transport inequalities to derive below-top-order energy estimates for the eikonal function quantities.

\begin{proposition}[Transport inequalities for the eikonal function quantities]
	\label{P:POINTWISETRANSPORTINEQUALITIESFOREIKFUNCTIONQUANTITIES} 
	The following pointwise estimates hold on $\characteristicdiamondtwoarg{[\leftubar,\ubarboot)}{[\moreinterestingu_1,\moreinterestingu_2]}$:
	\begin{subequations}
		\begin{align}
			|\Lunit \upmu| 
			& 
			\lesssim 
			\left|\comder \velocityarray\right| + \left|(\vortrenormalized,\GradEnt)\right|, 
			&&	\label{E:LMUPOINTWISE} 
			\\
			|\Lunit \tandersmall^N \upmu|, 
			\,
			|\tandersmall^N \Lunit \upmu| 
			& 
			\lesssim  \left| \comdersmall^{[1,N+1];1} \velocityarray\right| + \upmu \left| \tander^{[1,N+1]}\velocityarray\right| +  \left|\tander^{\le N} (\vortrenormalized,\GradEnt)\right|,	&& \mbox{if } 1 \leq N \leq \Ntop, 			\label{E:LTANGENTIALMUPOINTWISE} \\
			& \ \ +	|\tander^{[1,N]}\controlvars| + \fundbootsmall^{1/2} \left| \tandersmall^{[1,N]} \badcontrolvars\right|  \notag
			\\
			|\Lunit \tander^N \Lsmall^i|, 
			\,
			|\tander^N \Lunit \Lsmall^i| 
			& 
			\lesssim 
 \left| \tander^{[1,N+1]} \velocityarray\right| +  \left|\tander^{\le N} (\vortrenormalized,\GradEnt)\right| 
			+ 
			\auxbootsmall 
			|\tander^{[1,N]}\controlvars|, 
			\label{E:LUNITTANGENTIALDERIVATIVESOFLUNITIPOINTWISE} 
			&& \mbox{if } 0 \leq N \leq \Ntop, \\
			|\Lunit \tander^{N-1} \mytr_{\gtorus} \upchi|, 
			\,
			|\tander^{N-1} \Lunit \mytr_{\gtorus}\upchi| 
			& 
			\lesssim 
 \left| \tander^{[1,N+1]} \velocityarray\right| +  \left|\tander^{\le N} (\vortrenormalized,\GradEnt)\right|
 			+ 
			\auxbootsmall 
			|\tander^{[1,N]}\controlvars|,  
			\label{E:LTANGENTIALTRCHIPOINTWISE} 
			&& \mbox{if } 1 \leq N \leq \Ntop, \\
			|\angLie_{\Lunit} \angLie_{\tander}^{N-1} \upchi|_{\gtorus}, 
			\,
			|\angLie_{\tander}^{N-1}\angLie_{\Lunit} \upchi|_{\gtorus} 
			& 
			\lesssim
			 \left| \tander^{[1,N+1]} \velocityarray\right| +  \left|\tander^{\le N} (\vortrenormalized,\GradEnt)\right| 
 			+ 
			\auxbootsmall |\tander^{[1,N]}\controlvars|,
			&& \mbox{if } 1 \leq N \leq \Ntop,
			\label{E:ANGLIELTANGENTIALCHIPOINTWISE} 
			\\
			|\Lunit \comder^{N;1} \Lsmall^i|, 
			\,
			|\comder^{N;1} \Lunit \Lsmall^i| 
			& 
			\lesssim 
			 \left| \comdersmall^{[1,N+1];1} \velocityarray\right| + \upmu \left| \tander^{[1,N+1]}\velocityarray\right| +  \left|\tander^{\le N} (\vortrenormalized,\GradEnt)\right|,	&&\mbox{if } 1 \leq N \leq \Ntop.,	\label{E:LZLSMALLPOINTWISE}	\\
			& \ \ +	|\tander^{[1,N]}\controlvars| + \fundbootsmall^{1/2} \left| \tandersmall^{[1,N]} \badcontrolvars\right|  \notag 
			\\
			|\Lunit \comder^{N-1;1} \mytr_{\gtorus} \upchi|, 
			\,
			|\comder^{N-1;1} \Lunit \mytr_{\gtorus}\upchi| 
			& 
			\lesssim 
			 \left| \comdersmall^{[1,N+1];1} \velocityarray\right| + \upmu \left| \tander^{[1,N+1]}\velocityarray\right| +  \left|\tander^{\le N} (\vortrenormalized,\GradEnt)\right|,	&&\mbox{if } 2 \leq N \leq \Ntop.,	\label{E:LZTRCHIPOINTWISE}	\\
			& \ \ +	|\tander^{[1,N]}\controlvars| + \fundbootsmall^{1/2} \left| \tandersmall^{[1,N]} \badcontrolvars\right|  \notag 
			\\
			|\angLie_{\Lunit} \angLie_{\comder}^{N-1;1} \upchi|_{\gtorus}, 
			\,
			|\angLie_{\comder}^{N-1;1}\angLie_{\Lunit} \upchi|_{\gtorus} 
			& 
			\lesssim 
			 \left| \comdersmall^{[1,N+1];1} \velocityarray\right| + \upmu \left| \tander^{[1,N+1]}\velocityarray\right| +  \left|\tander^{\le N} (\vortrenormalized,\GradEnt)\right|,	&&\mbox{if } 2 \leq N \leq \Ntop.,	\label{E:ANGLIELZCHIPOINTWISE}	\\
			& \ \ +	|\tander^{[1,N]}\controlvars| + \fundbootsmall^{1/2} \left| \tandersmall^{[1,N]} \badcontrolvars\right|  \notag
		\end{align}
	\end{subequations}
\end{proposition}

\begin{proof}
The proposition was proved in \cite{jLjS2018}*{Prop.\,8.13} and \cite{abbrescia2022emergence}*{Prop.\,14.7}] with the RHS of the inequalities featuring $|\tander^{[1,N+1]}\wavearray|$ and $|\comdersmall^{N+1;1} \wavearray|$. The desired estimates follow from using \eqref{E:HIGHERORDERTANGENTIALDERIVATIVESOFDENSITY}--\eqref{E:HIGHERORDERMIXEDDERIVATIVESOFDENSITY} to bound the density in $\wavearray = (\LogDensity,v^1,v^2,v^3,s)$, as well as using the straightforward identity $\comder s = \smoothfunction(\controlvars,\badcontrolvars)  S$. 

\end{proof}

\subsection{Pointwise commutator estimates for $\upchi$ tied to a Codazzi-type identity}
\label{SS:CODAZZICOMMUTATOR} 

\begin{lemma}[Codazzi-type identity for $\upchi$;\cite{abbrescia2022emergence}*{Lemma\,14.8}]
	\label{L:CODAZZITYPEIDENTITY}
	There exist smooth functions, all schematically denoted by ``$\smoothfunction$,''
	such that the following identity holds:
	\begin{align} \label{E:CODAZZITYPEIDENTITY}
		\angdiv \upchi
		-
		\angrmD \mytr_{\gtorus} \upchi
		& 
		= 
		\smoothfunction(\tander^{\leq 1} \controlvars,\angrmD \vec{x}) \tander \controlvars
		+
		\smoothfunction(\controlvars,\angrmD \vec{x}) \tander^2 \wavearray.
	\end{align}
\end{lemma}

\begin{lemma}[Pointwise commutator estimates for $\upchi$ tied to a Codazzi-type identity; \cite{abbrescia2022emergence}*{Lemma\,14.9}]
	\label{L:CODAZZICOMMUTATORESTIMATES}
	Let $1 \leq N \leq \Ntop$. Then the following pointwise estimates 
	hold on $\characteristicdiamondtwoarg{[\leftubar,\ubarboot)}{[\moreinterestingu_1,\moreinterestingu_2]}$,
	where on LHS~\eqref{E:CODAZZICOMMUTATORESTIMATES},
	$\tander^{N-1}$ denotes the same order $N-1$ string of commutation vectorfields
	in each of the two terms: 
	\begin{align} \label{E:CODAZZICOMMUTATORESTIMATES} 
		\left| 
		\angdiv \angLie_{\tander}^{N-1} \upchi 
		- 
		\angrmD \tander^{N-1} \mytr_{\gtorus} \upchi 
		\right|_{\gtorus} 
		& 
		\lesssim \left| \tander^{[1,N+1]} \velocityarray\right| +  \left|\tander^{\le N} (\vortrenormalized,\GradEnt)\right| + \fundbootsmall \left| \tander^{[1,N]} \controlvars\right|  
	\end{align}
	
\end{lemma}
\begin{proof}
The lemma was proved in  \cite{abbrescia2022emergence}*{Lemma\,14.9} with the RHS of the inequalities featuring $|\tander^{[1,N+1]}\wavearray|$. The desired estimate follows from using \eqref{E:HIGHERORDERTANGENTIALDERIVATIVESOFDENSITY} to bound the density in $\wavearray = (\LogDensity,v^1,v^2,v^3,s)$, as well as using the straightforward identity $\tander s = \smoothfunction(\controlvars)  S$.

\end{proof}

\subsection{Pointwise estimates for the inhomogeneous terms in the commuted equations}
\label{SS:POINTWISEESTIMATESFORINHOMOGENEOUSTERMSINCOMMUTEDEQUATIONS}

\subsubsection{Pointwise estimates for the derivatives of the null forms}
\label{SSS:POINTWISEESTIMATESFORDERIVATIVESOFNULLFORMS}

\begin{lemma}[Pointwise estimates for the derivatives of the null forms]
	\label{L:POINTWISEESTIMATESFORDERIVATIVESOFNULLFORMS}
	Let $N \leq \Ntop$. The $\tander^N$-derivatives of the product of
	$\upmu$ and the terms defined in \eqref{E:TRANSPORTVORTVORTMAINTERMS}--\eqref{E:DIVENTROPYGRADIENTNULLFORM} 
	satisfy the following pointwise estimates
	on $\characteristicdiamondtwoarg{[\leftubar,\ubarboot)}{[\moreinterestingu_1,\moreinterestingu_2]}$:
	\begin{subequations}
		\begin{align} 
			\begin{split} \label{E:POINTWISEESIMATESFORDERIVATIVESOFNULLFORMSTRUCTUREMODIFIEDFLUIDVARIABLES}
				\left|\tander^N (\upmu \mainnullform_{(\VortVort)}^i) \right|, 
				\,
				\left|\tander^N (\upmu \mainnullform_{(\DivGradEnt)}) \right| 
				&
				\lesssim
				\left|
				\tander^{\leq N+1} (\vortrenormalized,\GradEnt)
				\right|
				\\
				& \ \
				+
				\fundbootsmall
				\left|
				\muX \tander^{[1,N]} \velocityarray
				\right| 
				+  
				\fundbootsmall
				\left|\tander^{[1,N+1]} \velocityarray \right| 
				+ 
				\fundbootsmall
				\left|\tandersmall^{[1,N]} \badcontrolvars \right|,
			\end{split}
			\\
			\left|\tander^N (\upmu \nullform_{(v)}^i) \right|
			& 
			\lesssim 
			\fundbootsmall
			\left|
			\muX \tander^{[1,N]} \velocityarray
			\right| 
			+  
			\left|
			\tander^{[1,N+1]} \velocityarray 
			\right| 
			+
			\left|\tander^{\le N} (\vortrenormalized,\GradEnt)\right| 
			+ 
			\fundbootsmall
			\left|
			\tandersmall^{[1,N]} \badcontrolvars 
			\right|,
			\label{E:POINTWISEESIMATESFORDERIVATIVESOFNULLFORMSTRUCTUREWAVEVARIABLES} 
			\\
			\begin{split} \label{E:POINTWISEESIMATESFORDERIVATIVESOFNULLFORMSTRUCTURETRANSPORTVARIABLES}
				\left| 
				\tander^N (\upmu \nullform_{(\VortVort)}^i) 
				\right|, 
				\, 
				\left|
				\tander^N (\upmu \nullform_{(\DivGradEnt)}) 
				\right|
				& 
				\lesssim 
				\fundbootsmall
				\left|
				\muX \tander^{[1,N]} \velocityarray
				\right| 
				+  
				\fundbootsmall
				\left|
				\tander^{[1,N+1]} \velocityarray 
				\right| 
				+ 
				\fundbootsmall
				\left|
				\tander^{\leq N} (\vortrenormalized,\GradEnt)
				\right| 
				+				
				\fundbootsmall
				\left|
				\tandersmall^{[1,N]} \badcontrolvars 
				\right|.
			\end{split}
		\end{align}
	\end{subequations}

\end{lemma}

\begin{proof}

The lemma was proved in  \cite{abbrescia2022emergence}*{Lemma\,14.10} with the RHS of the inequalities featuring $|\tander^{[1,N+1]}\wavearray|$ and $|\muX \tander^{[1,N]}\wavearray|$. The desired estimate follows from using \eqref{E:HIGHERORDERTANGENTIALDERIVATIVESOFDENSITY}--\eqref{E:HIGHERORDERMIXEDDERIVATIVESOFDENSITY} to bound the density in $\wavearray = (\LogDensity,v^1,v^2,v^3,s)$, as well as using the straightforward identity $\comder s = \smoothfunction(\controlvars,\badcontrolvars)  S$.

\end{proof}

\subsubsection{Pointwise estimates for the derivatives of the linear inhomogeneous terms}
\label{SSS:POINTWISEESTIMATESFORDERIVATIVESOFLINEARINHOMOGENEOUSTERMS}

\begin{lemma}[Pointwise estimates for the derivatives of the linear inhomogeneous terms]
\label{L:POINTWISEESTIMATESFORDERIVATIVESOFLINEARINHOMOGENEOUSTERMS}
Let $N \leq \Ntop$ and let $\tander^N \in \mathfrak{P}^{(N)}$. Consider the product of $\upmu$ and the terms defined in \eqref{E:VELOCITYILINEARORBETTER}, \eqref{E:SPECIFICVORTICITYLINEARORBETTER}--\eqref{E:RENORMALIZEDVORTICITYCURLLINEARORBETTER}. Then the following pointwise estimate holds on $\characteristicdiamondtwoarg{[\leftubar,\ubarboot)}{[\moreinterestingu_1,\moreinterestingu_2]}$:
	\begin{align} \label{E:POINTWISEESTIMATESFORALLLINEARTERMS}
		\left| \tander^N \left( \upmu \mathfrak{L}_{(v)}^i, \upmu \mathfrak{L}_{(\vortrenormalized)}^i,
				\,
				\upmu \mathfrak{L}_{(\GradEnt)}^i,
				\,
				\upmu \mathfrak{L}_{(\Flatdiv \vortrenormalized)},
				\,
				\upmu \mathfrak{L}_{(\VortVort)}^i
				\right) \right| & \lesssim \left| \tander^{\le N} (\vortrenormalized,\GradEnt)\right| + \fundbootsmall \left|\muX \tander^{[1,N]} \velocityarray \right| + \fundbootsmall \left| \tander^{[1,N+1]}\velocityarray\right| + \fundbootsmall \left| \tandersmall^{[1,N]} \badcontrolvars\right|.
	\end{align}

\end{lemma}

\begin{proof}
The lemma was proved in  \cite{abbrescia2022emergence}*{Lemma\,14.11} with the RHS of the inequalities featuring $|\tander^{[1,N+1]}\wavearray|$ and $|\muX \tander^{[1,N]}\wavearray|$. The desired estimate follows from using \eqref{E:HIGHERORDERTANGENTIALDERIVATIVESOFDENSITY}--\eqref{E:HIGHERORDERMIXEDDERIVATIVESOFDENSITY} to bound the density in $\wavearray = (\LogDensity,v^1,v^2,v^3,s)$, as well as using the straightforward identity $\comder s = \smoothfunction(\controlvars,\badcontrolvars)  S$.  

\end{proof}

\section{Trilinear structures in the error terms that drive the blow-up rate} \label{S:TRILINEARSTRUCTUREINTERMSTHATDRIVETHEBLOWUPRATE}

In this section we use the results of Sect.\,\ref{S:DENSITYDERIVATIVESINTERMSOFOTHERS} to reveal new trilinear structures for the crucial products driving the blow-up rates in the energy estimates.

\begin{proposition}[Trilinear structure in the error terms that drive the blow-up rate] \label{P:TRILINEARSTRUCTUREINTERMSTHATDRIVETHEBLOWUPRATE}
Let $\tander^N \in \mathfrak{P}^{(N)}$ denote any top order commutator of the form  $\tanderY^{\Ntop-1} \Lunit$ or $\tanderY^{\Ntop}$ (i.e., one with at least two factors of $\Lunit$ or a single $\Lunit$ that does not act first). Then the following estimates hold:
	\begin{align}
		\begin{split}
			&  \left| \frac{1}{2} (\muX v^1) \vec{G}_{LL}\diamond \muX \tander^{\Ntop} \wavearray - (\Lunit \upmu) \muX \tander^{\Ntop} v^1\right|  \\
			& \ \ \lesssim  \auxbootsmall \left| \comdersmall^{[1,\Ntop+1];1} \velocityarray\right| + \fundbootsmall \upmu \left| \tander^{[1,\Ntop+1]}\velocityarray\right|  + \upmu \left| \tander^{[1,\Ntop]}(\vortrenormalized,\GradEnt)\right|  \\
			& \ \ \ \ +  \left| \tander^{[1,\Ntop]}\controlvars\right| +  \left|\tandersmall^{[1,\Ntop]}\badcontrolvars\right| \label{E:TRILINEARSTRUCTUREFORMIXEDTOPORDERLMULTIPLIERTERM} 
		\end{split}
	\end{align}
	
In addition, the following estimate holds:
	\begin{align}
		\begin{split}
			& \left| \frac{1}{2} (\muX v^1) \vec{G}_{LL}\diamond \angLap\tanderY^{\Ntop-1} \wavearray - (\Lunit \upmu) \angLap \tanderY^{\Ntop-1}  v^1\right| \\
			& \ \ \lesssim \left| \Lunit \tander^{[1,\Ntop]}\velocityarray\right| +  \auxbootsmall \left| \tander^{[1,\Ntop+1]} \velocityarray\right|  + \left| \tander^{[1,\Ntop]}(\vortrenormalized,\GradEnt)\right| + \fundbootsmall \left|\tander^{[1,\Ntop]}\controlvars\right|. \label{E:TRILINEARSTRUCTUREFORPUREANGULARLMULTIPLIERTERM} 
		\end{split}
	\end{align}\end{proposition}

\begin{proof}

We first prove \eqref{E:TRILINEARSTRUCTUREFORMIXEDTOPORDERLMULTIPLIERTERM}. Using \eqref{E:MUTRANSPORT}, we have:
	\begin{align}
		\begin{split} \label{E:TRILINEARSTRUCTUREFORMIXEDTOPORDERLMULTIPLIERTERMSTEP1}
			\frac{1}{2} (\muX v^1) \vec{G}_{LL}\diamond \muX \tander^{\Ntop} \wavearray & = \frac{1}{2} (\muX v^1) G_{LL}^0 \muX \tander^{\Ntop}  \LogDensity + \frac{1}{2} \sum_{i=1}^3 (\muX v^1)G_{LL}^i \muX \tander^{\Ntop}  v^i + \frac{1}{2} (\muX v^1)G_{LL}^4 \muX \tander^{\Ntop}  \Ent \\
			& =  \frac{1}{2} G_{LL}^0\left(\muX v^1 \muX \tander^{\Ntop}  \LogDensity -( \muX \LogDensity )\muX \tander^{\Ntop}  v^1\right) + (\Lunit \upmu) \muX \tander^{\Ntop}   v^1 \\
			& \ \ - \frac{1}{2} \sum_{A=2}^3 (\muX v^A) G_{LL}^A \muX \tander^{\Ntop}  v^1 - \frac{1}{2} (\muX s) G_{LL}^4 \muX \tander^{\Ntop}  v^1 \\
			& \ \ + \sum_{A=2}^3 (\muX v^1)G_{LL}^A \muX \tander^{\Ntop} v^A + \frac{1}{2} (\muX v^1) G_{LL}^4 \muX \tander^{\Ntop}  \Ent
		\end{split}
	\end{align}
Note that bringing over the $(\Lunit \upmu) \muX  \tander^{\Ntop} v^1$ term on RHS\,\eqref{E:TRILINEARSTRUCTUREFORMIXEDTOPORDERLMULTIPLIERTERMSTEP1} over to the LHS yields the expression present in LHS \eqref{E:TRILINEARSTRUCTUREFORMIXEDTOPORDERLMULTIPLIERTERM}. Next, using the bootstrap assumptions (in particular \eqref{BA:AUXWAVEARRAYPARTIALLINFINITY}), \eqref{E:GLLAEXPRESSION}, Prop.\,\ref{P:SCHEMATICSTRUCTUREOFVARIOUSTENSORSINTERMSOFCONTROLVARS}, $\muX s = \upmu \Speed^{2} X_i S^i$, and \eqref{E:SMOOTHANGULARLAPLACIANPOINTWISEBOUNDEDBYCOMMUTATORVECTORFIELDS}, the second and third line on RHS\,\eqref{E:TRILINEARSTRUCTUREFORMIXEDTOPORDERLMULTIPLIERTERMSTEP1} are $\lesssim$ RHS\,\eqref{E:TRILINEARSTRUCTUREFORMIXEDTOPORDERLMULTIPLIERTERM}.

We now analyze the factor in parenthesis on RHS\,\eqref{E:TRILINEARSTRUCTUREFORPUREANGULARLMULTIPLIERTERMSTEP1}. Next, we use now \eqref{E:MUXRHOINTERMSOFVANDENT} to write:
	\begin{align}
		\begin{split} \label{E:TRILINEARSTRUCTUREFORMIXEDTOPORDERLMULTIPLIERTERMSTEP3}
			& (\muX v^1)\muX \tander^{\Ntop} \LogDensity - (\muX\LogDensity)\muX \tander^{\Ntop}v^1  = - \muX v^1 \left( X_i \muX \tander^{\Ntop} v^i\right) + (X_i \muX v^i)\muX \tander^{\Ntop} v^1 \\
			& \ \  + (\muX v^1) [\muX, \tander^{\Ntop}]\LogDensity - (\muX v^1) [X_i,\tander^{\Ntop}] \muX v^i - (\muX v^1)  X_i [\tander^{\Ntop},\muX]v^i \\
			& \ \ - (\muX v^1) \tander^{\Ntop} \left( \upmu X_i \Lunit v^i + \upmu  \exp(-\LogDensity)\frac{p_{;\Ent}}{\overline{\varrho}} \GradEnt^i X_i\right) + (\upmu X_i \Lunit v^i) \muX \tander^{\Ntop} v^1 + \upmu \exp(-\LogDensity) \frac{p_{;\Ent}}{\overline{\varrho}} (S^i X_i) \muX \tander^{\Ntop} v^1.
		\end{split}
	\end{align}
Using \eqref{E:HIGHERORDERTANGENTIALDERIVATIVESOFDENSITY}, commutator estimate \eqref{E:COMMUTATOROFMUXANDTANGENTIALCOMMUTATORS}, and the bootstrap assumptions, the second line on RHS\,\eqref{E:TRILINEARSTRUCTUREFORMIXEDTOPORDERLMULTIPLIERTERMSTEP3} is $\lesssim$ the $|\tander^{[1,\Ntop]}\controlvars| + \auxbootsmall |\tander^{[1,\Ntop]} \badcontrolvars|$  terms featured on RHS\,\eqref{E:TRILINEARSTRUCTUREFORMIXEDTOPORDERLMULTIPLIERTERM}.

We then observe a \emph{crucial} cancellation on the first line of RHS\,\eqref{E:TRILINEARSTRUCTUREFORMIXEDTOPORDERLMULTIPLIERTERMSTEP3} when $i=1$:
\begin{align} \label{E:TRILINEARSTRUCTUREFORMIXEDTOPORDERLMULTIPLIERTERMSTEP4}
	- \muX v^1 \left( X_i \muX \tander^{\Ntop} v^i\right) + (X_i \muX v^i)\muX \tander^{\Ntop} v^1 = - \muX v^1 X_A \muX \tander^{\Ntop} v^A + (X_A \muX v^A)\muX \tander^{\Ntop}  v^1. 
\end{align}
Since $X_A = \Speed^{-2} X^A = \Speed^{-2} \Xsmall^A$, we conclude the proof of \eqref{E:TRILINEARSTRUCTUREFORMIXEDTOPORDERLMULTIPLIERTERM} using \eqref{E:SCHEMATICSTRUCTUREOFXSMALL},  \eqref{E:TRILINEARSTRUCTUREFORPUREANGULARLMULTIPLIERTERMSTEP1}--\eqref{E:TRILINEARSTRUCTUREFORMIXEDTOPORDERLMULTIPLIERTERMSTEP4}.

We now prove \eqref{E:TRILINEARSTRUCTUREFORPUREANGULARLMULTIPLIERTERM}. First, using identical arguments used to prove \eqref{E:TRILINEARSTRUCTUREFORMIXEDTOPORDERLMULTIPLIERTERM} with $\muX \tander^{\Ntop}$ replaced by $\angLap \tanderY^{\Ntop-1}$, we have that: 
	\begin{align} \label{E:TRILINEARSTRUCTUREFORPUREANGULARLMULTIPLIERTERMSTEP1}
		\left| \frac{1}{2}( \muX v^1) \vec{G}_{\Lunit \Lunit} \diamond \angLap \tanderY^{\Ntop-1} \wavearray - (\Lunit \upmu) \angLap \tanderY^{\Ntop-1} v^1- \frac{1}{2} G_{\Lunit \Lunit}^0 \left((\muX v^1)\angLap \tanderY^{\Ntop-1} \LogDensity - (\muX\LogDensity)\angLap \tanderY^{\Ntop-1} v^1 \right)\right| \lesssim \text{RHS\,\eqref{E:TRILINEARSTRUCTUREFORMIXEDTOPORDERLMULTIPLIERTERM}}.
	\end{align}Setting $W = \Yvf{A}$ in \eqref{E:TANGENTIALDERIVATIVESOFRHOINTERMSOFVANDVORTANDENT} and then differentiating the resulting identity with respect to $\angLap \tanderY^{N-2}$, the same arguments used to prove \eqref{E:TRILINEARSTRUCTUREFORMIXEDTOPORDERLMULTIPLIERTERM} imply: 
\begin{align}
	\begin{split} \label{E:TRILINEARSTRUCTUREFORPUREANGULARLMULTIPLIERTERMSTEP2}
		\left| (\muX v^1) \angLap \tanderY^{\Ntop-1}  \LogDensity - (\muX \LogDensity) \angLap \tanderY^{\Ntop-1}  v^1 - \left\{ -  (\muX v^1)\left( X_i \angLap \tander^{\Ntop-1}  v^i\right) + \left( X_i \muX v^i\right) \angLap \tander^{\Ntop-1} v^1\right\} \right| \lesssim \text{RHS\,\eqref{E:TRILINEARSTRUCTUREFORPUREANGULARLMULTIPLIERTERM}}.
	\end{split}
\end{align}
We then observe a \emph{crucial} cancellation in the terms in the braces in \eqref{E:TRILINEARSTRUCTUREFORPUREANGULARLMULTIPLIERTERMSTEP2} when $i=1$: 
\begin{align} \label{E:TRILINEARSTRUCTUREFORPUREANGULARLMULTIPLIERTERMSTEP3}
	 -  (\muX v^1)\left( X_i \angLap \tander^{\Ntop-1}  v^i\right) + \left( X_i \muX v^i\right) \angLap \tander^{\Ntop-1} v^1 = -  (\muX v^1)\left( X_A \angLap \tander^{\Ntop-1}  v^A\right) + \left( X_A \muX v^A\right) \angLap \tander^{\Ntop-1} v^1\end{align}
Since $X_A = \Speed^{-2} X^A = \Speed^{-2} \Xsmall^A$, we conclude the proof of the second inequality in \eqref{E:TRILINEARSTRUCTUREFORPUREANGULARLMULTIPLIERTERM} using \eqref{E:SCHEMATICSTRUCTUREOFXSMALL} and the bootstrap assumptions.

\end{proof}

\section{Control of the flow map of  \texorpdfstring{$\newL$}{L}}
\label{S:ESTIMATESTIEDTOTHEFLOWMAPOFUBARNORMALIZEDNULLVECTORFIELD}
We continue to work under the assumptions of Sect.\,\ref{SS:SILENTFACTS}.
In this section, we derive various estimates tied to the
flow map of the vectorfield $\newL$.
We then use these results to derive preliminary
estimates for solutions $f$ to transport equations of the form
$\newL f = F$.

\subsection{Basic properties of the flow map of $\newL$ and estimates for solutions to $\newL f = F$}
\label{SS:PROPERTIESOFFLOWMAPOFNEWL}

\begin{lemma}[Basic properties of the flow map of $\newL$]
	\label{L:PROPERTIESOFFLOWMAPOFNEWL}
	Let $\newL$ be the null vectorfield defined in \eqref{E:DOUBLENULLEIKONALFUNCTIONNORMALIZEDNULLVECTORFIELDS},
	and let $\FlowmapnewL$ be the $\leftubar$-normalized flow map of
	$\newL$ with respect to the double-null coordinates
	$(\ubar,u,x^2,x^3)$, i.e., the solution to
	the following initial value problem:
	\begin{align} \label{E:FLOWMAPOFNEWLINDOUBLENULL}
		\nullgeop{\ubar}
		\FlowmapnewL(\ubar,u,x^2,x^3)
		& = 
		\newL \circ \FlowmapnewL(\ubar,u,x^2,x^3),
		&
		\FlowmapnewL(\leftubar,u,x^2,x^3)
		& = (\timefunction_0,u,x^2,x^3).
	\end{align}
	Then for $A = 2,3$ there exist functions 
	$\FlowmapnewLarg{A}
	: [\leftubar,\ubarboot] \times [\moreinterestingu_1,\moreinterestingu_2] \times \mathbb{T}^2
	\rightarrow \mathbb{T}$
	such that:
	\begin{align} \label{E:FORMOFNEWLGENERATORFLOWMAP}
	\FlowmapnewL(\ubar,u,x^2,x^3)
	& 
	=
	\left(
		\ubar,u,
		\FlowmapnewLarg{2}(\ubar,u,x^2,x^3),
		\FlowmapnewLarg{3}(\ubar,u,x^2,x^3)
	\right).
	\end{align}
	Moreover, $\FlowmapnewL$
	is a $C^{0,1}$ Lipeomorphism
	from $[\leftubar,\ubarboot] \times [\moreinterestingu_1,\moreinterestingu_2] \times \mathbb{T}^2$
	onto $[\leftubar,\ubarboot] \times [\moreinterestingu_1,\moreinterestingu_2] \times \mathbb{T}^2$
	satisfying:
	\begin{align} \label{E:C01DOUBLENULLBOUNDFORFLOWMAPOFNEWL}
		\left\|
			\FlowmapnewL
			-
			\mbox{\upshape I}
		\right\|_{C_{\textnormal{null}}^{0,1}([\leftubar,\ubarboot] \times [\moreinterestingu_1,\moreinterestingu_2] \times \mathbb{T}^2)}
		& \lesssim \auxbootsmall,
		& 
		\left\|
			\InverseFlowmapnewL
			-
			\mbox{\upshape I}
		\right\|_{C_{\textnormal{null}}^{0,1}([\leftubar,\ubarboot] \times [\moreinterestingu_1,\moreinterestingu_2] \times \mathbb{T}^2)}
		& \lesssim \auxbootsmall,
	\end{align}
	where $\mbox{\upshape I}(\ubar,u,x^2,x^3) \eqdef (\ubar,u,x^2,x^3)$
	is the identity map on double-null coordinate space and $\InverseFlowmapnewL$ 
	is the inverse function of $\FlowmapnewL$.
	In particular, for each fixed $(\ubar,u) \in[\leftubar,\ubarboot] \times [\moreinterestingu_1,\moreinterestingu_2]$,
	the map 
	$(x^2,x^3) 
	\mapsto 
	\left(\FlowmapnewLarg{2}(\ubar,u,x^2,x^3),\FlowmapnewLarg{3}(\ubar,u,x^2,x^3) \right)$
	is a $C^{0,1}$ Lipeomorphism from $\doublenulltoritwoarg{\leftubar}{u}$ onto $\doublenulltoritwoarg{\ubar}{u}$.
	
\end{lemma}

\begin{proof}
	\eqref{E:FORMOFNEWLGENERATORFLOWMAP} is a trivial consequence of \eqref{E:FLOWMAPOFNEWLINDOUBLENULL} and the identities $\newL \ubar =  1$ and $\newL u = 0$.
	
	To prove \eqref{E:C01DOUBLENULLBOUNDFORFLOWMAPOFNEWL},
	we first note that
	the functions $\left(\FlowmapnewLarg{2} - x^2,\FlowmapnewLarg{3} - x^3 \right)$
	solve the transport system:
	\begin{align} \label{E:NEWLFLOWMAPODETORUSCOMPONENTS}
	\nullgeop{\ubar}
	[\FlowmapnewLarg{A} - x^A](\ubar,u,x^2,x^3)
	& 
	= 
	\newL^A\left(\ubar,u,\FlowmapnewLarg{2},\FlowmapnewLarg{3}\right),
	&
	&(A=2,3)
	\end{align}
	with vanishing data at ingoing null time $\leftubar$.
	Next,
	we use definition \eqref{E:DOUBLENULLEIKONALFUNCTIONNORMALIZEDNULLVECTORFIELDS},
	Lemma\,\ref{L:RELATIONSHIPBETWEENGEOMETRICPARTIALDERIVATIVESANDDOUBLENULL},
	the bootstrap assumptions,
	Lemma\,\ref{L:PROPERTIESANDDIFFEOMORPHICEXTENSIONOFDOUBLENULLCOORDINATES},
	and Lemma\,\ref{L:CONTINUOUSEXTNESION}
	to deduce that:
	\begin{align} \label{E:NEWL2NEWL3DOUBLENULLC01BOUNDS}
	\left\| 
		\left(\newL^2,\newL^3 \right)
	\right\|_{C_{\textnormal{null}}^{0,1}([\leftubar,\ubarboot] \times [\moreinterestingu_1,\moreinterestingu_2] \times \mathbb{T}^2)}
	& 
	\lesssim \auxbootsmall.
	\end{align}
	Hence, commuting \eqref{E:NEWLFLOWMAPODETORUSCOMPONENTS} 
	up to one time with the adapted rough geometric coordinate partial derivatives,
	using \eqref{E:NEWL2NEWL3DOUBLENULLC01BOUNDS},
	and integrating with respect to ingoing null time,
	we find that for $\ubar \in [\leftubar,\ubarboot]$,
	we have:
	\begin{align} 
	\begin{split}\label{E:READYC01DOUBLENULLBOUNDFORFLOWMAPOFNEWL}
		\max_{A=2,3}
		\left\| 
			\FlowmapnewLarg{A} - x^A
		\right\|_{C_{\textnormal{null}}^{0,1}([\leftubar,\ubar] \times [\moreinterestingu_1,\moreinterestingu_2] \times \mathbb{T}^2)}
		& \leq C \auxbootsmall
			\\
		& \ \
			+
			C  \auxbootsmall
			\int_{\leftubar}^{\ubar}
				\max_{A=2,3}
				\left\| 
					\FlowmapnewLarg{A} - x^A
				\right\|_{C_{\textnormal{rough}}^{0,1}([\leftubar,\ubar'] \times [\moreinterestingu_1,\moreinterestingu_2] \times \mathbb{T}^2)}
			\, \rmd \ubar'.
	\end{split}
	\end{align}
	From \eqref{E:READYC01DOUBLENULLBOUNDFORFLOWMAPOFNEWL} and Gr\"{o}nwall's inequality,
	we find that
	$
	\max_{A=2,3}
		\left\| 
			\FlowmapnewLarg{A} - x^A
		\right\|_{C_{\textnormal{rough}}^{1,0}([\leftubar,\ubarboot] \times [\moreinterestingu_1,\moreinterestingu_2] \times \mathbb{T}^2)}
	\leq C \auxbootsmall
	$.
	From this bound and \eqref{E:FORMOFNEWLGENERATORFLOWMAP},
	we conclude the first bound stated in \eqref{E:C01DOUBLENULLBOUNDFORFLOWMAPOFNEWL}.
	From this bound and the inverse function theorem for Lipschitz functions, we conclude that
	$\FlowmapnewL$ is a $C_{\textnormal{null}}^{0,1}$ Lipeomorphism from
	$
	[\leftubar,\ubarboot] \times [\moreinterestingu_1,\moreinterestingu_2] \times \mathbb{T}^2
	$
	onto itself ($\FlowmapnewL$ is a global Lipeomorphism since it is close to the identity map).
	The second bound in \eqref{E:C01DOUBLENULLBOUNDFORFLOWMAPOFNEWL}
	follows from the identity
	$
	\InverseFlowmapnewL \circ \FlowmapnewL = \mbox{\upshape I}
	$.
\end{proof}
The following corollary is an immediate consequence of the preceding lemma and the Fundamental Theorem of Calculus.
\begin{corollary}[Pointwise identities and $L^\infty$ estimates tied to the integral curves of $\newL$] \label{C:TRANSPORTESTIMATESFORNEWLFEQUALSSOURCE}
Let $f$ be a function of the double-null coordinates on 
$[\leftubar,\ubarboot) \times [\moreinterestingu_1,\moreinterestingu_2] \times \mathbb{T}^2$,
let $\FlowmapnewL$ be the $\leftubar$-normalized flow map of $\newL$ from 
Lemma\,\ref{L:PROPERTIESOFFLOWMAPOFNEWL}.
Then relative to the double-null coordinates, the following identity holds for any 
$\leftubar \le \ubar_1 \le \ubar_2 < \ubarboot$
and any $(u,x^2,x^3) \in [\moreinterestingu_1,\moreinterestingu_2] \times \mathbb{T}^2$:
\begin{align} \label{E:TRANSPORTIDENTITYALONGLROUGHINTEGRALCURVES}
f \circ \FlowmapnewL(\ubar_2,u,x^2,x^3)
& = 
	f \circ \FlowmapnewL(\ubar_1,u,x^2,x^3)
	+ 
	\int_{\timefunction' = \ubar_1}^{\ubar_2} (\newL f) 
	\circ 
	\FlowmapnewL(\ubar',u,x^2,x^3) 
	\, \rmd \ubar'.
\end{align}
In addition, we have the following $L^{\infty}$ estimates:
\begin{subequations}
	\begin{align} 
		\| f \|_{L^{\infty}\left(\doublenulltoritwoarg{\ubar_2}{u}\right)} 
		& \leq 
		\| f\|_{L^{\infty}\left(\doublenulltoritwoarg{\ubar_1}{u}\right)} 
		+
		|\leftubar| \sup_{\ubar' \in [\ubar_1,\ubar_2]} 
		\| \newL f \|_{L^{\infty}\left(\doublenulltoritwoarg{\ubar'}{u}\right)},
		\label{E:LINFINITYTRANSPORTNEWLESTIMATE} 
	\end{align}
	\begin{align}
		\esssup_{(x^2,x^3) \in \T^2} \left| f \circ \FlowmapnewL(\ubar_2, u,x^2,x^3) -  f \circ \FlowmapnewL(\ubar_1, u,x^2,x^3) \right| \le |\leftubar| \sup_{\ubar' \in [\ubar_1,\ubar_2]} 
		\| \newL f \|_{L^{\infty}\left(\doublenulltoritwoarg{\ubar'}{u}\right)}.
		\label{E:ESSENTIALSUPPREMUMTRANSPORTNEWLLFESTIMATE} 
	\end{align}
\end{subequations}\end{corollary}

\subsection{Estimate for $\mydet \gnulltori$} 
\label{SS:DOUBLENULLTORIVOLUMEFORMCOMPARISONESTIMATE}
In the next lemma, we control the factor $\mydet \gnulltori(\ubar,u,x^2,x^3)$
featured in the area form $\voldoublenulltori$ (see \eqref{E:AREAFORMDOUBLENULLTORUS}) 
on the rough tori $\doublenulltoritwoarg{\ubar}{u}$.

\begin{lemma}[Estimate for $\mydet \gnulltori$] 
\label{L:DOUBLENULLTORIVOLUMEFORMCOMPARISONESTIMATE}
Recall that $\gnulltori$ is the first fundamental form of $\doublenulltoritwoarg{\ubar}{u}$
(see Def.\,\ref{D:FIRSTFUNDAMENTALFORMSOFDOUBLENULLTORI}) and that $\FlowmapnewL$ is the $\leftubar$-normalized flow map of
	$\newL$ with respect to the double-null coordinates
	$(\ubar,u,x^2,x^3)$ from Lemma\,\ref{L:PROPERTIESOFFLOWMAPOFNEWL}. 
Then for every $\ubar_1, \ubar_2 \in [\leftubar,\ubarboot]$
and every $(u,x^2,x^3) \in [\moreinterestingu_1,\moreinterestingu_2] \times \mathbb{T}^2$,
the following estimates hold,
where $\mydet \gnulltori$ is evaluated relative to the double-null coordinates
via the formula \eqref{E:GNULLTORICOMPONENTS}:
\begin{align} \label{E:ESTIMATEFORNULLTORIVOLUMEFORMDETERMINANT}
\mydet \gnulltori(\ubar_2,u,x^2,x^3)
& = 
\left\lbrace
	1 + \mathcal{O}(\auxbootsmall)
\right\rbrace
\mydet \gtorusroughfirstfund(\ubar_1,u,x^2,x^3)
=
1 + \mathcal{O}_{\mydiam}(\mathring{\upalpha}).
\end{align}\end{lemma}

\begin{proof}
First, using the expression of the double-null metric components 
$\gnulltori\left(\nullgeop{x^A},\nullgeop{x^B}\right)$ 
in \eqref{E:GNULLTORICOMPONENTS}, we compute:
\begin{align} \label{E:IDENTITYNEEDEDFORVOLFORMESTIMATEDOUBLENULLTORI}
\mydet \gnulltori 
& = 
\left(\frac{\Lunit \ubar}{\geop{t}\ubar}\right)^2 
\mydet \gtorus.
\end{align}
Since $\gtorus_{AB} = \gfour(\geop{x^A},\geop{x^B}) = \gfour(\p_A - \Yvfsmallcoeff{A}\p_1,\p_B - \Yvfsmallcoeff{B}\p_1)$, it follows from \eqref{E:SPLITMETRICINTOMINKOWSKIANDREMAINDERPART}--\eqref{E:METRICPERTURBATIONVANISHESATTRIVIALPSISOLUTION}, \eqref{E:DATAEPSILONISSMALLERTHANBOOTSTRAPEPSILONSMALLERTHANSQUAREOFDATAALPHA}, and Props.\,\ref{P:SCHEMATICSTRUCTUREOFVARIOUSTENSORSINTERMSOFCONTROLVARS} that $\mydet \gtorus = 1 +  \mathcal{O}_{\mydiam}(\mathring{\upalpha})$. 

Using \eqref{E:IDENTITYNEEDEDFORVOLFORMESTIMATEDOUBLENULLTORI},
our assumptions on the data from Sect.\,\ref{SSS:QUANTITATIVEASSUMPTIONSONDATAAWAYFROMSYMMETRY}, we deduce that
$\mydet \gnulltori(\leftubar,u,x^2,x^3)
=
1 + \mathcal{O}_{\mydiam}(\mathring{\upalpha})
$.
Moreover, using the identity $\Lunit \ubar = \geop{t} \ubar + L^A \geop{x^A} \ubar$, \eqref{E:NULLCOORDINATEPARTIALINGOINGEIKONALINTERMSOFGEOMETRICVF}, Props.\,\ref{P:SCHEMATICSTRUCTUREOFVARIOUSTENSORSINTERMSOFCONTROLVARS} and\,\ref{P:IMPROVEMENTOFAUXILIARYBOOTSTRAP}, \eqref{E:SMALLC01ESTIMATESFORUBAR}, we deduce that
$
\nullgeop{\ubar} \mydet \gnulltori
=
\mathcal{O}(\auxbootsmall)
$
which, in view of the mean value theorem,
yields that for any 
$(\ubar,u,x^2,x^3) \in [\leftubar,\ubarboot] \times [\moreinterestingu_1,\moreinterestingu_2] \times \mathbb{T}^2$,
we have
$
\mydet \gnulltori(\ubar,u,x^2,x^3)
=
\mydet \gnulltori(\leftubar,u,x^2,x^3)
+
\mathcal{O}(\auxbootsmall)
=
\left\lbrace
	1 + \mathcal{O}(\auxbootsmall)
\right\rbrace
\mydet \gnulltori(\leftubar,u,x^2,x^3)
$.
In total, these estimates imply \eqref{E:ESTIMATEFORNULLTORIVOLUMEFORMDETERMINANT}.

\end{proof}

\section{Control of the flow map of the ingoing null vectorfield}
\label{S:PROPERTIESOFTHEFLOWMAPOFNEWUL}
We continue to work under the assumptions of Sect.\,\ref{SS:SILENTFACTS}. In this short section, we derive useful properties of the flow map of $\newuL$. These properties will be used in Sect.\,\ref{S:SHARPCONTROLOFMUANDPROPERTIESOFCHOVGEOTOCARTESIAN} to construct a neighborhood of $\datahypfortimefunctiontwoarg{0}{[\leftubar,\ubarboot]}$ on which we have sharp control over $\Lunit \upmu$, see \eqref{E:SMALLNEIGHBORHOOD}--\eqref{E:OUTSIDEOFSMALLNEIGHBORHOODLOWERBOUNDFORMU}. We begin with the following lemma, which in particular, provides strict improvements of the bootstrap assumptions 

\begin{lemma}[Basic properties of the flow map of $\newuL$] \label{L:PROPERTIESOFFLOWMAPOFNEWUL}

Let $\newuL$ be the null vectorfield defined in \eqref{E:DOUBLENULLEIKONALFUNCTIONNORMALIZEDNULLVECTORFIELDS}, and recall that $\flowmapofnewularg{\Delta s}$ is the flow map of $\newuL$ (see Theorem\,\ref{T:CONSTRUCTIONOFTHEINGOINGEIKONALFUNCTIONCONTRACTIONMAPPING}). 

The map $\underline{F}$ from \eqref{E:BADEFINITIONOFFLOWOUTMAP} extends to a $C^{0,1}$-Lipeomorphism on the closure $\mbox{\upshape cl}\left(\underline{\mathscr{F}}\right)$. Moreover, $\mbox{\upshape cl}\left(\underline{\mathscr{F}}\right)$ is given by \eqref{E:BADEFINITIONOFFLOWOUTDOMAINFORUBAR} with $\domainforembeddingdatahypfortimefunctiontwoarg{0}{[-\ubarboot,-\leftubar]}$ in place of $\domainforembeddingdatahypfortimefunctiontwoarg{0}{(-\ubarboot,-\leftubar]}$, and $\underline{F}\left(\mbox{\upshape cl}\left(\underline{\mathscr{F}}\right)\right) = \characteristicdiamondtwoarg{[\leftubar,\ubarboot]}{[\moreinterestingu_1,\moreinterestingu_2]}$. In addition, the following estimates hold: 
		\begin{subequations}
			\begin{align}
				\left\| \underline{F} \right\|_{C^{0,1}\left(\mbox{\upshape cl}\left(\underline{\mathscr{F}}\right)\right)} & \le C, \label{E:FLOWMAPOFNEWULC01QUANTIFIED} \\
				\left\| \underline{F}^{-1} \right\|_{C^{0,1}_{\textnormal{geo}}\left( \characteristicdiamondtwoarg{[\leftubar,\ubarboot]}{[\moreinterestingu_1,\moreinterestingu_2]}\right)} & \le C, \label{E:INVERSEFLOWMAPOFNEWULC01QUANTIFIED}
			\end{align}
		\end{subequations}				

Moreover, for each $\ubar \in [\leftubar,\ubarboot]$ and each pair $u_1,u_2 \in [\moreinterestingu_1,\moreinterestingu_2]$, $\flowmapofnewularg{u_2 - u_1}$ is a $C^{0,1}_{\textnormal{geo}}$ Lipeomorphism from the double-null torus $\doublenulltoritwoarg{\ubar}{u_1}$ onto the double-null torus $\doublenulltoritwoarg{\ubar}{u_2}$. In particular, the integral curves of $\flowmapofnewularg{\Delta s}$ thread $\ingoingcharacteristicsurfacetwoarg{\ubar}{[\moreinterestingu_1,\moreinterestingu_2]}$. Moreover, for every $\ubar \in [\leftubar,\ubarboot]$, each integral curve of $\newuL$ passes through precisely one point on the $\newuL \upmu$-adapted torus $\argubarnewultorus{\ubar}$ defined in \eqref{E:DEFOFNEWULADAPTEDTORI}. Finally, for every  $\ubar \in [\leftubar,\ubarboot]$, each integral curve of $\newuL$ also passes through precisely one point on the $\upmu$-adapted torus $\twoargmumuxtorus{-\ubar}{0}$ defined in \eqref{E:MUXMUTORI}.

\end{lemma}

\begin{proof}

From \eqref{E:CLOSEDVERSIONC11BOUNDFORCHOVDOUBLENULLTOGEO}, it follows that $\newuL$ is a $C^{0,1}_{\textnormal{geo}}$ vectorfield on $\characteristicdiamondtwoarg{[\leftubar,\ubarboot]}{[\moreinterestingu_1,\moreinterestingu_2]}$. By the standard theory of flows for Lipschitz vectorfields, see e.g. \cite{rampazzo2007frobenius}, estimates \eqref{E:FLOWMAPOFNEWULC01QUANTIFIED}--\eqref{E:INVERSEFLOWMAPOFNEWULC01QUANTIFIED} hold.

That $\flowmapofnewularg{u_2- u_1} \colon \doublenulltoritwoarg{\ubar}{u_1} \to \doublenulltoritwoarg{\ubar}{u_2}$ is a $C^{0,1}_{\textnormal{geo}}$ Lipeomorhism is a consequence of $\newuL \ubar = 0, \, \newuL u = 1$ and the already proven properties of $\underline{F}$. 

From \eqref{E:BAUBARNEWULMUTORISTRUCTURE},
	\eqref{E:BAUBARNEWULMUTORILOCATION}, \eqref{E:LASTSLICELEVELSETSTRUCTUREANDLOCATIONOFMIN},
	and the facts that $\newuL \ubar = 0$
	and $\newuL u = 1$,
	it follows that for each fixed $\ubar \in [\leftubar,\ubarboot]$,
	every integral curve of
	$\newuL$
	in $ \ingoingcharacteristicsurfacetwoarg{\ubar}{[\moreinterestingu_1 ,\moreinterestingu_2]}$
	must intersect 
	$\argubarnewultorus{\ubar}$
	at one or more points in $\ingoingcharacteristicsurfacetwoarg{\ubar}{[\moreinterestingu_1, \frac{3}{4}\moreinterestingu_2 + \frac{1}{2} \interestingu]}$. Recalling that 
	$\newuL \upmu|_{\argubarnewultorus{\ubar}} = 0$,
	and using the transversal convexity bootstrap assumption \eqref{E:BAMUTRANSVERSALCONVEXITY} for
	$\newuL\newuL \upmu$, which now holds all the way up to the closure $\characteristicdiamondtwoarg{[\leftubar,\ubarboot]}{[\moreinterestingu_1,\moreinterestingu_2]}$ thanks to Lemma\,\ref{L:CONTINUOUSEXTNESION},
	we see that the intersection occurs at a unique point. 
	
The statement for the $(\upmu,\muX \upmu)$-adapted torus  $\twoargmumuxtorus{-\ubar}{0}$ follows from the fact that $\characteristicdiamondtwoarg{[\leftubar,\ubarboot]}{[\moreinterestingu_1,\moreinterestingu_2]}$ is a flow-out from $\datahypfortimefunctiontwoarg{0}{[\leftubar,\ubarboot]}$, see \eqref{E:BADEFINITIONOFFLOWOUTMAP}--\eqref{E:BAINVERSEFLOWMAPOFNEWULW1INFTYFINITE} and the respective estimates \eqref{E:FLOWMAPOFNEWULC01QUANTIFIED}--\eqref{E:INVERSEFLOWMAPOFNEWULC01QUANTIFIED} that now hold up to the closure.
\end{proof}

\section{The main $L^\infty$ estimates} \label{S:MAINLINFINITYESTIMATES}
We continue to work under the assumptions of Sect.\,\ref{SS:SILENTFACTS}.
In this section, we derive $L^{\infty}$ estimates for the wave and transport variables, as well as the eikonal function quantities 
that in particular yield improvements of the auxiliary bootstrap assumptions stated in Sect.\,\ref{SSS:AUXBOOTSTRAP}. We also derive pointwise estimates for various commutations used in our $L^2$ analysis. 

\subsection{$L^{\infty}$ estimates and improvement of the auxiliary bootstrap assumptions} \label{SS:LINFINITYFLUIDANDEIKONALANDIMPROVEMENTOFAUX}

\begin{proposition}[$L^{\infty}$ estimates and improvement of the auxiliary bootstrap assumptions] 
\label{P:IMPROVEMENTOFAUXILIARYBOOTSTRAP}
Under the parameter-size and initial data assumptions 
of Sects.\,\ref{SS:PARAMETERSIZEASSUMPTIONS}
and\,\ref{SS:ASSUMPTIONSONDATA}
and the bootstrap assumptions of 
Sects.\,\ref{SS:BOOTSTRAPSCAFFOLDING} and\,\ref{SS:MAINQUANTITATIVEBOOTSTRAPASSUMPTIONS},
the following estimates hold for
$(\ubar,u) \in [\leftubar,\ubarboot) \times [\moreinterestingu_1,\moreinterestingu_2]$
(where we recall that $\wavearray$ and $\wavearraypartial$ are defined in Def.\,\ref{D:ARRAYSOFWAVEVARIABLES},
that $\newL$ is defined in \eqref{E:DOUBLENULLEIKONALFUNCTIONNORMALIZEDNULLVECTORFIELDS},
and that in Sect.\,\ref{S:COMMUTATIONFORMULASIDENTITIESANDSTRINGSOFCOMMUTATORS},
we introduced notation for strings of commutation vectorfields).

\noindent \underline{\textbf{$L^{\infty}$ estimates for small quantities}}.\begin{align}  
	\left\| \RRiemann \right\|_{L^{\infty}\left(\doublenulltoritwoarg{\ubar}{u}\right)} 
	& \leq \mathring{\upalpha} 
		+ 
		C \fundbootsmall, 
		\label{E:LINFINITYIMPROVEMENTAUXR+} 
		\\
	\left\| \wavearraypartial \right\|_{L^{\infty}\left(\doublenulltoritwoarg{\ubar}{u}\right)} 
	& \leq C \fundbootsmall,
		 \label{E:IMPROVEAUXWAVEARRAYPARTIALLINFINITY} 
					\\
\begin{split} \label{E:LINFINITYIMPROVEMENTAUXWAVEARRAY} 
	\left\| \newL \tander^{\leq \Ntop-8} \muX \wavearray \right\|_{L^{\infty}\left(\doublenulltoritwoarg{\ubar}{u}\right)},
		\,
	\left\| \comdersmall^{[1,\Ntop-7];1} \wavearray \right\|_{L^{\infty}\left(\doublenulltoritwoarg{\ubar}{u}\right)},
	&
			\\
	\left\| \newL \tander^{\leq 4} \muX \muX \wavearray \right\|_{L^{\infty}\left(\doublenulltoritwoarg{\ubar}{u}\right)},
		\,
	\left\| \comdersmall^{[1,6];2} \wavearray \right\|_{L^{\infty}\left(\doublenulltoritwoarg{\ubar}{u}\right)},
		& \\
	\left\| \newL \tander^{\leq 2} \muX \muX \muX \wavearray \right\|_{L^{\infty}\left(\doublenulltoritwoarg{\ubar}{u}\right)},
		\,
	\left\| \comdersmall^{[1,5];3} \wavearray \right\|_{L^{\infty}\left(\doublenulltoritwoarg{\ubar}{u}\right)},
	& \\
	\left\| \newL \muX \muX \muX \muX \wavearray \right\|_{L^{\infty}\left(\doublenulltoritwoarg{\ubar}{u}\right)}
	& \leq 
	C \fundbootsmall,
	\end{split}
		\\
	\begin{split} \label{E:LINFINITYIMPROVEMENTAUXMUSMALL}
	\left\| \newL \tander^{[1,\Ntop-8]} \upmu \right\|_{L^{\infty}\left(\doublenulltoritwoarg{\ubar}{u}\right)}, 
		\,
	\left\| \tander_*^{[1,\Ntop-8]} \upmu \right\|_{L^{\infty}\left(\doublenulltoritwoarg{\ubar}{u}\right)},
	& 
	 \\
	\left\| \newL \comdersmall^{[1,5];1}  \upmu \right\|_{L^{\infty}\left(\doublenulltoritwoarg{\ubar}{u}\right)}, 
		\,
	\left\| \comderdoublesmall^{[1,5];1} \upmu \right\|_{L^{\infty}\left(\doublenulltoritwoarg{\ubar}{u}\right)},
	&	
		\\
	\left\| \newL \comdersmall^{[1,4];2} \upmu \right\|_{L^{\infty}\left(\doublenulltoritwoarg{\ubar}{u}\right)}, 
		\,
	\left\| \comderdoublesmall^{[1,4];2} \upmu \right\|_{L^{\infty}\left(\doublenulltoritwoarg{\ubar}{u}\right)}
	& \leq C \fundbootsmall,
	\end{split} 
		\\
	\left\| \Lsmall^1 \right\|_{L^{\infty}\left(\doublenulltoritwoarg{\ubar}{u}\right)} 
	& 
	\leq \mathring{\upalpha}
			+
			C \fundbootsmall,
		\label{E:LINFINITYIMPROVEMENTAUXL1SMALL} 
		\\
	\left\| \Lsmall^A \right\|_{L^{\infty}\left(\doublenulltoritwoarg{\ubar}{u}\right)}
	\label{E:LINFINITYIMPROVEMENTAUXLASMALL} 
	& 
	\leq C \fundbootsmall
				\\
	\begin{split} \label{E:LINFINITYIMPROVEMENTAUXMIXECTANGENTIALTRANSVERSALDERIVATIVESLISMALL}
	\left\| \newL \tander^{\leq \Ntop-7} \Lsmall^i \right\|_{L^{\infty}\left(\doublenulltoritwoarg{\ubar}{u}\right)}, 
			\,
	\left\| \tander^{[1,\Ntop-7]} \Lsmall^i \right\|_{L^{\infty}\left(\doublenulltoritwoarg{\ubar}{u}\right)}, 
			\\
	\left\| \newL \comder^{[1,\Ntop-8];1} \Lsmall^i \right\|_{L^{\infty}\left(\doublenulltoritwoarg{\ubar}{u}\right)}, 
		\,
	\left\| \comdersmall^{[1,\Ntop-8];1} \Lsmall^i \right\|_{L^{\infty}\left(\doublenulltoritwoarg{\ubar}{u}\right)}, 
			\\
	\left\| \newL \comder^{[1,5];2}\Lsmall^i \right\|_{L^{\infty}\left(\doublenulltoritwoarg{\ubar}{u}\right)}, 
		\,
	\left\| \comdersmall^{[1,5];2} \Lsmall^i \right\|_{L^{\infty}\left(\doublenulltoritwoarg{\ubar}{u}\right)},
		\\
	\left\| \newL \comder^{[1,4];3}\Lsmall^i \right\|_{L^{\infty}\left(\doublenulltoritwoarg{\ubar}{u}\right)}, 
		\,
	\left\| \comdersmall^{[1,4];3} \Lsmall^i \right\|_{L^{\infty}\left(\doublenulltoritwoarg{\ubar}{u}\right)}
	& 
	\leq C \fundbootsmall.
\end{split}
\end{align}\medskip

\noindent \underline{\textbf{$L^{\infty}$ estimates tied to pure transversal derivatives}}. 
\begin{align}
	\left 
		\| \muX^M \RRiemann 
	\right\|_{L^{\infty}\left(\doublenulltoritwoarg{\ubar}{u}\right)} 
	& 
	\leq 
	\mathring{\updelta}
	+ 
	C \fundbootsmall, 
	&
	1 \leq M \leq 4,
		\label{E:LINFINITYIMPROVEMENTAUXTRANSVERSALPDERIVATIVESRRIEMANNLARGE}  
			\\
	\left 
		\| \muX^M \wavearraypartial 
	\right\|_{L^{\infty}\left(\doublenulltoritwoarg{\ubar}{u}\right)} 
	& 
	\leq 
	C \fundbootsmall, 
	&
	1 \leq M \leq 4,
		\label{E:LINFINITYIMPROVEMENTAUXTRANSVERSALPDERIVATIVESPARTIALWAVEARRAYSMALL}  
			\\
	\left\| 
		\muX^M \upmu 
	\right\|_{L^{\infty}\left(\doublenulltoritwoarg{\ubar}{u}\right)} 
	& 
	\leq 
	1 +
		\frac{3}{2\blowupdeltadoublenull}
		\left\|
				\muX^M
				\left\lbrace
				G_{\Lunit\Lunit}^0
				\muX \RRiemann 
				\right\rbrace
			\right\|_{L^{\infty}\left(\doublenulltoritwoarg{\leftubar}{u}\right)}
	+ 
	C \fundbootsmall, 
		\label{E:LINFINITYIMPROVEMENTAUXMUANDTRANSVERSALDERIVATIVES} 
	&
	0 \leq M \leq 3,
		\\
	\left\|
		\newL \muX^M \upmu 
	\right\|_{L^{\infty}\left(\doublenulltoritwoarg{\ubar}{u}\right)}
	& 
	\leq 
	\frac{2}{3 \blowupdeltadoublenull} 
			\left\|
				\muX^M
				\left\lbrace
				G_{\Lunit\Lunit}^0
				\muX \RRiemann 
				\right\rbrace
			\right\|_{L^{\infty}\left(\doublenulltoritwoarg{\leftubar}{u}\right)}
	+ 
	C \fundbootsmall, 
		\label{E:LINFINITYIMPROVEMENTAUXNULLDERIVATIVEMUANDTRANSVERSALDERIVATIVES}  
	&
	0 \leq M \leq 3,
		\\
	\left\| 
		\muX^M \Lsmall^1
	\right\|_{L^{\infty}\left(\doublenulltoritwoarg{\ubar}{u}\right)} 
	& 
	\leq	 
	\mathring{\updelta}
	+ 
	C \fundbootsmall, 
		\label{E:LINFINITYESTIMATESFORMUXDERIVATIVESOFL1}
	&
	1 \leq M \leq 3,
		\\
	\left\| 
		\muX^M \Lsmall^A
	\right\|_{L^{\infty}\left(\doublenulltoritwoarg{\ubar}{u}\right)} 
	& 
	\leq	 
	C \fundbootsmall, 
		\label{E:LINFINITYESTIMATESFORMUXDERIVATIVESOFLA}
	&
	1 \leq M \leq 3.
\end{align}\end{proposition}

\begin{proof} 
The same proof as in \cite{abbrescia2022emergence}*{Prop.\,18.1} holds with $\Lrough \eqdef \frac{1}{\Lunit \timefunction}\Lunit$ in place of $\newL$, \cite{abbrescia2022emergence}*{Lemmas\,17.1 and 17.3} in place of Lemma\,\ref{L:PROPERTIESOFFLOWMAPOFNEWL} and Cor.\,\ref{C:TRANSPORTESTIMATESFORNEWLFEQUALSSOURCE}, and \cite{abbrescia2022emergence}*{(16.12b)} in place of \eqref{E:BALDERIVATIVEOFUBARAPPROXIMATELYUNITY}, respectively. 
\end{proof}

The following corollary is an immediate consequence of the fact that we have improved the auxiliary bootstrap assumptions. 

\begin{corollary}[$\mathring{\upalpha}^{1/2}$ and $\auxbootsmall$ can be replaced by $C\mathring{\upalpha}$ and  $C \fundbootsmall$] 
\label{C:IMPROVEAUX}
All prior inequalities whose RHS feature an explicit factor of $\mathring{\upalpha}^{1/2}, \auxbootsmall$ remain true with $\mathring{\upalpha}^{1/2},\, \auxbootsmall$ respectively replaced by $C \mathring{\upalpha},\, C \fundbootsmall$.
\end{corollary}

\subsection{$L^\infty$ commutator estimates}

In this section we provide several commutator estimates that will be used in our $L^2$ analysis.

\begin{lemma}[Commuting the angular Laplacian with strings of geometric vectorfields \cite{jSgHjLwW2016}*{Cor 8.9}] 
Let $1 \leq N \leq \Ntop$.
Then the following commutator estimates hold on $\characteristicdiamondtwoarg{[\leftubar,\ubarboot)}{[\moreinterestingu_1,\moreinterestingu_2]}$:
	\begin{align}
		\left| [\angLap, \tander^{N-2}] \wavearray \right| & \lesssim \left| \tander^{[1,N]}\wavearray \right| + \left| \tander^{[1,N-1]} \controlvars\right| \label{E:COMMINGANGULARLAPWITHTANGENTIAL} 
	\end{align}
\end{lemma}

\begin{lemma}[Commuting $\upmu$-weighted Cartesian derivatives with strings of geometric vectorfields \cite{abbrescia2022emergence}*{Lemma 24.2}] 
\label{L:POINTWISEESTIMATESFORMUWEIGHTEDCARTESIANCOMMUTATORWITHALLWAVEVARS} 
Let $1 \leq N \leq \Ntop$.
Then the following commutator estimates hold on $\characteristicdiamondtwoarg{[\leftubar,\ubarboot)}{[\moreinterestingu_1,\moreinterestingu_2]}$:
\begin{align}
\begin{split} \label{E:POINTWISEESTIMATESFORMUWEIGHTEDCARTESIANCOMMUTATORWITHALLWAVEVARS}
	& 
	\left| [\upmu \partial_i, \tander^N] (\vortrenormalized,\GradEnt) \right|, 
		\, 
	\left| [\upmu \Flatcurl, \tander^N] (\vortrenormalized,\GradEnt) \right|, 
		\, 
	\left| [\upmu \Flatdiv, \tander^N] (\vortrenormalized,\GradEnt) \right|  
		\\
	&  
	\lesssim 
	\left| \tander^{\leq N} (\vortrenormalized, \GradEnt) \right| 
	+ 
	\left| \muX \tander^{\leq N-1}(\vortrenormalized,\GradEnt)\right| 
	+  
	\fundbootsmall |\muX \tander^{[1,N-1]} \wavearray| 
	+
	\fundbootsmall \upmu |\tander^N \wavearray| 
	+ 
	\fundbootsmall |\tandersmall^{[1,N]} \badcontrolvars|. 
\end{split}
\end{align}
\end{lemma}

\section{Sharp control of $\upmu$, the properties of $\Upsilon$, and pointwise estimates of the double-null acoustic geometry}
\label{S:SHARPCONTROLOFMUANDPROPERTIESOFCHOVGEOTOCARTESIAN}

We continue to work under the assumptions of Sect.\,\ref{SS:SILENTFACTS}.
In this section, we derive sharp control of $\upmu$ and its derivatives, as well as strict improvements 
of the bootstrap assumptions from 
Sects.\,\ref{SSS:BAFORINVERSEFOLIATIONDENSITY}--\ref{SSS:BOOTSTRAPINGOINGEIKONALFUNCTION}. Because our $L^2$-type energies feature $\weight$ weights, 
these sharp estimates will reveal a relationship between $\upmu$ and $\weight$ which plays a fundamental role in our proof of the energy estimates.
Next, in Prop.\,\ref{P:HOMEOMORPHICANDDIFFEOMORPHICEXTENSIONOFCARTESIANCOORDINATES}, 
we derive homeomorphism and diffeomorphism properties 
of the change of variables map $\Upsilon$ from geometric coordinate to Cartesian coordinates,
which are crucial for understanding the structure of the singular boundary. 
The results of Prop.\,\ref{P:HOMEOMORPHICANDDIFFEOMORPHICEXTENSIONOFCARTESIANCOORDINATES}
also yield strict improvements of the bootstrap assumptions for
$\Upsilon$ stated in Sect.\,\ref{SSS:BACHOVMAPS}.

\subsection{Sharp control of $\upmu$ and its derivatives}
\label{SS:SHARPCONTROLOFMUANDDERIVATIVES}

\begin{proposition}[Sharp control of $\upmu$ and its derivatives]
\label{P:SHARPCONTROLOFMUANDDERIVATIVES}
The following estimates hold.

\medskip

\noindent \underline{\textbf{Minima of $\upmu$ occur precisely on $\argubarnewultorus{\ubar}$.}} For each $\ubar \in [\leftubar,\ubarboot]$, we have:
	\begin{align} 
		\min_{\ingoingcharacteristicsurfacetwoarg{\ubar}{[\moreinterestingu_1,\moreinterestingu_2]}} \upmu  =  \min_{\argubarnewultorus{\ubar}} \upmu \ge \frac{1}{1.01} \weight\big|_{\argubarnewultorus{\ubar}}  = \frac{-1}{1.01} \ubar.\label{E:KEYESTIMATECONTROLLINGINVERSEMUBYINVERSEWEIGHT} 
	\end{align}
	
\medskip		
		
\noindent \underline{\textbf{$\upmu$ is large when $|u| \geq \frac{3}{4}\interestingu$}}.
	The following lower bound holds, where $\boringregionmupositive > 0$ is the constant appearing in 
	\eqref{E:DATAMUISLARGEINBORINGREGION}:
	\begin{align} \label{E:MUISLARGEINBORINGREGION}
			\min_{\characteristicdiamondtwoarg{[\leftubar,\ubarboot]}{[\moreinterestingu_1,\moreinterestingu_2]} \setminus \characteristicdiamondtwoarg{[\leftubar,\ubarboot]}{[-\frac{3}{4}\interestingu,\frac{3}{4}\interestingu]}} \upmu
			& \geq \frac{\boringregionmupositive}{2}.
	\end{align}\medskip

\noindent \underline{\textbf{Transversal convexity of $\upmu$ and its consequences}}.
The following estimates hold:

	\begin{align} 
			\frac{\secondtransversalderivativemulowerbound}{5}
		\leq 
		\min_{\characteristicdiamondtwoarg{[\leftubar,\ubarboot]}{[\moreinterestingu_1,\moreinterestingu_2]}}
				\Big\lbrace
				\newuL \newuL \upmu, \newuL \muX \upmu, 
				\muX \muX \upmu \Big\rbrace &
			\leq 
			\max_{\characteristicdiamondtwoarg{[\leftubar,\ubarboot]}{[\moreinterestingu_1,\moreinterestingu_2]}}
			\Big\lbrace
				\newuL \newuL \upmu, \newuL \muX \upmu, \muX \muX \upmu
				\Big\rbrace
				\leq 
			\frac{5}{\secondtransversalderivativemulowerbound}, \label{E:MUTRANSVERSALCONVEXITY} \\
			\frac{\secondtransversalderivativemulowerbound}{5} \le  \essinf_{\characteristicdiamondtwoarg{[\leftubar,\ubarboot]}{[\moreinterestingu_1,\moreinterestingu_2]}} \nullgeop{u}\newuL \upmu & \le   \esssup_{\characteristicdiamondtwoarg{[\leftubar,\ubarboot]}{[\moreinterestingu_1,\moreinterestingu_2]}}  \nullgeop{u}\newuL \upmu \le \frac{5}{\secondtransversalderivativemulowerbound}. \label{E:MULESSREGULARTRANSVERSALCONVEXITY}
	\end{align}
	
	\begin{align} 
		\min_{\characteristicdiamondtwoarg{[\leftubar,\ubarboot]}{[\moreinterestingu_1,\moreinterestingu_2]}
			\backslash
			\characteristicdiamondtwoarg{[\leftubar,\ubarboot]}{[-\frac{1}{2} \interestingu,\frac{1}{2}\moreinterestingu_2 + \frac{1}{8} \interestingu]}} 
			|\newuL \upmu | 
			& \geq \frac{\secondtransversalderivativemulowerbound \interestingu}{8}.
			\label{E:REGIONWHERENEWULMULEVELSETISNOTLOCATED}
	\end{align}

	Moreover, the following pointwise estimates hold on
$\characteristicdiamondtwoarg{[\leftubar,\ubarboot]}{[\moreinterestingu_1,\moreinterestingu_2]}$:
\begin{align} \label{E:NEWULMUUBOUNDEDBYSQRTMU}
	|\newuL \upmu|	& \leq 
	C \sqrt{\upmu}.
	\end{align}

\medskip

\noindent \underline{\textbf{Location of $\newuLmulevelset^{[\leftubar,\ubarboot]}$ 
and $\argubarnewultorus{\ubar}$}}.
With
$\newuLmulevelset^{[\leftubar,\ubarboot]}$ 
and $\argubarnewultorus{\ubar}$
denoting the sets defined in Def.\,\ref{D:LEVELSETOFNEWULMUANDNEWULADAPTEDTORI},
we have:
\begin{subequations}
\begin{align} \label{E:NEWULMULEVELSETIMPROVELOCATION}
			\newuLmulevelset^{[\leftubar,\ubarboot]}
			& 
			\subset
			\characteristicdiamondtwoarg{[\leftubar,\ubarboot]}{[-\frac{1}{2}\interestingu,\frac{1}{2}\moreinterestingu_2 + \frac{1}{4}\interestingu]},
				\\
		\textnormal{for each $ \ubar \in [\leftubar,\ubarboot]$,} \qquad
		\argubarnewultorus{\ubar}
		& 
		\subset
		\ingoingcharacteristicsurfacetwoarg{\ubar}{[-\frac{1}{2}\interestingu, \frac{1}{2}\moreinterestingu_2 + \frac{1}{4} \interestingu]}.
			\label{E:IMPROVEDLEVELSETSTRUCTUREANDLOCATIONOFMINONNEWULMUTORI} 
\end{align}
\end{subequations}	
Moreover, with $\CHOVdoublenulltoubarnewulmu^{-1}$
denoting the inverse function of the function
$\CHOVdoublenulltoubarnewulmu$ defined in Def.\,\ref{D:CHOVINVOLVINGTHEDOUBLENULLCOORD}, we have:
\begin{align} \label{E:INTERMSOFPHICHOVMAPIMPROVEDLEVELSETSTRUCTUREANDLOCATIONOFMIN} 
\textnormal{For each $\ubar \in [\leftubar,\ubarboot]$,} \qquad
	\CHOVdoublenulltoubarnewulmu^{-1}
	\left(\lbrace \ubar \rbrace \times \lbrace 0 \rbrace \times \mathbb{T}^2 \right)
	\subset
	\lbrace \ubar \rbrace \times \left[-\frac{1}{2} \interestingu, \frac{1}{2}\moreinterestingu_2 + \frac{1}{4} \interestingu\right] \times \mathbb{T}^2.
\end{align}	
\medskip

\noindent \underline{\textbf{Control of null and almost null derivatives of $\upmu$:}} The following estimates hold.
\begin{subequations}
	\begin{align}
		-\frac{19}{16} \blowupdeltadoublenull \le \min_{\characteristicdiamondtwoarg{[\leftubar,\ubarboot]}{[\moreinterestingu_1,\moreinterestingu_2]}} \Lunit \upmu & \le \max_{\characteristicdiamondtwoarg{[\leftubar,\ubarboot]}{[\moreinterestingu_1,\moreinterestingu_2]}} \Lunit \upmu \le - \frac{13}{16} \blowupdeltadoublenull, \label{E:BOUNDSONLMUINTERESTINGREGION} \\
		-\frac{19}{16} \blowupdeltadoublenull \le \min_{\characteristicdiamondtwoarg{[\leftubar,\ubarboot]}{[\moreinterestingu_1,\moreinterestingu_2]}} \geop{t} \upmu & \le \max_{\characteristicdiamondtwoarg{[\leftubar,\ubarboot]}{[\moreinterestingu_1,\moreinterestingu_2]}} \Lunit \upmu \le - \frac{13}{16} \blowupdeltadoublenull \label{E:BOUNDSONGEOMETRICTDERIVATIVEMUINTERESTINGREGION} \\
		- \frac{117}{112}   \le \min_{\characteristicdiamondtwoarg{[\leftubar,\ubarboot]}{[\moreinterestingu_1,\moreinterestingu_2]}} \newL \upmu & \le \max_{\characteristicdiamondtwoarg{[\leftubar,\ubarboot]}{[\moreinterestingu_1,\moreinterestingu_2]}} \newL \upmu \le - \frac{133}{144}  \label{E:BOUNDSONNEWLMUINTERESTINGREGION}
	\end{align}
\end{subequations}\medskip

\noindent \underline{\textbf{Control of null derivatives of $\ubar$:}} The following estimates hold for $\Lunit \ubar = \ReciprocalLunitAppliedtoTimeFunction^{-1}$:
\begin{align}
	\frac{7}{9} \blowupdeltadoublenull \le \min_{\characteristicdiamondtwoarg{[\leftubar,\ubarboot]}{[\moreinterestingu_1,\moreinterestingu_2]}} \Lunit \ubar  = \ReciprocalLunitAppliedtoTimeFunction\le  \max_{\characteristicdiamondtwoarg{[\leftubar,\ubarboot]}{[\moreinterestingu_1,\moreinterestingu_2]}} \Lunit \ubar \le \frac{9}{7}\blowupdeltadoublenull.\label{E:BOUNDSOFLUBARINTERESTINGREGION}
\end{align}\medskip
	
	\noindent \underline{\textbf{A lower bound tied to the blowup in the interesting region}}.
	The following lower bounds holds:
	\begin{align} \label{E:LOWERBOUNDONMAGNITUDEOFXRPLUS}
	\min_{\characteristicdiamondtwoarg{[\leftubar,\ubarboot]}{[\moreinterestingu_1,\moreinterestingu_2]}}
	 \upmu |X \RRiemann|
	& \geq
		\frac{\blowupdeltadoublenull}{|\bar{\Speed}_{;\LogDensity} + 1|}, 
	\end{align}
	where $\blowupdeltadoublenull > 0$ is the data-parameter from \eqref{E:DELTASTARDOUBLENULLDEF}, and the vectorfield $X$ has Euclidean length satisfying
	$\sqrt{\sum_{a=1}^3 (X^a)^2} = 1 + \mathcal{O}(\mathring{\upalpha})$.

	\medskip
	
	\noindent \underline{\textbf{Especially sharp control in a small neighborhood}}.
	Recall that $\flowmapofnewularg{\Delta s}$ is the flow map of $\newuL$
	(see Theorem\,\ref{T:CONSTRUCTIONOFTHEINGOINGEIKONALFUNCTIONCONTRACTIONMAPPING}).
	There exists a small constant $\underline{\Delta S}_0$ depending on the modulus of continuity of the data,
	a neighborhood (in the geometric coordinate topology of $\characteristicdiamondtwoarg{[\leftubar,\ubarboot]}{[\moreinterestingu_1,\moreinterestingu_2]}$)
	$\smallneighborhoodofcreasearg{[\leftubar,\ubarboot]}$ of 
	$\newuLmulevelset^{[\leftubar,\ubarboot]}$
	of the form:
	\begin{align} \label{E:SMALLNEIGHBORHOOD}
		\smallneighborhoodofcreasearg{[\leftubar,\ubarboot]}
		=
		\flowmapofnewularg{(-\Delta S,\Delta S)}
		\left(\newuLmulevelset^{[\leftubar,\ubarboot]} \right)
		\eqdef 
		\bigcup_{\Delta s' \in (-\Delta S,\Delta S)}
		\flowmapofnewularg{\Delta s'}
		\left( \newuLmulevelset^{[\leftubar,\ubarboot]} \right)
	\end{align}
	such that:
	\begin{align} \label{E:SMALLNEIGHBORHOODCONTAINEDININTERESTINGREGION}
	\smallneighborhoodofcreasearg{[\leftubar,\ubarboot]}
	&
	\subset 
	\characteristicdiamondtwoarg{[\leftubar,\ubarboot]}{[\moreinterestingu_1,\moreinterestingu_2]},
	\end{align}
	and a constant $\awayfromsmallneightborhoodnmupositive > 0$ defined by:
	\begin{align} \label{E:MU2DEF}
		\awayfromsmallneightborhoodnmupositive
		\eqdef \min \left\lbrace \frac{\secondtransversalderivativemulowerbound}{5} (\underline{\Delta S}_0)^2, \frac{\boringregionmupositive}{2} \right\rbrace
	\end{align}
	(where $\boringregionmupositive$ is as in \eqref{E:DATAMUISLARGEINBORINGREGION} and
	$\secondtransversalderivativemulowerbound$ is as in \eqref{E:MUTRANSVERSALCONVEXITY})
	such that the following estimates hold:
	\begin{align} \label{E:NEWLMUISALMOSTMINUSONEINSMALLNEIGHBORHOOD}
	- 1.02 
	&
	\leq
	\min_{\smallneighborhoodofcreasearg{[\leftubar,\ubarboot]}} \newL \upmu 
	\leq 
	\max_{\smallneighborhoodofcreasearg{[\leftubar,\ubarboot]}} \newL \upmu 
	\leq 
	-0.98,
		\\
	\awayfromsmallneightborhoodnmupositive
	& 
	\leq
	\min_{
	\characteristicdiamondtwoarg{[\leftubar,\ubarboot]}{[\moreinterestingu_1,\moreinterestingu_2]}
	\backslash
	\smallneighborhoodofcreasearg{[\leftubar,\ubarboot]}} \upmu.
		\label{E:OUTSIDEOFSMALLNEIGHBORHOODLOWERBOUNDFORMU}
	\end{align}\end{proposition}

\begin{proof}
\begin{remark}[Silent use of Lemma\,\ref{L:CONTINUOUSEXTNESION}]
\label{R:SILENTUSEOFLEMMACONTINUOUSEXTNESION}
Throughout this proof, we sometimes silently use the continuous extension properties shown in Lemma\,\ref{L:CONTINUOUSEXTNESION}.
In particular, these properties allow us to extend various results that
were proved prior to the lemma on the half-open rough time interval $[\leftubar,\ubarboot)$
to the closed rough time interval
$[\leftubar,\ubarboot]$.
\end{remark}

\medskip
	
	\noindent \textbf{Proof of \eqref{E:MUTRANSVERSALCONVEXITY}--\eqref{E:MULESSREGULARTRANSVERSALCONVEXITY}:}
	To prove \eqref{E:MUTRANSVERSALCONVEXITY} for $\newuL \newuL \upmu$,
	we use first use \eqref{E:DOUBLENULLEIKONALFUNCTIONNORMALIZEDNULLVECTORFIELDS}, 
	Lemmas\,\ref{L:COMMUTATORSTOCOORDINATES},\,\ref{L:RELATIONSHIPBETWEENNULLCOORDINATEPARTIALSANDDOUBLENULLVECTORFIELDS},\,\ref{L:IDENTITIESFORDOUBLENULLFRAMEANDACOUSTICGEOMETRY}, and\,\ref{L:CONTINUOUSEXTNESION},
	and the estimate
	$\frac{1}{\Lunit \ubar} \approx 1$ (see \eqref{E:CLOSEDVERSIONLUBARESTIMATECHOVGEOTODOUBLENULL})
	to deduce that for $(\ubar,u) \in [\leftubar,\ubarboot] \times [\moreinterestingu_1,\moreinterestingu_1]$,
	we have
	$\| 
		\newL \newuL \newuL \upmu 
	\|_{L^{\infty}\left(\doublenulltoritwoarg{\ubar}{u}\right)} \leq C$.
	From this bound,
	\eqref{E:ESSENTIALSUPPREMUMTRANSPORTNEWLLFESTIMATE},
	and the data-assumption
	\eqref{E:DATATASSUMPTIONMUTRANSVERSALCONVEXITY},
	we deduce that in the region $\characteristicdiamondtwoarg{[\leftubar,\ubarboot]}{[\moreinterestingu_1,\moreinterestingu_2]}$,
	we have 
	$
	\frac{\secondtransversalderivativemulowerbound}{4}
	-
	C |\leftubar|
	\leq
	\newuL \newuL \upmu
	\leq
	\frac{4}{\secondtransversalderivativemulowerbound}
	+
	C |\leftubar|
	$,
	which, for $\leftubar$ sufficiently small,
	implies \eqref{E:MUTRANSVERSALCONVEXITY}
	for $\newuL \newuL \upmu$
	(recall that, as we highlighted in Sect.\,\ref{SS:CONVENTIONSFORCONSTANTS}, the constants $C$
	can be chosen to be independent of $\leftubar$).
	The remaining estimates in \eqref{E:MUTRANSVERSALCONVEXITY}
	follow from a nearly identical argument,
	where we use the identity
	\eqref{E:DOUBLENULLUDERIVATIVEINTERMSOFGEOMETRICVECTORFIELDSANDDERIVATIVESOFUBAR}, \eqref{E:SMALLC01ESTIMATESFORUBAR}, and Cor.\,\ref{C:PRELIMINARYESTIMATESFORTHEDOUBLENULLACOUSTICSCALARS}
	when deriving the estimates involving $\nullgeop{u} \newuL \upmu$.

\medskip
	
	\noindent \textbf{Proof of \eqref{E:MUISLARGEINBORINGREGION}, \eqref{E:REGIONWHERENEWULMULEVELSETISNOTLOCATED}, \eqref{E:BOUNDSONLMUINTERESTINGREGION}, \eqref{E:BOUNDSONGEOMETRICTDERIVATIVEMUINTERESTINGREGION}:} To prove \eqref{E:REGIONWHERENEWULMULEVELSETISNOTLOCATED}, we first argue as in the proof of \eqref{E:MUTRANSVERSALCONVEXITY}--\eqref{E:MULESSREGULARTRANSVERSALCONVEXITY} to deduce that $\| 
		\newL  \newuL \upmu 
	\|_{L^{\infty}\left(\doublenulltoritwoarg{\ubar}{u}\right)} \leq C$. From this bound, the estimate \eqref{E:ESSENTIALSUPPREMUMTRANSPORTNEWLLFESTIMATE}, and the data-assumption \eqref{E:DATATASSUMPTIONREGIONWHERENEWULMULEVELSETISNOTLOCATED}, we find that in the spacetime region $\characteristicdiamondtwoarg{[\leftubar,\ubarboot]}{[\moreinterestingu_1,\moreinterestingu_2]}
			\backslash
			\characteristicdiamondtwoarg{[\leftubar,\ubarboot]}{[-\frac{1}{2} \interestingu,\frac{1}{2}\moreinterestingu_2 + \frac{1}{8} \interestingu]}$, we have $|\newuL \upmu| \ge  \frac{\secondtransversalderivativemulowerbound \interestingu}{4}  - C |\leftubar|$. Assuming $|\leftubar|$ is sufficiently small, we conclude \eqref{E:REGIONWHERENEWULMULEVELSETISNOTLOCATED}. 
			
			The estimates \eqref{E:MUISLARGEINBORINGREGION} and \eqref{E:BOUNDSONLMUINTERESTINGREGION} follow from similar arguments based on the data assumptions \eqref{E:DATAMUISLARGEINBORINGREGION}--\eqref{E:DATATASSUMPTIONLMUQUANTITATIVENEGATIVITY}. We also conclude \eqref{E:BOUNDSONGEOMETRICTDERIVATIVEMUINTERESTINGREGION} by combining these arguments with the relation $\geop{t} \upmu = \Lunit \upmu + \mathcal{O}(\fundbootsmall)$, which follows from the identity $\geop{t}\upmu = \Lunit \upmu - L^A \geop{x^A} \upmu$, Lemma\,\ref{L:COMMUTATORSTOCOORDINATES},
	Prop.\,\ref{P:SCHEMATICSTRUCTUREOFVARIOUSTENSORSINTERMSOFCONTROLVARS},
	and the estimates of Prop.\,\ref{P:IMPROVEMENTOFAUXILIARYBOOTSTRAP}. 
	
\medskip
	
	\noindent \textbf{Proof of \eqref{E:BOUNDSOFLUBARINTERESTINGREGION}:} A much simpler, but similar, proof of \eqref{E:SHARPUNIFORMESTIMATEFORLUBARINFLOWOUTREGION} holds. We omit the proof, but remind the reader that it entails integrating \eqref{E:QUASILINEAREQUATIONSATISFIEDBYLUBAR} along integral curves of $\newuL$ and using Lemma\,\ref{L:PROPERTIESOFFLOWMAPOFNEWUL}. 
	
\medskip
\noindent \textbf{Proof of \eqref{E:BOUNDSONNEWLMUINTERESTINGREGION}:} This follows from \eqref{E:BOUNDSONLMUINTERESTINGREGION},\eqref{E:BOUNDSOFLUBARINTERESTINGREGION}, and the identity $\newL \upmu = \frac{1}{\Lunit \ubar} \Lunit \upmu$, which follows from \eqref{E:LUNITANDULUNITAPPLIEDTOEIKONALANDCARTESIANTIME} and \eqref{E:RELATIONBETWEENCARTESIANNORMALIZEDNULLVECTORFIELDSANDEIKONALFUNCTIONORMALIZEDNULLVECTORFIELDS}. 

\medskip	
	
	\noindent \textbf{Proof of \eqref{E:LOWERBOUNDONMAGNITUDEOFXRPLUS}}: We use \eqref{E:MUTRANSPORT},
\eqref{E:IDENTITYFORMAINTERMDRIVINGTHESHOCK},
Prop.\,\ref{P:SCHEMATICSTRUCTUREOFVARIOUSTENSORSINTERMSOFCONTROLVARS},
the estimates of Prop.\,\ref{P:IMPROVEMENTOFAUXILIARYBOOTSTRAP},
and
\eqref{E:BOUNDSONLMUINTERESTINGREGION}
to deduce that in $\twoargMrough{[\timefunction_0,\timefunctionboot],[-\interestingu,\interestingu]}{\muxmulevelsetvalue}$,
we have the estimate
$
\frac{1}{2}
\Speed^{-1}(\Speed^{-1} \Speed_{;\LogDensity} + 1)
|\muX \RRiemann| 
\geq 
\frac{13}{16} \blowupdeltadoublenull 
+ 
\mathcal{O}(\fundbootsmall) 
\geq 
\frac{3}{4} \blowupdeltadoublenull$.
Also Taylor expanding
$\Speed^{-1}(\Speed^{-1} \Speed_{;\LogDensity} + 1)$
around the background solution $\wavearray = 0$
and using 
\eqref{E:NONDEGENCONDITION}
and
\eqref{E:BACKGROUNDSOUNDSPEEDISUNITY},
we deduce the estimate
$
\Speed^{-1}(\Speed^{-1} \Speed_{;\LogDensity} + 1)
|\muX \RRiemann|
=
\left\lbrace
	1 + \mathcal{O}(\mathring{\upalpha})
\right\rbrace
|\bar{\Speed}_{;\LogDensity} + 1|
|\muX \RRiemann|
$. 
Combining these two estimates, we conclude \eqref{E:LOWERBOUNDONMAGNITUDEOFXRPLUS}.

The fact that $\sqrt{\sum_{a=1}^3 (X^a)^2} = 1 + \mathcal{O}(\mathring{\upalpha})$
follows from
\eqref{E:XSMALL},
\eqref{E:SCHEMATICSTRUCTUREOFXSMALL},
and the estimates of Prop.\,\ref{P:IMPROVEMENTOFAUXILIARYBOOTSTRAP}.

\medskip
\noindent \textbf{Proof of \eqref{E:NEWULMULEVELSETIMPROVELOCATION}--\eqref{E:IMPROVEDLEVELSETSTRUCTUREANDLOCATIONOFMINONNEWULMUTORI} and \eqref{E:INTERMSOFPHICHOVMAPIMPROVEDLEVELSETSTRUCTUREANDLOCATIONOFMIN}:} Since $\newuL \upmu = 0$ on $\newuLmulevelset^{[\leftubar,\ubarboot]}$, \eqref{E:NEWULMULEVELSETIMPROVELOCATION} follows from \eqref{E:REGIONWHERENEWULMULEVELSETISNOTLOCATED}. This also proves \eqref{E:IMPROVEDLEVELSETSTRUCTUREANDLOCATIONOFMINONNEWULMUTORI} by \eqref{E:LEVELSETOFNEWULMUFOLIATED}. 

\eqref{E:INTERMSOFPHICHOVMAPIMPROVEDLEVELSETSTRUCTUREANDLOCATIONOFMIN} follows then from \eqref{E:PHIINVERSEIMAGEOFTORUSCROSSMUINTERVALISTORUSCROSSINTERVALCONTAINEDININTERESTINGREGION} and the respective definitions \eqref{E:CHOVFROMGEOTODOUBLENULLCOORDINATES}, \eqref{E:CHOVFROMDOUBLENULLTOUBARNEWULMUCOORDINATES} of the change of variables maps $\CHOVgeotodoublenull$ and $\CHOVdoublenulltoubarnewulmu$. 

\medskip
\noindent \textbf{Proof of \eqref{E:KEYESTIMATECONTROLLINGINVERSEMUBYINVERSEWEIGHT}:} We first claim that the minimum of $\upmu$ along $\ingoingcharacteristicsurfacetwoarg{\ubar}{[\moreinterestingu_1,\moreinterestingu_2]}$ occurs precisely on $\argubarnewultorus{\ubar}$. In particular, this will imply the first identity in \eqref{E:KEYESTIMATECONTROLLINGINVERSEMUBYINVERSEWEIGHT}. From $\ubar|_{\datahypfortimefunctiontwoarg{0}{[\leftubar,\ubarboot]}} 
	= 
	- 
	\upmu|_{\datahypfortimefunctiontwoarg{0}{[\leftubar,\ubarboot]}}$, it follows from
	\eqref{E:MUISLARGEINBORINGREGION},
	\eqref{E:NEWULMULEVELSETIMPROVELOCATION}, and \eqref{E:MU1BIGGERTHANU0}
	that for each fixed $\ubar \in [\leftubar,\ubarboot]$,
	$\min_{\ingoingcharacteristicsurfacetwoarg{\ubar}{[\moreinterestingu_1,\moreinterestingu_2]}} \upmu$ is achieved only in
	the subset $\ingoingcharacteristicsurfacetwoarg{\ubar}{[-\frac{3}{4}\interestingu,\frac{3}{4}\interestingu]}$,
	which is interior to $\ingoingcharacteristicsurfacetwoarg{\ubar}{[\moreinterestingu_1,\moreinterestingu_2]}$ by \eqref{E:MOSTNEGATIVEMOREINTERESTINGVALUEOFU}--\eqref{E:MOSTPOSITIVEMOREINTERESTINGVALUEOFU}. Hence, since $\newuL$ is tangent to $\ingoingcharacteristicsurfacetwoarg{\ubar}{[\moreinterestingu_1,\moreinterestingu_2]}$, it must be that $\newuL \upmu = 0$ at the minima. Considering also \eqref{E:REGIONWHERENEWULMULEVELSETISNOTLOCATED}, we see that the minima of $\upmu$ in must belong to $\newuLmulevelset^{[\leftubar,\ubarboot]} \cap \ingoingcharacteristicsurfacetwoarg{\ubar}{[\moreinterestingu_1,\moreinterestingu_2]}$. This proves the claim and the first identity in \eqref{E:KEYESTIMATECONTROLLINGINVERSEMUBYINVERSEWEIGHT}. 

To prove the inequality in \eqref{E:KEYESTIMATECONTROLLINGINVERSEMUBYINVERSEWEIGHT}, consider any point $q \in \argubarnewultorus{\ubar}$. By Lemma\,\ref{L:PROPERTIESOFFLOWMAPOFNEWUL}, there exists a unique integral curve of $\newuL$ that joins $q$ to a point $q_0$ on $\twoargmumuxtorus{-\ubar}{0}$. Let $\underline{s} = \underline{s}(u)$ denote this integral curve, parametrized by the outgoing eikonal function (recall $\newuL u = 1$), and let $u_0$ and $u$ respectively denote the eikonal function values corresponding to $q_0$ and $q$. Next we note that the same proof of \eqref{E:NEWULMUISNEGATIVEONDATAHYPERSURFACE} implies:
	\begin{align}
		\newuL \upmu \big|_{\twoargmumuxtorus{-\ubar}{0}} = \{ \tfrac{1}{2} \Lunit \upmu + \mathcal{O}(\initialsmall)\}(-\ubar) < 0. \label{E:ESITMATEFORNEWULMUONPRIMALTORI}
	\end{align} This, $\newuL \upmu(q) = 0$, and the transversal convexity estimate \eqref{E:MUTRANSVERSALCONVEXITY} in particular implies $u > u_0$ and that the minimum of $\newuL \upmu \circ \underline{s}(u)$ occurs at $u = u_0$. By these considerations, the identity $\upmu(q_0) = -\ubar$, the Fundamental Theorem of Calculus, and estimates \eqref{E:MUTRANSVERSALCONVEXITY}, \eqref{E:BOUNDSONLMUINTERESTINGREGION}, \eqref{E:ESITMATEFORNEWULMUONPRIMALTORI}, we have:
\begin{subequations}
	\begin{align}
		\upmu(q)  & \ge - \ubar + C_1 \ubar(u - u_0), \label{E:ESTIMATEFORMUALONGINTEGRALCURVEOFNEWUL} \\
		0 & \ge C_2 \ubar + C_3 (u - u_0),\label{E:ESTIMATEFORNEWULMUALONGINTEGRALCURVEOFNEWUL}
	\end{align}
\end{subequations}
for some constants $C_1,C_2,C_3$ depending on $\blowupdelta$ and $\secondtransversalderivativemulowerbound$. Multiplying \eqref{E:ESTIMATEFORNEWULMUALONGINTEGRALCURVEOFNEWUL} by $\ubar$ and inserting the resulting bound $- \frac{C_2}{C_3} (\ubar)^2 \le \ubar (u - u_0)$ into \eqref{E:ESTIMATEFORMUALONGINTEGRALCURVEOFNEWUL}, the desired inequality in \eqref{E:KEYESTIMATECONTROLLINGINVERSEMUBYINVERSEWEIGHT} follows upon making $|\leftubar|$ sufficiently small. The last identity in \eqref{E:KEYESTIMATECONTROLLINGINVERSEMUBYINVERSEWEIGHT} follows easily from the definition \eqref{E:BLOWUPWEIGHT}. 

\medskip 
	
\noindent \textbf{Proof of \eqref{E:NEWULMUUBOUNDEDBYSQRTMU}:} Let $q \in \ingoingcharacteristicsurfacetwoarg{\ubar}{[\moreinterestingu_1,\moreinterestingu_2]}$. By Lemma \,\ref{L:PROPERTIESOFFLOWMAPOFNEWUL},
 there exists a unique integral curve of $\newuL$ that joins $q$ to a point $q_0$ on $\argubarnewultorus{\ubar}$. Let $\underline{s} = \underline{s}(u)$ denote this integral curve, parametrized by the outgoing eikonal function (recall $\newuL u = 1$), and let $u_0$ and $u$ respectively denote the eikonal function values corresponding to $q_0$ and $q$. That is, $\underline{s}(u) = q$, $\underline{s}(u_0) = q_0$, and $\newuL \upmu\circ \underline{s}(u_0) = 0$. Using \eqref{E:MUTRANSVERSALCONVEXITY}, the mean value theorem, and Taylor's theorem, we see that $|\newuL \upmu \circ \underline{s}(u) | \le \tfrac{5}{\secondtransversalderivativemulowerbound} |u - u_0|$ and $\upmu \circ \underline{s}(u) \ge \upmu  \circ \underline{s}(u_0)  + \tfrac{\secondtransversalderivativemulowerbound}{10}(u - u_0)^2 \ge  \tfrac{\secondtransversalderivativemulowerbound}{10}(u - u_0)^2$. Combining the above results, we find that at $q$, we have $|\newuL \upmu| \le \tfrac{5}{\secondtransversalderivativemulowerbound} \sqrt{\tfrac{10}{\secondtransversalderivativemulowerbound}} \sqrt{\upmu}$, which yields \eqref{E:NEWULMUUBOUNDEDBYSQRTMU}.
 
 \medskip

 \noindent \textbf{Proof of \eqref{E:NEWLMUISALMOSTMINUSONEINSMALLNEIGHBORHOOD}:} Recalling that $\flowmapofnewularg{\Delta s}$ is the flow map of $\newuL$,
for real numbers $\Delta s$ small and positive, we consider the set
$
\flowmapofnewularg{(-\Delta s,\Delta s)}
\left(\datahypfortimefunctiontwoarg{0}{[\leftubar,\ubarboot]} \right)
\eqdef
\bigcup_{\Delta s' \in (-\Delta s,\Delta s)}
\flowmapofnewularg{\Delta s'}
\left(\datahypfortimefunctiontwoarg{0}{[\leftubar,\ubarboot]} \right)$,
which by Lemma\,\ref{L:PROPERTIESOFFLOWMAPOFNEWUL} 
and \eqref{E:MUXMUKAPPALEVELSETLOCATION}
is a neighborhood of $\datahypfortimefunctiontwoarg{0}{[\leftubar,\ubarboot]}$
in $\characteristicdiamondtwoarg{[\leftubar,\ubarboot]}{[\moreinterestingu_1,\moreinterestingu_2]}$.
Next, using
Lemma\,\ref{L:COMMUTATORSTOCOORDINATES}, $\newL \eqdef \frac{1}{\Lunit \ubar}\Lunit$, \eqref{E:QUASILINEAREQUATIONSATISFIEDBYLUBAR}, 
Prop.\,\ref{P:SCHEMATICSTRUCTUREOFVARIOUSTENSORSINTERMSOFCONTROLVARS},
the estimates of Lemma\,\ref{L:PROPERTIESANDDIFFEOMORPHICEXTENSIONOFDOUBLENULLCOORDINATES}
and
Prop.\,\ref{P:IMPROVEMENTOFAUXILIARYBOOTSTRAP}, 
and \eqref{E:BOUNDSONLMUINTERESTINGREGION},
we deduce that
$
|\newuL \newL \upmu|
\leq C 
$.
Moreover, by \eqref{E:LUBARISLIKEMINUSLMUALONGDATASURFACE}, we can use the mean value theorem to deduce that in $\flowmapofnewularg{(-\Delta s,\Delta s)}
\left(\datahypfortimefunctiontwoarg{0}{[\leftubar,\ubarboot]} \right)$, the following estimate holds:
$
- 1.01 - C \Delta s
\leq
\newL \upmu \circ \flowmapofnewularg{\Delta s}
\leq - 0.99
+ C \Delta s$.
Choosing and fixing a value of $\Delta s$, which we denote by
$\Delta S_0$, to be sufficiently small (and positive) such that $C \Delta S_0 < .01$,
we arrive at \eqref{E:NEWLMUISALMOSTMINUSONEINSMALLNEIGHBORHOOD}
with
$\flowmapofnewularg{(-\Delta S_0,\Delta S_0)}
\left(\datahypfortimefunctiontwoarg{0}{[\leftubar,\ubarboot]} \right)
$ in place of $\smallneighborhoodofcreasearg{[\leftubar,\ubarboot]}$. We highlight the following fact: $\Delta S_0$ depends on the modulus of continuity of $\newL \upmu$ and is therefore controllable from the initial data. Note now that the proof of \eqref{E:KEYESTIMATECONTROLLINGINVERSEMUBYINVERSEWEIGHT}, specifically \eqref{E:ESTIMATEFORNEWULMUALONGINTEGRALCURVEOFNEWUL} and the preceding paragraph, shows that the maximum $u$-distance between $\twoargmumuxtorus{-\ubar}{0}$ and $\argubarnewultorus{\ubar}$ is comparable to $|\ubar|$. Consequently, we may assume that $\newuLmulevelset^{[\leftubar,\ubarboot]} \Subset \flowmapofnewularg{(-\Delta S_0,\Delta S_0)}
\left(\datahypfortimefunctiontwoarg{0}{[\leftubar,\ubarboot]} \right)
$ in place of $\smallneighborhoodofcreasearg{[\leftubar,\ubarboot]}$ by choosing $|\leftubar|$ sufficiently small compared to $\Delta S_0$, which was already revealed to depend only on the data. In particular, by setting 
	\[ \underline{\Delta S}_0 \eqdef \frac{1}{2} \min_{\substack{\ubar \in [\leftubar, \ubarboot] \\ (x^2,x^3) \in \T^2}}  \left\{\Cartesiantisafunctiononmumxtoriarg{\ubar}{0}(x^2,x^3) + \Delta S_0\right\}  - \Eikonalisafunctiononubarnewultori_{\ubar}(x^2,x^3) ,\] 
It follows that \eqref{E:NEWLMUISALMOSTMINUSONEINSMALLNEIGHBORHOOD} holds with $\smallneighborhoodofcreasearg{[\leftubar,\ubarboot]} \eqdef \flowmapofnewularg{(-\underline{\Delta S}_0,\underline{\Delta S}_0)} \left( \newuLmulevelset^{[\leftubar,\ubarboot]}\right)$.

\medskip

\noindent \textbf{Proof of \eqref{E:OUTSIDEOFSMALLNEIGHBORHOODLOWERBOUNDFORMU}:} Recall that
$
\newuLmulevelset^{[\leftubar,\ubarboot]}
=
\bigcup_{\ubar \in [\leftubar,\ubarboot]}
\argubarnewultorus{\ubar}
$
(see \eqref{E:LEVELSETOFNEWULMUFOLIATED})
and that for $\ubar \in [\leftubar,\ubarboot]$,
we have
$
\upmu|_{\argubarnewultorus{\ubar}}
\ge  \frac{-1}{1.01}\ubar 
$ by \eqref{E:KEYESTIMATECONTROLLINGINVERSEMUBYINVERSEWEIGHT}. For $\ubar \in [\leftubar,\ubarboot]$ and 
$q_0 \in \argubarnewultorus{\ubar}$,
the arguments we used in the proofs of \eqref{E:KEYESTIMATECONTROLLINGINVERSEMUBYINVERSEWEIGHT} and \eqref{E:NEWULMUUBOUNDEDBYSQRTMU} imply that
$
	\frac{\secondtransversalderivativemulowerbound}{5} (\Delta s)^2 
	\leq
	-
	\frac{1}{1.01} \ubar
	+
	\frac{\secondtransversalderivativemulowerbound}{5} (\Delta s)^2
	\leq
	\upmu
	\circ\flowmapofnewularg{\Delta s}(q_0) 
$.
In particular, by \eqref{E:LEVELSETOFNEWULMUFOLIATED},
then
$ \displaystyle
	\frac{\secondtransversalderivativemulowerbound}{5} (\underline{\Delta S}_0)^2 
	\leq
	\min_{\flowmapofnewularg{(-\underline{\Delta S}_0,\underline{\Delta S}_0)}
\left(\newuLmulevelset^{[\leftubar,\ubarboot]} \right)}
	\upmu
$.
From these observations,
the estimate \eqref{E:MUISLARGEINBORINGREGION},
and definition \eqref{E:MU2DEF},
we see that
\eqref{E:OUTSIDEOFSMALLNEIGHBORHOODLOWERBOUNDFORMU}
holds with $\underline{\Delta S}_0$ defined to be the small constant fixed in the proof of
\eqref{E:NEWLMUISALMOSTMINUSONEINSMALLNEIGHBORHOOD}.

 \end{proof}

\subsection{Homeomorphism and diffeomorphism properties of $\Upsilon$}
\label{SS:PROPERTIESOFCHOVFROMGEOMETRICTOCARTESIAN}
Our main goal in this section is to reveal the homeomorphism and diffeomorphism properties of
the change of variables map
$\Upsilon(t,u,x^2,x^3) = (t,x^1,x^2,x^3)$.
We start with the following monotonicity lemma, which plays an important role in
controlling $\Upsilon$.

\begin{lemma}[Monotonicity of $x^1$]
	\label{L:MONOTONICITYOFCARTESIANX1}
	The following identity holds:
		\begin{align}
			\nullgeop{u}x^1 \eqdef \left\{ \ReciprocaluLunitAppliedtoTimeFunction L^1 + \upmu X^1 + \MagnitueofinnerproductofnewLandnewuL  \frac{\angD^A \ubar X^A}{X^1} \right\} + \left\{ L^A + \upmu X^A - \MagnitueofinnerproductofnewLandnewuL \angD^A \ubar\right\} \left\{ \frac{X^A}{X^1} + L^1 \frac{\geop{x^A}\ubar}{\geop{t} \ubar} + \frac{L^B X^B \geop{x^A} \ubar}{X^1 \geop{t} \ubar}\right\}. \label{E:DOUBLENULLUDERIVATIVEOFCATERSIANX1}
		\end{align}
	Moreover, the following estimate holds on $\characteristicdiamondtwoarg{[\leftubar,\ubarboot]}{[\moreinterestingu_1,\moreinterestingu_2]}$:
		\begin{align}
			\nullgeop{u} x^1 = - \frac{1}{2}\upmu \left\{ 1 + \mathcal{O}_{\mydiam}(\mathring{\upalpha})
				\right\}.   \label{E:X1MONOTONOCITYESTIMATE}
		\end{align}
	Finally, for every fixed $[\ubar,x^2,x^3] \in [\leftubar,\ubarboot]\times\T^2$, the map $\ubar \mapsto x^1(\ubar,u,x^2,x^3)$ is strictly decreasing on $[\moreinterestingu_1,\moreinterestingu_2]$. 
	
\end{lemma}

\begin{proof}
\eqref{E:DOUBLENULLUDERIVATIVEOFCATERSIANX1} follows from \eqref{E:GEOMETRICVECTORFIELDSINTERMSOFCARTESIANONES}, \eqref{E:CONVENIENTIDENTITYFORNEWUL}, and \eqref{E:DOUBLENULLUDERIVATIVEINTERMSOFGEOMETRICVECTORFIELDSANDDERIVATIVESOFUBAR}. 

\eqref{E:X1MONOTONOCITYESTIMATE} then follows from \eqref{E:DOUBLENULLUDERIVATIVEOFCATERSIANX1}, \eqref{E:SHARPESTIMATEFORRATIOOFRATIOOFNULLGEOSICINNERPRODUCTANDFOLIATIONDENSITYANDINGOINGMU},  
	Prop.\,\ref{P:SCHEMATICSTRUCTUREOFVARIOUSTENSORSINTERMSOFCONTROLVARS},
	the estimates of Lemma\,\ref{L:PROPERTIESANDDIFFEOMORPHICEXTENSIONOFDOUBLENULLCOORDINATES},
	Prop.\,\ref{P:IMPROVEMENTOFAUXILIARYBOOTSTRAP},
	and Cor.\,\ref{C:IMPROVEAUX},
	and \eqref{E:DATAEPSILONISSMALLERTHANBOOTSTRAPEPSILONSMALLERTHANSQUAREOFDATAALPHA},
	which in particular imply that
	$X^1 = -1 + \Xsmall^1 = - 1 + \mathcal{O}_{\mydiam}(\mathring{\upalpha})$,
	$\Lunit^1 = 1 + \Lsmall^1 = 1 + \mathcal{O}_{\mydiam}(\mathring{\upalpha})$,
	and 
	$X^A, \Lunit^A = \mathcal{O}(\fundbootsmall) = \mathcal{O}_{\mydiam}(\mathring{\upalpha})$.

		To prove the monotonicity of $x^1$, we note that
	\eqref{E:X1MONOTONOCITYESTIMATE} and Prop.\,\ref{P:SHARPCONTROLOFMUANDDERIVATIVES}
	imply that $\nullgeop{u} x^1 < 0$,
	except in the case $\ubarboot = 0$,
	where $\nullgeop{u} x^1$ vanishes precisely along the torus
	$ \argubarnewultorus{0} = \twoargmumuxtorus{0}{0}$,
	which is contained
	in $\ingoingcharacteristicsurfacetwoarg{0}{[-\frac{1}{2}\interestingu,\frac{1}{2}\interestingu]}$.
	Moreover, $\newuL|_{ \argubarnewultorus{0} } = \nullgeop{u}|_{ \argubarnewultorus{0} } = \muX|_{ \argubarnewultorus{0} }$ by \eqref{E:NULLCOORDINATEPARTIALUINTERMSOFGEOMETRICVF}, \eqref{E:NEWULINTERMSOFNULLCOORDINATEPARTIALS}, \eqref{E:MUXUBARVANISHESATTHECREASE}, and \eqref{E:NULLACOUSTICSCALARSTHATAREVANISHINGLIKEMUATTHECREASE}. Hence, \eqref{E:MUTRANSVERSALCONVEXITY} implies that $\nullgeop{u} \newuL \upmu |_{ \argubarnewultorus{0} }>0$. From this, the graphical structure of the crease in double-null coordinates \eqref{E:LASTSLICETORIAREGRAPHSABOVEFLATTORIINGEOMETRICCOORDINATES}, we see that every integral curve of $\nullgeop{u}$ in $\ingoingcharacteristicsurfacetwoarg{0}{[\moreinterestingu_1,\moreinterestingu_2]}$ intersects $ \argubarnewultorus{0}$ in precisely one point. In total, we have shown that for every fixed
	$(\ubar,x^2,x^3) \in [\leftubar,\ubarboot] \times \mathbb{T}^2$,
	the map $u \rightarrow x^1(\ubar,u,x^2,x^3)$ on the domain $[\moreinterestingu_1,\moreinterestingu_2]$
	has a negative derivative, except at possibly a single point
	in $[-\frac{1}{2} \interestingu,\frac{1}{2} \interestingu]$.
	From this fact, we conclude that the map is strictly decreasing,
	as desired. 

\end{proof}

\begin{proposition}[Homeomorphism and diffeomorphism properties of $\Upsilon$]
		\label{P:HOMEOMORPHICANDDIFFEOMORPHICEXTENSIONOFCARTESIANCOORDINATES}
		The change of variables map $\Upsilon(t,u,x^2,x^3) = (t,x^1,x^2,x^3)$
		is \textbf{injective} on the compact set
		$\characteristicdiamondtwoarg{[\leftubar,\ubarboot]}{[\moreinterestingu_1,\moreinterestingu_2]}$
		and satisfies: 
		\begin{align} \label{E:C31BOUNDFORCHOVFROMGEOTOCARTESIAN}
		\| \Upsilon \|_{C_{\textnormal{geo}}^{3,1}\left(\characteristicdiamondtwoarg{[\leftubar,\ubarboot]}{[\moreinterestingu_1,\moreinterestingu_2]} \right)}
		\leq C.
		\end{align}
		In particular, $\Upsilon$ is a homeomorphism from 
		$\characteristicdiamondtwoarg{[\leftubar,\ubarboot]}{[\moreinterestingu_1,\moreinterestingu_2]}$
		onto its image.
		
		Moreover, with $d_{\textnormal{geo}} \Upsilon$ denoting the Jacobian matrix of $\Upsilon$,
		we have:
		\begin{align} \label{E:GEOTOCARTESIANJACOBIANDETERMINANTESTIMATE}
			\mydet d_{\textnormal{geo}} \Upsilon
			& = 
				  \frac{\upmu  X^1}{1 - \Yvfsmallcoeff{A} X^A} 
				\approx
				- 
				\upmu.
		\end{align}
		
		In addition, if $\ubarboot < 0$, then $\Upsilon$ is a diffeomorphism from
		$\characteristicdiamondtwoarg{[\leftubar,\ubarboot]}{[\moreinterestingu_1,\moreinterestingu_2]}$ onto its image.
		
\end{proposition}
\begin{proof}
	We already proved the bound \eqref{E:C31BOUNDFORCHOVFROMGEOTOCARTESIAN} in Lemma\,\ref{L:CONTINUOUSEXTNESION}.
	
	Next, we use
	\eqref{E:GEOMETRICVECTORFIELDSINTERMSOFCARTESIANONES},
	Prop.\,\ref{P:SCHEMATICSTRUCTUREOFVARIOUSTENSORSINTERMSOFCONTROLVARS},
	and
	Prop.\,\ref{P:IMPROVEMENTOFAUXILIARYBOOTSTRAP}
	to compute that
	$
	\frac{\partial \Upsilon(t,u,x^2,x^3)}{\partial (t,u,x^2,x^3)}
	= \begin{pmatrix}
			1 & 0 & 0 & 0  
				\\
			\Lunit^1 + * &   \frac{\upmu  X^1}{1 - \Yvfsmallcoeff{A} X^A} & * & * 
				\\
			0 & 0 & 1 & 0 				
				\\
			0 & 0 & 0 & 1
		\end{pmatrix}
	$,
	where here and in the rest of the proof,
	``$*$'' denotes any quantity that is pointwise bounded in magnitude by $\mathcal{O}(\mathring{\upalpha})$.
	Thus, $\mydet \frac{\partial \Upsilon(t,u,x^2,x^3)}{\partial (t,u,x^2,x^3)} = \upmu X^1 + *$,
	and therefore, using \eqref{E:XSMALL},
	Prop.\,\ref{P:SCHEMATICSTRUCTUREOFVARIOUSTENSORSINTERMSOFCONTROLVARS},
	and
	Prop.\,\ref{P:IMPROVEMENTOFAUXILIARYBOOTSTRAP},
	we compute that
	$\mydet \frac{\partial \Upsilon(t,u,x^2,x^3)}{\partial (t,u,x^2,x^3)} 
	=
	- \left\lbrace 1 + \mathcal{O}(\mathring{\upalpha}) \right\rbrace \upmu$,
	which yields \eqref{E:GEOTOCARTESIANJACOBIANDETERMINANTESTIMATE}.

	If $\ubarboot < 0$, then by Prop.\,\ref{P:SHARPCONTROLOFMUANDDERIVATIVES},
	$\upmu$ is uniformly positive on $\characteristicdiamondtwoarg{[\leftubar,\ubarboot]}{[\moreinterestingu_1,\moreinterestingu_2]}$,
	and from \eqref{E:GEOTOCARTESIANJACOBIANDETERMINANTESTIMATE} and the inverse function
	theorem, we see that
	$\Upsilon$ is a local diffeomorphism on $\characteristicdiamondtwoarg{[\leftubar,\ubarboot]}{[\moreinterestingu_1,\moreinterestingu_2]}$.
	Thus, to complete the proof, we need to show that $\Upsilon$ is injective on
	$\characteristicdiamondtwoarg{[\leftubar,\ubarboot]}{[\moreinterestingu_1,\moreinterestingu_2]}$,
	even if $\ubarboot = 0$. We will achieve this by
	proving the injectivity of the map
	$(\ubar,u,x^2,x^3) \rightarrow (\ubar,x^1,x^2,x^3)$
	on the domain
	$[\leftubar,\ubarboot] \times [\moreinterestingu_1,\moreinterestingu_2] \times \mathbb{T}^2$,
	and then the injectivity of the map 
	$(\ubar,x^1,x^2,x^3) \rightarrow (t,x^1,x^2,x^3)$;
	since the composition of two injective functions is injective, this 
	would finish the proof.

	To proceed, we first note that the injectivity of the map 
	$(\ubar,u,x^2,x^3) \rightarrow (\ubar,x^1,x^2,x^3)$
	on
	$[\leftubar,\ubarboot] \times [\moreinterestingu_1,\moreinterestingu_2] \times \mathbb{T}^2$
	follows from the monotonicity of the map $u \rightarrow x^1(\ubar,u,x^2,x^3)$
	guaranteed by Lemma\,\ref{L:MONOTONICITYOFCARTESIANX1}.
	For use below, we also note that by \eqref{E:LSMALLDEF}, 
	\eqref{E:NULLCOORDINATEPARTIALINGOINGEIKONALINTERMSOFGEOMETRICVF},
	Lemma\,\ref{L:COMMUTATORSTOCOORDINATES},
	\eqref{E:CLOSEDVERSIONKEYJACOBIANDETERMINANTESTIMATECHOVGEOTODOUBLENULL},
	and the estimates of Prop.\,\ref{P:IMPROVEMENTOFAUXILIARYBOOTSTRAP},
	we have:
	\begin{align} \label{E:DOUBLENULLUBARDERIVATIVEOFCARTESIANX1ISAPPROXIMATELYUNITY}
	\nullgeop{\ubar} x^1
	\approx \geop{t} x^1 
	& = 
	\Lunit x^1 
	-
	\Lunit^A \geop{x^A} x^1
	= 1
	+ 
	\Lsmall^1
	-
	\Lunit^A \geop{x^A} x^1
	\approx 1.
	\end{align}
	
	Now for each fixed $(x^2,x^3) \in \mathbb{T}^2$,
	let $\mathcal{I}_{x^2,x^3}$ denote the image of the set
	$[\leftubar,\ubarboot] \times [\moreinterestingu_1,\moreinterestingu_2]\times \lbrace (x^2,x^3) \rbrace$
	under the map $(\ubar,u,x^2,x^3) \rightarrow (\ubar,x^1,x^2,x^3)$.
	The arguments given above, including the monotonicity guaranteed by \eqref{E:DOUBLENULLUBARDERIVATIVEOFCARTESIANX1ISAPPROXIMATELYUNITY}, 
	imply that for each fixed $(u,x^2,x^3) \in [\moreinterestingu_1,\moreinterestingu_2]\times \mathbb{T}^2$,
	the map $\ubar \rightarrow x^1(\ubar,u,x^2,x^3)$
	is strictly increasing on $[\leftubar,\ubarboot]$, and that 
	for each fixed $(\ubar,x^2,x^3) \in [\leftubar,\ubarboot] \times \mathbb{T}^2$,
	the map
	$u \rightarrow x^1(\ubar,u,x^2,x^3)$ is strictly decreasing on $ [\moreinterestingu_1,\moreinterestingu_2]$.
	It follows that there exist scalar functions $\ubar \rightarrow \mathfrak{a}_{x^2,x^3}(\ubar)$
	and
	$\ubar \rightarrow \mathfrak{b}_{x^2,x^3}(\ubar)$
	on $[\leftubar,\ubarboot]$
	such that
	$\mathcal{I}_{x^2,x^3}
	=
	\lbrace (\ubar,x^1,x^2,x^3) \ | \ 
	\ubar \in [\leftubar,\ubarboot],
		\,
	\mathfrak{a}_{x^2,x^3}(\ubar) \leq x^1 \leq \mathfrak{b}_{x^2,x^3}(\ubar) \rbrace
	$,
	where $\mathfrak{a}_{x^2,x^3}(\cdot)$ and $\mathfrak{b}_{x^2,x^3}(\cdot)$ are $C^1$ functions of $\ubar$
	such that $\mathfrak{a}_{x^2,x^3}(\ubar) < \mathfrak{b}_{x^2,x^3}(\ubar)$ and
	$
	\frac{d}{d \ubar} \mathfrak{a}_{x^2,x^3}, 
		\,
	\frac{d}{d \ubar} \mathfrak{b}_{x^2,x^3}
	\approx 1
	$.

	To complete the proof, it remains for us to show that the map
	$(\ubar,x^1,x^2,x^3) \rightarrow (t,x^1,x^2,x^3)$
	is injective. Let $\mathcal{I}_{x^2,x^3}$ be the set from the previous paragraph.
	It suffices to show that for each fixed $(x^2,x^3) \in \mathbb{T}^2$,
	any two distinct points in $\mathcal{I}_{x^2,x^3}$ with the same $x^1$ coordinate
	must be mapped to distinct points under the map
	$(\ubar,x^1,x^2,x^3) \rightarrow (t,x^1,x^2,x^3)$.
	The structure of $\mathcal{I}_{x^2,x^3}$ revealed in the previous paragraph
	shows that for any two distinct points in $\mathcal{I}_{x^2,x^3}$ with the same $x^1$ coordinate,
	the straight line segment joining them
	(along which $x^1,x^2,x^3$ are constant and $\ubar$ varies)
	is contained in $\mathcal{I}_{x^2,x^3}$
	(this can be thought of as the vertical convexity of $\mathcal{I}_{x^2,x^3}$).
	Hence, to complete the proof, it suffices for us to show that the partial derivative of $t$
	with respect to $\ubar$ in the coordinate system
	$(\ubar,x^1,x^2,x^3)$ 
	is positive, except possibly when $\ubar = 0$.
	To proceed, we use 
	\eqref{E:LSMALLDEF},
	\eqref{E:XSMALL},
	\eqref{E:CARTESIANPARTIALTTOCOMMUTATORS}, 
	Lemma\,\ref{L:COMMUTATORSTOCOORDINATES}, 
	Prop.\,\ref{P:SCHEMATICSTRUCTUREOFVARIOUSTENSORSINTERMSOFCONTROLVARS}, \eqref{E:EIKONALEQUATIONWITHOUTTHEMUWEIGHT}, 
	and the estimates of 
	Lemma\,\ref{L:PROPERTIESANDDIFFEOMORPHICEXTENSIONOFDOUBLENULLCOORDINATES},
	and
	Prop.\,\ref{P:IMPROVEMENTOFAUXILIARYBOOTSTRAP},
	\eqref{E:BOUNDSONLMUINTERESTINGREGION},
	and \eqref{E:BOUNDSOFLUBARINTERESTINGREGION}
	to compute that
	$\partial_t = \Lunit + (1 + *) X + * \Yvf{2} + * \Yvf{3}$
	and that
	$\partial_t \ubar 
	\approx
	 \Lunit \ubar \approx 1
	$. We have therefore shown that the map $(\ubar ,x^1,x^2,x^3) \rightarrow (t,x^1,x^2,x^3)$
	is injective, as desired.

\end{proof}

\subsection{Pointwise estimates tied to the double-null geometry}
\label{SS:POINTWISEESTIMATESTIEDTODOUBLENULLACOUSTICGEOMETRY}
In this short section, we provide a lemma where we derive some useful pointwise estimates regarding the double-null charactacteristic geometry. 
Some of the estimates in the lemma,
such as the estimates
\eqref{E:POINTWISEBOUNDFORDOUBLENULLOROIDALTRACEOFDEFORMATIONTENSOROFNEWLVECTORFIELD}--\eqref{E:POINTWISEBOUNDFORDOUBLENULLTOROIDALTRACEOFDEFORMATIONTENSOROFNEWULVECTORFIELD}
for the deformation tensors of
$\newL$
and
$\newuL$, are necessary to estimate some identities related to integration by parts, see \eqref{E:IDENTITYUBARDERIVATIVEOFDOUBLENULLORUSINTEGRAL}--\eqref{E:IDENTITYEIKONALDERIVATIVEOFDOUBLENULLTORUSINTEGRAL}

\begin{lemma}[Pointwise estimates tied to the double-null acoustic geometry]
\label{L:POINTWISEESTIMATESTIEDTODOUBLENULLACOUSTICGEOMETRY}
Let 
$\deform{\newL}, \, \deform{\newuL}$ be the deformation tensors of the vectorfields 
$\newL, \, \newuL$ defined in \eqref{E:DOUBLENULLEIKONALFUNCTIONNORMALIZEDNULLVECTORFIELDS},
let $\mytr_{\gnulltori} \deform{\newL}$
and
$\mytr_{\gnulltori} \deform{\newuL}$ 
respectively denote their traces with respect to $\gnulltori$,
and let $\nullangDiv \ToriTangentVectorfieldAssociatedToDoubleNullFolliations$ be the vectorfield $\ToriTangentVectorfieldAssociatedToDoubleNullFolliations$
defined in \eqref{E:TORITANGENTVECTORFIELDASSOCIATEDTODOUBLENULLFRAME}.
Then the following pointwise estimates hold on 
$\characteristicdiamondtwoarg{[\leftubar,\ubarboot)}{[\moreinterestingu_1,\moreinterestingu_2]}$:
	\begin{align} 
			\sum_{\alpha = 0,1,2,3} |\newuL^\alpha| & \lesssim \upmu, \label{E:POINTWISEESTIMATEFORCARTESIANCOMPONENTSOFNEWUL}
		\\
			\sum_{\alpha = 0,1,2,3} |\ToriTangentVectorfieldAssociatedToDoubleNullFolliations^\alpha| & \lesssim \fundbootsmall, \label{E:POINTWISEESTIMATEFORCARTESIANCOMPONENTSOFDOUBLENULLTORITANGENTVECTORFIELD}
		\\
			|\mytr_{\gnulltori} \deform{\newL}|
			& \lesssim 1,
			\label{E:POINTWISEBOUNDFORDOUBLENULLOROIDALTRACEOFDEFORMATIONTENSOROFNEWLVECTORFIELD}
		\\
			|\mytr_{\gnulltori} \deform{\newuL}|
			& \lesssim 1.
			\label{E:POINTWISEBOUNDFORDOUBLENULLTOROIDALTRACEOFDEFORMATIONTENSOROFNEWULVECTORFIELD}
	\end{align}
 \end{lemma}
 
 \begin{proof} \hfill

\medskip

\noindent \textbf{Proof of \eqref{E:POINTWISEESTIMATEFORCARTESIANCOMPONENTSOFNEWUL}--\eqref{E:POINTWISEESTIMATEFORCARTESIANCOMPONENTSOFDOUBLENULLTORITANGENTVECTORFIELD}:} We first decompose the Cartesian components of $\newuL$ using \eqref{E:CONVENIENTIDENTITYFORULUNIT} as $\newuL^\alpha = \ReciprocaluLunitAppliedtoTimeFunction L^\alpha + \muX^\alpha - \MagnitueofinnerproductofnewLandnewuL \angD^\alpha \ubar \left(\geop{x^A}\right)^\alpha$. Using Prop.\,\ref{P:SCHEMATICSTRUCTUREOFVARIOUSTENSORSINTERMSOFCONTROLVARS},  and \eqref{E:MAGNITUDEOFINNERPRODUCTOFNEWLANDNEWULAPPROXIMATELYMU}--\eqref{E:RECIPROCALULUNITAPPLIEDTOTIMEFUNCTIONAPPROXIMATELYMU}, the desired estimate \eqref{E:POINTWISEESTIMATEFORCARTESIANCOMPONENTSOFNEWUL} follows. \eqref{E:POINTWISEESTIMATEFORCARTESIANCOMPONENTSOFDOUBLENULLTORITANGENTVECTORFIELD} follows from similar arguments using \eqref{E:TORITANGENTVECTORFIELDASSOCIATEDTODOUBLENULLFRAME}, \eqref{E:NULLCOORDINATEPARTIALXAINTERMSOFGEOMETRICVF}, and \eqref{E:SMALLC01ESTIMATESFORUBAR}.

\medskip
 
\noindent \textbf{Proof of \eqref{E:POINTWISEBOUNDFORDOUBLENULLOROIDALTRACEOFDEFORMATIONTENSOROFNEWLVECTORFIELD}:}

First, using \eqref{E:GNULLTORIINVERSECOMPONENDS} and \eqref{E:LIEGNULLTORUS}, we compute that relative to the double-null coordinates $(\ubar,u,x^2,x^3)$, we have $\mytr_{\gnulltori} \deform{\newuL} = (\gnulltori^{-1})(\mathrm{d} x^A, \mathrm{d} x^B) \newuL \gnulltori\left(\nullgeop{x^A},\nullgeop{x^B}\right) + 2 \nullgeop{x^A} \newuL^A$. We also compute $\newuL^A = \ReciprocaluLunitAppliedtoTimeFunction L^A + \upmu X^A - \MagnitueofinnerproductofnewLandnewuL \angD^A \ubar$ using \eqref{E:CONVENIENTIDENTITYFORULUNIT}. Next, we use \eqref{E:CHOVCOEFFICIENTSSMOOTHANGULARDERIVATIVESINTERMSOFDOUBLENULLONESANDL} and \eqref{E:TORITANGENTVECTORFIELDASSOCIATEDTODOUBLENULLFRAME} to write $\ToriTangentVectorfieldAssociatedToDoubleNullFolliations^A = (\gtorus^{-1})^{AB}  \geop{x^B} \ubar - \frac{|\angD \ubar|^2}{\Lunit \ubar} L^A$. Next we also use \eqref{E:CHOVCOEFFICIENTSSMOOTHANGULARDERIVATIVESINTERMSOFDOUBLENULLONESANDL}, \eqref{E:SMOOTHTORUSINVERSEFIRSTFUNDCOMPONENTSINTERMSOFDOUBLENULLTORUSINVERSEFIRSTFUNDCOMPONENTS}, and Props.\,\ref{P:SCHEMATICSTRUCTUREOFVARIOUSTENSORSINTERMSOFCONTROLVARS} and\,\ref{P:IMPROVEMENTOFAUXILIARYBOOTSTRAP}
to deduce the schematic identities
$\gnulltori\left(\nullgeop{x^A},\nullgeop{x^B}\right)
=
\smoothfunction
\left(\controlvars,
	\frac{1}{\geop{t}\ubar}, 
	\geop{x^2}\ubar, \geop{x^3} \ubar
\right)
$
and
$(\gtorusroughinversefirstfund)(\rmd x^A,\rmd x^B) 
=
\smoothfunction
\left(\controlvars,
	\frac{1}{\geop{t}\ubar}, 
	\geop{x^2}\ubar, \geop{x^3} \ubar
\right)$.
Considering these identities, \eqref{E:CONVENIENTIDENTITYFORULUNIT}, we deduce the following schematic identity:
	\begin{align}
		\mytr_{\gnulltori} \deform{\newuL}  = \smoothfunction\left( \comder \controlvars, \frac{1}{\geop{t}\ubar}, \tander^{[1,2]}\ubar, \newuL \tander \ubar\right). \label{E:SCHEMATICIDENTITYFORDOUBLENULLTORUSTRACEOFDEFORMATIONTENSOROFNEWUL}
	\end{align}
From \eqref{E:SCHEMATICIDENTITYFORDOUBLENULLTORUSTRACEOFDEFORMATIONTENSOROFNEWUL}, Lemma\,\ref{L:COMMUTATORSTOCOORDINATES}, Props.\,\ref{P:SCHEMATICSTRUCTUREOFVARIOUSTENSORSINTERMSOFCONTROLVARS} and\,\ref{P:IMPROVEMENTOFAUXILIARYBOOTSTRAP}, the estimates \eqref{E:ALLDERIVATIVESLINFTYESTIMATESFORUBARANDDERIVATIVES} and \eqref{E:CLOSEDVERSIONKEYJACOBIANDETERMINANTESTIMATECHOVGEOTODOUBLENULL}, and explicitly bounding the RHS\,\eqref{E:QUASILINEAREQUATIONSATISFIEDBYLUBAR},\eqref{E:QUASILINEAREQUATIONSATISFIEDBYSMOOTHANGULARDERIVATIVEOFUBAR} to estimate $\newuL \tander \ubar$, we arrive to estimate \eqref{E:POINTWISEBOUNDFORDOUBLENULLOROIDALTRACEOFDEFORMATIONTENSOROFNEWLVECTORFIELD}.

\medskip
 
\noindent \textbf{Proof of \eqref{E:POINTWISEBOUNDFORDOUBLENULLTOROIDALTRACEOFDEFORMATIONTENSOROFNEWULVECTORFIELD}:} 

Considering the identity $\newL = \frac{1}{\Lunit \ubar} \Lunit$ and arguing as in the proof of \eqref{E:SCHEMATICIDENTITYFORDOUBLENULLTORUSTRACEOFDEFORMATIONTENSOROFNEWUL}, we find that:
	\begin{align}
		\mytr_{\gnulltori} \deform{\newL}  = \smoothfunction\left( \tander^{\le 1} \controlvars, \frac{1}{\Lunit\ubar} \frac{1}{\geop{t}\ubar}, \tander^{[1,2]}\ubar \right). \label{E:SCHEMATICIDENTITYFORDOUBLENULLTORUSTRACEOFDEFORMATIONTENSOROFNEWL}
	\end{align}
	From \eqref{E:SCHEMATICIDENTITYFORDOUBLENULLTORUSTRACEOFDEFORMATIONTENSOROFNEWL}, Lemma\,\ref{L:COMMUTATORSTOCOORDINATES}, Props.\,\ref{P:SCHEMATICSTRUCTUREOFVARIOUSTENSORSINTERMSOFCONTROLVARS} and\,\ref{P:IMPROVEMENTOFAUXILIARYBOOTSTRAP}, the estimates \eqref{E:ALLDERIVATIVESLINFTYESTIMATESFORUBARANDDERIVATIVES}, \eqref{E:CLOSEDVERSIONKEYJACOBIANDETERMINANTESTIMATECHOVGEOTODOUBLENULL}, and \eqref{E:BOUNDSOFLUBARINTERESTINGREGION}, we arrive to the desired estimate \eqref{E:POINTWISEBOUNDFORDOUBLENULLTOROIDALTRACEOFDEFORMATIONTENSOROFNEWULVECTORFIELD}. 
 \end{proof}

 \section{Elliptic estimates for the characteristic geometry on the double-null tori} \label{S:ELLIPTICESTIMATESFORTHEACOUSTICGEOMETRYONTHEDOUBLENULLTORI}

We continue to work under the assumptions of Sect.\,\ref{SS:SILENTFACTS}.
In this section, we derive elliptic 
$L^2$ estimates for symmetric $\binom{0}{2}$-type tensorfields that are tangent to the geometric tori $\ell_{t,u}$. 
A key point is that the geometric tori are not adapted to the double-null foliations and hence the estimates will be derived on the  \emph{double-null} tori 
$\doublenulltoritwoarg{\ubar}{u}$. The main difficulty that arises from this is that the Gauss curvature of the double-null tori, which appears in the elliptic estimates and depends on \emph{three} derivatives of $\ubar$, is seemingly uncontrollable from the $C^{1,1}$ estimates derived in Lemma\,\ref{L:PROPERTIESANDDIFFEOMORPHICEXTENSIONOFDOUBLENULLCOORDINATES}. In Lemma\,\ref{L:GAUSSCURVATUREOFDOUBLENULLTORILINFINITYESTIMATE}, we will prove that all third order derivatives of $\ubar$ present in the Gauss curvature \emph{completely cancel}. Although similar cancellations hold in our prior work \cite{abbrescia2022emergence}, we did not have the need to derive them there because the rough time function was of regularity class $C^{2,1}$.

To obtain the elliptic estimates, we decompose
symmetric type $\binom{0}{2}$ $\ell_{t,u}$-tangent tensorfields $\upxi$ into 
a main piece that is tangent to $\doublenulltoritwoarg{\ubar}{u}$,
which we control with elliptic estimates on the double-null tori (see Lemma\,\ref{L:STANDARDELLIPTICONDOUBLENULLTORI}),
and error terms, which we must control with separate (easier) arguments.
Our primary application will be to apply the main elliptic estimate \eqref{E:ELLIPTICESTIMATECHI}
with $\upxi = \angLie_{\tander}^{\Ntop-1} \upchi$,
which will yield $L^2$-control of $\angLie_{\tander}^{\Ntop} \upchi$;
see the proof of \eqref{E:PRELIMINARYTOPORDERLIECHI}.
When $N = \Ntop$, these $L^2$ estimates are crucial for 
avoiding the loss of a derivative in the terms 
on RHSs \eqref{E:TOPCOMMUTEDWAVELFIRSTTHENALLYS}--\eqref{E:TOPCOMMUTEDWAVEALLYS}
that explicitly depend on the order $N$ derivatives $\mytr_{\gtorus}\upchi$.

We carried out the analogous tensorial decompositions described in the previous paragraph in our prior work \cite{abbrescia2022emergence}*{Sect.\,29}, with the rough tori  $\twoargroughtori{\timefunction,u}{0}$ in place of the double-null tori $\doublenulltoritwoarg{\ubar}{u}$. Since the calculations are nearly identical, and the main elliptic identity is well known, we are quite terse and focus instead on the cancellations in the Gauss curvature estimates. We clarify that the main geometric structure needed is that \eqref{E:SMOOTHANGULARDERIVATIVESINTERMSOFDOUBLENULLONESANDL} implies that the double-null tori coordinate vectorfields $\nullgeop{x^A}$ can be schematically decomposed into $(1 + \mathcal{O}(\fundbootsmall))\geop{x^A} + \mathcal{O}(\fundbootsmall)L$, and the same schematic decomposition occurred in \cite{abbrescia2022emergence}: $\roughgeop{x^A} = (1 + \mathcal{O}(\fundbootsmall))\geop{x^A} + \mathcal{O}(\fundbootsmall)L$.

\subsection{Statement of the main elliptic estimates}
\label{SS:ELLIPTICESTIMATESFORACOUSTICGEOMETRY}
We now state Prop.\,\ref{P:ELLIPTICESTIMATESONDOUBLENULLTORI},
which yields the main elliptic estimates.
Its proof is located in Sect.\,\ref{SS:PROOFOFP:ELLIPTICESTIMATESONDOUBLENULLTORI}.

\begin{proposition}[The main elliptic estimates for symmetric type $\binom{0}{2}$ $\ell_{t,u}$-tangent tensorfields] 
\label{P:ELLIPTICESTIMATESONDOUBLENULLTORI}
Let $\upxi$ be a symmetric type $\binom{0}{2}$ 
$\ell_{t,u}$-tangent tensorfield. 
Then the following estimate holds
for $(\ubar,u) \in [\leftubar,\ubarboot) \times [\moreinterestingu_1,\moreinterestingu_2]$:
\begin{align}
\begin{split} \label{E:ELLIPTICESTIMATECHI}
\sum_{\Singletan \in \{\Lunit,\Yvf{2},\Yvf{3}\}} 
  \int_{\characteristicdiamondtwoarg{[\leftubar,\ubar]}{[\moreinterestingu_1,u]}} 	
 	\upmu^2 |\angLie_{\Singletan} \upxi|_{\gtorus}^2 \weight^{\blowupratetoporderwave}
\, \voldiamond 
& 
\leq 
C \int_{\characteristicdiamondtwoarg{[\leftubar,\ubar]}{[\moreinterestingu_1,u]}} 	
  	\upmu^2 |\angLie_{\Lunit} \upxi|_{\gtorus}^2 \weight^{\blowupratetoporderwave}
\, \voldiamond 
+ 
C 
  \int_{\characteristicdiamondtwoarg{[\leftubar,\ubar]}{[\moreinterestingu_1,u]}} 	
	\upmu^2 |\angdiv \upxi|_{\gtorus}^2 \weight^{\blowupratetoporderwave}
\, \voldiamond 
	\\
& \ \
+ 
C 
\sum_{A=2,3} 
 \int_{\characteristicdiamondtwoarg{[\leftubar,\ubar]}{[\moreinterestingu_1,u]}} 	
 	\upmu^2 (\Yvf{A} \mytr_{\gtorus}\upxi)^2 \weight^{\blowupratetoporderwave}
\, \voldiamond \\
& \ \ 
+ 
C 
\fundbootsmall 
 \int_{\characteristicdiamondtwoarg{[\leftubar,\ubar]}{[\moreinterestingu_1,u]}} 	
 	|\upxi|_{\gtorus}^2 \weight^{\blowupratetoporderwave}
\, \voldiamond. 
\end{split}
\end{align}
\end{proposition}

\subsection{Differential operator pointwise identities and comparison estimates used in the proof of the elliptic estimates}
\label{SS:DIFFOPERATORSPOINTWISEIDENTITIESANDCOMPARISONNEEDEDFORELLIPTICESTIMATESONDOUBLENULLTORI}

\begin{lemma}[Identities involving contractions with $\gtorus$ and $\gnulltori$]
\label{L:SMOOTHANDDOUBLENULLFIRSTFUNDFORMSAGREEONPUTANGENTENSORS}
Let $\upeta$ be a $\nullhyparg{u}$-tangent tensorfield.
Recall that $|\cdot|_{\gtorus}$ and $|\cdot|_{\gnulltori}$
are defined in \eqref{E:SQUAREPOINTWISESEMINORMWITHRESPECTTOFIRSTFUNDOFSMOOTHTORI} and
\eqref{E:SQUAREPOINTWISESEMINORMWITHRESPECTTOFIRSTFUNDOFDOUBLENULLTORI}
respectively. Then the following identity holds:
\begin{align} \label{E:SEMINNORMSOFPUTANGENTTENSORSINDEPENDENTOFMETRIC}
|\upeta|_{\gtorus}  
& 
= |\upeta|_{\gnulltori}. 
\end{align}In addition, if $\upxi$ is a symmetric type $\binom{0}{2}$ $\ell_{t,u}$-tangent tensorfield, let $\breve{\upxi} \eqdef \nulltorusproject \upxi$ denote its $\gfour$-orthogonal projection onto $\doublenulltoritwoarg{\ubar}{u}$ defined in \eqref{E:TENSORFIELDPROJECTIONTONULLTORIDEFININGEQUATION}. Then we have: 
\begin{align}
|\upxi|_{\gtorus} & = \left| \breve{\upxi}\right|_{\gnulltori}. 
\label{E:NORMSOFDOUBLENULLTORICOMPONENTSANDFLATTORIAREEQUAL}
\end{align}In addition, if $\upeta$ is a symmetric type $\binom{0}{2}$ $\nullhyparg{u}$-tangent tensorfield, then with $\mytr_{\gtorus}$ and $\mytr_{\gnulltori}$ as in
\eqref{E:TRACES}
and 
\eqref{E:TRACEOFNULLTORITANGENT02TENSORS}, we have:
\begin{align} \label{E:TRACESOFPUTANGENTTENSORINDEPENDENTOFMETRIC} 
\mytr_{\gtorus} \upeta 
& = 
\mytr_{\gnulltori} \upeta. 
\end{align}Finally, if $\upxi$ and $\breve{\upxi}$ are as above then:
\begin{align} \label{E:TRACESOFDOUBLENULLTORICOMPONENTSANDFLATTORIAREEQUAL} 
\mytr_{\gtorus} \upxi 
& = 
\mytr_{\gnulltori} \breve{\upxi}. 
\end{align}
	
\end{lemma}

\begin{lemma}[Differential operator pointwise comparison estimates needed for the elliptic estimates]
\label{L:DIFFOPERATORSPOINTWISECOMPARISONNEEDEDFORELLIPTICESTIMATESONDOUBLENULLTORI}
Let $\varphi$ be a scalar function, 
let $\angrmd \varphi$ be the $\ell_{t,u}$-tangent one-form from Def.\,\ref{D:ANGULARDIFFERENTIAL},
and let $\nullangrmD \varphi$ be the $\doublenulltoritwoarg{\ubar}{u}$-tangent one-form
from Def.\,\ref{D:DOUBLENULLORTHOGONALPROJECTIONANDTANGENCY}.
Let $\upxi$ be a symmetric type $\binom{0}{2}$ $\ell_{t,u}$-tangent tensorfield,
and let $\breve{\upxi} \eqdef \nulltorusproject \upxi$ denote its $\gfour$-orthogonal projection onto $\doublenulltoritwoarg{\ubar}{u}$ defined in \eqref{E:TENSORFIELDPROJECTIONTONULLTORIDEFININGEQUATION}.
Then the following pointwise estimates hold
on $ \characteristicdiamondtwoarg{[\leftubar,\ubarboot)}{[\moreinterestingu_1,\moreinterestingu_2]}$: 
\begin{subequations}
\begin{align}
	|\nullangrmD \varphi|_{\gnulltori} 
	&
	=
	\left\lbrace
		1 + \mathcal{O}(\fundbootsmall)
	\right\rbrace
	|\angrmd \varphi|_{\gtorus}
	+ 
	\mathcal{O}(\fundbootsmall)
	|\Lunit \varphi|,
		\label{E:DOUBLENULLTORUSDIFFERENTIALAPPROXIMATEDBYSMOOTHTORUSDIFFERENTIAL} 
			\\
	|\nullangDiv \breve{\upxi}|_{\gnulltori}
	& =
	\left\lbrace
		1 + \mathcal{O}(\fundbootsmall)
	\right\rbrace
	|\angdiv \upxi|_{\gtorus}
	+
	\mathcal{O}(\fundbootsmall)
	|\angLie_{\Lunit} \upxi|_{\gtorus}
	+
	\mathcal{O}(\fundbootsmall)
	|\upxi|_{\gtorus},
		\label{E:DOUBLENULLTORUSDIVERGENCEAPPROXIMATEDBYSMOOTHTORUSDIVERGENCE} 
			\\
	|\nullangD \breve{\upxi}|_{\gnulltori}
	& =
	\left\lbrace
		1 + \mathcal{O}(\fundbootsmall)
	\right\rbrace
	|\angD \upxi|_{\gtorus}
	+
	\mathcal{O}(\fundbootsmall)
	|\angLie_{\Lunit} \upxi|_{\gtorus}
	+
	\mathcal{O}(\fundbootsmall)
	|\upxi|_{\gtorus}.
	\label{E:DOUBLENULLTORUSCOVARIANTDERIVATIVEAPPROXIMATEDBYSMOOTHTORUSCOVARIANTDERIVATIVE}
\end{align}
\end{subequations}\end{lemma}

\begin{proof}[Sketch of proof of Lemmas\,\ref{L:SMOOTHANDDOUBLENULLFIRSTFUNDFORMSAGREEONPUTANGENTENSORS}--\ref{L:DIFFOPERATORSPOINTWISECOMPARISONNEEDEDFORELLIPTICESTIMATESONDOUBLENULLTORI}]

The identical proofs of \cite{abbrescia2022emergence}*{Lemmas\,29.8 and 29.9} with $\twoargroughtori{\timefunction,u}{0}$ in place of the double-null tori $\doublenulltoritwoarg{\ubar}{u}$ and $\roughgeop{x^A}$ in place of $\nullgeop{x^A}$ hold. For the reader's convenience, we give a brief sketch. 

It is convenient to define $\{e_A\}_{A = 2,3}$ as the orthonormal frame on the geometric torus $\ell_{t,u}$ 
obtained from applying the Gram--Schmidt process to the geometric coordinate partial derivative vectorfields
$\left\{\geop{x^A} \right\}_{A = 2,3}$ with respect to $\gtorus$, starting with 
$e_2 \eqdef \frac{1}{\sqrt{\gtorus(\geop{x^2},\geop{x^2})}}\geop{x^2}$.
Similarly, $\left\{f_A \right\}_{A = 2,3}$ is defined to be 
the orthonormal frame on the double-null torus $\doublenulltoritwoarg{\ubar}{u}$ 
obtained from applying the Gram--Schmidt process to the adapted rough geometric coordinate partial derivative vectorfields
$\left\{\nullgeop{x^A} \right\}_{A = 2,3}$ with respect to $\gnulltori$,
starting with 
$f_2 \eqdef \frac{1}{\sqrt{\gnulltori(\nullgeop{x^2},\nullgeop{x^2})}}\nullgeop{x^2}$.

From the above, one can prove that there are $2\times 2$ orthogonal matrix-valued functions $\COVframe$ with components $\{\COVframe_{AB}\}_{A,B = 2,3}$ 
and scalar functions $\{\COVL_A \}_{A = 2,3}$ such that $e_A = 
\COVframe_{AB} f_B 
+ 
\COVL_A \Lunit$. Moreover, one can prove the estimates $|\COVframe_{AB}| < 1$ and $|\COVL| \lesssim \fundbootsmall$. 

Denoting the $\binom{2}{0}$ $\gfour$-duals of $\upxi$ and $\breve{\upxi}$ as $\upxi^{\#\#}$ and $\breve{\upxi}^{\#\#}$, straightforward multi-linear algebra implies that the difference $\breve{\upxi}^{\#\#} - \upxi^{\#\#}$ can be expressed as a linear combination of $\Lunit \otimes e_A, \, e_A \otimes \Lunit$, and $\Lunit\otimes\Lunit$ with small coefficients. The rest of the proof follows from straightforward calculations based on the identities $\upeta(\Lunit,\cdot) = \upxi(\Lunit,\cdot) = 0$, the Lie differentiation identity
$
	\Lie_{\Lunit} \upeta_{\alpha \beta}
	= 
	\Dfour_{\Lunit} \upeta_{\alpha \beta}
	+
	\upeta_{\kappa \beta} \Dfour_{\alpha} \Lunit^{\kappa}
	+
	\upeta_{\alpha \kappa} \Dfour_{\beta} \Lunit^{\kappa}
$, and $\Dfour_{e_A} \Lunit = \upchi(e_A,e_C) e_C - \upzeta_A \Lunit$ with  $\upzeta_A$ as in \eqref{E:TORISONTENSORFIELD}.

\end{proof}

\subsection{Standard elliptic estimates for symmetric type $\binom{0}{2}$ tensorfields on the double-null tori}
\label{SS:STANDARDELLIPTICONDOUBLENULLTORI}
In the next lemma, we provide standard elliptic estimates for
symmetric type $\binom{0}{2}$ $\doublenulltoritwoarg{\ubar}{u}$-tangent tensorfields.
Its proof is located in Sect.\,\ref{SSS:PROOFOFL:STANDARDELLIPTICONDOUBLENULLTORI}.

\begin{lemma}[Standard elliptic estimates for symmetric type $\binom{0}{2}$ tensorfields on the double-null tori]
\label{L:STANDARDELLIPTICONDOUBLENULLTORI}
Let $\Upxi$ be a symmetric type $\binom{0}{2}$ $\doublenulltoritwoarg{\ubar}{u}$-tangent tensorfield.
Then the following estimate
holds for $(\ubar,u) \in [\leftubar,\ubarboot) \times [\moreinterestingu_1,\moreinterestingu_2]$:
\begin{align} \label{E:STANDARDELLIPTICONDOUBLENULLTORI}
\int_{\doublenulltoritwoarg{\ubar}{u}} 
	\upmu^2 |\nullangD \Upxi|_{\gnulltori}^2 
\, \voldoublenulltori
& \leq 
6 
\int_{\doublenulltoritwoarg{\ubar}{u}}
	\upmu^2 |\nullangDiv \Upxi|_{\gnulltori}^2 
\, \voldoublenulltori
+ 
3 
\int_{\doublenulltoritwoarg{\ubar}{u}} 
	\upmu^2 |\nullangrmD \mytr_{\gnulltori} \Upxi|_{\gnulltori}^2
\, \voldoublenulltori
+ 
C 
\fundbootsmall 
\int_{\doublenulltoritwoarg{\ubar}{u}}
	|\Upxi|_{\gnulltori}^2
\, \voldoublenulltori.
\end{align}\end{lemma}

\subsubsection{The Gauss curvature of $\gnulltori$}
\label{SSS:GAUSSCURVATUREOFROUGHTORI}
To prove Lemma\,\ref{L:STANDARDELLIPTICONDOUBLENULLTORI},
we use the following Gauss curvature estimate.

\begin{lemma}[$L^{\infty}$ estimate for the Gauss curvature of $\gnulltori$] 
\label{L:GAUSSCURVATUREOFDOUBLENULLTORILINFINITYESTIMATE}
Recall that $\Gaussnulltori$ denotes the Gauss curvature of $(\doublenulltoritwoarg{\ubar}{u},\gtorusroughfirstfund)$
(see Sect.\,\ref{SS:CURVATURETENSORS}).
Then the following estimate holds for
$(\ubar,u) \in [\leftubar,\ubarboot] \times [\moreinterestingu_1,\moreinterestingu_2]$:
\begin{align} \label{E:LINFINITYBOUNDFORDOUBLENULLGAUSSCURVATURE}
\left\| \Gaussnulltori \right \|_{L^{\infty}\left(\doublenulltoritwoarg{\ubar}{u}\right)} 
& \leq C \fundbootsmall.  
\end{align}
\end{lemma}

\begin{proof}
Throughout this proof, for convenience, we denote $\gnulltori_{AB} \eqdef \gnulltori(\nullgeop{x^A},\nullgeop{x^B})$ and $(\gnulltori^{-1})^{AB} \eqdef \gnulltori^{-1}(\rmD x^A,\rmD x^B)$. 

From \eqref{E:RIEMANNCURVATURETENSOROFDOUBLENULLTORUS}, the convention $\nullangD_{XY}^2 Z= X^\alpha Y^\beta \nullangD_\alpha \nullangD_\beta Z$, and $[\nullgeop{x^A},\nullgeop{x^B}] = 0$, we have:
	\begin{align}
		\begin{split} \label{E:GAUSSCURVATUREESTIMATEINTERMEDIATESTEP1}
			\Riemnulltori\left(\nullgeop{x^A},\nullgeop{x^B},\nullgeop{x^C},\nullgeop{x^D}\right) & = \gnulltori\left( - \nullangD_A \nullangD_B \nullgeop{x^C}  + \nullangD_B \nullangD_A \nullgeop{x^C},\nullgeop{x^D}\right) \\
			& = \gnulltori_{ED} \left( \nullgeop{x^B} \nullangChristoff_{AC}^E - \nullgeop{x^A} \nullangChristoff_{BC}^E + \nullangChristoff_{BF}^E \nullangChristoff_{AC^F} - \nullangChristoff_{AF}^E \nullangChristoff_{BC}^F\right),
		\end{split}
	\end{align}
where $\nullangChristoff_{AB}^C = \frac{1}{2} (\gnulltori^{-1})^{CD} \left\{ \nullgeop{x^A} \gnulltori_{DB} + \nullgeop{x^B} \gnulltori_{DA} - \nullgeop{x^D} \gnulltori_{AB}\right\}$ are the Christoffel symbols of $\nullangD$ expressed in coordinates $(x^2,x^3)$ on $\doublenulltoritwoarg{\ubar}{u}$. From \eqref{E:NULLCOORDINATEPARTIALXAINTERMSOFGEOMETRICVF}, \eqref{E:GNULLTORICOMPONENTS}, and Prop.\,\ref{P:SCHEMATICSTRUCTUREOFVARIOUSTENSORSINTERMSOFCONTROLVARS}, we have the following schematic structure for the Christoffel symbols: 
	\begin{align}
		\nullangChristoff_{AB}^C = \smoothfunction\left(\tander \controlvars, \frac{1}{\Lunit \ubar}, \frac{1}{\geop{t} \ubar}, \tander \ubar\right) \cdot \left(\tander \controlvars, \tander^{[1,2]} \ubar\right).  \label{E:GAUSSCURVATUREESTIMATEINTERMEDIATESTEP2}
	\end{align}
From the results of Lemma\,\ref{L:PROPERTIESANDDIFFEOMORPHICEXTENSIONOFDOUBLENULLCOORDINATES}, including the estimates \eqref{E:CLOSEDVERSIONKEYJACOBIANDETERMINANTESTIMATECHOVGEOTODOUBLENULL}--\eqref{E:CLOSEDVERSIONLUBARESTIMATECHOVGEOTODOUBLENULL}, \eqref{E:CLOSEDVERSIONC11BOUNDFORCHOVDOUBLENULLTOGEO}, \eqref{E:SMALLC01ESTIMATESFORUBAR}, the bootstrap assumptions, Cor\,\ref{C:IMPROVEAUX}, and Rademacher's theorem, we have that $\left\| \nullangChristoff_{AB}^C \right\|_{L^{\infty}\left(\doublenulltoritwoarg{\ubar}{u}\right)} \lesssim \fundbootsmall$. By similar arguments, the contraction of the quadratic $\nullangChristoff^2$ terms on RHS\,\eqref{E:GAUSSCURVATUREESTIMATEINTERMEDIATESTEP1} against $(\gnulltori^{-1})^{AC}(\gnulltori^{-1})^{BD}$ satisfies the same $L^\infty$ estimate. Using definitions \eqref{E:RIEMANNCURVATURETENSOROFDOUBLENULLTORUS}--\eqref{E:2DGAUSSCURVATUREISTWICESCALARCURVATURE} of the double-null curvature, estimate \eqref{E:LINFINITYBOUNDFORDOUBLENULLGAUSSCURVATURE} will follow once we bound the top-order terms $\breve{\partial} \Upgamma$ in \eqref{E:GAUSSCURVATUREESTIMATEINTERMEDIATESTEP1}. 

Using the Leibniz rule and the definition of the Christoffel symbols, straightforward calculations imply: 
	\begin{align}
		\begin{split} \label{E:GAUSSCURVATUREESTIMATEINTERMEDIATESTEP3}
			\gnulltori_{ED} \left( \nullgeop{x^B} \nullangChristoff_{AC}^E - \nullgeop{x^A} \nullangChristoff_{BC}^E\right) & =  \nullgeop{x^A}\nullgeop{x^D} \gnulltori_{BC} + \nullgeop{x^B} \nullgeop{x^C} \gnulltori_{DA} - \nullgeop{x^B} \nullgeop{x^D} \gnulltori_{AC} - \nullgeop{x^A} \nullgeop{x^C} \gnulltori_{DB} \\
			& \ \ + \smoothfunction\left(\tander \controlvars, \frac{1}{\Lunit \ubar}, \frac{1}{\geop{t} \ubar}, \tander \ubar\right) \cdot \left(\tander \controlvars, \tander^{[1,2]} \ubar\right).
		\end{split}
	\end{align}
Since the $L^\infty$ norm of the fifth terms on RHS\,\eqref{E:GAUSSCURVATUREESTIMATEINTERMEDIATESTEP3} has already been estimated, we focus on the first four terms. Inserting the expression for $\gnulltori_{AB}$ from  \eqref{E:DOUBLENULLFIRSTFUNDCOMPONENTSINTERMSOFSMOOTHTORUSCOMPONENTS} into \eqref{E:GAUSSCURVATUREESTIMATEINTERMEDIATESTEP3}, we see that there are three different possibilities: \textbf{I)} both partial derivatives fall on the \emph{same} coefficient matrix $\gtorusdoublenullCOV^{-1}$; \textbf{II)} both angular partial derivatives $\breve{\partial}^2$ fall on $\gnulltori$; \textbf{III)} or the two angular partial derivatives get split up amongst the three factors $\gtorusdoublenullCOV^{-1} \gtorusdoublenullCOV^{-1} \gtorus$ on RHS\,\eqref{E:DOUBLENULLFIRSTFUNDCOMPONENTSINTERMSOFSMOOTHTORUSCOMPONENTS} by the Leibniz rule. Identity \eqref{E:INVERSECHOVCOEFFICIENTSSMOOTHANGULARDERIVATIVESINTERMSOFDOUBLENULLONESANDL} proves that cases \textbf{II)}--\textbf{III)} result in error terms that are precisely of the form $ \smoothfunction\left(\tander \controlvars, \frac{1}{\Lunit \ubar}, \frac{1}{\geop{t} \ubar}, \tander \ubar\right) \cdot \left(\tander \controlvars, \tander^{[1,2]} \ubar\right)$, whose $L^\infty$ norm we have already estimated. 

The remaining error terms resulting from case \textbf{I)} are:
	\begin{align}
		\begin{split} \label{E:GAUSSCURVATUREESTIMATEINTERMEDIATESTEP4}
			&  \left(\nullgeop{x^A}\nullgeop{x^D} \big( (\gtorusdoublenullCOV^{-1})_B^I\big) \right) (\gtorusdoublenullCOV^{-1})_C^J \gtorus_{IJ} +  \left(\nullgeop{x^A}\nullgeop{x^D} \big( (\gtorusdoublenullCOV^{-1})_C^J \big)\right)  (\gtorusdoublenullCOV^{-1})_B^I  \gtorus_{IJ} \\
			 + & \left( \nullgeop{x^B} \nullgeop{x^C} \big( (\gtorusdoublenullCOV^{-1})_D^I \big) \right) (\gtorusdoublenullCOV^{-1})_A^J \gtorus_{IJ}  + \left( \nullgeop{x^B} \nullgeop{x^C} \big( (\gtorusdoublenullCOV^{-1})_A^J \big) \right) (\gtorusdoublenullCOV^{-1})_D^I  \gtorus_{IJ} \\
			 - & \left( \nullgeop{x^B} \nullgeop{x^D} \big( (\gtorusdoublenullCOV^{-1})_A^I \big) \right) (\gtorusdoublenullCOV^{-1})_C^J \gtorus_{IJ} -  \left( \nullgeop{x^B} \nullgeop{x^D} \big(  (\gtorusdoublenullCOV^{-1})_C^J \big) \right) (\gtorusdoublenullCOV^{-1})_A^I \gtorus_{IJ} \\
			 - & \left(  \nullgeop{x^A} \nullgeop{x^C} \big( (\gtorusdoublenullCOV^{-1})_D^I \big) \right) (\gtorusdoublenullCOV^{-1})_B^J \gtorus_{IJ} - \left(  \nullgeop{x^A} \nullgeop{x^C} \big( (\gtorusdoublenullCOV^{-1})_B^J \big) \right) (\gtorusdoublenullCOV^{-1})_D^I   \gtorus_{IJ}. 
		\end{split}
	\end{align}
From \eqref{E:INVERSECHOVCOEFFICIENTSSMOOTHANGULARDERIVATIVESINTERMSOFDOUBLENULLONESANDL} and the Leibniz rule, we see that:
	\begin{align} \label{E:GAUSSCURVATUREESTIMATEINTERMEDIATESTEP5}
		\nullgeop{x^D} (\gtorusdoublenullCOV^{-1})_B^I = \nullgeop{x^D} \left( \frac{\geop{x^B} \ubar}{\geop{t} \ubar}\right) L^I + \frac{\geop{x^B} \ubar}{\geop{t} \ubar} \nullgeop{x^D} L^I. 
	\end{align}
From the same reasoning used in the prior two paragraphs, all of the error terms generated by the $ \frac{\geop{x^B} \ubar}{\geop{t} \ubar} \nullgeop{x^D} L^I$ in \eqref{E:GAUSSCURVATUREESTIMATEINTERMEDIATESTEP4} are of the form $ \smoothfunction\left(\tander \controlvars, \frac{1}{\Lunit \ubar}, \frac{1}{\geop{t} \ubar}, \tander \ubar\right) \cdot \left(\tander \controlvars, \tander^{[1,2]} \ubar\right)$. Expanding $\nullgeop{x^D}$ using the first identity in \eqref{E:NULLCOORDINATEPARTIALXAINTERMSOFGEOMETRICVF} when differentiating the first term on RHS\,\eqref{E:GAUSSCURVATUREESTIMATEINTERMEDIATESTEP5}, the remaining error terms are:
	\begin{align}
		\begin{split} \label{E:GAUSSCURVATUREESTIMATEINTERMEDIATESTEP6}
			& \nullgeop{x^A} \left(  \frac{ \twogeop{x^D}{x^B} \ubar }{\geop{t} \ubar} - \frac{\geop{x^B} \ubar \twogeop{x^D}{t} \ubar}{(\geop{t} \ubar)^2} - \frac{\geop{x^D} \ubar \twogeop{t}{x^B} \ubar}{(\geop{t} \ubar)^2} + \frac{\geop{x^D} \ubar \geop{x^B}\ubar \twogeop{t}{t} \ubar}{(\geop{t} \ubar)^3}\right) \Lunit^I(\gtorusdoublenullCOV^{-1})_C^J \gtorus_{IJ} \\
			+ & \nullgeop{x^A} \left(  \frac{ \twogeop{x^D}{x^C} \ubar }{\geop{t} \ubar} - \frac{\geop{x^C} \ubar \twogeop{x^D}{t} \ubar}{(\geop{t} \ubar)^2} - \frac{\geop{x^D} \ubar \twogeop{t}{x^C} \ubar}{(\geop{t} \ubar)^2} + \frac{\geop{x^D} \ubar \geop{x^C}\ubar \twogeop{t}{t} \ubar}{(\geop{t} \ubar)^3}\right) \Lunit^J(\gtorusdoublenullCOV^{-1})_B^I \gtorus_{IJ} \\
			+ &  \nullgeop{x^B}\left( \frac{\twogeop{x^C}{x^D}\ubar}{\geop{t}\ubar} - \frac{\geop{x^D}\ubar\twogeop{x^C}{t}\ubar}{(\geop{t}\ubar)^2} - \frac{\geop{x^C} \ubar \twogeop{t}{x^D}\ubar}{(\geop{t}\ubar)^2} + \frac{\geop{x^C}\ubar \geop{x^D}\ubar \twogeop{t}{t}\ubar}{(\geop{t}\ubar)^3}\right) \Lunit^I(\gtorusdoublenullCOV^{-1})_A^J \gtorus_{IJ}  \\
			+ &  \nullgeop{x^B}\left( \frac{\twogeop{x^C}{x^A}\ubar}{\geop{t}\ubar} - \frac{\geop{x^A}\ubar\twogeop{x^C}{t}\ubar}{(\geop{t}\ubar)^2} - \frac{\geop{x^C} \ubar\twogeop{t}{x^A}\ubar}{(\geop{t}\ubar)^2} + \frac{\geop{x^C}\ubar \geop{x^A}\ubar \twogeop{t}{t}\ubar}{(\geop{t}\ubar)^3}\right) \Lunit^J(\gtorusdoublenullCOV^{-1})_D^I \gtorus_{IJ} \\
			- & \nullgeop{x^B} \left( \frac{\twogeop{x^D}{x^A} \ubar}{\geop{t}\ubar} - \frac{\geop{x^A}\ubar \twogeop{x^D}{t}\ubar}{(\geop{t}\ubar)^2} - \frac{\geop{x^D}\ubar \twogeop{t}{x^A}\ubar}{(\geop{t}\ubar)^2} + \frac{\geop{x^D}\ubar \geop{x^A}\ubar \twogeop{t}{t}\ubar}{(\geop{t}\ubar)}\right) \Lunit^I(\gtorusdoublenullCOV^{-1})_C^J \gtorus_{IJ}  \\
			- & \nullgeop{x^B} \left( \frac{\twogeop{x^D}{x^C} \ubar}{\geop{t}\ubar} - \frac{\geop{x^C}\ubar \twogeop{x^D}{t}\ubar}{(\geop{t}\ubar)^2} - \frac{\geop{x^D}\ubar \twogeop{t}{x^C}\ubar}{(\geop{t}\ubar)^2} + \frac{\geop{x^D}\ubar \geop{x^C}\ubar \twogeop{t}{t}\ubar}{(\geop{t}\ubar)}\right) \Lunit^J(\gtorusdoublenullCOV^{-1})_A^I \gtorus_{IJ} \\
			- & \nullgeop{x^A} \left( \frac{\twogeop{x^C}{x^D}\ubar}{\geop{t}\ubar} - \frac{\geop{x^D}\ubar \twogeop{x^C}{t}\ubar}{(\geop{t}\ubar)^2} - \frac{\geop{x^C}\ubar \twogeop{t}{x^D}\ubar}{(\geop{t}\ubar)^2} + \frac{\geop{x^C}\ubar\geop{x^D}\ubar \twogeop{t}{t}\ubar}{(\geop{t}\ubar)^3} \right) \Lunit^I(\gtorusdoublenullCOV^{-1})_B^J \gtorus_{IJ} \\
			- & \nullgeop{x^A} \left( \frac{\twogeop{x^C}{x^B}\ubar}{\geop{t}\ubar} - \frac{\geop{x^B}\ubar \twogeop{x^C}{t}\ubar}{(\geop{t}\ubar)^2} - \frac{\geop{x^C}\ubar \twogeop{t}{x^B}\ubar}{(\geop{t}\ubar)^2} + \frac{\geop{x^C}\ubar\geop{x^B}\ubar \twogeop{t}{t}\ubar}{(\geop{t}\ubar)^3} \right) \Lunit^J(\gtorusdoublenullCOV^{-1})_D^I \gtorus_{IJ}
		\end{split}
	\end{align}
Applying the Leibniz rule one last time in \eqref{E:GAUSSCURVATUREESTIMATEINTERMEDIATESTEP6}, all of the error terms generated when the double-null coordinate derivatives $\nullgeop{x^A}$ fall on a \emph{first order} derivative $\geop{x^B}\ubar, \, \geop{t}\ubar$ are of the form $ \smoothfunction\left(\tander \controlvars, \frac{1}{\Lunit \ubar}, \frac{1}{\geop{t} \ubar}, \tander \ubar\right) \cdot \left(\tander \controlvars, \tander^{[1,2]} \ubar\right)$.

We now consider the final cases, which is when the double-null coordinate derivatives $\nullgeop{x^A}$ fall on a \emph{second order} derivatives $\twogeop{x^B}{x^C} \ubar, \, \twogeop{x^B}{t} \ubar, \twogeop{t}{t}\ubar$ in \eqref{E:GAUSSCURVATUREESTIMATEINTERMEDIATESTEP6}. Using the identity $\Lunit^J(\gtorusdoublenullCOV^{-1})_B^I = \Lunit^I(\gtorusdoublenullCOV^{-1})_B^J$, which follows easily from  \eqref{E:INVERSECHOVCOEFFICIENTSSMOOTHANGULARDERIVATIVESINTERMSOFDOUBLENULLONESANDL}, and $[\geop{x^A},\geop{x^B}] = [\geop{x^A},\geop{t}] = 0$, straightforward calculations show that the second and third lines of \eqref{E:GAUSSCURVATUREESTIMATEINTERMEDIATESTEP6} cancel the sixth and seventh lines of \eqref{E:GAUSSCURVATUREESTIMATEINTERMEDIATESTEP6}. 

This leaves only:
	\begin{align}
		\begin{split} \label{E:GAUSSCURVATUREESTIMATEINTERMEDIATESTEP7}
			& \nullgeop{x^A} \left(  \frac{ \twogeop{x^D}{x^B} \ubar }{\geop{t} \ubar} - \frac{\geop{x^B} \ubar \twogeop{x^D}{t} \ubar}{(\geop{t} \ubar)^2} - \frac{\geop{x^D} \ubar \twogeop{t}{x^B} \ubar}{(\geop{t} \ubar)^2} + \frac{\geop{x^D} \ubar \geop{x^B}\ubar \twogeop{t}{t} \ubar}{(\geop{t} \ubar)^3}\right) \Lunit^I(\gtorusdoublenullCOV^{-1})_C^J \gtorus_{IJ} \\
			+ &  \nullgeop{x^B}\left( \frac{\twogeop{x^C}{x^A}\ubar}{\geop{t}\ubar} - \frac{\geop{x^A}\ubar\twogeop{x^C}{t}\ubar}{(\geop{t}\ubar)^2} - \frac{\geop{x^C} \ubar\twogeop{t}{x^A}\ubar}{(\geop{t}\ubar)^2} + \frac{\geop{x^C}\ubar \geop{x^A}\ubar \twogeop{t}{t}\ubar}{(\geop{t}\ubar)^3}\right) \Lunit^J(\gtorusdoublenullCOV^{-1})_D^I \gtorus_{IJ} \\
			- & \nullgeop{x^B} \left( \frac{\twogeop{x^D}{x^A} \ubar}{\geop{t}\ubar} - \frac{\geop{x^A}\ubar \twogeop{x^D}{t}\ubar}{(\geop{t}\ubar)^2} - \frac{\geop{x^D}\ubar \twogeop{t}{x^A}\ubar}{(\geop{t}\ubar)^2} + \frac{\geop{x^D}\ubar \geop{x^A}\ubar \twogeop{t}{t}\ubar}{(\geop{t}\ubar)}\right) \Lunit^I(\gtorusdoublenullCOV^{-1})_C^J \gtorus_{IJ}  \\
			- & \nullgeop{x^A} \left( \frac{\twogeop{x^C}{x^B}\ubar}{\geop{t}\ubar} - \frac{\geop{x^B}\ubar \twogeop{x^C}{t}\ubar}{(\geop{t}\ubar)^2} - \frac{\geop{x^C}\ubar \twogeop{t}{x^B}\ubar}{(\geop{t}\ubar)^2} + \frac{\geop{x^C}\ubar\geop{x^B}\ubar \twogeop{t}{t}\ubar}{(\geop{t}\ubar)^3} \right) \Lunit^J(\gtorusdoublenullCOV^{-1})_D^I \gtorus_{IJ}
		\end{split}
	\end{align}			
Expanding $\nullgeop{x^A}$ using the first identity in \eqref{E:NULLCOORDINATEPARTIALXAINTERMSOFGEOMETRICVF}, 	straightforward but tedious calculations show that every term in \eqref{E:GAUSSCURVATUREESTIMATEINTERMEDIATESTEP7} cancels.

We have thus proven that the $L^{\infty}\left(\doublenulltoritwoarg{\ubar}{u}\right)$ norm of the contracted Riemann curvature $\Riemnulltori\left(\nullgeop{x^A},\nullgeop{x^B},\nullgeop{x^C},\nullgeop{x^D}\right)$ is $\lesssim \fundbootsmall$. Contracting  against $(\gnulltori^{-1})^{AC}(\gnulltori^{-1})^{BD}$, we conclude \eqref{E:LINFINITYBOUNDFORDOUBLENULLGAUSSCURVATURE}.

\end{proof}

\subsubsection{Proof of Lemma\,\ref{L:STANDARDELLIPTICONDOUBLENULLTORI}}
\label{SSS:PROOFOFL:STANDARDELLIPTICONDOUBLENULLTORI}
In this proof  we will use capital Latin indices to denote
the components of $\doublenulltoritwoarg{\ubar}{u}$-tangent tensorfields
with respect to the frame $\left\lbrace \nullgeop{x^A} \right\rbrace_{A=2,3}$
and co-frame $\lbrace \nullangrmD x^A \rbrace_{A=2,3}$ on $\doublenulltoritwoarg{\ubar}{u}$,
and we raise and lower indices with $\gnulltori^{-1}$ and $\gnulltori$.
For example, 
$\Upxi 
= 
\Upxi\left(\nullgeop{x^A},\nullgeop{x^B} \right)
\nullangrmD x^A \otimes \nullangrmD x^B
\eqdef \Upxi_{AB}
\nullangrmD x^A \otimes \nullangrmD x^B$ and $\Upxi^A_B = (\gnulltori^{-1})^{AC} \Upxi_{CB}$.
We start by defining $\breve{I} = \breve{I}^A \nullgeop{x^A}$ to be the
$\doublenulltoritwoarg{\ubar}{u}$-tangent vectorfield
with the following components relative to the coordinates $(x^2,x^3)$ on
$\doublenulltoritwoarg{\ubar}{u}$:
\begin{align} \label{E:CURRENTFORELLIPTICIDENTITYONDOUBLENULLTORIMANIFOLD}
	\breve{I}^A
	& 
	\eqdef  
	\upmu^2 \Upxi_{ B  C} \nullangD^{B} \Upxi^{AC} 
	- 
	\upmu^2 \Upxi^{AB} (\nullangDiv \Upxi)_B.
\end{align}
Next, with $\Gaussnulltori$ the Gauss curvature of
$\doublenulltoritwoarg{\ubar}{u}$,
we note the following standard identity,
which follows from the symmetry of $\Upxi$
(see \cite{jS2016b}*{Lemma~18.9} for the main ideas of the proof,
where we note that only trace-free tensorfields were handled in
\cite{jS2016b}*{Lemma~18.9} and thus RHSs\,\eqref{E:STANDARDELLIPTICONDOUBLENULLTORI} and \eqref{E:ELLIPTICIDENTITYONDOUBLENULLTORIMANIFOLD} feature
additional $\mytr_{\gnulltori} \Upxi$-dependent terms compared to \cite{jS2016b}*{Lemma~18.9}):
\begin{align}
\begin{split} \label{E:ELLIPTICIDENTITYONDOUBLENULLTORIMANIFOLD}
	\upmu^2 |\nullangD \Upxi|_{\gnulltori}^2 
	+ 
	2 \upmu^2 \Gaussnulltori |\Upxi|_{\gnulltori}^2
	& 
	= 
	2 \upmu^2 |\nullangDiv \Upxi|_{\gnulltori}^2 
	+ 
	\upmu^2 \Gaussnulltori (\mytr_{\gnulltori} \Upxi)^2 
	+
	\upmu^2 |\nullangrmD \mytr_{\gnulltori} \Upxi|_{\gnulltori}^2  
		\\
& \ \ 
		- 
		2 \upmu^2 
		\gnulltori^{-1}(\nullangDiv \Upxi,\nullangrmD 
		\mytr_{\gnulltori} \Upxi) 
		+ 
		2 
		\upmu  
		\Upxi^{AB} 
		\left(\nullgeop{x^A} \upmu \right)
		(\nullangDiv \Upxi)_{\widetilde B} 
	\\
& \ \
- 
2 \upmu 
\Upxi_{B C} 
\left(\nullgeop{x^A} \upmu\right) 
\nullangD^B \Upxi^{AC} 
+ 
\nullangDiv \breve{I}.
\end{split}
\end{align}
We then integrate \eqref{E:ELLIPTICIDENTITYONDOUBLENULLTORIMANIFOLD} over $\doublenulltoritwoarg{\ubar}{u}$
with respect to the area form $\voldoublenulltori$ defined in \eqref{E:AREAFORMDOUBLENULLTORUS}
and note that the integral of the perfect divergence term $\nullangDiv \breve{I}$ vanishes.
Next, we use the $\gnulltori$-Cauchy--Schwarz inequality and
Young's inequality to pointwise bound the three cross-terms on RHS\,\eqref{E:ELLIPTICIDENTITYONDOUBLENULLTORIMANIFOLD}
as follows:
\begin{align}
2 \left| 
	\upmu^2 \gnulltori^{-1}(\nullangDiv \Upxi, \nullangrmD \mytr_{\gnulltori} \Upxi) 
\right| 
& \leq 
	\upmu^2 
	\left|\nullangDiv \Upxi \right|_{\gnulltori}^2 
+ 
\upmu^2 
\left|\nullangrmD \mytr_{\gnulltori} \Upxi \right|_{\gnulltori}^2, 
		\label{E:ROUGHTORIELLIPTICELEMENTARYINEQUALITY1} 
			\\
2 \left| 
	\upmu \Upxi^{AB} \left(\nullgeop{x^A} \upmu \right)(\nullangDiv \Upxi)_{\widetilde B} 
\right| 
& \leq 
\upmu^2 
\left|\nullangDiv \Upxi \right|_{\gnulltori}^2 
+ 
\left|\nullangrmD \upmu \right|_{\gnulltori}^2 
\left|\Upxi \right|_{\gnulltori}^2, 
	\label{E:ROUGHTORIELLIPTICELEMENTARYINEQUALITY2}  
		\\
2 \left| \upmu \Upxi_{BC}  
\left(\nullgeop{x^A} \upmu \right) \nullangD^B \Upxi^{AC}\right| 
& 
\leq 
\frac{1}{3} 
\upmu^2 
|\nullangD \Upxi|_{\gnulltori}^2 
+ 
3 
\left|\nullangrmD \upmu \right|_{\gnulltori}^2 
\left|\Upxi \right|_{\gnulltori}^2.
 \label{E:ROUGHTORIELLIPTICELEMENTARYINEQUALITY3} 
\end{align}
Next,
we use \eqref{E:ANGDFPOINTWISEBOUNDEDBYCOMMUTATORVECTORFIELDS}, \eqref{E:POINTWISEDOULBENULLTORIDIFFERENTIALINORMINTERMSOFSMOOTHTORIDIFFERENTIALNORMANDL},
the bootstrap assumptions,
and Cor\,\ref{C:IMPROVEAUX}
to deduce:
\begin{align} \label{E:POINTWISEBOUNDDOUBLENULLGRADIENTOFMU}
\left|\nullangrmD \upmu \right|_{\gnulltori}
&
\lesssim 
|\rmd \upmu|_{\gtorus}
+
\fundbootsmall
|\Lunit \upmu|
\lesssim
|\tanderY \upmu|
+
\fundbootsmall
|\Lunit \upmu|
\lesssim
\fundbootsmall.
\end{align}
Using \eqref{E:POINTWISEBOUNDDOUBLENULLGRADIENTOFMU} to control the relevant factors on
RHSs\,\eqref{E:ROUGHTORIELLIPTICELEMENTARYINEQUALITY2}--\eqref{E:ROUGHTORIELLIPTICELEMENTARYINEQUALITY3},
using 
\eqref{E:LINFINITYBOUNDFORDOUBLENULLGAUSSCURVATURE}
to bound the factors of
$\Gaussnulltori$ in \eqref{E:ELLIPTICIDENTITYONDOUBLENULLTORIMANIFOLD},
using the inequality 
$|\mytr_{\gnulltori} \Upxi| 
\lesssim 
|\Upxi|_{\gnulltori}$,
and using $|\upmu| \lesssim 1$ (which follows from the bootstrap assumptions),
we conclude \eqref{E:STANDARDELLIPTICONDOUBLENULLTORI}.

\hfill $\qed$

\subsection{Proof of Prop.\,\ref{P:ELLIPTICESTIMATESONDOUBLENULLTORI}}
\label{SS:PROOFOFP:ELLIPTICESTIMATESONDOUBLENULLTORI}
We now prove Prop.\,\ref{P:ELLIPTICESTIMATESONDOUBLENULLTORI}. Let $\Upxi$ be a symmetric type $\binom{0}{2}$ $\doublenulltoritwoarg{\ubar}{u}$-tangent tensorfield. Then since $\nullgeop{x^A}  \weight = \nullgeop{x^A}(-\ubar) = 0$, the following weighted version of \eqref{E:STANDARDELLIPTICONDOUBLENULLTORI} holds:
\begin{align} \label{E:ELLIPTICESTIMATESONDOUBLENULLTORISTEP1}
\int_{\doublenulltoritwoarg{\ubar}{u}} 
	\upmu^2 |\nullangD \Upxi|_{\gnulltori}^2 \weight^{\blowupratetoporderwave}
\, \voldoublenulltori
& \leq 
6 
\int_{\doublenulltoritwoarg{\ubar}{u}}
	\upmu^2 |\nullangDiv \Upxi|_{\gnulltori}^2 \weight^{\blowupratetoporderwave}
\, \voldoublenulltori
+ 
3 
\int_{\doublenulltoritwoarg{\ubar}{u}} 
	\upmu^2 |\nullangrmD \mytr_{\gnulltori} \Upxi|_{\gnulltori}^2 \weight^{\blowupratetoporderwave}
\, \voldoublenulltori
+ 
C 
\fundbootsmall 
\int_{\doublenulltoritwoarg{\ubar}{u}}
	|\Upxi|_{\gnulltori}^2 \weight^{\blowupratetoporderwave}
\, \voldoublenulltori.
\end{align}Let $\upxi$ be a symmetric type $\binom{0}{2}$ $\ell_{t,u}$-tangent tensorfield, and let $\breve{\upxi} \eqdef \nulltorusproject \upxi$ denote its $\gfour$-orthogonal projection onto $\doublenulltoritwoarg{\ubar}{u}$ defined in \eqref{E:TENSORFIELDPROJECTIONTONULLTORIDEFININGEQUATION}. Since the term
$
 \int_{\characteristicdiamondtwoarg{[\leftubar,\ubar]}{[\moreinterestingu_1,u]}} 	
	\upmu^2 |\angLie_{\Lunit} \upxi|_{\gtorus}^2 \weight^{\blowupratetoporderwave}
\, \voldiamond 
$
on LHS\,\eqref{E:ELLIPTICESTIMATECHI} 
is manifestly bounded by RHS\,\eqref{E:ELLIPTICESTIMATECHI}, 
we only have to show that for $A=2,3$,
the term
$
 \int_{\characteristicdiamondtwoarg{[\leftubar,\ubar]}{[\moreinterestingu_1,u]}} 	
	\upmu^2 |\angLie_{\Yvf{A}} \upxi|_{\gtorus}^2 \weight^{\blowupratetoporderwave}
\, \voldiamond 
$
on LHS\,\eqref{E:ELLIPTICESTIMATECHI} is $\leq \mbox{RHS\,\eqref{E:ELLIPTICESTIMATECHI}}$. 
To proceed, we consider the inequality \eqref{E:ELLIPTICESTIMATESONDOUBLENULLTORISTEP1} 
with $\breve{\upxi}$ in the role of $\Upxi$.
Integrating the inequality with respect to $u'$ and $\ubar'$, as well as using  Lemmas\,\ref{L:SMOOTHANDDOUBLENULLFIRSTFUNDFORMSAGREEONPUTANGENTENSORS}
and\,\ref{L:DIFFOPERATORSPOINTWISECOMPARISONNEEDEDFORELLIPTICESTIMATESONDOUBLENULLTORI},
we deduce, in view of definition \eqref{E:VOLUMEFORMINGOINGNULLHYPERSURFACE}, that:
\begin{align}
\begin{split} \label{E:ELLIPTICESTIMATESONDOUBLENULLTORISTEP2}
 \int_{\characteristicdiamondtwoarg{[\leftubar,\ubar]}{[\moreinterestingu_1,u]}} 	
	\upmu^2 |\angD \upxi|_{\gtorus}^2 \weight^{\blowupratetoporderwave}
\, \voldiamond 
& \lesssim
 \int_{\characteristicdiamondtwoarg{[\leftubar,\ubar]}{[\moreinterestingu_1,u]}} 	
	\upmu^2 |\angdiv \upxi|_{\gtorus}^2 \weight^{\blowupratetoporderwave}
\, \voldiamond
+ 
 \int_{\characteristicdiamondtwoarg{[\leftubar,\ubar]}{[\moreinterestingu_1,u]}} 	
	\upmu^2 |\angD \mytr_{\gtorus} \upxi|_{\gtorus}^2\weight^{\blowupratetoporderwave}
\, \voldiamond
	\\
& \ \
+ 
\fundbootsmall
 \int_{\characteristicdiamondtwoarg{[\leftubar,\ubar]}{[\moreinterestingu_1,u]}} 	
	\upmu^2 |\Lunit \mytr_{\gtorus} \upxi|^2\weight^{\blowupratetoporderwave}
\, \voldiamond
+ 
\fundbootsmall
 \int_{\characteristicdiamondtwoarg{[\leftubar,\ubar]}{[\moreinterestingu_1,u]}} 	
	\upmu^2 |\angLie_{\Lunit} \upxi|_{\gtorus}^2\weight^{\blowupratetoporderwave}
\, \voldiamond \\
& \ \ 
+ 
\fundbootsmall 
  \int_{\characteristicdiamondtwoarg{[\leftubar,\ubar]}{[\moreinterestingu_1,u]}} 	
	|\upxi|_{\gtorus}^2\weight^{\blowupratetoporderwave}
\, \voldiamond.
\end{split}
\end{align}
Next, using the Leibniz rule, 
\eqref{E:TANDERGANDCHIESTIMATE},
and the bootstrap assumptions, we find that
$|\Lunit \mytr_{\gtorus} \upxi|
\lesssim 
|\angLie_{\Lunit} \upxi|_{\gtorus}
+
|\angLie_{\Lunit} \gtorus^{-1}|_{\gtorus}
|\upxi|_{\gtorus}
\lesssim 
|\angLie_{\Lunit} \upxi|_{\gtorus}
+
|\upxi|_{\gtorus}
$
and thus the third integral
$
\fundbootsmall
 \int_{\characteristicdiamondtwoarg{[\leftubar,\ubar]}{[\moreinterestingu_1,u]}} 	
	\upmu^2 |\Lunit \mytr_{\gtorus} \upxi|^2 \weight^{\blowupratetoporderwave}
\, \voldiamond
$
on RHS\,\eqref{E:ELLIPTICESTIMATESONDOUBLENULLTORISTEP2} 
is
bounded by the last two integrals on RHS\,\eqref{E:ELLIPTICESTIMATESONDOUBLENULLTORISTEP2}.
Next, we use the torsion-free property of the connection $\angD$ 
and the $\gtorus$-Cauchy--Schwarz inequality to deduce the pointwise estimate
$
|\angLie_{\Yvf{A}} \upxi|_{\gtorus}
\leq 
|\angDarg{\Yvf{A}} \upxi|_{\gtorus}
+
2 |\upxi|_{\gtorus}|\angD \Yvf{A}|_{\gtorus}
$.
Also using 
\eqref{E:POINTWISESEMINORMOFYVECTORFIELDS},
\eqref{E:SMOOTHTORUSNORMCOMPARBLETOTANGENTIALCONTRACTIONS},
the schematic identity and estimate $|\gfour(\Dfour_{\Yvf{A}} \Yvf{B},\Yvf{C})| = | \smoothfunction(\controlvars) \tanderY \controlvars | \lesssim \fundbootsmall$, we find that
$
|\angLie_{\Yvf{A}} \upxi|_{\gtorus}
\lesssim
|\angD \upxi|_{\gtorus}
+
\fundbootsmall
|\upxi|_{\gtorus}
$.
From this bound,
\eqref{E:ELLIPTICESTIMATESONDOUBLENULLTORISTEP2},
the estimates proved above,
and the pointwise bound 
$|\angD \mytr_{\gtorus} \upxi|_{\gtorus} 
\lesssim
\sum_{A=2,3} 
|\Yvf{A} \mytr_{\gtorus} \upxi|_{\gtorus}$,
which follows from \eqref{E:SMOOTHTORUSNORMCOMPARBLETOTANGENTIALCONTRACTIONS},
we conclude that
$
 \int_{\characteristicdiamondtwoarg{[\leftubar,\ubar]}{[\moreinterestingu_1,u]}} 	
	\upmu^2 |\angLie_{\Yvf{A}} \upxi|_{\gtorus}^2 \weight^{\blowupratetoporderwave}
\, \voldiamond 
\lesssim
\mbox{RHS\,\eqref{E:ELLIPTICESTIMATECHI}}
$.
We have therefore proved \eqref{E:ELLIPTICESTIMATECHI},
which completes the proof of Prop.\,\ref{P:ELLIPTICESTIMATESONDOUBLENULLTORI}.

\hfill $\qed$

 \section{The divergence theorem and general integration by parts identities} \label{S:DIVERGENCETHEOREMANDINTEGRATIONBYPARTSIDENTITIES}
 
We will encounter numerous error integrals in our $L^2$ analysis that, at first glance, appear to suffer from a loss of derivatives. In Sect.\,\ref{S:WAVEANDACOUSTICGEOMETRYAPRIORIESTIMATES} we will see that all of these potentially harmful integrals exhibit a crucial \emph{divergence} property that allows one to integrate away the derivative-losing terms. To that end, we develop a comprehensive study of the divergence theorem in this section. This is largely standard -- the nontrivial aspect has to do with proving that all boundary integrals on $\ingoingcharacteristicsurfacetwoarg{\ubar}{[\moreinterestingu_1,u]}$ generated from the divergence theorem \emph{only feature $\ingoingcharacteristicsurfacetwoarg{\ubar}{[\moreinterestingu_1,u]}$-tangential derivatives}.

\subsection{Divergence theorem identity} \label{SS:DIVERGENCETHEOREM}
In this section we prove the main divergence identity for spacetime vectorfields on the characteristic diamond $\characteristicdiamondtwoarg{[\leftubar,\ubar)}{[\moreinterestingu_1,u]}$.
\begin{lemma}[Divergence identity on the characteristic diamond] \label{L:ARBITRARYDIVERGENCETHEOREM}
	Let $\mathcal{J}$ be a spacetime vectorfield on  $\characteristicdiamondtwoarg{[\leftubar,\ubar)}{[\moreinterestingu_1,u]}$, and let $\mathcal{J} = \mathcal{J}^{\ubar} \nullgeop{\ubar} + \mathcal{J}^u \nullgeop{u} + \mathcal{J}^A \nullgeop{x^A}$ be its decomposition relative to the double-null coordinate vectorfields. Then the following identity holds:
	\begin{equation} \label{E:DIVERGENCEIDENTITYSPACETIMEVECTORFIELD}
		\int_{\ingoingcharacteristicsurfacetwoarg{\ubar}{[\moreinterestingu_1,u]}} \mathcal{J}^{\ubar} \MagnitueofinnerproductofnewLandnewuL \,  \volingoingnullhypersurface + \int_{\outgoingcharacteristicsurfacetwoarg{u}{[\leftubar,\ubar)}}  \mathcal{J}^{u} \MagnitueofinnerproductofnewLandnewuL \, \voloutgoingnullhypersurface  =	\int_{\ingoingcharacteristicsurfacetwoarg{\leftubar}{[\moreinterestingu_1,u]}} \mathcal{J}^{\ubar} \MagnitueofinnerproductofnewLandnewuL \, \volingoingnullhypersurface + \int_{\outgoingcharacteristicsurfacetwoarg{\moreinterestingu_1}{[\leftubar,\ubar)}}  \mathcal{J}^{u} \MagnitueofinnerproductofnewLandnewuL \, \voloutgoingnullhypersurface + \int_{\characteristicdiamondtwoarg{[\leftubar,\ubar)}{[\moreinterestingu_1,u]}} \Dfour_\alpha \mathcal{J}^\alpha \MagnitueofinnerproductofnewLandnewuL \, \voldiamond.
	\end{equation}
\end{lemma}

\begin{proof}
	By combining the standard identity for the divergence of a vectorfield expressed relative to the double-null coordinates 
	with the identity
	\eqref{E:DETACOUSTICALMETRICINDOUBLENULLCOORDS}, we have:
	\begin{align}
		\begin{split}  \label{E:IDFORDETWEIGHTEDCOVARIANTDIVERGENCEOFSPACETIMEVECTORFIELD}
			\sqrt{|\mbox{\upshape det} \gfour|} \, \Dfour_\alpha \mathcal{J}^\alpha
			& = 
			\nullgeop{\ubar} 
			\left( \sqrt{|\mbox{\upshape det} \gfour|} \mathcal{J}^{\ubar} \right) 
			+ 
			\nullgeop{u} \left( \sqrt{|\mbox{\upshape det} \gfour|} \mathcal{J}^u\right) 
			+ 
			\nullgeop{x^A} \left( \sqrt{|\mbox{\upshape det} \gfour|} \mathcal{J}^A \right) 
			\\
			& = 
			\nullgeop{\ubar} 
			\left( \mathcal{J}^{\ubar} \MagnitueofinnerproductofnewLandnewuL
			\sqrt{ \mbox{\upshape det} \gnulltori} 
			\right) 
			+ 
			\nullgeop{u} 
			\left( \mathcal{J}^u \MagnitueofinnerproductofnewLandnewuL \sqrt{ \mbox{\upshape det} \gnulltori} 
			\right) 
			+ 
			\nullgeop{x^A} 
			\left( 
			\sqrt{|\mbox{\upshape det} \gfour|} \mathcal{J}^A 
			\right).
		\end{split}
	\end{align}
	Integrating \eqref{E:IDFORDETWEIGHTEDCOVARIANTDIVERGENCEOFSPACETIMEVECTORFIELD} over 
	$\characteristicdiamondtwoarg{[\leftubar,\ubar)}{[\moreinterestingu_1,u]}$
	with respect to $\mathrm{d} x^2 \, \mathrm{d}x^3 \mathrm{d} u' \mathrm{d} \ubar'$, 
	using \eqref{E:VOLFORMACOUSTICALMETRICDOUBLENULLCOORDS} to relate the canonical volume form $\volcanonical_{\gfour}$ to $\voldiamond$, applying Fubini's theorem,
	and using that the integral of the last term on RHS\,\eqref{E:IDFORDETWEIGHTEDCOVARIANTDIVERGENCEOFSPACETIMEVECTORFIELD} 
	over $\T^2$ vanishes (since $\T^2$ is a closed manifold), we deduce \eqref{E:DIVERGENCEIDENTITYSPACETIMEVECTORFIELD}.
\end{proof}
We complement Lemma\,\ref{L:ARBITRARYDIVERGENCETHEOREM} with the following lemma, which provides a useful expression for the covariant divergence of a spacetime vectorfield.

\begin{lemma}[Covariant divergence identity for spacetime vectorfields] 
	\label{L:COVARIANTDIVERGENCEOFSPACETIMEVECTORFIELDINTERMSOFRESCALEDFRAME}
	Let $\mathcal{J}$ be a spacetime vectorfield.
	Consider the decomposition $\upmu \mathcal{J} 
	= 
	- \upmu \mathcal{J}_{\Lunit} \Lunit 
	- 
	\mathcal{J}_{\muX} \Lunit 
	- 
	\mathcal{J}_{\Lunit} \muX 
	+ 
	\upmu \smoothtorusproject \mathcal{J}$ of $\upmu \mathcal{J}$ 
	afforded by Lemma\,\ref{L:BASICPROPERTIESOFVECTORFIELDS} ,
	where $\mathcal{J}_{\Lunit} = \gfour(\mathcal{J},\Lunit)$, $\mathcal{J}_{\muX} = \gfour(\mathcal{J},\muX)$, 
	and $\smoothtorusproject \mathcal{J}$ is the $\ell_{t,u}$-projection of $\mathcal{J}$ (see Def.\,\ref{D:PROJECTIONTENSORFIELDSANDTANGENCYTOHYPERSURFACES}).
	Then the following identity holds:
	\begin{align} \label{E:COVARIANTDIVERGENCEOFSPACETIMEVECTORFIELDINTERMSOFRESCALEDFRAME}
		\upmu \Dfour_\alpha \mathcal{J}^{\alpha}= - L(\upmu \mathcal{J}_{\Lunit}) - L(\mathcal{J}_{\muX}) - \muX (\mathcal{J}_{\Lunit}) + \angdiv(\upmu \smoothtorusproject \mathcal{J}) - \upmu \mytr_{\gtorus} \angk \mathcal{J}_{\Lunit} - \mytr_{\gtorus} \upchi \mathcal{J}_{\muX}.
	\end{align}
	
	Consequently, the following identities hold:
	\begin{subequations}
		\begin{align}
			\Dfour_\alpha \Lunit^\alpha & = \frac{1}{\upmu} \Lunit \upmu + \mytr_{\gtorus} \upchi, \label{E:SPACETIMEDIVERGENCEOFL}, \\
			\Dfour_\alpha \Yvf{A}^\alpha & = \frac{1}{\upmu} \Yvf{A} \upmu + \frac{1}{2} \mytr_{\gtorus} \, \angdeform{\Yvf{A}} \label{E:SPACETIMEDIVERGENCEOFYA}, \\
			\Dfour_\alpha \muX^\alpha & =  \frac{1}{\upmu} \muX \upmu + \upmu \mytr_{\gtorus} \angk - \upmu \mytr_{\gtorus}\upchi, \label{E:SPACETIMEDIVERGENCEOFMUX} \\
	 \Dfour_\alpha \Transport^\alpha & = \frac{1}{\upmu} \Lunit \upmu +  \mytr_{\gtorus} \angk \label{E:SPACETIMEDIFERGENCEOFB}
		\end{align}
	\end{subequations}
	
Moreover, for the multiplier vectorfield $\multipliervectorfield$ from \eqref{E:MULTIPLIERVECTORFIELD}, we have: 
	\begin{align} \label{E:SPACETIMEDIVERGENCEOFMULTIPLIER}
		\Dfour_\alpha \multipliervectorfield^\alpha = \frac{1}{\upmu} \multipliervectorfield \upmu + 2  \Lunit \upmu  + 2 \upmu \mytr_{\gtorus}\angk + \mytr_{\gtorus}\upchi.
	\end{align}In addition, the following identities hold for scalar functions $\upupsilon$ and $\zeta$:
\begin{align}\label{E:SPACETIMEDIVERGENCEOFANGULARGRADIENT}
		\upmu \Dfour_\alpha \left( \frac{1}{\upmu} \upupsilon (\gtorus^{-1})^{\alpha\beta} \angD_\beta \zeta\right) = \gtorus( \angD \upupsilon, \angD \zeta) + \upupsilon \angLap \zeta.
	\end{align}Finally, we have the following identity:
	\begin{align} 
		\begin{split} \label{E:SPACETIMEDIVERGENCEOFMUWEIGHTEDCARTESIANSPATIALDERIVATIVES}
			\Dfour_\alpha (\upmu \p_i)^\alpha & \eqdef \frac{1}{\upmu} \muX \left( \frac{X^i}{\Speed^2} \upmu\right) + \frac{1}{\upmu} \Yvf{A} \left( \upmu \angularcoeffinmuweightedspatialcartesian_i^A\right) + \frac{1}{2} \angularcoeffinmuweightedspatialcartesian_i^A \, \mytr_{\gtorus} \angdeform{\Yvf{A}} \\
			& \ \ + \frac{\upmu X^i}{\Speed^2}  \mytr_{\gtorus} \angk - \frac{\upmu X^i}{\Speed^2} \mytr_{\gtorus}\upchi,
		\end{split}
	\end{align}
where $\angularcoeffinmuweightedspatialcartesian_i^A$ are the $3\times 2$ matrices from \eqref{E:ANGULARCOEFFINMUWEIGHTEDCARTESIAN}	

\end{lemma}

\begin{proof}
	The same proof of \cite[Lemma 4.3]{jSgHjLwW2016} holds for \eqref{E:COVARIANTDIVERGENCEOFSPACETIMEVECTORFIELDINTERMSOFRESCALEDFRAME} with minor modifications to account for the third space dimension.
	
	Identities \eqref{E:SPACETIMEDIVERGENCEOFL}--\eqref{E:SPACETIMEDIFERGENCEOFB} follow from straightforward applications of Lemma\,\ref{L:BASICPROPERTIESOFVECTORFIELDS} and \eqref{E:COVARIANTDIVERGENCEOFSPACETIMEVECTORFIELDINTERMSOFRESCALEDFRAME}. 
	
	Identity \eqref{E:SPACETIMEDIVERGENCEOFMULTIPLIER} follows from straightforward applications of the Leibniz rule and \eqref{E:SPACETIMEDIVERGENCEOFL},\eqref{E:SPACETIMEDIVERGENCEOFMUX}. 
	
	Identity \eqref{E:SPACETIMEDIVERGENCEOFANGULARGRADIENT} follows from \eqref{E:COVARIANTDIVERGENCEOFSPACETIMEVECTORFIELDINTERMSOFRESCALEDFRAME}, the angular Leibniz rule, and the fact that $ (\gtorus^{-1})^{\alpha\beta} \angrmD_\beta \zeta$ is $\gfour$-orthogonal to $\Lunit$ and $\muX$. 
	
	 Identity \eqref{E:SPACETIMEDIVERGENCEOFMUWEIGHTEDCARTESIANSPATIALDERIVATIVES}  follows from straightforward applications of Lemma\,\ref{L:BASICPROPERTIESOFVECTORFIELDS}, \eqref{E:EXPRESSIONFORMUWEIGHTEDSPATIALCARTESIANDERIVATIVESWRTCOMMUTATORS}, and \eqref{E:COVARIANTDIVERGENCEOFSPACETIMEVECTORFIELDINTERMSOFRESCALEDFRAME}. 
	
\end{proof}

\subsection{Integrating by parts a single differential operator} \label{SS:IBPASINGLEDIFFERENTIALOPERATOR}

The following lemma provides some preliminary integration by parts identities involving a single differential operator. 

\begin{lemma}[Preliminary integration by parts identities involving a single differential operator]  \label{L:IBPASINGLEDIFFERENTIALOPERATOR}
Let $(\ubar,u) \in [\leftubar,\ubarboot)\times [\moreinterestingu_1,\moreinterestingu_2]$, and let $\upupsilon$ and $\zeta$ be scalar functions on $\characteristicdiamondtwoarg{[\leftubar,\ubar)}{[\moreinterestingu_1,u]}$. Then the following integration by parts identities hold, where $A = 2,3$ in \eqref{E:PRELIMINARYINTEGRATIONBYPARTSWITHY}:
	\begin{subequations} 
		\begin{align}
			\begin{split} \label{E:PRELIMINARYINTEGRATIONBYPARTSWITHL}
				\int_{\characteristicdiamondtwoarg{[\leftubar,\ubar)}{[\moreinterestingu_1,u]}} (\Lunit \upupsilon) \zeta \ReciprocalLunitAppliedtoTimeFunction \, \voldiamond & = - \int_{\characteristicdiamondtwoarg{[\leftubar,\ubar)}{[\moreinterestingu_1,u]}} \upupsilon (\Lunit \zeta) \ReciprocalLunitAppliedtoTimeFunction \, \voldiamond \\
				& \ \
				- \int_{\characteristicdiamondtwoarg{[\leftubar,\ubar)}{[\moreinterestingu_1,u]}} \mytr_{\gtorus} \upchi \upupsilon \zeta \ReciprocalLunitAppliedtoTimeFunction \, \voldiamond\\
				& \ \ 
				+ \int_{\ingoingcharacteristicsurfacetwoarg{\ubar}{[\moreinterestingu_1,u]}} \upupsilon \zeta \, \volingoingnullhypersurface - \int_{\ingoingcharacteristicsurfacetwoarg{\leftubar}{[\moreinterestingu_1,u]}} \upupsilon \zeta \, \volingoingnullhypersurface,
			\end{split} \\
			\begin{split} \label{E:PRELIMINARYINTEGRATIONBYPARTSWITHY}
				\int_{\characteristicdiamondtwoarg{[\leftubar,\ubar)}{[\moreinterestingu_1,u]}} (\Yvf{A} \upupsilon) \zeta \ReciprocalLunitAppliedtoTimeFunction \, \voldiamond & = - \int_{\characteristicdiamondtwoarg{[\leftubar,\ubar)}{[\moreinterestingu_1,u]}} \upupsilon (\Yvf{A} \zeta) \ReciprocalLunitAppliedtoTimeFunction \, \voldiamond \\
				& \ \
				- \frac{1}{2} \int_{\characteristicdiamondtwoarg{[\leftubar,\ubar)}{[\moreinterestingu_1,u]}}\mytr_{\gtorus} \, \angdeform{\Yvf{A}} \upupsilon \zeta \ReciprocalLunitAppliedtoTimeFunction \, \voldiamond\\
				& \ \ 
				+ \int_{\ingoingcharacteristicsurfacetwoarg{\ubar}{[\moreinterestingu_1,u]}} \upupsilon \zeta (\Yvf{A} \ubar ) \ReciprocalLunitAppliedtoTimeFunction \, \volingoingnullhypersurface - \int_{\ingoingcharacteristicsurfacetwoarg{\leftubar}{[\moreinterestingu_1,u]}} \upupsilon \zeta ( \Yvf{A} \ubar )\ReciprocalLunitAppliedtoTimeFunction \, \volingoingnullhypersurface,
			\end{split} \\
			\begin{split} \label{E:PRELIMINARYINTEGRATIONBYPARTSWITHMUX}
				\int_{\characteristicdiamondtwoarg{[\leftubar,\ubar)}{[\moreinterestingu_1,u]}} (\muX \upupsilon) \zeta \ReciprocalLunitAppliedtoTimeFunction \, \voldiamond & = - \int_{\characteristicdiamondtwoarg{[\leftubar,\ubar)}{[\moreinterestingu_1,u]}} \upupsilon (\muX \zeta) \ReciprocalLunitAppliedtoTimeFunction \, \voldiamond \\
				& \ \
				- \int_{\characteristicdiamondtwoarg{[\leftubar,\ubar)}{[\moreinterestingu_1,u]}}\left\{ \upmu \mytr_{\gtorus} \angk - \upmu \mytr_{\gtorus} \upchi\right\} \upupsilon \zeta \ReciprocalLunitAppliedtoTimeFunction \, \voldiamond\\
				& \ \ 
				+ \int_{\ingoingcharacteristicsurfacetwoarg{\ubar}{[\moreinterestingu_1,u]}} \upupsilon \zeta (\muX \ubar ) \ReciprocalLunitAppliedtoTimeFunction \, \volingoingnullhypersurface - \int_{\ingoingcharacteristicsurfacetwoarg{\leftubar}{[\moreinterestingu_1,u]}} \upupsilon \zeta ( \muX \ubar )\ReciprocalLunitAppliedtoTimeFunction \, \volingoingnullhypersurface \\
				& \ \ 
				+ \int_{\outgoingcharacteristicsurfacetwoarg{u}{[\leftubar,\ubar)}} \upupsilon \zeta \ReciprocalLunitAppliedtoTimeFunction \, \voloutgoingnullhypersurface - \int_{\outgoingcharacteristicsurfacetwoarg{\moreinterestingu_1}{[\leftubar,\ubar)}} \upupsilon \zeta \ReciprocalLunitAppliedtoTimeFunction \, \voloutgoingnullhypersurface,
			\end{split} \\
			\begin{split} \label{E:PRELIMINARYINTEGRATIONBYPARTSWITHMB}
				\int_{\characteristicdiamondtwoarg{[\leftubar,\ubar)}{[\moreinterestingu_1,u]}} (\upmu \Transport \upupsilon) \zeta \ReciprocalLunitAppliedtoTimeFunction \, \voldiamond & = - \int_{\characteristicdiamondtwoarg{[\leftubar,\ubar)}{[\moreinterestingu_1,u]}} \upupsilon (\upmu \Transport \zeta) \ReciprocalLunitAppliedtoTimeFunction \, \voldiamond \\
				& \ \
				- \int_{\characteristicdiamondtwoarg{[\leftubar,\ubar)}{[\moreinterestingu_1,u]}}  \upmu (\mytr_{\gtorus} \angk )\upupsilon \zeta \ReciprocalLunitAppliedtoTimeFunction \, \voldiamond\\
				& \ \ 
				+ \int_{\ingoingcharacteristicsurfacetwoarg{\ubar}{[\moreinterestingu_1,u]}}\frac{\upmu}{\ingoingmu} \upupsilon \zeta  \ReciprocalLunitAppliedtoTimeFunction \, \volingoingnullhypersurface - \int_{\ingoingcharacteristicsurfacetwoarg{\leftubar}{[\moreinterestingu_1,u]}} \frac{\upmu}{\ingoingmu} \upupsilon \zeta \ReciprocalLunitAppliedtoTimeFunction \, \volingoingnullhypersurface \\
				& \ \ 
				+ \int_{\outgoingcharacteristicsurfacetwoarg{u}{[\leftubar,\ubar)}} \upupsilon \zeta \ReciprocalLunitAppliedtoTimeFunction \, \voloutgoingnullhypersurface - \int_{\outgoingcharacteristicsurfacetwoarg{\moreinterestingu_1}{[\leftubar,\ubar)}} \upupsilon \zeta \ReciprocalLunitAppliedtoTimeFunction \, \voloutgoingnullhypersurface,
			\end{split}
		\end{align}
	\end{subequations}
\end{lemma}

\begin{proof}
	The identity \eqref{E:PRELIMINARYINTEGRATIONBYPARTSWITHL} follows from \eqref{E:DIVERGENCEIDENTITYSPACETIMEVECTORFIELD} with $\mathcal{J} \eqdef \frac{1}{\upmu} \upupsilon \zeta \Lunit$. Indeed, \eqref{E:RATIOOFNULLGEOSICINNERPRODUCTANDFOLIATIONDENSITY} and \eqref{E:LUNITANDULUNITAPPLIEDTOEIKONALANDCARTESIANTIME} imply $\mathcal{J}^u = 0$, $\mathcal{J}^{\ubar} \MagnitueofinnerproductofnewLandnewuL = \upupsilon \zeta$. Next, using  \eqref{E:SPACETIMEDIVERGENCEOFL}, we have:
	\begin{align} 
		\begin{split} \label{E:CANCELLATIONOFSINGULARMUINDIVERGENCEIDENTITIES}
			\Dfour_\alpha \mathcal{J}^\alpha \MagnitueofinnerproductofnewLandnewuL & = (\Lunit \upupsilon) \zeta \ReciprocalLunitAppliedtoTimeFunction + \upupsilon (\Lunit  \zeta)\ReciprocalLunitAppliedtoTimeFunction - \frac{\Lunit  \upmu}{\upmu} \upupsilon \zeta  \ReciprocalLunitAppliedtoTimeFunction+ \frac{\Lunit  \upmu}{\upmu} \upupsilon \zeta \ReciprocalLunitAppliedtoTimeFunction + (\mytr_{\gtorus} \upchi) \upupsilon \zeta \ReciprocalLunitAppliedtoTimeFunction \\
			& =  (\Lunit  \upupsilon) \zeta\ReciprocalLunitAppliedtoTimeFunction + \upupsilon (\Lunit  \zeta) \ReciprocalLunitAppliedtoTimeFunction + (\mytr_{\gtorus} \upchi) \upupsilon \zeta \ReciprocalLunitAppliedtoTimeFunction.
		\end{split}
	\end{align}
	
	The identity \eqref{E:PRELIMINARYINTEGRATIONBYPARTSWITHY} follows from similar arguments applied to $\mathcal{J} \eqdef \frac{1}{\upmu} \upupsilon \zeta \Yvf{A}$. We omit the details, noting only that applying \eqref{E:SPACETIMEDIVERGENCEOFYA} yields a cancellation of $\frac{\Yvf{A}\upmu}{\upmu}$ similar to what happened to $\frac{\Lunit \upmu}{\upmu}$ in \eqref{E:CANCELLATIONOFSINGULARMUINDIVERGENCEIDENTITIES}.	
	
	The identity \eqref{E:PRELIMINARYINTEGRATIONBYPARTSWITHMUX} follows from similar arguments applied to $\mathcal{J} \eqdef \frac{1}{\upmu} \upupsilon \zeta \muX$ and using \eqref{E:SPACETIMEDIVERGENCEOFMUX}. We omit the details, and only highlight the presence of the outgoing characteristic boundary terms on RHS\,\eqref{E:PRELIMINARYINTEGRATIONBYPARTSWITHMUX}. It is easy to compute these using \eqref{E:DIVERGENCEIDENTITYSPACETIMEVECTORFIELD} and $\muX u = 1$. 
	
	The identity \eqref{E:PRELIMINARYINTEGRATIONBYPARTSWITHMB} follows from similar arguments applied to $\mathcal{J} \eqdef \upupsilon \zeta \Transport$ and using \eqref{E:SPACETIMEDIFERGENCEOFB}. Indeed, \eqref{E:LFOLIATIONDENSITY}--\eqref{E:ULFOLIATIONDENSITY} implies $\mathcal{J}^{\ubar}  \MagnitueofinnerproductofnewLandnewuL = \frac{\upmu}{\ingoingmu} \upupsilon \zeta \ReciprocalLunitAppliedtoTimeFunction$ and $\mathcal{J}^u  \MagnitueofinnerproductofnewLandnewuL = \upupsilon \zeta \ReciprocalLunitAppliedtoTimeFunction $. We omit the rest of the details.
	
\end{proof}

\subsection{Swapping two differential operators via integration by parts} \label{SS:IBPTWODIFFERENTIALOPERATORS}

In our analysis we will encounter integrals that schematically look like $\int_{\characteristicdiamondtwoarg{[\leftubar,\ubar)}{[\moreinterestingu_1,u]}} P_1 \upupsilon P_2 \zeta$ and are uncontrollable at a consistent level of regularity. Here, the volume form is suppressed for convenience, and $P_1, P_2$ are first-order differential operators that are \emph{both} transverse to the ingoing characteristics, e.g., $\gfour(P_2,\newuL),\, \gfour(P_1,\newuL) \not\equiv 0$. In our settings, it turns out that the bilinear product $\int_{\characteristicdiamondtwoarg{[\leftubar,\ubar)}{[\moreinterestingu_1,u]}} P_2 \upupsilon P_1 \zeta$, with the differential operators switched, features a perfect divergence structure, and hence the regularity issues can be integrated away. The reason this is possible is because, in our contexts, $\upupsilon$ and $\zeta$ will satisfy PDEs of \emph{different} characteristic speeds, e.g. suppose $\upmu \Transport \upupsilon = \cdots$ and $\upmu \Box_{\gfour} \zeta = \cdots$. The point of the following lemma is then to reveal that the transversality assumptions on $P_1, P_2$ imply that the boundary integrals on $\ingoingcharacteristicsurfacetwoarg{\ubar}{[\moreinterestingu_1,u]}$ generated from this procedure will only feature $\newuL \zeta$ or $\nullgeop{x^A} \zeta$. \emph{These are the only derivatives that we can control on the ingoing characteristics for the wave variables}, see Lemma\,\ref{L:COERCIVENESSOFL2CONTROLLINGQUANITIESUNIFIEDCOMMUTATOR}. Similar comments apply when $P_1, \, P_2$ are transverse to $\ingoingcharacteristicsurfacetwoarg{u}{[\leftubar,\ubar)}$.

\begin{lemma}[Integrating by parts two differential operators]\label{L:IBPTWODIFFERENTIALOPERATORS}
Let $P_1, P_2$ be two first order differential operators on $\characteristicdiamondtwoarg{[\leftubar,\ubar)}{[\moreinterestingu_1,u]}$ that are transverse to $\ingoingcharacteristicsurfacetwoarg{\ubar}{[\moreinterestingu_1,u]}$. That is, $P_1, P_2 \in \textnormal{span}\{\newuL, L, \nullgeop{x^2},\nullgeop{x^3}\}$ and there exist scalar functions $\underline{\alpha}_i,\underline{\beta}_i,\underline{\gamma}_i^B$ with $\underline{\beta}_i \not\equiv 0$ and $i = 1,2$ such that:
	\begin{subequations} 
		\begin{align} 
			P_1 & = \underline{\alpha}_1 \newuL + \underline{\beta}_1 L + \underline{\gamma}_1^A \nullgeop{x^A}, \label{E:DECOMPOSITIONOFIBPP1INDOUBLENULLFRAME} \\
			P_2 & = \underline{\alpha}_2 \newuL + \underline{\beta}_2 L + \underline{\gamma}_2^A \nullgeop{x^A}. \label{E:DECOMPOSITIONOFIBPP2INDOUBLENULLFRAME}
		\end{align}
	\end{subequations}
Then the following identity holds for scalar function $\upupsilon$ and $\zeta$: 
\begin{align}
	\begin{split} \label{E:IBPTWODIFFERENTIALOPERATORS}
		 \int_{\characteristicdiamondtwoarg{[\leftubar,\ubar)}{[\moreinterestingu_1,u]}} (P_1 \upupsilon) (P_2 \zeta)  \ReciprocalLunitAppliedtoTimeFunction
		\, \voldiamond
		&  = \int_{\characteristicdiamondtwoarg{[\leftubar,\ubar)}{[\moreinterestingu_1,u]}} (P_2 \upupsilon) (P_1 \zeta) \ReciprocalLunitAppliedtoTimeFunction
		\, \voldiamond +  \int_{\characteristicdiamondtwoarg{[\leftubar,\ubar)}{[\moreinterestingu_1,u]}}  \{ \ErrorIBP[ \upupsilon, \pmb{\p} \zeta]\} \ReciprocalLunitAppliedtoTimeFunction
		\, \voldiamond \\
		& \ \ -  \int_{\ingoingcharacteristicsurfacetwoarg{\leftubar}{[\moreinterestingu_1,u]}}  \left\{  \upupsilon (\underline{\alpha}_2 \newuL \zeta + \underline{\gamma}_2^A \nullgeop{x^A} \zeta)\underline{\beta}_1  - \upupsilon (\underline{\alpha}_1 \newuL \zeta + \underline{\gamma}_1^A \nullgeop{x^A} \zeta) \underline{\beta}_2\right\} \, \volingoingnullhypersurface  \\
		& \ \ +  \int_{\ingoingcharacteristicsurfacetwoarg{\ubar}{[\moreinterestingu_1,u]}} \left\{  \upupsilon (\underline{\alpha}_2 \newuL \zeta + \underline{\gamma}_2^A \nullgeop{x^A} \zeta)\underline{\beta}_1  - \upupsilon (\underline{\alpha}_1 \newuL \zeta + \underline{\gamma}_1^A \nullgeop{x^A} \zeta) \underline{\beta}_2\right\} \, \volingoingnullhypersurface \\
		& \ \  -  \int_{\outgoingcharacteristicsurfacetwoarg{\moreinterestingu_1}{[\leftubar,\ubar)}} \left\{ \upupsilon (P_2 \zeta)\underline{\alpha}_1  - \upupsilon (P_1 \zeta) \underline{\alpha}_2 \right\} \ReciprocalLunitAppliedtoTimeFunction \, \voloutgoingnullhypersurface \\
		& \ \ +  \int_{\outgoingcharacteristicsurfacetwoarg{u}{[\leftubar,\ubar)}}\left\{ \upupsilon (P_2 \zeta)\underline{\alpha}_1  - \upupsilon (P_1 \zeta) \underline{\alpha}_2 \right\} \ReciprocalLunitAppliedtoTimeFunction \, \voloutgoingnullhypersurface
	\end{split}
\end{align}
where:
\begin{align}
	\begin{split} \label{E:IBPTWODIFFERENTIALOPERATORSFIRSTERRORTYPE}
		\ErrorIBP[ \upupsilon, \pmb{\p} \zeta]  \eqdef   \upupsilon   [P_2,P_1]\zeta + \frac{P_1\upmu}{\upmu} \upupsilon P_2 \zeta -  \upupsilon (P_2 \zeta) \Dfour_\alpha P_1^\alpha - \frac{P_2 \upmu}{\upmu} \upupsilon P_1 \zeta + \upupsilon (P_1 \zeta) \Dfour_\alpha P_2^\alpha.
	\end{split}
\end{align}
	
Assume now that both $P_1, \, P_2$ are also transverse to $\outgoingcharacteristicsurfacetwoarg{u}{[\leftubar,\ubar)}$. That is, when expressed with respect to the smooth commutator frame,
	\begin{subequations} 
		\begin{align} 
			P_1 & = \underline{\alpha}_1 \muX + \beta_1 L + \gamma_1^A \Yvf{A}, \label{E:DECOMPOSITIONOFIBPP1INCOMMUTATORFRAME} \\
			P_2 & = \underline{\alpha}_2 \muX + \beta_2 L + \gamma_2^A \Yvf{A}. \label{E:DECOMPOSITIONOFIBPP2INCOMMUTATORFRAME}
		\end{align}
	\end{subequations}
where $\underline{\alpha}_i $ are the same constants as in \eqref{E:DECOMPOSITIONOFIBPP1INDOUBLENULLFRAME}--\eqref{E:DECOMPOSITIONOFIBPP1INDOUBLENULLFRAME} and now satisfy $\underline{\alpha}_i \not\equiv 0$. Then the boundary integrals on the outgoing characteristic surfaces only feature $\outgoingcharacteristicsurfacetwoarg{u}{[\leftubar,\ubar)}$-tangential derivatives. That is, the corresponding integrand on RHS\,\eqref{E:IBPTWODIFFERENTIALOPERATORS} is precisely: 
\begin{align} \label{E:GOODTANGENTIALDERIVATIVESONOUTGOINGCHARACTERISTICSAFTERTWOIBP}
	\upupsilon (P_2 \zeta)\underline{\alpha}_1  - \upupsilon (P_1 \zeta) \underline{\alpha}_2 = \upupsilon \left\{ \beta_2  \Lunit \zeta + \gamma_2^A\Yvf{A} \zeta \right\} \underline{\alpha}_1 - \upupsilon \left\{ \beta_1  \Lunit \zeta + \gamma_1^A\Yvf{A} \zeta \right\} \underline{\alpha}_2.
\end{align}
\end{lemma}

\begin{proof}
We first integrate by parts off the $P_1$ from $P_1 \upupsilon$ on LHS \eqref{E:IBPTWODIFFERENTIALOPERATORS}. More specifically, we apply the divergence theorem identity \eqref{E:IDFORDETWEIGHTEDCOVARIANTDIVERGENCEOFSPACETIMEVECTORFIELD} with $\mathcal{J}_1 = \frac{1}{\upmu} \upupsilon (P_2 \zeta) P_1$. To analyze the boundary terms generated from this, we use \eqref{E:RATIOOFNULLGEOSICINNERPRODUCTANDFOLIATIONDENSITY}, \eqref{E:LUNITANDULUNITAPPLIEDTOEIKONALANDCARTESIANTIME}, and \eqref{E:DECOMPOSITIONOFIBPP1INDOUBLENULLFRAME} to see that $\mathcal{J}_1^{\ubar} \MagnitueofinnerproductofnewLandnewuL = \upupsilon (P_2 \zeta) \beta_1$ and $\mathcal{J}_1^u \MagnitueofinnerproductofnewLandnewuL = \upupsilon (P_2 \zeta)\alpha_1 \ReciprocalLunitAppliedtoTimeFunction$. Applying \eqref{E:IDFORDETWEIGHTEDCOVARIANTDIVERGENCEOFSPACETIMEVECTORFIELD} then proves that LHS \eqref{E:IBPTWODIFFERENTIALOPERATORS} is equal to:
\begin{align}
	\begin{split} \label{E:STEP1IBPTWODIFFERENTIALOPERATORS}
		& \int_{\characteristicdiamondtwoarg{[\leftubar,\ubar)}{[\moreinterestingu_1,u]}} \left\{ - \upupsilon (P_2 P_1 \zeta) + \upupsilon [P_2, P_1] \zeta + \frac{P_1 \upmu}{\upmu} \upupsilon (P_2 \zeta) -  \upupsilon (P_2 \zeta) \Dfour_\alpha P_1^\alpha\right\} \ReciprocalLunitAppliedtoTimeFunction \, \voldiamond \\
		& \ \ -  \int_{\ingoingcharacteristicsurfacetwoarg{\leftubar}{[\moreinterestingu_1,u]}}  \upupsilon (P_2 \zeta) \beta_1  \, \volingoingnullhypersurface +  \int_{\ingoingcharacteristicsurfacetwoarg{\ubar}{[\moreinterestingu_1,u]}} \upupsilon (P_2 \zeta) \beta_1  \, \volingoingnullhypersurface -  \int_{\outgoingcharacteristicsurfacetwoarg{\moreinterestingu_1}{[\leftubar,\ubar)}} \upupsilon (P_2 \zeta)\alpha_1 \ReciprocalLunitAppliedtoTimeFunction   \, \voloutgoingnullhypersurface +  \int_{\outgoingcharacteristicsurfacetwoarg{u}{[\leftubar,\ubar)}} \upupsilon (P_2 \zeta)\alpha_1 \ReciprocalLunitAppliedtoTimeFunction  \, \voloutgoingnullhypersurface
	\end{split}
\end{align}
Applying \eqref{E:IDFORDETWEIGHTEDCOVARIANTDIVERGENCEOFSPACETIMEVECTORFIELD} again, but this time for $\mathcal{J}_2 = \frac{1}{\upmu} \upupsilon (P_1 \zeta)P_2$, the first spacetime integral in \eqref{E:STEP1IBPTWODIFFERENTIALOPERATORS} is precisely equal to:
\begin{align}
	\begin{split} \label{E:STEP2IBPTWODIFFERENTIALOPERATORS}
		& \int_{\characteristicdiamondtwoarg{[\leftubar,\ubar)}{[\moreinterestingu_1,u]}} \left\{ P_2\upupsilon (P_1 \zeta) - \frac{P_2 \upmu}{\upmu} \upupsilon (P_1 \zeta) +  \upupsilon (P_1 \zeta) \Dfour_\alpha P_2^\alpha\right\} \ReciprocalLunitAppliedtoTimeFunction \, \voldiamond \\
		& \ \ +  \int_{\ingoingcharacteristicsurfacetwoarg{\leftubar}{[\moreinterestingu_1,u]}}  \upupsilon (P_1 \zeta) \beta_2  \, \volingoingnullhypersurface -  \int_{\ingoingcharacteristicsurfacetwoarg{\ubar}{[\moreinterestingu_1,u]}} \upupsilon (P_1 \zeta) \beta_2  \, \volingoingnullhypersurface +  \int_{\outgoingcharacteristicsurfacetwoarg{\moreinterestingu_1}{[\leftubar,\ubar)}} \upupsilon (P_1 \zeta)\alpha_2 \ReciprocalLunitAppliedtoTimeFunction   \, \voloutgoingnullhypersurface -  \int_{\outgoingcharacteristicsurfacetwoarg{u}{[\leftubar,\ubar)}} \upupsilon (P_1 \zeta)\alpha_2 \ReciprocalLunitAppliedtoTimeFunction  \, \voloutgoingnullhypersurface
	\end{split}
\end{align}
Next we observe that the integrands on the ingoing characteristic surfaces generated from the two applications of the divergence theorems above are precisely $\upupsilon (P_2 \zeta) \beta_1 -  \upupsilon (P_1 \zeta) \beta_2$. Using \eqref{E:DECOMPOSITIONOFIBPP1INDOUBLENULLFRAME}--\eqref{E:DECOMPOSITIONOFIBPP2INDOUBLENULLFRAME}, we see \emph{a complete cancellation of all $\Lunit$-derivatives.} In particular, this difference is equal to $\upupsilon (\underline{\alpha}_2 \newuL \zeta + \underline{\gamma}_2^A \nullgeop{x^A} \zeta)\beta_1  - \upupsilon (\underline{\alpha}_1 \newuL \zeta + \underline{\gamma}_1^A \nullgeop{x^A} \zeta) \beta_2$. Adding together \eqref{E:STEP1IBPTWODIFFERENTIALOPERATORS}--\eqref{E:STEP2IBPTWODIFFERENTIALOPERATORS} and applying these observations concludes the proof of \eqref{E:IBPTWODIFFERENTIALOPERATORS}.

Now suppose both $P_1,P_2$ are also transverse to $\outgoingcharacteristicsurfacetwoarg{u}{[\leftubar,\ubar)}$, that is $\underline{\alpha}_i \not\equiv 0$ in \eqref{E:DECOMPOSITIONOFIBPP1INDOUBLENULLFRAME}--\eqref{E:DECOMPOSITIONOFIBPP1INDOUBLENULLFRAME}. Then since $P_1 u = \underline{\alpha_1}$ and $P_2 u = \underline{\alpha_2}$, the $\muX$-derivatives in the integrands on the outgoing characteristic boundary integrals cancel: $ \upupsilon (P_2 \zeta)\alpha_1  - \upupsilon (P_1 \zeta)\alpha_2  = \upupsilon (\beta_2 \Lunit \zeta + \gamma_2^A \Yvf{A} \zeta)\underline{\alpha}_1  - \upupsilon (\beta_1 \Lunit  \zeta + \gamma_1^A \Yvf{A} \zeta) \underline{\alpha}_2$. The identity \eqref{E:GOODTANGENTIALDERIVATIVESONOUTGOINGCHARACTERISTICSAFTERTWOIBP} then follows easily from \eqref{E:SMOOTHANGULARDERIVATIVESINTERMSOFDOUBLENULLONESANDL}. 

\end{proof}

\section{Integration by parts identities needed for the acoustic geometry} \label{S:IBPIDENTITIESNEEDFORTOPORDERENERGYESTIMATES}

\subsection{Integration by parts involving $\Lunit$} \label{SS:IBPIDENDITIESINVOLVINGL}

The identity in the following lemma will be used to control top-order wave equation error terms that are tied to integrating by parts derivative-losing acoustic geometry terms with an $\Lunit$.

\begin{lemma}[The key top-order integration by parts in $\Lunit$ identity tied to the acoustic geometry]
	\label{L:KEYIBPIDENTIFYFORWAVEEQUATIONENERGYESTIMATESINVOLVINGL} 
	Let $N = \Ntop$, $(\ubar,u) \in [\leftubar,\ubarboot)\times[\moreinterestingu_1,\moreinterestingu_2]$,
	and let $\varphi$ and $\upeta$ be scalar functions on $\characteristicdiamondtwoarg{[\leftubar,\ubar)}{[\moreinterestingu_1,u]}$. Let $ \tander^N \in \mathfrak{P}^{(N)}$, where $\mathfrak{P}^{(N)}$ denotes the set of top-order commutator operators from Def.\,\ref{D:STRINGSOFCOMMUTATIONVECTORFIELDS}.
	Then the following integration by parts identity holds:
	\begin{align}
		\begin{split} \label{E:KEYIBPIDENTIFYFORWAVEEQUATIONENERGYESTIMATES}
			&
			\int_{\characteristicdiamondtwoarg{[\leftubar,\ubar)}{[\moreinterestingu_1,u]}}
			 (\muX \varphi) 
			(\Lunit  \tander^N \varphi) 
			\Yvf{A} \upeta \, \weight^{\blowupratetoporderwave} \ReciprocalLunitAppliedtoTimeFunction
			\, \voldiamond 
			\\
			& 
			= 
			\int_{\characteristicdiamondtwoarg{[\leftubar,\ubar)}{[\moreinterestingu_1,u]}}  
			(\muX \varphi) 
			(\Yvf{A} \tander^N \varphi) 
			\Lunit \upeta \, \weight^{\blowupratetoporderwave} \ReciprocalLunitAppliedtoTimeFunction
			\, \voldiamond \\
			& \ \
			- \int_{\ingoingcharacteristicsurfacetwoarg{\ubar}{[\moreinterestingu_1,u]}} 
			( \muX \varphi) \left\{ \frac{1}{\Speed^2} (\gtorus^{-1})^{AB} \gtorusdoublenullCOV_B^C \nullgeop{x^C}\tander^N \varphi \right\} \upeta \, \weight^{\blowupratetoporderwave}
			\, \volingoingnullhypersurface 
			\\
			&  \ \ 
			+ 
			\int_{\ingoingcharacteristicsurfacetwoarg{\leftubar}{[\moreinterestingu_1,u]}} 
			(\muX \varphi) \left\{ \frac{1}{\Speed^2} (\gtorus^{-1})^{AB} \gtorusdoublenullCOV_B^C \nullgeop{x^C} \tander^N  \varphi \right\} \upeta \, \weight^{\blowupratetoporderwave}
			\, \volingoingnullhypersurface 
			\\
			& \ \ 
			+ 
			\int_{\characteristicdiamondtwoarg{[\leftubar,\ubar)}{[\moreinterestingu_1,u]}}
			\left\{ \ErrorIBP_1^{AG;L}[\tander^N \varphi;\upeta;\Yvf{A}]+ \ErrorIBP_2^{AG;L}[\tander^N \varphi;\upeta;\Yvf{A}]  \right\}\ReciprocalLunitAppliedtoTimeFunction
			\, \voldiamond,
		\end{split}
	\end{align} 
	where:  
	\begin{align} 
		\begin{split} \label{E:ERRORTERM1KEYIBPIDENTIFYFORWAVEEQUATIONENERGYESTIMATES}
			\ErrorIBP_1^{AG;L}[\tander^N \varphi;\upeta;\Yvf{A}]
			&  
			\eqdef  \weight^{\blowupratetoporderwave} \Big\{
			(\muX \varphi) 
			([\Lunit, \Yvf{A}] \tander^N \varphi) 
			\upeta  
			\\
			& \ \
			+ 
			(\Lunit \muX \varphi) 
			(\Yvf{A} \tander^N \varphi)
			\upeta 
			- 
			(\Yvf{A} \muX \varphi) 
			(\Lunit \tander^N \varphi)
			\upeta  \\
			& \ \ 
			\left. - \frac{1}{2} (\muX \varphi) (\Lunit \tander^N \varphi)\upeta \mytr_{\gtorus} \, \angdeform{\Yvf{A}} +  (\muX \varphi) (\Yvf{A} \tander^N \varphi )\upeta \mytr_{\gtorus} \upchi \right\} ,
		\end{split} \\
		\begin{split}
			\ErrorIBP_2^{AG;L}[\tander^N \varphi;\upeta;\Yvf{A}] & \eqdef - \blowupratetoporderwave \, \weight^{\blowupratetoporderwave - 1} (\Yvf{A} \weight)(\muX \varphi) (\Lunit \tander^N \varphi) \upeta  \\
			& \ \  +  \blowupratetoporderwave \, \weight^{\blowupratetoporderwave - 1} (\newL  \weight) (\muX \varphi) (\Yvf{A} \tander^N \varphi) \upeta \label{E:ERRORTERM2KEYIBPIDENTIFYFORWAVEEQUATIONENERGYESTIMATES}
		\end{split}
	\end{align}
	
\end{lemma}

\begin{proof}
	We begin by applying \eqref{E:IBPTWODIFFERENTIALOPERATORS} with $\upupsilon = \upeta$, $\zeta =  (\muX\varphi) \Lunit \tander^N \varphi \, \weight^{\blowupratetoporderwave}$, $L = P_1$, and $\Yvf{A} = P_2$. Tedious yet straightforward applications of the Leibniz rule applied to the error integral $ \int_{\characteristicdiamondtwoarg{[\leftubar,\ubar)}{[\moreinterestingu_1,u]}}  \{ \ErrorIBP[ \upupsilon, \pmb{\p} \zeta]\} \ReciprocalLunitAppliedtoTimeFunction
		\, \voldiamond$ on RHS\,\eqref{E:IBPTWODIFFERENTIALOPERATORS} produce the error integrals $\int_{\characteristicdiamondtwoarg{[\leftubar,\ubar)}{[\moreinterestingu_1,u]}}
			\left\{ \ErrorIBP^{AG;L}_1[\tander^N \varphi;\upeta;\Yvf{A}]+ \ErrorIBP^{AG;L}_2[\tander^N \varphi;\upeta;\Yvf{A}]  \right\} \ReciprocalLunitAppliedtoTimeFunction
			\, \voldiamond$ presend on RHS\,\eqref{E:KEYIBPIDENTIFYFORWAVEEQUATIONENERGYESTIMATES}. We omit most of these details and focus only the most delicate steps. Firstly, the terms in the boundary integrals follow from \eqref{E:YCOMMUTATORINTERMSOFNULLCOORDINATEVECTORFIELDSANDL} and \eqref{E:IBPTWODIFFERENTIALOPERATORS}; note that all of the $\Lunit$ derivatives cancel. Secondly,
inserting the expressions for $\Dfour_\alpha \Lunit^\alpha$ and $\Dfour_\alpha \Yvf{A}^\alpha$ from \eqref{E:SPACETIMEDIVERGENCEOFL}--\eqref{E:SPACETIMEDIVERGENCEOFYA} into \eqref{E:IBPTWODIFFERENTIALOPERATORSFIRSTERRORTYPE}, respectively, we see a total cancellation of the terms featuring singular $\upmu^{-1}$-coefficients.

We also point out that the outgoing characteristic boundary integrals from \eqref{E:IBPTWODIFFERENTIALOPERATORS} are not present in \eqref{E:KEYIBPIDENTIFYFORWAVEEQUATIONENERGYESTIMATES} because $\alpha_1 = \Lunit u = \Yvf{A} u = \alpha_2 = 0$. 	
	
\end{proof}

\section{Integration by parts identities needed for the derivative-losing specific vorticity} \label{S:IBPIDENTITIESNEEDEDFORTOPORDERVORTICITY}

The integration by parts identity in the following lemma is used to find a perfect divergence structure in the nonlinear interaction between vorticity and the velocity variables. In particular, it will be used to avoid derivative-loss at top order.

\begin{lemma}[The key top-order integration by parts identity tied to vorticity] \label{L:IBPIDENTITYNEEDEDFORTOPORDERVORTICITY}
For $\Sigma_t$-tangent vectorfields $\SigmatTan_1, \, \SigmatTan_2$, denote their Cartesian inner product by $\SigmatTan_1\cdot\SigmatTan_2 = \sum_{i=1}^3\SigmatTan_1^i\cdot\SigmatTan_2^i$. Then the following integration by parts identity holds: 
\begin{align}
	\begin{split} \label{E:IBPIDENTITYNEEDEDFORTOPORDERVORTICITY}
		 & \int_{\characteristicdiamondtwoarg{[\leftubar,\ubar)}{[\moreinterestingu_1,u]}} \Speed^2 \exp(\LogDensity) \upmu \left\{ \Flatcurl(\tander^N \vortrenormalized)\right\}  \cdot (\multipliervectorfield \tander^N v) \weight^{\blowupratetoporderwave} \ReciprocalLunitAppliedtoTimeFunction \, \voldiamond \\
		 & = - \int_{\characteristicdiamondtwoarg{[\leftubar,\ubar)}{[\moreinterestingu_1,u]}} \Speed^2  \exp(2\LogDensity) (\multipliervectorfield \tander^N \vortrenormalized)\cdot (\upmu \tander^N \vortrenormalized) \weight^{\blowupratetoporderwave} \ReciprocalLunitAppliedtoTimeFunction \, \voldiamond \\ 
		 & \ \ + \int_{\characteristicdiamondtwoarg{[\leftubar,\ubar)}{[\moreinterestingu_1,u]}} \left\{\Speed^2\exp(\LogDensity)  (\multipliervectorfield \tander^N \vortrenormalized)\cdot [\upmu,\tander^N] \Flatcurl v + \Speed^2 \exp(\LogDensity) (\multipliervectorfield \tander^N \vortrenormalized)\cdot [\tander^N,\upmu \Flatcurl] v \right\}\weight^{\blowupratetoporderwave} \ReciprocalLunitAppliedtoTimeFunction \, \voldiamond \\
		 & \ \  -  \int_{\characteristicdiamondtwoarg{[\leftubar,\ubar)}{[\moreinterestingu_1,u]}} \Speed^2\exp(\LogDensity) (\multipliervectorfield \tander^N \vortrenormalized^i)  \left\{\sum_{\substack{ \tander^{N_1}\tander^{N_2} = \tander^N \\ N_1 \le N-1}} \upmu (\tander^{N_1} \vortrenormalized^i) \tander^{N_2}(\exp(\LogDensity)) \right\} \weight^{\blowupratetoporderwave} \ReciprocalLunitAppliedtoTimeFunction \, \voldiamond \\
		 & \ \  + 
			\int_{\characteristicdiamondtwoarg{[\leftubar,\ubar)}{[\moreinterestingu_1,u]}}
			\left\{ \ErrorIBP_1^{\vortrenormalized;\text{top}}[\tander^N] + \ErrorIBP_2^{\vortrenormalized;\text{top}}[\tander^N]  \right\}\ReciprocalLunitAppliedtoTimeFunction
			\, \voldiamond, \\
		& \ \  -  \int_{\ingoingcharacteristicsurfacetwoarg{\leftubar}{[\moreinterestingu_1,u]}} \left\{ \Speed^2 \exp(\LogDensity)(\tander^N\vortrenormalized^b) \left( 2 \newuL \tander^N v^i +2 \MagnitueofinnerproductofnewLandnewuL \ToriTangentVectorfieldAssociatedToDoubleNullFolliations \tander^N v^i\right) \upepsilon_{iab}\left\{\left(\frac{\MagnitueofinnerproductofnewLandnewuL}{\ingoingmu} - \upmu\right)\frac{X^a}{\Speed^2} + \upmu \angularcoeffinmuweightedspatialcartesian_a^A \frac{\Yvf{A} \ubar}{\Lunit \ubar}\right\}  \right. \\
		& \ \ \ \ \ \ \ \ \ \ \left. - \Speed^2  \exp(\LogDensity)(\tander^N\vortrenormalized^b) \upepsilon_{iab}  \left(\frac{X^a}{\Speed^2} \newuL v^i  + \frac{X^a}{\Speed^2} \MagnitueofinnerproductofnewLandnewuL \ToriTangentVectorfieldAssociatedToDoubleNullFolliations \tander^N v^i + \frac{\upmu}{\Speed^2}\angularcoeffinmuweightedspatialcartesian_a^A (\gtorus^{-1})^{AB}  \gtorusdoublenullCOV_B^C \nullgeop{x^C} \tander^N v^i\right) \left(1 + \frac{2 \MagnitueofinnerproductofnewLandnewuL}{\ingoingmu}\right) \right\} \weight^{\blowupratetoporderwave} \volingoingnullhypersurface \\
		& \ \ +  \int_{\ingoingcharacteristicsurfacetwoarg{\ubar}{[\moreinterestingu_1,u]}} \left\{ \Speed^2 \exp(\LogDensity)(\tander^N\vortrenormalized^b) \left( 2 \newuL \tander^N v^i + 2\MagnitueofinnerproductofnewLandnewuL \ToriTangentVectorfieldAssociatedToDoubleNullFolliations \tander^N v^i\right) \upepsilon_{iab}\left\{\left(\frac{\MagnitueofinnerproductofnewLandnewuL}{\ingoingmu} - \upmu\right)\frac{X^a}{\Speed^2} + \upmu \angularcoeffinmuweightedspatialcartesian_a^A \frac{\Yvf{A} \ubar}{\Lunit \ubar}\right\}  \right. \\
		& \ \ \ \ \ \ \ \ \ \  \left. - \Speed^2  \exp(\LogDensity)(\tander^N\vortrenormalized^b) \upepsilon_{iab}  \left(\frac{X^a}{\Speed^2} \newuL v^i  + \frac{X^a}{\Speed^2} \MagnitueofinnerproductofnewLandnewuL \ToriTangentVectorfieldAssociatedToDoubleNullFolliations \tander^N v^i + \frac{\upmu}{\Speed^2}\angularcoeffinmuweightedspatialcartesian_a^A (\gtorus^{-1})^{AB}  \gtorusdoublenullCOV_B^C \nullgeop{x^C} \tander^N v^i\right) \left(1 + \frac{2 \MagnitueofinnerproductofnewLandnewuL}{\ingoingmu}\right) \right\} \weight^{\blowupratetoporderwave} \volingoingnullhypersurface  \\
		& \ \ -  \int_{\outgoingcharacteristicsurfacetwoarg{\moreinterestingu_1}{[\leftubar,\ubar)}}  \left\{  \exp(\LogDensity)(\tander^N\vortrenormalized^b) ( (1 + 2\upmu) \Lunit \tander^N v^i) \upepsilon_{iab}X^a - 2 \exp(\LogDensity) \Speed^2 \left(\tander^N \vortrenormalized^b\right) \left( \upepsilon_{iab} \upmu \angularcoeffinmuweightedspatialcartesian_a^A \Yvf{A} \tander^N v^i\right)  \right\} \weight^{\blowupratetoporderwave} \ReciprocalLunitAppliedtoTimeFunction \, \voloutgoingnullhypersurface  \\
		& \ \ +  \int_{\outgoingcharacteristicsurfacetwoarg{u}{[\leftubar,\ubar)}}  \left\{  \exp(\LogDensity)(\tander^N\vortrenormalized^b) ( (1 + 2\upmu) \Lunit \tander^N v^i) \upepsilon_{iab}X^a - 2\exp(\LogDensity) \Speed^2 \left(\tander^N \vortrenormalized^b\right) \left( \upepsilon_{iab} \upmu \angularcoeffinmuweightedspatialcartesian_a^A \Yvf{A} \tander^N v^i\right)  \right\} \weight^{\blowupratetoporderwave} \ReciprocalLunitAppliedtoTimeFunction \, \voloutgoingnullhypersurface  
	\end{split}
\end{align}
where:
\begin{subequations}
	\begin{align}
		\begin{split} \label{E:ERRORTERM1IBPIDENTITYNEEDEDFORTOPORDERVORTICITY}
			\ErrorIBP_1^{\vortrenormalized;\text{top}}[\tander^N] & \eqdef  \weight^{\blowupratetoporderwave} \left\{ \Speed^2 \exp(\LogDensity) \left(\tander^N\vortrenormalized^b\right) [\multipliervectorfield, \upepsilon_{iab} \upmu \p_a] \tander^N v^i  \right. \\
			& \ \  - 2  \exp(\LogDensity) \Speed \left( \Speed_{;\LogDensity} \upepsilon_{iab} \upmu (\p_a \LogDensity) + \Speed_{;\Ent} \upepsilon_{iab} \upmu (\p_a \Ent) \right) \left(\tander^N \vortrenormalized^b\right) \multipliervectorfield  \tander^N v^i   \\
			& \ \ - \Speed^2 \left( \upepsilon_{iab} \p_a (\exp(\LogDensity))\right)  \left(\tander^N \vortrenormalized^b\right) \multipliervectorfield  \tander^N v^i \\
			& \ \ - \Speed^2 \exp(\LogDensity) \upepsilon_{iab}\left(  \muX \left( \frac{X^a}{\Speed^2}\right) +\Yvf{A} \angularcoeffinmuweightedspatialcartesian_a^A + \frac{1}{2} \angularcoeffinmuweightedspatialcartesian_a^A \, \mytr_{\gtorus} \angdeform{\Yvf{A}} + \frac{\upmu X^a}{\Speed^2}  \mytr_{\gtorus} \angk - \frac{\upmu X^a}{\Speed^2} \mytr_{\gtorus}\upchi\right)  \left(\tander^N \vortrenormalized^b\right) \multipliervectorfield  \tander^N v^i \\
			&  \ \ + 2 \Speed \left( \Speed_{;\LogDensity}(\multipliervectorfield \LogDensity) + \Speed_{;\Ent} (\multipliervectorfield \Ent) \right)  \exp(\LogDensity)\left(\tander^N \vortrenormalized^b\right) \upepsilon_{iab}\upmu \p_a (\tander^N v^i) \\
			& \ \ + \Speed^2 \left(\multipliervectorfield \exp(\LogDensity)\right)  \left(\tander^N \vortrenormalized^b\right) \upepsilon_{iab}\upmu \p_a (\tander^N v^i) \\
			& \ \  \left. + \Speed^2 \exp(\LogDensity) \left(  2  \Lunit \upmu  + 2 \upmu \mytr_{\gtorus}\angk + \mytr_{\gtorus}\upchi \right) \left(\tander^N \vortrenormalized^b\right) \upepsilon_{iab}\upmu \p_a (\tander^N v^i) \right\},
		\end{split} \\
		\begin{split} \label{E:ERRORTERM2IBPIDENTITYNEEDEDFORTOPORDERVORTICITY}
			\ErrorIBP_2^{\vortrenormalized;\text{top}}[\tander^N] & \eqdef - \blowupratetoporderwave\weight^{\blowupratetoporderwave-1} \upepsilon_{iab}(\upmu \p_a \weight)\Speed^2 \exp(\LogDensity) \left(\tander^N \vortrenormalized^b\right) \multipliervectorfield  \tander^N v^i \\
			& \ \  + \blowupratetoporderwave\weight^{\blowupratetoporderwave-1}( \multipliervectorfield \weight)\Speed^2 \exp(\LogDensity) \left(\tander^N \vortrenormalized^b\right) \upepsilon_{iab}\upmu \p_a (\tander^N v^i) .
		\end{split}
	\end{align}
\end{subequations}
\end{lemma}

\begin{proof}
We first expand the Euclidean inner product and curl operator in the integrand of LHS \eqref{E:IBPIDENTITYNEEDEDFORTOPORDERVORTICITY} as $\Speed^2\exp(\LogDensity) \upmu \upepsilon_{iab}\p_a (\tander^N \vortrenormalized^b) ( \multipliervectorfield \tander^N v^i) \weight^{\blowupratetoporderwave}$. 
For each $i,a,b = 1,2,3$, we consider the integration by parts identity \eqref{E:IBPTWODIFFERENTIALOPERATORS} with $P_1 = \upmu \upepsilon_{iab}\p_a$, $P_2 = \multipliervectorfield$, $\upupsilon =  \Speed^2 \left(\tander^N \vort^b\right)  \weight^{\blowupratetoporderwave}$, and $\zeta = \tander^N v^i$ and use Lemma\,\ref{L:MOREUSEFULFRAMEDECOMPOSITIONS}, in particular \eqref{E:MULTIPLIERINTERMSOFDOUBLENULLFRAME} and \eqref{E:MUWEIGHTEDCARTESIANINTERMSOFDOUBLENULLFRAME}. 

We omit most of the straightforward, yet tedious details needed for the argument in the previous paragraph. However, we do make some clarifying remarks. 

Firstly, the $\frac{1}{\upmu}$-singular terms in RHS\,\eqref{E:IBPTWODIFFERENTIALOPERATORSFIRSTERRORTYPE} are all cancelled by the $\Dfour_\alpha P_1^\alpha$ and $\Dfour_\alpha P_2^\alpha$ terms, see \eqref{E:SPACETIMEDIVERGENCEOFMULTIPLIER}--\eqref{E:SPACETIMEDIVERGENCEOFMUWEIGHTEDCARTESIANSPATIALDERIVATIVES}.

Secondly, we now explain the presence of the first three spacetime integrals on RHS\,\eqref{E:IBPIDENTITYNEEDEDFORTOPORDERVORTICITY}. In applying the procedure from the above discussions, we come to the following integrand: $\Speed^2 \exp(\LogDensity) \multipliervectorfield(\tander^N \vortrenormalized^b) \upepsilon_{iab} \upmu \p_a \tander^N v^i$. From the antisymmetry of $\upepsilon_{iab}$, this is precisely equal to $- \Speed^2 \exp(\LogDensity) \multipliervectorfield(\tander^N \vortrenormalized) \cdot \upmu \Flatcurl \tander^N v$. By differentiating \eqref{E:SPECIFICVORTICITYDEF} with respect to $\tander^N$ and then multiplying by $\upmu$, straightforward commutations justify the first three spacetime integrals on RHS\,\eqref{E:IBPIDENTITYNEEDEDFORTOPORDERVORTICITY}. 

Finally, the precise form of all of the boundary integrals on ingoing and outgoing characteristic surfaces in \eqref{E:IBPIDENTITYNEEDEDFORTOPORDERVORTICITY} follows from \eqref{E:MULTIPLIERINTERMSOFDOUBLENULLFRAME}, \eqref{E:MUWEIGHTEDCARTESIANINTERMSOFDOUBLENULLFRAME}--\eqref{E:EXPRESSIONFORMUWEIGHTEDSPATIALCARTESIANDERIVATIVESWRTCOMMUTATORS}, and the precise form of the boundary terms in \eqref{E:IBPTWODIFFERENTIALOPERATORS} and \eqref{E:GOODTANGENTIALDERIVATIVESONOUTGOINGCHARACTERISTICSAFTERTWOIBP}.

\end{proof}

\begin{lemma}[Perfect divergence structure for derivative-losing vorticity] \label{L:PERFECTDIVERGENCESTRUCTUREFOROVERDIFFERENTIATEDVORTICITY}
For $\Sigma_t$-tangent vectorfields $\SigmatTan_1, \, \SigmatTan_2$, denote their Cartesian inner product by $\SigmatTan_1\cdot\SigmatTan_2 = \sum_{i=1}^3\SigmatTan_1^i\cdot\SigmatTan_2^i$. Then the following identity holds:
	\begin{align}
		\begin{split} \label{E:PERFECTDIVERGENCESTRUCTUREFOROVERDIFFERENTIATEDVORTICITY}
		&  - \int_{\characteristicdiamondtwoarg{[\leftubar,\ubar)}{[\moreinterestingu_1,u]}} \Speed^2 \exp(2\LogDensity) (\multipliervectorfield \tander^N \vortrenormalized)\cdot (\upmu \tander^N \vortrenormalized) \weight^{\blowupratetoporderwave} \ReciprocalLunitAppliedtoTimeFunction \, \voldiamond 
		= - 2  \int_{\characteristicdiamondtwoarg{[\leftubar,\ubar)}{[\moreinterestingu_1,u]}} \Speed^2 \exp(2\LogDensity)  (\upmu \Transport \tander^N \vortrenormalized)\cdot (\upmu \tander^N \vortrenormalized) \weight^{\blowupratetoporderwave} \ReciprocalLunitAppliedtoTimeFunction \, \voldiamond \\
		 & \ \ - \frac{1}{2} \int_{\ingoingcharacteristicsurfacetwoarg{\ubar}{[\moreinterestingu_1,u]}}  \upmu \Speed^2 \exp(2\LogDensity) |\tander^N \vortrenormalized|^2 \weight^{\blowupratetoporderwave}  \volingoingnullhypersurface +  \frac{1}{2} \int_{\ingoingcharacteristicsurfacetwoarg{\leftubar}{[\moreinterestingu_1,u]}}  \upmu \Speed^2 \exp(2\LogDensity) |\tander^N \vortrenormalized|^2 \weight^{\blowupratetoporderwave}  \volingoingnullhypersurface \\
		 & \ \ +  \int_{\characteristicdiamondtwoarg{[\leftubar,\ubar)}{[\moreinterestingu_1,u]}} \upmu  \exp(2\LogDensity)\Speed (\Speed_{;\LogDensity} \Lunit \LogDensity + \Speed_{;\Ent} \Lunit \Ent) |\tander^N \vortrenormalized|^2 \weight^{\blowupratetoporderwave} \ReciprocalLunitAppliedtoTimeFunction \, \voldiamond \\
		 & \ \ +  \int_{\characteristicdiamondtwoarg{[\leftubar,\ubar)}{[\moreinterestingu_1,u]}} \upmu \Lunit\left( \exp(2\LogDensity)\right)\Speed^2 |\tander^N \vortrenormalized|^2 \weight^{\blowupratetoporderwave} \ReciprocalLunitAppliedtoTimeFunction \, \voldiamond + \frac{1}{2}\int_{\characteristicdiamondtwoarg{[\leftubar,\ubar)}{[\moreinterestingu_1,u]}} (\Lunit \upmu) \Speed^2 \exp(2\LogDensity) |\tander^N \vortrenormalized|^2 \weight^{\blowupratetoporderwave} \ReciprocalLunitAppliedtoTimeFunction \, \voldiamond \\
		 & \ \ + \frac{\blowupratetoporderwave}{2}  \int_{\characteristicdiamondtwoarg{[\leftubar,\ubar)}{[\moreinterestingu_1,u]}}  \upmu \Speed^2 \exp(2\LogDensity)  |\tander^N \vortrenormalized|^2  (\newL \weight) \weight^{\blowupratetoporderwave-1} \, \voldiamond \\
		 & \ \ + \frac{1}{2} \int_{\characteristicdiamondtwoarg{[\leftubar,\ubar)}{[\moreinterestingu_1,u]}}   \mytr_{\gtorus}\upchi \upmu \Speed^2 \exp(2\LogDensity)  |\tander^N \vortrenormalized|^2 \weight^{\blowupratetoporderwave} \ReciprocalLunitAppliedtoTimeFunction \, \voldiamond 
		\end{split}
	\end{align}\end{lemma}

\begin{proof}
We first use \eqref{E:MULTIPLIERVECTORFIELD} to express the Euclidean product on LHS \eqref{E:PERFECTDIVERGENCESTRUCTUREFOROVERDIFFERENTIATEDVORTICITY} as $ (\Lunit \tander^N \vortrenormalized) \cdot (\upmu \tander^N \vortrenormalized) + 2 (\upmu \Transport \tander^N \vortrenormalized)\cdot (\upmu \tander^N \vortrenormalized)$. The integral genenated by the term featuring $\Transport$ is directly present on RHS\,\eqref{E:PERFECTDIVERGENCESTRUCTUREFOROVERDIFFERENTIATEDVORTICITY}.

The desired identity \eqref{E:PERFECTDIVERGENCESTRUCTUREFOROVERDIFFERENTIATEDVORTICITY} then follows from applying \eqref{E:PRELIMINARYINTEGRATIONBYPARTSWITHL} to the integrand featuring $\Lunit$ with $\upupsilon = \frac{1}{2} \upmu \Speed^2\exp(2\LogDensity) |\tander^N \vortrenormalized|^2\weight^{\blowupratetoporderwave}$ and $\zeta = 1$. We omit the straightforward, yet tedious details. 
\end{proof}

\begin{lemma}[Harmless derivative-losing integration by parts identities tied to vorticity] \label{L:HARMLESSOVERDIFFERENTIATEDIBPFORVORTICITY}
For $\Sigma_t$-tangent vectorfields $\SigmatTan_1, \, \SigmatTan_2$, denote their Cartesian inner product by $\SigmatTan_1\cdot\SigmatTan_2 = \sum_{i=1}^3\SigmatTan_1^i\cdot\SigmatTan_2^i$. Then the following identities hold:
	\begin{align}
		\begin{split} \label{E:HARMLESSOVERDIFFERENTIATEDIBPFORVORTICITY}
			&  \int_{\characteristicdiamondtwoarg{[\leftubar,\ubar)}{[\moreinterestingu_1,u]}} \left\{\Speed^2 \exp(\LogDensity) (\multipliervectorfield \tander^N \vortrenormalized)\cdot [\upmu,\tander^N] (\Flatcurl v) + \Speed^2\exp(\LogDensity) (\multipliervectorfield \tander^N \vortrenormalized)\cdot [\tander^N,\upmu \Flatcurl] v \right\}\weight^{\blowupratetoporderwave} \ReciprocalLunitAppliedtoTimeFunction \, \voldiamond \\
			&  \ \ \ \ - \int_{\characteristicdiamondtwoarg{[\leftubar,\ubar)}{[\moreinterestingu_1,u]}} \Speed^2 \exp(\LogDensity) (\multipliervectorfield \tander^N \vortrenormalized^i ) \left\{\sum_{\substack{ \tander^{N_1}\tander^{N_2} = \tander^N \\ N_1 \le N-1}} \upmu (\tander^{N_1} \vortrenormalized^i) \tander^{N_2}(\exp(\LogDensity)) \right\} \weight^{\blowupratetoporderwave} \ReciprocalLunitAppliedtoTimeFunction \, \voldiamond \\
			& =  2 \int_{\characteristicdiamondtwoarg{[\leftubar,\ubar)}{[\moreinterestingu_1,u]}} \left\{\Speed^2\exp(\LogDensity)  (\upmu \Transport \tander^N \vortrenormalized)\cdot [\upmu,\tander^N] (\Flatcurl v) + \Speed^2\exp(\LogDensity) (\upmu \Transport  \tander^N \vortrenormalized)\cdot [\tander^N,\upmu \Flatcurl] v \right\}\weight^{\blowupratetoporderwave} \ReciprocalLunitAppliedtoTimeFunction \, \voldiamond \\
			& \ \ - \int_{\characteristicdiamondtwoarg{[\leftubar,\ubar)}{[\moreinterestingu_1,u]}} \left\{\Speed^2\exp(\LogDensity) (\tander^N \vortrenormalized) \cdot \Lunit \left(  [\upmu,\tander^N] (\Flatcurl v) \right) + \Speed^2\exp(\LogDensity)  (\tander^N \vortrenormalized) \cdot \Lunit \left( [\tander^N,\upmu \Flatcurl ] v\right) \right\} \weight^{\blowupratetoporderwave} \ReciprocalLunitAppliedtoTimeFunction \, \voldiamond \\
			& \ \ - 2 \int_{\characteristicdiamondtwoarg{[\leftubar,\ubar)}{[\moreinterestingu_1,u]}} \Speed^2 \exp(\LogDensity) (\upmu \Transport \tander^N \vortrenormalized^i ) \left\{\sum_{\substack{ \tander^{N_1}\tander^{N_2} = \tander^N \\ N_1 \le N-1}}  \upmu (\tander^{N_1} \vortrenormalized^i) \tander^{N_2}(\exp(\LogDensity)) \right\} \weight^{\blowupratetoporderwave} \ReciprocalLunitAppliedtoTimeFunction \, \voldiamond \\			
			& \ \  + \int_{\characteristicdiamondtwoarg{[\leftubar,\ubar)}{[\moreinterestingu_1,u]}} \Speed^2 \exp(\LogDensity) ( \tander^N \vortrenormalized^i ) \left\{\sum_{\substack{ \tander^{N_1}\tander^{N_2} = \tander^N \\ N_1 \le N-1}} \Lunit \left(\upmu (\tander^{N_1} \vortrenormalized^i) \tander^{N_2}(\exp(\LogDensity))\right) \right\} \weight^{\blowupratetoporderwave} \ReciprocalLunitAppliedtoTimeFunction \, \voldiamond \\
			& \ \ + \int_{\characteristicdiamondtwoarg{[\leftubar,\ubar)}{[\moreinterestingu_1,u]}}
			\left\{ \ErrorIBP_1^{\vortrenormalized;\text{below-top}}[\tander^N] + \ErrorIBP_2^{\vortrenormalized;\text{below-top}}[\tander^N]  \right\}\ReciprocalLunitAppliedtoTimeFunction
			\, \voldiamond \\
			& \ \ + \int_{\ingoingcharacteristicsurfacetwoarg{\ubar}{[\moreinterestingu_1,u]}} \left\{\Speed^2\exp(\LogDensity)  (\tander^N \vortrenormalized) \cdot  [\upmu,\tander^N]\vortrenormalized  + \Speed^2(\tander^N \vortrenormalized) \cdot [\tander^N, \upmu \Flatcurl] v \right\}  \weight^{\blowupratetoporderwave} \, \volingoingnullhypersurface \\
			& \ \ -  \int_{\ingoingcharacteristicsurfacetwoarg{\ubar}{[\moreinterestingu_1,u]}}  \Speed^2 \exp(\LogDensity) ( \tander^N \vortrenormalized^i ) \left\{\sum_{\substack{ \tander^{N_1}\tander^{N_2} = \tander^N \\ N_1 \le N-1}}  \upmu (\tander^{N_1} \vortrenormalized^i) \tander^{N_2}(\exp(\LogDensity)) \right\} \weight^{\blowupratetoporderwave} \, \volingoingnullhypersurface \\
			&  \ \ - \int_{\ingoingcharacteristicsurfacetwoarg{\leftubar}{[\moreinterestingu_1,u]}} \left\{\Speed^2\exp(\LogDensity)  (\tander^N \vortrenormalized) \cdot  [\upmu,\tander^N] \vortrenormalized  + \Speed^2\exp(\LogDensity) (\tander^N \vortrenormalized) \cdot [\tander^N, \upmu \Flatcurl] v \right\}  \weight^{\blowupratetoporderwave} \, \volingoingnullhypersurface \\
			& \ \ +  \int_{\ingoingcharacteristicsurfacetwoarg{\leftubar}{[\moreinterestingu_1,u]}}  \Speed^2 \exp(\LogDensity) ( \tander^N \vortrenormalized^i ) \left\{\sum_{\substack{ \tander^{N_1}\tander^{N_2} = \tander^N \\ N_1 \le N-1}}  \upmu (\tander^{N_1} \vortrenormalized^i) \tander^{N_2}(\exp(\LogDensity)) \right\} \weight^{\blowupratetoporderwave} \, \volingoingnullhypersurface ,
		\end{split}
	\end{align}
where: 
\begin{subequations}
	\begin{align}
		\begin{split} \label{E:ERROR1HARMLESSOVERDIFFERENTIATEDIBPFORVORTICITY}
			\ErrorIBP_1^{\vort;\text{below-top}}[\tander^N] & = \weight^{\blowupratetoporderwave} \Bigg\{  \Speed^2 \left( \Lunit(\exp(\LogDensity))\right) \left(( \tander^N \vortrenormalized)\cdot   [\upmu,\tander^N] (\Flatcurl v) +  (\tander^N \vortrenormalized) \cdot [\tander^N, \upmu \Flatcurl] v\right)  \\
			& \ \ 	-2\Speed (\Speed_{;\LogDensity} \Lunit \LogDensity + \Speed_{;\Ent}\Lunit \Ent) \left(  (\tander^N \vortrenormalized) \cdot  [\upmu,\tander^N] (\Flatcurl v) +  (\tander^N \vortrenormalized) \cdot [\tander^N, \upmu \Flatcurl] v\right)  \\
			& \ \ -  \mytr_{\gtorus}\upchi \left(\Speed^2 \exp(\LogDensity) (\tander^N\vortrenormalized) \cdot  [\upmu,\tander^N] \vortrenormalized  + \Speed^2\exp(\LogDensity)(\tander^N \vortrenormalized) \cdot [\tander^N, \upmu \Flatcurl] v\right) \\
			& \ \ + (\mytr_{\gtorus}\upchi) \Speed^2 \exp(\LogDensity) ( \tander^N \vortrenormalized^i ) \left\{\sum_{\substack{ \tander^{N_1}\tander^{N_2} = \tander^N \\ N_1 \le N-1}}  \upmu (\tander^{N_1} \vortrenormalized^i) \tander^{N_2}(\exp(\LogDensity))\right\} \Bigg\} 
		\end{split} \\
		\begin{split} \label{E:ERROR2HARMLESSOVERDIFFERENTIATEDIBPFORVORTICITY} 
			 \ErrorIBP_2^{\vortrenormalized;\text{below-top}}[\tander^N] & = -\blowupratetoporderwave \weight^{\blowupratetoporderwave -1}  (\Lunit \weight) \Bigg\{ \Speed^2 \exp(\LogDensity) (\tander^N \vortrenormalized) \cdot  [\upmu,\tander^N] \Flatcurl v  + \Speed^2\exp(\LogDensity)(\tander^N \vortrenormalized) \cdot [\tander^N, \upmu \Flatcurl] v  \\
			 & \ \ - \Speed^2 \exp(\LogDensity) ( \tander^N \vortrenormalized^i ) \left\{\sum_{\substack{ \tander^{N_1}\tander^{N_2} = \tander^N \\ N_1 \le N-1}}  \upmu (\tander^{N_1} \vortrenormalized^i) \tander^{N_2}(\exp(\LogDensity))\right\} \Bigg\}
		\end{split}
	\end{align}
\end{subequations}\end{lemma}

\begin{proof}
We first apply \eqref{E:MULTIPLIERVECTORFIELD} to LHS \eqref{E:HARMLESSOVERDIFFERENTIATEDIBPFORVORTICITY} and keep the resulting integrals featuring $\Transport$ unchanged. For the integrals featuring $\Lunit$, we apply \eqref{E:PRELIMINARYINTEGRATIONBYPARTSWITHL} to integrate by parts off the $\Lunit$ from the $\tander^N \vortrenormalized$. We omit the straightforward, yet tedious details.

\end{proof}

\section{Basic ingredients in the $L^2$ analysis} \label{S:BASICINGREDIENTSINL2ANALYSIS}

We continue to work under the assumptions of Sect.\,\ref{SS:SILENTFACTS}.
In this section, we establish some preliminary ingredients that we will use when we derive energy estimates.
In Sect.\,\ref{SS:DIFFERENTIALANDINTEGRALIDENTITIESINVOLVINGDOUBLENULLTORI}, we derive some differential and integral identities involving the double-null tori
$\doublenulltoritwoarg{\ubar}{u}$. 
In Sect.\,\ref{SS:SOBOLEVEMBEDDINGANDFTCESTIMATESONDOUBLENULLTORI}, we establish some basic Sobolev embedding estimates 
and fundamental theorem of calculus-type estimates on the double-null tori. 
In Sect.\,\ref{S:IBPIDENTITIESNEEDFORTOPORDERENERGYESTIMATES}, 
we use the identities from Sect.\,\ref{S:DIVERGENCETHEOREMANDINTEGRATIONBYPARTSIDENTITIES} 
to prove various integration by parts identities 
that will play a key role in the most singular energy estimates for the wave variables. 
Next, in Sect.\,\ref{SS:FUNDAMENTALENERGYIDENTITIES},
we use the vectorfield multiplier method to construct 
the building block $L^2$-based spacetime energies and null-fluxes that we will use to control the wave variables $\wavearray$.
Furthermore, in Prop.\,\ref{P:ENERGYNULLFLUXINTEGRALIDENTITIES}, 
we establish the fundamental energy-null-flux
integral identities that we exploit in our $L^2$ analysis.

In Sect.\,\ref{SS:FUNDAMENTALL2CONTROLLINGWAVEQUANTITIES}, we use the building block energies and null fluxes to 
define the quantities that we will use to control the solution in $\weight$-weighted $L^2$ spaces.
Of particular interest are the spacetime integrals 
 $\spacetimecoercive{q}$ 
and
$\newspacetimecoercive{q}$
defined in
\eqref{E:SPACETIMECOERCIVEINTEGRALS} and \eqref{E:NEWSPACETIMECOERCIVEINTEGRALS}
respectively. These spacetime integrals
appear on the left-hand side of our energy identity \eqref{E:ENERGYNULLFLUXINTEGRALIDENTITIESWAVE},
the former is fundamental for controlling error integrals that involve the quantities $\angrmd \tander^N \Psi$ in a way similar to previous shock formation results \cite{dC2007, gHsKjSwW2016, jLjS2018,abbrescia2022emergence, LukSpeck2024stability}. The latter is a new spacetime coercive term that plays the central role in the weighted energy hierarchies of Sect.\,\ref{S:SATEMENTOFAPRIORIL2ESTIMATES}.
Finally, in Sect.\,\ref{SS:COERCIVENESSFUNDAMENTALL2CONTROLLINGQUANTITIES},
we exhibit the key coerciveness properties of our $L^2$-controlling quantities
with respect to the weighted $L^2$ norms from Sect.\,\ref{SSS:GEOMETRICL2NORMS}.

\subsection{Differential and integral identities involving $\doublenulltoritwoarg{\ubar}{u}$} 

\label{SS:DIFFERENTIALANDINTEGRALIDENTITIESINVOLVINGDOUBLENULLTORI}
The following lemma, though standard, plays an important role in our proof of the 
spacetime energy--null-flux identities for the wave variables
(see Prop.\,\ref{P:ENERGYNULLFLUXINTEGRALIDENTITIES}).

\begin{lemma}[Differential and integral identities involving $\doublenulltoritwoarg{\ubar}{u}$] 
	\label{L:DIFFERENTIALANDINTEGRALIDENTITIESONROUGHTORI}
	Let $f$ be a scalar function,
	and for 
	$Z \in \{\newuL, \newL\}$, let $\mytr_{\gnulltori} \deform{Z}$ be the $\gnulltori$-trace 
	(see Def.\,\ref{D:TRACEOFDOUBLENULLTORITANGENT02TENSORS})
	of the deformation tensor $\deform{Z}$ of $Z$.
	Then the following integral identities hold:
	\begin{subequations}
		\begin{align}
			\nullgeop{\ubar} 
			\left( 
			\int_{\doublenulltoritwoarg{\ubar}{u}} 
			f 
			\, \voldoublenulltori
			\right) 
			& = 
			\ \int_{\doublenulltoritwoarg{\ubar}{u}} 
			\left\lbrace 
			\newL f 
			+ 
			\frac{1}{2} f \mytr_{\gnulltori} \deform{\newL}
			\right\rbrace \, \voldoublenulltori,
			\label{E:IDENTITYUBARDERIVATIVEOFDOUBLENULLORUSINTEGRAL} 
			\\
			\nullgeop{u} \left(  \int_{\doublenulltoritwoarg{\ubar}{u}} f \, \voldoublenulltori\right) & =  \int_{\doublenulltoritwoarg{\ubar}{u}} 
			\left\lbrace 
			\newuL f 
			+ 
			\frac{1}{2} f \mytr_{\gtorusroughfirstfund} \deform{\newuL}
			\right\rbrace \, \voldoublenulltori.
			\label{E:IDENTITYEIKONALDERIVATIVEOFDOUBLENULLTORUSINTEGRAL}
		\end{align}
	\end{subequations}
	
	Moreover, for $u_1 \leq u_2$, we have:
	\begin{align} \label{E:NEWULINTEGRALONDOUBLENULLTORI}
		\int_{\ingoingcharacteristicsurfacetwoarg{\ubar}{[u_1,u_2]}} 
		\left\lbrace
		\newuL f 
		+ 
		\frac{1}{2} f \mytr_{\gnulltori} \deform{\newuL}
		\right\rbrace
		\, \volingoingnullhypersurface
		& =
		\int_{\doublenulltoritwoarg{\ubar}{u_2}} 
		f 
		\, \voldoublenulltori
		-
		\int_{\doublenulltoritwoarg{\ubar}{u_1}} 
		f 
		\, \voldoublenulltori.
	\end{align}
	
\end{lemma}

\begin{proof}
Since $\newL, \newuL$ is $\gfour$-orthogonal to 
$\doublenulltoritwoarg{\ubar}{u}$ and the identities $\newuL \ubar = \newL u = 0$ and $\newuL u = \newL \ubar = 1$ are satisfied, the same proofs of \cite{abbrescia2022emergence}*{Lemmas 6.3, 10.5}.  
hold.

\end{proof}

\subsection{Sobolev embedding and fundamental theorem of calculus-type estimates on the double-null tori} \label{SS:SOBOLEVEMBEDDINGANDFTCESTIMATESONDOUBLENULLTORI} In this section, 
on the double-null tori $\doublenulltoritwoarg{\ubar}{u}$,
we derive $L^{\infty}$ Sobolev embedding estimates as well as simple $L^2$ estimates
that rely on the fundamental theorem of calculus.

\begin{lemma}[Sobolev embedding and fundamental theorem of calculus-type 
estimates on $\doublenulltoritwoarg{\ubar}{u}$] 
\label{L:SOBOLEVEMBEDDINGANDFTCL2ESTIMATESONROUGHTORI}
Let $f$ be a scalar function on $\characteristicdiamondtwoarg{[\leftubar,\ubarboot)}{[\moreinterestingu_1,\moreinterestingu_2]}$. 
Then the following estimates hold for $(\ubar,u) \in [\leftubar,\ubarboot)\times [\moreinterestingu_1,\moreinterestingu_2]$:
\begin{align}
	\| \tander^N f \|_{L^2\left(\doublenulltoritwoarg{\ubar}{u}\right)}^2 
	& 
	\lesssim 
	\| \tander^N f  \|_{L^2\left(\doublenulltoritwoarg{\leftubar}{u}\right)}^2 
	+ 
	\int_{\outgoingcharacteristicsurfacetwoarg{u}{[\leftubar,\ubarboot)}} 
		\frac{1}{\Lunit \ubar} |\Lunit \tander^N f|^2 
	\, \voloutgoingnullhypersurface,
		\label{E:DOUBLENULLTORUSINNULLHYPERSURFACEL2FUNDAMENTALTHEOREMOFCALCULUSESTIMATE} 
\end{align}
\begin{subequations}
\begin{align}
	\| f \|_{L^{\infty}\left(\doublenulltoritwoarg{\ubar}{u}\right)} 
	& \lesssim 
	\| \tander^{\leq 2} f \|_{L^2\left(\doublenulltoritwoarg{\ubar}{u}\right)},
		\label{E:H2LINFINITYSOBOLEVEMBEDDINGDOUBLENULLTORUS} 
			\\
	\| f\|_{L^{\infty}\left(\doublenulltoritwoarg{\ubar}{u}\right)}^2 
	& \lesssim 
	\| 
		\tander^{\leq 2} f  
	\|_{L^2\left(\doublenulltoritwoarg{\leftubar}{u}\right)}^2 
		+ 
		\int_{\outgoingcharacteristicsurfacetwoarg{u}{[\leftubar,\ubarboot]}} \frac{1}{\Lunit \ubar} |\Lunit \tander^{\leq 2} f|^2 \, 		\voloutgoingnullhypersurface.  
	\label{E:H2LINFINITYFUNDAMENTALTHEOREMOFCALCULUSPLUSSOBOLEVEMBEDDINGONDOUBLENULLTORUS}
\end{align}
\end{subequations}
\end{lemma}

\begin{proof}
Using \eqref{E:IDENTITYUBARDERIVATIVEOFDOUBLENULLORUSINTEGRAL} with $f^2$ in place of $f$, the identity $\newL = \frac{1}{\Lunit \ubar} \Lunit$, the pointwise estimate \eqref{E:POINTWISEBOUNDFORDOUBLENULLOROIDALTRACEOFDEFORMATIONTENSOROFNEWLVECTORFIELD},
and Young's inequality, we deduce
$
\left|
\nullgeop{\ubar} 
\| \tander^N f \|_{L^2\left(\doublenulltoritwoarg{\ubar}{u}\right)}^2 
\right|
\leq 
\left\|
	\frac{1}{\sqrt{\Lunit \ubar}} \Lunit \tander^N f 
\right\|_{L^2(\doublenulltoritwoarg{\ubar}{u})}^2 
+ 
C \| \tander^N f\|_{L^2\left(\doublenulltoritwoarg{\ubar}{u}\right)}^2
$. 
Integrating this inequality with respect to $\ubar$, 
applying Gr\"{o}nwall's inequality, 
and also using the identity 
$\int_{\ubar' = \leftubar}^{\ubar} 
	\left\| 
		\frac{1}{\sqrt{\Lunit \ubar'}} \tander^N f
	\right\|_{L^2\left(\doublenulltoritwoarg{\ubar'}{u}\right)}^2 
\, \rmd\ubar' 
= \left\| 
	\frac{1}{\sqrt{\Lunit \ubar}} \Lunit \tander^N f 
	\right \|_{L^2\left(\outgoingcharacteristicsurfacetwoarg{u}{[\leftubar,\ubarboot)}\right)}^2
$, we conclude the inequality \eqref{E:DOUBLENULLTORUSINNULLHYPERSURFACEL2FUNDAMENTALTHEOREMOFCALCULUSESTIMATE}.


We now prove 
\eqref{E:H2LINFINITYSOBOLEVEMBEDDINGDOUBLENULLTORUS}--\eqref{E:H2LINFINITYFUNDAMENTALTHEOREMOFCALCULUSPLUSSOBOLEVEMBEDDINGONDOUBLENULLTORUS}. 
We begin with the following estimate, which holds at fixed $(\ubar,u)$ by virtue of
the standard Sobolev embedding result $H^2(\T^2) \hookrightarrow L^{\infty}(\T^2)$:
\begin{align} \label{E:STANDARDH2LINFINITYSOBOLEVEMBEDDINGONT2}
	\| f \|_{L^{\infty}\left(\doublenulltoritwoarg{\ubar}{u}\right)} 
	& \lesssim 
	\sum_{I + J \leq 2} 
		\left\lbrace 
			\int_{\T^2} 
				\left| 
					\left(\nullgeop{x^2}\right)^I \left(\nullgeop{x^3}\right)^J 
					f(\ubar,u,x^2,x^3)
				\right|^2 
			\, 
			\rmd x^2 \rmd x^3
		\right\rbrace^{1/2}.
\end{align}
Using \eqref{E:NULLCOORDINATEPARTIALXAINTERMSOFGEOMETRICVF},
Lemma\,\ref{L:LINFTYESTIMATESFORINGOINGEIKONALFUNCTIONANDDERIVATIVES} to estimate derivatives of 
$\ubar$ in \eqref{E:NULLCOORDINATEPARTIALXAINTERMSOFGEOMETRICVF}, 
the identities \eqref{E:GEOP2TOCOMMUTATORS}--\eqref{E:GEOP3TOCOMMUTATORS},
Prop.\,\ref{P:SCHEMATICSTRUCTUREOFVARIOUSTENSORSINTERMSOFCONTROLVARS},
the bootstrap assumptions,
and the area form estimate 
\eqref{E:ESTIMATEFORNULLTORIVOLUMEFORMDETERMINANT}, 
we deduce, in view of definitions \eqref{E:AREAFORMDOUBLENULLTORUS} and \eqref{E:DOUBLENULLTORUSINTEGRAL},
that $\mbox{RHS\,\eqref{E:STANDARDH2LINFINITYSOBOLEVEMBEDDINGONT2}}
\lesssim 
\| \tander^{\leq 2} f \|_{L^2\left(\doublenulltoritwoarg{\ubar}{u}\right)}
$.
We have therefore proved \eqref{E:H2LINFINITYSOBOLEVEMBEDDINGDOUBLENULLTORUS}. 
The estimate 
\eqref{E:H2LINFINITYFUNDAMENTALTHEOREMOFCALCULUSPLUSSOBOLEVEMBEDDINGONDOUBLENULLTORUS} then
then follows from 
\eqref{E:H2LINFINITYSOBOLEVEMBEDDINGDOUBLENULLTORUS}
and
\eqref{E:DOUBLENULLTORUSINNULLHYPERSURFACEL2FUNDAMENTALTHEOREMOFCALCULUSESTIMATE} with $N = 0,1,2$.  

\end{proof}

\subsection{Fundamental energy identities} \label{SS:FUNDAMENTALENERGYIDENTITIES}

In this section, 
we derive the fundamental integral identities used to derive estimates.

\subsubsection{Energy-momentum tensor, energy currents, and the multiplier vectorfield for the wave variables}
\label{SSS:ENERGYMOMENTUMTENSORANDCURRENTSFORWAVEEQUATIONS}
To derive energy identities for the wave variables $\wavearray$, which solve the quasilinear wave equations
\eqref{E:COVARIANTWAVEEQUATIONSWAVEVARIABLES},
we rely on the well-known multiplier method, which we now introduce.

Let $f$ be a scalar function (in our applications, $f$ will be some vectorfield derivative of one of the wave variables). 
The energy-momentum tensor associated to $f$ is the following symmetric type $\binom{0}{2}$ tensorfield,
where we recall that $\Dfour$ is the Levi-Civita connection of the spacetime metric $\gfour$:
\begin{align} \label{E:DEFOFENERGYMOMENTUM} 
	\enmomem_{\alpha \beta} 
	&	= 
	\enmomem_{\alpha \beta}[f] 
	\eqdef 
	(\Dfour_\alpha f) 
	\Dfour_\beta f 
	-
	\frac{1}{2} 
	\gfour_{\alpha \beta} 
	(\gfour^{-1})^{\kappa \lambda}
	(\Dfour_{\kappa} f)
	\Dfour_{\lambda} f.
\end{align}Given any scalar function $f$ and any multiplier vectorfield $Z$, 
we define the corresponding energy current vectorfield as follows:
\begin{align}  \label{E:WAVEEQUATIONENERGYCURRENT}
	\Jenarg{Z}{\alpha}[f] 
	& 
	\eqdef \enmomem^{\alpha \beta}[f] Z_{\beta}.
\end{align}
Recall that the deformation tensor of $Z$ is given by \eqref{E:DEFORMATIONTENSORDEF}.

The starting point for our derivation of our $L^2$-type integral identities 
for solutions of covariant wave equations is the following identity \emph{for the weighted current $ \, \weight^{q}\, \Jen{Z}[f]$} for some $M \in [0,\infty)$,
which follows easily from the definitions:
\begin{align} \label{E:DIVERGENCEOFWAVEENERGYCURRENT}
	\Dfour_\alpha \left( \, \weight^{q} \Jenarg{Z}{\alpha}[f]\right) 
	& =  
	\, \weight^{q} (\square_{\gfour(\wavearray)} f) Z f 
	+  
	\frac{1}{2}  \, \weight^{q} \enmomem^{\alpha \beta} \deformarg{Z}{\alpha}{\beta} + q  \, \weight^{q-1} \, \Jenarg{Z}{\alpha}[f]\p_\alpha  \weight.
\end{align}
We highlight already the last term on the RHS\,\eqref{E:DIVERGENCEOFWAVEENERGYCURRENT}: $\Jenarg{Z}{\alpha}[f] \p_\alpha \weight$. This term was not present in our previous work \cite{abbrescia2022emergence} and plays a crucial role in our $L^2$ analysis. We analyze $\Jenarg{Z}{\alpha}[f] \p_\alpha \weight$ in Sect.\,\ref{SSS:NEWBULKTERMPRODUCEDBYTHEWEIGHTS}.

In order to obtain wave equation energy estimates that are sufficient to allow us
to track the solution up to the singular boundary, 
we use the following multiplier vectorfield defined in \eqref{E:MULTIPLIERVECTORFIELD}.

\begin{remark}[Weighted multipliers] \label{R:WEIGHTEDMULTIPLIER} In the applications to come, we will use \eqref{E:DIVERGENCEOFWAVEENERGYCURRENT} with $Z \eqdef \multipliervectorfield$. We note that this is equivalent to using the \emph{unweighted} current $\Jenarg{Z}{\alpha}[f]$, and instead using the \emph{weighted multiplier} $\, \weight^{q} \multipliervectorfield$.
\end{remark}

\subsubsection{Decomposition of components of the energy-momentum tensor}
\label{SSS:DECOMPOSITIONOFENERGYMOMENTUMTENSORCOMPOENTNS}
We start with the following lemma, 
which yields identities for various components of the energy-momentum tensor. 

\begin{lemma}[Decomposition of various components of the energy-momentum tensor]
	\label{L:DECOMPOSITIONOFCOMPOENTNSOFENERGYMOMENTUMTENSOR} 
	Let $f$ be a scalar function, and let $\enmomem$ be the corresponding energy-momentum tensor
	defined in \eqref{E:DEFOFENERGYMOMENTUM}.
	Then the following identities hold:
	\begin{subequations}
		\begin{align}
			\enmomem[f](\Lunit,\Lunit) & = (\Lunit f)^2, 
			\label{E:QLL} 
			\\
			\enmomem[f](X,\Lunit) 
			& = 
			- \frac{1}{2} (\Lunit f)^2 + \frac{1}{2} | \angD f|_{\gtorus}^2, 
			\label{E:QLX} 
			\\
			\enmomem[f](X,X) 
			& = \frac{1}{2} (\Lunit f)^2 + (\Lunit f) Xf + (X f)^2 - \frac{1}{2} |\angD f|_{\gtorus}^2. \label{E:QXX}
		\end{align}
	\end{subequations}
	Moreover, the following identities also hold: 
	\begin{subequations}
		\begin{align}
			\enmomem[f](\multipliervectorfield,\Lunit) 
			& 
			= (1 + \upmu)(\Lunit f)^2 + \upmu |\angD f|^2_{\gtorus}, 
			\label{E:QBREVETL} 
			\\
			\enmomem[f](\multipliervectorfield,\newuL) & = 2  (\newuL f)^2 + 2  \MagnitueofinnerproductofnewLandnewuL \ToriTangentVectorfieldAssociatedToDoubleNullFolliations f (\newuL f) + \frac{1}{2} \left( \upmu + 2 \upmu \ReciprocaluLunitAppliedtoTimeFunction\right) | \nullangrmD f|^2_{\gnulltori}, \label{E:QBREVETNEWUL}
		\end{align}
	\end{subequations}
	where $\ToriTangentVectorfieldAssociatedToDoubleNullFolliations$ is as in \eqref{E:TORITANGENTVECTORFIELDASSOCIATEDTODOUBLENULLFRAME}.  
\end{lemma}
\begin{proof}
	Identities \eqref{E:QLL}--\eqref{E:QBREVETL} were proven in \cite{abbrescia2022emergence}*{Lemma 21.13}; we omit the details, noting only that the proofs are based on Lemmas\,\ref{L:BASICPROPERTIESOFVECTORFIELDS} and\,\ref{L:USEFULIDENTITIESFORFIRSTFUNDAMENTALFORM}.
	
	We now prove \eqref{E:QBREVETL}. Contracting \eqref{E:DEFOFENERGYMOMENTUM} with $\multipliervectorfield^\alpha \newuL^\beta$ and using \eqref{E:ACOUSTICALINVERSEMETRICINTERMSOFNULLVECTORFIELDSANDDOUBLENULLSPHEREINVERSEFIRSTFUND}, we find: 
	\begin{equation} \label{E:ENMOMEMCONTRACTMULTIPLIERANDNEWUL}
		\enmomem[f](\multipliervectorfield,\newuL) = \multipliervectorfield f (\newuL f) - \frac{1}{2} \gfour(\multipliervectorfield, \newuL) \left\{ -\frac{2}{\MagnitueofinnerproductofnewLandnewuL} \newL f (\newuL f) + |\nullangrmD f|^2_{\gnulltori}\right\}.
	\end{equation}
	Next, we claim that the following identities hold:
	\begin{subequations}
		\begin{align}
			\multipliervectorfield f (\newuL f) & = \left( \frac{1}{\ReciprocalLunitAppliedtoTimeFunction} + 2 \frac{\upmu}{\ingoingmu}\right) \newL f(\newuL f) + 2 (\newuL f)^2 + 2  \MagnitueofinnerproductofnewLandnewuL \ToriTangentVectorfieldAssociatedToDoubleNullFolliations f (\newuL f) 	\label{E:PRODUCTOFMULTIPLIERVECTORFIELDANDNEWULACTINGONSCALARS}  \\
			\gfour(\multipliervectorfield,\newuL) & = - \frac{\MagnitueofinnerproductofnewLandnewuL}{\ReciprocalLunitAppliedtoTimeFunction} - 2 \upmu \ReciprocaluLunitAppliedtoTimeFunction. \label{E:GFOURINNERPRODUCTMULTIPLIERVECTORFIELDANDNEWUL}
		\end{align}
	\end{subequations}
	Straightforward algebra using \eqref{E:ENMOMEMCONTRACTMULTIPLIERANDNEWUL}--\eqref{E:GFOURINNERPRODUCTMULTIPLIERVECTORFIELDANDNEWUL} proves \eqref{E:QBREVETNEWUL}. We now prove \eqref{E:PRODUCTOFMULTIPLIERVECTORFIELDANDNEWULACTINGONSCALARS}--\eqref{E:GFOURINNERPRODUCTMULTIPLIERVECTORFIELDANDNEWUL}. The first identity  \eqref{E:PRODUCTOFMULTIPLIERVECTORFIELDANDNEWULACTINGONSCALARS} follows from multiplying \eqref{E:NEWULWITHRESPECTTOLOLDULUNIT} by $\upmu$, solving for $\upmu \Transport$, and inserting the resulting expression into  \eqref{E:MULTIPLIERVECTORFIELD}.

Next, we express $\multipliervectorfield = L + 2 \upmu \Transport = \frac{1}{\ReciprocalLunitAppliedtoTimeFunction} \newL + 2 \upmu \Transport$. Identities \eqref{E:INNERPRODUCTOFEIKONALFUNCTIONNORMALIZEDNULLVECTORFIELDS} and \eqref{E:RELATIONBETWEENCARTESIANNORMALIZEDNULLVECTORFIELDSANDEIKONALFUNCTIONORMALIZEDNULLVECTORFIELDS}--\eqref{E:EQUIVALENTLUNITANDULUNITAPPLIEDTOEIKONALANDCARTESIANTIME} conclude the proof of \eqref{E:GFOURINNERPRODUCTMULTIPLIERVECTORFIELDANDNEWUL}.
\end{proof}

\subsubsection{The new bulk term produced by the weights}
\label{SSS:NEWBULKTERMPRODUCEDBYTHEWEIGHTS}

In this section we analyze the new bulk term produced by incorporating weights with multiplier method, namely, the last term on RHS\,\eqref{E:DIVERGENCEOFWAVEENERGYCURRENT}. The main result is provided by the following lemma.

\begin{lemma}[The new bulk terms produced by the weights] \label{L:THENEWBULKTERMPRODUCEDBYTHEWEIGHTS}
Let $f$ be a scalar function, let $\multipliervectorfield$ be the multiplier vectorfield defined in \eqref{E:MULTIPLIERVECTORFIELD}, and let $\Jenarg{\multipliervectorfield}{\alpha}[f]$ be its corresponding energy current vectorfield. Then the following identity holds:
\begin{equation} \label{E:THENEWBULKTERMPRODUCEDBYTHEWEIGHTS}
	\Jenarg{\multipliervectorfield}{\alpha}[f] \p_\alpha \weight  =  - \frac{1}{\MagnitueofinnerproductofnewLandnewuL} \newL\weight \left\{2  (\newuL f)^2 + 2  \MagnitueofinnerproductofnewLandnewuL \ToriTangentVectorfieldAssociatedToDoubleNullFolliations f (\newuL f) + \frac{1}{2} \left( \upmu + 2 \upmu \ReciprocaluLunitAppliedtoTimeFunction\right) | \nullangrmD f|^2_{\gnulltori}\right\}.
\end{equation}

\end{lemma}

\begin{proof}
Since $\weight$ is proportional to $\ubar$ by \eqref{E:BLOWUPWEIGHT}, it follows from \eqref{E:ACOUSTICALINVERSEMETRICINTERMSOFNULLVECTORFIELDSANDDOUBLENULLSPHEREINVERSEFIRSTFUND} that $\gfour^{\alpha\delta} \p_\alpha \weight = -\frac{1}{\MagnitueofinnerproductofnewLandnewuL} \newuL^{\delta} \newL \weight$. The desired identity \eqref{E:THENEWBULKTERMPRODUCEDBYTHEWEIGHTS} then follows from  \eqref{E:WAVEEQUATIONENERGYCURRENT} with $Z = \multipliervectorfield$ and \eqref{E:QBREVETNEWUL}.

\end{proof}

\subsubsection{Building block spacetime energy--null-fluxes} \label{SSS:BUILDINGBLOCKENERGIESANDNULLFLUXES} We are now ready to define our building block null-fluxes for the wave variables. We also define spacetime integrals that play a crucial role in our energy estimates; see Def.\,\ref{D:SPACETIMECOERCIVEINTEGRALS}. These building blocks will also feature the weights $\, \weight^{q}$ for some $q \in [0,\infty)$, which are crucial for our a priori estimates. We recall that, at the top derivative order, the weights for the wave variables $\wavearray$ will be $\blowupratetoporderwave$.

\begin{definition}[Null-fluxes for the wave and transport variables] \label{D:NULLFLUXESFORTHEWAVEANDTRANSPORTVARIABLES} Let $f$ be a scalar function on $\characteristicdiamondtwoarg{[\leftubar,\ubar)}{[\moreinterestingu_1,u]}$ and let $q \in [0,\infty)$. We respectively define the {wave} ingoing and outgoing $\weight$-weighted null-fluxes associated to $f$ as follows:
	\begin{subequations}
		\begin{align}
			\ingoingfluxwave{q}[f](\ubar,u) & \eqdef \int_{\ingoingcharacteristicsurfacetwoarg{\ubar}{[\moreinterestingu_1,u]}} \enmomem[f](\multipliervectorfield,\newuL) \, \weight^{q}   \volingoingnullhypersurface, \label{E:WAVEINGOINGNULLFLUX} \\
			\outgoingfluxwave{q}[f] (\ubar,u) & \eqdef \int_{\outgoingcharacteristicsurfacetwoarg{u}{[\leftubar,\ubar)}} \enmomem[f](\multipliervectorfield,\Lunit) \, \weight^{q}  \ReciprocalLunitAppliedtoTimeFunction \voloutgoingnullhypersurface. \label{E:WAVEOUTGOINGNULLFLUX}
		\end{align}
	\end{subequations}
	
We respectively define the \emph{transport} ingoing and outgoing $\weight$-weighted null-fluxes associated to $f$ as follows:
	\begin{subequations}
		\begin{align}
			\ingoingfluxtransport{q}[f](\ubar,u) & \eqdef \int_{\ingoingcharacteristicsurfacetwoarg{\ubar}{[\moreinterestingu_1,u]}} f^2  \, \weight^{q} \ReciprocaluLunitAppliedtoTimeFunction \volingoingnullhypersurface, \label{E:TRANSPORTINGOINGNULLFLUX} \\
			\outgoingfluxtransport{q}[f] (\ubar,u) & \eqdef \int_{\outgoingcharacteristicsurfacetwoarg{u}{[\leftubar,\ubar)	}} f^2  \, \weight^{q} \ReciprocalLunitAppliedtoTimeFunction \voloutgoingnullhypersurface. \label{E:TRANSPORTOUTGOINGNULLFLUX}
		\end{align}
	\end{subequations}
	
\end{definition}

In the next definition we introduce \emph{two types} of spacetime integrals that are of crucial importance in our energy estimates. The first type dates back to Christodoulou's monumental work on shock formation \cite{dC2007} and an appropriate version of it is present in our prior work \cite{abbrescia2022emergence}. The second type is a new coercive integral that manifests itself in the \emph{weighted} $L^2$-energy approach.

\begin{definition}[Coercive spacetime integrals] \label{D:SPACETIMECOERCIVEINTEGRALS} 
	Let $\upxi$ be a scalar function or an $\ell_{t,u}$-tangent one form defined on $\characteristicdiamondtwoarg{[\leftubar,\ubar)}{[\moreinterestingu_1,u]}$. Then we define the following spacetime integral:
	\begin{equation}
		\spacetimecoercive{q}[\upxi](\ubar,u) \eqdef - \int_{\characteristicdiamondtwoarg{[\leftubar,\ubar)}{[\moreinterestingu_1,u]}}  ( \newL \upmu) |\upxi|^2 \mathbf{1}_{\characteristicdiamondtwoarg{[\leftubar,\ubar)}{[\moreinterestingu_1,u]} \cap \smallneighborhoodofcreasearg{[\leftubar,\ubarboot]}}
 \, \weight^{q}  \voldiamond, \label{E:SPACETIMECOERCIVEINTEGRALS}
	\end{equation}
	where $\smallneighborhoodofcreasearg{[\leftubar,\ubarboot]}$ is the spacetime region from \eqref{E:SMALLNEIGHBORHOOD}. In the case that $\upxi$ is an $\ell_{t,u}$-tangent one form, $|\xi|^2$ is understood to be $|\xi|_{\gtorus}^2$ in \eqref{E:SPACETIMECOERCIVEINTEGRALS}.

	Next, let $\upxi$ be a scalar function or a $\doublenulltoritwoarg{\ubar}{u}$-tangent one form  defined on $\characteristicdiamondtwoarg{[\leftubar,\ubar)}{[\moreinterestingu_1,u]}$.  Then we define the following spacetime integral:
	\begin{equation}
		\newspacetimecoercive{q}[\upxi](\ubar,u) \eqdef - \int_{\characteristicdiamondtwoarg{[\leftubar,\ubar)}{[\moreinterestingu_1,u]}}   \newL \weight \, |\upxi|^2 \, \weight^{q}  \voldiamond. \label{E:NEWSPACETIMECOERCIVEINTEGRALS}
	\end{equation}
	In the case that $\upxi$ is an $\doublenulltoritwoarg{\ubar}{u}$-tangent one form , $|\upxi|^2$ is understood to be $|\upxi|_{\gnulltori}^2$ in \eqref{E:NEWSPACETIMECOERCIVEINTEGRALS}.
\end{definition}

\begin{remark}[Missing factor of $\frac{1}{2}$ in coercive spacetime integral] The analog of our spacetime integral \eqref{E:SPACETIMECOERCIVEINTEGRALS} in our prior work \cite{abbrescia2022emergence} featured an additional factor of $\frac{1}{2}$. We have dropped this factor and have chosen to instead keep the coefficients in the energy identity, see the third integral in \eqref{E:ENERGYNULLFLUXINTEGRALIDENTITIESWAVE}.
\end{remark}

\subsubsection{The fundamental spacetime energy--null-flux integral identities} 
\label{SSS:ENERGYNULLFLUXINTEGRALIDENTITIES}
In the next proposition, 
we provide the fundamental null-flux integral identities that form the foundation
of our hyperbolic energy estimates for the wave variables. 

\begin{proposition}[Fundamental spacetime energy--null-flux identities] \label{P:ENERGYNULLFLUXINTEGRALIDENTITIES}
	\hfill
	
	\noindent \underline{\textbf{Wave equation spacetime energy--null-flux identity}}.
	Suppose that on $\characteristicdiamondtwoarg{[\leftubar,\ubar)}{[\moreinterestingu_1,u]}$, the scalar function $f$ is the solution to the inhomogeneous covariant wave equation $\upmu \Box_{\gfour} f = \mathfrak{G}$. Then the following identity holds for $q \in [0,\infty)$:
	\begin{align}
		\begin{split} \label{E:ENERGYNULLFLUXINTEGRALIDENTITIESWAVE}
			& \ingoingfluxwave{q}[f](\ubar,u) + \outgoingfluxwave{q}[f](\ubar,u) + \frac{1}{2} \spacetimecoercive{q}[ \angrmD f](\ubar,u)   			+ \frac{q}{2} \newspacetimecoercive{q-1}\left[ \sqrt{\upmu + 2\ReciprocaluLunitAppliedtoTimeFunction \upmu} \, \nullangrmD f\right](\ubar,u) + 2 q \newspacetimecoercive{q-1}\left[ \newuL f\right](\ubar,u)  \\
			& = \ingoingfluxwave{q}[f](\leftubar,u) + \outgoingfluxwave{q}[f](\ubar,\moreinterestingu_1) - \int_{\characteristicdiamondtwoarg{[\leftubar,\ubar)}{[\moreinterestingu_1,u]}} \left\{ \Lunit f + 2 \upmu \Transport f \right\} \mathfrak{G} \, \weight^{q} \ReciprocalLunitAppliedtoTimeFunction \voldiamond  \\
			& \ \
			+ \int_{\characteristicdiamondtwoarg{[\leftubar,\ubar)}{[\moreinterestingu_1,u]}}  
			{^{(\multipliervectorfield)}\mathfrak{B}}[f]\, \weight^{q} \ReciprocalLunitAppliedtoTimeFunction \voldiamond  + 2q \int_{\characteristicdiamondtwoarg{[\leftubar,\ubar)}{[\moreinterestingu_1,u]}}  \MagnitueofinnerproductofnewLandnewuL  (\ToriTangentVectorfieldAssociatedToDoubleNullFolliations f )(\newuL f) (\newL \weight) \, \weight^{q-1} \voldiamond
		\end{split}
	\end{align}
	The error term ${^{(\multipliervectorfield)}\mathfrak{B}}[f]$ appearing 
	in the last integral on RHS\,\eqref{E:ENERGYNULLFLUXINTEGRALIDENTITIESWAVE}
	can be decomposed as follows: 
	\begin{align}  \label{E:WAVEENERGYIDENTITYBULKERRORTERM}
		{^{(\multipliervectorfield)}\mathfrak{B}}[f] 
		& \eqdef \frac{1}{2} \left| \angrmD f\right|^2_{\gtorus} \mathbf{1}_{\characteristicdiamondtwoarg{[\leftubar,\ubar)}{[\moreinterestingu_1,u]} \setminus \smallneighborhoodofcreasearg{[\leftubar,\ubarboot]}}  +
		\sum_{i=1}^6 {^{(\multipliervectorfield)}\mathfrak{B}_{(i)}[f]},
	\end{align}
	where:
	\begin{subequations}
		\begin{align}
			{^{(\multipliervectorfield)}\mathfrak{B}_{(1)}[f]} 
			& 
			\eqdef (\Lunit f)^2 
			\left\lbrace 
			- \frac{1}{2} \Lunit \upmu + \muX \upmu
			- 
			\frac{1}{2} \upmu \mytr_{\gtorus} \upchi 
			- 
			 \upmu^2 \mytr_{\gtorus} \angktan 
			- 
			 \mytr_{\gtorus} \angktrans
			\right\rbrace,  
			\label{E:WAVEENERGYIDENTITYBULKERRORTERM1} 
			\\
			{^{(\multipliervectorfield)}\mathfrak{B}_{(2)}[f]} 
			& 
			\eqdef 
			- (\Lunit f) (\muX f) 
			\left\lbrace  
			\mytr_{\gtorus} \upchi 
			+ 
			2 \upmu \mytr_{\gtorus} \angktan 
			+ 
			2 \mytr_{\gtorus} \angktrans 
			\right\rbrace, 
			\label{E:WAVEENERGYIDENTITYBULKERRORTERM2} 
			\\
			{^{(\multipliervectorfield)}\mathfrak{B}_{(3)}[f]} 
			& \eqdef  
			|\angrmD f|_{\gtorus}^2 
			\left\lbrace \muX \upmu + 2 \upmu  \Lunit \upmu +			+ 
			\frac{1}{2}
			\upmu \mytr_{\gtorus}\upchi 
			+ 
			\upmu^2 \mytr_{\gtorus} \angktan 
			+ 
			\upmu \mytr_{\gtorus} \angktrans 
			\right\rbrace, 
			\label{E:WAVEENERGYIDENTITYBULKERRORTERM3}
			\\
			{^{(\multipliervectorfield)}\mathfrak{B}_{(4)}[f]} 
			& 
			\eqdef 
			(\Lunit f) (\angrmD^{\#} f)
			\cdot
			\left\lbrace 
			(1 - 2  \upmu) \angrmD \upmu 
			+ 
			2 \upmu \zetatan 
			+ 
			2 \zetatrans 
			\right\rbrace, 
			\label{E:WAVEENERGYIDENTITYBULKERRORTERM4} 
			\\
			{^{(\multipliervectorfield)}\mathfrak{B}_{(5)}[f]} 
			& \eqdef 
			-2 
			(\muX f) 
			(\angrmD^{\#} f)
			\cdot 
			\left\lbrace 
			\angrmD \upmu 
			+ 
			2 \upmu \zetatan 
			+ 
			2 \zetatrans 
			\right\rbrace, \label{E:WAVEENERGYIDENTITYBULKERRORTERM5} 
			\\
			{^{(\multipliervectorfield)}\mathfrak{B}_{(6)}[f]} 
			& \eqdef 
			- 
			\upmu \angrmD^{\#} f \otimes \angrmD^{\#} f 
			\cdot 
			\left\lbrace 
			\upchi + 2 \upmu \angktan + 2 \angktrans 
			\right\rbrace, \label{E:WAVEENERGYIDENTITYBULKERRORTERM6}
		\end{align}
	\end{subequations} 
	In \eqref{E:WAVEENERGYIDENTITYBULKERRORTERM1}--\eqref{E:WAVEENERGYIDENTITYBULKERRORTERM6},
	the $\ell_{t,u}$-tangent tensorfields $\upchi$,
	$\angktan$, $\angktrans$, $\zetatan$, and $\zetatrans$
	are as in Lemma\,\ref{L:USEFULIDENTITIESANDDECOMPOSITIONSFORSECONDFUNDAMENTALFORMSANDTORSION}.

	\hfill
	
	\noindent \underline{\textbf{Transport equation spacetime energy--null-flux identity}}.	
	Suppose that on $\characteristicdiamondtwoarg{[\leftubar,\ubar)}{[\moreinterestingu_1,u]}$, the scalar function $f$ is the solution to the inhomogeneous transport equation $\upmu \Transport f = \mathfrak{G}$. Then the following identity holds for $q \in [0,\infty)$:
	\begin{align} 
		\begin{split} \label{E:ENERGYNULLFLUXINTEGRALIDENTITIESTRANSPORT} 
			& \ingoingfluxtransport{q}[f](\ubar,u) + \outgoingfluxtransport{q}[f](\ubar,u) + q\newspacetimecoercive{q-1}[\sqrt{\ReciprocaluLunitAppliedtoTimeFunction} f](\ubar,u)  \\
			&
			= \ingoingfluxtransport{q}[f](\leftubar,u) + \outgoingfluxtransport{q}[f](\ubar,\moreinterestingu_1) 
			\\
			& \ \ +  2 \int_{\characteristicdiamondtwoarg{[\leftubar,\ubar)}{[\moreinterestingu_1,u]}} f \mathfrak{G} \, \, \weight^{q} \ReciprocalLunitAppliedtoTimeFunction \voldiamond    + \int_{\characteristicdiamondtwoarg{[\leftubar,\ubar)}{[\moreinterestingu_1,u]}} \left\{\Lunit \upmu + \upmu  \mytr_{\gtorus} \angk \right\}
			f^2 \, \weight^{q} \ReciprocalLunitAppliedtoTimeFunction \voldiamond . 
		\end{split}
	\end{align}
	\end{proposition}

		\begin{proof}
		\textbf{\underline{Proof of \eqref{E:ENERGYNULLFLUXINTEGRALIDENTITIESWAVE}:}}
		
		Throughout, we often use the abbreviated notation $\mathbf{J} \eqdef \, \Jen{\multipliervectorfield}[f]$ to denote the energy current defined in \eqref{E:WAVEEQUATIONENERGYCURRENT}, where $Z \eqdef \multipliervectorfield$ and $\multipliervectorfield$ is defined in \eqref{E:MULTIPLIERVECTORFIELD}. To proceed, we write $\mathbf{J}$ in terms of the double-null coordinate partial derivative vectorfields: $\mathbf{J} = \mathbf{J}^{\ubar} \nullgeop{\ubar} + \mathbf{J}^u \nullgeop{u} + \mathbf{J}^A \nullgeop{x^A}$. Next, using \eqref{E:RATIOOFNULLGEOSICINNERPRODUCTANDFOLIATIONDENSITY} and \eqref{E:DOUBLENULLCARTESIANTIMENORMALIZEDNULLVECTORFIELDS}--\eqref{E:DOUBLENULLEIKONALFUNCTIONNORMALIZEDNULLVECTORFIELDS}, we find that:
		\begin{subequations}
			\begin{align}
				\mathbf{J}^{\ubar} & = - \newuL_\alpha \Jenarg{\multipliervectorfield}{\alpha}[f] \frac{1}{\MagnitueofinnerproductofnewLandnewuL} = - \enmomem[f](\multipliervectorfield,\newuL) \frac{1}{\MagnitueofinnerproductofnewLandnewuL}, \label{E:UBARCOORDINATEOFWAVEENERGYCURRENT} \\
				\mathbf{J}^{u} & = - \Lunit_\alpha \Jenarg{\multipliervectorfield}{\alpha}[f] \frac{\ReciprocalLunitAppliedtoTimeFunction}{\MagnitueofinnerproductofnewLandnewuL} = - \enmomem[f](\multipliervectorfield,\Lunit) \frac{\ReciprocalLunitAppliedtoTimeFunction}{\MagnitueofinnerproductofnewLandnewuL}. \label{E:UCOORDINATEOFWAVEENERGYCURRENT}
			\end{align}
		\end{subequations}
		We now use the definitions of the null-fluxes \eqref{E:WAVEINGOINGNULLFLUX}--\eqref{E:WAVEOUTGOINGNULLFLUX}, identities \eqref{E:UBARCOORDINATEOFWAVEENERGYCURRENT}--\eqref{E:UCOORDINATEOFWAVEENERGYCURRENT}, the divergence theorem \eqref{E:DIVERGENCEIDENTITYSPACETIMEVECTORFIELD} (applied to $\mathcal{J} = \, \weight^{q} \mathbf{J}$), and \eqref{E:DIVERGENCEOFWAVEENERGYCURRENT} to derive:
		\begin{align}
			\begin{split} \label{E:MAINDIVERGENCEIDENTITYFORWAVEENERGYNULLFLUXESINTERMEDIATESTEP1}
				\ingoingfluxwave{q}[f](\ubar,u) + \outgoingfluxwave{q}[f](\ubar,u) 
				& = \ingoingfluxwave{q}[f](\leftubar,u) + \outgoingfluxwave{q}[f](\ubar,\moreinterestingu_1)  \\
				& \ \ 
				- \int_{\characteristicdiamondtwoarg{[\leftubar,\ubar)}{[\moreinterestingu_1,u]}}  (\Box_{\gfour} f) \multipliervectorfield f \MagnitueofinnerproductofnewLandnewuL  \weight^{q}  \, \voldiamond  \\
				& \ \ 
				-  \int_{\characteristicdiamondtwoarg{[\leftubar,\ubar)}{[\moreinterestingu_1,u]}}
				\frac{1}{2}   \enmomem^{\alpha \beta} \deformarg{\multipliervectorfield}{\alpha}{\beta} \MagnitueofinnerproductofnewLandnewuL  \weight^{q}  \, \voldiamond \\
				& \ \ - q  \int_{\characteristicdiamondtwoarg{[\leftubar,\ubar)}{[\moreinterestingu_1,u]}} \Jenarg{\multipliervectorfield}{\alpha}[f]\p_\alpha\weight
				\MagnitueofinnerproductofnewLandnewuL  \weight^{q-1}\, \voldiamond .
			\end{split}
		\end{align}
		Using $\MagnitueofinnerproductofnewLandnewuL  = \upmu \ReciprocalLunitAppliedtoTimeFunction$ and \eqref{E:MULTIPLIERVECTORFIELD}, we see that the first spacetime integral on RHS\,\eqref{E:MAINDIVERGENCEIDENTITYFORWAVEENERGYNULLFLUXESINTERMEDIATESTEP1} is manifestly present as the first spacetime integral on the RHS of \eqref{E:ENERGYNULLFLUXINTEGRALIDENTITIESWAVE} after substituting $\upmu \Box_{\gfour} f = \mathfrak{G}$ and $\multipliervectorfield = \Lunit + 2\upmu \Transport$.

Next, using 
		definitions\,\eqref{E:DEFOFENERGYMOMENTUM},
		\eqref{E:DEFORMATIONTENSORDEF},
		and \eqref{E:MULTIPLIERVECTORFIELD},
		the identity 
		$(\gfour^{-1})^{\alpha \delta} 
		= 
		- 
		\Lunit^{\alpha} \Lunit^{\delta} 
		- 
		\Lunit^{\alpha} X^{\delta} 
		- X^{\alpha}\Lunit^{\delta} 
		+ 
		(\gtorus^{-1})^{\alpha \delta}
		$
		which follows from \eqref{E:SMOOTHTORUSMETRICINTERMSOFSIGMATMETRICANDX},
		and the analogous identity for $(\gfour^{-1})^{\beta \sigma}$, tedious but straightforward calculations using Lemma\,\ref{L:DECOMPOSITIONOFCOMPOENTNSOFENERGYMOMENTUMTENSOR} imply:
		\begin{align}
			\begin{split} \label{E:MAINDIVERGENCEIDENTITYFORWAVEENERGYNULLFLUXESINTERMEDIATESTEP2}
				-\frac{1}{2} \upmu \enmomem^{\alpha \beta} \deformarg{\multipliervectorfield}{\alpha}{\beta}  & =  (\Lunit f)^2 
				\left\lbrace 
				- \frac{1}{2} \Lunit \upmu 
				+ 
				\muX \upmu 
				- 
				\frac{1}{2} \upmu \mytr_{\gtorus} \upchi 
				-
				\upmu^2 \mytr_{\gtorus} \angktan 
				-
				\upmu \mytr_{\gtorus} \angktrans
				\right\rbrace   \\
				&
				\ \ 
				- (\Lunit f) (\muX f) 
				\left\lbrace 
				\mytr_{\gtorus} \upchi 
				+ 
				2 \upmu \mytr_{\gtorus} \angktan 
				+ 
				2 \mytr_{\gtorus} \angktrans 
				\right\rbrace \\
				& 
				\ \ 
				+
				|\angrmD f|_{\gtorus}^2 
				\left\lbrace 
				\frac{1}{2} L\upmu + \muX \upmu
				+ 
				2 \upmu  \Lunit \upmu 
				+ 
				\frac{1}{2}
				\upmu \mytr_{\gtorus}\upchi 
				+ 
				\upmu^2 \mytr_{\gtorus} \angktan 
				+ 
				\upmu \mytr_{\gtorus} \angktrans 
				\right\rbrace, 
				\\
				&
				\ \ +
				(\Lunit f) (\angrmD^{\#} f)
				\cdot
				\left\lbrace 
				(1 - 2 \upmu) \angrmD \upmu 
				+ 
				2 \upmu \zetatan 
				+ 
				2 \zetatrans 
				\right\rbrace, 
				\\
				& \ \ 
				-2 
				(\muX f) 
				(\angrmD^{\#} f)
				\cdot 
				\left\lbrace 
				\angrmD \upmu 
				+ 
				2 \upmu \zetatan 
				+ 
				2 \zetatrans 
				\right\rbrace \\
				& \ \
				- 
				\upmu \angrmD^{\#} f \otimes \angrmD^{\#} f 
				\cdot 
				\left\lbrace 
				\upchi + 2 \upmu \angktan + 2 \angktrans 
				\right\rbrace. 
			\end{split}
		\end{align}
		
We direct the readers to  \cite{jSgHjLwW2016}*{Lemma~3.3} for more details of the proof of \eqref{E:MAINDIVERGENCEIDENTITYFORWAVEENERGYNULLFLUXESINTERMEDIATESTEP2}. 
		Next, we focus on the integral generated by the $\frac{1}{2} |\angrmD f|_{\gtorus}^2 \Lunit \upmu$ term present in the third line of RHS\,\eqref{E:MAINDIVERGENCEIDENTITYFORWAVEENERGYNULLFLUXESINTERMEDIATESTEP2}. Decomposing the integral using: 
	\begin{align}
		1
		= 
		\mathbf{1}_{\characteristicdiamondtwoarg{[\leftubar,\ubar)}{[\moreinterestingu_1,u]} \cap \smallneighborhoodofcreasearg{[\leftubar,\ubarboot]}} +
		\mathbf{1}_{\characteristicdiamondtwoarg{[\leftubar,\ubar)}{[\moreinterestingu_1,u]} \setminus \smallneighborhoodofcreasearg{[\leftubar,\ubarboot]}} \label{E:MAINDIVERGENCEIDENTITYFORWAVEENERGYNULLFLUXESINTERMEDIATESTEP3}
	\end{align}
and moving the integral generated by the $\frac{1}{2} |\angrmD f|_{\gtorus}^2 \Lunit \upmu \mathbf{1}_{\characteristicdiamondtwoarg{[\leftubar,\ubar)}{[\moreinterestingu_1,u]} \cap \smallneighborhoodofcreasearg{[\leftubar,\ubarboot]}}$ term over to the LHS of \eqref{E:MAINDIVERGENCEIDENTITYFORWAVEENERGYNULLFLUXESINTERMEDIATESTEP1} produces the $\frac{1}{2}\spacetimecoercive{q}[ \angrmD f](\ubar,u)$ term present on LHS\,\eqref{E:ENERGYNULLFLUXINTEGRALIDENTITIESWAVE}, see \eqref{E:SPACETIMECOERCIVEINTEGRALS}. We clarify that the LHS \eqref{E:MAINDIVERGENCEIDENTITYFORWAVEENERGYNULLFLUXESINTERMEDIATESTEP2} features an $\upmu$ while the penultimate integral on RHS\,\eqref{E:MAINDIVERGENCEIDENTITYFORWAVEENERGYNULLFLUXESINTERMEDIATESTEP1} features a $\MagnitueofinnerproductofnewLandnewuL = \upmu\ReciprocalLunitAppliedtoTimeFunction$. The remaining factor of $\ReciprocalLunitAppliedtoTimeFunction$ is multiplied by $\Lunit \upmu$ to produce $\newL \upmu$ present in \eqref{E:SPACETIMECOERCIVEINTEGRALS}.
The remaining integral generated by $\frac{1}{2} |\angrmD f|_{\gtorus}^2 \Lunit \upmu \mathbf{1}_{\characteristicdiamondtwoarg{[\leftubar,\ubar)}{[\moreinterestingu_1,u]} \cap \smallneighborhoodofcreasearg{[\leftubar,\ubarboot]}}$ corresponds to the first term present on the RHS\,\eqref{E:WAVEENERGYIDENTITYBULKERRORTERM}.
	
	The remaining terms in \eqref{E:MAINDIVERGENCEIDENTITYFORWAVEENERGYNULLFLUXESINTERMEDIATESTEP2} are present in \eqref{E:WAVEENERGYIDENTITYBULKERRORTERM}--\eqref{E:WAVEENERGYIDENTITYBULKERRORTERM6}. We remark that in the last two products
$
-
\upmu^2 \mytr_{\gtorus} \angktan 
- 
\upmu \mytr_{\gtorus} \angktrans
$  
on RHS\,\eqref{E:WAVEENERGYIDENTITYBULKERRORTERM1},
we have added a factor of $\upmu$ that was mistakenly omitted from
\cite{jSgHjLwW2016}*{Equation~(3.14a)}; this factor will be negligible in the context of our estimates.

We now analyze the last integral present in the RHS\,\eqref{E:MAINDIVERGENCEIDENTITYFORWAVEENERGYNULLFLUXESINTERMEDIATESTEP1}. Using \eqref{E:THENEWBULKTERMPRODUCEDBYTHEWEIGHTS}, we see that this integral is equal to: 
\begin{equation} \label{E:MAINDIVERGENCEIDENTITYFORWAVEENERGYNULLFLUXESINTERMEDIATESTEP4}
	- \frac{q}{2} \newspacetimecoercive{q-1}[ \sqrt{1 + 2 \upmu\ReciprocaluLunitAppliedtoTimeFunction} \nullangrmD f](\ubar,u) - 2 q \newspacetimecoercive{q-1}[\newuL f](\ubar,u) + 2q \int_{\characteristicdiamondtwoarg{[\leftubar,\ubar)}{[\moreinterestingu_1,u]} }\MagnitueofinnerproductofnewLandnewuL \ToriTangentVectorfieldAssociatedToDoubleNullFolliations f (\newuL f) (\newL \weight) \, \weight^{q-1} \voldiamond,
\end{equation} 
see \eqref{E:NEWSPACETIMECOERCIVEINTEGRALS}.  
Combining \eqref{E:MAINDIVERGENCEIDENTITYFORWAVEENERGYNULLFLUXESINTERMEDIATESTEP1}--\eqref{E:MAINDIVERGENCEIDENTITYFORWAVEENERGYNULLFLUXESINTERMEDIATESTEP4} concludes the proof of \eqref{E:ENERGYNULLFLUXINTEGRALIDENTITIESWAVE}.

\noindent \textbf{\underline{Proof of \eqref{E:ENERGYNULLFLUXINTEGRALIDENTITIESTRANSPORT}:}} 

Consider the vectorfield $\mathbf{J} = f^2 \Transport$. Using \eqref{E:LFOLIATIONDENSITY}--\eqref{E:ULFOLIATIONDENSITY} and \eqref{E:EQUIVALENTLUNITANDULUNITAPPLIEDTOEIKONALANDCARTESIANTIME}, it follows that   and the expansion of $\mathbf{J}$ in terms of the double-null coordinate derivative vectorfields is $\mathbf{J} = \frac{1}{\ingoingmu} f^2 \nullgeop{\ubar} + \frac{1}{\upmu}f^2 \nullgeop{u} + \mathbf{J}^A \nullgeop{x^A}$. Next, using $\MagnitueofinnerproductofnewLandnewuL = \upmu \ReciprocalLunitAppliedtoTimeFunction = \ingoingmu \ReciprocaluLunitAppliedtoTimeFunction$, the divergence identity \eqref{E:DIVERGENCEIDENTITYSPACETIMEVECTORFIELD} applied to the vectorfield $\mathcal{J} \eqdef \, \weight^{q} \mathbf{J}$, and the definitions of the null-fluxes \eqref{E:TRANSPORTINGOINGNULLFLUX}--\eqref{E:TRANSPORTOUTGOINGNULLFLUX}, we see:
		\begin{align}
			\begin{split} \label{E:MAINDIVERGENCEIDENTITYFORTRANSPORTENERGYNULLFLUXESINTERMEDIATESTEP1}
				& \ingoingfluxtransport{q}[f](\ubar,u) + \outgoingfluxtransport{q}[f](\ubar,u)  \\
				& \ \  = \ingoingfluxtransport{q}[f](\leftubar,u) + \outgoingfluxtransport{q}[f](\ubar,\moreinterestingu_1)  \\
				& \ \ 
				+ 2 \int_{\characteristicdiamondtwoarg{[\leftubar,\ubar)}{[\moreinterestingu_1,u]}} \, \weight^{q} \mathfrak{G} f  \ReciprocalLunitAppliedtoTimeFunction \, \voldiamond + \int_{\characteristicdiamondtwoarg{[\leftubar,\ubar)}{[\moreinterestingu_1,u]}} \, \weight^{q} f^2 \upmu \Dfour_\alpha \Transport^\alpha \ReciprocalLunitAppliedtoTimeFunction \, \voldiamond  \\
				& \ \
				+ q \int_{\characteristicdiamondtwoarg{[\leftubar,\ubar)}{[\moreinterestingu_1,u]}} \, \weight^{q-1} (\MagnitueofinnerproductofnewLandnewuL \Transport \weight)f^2  \, \voldiamond
			\end{split}
		\end{align}
		Using the covariant divergence identity $\upmu \Dfour_\alpha \Transport^\alpha = \Lunit \upmu + \upmu\mytr_{\gtorus} \angk$ (which follows from \eqref{E:COVARIANTDIVERGENCEOFSPACETIMEVECTORFIELDINTERMSOFRESCALEDFRAME}), the desired identity \eqref{E:ENERGYNULLFLUXINTEGRALIDENTITIESTRANSPORT} follows from the identity $\MagnitueofinnerproductofnewLandnewuL \Transport \weight = - \frac{\MagnitueofinnerproductofnewLandnewuL}{\ingoingmu} = -  \ReciprocaluLunitAppliedtoTimeFunction$ in \eqref{E:MAINDIVERGENCEIDENTITYFORTRANSPORTENERGYNULLFLUXESINTERMEDIATESTEP1}, bringing over the integral generated by the $\ReciprocaluLunitAppliedtoTimeFunction$ terms to the LHS of the equality, resulting in the $\newspacetimecoercive{q-1}[\sqrt{\ReciprocaluLunitAppliedtoTimeFunction} f]$ term on the LHS \eqref{E:ENERGYNULLFLUXINTEGRALIDENTITIESTRANSPORT}.

\end{proof}

\subsection{Energy identities for the top-order acoustic geometry variables} \label{SS:ENERGYIDENTITIESFORTOPORDERACOUSTICQUANTITIES}

In this section, we derive spacetime energy--null flux identities that are the used to obtain $L^2$-control of the acoustic geometry at the top derivative level. These include the fully modified and partially modified quantities introduced in Sect.\,\ref{S:CONSTRUCTIONOFMODIFIEDQUANTITIES}, as well as the top-order term $\tander^{\Ntop} \mytr_{\gtorus} \upchi$, whose energy is weaker by a power of $\weight$ than that of the fully modified quantities. These will provide $L^2$-control of the acoustic geometry at the top derivative level, and are necessary in view of our weighted energy estimate approach. 

\begin{definition}[The function $\sumofbrevexvisquared$] \label{D:FACTORINTOPACOUSTICENERGYNEEDEDFORSHARPCONSTANTS}
We define $\sumofbrevexvisquared$ to be the scalar function:
	\begin{align}
		\sumofbrevexvisquared \eqdef \sqrt{(\muX v^1)^2 + (\muX v^2)^2 + (\muX v^3)^2}. \label{E:FACTORINTOPACOUSTICENERGYNEEDEDFORSHARPCONSTANTS}
	\end{align}
	
\end{definition}

\begin{proposition}[Spacetime energy--null flux identities for the top-order characteristic geometry] \label{P:ENERGYIDENTITIESFORTOPORDERACOUSTICQUANTITIES}
	Let $N = \Ntop$, let $\| \cdot \|_{L^2_{q}\left(\ingoingcharacteristicsurfacetwoarg{\ubar}{[\moreinterestingu_1,\moreinterestingu_2]}\right)}$ denote the $\, \weight^{q}$-weighted $L^2$-norm on the ingoing null hypersurfaces of \eqref{E:GEOMETRICL2NORMSINGOINGNULLHYPERSURFACESANDSPACETIMEREGIONS}, let $\spacetimecoercive{q}, \, \newspacetimecoercive{q}$ denote the coercive spacetime integrals of \eqref{E:SPACETIMECOERCIVEINTEGRALS}--\eqref{E:SPACETIMECOERCIVEINTEGRALS}, and let $\sumofbrevexvisquared$ denote the function from \eqref{E:FACTORINTOPACOUSTICENERGYNEEDEDFORSHARPCONSTANTS}. 
	
	\noindent \textbf{\underline{Identities for $\fullymodquant{\tander^N}$:}} Let $\fullymodquant{\tander^N}$ be the fully modified quantity of Def.\,\ref{D:FULLYANDPARTIALLYMODIFIEDQUANTITIES}. Let
$\tander^N \in \mathfrak{P}^{(N)}$ be the set of order $N$ $\nullhyparg{u}$-tangential commutator operators from
Def.\,\ref{D:STRINGSOFCOMMUTATIONVECTORFIELDS}. Then for any $(\ubar,u) \in [\leftubar,\ubarboot)\times [\moreinterestingu_1,\moreinterestingu_2]$, the following identity holds: 
	\begin{align}
		\begin{split} \label{E:SPACETIMEENERGYNULLFLUXIDENTITYFORFULLYMODIFIEDQUANTITY}
			& \left\| \ReciprocalLunitAppliedtoTimeFunction \sumofbrevexvisquared\, \fullymodquant{\tander^N}  \right\|^2_{L^2_{\blowupratetoporderwave}\left( \ingoingcharacteristicsurfacetwoarg{\ubar}{[\moreinterestingu_1,u]}\right)} + \blowupratetoporderwave \newspacetimecoercive{\blowupratetoporderwave-1}\left[\ReciprocalLunitAppliedtoTimeFunction \sumofbrevexvisquared \fullymodquant{\tander^N}  \right](\ubar,u) \\
			& =  \left\| \ReciprocalLunitAppliedtoTimeFunction \sumofbrevexvisquared\, \fullymodquant{\tander^N}  \right\|^2_{L^2_{\blowupratetoporderwave}\left( \ingoingcharacteristicsurfacetwoarg{\leftubar}{[\moreinterestingu_1,u]}\right)} 
			+ \int_{ \characteristicdiamondtwoarg{[\leftubar,\ubar)}{[\moreinterestingu_1,u]}} \Big\{ 2 \ReciprocalLunitAppliedtoTimeFunction^2 \sumofbrevexvisquared^2 \,  \fullymodquant{\tander^N}  \Lunit \left( \fullymodquant{\tander^N}\right)  \\
			& \ \ 
			+ \left. 2 \ReciprocalLunitAppliedtoTimeFunction^2 \sumofbrevexvisquared \left( \Lunit \sumofbrevexvisquared\right) \left( \fullymodquant{\tander^N}\right)^2  + 2 \ReciprocalLunitAppliedtoTimeFunction (\Lunit \ReciprocalLunitAppliedtoTimeFunction)  \sumofbrevexvisquared^2 \left( \fullymodquant{\tander^N}\right)^2  +  \ReciprocalLunitAppliedtoTimeFunction^2 \sumofbrevexvisquared^2 \left( \fullymodquant{\tander^N}\right)^2    \mytr_{\gtorus} \upchi \right\} \, \weight^{\blowupratetoporderwave} \ReciprocalLunitAppliedtoTimeFunction \, \voldiamond.
		\end{split}
	\end{align}
	
	 Let $\tander^N \in \mathfrak{P}^{(N)}$ be the set of order $N$ $\nullhyparg{u}$-tangential commutator operators from
Def.\,\ref{D:STRINGSOFCOMMUTATIONVECTORFIELDS}. Then, the following identity holds:
	\begin{align}
		\begin{split} \label{E:SPACETIMEENERGYNULLFLUXIDENTITYFORIMPRECISETOPORDERCHI}
			& \left\|  \fullymodquant{\tander^N} \right\|^2_{L^2_{\blowupratetoporderwave+1}\left( \ingoingcharacteristicsurfacetwoarg{\ubar}{[\moreinterestingu_1,u]}\right)} + (\blowupratetoporderwave+1)\newspacetimecoercive{\blowupratetoporderwave }\left[ \fullymodquant{\tander^N} \right](\ubar,u)  \\
			& =  \left\| \fullymodquant{\tander^N} \right\|^2_{L^2_{\blowupratetoporderwave+1}\left( \ingoingcharacteristicsurfacetwoarg{\leftubar}{[\moreinterestingu_1,u]}\right)}+ \int_{ \characteristicdiamondtwoarg{[\leftubar,\ubar)}{[\moreinterestingu_1,u]}} \left\{  2 \,  \fullymodquant{\tander^N}  \Lunit \left( \fullymodquant{\tander^N}\right) + \left( \fullymodquant{\tander^N}\right)^2   \mytr_{\gtorus} \upchi  \right\} \weight^{\blowupratetoporderwave+1} \ReciprocalLunitAppliedtoTimeFunction \, \voldiamond.
		\end{split}
	\end{align} 

	\noindent \textbf{\underline{Identity for $\partialmodquant{\tanderY^{N-1}}$:}} Let $\partialmodquant{\tanderY^{N-1}}$ be the partially modified quantity of Def.\,\ref{D:FULLYANDPARTIALLYMODIFIEDQUANTITIES}. Then for any $(\ubar,u) \in [\leftubar,\ubarboot)\times [\moreinterestingu_1,\moreinterestingu_2]$, the following identity holds: 
	\begin{align}
		\begin{split} \label{E:SPACETIMEENERGYNULLFLUXIDENTITYFORPARTIALLYMODIFIEDQUANTITY}
			& \left\| \sumofbrevexvisquared \, \partialmodquant{\tanderY^{N-1}} \right\|^2_{L^2_{\blowupratetoporderwave-1}\left( \ingoingcharacteristicsurfacetwoarg{\ubar}{[\moreinterestingu_1,u]}\right)} +(\blowupratetoporderwave-1)\newspacetimecoercive{\blowupratetoporderwave-2}\left[\sumofbrevexvisquared \, \partialmodquant{\tanderY^{N-1}}\right](\ubar,u) \\
			&  =   \left\|  \sumofbrevexvisquared \, \partialmodquant{\tanderY^{N-1}} \right\|^2_{L^2_{\blowupratetoporderwave-1}\left( \ingoingcharacteristicsurfacetwoarg{\leftubar}{[\moreinterestingu_1,u]}\right)} 
			+ \int_{ \characteristicdiamondtwoarg{[\leftubar,\ubar)}{[\moreinterestingu_1,u]}} \Big\{ 2  \sumofbrevexvisquared^2  \,   \partialmodquant{\tanderY^{N-1}} \Lunit \left( \partialmodquant{\tanderY^{N-1}}\right)   \\
			& \ \ 
			+ \left. 2   \sumofbrevexvisquared \left( \Lunit \sumofbrevexvisquared\right) \left(  \partialmodquant{\tanderY^{N-1}}  \right)^2  +  \sumofbrevexvisquared^2 \left( \partialmodquant{\tanderY^{N-1}} \right)^2   \mytr_{\gtorus} \upchi \right\} \, \weight^{\blowupratetoporderwave-1} \ReciprocalLunitAppliedtoTimeFunction \, \voldiamond.
		\end{split}
	\end{align}		
	%
	

\end{proposition}

\begin{proof}
	We consider the divergence theorem identity \eqref{E:DIVERGENCEIDENTITYSPACETIMEVECTORFIELD} with $\mathcal{J} = \frac{1}{\upmu}\weight^{\blowupratetoporderwave} \ReciprocalLunitAppliedtoTimeFunction^2\sumofbrevexvisquared^2 \left( \fullymodquant{\tanderY^N}\right)^2 \ReciprocalLunitAppliedtoTimeFunction^2 L$.  Using $\MagnitueofinnerproductofnewLandnewuL = \upmu \ReciprocalLunitAppliedtoTimeFunction$ and $\Lunit \ubar = \frac{1}{\ReciprocalLunitAppliedtoTimeFunction}$, it follows that: 
	\begin{equation} \label{E:SPACETIMEENERGYFORFULLYMODIFIEDQUANTITYINTERMEDIATESTEP1}
		\mathcal{J}^{\ubar} \MagnitueofinnerproductofnewLandnewuL = \weight^{\blowupratetoporderwave} \sumofbrevexvisquared^2  \left( \fullymodquant{\tanderY^N}\right)^2 \ReciprocalLunitAppliedtoTimeFunction^2.
	\end{equation}
	The Leibniz rule and $\newL =  \ReciprocalLunitAppliedtoTimeFunction \Lunit$ imply:
	\begin{align}
		\begin{split} \label{E:SPACETIMEENERGYFORFULLYMODIFIEDQUANTITYINTERMEDIATESTEP2}
			\Dfour_\alpha \mathcal{J}^\alpha \MagnitueofinnerproductofnewLandnewuL & =  \, \weight^{\blowupratetoporderwave}\ReciprocalLunitAppliedtoTimeFunction \left\{ 2 \sumofbrevexvisquared^2 \Lunit \left( \fullymodquant{\tanderY^N}\right) \,  \fullymodquant{\tander^N}  
			+ 2 \sumofbrevexvisquared \left( \Lunit \sumofbrevexvisquared\right) \left( \fullymodquant{\tanderY^N}\right)^2 
			+ 2  \sumofbrevexvisquared^2 \left( \fullymodquant{\tanderY^N}\right)^2 \ReciprocalLunitAppliedtoTimeFunction (\Lunit \ReciprocalLunitAppliedtoTimeFunction) \right. \\
			& \ \ \left.  - \frac{\Lunit \upmu}{\upmu}\sumofbrevexvisquared^2  \left( \fullymodquant{\tanderY^N}\right)^2 \right\} 
			+ \blowupratetoporderwave \weight^{\blowupratetoporderwave - 1} \newL \weight  \sumofbrevexvisquared^2 \left( \fullymodquant{\tanderY^N}\right)^2 \ReciprocalLunitAppliedtoTimeFunction^2 + \MagnitueofinnerproductofnewLandnewuL \weight^{\blowupratetoporderwave}   \sumofbrevexvisquared^2 \left( \fullymodquant{\tanderY^N}\right)^2 \ReciprocalLunitAppliedtoTimeFunction^2 \left\{  \Dfour_\alpha \Lunit^\alpha \right\}
		\end{split}
	\end{align}
	Next we note $ \Dfour_\alpha \Lunit^\alpha = \frac{ \Lunit \upmu}{\upmu} + \mytr_{\gtorus} \upchi$, which follows from \eqref{E:COVARIANTDIVERGENCEOFSPACETIMEVECTORFIELDINTERMSOFRESCALEDFRAME}. The desired identity \eqref{E:SPACETIMEENERGYNULLFLUXIDENTITYFORFULLYMODIFIEDQUANTITY} follows from inserting \eqref{E:SPACETIMEENERGYFORFULLYMODIFIEDQUANTITYINTERMEDIATESTEP1}--\eqref{E:SPACETIMEENERGYFORFULLYMODIFIEDQUANTITYINTERMEDIATESTEP2} into \eqref{E:DIVERGENCEIDENTITYSPACETIMEVECTORFIELD}, noting that the $\frac{\Lunit \upmu}{\upmu}$ term on the second line of \eqref{E:SPACETIMEENERGYFORFULLYMODIFIEDQUANTITYINTERMEDIATESTEP2} \emph{exactly cancels} the resulting $\frac{\Lunit \upmu}{\upmu}$ term from $\Dfour_\alpha \Lunit^\alpha$, and bringing the integral generated by the $\newL \weight$ term on RHS\,\eqref{E:SPACETIMEENERGYFORFULLYMODIFIEDQUANTITYINTERMEDIATESTEP2} over to the LHS (see \eqref{E:NEWSPACETIMECOERCIVEINTEGRALS}).
	
	Similarly, \eqref{E:SPACETIMEENERGYNULLFLUXIDENTITYFORIMPRECISETOPORDERCHI} follows from similar arguments applied to $\mathcal{J} = \frac{1}{\upmu}  \weight^{\blowupratetoporderwave+1} \left( \fullymodquant{\tander^N}\right)^2 \Lunit$.  We omit the details.
	
	The analogous identity \eqref{E:SPACETIMEENERGYNULLFLUXIDENTITYFORPARTIALLYMODIFIEDQUANTITY} for $\partialmodquant{\tanderY^{N-1}}$ follows from similar arguments applied to the spacetime vectorfield $\mathcal{J} = \frac{1}{\upmu} \weight^{\blowupratetoporderwave-1}  \sumofbrevexvisquared^2 \left( \partialmodquant{\tanderY^{N-1}}\right)^2  \Lunit$.

\end{proof}

\subsection{Energy identities for below-top-order characteristic geometry variables} \label{SS:ENERGYIDENTITIESFORBELOWTOPORDERACOUSTICGEOMETRYVARIABLES}
In this section, we derive analogous energy identities from Sect.\,\ref{SS:ENERGYIDENTITIESFORTOPORDERACOUSTICQUANTITIES} for below-top-order derivatives of the eikonal function quantities $\upmu, \Lunit^i, \upchi$, and $\mytr_{\gtorus} \upchi$. The $L^2$-estimates derived from these identities \emph{involve a loss of one order of differentiability} relative to $\Psi$, see Sect.\,\ref{SS:PRELIMINARYL2ESTIMATESFORBELOWTOPORDERDERIVATIVESOFACOUSTICGEOMETRY}. The same identity can be used to derive a derivative-losing estimate for $P^N\Psi$, see Sect.\,\ref{SS:WAVEVARIABLEL2ESTIAMTESTHATLOSEONEDERIVATIVE}.

\begin{proposition}[Spacetime energy--null flux identities for the below-top-order characteristic geometry] \label{P:ENERGYIDENTITIESFOREBLOWTOPORDERACOUSTICQUANTITIES}
	Let $N \le \Ntop$, let $\| \cdot \|_{L^2_{q}\left(\ingoingcharacteristicsurfacetwoarg{\ubar}{[\moreinterestingu_1,\moreinterestingu_2]}\right)}$ denote the $\, \weight^{q}$-weighted $L^2$-norm on the ingoing null hypersurfaces of \eqref{E:GEOMETRICL2NORMSINGOINGNULLHYPERSURFACESANDSPACETIMEREGIONS}, and let $\spacetimecoercive{q}$ denote the coercive spacetime integral of \eqref{E:NEWSPACETIMECOERCIVEINTEGRALS}. For 
	\[\varphi \in \{\tander^N\Psi, \, \tandersmall^N \upmu, \, \tander^N \Lunit^i, \, \tander^{N-1} \upchi, \, |\angLie_{\tander}^{N-1} \upchi|_{\gtorus}, \, \comdersmall^{N;1} \Lunit^i, \, \comder^{N-1;1} \mytr_{\gtorus} \upchi,\, |\angLie_{\comder}^{N-1;1} \upchi|_{\gtorus}\},\] the following spacetime energy--null-flux identity holds: 
	\begin{align}
		\begin{split} \label{E:SPACETIMEENERGYNULLFLUXIDENTITYFORBELOWTOPORDERACOUSTIC}
			& \left\|  \varphi \right\|^2_{L^2_{q}\left( \ingoingcharacteristicsurfacetwoarg{\ubar}{[\moreinterestingu_1,u]}\right)} +  q\newspacetimecoercive{q-1}[\varphi](\ubar,u) \\
			& =  \left\|  \varphi \right\|^2_{L^2_{q}\left( \ingoingcharacteristicsurfacetwoarg{\leftubar}{[\moreinterestingu_1,u]}\right)} + \int_{ \characteristicdiamondtwoarg{[\leftubar,\ubar)}{[\moreinterestingu_1,u]}} \left\{ 2 \varphi \Lunit \varphi  + \varphi^2   \mytr_{\gtorus} \upchi \right\} \, \weight^{q} \ReciprocalLunitAppliedtoTimeFunction \, \voldiamond.		\end{split}
	\end{align}
	
\end{proposition}
\begin{proof}
	The result follows from similar arguments used to prove \eqref{E:SPACETIMEENERGYNULLFLUXIDENTITYFORFULLYMODIFIEDQUANTITY}--\eqref{E:SPACETIMEENERGYNULLFLUXIDENTITYFORPARTIALLYMODIFIEDQUANTITY} with $\mathcal{J} = \, \weight^{q-1} \varphi^2 L$. We omit the details.
\end{proof}

\subsection{The fundamental $L^2$-velocity-controlling quantities}
\label{SS:FUNDAMENTALL2CONTROLLINGWAVEQUANTITIES}
In this section, we use the building block quantities
from Sect.\,\ref{SSS:BUILDINGBLOCKENERGIESANDNULLFLUXES} to 
construct the ``fundamental $L^2$-velocity-controlling quantities'' that we use to control
$\velocityarray$, and their derivatives in $L^2$ in various regions.

\begin{definition}[The fundamental $L^2$-velocity-controlling quantities] 
	\label{D:MAINCOERCIVE} 
	Let $\blowuprateofwaveWRTnumberofcommutations(N)$ denote the powers of $\weight$ from Def.\,\ref{D:NOTATIONFORPOWERSOFWEIGHT}.  Let $1 \le N \le \Ntop$ be an integer. Then in terms of the null-flux quantities of Def.\,\ref{D:NULLFLUXESFORTHEWAVEANDTRANSPORTVARIABLES} 
	and the differentiated conventions of Def.\,\ref{D:STRINGSOFCOMMUTATIONVECTORFIELDS}, 
	we define the following $L^2$-controlling quantities:
	\begin{subequations}
			\begin{align}  
				\begin{split}\label{E:UNIFIEDCOMMUTATORVELOCITYNULLFLUXL2CONTROLLINGQUANTITY} 
					\totalfluxcontrolvelocity_N(\ubar,u) 
					& 
					\eqdef 
					\max_{\tander^N \in \mathcal{P}^{(N)}}
					\sup_{\substack{\ubar' \in [\leftubar,\ubar] 
							\\
							u' \in [\moreinterestingu_1,u]}} 					
					\sum_{i=1}^3			
					\left\lbrace 
					 \ingoingfluxwave{\blowuprateofwaveWRTnumberofcommutations(N)}\left[\tander^N v^i\right](\ubar',u') + 
					\outgoingfluxwave{\blowuprateofwaveWRTnumberofcommutations(N)}\left[\tander^N v^i\right](\ubar',u')
					\right\rbrace,  
				\end{split}
				\\
				\strongangularcontrolvelocity_N(\ubar,u) 
				& 
				\eqdef 
					\max_{\tander^N \in \mathcal{P}^{(N)}} \sum_{i=1}^3
						\spacetimecoercive{\blowuprateofwaveWRTnumberofcommutations(N)}\left[\angrmD \tander^N v^i\right] (\ubar,u), 
				\label{E:UNIFIEDCOMMUTATORSPACETIMESTRONGANGULARVELOCITYCONTROL}  \\			
				\weakspacetimeangularvelocitycontrol_N(\ubar,u) & \eqdef 
					\begin{cases}
						\max_{\tander^N \in \mathcal{P}^{(N)}} \sum_{i=1}^3
						\newspacetimecoercive{\blowuprateofwaveWRTnumberofcommutations(N) - 1}\left[\sqrt{\upmu + 2 \ReciprocaluLunitAppliedtoTimeFunction \upmu} \nullangrmD \tander^N v^i\right] (\ubar,u),  & \blowuprateofwaveWRTnumberofcommutations(N) > 0, \\
						0 & \blowuprateofwaveWRTnumberofcommutations(N) = 0, 
					\end{cases}
					\label{E:UNIFIEDCOMMUTATORSPACETIMEL2ANGULARCONTROLLINGQUANTITY} \\	
				\weakspacetimetransvelocitycontrol_N(\ubar,u) 
				& 
				\eqdef
					\begin{cases}
						\max_{\tander^N \in \mathcal{P}^{(N)}} 
						\sum_{i=1}^3
						\newspacetimecoercive{\blowuprateofwaveWRTnumberofcommutations(N) - 1}\left[ \newuL \tander^N v^i\right] (\ubar,u), & \blowuprateofwaveWRTnumberofcommutations(N) > 0 \\
						0 & \blowuprateofwaveWRTnumberofcommutations(N) = 0,
					\end{cases}   \label{E:UNIFIEDCOMMUTATORVELOCITYSPACETIMEL2TRANSCONTROLLINGQUANTITY}
				\\		
				\totalcontrolwave_N(\ubar,u)
				&				
				\eqdef 
				\totalfluxcontrolvelocity_N(\ubar,u) +
				\strongangularcontrolvelocity_N(\ubar,u) + \weakspacetimetransvelocitycontrol_N(\ubar,u)  + \weakspacetimeangularvelocitycontrol_N(\ubar,u),
				\label{E:VELOCITYFLUXANDSTRONGANGULARL2CONTROLLINGQUANTITY}
			\end{align}
		\end{subequations}
\end{definition}

\begin{definition}[Summed $L^2$-controlling quantities] 
	\label{D:SUMMEDL2CONTROLLINGQUANTITIES}
	For positive integers $N_1 < N_2$ and non-negative integers $N$, 
	we define the following summed $L^2$-controlling quantities:
	\begin{align} \label{D:ENERGYSUMCONVENTIONS}
		\totalfluxcontrolvelocity_{[N_1,N_2]}(\ubar,u) 
		& 
		\eqdef \sum_{M = N_1}^{N_2} \totalfluxcontrolvelocity_M(\ubar,u),	\end{align}
	and similarly for the other controlling quantities. 
\end{definition}

\begin{remark}[Notation change from \cite{abbrescia2022emergence}] 

We note that in our companion work \cite{abbrescia2022emergence}, we used $\mathbb{Q}_N$ to denote the $L^2$-controlling quantities on the outgoing null hypersurfaces (as well as on the rough foliations). The reason we have changed this to $\totalfluxcontrolvelocity_N$ is because we will only derive energy estimates for the \emph{velocity} in the current paper and wanted to have a visual reminder for the reader, hence the ``$\totalfluxcontrolvelocity$''. Moreover, we stress: \emph{the $L^2$-controlling quantities on the outgoing null hypersurfaces for the \textbf{specific vorticity} in \cite{abbrescia2022emergence} is denoted by $\mathbb{V}_N$}. In the current work, we denote these by $\fluxcontrolVort_N$, see Def.\,\ref{D:FUNDAMENTALTRANSPORTL2CONTROLLINGQUANTITIES}.
\end{remark}

\subsection{The fundamental $L^2$-transport-controlling quantities} \label{SS:FUNDAMENTALTRANSPORTL2CONTROLLINGQUANTITIES} We now use the building block quantities from Sect.\,\ref{SSS:BUILDINGBLOCKENERGIESANDNULLFLUXES} to construct the ``fundamental $L^2$-transport-controlling quantities'' that we use to control $(\vortrenormalized,\GradEnt)$ and their derivatives in $L^2_q$ in various regions.

\begin{definition}[The fundamental $L^2$-transport-controlling quantities] \label{D:FUNDAMENTALTRANSPORTL2CONTROLLINGQUANTITIES} 
Let $\blowuprateoftransportWRTnumberofcommutations(N)$ denote the power of $\weight$ from Def.\,\ref{D:NOTATIONFORPOWERSOFWEIGHT}. Then in terms of the null-flux quantities of Def.\,\ref{D:NULLFLUXESFORTHEWAVEANDTRANSPORTVARIABLES} 
	and the differentiated conventions of Def.\,\ref{D:STRINGSOFCOMMUTATIONVECTORFIELDS}, 
	we define the following $L^2$-transport-controlling quantities:
	\begin{itemize}
		\item \underline{\textbf{Vorticity controlling quantities}}.

			\begin{subequations}
				\begin{align} 
					\fluxcontrolVort_N(\ubar,u) 
					& 
					\eqdef 
					\max_{\substack{\tander^N \in \mathfrak{P}^{(N)} \\ \vortrenormalized \in \{ \vortrenormalized^1,\vortrenormalized^2,\vortrenormalized^3\}}} 
						\sup_{\substack{\ubar' \in [\leftubar,\ubar]
						\\
						u' \in [\moreinterestingu_1,u]}} 
						\left\lbrace   \ingoingfluxtransport{\blowuprateoftransportWRTnumberofcommutations(N)}[\tander^N  \vortrenormalized](\ubar',u')  + \outgoingfluxtransport{\blowuprateoftransportWRTnumberofcommutations(N)}[\tander^N  \vortrenormalized](\ubar',u') \right\}, \label{E:FLUXVORTICITYL2CONTROLLINGQUANTITIES} \\
					\bulkcontrolVort_N(\ubar,u) 
					& 
					\eqdef 
						\begin{cases}
							\max_{\substack{\tander^N \in \mathfrak{P}^{(N)} \\ \vortrenormalized \in \{ \vortrenormalized^1,\vortrenormalized^2,\vortrenormalized^3\}}} 
							\newspacetimecoercive{\blowuprateoftransportWRTnumberofcommutations(N) - 1}[\sqrt{\ReciprocaluLunitAppliedtoTimeFunction}\, \tander^N  \vortrenormalized ] (\ubar,u), & \blowuprateoftransportWRTnumberofcommutations(N) > 0 \\
							0 & \blowuprateoftransportWRTnumberofcommutations(N)  = 0,
						\end{cases}
						\label{E:NEWSPACETIMEVORTICITYL2CONTROLLINGQUANTITY} \\
					\toricontrolVort_N(\ubar,u) 
					& \eqdef \max_{\substack{\tander^N \in \mathfrak{P}^{(N)} \\ \vortrenormalized \in \{ \vortrenormalized^1,\vortrenormalized^2,\vortrenormalized^3\}}} 
					\begin{cases} 
						\left\| \frac{1}{\sqrt{\upmu}} \tander^N \vortrenormalized \right\|_{L^2_{\blowupratetoporderacoustic}(\doublenulltoritwoarg{\ubar}{u})}^2 & N = \Ntop, \\
						\left\|  \tander^N \vortrenormalized \right\|_{L^2_{\blowuprateofacousticgeoWRTnumberofcommutations(N)}(\doublenulltoritwoarg{\ubar}{u})}^2 & N \le \Ntop -1.
					\end{cases} \label{E:NULLTORIVORTICITYL2CONTROLLINGQUANTITY}
				\end{align}
			\end{subequations}			
		\item \underline{\textbf{Gradient entropy controlling quantities}}.			
			\begin{subequations}
				\begin{align}
					\fluxcontrolGradEnt_N(\ubar,u) 
					& 
					\eqdef 
					\max_{\substack{\tander^N \in \mathfrak{P}^{(N)} \\ \GradEnt \in \{ \GradEnt^1,\GradEnt^2,\GradEnt^3\}}} 
						\sup_{\substack{\ubar' \in [\leftubar,\ubar]
						\\
						u' \in [\moreinterestingu_1,u]}} 
						\left\lbrace   \ingoingfluxtransport{\blowuprateoftopordertransportWRTnumberofcommutations(N)}[\tander^N  \GradEnt](\ubar',u')  + \outgoingfluxtransport{\blowuprateoftopordertransportWRTnumberofcommutations(N)}[\tander^N  \GradEnt](\ubar',u') \right\}, \label{E:FLUXGRADENTL2CONTROLLINGQUANTITIES} \\
					\bulkcontrolGradEnt_N(\ubar,u) 
					& 
					\eqdef 
						\begin{cases}
							\max_{\substack{\tander^N \in \mathfrak{P}^{(N)} \\ \GradEnt \in \{ \GradEnt^1,\GradEnt^2,\GradEnt^3\}}} 
						\newspacetimecoercive{\blowuprateoftopordertransportWRTnumberofcommutations(N) - 1}[\sqrt{\ReciprocaluLunitAppliedtoTimeFunction} \, \tander^N  \GradEnt] (\ubar,u),  & \blowuprateoftransportWRTnumberofcommutations(N) >0 \\
						0 & \blowuprateoftransportWRTnumberofcommutations(N)  = 0,
						\end{cases}
						\label{E:GRADENTNEWSPACETIMEL2CONTROLLINGQUANTITY} \\
					\toricontrolGradEnt_N(\ubar,u) 
					& \eqdef \max_{\substack{\tander^N \in \mathfrak{P}^{(N)} \\ \GradEnt \in \{ \GradEnt^1,\GradEnt^2,\GradEnt^3\}}} 
					\begin{cases} 
						\left\| \frac{1}{\sqrt{\upmu}} \tander^N \GradEnt \right\|_{L^2_{\blowupratetoporderacoustic}(\doublenulltoritwoarg{\ubar}{u})}^2 & N = \Ntop, \\
						\left\|  \tander^N \GradEnt \right\|_{L^2_{\blowuprateofacousticgeoWRTnumberofcommutations(N)}(\doublenulltoritwoarg{\ubar}{u})}^2 & N \le \Ntop -1.
					\end{cases} \label{E:NULLTORIGRADENTL2CONTROLLINGQUANTITY}					
				\end{align}
			\end{subequations}
		\item \underline{\textbf{Modified top order vorticity controlling quantities}}.			
			\begin{subequations}
				\begin{align}
					\fluxcontrolVortVort_N(\ubar,u) 
					& 
					\eqdef 
					\max_{\substack{\tander^N \in \mathfrak{P}^{(N)} \\ \VortVort \in \{ \VortVort^1,\VortVort^2,\VortVort^3\}}} 
						\sup_{\substack{\ubar' \in [\leftubar,\ubar]
						\\
						u' \in [\moreinterestingu_1,u]}} 
						\left\lbrace   \ingoingfluxtransport{\blowuprateofacousticgeoWRTnumberofcommutations(N)}[\tander^N  \VortVort](\ubar',u')  + \outgoingfluxtransport{\blowuprateofacousticgeoWRTnumberofcommutations(N)}[\tander^N  \VortVort](\ubar',u') \right\}, \label{E:FLUXVORTVORTL2CONTROLLINGQUANTITIES} \\
					\bulkcontrolVortVort_N(\ubar,u) 
					& 
					\eqdef 
						\begin{cases}
							\max_{\substack{\tander^N \in \mathfrak{P}^{(N)} \\ \VortVort \in \{ \VortVort^1,\VortVort^2,\VortVort^3\}}} 
						\newspacetimecoercive{\blowuprateofacousticgeoWRTnumberofcommutations(N) - 1}[\sqrt{\ReciprocaluLunitAppliedtoTimeFunction} \, \tander^N  \VortVort] (\ubar,u),  & \blowuprateofacousticgeoWRTnumberofcommutations(N) > 0, \\
					0 & \blowuprateofacousticgeoWRTnumberofcommutations(N) = 0,
						\end{cases} \label{E:VORTVORTNEWSPACETIMEL2CONTROLLINGQUANTITY} \\
					\toricontrolVortVort_N(\ubar,u) 
					& \eqdef \max_{\substack{\tander^N \in \mathfrak{P}^{(N)} \\ \VortVort \in \{ \VortVort^1,\VortVort^2,\VortVort^3\}}} 
						\left\|  \tander^N \VortVort \right\|_{L^2_{\blowuprateofacousticgeoWRTnumberofcommutations(N+1)}(\doublenulltoritwoarg{\ubar}{u})}^2. \label{E:NULLTORIVORTVORTL2CONTROLLINGQUANTITY}
				\end{align}
			\end{subequations}
		\item \underline{\textbf{Modified top order entropy controlling quantities}}.			
			\begin{subequations}
				\begin{align}
					\fluxcontrolDivGradEnt_N(\ubar,u) 
					& 
					\eqdef 
					\max_{\tander^N \in \mathfrak{P}^{(N)}} 
						\sup_{\substack{\ubar' \in [\leftubar,\ubar]
						\\
						u' \in [\moreinterestingu_1,u]}} 
						\left\lbrace   \ingoingfluxtransport{\blowuprateofacousticgeoWRTnumberofcommutations(N)}[\tander^N  \DivGradEnt](\ubar',u')  + \outgoingfluxtransport{\blowuprateofacousticgeoWRTnumberofcommutations(N)}[\tander^N  \DivGradEnt](\ubar',u') \right\}, \label{E:FLUXDIVGRADENTL2CONTROLLINGQUANTITIES} \\
					\bulkcontrolDivGradEnt_N(\ubar,u) 
					& 
					\eqdef 
						\begin{cases}
							\max_{\tander^N \in \mathfrak{P}^{(N)}}
						\newspacetimecoercive{\blowuprateofacousticgeoWRTnumberofcommutations(N) - 1}[\sqrt{\ReciprocaluLunitAppliedtoTimeFunction} \, \tander^N  \DivGradEnt] (\ubar,u),  & \blowuprateofacousticgeoWRTnumberofcommutations(N) > 0, \\
						0 & \blowuprateofacousticgeoWRTnumberofcommutations(N) = 0, 
					\end{cases} \label{E:DIVGRADENTNEWSPACETIMEL2CONTROLLINGQUANTITY} \\
					\toricontrolDivGradEnt_N(\ubar,u) 
					& \eqdef \max_{\tander^N \in \mathfrak{P}^{(N)}}
						\left\|  \tander^N \DivGradEnt \right\|_{L^2_{\blowuprateofacousticgeoWRTnumberofcommutations(N+1)}(\doublenulltoritwoarg{\ubar}{u})}^2. \label{E:NULLTORIDIVGRADENTL2CONTROLLINGQUANTITY}					
				\end{align}
			\end{subequations}
\end{itemize}

\end{definition}

\subsection{The fundamental $L^2$-acoustic geometry-controlling quantities}
\label{SS:FUNDAMENTALL2CONTROLLINGACOUSTICQUANTITIES}

In this section we define the analogs of the controlling quantities from Sect.\,\ref{SS:FUNDAMENTALL2CONTROLLINGWAVEQUANTITIES} for the acoustic geometry variables $\fullymodquant{\tander^N} $ and $\partialmodquant{\tander^N}$. Due to the derivative losing estimates of Sect.\,\ref{S:PRELIMINARYL2ESTIMATESFORBELOWTOPORDERDERIVATIVESOFACOUSTICGEOMETRYANDDERIVATIVELOSING}, we will only need these controlling quantities at the top-order.

\begin{definition}[The fundamental $L^2$-acoustic geometry-controlling quantities] \label{D:FUNDAMENTALL2CONTROLLINGACOUSTICQUANTITIES}
Let $N = \Ntop$ and $\blowuprateofwaveWRTnumberofcommutations(N) = \blowupratetoporderwave$  denote the powers of $\weight$ from Def.\,\ref{D:NOTATIONFORPOWERSOFWEIGHT}. Let $\sumofbrevexvisquared$ denote the function from \eqref{E:FACTORINTOPACOUSTICENERGYNEEDEDFORSHARPCONSTANTS}. Then in terms of the differentiated conventions of Def.\,\ref{D:STRINGSOFCOMMUTATIONVECTORFIELDS}, 
we define the following $L^2$-controlling quantities:
	\begin{itemize}[leftmargin=*]
		\item \underline{\textbf{Precise fully modified controlling quantities}}.
			\begin{subequations} 
				\begin{align}
					\totalfluxcontrolprecisefullymodifiedchi(\ubar,u) & \eqdef \max_{\tander^{N} \in {\mathfrak{P}}^{(N)}} \sup_{\substack{\ubar' \in [\leftubar,\ubar] \\ u' \in [\moreinterestingu_1,u]}}  \left\|\ReciprocalLunitAppliedtoTimeFunction \sumofbrevexvisquared \, \fullymodquant{\tander^{N}}  \right\|^2_{L^2_{\blowupratetoporderwave}\left( \ingoingcharacteristicsurfacetwoarg{\ubar}{[\moreinterestingu_1,u]}\right)},  \label{E:PRECISEFULLYMODQUANTFLUXTOTALCONTROL} \\
					\totalbulkcontrolprecisefullymodifiedchi(\ubar,u) & \eqdef    \max_{\tander^{N} \in {\mathfrak{P}}^{(N)}} \newspacetimecoercive{\blowupratetoporderwave-1}\left[ \ReciprocalLunitAppliedtoTimeFunction \sumofbrevexvisquared \, \fullymodquant{\tander^{N}}  \right](\ubar,u).\label{E:PRECISEFULLYMODQUANTBULKTOTALCONTROL} 
				\end{align}
			\end{subequations} 
		\item  \underline{\textbf{Imprecise fully modified controlling quantities}}.
			\begin{subequations}
				\begin{align} 
					\totalfluxcontrolimprecisefullymodifiedchi(\ubar,u) & \eqdef \max_{\tander^{N}  \in {\mathfrak{P}}^{(N)}} \sup_{\substack{\ubar' \in [\leftubar,\ubar] \\ u' \in [\moreinterestingu_1,u]}}  \left\| \fullymodquant{\tander^{N} } \right\|^2_{L^2_{\blowupratetoporderwave + 1}\left( \ingoingcharacteristicsurfacetwoarg{\ubar}{[\moreinterestingu_1,u]}\right)},  \label{E:IMPRECISEFULLYMODQUANTFLUXOTALCONTROL} \\
					\totalbulkcontrolimprecisefullymodifiedchi(\ubar,u) & \eqdef  \max_{\tander^{N}  \in {\mathfrak{P}}^{(N)}} \newspacetimecoercive{\blowupratetoporderwave}\left[ \fullymodquant{\tander^{N}  } \right](\ubar,u). \label{E:IMPRECISEFULLYMODQUANTBULKTOTALCONTROL}
				\end{align} 
			\end{subequations}
		\item \underline{\textbf{Partially modified controlling quantities}}.				
			\begin{subequations}
				\begin{align}				
				 	\totalfluxcontrolpartialmodifiedchi(\ubar,u) & \eqdef \max_{\tanderY^{N-1} \in {\mathfrak{Y}}^{(N-1)} } \left\| \sumofbrevexvisquared \, \partialmodquant{\tanderY^{N-1}} \right\|^2_{L^2_{\blowupratetoporderwave - 1}\left( \ingoingcharacteristicsurfacetwoarg{\ubar}{[\moreinterestingu_1,u]}\right)} ,\label{E:PARTIALMODQUANTFLUXTOTALCONTROL} \\
					\totalbulkcontrolpartialmodifiedchi(\ubar,u) & \eqdef \max_{\belowtopordercomforpartialmodifiedchi^{N-1} \in {\mathfrak{Y}}^{(N-1)} }\newspacetimecoercive{\blowupratetoporderwave - 2}\left[\sumofbrevexvisquared \, \partialmodquant{\tanderY^{N-1}}\right](\ubar,u). \label{E:PARTIALMODQUANTBULKTOTALCONTROL}
				\end{align}
			\end{subequations} 
	\end{itemize}

\end{definition}

\subsection{The master controlling quantity} \label{SS:MASTERCONTROLLINGQUANTITY}
At the top order, we will develop a Gr\"onwall argument to control the energies of both the wave and acoustic geometry controlling quantities simultaneously. We define below the master controlling quantities for this purpose. We do not define an analogous master controlling quantity for below-top-order $N \le \Ntop -1$ energies as they only involve the \emph{wave}.

\begin{definition}[The master controlling quantity] \label{D:MASTERCONTROLLINGQUANTITY}
	For $(\ubar,u) \in [\leftubar,\ubarboot)\times[\moreinterestingu_1,\moreinterestingu_2]$, we define the \emph{top-order master controlling quantity} 
	by:
	\begin{align}
		\begin{split}\label{E:MASTERCONTROLLINGQUANTITY}
			 \mastercontroltop(\ubar,u) & =  \totalfluxcontrolvelocity_N(\ubar,u) +
				\frac{1}{2} \strongangularcontrolvelocity_N(\ubar,u) + 2\blowupratetoporderwave \weakspacetimetransvelocitycontrol_{\Ntop}(\ubar,u)  +\frac{\blowupratetoporderwave}{2}  \weakspacetimeangularvelocitycontrol_{\Ntop}(\ubar,u) \\
				& \ \  + \totalfluxcontrolprecisefullymodifiedchi(\ubar,u) + \blowupratetoporderwave \totalbulkcontrolprecisefullymodifiedchi(\ubar,u)  \\
			&  \ \ + \totalfluxcontrolimprecisefullymodifiedchi(\ubar,u) + (\blowupratetoporderwave + 1) \totalbulkcontrolimprecisefullymodifiedchi(\ubar,u) \\
			& \ \ + \totalfluxcontrolpartialmodifiedchi(\ubar,u) +(\blowupratetoporderwave - 1)\totalbulkcontrolpartialmodifiedchi(\ubar,u).
		\end{split} 
	\end{align}
	
\end{definition} 

\begin{remark}[The role of the coefficients in the master controlling quantity] \label{R:ROLEOFBLOWUPRATEINMASTERCONTROL} We highlight the fact that the coefficients featuring $\blowupratetoporderwave$ in \eqref{E:MASTERCONTROLLINGQUANTITY} are precisely the ones which accompany the new spacetime bulk energies in the fundamental energy inequalities, e.g. \eqref{E:ENERGYNULLFLUXINTEGRALIDENTITIESWAVE}. We keep track of these universal constants as they will be used to absorb harmful error integrals present in the energy estimates featuring both large and small coefficients.

\end{remark}

\subsection{The coerciveness of the fundamental $L^2$-controlling quantities}
\label{SS:COERCIVENESSFUNDAMENTALL2CONTROLLINGQUANTITIES}
In this section, we exhibit the coerciveness properties of the $L^2$ controlling quantities
from Def.\,\ref{D:MAINCOERCIVE}.

\subsubsection{Coerciveness of the spacetime integrals $\spacetimecoercive{q}[\upxi](\ubar,u)$ 
and
$\newspacetimecoercive{q}[\upxi](\ubar,u)$}
\label{SSS:COERCIVENESSOFSPACETIMEINTEGRALS}

We begin with the spacetime bulk integrals of Def.\,\ref{D:SPACETIMECOERCIVEINTEGRALS}. Recall that these spacetime integrals
appear on the left-hand side of our energy identities \eqref{E:ENERGYNULLFLUXINTEGRALIDENTITIESWAVE}--\eqref{E:ENERGYNULLFLUXINTEGRALIDENTITIESTRANSPORT}
for the velocity and transport variables.
  
  \begin{lemma}[Coerciveness of the spacetime integrals $\spacetimecoercive{q}$ and $\newspacetimecoercive{q}$] \label{L:COERCIVITYOFSPACETIMEL2INTEGRALS} Let $\xi$ be a scalar function or an $\ell_{t,u}$-tangent one form defined on $\characteristicdiamondtwoarg{[\leftubar,\ubarboot)}{[\moreinterestingu_1,u]}$. Then the following estimate holds:
\begin{align} 
	\spacetimecoercive{q}[\xi](\ubar,u) \ge 0.92 \left\| \xi \right\|^2_{L^2_q\left(\characteristicdiamondtwoarg{[\leftubar,\ubar)}{[\moreinterestingu_1,u]} \cap \smallneighborhoodofcreasearg{[\leftubar,\ubarboot]}\right)}. \label{E:COERCIVITYOLDSPACETIMETERM}
\end{align}
Next, let $\xi$ be a scalar function or a $\doublenulltoritwoarg{\ubar}{u}$-tangent one form  defined on $\characteristicdiamondtwoarg{[\leftubar,\ubarboot)}{[\moreinterestingu_1,u]}$. Then the following identity holds:
\begin{align}
	\newspacetimecoercive{q}[\xi](\ubar,u)  =  \left\| \xi \right\|^2_{L^2_q\left(\characteristicdiamondtwoarg{[\leftubar,\ubarboot)}{[\moreinterestingu_1,\moreinterestingu_2]}\right)}. \label{E:COERCIVITYNEWSPACETIMETERM}
\end{align}
\end{lemma}

\begin{proof}

The desired estimate \eqref{E:COERCIVITYOLDSPACETIMETERM} follows immediately from the definition \eqref{E:SPACETIMECOERCIVEINTEGRALS} and \eqref{E:BOUNDSONNEWLMUINTERESTINGREGION}. The identity \eqref{E:COERCIVITYNEWSPACETIMETERM} follows from $\newL \ubar = 1$ and the definition of the weight \eqref{E:BLOWUPWEIGHT}.
\end{proof}

\subsubsection{Coercivity of top-order null fluxes and spacetime integrals for the velocity}
\label{SSS:COERCIVENESSOFTOPORDERL2CONTROLLINGQUANTITIES}

In the next lemma,
we exhibit the coerciveness of
the $L^2$-velocity-controlling quantities from Def.\,\ref{D:MAINCOERCIVE} at top-order.
In particular, we exhibit their $L^2$-coerciveness properties on both the ingoing and outgoing the characteristic surfaces
$\ingoingcharacteristicsurfacetwoarg{\ubar}{[\moreinterestingu_1,u]}$,
$\outgoingcharacteristicsurfacetwoarg{u}{[\leftubar,\ubar)}$ as well as the characteristic diamond $\characteristicdiamondtwoarg{[\leftubar,\ubarboot)}{[\moreinterestingu_1,\moreinterestingu_2]}$.

\begin{lemma}[The coercivity of the velocity $L^2$-controlling quantities]
	\label{L:COERCIVENESSOFL2CONTROLLINGQUANITIESUNIFIEDCOMMUTATOR} 
Let $\totalfluxcontrolvelocity_N(\ubar,u), \dots, \weakspacetimetransvelocitycontrol_N(\ubar,u)$ be the velocity $L^2$-controlling quantities from Def.\,\eqref{E:UNIFIEDCOMMUTATORVELOCITYNULLFLUXL2CONTROLLINGQUANTITY}--\eqref{E:UNIFIEDCOMMUTATORVELOCITYSPACETIMEL2TRANSCONTROLLINGQUANTITY},
let $(\ubar,u) \in [\leftubar,\ubarboot)\times [\moreinterestingu_1,\moreinterestingu_2]$,
and let $(\ubar',u') \in [\leftubar,\ubar] \times [\moreinterestingu_1,u]$.
Then for $1 \leq N \leq \Ntop$,
the following lower bounds hold,
where $\mathfrak{P}^{(N)}$ is the set of order $N$ commutator operators from
Def.\,\ref{D:NOTATIONFORPOWERSOFWEIGHT}, $L^2_q$ are the weighted norms from Def.\,\ref{D:GEOMETRICL2NORMS}, and $\blowuprateofwaveWRTnumberofcommutations(N)$ is as in Def.\,\ref{D:NOTATIONFORPOWERSOFWEIGHT}:
\begin{subequations}
	\begin{align}
		\begin{split} \label{E:COERCIVENESSOFNULLFLUXCONTROLWAVE}
			\totalfluxcontrolvelocity_N(\ubar,u)  & \geq \max_{\tander^N \in \mathfrak{P}^{(N)}} \left\{ \sum_{i=1}^3 \frac{1}{2} \left\|  
			\sqrt{\upmu} \left| \nullangrmD \tander^N v^i \right|_{\gnulltori} 
			\right\|^2_{L_{\blowuprateofwaveWRTnumberofcommutations(N)}^2\left(\ingoingcharacteristicsurfacetwoarg{\ubar'}{[\moreinterestingu,u']}\right)}, 
			 \sum_{i=1}^3	
			\left\| 
			\newuL \tander^N v 
			\right\|_{L_{\blowuprateofwaveWRTnumberofcommutations(N)}^2\left(\ingoingcharacteristicsurfacetwoarg{\ubar'}{[\moreinterestingu,u']}\right)}^2, \right.
			\\
			& \ \ \ \ \ \ \ \ \ \ \ \ \ \ \ 
			\left. 
			 \sum_{i=1}^3 \left\| 
			\frac{1}{\sqrt{\Lunit \ubar}} \Lunit \tander^N v^i
			\right\|_{L_{\blowuprateofwaveWRTnumberofcommutations(N)}^2\left(\outgoingcharacteristicsurfacetwoarg{u'}{[\leftubar,\ubar')}\right)}^2, 
			\,  \sum_{i=1}^3
			\left\| 
			\frac{\upmu}{\sqrt{\Lunit \ubar}}  
			\left| \angrmd \tander^N v^i \right|_{\gtorus} 
			\right\|_{L_{\blowuprateofwaveWRTnumberofcommutations(N)}^2\left(\outgoingcharacteristicsurfacetwoarg{u'}{[\leftubar,\ubar')}\right)}^2
			\right\rbrace,
		\end{split}	\\
		\begin{split} \label{E:UNIFIEDCOMMUTATORMIXANGULARSPACEITMECOERCIVITYWITHANGULARLAPLACIANINTHEINTEGRAND} 
			\strongangularcontrolvelocity_N(\ubar,u) & \ge  \max_{\tander^N \in \mathfrak{P}^{(N)}}  \sum_{i=1}^3  0.92 \left\| \angrmD  \tander^N v^i \right\|^2_{L^2_{\blowuprateofwaveWRTnumberofcommutations(N)}\left( \characteristicdiamondtwoarg{[\leftubar,\ubar)}{[\moreinterestingu_1,u]} \cap \smallneighborhoodofcreasearg{[\leftubar,\ubarboot]}\right)} 
		\end{split} \\ 
		\begin{split} \label{E:WAVENEWSPACETIMENULLANGULARL2CONTROLLINGQUANTITY} 
			\weakspacetimeangularvelocitycontrol_N(\ubar,u) & \ge  \max_{\tander^N \in \mathfrak{P}^{(N)}}  \sum_{i=1}^3   \left\| \sqrt{\upmu} \left| \nullangrmD   \tander^N v^i \right|_{\gnulltori} \right\|^2_{L^2_{\blowuprateofwaveWRTnumberofcommutations(N)-1}\left(  \characteristicdiamondtwoarg{[\leftubar,\ubar')}{[\moreinterestingu_1,u']}\right)} 
		\end{split}	\\
		\begin{split} \label{E:WAVENEWSPACETIMEL2TRANSVERSALCONTROLLINGQUANTITY} 
		\weakspacetimetransvelocitycontrol_N(\ubar,u)  & \ge \max_{\tander^N \in \mathfrak{P}^{(N)}}   \sum_{i=1}^3 \left\|  \newuL \tander^N v^i \right\|^2_{L^2_{\blowuprateofwaveWRTnumberofcommutations(N)-1}\left(  \characteristicdiamondtwoarg{[\leftubar,\ubar')}{[\moreinterestingu_1,u']}\right)} 
		\end{split}		
	\end{align}
\end{subequations}
\end{lemma}

\begin{proof}

We first prove \eqref{E:COERCIVENESSOFNULLFLUXCONTROLWAVE}. We fix any $(\ubar,u) \in [\leftubar,\ubarboot)\times [\moreinterestingu_1,\moreinterestingu_2]$,
$(\ubar',u') \in [\leftubar,\ubar] \times [\moreinterestingu_1,u]$,
and
$\tander^N \in \mathfrak{P}^{(N)}$.
First, using \eqref{E:LUNITANDULUNITAPPLIEDTOEIKONALANDCARTESIANTIME}, \eqref{E:WAVEOUTGOINGNULLFLUX}
and \eqref{E:QBREVETL} with $f \eqdef \tander^N v^i$, summing over $i = 1,2,3$
\eqref{E:GEOMETRICL2NORMSTORIANDOUTGOINGNULLHYPERSURFACES},
and \eqref{E:UNIFIEDCOMMUTATORVELOCITYNULLFLUXL2CONTROLLINGQUANTITY},
we find that
$ \sum_{i=1}^3
\left\| \frac{1}{\sqrt{\Lunit \ubar}} \Lunit \tander^N v^i
	\right\|_{L_{\blowuprateofwaveWRTnumberofcommutations(N)}^2\left(\outgoingcharacteristicsurfacetwoarg{u'}{[\leftubar,\ubar')}\right)}^2 \leq
\totalfluxcontrolvelocity_N(\ubar,u)
$
and
$ \sum_{i=1}^3
\left\| \frac{\sqrt{\upmu}}{\sqrt{\Lunit \ubar}}  
		\left| \angrmd\tander^N v^i \right|_{\gtorus} 
	\right\|_{L_{\blowuprateofwaveWRTnumberofcommutations(N)}^2\left(\outgoingcharacteristicsurfacetwoarg{u'}{[\leftubar,\ubar')}\right)}^2 
\leq
\totalfluxcontrolvelocity_N(\ubar,u)
$
as desired. 

Next, from \eqref{E:SMOOTHTORUSINVERSEFIRSTFUNDCOMPONENTSINTERMSOFDOUBLENULLTORUSINVERSEFIRSTFUNDCOMPONENTS}, \eqref{E:TORITANGENTVECTORFIELDASSOCIATEDTODOUBLENULLFRAME}, we compute: 
\begin{align} \label{E:EXPRESSIONFORNULLTORIVECTORFIELDACTINGONSCALAR}
	\ToriTangentVectorfieldAssociatedToDoubleNullFolliations f = \gnulltori^{-1}\left(\mathrm{d} x^{A},\mathrm{d} x^B\right) (\gtorusdoublenullCOV^{-1})_B^C \geop{x^C} \ubar \nullgeop{x^A} f.
\end{align}
Given any scalar function $f$ on $\doublenulltoritwoarg{\ubar}{u}$, \eqref{E:EXPRESSIONFORNULLTORIVECTORFIELDACTINGONSCALAR},  \eqref{E:SMALLC01ESTIMATESFORUBAR}, the Cauchy-Schwarz inequality, and the identity $\nullangrmD f \cdot \nullgeop{x^A} = \nullgeop{x^A} f$,  we compute $|\ToriTangentVectorfieldAssociatedToDoubleNullFolliations f| \lesssim \fundbootsmall | \nullangrmD f|_{\gnulltori}$. Next we re-express the first two terms on RHS\,\eqref{E:QBREVETNEWUL} as $(\newuL f)^2 + (\newuL f + \MagnitueofinnerproductofnewLandnewuL \ToriTangentVectorfieldAssociatedToDoubleNullFolliations f)^2 - \MagnitueofinnerproductofnewLandnewuL^2 (\ToriTangentVectorfieldAssociatedToDoubleNullFolliations f)^2$. On the other hand, from  \eqref{E:MAGNITUDEOFINNERPRODUCTOFNEWLANDNEWULAPPROXIMATELYMU} and \eqref{E:SHARPESTIMATEFORRATIOOFRATIOOFNULLGEOSICINNERPRODUCTANDFOLIATIONDENSITYANDINGOINGMU}, we see that $\upmu \ReciprocaluLunitAppliedtoTimeFunction |\nullangrmD f|_{\gnulltori}^2 - \MagnitueofinnerproductofnewLandnewuL^2 (\ToriTangentVectorfieldAssociatedToDoubleNullFolliations f)^2 \ge \frac{1}{4} \upmu^2 |\nullangrmD f|_{\gnulltori}^2$. Altogether, we have shown that:
	\begin{align}
		\enmomem[f](\multipliervectorfield,\newuL) \ge (\newuL f)^2 + (\newuL f + \MagnitueofinnerproductofnewLandnewuL \ToriTangentVectorfieldAssociatedToDoubleNullFolliations f)^2 + \frac{1}{4}\left( 2\upmu + \upmu^2\right) |\nullangrmD f|_{\gnulltori}^2. \label{E:COERCIVITYFORQBREVETNEWUL}
	\end{align}
From \eqref{E:GEOMETRICL2NORMSINGOINGNULLHYPERSURFACESANDSPACETIMEREGIONS}, \eqref{E:WAVEINGOINGNULLFLUX},  \eqref{E:UNIFIEDCOMMUTATORVELOCITYNULLFLUXL2CONTROLLINGQUANTITY}, and the now proved \eqref{E:COERCIVITYFORQBREVETNEWUL}, the desired bounds $\frac{1}{2}\sum_{i=1}^3
		\left\|  
			\sqrt{\upmu} \left| \nullangrmD \tander^N v^i \right|_{\gnulltori} 
		\right\|_{L_{\blowuprateofwaveWRTnumberofcommutations(N)}^2\left(\ingoingcharacteristicsurfacetwoarg{\ubar'}{[\moreinterestingu,u']}\right)}^2 \le \totalfluxcontrolvelocity_N(\ubar,u)$ and $ \sum_{i=1}^3
	\left\| 
		\newuL \tander^N v^i 
	\right\|_{L_{\blowuprateofwaveWRTnumberofcommutations(N)}^2\left(\ingoingcharacteristicsurfacetwoarg{\ubar'}{[\moreinterestingu,u']}\right)}^2 \le  \totalfluxcontrolvelocity_N(\ubar,u)$ hold. In total, we have proved \eqref{E:COERCIVENESSOFNULLFLUXCONTROLWAVE}

Estimates \eqref{E:UNIFIEDCOMMUTATORMIXANGULARSPACEITMECOERCIVITYWITHANGULARLAPLACIANINTHEINTEGRAND}--\eqref{E:WAVENEWSPACETIMEL2TRANSVERSALCONTROLLINGQUANTITY} follow immediately from Def.\,\ref{D:SPACETIMECOERCIVEINTEGRALS} and Lemma\,\ref{L:COERCIVITYOFSPACETIMEL2INTEGRALS}.

\end{proof}

\subsubsection{Coerciveness estimates of null fluxes and spacetime integrals for the transport variables}
\label{SSS:COERCIVENESSOFL2VORTICITYANDENTROPYCONTROLLINGQUANTITIES}
In the next lemma,
we exhibit the coercivity of
the $L^2$-transport-controlling quantities from Def.\,\ref{D:FUNDAMENTALTRANSPORTL2CONTROLLINGQUANTITIES}. The proof follows immediately from the definition; we omit the proof.

\begin{lemma}[The coercivity of the fundamental $L^2$-controlling quantities for the vorticity and gradient entropy]
	\label{L:COERCIVENESSOFL2CONTROLLINGQUANITIESFORVORTICITYANDENTROPY}
Let $\fluxcontrolVort_N(\ubar,u), \dots, \bulkcontrolDivGradEnt_N(\ubar,u)$ be the $L^2$-transport-controlling quantities from Def.\,\ref{D:FUNDAMENTALTRANSPORTL2CONTROLLINGQUANTITIES},
let $(\ubar,u) \in [\leftubar,\ubarboot)\times [\moreinterestingu_1,\moreinterestingu_2]$,
and let $(\ubar',u') \in [\leftubar,\ubar] \times [\moreinterestingu_1,u]$.
Then for $1 \leq N \leq \Ntop$,
the following lower bounds hold,
where $\mathfrak{P}^{(N)}$ is the set of order $N$ $\nullhyparg{u}$-tangential commutator operators from
Def.\,\ref{D:STRINGSOFCOMMUTATIONVECTORFIELDS}, $L^2_q$ are the weighted norms from Def.\,\ref{D:GEOMETRICL2NORMS}, and $\blowuprateoftransportWRTnumberofcommutations(N), \blowuprateofacousticgeoWRTnumberofcommutations(N)$ are as in Def.\,\ref{D:NOTATIONFORPOWERSOFWEIGHT}:
\begin{subequations}
	\begin{align}
		\fluxcontrolVort_N(\ubar,u) 
		& 
		\geq 
		 \max_{\substack{\tander^N \in \mathfrak{P}^{(N)} \\ \vortrenormalized \in \{ \vortrenormalized^1,\vortrenormalized^2,\vortrenormalized^3\}}} 
		\left\lbrace 
		\left\| 
		\sqrt{\ReciprocaluLunitAppliedtoTimeFunction} \tander^N \vortrenormalized 
		\right\|_{L_{\blowuprateoftransportWRTnumberofcommutations(N)}^2\left(\ingoingcharacteristicsurfacetwoarg{\ubar'}{[\moreinterestingu,u']}\right)}^2, 
		\,
		\left\| 
		\frac{1}{\sqrt{\Lunit\ubar}}  \tander^N \vortrenormalized 
		\right\|_{L_{\blowuprateoftransportWRTnumberofcommutations(N)}^2\left(\outgoingcharacteristicsurfacetwoarg{u'}{[\leftubar,\ubar')}\right)}^2 
		\right\rbrace, 
		\label{E:COERCIVENESSOFCONTROLVORT} 
		\\
		\bulkcontrolVort_N(\ubar,u) & = 	 \max_{\substack{\tander^N \in \mathfrak{P}^{(N)} \\ \vortrenormalized \in \{ \vortrenormalized^1,\vortrenormalized^2,\vortrenormalized^3\}}} \left\| \sqrt{\ReciprocaluLunitAppliedtoTimeFunction} \tander^N \vortrenormalized \right\|_{L_{\blowuprateoftransportWRTnumberofcommutations(N)-1}^2\left(\characteristicdiamondtwoarg{[\leftubar,\ubar')}{[\moreinterestingu_1,u']}\right)}^2, \label{E:COERCIVENESSOFNEWBULKVORTICITY} \\
		\fluxcontrolGradEnt_N(\timefunction,u) 
		& \geq 
		 \max_{\substack{\tander^N \in \mathfrak{P}^{(N)} \\ \GradEnt \in \{ \GradEnt^1,\GradEnt^2,\GradEnt^3\}}}
		\left\lbrace 
		\left\| 
		\sqrt{\ReciprocaluLunitAppliedtoTimeFunction} \tander^N  \tander^N \GradEnt 
		\right\|_{L_{\blowuprateoftransportWRTnumberofcommutations(N)}^2\left(\ingoingcharacteristicsurfacetwoarg{\ubar}{[\moreinterestingu,u]}\right)}^2, 
		\,
		\left\| 
		\frac{1}{\sqrt{\Lunit\ubar}}  \tander^N \GradEnt 
		\right\|_{L_{\blowuprateoftransportWRTnumberofcommutations(N)}^2\left(\outgoingcharacteristicsurfacetwoarg{u}{[\leftubar,\ubar)}\right)}^2
		\right\rbrace, 
		\label{E:COERCIVENESSOFCONTROLGRADENT} \\
		\bulkcontrolGradEnt_N(\ubar,u) & = \max_{\substack{\tander^N \in \mathfrak{P}^{(N)} \\ \GradEnt \in \{ \GradEnt^1,\GradEnt^2,\GradEnt^3\}}}	\left\| \sqrt{\ReciprocaluLunitAppliedtoTimeFunction} \tander^N \GradEnt \right\|_{L_{\blowuprateoftransportWRTnumberofcommutations(N)-1}^2\left(\characteristicdiamondtwoarg{[\leftubar,\ubar')}{[\moreinterestingu_1,u']}\right)}^2, \label{E:COERCIVENESSOFNEWBULKGRADENT} \\
		\fluxcontrolVortVort_N(\ubar,u) 
		& 
		\geq 
		 \max_{\substack{\tander^N \in \mathfrak{P}^{(N)} \\ \VortVort \in \{ \VortVort^1,\VortVort^2,\VortVort^3\}}} 
		\left\lbrace 
		\left\| 
		\sqrt{\ReciprocaluLunitAppliedtoTimeFunction} \tander^N \VortVort 
		\right\|_{L_{\blowuprateofacousticgeoWRTnumberofcommutations(N)}^2\left(\ingoingcharacteristicsurfacetwoarg{\ubar'}{[\moreinterestingu,u']}\right)}^2, 
		\,
		\left\| 
		\frac{1}{\sqrt{\Lunit\ubar}}  \tander^N \VortVort 
		\right\|_{L_{\blowuprateofacousticgeoWRTnumberofcommutations(N)}^2\left(\outgoingcharacteristicsurfacetwoarg{u'}{[\leftubar,\ubar')}\right)}^2 
		\right\rbrace, 
		\label{E:COERCIVENESSOFCONTROLVORTVORT} 
		\\
		\bulkcontrolVortVort_N(\ubar,u) & = 	 \max_{\substack{\tander^N \in \mathfrak{P}^{(N)} \\ \VortVort \in \{ \VortVort^1,\VortVort^2,\VortVort^3\}}} \left\| \sqrt{\ReciprocaluLunitAppliedtoTimeFunction} \tander^N \VortVort \right\|_{L_{\blowuprateofacousticgeoWRTnumberofcommutations(N)-1}^2\left(\characteristicdiamondtwoarg{[\leftubar,\ubar')}{[\moreinterestingu_1,u']}\right)}^2, \label{E:COERCIVENESSOFNEWBULKVORTVORT} \\	
		\fluxcontrolDivGradEnt_N(\ubar,u) 
		& 
		\geq 
		 \max_{\tander^N \in \mathfrak{P}^{(N)}}
		\left\lbrace 
		\left\| 
		\sqrt{\ReciprocaluLunitAppliedtoTimeFunction} \tander^N \DivGradEnt 
		\right\|_{L_{\blowuprateofacousticgeoWRTnumberofcommutations(N)}^2\left(\ingoingcharacteristicsurfacetwoarg{\ubar'}{[\moreinterestingu,u']}\right)}^2, 
		\,
		\left\| 
		\frac{1}{\sqrt{\Lunit\ubar}}  \tander^N \DivGradEnt 
		\right\|_{L_{\blowuprateofacousticgeoWRTnumberofcommutations(N)}^2\left(\outgoingcharacteristicsurfacetwoarg{u'}{[\leftubar,\ubar')}\right)}^2 
		\right\rbrace, 
		\label{E:COERCIVENESSOFCONTROLDIVGRADENT} 
		\\
		\bulkcontrolDivGradEnt_N(\ubar,u) & = 	 \max_{\tander^N \in \mathfrak{P}^{(N)}} \left\| \sqrt{\ReciprocaluLunitAppliedtoTimeFunction} \tander^N \DivGradEnt \right\|_{L_{\blowuprateofacousticgeoWRTnumberofcommutations(N)-1}^2\left(\characteristicdiamondtwoarg{[\leftubar,\ubar')}{[\moreinterestingu_1,u']}\right)}^2. \label{E:COERCIVENESSOFNEWBULKDIVGRADENT}	
	\end{align}
\end{subequations} 
\end{lemma}

\subsection{Estimates for the tangential derivatives of wave variables that allow a loss of one derivative}

In this short section we use the results of Lemma \,\ref{L:COERCIVENESSOFL2CONTROLLINGQUANITIESUNIFIEDCOMMUTATOR} to derive $L^2$ estimates for the velocity variables that feature a loss of one derivative. Importantly, however, this does \emph{not} come at a loss of $\weight$-weights. The loss of one derivative also allows us to control \emph{tangential} derivatives on the \emph{ingoing} characteristics, which is not possible at top-order without a loss.

\begin{lemma}[$L^2$ estimates for the tangential derivatives of wave variables that allow a loss of one derivative] \label{L:DERIVATIVELOSINGL2ESTIMATEFORTANGENTIALDERIVATIVES}

Let $1 \leq N \leq \Ntop$, fix $(\ubar,u) \in [\leftubar,\ubarboot)\times[\moreinterestingu_1,\moreinterestingu_2]$,  let $\blowuprateofwaveWRTnumberofcommutations(N)$ denote the powers of $\weight$ from Def.\,\ref{D:NOTATIONFORPOWERSOFWEIGHT}, and let $\totalfluxcontrolvelocity_N(\ubar,u)$ be the $L^2$-controlling quantities from Def.\,\ref{D:MAINCOERCIVE}. Then the following estimates hold for $(\ubar',u') \in [\leftubar,\ubar]\times[\moreinterestingu_1,u]$:
\begin{align}
	\left\| 
		\tander^N \velocityarray 
	\right\|_{L_{\blowuprateofwaveWRTnumberofcommutations(N)}^2\left(\ingoingcharacteristicsurfacetwoarg{\ubar'}{[\moreinterestingu,u']}\right)}^2 & \lesssim  \initialsmalldoublenull^2 + \totalfluxcontrolvelocity_N(\ubar,u).\label{E:DERIVATIVELOSINGL2ESTIMATEFORTANGENTIALDERIVATIVES}
\end{align}
\end{lemma}

\begin{proof}
Let $v \in \velocityarray = \{v^1,v^2,v^3\}$. Applying \eqref{E:SPACETIMEENERGYNULLFLUXIDENTITYFORBELOWTOPORDERACOUSTIC} to $\varphi = \tander^N v$, using the positivity of the term $\newspacetimecoercive{q-1}[\tander^N v](\ubar',u')$ revealed in \eqref{E:COERCIVITYNEWSPACETIMETERM}, bounding $|\mytr_{\gtorus}\upchi| \lesssim \fundbootsmall$ (which itself follows from \eqref{E:CHIEXPRESSSIONINTERMSOFDERIVATIVESOFLUNITI}, 
Prop.\,\ref{P:SCHEMATICSTRUCTUREOFVARIOUSTENSORSINTERMSOFCONTROLVARS},
the bootstrap assumptions,
and Prop.\,\ref{P:IMPROVEMENTOFAUXILIARYBOOTSTRAP}), \eqref{E:BOUNDSOFLUBARINTERESTINGREGION}, Young's inequality, and Gr\"onwall's inequality, we find:
\begin{align}
	\begin{split} \label{E:DERIVATIVELOSINGESTIMATEFORTANGENTIALWAVESTEP1}
		\int_{\ingoingcharacteristicsurfacetwoarg{\ubar'}{[\moreinterestingu_1,u']}} \left( \tander^N v \right)^2 \weight^{\blowuprateofwaveWRTnumberofcommutations(N)}   \volingoingnullhypersurface & \lesssim \int_{\ingoingcharacteristicsurfacetwoarg{\leftubar}{[\moreinterestingu_1,u']}} \left( \tander^N v \right)^2 \weight^{\blowuprateofwaveWRTnumberofcommutations(N)} \volingoingnullhypersurface + \int_{\characteristicdiamondtwoarg{[\leftubar,\ubar']}{[\moreinterestingu_1,u']}} \frac{1}{\Lunit \ubar} \left( \Lunit \tander^N v \right)^2\weight^{\blowuprateofwaveWRTnumberofcommutations(N)} \voldiamond \\
		&  \lesssim \int_{u'' = \moreinterestingu_1}^{u'} \left\| \tander^{N} \velocityarray 
	\right\|_{L^2_{\blowuprateofwaveWRTnumberofcommutations(N)}(\doublenulltoritwoarg{\leftubar}{u''})}^2 \, \rmd u'' + \int_{u'' = \moreinterestingu_1}^{u'} \int_{\outgoingcharacteristicsurfacetwoarg{u''}{[\leftubar,\ubar']}}  \frac{1}{\Lunit \ubar} \left( \Lunit \tander^N \velocityarray \right)^2\weight^{\blowuprateofwaveWRTnumberofcommutations(N)} \volingoingnullhypersurface  \rmd u'' \\
	&  \lesssim \initialsmalldoublenull^2 + \totalfluxcontrolvelocity_N(\ubar,u).
	\end{split}
\end{align}
We clarify that the last inequality in \eqref{E:DERIVATIVELOSINGESTIMATEFORTANGENTIALWAVESTEP1} follows from \eqref{E:SMALLDATAOFPSIONINITIALDOUBLENULLTORI}, \eqref{E:COERCIVENESSOFNULLFLUXCONTROLWAVE}, and the monotonicity of $\totalfluxcontrolvelocity_N(\ubar,u)$.

\end{proof} 

\subsection{Estimates for $\muX$ derivatives of the velocity} \label{SS:SPACETIMEESTIMATEFORMUXDERIVATIVESOFWAVEVARIABLES}
In this short section we use the results of Lemma\,\ref{L:COERCIVENESSOFL2CONTROLLINGQUANITIESUNIFIEDCOMMUTATOR} to derive $L^2$ estimates for the velocity that feature a transversal $\muX$ differentiation. This is necessary because $\muX$ is part of our multiplier method, see \eqref{E:MULTIPLIERVECTORFIELD}, whereas the controlling quantities are adapted to the null frame. 

\begin{lemma}[$L^2$ estimates for $\muX$ derivatives of wave variables] \label{L:L2ESTIMATESFORMUXTRANSVERSALDERIVATIVESOFWAVEVARIABLES}

Let $1 \leq N \leq \Ntop$, fix $(\ubar,u) \in [\leftubar,\ubarboot)\times[\moreinterestingu_1,\moreinterestingu_2]$,  let $\blowuprateofwaveWRTnumberofcommutations(N)$ denote the powers of $\weight$ from Def.\,\ref{D:NOTATIONFORPOWERSOFWEIGHT}, and let $\totalfluxcontrolvelocity_N(\ubar,u)$ be the $L^2$-controlling quantities from Def.\,\ref{D:MAINCOERCIVE}. Then the following estimates hold for $(\ubar,u) \in [\leftubar,\ubarboot)\times[\moreinterestingu_1,\moreinterestingu_2]$:
\begin{align}
	\left\| 
		\muX \tander^N \velocityarray 
	\right\|_{L_{\blowuprateofwaveWRTnumberofcommutations(N)}^2\left(\characteristicdiamondtwoarg{[\leftubar,\ubar]}{[\moreinterestingu_1,u]}\right)}^2 & \le  \int_{\ubar' = \leftubar}^{\ubar} \totalfluxcontrolvelocity_N(\ubar',u) \, \rmd \ubar' + C\int_{\moreinterestingu_1}^u \totalfluxcontrolvelocity_N(\ubar,u') \, \rmd u'. \label{E:L2ESTIMATESFORMUXTRANSVERSALDERIVATIVESOFWAVEVARIABLES}
\end{align}
\end{lemma}

\begin{proof}
The result follows from the iterated integral identities \eqref{E:OUTGOINGNULLHYPINTEGRAL}--\eqref{E:DIAMONDSPACETIMEINTEGRAL}, using \eqref{E:POINTWISESTIMATEFORBREVEXINTERMSOFNEWUL} and Cor\,\ref{C:IMPROVEAUX} to bound the integrand on LHS \eqref{E:L2ESTIMATESFORMUXTRANSVERSALDERIVATIVESOFWAVEVARIABLES}, the bootstrap assumptions, and \eqref{E:COERCIVENESSOFNULLFLUXCONTROLWAVE}.
\end{proof}

\section{The elliptic-hyperbolic integral identities}   \label{S:ELLIPTICHYPERBOLICIDENTITIES}

\subsection{Basic geometric constructions and definitions}
\label{SS:ELLIPTICHYPERBOLICIDENTITIESGEOMETRICCONSTRUCTIONSANDDEFINITIONS}
In this section, we define some basic geometric objects that play a role
in our derivation of the localized integral identities. 

\begin{definition}[Projection onto $\nullhyparg{u}$]  
	\label{D:CHARACTERISTICHYPERSURFACEPROJECTIONTENSORFIELD} \hfill
	\begin{enumerate}
		\item We define the type $\binom{1}{1}$ projection tensorfield $\Nullhypersurfaceproject$
		onto the characteristic hypersurfaces $\nullhyparg{u}$ as follows,
		where $\updelta_{\beta}^{\ \alpha}$ denotes the Kronecker delta and $\gfour(\Lunit,\uLunit) = - \MagnitueofinnerproductofLunitanduLunit$ is given as in \eqref{E:DEFOFINNERPRODUCTOFLUNITANDULUNIT}: 
		\begin{align}
			\Nullhypersurfaceproject_{\beta}^{\ \alpha}
			& 
			\eqdef 
			\updelta_{\beta}^{\ \alpha}
			+
			\frac{1}{\MagnitueofinnerproductofLunitanduLunit}
			\uLunit^{\alpha}
			\Lunit_{\beta}  \label{E:NULLHYPERSURFACEPROJECTIONINTERMSOFUNITL} \\
			& = 	\updelta_{\beta}^{\ \alpha}
			+
			\frac{1}{\MagnitueofinnerproductofnewLandnewuL}
			\newuL^{\alpha}
			\newL_{\beta} \label{E:NULLHYPERSURFACEPROJECTIONINTERMSOFNEWL}
		\end{align}
		\item Given any type $\binom{m}{n}$ spacetime tensorfield $\upxi$, 
		we define its $\nullhyparg{u}$-projection $\Nullhypersurfaceproject \upxi$ 
		as follows:
		\begin{align}
			(\Nullhypersurfaceproject \upxi)_{\beta_1 \cdots \beta_n}^{\alpha_1 \cdots \alpha_m}
			& 
			\eqdef 
			\Nullhypersurfaceproject_{\widetilde{\alpha}_1}^{\ \alpha_1}
			\cdots 
			\Nullhypersurfaceproject_{\widetilde{\alpha}_m}^{\ \alpha_m}
			(\gfour^{-1})^{\widetilde{\beta}_1 \widetilde{\gamma}_1}
			\gfour_{\beta_1 \gamma_1}
			\Nullhypersurfaceproject_{\widetilde{\beta}_1}^{\ \gamma_1}
			\cdots 
			(\gfour^{-1})^{\widetilde{\beta}_n \widetilde{\gamma}_n}
			\gfour_{\beta_n \gamma_n}
			\Nullhypersurfaceproject_{\widetilde{\beta}_n}^{\ \gamma_n}
			\upxi_{\widetilde{\gamma}_1 \cdots \widetilde{\gamma}_n}^{\widetilde{\alpha}_1 \cdots \widetilde{\alpha}_m}.
			\label{E:PROJECTIONOFTENSORONTONULLHYPERSURFACE} 
		\end{align}
	\end{enumerate}
\end{definition}

\begin{remark}[Lack of symmetry]
	\label{R:LACKOFSYMMETRYNULLHYPORJECTION}
	The type $\binom{0}{2}$ tensorfield 
	$
	\gfour_{\alpha \widetilde{\alpha}}
	\Nullhypersurfaceproject_{\beta}^{\ \widetilde{\alpha}}
	$ is not symmetric. 
	This is the reason that in equation \eqref{E:PROJECTIONOFTENSORONTONULLHYPERSURFACE},
	there are factors of $\gfour^{-1}$ and $\gfour$ and we
	were careful about the placement of the indices on $\Nullhypersurfaceproject$
	that are contracted against the lower indices of $\upxi$,
	unlike in equations
	\eqref{E:PROJECTIONOFTENSORONTOCARTESIANSIGMAT},
	\eqref{E:PROJECTIONOFTENSORONTOFLATTORUS},
	and
	\eqref{E:TENSORFIELDPROJECTIONTONULLTORIDEFININGEQUATION}.
\end{remark}

\subsubsection{The Cartesian gradient of $\upxi$ and $|Z \SigmatTan|_g$}
\label{SSS:DERIVATIVEOFATENSORFIELD}

\begin{definition}[The Cartesian gradient of $\upxi$ and $|Z \SigmatTan|_g$]
	\label{D:DERIVATIVEOFATENSORFIELD}
	\item 
	Given a type $\binom{m}{n}$ spacetime tensorfield $\upxi$, 
	we define its \emph{Cartesian gradient} 
	$\pmb{\partial} \upxi$ to be the type $\binom{m}{n+1}$ spacetime tensorfield with Cartesian components:
	\begin{align} \label{E:CARTESIANGRADIENTOFTENSORFIELD}
		(\pmb{\partial} \upxi)_{\beta_1 \beta_2 \cdots \beta_{n+1}}^{\alpha_1 \cdots \alpha_m}
		& \eqdef
		\p_{\beta_1} \upxi_{\beta_2 \cdots \beta_{n+1}}^{\alpha_1 \cdots \alpha_m}.
	\end{align}
	\item Let $\SigmatTan$ be a $\Sigma_t$-tangent vectorfield and let $Z$ be a spacetime vectorfield.
	Relative to the Cartesian coordinates, 
	we define (see Remark~\ref{R:OMITTINGZEROCOMPONENTINCARTESIANCOORDINATES}) $|Z \SigmatTan|_g \geq 0$
	as follows:
	\begin{align} \label{E:SIGMATPOINTWISENOROMOFZDERIVATIVEOFSIGMATTANGENTVECTORFIELD}
		|Z \SigmatTan|_g^2 
		& \eqdef g_{ab} (Z \SigmatTan^a) Z \SigmatTan^b,
	\end{align}
	where as usual $Z \SigmatTan^a = Z^{\alpha} \partial_{\alpha} \SigmatTan^a$.
\end{definition}

Note that $|Z \SigmatTan|_g$ is the $|\cdot|_g$ norm of the $\Sigma_t$-tangent vectorfield with the Cartesian spatial components
$Z \SigmatTan^i$, $i=1,2,3$. We also note that, viewing $Z \SigmatTan$ as a spacetime vectorfield with $Z\SigmatTan^0 =0$, we have:
	\begin{align} \label{E:SPACETIMENORMOFSIGMATTTANVECTORFIELDISTHESAMEASSIGMATINDUCEDNORM}
		|Z\SigmatTan|_{\gfour} = |Z \SigmatTan|_g .
	\end{align}

\subsection{The coercive elliptic-hyperbolic quadratic form and its coercivity}
\label{SS:NULLHYPERSURFACEADAPTEDCOERCIVEQUADRATICFORM}

\subsubsection{The coercive elliptic-hyperbolic quadratic form}
\label{SSS:DEFINITIONOFELLIPTICHYPERBOLICQUADRATICFORM}
In the next definition, we introduce the solution-adapted quadratic form $\ellipticCoerciveQuadratic$ 
that we will use to control the top-order derivatives of the specific vorticity and entropy gradient. 

\begin{definition}[The coercive elliptic-hyperbolic quadratic form]
	\label{D:NULLHYPERSURFACEADAPTEDCOERCIVEQUADRATICFORM}
	Let $\Nullhypersurfaceproject_{\beta}^{\ \alpha}$ be as in Def.\,\ref{D:CHARACTERISTICHYPERSURFACEPROJECTIONTENSORFIELD} and let $\newuL^\perp\eqdef \ReciprocaluLunitAppliedtoTimeFunction \uLunit^\perp$ and $\newL^\perp \eqdef \ReciprocalLunitAppliedtoTimeFunction \Lunit^\perp$, where $\uLunit^\perp$ and $\Lunit^\perp$ are as in Lemma\,\ref{L:BINTERMSOFDOUBLENULLFRAME}. Then we define
	$\ellipticCoerciveQuadratic[\pmb{\partial} \SigmatTan,\pmb{\partial} \SigmatTan]$
	to be the following quadratic form associated to $\SigmatTan$:
	\begin{align} 	
		\begin{split} \label{E:NULLHYPERSURFACEADAPTEDCOERCIVEQUADRATICFORM}
			\ellipticCoerciveQuadratic[\pmb{\partial} \SigmatTan,\pmb{\partial} \SigmatTan]
			&
			\eqdef 
			\Nullhypersurfaceproject_{\beta}^{\ \alpha} \Nullhypersurfaceproject_{\gamma}^{\ \delta} (\p_\alpha \SigmatTan^\gamma) \p_\kappa \SigmatTan^\nu
			(\gfour^{-1})^{\beta \kappa} \gfour_{\delta \nu}  + \left| \Transport \SigmatTan\right|^2_{g} \\
			& \ \ + \frac{1}{\MagnitueofinnerproductofnewLandnewuL^2} (\Transport \SigmatTan^a) \newuL^\perp_a (\newL^\perp \SigmatTan^b)\newL^\perp_b + \frac{1}{\MagnitueofinnerproductofnewLandnewuL^2} (\newL \SigmatTan^a) \newL^\perp_a (\Transport \SigmatTan^b)\newuL^\perp_b \\
			& \ \ - \frac{1}{\MagnitueofinnerproductofnewLandnewuL^2} \left( \frac{3-\sqrt{6}}{6} \newuL^\perp_a \newL^\perp \SigmatTan^a\right)^2 - \frac{1}{\MagnitueofinnerproductofnewLandnewuL^2} \left( \frac{3-\sqrt{6}}{6} \newL^\perp_a \newuL^\perp \SigmatTan^a\right)^2 - \frac{1}{6} (\nullangpartial_a \SigmatTan^a)^2  \\
		\end{split}
	\end{align}
\end{definition}

\subsubsection{The coercivity of the elliptic-hyperbolic quadratic form}
\label{SSS:COERCIVITYOFELLIPTICHYPERBOLICQUADRATICFORM}
In the next lemma, we exhibit the coercivity of 
$\ellipticCoerciveQuadratic[\pmb{\partial} \SigmatTan,\pmb{\partial} \SigmatTan]$.

\begin{lemma}[Coercivity of {$\ellipticCoerciveQuadratic[\pmb{\partial} \SigmatTan,\pmb{\partial} \SigmatTan]$}]
	\label{L:COERCIVENESSOFELLITPICHYPERBOLICQUADRATICFORM}
	On $\characteristicdiamondtwoarg{[\leftubar,\ubar)}{[\moreinterestingu_1,u]}$,
	the quadratic form
	$\ellipticCoerciveQuadratic$ from Def.~\ref{D:NULLHYPERSURFACEADAPTEDCOERCIVEQUADRATICFORM} satisfies the following estimates: 
		\begin{subequations}
			\begin{align}
				\frac{1}{5} \sum_{a=1}^3 |\p_a \SigmatTan|^2_g & \le \ellipticCoerciveQuadratic[\pmb{\partial} \SigmatTan,\pmb{\partial} \SigmatTan] \label{E:QUADRATICFORMCOERCIVEINALLCARTESIANSPATIALDERIVATIVES} \\
				\frac{1}{5} |\Transport \SigmatTan|^2_{g} &  \le \ellipticCoerciveQuadratic[\pmb{\partial} \SigmatTan,\pmb{\partial} \SigmatTan], \label{E:QUADRATICFORMCOERCIVEINBDERIVATIVES} \\
				\frac{1}{6}  |\pmb{\p} \SigmatTan|_{\gnulltori}^2  & \le  \ellipticCoerciveQuadratic[\pmb{\partial} \SigmatTan,\pmb{\partial} \SigmatTan], \label{E:QUADRATICFORMCOERCIVEINNULLANGULAR} 
			\end{align}
		\end{subequations}
	where
	$|\pmb{\p} \SigmatTan|_{\gnulltori}^2 = (\gnulltori^{-1})^{\alpha_1 \alpha_2} \gnulltori_{\beta_1\beta_2} (\p_{\alpha_1} \SigmatTan^{\beta_1})\p_{\alpha_2} \SigmatTan^{\beta_2}$.
	
	In particular, $\ellipticCoerciveQuadratic$ is quantitatively positive definite 
	on the space of Cartesian gradients of $\Sigma_t$-tangent vectorfields $\SigmatTan$
	in the following sense, 
	\begin{align} 
		\ellipticCoerciveQuadratic[\pmb{\partial} \SigmatTan,\pmb{\partial} \SigmatTan]
		\approx 
		\sum_{\alpha = 0}^3 |\partial_{\alpha} \SigmatTan|_g^2. \label{E:COERCIVENESSOFELLIPTICHYPERBOLICQUADRATICFORM}
	\end{align}
\end{lemma}

\begin{proof}
	We start by focusing on the first quadratic terms on the RHS\,\eqref{E:NULLHYPERSURFACEADAPTEDCOERCIVEQUADRATICFORM}, namely $\Nullhypersurfaceproject_{\beta}^{\ \alpha} \Nullhypersurfaceproject_{\gamma}^{\ \delta} (\p_\alpha \SigmatTan^\gamma) \p_\kappa \SigmatTan^\nu (\gfour^{-1})^{\beta \kappa} \gfour_{\delta \nu} $. Using \eqref{E:NULLHYPERSURFACEPROJECTIONINTERMSOFNEWL}, we have the following identity: 
	\begin{align}
		\begin{split} \label{E:NULLHYPERSURFACEADAPTEDCOERCIVEQUADRATICFORMSTEP1}
			\Nullhypersurfaceproject_{\beta}^{\ \alpha} \Nullhypersurfaceproject_{\gamma}^{\ \delta} (\p_\alpha \SigmatTan^\gamma) \p_\kappa \SigmatTan^\nu (\gfour^{-1})^{\beta \kappa} \gfour_{\delta \nu} & = (\gfour^{-1})^{\alpha\beta} \gfour_{\gamma\delta} \p_\alpha \SigmatTan^\gamma \p_\beta \SigmatTan^\delta + \frac{1}{\MagnitueofinnerproductofnewLandnewuL}\gfour_{\gamma \delta} (\newL \SigmatTan^\gamma) \newuL \SigmatTan^\delta  \\
			& \ \ + \frac{1}{\MagnitueofinnerproductofnewLandnewuL} (\gfour^{-1})^{\alpha\beta} (\p_\alpha \SigmatTan^\gamma) \newL_\gamma (\p_\beta \SigmatTan^\delta)\newuL_\delta + \frac{1}{\MagnitueofinnerproductofnewLandnewuL^2} (\newuL \SigmatTan^\alpha) \newL_\alpha (\newL \SigmatTan^\beta) \newuL_\beta
		\end{split}
	\end{align}
Using \eqref{E:ACOUSTICALINVERSEMETRICINTERMSOFNULLVECTORFIELDSANDDOUBLENULLSPHEREINVERSEFIRSTFUND} and the symmetry of $\gfour$ and $\gnulltori^{-1}$, the sum of the first two terms on RHS\,\eqref{E:NULLHYPERSURFACEADAPTEDCOERCIVEQUADRATICFORMSTEP1} are:
	\begin{align}
		\begin{split} \label{E:NULLHYPERSURFACEADAPTEDCOERCIVEQUADRATICFORMSTEP2}
			(\gfour^{-1})^{\alpha\beta} \gfour_{\gamma\delta} \p_\alpha \SigmatTan^\gamma \p_\beta \SigmatTan^\delta + \frac{1}{\MagnitueofinnerproductofnewLandnewuL}\gfour_{\gamma\delta} (\newuL \SigmatTan^\gamma) \newuL \SigmatTan^\delta & = \frac{1}{2} (\gfour^{-1})^{\alpha\beta} \gfour_{\gamma\delta} \p_\alpha \SigmatTan^\gamma \p_\beta \SigmatTan^\delta + \frac{1}{2}  \gfour_{\gamma\delta} (\gnulltori^{-1})^{\alpha\beta}( \p_\alpha \SigmatTan^\gamma) \p_\beta\SigmatTan^\delta \\
			& =  \ \  \frac{1}{2} (\gfour^{-1})^{\alpha\beta} \gfour_{\gamma\delta} \p_\alpha \SigmatTan^\gamma \p_\beta \SigmatTan^\delta + \frac{1}{2} |\pmb{\p} \SigmatTan|_{\gnulltori}^2 \\
			& \ \ -\frac{1}{\MagnitueofinnerproductofnewLandnewuL} \newuL_\gamma \newL_\delta (\gnulltori^{-1})^{\alpha\beta}( \p_\alpha \SigmatTan^\gamma) \p_\beta\SigmatTan^\delta 
		\end{split}
	\end{align}
Using again \eqref{E:ACOUSTICALINVERSEMETRICINTERMSOFNULLVECTORFIELDSANDDOUBLENULLSPHEREINVERSEFIRSTFUND} again, the last two terms on RHS\,\eqref{E:NULLHYPERSURFACEADAPTEDCOERCIVEQUADRATICFORMSTEP1} are given by:  
 \begin{align}
		\begin{split} \label{E:NULLHYPERSURFACEADAPTEDCOERCIVEQUADRATICFORMSTEP3}
			&  \frac{1}{\MagnitueofinnerproductofnewLandnewuL} (\gfour^{-1})^{\alpha\beta} (\p_\alpha \SigmatTan^\gamma) \newL_\gamma (\p_\beta \SigmatTan^\delta)\newuL_\delta + \frac{1}{\MagnitueofinnerproductofnewLandnewuL^2} (\newuL \SigmatTan^\alpha) \newL_\alpha (\newL \SigmatTan^\beta) \newuL_\beta \\
			& \ \  = - \frac{1}{\MagnitueofinnerproductofnewLandnewuL^2}  (\newuL \SigmatTan^\alpha) \newL_\alpha (\newL \SigmatTan^\beta) \newuL_\beta 
			- \frac{1}{\MagnitueofinnerproductofnewLandnewuL^2}  (\newuL \SigmatTan^\alpha) \newuL_\alpha (\newL \SigmatTan^\beta) \newL_\beta \\
			 & \ \ \ \  +  \frac{1}{\MagnitueofinnerproductofnewLandnewuL} (\gnulltori^{-1})^{\alpha\beta} (\p_\alpha \SigmatTan^\gamma) \newL_\gamma (\p_\beta \SigmatTan^\delta)\newuL_\delta 
			 + \frac{1}{\MagnitueofinnerproductofnewLandnewuL^2} (\newuL \SigmatTan^\alpha) \newL_\alpha (\newL \SigmatTan^\beta) \newuL_\beta \\
			 &  \ \ = - \frac{1}{\MagnitueofinnerproductofnewLandnewuL^2}  (\newuL \SigmatTan^\alpha) \newuL_\alpha (\newL \SigmatTan^\beta) \newL_\beta 
			 +  \frac{1}{\MagnitueofinnerproductofnewLandnewuL} (\gnulltori^{-1})^{\alpha\beta} (\p_\alpha \SigmatTan^\gamma) \newL_\gamma (\p_\beta \SigmatTan^\delta)\newuL_\delta, 
		\end{split}
	\end{align}
where we note the cancellation of the $ \frac{1}{\MagnitueofinnerproductofnewLandnewuL^2} (\newuL \SigmatTan^\alpha) \newL_\alpha (\newL \SigmatTan^\beta) \newuL_\beta$ terms. Next we note that the second term on RHS\,\eqref{E:NULLHYPERSURFACEADAPTEDCOERCIVEQUADRATICFORMSTEP3} completely cancels the last term on RHS\,\eqref{E:NULLHYPERSURFACEADAPTEDCOERCIVEQUADRATICFORMSTEP2}. Using \eqref{E:FIRSTFUNDOFSIGMATINTERMSOFACOUSTICALMETRICANDTRANSPORT} and the $\Sigma_t$-tangency of $\SigmatTan$, in total with \eqref{E:NULLHYPERSURFACEADAPTEDCOERCIVEQUADRATICFORMSTEP1}--\eqref{E:NULLHYPERSURFACEADAPTEDCOERCIVEQUADRATICFORMSTEP3} we have shown the following:
	\begin{align}
		\begin{split} \label{E:NULLHYPERSURFACEADAPTEDCOERCIVEQUADRATICFORMSTEP4}
			& \Nullhypersurfaceproject_{\beta}^{\ \alpha} \Nullhypersurfaceproject_{\gamma}^{\ \delta} (\p_\alpha \SigmatTan^\gamma) \p_\kappa \SigmatTan^\nu (\gfour^{-1})^{\beta \kappa} \gfour_{\delta \nu} +  |\Transport \SigmatTan|_g^2 \\
			& \ \  =  \frac{1}{2} |\p \SigmatTan|_g^2 + \frac{1}{2} |\Transport \SigmatTan|_g^2  - \frac{1}{\MagnitueofinnerproductofnewLandnewuL^2}  (\newuL \SigmatTan^\alpha) \newuL_\alpha (\newL \SigmatTan^\beta) \newL_\beta 
			+ \frac{1}{2} |\pmb{\p} \SigmatTan|_{\gnulltori}^2 .
		\end{split}
	\end{align}
	  
We now focus on the third term on RHS\,\eqref{E:NULLHYPERSURFACEADAPTEDCOERCIVEQUADRATICFORMSTEP4}. Using \eqref{E:LUNITINTERMSOFBANDSIGMATTAN}--\eqref{E:ULUNITINTERMSOFBANDSIGMATTAN} and  $\newuL^\perp\eqdef \ReciprocaluLunitAppliedtoTimeFunction \uLunit^\perp$ and $\newL^\perp \eqdef \ReciprocalLunitAppliedtoTimeFunction \Lunit^\perp$, it follows that: 
	\begin{align}
		\begin{split} \label{E:NULLHYPERSURFACEADAPTEDCOERCIVEQUADRATICFORMSTEP5}
			- \frac{1}{\MagnitueofinnerproductofnewLandnewuL^2}  (\newuL \SigmatTan^\alpha) \newuL_\alpha (\newL \SigmatTan^\beta) \newL_\beta & = -\frac{\ReciprocaluLunitAppliedtoTimeFunction^2 \ReciprocalLunitAppliedtoTimeFunction^2 }{\MagnitueofinnerproductofnewLandnewuL^2} \uLunit^\perp_a \Lunit^\perp_b (\Transport \SigmatTan^a)\Transport \SigmatTan^b - \frac{\ReciprocaluLunitAppliedtoTimeFunction^2 \ReciprocalLunitAppliedtoTimeFunction^2 }{\MagnitueofinnerproductofnewLandnewuL^2} \uLunit^\perp_a \Lunit^\perp_b (\uLunit^\perp \SigmatTan^a)\Lunit^\perp \SigmatTan^b \\
			& \ \ - \frac{1}{\MagnitueofinnerproductofnewLandnewuL^2} (\Transport \SigmatTan^a) \newuL^\perp_a (\newL^\perp \SigmatTan^b)\newL^\perp_b - \frac{1}{\MagnitueofinnerproductofnewLandnewuL^2} (\newL \SigmatTan^a) \newL^\perp_a (\Transport \SigmatTan^b)\newuL^\perp_b
		\end{split}
	\end{align}
We note that the last two terms on RHS\,\eqref{E:NULLHYPERSURFACEADAPTEDCOERCIVEQUADRATICFORMSTEP5} completely cancel the third and fourth terms in the definition of the quadratic form $\ellipticCoerciveQuadratic[\pmb{\partial} \SigmatTan,\pmb{\partial} \SigmatTan]$ in \eqref{E:NULLHYPERSURFACEADAPTEDCOERCIVEQUADRATICFORM}.  In total, this cancelation and identities \eqref{E:NULLHYPERSURFACEADAPTEDCOERCIVEQUADRATICFORMSTEP4}--\eqref{E:NULLHYPERSURFACEADAPTEDCOERCIVEQUADRATICFORMSTEP5} imply: 
\begin{align}
		\begin{split} \label{E:NULLHYPERSURFACEADAPTEDCOERCIVEQUADRATICFORMSTEP6}
			\ellipticCoerciveQuadratic[\pmb{\partial} \SigmatTan,\pmb{\partial} \SigmatTan] & =  \frac{1}{2} |\p \SigmatTan|_g^2  - \frac{\ReciprocaluLunitAppliedtoTimeFunction^2 \ReciprocalLunitAppliedtoTimeFunction^2 }{\MagnitueofinnerproductofnewLandnewuL^2} \uLunit^\perp_a \Lunit^\perp_b (\uLunit^\perp \SigmatTan^a)\Lunit^\perp \SigmatTan^b \\
			& - \frac{1}{\MagnitueofinnerproductofnewLandnewuL^2} \left( \frac{3-\sqrt{6}}{6} \newuL^\perp_a \newL^\perp \SigmatTan^a\right)^2 - \frac{1}{\MagnitueofinnerproductofnewLandnewuL^2} \left( \frac{3-\sqrt{6}}{6} \newL^\perp_a \newuL^\perp \SigmatTan^a\right)^2 \\
			& \ \ + \frac{1}{2} |\Transport \SigmatTan|_g^2 -\frac{\ReciprocaluLunitAppliedtoTimeFunction^2 \ReciprocalLunitAppliedtoTimeFunction^2 }{\MagnitueofinnerproductofnewLandnewuL^2} \uLunit^\perp_a \Lunit^\perp_b (\Transport \SigmatTan^a)\Transport \SigmatTan^b  \\
			& \ \ + \frac{1}{2} |\pmb{\p} \SigmatTan|_{\gnulltori}^2 - \frac{1}{6} (\nullangpartial_a \SigmatTan^a)^2  \\
		\end{split}
	\end{align}
By \eqref{E:SIZEOFLUBITPERPANDULUNITPERPANDTHEIRINNERPRODUCT}, the fact that $\Lunit^\perp$ and $\uLunit^\perp$ are $\Sigma_t$-tangent, the $g = \gfour + \Transport_\flat\otimes\Transport_\flat$ Cauchy-Schwarz inequality, we have that:
	\begin{align} \label{E:NULLHYPERSURFACEADAPTEDCOERCIVEQUADRATICFORMSTEP7}
		\left| \uLunit^\perp_a \Lunit^\perp \SigmatTan\right|, \left| \Lunit^\perp_a \uLunit^\perp \SigmatTan\right|, \left| \uLunit^\perp_a \uLunit^\perp \SigmatTan\right|, \left| \Lunit^\perp_a \Lunit^\perp \SigmatTan\right| \le \left| \p \SigmatTan\right|_g.
	\end{align}
From \eqref{E:RATIOOFNULLGEOSICINNERPRODUCTANDFOLIATIONDENSITY} and estimate \eqref{E:SHARPESTIMATEFORRATIOOFRATIOOFNULLGEOSICINNERPRODUCTANDFOLIATIONDENSITYANDINGOINGMU},  we have that $\frac{\ReciprocaluLunitAppliedtoTimeFunction^2 \ReciprocalLunitAppliedtoTimeFunction^2}{\MagnitueofinnerproductofnewLandnewuL^2}  = \frac{1}{4} + \mathcal{O}(\fundbootsmall)$. Hence, by piciking $\epsilon > 0$ sufficiently small, we can then make: 
	\begin{align} \label{E:NULLHYPERSURFACEADAPTEDCOERCIVEQUADRATICFORMSTEP8}
		\frac{1}{2} - \frac{\ReciprocaluLunitAppliedtoTimeFunction^2 \ReciprocalLunitAppliedtoTimeFunction^2}{\MagnitueofinnerproductofnewLandnewuL^2} - 2 \frac{\ReciprocaluLunitAppliedtoTimeFunction^2 \ReciprocalLunitAppliedtoTimeFunction^2}{\MagnitueofinnerproductofnewLandnewuL^2}\left( \frac{3-\sqrt{6}}{6}\right)^2 \ge \frac{1}{5}.
	\end{align}
Estimates \eqref{E:NULLHYPERSURFACEADAPTEDCOERCIVEQUADRATICFORMSTEP7}--\eqref{E:NULLHYPERSURFACEADAPTEDCOERCIVEQUADRATICFORMSTEP8} prove that the first two lines on RHS\,\eqref{E:NULLHYPERSURFACEADAPTEDCOERCIVEQUADRATICFORMSTEP6} are $\ge \frac{1}{5} |\p \SigmatTan|_g^2$.  Nearly identical arguments prove that the third line on RHS\,\eqref{E:NULLHYPERSURFACEADAPTEDCOERCIVEQUADRATICFORMSTEP6} is $\ge \frac{1}{5} |\Transport \SigmatTan|_g^2$. 

Finally, since $\SigmatTan^0 = 0$ and $|\nulltorusproject|_{\gnulltori}^2 = 2$ where $\nulltorusproject$ is as in  \eqref{E:DOUBLENULLTORIPROJECTIONDEFININGEQUATION}, by the $\gnulltori$ Cauchy-Schwarz inequality, we have that: 
	\begin{align} \label{E:NULLHYPERSURFACEADAPTEDCOERCIVEQUADRATICFORMSTEP9}
		\frac{1}{6} \left( \nullangpartial_a \SigmatTan^a\right)^2 = \frac{1}{6} \left( \nullangpartial_\alpha \SigmatTan^\beta\right)^2 = \left( \nulltorusproject_\beta^{\ \alpha} \p_\alpha \SigmatTan^\beta\right)^2 \le \frac{1}{3} |\pmb{\p} \SigmatTan|_{\gnulltori}^2.
	\end{align}
This proves that the last line on RHS\,\eqref{E:NULLHYPERSURFACEADAPTEDCOERCIVEQUADRATICFORMSTEP6} is $\ge \frac{1}{6} |\pmb{\p} \SigmatTan|_{\gnulltori}^2$. This concludes the proof of \eqref{E:QUADRATICFORMCOERCIVEINNULLANGULAR} and the lemma.
	
\end{proof}

\subsection{The characteristic currents}
\label{SS:CHARACTERISTICCURRENT}
The $\nullhyparg{u}$-tangent vectorfield $\ehcurrent[\SigmatTan,\pmb{\p}\SigmatTan]$ in the next definition plays a key role in our analysis. 
In our proof of Prop.~\ref{P:INTEGRALIDENTITYFORELLIPTICHYPERBOLICCURRENT}, we
use them for bookkeeping when integrating by parts.
We sometimes refer to the $\ehcurrent[\SigmatTan,\pmb{\p}\SigmatTan]$ 
as ``characteristic currents'' since they are tangent to $\nullhyparg{u}$, 
or ``elliptic-hyperbolic currents'' since they are the basic ingredient for the elliptic-hyperbolic integral identities. 
In our prior work \cite{abbrescia2022emergence}, we used related -- but distinct --
$\hypthreearg{\timefunction}{[- \rightu,u]}{\upkappa}$-tangent currents 
to derive elliptic-hyperbolic integral identities.
Compared to the currents in \cite{abbrescia2022emergence}, 
the ones featured in the next definition are better adapted to the structure of the singularity in the sense
that they do not generate any critical-strength error terms in our top-order $L^2$ estimates.
Moreover, when we integrate by parts over the spacetime region $\characteristicdiamondtwoarg{[\leftubar,\ubar)}{[\moreinterestingu_1,u]}$,
the $\nullhyparg{u}$-tangency of the $\ehcurrent[\SigmatTan,\pmb{\p}\SigmatTan]$
allows us to avoid boundary integrals along $\nullhyparg{u}$. This is important because
some of the acoustic geometry error terms (such as the top-order derivatives of $\upmu$)
do not have sufficiently regularity to be controlled in $L^2$ along $\nullhyparg{u}$.

We find it convenient for our upcoming analysis to extract from the characteristic current into a  ``flat part'' $\flatcurrent$. We point out that the usual Hodge elliptic identities on for vectorfields on $\R^3$ follows from the divergence theorem applied to the $\flatcurrent$.

\subsubsection{The definition of the $\nullhyparg{u}$-tangent characteristic and flat current}
\label{SSS:DEFOFCHARCURRENT}

\begin{definition}[The $\nullhyparg{u}$-tangent elliptic-hyperbolic characteristic current]
	\label{D:PUTANGENTELLIPTICHYPERBOLICCURRENT}
	Let $\SigmatTan$ be a $\Sigma_t$-tangent vectorfield.
	We define the characteristic current
	to be the vectorfield $\ehcurrent = \ehcurrent[\SigmatTan,\pmb{\partial} \SigmatTan]$
	with the following components, $(\alpha = 0,1,2,3)$:
	\begin{align} \label{E:PUTANGENTELLIPTICHYPERBOLICCURRENT}
		\ehcurrent^{\alpha}[\SigmatTan,\pmb{\partial} \SigmatTan]
		& \eqdef 
		\SigmatTan^{\gamma} 
		\Nullhypersurfaceproject_{\gamma}^{\ \beta} 
		\Nullhypersurfaceproject_{\kappa}^{\ \alpha} \partial_{\beta} \SigmatTan^{\kappa}
		-
		\SigmatTan^{\gamma} \Nullhypersurfaceproject_{\gamma}^{\ \alpha} 
		\Nullhypersurfaceproject_{\beta}^{\ \kappa}
		\partial_{\kappa} \SigmatTan^{\beta}.
	\end{align}
\end{definition}

\begin{definition}[The flat elliptic-hyperbolic current]
	\label{D:FLATELLIPTICHYPERBOLICCURRENT}
	Let $\SigmatTan$ be a $\Sigma_t$-tangent vectorfield.
	We define the characteristic current
	to be the vectorfield $\flatcurrent = \flatcurrent[\SigmatTan,\pmb{\partial} \SigmatTan]$
	with the following components, $(\alpha = 0,1,2,3)$:
	\begin{align} \label{E:FLATELLIPTICHYPERBOLICCURRENT}
		\flatcurrent^{\alpha}[\SigmatTan,\pmb{\partial} \SigmatTan]
		& \eqdef 
		\SigmatTan^{\beta} 
		\partial_{\beta} \SigmatTan^{\alpha}
		-
		\SigmatTan^{\alpha} 
		\partial_{\beta} \SigmatTan^{\beta}.
	\end{align}
\end{definition}

\begin{remark}[{$\ehcurrent^{\alpha}[\SigmatTan,\pmb{\partial} \SigmatTan]$} is $\nullhyparg{u}$-tangent]
	\label{R:CHARACTERISTICCURRENTISPUTANGENT}
	Since for any vectorfield $Z$,
	the vectorfield $\Nullhypersurfaceproject_{\beta}^{\ \alpha} Z^{\beta}$
	is $\nullhyparg{u}$-tangent,
	it follows from \eqref{E:PUTANGENTELLIPTICHYPERBOLICCURRENT} 
	that indeed, $\ehcurrent^{\alpha}[\SigmatTan,\pmb{\partial} \SigmatTan]$ 
	is $\nullhyparg{u}$-tangent.
\end{remark}

\begin{remark}[The flat current is $\Sigma_t$-tangent] Since $\SigmatTan$ is $\Sigma_t$-tangent, it follow that $\flatcurrent^0[\SigmatTan,\pmb{\partial}\SigmatTan] = 0$ and can be expressed as $\flatcurrent[\SigmatTan,\pmb{\p}\SigmatTan] = \SigmatTan^b \p_b \SigmatTan^a - \SigmatTan^a \p_b \SigmatTan^b$. However, it is very convenient to view $\flatcurrent$ as both a $\Sigma_t$-tangent vectorfield and as a spacetime vectorfield with $\flatcurrent^0 = 0$.
\end{remark}

\begin{lemma}[Relationship between the characteristic and flat currents] \label{L:RELATIONSHIPBETWEENFLATANDCHARACTERISTICCURRENTS}

The following identity holds:
	\begin{align}
		\begin{split} \label{E:RELATIONSHIPBETWEENFLATANDCHARACTERISTICCURRENTS}
			\ehcurrent^\alpha[\SigmatTan,\pmb{\p}\SigmatTan] & = \flatcurrent^\alpha[\SigmatTan,\pmb{\p}\SigmatTan] + \frac{1}{\MagnitueofinnerproductofnewLandnewuL} \newuL^\alpha \newL_\gamma \SigmatTan^\beta \p_\beta \SigmatTan^\gamma +  \frac{1}{\MagnitueofinnerproductofnewLandnewuL} \SigmatTan^\beta \newL_\beta \newuL \SigmatTan^\alpha \\
			& \ \ -  \frac{1}{\MagnitueofinnerproductofnewLandnewuL} \SigmatTan^\alpha \newL_\beta \newuL \SigmatTan^\beta -  \frac{1}{\MagnitueofinnerproductofnewLandnewuL} \newuL^\alpha \SigmatTan^\beta \newL_\beta \p_\gamma \SigmatTan^\gamma.
		\end{split}
	\end{align}
\end{lemma}

\begin{proof}
	We insert the definition of $\Nullhypersurfaceproject_{\beta}^{\ \alpha}$ \eqref{E:NULLHYPERSURFACEPROJECTIONINTERMSOFNEWL} into the RHS\,\eqref{E:PUTANGENTELLIPTICHYPERBOLICCURRENT} with the appropriate indeces relabeled. The desired identity then follows from straightforward calculations, which depends on a crucial cancellation of the cross terms $\frac{1}{\MagnitueofinnerproductofnewLandnewuL^2} \SigmatTan^\beta \newL_\beta \newuL^\alpha \newL_\gamma \newuL \SigmatTan^\gamma$.
\end{proof}

\subsection{The covariant divergence identity satisfied by the elliptic-hyperbolic current}
\label{SSS:COVARIANTDIVERGENCEIDENTITYSATISFIEDBYCHARCURRENT}
In the next lemma, we provide the main covariant divergence identity
satisfied by the current $\ehcurrent^{\alpha}[\SigmatTan,\pmb{\partial} \SigmatTan]$ 
from Def.~\ref{D:PUTANGENTELLIPTICHYPERBOLICCURRENT}.
The identity forms the starting point for the divergence-theorem based proof of Prop.~\ref{P:INTEGRALIDENTITYFORELLIPTICHYPERBOLICCURRENT}.

\begin{lemma}[Covariant divergence identity for the elliptic-hyperbolic current ]
	\label{L:COVARIANTDIVERGENCEIDENTITYFORELLIPTICHYPERBOLICCURRENT}
	Let 
	$\SigmatTan$ be a $\Sigma_t$-tagnent vectorfield,
	let $\ellipticCoerciveQuadratic[\pmb{\partial} \SigmatTan,\pmb{\partial} \SigmatTan]$
	be the quadratic form defined by \eqref{E:NULLHYPERSURFACEADAPTEDCOERCIVEQUADRATICFORM},
	and let $\widetilde \weight$ be an arbitrary ``weight function,'' not necessarily the weight $\weight$ from \eqref{E:BLOWUPWEIGHT}. 
	Then the following identity holds relative to the Cartesian coordinates, 
	where $\Dfour$ is the Levi-Civita connection of $\gfour$:
	\begin{align}	
		\begin{split} \label{E:COVARIANTDIVERGENCEIDENTITYFORELLIPTICHYPERBOLICCURRENT}
			\widetilde \weight \ellipticCoerciveQuadratic[\pmb{\partial} \SigmatTan,\pmb{\partial} \SigmatTan]
			& = \Dfour_{\alpha} (\widetilde \weight \ehcurrent^{\alpha}[\SigmatTan,\pmb{\partial} \SigmatTan])
			+
			\widetilde \weight
			\mathfrak{J}_{(\textnormal{Antisymmetric)}}[\pmb{\partial} \SigmatTan,\pmb{\partial} \SigmatTan]
			+
			\widetilde \weight
			\mathfrak{J}_{(\textnormal{Div})}[\pmb{\partial} \SigmatTan,\pmb{\partial} \SigmatTan]
			\\
			& \ \
			+
			\mathfrak{J}_{(\pmb{\partial} \widetilde \weight)}[\SigmatTan,\pmb{\partial} \SigmatTan]
			+
			\widetilde \weight
			\mathfrak{J}_{(\textnormal{Absorb}-1)}[\SigmatTan,\pmb{\partial} \SigmatTan]
			+
			\widetilde \weight
			\mathfrak{J}_{(\textnormal{Absorb}-2)}[\SigmatTan,\pmb{\partial} \SigmatTan]
			\\
			& \ \
			+
			\widetilde \weight
			\mathfrak{J}_{(\textnormal{Material})}[\pmb{\partial} \SigmatTan,\pmb{\partial} \SigmatTan]
			+
			\widetilde \weight
			\mathfrak{J}_{(\textnormal{Null Geometry})}[\SigmatTan,\pmb{\partial} \SigmatTan],
		\end{split}
	\end{align}
	where with 
	$(\mathrm{d} \SigmatTan_{\flat})_{\alpha \beta} 
	\eqdef \partial_{\alpha} \SigmatTan_{\beta} 
	- 
	\partial_{\beta} \SigmatTan_{\alpha}
	$,
	we have:
	\begin{subequations}
		\begin{align} \label{E:ANTISYMMETRICNULLCURRENTSPACETIMERRORTERM}
			\mathfrak{J}_{(\textnormal{Antisymmetric)}}[\pmb{\partial} \SigmatTan,\pmb{\partial} \SigmatTan]
			& \eqdef
			 \frac{1}{2} \Nullhypersurfaceproject_{\beta}^{\ \alpha} \Nullhypersurfaceproject_{\gamma}^{\ \delta} (\mathrm{d} \SigmatTan_\flat)_{\alpha\gamma} (\mathrm{d} \SigmatTan_\flat)_{\kappa\delta} (\gfour^{-1})^{\beta\kappa}  (\gfour^{-1})^{\gamma\nu}
			\\
			\mathfrak{J}_{(\textnormal{Div})}[\pmb{\partial} \SigmatTan,\pmb{\partial} \SigmatTan]
			& \eqdef
			\left( \frac{6-\sqrt{6}}{6}\right)^2 (\p_a \SigmatTan^a)^2,
			\label{E:DIVERGENCENULLCURRENTSPACETIMERRORTERM} 
			\\
			\mathfrak{J}_{(\pmb{\partial} \widetilde \weight)}[\SigmatTan,\pmb{\partial} \SigmatTan]
			& \eqdef
			-
			\ehcurrent^{\alpha}[\SigmatTan,\pmb{\partial} \SigmatTan] \partial_{\alpha} \widetilde \weight,
			\label{E:DERIVATIVEOFWEIGHTNULLCURRENTSPACETIMERRORTERM} 
			\\
			\begin{split} 		\label{E:ABSORBABLEANTISYMMETRICANDDIVERGENCENULLCURRENTSPACETIMERRORTERM} 
				\mathfrak{J}_{(\textnormal{Absorb}-1)}[\SigmatTan,\pmb{\partial} \SigmatTan]
				& \eqdef
				 \Nullhypersurfaceproject_{\beta}^{\ \alpha} \Nullhypersurfaceproject_{\gamma}^{\ \delta}  \SigmatTan^{\nu} \p_\alpha \gfour_{\nu \gamma'} (\mathrm{d} \SigmatTan_\flat)_{\delta\kappa}  (\gfour^{-1})^{\beta\kappa} (\gfour^{-1})^{\gamma \gamma'}
				\\
				& \ \
				- 2 \left\{  \frac{6-\sqrt{6}}{6}\p_a \SigmatTan^a + \frac{\sqrt{6}}{6} \nullangpartial_a \SigmatTan^a + \frac{\ReciprocaluLunitAppliedtoTimeFunction \ReciprocalLunitAppliedtoTimeFunction}{2\MagnitueofinnerproductofnewLandnewuL} \left(\frac{3-\sqrt{6}}{6}\right) \Lunit^\perp_a \uLunit^\perp \SigmatTan^a 
			+ \frac{\ReciprocaluLunitAppliedtoTimeFunction \ReciprocalLunitAppliedtoTimeFunction}{2\MagnitueofinnerproductofnewLandnewuL} \left(\frac{3-\sqrt{6}}{6}\right) \uLunit^\perp_a \Lunit^\perp \SigmatTan^a \right\} \\
				& \qquad \qquad \qquad \times \left\{  \frac{\ReciprocaluLunitAppliedtoTimeFunction \ReciprocalLunitAppliedtoTimeFunction}{2\MagnitueofinnerproductofnewLandnewuL} \Lunit_a^\perp  \Transport \SigmatTan^a 
				- \frac{\ReciprocaluLunitAppliedtoTimeFunction \ReciprocalLunitAppliedtoTimeFunction}{2\MagnitueofinnerproductofnewLandnewuL} \uLunit_a^\perp \Transport \SigmatTan^a + \frac{\ReciprocaluLunitAppliedtoTimeFunction \ReciprocalLunitAppliedtoTimeFunction}{2\MagnitueofinnerproductofnewLandnewuL} \upepsilon_{bak} \Lunit_a^\perp (\uLunit^\perp)^b (\Flatcurl \SigmatTan)^k\right\} \\
				& - \left\{  \frac{\ReciprocaluLunitAppliedtoTimeFunction \ReciprocalLunitAppliedtoTimeFunction}{2\MagnitueofinnerproductofnewLandnewuL} \Lunit_a^\perp  \Transport \SigmatTan^a 
			- \frac{\ReciprocaluLunitAppliedtoTimeFunction \ReciprocalLunitAppliedtoTimeFunction}{2\MagnitueofinnerproductofnewLandnewuL} \uLunit_a^\perp \Transport \SigmatTan^a + \frac{\ReciprocaluLunitAppliedtoTimeFunction \ReciprocalLunitAppliedtoTimeFunction}{2\MagnitueofinnerproductofnewLandnewuL} \upepsilon_{bak} \Lunit_a^\perp (\uLunit^\perp)^b (\Flatcurl \SigmatTan)^k\right\}^2
			\end{split}
			\\
			\begin{split}
				\mathfrak{J}_{(\textnormal{Absorb}-2)}[\SigmatTan,\pmb{\partial} \SigmatTan]
				& \eqdef
				-
				\Chfour_{\alpha \ \beta}^{\ \alpha} 
				\ehcurrent^{\beta}[\SigmatTan,\pmb{\partial} \SigmatTan]
				-
				\Nullhypersurfaceproject_{\beta}^{\ \alpha} \Nullhypersurfaceproject_{\gamma}^{\ \delta} (\p_\alpha \SigmatTan^\gamma)  \SigmatTan^\nu
			(\gfour^{-1})^{\beta \kappa} \p_\kappa \gfour_{\delta \nu} \\
				& \ \ 
				+
				\Nullhypersurfaceproject_{\beta}^{\ \alpha} \Nullhypersurfaceproject_{\gamma}^{\ \delta} (\p_\alpha \SigmatTan^{\gamma}) \SigmatTan^\kappa  (\gfour^{-1})^{\beta \nu} \p_\delta \gfour_{\kappa \nu}
				\label{E:EASYABSORBABLENULLCURRENTSPACETIMERRORTERM} 
			\end{split}
			\\
			\mathfrak{J}_{(\textnormal{Material})}[\pmb{\partial} \SigmatTan,\pmb{\partial} \SigmatTan] 
			& \eqdef
			|\Transport \SigmatTan|_g^2 + \frac{1}{\MagnitueofinnerproductofnewLandnewuL^2} (\Transport \SigmatTan^a) \newuL^\perp_a (\newL^\perp \SigmatTan^b)\newL^\perp_b + \frac{1}{\MagnitueofinnerproductofnewLandnewuL^2} (\newL \SigmatTan^a) \newL^\perp_a (\Transport \SigmatTan^b)\newuL^\perp_b,
			\label{E:QUADRATICINMATERIALDERIVATIVESNULLCURRENTSPACETIMERRORTERM} 
			\\
			\mathfrak{J}_{(\textnormal{Null Geometry})}[\SigmatTan,\pmb{\partial} \SigmatTan]
			& \eqdef
			\SigmatTan^{\gamma} 
			\left\lbrace
			\partial_{\alpha}
			\left(
			\Nullhypersurfaceproject_{\gamma}^{\ \alpha} 
			\Nullhypersurfaceproject_{\kappa}^{\ \delta} 
			\right)
			\right\rbrace
			\partial_{\delta} \SigmatTan^{\kappa}
			-
			\SigmatTan^{\gamma} 
			\left\lbrace
			\partial_{\alpha}
			\left(
			\Nullhypersurfaceproject_{\gamma}^{\ \delta} 
			\Nullhypersurfaceproject_{\kappa}^{\ \alpha} 
			\right)
			\right\rbrace
			\partial_{\delta} \SigmatTan^{\kappa},
			\label{E:DERIVATIVESOFNULLGEOMETRYNULLCURRENTSPACETIMERRORTERM}
		\end{align}
	\end{subequations}
	and on RHS\,\eqref{E:EASYABSORBABLENULLCURRENTSPACETIMERRORTERM},
	$\Chfour_{\alpha \ \beta}^{\ \gamma}  
	\eqdef \frac{1}{2} (\gfour^{-1})^{\gamma \delta} 
	\left(\partial_{\alpha} \gfour_{\delta \beta} 
	+ 
	\partial_{\beta} \gfour_{\alpha \delta} 
	-
	\partial_{\delta} \gfour_{\alpha \beta}
	\right)$
	are the Cartesian Christoffel symbols of $\gfour$.

\end{lemma}

\begin{proof}

	We begin by expanding $\Dfour_\alpha \ehcurrent^\alpha[\SigmatTan,\pmb{\p}\SigmatTan]$ relative to Cartesian coordinate partial derivatives:
	\begin{align} \label{E:COVARIANTDIVERGENCEIDENTITYFORELLIPTICHYPERBOLICCURRENTSTEP1}
		\begin{split}
			\Dfour_\alpha \ehcurrent^\alpha[\SigmatTan,\pmb{\p}\SigmatTan] & = 	
			\Nullhypersurfaceproject_{\beta}^{\ \alpha} \Nullhypersurfaceproject_{\gamma}^{\ \delta} (\p_\alpha \SigmatTan^{\gamma} )\partial_{\delta} \SigmatTan^{\beta}
			-
			\left(\Nullhypersurfaceproject_{\beta}^{\ \alpha} 
			\partial_{\alpha} \SigmatTan^{\beta}\right)^2 + \Chfour_{\alpha \ \beta}^{\ \alpha}  \ehcurrent^\beta[\SigmatTan,\pmb{\partial}\SigmatTan] \\
			& \ \ -
			\SigmatTan^{\gamma} 
			\left\lbrace
			\partial_{\alpha}
			\left(
			\Nullhypersurfaceproject_{\gamma}^{\ \alpha} 
			\Nullhypersurfaceproject_{\kappa}^{\ \delta} 
			\right)
			\right\rbrace
			\partial_{\delta} \SigmatTan^{\kappa}
			+
			\SigmatTan^{\gamma} 
			\left\lbrace
			\partial_{\alpha}
			\left(
			\Nullhypersurfaceproject_{\gamma}^{\ \delta} 
			\Nullhypersurfaceproject_{\kappa}^{\ \alpha} 
			\right)
			\right\rbrace
			\partial_{\delta} \SigmatTan^{\kappa},
		\end{split}
	\end{align}
Let us on the first derivative quadratic term on the RHS\,\eqref{E:COVARIANTDIVERGENCEIDENTITYFORELLIPTICHYPERBOLICCURRENTSTEP1}. Writing $\p_\delta \SigmatTan^\beta = \p_\delta ((\gfour^{-1})^{\beta\kappa}\SigmatTan_\kappa)$, we have: 
	\begin{align}
		\begin{split} \label{E:COVARIANTDIVERGENCEIDENTITYFORELLIPTICHYPERBOLICCURRENTSTEP2}
			\Nullhypersurfaceproject_{\beta}^{\ \alpha} \Nullhypersurfaceproject_{\gamma}^{\ \delta} (\p_\alpha \SigmatTan^{\gamma} )\partial_{\delta} \SigmatTan^{\beta} & = \Nullhypersurfaceproject_{\beta}^{\ \alpha} \Nullhypersurfaceproject_{\gamma}^{\ \delta} (\p_\alpha \SigmatTan^{\gamma} )(\partial_{\kappa} \SigmatTan_\delta) (\gfour^{-1})^{\beta\kappa} +  \Nullhypersurfaceproject_{\beta}^{\ \alpha} \Nullhypersurfaceproject_{\gamma}^{\ \delta} (\p_\alpha \SigmatTan^{\gamma} )\left( \p_\delta \SigmatTan_\kappa - \partial_{\kappa} \SigmatTan_\delta\right) (\gfour^{-1})^{\beta\kappa} \\
			& \ \ - \Nullhypersurfaceproject_{\beta}^{\ \alpha} \Nullhypersurfaceproject_{\gamma}^{\ \delta} (\p_\alpha \SigmatTan^{\gamma}) \SigmatTan^\kappa  (\gfour^{-1})^{\beta \nu} \p_\delta \gfour_{\kappa \nu}.
		\end{split}
	\end{align}
	
Denoting $(\mathrm{d} \SigmatTan_\flat)_{\delta\kappa} \eqdef \left( \p_\delta \SigmatTan_\kappa - \partial_{\kappa} \SigmatTan_\delta\right)$ and writing $\p_\alpha \SigmatTan^\gamma = \p_\alpha((\gfour^{-1})^{\gamma\nu} \SigmatTan_{\nu})$, the second term on RHS is equal to: 
	\begin{align}
		\begin{split} \label{E:COVARIANTDIVERGENCEIDENTITYFORELLIPTICHYPERBOLICCURRENTSTEP3}
			& \Nullhypersurfaceproject_{\beta}^{\ \alpha} \Nullhypersurfaceproject_{\gamma}^{\ \delta} (\p_\alpha \SigmatTan^{\gamma} )\left( \p_\delta \SigmatTan_\kappa - \partial_{\kappa} \SigmatTan_\delta\right) (\gfour^{-1})^{\beta\kappa}  \\
			& \ \ =  \Nullhypersurfaceproject_{\beta}^{\ \alpha} \Nullhypersurfaceproject_{\gamma}^{\ \delta} (\p_\alpha \SigmatTan_{\nu}) (\mathrm{d} \SigmatTan_\flat)_{\delta\kappa}  (\gfour^{-1})^{\beta\kappa}  (\gfour^{-1})^{\gamma\nu} -   \Nullhypersurfaceproject_{\beta}^{\ \alpha} \Nullhypersurfaceproject_{\gamma}^{\ \delta}  \SigmatTan^{\nu} \p_\alpha \gfour_{\nu \gamma'} (\mathrm{d} \SigmatTan_\flat)_{\delta\kappa}  (\gfour^{-1})^{\beta\kappa} (\gfour^{-1})^{\gamma \gamma'} \\ 
			& \ \ = - \frac{1}{2} \Nullhypersurfaceproject_{\beta}^{\ \alpha} \Nullhypersurfaceproject_{\gamma}^{\ \delta} (\mathrm{d} \SigmatTan_\flat)_{\alpha\gamma} (\mathrm{d} \SigmatTan_\flat)_{\kappa\delta} (\gfour^{-1})^{\beta\kappa}  (\gfour^{-1})^{\gamma\nu}  -   \Nullhypersurfaceproject_{\beta}^{\ \alpha} \Nullhypersurfaceproject_{\gamma}^{\ \delta}  \SigmatTan^{\nu} \p_\alpha \gfour_{\nu \gamma'} (\mathrm{d} \SigmatTan_\flat)_{\delta\kappa}  (\gfour^{-1})^{\beta\kappa} (\gfour^{-1})^{\gamma \gamma'}.
		\end{split}
	\end{align}Coming back to the first term on RHS\,\eqref{E:COVARIANTDIVERGENCEIDENTITYFORELLIPTICHYPERBOLICCURRENTSTEP2}, we have:
	\begin{align}
		\begin{split} \label{E:COVARIANTDIVERGENCEIDENTITYFORELLIPTICHYPERBOLICCURRENTSTEP4}
			 \Nullhypersurfaceproject_{\beta}^{\ \alpha} \Nullhypersurfaceproject_{\gamma}^{\ \delta} (\p_\alpha \SigmatTan^{\gamma} )(\partial_{\kappa} \SigmatTan_\delta) (\gfour^{-1})^{\beta\kappa} & =  \Nullhypersurfaceproject_{\beta}^{\ \alpha} \Nullhypersurfaceproject_{\gamma}^{\ \delta} (\p_\alpha \SigmatTan^{\gamma} )\p_\kappa(\SigmatTan^\nu \gfour_{\delta\nu})(\gfour^{-1})^{\beta\kappa} \\
			 & = \Nullhypersurfaceproject_{\beta}^{\ \alpha} \Nullhypersurfaceproject_{\gamma}^{\ \delta} (\p_\alpha \SigmatTan^\gamma) \p_\kappa \SigmatTan^\nu
			(\gfour^{-1})^{\beta \kappa} \gfour_{\delta \nu} + \Nullhypersurfaceproject_{\beta}^{\ \alpha} \Nullhypersurfaceproject_{\gamma}^{\ \delta} (\p_\alpha \SigmatTan^\gamma)  \SigmatTan^\nu
			(\gfour^{-1})^{\beta \kappa} \p_\kappa \gfour_{\delta \nu}.
		\end{split}
	\end{align}Next, we focus on the perfect divergence term in RHS\,\eqref{E:COVARIANTDIVERGENCEIDENTITYFORELLIPTICHYPERBOLICCURRENTSTEP1}, namely $ \left(\Nullhypersurfaceproject_{\beta}^{\ \alpha}  \partial_{\alpha} \SigmatTan^{\beta}\right)^2$. Using the $\Sigma_t$-tangency of $\SigmatTan$, we write:
	\begin{align} 
		\Nullhypersurfaceproject_{\beta}^{\ \alpha}  \partial_{\alpha} \SigmatTan^{\beta} & = \p_a \SigmatTan^a + \frac{1}{2\MagnitueofinnerproductofnewLandnewuL} \newL^\perp_a \newuL \SigmatTan^a +  \frac{1}{2\MagnitueofinnerproductofnewLandnewuL} \newL^\perp_a \newuL \SigmatTan^a. \label{E:COVARIANTDIVERGENCEIDENTITYFORELLIPTICHYPERBOLICCURRENTSTEP5}
	\end{align}		 
Using \eqref{E:ULUNITINTERMSOFBANDSIGMATTAN}, we express the last two terms on RHS\,\eqref{E:COVARIANTDIVERGENCEIDENTITYFORELLIPTICHYPERBOLICCURRENTSTEP5} as $\frac{\ReciprocaluLunitAppliedtoTimeFunction \ReciprocalLunitAppliedtoTimeFunction}{2\MagnitueofinnerproductofnewLandnewuL} \newL^\perp_a \newuL \SigmatTan^a = \frac{\ReciprocaluLunitAppliedtoTimeFunction \ReciprocalLunitAppliedtoTimeFunction}{2\MagnitueofinnerproductofnewLandnewuL} \Lunit_a^\perp \uLunit^\perp \SigmatTan^a + \frac{\ReciprocaluLunitAppliedtoTimeFunction \ReciprocalLunitAppliedtoTimeFunction}{2\MagnitueofinnerproductofnewLandnewuL} \Lunit_a^\perp  \Transport \SigmatTan^a$.  Since $\uLunit^\perp$ and $\Lunit^\perp$ are $\Sigma_t$-tangent, the first of these terms satisfies the following:
	\begin{align}
		\begin{split} \label{E:COVARIANTDIVERGENCEIDENTITYFORELLIPTICHYPERBOLICCURRENTSTEP6}
			 \frac{\ReciprocaluLunitAppliedtoTimeFunction \ReciprocalLunitAppliedtoTimeFunction}{2\MagnitueofinnerproductofnewLandnewuL} \Lunit_a^\perp \uLunit^\perp \SigmatTan^a  & =  \frac{\ReciprocaluLunitAppliedtoTimeFunction \ReciprocalLunitAppliedtoTimeFunction}{2\MagnitueofinnerproductofnewLandnewuL} \Lunit_a^\perp (\uLunit^\perp)^b \p_b \SigmatTan^a \\
			 & =  \frac{\ReciprocaluLunitAppliedtoTimeFunction \ReciprocalLunitAppliedtoTimeFunction}{2\MagnitueofinnerproductofnewLandnewuL} \Lunit_a^\perp (\uLunit^\perp)^b \p_a \SigmatTan^b +  \frac{\ReciprocaluLunitAppliedtoTimeFunction \ReciprocalLunitAppliedtoTimeFunction}{2\MagnitueofinnerproductofnewLandnewuL} \upepsilon_{bak} \Lunit_a^\perp (\uLunit^\perp)^b (\Flatcurl \SigmatTan)^k \\
			 & = \frac{\ReciprocaluLunitAppliedtoTimeFunction \ReciprocalLunitAppliedtoTimeFunction}{2\MagnitueofinnerproductofnewLandnewuL} (\Lunit^\perp)^a \uLunit^\perp_b \p_a \SigmatTan^b +  \frac{\ReciprocaluLunitAppliedtoTimeFunction \ReciprocalLunitAppliedtoTimeFunction}{2\MagnitueofinnerproductofnewLandnewuL} \upepsilon_{bak} \Lunit_a^\perp (\uLunit^\perp)^b (\Flatcurl \SigmatTan)^k \\
			 & = \frac{\ReciprocaluLunitAppliedtoTimeFunction \ReciprocalLunitAppliedtoTimeFunction}{2\MagnitueofinnerproductofnewLandnewuL} \uLunit_a^\perp \Lunit^\perp \SigmatTan^a +  \frac{\ReciprocaluLunitAppliedtoTimeFunction \ReciprocalLunitAppliedtoTimeFunction}{2\MagnitueofinnerproductofnewLandnewuL} \upepsilon_{bak} \Lunit_a^\perp (\uLunit^\perp)^b (\Flatcurl \SigmatTan)^k \\
			 & = \frac{\ReciprocaluLunitAppliedtoTimeFunction \ReciprocalLunitAppliedtoTimeFunction}{2\MagnitueofinnerproductofnewLandnewuL} \uLunit_a^\perp \Lunit \SigmatTan^a - \frac{\ReciprocaluLunitAppliedtoTimeFunction \ReciprocalLunitAppliedtoTimeFunction}{2\MagnitueofinnerproductofnewLandnewuL} \uLunit_a^\perp \Transport \SigmatTan^a  +  \frac{\ReciprocaluLunitAppliedtoTimeFunction \ReciprocalLunitAppliedtoTimeFunction}{2\MagnitueofinnerproductofnewLandnewuL} \upepsilon_{bak} \Lunit_a^\perp (\uLunit^\perp)^b (\Flatcurl \SigmatTan)^k
		\end{split}
	\end{align}
Combining \eqref{E:COVARIANTDIVERGENCEIDENTITYFORELLIPTICHYPERBOLICCURRENTSTEP5}--\eqref{E:COVARIANTDIVERGENCEIDENTITYFORELLIPTICHYPERBOLICCURRENTSTEP6}, we have thus proven that:
	\begin{align}
		\begin{split} \label{E:COVARIANTDIVERGENCEIDENTITYFORELLIPTICHYPERBOLICCURRENTSTEP7}
			\Nullhypersurfaceproject_{\beta}^{\ \alpha}  \partial_{\alpha} \SigmatTan^{\beta} & = \p_a \SigmatTan^a +
			\frac{\ReciprocaluLunitAppliedtoTimeFunction \ReciprocalLunitAppliedtoTimeFunction}{2\MagnitueofinnerproductofnewLandnewuL}  \Lunit^\perp_a \uLunit^\perp \SigmatTan^a 
			+ \frac{\ReciprocaluLunitAppliedtoTimeFunction \ReciprocalLunitAppliedtoTimeFunction}{2\MagnitueofinnerproductofnewLandnewuL} \uLunit^\perp_a \Lunit \SigmatTan^a 
			\\
			& \ \ + \frac{\ReciprocaluLunitAppliedtoTimeFunction \ReciprocalLunitAppliedtoTimeFunction}{2\MagnitueofinnerproductofnewLandnewuL} \Lunit_a^\perp  \Transport \SigmatTan^a 
			- \frac{\ReciprocaluLunitAppliedtoTimeFunction \ReciprocalLunitAppliedtoTimeFunction}{2\MagnitueofinnerproductofnewLandnewuL} \uLunit_a^\perp \Transport \SigmatTan^a + \frac{\ReciprocaluLunitAppliedtoTimeFunction \ReciprocalLunitAppliedtoTimeFunction}{2\MagnitueofinnerproductofnewLandnewuL} \upepsilon_{bak} \Lunit_a^\perp (\uLunit^\perp)^b (\Flatcurl \SigmatTan)^k.
		\end{split}
	\end{align}
Since $\uLunit^\perp_a = \uLunit_a$ and $\Lunit_a = \Lunit_a^\perp$, using \eqref{E:DEFOFINNERPRODUCTOFLUNITANDULUNIT} and \eqref{E:NEWULWITHRESPECTTOLOLDULUNIT}, the first line in RHS\,\eqref{E:COVARIANTDIVERGENCEIDENTITYFORELLIPTICHYPERBOLICCURRENTSTEP7} is equal to: 
	\begin{align}
		\begin{split} \label{E:COVARIANTDIVERGENCEIDENTITYFORELLIPTICHYPERBOLICCURRENTSTEP8}
			 \p_a \SigmatTan^a +
			\frac{\ReciprocaluLunitAppliedtoTimeFunction \ReciprocalLunitAppliedtoTimeFunction}{2\MagnitueofinnerproductofnewLandnewuL}  \Lunit^\perp_a \uLunit^\perp \SigmatTan^a 
			+ \frac{\ReciprocaluLunitAppliedtoTimeFunction \ReciprocalLunitAppliedtoTimeFunction}{2\MagnitueofinnerproductofnewLandnewuL} \uLunit^\perp_a \Lunit \SigmatTan^a & =  \frac{6-\sqrt{6}}{6}\p_a \SigmatTan^a + \frac{\sqrt{6}}{6} \p_a \SigmatTan^a
			+ \frac{\ReciprocaluLunitAppliedtoTimeFunction \ReciprocalLunitAppliedtoTimeFunction}{2\MagnitueofinnerproductofnewLandnewuL}  \Lunit^\perp_a \uLunit^\perp \SigmatTan^a 
			+ \frac{\ReciprocaluLunitAppliedtoTimeFunction \ReciprocalLunitAppliedtoTimeFunction}{2\MagnitueofinnerproductofnewLandnewuL} \uLunit^\perp_a \Lunit \SigmatTan^a \\
			& =  \frac{6-\sqrt{6}}{6}\p_a \SigmatTan^a + \frac{\sqrt{6}}{6} \nullangpartial_a \SigmatTan^a \\
			& \ \ + \frac{\ReciprocaluLunitAppliedtoTimeFunction \ReciprocalLunitAppliedtoTimeFunction}{\MagnitueofinnerproductofnewLandnewuL} \left(\frac{3-\sqrt{6}}{6}\right) \Lunit^\perp_a \uLunit^\perp \SigmatTan^a 
			+ \frac{\ReciprocaluLunitAppliedtoTimeFunction \ReciprocalLunitAppliedtoTimeFunction}{\MagnitueofinnerproductofnewLandnewuL} \left(\frac{3-\sqrt{6}}{6}\right) \uLunit^\perp_a \Lunit^\perp \SigmatTan^a 
		\end{split}
	\end{align}
Combining \eqref{E:COVARIANTDIVERGENCEIDENTITYFORELLIPTICHYPERBOLICCURRENTSTEP7}--\eqref{E:COVARIANTDIVERGENCEIDENTITYFORELLIPTICHYPERBOLICCURRENTSTEP8}, we have that: 
	\begin{align}
		\begin{split} \label{E:COVARIANTDIVERGENCEIDENTITYFORELLIPTICHYPERBOLICCURRENTSTEP9}
			\left( \Nullhypersurfaceproject_{\beta}^{\ \alpha}  \partial_{\alpha} \SigmatTan^{\beta} \right)^2 
			& = \left( \frac{6-\sqrt{6}}{6}\right)^2 (\p_a \SigmatTan^a)^2 + \frac{1}{6} ( \nullangpartial_a \SigmatTan^a)^2 
			+  \frac{\ReciprocaluLunitAppliedtoTimeFunction^2 \ReciprocalLunitAppliedtoTimeFunction^2}{\MagnitueofinnerproductofnewLandnewuL^2} \left( \frac{3-\sqrt{6}}{6} \uLunit^\perp_a \Lunit^\perp \SigmatTan^a\right)^2
			+ \frac{\ReciprocaluLunitAppliedtoTimeFunction^2 \ReciprocalLunitAppliedtoTimeFunction^2}{\MagnitueofinnerproductofnewLandnewuL^2}  \left( \frac{3-\sqrt{6}}{6} \Lunit^\perp_a \uLunit^\perp \SigmatTan^a\right)^2 \\
			& \ \ + 2 \left\{  \frac{6-\sqrt{6}}{6}\p_a \SigmatTan^a + \frac{\sqrt{6}}{6} \nullangpartial_a \SigmatTan^a + \frac{\ReciprocaluLunitAppliedtoTimeFunction \ReciprocalLunitAppliedtoTimeFunction}{\MagnitueofinnerproductofnewLandnewuL} \left(\frac{3-\sqrt{6}}{6}\right) \Lunit^\perp_a \uLunit^\perp \SigmatTan^a 
			+ \frac{\ReciprocaluLunitAppliedtoTimeFunction \ReciprocalLunitAppliedtoTimeFunction}{\MagnitueofinnerproductofnewLandnewuL} \left(\frac{3-\sqrt{6}}{6}\right) \uLunit^\perp_a \Lunit^\perp \SigmatTan^a \right\} \\
			& \qquad \qquad \qquad \times \left\{  \frac{\ReciprocaluLunitAppliedtoTimeFunction \ReciprocalLunitAppliedtoTimeFunction}{2\MagnitueofinnerproductofnewLandnewuL} \Lunit_a^\perp  \Transport \SigmatTan^a 
			- \frac{\ReciprocaluLunitAppliedtoTimeFunction \ReciprocalLunitAppliedtoTimeFunction}{2\MagnitueofinnerproductofnewLandnewuL} \uLunit_a^\perp \Transport \SigmatTan^a + \frac{\ReciprocaluLunitAppliedtoTimeFunction \ReciprocalLunitAppliedtoTimeFunction}{2\MagnitueofinnerproductofnewLandnewuL} \upepsilon_{bak} \Lunit_a^\perp (\uLunit^\perp)^b (\Flatcurl \SigmatTan)^k\right\} \\
			& \ \ +  \left\{  \frac{\ReciprocaluLunitAppliedtoTimeFunction \ReciprocalLunitAppliedtoTimeFunction}{2\MagnitueofinnerproductofnewLandnewuL} \Lunit_a^\perp  \Transport \SigmatTan^a 
			- \frac{\ReciprocaluLunitAppliedtoTimeFunction \ReciprocalLunitAppliedtoTimeFunction}{2\MagnitueofinnerproductofnewLandnewuL} \uLunit_a^\perp \Transport \SigmatTan^a + \frac{\ReciprocaluLunitAppliedtoTimeFunction \ReciprocalLunitAppliedtoTimeFunction}{2\MagnitueofinnerproductofnewLandnewuL} \upepsilon_{bak} \Lunit_a^\perp (\uLunit^\perp)^b (\Flatcurl \SigmatTan)^k\right\}^2.
		\end{split}
	\end{align}
Combining \eqref{E:COVARIANTDIVERGENCEIDENTITYFORELLIPTICHYPERBOLICCURRENTSTEP1}--\eqref{E:COVARIANTDIVERGENCEIDENTITYFORELLIPTICHYPERBOLICCURRENTSTEP9}, we have thus proven that: 
\begin{align}
		\begin{split} \label{E:COVARIANTDIVERGENCEIDENTITYFORELLIPTICHYPERBOLICCURRENTSTEP10}
			&  \Nullhypersurfaceproject_{\beta}^{\ \alpha} \Nullhypersurfaceproject_{\gamma}^{\ \delta} (\p_\alpha \SigmatTan^\gamma) \p_\kappa \SigmatTan^\nu	(\gfour^{-1})^{\beta \kappa} \gfour_{\delta \nu} 
			 - \frac{1}{\MagnitueofinnerproductofnewLandnewuL^2} \left( \frac{3-\sqrt{6}}{6} \newuL^\perp_a \newL^\perp \SigmatTan^a\right)^2 - \frac{1}{\MagnitueofinnerproductofnewLandnewuL^2} \left( \frac{3-\sqrt{6}}{6} \newL^\perp_a \newuL^\perp \SigmatTan^a\right)^2 - \frac{1}{6} (\nullangpartial_a \SigmatTan^a)^2 \\
			 & \ \ = \Dfour_\alpha \ehcurrent^\alpha[\SigmatTan,\pmb{\p}\SigmatTan]\\
			 & 
			 \ \ \ \ 
			  + \frac{1}{2} \Nullhypersurfaceproject_{\beta}^{\ \alpha} \Nullhypersurfaceproject_{\gamma}^{\ \delta} (\mathrm{d} \SigmatTan_\flat)_{\alpha\gamma} (\mathrm{d} \SigmatTan_\flat)_{\kappa\delta} (\gfour^{-1})^{\beta\kappa}  (\gfour^{-1})^{\gamma\nu} 
			  +  \Nullhypersurfaceproject_{\beta}^{\ \alpha} \Nullhypersurfaceproject_{\gamma}^{\ \delta}  \SigmatTan^{\nu} \p_\alpha \gfour_{\nu \gamma'} (\mathrm{d} \SigmatTan_\flat)_{\delta\kappa}  (\gfour^{-1})^{\beta\kappa} (\gfour^{-1})^{\gamma \gamma'} \\
			  & \ \ \ \ 
				- 2 \left\{  \frac{6-\sqrt{6}}{6}\p_a \SigmatTan^a + \frac{\sqrt{6}}{6} \nullangpartial_a \SigmatTan^a + \frac{\ReciprocaluLunitAppliedtoTimeFunction \ReciprocalLunitAppliedtoTimeFunction}{2\MagnitueofinnerproductofnewLandnewuL} \left(\frac{3-\sqrt{6}}{6}\right) \Lunit^\perp_a \uLunit^\perp \SigmatTan^a 
			+ \frac{\ReciprocaluLunitAppliedtoTimeFunction \ReciprocalLunitAppliedtoTimeFunction}{2\MagnitueofinnerproductofnewLandnewuL} \left(\frac{3-\sqrt{6}}{6}\right) \uLunit^\perp_a \Lunit^\perp \SigmatTan^a \right\} \\
				& \qquad \qquad \qquad \times \left\{  \frac{\ReciprocaluLunitAppliedtoTimeFunction \ReciprocalLunitAppliedtoTimeFunction}{2\MagnitueofinnerproductofnewLandnewuL} \Lunit_a^\perp  \Transport \SigmatTan^a 
				- \frac{\ReciprocaluLunitAppliedtoTimeFunction \ReciprocalLunitAppliedtoTimeFunction}{2\MagnitueofinnerproductofnewLandnewuL} \uLunit_a^\perp \Transport \SigmatTan^a + \frac{\ReciprocaluLunitAppliedtoTimeFunction \ReciprocalLunitAppliedtoTimeFunction}{2\MagnitueofinnerproductofnewLandnewuL} \upepsilon_{bak} \Lunit_a^\perp (\uLunit^\perp)^b (\Flatcurl \SigmatTan)^k\right\} \\
			& \ \ \ \ 
			- \left\{  \frac{\ReciprocaluLunitAppliedtoTimeFunction \ReciprocalLunitAppliedtoTimeFunction}{2\MagnitueofinnerproductofnewLandnewuL} \Lunit_a^\perp  \Transport \SigmatTan^a 
			- \frac{\ReciprocaluLunitAppliedtoTimeFunction \ReciprocalLunitAppliedtoTimeFunction}{2\MagnitueofinnerproductofnewLandnewuL} \uLunit_a^\perp \Transport \SigmatTan^a + \frac{\ReciprocaluLunitAppliedtoTimeFunction \ReciprocalLunitAppliedtoTimeFunction}{2\MagnitueofinnerproductofnewLandnewuL} \upepsilon_{bak} \Lunit_a^\perp (\uLunit^\perp)^b (\Flatcurl \SigmatTan)^k\right\}^2 \\
			  &  \ \ \ \ 
			  + \Nullhypersurfaceproject_{\beta}^{\ \alpha} \Nullhypersurfaceproject_{\gamma}^{\ \delta} (\p_\alpha \SigmatTan^{\gamma}) \SigmatTan^\kappa  (\gfour^{-1})^{\beta \nu} \p_\delta \gfour_{\kappa \nu} - \Nullhypersurfaceproject_{\beta}^{\ \alpha} \Nullhypersurfaceproject_{\gamma}^{\ \delta} (\p_\alpha \SigmatTan^\gamma)  \SigmatTan^\nu
			(\gfour^{-1})^{\beta \kappa} \p_\kappa \gfour_{\delta \nu} -
				\Chfour_{\alpha \ \beta}^{\ \alpha} 
				\ehcurrent^{\beta}[\SigmatTan,\pmb{\partial} \SigmatTan] \\
			 &
			 \ \ \ \
			 + 
			\SigmatTan^{\gamma} 
			\left\lbrace
			\partial_{\alpha}
			\left(
			\Nullhypersurfaceproject_{\gamma}^{\ \alpha} 
			\Nullhypersurfaceproject_{\kappa}^{\ \delta} 
			\right)
			\right\rbrace
			\partial_{\delta} \SigmatTan^{\kappa}
			-
			\SigmatTan^{\gamma} 
			\left\lbrace
			\partial_{\alpha}
			\left(
			\Nullhypersurfaceproject_{\gamma}^{\ \delta} 
			\Nullhypersurfaceproject_{\kappa}^{\ \alpha} 
			\right)
			\right\rbrace
			\partial_{\delta} \SigmatTan^{\kappa}.
		\end{split}
	\end{align}
			 
Taking into account the defining expression for $\ellipticCoerciveQuadratic[\pmb{\partial} \SigmatTan,\pmb{\partial} \SigmatTan]$ in \eqref{E:NULLHYPERSURFACEADAPTEDCOERCIVEQUADRATICFORM}, the proof of the lemma for $\widetilde{\weight}$ follows from adding $|\Transport \SigmatTan|_g^2 + \frac{1}{\MagnitueofinnerproductofnewLandnewuL^2} (\Transport \SigmatTan^a) \newuL^\perp_a (\newL^\perp \SigmatTan^b)\newL^\perp_b + \frac{1}{\MagnitueofinnerproductofnewLandnewuL^2} (\newL \SigmatTan^a) \newL^\perp_a (\Transport \SigmatTan^b)\newuL^\perp_b$ to both sides of \eqref{E:COVARIANTDIVERGENCEIDENTITYFORELLIPTICHYPERBOLICCURRENTSTEP10}.

The identity \eqref{E:COVARIANTDIVERGENCEIDENTITYFORELLIPTICHYPERBOLICCURRENT} incorporating an arbitrary weight $\widetilde{\weight}$ follows from multiplying the proved identity by $\widetilde{\weight}$ and differentiating by parts $\widetilde{\weight} \Dfour_\alpha \ehcurrent^\alpha[\SigmatTan,\pmb{\p}\SigmatTan] = \Dfour_\alpha \left( \widetilde{\weight} \Dfour_\alpha \ehcurrent^\alpha[\SigmatTan,\pmb{\p}\SigmatTan]\right) - \ehcurrent^\alpha[\SigmatTan,\pmb{\p}\SigmatTan]\p_\alpha\widetilde{\weight}$.

\end{proof}

\subsection{Geometric structure in the exterior derivatives $\rm{d} \tander^N \vortrenormalized$ and $\rm{d} \tander^N \GradEnt$ and higher order commutations of the fluid variables} \label{SS:GEOMETRICSTRUCTUREOFEXTERIORDERIVATIVES}

Since $\ehcurrent$ is $\nullhyparg{u}$-tangent, applying the version of the divergence theorem on the characteristic diamond from Lemma\,\ref{L:ARBITRARYDIVERGENCETHEOREM} will  only introduce boundary terms on $\ingoingcharacteristicsurfacetwoarg{\ubar}{[\moreinterestingu_1,u]}$ and $\ingoingcharacteristicsurfacetwoarg{\leftubar}{[\moreinterestingu_1,u]}$, which feature integrands with the following contraction: $\ehcurrent^\alpha[\SigmatTan,\pmb{\partial}\SigmatTan] \newuL_\alpha$ where $\pmb{\p}\SigmatTan = \pmb{\p}\tander^N\vortrenormalized$ or $\pmb{\p}\tander^N \GradEnt$. The following lemma will be crucial to prove that all top order derivatives of the velocity $\tander^{\Ntop+1}v$ featured on these boundary terms involve only $\ingoingcharacteristicsurfacetwoarg{\ubar}{[\moreinterestingu_1,u]}$-tangential derivatives as the last differential operator (e.g. $\newuL \tander^N v$).

\begin{lemma}[Sharp decomposition of  $\rm{d} \tander^N \vortrenormalized$ and $\rm{d} \tander^N \GradEnt$] \label{L:SHARPDECOMPOSITIONOFTOPORDEREXTERIORDERIVATIVEOFVORTANDGRADENT}
For a $\Sigma_t$-tangent vectorfield $\SigmatTan$, let $(\mathrm{d} \SigmatTan_\flat)_{\alpha\beta} = \partial_\alpha \SigmatTan_\beta - \partial_\beta \SigmatTan_\alpha$ denote its spacetime exterior derivative. Let $\upepsilon_{\alpha\beta\gamma\delta}$ denote the fully anti-symmetric spacetime tensor normalized by $\upepsilon_{0123} = 1$. Then, the following identities hold:
	\begin{subequations}
		\begin{align}
			\begin{split} \label{E:SHARPDECOMPOSITIONOFTOPORDEREXTERIORDERIVATIVEOFVORT}
				(\rmd[(\tander^N \vortrenormalized)_\flat])_{\alpha\beta} & \eqdef \p_\alpha (\tander^N \vortrenormalized)_\beta  -\p_\beta(\tander^N\vortrenormalized)_\alpha \\
				& = 2 (\partial_\beta \ln \Speed)(\tander^N \vortrenormalized)_\alpha - 2 (\partial_\alpha \ln \Speed)(\tander^N \vortrenormalized)_\beta  - \Transport_\alpha (\tander^N \vortrenormalized)_a \p_\beta v^a + \Transport_\beta (\tander^N \vortrenormalized)_a \p_\alpha v^a \\
				& \ \  - \Transport_\alpha ( \vortrenormalized)_a \p_\beta \tander^N v^a + \Transport_\beta ( \vortrenormalized)_a \p_\alpha \tander^N v^a + \Speed^{-2} \exp(\LogDensity) \upepsilon_{\alpha\beta\gamma\delta} \Transport^\gamma \tander^N \VortVort^\delta \\ 
				& \ \ + \Speed^{-4} \exp(-2\LogDensity) \frac{p_{;\Ent}}{\overline{\varrho}} \upepsilon_{\alpha\beta\gamma\delta} \left\{ -(\Transport \tander^N v^\gamma)\GradEnt^\delta + \Transport^\gamma \GradEnt^\delta (\p_a \tander^N v^a) - \Transport^\gamma S^a\p_a \tander^N v^\delta\right\} \\
				& \ \ + \frac{1}{\upmu} \left\{ \updelta_\alpha^0 \gfour_{\beta d} - \updelta_\beta^0 \gfour_{\alpha d}\right\}\, ^{(\tander^N)}\ehcurrent_{(\vortrenormalized;1)}^d  + \frac{1}{\upmu} \Speed^{-2} \upepsilon_{\alpha\beta\gamma d} \Transport^\gamma \,^{(\tander^N)}\ehcurrent_{(\vortrenormalized;2)}^d,
			\end{split} \\
			\begin{split} \label{E:SHARPDECOMPOSITIONOFTOPORDEREXTERIORDERIVATIVEOFGRADENT}
				(\rmd[(\tander^N \GradEnt)_\flat])_{\alpha\beta} & \eqdef \p_\alpha (\tander^N \GradEnt)_\beta  -\p_\beta(\tander^N\GradEnt)_\alpha \\			
				& =  2 (\partial_\beta \ln \Speed)(\tander^N \GradEnt)_\alpha - 2 (\partial_\alpha \ln \Speed)(\tander^N \vortrenormalized)_\beta  - \Transport_\alpha (\tander^N \GradEnt)_a \p_\beta v^a + \Transport_\beta (\tander^N \GradEnt)_a \p_\alpha v^a \\
				& \ \  + \Transport_\alpha ( \GradEnt)_a \p_\beta \tander^N v^a - \Transport_\beta ( \GradEnt)_a \p_\alpha \tander^N v^a \\
				& \ \ + \frac{1}{\upmu} \left\{ \updelta_\alpha^0 \gfour_{\beta d} - \updelta_\beta^0 \gfour_{\alpha d}\right\}\, ^{(\tander^N)}\ehcurrent_{(\GradEnt;1)}^d  + \frac{1}{\upmu} \Speed^{-2} \upepsilon_{\alpha\beta\gamma d} \Transport^\gamma \, ^{(\tander^N)}\ehcurrent_{(\GradEnt;2)}^d, 
			\end{split} 				
		\end{align}
	\end{subequations}
where:  
	\begin{subequations}
		\begin{align}
			\begin{split} \label{E:SHARPDECOMPOSITIONOFTOPORDEREXTERIORDERIVATIVEOFVORTFIRSTERRORTERM}
				^{(\tander^N)}\ehcurrent_{(\vortrenormalized;1)}^d & = [\upmu\Transport,\tander^N] \vortrenormalized^d - \vortrenormalized^a[\upmu \p_d,\tander^N]v^a + \sum_{\substack{ \tander^{N_1}\tander^{N_2} = \tander^N \\ N_2 \le N - 1}} (\tander^{N_1} \vortrenormalized^a) \tander^{N_2} (\upmu \p_d v^a) \\
				& \ \ - \exp(-2\LogDensity) \Speed^{-2} \frac{p_{;\Ent}}{\overline{\varrho}} \upepsilon_{dab} \left([ \upmu \Transport,\tander^N]v^a\right) \GradEnt^b \\
				& \ \ + \sum_{\substack{ \tander^{N_1}\tander^{N_2}\tander^{N_3} = \tander^N \\ N_2 \le N - 1}} \upepsilon_{dab} \left(\tander^{N_1} \left\{ \exp(-2\LogDensity)\Speed^{-2} \frac{p_{;\Ent}}{\overline{\varrho}} \right\}\right)\left( \tander^{N_2} (\upmu \Transport v^a) \right) \tander^{N_3} \GradEnt^b, 
			\end{split} \\
			\begin{split} \label{E:SHARPDECOMPOSITIONOFTOPORDEREXTERIORDERIVATIVEOFVORTSECONDERRORTERM}
				^{(\tander^N)}\ehcurrent_{(\vortrenormalized;2)}^d & = \upepsilon_{dab}[\upmu \p_a,\tander^N]\vortrenormalized^b +  \sum_{\substack{ \tander^{N_1}\tander^{N_2}\tander^{N_3} = \tander^N \\ N_3 \le N - 1}} \left( \tander^{N_1} \upmu\right) \left(\tander^{N_2}\exp(\LogDensity)\right) \tander^{N_3} \VortVort \\
				&  \ \ - \exp(-2\LogDensity) \Speed^{-2} \frac{p_{;\Ent}}{\overline{\varrho}} \left( [\upmu\p_a,\tander^N] v^a\right) \GradEnt^d + \exp(-2\LogDensity) \Speed^{-2} \frac{p_{;\Ent}}{\overline{\varrho}} \GradEnt^a [\upmu\p_a,\tander^N] v^d  \\
				& \ \ +  \sum_{\substack{\tander^{N_1}\tander^{N_2}\tander^{N_3} = \tander^N \\ N_2 \le N - 1}} \left( \tander^{N_1} \left\{ \exp(-2\LogDensity) \Speed^{-2} \frac{p_{;\Ent}}{\overline{\varrho}} \right\}\right) \left( \tander^{N_2} \left\{ \upmu \p_a v^a\right\} \right) \tander^{N_3} \GradEnt^d \\
				& \ \ - \sum_{\substack{\tander^{N_1}\tander^{N_2}\tander^{N_3} = \tander^N \\ N_3 \le N - 1}} \left( \tander^{N_1} \left\{ \exp(-2\LogDensity) \Speed^{-2} \frac{p_{;\Ent}}{\overline{\varrho}} \right\}\right) \left( \tander^{N_2} \GradEnt^a\right) \tander^{N_3} \left( \upmu \p_a v^d\right),
			\end{split}
		\end{align}
	\end{subequations}
and: 
	\begin{subequations}
		\begin{align} 
			^{(\tander^N)}\ehcurrent_{(\GradEnt;1)}^d & = [\upmu \Transport,\tander^N] \GradEnt^d + S^a[\upmu \p_d ,\tander^N]v^a - \sum_{\substack{ \tander^{N_1}\tander^{N_2} = \tander^N \\ N_2 \le N - 1}} \left(\tander^{N_1} \GradEnt^a\right) \tander^{N_2} \left( \upmu \p_d v^a\right), \label{E:SHARPDECOMPOSITIONOFTOPORDEREXTERIORDERIVATIVEOFGRADENTFIRSTERRORTERM} \\
			^{(\tander^N)}\ehcurrent_{(\GradEnt;2)}^d & = \upepsilon_{dab} [\upmu\p_a,\tander^N] \GradEnt^b. \label{E:SHARPDECOMPOSITIONOFTOPORDEREXTERIORDERIVATIVEOFGRADENTSECONDERRORTERM}
		\end{align}
	\end{subequations}
\end{lemma}

\begin{proof}
This was proved in \cite{abbrescia2025remarkable}*{Lemma 9.7} with $\upmu$ taking the role of $\weight_{(2)}$. We note that the factors of $\Transport_\alpha$ present in the RHS\,\eqref{E:SHARPDECOMPOSITIONOFTOPORDEREXTERIORDERIVATIVEOFVORT}--\eqref{E:SHARPDECOMPOSITIONOFTOPORDEREXTERIORDERIVATIVEOFGRADENT} were expressed as $-\updelta_\alpha^0$ in \cite{abbrescia2025remarkable}*{Lemma 9.7}.

\end{proof}

\subsection{Top order commutation identities for first order fluid equations and modified fluid variables} \label{SS:STRUCTUREOFHIGHERORDERFLUIDCOMMUTATIONS}

In order to properly exploit the geometric structures of Lemma\,\ref{L:RELATIONSHIPBETWEENFLATANDCHARACTERISTICCURRENTS}, we will need to carefully commute the first order formulation of the compressible Euler equations \eqref{E:FIRSTORDERFORMULATION} with $\tander^N$ as well as the definition of the modified fluid variables.

\begin{lemma}[Top order commutation identities in the first order formulation] \label{E:COMMUTEDFIRSTORDERFORMULATION}
The following identities hold for $i = 1,2,3$:
	\begin{subequations}
		\begin{align}
			\upmu \Transport \tander^N \LogDensity & = -\upmu \p_a \tander^N v^a + \, ^{(\tander^N)}\ehcurrent_{(\LogDensity)},  \label{E:TOPORDERCOMMUTEDTRANSPORTEQUATIONFORDENSITY} \\ 
			\upmu \Transport \tander^N v^i & = - \upmu \Speed^2 \p_i \tander^N \LogDensity + \, ^{(\tander^N)}\ehcurrent_{(v;1)}^i,  \label{E:TOPORDERCOMMUTEDTRANSPORTEQUATIONFORVELOCITY} \\
			\upmu \Transport \tander^N \vortrenormalized^i & = \upmu \vortrenormalized^a \p_i \tander^N v^a - \upmu \exp(-2\LogDensity) \Speed^{-2} \frac{p_{;\Ent}}{\overline{\varrho}} \upepsilon_{iab}\left(\Transport \tander^N v^a\right) \GradEnt^b + 	\, ^{(\tander^N)}\ehcurrent_{(\vortrenormalized;1)}^i,  \label{E:TOPORDERCOMMUTEDTRANSPORTEQUATIONFORVORTICITY}\\
			\upmu \Transport \tander^N \GradEnt^i & = - \upmu \GradEnt^a \p_i \tander^N v^a + \, ^{(\tander^N)}\ehcurrent_{(\GradEnt;1)}^i, \label{E:TOPORDERCOMMUTEDTRANSPORTEQUATIONFORGRADENT}
		\end{align}
	\end{subequations}
where $^{(\tander^N)}\ehcurrent_{(\vortrenormalized;1)}^i$ and $^{(\tander^N)}\ehcurrent_{(\GradEnt;1)}^i$ are as in \eqref{E:SHARPDECOMPOSITIONOFTOPORDEREXTERIORDERIVATIVEOFVORTFIRSTERRORTERM},\eqref{E:SHARPDECOMPOSITIONOFTOPORDEREXTERIORDERIVATIVEOFGRADENTFIRSTERRORTERM}, respectively, and: 
\begin{subequations}
		\begin{align}
			 ^{(\tander^N)}\ehcurrent_{(\LogDensity)} & = [\upmu\Transport,\tander^N]\LogDensity + [\upmu\p_a,\tander^N]v^a, \label{E:TOPORDERCOMMUTEDTRANSPORTEQUATIONFORDENSITYERRORTERM} \\
			 \begin{split} \label{E:TOPORDERCOMMUTEDTRANSPORTEQUATIONFORVELOCITYERRORTERM}
			 	^{(\tander^N)}\ehcurrent_{(v;1)}^i & = [\upmu \Transport ,\tander^N ] v^i + \Speed^2 [\upmu \p_i,\tander^N]\LogDensity - \sum_{\substack{ \tander^{N_1}\tander^{N_2} =\tander^N \\ N_2 \le N-1}} \left( \tander^{N_1}(\Speed^2)\right) \tander^{N_2}\left( \upmu \p_i \LogDensity\right) \\
				 & \ \ - \tander^N \left\{ \upmu \exp(-\LogDensity) \frac{p_{;\Ent}}{\overline{\varrho}} \GradEnt^i\right\}.
			\end{split}
		\end{align}
	\end{subequations}
In addition, the following identities hold: 
\begin{subequations}
		\begin{align}
			\upmu \Flatcurl \left( \tander^N v\right)^i & = \upmu \exp(\LogDensity)  \tander^N \vortrenormalized^i + \, ^{(\tander^N)}\ehcurrent_{(v;2)}^i \label{E:TOPORDERCOMMUTEDCURLOFV} \\
			\upmu \Flatdiv \left( \tander^N\vortrenormalized\right) & = - \vortrenormalized^a \upmu \p_a \tander^N \LogDensity + \, ^{(\tander^N)}\ehcurrent_{(\vortrenormalized;3)}, \label{E:TOPORDERCOMMUTEDDIVERGENCEOFVORT} \\			
			\begin{split} \label{E:TOPORDERCOMMUTEDCURLOFCURL} 
				\upmu \left( \Flatcurl \tander^N \vortrenormalized\right)^i & = \upmu \exp(\LogDensity) \tander^N\VortVort^i + \upmu \exp(-2\LogDensity) \Speed^{-2} \frac{p_{;\Ent}}{\overline{\varrho}} \left( \p_a \tander^N v^a\right) \GradEnt^i \\
				& - \upmu \exp(-2\LogDensity) \Speed^{-2} \frac{p_{;\Ent}}{\overline{\varrho}} \GradEnt^a \p_a \tander^N v^i + \, ^{(\tander^N)}\ehcurrent_{(\vortrenormalized;2)}^i ,
			\end{split} \\
			\begin{split}  \label{E:TOPORDERCOMMUTEDDIVOFGRADENT}
				\upmu \Flatdiv\left( \tander^N \GradEnt\right) & = \upmu \exp(2\LogDensity) \tander^N \DivGradEnt + \upmu S^a\p_a \tander^N \LogDensity + \, ^{(\tander^N)}\ehcurrent_{(\GradEnt;3)},
			\end{split}
		\end{align}
	\end{subequations}
where  $^{(\tander^N)}\ehcurrent_{(\vortrenormalized;2)}^i$ is as in \eqref{E:SHARPDECOMPOSITIONOFTOPORDEREXTERIORDERIVATIVEOFVORTSECONDERRORTERM} and: 
	\begin{align}
			 \, ^{(\tander^N)}\ehcurrent_{(v;2)}^i & = [\upmu\Flatcurl,\tander^N] v^i  + [\tander^N,\upmu] v^i + \sum_{\substack{\tander^{N_1}\tander^{N_2} = \tander^N \\ N_1 \le N-1}} \upmu (\tander^{N_1} \vortrenormalized^i) \tander^{N_2}(\exp(\LogDensity)) \label{E:TOPORDERCURLOFVERRORTERM} \\
			^{(\tander^N)}\ehcurrent_{(\vortrenormalized;3)} & =  [\upmu \p_a,\tander^N]\vortrenormalized^a + \vortrenormalized^a [\upmu\p_a,\tander^N]\LogDensity  - \sum_{\substack{\tander^{N_1}\tander^{N_2}=\tander^N \\ N_2 \le N - 1}} \left(\tander^{N_1}\vortrenormalized^a\right) \tander^{N_2}\left(\upmu \p_a \LogDensity\right)\label{E:TOPORDERCOMMUTEDDIVERGENCEOFVORTERRORTERM} \\
		\begin{split}
			^{(\tander^N)}\ehcurrent_{(\GradEnt;3)} & = [\upmu \p_a, \tander^N] \GradEnt^a + \sum_{\substack{\tander^{N_1}\tander^{N_2}\tander^{N_3}=\tander^N \\ N_3 \le N-1}} \left( \tander^{N_1}\upmu\right) \left(\tander^{N_2} \exp(2\LogDensity)\right) \tander^{N_3} \DivGradEnt \\
			& \ \ - S^a[\upmu\p_a,\tander^N ]\LogDensity + \sum_{\substack{\tander^{N_1}\tander^{N_2}=\tander^N \\ N_2 \le N - 1}} \left( \tander^{N_1}\GradEnt^a\right)\tander^{N_2}\left(\upmu\p_a \LogDensity\right). 
			\label{E:SHARPDECOMPOSITIONOFTOPORDERDIVGRADENTERRORTERM}
		\end{split}
	\end{align}
\end{lemma}

\begin{proof}
All identities except for \eqref{E:TOPORDERCOMMUTEDCURLOFV} were proved in \cite{abbrescia2025remarkable}*{Prop. 9.4} with $\upmu$ taking the role of $\weight_{(2)}$. The identity \eqref{E:TOPORDERCOMMUTEDCURLOFV} follows immediately from differentiating the definition of the specific vorticity \eqref{E:SPECIFICVORTICITYDEF}, the Leibniz rule, and straightforward commutations.

\end{proof}

\subsection{Contraction of the characteristic current with $\newuL$} \label{SS:CONTRACTIONOFCHARACTERISTICCURRENTSANDNEWUL}

In this section we use the results of Sects.\,\ref{SS:GEOMETRICSTRUCTUREOFEXTERIORDERIVATIVES}--\ref{SS:STRUCTUREOFHIGHERORDERFLUIDCOMMUTATIONS} to prove that the contraction of the characteristic currents $\ehcurrent^\alpha[\tander^N \vortrenormalized,\pmb{\p} \tander^N \vortrenormalized]$ and $\ehcurrent^\alpha[\tander^N \GradEnt ,\pmb{\p}\tander^N \GradEnt]$ against $\newuL$ produces only top order derivatives of the velocity with $\ingoingcharacteristicsurfacetwoarg{\ubar}{[\moreinterestingu_1,u]}$-tangential derivatives as the last differential operator (e.g. $\newuL \tander^N v$).

\subsubsection{Contraction of the flat current with $\newuL$} \label{SSS:CONTRACTIONOFFLATCURRENTSANDNEWUL}
We begin with the following two lemmas, which reveal the hidden structure between the contraction of $\newuL_\alpha$ with the flat current for top order vorticity and gradient entropy. Crucially, the perfect $\newuL$ derivatives in the lemma feature a truly remarkable overall positive sign.

\begin{lemma}[Remarkable structure in the contraction of the flat current with $\newuL$ for top-order vorticity] \label{L:CONTRACTIONOFFLATCURRENTFORVORTANDNEWUL}
Let $\flatcurrent[\SigmatTan,\pmb{\p}\SigmatTan]$ denote the flat current from \eqref{E:FLATELLIPTICHYPERBOLICCURRENT}, where $\SigmatTan$ is a $\Sigma_t$-tangent vectorfield. Then, the following identities hold:
	\begin{align}
		\begin{split} \label{E:CONTRACTIONOFFLATCURRENTFORVORTANDNEWUL}
			\flatcurrent^\alpha [\tander^N \vortrenormalized,\pmb{\p}\tander^N\vortrenormalized] \newuL_\alpha & = \frac{1}{2} \newuL \left( |\tander^N \vortrenormalized|_{\gfour}^2\right)+  \, ^{(\tander^N \vortrenormalized)}\Currentboundaryerrorhavetocontrolprincipal^{(\textnormal{Flat})} + \, ^{(\tander^N \vortrenormalized)}\Currentboundaryerrorhavetocontrollowerorder^{(\textnormal{Flat})},
		\end{split}
	\end{align}
	where:
	\begin{align}
		\begin{split} \label{E:PRINCIPALERRORTERMHAVETOCONTROLKEYIDPUTANGENTCURRENTCONTRACTEDAGAINSTVECTORFIELDFORVORTICITY}
			\, ^{(\tander^N \vortrenormalized)}\Currentboundaryerrorhavetocontrolprincipal^{(\textnormal{Flat})} & =  \Speed^{-2} \exp(\LogDensity) \upepsilon_{\beta\gamma\delta\kappa} (\tander^N \vortrenormalized^\beta )\newuL^\gamma \Transport^\delta \tander^N \VortVort^\kappa \\
			& \ \ + \frac{1}{\MagnitueofinnerproductofnewLandnewuL} \Speed^{-4} \exp(-2\LogDensity) \frac{p_{;\Ent}}{\overline{\varrho}} \upepsilon_{\beta\gamma\delta i}  \Transport^\delta S^a \newuL^\gamma (\tander^N \vortrenormalized^\beta) \newL_i \newuL  \tander^N v^a \\
			& \ \ 	- \Speed^{-4} \exp(-2\LogDensity) \frac{p_{;\Ent}}{\overline{\varrho}} \upepsilon_{\beta\gamma\delta i}    \Transport^\delta S^a \newuL^\gamma  (\tander^N \vortrenormalized^\beta) \nullangpartial_i\tander^N v^a \\
			& \ \   + \frac{1}{\ingoingmu} (\tander^N \vortrenormalized^\beta) (\vortrenormalized_a) \newL_\beta \newuL \tander^N v^a   - \frac{1}{\upmu} (\tander^N \vortrenormalized^\alpha)\newuL_\alpha \vortrenormalized_a \newuL \tander^N v^a  \\
			& \ \ + (\tander^N \vortrenormalized^\alpha)\newuL_\alpha \vortrenormalized_a  \breve{\slashed{\Transport}} \tander^N v^a - \ReciprocaluLunitAppliedtoTimeFunction (\tander^N \vortrenormalized^\beta) (\vortrenormalized_a) \nullangpartial_\beta \tander^N v^a \\
			& \ \ - \frac{1}{\MagnitueofinnerproductofnewLandnewuL} \Speed^{-4} \exp(-2\LogDensity) \frac{p_{;\Ent}}{\overline{\varrho}}\upepsilon_{\beta\gamma \delta \kappa} (\tander^N \vortrenormalized^\beta)  (\newL^\delta \newuL \tander^N \LogDensity )\GradEnt^\kappa \newuL^\gamma \\
			& \ \ + \Speed^{-2} \exp(-2\LogDensity) \frac{p_{;\Ent}}{\overline{\varrho}} \upepsilon_{\beta\gamma i \kappa} (\tander^N \vortrenormalized^\beta) (\nullangpartial_i \tander^N \LogDensity )\GradEnt^\kappa \newuL^\gamma \\
			& \ \  -  \Speed^{-4} \exp(-2\LogDensity) \frac{p_{;\Ent}}{\overline{\varrho}} \upepsilon_{\beta\gamma\delta\kappa}   (\tander^N \vortrenormalized^\beta)  \Transport^\delta \GradEnt^\kappa (\breve{\slashed{\Transport}} \tander^N \LogDensity)\newuL^\gamma, 
		\end{split}
	\end{align}
	and:

	\begin{align}
		\begin{split} \label{E:LOWERORDERERRORTERMHAVETOCONTROLKEYIDPUTANGENTCURRENTCONTRACTEDAGAINSTVECTORFIELDFORVORTICITY}
			\, ^{(\tander^N \vortrenormalized)}\Currentboundaryerrorhavetocontrollowerorder^{(\textnormal{Flat})} & = (\tander^N \vortrenormalized^\alpha) (\tander^N \vortrenormalized^\beta) \newuL \gfour_{\alpha\beta} - (\tander^N \vortrenormalized^\beta) (\tander^N \vortrenormalized^\gamma) \newuL^\alpha \p_\beta \gfour_{\alpha\gamma} \\
			& \ \ 	+  2  (\newuL \ln \Speed)(\tander^N \vortrenormalized^\beta)(\tander^N \vortrenormalized)_\beta - 2 (\tander^N \vortrenormalized^\beta) (\partial_\beta \ln \Speed)(\tander^N \vortrenormalized)_\gamma \newuL^\gamma  \\
			& \ \  + (\tander^N \vortrenormalized^\beta) \Transport_\gamma \newuL^\gamma (\tander^N \vortrenormalized)_a \p_\beta v^a \\
			& \ \ + \frac{1}{\upmu} \left\{ \updelta_\beta^0 \gfour_{\gamma d} - \updelta_\gamma^0 \gfour_{\beta i}\right\}(\tander^N \vortrenormalized^\beta)  \newuL^\gamma \, ^{(\tander^N)} \ehcurrent_{(\vortrenormalized;1)}^i  + \frac{1}{\upmu} \Speed^{-2} \upepsilon_{\beta\gamma\delta i} (\tander^N \vortrenormalized^\beta)  \newuL^\gamma \Transport^\delta \,^{(\tander^N)}\ehcurrent_{(\vortrenormalized;2)}^i \\
			& \ \ + \frac{1}{\upmu} (\tander^N \vortrenormalized^\alpha)\newuL_\alpha \vortrenormalized_a  \, ^{(\tander^N)}\ehcurrent_{(v;1)}^a  - \frac{1}{\upmu} (\tander^N \vortrenormalized^\alpha)\newuL_\alpha \, ^{(\tander^N)}\ehcurrent_{(\vortrenormalized;3)} \\
			& \ \ + \frac{1}{\upmu}  \Speed^{-4} \exp(-2\LogDensity) \frac{p_{;\Ent}}{\overline{\varrho}} \upepsilon_{\beta\gamma\delta\kappa} \Transport^\delta \GradEnt^\kappa  (\tander^N \vortrenormalized^\beta)  \, ^{(\tander^N)}\ehcurrent_{(\LogDensity)}  \newuL^\gamma \\
			& \ \ - \frac{1}{\upmu} \Speed^{-4} \exp(-2\LogDensity) \frac{p_{;\Ent}}{\overline{\varrho}} \upepsilon_{\beta\gamma i \kappa}  (\tander^N \vortrenormalized^\beta)\, ^{(\tander^N)}\ehcurrent_{(v;1)}^i  \GradEnt^\kappa  \newuL^\gamma \\
			& \ \  -  \Speed^{-4} \exp(-\LogDensity) \frac{p_{;\Ent}}{\overline{\varrho}}\upepsilon_{\beta\gamma\delta i} \upepsilon_{aij} \Transport^\delta S^a  \newuL^\gamma (\tander^N \vortrenormalized^\beta)  \tander^N \vortrenormalized^j \\
			& \ \ - \frac{1}{\upmu} \Speed^{-4} \exp(-2\LogDensity) \frac{p_{;\Ent}}{\overline{\varrho}}  \upepsilon_{\beta\gamma\delta i} \upepsilon_{aij}  \Transport^\delta S^a  \newuL^\gamma( \tander^N \vortrenormalized^\beta ) \, ^{(\tander^N)}\ehcurrent_{(v;2)}^j
		\end{split}
	\end{align}

\end{lemma}
\begin{proof}
	Throughout this proof, we often silently view $\Sigma_t$-tangent vectorfields $\SigmatTan^a$ as spacetime vectorfields $\SigmatTan^\alpha$ where $\alpha = 0$.
	
	We begin by writing:
	\begin{align} \label{E:CONTRACTIONOFFLATCURRENTFORVORTANDNEWULSTEP1}
		\p_\beta \tander^N \vortrenormalized^\alpha = (\p_\beta (\tander^N \vortrenormalized)_\gamma) (\gfour^{-1})^{\alpha\gamma} - (\tander^N \vortrenormalized^\gamma) (\gfour^{-1})^{\alpha \delta}  \p_\beta \gfour_{\delta\gamma}.
	\end{align}
	 Applying \eqref{E:SHARPDECOMPOSITIONOFTOPORDEREXTERIORDERIVATIVEOFVORT} to switch the $\beta-\gamma$ indeces in the first term on RHS\,\eqref{E:CONTRACTIONOFFLATCURRENTFORVORTANDNEWULSTEP1}, inserting the resulting expression into the first term of the defining identity for the flat current \eqref{E:FLATELLIPTICHYPERBOLICCURRENT}, contracting with $\newuL_\alpha$, and relabling the indeces, we find:
	 \begin{align}
	 	\begin{split} \label{E:CONTRACTIONOFFLATCURRENTFORVORTANDNEWULSTEP2}
			\flatcurrent^\alpha [\tander^N \vortrenormalized,\pmb{\p}\tander^N\vortrenormalized] \newuL_\alpha & = (\tander^N \vortrenormalized^\alpha ) \newuL (\tander^N \vortrenormalized)_\alpha - (\tander^N \vortrenormalized^\beta) \Transport_\beta \newuL \tander^N v^a + (\tander^N \vortrenormalized^\beta) \Transport_\gamma \newuL^\gamma (\vortrenormalized)_a \p_\beta \tander^N v^a \\
			& \ \ + \Speed^{-4} \exp(-2\LogDensity) \frac{p_{;\Ent}}{\overline{\varrho}} \upepsilon_{\beta\gamma\delta\kappa} (\tander^N \vortrenormalized^\beta)  \left\{ -(\Transport \tander^N v^\delta)\GradEnt^\kappa + \Transport^\delta \GradEnt^\kappa (\p_a \tander^N v^a) - \Transport^\delta S^a\p_a \tander^N v^\kappa \right\} \newuL^\gamma \\
			& \ \ + \Speed^{-2} \exp(\LogDensity) \upepsilon_{\beta\gamma\delta\kappa} (\tander^N \vortrenormalized^\beta )\newuL^\gamma \Transport^\delta \tander^N \VortVort^\kappa - (\tander^N \vortrenormalized^\alpha) \p_\beta (\tander^N \vortrenormalized^\beta) \newuL_\alpha\\
			& \ \ 
			- (\tander^N \vortrenormalized^\beta) (\tander^N \vortrenormalized^\gamma) \newuL^\alpha \p_\beta \gfour_{\alpha\gamma} 
			+  2  (\newuL \ln \Speed)(\tander^N \vortrenormalized^\beta)(\tander^N \vortrenormalized)_\beta - 2 (\tander^N \vortrenormalized^\beta) (\partial_\beta \ln \Speed)(\tander^N \vortrenormalized)_\gamma \newuL^\gamma  \\
			& \ \ - (\tander^N \vortrenormalized^\beta)   \Transport_\beta (\tander^N \vortrenormalized)_a \newuL  v^a + (\tander^N \vortrenormalized^\beta) \Transport_\gamma \newuL^\gamma (\tander^N \vortrenormalized)_a \p_\beta v^a \\
			& \ \ + \frac{1}{\upmu} \left\{ \updelta_\beta^0 \gfour_{\gamma d} - \updelta_\gamma^0 \gfour_{\beta i}\right\}(\tander^N \vortrenormalized^\beta)  \newuL^\gamma \, ^{(\tander^N)} \ehcurrent_{(\vortrenormalized;1)}^i  + \frac{1}{\upmu} \Speed^{-2} \upepsilon_{\beta\gamma\delta i} (\tander^N \vortrenormalized^\beta)  \newuL^\gamma \Transport^\delta \,^{(\tander^N)}\ehcurrent_{(\vortrenormalized;2)}^i,
		\end{split}
	\end{align}
	where $\, ^{(\tander^N)} \ehcurrent_{(\vortrenormalized;1)}^i$ and $\,^{(\tander^N)}\ehcurrent_{(\vortrenormalized;2)}^i$ are as in \eqref{E:SHARPDECOMPOSITIONOFTOPORDEREXTERIORDERIVATIVEOFVORTFIRSTERRORTERM}--\eqref{E:SHARPDECOMPOSITIONOFTOPORDEREXTERIORDERIVATIVEOFVORTSECONDERRORTERM}. 
	Writing the first term on RHS\,\eqref{E:CONTRACTIONOFFLATCURRENTFORVORTANDNEWULSTEP2} as:
	\begin{align} \label{E:CONTRACTIONOFFLATCURRENTFORVORTANDNEWULSTEP3}
		(\tander^N \vortrenormalized^\alpha ) \newuL (\tander^N \vortrenormalized)_\alpha  & = \frac{1}{2} \newuL \left( |\tander^N \vortrenormalized|_{\gfour}^2\right) + (\tander^N \vortrenormalized^\alpha) (\tander^N \vortrenormalized^\beta) \newuL \gfour_{\alpha\beta},
	\end{align}
	the last term on RHS\,\eqref{E:CONTRACTIONOFFLATCURRENTFORVORTANDNEWULSTEP3} and the last three lines on the RHS\,\eqref{E:CONTRACTIONOFFLATCURRENTFORVORTANDNEWULSTEP2} (except for $- (\tander^N \vortrenormalized^\beta)   \Transport_\beta (\tander^N \vortrenormalized)_a \newuL  v^a = 0$, which follows since $\tander^N \vortrenormalized$ is $\Sigma_t$-tangent) are present on the RHS 
\eqref{E:LOWERORDERERRORTERMHAVETOCONTROLKEYIDPUTANGENTCURRENTCONTRACTEDAGAINSTVECTORFIELDFORVORTICITY}. The first term in \eqref{E:CONTRACTIONOFFLATCURRENTFORVORTANDNEWULSTEP3} is present in \eqref{E:CONTRACTIONOFFLATCURRENTFORVORTANDNEWUL}. Similarly, we note that the second term on RHS\,\eqref{E:CONTRACTIONOFFLATCURRENTFORVORTANDNEWULSTEP2} vanishes as $\tander^N \vortrenormalized$ is $\Sigma_t$-tangent. 

Next, we analyze the sum of the third term on first line and the second term in the third line in the RHS\,\eqref{E:CONTRACTIONOFFLATCURRENTFORVORTANDNEWULSTEP2}:
\begin{align}
		 (\tander^N \vortrenormalized^\beta) \Transport_\gamma \newuL^\gamma (\vortrenormalized)_a \p_\beta \tander^N v^a  - (\tander^N \vortrenormalized^\alpha) \p_\beta (\tander^N \vortrenormalized^\beta) \newuL_\alpha. \label{E:CONTRACTIONOFFLATCURRENTFORVORTANDNEWULSTEP4}
	\end{align}
We first use \eqref{E:RELATIONBETWEENCARTESIANNORMALIZEDNULLVECTORFIELDSANDEIKONALFUNCTIONORMALIZEDNULLVECTORFIELDS}--\eqref{E:EQUIVALENTLUNITANDULUNITAPPLIEDTOEIKONALANDCARTESIANTIME} and \eqref{E:DECOMPOOFPARTIALALPHAINTODOUBLENULLFRAME} to write the first term in \eqref{E:CONTRACTIONOFFLATCURRENTFORVORTANDNEWULSTEP4} as: 
	\begin{align}
		\frac{\ReciprocaluLunitAppliedtoTimeFunction}{\MagnitueofinnerproductofnewLandnewuL} (\tander^N \vortrenormalized^\beta) (\vortrenormalized_a) \newL_\beta \newuL \tander^N v^a + \frac{\ReciprocaluLunitAppliedtoTimeFunction}{\MagnitueofinnerproductofnewLandnewuL}  (\tander^N \vortrenormalized^\beta) (\vortrenormalized_a) \newuL_\beta \newL \tander^N v^a - \ReciprocaluLunitAppliedtoTimeFunction (\tander^N \vortrenormalized^\beta) (\vortrenormalized_a) \nullangpartial_\beta \tander^N v^a. \label{E:CONTRACTIONOFFLATCURRENTFORVORTANDNEWULSTEP5} 
	\end{align}
Next, since $\tander^N \vortrenormalized$ is $\Sigma_t$-tangent and so $\p_\beta (\tander^N \vortrenormalized^\beta) = \p_a (\tander^N \vortrenormalized^a)$, we express the second factor on RHS\,\eqref{E:CONTRACTIONOFFLATCURRENTFORVORTANDNEWULSTEP4} using \eqref{E:TOPORDERCOMMUTEDDIVERGENCEOFVORT} as $(\tander^N \vortrenormalized^\alpha) (\vortrenormalized^a \p_a \tander^N \LogDensity - \frac{1}{\upmu}  \, ^{(\tander^N)}\ehcurrent_{(\vortrenormalized;3)})\newuL_\alpha $, where $ \, ^{(\tander^N)}\ehcurrent_{(\vortrenormalized;3)}$ is as in \eqref{E:TOPORDERCOMMUTEDDIVERGENCEOFVORTERRORTERM}. We then use $\Speed^{-2} \vortrenormalized^a = \vortrenormalized_a$ (see \eqref{E:TOPORDERCOMMUTEDTRANSPORTEQUATIONFORVELOCITY}) and use \eqref{E:BINTERMSOFDOUBLENULLFRAME}, \eqref{E:TOPORDERCOMMUTEDTRANSPORTEQUATIONFORVELOCITY} to write: 
	\begin{align}
		\begin{split} \label{E:CONTRACTIONOFFLATCURRENTFORVORTANDNEWULSTEP6}
			- (\tander^N \vortrenormalized^\alpha) \p_\beta (\tander^N \vortrenormalized^\beta) \newuL_\alpha & = (\tander^N \vortrenormalized^\alpha)\newuL_\alpha \vortrenormalized^a \p_a \tander^N \LogDensity - \frac{1}{\upmu} (\tander^N \vortrenormalized^\alpha)\newuL_\alpha \, ^{(\tander^N)}\ehcurrent_{(\vortrenormalized;3)}\\
			& = - (\tander^N \vortrenormalized^\alpha)\newuL_\alpha \vortrenormalized_a \Transport \tander^N v^a + \frac{1}{\upmu} (\tander^N \vortrenormalized^\alpha)\newuL_\alpha \vortrenormalized_a  \, ^{(\tander^N)}\ehcurrent_{(v;1)}^a  - \frac{1}{\upmu} (\tander^N \vortrenormalized^\alpha)\newuL_\alpha \, ^{(\tander^N)}\ehcurrent_{(\vortrenormalized;3)} \\
			& =  - \frac{1}{\upmu} (\tander^N \vortrenormalized^\alpha)\newuL_\alpha \vortrenormalized_a \newuL \tander^N v^a -  \frac{1}{\ingoingmu} (\tander^N \vortrenormalized^\alpha)\newuL_\alpha \vortrenormalized_a \newL \tander^N v^a + (\tander^N \vortrenormalized^\alpha)\newuL_\alpha \vortrenormalized_a  \breve{\slashed{\Transport}} \tander^N v^a \\
			&  \ \ + \frac{1}{\upmu} (\tander^N \vortrenormalized^\alpha)\newuL_\alpha \vortrenormalized_a  \, ^{(\tander^N)}\ehcurrent_{(v;1)}^a  - \frac{1}{\upmu} (\tander^N \vortrenormalized^\alpha)\newuL_\alpha \, ^{(\tander^N)}\ehcurrent_{(\vortrenormalized;3)}.
		\end{split}
	\end{align}
Since $\frac{\ReciprocaluLunitAppliedtoTimeFunction}{\MagnitueofinnerproductofnewLandnewuL} = \frac{1}{\ingoingmu}$ (see \eqref{E:RATIOOFNULLGEOSICINNERPRODUCTANDFOLIATIONDENSITY}), we find that there is a complete cancellation of the terms featuring $\newL \tander^N v^a$ in the sum of \eqref{E:CONTRACTIONOFFLATCURRENTFORVORTANDNEWULSTEP5}--\eqref{E:CONTRACTIONOFFLATCURRENTFORVORTANDNEWULSTEP6}, and in total we have proved: 
	\begin{align}
		\begin{split} \label{E:CONTRACTIONOFFLATCURRENTFORVORTANDNEWULSTEP7}
			(\tander^N \vortrenormalized^\beta) \Transport_\gamma \newuL^\gamma (\vortrenormalized)_a \p_\beta \tander^N v^a  - (\tander^N \vortrenormalized^\alpha) \p_\beta (\tander^N \vortrenormalized^\beta) \newuL_\alpha 
			& = \frac{1}{\ingoingmu} (\tander^N \vortrenormalized^\beta) (\vortrenormalized_a) \newL_\beta \newuL \tander^N v^a   - \frac{1}{\upmu} (\tander^N \vortrenormalized^\alpha)\newuL_\alpha \vortrenormalized_a \newuL \tander^N v^a  \\
			& \ \ + (\tander^N \vortrenormalized^\alpha)\newuL_\alpha \vortrenormalized_a  \breve{\slashed{\Transport}} \tander^N v^a - \ReciprocaluLunitAppliedtoTimeFunction (\tander^N \vortrenormalized^\beta) (\vortrenormalized_a) \nullangpartial_\beta \tander^N v^a \\
			&  \ \ + \frac{1}{\upmu} (\tander^N \vortrenormalized^\alpha)\newuL_\alpha \vortrenormalized_a  \, ^{(\tander^N)}\ehcurrent_{(v;1)}^a  - \frac{1}{\upmu} (\tander^N \vortrenormalized^\alpha)\newuL_\alpha \, ^{(\tander^N)}\ehcurrent_{(\vortrenormalized;3)}.
		\end{split}
	\end{align}The first four terms and the last two terms on RHS\,\eqref{E:CONTRACTIONOFFLATCURRENTFORVORTANDNEWULSTEP7} are present in RHS\,\eqref{E:PRINCIPALERRORTERMHAVETOCONTROLKEYIDPUTANGENTCURRENTCONTRACTEDAGAINSTVECTORFIELDFORVORTICITY} and \eqref{E:LOWERORDERERRORTERMHAVETOCONTROLKEYIDPUTANGENTCURRENTCONTRACTEDAGAINSTVECTORFIELDFORVORTICITY}, respectively. This concludes the identification of the principal terms generated from the sum of the third term on first line and the second term in the third line in the RHS\,\eqref{E:CONTRACTIONOFFLATCURRENTFORVORTANDNEWULSTEP2}.
	
We now consider the second line on RHS\,\eqref{E:CONTRACTIONOFFLATCURRENTFORVORTANDNEWULSTEP2}, specifically the first two terms (and drop the factors of $ \Speed^{-4} \exp(-2\LogDensity) \frac{p_{;\Ent}}{\overline{\varrho}}\tander^N \vortrenormalized^\beta $ for now):
	\begin{align} \label{E:CONTRACTIONOFFLATCURRENTFORVORTANDNEWULSTEP8}
		-\upepsilon_{\beta\gamma\delta\kappa}(\Transport \tander^N v^\delta)\GradEnt^\kappa \newuL^\gamma +\upepsilon_{\beta\gamma\delta\kappa} \Transport^\delta \GradEnt^\kappa (\p_a \tander^N v^a) \newuL^\gamma.
	\end{align}
Since $\Transport \tander^N v$ is $\Sigma_t$-tangent, we use \eqref{E:TOPORDERCOMMUTEDTRANSPORTEQUATIONFORVELOCITY} to write the first term as $\upepsilon_{\beta\gamma i \kappa} \Speed^2 (\p_i \tander^N \LogDensity )\GradEnt^\kappa \newuL^\gamma - \frac{1}{\upmu} \upepsilon_{\beta\gamma i \kappa} \, ^{(\tander^N)}\ehcurrent_{(v;1)}^i \GradEnt^\kappa \newuL^\gamma$. We then express first term featuring $\partial_i \tander^N\LogDensity$ in terms of the double-null frame using \eqref{E:DECOMPOOFPARTIALALPHAINTODOUBLENULLFRAME} to find: 
	\begin{align} 
		\begin{split} \label{E:CONTRACTIONOFFLATCURRENTFORVORTANDNEWULSTEP9}
			-\upepsilon_{\beta\gamma\delta\kappa}(\Transport \tander^N v^\delta)\GradEnt^\kappa \newuL^\gamma  & = \upepsilon_{\beta\gamma i \kappa} \Speed^2 (\p_i \tander^N \LogDensity )\GradEnt^\kappa \newuL^\gamma - \frac{1}{\upmu} \upepsilon_{\beta\gamma i \kappa} \, ^{(\tander^N)}\ehcurrent_{(v;1)}^i \GradEnt^\kappa \newuL^\gamma \\
			& = - \frac{1}{\MagnitueofinnerproductofnewLandnewuL} \upepsilon_{\beta\gamma i \kappa} \Speed^2 (\newL_i \newuL \tander^N \LogDensity )\GradEnt^\kappa \newuL^\gamma - \frac{1}{\MagnitueofinnerproductofnewLandnewuL}  \upepsilon_{\beta\gamma i \kappa} \Speed^2 (\newuL_i \newL \tander^N \LogDensity )\GradEnt^\kappa \newuL^\gamma \\
			& \ \ +  \upepsilon_{\beta\gamma i \kappa} \Speed^2 (\nullangpartial_i \tander^N \LogDensity )\GradEnt^\kappa \newuL^\gamma - \frac{1}{\upmu} \upepsilon_{\beta\gamma i \kappa} \, ^{(\tander^N)}\ehcurrent_{(v;1)}^i \GradEnt^\kappa \newuL^\gamma \\
		\end{split}
	\end{align}
Let us momentarily focus on the term featuring $\newuL_i \newL \tander^N\LogDensity$. Notice now that by \eqref{E:RELATIONBETWEENCARTESIANNORMALIZEDNULLVECTORFIELDSANDEIKONALFUNCTIONORMALIZEDNULLVECTORFIELDS}, \eqref{E:ULUNITINTERMSOFBANDSIGMATTAN}, and $\Transport_\alpha = - \updelta_\alpha^0$, we may write $\newuL_i = \ReciprocaluLunitAppliedtoTimeFunction (\uLunit^\perp)_i$. Since $\uLunit^\perp$ is $\Sigma_t$-tangent, it also follows that $\upepsilon_{\beta\gamma i \kappa} \Speed^2 (\uLunit^\perp)_i = \upepsilon_{\beta\gamma i \kappa}  (\uLunit^\perp)^i = \upepsilon_{\beta\gamma \delta  \kappa} (\uLunit^\perp)^\delta$. Using \eqref{E:ULUNITINTERMSOFBANDSIGMATTAN} once again, and a symmetric argument for the term featuring $\newL_i \newuL \tander^N \LogDensity$, we have that: 
\begin{align} 
		\begin{split} \label{E:CONTRACTIONOFFLATCURRENTFORVORTANDNEWULSTEP10}
			-\upepsilon_{\beta\gamma\delta\kappa}(\Transport \tander^N v^\delta)\GradEnt^\kappa \newuL^\gamma  & = - \frac{1}{\MagnitueofinnerproductofnewLandnewuL} \upepsilon_{\beta\gamma \delta \kappa}  (\newL^\delta \newuL \tander^N \LogDensity )\GradEnt^\kappa \newuL^\gamma 	+ \frac{\ReciprocalLunitAppliedtoTimeFunction}{\MagnitueofinnerproductofnewLandnewuL} \upepsilon_{\beta\gamma \delta \kappa}  (\Transport^\delta \newuL \tander^N \LogDensity )\GradEnt^\kappa \newuL^\gamma \\
			& \ \  - \frac{1}{\MagnitueofinnerproductofnewLandnewuL} \upepsilon_{\beta\gamma \delta \kappa}  (\newuL^\delta \newL \tander^N \LogDensity )\GradEnt^\kappa \newuL^\gamma 	+ \frac{\ReciprocaluLunitAppliedtoTimeFunction}{\MagnitueofinnerproductofnewLandnewuL} \upepsilon_{\beta\gamma \delta \kappa}  (\Transport^\delta \newL \tander^N \LogDensity )\GradEnt^\kappa \newuL^\gamma \\
			& \ \ +  \upepsilon_{\beta\gamma i \kappa} \Speed^2 (\nullangpartial_i \tander^N \LogDensity )\GradEnt^\kappa \newuL^\gamma - \frac{1}{\upmu} \upepsilon_{\beta\gamma i \kappa} \, ^{(\tander^N)}\ehcurrent_{(v;1)}^i \GradEnt^\kappa \newuL^\gamma 
		\end{split}
	\end{align}
We now observe the miraculous cancellation of the third term on RHS\,\eqref{E:CONTRACTIONOFFLATCURRENTFORVORTANDNEWULSTEP10} due to the anty-sytmmetry of $\upepsilon_{\beta\gamma\delta\kappa}$: 
	\begin{align}	\label{E:CONTRACTIONOFFLATCURRENTFORVORTANDNEWULSTEP11}
		- \frac{1}{\MagnitueofinnerproductofnewLandnewuL} \upepsilon_{\beta\gamma \delta \kappa}  (\newuL^\delta \newL \tander^N \LogDensity )\GradEnt^\kappa \newuL^\gamma  = 0. 
	\end{align}We now come to the second term on RHS\,\eqref{E:CONTRACTIONOFFLATCURRENTFORVORTANDNEWULSTEP8}. Using \eqref{E:TOPORDERCOMMUTEDTRANSPORTEQUATIONFORDENSITY}, this term is equivalent to: 
\begin{align} \label{E:CONTRACTIONOFFLATCURRENTFORVORTANDNEWULSTEP12}
		\upepsilon_{\beta\gamma\delta\kappa} \Transport^\delta \GradEnt^\kappa (\p_a \tander^N v^a) \newuL^\gamma & = -\upepsilon_{\beta\gamma\delta\kappa} \Transport^\delta \GradEnt^\kappa (\Transport \tander^N\LogDensity) \newuL^\gamma + \frac{1}{\upmu} \upepsilon_{\beta\gamma\delta\kappa} \Transport^\delta \GradEnt^\kappa  \, ^{(\tander^N)}\ehcurrent_{(\LogDensity)}  \newuL^\gamma,  
	\end{align}
where $ \, ^{(\tander^N)}\ehcurrent_{(\LogDensity)}$ is as in \eqref{E:TOPORDERCOMMUTEDTRANSPORTEQUATIONFORDENSITY}. Using \eqref{E:BINTERMSOFDOUBLENULLFRAME}, it follows that: 
	\begin{align}
		\begin{split} \label{E:CONTRACTIONOFFLATCURRENTFORVORTANDNEWULSTEP13}
			\upepsilon_{\beta\gamma\delta\kappa} \Transport^\delta \GradEnt^\kappa (\p_a \tander^N v^a) \newuL^\gamma & =  -\frac{1}{\ingoingmu} \upepsilon_{\beta\gamma\delta\kappa}  \Transport^\delta \GradEnt^\kappa (\newL \tander^N \LogDensity)\newuL^\gamma  -\frac{1}{\upmu} \upepsilon_{\beta\gamma\delta\kappa}  \Transport^\delta \GradEnt^\kappa (\newuL \tander^N \LogDensity)\newuL^\gamma  \\
			&  \ \ -  \upepsilon_{\beta\gamma\delta\kappa}  \Transport^\delta \GradEnt^\kappa (\breve{\slashed{\Transport}} \tander^N \LogDensity)\newuL^\gamma  + \frac{1}{\upmu} \upepsilon_{\beta\gamma\delta\kappa} \Transport^\delta \GradEnt^\kappa  \, ^{(\tander^N)}\ehcurrent_{(\LogDensity)}  \newuL^\gamma. 
		\end{split}
	\end{align}
Since $\frac{\ReciprocaluLunitAppliedtoTimeFunction}{\MagnitueofinnerproductofnewLandnewuL} = \frac{1}{\ingoingmu}$, it follows that there is a miraculous cancellation of the terms featuring $\Transport^\delta \newL \tander^N \LogDensity$ in \eqref{E:CONTRACTIONOFFLATCURRENTFORVORTANDNEWULSTEP10} and \eqref{E:CONTRACTIONOFFLATCURRENTFORVORTANDNEWULSTEP13}. Although not crucial, there is also a cancellation of the terms featuring $\Transport^\delta \newuL \tander^N \LogDensity$. In total, by summing 
		
In total, from \eqref{E:CONTRACTIONOFFLATCURRENTFORVORTANDNEWULSTEP8}--\eqref{E:CONTRACTIONOFFLATCURRENTFORVORTANDNEWULSTEP13}, we have the following identity:
\begin{align} \label{E:CONTRACTIONOFFLATCURRENTFORVORTANDNEWULSTEP14}
		\begin{split}
			-\upepsilon_{\beta\gamma\delta\kappa}(\Transport \tander^N v^\delta)\GradEnt^\kappa \newuL^\gamma +\upepsilon_{\beta\gamma\delta\kappa} \Transport^\delta \GradEnt^\kappa (\p_a \tander^N v^a) \newuL^\gamma & =  - \frac{1}{\MagnitueofinnerproductofnewLandnewuL} \upepsilon_{\beta\gamma \delta \kappa}  (\newL^\delta \newuL \tander^N \LogDensity )\GradEnt^\kappa \newuL^\gamma +  \upepsilon_{\beta\gamma i \kappa} \Speed^2 (\nullangpartial_i \tander^N \LogDensity )\GradEnt^\kappa \newuL^\gamma \\
			& \ \  -  \upepsilon_{\beta\gamma\delta\kappa}  \Transport^\delta \GradEnt^\kappa (\breve{\slashed{\Transport}} \tander^N \LogDensity)\newuL^\gamma  + \frac{1}{\upmu} \upepsilon_{\beta\gamma\delta\kappa} \Transport^\delta \GradEnt^\kappa  \, ^{(\tander^N)}\ehcurrent_{(\LogDensity)}  \newuL^\gamma \\
			& \ \ - \frac{1}{\upmu} \upepsilon_{\beta\gamma i \kappa} \, ^{(\tander^N)}\ehcurrent_{(v;1)}^i \GradEnt^\kappa \newuL^\gamma,
		\end{split}
	\end{align}
where the first three terms on RHS\,\eqref{E:CONTRACTIONOFFLATCURRENTFORVORTANDNEWULSTEP14} feature only $\ingoingcharacteristicsurfacetwoarg{\ubar}{[\moreinterestingu_1,u]}$-tangential derivatives as the last differential operator acting on $\tander^N v^a$ or $\tander^N \LogDensity$. Taking into account the factors of  $ \Speed^{-4} \exp(-2\LogDensity) \frac{p_{;\Ent}}{\overline{\varrho}}\tander^N \vortrenormalized^\beta $ present  on the second line of RHS\,\eqref{E:CONTRACTIONOFFLATCURRENTFORVORTANDNEWULSTEP2}, all of the terms in \eqref{E:CONTRACTIONOFFLATCURRENTFORVORTANDNEWULSTEP14} are featured in \eqref{E:PRINCIPALERRORTERMHAVETOCONTROLKEYIDPUTANGENTCURRENTCONTRACTEDAGAINSTVECTORFIELDFORVORTICITY}--\eqref{E:LOWERORDERERRORTERMHAVETOCONTROLKEYIDPUTANGENTCURRENTCONTRACTEDAGAINSTVECTORFIELDFORVORTICITY}.

The last remaining term in \eqref{E:CONTRACTIONOFFLATCURRENTFORVORTANDNEWULSTEP2} is $ -\upepsilon_{\beta\gamma\delta\kappa}  \Transport^\delta S^a\p_a \tander^N v^\kappa \newuL^\gamma = - \upepsilon_{\beta\gamma\delta i}\Transport^\delta S^a \p_a \tander^N v^i$, where we have dropped momentarily the factors of $ \Speed^{-4} \exp(-2\LogDensity) \frac{p_{;\Ent}}{\overline{\varrho}}\tander^N \vortrenormalized^\beta $. Writing $\p_a \tander^N v^i = \p_i \tander^N v^a + \upepsilon_{aij} \Flatcurl(\tander^N v)^j$ and using \eqref{E:TOPORDERCOMMUTEDCURLOFV}, we have that this term satisfies:
	\begin{align}
		\begin{split} \label{E:CONTRACTIONOFFLATCURRENTFORVORTANDNEWULSTEP15}
			 -\upepsilon_{\beta\gamma\delta\kappa}  \Transport^\delta S^a\p_a \tander^N v^\kappa \newuL^\gamma & = - \upepsilon_{\beta\gamma\delta i}  \Transport^\delta S^a \newuL^\gamma  \p_i \tander^N v^a -\upepsilon_{\beta\gamma\delta i} \upepsilon_{aij}  \Transport^\delta S^a  \newuL^\gamma  \Flatcurl(\tander^N v)^j  \\
			 & = - \upepsilon_{\beta\gamma\delta i}  \Transport^\delta S^a \newuL^\gamma  \p_i \tander^N v^a - \upepsilon_{\beta\gamma\delta i} \upepsilon_{aij}  \exp(\LogDensity) \Transport^\delta S^a  \newuL^\gamma \tander^N \vortrenormalized^j \\
			 & \ \ - \frac{1}{\upmu}  \upepsilon_{\beta\gamma\delta i} \upepsilon_{aij}  \Transport^\delta S^a  \newuL^\gamma \, ^{(\tander^N)}\ehcurrent_{(v;2)}^i
		\end{split}
	\end{align}
The last two terms on RHS\,\eqref{E:CONTRACTIONOFFLATCURRENTFORVORTANDNEWULSTEP15} are present on RHS\,\eqref{E:LOWERORDERERRORTERMHAVETOCONTROLKEYIDPUTANGENTCURRENTCONTRACTEDAGAINSTVECTORFIELDFORVORTICITY} upon taking to account the additional factors of $\Speed^{-4} \exp(-2\LogDensity) \frac{p_{;\Ent}}{\overline{\varrho}}\tander^N \vortrenormalized^\beta$. 

We now focus on the first term on RHS\,\eqref{E:CONTRACTIONOFFLATCURRENTFORVORTANDNEWULSTEP15}. Using \eqref{E:DECOMPOOFPARTIALALPHAINTODOUBLENULLFRAME}, this can be expressed as:
	\begin{align}
		\begin{split} \label{E:CONTRACTIONOFFLATCURRENTFORVORTANDNEWULSTEP16}
			- \upepsilon_{\beta\gamma\delta i}  \Transport^\delta S^a \newuL^\gamma  \p_i \tander^N v^a 
			& = \frac{1}{\MagnitueofinnerproductofnewLandnewuL} \upepsilon_{\beta\gamma\delta i}  \Transport^\delta S^a \newuL^\gamma \newuL_i \newL  \tander^N v^a
			+ \frac{1}{\MagnitueofinnerproductofnewLandnewuL} \upepsilon_{\beta\gamma\delta i}  \Transport^\delta S^a \newuL^\gamma \newL_i \newuL  \tander^N v^a \\
			& \ \ -  \upepsilon_{\beta\gamma\delta i}  \Transport^\delta S^a \newuL^\gamma \nullangpartial_i\tander^N v^a.
		\end{split}
	\end{align} 
Next we note the following series of identities:
	\begin{align} \label{E:CONTRACTIONOFFLATCURRENTFORVORTANDNEWULSTEP17} 
		\upepsilon_{\beta\gamma\delta i} \newuL_i =\ReciprocaluLunitAppliedtoTimeFunction \upepsilon_{\beta\gamma\delta i} (\uLunit^\perp)_i =  \Speed^{-2} \ReciprocaluLunitAppliedtoTimeFunction \upepsilon_{\beta\gamma\delta i} (\uLunit^\perp)^i = \Speed^{-2}\ReciprocaluLunitAppliedtoTimeFunction  \upepsilon_{\beta\gamma\delta \kappa} (\uLunit^\perp)^\kappa,
	\end{align}
which follow from $\newuL = \ReciprocaluLunitAppliedtoTimeFunction\uLunit$, the decomposition $\uLunit = \Transport + \uLunit^\perp$ where $\uLunit^\perp$ is $\Sigma_t$-tangent (recall \eqref{E:ULUNITINTERMSOFBANDSIGMATTAN}) and $\Transport_\alpha = - \updelta_\alpha^0$. Inserting \eqref{E:CONTRACTIONOFFLATCURRENTFORVORTANDNEWULSTEP17} into \eqref{E:CONTRACTIONOFFLATCURRENTFORVORTANDNEWULSTEP16} and using \eqref{E:ULUNITINTERMSOFBANDSIGMATTAN} once again, we have that the first term on RHS\,\eqref{E:CONTRACTIONOFFLATCURRENTFORVORTANDNEWULSTEP16} completely vanishes thanks to the anti-symmetric property of $\upepsilon$:
	\begin{align}
		\upepsilon_{\beta\gamma\delta \kappa} \Transport^\delta \uLunit^\gamma \uLunit^\kappa =  \upepsilon_{\beta\gamma\delta \kappa} \Transport^\delta \uLunit^\gamma \Transport^\kappa =  0. \label{E:CONTRACTIONOFFLATCURRENTFORVORTANDNEWULSTEP18}
	\end{align}
We note that the remaining two terms on RHS\,\eqref{E:CONTRACTIONOFFLATCURRENTFORVORTANDNEWULSTEP16} feature only $\ingoingcharacteristicsurfacetwoarg{\ubar}{[\moreinterestingu_1,u]}$-tangential derivatives as the last differential operator acting on $\tander^N v^a$ and appear on the RHS\,\eqref{E:PRINCIPALERRORTERMHAVETOCONTROLKEYIDPUTANGENTCURRENTCONTRACTEDAGAINSTVECTORFIELDFORVORTICITY} after accounting for the factors of $\Speed^{-4} \exp(-2\LogDensity) \frac{p_{;\Ent}}{\overline{\varrho}}\tander^N \vortrenormalized^\beta$.

\end{proof}

\begin{lemma}[Remarkable structure in the contraction of the flat current with $\newuL$ for top-order gradient entropy] \label{L:CONTRACTIONOFFLATCURRENTFORGRADENTANDNEWUL}
Let $\flatcurrent[\SigmatTan,\pmb{\p}\SigmatTan]$ denote the flat current from \eqref{E:FLATELLIPTICHYPERBOLICCURRENT}, where $\SigmatTan$ is a $\Sigma_t$-tangent vectorfield. Then, the following identities hold:
	\begin{align}
		\begin{split} \label{E:CONTRACTIONOFFLATCURRENTFORGRADENTANDNEWUL}
			\flatcurrent^\alpha [\tander^N \GradEnt,\pmb{\p}\tander^N\GradEnt] \newuL_\alpha & = \frac{1}{2} \newuL \left( |\tander^N \GradEnt|_{\gfour}^2\right)+  \, ^{(\tander^N \GradEnt)}\Currentboundaryerrorhavetocontrolprincipal^{(\textnormal{Flat})} + \, ^{(\tander^N \GradEnt)}\Currentboundaryerrorhavetocontrollowerorder^{(\textnormal{Flat})}
		\end{split}
	\end{align}
	
	\begin{align}
		\begin{split} \label{E:PRINCIPALERRORTERMHAVETOCONTROLKEYIDPUTANGENTCURRENTCONTRACTEDAGAINSTVECTORFIELDFORGRADENT}
			\, ^{(\tander^N \GradEnt)}\Currentboundaryerrorhavetocontrolprincipal^{(\textnormal{Flat})} & = - \exp(2\LogDensity)( \tander^N \GradEnt^\alpha) \newuL _\alpha \tander^N \DivGradEnt \\
			& \ \  -  \frac{\ReciprocaluLunitAppliedtoTimeFunction}{\MagnitueofinnerproductofnewLandnewuL} (\tander^N \GradEnt^\beta) \GradEnt_a \newL_\beta \newuL \tander^N v^a +  \ReciprocaluLunitAppliedtoTimeFunction (\tander^N \GradEnt^\beta) \GradEnt_a \nullangpartial_\beta \tander^N v^a  \\
			& \ \   + \frac{1}{\upmu}(\tander^N \GradEnt^\alpha) \newuL_\alpha \GradEnt_a \newuL \tander^N v^a -   (\tander^N \GradEnt^\alpha) \newuL_\alpha \GradEnt_a  \breve{\slashed{\Transport}} \tander^N v^a.
		\end{split}
	\end{align}	\begin{align}
		\begin{split} \label{E:LOWERORDERERRORTERMHAVETOCONTROLKEYIDPUTANGENTCURRENTCONTRACTEDAGAINSTVECTORFIELDFORGRADENT}
			\, ^{(\tander^N \GradEnt)}\Currentboundaryerrorhavetocontrollowerorder^{(\textnormal{Flat})} & = (\tander^N \GradEnt^\alpha) (\tander^N \GradEnt^\beta) \newuL \gfour_{\alpha\beta}	- (\tander^N \GradEnt^\beta) (\tander^N \GradEnt^\gamma) \newuL^\alpha \p_\beta \gfour_{\alpha\gamma} 
			\\
			& \ \ +  2  (\newuL \ln \Speed)(\tander^N \GradEnt^\beta)(\tander^N \GradEnt)_\beta - 2 (\tander^N \GradEnt^\beta) (\partial_\beta \ln \Speed)(\tander^N \GradEnt)_\gamma \newuL^\gamma  \\
			& \ \ + (\tander^N \GradEnt^\beta) \Transport_\gamma \newuL^\gamma (\tander^N \GradEnt)_a \p_\beta v^a \\
			& \ \ + \frac{1}{\upmu} \left\{ \updelta_\beta^0 \gfour_{\gamma d} - \updelta_\gamma^0 \gfour_{\beta i}\right\}(\tander^N \GradEnt^\beta)  \newuL^\gamma \, ^{(\tander^N)} \ehcurrent_{(\GradEnt;1)}^i  + \frac{1}{\upmu} \Speed^{-2} \upepsilon_{\beta\gamma\delta i} (\tander^N \GradEnt^\beta)  \newuL^\gamma \Transport^\delta \,^{(\tander^N)}\ehcurrent_{(\GradEnt;2)}^i \\
			& \ \ - \frac{1}{\upmu} ( \tander^N \GradEnt^\alpha) \newuL _\alpha \, ^{(\tander^N)}\ehcurrent_{(\GradEnt;3)} - \frac{1}{\upmu} (\tander^N \GradEnt^\alpha) \newuL_\alpha S_a\, ^{(\tander^N)}\ehcurrent_{(v;1)}^a.
		\end{split}
	\end{align}

\end{lemma}

\begin{proof}
	Throughout this proof, we often silently view $\Sigma_t$-tangent vectorfields $\SigmatTan^a$ as spacetime vectorfields $\SigmatTan^\alpha$ where $\alpha = 0$.
	
	We begin by writing: 
	\begin{align} \label{E:CONTRACTIONOFFLATCURRENTFORGRADENTANDNEWULSTEP1}
		\p_\beta \tander^N \GradEnt^\alpha = (\p_\beta (\tander^N \GradEnt)_\gamma) (\gfour^{-1})^{\alpha\gamma} - (\tander^N \GradEnt^\gamma) (\gfour^{-1})^{\alpha \delta}  \p_\beta \gfour_{\delta\gamma}.
	\end{align}
	 Applying \eqref{E:SHARPDECOMPOSITIONOFTOPORDEREXTERIORDERIVATIVEOFVORT} to switch the $\beta-\gamma$ indeces in the first term on RHS\,\eqref{E:CONTRACTIONOFFLATCURRENTFORVORTANDNEWULSTEP1}, inserting the resulting expression into the first term of the defining identity for the flat current \eqref{E:FLATELLIPTICHYPERBOLICCURRENT}, contracting with $\newuL_\alpha$, and relabling the indeces, we find:
	 \begin{align}
	 	\begin{split} \label{E:CONTRACTIONOFFLATCURRENTFORGRADENTANDNEWULSTEP2}
			\flatcurrent^\alpha [\tander^N \GradEnt,\pmb{\p}\tander^N\GradEnt] \newuL_\alpha & = (\tander^N \GradEnt^\alpha ) \newuL (\tander^N \GradEnt)_\alpha + (\tander^N \GradEnt^\beta) \Transport_\beta \newuL \tander^N v^a \\
			& \ \ - (\tander^N \GradEnt^\beta) \Transport_\gamma \newuL^\gamma (\GradEnt)_a \p_\beta \tander^N v^a - (\tander^N \GradEnt^\alpha) \p_\beta (\tander^N \GradEnt^\beta) \newuL_\alpha\\
			& \ \ 
			- (\tander^N \GradEnt^\beta) (\tander^N \GradEnt^\gamma) \newuL^\alpha \p_\beta \gfour_{\alpha\gamma} 
			+  2  (\newuL \ln \Speed)(\tander^N \GradEnt^\beta)(\tander^N \GradEnt)_\beta - 2 (\tander^N \GradEnt^\beta) (\partial_\beta \ln \Speed)(\tander^N \GradEnt)_\gamma \newuL^\gamma  \\
			& \ \ - (\tander^N \GradEnt^\beta)   \Transport_\beta (\tander^N \GradEnt)_a \newuL  v^a + (\tander^N \GradEnt^\beta) \Transport_\gamma \newuL^\gamma (\tander^N \GradEnt)_a \p_\beta v^a \\
			& \ \ + \frac{1}{\upmu} \left\{ \updelta_\beta^0 \gfour_{\gamma d} - \updelta_\gamma^0 \gfour_{\beta i}\right\}(\tander^N \GradEnt^\beta)  \newuL^\gamma \, ^{(\tander^N)} \ehcurrent_{(\GradEnt;1)}^i  + \frac{1}{\upmu} \Speed^{-2} \upepsilon_{\beta\gamma\delta i} (\tander^N \GradEnt^\beta)  \newuL^\gamma \Transport^\delta \,^{(\tander^N)}\ehcurrent_{(\GradEnt;2)}^i,
		\end{split}
	\end{align}
where $ \, ^{(\tander^N)} \ehcurrent_{(\GradEnt;1)}^i$ and $ \,^{(\tander^N)}\ehcurrent_{(\GradEnt;2)}^i$ are as in \eqref{E:SHARPDECOMPOSITIONOFTOPORDEREXTERIORDERIVATIVEOFGRADENTFIRSTERRORTERM}--\eqref{E:SHARPDECOMPOSITIONOFTOPORDEREXTERIORDERIVATIVEOFGRADENTSECONDERRORTERM}

Writing the first term on RHS\,\eqref{E:CONTRACTIONOFFLATCURRENTFORGRADENTANDNEWULSTEP2} as:
	\begin{align} \label{E:CONTRACTIONOFFLATCURRENTFORGRADENTANDNEWULSTEP3}
		(\tander^N \GradEnt^\alpha ) \newuL (\tander^N \GradEnt)_\alpha  & = \frac{1}{2} \newuL \left( |\tander^N \GradEnt|_{\gfour}^2\right) + (\tander^N \GradEnt^\alpha) (\tander^N \GradEnt^\beta) \newuL \gfour_{\alpha\beta},
	\end{align}
	the last term on RHS\,\eqref{E:CONTRACTIONOFFLATCURRENTFORGRADENTANDNEWULSTEP3} and the last three lines on the RHS\,\eqref{E:CONTRACTIONOFFLATCURRENTFORGRADENTANDNEWULSTEP2} (except for $(\tander^N \GradEnt^\beta)   \Transport_\beta (\tander^N \GradEnt)_a \newuL  v^a = 0$, which follows since $\tander^N \GradEnt$ is $\Sigma_t$-tangent) are present on the RHS 
\eqref{E:LOWERORDERERRORTERMHAVETOCONTROLKEYIDPUTANGENTCURRENTCONTRACTEDAGAINSTVECTORFIELDFORGRADENT}. The first term in \eqref{E:CONTRACTIONOFFLATCURRENTFORGRADENTANDNEWULSTEP3} is present in \eqref{E:CONTRACTIONOFFLATCURRENTFORGRADENTANDNEWUL}. Similarly, we note that the second term on RHS\,\eqref{E:CONTRACTIONOFFLATCURRENTFORVORTANDNEWULSTEP2} vanishes as $\tander^N \vortrenormalized$ is $\Sigma_t$-tangent. 

We now focus on the fourth term on RHS\,\eqref{E:CONTRACTIONOFFLATCURRENTFORGRADENTANDNEWULSTEP2}. Since $\tander^N \GradEnt$ is $\Sigma_t$-tangent, we may write $\p_\beta \tander^N \GradEnt^\beta = \p_a \tander^N \GradEnt^a$. Inserting the identity to \eqref{E:TOPORDERCOMMUTEDDIVOFGRADENT}, we have that the remaining terms on RHS\,\eqref{E:CONTRACTIONOFFLATCURRENTFORGRADENTANDNEWULSTEP2} are: 
\begin{align} \label{E:CONTRACTIONOFFLATCURRENTFORGRADENTANDNEWULSTEP4}
		- (\tander^N \GradEnt^\beta) \Transport_\gamma \newuL^\gamma (\GradEnt)_a \p_\beta \tander^N v^a - \exp(2\LogDensity)( \tander^N \GradEnt^\alpha) \newuL _\alpha \tander^N \DivGradEnt - ( \tander^N \GradEnt^\alpha) \newuL _\alpha S^a\p_a\tander^N \LogDensity - \frac{1}{\upmu} ( \tander^N \GradEnt^\alpha) \newuL _\alpha \, ^{(\tander^N)}\ehcurrent_{(\GradEnt;3)},
	\end{align}
where $\, ^{(\tander^N)}\ehcurrent_{(\GradEnt;3)}$ is as in \eqref{E:SHARPDECOMPOSITIONOFTOPORDERDIVGRADENTERRORTERM}. The second and last terms on RHS\,\eqref{E:CONTRACTIONOFFLATCURRENTFORGRADENTANDNEWULSTEP4} is present in RHS\,\eqref{E:PRINCIPALERRORTERMHAVETOCONTROLKEYIDPUTANGENTCURRENTCONTRACTEDAGAINSTVECTORFIELDFORGRADENT}--\eqref{E:LOWERORDERERRORTERMHAVETOCONTROLKEYIDPUTANGENTCURRENTCONTRACTEDAGAINSTVECTORFIELDFORGRADENT}, respectively.

We hence focus on the remaining terms: 
\begin{align} \label{E:CONTRACTIONOFFLATCURRENTFORGRADENTANDNEWULSTEP5}
		- (\tander^N \GradEnt^\beta) \Transport_\gamma \newuL^\gamma (\GradEnt)_a \p_\beta \tander^N v^a - ( \tander^N \GradEnt^\alpha) \newuL _\alpha S^a\p_a\tander^N \LogDensity.
	\end{align}
Since $\Transport_\gamma \newuL^\gamma = - \ReciprocaluLunitAppliedtoTimeFunction$ by \eqref{E:EQUIVALENTLUNITANDULUNITAPPLIEDTOEIKONALANDCARTESIANTIME}, we use \eqref{E:TOPORDERCOMMUTEDTRANSPORTEQUATIONFORVELOCITY} to express \eqref{E:CONTRACTIONOFFLATCURRENTFORGRADENTANDNEWULSTEP5} as:
\begin{align}
		\ReciprocaluLunitAppliedtoTimeFunction (\tander^N \GradEnt^\beta) S_a \p_\beta \tander^N v^a +\Speed^{-2} (\tander^N \GradEnt^\alpha) \newuL_\alpha S^a \Transport \tander^N v^a - \frac{1}{\upmu} \Speed^{-2} (\tander^N \GradEnt^\alpha) \newuL_\alpha S^a\, ^{(\tander^N)}\ehcurrent_{(v;1)}^a,  \label{E:CONTRACTIONOFFLATCURRENTFORGRADENTANDNEWULSTEP6}
	\end{align}
where $ \, ^{(\tander^N)}\ehcurrent_{(v;1)}^a$ is as in \eqref{E:TOPORDERCOMMUTEDTRANSPORTEQUATIONFORVELOCITYERRORTERM}. Since $\Speed^{-2} \GradEnt^a = \GradEnt_a$, the last term on RHS\,\eqref{E:CONTRACTIONOFFLATCURRENTFORGRADENTANDNEWULSTEP6} is present in RHS\,\eqref{E:LOWERORDERERRORTERMHAVETOCONTROLKEYIDPUTANGENTCURRENTCONTRACTEDAGAINSTVECTORFIELDFORGRADENT}.
	 
Once again using $\Speed^{-2} \GradEnt^a = \GradEnt_a$, identities \eqref{E:DECOMPOOFPARTIALALPHAINTODOUBLENULLFRAME} and \eqref{E:BINTERMSOFDOUBLENULLFRAME} imply that the first two terms on RHS\,\eqref{E:CONTRACTIONOFFLATCURRENTFORGRADENTANDNEWULSTEP6}, we have:
\begin{align}
		\begin{split} \label{E:CONTRACTIONOFFLATCURRENTFORGRADENTANDNEWULSTEP7}
			& \ReciprocaluLunitAppliedtoTimeFunction (\tander^N \GradEnt^\beta) S_a \p_\beta \tander^N v^a +\Speed^{-2} (\tander^N \GradEnt^\alpha) \newuL_\alpha S^a \Transport \tander^N v^a \\
			& \ \  = - \frac{\ReciprocaluLunitAppliedtoTimeFunction}{\MagnitueofinnerproductofnewLandnewuL} (\tander^N \GradEnt^\beta) \GradEnt_a \newuL_\beta \newL \tander^N v^a -  \frac{\ReciprocaluLunitAppliedtoTimeFunction}{\MagnitueofinnerproductofnewLandnewuL} (\tander^N \GradEnt^\beta) \GradEnt_a \newL_\beta \newuL \tander^N v^a +  \ReciprocaluLunitAppliedtoTimeFunction (\tander^N \GradEnt^\beta) \GradEnt_a \nullangpartial_\beta \tander^N v^a  \\
			& \ \ \ \ + \frac{1}{\ingoingmu}(\tander^N \GradEnt^\alpha) \newuL_\alpha \GradEnt_a \newL \tander^N v^a + \frac{1}{\upmu}(\tander^N \GradEnt^\alpha) \newuL_\alpha \GradEnt_a \newuL \tander^N v^a -   (\tander^N \GradEnt^\alpha) \newuL_\alpha \GradEnt_a  \breve{\slashed{\Transport}} \tander^N v^a.
		\end{split}
	\end{align}
Since $ \frac{\ReciprocaluLunitAppliedtoTimeFunction}{\MagnitueofinnerproductofnewLandnewuL} = \frac{1}{\ingoingmu}$ by \eqref{E:RATIOOFNULLGEOSICINNERPRODUCTANDFOLIATIONDENSITY}, we see a complete cancellation of the two terms on RHS\,\eqref{E:CONTRACTIONOFFLATCURRENTFORGRADENTANDNEWULSTEP7} featuring $\newL \tander^N v^a$. The remaining terms on RHS\,\eqref{E:CONTRACTIONOFFLATCURRENTFORGRADENTANDNEWULSTEP7} all  feature only $\ingoingcharacteristicsurfacetwoarg{\ubar}{[\moreinterestingu_1,u]}$-tangential derivatives as the last differential operator acting on $\tander^N v^a$ and are present on RHS\,\eqref{E:PRINCIPALERRORTERMHAVETOCONTROLKEYIDPUTANGENTCURRENTCONTRACTEDAGAINSTVECTORFIELDFORGRADENT}
\end{proof}

\subsubsection{Contraction of the non-flat component of the characteristic current with $\newuL$} \label{SSS:CONTRATIONOFNONFLATPARTSOFCHARACTERISTICCURRENTWITHNEWUL}

In this section we compute the contraction of the non-flat components of the characteristic current with $\newuL$, namely $(\ehcurrent[\SigmatTan,\pmb{\p}\SigmatTan] - \flatcurrent[\SigmatTan,\pmb{\p}\SigmatTan])^\alpha \newuL_\alpha$.

	\begin{lemma}[Preliminary structure for the non-flat component of $\ehcurrent$ contracted with $\newuL$] \label{L:PRELIMINARYSTRUCTUREFORNONFLATCURRENTCONTRACTEDWITHNEWUL}
	The following identity holds for $\Sigma_t$-tangent vectorfields $\SigmatTan$:
		\begin{align}
			(\ehcurrent[\SigmatTan,\pmb{\p}\SigmatTan] - \flatcurrent[\SigmatTan,\pmb{\p}\SigmatTan])^\alpha \newuL_\alpha & =  \frac{\ingoingmu}{\MagnitueofinnerproductofnewLandnewuL} \SigmatTan^\alpha \newuL_\alpha (\newuL \SigmatTan^\beta)  \breve{\slashed{\Transport}}_\beta - \frac{\ingoingmu}{\MagnitueofinnerproductofnewLandnewuL} \SigmatTan^\beta \breve{\slashed{\Transport}}_\beta (\newuL \SigmatTan^\alpha)\newuL_\alpha. \label{E:PRELIMINARYSTRUCTUREFORNONFLATCURRENTCONTRACTEDWITHNEWUL}
		\end{align} 

	\end{lemma}
	\begin{proof}
		The following identity follows immediately from $\newuL^\alpha \newuL_\alpha = 0$:
		\begin{align} \label{E:PRELIMINARYSTRUCTUREFORNONFLATCURRENTCONTRACTEDWITHNEWULSTEP1}
			(\ehcurrent[\SigmatTan,\pmb{\p}\SigmatTan] - \flatcurrent[\SigmatTan,\pmb{\p}\SigmatTan])^\alpha \newuL_\alpha & =  \frac{1}{\MagnitueofinnerproductofnewLandnewuL} \SigmatTan^\beta \newL_\beta (\newuL \SigmatTan^\alpha)\newuL_\alpha -  \frac{1}{\MagnitueofinnerproductofnewLandnewuL} \SigmatTan^\alpha \newuL_\alpha  (\newuL \SigmatTan^\beta )\newL_\beta.
		\end{align}
		By the $\Sigma_t$-tangency of $\SigmatTan$, $0 = \SigmatTan^\alpha \Transport_\alpha = (\newuL \SigmatTan^\alpha)\Transport_\alpha$ and \eqref{E:BINTERMSOFDOUBLENULLFRAME}, \eqref{E:TORITANGENTVECTORFIELDASSOCIATEDTODOUBLENULLISPROJECTIONOFTRANSPORT} imply that RHS\,\eqref{E:PRELIMINARYSTRUCTUREFORNONFLATCURRENTCONTRACTEDWITHNEWULSTEP1} is equal to: 
		\begin{align}
			\begin{split} \label{E:PRELIMINARYSTRUCTUREFORNONFLATCURRENTCONTRACTEDWITHNEWULSTEP2}
				- \frac{\ingoingmu}{\MagnitueofinnerproductofnewLandnewuL} \SigmatTan^\beta \newuL_\beta (\newuL \SigmatTan^\alpha)\newuL_\alpha - \frac{\ingoingmu}{\MagnitueofinnerproductofnewLandnewuL} \SigmatTan^\beta \breve{\slashed{\Transport}}_\beta (\newuL \SigmatTan^\alpha)\newuL_\alpha & + \frac{\ingoingmu}{\MagnitueofinnerproductofnewLandnewuL} \SigmatTan^\alpha \newuL_\alpha (\newuL \SigmatTan^\beta) \newuL_\beta +  \frac{\ingoingmu}{\MagnitueofinnerproductofnewLandnewuL} \SigmatTan^\alpha \newuL_\alpha (\newuL \SigmatTan^\beta)  \breve{\slashed{\Transport}}_\beta \\
				& =  \frac{\ingoingmu}{\MagnitueofinnerproductofnewLandnewuL} \SigmatTan^\alpha \newuL_\alpha (\newuL \SigmatTan^\beta)  \breve{\slashed{\Transport}}_\beta - \frac{\ingoingmu}{\MagnitueofinnerproductofnewLandnewuL} \SigmatTan^\beta \breve{\slashed{\Transport}}_\beta (\newuL \SigmatTan^\alpha)\newuL_\alpha,
			\end{split}
		\end{align}
		as desired.
	\end{proof}
	
	In addition to several integration by parts arguments done in the forthcoming sections, the proof that the RHS\,\eqref{E:PRELIMINARYSTRUCTUREFORNONFLATCURRENTCONTRACTEDWITHNEWUL} will only feature $\ingoingcharacteristicsurfacetwoarg{\ubar}{[\moreinterestingu_1,u]}$-tangential derivatives as the last differential operator acting on $\tander^N v^a$ depends on a careful analysis of the following crucial factor $\SigmatTan^\alpha \newuL_\alpha (\newuL \SigmatTan^\beta) \breve{\slashed{\Transport}}_\beta$ for $\SigmatTan = \tander^N \vortrenormalized$ and $\SigmatTan = \tander^N\GradEnt$ revealed in Lemma\,\ref{L:PRELIMINARYSTRUCTUREFORNONFLATCURRENTCONTRACTEDWITHNEWUL}
	
	\begin{lemma}[Remarkable structure in the $\SigmatTan^\alpha \newuL_\alpha (\newuL \SigmatTan^\beta) \breve{\slashed{\Transport}}_\beta$ error terms] \label{L:MIRACULOUSSTRUCTUREINTHEVORTICITYANDENTROPYBOUNDARYERRORTERMS}
		The following identity holds for the top-order vorticity:
			\begin{align}
				\begin{split} \label{E:MIRACULOUSSTRUCTUREINTHEVORTICITYBOUNDARYERRORTERM}
					(\tander^N \vortrenormalized^\alpha) \newuL_\alpha (\newuL \tander^N \vortrenormalized^\beta) \breve{\slashed{\Transport}}_\beta & = \nullangDiv\left( \frac{1}{2} (\tander^N \vortrenormalized^\alpha \newuL_\alpha)^2 \breve{\slashed{\Transport}}\right) + \,^{(\tander^N \vortrenormalized)}\Currentboundaryerrorhavetocontrolprincipal^{(\textnormal{Non-flat})} + \,^{(\tander^N \vortrenormalized)}\Currentboundaryerrorhavetocontrollowerorder^{(\textnormal{Non-flat})}, 
				\end{split}
			\end{align}
		where: 
			\begin{align}
				\begin{split} \label{E:MIRACULOUSSTRUCTUREINTHEVORTICITYBOUNDARYERRORTERMPRINCIPAL}
					^{(\tander^N \vortrenormalized)}\Currentboundaryerrorhavetocontrolprincipal^{(\textnormal{Non-flat})} & =  \Speed^{-2} \exp(\LogDensity) \upepsilon_{\gamma\beta\delta\kappa} (\tander^N \vortrenormalized^\alpha) \newuL_\alpha \newuL^\gamma \breve{\slashed{\Transport}}^\beta \Transport^\delta \tander^N \VortVort^\kappa \\
					& \ \ -(\tander^N \vortrenormalized^\alpha)\newuL_\alpha  \newuL^\gamma \Transport_\gamma ( \vortrenormalized)_a \breve{\slashed{\Transport}} \tander^N v^a + (\tander^N \vortrenormalized^\alpha)\newuL_\alpha \breve{\slashed{\Transport}}^\beta \Transport_\beta  ( \vortrenormalized)_a \newuL \tander^N v^a \\
				 	& \ \ + \frac{1}{\MagnitueofinnerproductofnewLandnewuL} \Speed^{-4} \exp(-2\LogDensity) \frac{p_{;\Ent}}{\overline{\varrho}}\upepsilon_{\beta\gamma \delta \kappa} (\tander^N \vortrenormalized^\alpha) \newuL_\alpha \GradEnt^\kappa \newuL^\gamma \breve{\slashed{\Transport}}^\beta  \newL^\delta \newuL \tander^N \LogDensity \\
					& \ \ + \Speed^{-2} \exp(-2\LogDensity) \frac{p_{;\Ent}}{\overline{\varrho}} \upepsilon_{\beta\gamma i \kappa} (\tander^N \vortrenormalized^\alpha) \newuL_\alpha \GradEnt^\kappa \newuL^\gamma \breve{\slashed{\Transport}}^\beta \nullangpartial_i \tander^N \LogDensity  \\
					& \ \ - \Speed^{-4} \exp(-2\LogDensity) \frac{p_{;\Ent}}{\overline{\varrho}} \upepsilon_{\beta\gamma\delta\kappa} (\tander^N \vortrenormalized^\alpha) \newuL_\alpha \Transport^\delta  \GradEnt^\kappa \newuL^\gamma \breve{\slashed{\Transport}}^\beta  \breve{\slashed{\Transport}} \tander^N \LogDensity \\
			& \ \ 	+ \frac{1}{\MagnitueofinnerproductofnewLandnewuL}\Speed^{-4} \exp(-2\LogDensity) \frac{p_{;\Ent}}{\overline{\varrho}} \upepsilon_{\beta\gamma\delta i} (\tander^N \vortrenormalized^\alpha) \newuL_\alpha \Transport^\delta S^a \newuL^\gamma \breve{\slashed{\Transport}}^\beta \newL_i \newuL  \tander^N v^a  \\
			& \ \ -  \Speed^{-4} \exp(-2\LogDensity) \frac{p_{;\Ent}}{\overline{\varrho}}\upepsilon_{\beta\gamma\delta i} (\tander^N \vortrenormalized^\alpha) \newuL_\alpha \Transport^\delta S^a \newuL^\gamma \breve{\slashed{\Transport}}^\beta \nullangpartial_i\tander^N v^a	
				\end{split}
			\end{align}
		and: 
			\begin{align}
				\begin{split} \label{E:MIRACULOUSSTRUCTUREINTHEVORTICITYBOUNDARYERRORTERMLOWERORDER}
					^{(\tander^N \vortrenormalized)}\Currentboundaryerrorhavetocontrollowerorder^{(\textnormal{Non-flat})} & =  - (\tander^N \vortrenormalized^\alpha) \newuL_\alpha(\tander^N \vortrenormalized)^\beta \breve{\slashed{\Transport}}^\delta \newuL \gfour_{\beta\delta}	\\
					& \ \   + 2(\tander^N \vortrenormalized^\alpha)\newuL_\alpha \newuL^\gamma (\tander^N\vortrenormalized)_\gamma (\breve{\slashed{\Transport}} \ln \Speed)- 2 (\tander^N \vortrenormalized^\alpha)\newuL_\alpha (\newuL \ln \Speed)\breve{\slashed{\Transport}}^\beta (\tander^N \vortrenormalized)_\beta \\
					& \ \   - (\tander^N \vortrenormalized^\alpha)\newuL_\alpha  \newuL^\gamma \Transport_\gamma (\tander^N \vortrenormalized)_a \breve{\slashed{\Transport}} v^a + (\tander^N \vortrenormalized^\alpha)\newuL_\alpha \breve{\slashed{\Transport}}^\beta \Transport_\beta (\tander^N \vortrenormalized)_a (\newuL v^a) \\
					& \ \ + \frac{1}{\upmu} \left\{ \updelta_\gamma^0 \gfour_{\beta i} - \updelta_\beta^0 \gfour_{\gamma i}\right\}(\tander^N \vortrenormalized^\alpha)\newuL_\alpha \newuL^\gamma \breve{\slashed{\Transport}}^\beta \, ^{(\tander^N)}\ehcurrent_{(\vortrenormalized;1)}^i  + \frac{1}{\upmu} \Speed^{-2} \upepsilon_{\gamma\beta\delta i} (\tander^N \vortrenormalized^\alpha)\newuL_\alpha \newuL^\gamma \breve{\slashed{\Transport}}^\beta \Transport^\delta \, ^{(\tander^N)}\ehcurrent_{(\vortrenormalized;2)}^i	\\
					& \ \ - (\tander^N \vortrenormalized^\alpha \newuL_\alpha)^2 \nullangDiv \breve{\slashed{\Transport}} - (\tander^N \vortrenormalized^\alpha) \newuL_\alpha (\tander^N \vortrenormalized)_\gamma \breve{\slashed{\Transport}} \newuL^\gamma	\\
					& \ \  + \frac{1}{\upmu}\Speed^{-4} \exp(-2\LogDensity) \frac{p_{;\Ent}}{\overline{\varrho}} \upepsilon_{\beta\gamma\delta\kappa} (\tander^N \vortrenormalized^\alpha) \newuL_\alpha \Transport^\delta \GradEnt^\kappa  \, ^{(\tander^N)}\ehcurrent_{(\LogDensity)}  \newuL^\gamma \\
					& \ \  - \frac{1}{\upmu}\Speed^{-4} \exp(-2\LogDensity) \frac{p_{;\Ent}}{\overline{\varrho}} \upepsilon_{\beta\gamma i \kappa} \, ^{(\tander^N)}\ehcurrent_{(v;1)}^i \GradEnt^\kappa \newuL^\gamma \\
					& \ \ - \Speed^{-4} \exp(-\LogDensity) \frac{p_{;\Ent}}{\overline{\varrho}}\upepsilon_{\beta\gamma\delta i} \upepsilon_{aij}  (\tander^N \vortrenormalized^\alpha) \newuL_\alpha  \Transport^\delta S^a  \newuL^\gamma \breve{\slashed{\Transport}}^\beta  \tander^N \vortrenormalized^j  \\
					& \ \   - \frac{1}{\upmu}  \Speed^{-4} \exp(-2\LogDensity) \frac{p_{;\Ent}}{\overline{\varrho}} \upepsilon_{\beta\gamma\delta i} \upepsilon_{aij}  (\tander^N \vortrenormalized^\alpha) \newuL_\alpha \Transport^\delta S^a  \newuL^\gamma \breve{\slashed{\Transport}}^\beta  \, ^{(\tander^N)}\ehcurrent_{(v;2)}^j.					
				\end{split}
			\end{align}Similarly, the following identities hold for the top-order gradient entropy:
	\begin{align}
				\begin{split} \label{E:MIRACULOUSSTRUCTUREINTHEGRADIENTENTROPYBOUNDARYERRORTERM}
					(\tander^N \GradEnt^\alpha) \newuL_\alpha (\newuL \tander^N \GradEnt^\beta) \breve{\slashed{\Transport}}_\beta & = \nullangDiv\left( \frac{1}{2} (\tander^N \GradEnt^\alpha \newuL_\alpha)^2 \breve{\slashed{\Transport}}\right) + \,^{(\tander^N \GradEnt)}\Currentboundaryerrorhavetocontrolprincipal^{(\textnormal{Non-flat})} + \,^{(\tander^N \GradEnt)}\Currentboundaryerrorhavetocontrollowerorder^{(\textnormal{Non-flat})}, 
				\end{split}
			\end{align}
		where: 
			\begin{align}
				\begin{split} \label{E:MIRACULOUSSTRUCTUREINTHEGRADIENTENTROPYBOUNDARYERRORTERMPRINCIPAL}
					^{(\tander^N \GradEnt)}\Currentboundaryerrorhavetocontrolprincipal^{(\textnormal{Non-flat})} & = (\tander^N \GradEnt^\alpha)\newuL_\alpha  \newuL^\gamma \Transport_\gamma ( \GradEnt)_a \breve{\slashed{\Transport}} \tander^N v^a -  (\tander^N \GradEnt^\alpha)\newuL_\alpha \breve{\slashed{\Transport}}^\beta \Transport_\beta  ( \GradEnt)_a \newuL \tander^N v^a 
				\end{split}
			\end{align}
		and: 
			\begin{align}
				\begin{split} \label{E:MIRACULOUSSTRUCTUREINTHEGRADIENTENTROPYBOUNDARYERRORTERMLOWERORDER}
					^{(\tander^N \GradEnt)}\Currentboundaryerrorhavetocontrollowerorder^{(\textnormal{Non-flat})} & = - (\tander^N \GradEnt^\alpha) \newuL_\alpha(\tander^N \GradEnt)^\beta \breve{\slashed{\Transport}}^\delta \newuL \gfour_{\beta\delta} \\
					& \ \ + 2(\tander^N \GradEnt^\alpha)\newuL_\alpha \newuL^\gamma (\tander^N\GradEnt)_\gamma (\breve{\slashed{\Transport}} \ln \Speed)- 2 (\tander^N \GradEnt^\alpha)\newuL_\alpha (\newuL \ln \Speed)\breve{\slashed{\Transport}}^\beta (\tander^N \GradEnt)_\beta \\
					& \ \  - (\tander^N \GradEnt^\alpha)\newuL_\alpha  \newuL^\gamma \Transport_\gamma (\tander^N \GradEnt)_a \breve{\slashed{\Transport}} v^a + (\tander^N \GradEnt^\alpha)\newuL_\alpha \breve{\slashed{\Transport}}^\beta \Transport_\beta (\tander^N \GradEnt)_a (\newuL v^a) \\
					& \ \ + \frac{1}{\upmu} \left\{ \updelta_\gamma^0 \gfour_{\beta i} - \updelta_\beta^0 \gfour_{\gamma i}\right\}(\tander^N \GradEnt^\alpha)\newuL_\alpha \newuL^\gamma \breve{\slashed{\Transport}}^\beta \, ^{(\tander^N)}\ehcurrent_{(\GradEnt;1)}^i  \\
					& \ \  + \frac{1}{\upmu} \Speed^{-2} \upepsilon_{\gamma\beta\delta i} (\tander^N \GradEnt^\alpha)\newuL_\alpha \newuL^\gamma \breve{\slashed{\Transport}}^\beta \Transport^\delta \, ^{(\tander^N)}\ehcurrent_{(\GradEnt;2)}^i \\
					& \ \ - (\tander^N \GradEnt^\alpha \newuL_\alpha)^2 \nullangDiv \breve{\slashed{\Transport}} - (\tander^N \GradEnt^\alpha) \newuL_\alpha (\tander^N \GradEnt)_\gamma \breve{\slashed{\Transport}} \newuL^\gamma.
				\end{split}
			\end{align}
			
	\end{lemma}

\begin{proof}

	Throughout this proof, we often silently view $\Sigma_t$-tangent vectorfields $\SigmatTan^a$ as spacetime vectorfields $\SigmatTan^\alpha$ where $\alpha = 0$.
	
\medskip

\noindent \underline{\textbf{Proof of the remarkable structure for top-order vorticity:}} We begin by writing LHS \eqref{E:MIRACULOUSSTRUCTUREINTHEVORTICITYBOUNDARYERRORTERM} as: 
	\begin{align}
		\begin{split}  \label{E:MIRACULOUSSTRUCTUREINTHEVORTICITYBOUNDARYERRORTERMSTEP1}
			(\tander^N \vortrenormalized^\alpha) \newuL_\alpha (\newuL \tander^N \vortrenormalized^\beta) \breve{\slashed{\Transport}}_\beta & = (\tander^N \vortrenormalized^\alpha)\newuL_\alpha \newuL^\gamma \p_\gamma (\gfour^{\beta\delta}(\tander \vortrenormalized)_\delta) \slashed{\Transport}_\beta \\
			& = (\tander^N \vortrenormalized^\alpha)\newuL_\alpha \newuL^\gamma (\p_\gamma (\tander^N \vortrenormalized)_\beta) \breve{\slashed{\Transport}}^\beta - (\tander^N \vortrenormalized^\alpha) \newuL_\alpha(\tander^N \vortrenormalized)^\beta \breve{\slashed{\Transport}}^\delta \newuL \gfour_{\beta\delta}.
		\end{split}
	\end{align}
The second term on RHS\,\eqref{E:MIRACULOUSSTRUCTUREINTHEVORTICITYBOUNDARYERRORTERMSTEP1} is present in RHS\,\eqref{E:MIRACULOUSSTRUCTUREINTHEVORTICITYBOUNDARYERRORTERMLOWERORDER}. 

We apply \eqref{E:SHARPDECOMPOSITIONOFTOPORDEREXTERIORDERIVATIVEOFVORT} to the first term RHS\,\eqref{E:MIRACULOUSSTRUCTUREINTHEVORTICITYBOUNDARYERRORTERMSTEP1} to find that:
	\begin{align}
		\begin{split}  \label{E:MIRACULOUSSTRUCTUREINTHEVORTICITYBOUNDARYERRORTERMSTEP2}
			& (\tander^N \vortrenormalized^\alpha)\newuL_\alpha \newuL^\gamma (\p_\gamma (\tander^N \vortrenormalized)_\beta) \breve{\slashed{\Transport}}^\beta \\
			& \ \  = (\tander^N \vortrenormalized^\alpha)\newuL_\alpha \newuL^\gamma \breve{\slashed{\Transport}} (\tander^N\vortrenormalized)_\gamma + \Speed^{-2} \exp(\LogDensity) \upepsilon_{\gamma\beta\delta\kappa} (\tander^N \vortrenormalized^\alpha) \newuL_\alpha \newuL^\gamma \breve{\slashed{\Transport}}^\beta \Transport^\delta \tander^N \VortVort^\kappa  \\
			& \ \ \ \   -(\tander^N \vortrenormalized^\alpha)\newuL_\alpha  \newuL^\gamma \Transport_\gamma ( \vortrenormalized)_a \breve{\slashed{\Transport}} \tander^N v^a + (\tander^N \vortrenormalized^\alpha)\newuL_\alpha \breve{\slashed{\Transport}}^\beta \Transport_\beta  ( \vortrenormalized)_a \newuL \tander^N v^a \\
			& \ \ \ \ + \Speed^{-4} \exp(-2\LogDensity) \frac{p_{;\Ent}}{\overline{\varrho}} (\tander^N \vortrenormalized^\alpha)\newuL_\alpha \upepsilon_{\gamma \beta\delta \kappa} \newuL^\gamma \breve{\slashed{\Transport}}^\beta  \left\{  -(\Transport \tander^N v^\delta)\GradEnt^\kappa + \Transport^\beta \GradEnt^\kappa (\p_a \tander^N v^a) - \Transport^\delta S^a\p_a \tander^N v^\kappa \right\} \\
			& \ \  \ \ + 2(\tander^N \vortrenormalized^\alpha)\newuL_\alpha \newuL^\gamma (\tander^N\vortrenormalized)_\gamma (\breve{\slashed{\Transport}} \ln \Speed)- 2 (\tander^N \vortrenormalized^\alpha)\newuL_\alpha (\newuL \ln \Speed)\breve{\slashed{\Transport}}^\beta (\tander^N \vortrenormalized)_\beta \\
			& \ \ \ \  - (\tander^N \vortrenormalized^\alpha)\newuL_\alpha  \newuL^\gamma \Transport_\gamma (\tander^N \vortrenormalized)_a \breve{\slashed{\Transport}} v^a + (\tander^N \vortrenormalized^\alpha)\newuL_\alpha \breve{\slashed{\Transport}}^\beta \Transport_\beta (\tander^N \vortrenormalized)_a (\newuL v^a) \\
			& \ \ \ \ + \frac{1}{\upmu} \left\{ \updelta_\gamma^0 \gfour_{\beta i} - \updelta_\beta^0 \gfour_{\gamma i}\right\}(\tander^N \vortrenormalized^\alpha)\newuL_\alpha \newuL^\gamma \breve{\slashed{\Transport}}^\beta \, ^{(\tander^N)}\ehcurrent_{(\vortrenormalized;1)}^i  + \frac{1}{\upmu} \Speed^{-2} \upepsilon_{\gamma\beta\delta i} (\tander^N \vortrenormalized^\alpha)\newuL_\alpha \newuL^\gamma \breve{\slashed{\Transport}}^\beta \Transport^\delta \, ^{(\tander^N)}\ehcurrent_{(\vortrenormalized;2)}^i, 
		\end{split}
	\end{align}
where $\, ^{(\tander^N)}\ehcurrent_{(\vortrenormalized;1)}^i$ and $\, ^{(\tander^N)}\ehcurrent_{(\vortrenormalized;2)}^i$ are as in \eqref{E:SHARPDECOMPOSITIONOFTOPORDEREXTERIORDERIVATIVEOFVORTFIRSTERRORTERM}--\eqref{E:SHARPDECOMPOSITIONOFTOPORDEREXTERIORDERIVATIVEOFVORTSECONDERRORTERM}. The last three lines on RHS\,\eqref{E:MIRACULOUSSTRUCTUREINTHEVORTICITYBOUNDARYERRORTERMSTEP2} are present in RHS\,\eqref{E:MIRACULOUSSTRUCTUREINTHEVORTICITYBOUNDARYERRORTERMLOWERORDER}. The term featuring $\tander^N \VortVort^\kappa$ as well as the second line on RHS\,\eqref{E:MIRACULOUSSTRUCTUREINTHEVORTICITYBOUNDARYERRORTERMSTEP2} features only top order derivatives of the velocity with $\ingoingcharacteristicsurfacetwoarg{\ubar}{[\moreinterestingu_1,u]}$-tangential derivatives as the last differential operator and are present in RHS\,\eqref{E:MIRACULOUSSTRUCTUREINTHEVORTICITYBOUNDARYERRORTERMPRINCIPAL}. We differentiate the first term on RHS\,\eqref{E:MIRACULOUSSTRUCTUREINTHEVORTICITYBOUNDARYERRORTERMSTEP2} as: 
	\begin{align}
		\begin{split} \label{E:MIRACULOUSSTRUCTUREINTHEVORTICITYBOUNDARYERRORTERMSTEP3}
			(\tander^N \vortrenormalized^\alpha)\newuL_\alpha \newuL^\gamma \breve{\slashed{\Transport}} (\tander^N\vortrenormalized)_\gamma  & = \nullangDiv\left( \frac{1}{2} (\tander^N \vortrenormalized^\alpha \newuL_\alpha)^2 \breve{\slashed{\Transport}}\right) - (\tander^N \vortrenormalized^\alpha \newuL_\alpha)^2 \nullangDiv \breve{\slashed{\Transport}} - (\tander^N \vortrenormalized^\alpha) \newuL_\alpha (\tander^N \vortrenormalized)_\gamma \breve{\slashed{\Transport}} \newuL^\gamma.
		\end{split}
	\end{align}
	The first term on RHS\,\eqref{E:MIRACULOUSSTRUCTUREINTHEVORTICITYBOUNDARYERRORTERMSTEP3} is present on RHS\,\eqref{E:MIRACULOUSSTRUCTUREINTHEVORTICITYBOUNDARYERRORTERM}, while the last two are present on RHS\,\eqref{E:MIRACULOUSSTRUCTUREINTHEVORTICITYBOUNDARYERRORTERMLOWERORDER}.
	
	It remains to handle the error terms featuring $\upepsilon_{\gamma\beta\delta\kappa}\newuL^\gamma \breve{\slashed{\Transport}}^\beta$ on RHS\,\eqref{E:MIRACULOUSSTRUCTUREINTHEVORTICITYBOUNDARYERRORTERMSTEP2}. We highlight that these terms were already analyzed in the large as  they were present on RHS\,\eqref{E:CONTRACTIONOFFLATCURRENTFORVORTANDNEWULSTEP2} with $\tander^N\vortrenormalized^\beta$ in the place of $\breve{\slashed{\Transport}}^\beta$. However, this contraction in the $\beta$ index \emph{played no role in revealing the structure of these terms in Lemma\,\ref{L:CONTRACTIONOFFLATCURRENTFORVORTANDNEWUL}.} Indeed, what was fundamental was that the contraction in the $\gamma$ index was with $\newuL^\gamma$. Hence, using the same analysis as \eqref{E:CONTRACTIONOFFLATCURRENTFORVORTANDNEWULSTEP8}--\eqref{E:CONTRACTIONOFFLATCURRENTFORVORTANDNEWULSTEP18} and momentarily dropping the factors of $\Speed^{-4} \exp(-2\LogDensity) \frac{p_{;\Ent}}{\overline{\varrho}}(\tander^N \vortrenormalized^\alpha) \newuL_\alpha$, it follows that: 
	\begin{align}
		\begin{split} \label{E:MIRACULOUSSTRUCTUREINTHEVORTICITYBOUNDARYERRORTERMSTEP4}
			& \upepsilon_{\gamma \beta\delta \kappa} \newuL^\gamma \breve{\slashed{\Transport}}^\beta  \left\{  -(\Transport \tander^N v^\delta)\GradEnt^\kappa + \Transport^\beta \GradEnt^\kappa (\p_a \tander^N v^a) - \Transport^\delta S^a\p_a \tander^N v^\kappa \right\} \\
			 & = \ \  \frac{1}{\MagnitueofinnerproductofnewLandnewuL} \upepsilon_{\beta\gamma \delta \kappa} \GradEnt^\kappa \newuL^\gamma \breve{\slashed{\Transport}}^\beta  \newL^\delta \newuL \tander^N \LogDensity +  \upepsilon_{\beta\gamma i \kappa} \Speed^2  \GradEnt^\kappa \newuL^\gamma \breve{\slashed{\Transport}}^\beta \nullangpartial_i \tander^N \LogDensity  -  \upepsilon_{\beta\gamma\delta\kappa}  \Transport^\delta  \GradEnt^\kappa \newuL^\gamma \breve{\slashed{\Transport}}^\beta  \breve{\slashed{\Transport}} \tander^N \LogDensity \\
			& \ \ 	+ \frac{1}{\MagnitueofinnerproductofnewLandnewuL} \upepsilon_{\beta\gamma\delta i}  \Transport^\delta S^a \newuL^\gamma \breve{\slashed{\Transport}}^\beta \newL_i \newuL  \tander^N v^a  -  \upepsilon_{\beta\gamma\delta i}  \Transport^\delta S^a \newuL^\gamma \breve{\slashed{\Transport}}^\beta \nullangpartial_i\tander^N v^a \\
			& \ \  + \frac{1}{\upmu} \upepsilon_{\beta\gamma\delta\kappa} \Transport^\delta \GradEnt^\kappa  \, ^{(\tander^N)}\ehcurrent_{(\LogDensity)}  \newuL^\gamma  - \frac{1}{\upmu} \upepsilon_{\beta\gamma i \kappa} \, ^{(\tander^N)}\ehcurrent_{(v;1)}^i \GradEnt^\kappa \newuL^\gamma \\
			& \ \ - \upepsilon_{\beta\gamma\delta i} \upepsilon_{aij}  \exp(\LogDensity) \Transport^\delta S^a  \newuL^\gamma \breve{\slashed{\Transport}}^\beta  \tander^N \vortrenormalized^j   - \frac{1}{\upmu}  \upepsilon_{\beta\gamma\delta i} \upepsilon_{aij}  \Transport^\delta S^a  \newuL^\gamma \breve{\slashed{\Transport}}^\beta  \, ^{(\tander^N)}\ehcurrent_{(v;2)}^i
		\end{split}
	\end{align}
Accounting for the factors of $\Speed^{-4} \exp(-2\LogDensity) \frac{p_{;\Ent}}{\overline{\varrho}}(\tander^N \vortrenormalized^\alpha) \newuL_\alpha$, the first two lines on RHS\,\eqref{E:MIRACULOUSSTRUCTUREINTHEVORTICITYBOUNDARYERRORTERMSTEP4} are present on RHS\,\eqref{E:MIRACULOUSSTRUCTUREINTHEVORTICITYBOUNDARYERRORTERMPRINCIPAL}, while the last two lines are present in \eqref{E:MIRACULOUSSTRUCTUREINTHEVORTICITYBOUNDARYERRORTERMLOWERORDER}.

\medskip

\noindent \underline{\textbf{Proof of the remarkable structure for top-order gradient entropy:}} We begin by writing LHS \eqref{E:MIRACULOUSSTRUCTUREINTHEGRADIENTENTROPYBOUNDARYERRORTERM} as: 
	\begin{align}
		\begin{split}  \label{E:MIRACULOUSSTRUCTUREINTHEGRADENTBOUNDARYERRORTERMSTEP1}
			(\tander^N \GradEnt^\alpha) \newuL_\alpha (\newuL \tander^N \GradEnt^\beta) \breve{\slashed{\Transport}}_\beta & = (\tander^N \GradEnt^\alpha)\newuL_\alpha \newuL^\gamma \p_\gamma (\gfour^{\beta\delta}(\tander \GradEnt)_\delta) \slashed{\Transport}_\beta \\
			& = (\tander^N \GradEnt^\alpha)\newuL_\alpha \newuL^\gamma (\p_\gamma (\tander^N \GradEnt)_\beta) \breve{\slashed{\Transport}}^\beta - (\tander^N \GradEnt^\alpha) \newuL_\alpha(\tander^N \GradEnt)^\beta \breve{\slashed{\Transport}}^\delta \newuL \gfour_{\beta\delta}.
		\end{split}
	\end{align}
The second term on RHS\,\eqref{E:MIRACULOUSSTRUCTUREINTHEVORTICITYBOUNDARYERRORTERMSTEP1} is present in RHS\,\eqref{E:MIRACULOUSSTRUCTUREINTHEVORTICITYBOUNDARYERRORTERMLOWERORDER}. 

We apply \eqref{E:SHARPDECOMPOSITIONOFTOPORDEREXTERIORDERIVATIVEOFGRADENT} to the first term RHS\,\eqref{E:MIRACULOUSSTRUCTUREINTHEGRADENTBOUNDARYERRORTERMSTEP1} to find: 
\begin{align}
		\begin{split}  \label{E:MIRACULOUSSTRUCTUREINTHEGRADENTBOUNDARYERRORTERMSTEP2}
			& (\tander^N \GradEnt^\alpha)\newuL_\alpha \newuL^\gamma (\p_\gamma (\tander^N \GradEnt)_\beta) \breve{\slashed{\Transport}}^\beta \\
			& \ \  = (\tander^N \GradEnt^\alpha)\newuL_\alpha \newuL^\gamma \breve{\slashed{\Transport}} (\tander^N\GradEnt)_\gamma  \\
			& \ \ \ \   + (\tander^N \GradEnt^\alpha)\newuL_\alpha  \newuL^\gamma \Transport_\gamma ( \GradEnt)_a \breve{\slashed{\Transport}} \tander^N v^a -  (\tander^N \GradEnt^\alpha)\newuL_\alpha \breve{\slashed{\Transport}}^\beta \Transport_\beta  ( \GradEnt)_a \newuL \tander^N v^a \\
			& \ \  \ \ + 2(\tander^N \GradEnt^\alpha)\newuL_\alpha \newuL^\gamma (\tander^N\GradEnt)_\gamma (\breve{\slashed{\Transport}} \ln \Speed)- 2 (\tander^N \GradEnt^\alpha)\newuL_\alpha (\newuL \ln \Speed)\breve{\slashed{\Transport}}^\beta (\tander^N \GradEnt)_\beta \\
			& \ \ \ \  - (\tander^N \GradEnt^\alpha)\newuL_\alpha  \newuL^\gamma \Transport_\gamma (\tander^N \GradEnt)_a \breve{\slashed{\Transport}} v^a + (\tander^N \GradEnt^\alpha)\newuL_\alpha \breve{\slashed{\Transport}}^\beta \Transport_\beta (\tander^N \GradEnt)_a (\newuL v^a) \\
			& \ \ \ \ + \frac{1}{\upmu} \left\{ \updelta_\gamma^0 \gfour_{\beta i} - \updelta_\beta^0 \gfour_{\gamma i}\right\}(\tander^N \GradEnt^\alpha)\newuL_\alpha \newuL^\gamma \breve{\slashed{\Transport}}^\beta \, ^{(\tander^N)}\ehcurrent_{(\GradEnt;1)}^i  + \frac{1}{\upmu} \Speed^{-2} \upepsilon_{\gamma\beta\delta i} (\tander^N \GradEnt^\alpha)\newuL_\alpha \newuL^\gamma \breve{\slashed{\Transport}}^\beta \Transport^\delta \, ^{(\tander^N)}\ehcurrent_{(\GradEnt;2)}^i.
		\end{split}
	\end{align}
The second line on RHS\,\eqref{E:MIRACULOUSSTRUCTUREINTHEGRADENTBOUNDARYERRORTERMSTEP2} is present in RHS\,\eqref{E:MIRACULOUSSTRUCTUREINTHEGRADIENTENTROPYBOUNDARYERRORTERMPRINCIPAL}, while the third, fourth, and fifth line on RHS\,\eqref{E:MIRACULOUSSTRUCTUREINTHEGRADENTBOUNDARYERRORTERMSTEP2} are present in \eqref{E:MIRACULOUSSTRUCTUREINTHEGRADIENTENTROPYBOUNDARYERRORTERMLOWERORDER}. We differentiate the first term on RHS\,\eqref{E:MIRACULOUSSTRUCTUREINTHEVORTICITYBOUNDARYERRORTERMSTEP2} as: 
	\begin{align}
		\begin{split} \label{E:MIRACULOUSSTRUCTUREINTHEGRADENTBOUNDARYERRORTERMSTEP3}
			(\tander^N \GradEnt^\alpha)\newuL_\alpha \newuL^\gamma \breve{\slashed{\Transport}} (\tander^N\GradEnt)_\gamma  & = \nullangDiv\left( \frac{1}{2} (\tander^N \GradEnt^\alpha \newuL_\alpha)^2 \breve{\slashed{\Transport}}\right) - (\tander^N \GradEnt^\alpha \newuL_\alpha)^2 \nullangDiv \breve{\slashed{\Transport}} - (\tander^N \GradEnt^\alpha) \newuL_\alpha (\tander^N \GradEnt)_\gamma \breve{\slashed{\Transport}} \newuL^\gamma.
		\end{split}
	\end{align}
	The first term on RHS\,\eqref{E:MIRACULOUSSTRUCTUREINTHEVORTICITYBOUNDARYERRORTERMSTEP3} is present on RHS\,\eqref{E:MIRACULOUSSTRUCTUREINTHEGRADIENTENTROPYBOUNDARYERRORTERM}, while the last two are present on RHS\,\eqref{E:MIRACULOUSSTRUCTUREINTHEGRADIENTENTROPYBOUNDARYERRORTERMLOWERORDER}.

\end{proof}

\subsubsection{The remarkable structure in the ingoing null-flux boundary terms in the elliptic hyperbolic identities} \label{SS:MIRACULOUSSTRUCTUREINTHEBOUNDARYTERM} Here we use the results of the previous subsection to reveal the remarkable structure of the boundary terms in the elliptic hyperbolic integral identities. Specifically, after an integration by parts, the co-dimension 2 integrals feature an overall \emph{good sign}. This is detailed in the following lemma:

\begin{lemma}[Remarkable structure in the $\ingoingcharacteristicsurfacetwoarg{\ubar_2}{[u_1,u_2]}$ boundary terms in the elliptic-hyperbolic integral identities] \label{L:MIRACULOUSSTRUCTUREINTHEBOUNDARYTERM} Let $\SigmatTan = \tander^N \vortrenormalized$ or $\tander^N \GradEnt$. Then, for any $u_1 \le u_2$, the following identity holds:
\begin{align}
		\begin{split} \label{E:MIRACULOUSSTRUCTUREINTHEBOUNDARYTERM}
			- \int_{\ingoingcharacteristicsurfacetwoarg{\ubar}{[u_1,u_2]}} \frac{1}{\upmu} \weight^{\blowupratetoporderacoustic} \newuL_\alpha \ehcurrent^\alpha \volingoingnullhypersurface  & = -\int_{\ingoingcharacteristicsurfacetwoarg{\ubar_1}{[u_1,u_2]}} \left\{ \, ^{(\SigmatTan)}\Currentboundaryerrorhavetocontrolprincipalweighted{\blowupratetoporderacoustic}  + \, ^{(\SigmatTan)}\Currentboundaryerrorhavetocontrollowerorderweighted{\blowupratetoporderacoustic}  \right\} \volingoingnullhypersurface \\
			& \ \ - \int_{\doublenulltoritwoarg{\ubar}{u_2}} \frac{1}{\upmu} \coercivequadraticformintoriboundaryterms[\SigmatTan,\SigmatTan]  \weight^{\blowupratetoporderacoustic} \, \voldoublenulltori 
			+ \int_{\doublenulltoritwoarg{\ubar}{u_1}} \frac{1}{\upmu}\coercivequadraticformintoriboundaryterms[\SigmatTan,\SigmatTan]   \weight^{\blowupratetoporderacoustic} \, \voldoublenulltori,
		\end{split}
	\end{align}
where:
	\begin{align} \label{E:TORIBOUNDARYTERMSFORELLIPTICHYPERBOLICIDENTITIES}
		\coercivequadraticformintoriboundaryterms[\SigmatTan,\SigmatTan]  \eqdef	 \frac{1}{2} |\SigmatTan|_{\gfour}^2 - \frac{\ingoingmu}{\MagnitueofinnerproductofnewLandnewuL}  \SigmatTan^\alpha \newuL_\alpha   \SigmatTan^\beta \breve{\slashed{\Transport}}_\beta 
	\end{align}
	and:
	\begin{subequations}
		\begin{align}
			\, ^{(\SigmatTan)}\Currentboundaryerrorhavetocontrolprincipalweighted{\blowupratetoporderacoustic} & = \left\{ \frac{1}{\upmu}  \, ^{(\SigmatTan)}\Currentboundaryerrorhavetocontrolprincipal^{(\textnormal{Flat})} 
			+ 2\frac{\ingoingmu}{\upmu \MagnitueofinnerproductofnewLandnewuL}   \,^{(\SigmatTan )}\Currentboundaryerrorhavetocontrolprincipal^{(\textnormal{Non-flat})} \right\}  \weight^{\blowupratetoporderacoustic} \label{E:PRINCIPALERRORTERMHAVETOCONTROLKEYIDPUTANGENTCURRENTCONTRACTEDAGAINSTVECTORFIELDFORVORTICITYWEIGHTED} \\
			\begin{split}
				\, ^{(\SigmatTan)}\Currentboundaryerrorhavetocontrollowerorderweighted{\blowupratetoporderacoustic} & =  \left\{\frac{1}{\upmu} \, ^{(\SigmatTan)}\Currentboundaryerrorhavetocontrollowerorder^{(\textnormal{Flat})} 
				+ 2 \frac{\ingoingmu}{\upmu \MagnitueofinnerproductofnewLandnewuL}  \,^{(\SigmatTan)}\Currentboundaryerrorhavetocontrollowerorder^{(\textnormal{Non-flat})}\right\} \weight^{\blowupratetoporderacoustic} \\
				& \ \  + \left\{\frac{\newuL \upmu}{2\upmu^2} |\SigmatTan|_{\gfour}^2
				+ \newuL\left( \frac{1}{\upmu}\right) \SigmatTan^\alpha \uLunit_\alpha \SigmatTan^\beta \breve{\slashed{\Transport}}_\beta 
				+ \breve{\slashed{\Transport}} \left( \frac{\ingoingmu}{\upmu \MagnitueofinnerproductofnewLandnewuL} \right) ( \SigmatTan^\alpha\newuL_\alpha)^2  \right. \\
				 & \ \ \ \ \left. + \frac{1}{\upmu}  \SigmatTan^\alpha \SigmatTan^\beta \newuL(\uLunit_\alpha \breve{\slashed{\Transport}}_\beta)
				 + \frac{1}{2 \upmu}  \SigmatTan^\alpha \uLunit_\alpha   \SigmatTan^\beta \breve{\slashed{\Transport}}_\beta \mytr_{\gnulltori} \deform{\newuL} 
				 -   \frac{1}{4\upmu} |\SigmatTan|_{\gfour}^2 \mytr_{\gnulltori} \deform{\newuL}  \right\} \weight^{\blowupratetoporderacoustic},\label{E:LOWERORDERERRORTERMHAVETOCONTROLKEYIDPUTANGENTCURRENTCONTRACTEDAGAINSTVECTORFIELDFORVORTICITYWEIGHTED}
			\end{split}
		\end{align}
	\end{subequations}
and $^{(\SigmatTan)}\Currentboundaryerrorhavetocontrolprincipal^{(\textnormal{Flat})} , \, ^{(\SigmatTan)}\Currentboundaryerrorhavetocontrolprincipal^{(\textnormal{Non-flat})}, \, ^{(\SigmatTan)}\Currentboundaryerrorhavetocontrollowerorder^{(\textnormal{Flat})}, \, ^{(\SigmatTan)}\Currentboundaryerrorhavetocontrollowerorder^{(\textnormal{Non-flat})}$ are as in Lemmas\,\ref{L:CONTRACTIONOFFLATCURRENTFORVORTANDNEWUL} and\,\ref{L:MIRACULOUSSTRUCTUREINTHEVORTICITYANDENTROPYBOUNDARYERRORTERMS} for $\SigmatTan = \tander^N \vortrenormalized$ and Lemmas\,\ref{L:CONTRACTIONOFFLATCURRENTFORGRADENTANDNEWUL} and\,\ref{L:MIRACULOUSSTRUCTUREINTHEVORTICITYANDENTROPYBOUNDARYERRORTERMS} for $\SigmatTan = \tander^N \GradEnt$.

Moreover, the quadratic form $\coercivequadraticformintoriboundaryterms[\SigmatTan,\SigmatTan]  $ in \eqref{E:TORIBOUNDARYTERMSFORELLIPTICHYPERBOLICIDENTITIES} is coercive and satisfies the following estimate:
	\begin{align} \label{E:KEYCOERCIVITYONTORIBOUNDARYTERMSFORELLIPTICHYPERBOLICIDENTITIES}
		\coercivequadraticformintoriboundaryterms[\SigmatTan,\SigmatTan]  \approx \frac{1}{4} |\SigmatTan|_g^2 + \frac{1}{2} |\SigmatTan|_{\gnulltori}^2.
	\end{align}
\end{lemma}

\begin{proof}

Let $\SigmatTan = \tander^N \vortrenormalized$ or $\tander^N \GradEnt$. Then by the decomposition of $\ehcurrent[\SigmatTan,\pmb{\p}\SigmatTan]$ in its flat and non-flat components \eqref{E:RELATIONSHIPBETWEENFLATANDCHARACTERISTICCURRENTS}, the results of Lemmas\,\ref{L:CONTRACTIONOFFLATCURRENTFORVORTANDNEWUL}--\ref{L:PRELIMINARYSTRUCTUREFORNONFLATCURRENTCONTRACTEDWITHNEWUL} imply that: 
\begin{align}
		\begin{split} \label{E:MIRACULOUSSTRUCTUREINTHEBOUNDARYTERMSTEP1}
			- \int_{\ingoingcharacteristicsurfacetwoarg{\ubar}{[u_1,u_2]}} \frac{1}{\upmu} \newuL_\alpha \ehcurrent^\alpha[\SigmatTan,\pmb{\p}\SigmatTan] \weight^{\blowupratetoporderacoustic} \, \volingoingnullhypersurface 
			 & = - \int_{\ingoingcharacteristicsurfacetwoarg{\ubar}{[u_1,u_2]}} \frac{1}{\upmu} \left\{ \frac{1}{2} \newuL \left( |\SigmatTan|_{\gfour}^2\right) 
			 +  \frac{\ingoingmu}{\MagnitueofinnerproductofnewLandnewuL} \SigmatTan^\alpha \newuL_\alpha (\newuL \SigmatTan^\beta)  \breve{\slashed{\Transport}}_\beta 
			 - \frac{\ingoingmu}{\MagnitueofinnerproductofnewLandnewuL} \SigmatTan^\beta \breve{\slashed{\Transport}}_\beta (\newuL \SigmatTan^\alpha)\newuL_\alpha \right. \\
			& \ \ \left. + \, ^{(\SigmatTan)}\Currentboundaryerrorhavetocontrolprincipal^{(\textnormal{Flat})} +  \, ^{(\SigmatTan)}\Currentboundaryerrorhavetocontrollowerorder^{(\textnormal{Flat})}  \right\} \weight^{\blowupratetoporderacoustic} \, \volingoingnullhypersurface.
		\end{split}
	\end{align}
We focus our attention first on the third term on RHS\,\eqref{E:MIRACULOUSSTRUCTUREINTHEBOUNDARYTERMSTEP1}. Using $\frac{\ingoingmu}{\MagnitueofinnerproductofnewLandnewuL} \newuL_\alpha = \uLunit_\alpha$ (which follows from \eqref{E:RATIOOFNULLGEOSICINNERPRODUCTANDFOLIATIONDENSITY} and \eqref{E:RELATIONBETWEENCARTESIANNORMALIZEDNULLVECTORFIELDSANDEIKONALFUNCTIONORMALIZEDNULLVECTORFIELDS}) Adding and subtracting $\frac{\ingoingmu}{2 \MagnitueofinnerproductofnewLandnewuL}  \SigmatTan^\alpha \newuL_\alpha   \SigmatTan^\beta \breve{\slashed{\Transport}}_\beta 
(\mytr_{\gnulltori} \deform{\newuL}) \weight^{\blowupratetoporderacoustic}$ and then integrating by parts using \eqref{E:NEWULINTEGRALONDOUBLENULLTORI}, we find that: 
\begin{align}
		\begin{split} \label{E:MIRACULOUSSTRUCTUREINTHEBOUNDARYTERMSTEP2}
			 \int_{\ingoingcharacteristicsurfacetwoarg{\ubar}{[u_1,u_2]}} \frac{1}{\upmu}  \frac{\ingoingmu}{\MagnitueofinnerproductofnewLandnewuL} \SigmatTan^\beta \breve{\slashed{\Transport}}_\beta (\newuL \SigmatTan^\alpha)\newuL_\alpha  \weight^{\blowupratetoporderacoustic} \, \volingoingnullhypersurface
			 & = -  \int_{\ingoingcharacteristicsurfacetwoarg{\ubar}{[u_1,u_2]}}   \left\{ \frac{1}{\upmu} \frac{\ingoingmu}{\MagnitueofinnerproductofnewLandnewuL} \SigmatTan^\alpha \newuL_\alpha (\newuL \SigmatTan^\beta)  \breve{\slashed{\Transport}}_\beta  + \newuL\left( \frac{1}{\upmu}\right) \SigmatTan^\alpha \uLunit_\alpha \SigmatTan^\beta \breve{\slashed{\Transport}}_\beta  \right. \\
			 & \ \ \left. + \frac{1}{\upmu}  \SigmatTan^\alpha \SigmatTan^\beta \newuL(\uLunit_\alpha \breve{\slashed{\Transport}}_\beta) + \frac{1}{2 \upmu}  \SigmatTan^\alpha \uLunit_\alpha   \SigmatTan^\beta \breve{\slashed{\Transport}}_\beta \mytr_{\gnulltori} \deform{\newuL}  \right\} \weight^{\blowupratetoporderacoustic} \, \volingoingnullhypersurface \\
			 & \ \ + \int_{\doublenulltoritwoarg{\ubar}{u_2}} \frac{\ingoingmu}{\upmu \MagnitueofinnerproductofnewLandnewuL}  \SigmatTan^\alpha \newuL_\alpha   \SigmatTan^\beta \breve{\slashed{\Transport}}_\beta  \weight^{\blowupratetoporderacoustic}  \, \voldoublenulltori - \int_{\doublenulltoritwoarg{\ubar}{u_1}} \frac{\ingoingmu}{\upmu \MagnitueofinnerproductofnewLandnewuL}  \SigmatTan^\alpha \newuL_\alpha   \SigmatTan^\beta \breve{\slashed{\Transport}}_\beta  \weight^{\blowupratetoporderacoustic}  \, \voldoublenulltori
		\end{split}
	\end{align} 
	
Next, adding and subtracting $\frac{1}{2 \upmu} |\SigmatTan|_{\gfour}^2 (\mytr_{\gnulltori} \deform{\newuL}) \weight^{\blowupratetoporderacoustic}$  to the first error integral on RHS\,\eqref{E:MIRACULOUSSTRUCTUREINTHEBOUNDARYTERMSTEP1} and integrating by parts with \eqref{E:NEWULINTEGRALONDOUBLENULLTORI}, we have that:
\begin{align}
		\begin{split} \label{E:MIRACULOUSSTRUCTUREINTHEBOUNDARYTERMSTEP3}
			- \int_{\ingoingcharacteristicsurfacetwoarg{\ubar}{[u_1,u_2]}} \frac{1}{2\upmu}  \newuL \left( |\SigmatTan|_{\gfour}^2\right)  \weight^{\blowupratetoporderacoustic} \, \volingoingnullhypersurface & = \int_{\ingoingcharacteristicsurfacetwoarg{\ubar}{[u_1,u_2]}} \left\{- \frac{\newuL \upmu}{2\upmu^2} |\SigmatTan|_{\gfour}^2 + \frac{1}{4\upmu} |\SigmatTan|_{\gfour}^2 \mytr_{\gnulltori} \deform{\newuL} \right\}  \weight^{\blowupratetoporderacoustic} \, \volingoingnullhypersurface \\
			& \ \ - \int_{\doublenulltoritwoarg{\ubar}{u_2}} \frac{1}{2\upmu} |\SigmatTan|_{\gfour}^2  \weight^{\blowupratetoporderacoustic} \, \voldoublenulltori
			+ \int_{\doublenulltoritwoarg{\ubar}{u_1}} \frac{1}{2\upmu} |\SigmatTan|_{\gfour}^2  \weight^{\blowupratetoporderacoustic} \, \voldoublenulltori.
		\end{split}
	\end{align}Combining \eqref{E:MIRACULOUSSTRUCTUREINTHEBOUNDARYTERMSTEP1}--\eqref{E:MIRACULOUSSTRUCTUREINTHEBOUNDARYTERMSTEP2},  we see that we have proved: 
\begin{align}
		\begin{split} \label{E:MIRACULOUSSTRUCTUREINTHEBOUNDARYTERMSTEP4}
			& - \int_{\ingoingcharacteristicsurfacetwoarg{\ubar}{[u_1,u_2]}} \frac{1}{\upmu} \newuL_\alpha \ehcurrent^\alpha[\SigmatTan,\pmb{\p}\SigmatTan] \weight^{\blowupratetoporderacoustic} \, \volingoingnullhypersurface \\
			&\ \  = \int_{\ingoingcharacteristicsurfacetwoarg{\ubar}{[u_1,u_2]}}  \left\{- \frac{\newuL \upmu}{2\upmu^2} |\SigmatTan|_{\gfour}^2
			- 2\frac{\ingoingmu}{\upmu \MagnitueofinnerproductofnewLandnewuL} \SigmatTan^\alpha \newuL_\alpha (\newuL \SigmatTan^\beta)  \breve{\slashed{\Transport}}_\beta  
			- \newuL\left( \frac{1}{\upmu}\right) \SigmatTan^\alpha \uLunit_\alpha \SigmatTan^\beta \breve{\slashed{\Transport}}_\beta  \right. \\
			 & \ \ \ \ \left. - \frac{1}{\upmu}  \SigmatTan^\alpha \SigmatTan^\beta \newuL(\uLunit_\alpha \breve{\slashed{\Transport}}_\beta)
			 - \frac{1}{2 \upmu}  \SigmatTan^\alpha \uLunit_\alpha   \SigmatTan^\beta \breve{\slashed{\Transport}}_\beta \mytr_{\gnulltori} \deform{\newuL} 
			 +   \frac{1}{4\upmu} |\SigmatTan|_{\gfour}^2 \mytr_{\gnulltori} \deform{\newuL}  \right\} \weight^{\blowupratetoporderacoustic} \, \volingoingnullhypersurface  \\
			 & \ \ \ \  -  \int_{\ingoingcharacteristicsurfacetwoarg{\ubar}{[u_1,u_2]}}  \left\{\, ^{(\SigmatTan)}\Currentboundaryerrorhavetocontrolprincipal^{(\textnormal{Flat})} +  \, ^{(\SigmatTan)}\Currentboundaryerrorhavetocontrollowerorder^{(\textnormal{Flat})}  \right\} \weight^{\blowupratetoporderacoustic} \, \volingoingnullhypersurface \\
			 & \ \ \ \ - \int_{\doublenulltoritwoarg{\ubar}{u_2}} \frac{1}{\upmu} \left\{ \frac{1}{2} |\SigmatTan|_{\gfour}^2 -\frac{\ingoingmu}{\MagnitueofinnerproductofnewLandnewuL}  \SigmatTan^\alpha \newuL_\alpha   \SigmatTan^\beta \breve{\slashed{\Transport}}_\beta \right\}  \weight^{\blowupratetoporderacoustic} \, \voldoublenulltori\\
			 & \ \ \ \ 
			+ \int_{\doublenulltoritwoarg{\ubar}{u_1}} \frac{1}{\upmu} \left\{ \frac{1}{2} |\SigmatTan|_{\gfour}^2 - \frac{\ingoingmu}{\MagnitueofinnerproductofnewLandnewuL}  \SigmatTan^\alpha \newuL_\alpha   \SigmatTan^\beta \breve{\slashed{\Transport}}_\beta \right\}  \weight^{\blowupratetoporderacoustic} \, \voldoublenulltori.
		\end{split}
	\end{align}

We now focus on the second term on RHS\,\eqref{E:MIRACULOUSSTRUCTUREINTHEBOUNDARYTERMSTEP4}, namely the boundary integral featuring $- 2\frac{\ingoingmu}{\upmu \MagnitueofinnerproductofnewLandnewuL} \SigmatTan^\alpha \newuL_\alpha (\newuL \SigmatTan^\beta)  \breve{\slashed{\Transport}}_\beta$. To handle the $\nullangDiv$ terms featured on RHS\,\eqref{E:MIRACULOUSSTRUCTUREINTHEVORTICITYBOUNDARYERRORTERM} and \eqref{E:MIRACULOUSSTRUCTUREINTHEGRADIENTENTROPYBOUNDARYERRORTERM}, we differentiate by parts:
	\begin{align} 
		- \frac{\ingoingmu}{\upmu \MagnitueofinnerproductofnewLandnewuL} \nullangDiv\left(  (\SigmatTan^\alpha \newuL_\alpha)^2 \breve{\slashed{\Transport}}\right) \weight^{\blowupratetoporderacoustic} = - \nullangDiv\left( \frac{\ingoingmu}{\upmu \MagnitueofinnerproductofnewLandnewuL} (\SigmatTan^\alpha \newuL_\alpha)^2 \weight^{\blowupratetoporderacoustic} \breve{\slashed{\Transport}}\right)+  \breve{\slashed{\Transport}} \left( \frac{\ingoingmu}{\upmu \MagnitueofinnerproductofnewLandnewuL} \right)( \SigmatTan^\alpha\newuL_\alpha)^2\weight^{\blowupratetoporderacoustic} \label{E:MIRACULOUSSTRUCTUREINTHEBOUNDARYTERMSTEP5}
	\end{align}
where the note that $\breve{\slashed{\Transport}} \weight = \breve{\slashed{\Transport}} (-\ubar) = 0$ because $\breve{\slashed{\Transport}}$ is $\doublenulltoritwoarg{\ubar}{u}$-tangent. As the integral over $\ingoingcharacteristicsurfacetwoarg{\ubar}{[u_1,u_2]}$ of the first term on RHS\,\eqref{E:MIRACULOUSSTRUCTUREINTHEBOUNDARYTERMSTEP5}  vanishes, we have that: 
	\begin{align}
		\begin{split} \label{E:MIRACULOUSSTRUCTUREINTHEBOUNDARYTERMSTEP6}
			-2 \int_{\ingoingcharacteristicsurfacetwoarg{\ubar}{[u_1,u_2]}} \frac{\ingoingmu}{\upmu \MagnitueofinnerproductofnewLandnewuL} \SigmatTan^\alpha \newuL_\alpha (\newuL \SigmatTan^\beta)  \breve{\slashed{\Transport}}_\beta  \weight^{\blowupratetoporderacoustic} \, \volingoingnullhypersurface 
			& =  \int_{\ingoingcharacteristicsurfacetwoarg{\ubar}{[u_1,u_2]}}  \left\{ \breve{\slashed{\Transport}} \left( \frac{\ingoingmu}{\upmu \MagnitueofinnerproductofnewLandnewuL} \right) ( \SigmatTan^\alpha\newuL_\alpha)^2  - 2\frac{\ingoingmu}{\upmu \MagnitueofinnerproductofnewLandnewuL}   \,^{(\SigmatTan )}\Currentboundaryerrorhavetocontrolprincipal^{(\textnormal{Non-flat})} -2 \frac{\ingoingmu}{\upmu \MagnitueofinnerproductofnewLandnewuL}  \,^{(\SigmatTan)}\Currentboundaryerrorhavetocontrollowerorder^{(\textnormal{Non-flat})}\right\} \weight^{\blowupratetoporderacoustic} \, \volingoingnullhypersurface. 
		\end{split}
	\end{align}
	
Inserting \eqref{E:MIRACULOUSSTRUCTUREINTHEBOUNDARYTERMSTEP6} into the second integral on RHS\,\eqref{E:MIRACULOUSSTRUCTUREINTHEBOUNDARYTERMSTEP4}, we conclude the proof of \eqref{E:MIRACULOUSSTRUCTUREINTHEBOUNDARYTERM}.

It remains to prove \eqref{E:KEYCOERCIVITYONTORIBOUNDARYTERMSFORELLIPTICHYPERBOLICIDENTITIES}. Since $\SigmatTan^\beta \Transport_\beta = 0$, \eqref{E:BINTERMSOFDOUBLENULLFRAME}--\eqref{E:TORITANGENTVECTORFIELDASSOCIATEDTODOUBLENULLISPROJECTIONOFTRANSPORT} implies that $-\frac{\ingoingmu}{\MagnitueofinnerproductofnewLandnewuL}   \SigmatTan^\beta \breve{\slashed{\Transport}}_\beta = \frac{\ingoingmu}{\upmu \MagnitueofinnerproductofnewLandnewuL}\SigmatTan^\beta\newuL_\beta + \frac{1}{\MagnitueofinnerproductofnewLandnewuL} \SigmatTan^\beta \newL_\beta $. Inserting this identity into RHS\,\eqref{E:TORIBOUNDARYTERMSFORELLIPTICHYPERBOLICIDENTITIES} and using \eqref{E:ACOUSTICALMETRICINTERMSOFNULLVECTORFIELDSANDDOUBLENULLSPHEREFIRSTFUND}, we have that: 
\begin{align}
		\begin{split}
			\coercivequadraticformintoriboundaryterms[\SigmatTan,\SigmatTan] & = - \frac{1}{\MagnitueofinnerproductofnewLandnewuL}\SigmatTan^\alpha \newuL_\alpha \SigmatTan^\beta \newL_\beta + \frac{1}{2} |\SigmatTan|_{\gnulltori}^2 +  \frac{\ingoingmu}{\upmu \MagnitueofinnerproductofnewLandnewuL}(\SigmatTan^\alpha\newuL_\beta)^2 + \frac{1}{\MagnitueofinnerproductofnewLandnewuL} \SigmatTan^\alpha \newuL_\alpha \SigmatTan^\beta \newL_\beta \\
			& =  \frac{\ingoingmu}{\upmu \MagnitueofinnerproductofnewLandnewuL}(\SigmatTan^\alpha\newuL_\alpha)^2 + \frac{1}{2} |\SigmatTan|_{\gnulltori}^2 
		\end{split}
	\end{align}\noindent Finally, using $ \frac{\ingoingmu}{\upmu \MagnitueofinnerproductofnewLandnewuL}(\SigmatTan^\alpha\newuL_\alpha)^2=  \frac{\ReciprocaluLunitAppliedtoTimeFunction}{ \upmu}(\SigmatTan^\alpha \uLunit_\alpha )^2 = \frac{\ReciprocaluLunitAppliedtoTimeFunction}{ \upmu}(\SigmatTan^a \uLunit_a^\perp )^2$ (which follows from the $\Sigma_t$-tangency of $\uLunit^\perp$), the second identity in \eqref{E:SIZEOFLUBITPERPANDULUNITPERPANDTHEIRINNERPRODUCT} and estimate \eqref{E:SHARPESTIMATEFORRATIOOFRATIOOFNULLGEOSICINNERPRODUCTANDFOLIATIONDENSITYANDINGOINGMU} prove \eqref{E:KEYCOERCIVITYONTORIBOUNDARYTERMSFORELLIPTICHYPERBOLICIDENTITIES} upon taking $\fundbootsmall$ sufficiently small.
	
\end{proof}

\subsection{The main elliptic-hyperbolic integral identity} \label{SS:P:INTEGRALIDENTITYFORELLIPTICHYPERBOLICCURRENT} We now state and prove the main elliptic-hyperbolic integral identity.

\begin{proposition}[The main elliptic-hyperbolic integral identity] \label{P:INTEGRALIDENTITYFORELLIPTICHYPERBOLICCURRENT}
	Let $\SigmatTan$ be a $\Sigma_t$-tangent vectorfield and let $\ellipticCoerciveQuadratic[\pmb{\partial}\SigmatTan,\pmb{\partial}\SigmatTan]$ be the quadratic form from Def.~\ref{D:NULLHYPERSURFACEADAPTEDCOERCIVEQUADRATICFORM}. Then for any $\ubar_1 \le \ubar_2$ and $u_1 \le u_2$, the following integral identity holds:
	\begin{align}
		\begin{split} \label{E:INTEGRALIDENTITYFORELLIPTICHYPERBOLICCURRENT}
			&  \int_{\characteristicdiamondtwoarg{[\ubar_1,\ubar_2)}{[u_1,u_2]}} \weight^{\blowupratetoporderacoustic} \ReciprocalLunitAppliedtoTimeFunction \, \ellipticCoerciveQuadratic[\pmb{\partial} \SigmatTan,\pmb{\partial} \SigmatTan]  \, \voldiamond + \int_{\doublenulltoritwoarg{\ubar_2}{u_2}}  \frac{1}{\upmu} \weight^{\blowupratetoporderacoustic} \coercivequadraticformintoriboundaryterms[\SigmatTan,\SigmatTan] \voldoublenulltori \\
			& 
			= \int_{\doublenulltoritwoarg{\ubar_2}{u_1}}  \frac{1}{\upmu} \weight^{\blowupratetoporderacoustic}  \coercivequadraticformintoriboundaryterms[\SigmatTan,\SigmatTan] \voldoublenulltori + \int_{\doublenulltoritwoarg{\ubar_1}{u_2}}  \frac{1}{\upmu} \weight^{\blowupratetoporderacoustic} \coercivequadraticformintoriboundaryterms[\SigmatTan,\SigmatTan] \voldoublenulltori - \int_{\doublenulltoritwoarg{\ubar_1}{u_1}}  \frac{1}{\upmu} \weight^{\blowupratetoporderacoustic} \coercivequadraticformintoriboundaryterms[\SigmatTan,\SigmatTan] \voldoublenulltori \\
			& \ \
			+ \int_{\ingoingcharacteristicsurfacetwoarg{\ubar_1}{[u_1,u_2]}} \left\{ \, ^{(\SigmatTan)}\Currentboundaryerrorhavetocontrolprincipalweighted{\blowupratetoporderacoustic}  + \, ^{(\SigmatTan)}\Currentboundaryerrorhavetocontrollowerorderweighted{\blowupratetoporderacoustic}  \right\} \volingoingnullhypersurface \\
			& \ \ 
			- 
			\int_{\ingoingcharacteristicsurfacetwoarg{\ubar_2}{[u_1,u_2]}} \left\{ \, ^{(\SigmatTan)}\Currentboundaryerrorhavetocontrolprincipalweighted{\blowupratetoporderacoustic}  + ^{(\SigmatTan)}\Currentboundaryerrorhavetocontrollowerorderweighted{\blowupratetoporderacoustic} \right\} \volingoingnullhypersurface + \int_{\characteristicdiamondtwoarg{[\ubar_1,\ubar_2)}{[-\rightu,\leftu]}} \weight^{\blowupratetoporderacoustic} \ReciprocalLunitAppliedtoTimeFunction \, \EllipticHyperbolicCurrentIntegralIdentityTotalSpacetimeErrorTerm[\SigmatTan,\pmb{\partial} \SigmatTan]  \, \voldiamond,
		\end{split}
	\end{align}
	where $\coercivequadraticformintoriboundaryterms[\SigmatTan,\SigmatTan]$ is given by \eqref{E:TORIBOUNDARYTERMSFORELLIPTICHYPERBOLICIDENTITIES}, $\, ^{(\SigmatTan)}\Currentboundaryerrorhavetocontrolprincipalweighted{\blowupratetoporderacoustic}
	$ and $^{(\SigmatTan)}\Currentboundaryerrorhavetocontrollowerorderweighted{\blowupratetoporderacoustic}$ are as in \eqref{E:PRINCIPALERRORTERMHAVETOCONTROLKEYIDPUTANGENTCURRENTCONTRACTEDAGAINSTVECTORFIELDFORVORTICITYWEIGHTED}--\eqref{E:LOWERORDERERRORTERMHAVETOCONTROLKEYIDPUTANGENTCURRENTCONTRACTEDAGAINSTVECTORFIELDFORVORTICITYWEIGHTED},
	\begin{multline}
		\label{E:ELLIPTICHYPERBOLICINTEGRALIDENTITYBULKERRORTERM}
		\EllipticHyperbolicCurrentIntegralIdentityTotalSpacetimeErrorTerm[\SigmatTan,\pmb{\partial} \SigmatTan]
		\eqdef
		\mathfrak{J}_{(\textnormal{Antisymmetric)}}[\pmb{\partial} \SigmatTan,\pmb{\partial} \SigmatTan]
		+
		\mathfrak{J}_{(\textnormal{Div})}[\pmb{\partial} \SigmatTan,\pmb{\partial} \SigmatTan]
		+
		\upmu \weight^{-\blowupratetoporderacoustic}  
		\mathfrak{J}_{(\pmb{\partial} \frac{1}{\upmu} \weight^{\blowupratetoporderacoustic})}[\SigmatTan,\pmb{\partial} \SigmatTan]
		\\
		+
		\mathfrak{J}_{(\textnormal{Absorb}-1)}[\SigmatTan,\pmb{\partial} \SigmatTan]
		+
		\mathfrak{J}_{(\textnormal{Absorb}-2)}[\SigmatTan,\pmb{\partial} \SigmatTan]
		+
		\mathfrak{J}_{(\textnormal{Material})}[\pmb{\partial} \SigmatTan,\pmb{\partial} \SigmatTan]
		+
		\mathfrak{J}_{(\textnormal{Null Geometry})}[\SigmatTan,\pmb{\partial} \SigmatTan]
	\end{multline}
	and the error terms 
	$\mathfrak{J}_{(\textnormal{Antisymmetric)}}[\pmb{\partial} \SigmatTan,\pmb{\partial} \SigmatTan], 
	\cdots, 
	\mathfrak{J}_{(\textnormal{Null Geometry})}[\SigmatTan,\pmb{\partial} \SigmatTan]$
	on RHS\,\eqref{E:ELLIPTICHYPERBOLICINTEGRALIDENTITYBULKERRORTERM}
	are defined in
	\eqref{E:ANTISYMMETRICNULLCURRENTSPACETIMERRORTERM}--\eqref{E:DERIVATIVESOFNULLGEOMETRYNULLCURRENTSPACETIMERRORTERM}.

\end{proposition}

\begin{proof}
	We main idea of the proof is to couple the two divergence identities  \eqref{E:DIVERGENCEIDENTITYSPACETIMEVECTORFIELD} (with $\mathcal{J} = \frac{1}{\upmu}\weight^{\blowupratetoporderacoustic} \ehcurrent[\SigmatTan,\pmb{\partial}]$) and \eqref{E:COVARIANTDIVERGENCEIDENTITYFORELLIPTICHYPERBOLICCURRENT} (with $\widetilde \weight = \frac{1}{\upmu} \weight^{\blowupratetoporderacoustic}$). Using $\mathcal{J}^u = 0$ since $\ehcurrent$ is $\nullhyp_u$-tangent, and $\mathcal{J}^{\ubar} \MagnitueofinnerproductofnewLandnewuL = - \frac{1}{\upmu} \weight^{\blowupratetoporderacoustic} \newuL_\alpha \ehcurrent^\alpha$ (see \eqref{E:DOUBLENULLEIKONALFUNCTIONNORMALIZEDNULLVECTORFIELDS}), these two divergence theorems imply: 
	\begin{align}
		\begin{split} \label{E:INTEGRALIDENTITYFORELLIPTICHYPERBOLICCURRENTINTERMEDIATESTEP1}
			& \int_{\characteristicdiamondtwoarg{[\ubar_1,\ubar_2)}{[u_1,u_2]}} \weight^{\blowupratetoporderacoustic} \ReciprocalLunitAppliedtoTimeFunction\, \ellipticCoerciveQuadratic[\pmb{\partial} \SigmatTan,\pmb{\partial} \SigmatTan]  \, \voldiamond \\
			& = - 
			\int_{\ingoingcharacteristicsurfacetwoarg{\ubar_2}{[u_1,u_2]}} \frac{1}{\upmu} \weight^{\blowupratetoporderacoustic} \newuL_\alpha \ehcurrent^\alpha \volingoingnullhypersurface  + 
			\int_{\ingoingcharacteristicsurfacetwoarg{\ubar_1}{[u_1,u_2]}} \frac{1}{\upmu} \weight^{\blowupratetoporderacoustic} \newuL_\alpha \ehcurrent^\alpha \volingoingnullhypersurface  \\
			& \ \ 
			+ \int_{\characteristicdiamondtwoarg{[\ubar_1,\ubar_2)}{[u_1,u_2]}}   \weight^{\blowupratetoporderacoustic} \ReciprocalLunitAppliedtoTimeFunction \, \EllipticHyperbolicCurrentIntegralIdentityTotalSpacetimeErrorTerm[\SigmatTan,\pmb{\partial} \SigmatTan] \, \voldiamond
		\end{split}
	\end{align}
	We emphasize that the signs of the first two integrals are tied to the negative sign in $\mathcal{J}^{\ubar} = - \frac{1}{\upmu} \weight^{\blowupratetoporderacoustic} \newuL_\alpha \ehcurrent^\alpha$. Applying Lemma\,\ref{L:MIRACULOUSSTRUCTUREINTHEBOUNDARYTERM} to the two ingoing flux boundary terms in \eqref{E:INTEGRALIDENTITYFORELLIPTICHYPERBOLICCURRENTINTERMEDIATESTEP1} concludes the proof of the proposition.
	\end{proof}


\section{Pointwise estimates for the error terms in the commuted wave equations} 
\label{S:POINTWISESTIMATESFORWAVEEQUATIONS}
We continue to work under the assumptions of Sect.\,\ref{SS:SILENTFACTS}.
In this section, we derive pointwise estimates for the error terms that arise
when we commute the wave equations \eqref{E:COVARIANTWAVEEQUATIONSWAVEVARIABLES} up to $\Ntop$ times, where we recall that $\Ntop$ is a fixed integer satisfying \eqref{E:NTOPLARGENESSASSUMPTION}.
More precisely,
for $v^i \in \{v^1,v^2,v^3\}$ 
and $1 \leq N \leq \Ntop$,
we derive pointwise estimates for the inhomogeneous term $ {^{( \tander^{\Ntop},v^i)} \mathfrak{G}}$ 
in the $\upmu$-weighted geometric wave equation
$\upmu \Box_{\gfour}  \tander^{\Ntop} v^i = {^{( \tander^{\Ntop},v^i)} \mathfrak{G}}$ satisfied by $ \tander^{\Ntop}v^i$. 
These pointwise estimates are a preliminary ingredient for the $L^2$ estimates that we derive later on.
Many of the terms appearing in $ {^{( \topordertancom^{\Ntop},v^i)} \mathfrak{G}}$ are harmless from the point of view of regularity 
and the strength of their singularity;  
the bulk of our effort goes towards the most difficult terms, 
which involve the top-order derivatives of the eikonal function
and which we handle by using the modified quantities from Def.\,\ref{D:FULLYANDPARTIALLYMODIFIEDQUANTITIES}.

\subsection{Identification of the most difficult error terms in the commuted wave equations}
\label{SS:MOSTDIFFICULTTERMSINCOMMUTEDWAVEEQUATIONS}
Most of the terms in the commuted wave equations are harmless from the point of view of regularity 
and the strength of their singularity. The next definition captures these ``harmless'' error terms.

\subsubsection{Harmless wave equation error terms}
\label{SSS:HARMLESSWAVEEQUATIONERRORTERMS}
\begin{definition}[Harmless wave equation error terms] 
\label{D:HARMLESSWAVE}
Let $1 \leq N \leq \Ntop$.
We define $\HarmlessWave{N}$ to be any term
that satisfies the following pointwise estimate 
on $\characteristicdiamondtwoarg{[\leftubar,\ubarboot)}{[\moreinterestingu_1,\moreinterestingu_2]}$:
	\begin{align} \label{E:HARMLESSWAVE}
		\left| \HarmlessWave{[1,N]} \right|
		&
		\lesssim \left|\comdersmall^{[1,N+1]; \le 1}\velocityarray\right| + \left|\tander^{\le N} (\vortrenormalized,\GradEnt)\right| + 
				\left| \comdersmall^{[1,N];\le 1} \controlvars \right| + \left| \tandersmall^{[1,N]}\badcontrolvars\right|.
	\end{align}
\end{definition}

The following simple lemma shows that $\HarmlessWave{\Ntop - 12}$ terms are small in 
the norm $\| \cdot \|_{L^{\infty}(\doublenulltoritwoarg{\ubar}{u})}$.

\begin{lemma}[$L^{\infty}$ estimates for $\HarmlessWave{\Ntop - 12}$]  
\label{L:LINFINITYESTIMATESFORHARMLESSWAVEERRORTERMS}
Let $\HarmlessWave{\Ntop - 12}$ be as in Def.\,\ref{D:HARMLESSWAVE}.
Then the following estimate holds for $(\ubar,u) \in [\leftubar,\ubarboot) \times [\moreinterestingu_1,\moreinterestingu_2]$:
\begin{align} \label{E:LINFINITYESTIMATESFORHARMLESSWAVEERRORTERMS}
	\left\| 
		\HarmlessWave{\Ntop - 12} 
	\right\|_{L^{\infty}(\doublenulltoritwoarg{\ubar}{u})} 
	&
	\lesssim 
	\fundbootsmall.
\end{align}
\end{lemma}

\begin{proof}
The estimate \eqref{E:LINFINITYESTIMATESFORHARMLESSWAVEERRORTERMS} follows from definition \eqref{E:HARMLESSWAVE} and
Prop.\,\ref{P:IMPROVEMENTOFAUXILIARYBOOTSTRAP}.
\end{proof}

\subsubsection{The most difficult error terms in the commuted wave equations}
\label{SSS:MOSTDIFFICULTWAVETERMS}
In the following proposition, we identify the most difficult error terms
in the commuted wave equations satisfied by the velocity $\velocityarray$.

\begin{proposition}[Identification of the most difficult error terms in the commuted wave equations] 
\label{P:MOSTDIFFICULTWAVETERMS} 
 Let $\velocityarray \eqdef (v^1,v^2,v^3)$ 
be solutions to the covariant wave equations \eqref{E:COVARIANTWAVEEQUATIONSWAVEVARIABLES}. We denote the product of $\upmu$ and the RHS of the covariant wave equation satisfied by 
$v^i$ by $\mathfrak{G}_{i}$, 
i.e., 
$\upmu \Box_{\gfour} v^i = \mathfrak{G}_{i}$. \begin{subequations}  
\begin{align} \label{E:TOPCOMMUTEDWAVELFIRSTTHENALLYS}
\upmu \Box_{\gfour} (\tanderY^{N-1} \Lunit v^i) 
& 
= 
\angrmd^{\sharp} v^i  \cdot \upmu \angrmd \tanderY^{N-1} \mytr_{\gtorus}\upchi 
+ 
\tanderY^{N-1} \Lunit \mathfrak{G}_{i} 
+ 
\HarmlessWave{N},	
	\\
\begin{split} \label{E:TOPCOMMUTEDWAVEALLYS}
\upmu \Box_{\gfour} (\tanderY^{N-1}\Yvf{A} v^i) 
& 
= 
(\muX v^i)  \tanderY^{N-1} \Yvf{A} \mytr_{\gtorus} \upchi 
+ 
(\Speed^{-2} X^A) \angrmd^{\sharp} v^i \cdot \upmu \angrmd \tanderY^{N-1} \mytr_{\gtorus}\upchi 
\\
& \ \
+ 
\tanderY^{N-1}\Yvf{A} \mathfrak{G}_{i} 
+ 
\HarmlessWave{N}. 
\end{split}
\end{align}Moreover, if $1 \leq N \leq \Ntop$ and $\tander^N$
denotes any order $N$ string of $\mathcal{P}_u$-tangent commutator
other than the ones appearing on LHSs~\eqref{E:TOPCOMMUTEDWAVELFIRSTTHENALLYS}--\eqref{E:TOPCOMMUTEDWAVEALLYS},
(i.e., if $\tander^N$ features at least two copies of $\Lunit$ or only a single $\Lunit$ but does not act first like it does in
\eqref{E:TOPCOMMUTEDWAVELFIRSTTHENALLYS}), 
then $\tander^N v^i$ obeys the following wave equation:
\begin{align} \label{E:TOPCOMMUTEDWAVENOTDIFFICULT}
\upmu \Box_{\gfour} (\tander^N v^i) 
	& = 
	\tander^N \mathfrak{G}_i
	+ 
	\HarmlessWave{N}.
\end{align}
\end{subequations}
\end{proposition}

\begin{proof}
In $2D$, a detailed proof of \eqref{E:TOPCOMMUTEDWAVELFIRSTTHENALLYS}--\eqref{E:TOPCOMMUTEDWAVENOTDIFFICULT} was provided in \cite{jLjS2018}*{Proposition~13.2}, 
except that all of the above-top-order vorticity or entropy involving terms have been soaked  into our definition of $\mathfrak{G}_{i}$. 
Only minor changes are needed to account for the third space dimension, so we omit the details here.  

\end{proof}

\subsubsection{Pointwise estimates for the inhomogeneous terms in the transport equations satisfied by the modified quantities} 
\label{SSS:POINTWISEESTIMATESFORINHOMOGENEOUSTERMSINTRANSPORTEQUATIONSFORMODIFIEDQUANTITIES}
We start with the following lemma, which provides pointwise estimates for 
the inhomogeneous terms in the transport equations satisfied by the fully modified and partially modified quantities.

\begin{lemma}[Pointwise estimates for inhomogeneous terms tied to the modified quantities]
\label{L:POINTWISEESTIMATESFORINHOMOGENEOUSTERMSINTRANSPORTEQUATIONSFORMODIFIEDQUANTITIES}
Let $N = \Ntop$.

\medskip
\noindent \underline{\textbf{Estimates tied to the fully modified quantities}}.
Let $\tander^N \in \mathfrak{P}^{(N)}$ where $\mathfrak{P}^{(N)}$
is the set of order $N$ $\ell_{t,u}$-tangential commutator operators from Def.\,\ref{D:STRINGSOFCOMMUTATIONVECTORFIELDS}. 
Let $\mathfrak{X}$ be the term defined in \eqref{E:MODIFIEDQUANTITYINHOM},
and let $\mathfrak{A}$ be the term appearing on RHS\,\eqref{E:RICCICONTRACTEDLANDL}
and enjoying the schematic structure \eqref{E:AINHOMRIC}.
Then the following pointwise estimates hold on 
$\characteristicdiamondtwoarg{[\leftubar,\ubarboot)}{[\moreinterestingu_1,\moreinterestingu_2]}$:
\begin{subequations}
\begin{align}
	\begin{split}
		\left| 
		(\muX v^1) \tander^N \mathfrak{X} + 2 (\Lunit \upmu) \muX \tander^N v^1\right|
		& 
		\lesssim 
		  \fundbootsmall \left| \newuL \tander^{[1,N]} \velocityarray\right| + \fundbootsmall \upmu \left| \tander^{[1,N+1]}\velocityarray\right|  \\
		& \ \ \ \ + \upmu \left| \tander^{\le N}(\vortrenormalized,\GradEnt)\right| +  \left| \tander^{[1,N]}\controlvars\right| +  \left|\tandersmall^{[1,N]}\badcontrolvars\right| \label{E:POINTWISESUMOFMODIFIEDQUANTITYINHOMANDGLL} 
	\end{split} \\
	\begin{split}
		\left|\tander^N \mathfrak{X} \right| 
		& \lesssim 
		\left|\newuL \tander^N\velocityarray \right| + \upmu  \left|\tander^{[1,N+1]} \velocityarray\right|
		+ 
		\left|\comdersmall^{[1,N];1}\velocityarray \right| \\
		& \ \ \ \ +  \left| \tander^{\le N}(\vortrenormalized,\GradEnt)\right| 
		+
		\left|\tander^{[1,N]}\controlvars\right| 
		+ 
		\left|\tandersmall^{[1,N]}\badcontrolvars \right|,
		\label{E:POINTWISEESTIMATETANGENTDERIVATIVESOFMODQUANTINHOM} 
	\end{split}  
			\\
	\begin{split}
		\label{E:POINTWISEESTIMATESMOOTHTORUSTERIVATIVESDERIVATIVESOFRICLLINHOM}
		|\tander^N \mathfrak{A}|
		& \lesssim \upmu\left| \tander^N(\VortVort,\DivGradEnt)\right| +  \left| \tander^{\le N-1}(\VortVort,\DivGradEnt)\right| + \left| \tander^{\le N}(\vortrenormalized,\GradEnt)\right| +  \left| \comdersmall^{[1,N+1];1} \velocityarray\right|  \\
		& \ \ \ \ 
		+ \upmu \left| \tander^{[1,N+1]}\velocityarray\right| +  \left|\tander^{\le N} (\vortrenormalized,\GradEnt)\right| + 	|\tander^{[1,N]} \controlvars | 
		+ 
		|\tandersmall^{[1,N]}\badcontrolvars|.
	\end{split}
\end{align}
\end{subequations}\medskip

\noindent \underline{\textbf{Estimates tied to the partially modified quantities}}.
If $\tanderY^{N-1} \in \mathfrak{Y}^{(N-1)}$ 
and $\widetilde{\mathfrak{X}}$, 
$\partialmodquantinhom{\tanderY^{N-1}}$,
${^{(\tanderY^{N-1})}\mathfrak{B}}$
are as defined in \eqref{E:PARTIALMODIFIEDQUANTITYINHOMZEROORDER},
\eqref{E:PARTIALMODIFIEDQUANTITYINHOM} 
(with $\tanderY^{N-1}$ in the role of $\tander^N$), 
and \eqref{E:PARTIALMODIFIEDQUANTITYINHOM} (with $\tanderY^{N-1}$ in the role of $\tander^{N-1}$) 
respectively, then the following pointwise estimates hold on 
$\characteristicdiamondtwoarg{[\leftubar,\ubarboot)}{[\moreinterestingu_1,\moreinterestingu_2]}$:
\begin{subequations}
\begin{align}
\left|\partialmodquantinhom{\tanderY^{N-1}} \right|
& \lesssim 
\left| \tander^{[1,N]} \velocityarray\right| +  \left|\tander^{\le N-1} (\vortrenormalized,\GradEnt)\right| + \fundbootsmall \left| \tander^{[1,N-1]} \controlvars\right|, 
\label{E:POINTWISEBELOWTOPORDERPARTIALMODQUANTINHOM} 
		\\
\left| \Lunit \partialmodquantinhom{\tanderY^{N-1}} \right|, 
	\,
\left| \Yvf{A} \partialmodquantinhom{\tanderY^{N-1}}\right| 
& 
\lesssim 
\left| \tander^{[1,N+1]} \velocityarray\right| +  \left|\tander^{\le N} (\vortrenormalized,\GradEnt)\right| + \fundbootsmall \left| \tander^{[1,N]} \controlvars\right|, 
\label{E:POINTWISELANDYDERIVATIVESOFTOPORDERPARTIALMODQUANTINHOM} 
		\\
\left|{^{(\tanderY^{N-1})}\mathfrak{B}} \right| 
& 
\lesssim 
\left |\tander^{[1,N]} \controlvars\right|. 
	\label{E:POINTWISETANGENTDERIVATIVEOFBINHOM}
\end{align}
\end{subequations}\end{lemma}

\begin{proof}
To prove \eqref{E:POINTWISESUMOFMODIFIEDQUANTITYINHOMANDGLL}, we write: 
	\begin{align}
		\begin{split} \label{E:POINTWISESUMOFMODIFIEDQUANTITYINHOMANDGLLSTEP1}
			\tander^N\mathfrak{X} & =  - \vec{G}_{\Lunit\Lunit} \diamond \muX \tander^N \wavearray  - \vec{G}_{\Lunit\Lunit} \diamond [\tander^N,\muX] \wavearray - \left[ \tander^N,\vec{G}_{\Lunit\Lunit}\right] \diamond \muX  \wavearray \\
			& \ \  + \tander^N \left\{ - \frac{1}{2} \upmu \mytr_{\gtorus} \angG \diamond \Lunit \wavearray - \frac{1}{2} \upmu \vec{G}_{\Lunit \Lunit} \diamond \Lunit \wavearray + \upmu \angG_{\Lunit}^{\#} \diamond \cdot \angrmD \wavearray\right\}
		\end{split}
	\end{align}
Multiplying both sides of \eqref{E:POINTWISESUMOFMODIFIEDQUANTITYINHOMANDGLLSTEP1} by $\muX v^1$ and adding $2 (\Lunit \upmu) \muX \tander^N v^1$ to both sides as well, the desired bound \eqref{E:POINTWISESUMOFMODIFIEDQUANTITYINHOMANDGLL} follows from commutator estimate \eqref{E:COMMUTATOROFMUXANDTANGENTIALCOMMUTATORS}, Prop.\,\ref{P:ENTROPYDERIVATIVESINTERMSOFOTHERS},  pointwise bound \eqref{E:ARBITRARYCOMMUTATORSTRINGESTIMATE}, the trilinear structure \eqref{E:TRILINEARSTRUCTUREFORMIXEDTOPORDERLMULTIPLIERTERM},  and Cor.\,\ref{C:IMPROVEAUX}.

All of the remaining estimates stated in the lemma are straightforward consequences of
Lemma\,\ref{L:ANGULARDIFFERENTIALCOMMUTESWITHANGLIE},
Prop.\,\ref{P:SCHEMATICSTRUCTUREOFVARIOUSTENSORSINTERMSOFCONTROLVARS},
Lemma\,\ref{L:SCHEMATICEXPRESSIONFORANGULARLAPLACIAN},
the commutator estimate \eqref{E:COMMUTATOROFTANGENTIALANDTANGENTIALCOMMUTATORS},
and the estimates of Props.\,\ref{P:ENTROPYDERIVATIVESINTERMSOFOTHERS} and\,\ref{P:IMPROVEMENTOFAUXILIARYBOOTSTRAP}.
\end{proof}

\section{Pointwise estimates for controlling the specific vorticity, entropy gradient, and modified fluid variables}
\label{S:POINTWISEESTIMATESFORCONTROLLINGSPECIFICVORTICITYANDENTROPYGRADIENT}
We continue to work under the assumptions of Sect.\,\ref{SS:SILENTFACTS}.
In this section, we derive a variety of pointwise estimates that we will use
in Sects.\,\ref{S:BELOWTOPORDERHYPERBOLICL2ESTIMATESFORSPECIFICVORTICITYANDENTROPYGRADIENT}--\eqref{S:TOPORDERELIPTICHYPERBOLICL2ESTIMATESFORSPECIFICVORTICITYANDENTROPYGRADIENT}, 
when we derive $L^2$ estimates 
for $\vortrenormalized$, $\GradEnt$, $\VortVort$, and $\DivGradEnt$ up to top-order.
Among the estimates we derive are pointwise estimates for the error terms in the
elliptic-hyperbolic integral identity \eqref{E:INTEGRALIDENTITYFORELLIPTICHYPERBOLICCURRENT}.

\subsection{Simple commutator estimates involving $\vort$ and $\GradEnt$}
\label{SSS:SIMPLECOMMUTATORLEMMAINVOLVINGTRANSPORTVARIABLE}
The following simple lemma will be used in our analysis. 

\begin{lemma}[Commuting $\upmu$-weighted Cartesian derivatives with the geometric vectorfields] 
\label{L:POINTWISEESTIMATESFORMUWEIGHTEDCARTESIANCOMMUTATOR} 
Let $1 \leq N \leq \Ntop$.
Then the following commutator estimates hold on $\characteristicdiamondtwoarg{[\leftubar,\ubarboot)}{[\moreinterestingu_1,\moreinterestingu_2]}$:
\begin{align}
\begin{split} \label{E:POINTWISEESTIMATESFORMUWEIGHTEDCARTESIANCOMMUTATOR}
	& 
	\left| [\upmu \partial_i, \tander^N] (\vortrenormalized,\GradEnt) \right|, 
		\, 
	\left| [\upmu \Flatcurl, \tander^N] (\vortrenormalized,\GradEnt) \right|, 
		\, 
	\left| [\upmu \Flatdiv, \tander^N] (\vortrenormalized,\GradEnt) \right|  
		\\
	&  
	\lesssim 
	\left| \tander^{\leq N} (\vortrenormalized, \GradEnt) \right| 
	+ 
	\left| \muX \tander^{\leq N-1}(\vortrenormalized,\GradEnt)\right| 
	+  
	\fundbootsmall |\muX \tander^{[1,N-1]} \velocityarray| 
	+
	\fundbootsmall \upmu |\tander^{[1,N]} \velocityarray| 
	+ 
	\fundbootsmall |\tandersmall^{[1,N]} \badcontrolvars|. 
\end{split}
\end{align}
\end{lemma}

\begin{proof}
The desired estimate follows from \eqref{E:POINTWISEESTIMATESFORMUWEIGHTEDCARTESIANCOMMUTATORWITHALLWAVEVARS} using \eqref{E:HIGHERORDERTANGENTIALDERIVATIVESOFDENSITY}--\eqref{E:HIGHERORDEREXACTLYONEMUXINDERIVATIVESOFDENSITY} to bound the density in $\wavearray = (\LogDensity,v^1,v^2,v^3,s)$, as well as using the straightforward identity $\comder s = \smoothfunction(\controlvars,\badcontrolvars)  S$.  
\end{proof}

\subsection{Pointwise estimates for $\vortrenormalized$, $\GradEnt$, $\VortVort$, $\DivGradEnt$, and their derivatives}
\label{SS:POINTWISEESTIMATESFORTRANSPORTVARIABLESANDTHEIRDERIVATIVES}
In this section, we derive
pointwise estimates for $\vortrenormalized$, $\GradEnt$, $\VortVort$, $\DivGradEnt$, and various derivatives
of these quantities.  We provide the main estimates in Prop.\,\ref{P:POINTWISESTIAMTESFORALLTHETRANSPORTVARIABLES}.

\subsubsection{The main pointwise estimates}
\label{SSS:POINTWISEESTIMATESFORTRANSPORTVARIABLESANDTHEIRDERIVATIVES}

\begin{proposition}[Pointwise estimates for $\vortrenormalized,\GradEnt,\VortVort,\DivGradEnt$, and their derivatives]
	\label{P:POINTWISESTIAMTESFORALLTHETRANSPORTVARIABLES}
	The following pointwise estimates hold on 
	$\characteristicdiamondtwoarg{[\leftubar,\ubarboot)}{[\moreinterestingu_1,\moreinterestingu_2]}$:
	
	\medskip 
	
	\noindent \underline{\textbf{Transport estimates}}.
	For $0 \leq N \leq \Ntop$, we have:	
	\begin{align}
		|\upmu \Transport \tander^N (\vortrenormalized,\GradEnt)|
		& \lesssim 
			|\tander^{\leq N} (\vortrenormalized,\GradEnt)|
			+
			\fundbootsmall |\newuL \tander^{[1,N]} \velocityarray|
			+
			\fundbootsmall \upmu |\tander^{[1,N + 1]} \velocityarray|
			+
			\fundbootsmall |\tandersmall^{[1,N]} \badcontrolvars|,
				\label{E:COMMUTEDTRANSPORTPOINTWISEESTIMATESFORSPECIFICVORTICITYANDENTROPYGRADIENT} 
	\end{align}
	
	\begin{subequations}
	\begin{align}
		|\upmu \Transport \tander^N \VortVort|
		& \lesssim 
			|\tander^{\leq N} \VortVort|
			+
			|\tander^{\leq N+1} (\vortrenormalized,\GradEnt)|
			+
			\fundbootsmall |\newuL \tander^{[1,N]} \velocityarray|
			+
			\fundbootsmall |\tander^{[1,N + 1]} \velocityarray|
			+
			\fundbootsmall |\tandersmall^{[1,N]} \badcontrolvars|,
			\label{E:COMMUTEDTRANSPORTPOINTWISEESTIMATESFORMODIFIEDCURLOFVORT}
				\\
		|\upmu \Transport \tander^N \DivGradEnt|
		& \lesssim 
			|\tander^{\leq N} \DivGradEnt|
			+
			|\tander^{\leq N+1} (\vortrenormalized,\GradEnt)|
			+
			\fundbootsmall |\newuL \tander^{[1,N]} \velocityarray|
			+
			\fundbootsmall | \tander^{[1,N + 1]} \velocityarray|
			+
			\fundbootsmall |\tandersmall^{[1,N]} \badcontrolvars|.
			\label{E:COMMUTEDTRANSPORTPOINTWISEESTIMATESFORMODIFIEDDIVERGENCEOFENTROPYGRADIENT}
	\end{align}
	\end{subequations}	
	\medskip
	
	\noindent \underline{\textbf{Algebraic estimates for transversal derivatives in terms of tangential derivatives}}.
	For $0 \leq N \leq \Ntop$, we have:	
	\begin{subequations}
	\begin{align}
			|\muX \tander^N \vortrenormalized|,
				\,
			|\tander^N \muX \vortrenormalized| 
			& \lesssim 
			\upmu |\Lunit \tander^N \vortrenormalized|
			+
			|\tander^{\leq N} (\vortrenormalized,\GradEnt)|
			+
			\fundbootsmall |\newuL \tander^{[1,N]} \velocityarray|
			+
			\fundbootsmall \upmu | \tander^{[1,N + 1]} \velocityarray|
			+
			\fundbootsmall |\tandersmall^{[1,N]} \badcontrolvars|,
				\label{E:COMMUTEDTRANSVERSALDERVIATVIESOFVORTICITYINTERMSOFTANGENTIAL} 
					\\
			|\muX \tander^N \GradEnt|,
				\,
			|\tander^N \muX \GradEnt|
			& \lesssim 
			\upmu |\Lunit \tander^N \GradEnt|
			+
			|\tander^{\leq N} (\vortrenormalized,\GradEnt)|
			+
			\fundbootsmall |\newuL \tander^{[1,N]} \velocityarray|
			+
			\fundbootsmall \upmu | \tander^{[1,N + 1]} \velocityarray|
			+
			\fundbootsmall |\tandersmall^{[1,N]} \badcontrolvars|.
				\label{E:COMMUTEDTRANSVERSALDERVIATVIESOFENTROPYGRADIENTINTERMSOFTANGENTIAL} 
		\end{align}
		\end{subequations}

	\medskip
	
	\noindent \underline{\textbf{Algebraic estimates for $(\Flatdiv \vortrenormalized, \Flatdiv \GradEnt)$ and 
	 $(\Flatcurl \vortrenormalized, \Flatcurl \GradEnt)$ in terms of $(\VortVort,\DivGradEnt)$}}.
	For $0 \leq N \leq \Ntop$, we have:	
	\begin{subequations}
	\begin{align}
	\begin{split} 	\label{E:COMMUTEDPOINTWISEEUCLIDEANCURLOFVORTICITY}
		\left|
			\Flatcurl \tander^N \vortrenormalized
		\right|
		& \lesssim 
			\left|
				\tander^N \VortVort 
			\right|
			+
			\frac{1}{\upmu}
			\left|
				\tander^{\leq N-1} \VortVort 
			\right|
			+
			\frac{1}{\upmu}
			\left|
				\tander^{\leq N} (\vortrenormalized,\GradEnt) 
			\right|
				\\
		& \ \
			+
			\frac{\fundbootsmall}{\upmu} 
			\left|	
				\newuL \tander^{[1,N]} \velocityarray 
			\right|
			+
			\fundbootsmall 
			\left|
				\tander^{[1,N + 1]} \velocityarray 
			\right|
			+
			\frac{\fundbootsmall}{\upmu}
			\left|
				\tandersmall^{[1,N]} \badcontrolvars 
			\right|,
			\end{split} 
				\\
		\begin{split} 	\label{E:COMMUTEDPOINTWISEEUCLIDEANEUCLIDEANDIVERGENCEOFENTROPYGRADIEENT} 
		\left|
			\Flatdiv \tander^N \GradEnt
		\right|
		& \lesssim	
			\left|
				\tander^N \DivGradEnt 
			\right|
			+
			\frac{1}{\upmu}
			\left|
				\tander^{\leq N-1} \DivGradEnt 
			\right|
			+
			\frac{1}{\upmu}
			\left|
				\tander^{\leq N} (\vortrenormalized,\GradEnt) 
			\right|
				\\
		& \ \
			+
			\frac{\fundbootsmall}{\upmu}
			\left|	
				\newuL \tander^{[1,N]} \velocityarray 
			\right|
			+
			\fundbootsmall 
			\left|
				\tander^{[1,N + 1]}\velocityarray 
			\right|
			+
			\frac{\fundbootsmall}{\upmu}
			\left|
				\tandersmall^{[1,N]} \badcontrolvars 
			\right|,
		\end{split}
				\\
		\left|
			\Flatdiv \tander^N \vortrenormalized
		\right|,
			\,
		\left|
			\Flatcurl \tander^N \GradEnt
		\right|
		& \lesssim 
			\frac{1}{\upmu}
			|\tander^{\leq N} (\vortrenormalized,\GradEnt)|
			+
			\frac{\fundbootsmall}{\upmu}
			|\newuL \tander^{[1,N]} \velocityarray|
			+
			\fundbootsmall | \tander^{[1,N + 1]} \velocityarray|
			+
			\frac{\fundbootsmall}{\upmu}
			|\tandersmall^{[1,N]} \badcontrolvars|.
				\label{E:COMMUTEDPOINTWISEEUCLIDEANDIVERGENCEOFVORTICITYANDEUCLIDEANCURLOFENTROPYGRADIENT}
	\end{align}
	\end{subequations}
	
	\medskip
	
	\noindent \underline{\textbf{Estimates for the exterior derivative of 
	$(\tander^N \vortrenormalized)_{\flat}$ and 
	$(\tander^N \GradEnt)_{\flat}$}}.
	The following estimates hold, where
	$(\rmd \tander^N \vortrenormalized)_{\flat}$ is the two-form with the Cartesian components components
	$\partial_{\alpha} (\tander^N \vortrenormalized)_{\beta}
	-
	\partial_{\beta} (\tander^N \vortrenormalized)_{\alpha}
	$
	(where 
	$(\tander^N \vortrenormalized)_{\alpha}
	= \gfour_{\alpha \gamma} \tander^N \vortrenormalized^{\gamma}
	$), and similarly for $\rmd (\tander^N \GradEnt)_{\flat}$:
	\begin{subequations}
	\begin{align}
	\begin{split} \label{E:POINTWISEBOUNDEXTERIORDERIVATIVEOFTOPORDERDERIVATIVESOFVORTICITY} 
		\sum_{\alpha,\beta = 0}^4 |\rmd ((\tander^N \vortrenormalized)_{\flat})_{\alpha\beta} |
		& \lesssim 
			\left|
				\tander^N \VortVort 
			\right|
			+
			\frac{1}{\upmu}
			\left|
				\tander^{\le N -1} \VortVort 
			\right|
			+
			\frac{1}{\upmu}
			\left|
				\tander^{\le N} (\vortrenormalized,\GradEnt) 
			\right|
				\\
		& \ \
			+
			\frac{\fundbootsmall}{\upmu}
			\left|	
				\newuL \tander^{[1,N]} \wavearray 
			\right|
			+
			\fundbootsmall 
			\left|
				\tander^{[1,N + 1]} \wavearray 
			\right|
			+
			\frac{\fundbootsmall}{\upmu}
			\left|
				\tandersmall^{[1,\Ntop]} \badcontrolvars 
			\right|,
		\end{split}		
					\\
		\sum_{\alpha,\beta = 0}^4 |\rmd ((\tander^N \GradEnt)_{\flat})_{\alpha\beta}| 
		& \lesssim 
			\frac{1}{\upmu}
			\left|
				\tander^{\le N} (\vortrenormalized,\GradEnt) 
			\right|
			+
			\frac{\fundbootsmall}{\upmu}
			\left|	
				\newuL \tander^{[1,N]} \wavearray 
			\right|
			+
			\fundbootsmall 
			\left|
				\tander^{[1,N + 1]}  \wavearray 
			\right|
			+
			\frac{\fundbootsmall}{\upmu}
			\left|
				\tandersmall^{[1,\Ntop]} \badcontrolvars 
			\right|.
				\label{E:POINTWISEBOUNDEXTERIORDERIVATIVEOFTOPORDERDERIVATIVESOFENTROPYGRADIENT} 
			\end{align}
		\end{subequations}
	
\end{proposition}

\begin{proof}
Let $0 \le N \le \Ntop$. Using  Lemma\,\ref{L:RELATIONSHIPBETWEENCARTESIANPARTIALDERIVATIVESANDSMOOTHGEOMETRICCOMMUTATORS}, Props.\,\ref{P:SCHEMATICSTRUCTUREOFVARIOUSTENSORSINTERMSOFCONTROLVARS},\,\ref{P:COMMUTATORESTIMATES} and\,\ref{P:IMPROVEMENTOFAUXILIARYBOOTSTRAP}, and the results of Sect.~\ref{SS:TRANSVERSALDERIVATIVESINTERMSOFTANGENTIALONES}, all of the estimates \eqref{E:COMMUTEDTRANSPORTPOINTWISEESTIMATESFORSPECIFICVORTICITYANDENTROPYGRADIENT}--\eqref{E:POINTWISEBOUNDEXTERIORDERIVATIVEOFTOPORDERDERIVATIVESOFENTROPYGRADIENT} were proven in \cite{abbrescia2022emergence}*{Prop.\,23.3} with $|\newuL \tander^{[1,N]} \velocityarray|$ and $|\tander^{[1,N + 1]} \velocityarray|$ replaced by $|\muX \tander^{[1,N]} \wavearray|$ and $|\tander^{N+1}\wavearray|$, respectively.\footnote{Actually, in \cite{abbrescia2022emergence}*{Prop.\,23.3}, we introduced a Riemannian metric on the space of type $\binom{0}{2}$-tensorfields and the proposition there featured $|\mathrm{d}(\tander^N(\vortrenormalized,\GradEnt)_\flat|_{\hfour}$. However, we proved in \cite{abbrescia2022emergence}*{Lemma.\,21.7} that: 
	\begin{align} |\upxi|_{\hfour} \approx \sum_{\substack{0 \le \alpha_1,\dots, \alpha_m \le 3 \\ 0 \le \beta_1,\dots,\beta_n \le 3}} \left|\xi_{\beta_1\cdots\beta_n}^{\alpha_1\cdots\alpha_m}\right| \label{E:HFOURINTERMSOFCARTESIANNORM}
	\end{align}
for any type $\binom{m}{n}$-tensorfield.} To conclude the proof, we use \eqref{E:HIGHERORDERTANGENTIALDERIVATIVESOFDENSITY}--\eqref{E:HIGHERORDEREXACTLYONEMUXINDERIVATIVESOFDENSITY} to bound the density in $\wavearray = (\LogDensity,v^1,v^2,v^3,s)$ (as well as using the straightforward identity $\comder s = \smoothfunction(\controlvars,\badcontrolvars)  S$) and then use \eqref{E:POINTWISESTIMATEFORBREVEXINTERMSOFNEWUL} to bound $|\muX \tander^{[1,N]} \velocityarray| \le |\newuL  \tander^{[1,N]} \velocityarray| + \upmu |\tander^{[1,N + 1]} \velocityarray|$. 
\end{proof}

\subsection{Pointwise estimates for the elliptic-hyperbolic integral identity error terms}
\label{SS:POINTWISEESTIMTAESFORELLIPTICHYPERBOLICIDENTITYERORTERMS} 
Recall that to prove our top-order $L^2$ estimates for the specific vorticity and entropy gradient,
we will rely on the elliptic-hyperbolic integral identity \eqref{E:INTEGRALIDENTITYFORELLIPTICHYPERBOLICCURRENT}
with $\tander^N \vortrenormalized$ and $\tander^N \GradEnt$
in the role of $\SigmatTan$.
In the next proposition, we derive pointwise estimates for the error terms appearing in the identity.

\begin{proposition}[Pointwise estimates for the elliptic-hyperbolic integral identity error terms]
	\label{P:POINTWISEESTIMTAESFORELLIPTICHYPERBOLICIDENTITYERORTERMS} 
	Let $\varsigma \in (0,1]$ and $N = \Ntop$.
	Then the error terms 
	appearing in the elliptic-hyperbolic integral identity \eqref{E:INTEGRALIDENTITYFORELLIPTICHYPERBOLICCURRENT}
	(with $\tander^N \vortrenormalized$ and $\tander^N \GradEnt$ in the role of $\SigmatTan$)
	satisfy the following pointwise estimates 
	on $\characteristicdiamondtwoarg{[\leftubar,\ubarboot)}{[\moreinterestingu_1,\moreinterestingu_2]}$,
	where the implicit constants are \textbf{independent} of $\varsigma$.
	
	\medskip 
	\noindent \underline{\textbf{Estimates for controlling bulk error integrals}}.
		\begin{subequations}
		\begin{align}
		\begin{split} \label{E:ANTISYMMETRICVORTICITYELLIPTICHYPERBOLICIDENTITYERORTERMPOINTWISE}	
		\left|
			\mathfrak{J}_{(\textnormal{Antisymmetric})}
			[\pmb{\partial} \tander^N \vortrenormalized,\pmb{\partial} \tander^N \vortrenormalized]
		\right|
		&
		\lesssim
			\left|
				\tander^N \VortVort 
			\right|^2
			+
			\frac{1}{\upmu^2}
			\left|
				\tander^{\le N -1} \VortVort 
			\right|^2				
			+
			\frac{1}{\upmu^2}
			\left|
				\tander^{\le N} (\vortrenormalized,\GradEnt) 
			\right|^2
				\\
		& \ \
			+
			\frac{\fundbootsmall^2}{\upmu^2}
			\left|	
				\newuL \tander^{[1,N]} \velocityarray 
			\right|^2
			+
			\fundbootsmall^2 
			\left|
				\tander^{[1,N+1]} \velocityarray 
			\right|^2
			+
			\frac{\fundbootsmall^2}{\upmu^2}
			\left|
				\tandersmall^{[1,N]} \badcontrolvars 
			\right|^2,
		\end{split}
			\\
		\begin{split} \label{E:ANTISYMMETRICENTROPYGRADIENTELLIPTICHYPERBOLICIDENTITYERORTERMPOINTWISE} 
		\left|
			\mathfrak{J}_{(\textnormal{Antisymmetric})}[\pmb{\partial} \tander^N\GradEnt,\pmb{\partial} \tander^N \GradEnt]
		\right|
		& \lesssim
			\frac{1}{\upmu^2}
			\left|
				\tander^{\le N} (\vortrenormalized,\GradEnt) 
			\right|^2
				\\
		&  
			\ \
			+
			\frac{\fundbootsmall^2}{\upmu^2}
			\left|	
				\newuL \tander^{[1,N]} \velocityarray 
			\right|^2
			+
			\fundbootsmall^2 
			\left|
				\tander^{[1,N+1]} \velocityarray 
			\right|^2
			+
			\frac{\fundbootsmall^2}{\upmu^2}
			\left|
				\tandersmall^{[1,N]} \badcontrolvars 
			\right|^2,				
	\end{split}
	\end{align}
\end{subequations}	
	
\begin{subequations}
	\begin{align}
	\begin{split} \label{E:DIVVORTICITYELLIPTICHYPERBOLICIDENTITYERORTERMPOINTWISE}	 
	\left|
			\mathfrak{J}_{(\textnormal{Div})}[\pmb{\partial} \tander^N \vortrenormalized,\pmb{\partial} \tander^N\vortrenormalized]
		\right| 			
	& \lesssim
			\frac{1}{\upmu^2}
			|\tander^{\le N} (\vortrenormalized,\GradEnt)|^2
				\\
	& \ \
			+
			\frac{\fundbootsmall^2}{\upmu^2}
			|\newuL \tander^{[1,N]} \velocityarray|^2
			+
			\fundbootsmall^2 |\tander^{[1,N+1]} \velocityarray|^2
			+
			\frac{\fundbootsmall^2}{\upmu^2}
			|\tandersmall^{[1,N]} \badcontrolvars|^2, 
	\end{split}
					\\
	\begin{split} \label{E:DIVENTROPYGRADIENTELLIPTICHYPERBOLICIDENTITYERORTERMPOINTWISE} 
		\left|
			\mathfrak{J}_{(\textnormal{Div})}[\pmb{\partial} \tander^N \GradEnt,\pmb{\partial} \tander^N\GradEnt]
		\right| 			
		& \lesssim
			\left|
				\tander^N \DivGradEnt 
			\right|^2
			+
			\frac{1}{\upmu^2}
			\left|
				\tander^{\le N -1} \DivGradEnt 
			\right|^2
			+
			\frac{1}{\upmu^2}
			|\tander^{\le N} (\vortrenormalized,\GradEnt)|^2
				\\
		& \ \
			+
			\frac{\fundbootsmall^2}{\upmu^2}
			|\newuL \tander^{[1,N]} \velocityarray|^2
			+
			\fundbootsmall^2 |\tander^{[1,N+1]} \velocityarray|^2
			+
			\frac{\fundbootsmall^2}{\upmu^2}
			|\tandersmall^{[1,N]} \badcontrolvars|^2,
	\end{split}
	\end{align}
\end{subequations}
	
\begin{subequations}
\begin{align}
		\upmu \weight^{-\blowupratetoporderacoustic}  
		\left|
			\mathfrak{J}_{(\pmb{\partial} \frac{1}{\upmu} \weight^{\blowupratetoporderacoustic})}
				[\tander^N \vortrenormalized,\pmb{\partial} \tander^N \vortrenormalized]
		\right|
		& \lesssim
			\varsigma 
			\ellipticCoerciveQuadratic[\pmb{\partial}\tander^N \vortrenormalized,\pmb{\partial} \tander^N			
			\vortrenormalized]
			+
			\frac{1}{\varsigma}
			\frac{1}{\weight^2}
			|\tander^{\le N} \vortrenormalized|^2,
				\label{E:DERIVATIVEOFONEOVERMUWEIGHTVORTICITYELLIPTICHYPERBOLICIDENTITYERORTERMPOINTWISE} 
					\\
		\upmu \weight^{-\blowupratetoporderacoustic}  
		\left|
			\mathfrak{J}_{(\pmb{\partial} \frac{1}{\upmu} \weight^{\blowupratetoporderacoustic})}[\tander^N \GradEnt,\pmb{\partial} \tander^N \GradEnt]
		\right|
		& \lesssim
			\varsigma 
			\ellipticCoerciveQuadratic[\pmb{\partial}\tander^N \GradEnt,\pmb{\partial} \tander^N			
			\GradEnt]
			+
			\frac{1}{\varsigma} 
			\frac{1}{\weight^2}
			|\tander^{\le N} \GradEnt|^2,	
			\label{E:DERIVATIVEOFONEOVERMUWEIGHTENTROPYGRADIENTELLIPTICHYPERBOLICIDENTITYERORTERMPOINTWISE} 
	\end{align} 
	\end{subequations}
	
\begin{subequations}
	\begin{align}
	\begin{split} \label{E:FIRSTABSORBTERMVORTICITYELLIPTICHYPERBOLICIDENTITYERORTERMPOINTWISE}  
	\left|
			\mathfrak{J}_{(\textnormal{Absorb-1})}[\tander^N \vortrenormalized,\pmb{\partial} \tander^N\vortrenormalized]
		\right| 
		& 
		\lesssim
			\varsigma 
			\ellipticCoerciveQuadratic[\pmb{\partial}\tander^N \vortrenormalized,\pmb{\partial} \tander^N\vortrenormalized]
					\\
		& \ \
			+
			\left(1 + \frac{1}{\varsigma} \right)
			\left|
				\tander^N \VortVort 
			\right|^2
			+
			\left(1 + \frac{1}{\varsigma} \right)
			\frac{1}{\upmu^2}
			\left|
				\tander^{\le N -1} \VortVort 
			\right|^2
					\\					
			& \ \
			+
			\left(1 + \frac{1}{\varsigma} \right)
			\frac{1}{\upmu^2}
			|\tander^{\le N} (\vortrenormalized,\GradEnt)|^2
				\\
		& \ \
			+
			\left(1 + \frac{1}{\varsigma} \right)
			\frac{\fundbootsmall^2}{\upmu^2}
			|\newuL \tander^{[1,N]} \velocityarray|^2
			+
			\left(1 + \frac{1}{\varsigma} \right)
			\fundbootsmall^2 |\tander^{[1,N+1]} \velocityarray|^2
				\\
		& \ \
			+
			\left(1 + \frac{1}{\varsigma} \right)
			\frac{\fundbootsmall^2}{\upmu^2}
			|\tandersmall^{[1,N]} \badcontrolvars|^2,
		\end{split}
			\\
		\begin{split} \label{E:FIRSTABSORBTERMENTROPYGRADIENTELLIPTICHYPERBOLICIDENTITYERORTERMPOINTWISE} 
		\left|
			\mathfrak{J}_{(\textnormal{Absorb-1})}[\tander^N \GradEnt,\pmb{\partial} \tander^N \GradEnt]
		\right|
		 &
		\lesssim
			\varsigma 
			\ellipticCoerciveQuadratic[\pmb{\partial}\tander^N \GradEnt,\pmb{\partial} \tander^N \GradEnt]
				\\
			& \ \
				+
			\left(1 + \frac{1}{\varsigma} \right)
			\frac{1}{\upmu^2}
			|\tander^{\le N} (\vortrenormalized,\GradEnt)|^2
				\\
		& \ \
			+
			\left(1 + \frac{1}{\varsigma} \right)
			\frac{\fundbootsmall^2}{\upmu^2}
			|\newuL \tander^{[1,N]} \velocityarray|^2
			+
			\left(1 + \frac{1}{\varsigma} \right)
			\fundbootsmall^2  | \tander^{[1,N+1]} \velocityarray|^2
				\\
		& \ \
			+
			\left(1 + \frac{1}{\varsigma} \right)
			\frac{\fundbootsmall^2}{\upmu^2}
			|\tandersmall^{[1,N]} \badcontrolvars|^2,
	\end{split}
	\end{align}
\end{subequations}	\begin{subequations}
		\begin{align}
		\left|
			\mathfrak{J}_{(\textnormal{Absorb-2})}[\tander^N \vortrenormalized,\pmb{\partial} \tander^N\vortrenormalized]
		\right|
		& \lesssim
			\varsigma 
			\ellipticCoerciveQuadratic[\pmb{\partial}\tander^N \vortrenormalized,\pmb{\partial} \tander^N\vortrenormalized]
			+
			\frac{1}{\varsigma} 
			\frac{1}{\upmu^2}
			|\tander^N \vortrenormalized|^2,
			 \label{E:SECONDABSORBTERMVORTICITYELLIPTICHYPERBOLICIDENTITYERORTERMPOINTWISE} 
				\\
		\left|
			\mathfrak{J}_{(\textnormal{Absorb-2})}[\tander^N \GradEnt,\pmb{\partial} \tander^N \GradEnt]
		\right|
		& \lesssim
			\varsigma 
			\ellipticCoerciveQuadratic[\pmb{\partial}\tander^N \GradEnt,\pmb{\partial} \tander^N \GradEnt]
			+
			\frac{1}{\varsigma} 
			\frac{1}{\upmu^2}
			|\tander^N \GradEnt|^2,
			\label{E:SECONDABSORBTERMENTROPYGRADIENTELLIPTICHYPERBOLICIDENTITYERORTERMPOINTWISE} 
	\end{align}
	\end{subequations}	\begin{subequations}
	\begin{align}
	\begin{split} \label{E:QUADRATICINMATERIALDERIVATIVETERMVORTICITYELLIPTICHYPERBOLICIDENTITYERORTERMPOINTWISE} 
		\left|
			\mathfrak{J}_{(\textnormal{Material})}[\tander^N \vortrenormalized,\pmb{\partial} \tander^N\vortrenormalized]
		\right| 
		 &
		\lesssim
			\frac{1}{\upmu^2}
			|\tander^{\le N} (\vortrenormalized,\GradEnt)|^2
					\\
			& \ \
				+
			\frac{\fundbootsmall^2}{\upmu^2}
			|\newuL \tander^{[1,N]} \velocityarray|^2
			+
			\fundbootsmall^2 | \tander^{[1,N+1]} \velocityarray|^2
			+
			\frac{\fundbootsmall^2}{\upmu^2}
			|\tandersmall^{[1,N]} \badcontrolvars|^2,
		\end{split} 	
						\\
		\begin{split} 
		\label{E:QUADRATICINMATERIALDERIVATIVETERMENTROPYGRADIENTELLIPTICHYPERBOLICIDENTITYERORTERMPOINTWISE} 
		\left|
			\mathfrak{J}_{(\textnormal{Material})}[\tander^N \GradEnt,\pmb{\partial} \tander^N \GradEnt]
		\right| 
		 &
		\lesssim
		\frac{1}{\upmu^2}
			|\tander^{\le N} (\vortrenormalized,\GradEnt)|^2
				\\
		& \ \
			+
			\frac{\fundbootsmall^2}{\upmu^2}
			|\newuL \tander^{[1,N]} \velocityarray|^2
			+
			\fundbootsmall^2 | \tander^{[1,N+1]} \velocityarray|^2
			+
			\frac{\fundbootsmall^2}{\upmu^2}
			|\tandersmall^{[1,N]} \badcontrolvars|^2,
	\end{split}
	\end{align}
\end{subequations}
	
	\begin{subequations}
	\begin{align}
		\left|
			\mathfrak{J}_{(\textnormal{Null Geometry})}[\tander^N \vortrenormalized,\pmb{\partial} \tander^N \vortrenormalized]
		\right|
		& \lesssim
			\varsigma 
			\ellipticCoerciveQuadratic[\pmb{\partial}\tander^N \vortrenormalized,\pmb{\partial} \tander^N\vortrenormalized]
			+
			\frac{1}{\varsigma} 
			\frac{1}{\upmu^2}
			|\tander^N \vortrenormalized|^2,
			 \label{E:NULLGEOMETRYVORTICITYELLIPTICHYPERBOLICIDENTITYERORTERMPOINTWISE}  
				\\
		\left|
			\mathfrak{J}_{(\textnormal{Null Geometry})}[\tander^N \GradEnt,\pmb{\partial} \tander^N \GradEnt]
		\right|
		& \lesssim
			\varsigma 
			\ellipticCoerciveQuadratic[\pmb{\partial}\tander^N \GradEnt,\pmb{\partial} \tander^N \GradEnt]
			+
			\frac{1}{\varsigma} 
			\frac{1}{\upmu^2}
			|\tander^N \GradEnt|^2.
			\label{E:NULLGEOMETRYENTROPYGRADIENTELLIPTICHYPERBOLICIDENTITYERORTERMPOINTWISE}
	\end{align}
	\end{subequations} 
	\medskip
	\noindent \underline{\textbf{Estimates for controlling boundary error integrals}}. 
		
	\begin{subequations}
	\begin{align}
	\begin{split} \label{E:PRINCIPALSIGAMTILDEELLIPTICHYPERBOLICVORTICITYPOINTWISE}	
		\weight^{- \blowupratetoporderacoustic} \left|
			\, ^{(\tander^N \vortrenormalized)}\Currentboundaryerrorhavetocontrolprincipalweighted{\blowupratetoporderacoustic} \right|
		& 
		\lesssim
				\upmu
			\left|
				\tander^N \VortVort 
			\right|^2
			+
			\frac{1}{\upmu}
			\left|
				\tander^{\le N} (\vortrenormalized,\GradEnt) 
			\right|^2
					\\
		& \ \
			+
			\frac{\fundbootsmall^2}{\upmu}
			\left|	
				\newuL \tander^{[1,N]} \velocityarray 
			\right|^2
			+
			\fundbootsmall^2 
				\upmu
			\left| \nullangrmD 
				\tander^{[1,N]} \velocityarray 
			\right|^2_{\gnulltori}
				\\
		& \ \
			+
			\frac{\fundbootsmall^2}{\upmu}
			\left|
				\comdersmall^{[1,N];\le 1} \controlvars 
			\right|^2,
	\end{split}
				\\
	\begin{split} \label{E:PRINCIPALSIGAMTILDEELLIPTICHYPERBOLICENTROPYGRADIENTPOINTWISE}
		\weight^{- \blowupratetoporderacoustic}  \left|
			\, ^{(\tander^N \GradEnt)}\Currentboundaryerrorhavetocontrolprincipalweighted{\blowupratetoporderacoustic} \right|
		& 
		\lesssim
				\upmu
			\left|
				\tander^N \DivGradEnt
			\right|^2
			+
			\frac{1}{\upmu}
			\left|
				\tander^{N} \GradEnt
			\right|^2
				 \\
		& \ \
			+
			\frac{\fundbootsmall^2}{\upmu}
			\left|	
				\newuL \tander^{N} \velocityarray 
			\right|^2
			+
			\fundbootsmall^2 
				\upmu
			\left| \nullangrmD
				\tander^{N} \velocityarray 
			\right|_{\gnulltori}^2
	\end{split} 
	\end{align}
	\end{subequations}
	
	\begin{subequations}
	\begin{align}
		\begin{split}
			\weight^{- \blowupratetoporderacoustic}  \left| 
			\, ^{(\tander^N \vortrenormalized)} \Currentboundaryerrorhavetocontrollowerorderweighted{\blowupratetoporderacoustic}
		\right|
			& \lesssim \frac{1}{\upmu}  |\tander^{\le N-1}(\VortVort,\DivGradEnt)|_g^2 +
			\frac{1}{\upmu^{3/2}} |\tander^{\le N}(\vortrenormalized,\GradEnt)|_g^2 + \frac{\fundbootsmall^2}{\upmu} \left|\newuL \tander^{[1,N-1]} \velocityarray\right|^2 + \frac{\fundbootsmall^2}{\upmu} \left|\tander^{[1,N]} \velocityarray\right|^2 \\
			& \ \  + \frac{\fundbootsmall^2}{\upmu} \left| \tandersmall^{[1,N-1]} \badcontrolvars\right|^2 + \frac{\fundbootsmall^2}{\upmu} \left| \comdersmall^{[1,N-1];\le 1} \controlvars\right|^2,
			\label{E:LOWERORDERSIGAMTILDEELLIPTICHYPERBOLICVORTICITYPOINTWISE}	
		\end{split} \\
		\begin{split}				
			\weight^{- \blowupratetoporderacoustic} 	\left| \, ^{(\tander^N \GradEnt)} \Currentboundaryerrorhavetocontrollowerorderweighted{\blowupratetoporderacoustic} \right|
			& \lesssim \frac{1}{\upmu}  |\tander^{\le N-1}(\VortVort,\DivGradEnt)|_g^2 +
			\frac{1}{\upmu^{3/2}} |\tander^{\le N}(\vortrenormalized,\GradEnt)|_g^2 + \frac{\fundbootsmall^2}{\upmu} \left|\newuL \tander^{[1,N-1]} \velocityarray\right|^2 + \frac{\fundbootsmall^2}{\upmu} \left|\tander^{[1,N]} \velocityarray\right|^2 \\
			& \ \  + \frac{\fundbootsmall^2}{\upmu} \left| \tandersmall^{[1,N-1]} \badcontrolvars\right|^2 + \frac{\fundbootsmall^2}{\upmu} \left| \comdersmall^{[1,N-1];\le 1} \controlvars\right|^2.
			\label{E:LOWERORDERSIGAMTILDEELLIPTICHYPERBOLIENTROPYGRADIENTPOINTWISE}	
		\end{split}
	\end{align}
	\end{subequations}

\end{proposition}

\begin{proof}

\medskip \hfill

\noindent \textbf{Proof of \eqref{E:ANTISYMMETRICVORTICITYELLIPTICHYPERBOLICIDENTITYERORTERMPOINTWISE}--\eqref{E:ANTISYMMETRICENTROPYGRADIENTELLIPTICHYPERBOLICIDENTITYERORTERMPOINTWISE}:}

Let $\SigmatTan = \tander^N \vortrenormalized$ or $\tander^N \GradEnt$. From \eqref{E:SCHEMATICSTRUCTUREOFMETRICETC} and the bootstrap assumptions, we may bound the RHS\,\eqref{E:ANTISYMMETRICNULLCURRENTSPACETIMERRORTERM} as: 
	\begin{align}	\label{E:ANTISYMMETRICVORTICITYORGRADENTELLIPTICHYPERBOLICIDENTITYERORTERMPOINTWISESTEP1}
		\left| \mathfrak{J}_{(\textnormal{Antisymmetric)}}[\pmb{\partial} \SigmatTan,\pmb{\partial} \SigmatTan] \right| \lesssim \left(\sum_{\alpha,\beta =0}^3 \left|(\mathrm{d}\SigmatTan_\flat)_{\alpha\beta}\right|\right)^2  \left(\sum_{\alpha,\beta=0}^3 \left|\Nullhypersurfaceproject_\beta^{\ \alpha}\right|\right)^2.
	\end{align}
From \eqref{E:DEFOFINNERPRODUCTOFLUNITANDULUNIT}, \eqref{E:CONVENIENTIDENTITYFORULUNIT} and Cor.\,\ref{C:PRELIMINARYESTIMATESFORTHEDOUBLENULLACOUSTICSCALARS}, we have that $\frac{1}{\MagnitueofinnerproductofLunitanduLunit} \approx 1$, $|\uLunit^\alpha| \lesssim 1$ and hence \eqref{E:NULLHYPERSURFACEPROJECTIONINTERMSOFUNITL} implies:
	\begin{align}
		\sum_{\alpha,\beta = 0}^3 |\Nullhypersurfaceproject_\beta^{\ \alpha}| \lesssim 1. \label{E:ANTISYMMETRICVORTICITYORGRADENTELLIPTICHYPERBOLICIDENTITYERORTERMPOINTWISESTEP2}
	\end{align}\noindent From \eqref{E:ANTISYMMETRICVORTICITYORGRADENTELLIPTICHYPERBOLICIDENTITYERORTERMPOINTWISESTEP2}, we may insert \eqref{E:SHARPDECOMPOSITIONOFTOPORDEREXTERIORDERIVATIVEOFVORT} (and the respective error terms \eqref{E:SHARPDECOMPOSITIONOFTOPORDEREXTERIORDERIVATIVEOFVORTFIRSTERRORTERM}--\eqref{E:SHARPDECOMPOSITIONOFTOPORDEREXTERIORDERIVATIVEOFVORTSECONDERRORTERM})   and \eqref{E:POINTWISEESTIMATESFORMUWEIGHTEDCARTESIANCOMMUTATOR}  for $\SigmatTan = \tander^N \vortrenormalized$ into the right hand side of \eqref{E:ANTISYMMETRICVORTICITYORGRADENTELLIPTICHYPERBOLICIDENTITYERORTERMPOINTWISESTEP1}. Straightforward applications of Young's inequality, Lemma\,\ref{L:RELATIONSHIPBETWEENCARTESIANPARTIALDERIVATIVESANDSMOOTHGEOMETRICCOMMUTATORS} (which implies $\pmb{\p} \tander^N \velocityarray = \frac{1}{\upmu} \smoothfunction(\controlvars) \muX \tander^N \velocityarray + \smoothfunction(\controlvars) \Lunit \tander^N \velocityarray + \smoothfunction^A(\controlvars) \Yvf{A} \tander^N \velocityarray$), the bound $|\p_\alpha \gfour_{\beta\gamma}| \lesssim |\pmb{\p}\wavearray| \lesssim \frac{1}{\upmu}$,
and \eqref{E:POINTWISESTIMATEFORBREVEXINTERMSOFNEWUL}, collectively imply the desired estimate \eqref{E:ANTISYMMETRICVORTICITYELLIPTICHYPERBOLICIDENTITYERORTERMPOINTWISE}. 

Similar arguments hold for \eqref{E:ANTISYMMETRICENTROPYGRADIENTELLIPTICHYPERBOLICIDENTITYERORTERMPOINTWISE} using \eqref{E:SHARPDECOMPOSITIONOFTOPORDEREXTERIORDERIVATIVEOFGRADENT} and \eqref{E:SHARPDECOMPOSITIONOFTOPORDEREXTERIORDERIVATIVEOFGRADENTFIRSTERRORTERM}--\eqref{E:SHARPDECOMPOSITIONOFTOPORDEREXTERIORDERIVATIVEOFGRADENTSECONDERRORTERM}. We omit the details and highlight only the fact that there is no $\tander^N\VortVort$ term on RHS\,\eqref{E:ANTISYMMETRICENTROPYGRADIENTELLIPTICHYPERBOLICIDENTITYERORTERMPOINTWISE} due to its lack in \eqref{E:SHARPDECOMPOSITIONOFTOPORDEREXTERIORDERIVATIVEOFGRADENT} and \eqref{E:SHARPDECOMPOSITIONOFTOPORDEREXTERIORDERIVATIVEOFGRADENTFIRSTERRORTERM}--\eqref{E:SHARPDECOMPOSITIONOFTOPORDEREXTERIORDERIVATIVEOFGRADENTSECONDERRORTERM}.
			
\medskip

\noindent \textbf{Proof of \eqref{E:DIVVORTICITYELLIPTICHYPERBOLICIDENTITYERORTERMPOINTWISE}--\eqref{E:DIVENTROPYGRADIENTELLIPTICHYPERBOLICIDENTITYERORTERMPOINTWISE}:}

To prove \eqref{E:DIVVORTICITYELLIPTICHYPERBOLICIDENTITYERORTERMPOINTWISE}, we insert \eqref{E:TOPORDERCOMMUTEDDIVERGENCEOFVORT} into the RHS\,\eqref{E:DIVERGENCENULLCURRENTSPACETIMERRORTERM}, use Prop.\,\ref{P:ENTROPYDERIVATIVESINTERMSOFOTHERS}, and similar arguments used to prove \eqref{E:ANTISYMMETRICVORTICITYELLIPTICHYPERBOLICIDENTITYERORTERMPOINTWISE}--\eqref{E:ANTISYMMETRICENTROPYGRADIENTELLIPTICHYPERBOLICIDENTITYERORTERMPOINTWISE} to complete the estimate. 

Nearly identical arguments imply \eqref{E:DIVENTROPYGRADIENTELLIPTICHYPERBOLICIDENTITYERORTERMPOINTWISE}, and we highlight only that the $\tander^N\DivGradEnt$ present on RHS\,\eqref{E:TOPORDERCOMMUTEDDIVOFGRADENT} produces the top-order terms on RHS\,\eqref{E:DIVENTROPYGRADIENTELLIPTICHYPERBOLICIDENTITYERORTERMPOINTWISE}. 

\medskip

\noindent \textbf{Proof of \eqref{E:DERIVATIVEOFONEOVERMUWEIGHTVORTICITYELLIPTICHYPERBOLICIDENTITYERORTERMPOINTWISE}--\eqref{E:DERIVATIVEOFONEOVERMUWEIGHTENTROPYGRADIENTELLIPTICHYPERBOLICIDENTITYERORTERMPOINTWISE},  \eqref{E:SECONDABSORBTERMVORTICITYELLIPTICHYPERBOLICIDENTITYERORTERMPOINTWISE}--\eqref{E:SECONDABSORBTERMENTROPYGRADIENTELLIPTICHYPERBOLICIDENTITYERORTERMPOINTWISE}, and \eqref{E:NULLGEOMETRYVORTICITYELLIPTICHYPERBOLICIDENTITYERORTERMPOINTWISE}--\eqref{E:NULLGEOMETRYENTROPYGRADIENTELLIPTICHYPERBOLICIDENTITYERORTERMPOINTWISE}:}

From the definition of \eqref{E:DERIVATIVEOFWEIGHTNULLCURRENTSPACETIMERRORTERM} and the Leibniz rule we deduce the identity: 
\begin{align}\label{E:FIRSTSTEPPOINTWISEESTIMATEFORDERIVATIVEOFONEOVERMUWEIGHTEHIDENTITYBULKTERM}
	\upmu \weight^{-\blowupratetoporderacoustic} \mathfrak{J}_{(\pmb{\partial} \frac{1}{\upmu} \weight^{\blowupratetoporderacoustic})} [\SigmatTan,\pmb{\partial} \SigmatTan] = \frac{1}{\upmu} \ehcurrent^{\alpha}[\SigmatTan,\pmb{\partial} \SigmatTan] \partial_{\alpha} \upmu - \frac{\blowupratetoporderacoustic}{\weight} \ehcurrent^{\alpha}[\SigmatTan,\pmb{\partial} \SigmatTan] \partial_{\alpha}\weight.
\end{align}
Using the tensorial decomposition of the characteristic current $\ehcurrent[\SigmatTan,\pmb{\p}\SigmatTan] = - \left(\ehcurrent^\alpha [\SigmatTan,\pmb{\p}\SigmatTan]X_\alpha\right) \Lunit + \smoothtorusproject \ehcurrent[\SigmatTan,\pmb{\p}\SigmatTan]$ where $ \smoothtorusproject \ehcurrent[\SigmatTan,\pmb{\p}\SigmatTan]$ is the $\ell_{t,u}$-projection of $\ehcurrent[\SigmatTan,\pmb{\p}\SigmatTan]$ (see Def.\,\ref{D:PROJECTIONTENSORFIELDSANDTANGENCYTOHYPERSURFACES} and Lemma\,\ref{L:COVARIANTDIVERGENCEOFSPACETIMEVECTORFIELDINTERMSOFRESCALEDFRAME}), and the $\gtorus$ Cauchy--Schwarz inequality,
we deduce the following pointwise bound: 
\begin{align} \label{E:FIRSTPOINTWISEBOUNDFORCHARCURRENTDERIVATIVEOF1OVERMUWEIGHTERRORTERM}
		 \frac{1}{\upmu} \left|\ehcurrent^{\alpha}[\SigmatTan,\pmb{\partial} \SigmatTan] \partial_{\alpha} \upmu\right|
		 \lesssim
		\frac{1}{\upmu} \left|\Lunit \upmu\right| \left|X_\alpha \ehcurrent^{\alpha}[\SigmatTan,\pmb{\partial} \SigmatTan]\right| 
		+ \frac{1}{\upmu}|\angD \upmu|_{\gtorus}\left|\smoothtorusproject \ehcurrent[\SigmatTan,\pmb{\partial} \SigmatTan]\right|_{\gtorus}.
	\end{align}\noindent To control the contraction $X_\alpha \ehcurrent^\alpha[\SigmatTan,\pmb{\p}\SigmatTan]$, we use the fact that $X$ is $\Sigma_t$-tangent, \eqref{E:ULUNITINTERMSOFBANDSIGMATTAN}--\eqref{E:SIZEOFLUBITPERPANDULUNITPERPANDTHEIRINNERPRODUCT}, and \eqref{E:MAGNITUDEOFINNERPRODUCTOFNEWLANDNEWULAPPROXIMATELYMU}--\eqref{E:RECIPROCALULUNITAPPLIEDTOTIMEFUNCTIONAPPROXIMATELYMU} to bound $\frac{1}{\MagnitueofinnerproductofnewLandnewuL}| X_\alpha \newuL^\alpha| \lesssim \frac{\ReciprocaluLunitAppliedtoTimeFunction}{\MagnitueofinnerproductofnewLandnewuL} |X|_g |\uLunit^\perp| \approx 1$. From this, the identity \eqref{E:RELATIONSHIPBETWEENFLATANDCHARACTERISTICCURRENTS}, and $|\uLunit^\alpha| \lesssim 1$, we find that $ \left|X_\alpha \ehcurrent^{\alpha}[\SigmatTan,\pmb{\partial} \SigmatTan]\right| \lesssim |\SigmatTan|_g\left(\sum_{\alpha=0}^3|\partial_\alpha \SigmatTan|_g\right)$. The second term on RHS\,\eqref{E:FIRSTPOINTWISEBOUNDFORCHARCURRENTDERIVATIVEOF1OVERMUWEIGHTERRORTERM} is similarly bounded using the last estimate in \eqref{E:SMOOTHTORUSNORMCOMPARBLETOTANGENTIALCONTRACTIONS}, Prop.\,\ref{P:SCHEMATICSTRUCTUREOFVARIOUSTENSORSINTERMSOFCONTROLVARS},
the estimates of Prop.\,\ref{P:IMPROVEMENTOFAUXILIARYBOOTSTRAP}. Altogether, from  \eqref{E:FIRSTPOINTWISEBOUNDFORCHARCURRENTDERIVATIVEOF1OVERMUWEIGHTERRORTERM}, \eqref{E:COERCIVENESSOFELLIPTICHYPERBOLICQUADRATICFORM}, and Young's inequality, we deduce the following pointwise bound:
	\begin{align} \label{E:SECONDPOINTWISEBOUNDFORCHARCURRENTDERIVATIVEOF1OVERMUWEIGHTERRORTERM}
		 \frac{1}{\upmu} \left|\ehcurrent^{\alpha}[\SigmatTan,\pmb{\partial} \SigmatTan] \partial_{\alpha} \upmu\right|
		& 
		\lesssim
		\frac{1}{\upmu}
		|\SigmatTan|_g
		\left(\sum_{\alpha=0}^3|\partial_\alpha \SigmatTan|_g\right) \lesssim \frac{\varsigma}{\Lunit \ubar} \ellipticCoerciveQuadratic[\pmb{\partial} \SigmatTan,\pmb{\partial} \SigmatTan] + \frac{1}{\varsigma \upmu^2} |\SigmatTan|_g^2.
	\end{align}
Using similar ideas, but instead relying on the bounds $|\Lunit \weight| = |\Lunit\ubar| \approx 1$ and $|\angD \weight|_{\gtorus} = |\angD \ubar|_{\gtorus} \lesssim \fundbootsmall$, which follow easily from Cor.~\ref{C:ELLTUPROJECTEDVERSIONOFCARTESIANPARTIALDERIVATIVES}, \eqref{E:BLOWUPWEIGHT}, the pointwise comparison estimate \eqref{E:SMOOTHTORUSNORMCOMPARBLETOTANGENTIALCONTRACTIONS}, Lemm~\ref{L:PROPERTIESANDDIFFEOMORPHICEXTENSIONOFDOUBLENULLCOORDINATES}, and Cor.~\ref{C:IMPROVEAUX}, we compute: 
\begin{align} \label{E:THIRDPOINTWISEBOUNDFORCHARCURRENTDERIVATIVEOF1OVERMUWEIGHTERRORTERM}
		 \frac{\blowupratetoporderacoustic}{\weight} \left|\ehcurrent^{\alpha}[\SigmatTan,\pmb{\partial} \SigmatTan] \partial_{\alpha} \upmu\right|
		& 
		\lesssim
		\frac{1}{\weight}
		|\SigmatTan|_g
		\left(\sum_{\alpha=0}^3|\partial_\alpha \SigmatTan|_g\right) 
		\lesssim \varsigma \ellipticCoerciveQuadratic[\pmb{\partial} \SigmatTan,\pmb{\partial} \SigmatTan] + \frac{1}{\varsigma \weight^2} |\SigmatTan|_g^2.
	\end{align}
Using \eqref{E:FIRSTSTEPPOINTWISEESTIMATEFORDERIVATIVEOFONEOVERMUWEIGHTEHIDENTITYBULKTERM}--\eqref{E:THIRDPOINTWISEBOUNDFORCHARCURRENTDERIVATIVEOF1OVERMUWEIGHTERRORTERM} with 
$\tander^N \vortrenormalized$
and
$\tander^N \GradEnt$
in the role of $\SigmatTan$ and \eqref{E:KEYESTIMATECONTROLLINGINVERSEMUBYINVERSEWEIGHT}, we deduce the desired bounds
\eqref{E:DERIVATIVEOFONEOVERMUWEIGHTVORTICITYELLIPTICHYPERBOLICIDENTITYERORTERMPOINTWISE}--\eqref{E:DERIVATIVEOFONEOVERMUWEIGHTENTROPYGRADIENTELLIPTICHYPERBOLICIDENTITYERORTERMPOINTWISE}.

Similar arguments using the pointwise bounds $|\Chfour_{\alpha \ \beta}^{\ \alpha}|, |\p_\alpha \gfour_{\beta\gamma}| \lesssim | \pmb{\p}\wavearray|\lesssim \frac{1}{\upmu}$ prove \eqref{E:SECONDABSORBTERMVORTICITYELLIPTICHYPERBOLICIDENTITYERORTERMPOINTWISE}--\eqref{E:SECONDABSORBTERMENTROPYGRADIENTELLIPTICHYPERBOLICIDENTITYERORTERMPOINTWISE}, we omit the details. 

Next, \eqref{E:NULLGEOMETRYVORTICITYELLIPTICHYPERBOLICIDENTITYERORTERMPOINTWISE}--\eqref{E:NULLGEOMETRYENTROPYGRADIENTELLIPTICHYPERBOLICIDENTITYERORTERMPOINTWISE}  hold using similar Young inequality type-arguments as well as the pointwise bound $|\p_\alpha \Nullhypersurfaceproject_\beta^{\ \gamma}| \lesssim |\pmb{\p} \wavearray| \lesssim \frac{1}{\upmu}$.

\medskip

\noindent \textbf{Proof of \eqref{E:FIRSTABSORBTERMVORTICITYELLIPTICHYPERBOLICIDENTITYERORTERMPOINTWISE}--\eqref{E:FIRSTABSORBTERMENTROPYGRADIENTELLIPTICHYPERBOLICIDENTITYERORTERMPOINTWISE}:}

From \eqref{E:ANTISYMMETRICVORTICITYORGRADENTELLIPTICHYPERBOLICIDENTITYERORTERMPOINTWISESTEP2} and $|\p_\alpha \gfour_{\nu\gamma'} |\lesssim \frac{1}{\upmu}$ we bound the first term on RHS\,\eqref{E:ABSORBABLEANTISYMMETRICANDDIVERGENCENULLCURRENTSPACETIMERRORTERM} using Young's inequality as:
\begin{align}
		\left|  \Nullhypersurfaceproject_{\beta}^{\ \alpha} \Nullhypersurfaceproject_{\gamma}^{\ \delta}  \SigmatTan^{\nu} \p_\alpha \gfour_{\nu \gamma'} (\mathrm{d} \SigmatTan_\flat)_{\delta\kappa}  (\gfour^{-1})^{\beta\kappa} (\gfour^{-1})^{\gamma \gamma'} \right|  \lesssim \varsigma \ellipticCoerciveQuadratic[\pmb{\partial} \SigmatTan,\pmb{\partial} \SigmatTan] + \frac{1}{\varsigma \upmu^2} |\SigmatTan|_g^2. \label{E:FIRSTABSORBTERMVORTICITYANDENTROPYELLIPTICHYPERBOLICIDENTITYERORTERMPOINTWISESTEP1}
	\end{align}We now focus on the second line on RHS\,\eqref{E:ABSORBABLEANTISYMMETRICANDDIVERGENCENULLCURRENTSPACETIMERRORTERM}.  First, note that by $\frac{\ReciprocaluLunitAppliedtoTimeFunction \ReciprocalLunitAppliedtoTimeFunction}{\MagnitueofinnerproductofnewLandnewuL} = \frac{\ReciprocaluLunitAppliedtoTimeFunction}{\upmu} \approx 1$ (see \eqref{E:RECIPROCALULUNITAPPLIEDTOTIMEFUNCTIONAPPROXIMATELYMU}), \eqref{E:LUNITINTERMSOFBANDSIGMATTAN}--\eqref{E:SIZEOFLUBITPERPANDULUNITPERPANDTHEIRINNERPRODUCT}, \eqref{E:QUADRATICFORMCOERCIVEINNULLANGULAR}, and \eqref{E:COERCIVENESSOFELLIPTICHYPERBOLICQUADRATICFORM}, we have that the first term in braces, namely $\left\{  \frac{6-\sqrt{6}}{6}\p_a \SigmatTan^a + \cdots +  \frac{\ReciprocaluLunitAppliedtoTimeFunction \ReciprocalLunitAppliedtoTimeFunction}{2\MagnitueofinnerproductofnewLandnewuL} \left(\frac{3-\sqrt{6}}{6}\right) \uLunit^\perp_a \Lunit^\perp \SigmatTan^a \right\}$, is bounded in magnitude by $C\sqrt{\ellipticCoerciveQuadratic[\pmb{\partial} \SigmatTan,\pmb{\partial} \SigmatTan]}$. When $\SigmatTan = \tander^N \vortrenormalized$, dividing  \eqref{E:COMMUTEDTRANSPORTPOINTWISEESTIMATESFORSPECIFICVORTICITYANDENTROPYGRADIENT} by $\upmu$ and inserting the resulting estimate for $\Transport \tander^N \vortrenormalized$, as well as inserting \eqref{E:COMMUTEDPOINTWISEEUCLIDEANCURLOFVORTICITY} for the $\Flatcurl \tander^N\vortrenormalized$, straightforward applications of Young's inequality imply that the second line is $\lesssim$ RHS\,\eqref{E:FIRSTABSORBTERMVORTICITYELLIPTICHYPERBOLICIDENTITYERORTERMPOINTWISE}.

Next, when $\SigmatTan = \tander^N \GradEnt$, we use Young's inequality after inserting \eqref{E:COMMUTEDTRANSPORTPOINTWISEESTIMATESFORSPECIFICVORTICITYANDENTROPYGRADIENT} and \eqref{E:COMMUTEDPOINTWISEEUCLIDEANDIVERGENCEOFVORTICITYANDEUCLIDEANCURLOFENTROPYGRADIENT} to deduce that the magnitude of the second line of \eqref{E:ABSORBABLEANTISYMMETRICANDDIVERGENCENULLCURRENTSPACETIMERRORTERM} for this case is $\lesssim$ RHS\,\eqref{E:FIRSTABSORBTERMENTROPYGRADIENTELLIPTICHYPERBOLICIDENTITYERORTERMPOINTWISE}.

Similar arguments from the previous two paragraphs also prove that the last line of RHS\,\eqref{E:ABSORBABLEANTISYMMETRICANDDIVERGENCENULLCURRENTSPACETIMERRORTERM} is $\lesssim$ \eqref{E:FIRSTABSORBTERMVORTICITYELLIPTICHYPERBOLICIDENTITYERORTERMPOINTWISE}--\eqref{E:FIRSTABSORBTERMENTROPYGRADIENTELLIPTICHYPERBOLICIDENTITYERORTERMPOINTWISE}, we omit the repetitive details.

\medskip

\noindent \textbf{Proof of \eqref{E:QUADRATICINMATERIALDERIVATIVETERMVORTICITYELLIPTICHYPERBOLICIDENTITYERORTERMPOINTWISE}--\eqref{E:QUADRATICINMATERIALDERIVATIVETERMENTROPYGRADIENTELLIPTICHYPERBOLICIDENTITYERORTERMPOINTWISE}:}

Dividing  \eqref{E:COMMUTEDTRANSPORTPOINTWISEESTIMATESFORSPECIFICVORTICITYANDENTROPYGRADIENT} by $\upmu$, it follows that $|\Transport(\tander^N(\vortrenormalized,\GradEnt))|_g^2 \lesssim$ RHS\,\eqref{E:QUADRATICINMATERIALDERIVATIVETERMVORTICITYELLIPTICHYPERBOLICIDENTITYERORTERMPOINTWISE}--\eqref{E:QUADRATICINMATERIALDERIVATIVETERMENTROPYGRADIENTELLIPTICHYPERBOLICIDENTITYERORTERMPOINTWISE}.

By the Cauchy--Schwarz inequality with respect to $g$, since $\newuL^\perp = \ReciprocaluLunitAppliedtoTimeFunction\uLunit^\perp$, \eqref{E:SIZEOFLUBITPERPANDULUNITPERPANDTHEIRINNERPRODUCT} implies that: 
\begin{align}
		\left| \frac{1}{\MagnitueofinnerproductofnewLandnewuL^2} (\Transport (\tander^N \vortrenormalized,\tander^N \GradEnt)^a) \newuL^\perp_a \right| \lesssim \frac{1}{\upmu} |\Transport \tander^N(\vortrenormalized,\GradEnt)|_g.\label{E:QUADRATICINMATERIALDERIVATIVETERMVORTICITYANDENTROPYELLIPTICHYPERBOLICIDENTITYERORTERMPOINTWISESTEP1}
	\end{align}
Inserting \eqref{E:QUADRATICINMATERIALDERIVATIVETERMVORTICITYANDENTROPYELLIPTICHYPERBOLICIDENTITYERORTERMPOINTWISESTEP1} and subsequently \eqref{E:COMMUTEDTRANSPORTPOINTWISEESTIMATESFORSPECIFICVORTICITYANDENTROPYGRADIENT} into the last two terms on RHS\,\eqref{E:QUADRATICINMATERIALDERIVATIVESNULLCURRENTSPACETIMERRORTERM}, and applying Young's inequality, we conclude the proof.

\medskip

\noindent \textbf{Proof of \eqref{E:PRINCIPALSIGAMTILDEELLIPTICHYPERBOLICVORTICITYPOINTWISE}--\eqref{E:PRINCIPALSIGAMTILDEELLIPTICHYPERBOLICENTROPYGRADIENTPOINTWISE}:}

Let $\SigmatTan = \tander^N\vortrenormalized$. Due to the factor of $\weight^{-\blowupratetoporderacoustic}$ in LHS \eqref{E:PRINCIPALSIGAMTILDEELLIPTICHYPERBOLICVORTICITYPOINTWISE}--\eqref{E:PRINCIPALSIGAMTILDEELLIPTICHYPERBOLICENTROPYGRADIENTPOINTWISE}, it suffices to bound $ \frac{1}{\upmu}  \, ^{(\SigmatTan)}\Currentboundaryerrorhavetocontrolprincipal^{(\textnormal{Flat})} + 2\frac{\ingoingmu}{\upmu \MagnitueofinnerproductofnewLandnewuL}   \,^{(\SigmatTan )}\Currentboundaryerrorhavetocontrolprincipal^{(\textnormal{Non-flat})} $, where $\, ^{(\SigmatTan)}\Currentboundaryerrorhavetocontrolprincipal^{(\textnormal{Flat})}$ is given in \eqref{E:PRINCIPALERRORTERMHAVETOCONTROLKEYIDPUTANGENTCURRENTCONTRACTEDAGAINSTVECTORFIELDFORVORTICITY}.

To bound  $\frac{1}{\upmu}  \, ^{(\SigmatTan)}\Currentboundaryerrorhavetocontrolprincipal^{(\textnormal{Flat})} $, we first focus on all of the terms in RHS\,\eqref{E:PRINCIPALERRORTERMHAVETOCONTROLKEYIDPUTANGENTCURRENTCONTRACTEDAGAINSTVECTORFIELDFORVORTICITY} featuring $\newuL \tander^N v^i$ or $\newuL \tander^N \LogDensity$ that \emph{also} feature factors of $\frac{1}{\upmu}$ or $\frac{1}{\MagnitueofinnerproductofnewLandnewuL}$. Observe that all of these terms involve contractions of $\newuL =\ReciprocaluLunitAppliedtoTimeFunction \uLunit $ which effectively cancel the singular factors by \eqref{E:MAGNITUDEOFINNERPRODUCTOFNEWLANDNEWULAPPROXIMATELYMU}--\eqref{E:RECIPROCALULUNITAPPLIEDTOTIMEFUNCTIONAPPROXIMATELYMU}. Moreover, all of these terms feature small factors of \emph{undifferentiated} $\GradEnt$ or $\vortrenormalized$, which are bounded in magnitude by $\fundbootsmall$. Hence, by Young's inequality, these are all bounded in magnitude by:
	\begin{align}
		\lesssim \frac{\fundbootsmall}{\upmu}  \left(\left| \newuL \tander^N \velocityarray \right| + \left| \newuL \tander^N \LogDensity\right|\right) \left|\tander^N \vortrenormalized\right| \lesssim \frac{\fundbootsmall^2}{\upmu}  \left(\left| \newuL \tander^N \velocityarray \right|^2 + \left| \newuL \tander^N \LogDensity\right|^2\right) + \frac{1}{\upmu} \left|\tander^N \vortrenormalized\right|^2.	\label{E:PRINCIPALSIGAMTILDEELLIPTICHYPERBOLICVORTICITYPOINTWISESTEP1}
	\end{align}
To bound the $\newuL \tander^N \LogDensity$ term, we use Prop.\,\ref{P:ENTROPYDERIVATIVESINTERMSOFOTHERS}. Specifically, since \eqref{E:LRHOINTERMSOFVANDENT} and \eqref{E:TANGENTIALDERIVATIVESOFRHOINTERMSOFVANDVORTANDENT} imply that any of $\{ \Lunit \LogDensity, \Yvf{2} \LogDensity, \Yvf{3} \LogDensity\}$ can be expressed schematically as $\smoothfunction(\controlvars) \tander \velocityarray + \smoothfunction(\controlvars) (\vortrenormalized,\GradEnt)$. Since $\newuL \tander^N = \newuL \tander^{N-1}\Singletan $ for some $\Singletan \in \{\Lunit,\Yvf{2},\Yvf{3}\}$, the RHS\,\eqref{E:PRINCIPALSIGAMTILDEELLIPTICHYPERBOLICVORTICITYPOINTWISESTEP1} $\lesssim$  \eqref{E:PRINCIPALSIGAMTILDEELLIPTICHYPERBOLICVORTICITYPOINTWISE}. The remaining term on RHS\,\eqref{E:PRINCIPALERRORTERMHAVETOCONTROLKEYIDPUTANGENTCURRENTCONTRACTEDAGAINSTVECTORFIELDFORVORTICITY} featuring $\newuL \tander^N v^i$, namely $\frac{1}{\ingoingmu} (\tander^N \vortrenormalized^\beta) (\vortrenormalized_a) \newL_\beta \newuL \tander^N v^a$ does not feature the regularizing contraction of $\newuL$. However, $\frac{1}{\ingoingmu} \approx 1$ and so this term is also bounded in magnitude by $\lesssim$ RHS\,\eqref{E:PRINCIPALSIGAMTILDEELLIPTICHYPERBOLICVORTICITYPOINTWISE}

The terms featuring null-angular directions such as $\nullangpartial_i \tander^N v^a$ or $\breve{\slashed{\Transport}}\tander^N v^a$ all feature the additional regular factor of $\ReciprocaluLunitAppliedtoTimeFunction \approx \upmu$ from the $\newuL$ contractions as well as the small factors of undifferentiated $\GradEnt$ or $\vortrenormalized$. We therefore bound these as:
	\begin{align} \label{E:PRINCIPALSIGAMTILDEELLIPTICHYPERBOLICVORTICITYPOINTWISESTEP2}
		\lesssim \fundbootsmall \frac{\ReciprocaluLunitAppliedtoTimeFunction}{\upmu}  \left|  \nullangrmD \tander^N v^a\right|_{\gnulltori} \left| \tander^N \vortrenormalized\right| \lesssim \fundbootsmall^2 \upmu \left| \nullangrmD \tander^N \velocityarray\right|^2 + \frac{1}{\upmu} \left| \tander^N \vortrenormalized\right|^2.
	\end{align}
The terms featuring $\nullangpartial_i \tander^N \LogDensity$ or $\breve{\slashed{\Transport}}\tander^N \LogDensity$ are also $\lesssim$ RHS\,\eqref{E:PRINCIPALSIGAMTILDEELLIPTICHYPERBOLICVORTICITYPOINTWISE} using similar arguments as in the previous paragraph.
		
The remaining term on RHS\,\eqref{E:PRINCIPALERRORTERMHAVETOCONTROLKEYIDPUTANGENTCURRENTCONTRACTEDAGAINSTVECTORFIELDFORVORTICITY} is the one featuring $\tander^N \VortVort^\kappa$. Since this also features the regular $\newuL$ contraction which cancels the $\frac{1}{\upmu}$, we bound these as: 
	\begin{align}
		\lesssim \upmu |\tander^N \VortVort|^2 + \frac{1}{\upmu}|\tander^N \vortrenormalized|^2. \label{E:PRINCIPALSIGAMTILDEELLIPTICHYPERBOLICVORTICITYPOINTWISESTEP3}
	\end{align}Strictly easier arguments hold in the case for $\SigmatTan = \tander^N \GradEnt$ using instead \eqref{E:CONTRACTIONOFFLATCURRENTFORGRADENTANDNEWUL}, whose terms all feature either the same regularizing contractions of $\newuL$ whenever there is a singular factor of $\frac{1}{\upmu}$ or $\frac{1}{\MagnitueofinnerproductofnewLandnewuL}$ when $\newuL \tander^N v^a$ is present, or feature an overall factor of $\ReciprocaluLunitAppliedtoTimeFunction \approx \upmu$ otherwise. We say the arguments are strictly easier because there are no top-order derivatives of $\LogDensity$ on RHS\,\eqref{E:CONTRACTIONOFFLATCURRENTFORGRADENTANDNEWUL}. This concludes our estimates of $ \frac{1}{\upmu}  \, ^{(\SigmatTan)}\Currentboundaryerrorhavetocontrolprincipal^{(\textnormal{Flat})} $ for both $\SigmatTan = \tander^N \vortrenormalized$ or $\tander^N \GradEnt$.

To control the remaining $2\frac{\ingoingmu}{\upmu \MagnitueofinnerproductofnewLandnewuL}   \,^{(\SigmatTan )}\Currentboundaryerrorhavetocontrolprincipal^{(\textnormal{Non-flat})} $ terms, we observe that the overall coefficient is of order $\frac{\ingoingmu}{\upmu \MagnitueofinnerproductofnewLandnewuL} \approx \frac{1}{\upmu^2}$. However, from \eqref{E:MIRACULOUSSTRUCTUREINTHEVORTICITYBOUNDARYERRORTERMPRINCIPAL} and \eqref{E:MIRACULOUSSTRUCTUREINTHEGRADIENTENTROPYBOUNDARYERRORTERMPRINCIPAL}, we see that these terms these terms all feature \emph{an additional contraction against another copy of $\newuL$} compared to \eqref{E:PRINCIPALERRORTERMHAVETOCONTROLKEYIDPUTANGENTCURRENTCONTRACTEDAGAINSTVECTORFIELDFORVORTICITY}--\eqref{E:PRINCIPALERRORTERMHAVETOCONTROLKEYIDPUTANGENTCURRENTCONTRACTEDAGAINSTVECTORFIELDFORGRADENT}. Hence, nearly identical arguments as in \eqref{E:PRINCIPALSIGAMTILDEELLIPTICHYPERBOLICVORTICITYPOINTWISESTEP1}--\eqref{E:PRINCIPALSIGAMTILDEELLIPTICHYPERBOLICVORTICITYPOINTWISESTEP3} hold; we omit the details.

\medskip

\noindent \textbf{Proof of \eqref{E:LOWERORDERSIGAMTILDEELLIPTICHYPERBOLICVORTICITYPOINTWISE}--\eqref{E:LOWERORDERSIGAMTILDEELLIPTICHYPERBOLIENTROPYGRADIENTPOINTWISE}:}

We begin by multiplying \eqref{E:LOWERORDERERRORTERMHAVETOCONTROLKEYIDPUTANGENTCURRENTCONTRACTEDAGAINSTVECTORFIELDFORVORTICITYWEIGHTED} by $\weight^{-\blowupratetoporderacoustic}$, and show first that $\frac{1}{\upmu} \, ^{(\SigmatTan)}\Currentboundaryerrorhavetocontrollowerorder^{(\textnormal{Flat})}  + 2 \frac{\ingoingmu}{\upmu \MagnitueofinnerproductofnewLandnewuL}  \,^{(\SigmatTan)}\Currentboundaryerrorhavetocontrollowerorder^{(\textnormal{Non-flat})}$ are $\lesssim$ RHS\,\eqref{E:LOWERORDERSIGAMTILDEELLIPTICHYPERBOLICVORTICITYPOINTWISE}--\eqref{E:LOWERORDERSIGAMTILDEELLIPTICHYPERBOLIENTROPYGRADIENTPOINTWISE}. 

First, note that by expressing $\upmu\p_\alpha$ in terms of the commutator vectorfields (see Lemma\,\ref{L:RELATIONSHIPBETWEENCARTESIANPARTIALDERIVATIVESANDSMOOTHGEOMETRICCOMMUTATORS}),
the commutator terms on RHS\,\eqref{E:LOWERORDERERRORTERMHAVETOCONTROLKEYIDPUTANGENTCURRENTCONTRACTEDAGAINSTVECTORFIELDFORVORTICITY},  \eqref{E:LOWERORDERERRORTERMHAVETOCONTROLKEYIDPUTANGENTCURRENTCONTRACTEDAGAINSTVECTORFIELDFORGRADENT},  \eqref{E:MIRACULOUSSTRUCTUREINTHEVORTICITYBOUNDARYERRORTERMLOWERORDER}, and \eqref{E:MIRACULOUSSTRUCTUREINTHEGRADIENTENTROPYBOUNDARYERRORTERMLOWERORDER} (e.g. $\vortrenormalized^a [\upmu \p_d,\tander^N] v^a$ in $\, ^{(\tander^N)}\ehcurrent_{(\vortrenormalized;1)}^d $, see \eqref{E:SHARPDECOMPOSITIONOFTOPORDEREXTERIORDERIVATIVEOFVORTFIRSTERRORTERM}) are schematically of the form  $\tander^{[1,N]} \velocityarray $ up to harmless error terms. These are in turn bounded by $|\nullangrmD \tander^{[1,N-1]} \velocityarray|_{\gnulltori} + |\Lunit \tander^{[1,N-1]} \velocityarray|$ by \eqref{E:POINTWISESMOOTHTORIDIFFERENTIALNORMINTERMSOFDOULBENULLTORIDIFFERENTIALINORMANDL}. All this to say, these lower-order commutator error terms do not exclusively feature derivatives of $v^i$ or $\LogDensity$ whose last derivative is tangent to $\ingoingcharacteristicsurfacetwoarg{\ubar}{[\moreinterestingu_1,u]}$-tangential derivatives as the last differential operator, as did the principal terms in the proof of \eqref{E:PRINCIPALSIGAMTILDEELLIPTICHYPERBOLICVORTICITYPOINTWISE}--\eqref{E:PRINCIPALSIGAMTILDEELLIPTICHYPERBOLICENTROPYGRADIENTPOINTWISE}. We highlight then the fact that \emph{all} such lower order terms on RHS\,\eqref{E:LOWERORDERERRORTERMHAVETOCONTROLKEYIDPUTANGENTCURRENTCONTRACTEDAGAINSTVECTORFIELDFORVORTICITY},  \eqref{E:LOWERORDERERRORTERMHAVETOCONTROLKEYIDPUTANGENTCURRENTCONTRACTEDAGAINSTVECTORFIELDFORGRADENT} which feature a singular $\frac{1}{\upmu}$ also feature a contraction against $\newuL = \ReciprocaluLunitAppliedtoTimeFunction\uLunit$, which cancel up to a constant. Moreover, \eqref{E:POINTWISEESTIMATESFORMUWEIGHTEDCARTESIANCOMMUTATOR} imples that all terms which are not purely quadratic in $|\tander^{\le N}(\vortrenormalized,\GradEnt)|$ (e.g. those terms featuring $ \,^{(\tander^N)}\ehcurrent_{(\vortrenormalized;1)}^i \,^{(\tander^N)}\ehcurrent_{(\vortrenormalized;2)}^i,  \,^{(\tander^N)}\ehcurrent_{(\vortrenormalized;3)}$ or $ \,^{(\tander^N)}\ehcurrent_{(\GradEnt;1)}^i \,^{(\tander^N)}\ehcurrent_{(\GradEnt;2)}^i,  \,^{(\tander^N)}\ehcurrent_{(\GradEnt;3)}$, which, for example, have the aforementioned terms of the form $|\tander^{\le N} (\vortrenormalized,\GradEnt)|\cdot | \tander^{[1,N]}\velocityarray|$)  all feature an $\mathcal{O}(\fundbootsmall)$ factor. Hence, we bound: 
\begin{align}
		\begin{split} \label{E:LOWERORDERSIGAMTILDEELLIPTICHYPERBOLICVORTICITYANDENTROPYPOINTWISESTEP1}
			\left| \frac{1}{\upmu}  \, ^{(\SigmatTan)}\Currentboundaryerrorhavetocontrollowerorder^{(\textnormal{Flat})} \right| 
			& \lesssim \frac{1}{\upmu}  |\tander^{\le N-1}(\VortVort,\DivGradEnt)|_g^2 +
			\frac{1}{\upmu} |\tander^{\le N}(\vortrenormalized,\GradEnt)|_g^2 + \frac{\fundbootsmall^2}{\upmu} \left|\newuL \tander^{[1,N-1]} \velocityarray\right|^2 + \frac{\fundbootsmall^2}{\upmu} \left|\tander^{[1,N]} \velocityarray\right|^2 \\
			& \ \  + \frac{\fundbootsmall^2}{\upmu} \left| \tandersmall^{[1,N-1]} \badcontrolvars\right|^2 + \frac{\fundbootsmall^2}{\upmu} \left| \comdersmall^{[1,N-1];\le 1} \controlvars\right|^2.
		\end{split}
	\end{align}Even though $2 \frac{\ingoingmu}{\upmu \MagnitueofinnerproductofnewLandnewuL}  \,^{(\SigmatTan)}\Currentboundaryerrorhavetocontrollowerorder^{(\textnormal{Non-flat})}$ features an overall singular factor $ \frac{\ingoingmu}{\upmu \MagnitueofinnerproductofnewLandnewuL} \approx \upmu^{-2}$, nearly identical arguments as above prove that $|2 \frac{\ingoingmu}{\upmu \MagnitueofinnerproductofnewLandnewuL}  \,^{(\SigmatTan)}\Currentboundaryerrorhavetocontrollowerorder^{(\textnormal{Non-flat})}| \lesssim$ RHS\,\eqref{E:LOWERORDERSIGAMTILDEELLIPTICHYPERBOLICVORTICITYANDENTROPYPOINTWISESTEP1}. This is because the RHS\,\eqref{E:MIRACULOUSSTRUCTUREINTHEVORTICITYBOUNDARYERRORTERMLOWERORDER}, and \eqref{E:MIRACULOUSSTRUCTUREINTHEGRADIENTENTROPYBOUNDARYERRORTERMLOWERORDER} all feature an additional contraction against $\newuL$, which adds an additional overall factor of $\upmu$. We omit the repetitive details.

Finally, it remains to prove that the terms in braces on RHS\,\eqref{E:LOWERORDERERRORTERMHAVETOCONTROLKEYIDPUTANGENTCURRENTCONTRACTEDAGAINSTVECTORFIELDFORVORTICITYWEIGHTED} are $\lesssim$ RHS  \eqref{E:LOWERORDERSIGAMTILDEELLIPTICHYPERBOLICVORTICITYPOINTWISE}--\eqref{E:LOWERORDERSIGAMTILDEELLIPTICHYPERBOLIENTROPYGRADIENTPOINTWISE}. From the transport equations of Lemma\,\ref{L:QUASILINEARTRANSPORTSYSTEMFORDERIVATIVESOFUBAR}, the pointwise estimates of Lemma\,\ref{L:PROPERTIESANDDIFFEOMORPHICEXTENSIONOFDOUBLENULLCOORDINATES}, and the bound \eqref{E:POINTWISEBOUNDFORDOUBLENULLTOROIDALTRACEOFDEFORMATIONTENSOROFNEWULVECTORFIELD}, it follows that the terms $ \frac{1}{\upmu}  \SigmatTan^\alpha \SigmatTan^\beta \newuL(\uLunit_\alpha \breve{\slashed{\Transport}}_\beta)
				 + \frac{1}{2 \upmu}  \SigmatTan^\alpha \uLunit_\alpha   \SigmatTan^\beta \breve{\slashed{\Transport}}_\beta \mytr_{\gnulltori} \deform{\newuL} 
				 -   \frac{1}{4\upmu} |\SigmatTan|_{\gfour}^2 \mytr_{\gnulltori} \deform{\newuL}$ are $\lesssim \frac{1}{\upmu} |\tander^{N}(\vortrenormalized,\GradEnt)|^2$.

Finally, it suffices to bound the terms $\frac{\newuL \upmu}{2\upmu^2} |\SigmatTan|_{\gfour}^2
				+ \newuL\left( \frac{1}{\upmu}\right) \SigmatTan^\alpha \uLunit_\alpha \SigmatTan^\beta \breve{\slashed{\Transport}}_\beta 
				+ \breve{\slashed{\Transport}} \left( \frac{\ingoingmu}{\upmu \MagnitueofinnerproductofnewLandnewuL} \right) ( \SigmatTan^\alpha\newuL_\alpha)^2$. Since the last term features two contractions of $\newuL_\alpha$ and $\breve{\slashed{\Transport}}\left( \frac{\ingoingmu}{\upmu \MagnitueofinnerproductofnewLandnewuL} \right) \approx \frac{1}{\upmu^3}$, this term is bounded by $\frac{1}{\upmu}  |\tander^{N}(\vortrenormalized,\GradEnt)|^2$. From \eqref{E:NEWULMUUBOUNDEDBYSQRTMU}, the first two terms are bounded by $\frac{1}{\upmu^{\frac{3}{2}}}  |\tander^{N}(\vortrenormalized,\GradEnt)|^2$. As $\frac{1}{\upmu} \lesssim \frac{1}{\upmu^{3/2}}$, we conclude the proof of \eqref{E:LOWERORDERSIGAMTILDEELLIPTICHYPERBOLICVORTICITYPOINTWISE}--\eqref{E:LOWERORDERSIGAMTILDEELLIPTICHYPERBOLIENTROPYGRADIENTPOINTWISE} and the proposition.

\end{proof}

\section{Statement of the a priori $L^2$ estimates and  
data-estimates for the $L^2$-controlling quantities, and 
bootstrap assumptions for the wave variable energies} \label{S:SATEMENTOFAPRIORIL2ESTIMATES}

We continue to work under the assumptions of Sect.\,\ref{SS:SILENTFACTS}. In Sect.\,\ref{SS:STATEMENTOFAPRIORIL2ESTIMATESFORMASTERCONTROL}, we state the a priori $L^2$ estimates for the master controlling quantities $\mastercontroltop(\ubar,u)$ defined in Sect.\,\ref{SS:MASTERCONTROLLINGQUANTITY}.
In Sect.\,\ref{SS:STATEMENTOFAPRIORIL2ESTIMATESFORWAVEVARIABLES}, we state the a priori energy--null-flux estimates
for the wave variables. We state the corresponding energy estimates for the acoustic geometry in Sect.\,\ref{SS:STATEMENTOFAPRIORIL2ESTIMATESFORACOUSTICGEOMETRY}.
The proofs take considerable effort and form the focus of the paper
through Sect.\,\ref{S:WAVEANDACOUSTICGEOMETRYAPRIORIESTIMATES}.
In Sect.\,\ref{SS:DATAESTIMATESFORL2CONTROLLINGQUANTITIES}, 
we show that along the data hypersurfaces 
$\ingoingcharacteristicsurfacetwoarg{\ubar}{[\moreinterestingu_1,u]}$ and $\outgoingcharacteristicsurfacetwoarg{\moreinterestingu_1}{[\leftubar,0]}$, the $L^2$-controlling quantities are bounded by $\lesssim \initialsmalldoublenull^2$,
i.e., the $L^2$-controlling quantities have small data. We emphasize that all of the estimates stated in this section feature \emph{bounded} energies, a distinct feature from all prior works on shocks where the mid-to-top-order energies featured singular weights, see for example the RHS of the inequalities \cite{abbrescia2022emergence}*{(25.1a)}. This is difference is attributed to the fact that we incorporated the weights into the definition of the energies and null-fluxes, see Sect.\,\ref{S:BASICINGREDIENTSINL2ANALYSIS}.


\subsection{Statement of the top order a priori $L^2$ estimates for the master controlling quantities}
\label{SS:STATEMENTOFAPRIORIL2ESTIMATESFORMASTERCONTROL}

In the following proposition, we state the main $L^2$ a priori estimate needed at the top order. It's proof is located in Sect.\,\ref{SSS:PROOFOFTOPORDERAPRIORIL2ESTIMATEMASTERCONTROLANDALLORDERWAVES}. 

\begin{proposition}[The main a priori estimates for $\mastercontroltop$] \label{P:TOPORDERAPRIORIL2ESTIMATESMASTERCONTROL} 
	Let $\mastercontroltop(\ubar,u)$ 
	be the top order $L^2$-master controlling quantity for the wave variables 
	$\wavearray$ and acoustic geometry, as defined in \eqref{E:MASTERCONTROLLINGQUANTITY}.
	Under the data-assumptions of Sect.\,\ref{S:ASSUMPTIONSONTHEDATA},
	the parameter size-assumptions of Sect.\,\ref{SS:PARAMETERSIZEASSUMPTIONS},
	and the bootstrap assumptions of Sect.\,\ref{S:BOOTSTRAPEVERYTHINGEXCEPTENERGIES},
	there exists a constant $C > 0$
	such that the following estimates hold for 
	$(\ubar,u) \in [\leftubar,\ubarboot) \times [\moreinterestingu_1,\moreinterestingu_2]$: 
\begin{align}
	\mastercontroltop(\ubar,u)
	& \leq C \initialsmalldoublenull^2.\label{E:MAINTOPORDERL2ESTIMATESMASTERCONTROLQUANTITY}
\end{align}
\end{proposition}

\subsection{Statement of the a priori $L^2$ estimates for the wave variables}
\label{SS:STATEMENTOFAPRIORIL2ESTIMATESFORWAVEVARIABLES}
In the following proposition, we state our main a priori energy estimates
for the wave variables. We stress that for $N = \Ntop$, \eqref{E:MAINWAVEENERGYESTIMATES} follows from \eqref{E:MAINTOPORDERL2ESTIMATESMASTERCONTROLQUANTITY} by virtue of \eqref{E:MASTERCONTROLLINGQUANTITY}. The estimates for $N \le \Ntop-1$ are in fact nonlinearly coupled to estimates of $\mastercontroltop(\ubar,u)$, and its proof is also located in Sect.\,\ref{SSS:PROOFOFTOPORDERAPRIORIL2ESTIMATEMASTERCONTROLANDALLORDERWAVES}.

\begin{proposition}[The main a priori estimates for $\totalcontrolwave_{[1,\Ntop]}$] \label{P:APRIORIL2ESTIMATESWAVEVARIABLES} 
	Let $\totalcontrolwave_{N}(\ubar,u)$ 
	be the $L^2$-controlling quantity for the wave variables 
	$\wavearray$, as defined in \eqref{E:VELOCITYFLUXANDSTRONGANGULARL2CONTROLLINGQUANTITY}.
	Under the data-assumptions of Sect.\,\ref{S:ASSUMPTIONSONTHEDATA},
	the parameter size-assumptions of Sect.\,\ref{SS:PARAMETERSIZEASSUMPTIONS},
	and the bootstrap assumptions of Sect.\,\ref{S:BOOTSTRAPEVERYTHINGEXCEPTENERGIES},
	there exists a constant $C > 0$
	such that the following estimates hold for 
	$(\ubar,u) \in [\leftubar,\ubarboot) \times [\moreinterestingu_1,\moreinterestingu_2]$: 
\begin{align}
\totalcontrolwave_{N}(\ubar,u)  
& \leq C \initialsmalldoublenull^2 
&& 
\mbox{if } 1 \le N \le \Ntop. 
	\label{E:MAINWAVEENERGYESTIMATES}
\end{align}
Consequently, by \eqref{E:BLOWUPRATEOFWAVEWRTNUMBEROFCOMMUTATORS} and \eqref{E:COERCIVENESSOFNULLFLUXCONTROLWAVE}, we have the following \underline{non-weighted} estimates  for $1 \le N \le \Ntop -4$:
\begin{align} \label{E:NONWEIGHTEDL2ESTIMATESFORWAVE}
	\left\| \sqrt{\upmu} \nullangrmD \tander^N \wavearray\right\|_{L^2_0\left(\ingoingcharacteristicsurfacetwoarg{\ubar}{[\moreinterestingu_1,u]}\right)}, \, \left\| \newuL \tander^N \wavearray\right\|_{L^2_0\left(\ingoingcharacteristicsurfacetwoarg{\ubar}{[\moreinterestingu_1,u]}\right)}, \, \left\| \Lunit \tander^N \wavearray\right\|_{L^2_0\left(\outgoingcharacteristicsurfacetwoarg{u}{[\leftubar,\ubar)}\right)}, \,  \left\| \sqrt{\upmu} \angrmd \tander^N \wavearray\right\|_{L^2_0\left(\outgoingcharacteristicsurfacetwoarg{u}{[\leftubar,\ubar)}\right)} \le C \initialsmalldoublenull
\end{align}\end{proposition}

\subsection{Statement of the a priori $L^2$ estimates for the acoustic geometry}
\label{SS:STATEMENTOFAPRIORIL2ESTIMATESFORACOUSTICGEOMETRY}
The following proposition provides our main a priori null-flux estimates
for the acoustic geometry. Its proof is located in Sect.\,\ref{SS:PROOFOFENERGYESTIMATESFORACOUSTICGEOMETRY}.

\begin{proposition}[The main a priori estimates for the acoustic geometry along the rough foliations] 
\label{P:APRIORIL2ESTIMATESACOUSTICGEOMETRY}
	Under the data-assumptions of Sect.\,\ref{S:ASSUMPTIONSONTHEDATA},
	the parameter size-assumptions of Sect.\,\ref{SS:PARAMETERSIZEASSUMPTIONS},
	and the bootstrap assumptions of Sect.\,\ref{S:BOOTSTRAPEVERYTHINGEXCEPTENERGIES},
	there exists a constant $C > 0$
	such that the following estimates hold for 
	$(\ubar,u) \in [\leftubar,\ubarboot) \times [\moreinterestingu_1,\moreinterestingu_2]$: 
	
	\begin{align}
			\begin{pmatrix}
				\left\| \tander^{\le N} \mytr_{\gtorus}\upchi \right\|_{L^2_{\blowuprateofacousticgeoWRTnumberofcommutations(N)}\left( \ingoingcharacteristicsurfacetwoarg{\ubar}{[\moreinterestingu_1,u]}\right)} \\
				\left\| \angLie_{\tander}^{\le N} \upchi \right\|_{L^2_{\blowuprateofacousticgeoWRTnumberofcommutations(N)}\left( \ingoingcharacteristicsurfacetwoarg{\ubar}{[\moreinterestingu_1,u]}\right)} \\
				\left\| \tander^{[1,N+1]} \Lsmall^i\right\|_{L^2_{\blowuprateofacousticgeoWRTnumberofcommutations(N)}\left( \ingoingcharacteristicsurfacetwoarg{\ubar}{[\moreinterestingu_1,u]}\right)}	\\
				\left\| \tandersmall^{[1,N+1]}  \upmu \right\|_{L^2_{\blowuprateofacousticgeoWRTnumberofcommutations(N)}\left( \ingoingcharacteristicsurfacetwoarg{\ubar}{[\moreinterestingu_1,u]}\right)}
			\end{pmatrix} & \le C \initialsmalldoublenull, &  1 \le N \le \Ntop -1.  \label{E:MAINACOUSTICBELOWTOPESTIMATES}
	\end{align}
	
	Moreover, the following estimates hold at the top order:
	\begin{align}
		\begin{split}
			\begin{pmatrix}
				\left\| \tander^{\Ntop} \mytr_{\gtorus}\upchi \right\|_{L^2_{\blowupratetoporderwave+2}\left(\characteristicdiamondtwoarg{[\leftubar,\ubar)}{[\moreinterestingu_1,u]}\right)} \\
				\left\| \angLie_{\tander}^{\Ntop} \upchi \right\|_{L^2_{\blowupratetoporderwave+2}\left(\characteristicdiamondtwoarg{[\leftubar,\ubar)}{[\moreinterestingu_1,u]}\right)} \\
			\end{pmatrix}
			\le C  \initialsmalldoublenull. \label{E:MAINACOUSTICTOPORDERESTIMATES}
		\end{split}
	\end{align}	
\end{proposition}

\subsection{Statement of the a priori $L^2$ estimates and  
data-estimates for the $L^2$-transport-controlling quantities} \label{S:SATEMENTOFAPRIORIL2ESTIMATESFORTRANSPORTVARIABLES2DEULER}

In this short subsection we provide the statement of the main  a priori $L^2$ estimates and  
data-estimates for the $L^2$-transport-controlling quantities.

\begin{proposition}[The main below-top-order a priori estimates for the transport variables and the modified fluid variables] \label{P:APRIORIL2ESTIMATESFORTRANSPORTVARIABLES}
	Let $N = \Ntop$ and let $\fluxcontrolVort_N(\ubar,u),  \dots, \bulkcontrolDivGradEnt_N(\ubar,u)$ be the $L^2$-transport-controlling quantities as defined in Def.\,\ref{D:FUNDAMENTALTRANSPORTL2CONTROLLINGQUANTITIES}. 
	Under the parameter size-assumptions of Sect.\,\ref{SS:PARAMETERSIZEASSUMPTIONS}, the data-assumptions of Sects.\,\ref{S:ASSUMPTIONSONTHEDATA}, 
	and the bootstrap assumptions of Sect.\,\ref{S:BOOTSTRAPEVERYTHINGEXCEPTENERGIES},
	there exists a constant $C > 0$
	such that the following estimates hold for 
	$(\ubar,u) \in [\leftubar,\ubarboot) \times [\moreinterestingu_1,\moreinterestingu_2]$: 
\begin{subequations} \label{E:APRIORIL2ESTIMATESFORTRANSPORTVARIABLES}
	\begin{align}
		\fluxcontrolVort_{\le N} (\ubar,u) + \bulkcontrolVort_{\le N} (\ubar,u)
		& \leq C \initialsmalldoublenull^2, \label{E:APRIORIL2ESTIMATESFORVORTICITYFLUX} \\
		\fluxcontrolGradEnt_{\le N} (\ubar,u) + \bulkcontrolGradEnt_{\le N} (\ubar,u)& \leq C \initialsmalldoublenull^2, \label{E:APRIORIL2ESTIMATESFORGRADENTFLUX} \\
		\fluxcontrolVortVort_{\le N-1} (\ubar,u) + \bulkcontrolVortVort_{\le N-1} (\ubar,u) & \leq C \initialsmalldoublenull^2, \label{E:APRIORIL2BELOWTOPESTIMATESFORVORTVORTFLUX} \\	
		\fluxcontrolDivGradEnt_{\le N-1} (\ubar,u) + \bulkcontrolDivGradEnt_{\le N-1} (\ubar,u) & \leq C \initialsmalldoublenull^2. \label{E:APRIORIL2BELOWTOPESTIMATESFORDIVGRADENTFLUX} 
	\end{align}
\end{subequations}	
In addition, let $\toricontrolVort,\toricontrolGradEnt,\toricontrolVortVort, \toricontrolDivGradEnt$ denote the $L^2$-controlling quantities on the double-null tori  defined in Def.\,\ref{D:FUNDAMENTALTRANSPORTL2CONTROLLINGQUANTITIES}. Then, the following a priori estimates hold:
\begin{subequations}
		\begin{align}
			\toricontrolVort_{\le N - 1}(\ubar,u), 
			\,
			\toricontrolGradEnt_{\le N - 1} 
			& 
			\leq C \initialsmall^2, 
			\label{E:MAINTORIL2VORTGRADENTBELOWTOPORDERBLOWUP} \\
			\toricontrolVortVort_{\le N - 2}(\ubar,u),
			\,
			\toricontrolDivGradEnt_{\le N - 2}(\ubar,u)
			& \leq C \initialsmall^2.
			\label{E:MAINTORIL2MODIFIEDFLUIDVARIABLESBELOWTOPORDERBLOWUP}
	\end{align}
\end{subequations}

\end{proposition}

The proof of Prop.\ref{P:APRIORIL2ESTIMATESFORTRANSPORTVARIABLES} will be given in Sect.\,\ref{SS:PROOFOFBELOWTOPORDERTRANSPORTENERGYESTIMATES} and\,\ref{SS:PROOFOFBELOWTOPORDERROUGHTORIENERGYESTIMATES}.

\begin{proposition}[The main top-order a priori $L^2$ estimates for the transport and modified fluid variables]
	\label{P:TOPORDERMODIFIEDAPRIORIL2ESTIMATESANDDOUBLENULLTORIESTIMATES} 
	Let $\fluxcontrolVortVort_N(\ubar,u),  \fluxcontrolDivGradEnt_N(\ubar,u), \bulkcontrolDivGradEnt_N(\ubar,u), \bulkcontrolDivGradEnt_N(\ubar,u)$ denote the $L^2$-controlling for the modified fluid variables from Def.\,\ref{D:FUNDAMENTALTRANSPORTL2CONTROLLINGQUANTITIES}. Under the data-assumptions of Sect.\,\ref{S:ASSUMPTIONSONTHEDATA},
	the parameter size-assumptions of Sect.\,\ref{SS:PARAMETERSIZEASSUMPTIONS},
	and the bootstrap assumptions of Sect.\,\ref{S:BOOTSTRAPEVERYTHINGEXCEPTENERGIES},
	there exists a constant $C > 0$
	such that the following estimates hold for 
	$(\ubar,u) \in [\leftubar,\ubarboot) \times [\moreinterestingu_1,\moreinterestingu_2]$: 
	
	\begin{subequations}
		\begin{align}
			\fluxcontrolVortVort_{\Ntop}(\ubar,u) + \bulkcontrolVortVort_{\Ntop}(\ubar,u) & \le C \initialsmall^2, \label{E:TOPORDERVORTVORTL2APRIORIESTIMATES} \\
			\fluxcontrolDivGradEnt_{\Ntop}(\ubar,u) + \bulkcontrolDivGradEnt_{\Ntop}(\ubar,u) & \le C \initialsmall^2. \label{E:TOPORDERDIVGRADENTL2APRIORIESTIMATES} 
		\end{align}
	\end{subequations}
	
	In addition, let $\toricontrolVort,\toricontrolGradEnt,\toricontrolVortVort, \toricontrolDivGradEnt$ denote the $L^2$-controlling quantities on the double-null tori  defined in Def.\,\ref{D:FUNDAMENTALTRANSPORTL2CONTROLLINGQUANTITIES}. Then, the following a priori estimates hold:

	\begin{subequations}
		\begin{align}
			\toricontrolVort_{\Ntop}(\ubar,u), 
			\,
			\toricontrolGradEnt_{\Ntop}(\ubar,u)
			& 
			\leq C \initialsmall^2, 
			\label{E:MAINTORIL2VORTGRADENTTOPORDER} \\
			\toricontrolVortVort_{\Ntop - 1}(\ubar,u),
			\,
			\toricontrolDivGradEnt_{\Ntop - 1}(\ubar,u)
			& \leq C \initialsmall^2.
			\label{E:MAINTORIL2MODIFIEDFLUIDVARIABLESTOPORDER}
		\end{align}
	\end{subequations}

\end{proposition}

The proof of Prop.\,\ref{P:TOPORDERMODIFIEDAPRIORIL2ESTIMATESANDDOUBLENULLTORIESTIMATES} will be given in Sect.\,\ref{S:TOPORDERELIPTICHYPERBOLICL2ESTIMATESFORSPECIFICVORTICITYANDENTROPYGRADIENT}.

\begin{remark}[The need for Prop.\,\ref{P:TOPORDERMODIFIEDAPRIORIL2ESTIMATESANDDOUBLENULLTORIESTIMATES}] We highly that the statement and conclusions of Prop.\,\ref{P:TOPORDERMODIFIEDAPRIORIL2ESTIMATESANDDOUBLENULLTORIESTIMATES} could have easily been included in Prop.\,\ref{P:APRIORIL2ESTIMATESFORTRANSPORTVARIABLES}. However, we chose to give those estimates their own proposition as their proofs require much more delicate arguments than Prop.\,\ref{P:APRIORIL2ESTIMATESFORTRANSPORTVARIABLES} that involve the elliptic-hyperbolic integral identities of Sect.\,\ref{S:ELLIPTICHYPERBOLICIDENTITIES}.
\end{remark}

\subsection{Data-estimates for the $L^2$-controlling quantities} 
\label{SS:DATAESTIMATESFORL2CONTROLLINGQUANTITIES}
In this section, we derive estimates for the data of the
$L^2$-controlling quantities
$\totalcontrolwave_{[1,\Ntop]}, \,  \totalfluxcontrolprecisefullymodifiedchi$, etc. 
We will use these data-estimates in our proofs of
Props.\,\ref{P:APRIORIL2ESTIMATESWAVEVARIABLES}--\ref{P:APRIORIL2ESTIMATESACOUSTICGEOMETRY}.

\begin{lemma}[The $L^2$-controlling quantities are initially small] 
\label{L:ALLL2CONTROLLINGQUANTITIESINITIALLYSMALL}
The following estimates hold for $(\ubar,u) \in [\leftubar,\ubarboot) \times [\moreinterestingu_1,\moreinterestingu_2]$:

	\begin{align}
		\totalcontrolwave_{[1,\Ntop]}(\leftubar,u) & \leq C\initialsmalldoublenull^2, & &  \totalcontrolwave_{[1,\Ntop]}(\ubar,\moreinterestingu)  \leq C \initialsmalldoublenull^2, \label{E:WAVEL2CONTROLLINGINITIALLYSMALL} 
	\end{align}
	\begin{subequations}
		\begin{align} 
			\fluxcontrolVort_{\le \Ntop} (\leftubar,u) & \le C \initialsmalldoublenull^2, & & & \fluxcontrolVort_{\le \Ntop} (\ubar,\moreinterestingu_1) & \le C \initialsmalldoublenull^2, \label{E:FLUXFORVORTICITYINITIALLSMALL} \\
			\fluxcontrolGradEnt_{\le \Ntop} (\leftubar,u) & \le C \initialsmalldoublenull^2, & & & \fluxcontrolGradEnt_{\le \Ntop} (\ubar,\moreinterestingu_1) & \le C \initialsmalldoublenull^2, \label{E:FLUXFORGRADENTINITIALLSMALL} \\
			\fluxcontrolVortVort_{\le \Ntop} (\leftubar,u) & \le C \initialsmalldoublenull^2, & & & \fluxcontrolVortVort_{\le \Ntop} (\ubar,\moreinterestingu_1) & \le C \initialsmalldoublenull^2, \label{E:FLUXFORVORTVORTNITIALLSMALL} \\
			\fluxcontrolDivGradEnt_{\le \Ntop} (\leftubar,u) & \le C \initialsmalldoublenull^2, & & & \fluxcontrolDivGradEnt_{\le \Ntop} (\ubar,\moreinterestingu_1) & \le C \initialsmalldoublenull^2, \label{E:FLUXFORDIVGRADENTTNITIALLSMALL} 
		\end{align}
	\end{subequations}
	\begin{subequations}
		\begin{align}
			\totalfluxcontrolprecisefullymodifiedchi(\leftubar,u) & \leq C \initialsmalldoublenull^2, \label{E:PRECISEFULLYMODQUANTL2CONTROLLINGINITIALLYSMALL} \\
			\totalfluxcontrolimprecisefullymodifiedchi(\leftubar,u) & \leq C\initialsmalldoublenull^2 \label{E:IMPRECISEFULLYMODQUANTCONTROLLINGINITIALLYSMALL} \\
			\totalfluxcontrolpartialmodifiedchi(\leftubar,u) & \leq C\initialsmalldoublenull^2. \label{E:PARTIALMODQUANT2CONTROLLINGINITIALLYSMALL}
		\end{align}
	\end{subequations}
\end{lemma}

\begin{proof}
	\eqref{E:WAVEL2CONTROLLINGINITIALLYSMALL} follows as a straightforward consequence of
	the data-assumptions  
\eqref{E:TRANSVERSALDERIVATIVEOFTANGENTIALL2NORMSOFWAVEVARIABLESSMALLALONGINITIALINGOINGNULLSURFACE}--\eqref{E:SMOOTHANGULARDERIVATIVEOFTANGENTIALL2NORMSOFWAVEVARIABLESSMALLALONGINITIALOUTGOINGCHARACTERISTICHYPERSURFACE},  
	definitions 
	\eqref{E:WAVEINGOINGNULLFLUX}--\eqref{E:WAVEOUTGOINGNULLFLUX},
	\eqref{E:SPACETIMECOERCIVEINTEGRALS},
	\eqref{E:UNIFIEDCOMMUTATORVELOCITYNULLFLUXL2CONTROLLINGQUANTITY}--\eqref{E:UNIFIEDCOMMUTATORSPACETIMESTRONGANGULARVELOCITYCONTROL} and \eqref{E:VELOCITYFLUXANDSTRONGANGULARL2CONTROLLINGQUANTITY}, 
	and the identities 
	\eqref{E:QBREVETL}--\eqref{E:QBREVETNEWUL}, and estimate \eqref{E:COERCIVENESSOFNULLFLUXCONTROLWAVE}.
	
	Estimates \eqref{E:FLUXFORVORTICITYINITIALLSMALL}--\eqref{E:FLUXFORDIVGRADENTTNITIALLSMALL} follow immediately from Lemma\,\ref{L:COERCIVENESSOFL2CONTROLLINGQUANITIESFORVORTICITYANDENTROPY} and estimate $\frac{1}{\Lunit \ubar} \approx 1$.
	
	Estimates \eqref{E:PRECISEFULLYMODQUANTL2CONTROLLINGINITIALLYSMALL}--\eqref{E:PARTIALMODQUANT2CONTROLLINGINITIALLYSMALL} follow immediately from \eqref{E:L2DATASSUMPTIONFORFULLYMODIFIEDCHITOPORDER}--\eqref{E:L2DATASSUMPTIONFORPARTIALLYMODIFIEDCHITOPORDER}.

\end{proof}

\subsection{Bootstrap assumptions for the $\totalcontrolwave_N$}
\label{SS:BOOTSTRAPASSUMPTIONSFORTHEWAVEENERGIES}
For future use in a continuation criteria argument,
we find it convenient to make bootstrap assumptions for the $L^2$-controlling quantities
for the wave variables.
Specifically, with
$\weakspacetimetransvelocitycontrol_N(\ubar,u)$
(see definition  \eqref{E:VELOCITYFLUXANDSTRONGANGULARL2CONTROLLINGQUANTITY}) 
denoting the $L^2$-controlling quantity for the wave variables $\wavearray$,
we assume that the following bootstrap assumptions hold for 
$(\ubar,u) \in [\leftubar,\ubarboot) \times [\moreinterestingu_1,\moreinterestingu_2]$, 
where $\fundbootsmall$ is the bootstrap parameter from Sect.\,\ref{SSS:FUNDAMENTALQUANTITATIVE}:
 \begin{subequations}
	\begin{align}
		\totalcontrolwave_{N}(\ubar,u)
		&  
		\leq \fundbootsmall 
		&& \mbox{if } 1 \le N \le \Ntop. 
		\label{E:MAINWAVEENERGYBOOTSTRAP}
	\end{align}
\end{subequations}\begin{remark}[The wave variable energy estimates improve the bootstrap assumptions]
	\label{R:WAVEENERGYESTIMATESYIELDSTRICTIMPROVEMENTSOVERWAVEENERGYBOOTSTRAPASSUMPTIONS}
	Note that when $\initialsmalldoublenull$ is sufficiently small,
	the estimates of Props.\,\ref{P:TOPORDERAPRIORIL2ESTIMATESMASTERCONTROL}--\ref{P:APRIORIL2ESTIMATESWAVEVARIABLES} 
	yield strict improvements over the bootstrap assumptions
	\eqref{E:MAINWAVEENERGYBOOTSTRAP}.
\end{remark}

\section{Preliminary below-top-order $L^2$ estimates for the acoustic geometry and a derivative-losing estimate} 
\label{S:PRELIMINARYL2ESTIMATESFORBELOWTOPORDERDERIVATIVESOFACOUSTICGEOMETRYANDDERIVATIVELOSING}
We continue to work under the assumptions of Sect.\,\ref{SS:SILENTFACTS}.
In this short section, 
we derive preliminary $L^2$ estimates for the below-top-order derivatives of the
eikonal function quantities $\upmu$, $\Lunit^i$, $\upchi$, and $\mytr_{\gtorus}\upchi$. 
We state the bounds in terms of the wave $L^2$-controlling quantities $\totalfluxcontrolvelocity_{[1,N]}(\ubar,u)$
from Def.\,\ref{D:MAINCOERCIVE} 
and the initial data-size-parameter $\initialsmalldoublenull$.
We provide the main estimates in Lemma\,\ref{L:PRELIMINARYBELOWTOPORDERL2ESTIMATESFOREIKONALFUINCTIONQUANTITIES}.
The estimates are rather straightforward 
consequences of the transport inequalities provided by
Props.\,\ref{P:POINTWISETRANSPORTINEQUALITIESFOREIKFUNCTIONQUANTITIES}
and\,\ref{P:ENERGYIDENTITIESFOREBLOWTOPORDERACOUSTICQUANTITIES}. 
In Cor\,\ref{C:NONSINGULARL2ESTIMATESFORWAVEVARIABLESTHATLOSEONEDERIVATIVE},
we derive derivative-losing $L^2$ estimates for $\wavearray$ that do not involve
any loss of the weights $\weight$. The absence of such a loss
is important for our proof that the $\weight$-weighted energies for the wave variables become less and less singular with respect to 
powers of $\weight$
as we descend below the top-order
(see Prop.\,\ref{P:APRIORIL2ESTIMATESWAVEVARIABLES}).

Most of the $L^2$ estimates we derive in this section lose one derivative.
In Prop.\,\ref{P:TOPORDERIMPRECISEL2ESTIMATEMUCHI}, we will derive complementary 
$L^2$ estimates for the top-order derivatives of $\upchi$ and $\mytr_{\gtorus}\upchi$.
Those estimates are much harder to prove because we cannot afford to lose 
any derivatives, which forces us to rely on the modified quantities from Sect.\,\ref{S:CONSTRUCTIONOFMODIFIEDQUANTITIES}
and elliptic estimates for the top-order derivatives of $\upchi$.

\subsection{Preliminary below-top-order $L^2$ estimates for the eikonal function quantities}
\label{SS:PRELIMINARYL2ESTIMATESFORBELOWTOPORDERDERIVATIVESOFACOUSTICGEOMETRY}

\begin{lemma}[Preliminary below-top-order $L^2$ estimates for the eikonal function quantities] \label{L:PRELIMINARYBELOWTOPORDERL2ESTIMATESFOREIKONALFUINCTIONQUANTITIES}
	Let $1 \le N \le \Ntop -1$ and consider any 
	\begin{subequations}
		\begin{align}
			\varphi_N  & \in \left\{ \tander^{[1,N+1]} \Lsmall^i, \tander^{\le N} \mytr_{\gtorus} \upchi, |\angLie_\tander^{\le N} \upchi|_{\gtorus}\}, \label{E:ARBITRARYCHOICEOFPRELIMINARYBELOWTOPORDEREIKONALFUINCTIONQUANTITIESNOMUX} \right\} \\
			 \widetilde{\varphi}_N & \in \left\{\tandersmall^{[1,N+1]} \upmu, \comdersmall^{[1,N+1];1} \Lsmall^i, \comder^{\le N;1} \mytr_{\gtorus}\upchi, |\angLie_{\comder}^{\le N;1} \upchi|_{\gtorus}\right\}. \label{E:ARBITRARYCHOICEOFPRELIMINARYBELOWTOPORDEREIKONALFUINCTIONQUANTITIESWITHMUX}
		\end{align}
	\end{subequations}	
	Let $\blowuprateofacousticgeoWRTnumberofcommutations(N)$ denote the blowup rate of the null geometry variables defined as in \eqref{E:BLOWUPRATEOFACOUSTICBELOWTOPWRTNUMBEROFCOMMUTATORS}. Then for sufficiently small $\varsigma_0 \in (0,1]$, the following estimates hold for $(\ubar,u) \in [\leftubar,\ubarboot) \times [\moreinterestingu_1,\moreinterestingu_2]$, where the implicit constants are \textbf{independent} of $\varsigma_0$: 
	\begin{align}
		\left\| \varphi_N  \right\|^2_{L^2_{\blowuprateofacousticgeoWRTnumberofcommutations(N)}\left( \ingoingcharacteristicsurfacetwoarg{\ubar}{[\moreinterestingu_1,u]}\right)}		  
		&    \lesssim \initialsmalldoublenull^2 + (1 + \varsigma_0^{-1}) \int_{u' = \moreinterestingu_1}^u \totalfluxcontrolvelocity_{[1,N+1]}(\ubar,u') \, \mathrm{d} u' + (1 + \varsigma_0^{-1}) \strongangularcontrolvelocity_{[1,N+1]}(\ubar,u) \\
		& \ \ + \int_{u' = \moreinterestingu_1}^u \left\{\fluxcontrolVort_{\le N+1}(\ubar,u')+ \fluxcontrolGradEnt_{\le N+1}(\ubar,u')\right\} \, \mathrm{d} u',   \label{E:PRELIMINARYBELOWTOPORDERL2ESTIMATESFOREIKONALFUINCTIONQUANTITIES} \\
		\begin{split} \label{E:PRELIMINARYBELOWTOPORDERL2ESTIMATESFOREIKONALFUINCTIONQUANTITIESWITHZ}
			\left\| \widetilde{\varphi}_N  \right\|^2_{L^2_{\blowuprateofacousticgeoWRTnumberofcommutations(N)}\left( \ingoingcharacteristicsurfacetwoarg{\ubar}{[\moreinterestingu_1,u]}\right)}		  
			&    \lesssim \initialsmalldoublenull^2 +  (1 + \varsigma_0^{-1}) \int_{\ubar' = \leftubar}^{\ubar} \totalfluxcontrolvelocity_{[1,N+1]}(\ubar',u) \, \mathrm{d} \ubar'   \\
			& \ \ +  (1 + \varsigma_0^{-1}) \int_{u' = \moreinterestingu_1}^u \totalfluxcontrolvelocity_{[1,N+1]}(\ubar,u') \, \mathrm{d} u' + (1 + \varsigma_0^{-1}) \strongangularcontrolvelocity_{[1,N+1]}(\ubar,u) \\
			& \ \ + \int_{u' = \moreinterestingu_1}^u \left\{\fluxcontrolVort_{\le N+1}(\ubar,u')+ \fluxcontrolGradEnt_{\le N+1}(\ubar,u')\right\} \, \mathrm{d} u'. 
		\end{split}
	\end{align} 
	
Moreover, whenever $\blowuprateofacousticgeoWRTnumberofcommutations(N-1) > 0$, we have: 
	\begin{align}
		\left\| \tander^N \Lsmall^i  \right\|^2_{L^2_{\blowuprateofacousticgeoWRTnumberofcommutations(N-1) -1 }\left(\characteristicdiamondtwoarg{[\leftubar,\ubar)}{[\moreinterestingu_1,u]}\right)}		  
		&    \lesssim \initialsmalldoublenull^2 + (1 + \varsigma_0^{-1}) \int_{u' = \moreinterestingu_1}^u \totalfluxcontrolvelocity_{[1,N+1]}(\ubar,u') \, \mathrm{d} u' + (1 + \varsigma_0^{-1}) \strongangularcontrolvelocity_{[1,N+1]}(\ubar,u) \\
		& \ \ + \int_{u' = \moreinterestingu_1}^u \left\{\fluxcontrolVort_{\le N+1}(\ubar,u')+ \fluxcontrolGradEnt_{\le N+1}(\ubar,u')\right\} \, \mathrm{d} u',   \label{E:PRELIMINARYJUSTBELOWTOPORDERL2SPACETIMEESTIMATESFORLSMALL} \\
		\begin{split} \label{E:PRELIMINARYJUSTBELOWTOPORDERL2SPACETIMEESTIMATESFORMU}
			\left\| \tandersmall^N \upmu  \right\|^2_{L^2_{\blowuprateofacousticgeoWRTnumberofcommutations(N-1) -1}\left( \characteristicdiamondtwoarg{[\leftubar,\ubar)}{[\moreinterestingu_1,u]}\right)}		  
			&    \lesssim \initialsmalldoublenull^2 +  (1 + \varsigma_0^{-1}) \int_{\ubar' = \leftubar}^{\ubar} \totalfluxcontrolvelocity_{[1,N+1]}(\ubar',u) \, \mathrm{d} \ubar'   \\
			& \ \ +  (1 + \varsigma_0^{-1}) \int_{u' = \moreinterestingu_1}^u \totalfluxcontrolvelocity_{[1,N+1]}(\ubar,u') \, \mathrm{d} u' + (1 + \varsigma_0^{-1}) \strongangularcontrolvelocity_{[1,N+1]}(\ubar,u) \\
			& \ \ + \int_{u' = \moreinterestingu_1}^u \left\{\fluxcontrolVort_{\le N+1}(\ubar,u')+ \fluxcontrolGradEnt_{\le N+1}(\ubar,u')\right\} \, \mathrm{d} u'. 
		\end{split}
	\end{align}

\end{lemma}

\begin{proof}  \hfill

\noindent \underline{\textbf{Preliminary pointwise estimates:}}

We first prove the following bounds:
	\begin{subequations}
		\begin{align}
			\left| \Lunit \left( |\angLie_{\tander}^{N-1} \upchi|_{\gtorus}\right)\right| & \lesssim |\angLie_{L}\angLie_{\tander}^{N-1} \upchi|_{\gtorus} + \fundbootsmall |\angLie_{\tander}^{N-1} \upchi|_{\gtorus}, \label{E:PRELIMINARYPOINTWISEFORNORMOFLIETANDERCHI} \\
			\left| \Lunit \left( |\angLie_{\comder}^{N-1;1} \upchi|_{\gtorus}\right)\right| & \lesssim |\angLie_{L}\angLie_{\comder}^{N-1;1} \upchi|_{\gtorus} + \fundbootsmall |\angLie_{\tander}^{N-1} \upchi|_{\gtorus},  \label{E:PRELIMINARYPOINTWISEFORNORMOFLIECOMDERCHI}
		\end{align}
	\end{subequations}
which will be used to derive estimates for the transport equation satisfied by $|\angLie_{\tander}^{N-1} \upchi|_{\gtorus}, \, |\angLie_{\comder}^{N-1;1} \upchi|_{\gtorus}$ along $\Lunit$. 
To see \eqref{E:PRELIMINARYPOINTWISEFORNORMOFLIETANDERCHI}, we use the Leibniz rule for $\ell_{t,u}$-projected Lie derivatives,
the identity $\angLie_{\Lunit} \gtorus^{-1} = - 2 \upchi^{\# \#}$ 
(which is a straightforward consequence of \eqref{E:NULLSECONDFUND}),
and the estimates 
$|\mytr_{\gtorus}\upchi|, \, |\upchi|_{\gtorus} \lesssim \fundbootsmall$, 
which themselves follow from \eqref{E:CHIEXPRESSSIONINTERMSOFDERIVATIVESOFLUNITI}, 
Prop.\,\ref{P:SCHEMATICSTRUCTUREOFVARIOUSTENSORSINTERMSOFCONTROLVARS},
the bootstrap assumptions,
and Prop.\,\ref{P:IMPROVEMENTOFAUXILIARYBOOTSTRAP}. The bound \eqref{E:PRELIMINARYPOINTWISEFORNORMOFLIECOMDERCHI} follows from similar arguments which we omit.

\medskip

\noindent \underline{\textbf{Transport equation estimates for pure tangential commutations:} \emph{ingoing characteristic integrals.}} 

	We fix an integer $1 \le N \le \Ntop -1$, a point $u \in [\moreinterestingu_1,\moreinterestingu_2]$, and define the following function of $\ubar \in [\leftubar,\ubarboot)$:
	\begin{align}
		p_N(\ubar) & \eqdef \sum_{N' = 0}^{N}   \int_{\ingoingcharacteristicsurfacetwoarg{\ubar}{[\moreinterestingu_1,u]}} \left\{ \left(\tander^{N'} \mytr_{\gtorus} \upchi\right)^2 + \left| \angLie_{\tander}^{N'} \upchi \right|_{\gtorus}^2 + \sum_{i=1}^3  \left( \tander^{N'+1} \Lsmall^i\right)^2  \right\} \weight^{\blowuprateofacousticgeoWRTnumberofcommutations(N')}  \, \volingoingnullhypersurface. \label{E:SQUAREL2NORMSOFBELOWTOPORDERACOUSTICVARIABLESWITHOUTMUXDIFFERENTIATION} 
	\end{align} 
Suppose that $\Ntop-3 \le N' \le \Ntop -1$, which implies $\blowuprateofacousticgeoWRTnumberofcommutations(N') > 0$. Applying the identity \eqref{E:SPACETIMEENERGYNULLFLUXIDENTITYFORBELOWTOPORDERACOUSTIC} for $F_{N'} \in \{ \tander^{N'+1} \Lsmall^i, \, \tander^{N'} \mytr_{\gtorus} \upchi, \, |\angLie_{\tander}^{N'} \upchi |_{\gtorus}\}$ with $q =\blowuprateofacousticgeoWRTnumberofcommutations(N')$, inserting the estimates \eqref{E:LUNITTANGENTIALDERIVATIVESOFLUNITIPOINTWISE}--\eqref{E:ANGLIELTANGENTIALCHIPOINTWISE} for $\Lunit \varphi$ along with Cor\,\ref{C:IMPROVEAUX}, using $\frac{1}{\Lunit \ubar} = \ReciprocalLunitAppliedtoTimeFunction \approx 1$, and applying Young's inequality, we find:
\begin{align} \label{E:SQUAREL2NORMSOFBELOWTOPORDERACOUSTICVARIABLESWITHOUTMUXDIFFERENTIATIONSTEP1}
	\begin{split}
		& \int_{\ingoingcharacteristicsurfacetwoarg{\ubar}{[\moreinterestingu_1,u]}} \left| F_{N'} \right|^2 \weight^{\blowuprateofacousticgeoWRTnumberofcommutations(N')} \volingoingnullhypersurface + \blowuprateofacousticgeoWRTnumberofcommutations(N') \newspacetimecoercive{\blowuprateofacousticgeoWRTnumberofcommutations(N')-1}[F_{N'}](\ubar,u) \\
		& \ \ \lesssim p_{N'} (\leftubar) + \fundbootsmall \int_{\ubar' = \leftubar}^{\ubar} p_{N'} (\ubar',u)\, \mathrm{d} \ubar + \int_{\characteristicdiamondtwoarg{[\leftubar,\ubar)}{[\moreinterestingu_1,u]}} \left| F_{N'}\right| \cdot \left| \tander^{[1,N'+2]} \velocityarray \right| \weight^{\blowuprateofacousticgeoWRTnumberofcommutations(N')} \ReciprocalLunitAppliedtoTimeFunction \, \voldiamond \\
		& \ \ \ \ + \int_{\characteristicdiamondtwoarg{[\leftubar,\ubar)}{[\moreinterestingu_1,u]}} \left| F_{N'}\right| \cdot \left| \tander^{\le N'+1}(\vortrenormalized,\GradEnt)\right| \weight^{\blowuprateofacousticgeoWRTnumberofcommutations(N')} \ReciprocalLunitAppliedtoTimeFunction \, \voldiamond.
	\end{split}
\end{align}
We now consider the case where $\tander^{N'+2}\velocityarray = \Lunit \tander^{N'+1} \velocityarray$ in the second-to-last integral in \eqref{E:SQUAREL2NORMSOFBELOWTOPORDERACOUSTICVARIABLESWITHOUTMUXDIFFERENTIATIONSTEP1}. Multiplying and dividing the integrand by $(\varsigma_0\weight)^{1/2}$ and applying Young's inequality, we find that:
	\begin{align}
		\begin{split} \label{E:SQUAREL2NORMSOFBELOWTOPORDERACOUSTICVARIABLESWITHOUTMUXDIFFERENTIATIONSTEP2}
			 \int_{\characteristicdiamondtwoarg{[\leftubar,\ubar)}{[\moreinterestingu_1,u]}} \left| F_{N'}\right| \cdot \left| \Lunit \tander^{N'+1} \velocityarray \right| \weight^{\blowuprateofacousticgeoWRTnumberofcommutations(N')} \ReciprocalLunitAppliedtoTimeFunction \, \voldiamond
			 & \lesssim \varsigma_0 	 \int_{\characteristicdiamondtwoarg{[\leftubar,\ubar)}{[\moreinterestingu_1,u]}} \left| F_{N'}\right|^2  \weight^{\blowuprateofacousticgeoWRTnumberofcommutations(N') - 1} \voldiamond \\
			 & \ \ + \varsigma_0^{-1}  \int_{\characteristicdiamondtwoarg{[\leftubar,\ubar)}{[\moreinterestingu_1,u]}} \left| \Lunit \tander^{N'+1} \velocityarray \right|^2 \weight^{\blowuprateofacousticgeoWRTnumberofcommutations(N')+1} \ReciprocalLunitAppliedtoTimeFunction \, \voldiamond.
		\end{split}
	\end{align}
If $\varsigma_0$ is sufficiently small, then \eqref{E:COERCIVITYNEWSPACETIMETERM} implies that the first integral on RHS\,\eqref{E:SQUAREL2NORMSOFBELOWTOPORDERACOUSTICVARIABLESWITHOUTMUXDIFFERENTIATIONSTEP2} may be absorbed by the spacetime integral $\newspacetimecoercive{\blowuprateofacousticgeoWRTnumberofcommutations(N')-1}[\varphi](\ubar,u)$ present on the LHS\,\eqref{E:SQUAREL2NORMSOFBELOWTOPORDERACOUSTICVARIABLESWITHOUTMUXDIFFERENTIATIONSTEP1}. Moreover, the second integral on RHS\,\eqref{E:SQUAREL2NORMSOFBELOWTOPORDERACOUSTICVARIABLESWITHOUTMUXDIFFERENTIATIONSTEP2} is $\lesssim \varsigma_0^{-1}  \int_{u' = \moreinterestingu_1}^u \totalfluxcontrolvelocity_{N'+1}(\ubar,u') \, \mathrm{d} u'$ by \eqref{E:COERCIVENESSOFNULLFLUXCONTROLWAVE} and the identity $\blowuprateofacousticgeoWRTnumberofcommutations(N') +1 = \blowuprateofwaveWRTnumberofcommutations(N'+1)$, which holds whenever $\blowuprateofacousticgeoWRTnumberofcommutations(N') > 0$ and $N' \le \Ntop -1$ (see \eqref{E:BLOWUPRATEOFWAVEWRTNUMBEROFCOMMUTATORS}--\eqref{E:BLOWUPRATEOFACOUSTICBELOWTOPWRTNUMBEROFCOMMUTATORS}). A similar argument as in \eqref{E:SQUAREL2NORMSOFBELOWTOPORDERACOUSTICVARIABLESWITHOUTMUXDIFFERENTIATIONSTEP2} can be used in the case where $\tander^{N'+2}\velocityarray = \Yvf{A} \tander^{N'+1} \velocityarray$ for some $A = 2,3$, but one instead bounds the resulting integral by $\varsigma_0^{-1}\strongangularcontrolvelocity_{N'}(\ubar,u)$ using \eqref{E:SMOOTHTORUSNORMCOMPARBLETOTANGENTIALCONTRACTIONS}, \eqref{E:UNIFIEDCOMMUTATORSPACETIMESTRONGANGULARVELOCITYCONTROL}, and \eqref{E:COERCIVITYOLDSPACETIMETERM}.

To bound the last error integral on RHS\,\eqref{E:SQUAREL2NORMSOFBELOWTOPORDERACOUSTICVARIABLESWITHOUTMUXDIFFERENTIATIONSTEP1}, we use the identity $\blowuprateofacousticgeoWRTnumberofcommutations(N') = \blowuprateoftransportWRTnumberofcommutations(N'+1)$ and Young's inequality:
	\begin{align} 
		\begin{split} \label{E:SQUAREL2NORMSOFBELOWTOPORDERACOUSTICVARIABLESWITHOUTMUXDIFFERENTIATIONSTEP2B}
			 \int_{\characteristicdiamondtwoarg{[\leftubar,\ubar)}{[\moreinterestingu_1,u]}} \left| F_{N'}\right| \cdot \left| \tander^{\le N' + 1}(\vortrenormalized,\GradEnt) \right| \weight^{\blowuprateofacousticgeoWRTnumberofcommutations(N')} \ReciprocalLunitAppliedtoTimeFunction \, \voldiamond 
			 & \lesssim \int_{\characteristicdiamondtwoarg{[\leftubar,\ubar)}{[\moreinterestingu_1,u]}} \left| F_{N'}\right|^2  \weight^{\blowuprateofacousticgeoWRTnumberofcommutations(N') } \voldiamond \\
			 & \ \  +   \int_{\characteristicdiamondtwoarg{[\leftubar,\ubar)}{[\moreinterestingu_1,u]}} \left| \tander^{\le N'+1}(\vortrenormalized,\GradEnt)\right|^2 \weight^{\blowuprateoftransportWRTnumberofcommutations(N'+1)} \ReciprocalLunitAppliedtoTimeFunction \, \voldiamond.
		\end{split}
	\end{align}
The last integral on \eqref{E:SQUAREL2NORMSOFBELOWTOPORDERACOUSTICVARIABLESWITHOUTMUXDIFFERENTIATIONSTEP2B} is easily seen to be $\lesssim \int_{u' = \moreinterestingu_1}^u \left\{\fluxcontrolVort_{\le N+1}(\ubar,u')+ \fluxcontrolGradEnt_{\le N+1}(\ubar,u')\right\} \, \mathrm{d} u'$.

Suppose now that $N' = \Ntop - 4$ so that  $0 = \blowuprateofacousticgeoWRTnumberofcommutations(N') =\blowuprateoftransportWRTnumberofcommutations(N'+1)$ but $0 < \blowuprateofwaveWRTnumberofcommutations(N'+1) = 0.13 < 1$, where the second inequality is a consequence of \eqref{E:BLOWUPRATEOFWAVEWRTNUMBEROFCOMMUTATORS}--\eqref{E:BLOWUPRATEOFACOUSTICBELOWTOPWRTNUMBEROFCOMMUTATORS}. Applying \eqref{E:SPACETIMEENERGYNULLFLUXIDENTITYFORBELOWTOPORDERACOUSTIC} as above but now with $q = 0$, we find that:
\begin{align} \label{E:SQUAREL2NORMSOFBELOWTOPORDERACOUSTICVARIABLESWITHOUTMUXDIFFERENTIATIONSTEP3}
		\begin{split}
			\int_{\ingoingcharacteristicsurfacetwoarg{\ubar}{[\moreinterestingu_1,u]}} \left| F_{N'} \right|^2  \volingoingnullhypersurface  
			& \lesssim  p_{N'} (\leftubar) + \fundbootsmall \int_{\ubar' = \leftubar}^{\ubar} p_{N'} (\ubar',u)\, \mathrm{d} \ubar + \int_{\characteristicdiamondtwoarg{[\leftubar,\ubar)}{[\moreinterestingu_1,u]}} \left| F_{N'}\right| \cdot \left| \tander^{[1,N'+2]} \velocityarray \right| \ReciprocalLunitAppliedtoTimeFunction \, \voldiamond. \\
		& \ \ \ \ \int_{\characteristicdiamondtwoarg{[\leftubar,\ubar)}{[\moreinterestingu_1,u]}} \left| F_{N'}\right|^2  \voldiamond   +   \int_{\characteristicdiamondtwoarg{[\leftubar,\ubar)}{[\moreinterestingu_1,u]}} \left| \tander^{\le N'+1}(\vortrenormalized,\GradEnt)\right|^2 \ReciprocalLunitAppliedtoTimeFunction \, \voldiamond.
		\end{split}
	\end{align}
Upon multiplying and dividing by $\weight^{\blowuprateofwaveWRTnumberofcommutations(N'+1)/2}$ in the second integral in RHS\,\eqref{E:SQUAREL2NORMSOFBELOWTOPORDERACOUSTICVARIABLESWITHOUTMUXDIFFERENTIATIONSTEP3}, Young's inequality and \eqref{E:BLOWUPWEIGHT} imply:
   \begin{align} \label{E:SQUAREL2NORMSOFBELOWTOPORDERACOUSTICVARIABLESWITHOUTMUXDIFFERENTIATIONSTEP4}
	\begin{split}
		  \int_{\characteristicdiamondtwoarg{[\leftubar,\ubar)}{[\moreinterestingu_1,u]}} \left| F_{N'}\right| \cdot \left| \tander^{[1,N'+1]} \velocityarray \right| \ReciprocalLunitAppliedtoTimeFunction \, \voldiamond &  \lesssim \int_{\ubar' = \leftubar}^{\ubar} |\ubar'|^{-\blowuprateofwaveWRTnumberofcommutations(N'+1)} \int_{\ingoingcharacteristicsurfacetwoarg{\ubar'}{[\moreinterestingu_1,u]}}    \left| F_{N'} \right|^2  \volingoingnullhypersurface \mathrm{d} \ubar' \\
		& \ \ +  \int_{\characteristicdiamondtwoarg{[\leftubar,\ubar)}{[\moreinterestingu_1,u]}} \left|\tander^{[1,N'+2]}\velocityarray \right|^2 \weight^{\blowuprateofwaveWRTnumberofcommutations(N'+1)}  \voldiamond,
	\end{split}
\end{align}
where we note that the first integral on RHS\,\eqref{E:SQUAREL2NORMSOFBELOWTOPORDERACOUSTICVARIABLESWITHOUTMUXDIFFERENTIATIONSTEP4} has an integrable power in $\ubar$. The second integral on RHS\,\eqref{E:SQUAREL2NORMSOFBELOWTOPORDERACOUSTICVARIABLESWITHOUTMUXDIFFERENTIATIONSTEP4} can once again be bounded by $ \int_{u' = \moreinterestingu_1}^u \totalfluxcontrolvelocity_{N'}(\ubar,u') \, \mathrm{d} u'$ (in the case $\tander^{N'+2} \velocityarray = \Lunit \tander^{N'+1}\velocityarray $) and $ \strongangularcontrolvelocity_{N'}(\ubar,u)$ (in the case $\tander^{N'+2} \velocityarray = \Yvf{A} \tander^{N'+1}\velocityarray$)  using \eqref{E:SMOOTHTORUSNORMCOMPARBLETOTANGENTIALCONTRACTIONS}, \eqref{E:UNIFIEDCOMMUTATORSPACETIMESTRONGANGULARVELOCITYCONTROL}, \eqref{E:COERCIVENESSOFNULLFLUXCONTROLWAVE}, \eqref{E:COERCIVENESSOFNULLFLUXCONTROLWAVE}, and \eqref{E:COERCIVITYOLDSPACETIMETERM}.

Combining \eqref{E:SQUAREL2NORMSOFBELOWTOPORDERACOUSTICVARIABLESWITHOUTMUXDIFFERENTIATIONSTEP3}--\eqref{E:SQUAREL2NORMSOFBELOWTOPORDERACOUSTICVARIABLESWITHOUTMUXDIFFERENTIATIONSTEP4} and applying Gr\"onwall's inequality we find that:
\begin{align} \label{E:SQUAREL2NORMSOFBELOWTOPORDERACOUSTICVARIABLESWITHOUTMUXDIFFERENTIATIONSTEP5}
		\begin{split}
			\int_{\ingoingcharacteristicsurfacetwoarg{\ubar}{[\moreinterestingu_1,u]}} \left| F_{N'} \right|^2  \volingoingnullhypersurface
			& \lesssim p_{N'} (\leftubar) + \fundbootsmall \int_{\ubar' = \leftubar}^{\ubar} p_{N'} (\ubar',u)\, \mathrm{d} \ubar  +  \int_{u' = \moreinterestingu_1}^u \totalfluxcontrolvelocity_{N'}(\ubar,u') \, \mathrm{d} u' + \strongangularcontrolvelocity_{N'}(\ubar,u) \\	
			& \ \ + \int_{u' = \moreinterestingu_1}^u \left\{\fluxcontrolVort_{\le N+1}(\ubar,u')+ \fluxcontrolGradEnt_{\le N+1}(\ubar,u')\right\} \, \mathrm{d} u',
		\end{split}
	\end{align}
The following estimate holds for the remaining case when $N' \le \Ntop -5$ so that  $0 = \blowuprateofacousticgeoWRTnumberofcommutations(N') = \blowuprateofwaveWRTnumberofcommutations(N')$:
	\begin{align} \label{E:SQUAREL2NORMSOFBELOWTOPORDERACOUSTICVARIABLESWITHOUTMUXDIFFERENTIATIONSTEP6}
		\begin{split}
			\int_{\ingoingcharacteristicsurfacetwoarg{\ubar}{[\moreinterestingu_1,u]}} \left| \varphi_{N'} \right|^2  \volingoingnullhypersurface & \lesssim p_{N'} (\leftubar) +  \int_{\ubar' = \leftubar}^{\ubar} p_{N'} (\ubar',u)\, \mathrm{d} \ubar  +   \int_{u' = \moreinterestingu_1}^u \totalfluxcontrolvelocity_{N'}(\ubar,u') \, \mathrm{d} u'  + \strongangularcontrolvelocity_{N'}(\ubar,u) \\
			& + \int_{u' = \moreinterestingu_1}^u \left\{\fluxcontrolVort_{\le N+1}(\ubar,u')+ \fluxcontrolGradEnt_{\le N+1}(\ubar,u')\right\} \, \mathrm{d} u'.
		\end{split}
	\end{align}
The proof of \eqref{E:SQUAREL2NORMSOFBELOWTOPORDERACOUSTICVARIABLESWITHOUTMUXDIFFERENTIATIONSTEP6} is similar to \eqref{E:SQUAREL2NORMSOFBELOWTOPORDERACOUSTICVARIABLESWITHOUTMUXDIFFERENTIATIONSTEP5} and we omit it, noting only that we do not need to worry about singular powers of $\weight = |\ubar|$. 

Putting together the results of \eqref{E:SQUAREL2NORMSOFBELOWTOPORDERACOUSTICVARIABLESWITHOUTMUXDIFFERENTIATIONSTEP1}--\eqref{E:SQUAREL2NORMSOFBELOWTOPORDERACOUSTICVARIABLESWITHOUTMUXDIFFERENTIATIONSTEP6}, summing over $F_{N'} \in \{ \tander^{N'+1} \Lsmall^i, \, \tander^{N'} \mytr_{\gtorus} \upchi, \, |\angLie_{\tander}^{N'} \upchi |_{\gtorus}\}$ and over $N'$ we conclude the estimate: 
\begin{align}
		\begin{split} \label{E:SQUAREL2NORMSOFBELOWTOPORDERACOUSTICVARIABLESWITHOUTMUXDIFFERENTIATIONSTEP7}
			p_N(\ubar) &  \lesssim p_{N} (\leftubar) +  (1 + \fundbootsmall) \int_{\ubar' = \leftubar}^{\ubar} p_{N} (\ubar',u)\, \mathrm{d} \ubar  + (1 + \varsigma_0^{-1})  \int_{u' = \moreinterestingu_1}^u \totalfluxcontrolvelocity_{[1,N]}(\ubar,u') \, \mathrm{d} u'  + (1 + \varsigma_0^{-1}) \strongangularcontrolvelocity_{[1,N]}(\ubar,u) \\
			& \ \  + \int_{u' = \moreinterestingu_1}^u \left\{\fluxcontrolVort_{\le N+1}(\ubar,u')+ \fluxcontrolGradEnt_{\le N+1}(\ubar,u')\right\} \, \mathrm{d} u'.
		\end{split}
	\end{align}
Applying Gr\"onwall's inequality to \eqref{E:SQUAREL2NORMSOFBELOWTOPORDERACOUSTICVARIABLESWITHOUTMUXDIFFERENTIATIONSTEP7} and estimating $p_N(\leftubar)$ using the data estimates of Sect.\,\ref{SSS:QUANTITATIVEASSUMPTIONSONDATAAWAYFROMSYMMETRY} concludes the proof of \eqref{E:PRELIMINARYBELOWTOPORDERL2ESTIMATESFOREIKONALFUINCTIONQUANTITIES} for $\varphi \in \{ \tander^{N+1} \Lsmall^i, \, \tander^{N} \mytr_{\gtorus} \upchi, \, |\angLie_{\tander}^{N} \upchi |_{\gtorus}\}$.

We note that this also proves \eqref{E:PRELIMINARYJUSTBELOWTOPORDERL2SPACETIMEESTIMATESFORLSMALL} by virtue of the $\newspacetimecoercive{\blowuprateofacousticgeoWRTnumberofcommutations(N')-1}[\varphi](\ubar,u)$ term on LHS \eqref{E:SQUAREL2NORMSOFBELOWTOPORDERACOUSTICVARIABLESWITHOUTMUXDIFFERENTIATIONSTEP1}, see \eqref{E:COERCIVITYNEWSPACETIMETERM}.

\medskip

\noindent \underline{\textbf{Transport equation estimates for commutations involving transversal derivatives:}}

We fix an integer $1 \le N \le \Ntop -1$, a point $u \in [\moreinterestingu_1,\moreinterestingu_2]$, and define the following function of $\ubar \in [\leftubar,\ubarboot)$:
\begin{align}
	\tilde{p}_N(\ubar) & \eqdef \sum_{N' = 0}^{N}   \int_{\ingoingcharacteristicsurfacetwoarg{\ubar}{[\moreinterestingu_1,u]}} \left\{ \left(\comder^{N';1} \mytr_{\gtorus} \upchi\right)^2 + \left| \angLie_{\comder}^{N';1} \upchi \right|_{\gtorus}^2 + \sum_{i=1}^3  \left( \comdersmall^{N'+1;1} \Lsmall^i\right)^2   + \left(\tandersmall^{N+1} \upmu\right)^2 \right\} \weight^{\blowuprateofacousticgeoWRTnumberofcommutations(N')}  \, \volingoingnullhypersurface. \label{E:SQUAREL2NORMSOFBELOWTOPORDERACOUSTICVARIABLESWITHMUXDIFFERENTIATION}
\end{align}
Then it can be shown that $\tilde{p}_N(\ubar)$ satisfies:
\begin{align}
		\begin{split} \label{E:SQUAREL2NORMSOFBELOWTOPORDERACOUSTICVARIABLESWITHMUXDIFFERENTIATIONSTEP1}
			\tilde{p}_N(\ubar)  & \lesssim \tilde{p}_{N} (\leftubar) +   (1 + \varsigma_0^{-1}) \int_{\ubar' = \leftubar}^{\ubar} \totalfluxcontrolvelocity_{[1,N]}(\ubar',u) \, \mathrm{d} \ubar'  +  (1 + \varsigma_0^{-1}) \int_{u' = \moreinterestingu_1}^u \totalfluxcontrolvelocity_{[1,N]}(\ubar,u') \, \mathrm{d} u'  + (1 + \varsigma_0^{-1}) \strongangularcontrolvelocity_{N}(\ubar,u) \\
			& \ \  + \int_{u' = \moreinterestingu_1}^u \left\{\fluxcontrolVort_{\le N+1}(\ubar,u')+ \fluxcontrolGradEnt_{\le N+1}(\ubar,u')\right\} \, \mathrm{d} u'.
		\end{split}
	\end{align}We omit most of the details of the proof of \eqref{E:SQUAREL2NORMSOFBELOWTOPORDERACOUSTICVARIABLESWITHMUXDIFFERENTIATIONSTEP1} as it is very similar to \eqref{E:SQUAREL2NORMSOFBELOWTOPORDERACOUSTICVARIABLESWITHOUTMUXDIFFERENTIATIONSTEP7}. The main difference is that upon using \eqref{E:LTANGENTIALMUPOINTWISE} and \eqref{E:LZLSMALLPOINTWISE}--\eqref{E:ANGLIELZCHIPOINTWISE}, one encounters a transversal derivative error term of the form $| \comdersmall^{[1,N+1];1}\wavearray |$. Using \eqref{E:L2ESTIMATESFORMUXTRANSVERSALDERIVATIVESOFWAVEVARIABLES}, the same proof as in \eqref{E:SQUAREL2NORMSOFBELOWTOPORDERACOUSTICVARIABLESWITHOUTMUXDIFFERENTIATIONSTEP7} follows. 

Finally, we note that in the process of deriving \eqref{E:SQUAREL2NORMSOFBELOWTOPORDERACOUSTICVARIABLESWITHMUXDIFFERENTIATIONSTEP1}, one similarly uses \eqref{E:SPACETIMEENERGYNULLFLUXIDENTITYFORBELOWTOPORDERACOUSTIC} and hence there are coercive bulk terms such as $\blowuprateofacousticgeoWRTnumberofcommutations(N') \newspacetimecoercive{\blowuprateofacousticgeoWRTnumberofcommutations(N')-1}[\tander^{N'} \upmu](\ubar,u)$ which we chose not to explicitly write on LHS \eqref{E:SQUAREL2NORMSOFBELOWTOPORDERACOUSTICVARIABLESWITHMUXDIFFERENTIATIONSTEP1}. The case where $N' = \Ntop$ proves \eqref{E:PRELIMINARYJUSTBELOWTOPORDERL2SPACETIMEESTIMATESFORMU}.

\end{proof}
\begin{remark}[The constant $\varsigma_0$] \label{R:SMALLNESSCONSTANTFORBELOWTOPORDERACOUSTICVARIABLES}
For the rest of the paper, we consider $\varsigma_0 \in (0,1]$ to be the fixed constant from Lemma\,\ref{L:PRELIMINARYBELOWTOPORDERL2ESTIMATESFOREIKONALFUINCTIONQUANTITIES}, with the \emph{only} exception being Lemma\,\ref{L:BELOWTOPORDERENERGYINTEGRALINEQUALITIESSPECIFICVORTICITYANDENTROPYGRADIENT} where it is potentially made smaller. It is important to make this distinction because large constants of the form $1 + \varsigma_0^{-1}$ feature in many of the $L^2$ a priori estimates, and in order to close the estimates, these large constants must be carefully compensated with a small coefficient, e.g. $\fundbootsmall(1 + \varsigma_0^{-1})$. 
\end{remark}

\subsection{$L^2$ estimates for $\velocityarray$ that lose one derivative}
\label{SS:WAVEVARIABLEL2ESTIAMTESTHATLOSEONEDERIVATIVE}

\begin{corollary}[$L^2$ estimates for $\velocityarray$ that lose one derivative]
\label{C:NONSINGULARL2ESTIMATESFORWAVEVARIABLESTHATLOSEONEDERIVATIVE}
Let $1 \leq N \leq \Ntop$,
and recall the vectorfield conventions established in Def.\,\ref{D:STRINGSOFCOMMUTATIONVECTORFIELDS}. Let $\blowuprateofwaveWRTnumberofcommutations(N)$ denote the blowup rate of the wave variables defined in \eqref{E:BLOWUPRATEOFWAVEWRTNUMBEROFCOMMUTATORS}. Let $\varsigma_0$ be the constant featured in Lemma\,\ref{L:PRELIMINARYBELOWTOPORDERL2ESTIMATESFOREIKONALFUINCTIONQUANTITIES} (see also Remark\,\ref{R:SMALLNESSCONSTANTFORBELOWTOPORDERACOUSTICVARIABLES}). Then the following estimates hold for $(\ubar,u) \in [\leftubar,\ubarboot) \times [\moreinterestingu_1,\moreinterestingu_2]$, where the implicit constants are \textbf{independent} of $\varsigma_0$:
\begin{align}
	\begin{split}
		\left\| 
			\comder^{N;1} \velocityarray \right\|_{L_{\blowuprateofwaveWRTnumberofcommutations(N)}^2\left(\ingoingcharacteristicsurfacetwoarg{\ubar}{[\moreinterestingu,u]}\right)}^2 & \lesssim  \initialsmalldoublenull^2 + \totalfluxcontrolvelocity_{[1,N]}(\ubar,u) +  (1 + \varsigma_0^{-1}) \int_{\ubar' = \leftubar}^{\ubar} \totalfluxcontrolvelocity_{[1,N-1]}(\ubar',u) \, \mathrm{d} \ubar'   \\
			& \ \ +  (1 + \varsigma_0^{-1}) \int_{u' = \moreinterestingu_1}^u \totalfluxcontrolvelocity_{[1,N-1]}(\ubar,u') \, \mathrm{d} u' + \fundbootsmall(1 + \varsigma_0^{-1}) \strongangularcontrolvelocity_{[1,N-1]}(\ubar,u) \\
			& \ \ + \int_{u' = \moreinterestingu_1}^u \left\{\fluxcontrolVort_{\le N-1}(\ubar,u')+ \fluxcontrolGradEnt_{\le N-1}(\ubar,u')\right\} \, \mathrm{d} u'.\label{E:DERIVATIVELOSINGCOERCIVITYFORWAVE}
	\end{split}
\end{align}\end{corollary} 

\begin{proof}

The starting point of the proof is applying Cor\,\ref{C:IMPROVEAUX} to \eqref{E:ARBITRARYCOMMUTATORSTRINGESTIMATE} with $N+1$ replaced by $N$, which we write as: 
	\begin{align} \label{E:ARBITRARYBELOWTOPCOMMUTATORSTRINGESTIMATE}
		\begin{split}
			\left|\comdersmall^{[1,N];1} \velocityarray \right| &  \lesssim \left|\newuL \tander^{[1,N-1]} \velocityarray \right| + \upmu \left| \nullangrmD \tander^{[1,N-1]} \velocityarray \right|_{\gnulltori} + \left| \tander^{[1,N]} \velocityarray \right|  + \fundbootsmall \left|\comdersmall^{[1,N-1];1} \controlvars \right| + \fundbootsmall \left| \tandersmall^{[1,N-1]} \badcontrolvars \right|.
		\end{split}
	\end{align}
We note that $\comder^{N;1} = \tander$ when $N = 1$ and then only the third term on RHS\,\eqref{E:ARBITRARYBELOWTOPCOMMUTATORSTRINGESTIMATE} is present. Taking $\left\| \cdot \right\|_{L_{\blowuprateofwaveWRTnumberofcommutations(N)}^2\left(\ingoingcharacteristicsurfacetwoarg{\ubar}{[\moreinterestingu,u]}\right)}^2 $ of \eqref{E:ARBITRARYBELOWTOPCOMMUTATORSTRINGESTIMATE}, using \eqref{E:COERCIVENESSOFNULLFLUXCONTROLWAVE} and the trivial bound $\weight^{\blowuprateofwaveWRTnumberofcommutations(N)} \lesssim \weight^{\blowuprateofwaveWRTnumberofcommutations(N-1)}$ to bound the first two resulting integrals by $ \totalfluxcontrolvelocity_{[1,N-1]}(\ubar,u)$, using \eqref{E:DERIVATIVELOSINGL2ESTIMATEFORTANGENTIALDERIVATIVES} to bound the third resulting integral, and using \eqref{E:PRELIMINARYBELOWTOPORDERL2ESTIMATESFOREIKONALFUINCTIONQUANTITIES}--\eqref{E:PRELIMINARYBELOWTOPORDERL2ESTIMATESFOREIKONALFUINCTIONQUANTITIESWITHZ} to bound the fourth and fifth resulting integrals, we conclude the proof.
\end{proof}


\section{Below-top-order hyperbolic $L^2$ estimates for the specific vorticity and entropy gradient}
\label{S:BELOWTOPORDERHYPERBOLICL2ESTIMATESFORSPECIFICVORTICITYANDENTROPYGRADIENT}
We continue to work under the assumptions of Sect.\,\ref{SS:SILENTFACTS}.
In this short section, 
we prove the below-top-order $L^2$ estimates for the specific vorticity and entropy gradient as well as the modified fluid variables. 
Specifically, we prove Prop.\,\ref{P:APRIORIL2ESTIMATESFORTRANSPORTVARIABLES}. We derive a preliminary energy integral inequality 
in Sect.\,\ref{SS:BELOWTOPORDERHYPERBOLICL2ESTIMATESFORSPECIFICVORTICITYANDENTROPYGRADIENT},
and we prove the final estimates in 
Sects.\,\ref{SS:PROOFOFBELOWTOPORDERTRANSPORTENERGYESTIMATES}--\ref{SS:PROOFOFBELOWTOPORDERROUGHTORIENERGYESTIMATES}. These estimates are relatively straightforward consequences of
the transport energy identity \eqref{E:ENERGYNULLFLUXINTEGRALIDENTITIESTRANSPORT}
and various pointwise estimates we have already established,
including the ones provided by Prop.\,\ref{P:POINTWISESTIAMTESFORALLTHETRANSPORTVARIABLES}.
In Sect.\,\ref{S:TOPORDERELIPTICHYPERBOLICL2ESTIMATESFORSPECIFICVORTICITYANDENTROPYGRADIENT},
we will prove the top-order estimates of Prop.\,\ref{P:TOPORDERMODIFIEDAPRIORIL2ESTIMATESANDDOUBLENULLTORIESTIMATES}. The proofs of these estimates are much more difficult because they rely on
the intricate elliptic-hyperbolic integral identity \eqref{E:INTEGRALIDENTITYFORELLIPTICHYPERBOLICCURRENT}.

In Sect.\,\ref{S:WAVEANDACOUSTICGEOMETRYAPRIORIESTIMATES}, we will use the estimates for 
$\vortrenormalized$,
$\GradEnt$,
$\VortVort$, and $\DivGradEnt$
that we derive in this section in our proof of the wave a priori estimates, 
which we stated as Prop.\,\ref{P:APRIORIL2ESTIMATESWAVEVARIABLES}.
Hence, we highlight that for the logic of the paper, it is important that 
\textbf{the estimates we derive in this section do not rely on the wave estimates of Prop.\,\ref{P:APRIORIL2ESTIMATESWAVEVARIABLES}};
our proofs there 
instead rely on the bootstrap assumptions
\eqref{E:MAINWAVEENERGYBOOTSTRAP}
for the wave energies, which are \emph{weaker} than the estimates that we stated  
in Prop.\,\ref{P:APRIORIL2ESTIMATESWAVEVARIABLES}.

\subsection{Integral inequalities for the below-top-order vorticity- and entropy gradient-controlling quantities}
\label{SS:BELOWTOPORDERHYPERBOLICL2ESTIMATESFORSPECIFICVORTICITYANDENTROPYGRADIENT}
We begin with the following preliminary lemma,
which provides integral inequalities
for the below-top-order vorticity- and entropy gradient-controlling quantities.

\begin{lemma}[Integral inequalities for the below-top-order vorticity- and entropy gradient-controlling quantities] 
\label{L:BELOWTOPORDERENERGYINTEGRALINEQUALITIESSPECIFICVORTICITYANDENTROPYGRADIENT}
	Let $0 \leq N \leq \Ntop$, and recall that the $L^2$-controlling quantities
	$\fluxcontrolVort_{ N}$, 
	$\fluxcontrolGradEnt_{ N}$,
	$\cdots$, are defined in
	Defs.\,\ref{D:MAINCOERCIVE} and\,\ref{D:SUMMEDL2CONTROLLINGQUANTITIES}.
	Then by making the constant $\varsigma_0$ from Lemma\,\ref{L:PRELIMINARYBELOWTOPORDERL2ESTIMATESFOREIKONALFUINCTIONQUANTITIES} (see also Remark\,\ref{R:SMALLNESSCONSTANTFORBELOWTOPORDERACOUSTICVARIABLES}) potentially even smaller, the following integral inequality holds for
	$(\ubar,u) \in [\leftubar,\ubarboot) \times [\moreinterestingu_1,\moreinterestingu_2]$:
	\begin{align}   \label{E:BELOWTOPORDERENERGYINTEGRALINEQUALITIESSPECIFICVORTICITYANDENTROPYGRADIENT}
		\begin{split}
			& \fluxcontrolVort_{N}(\ubar,u)
			+
			\fluxcontrolGradEnt_{N}(\ubar,u) + \blowuprateoftransportWRTnumberofcommutations(N) \left( \bulkcontrolVort_N(\ubar,u) + \bulkcontrolGradEnt_N(\ubar,u) \right)\\
			& \ \   \lesssim
			\initialsmall^2
				+ 
			\int_{\ubar' = \leftubar}^{\ubar} \frac{1}{|\ubar'|^{0.13}} 
			\left\lbrace 
				\fluxcontrolVort_{\leq N}(\ubar,u')
				+
				\fluxcontrolGradEnt_{\leq N}(\ubar,u')
			\right\rbrace
			\, \rmd \ubar' 
			\\	
			& \ \ \ \ +
			\int_{u' = \moreinterestingu_1}^u
				\left\lbrace 
					\fluxcontrolVort_{\leq N}(\ubar,u')
					+
					\fluxcontrolGradEnt_{\leq N}(\ubar,u')
				\right\rbrace
			\, \rmd u'
			+
			(1+\varsigma_0^{-1})\fundbootsmall^2 
			\int_{\ubar' = \leftubar}^{\ubar} \frac{1}{|\ubar'|^{0.13}} 
				\totalfluxcontrolvelocity_{[1,N]}(\ubar',u) 
			\, \rmd \ubar' \\
			& \ \ \ \ + (1 + \varsigma_0^{-1})\fundbootsmall^2 \int_{u' = \moreinterestingu_1}^u \totalfluxcontrolvelocity_{[1,N+1]}(\ubar,u') \, \mathrm{d} u' + (1+\varsigma_0^{-1}) \fundbootsmall^2 \strongangularcontrolvelocity_{[1,N]}(\ubar,u) \\
			& \ \ \ \ + \fundbootsmall^2 \weakspacetimetransvelocitycontrol_{[1,N]}(\ubar,u) + \fundbootsmall^2 \weakspacetimeangularvelocitycontrol_{[1,N]}(\ubar,u)
		\end{split}
	\end{align}

	Moreover, the following estimates hold for $0 \leq N \leq \Ntop - 1$:
	\begin{align} 
	\begin{split} \label{E:BELOWTOPORDERENERGYINTEGRALINEQUALITIESFORMODIFIEDFLUID}
		\fluxcontrolVortVort_{\leq N}(\ubar,u)
		+
		\fluxcontrolDivGradEnt_{\leq N}(\ubar,u)
		& \lesssim
		\initialsmall^2 +	\int_{\ubar' = \leftubar}^{\ubar}	\left\lbrace 
				\fluxcontrolVortVort_{\leq N}(\ubar,u')
				+
				\fluxcontrolDivGradEnt_{\leq N}(\ubar,u')
			\right\rbrace
			\, \rmd \ubar' 
		+
		\int_{u' = \moreinterestingu_1}^u
			\left\lbrace
				\;{\leq N}(\ubar,u')
				+
				\fluxcontrolDivGradEnt_{\leq N}(\ubar,u)
			\right\rbrace
		\, \rmd u'
			\\
	& \ \
		+
		\int_{u' = \moreinterestingu_1}^u 
			\left\lbrace
				\fluxcontrolVort_{\leq N+1}(\ubar,u')
				+
				\fluxcontrolGradEnt_{\leq N+1}(\ubar,u')
			\right\rbrace
		\, \rmd u'	
			+
		(1 + \varsigma_0^{-1})\fundbootsmall^2 
		\int_{\ubar' = \leftubar}^{\ubar}	\totalfluxcontrolvelocity_{[1,N]}(\ubar',u) 
		\, \rmd \ubar' \\
		& \ \ 
		+
		(1 + \varsigma_0^{-1})\fundbootsmall^2 
		\int_{u' = \moreinterestingu_1}^u 
			\totalfluxcontrolvelocity_{[1,N]}(\ubar,u')
		\, \rmd u' 
		+ 
		(1 + \varsigma_0^{-1})\fundbootsmall^2 \strongangularcontrolvelocity_{[1,N]}(\ubar,u).
	\end{split}
	\end{align}
	
\end{lemma}

\begin{proof}
	We first prove \eqref{E:BELOWTOPORDERENERGYINTEGRALINEQUALITIESSPECIFICVORTICITYANDENTROPYGRADIENT}.
	Let $N \in \{ 0,\dots, \Ntop\}$ and consider the transport energy identity \eqref{E:ENERGYNULLFLUXINTEGRALIDENTITIESTRANSPORT}
	with $(\tander^N \vortrenormalized,\tander^N \GradEnt)$ in the role of $f$ and $q = \blowuprateoftransportWRTnumberofcommutations(N)$ is as in Def.\,\ref{D:NOTATIONFORPOWERSOFWEIGHT}.
	We use the bootstrap assumptions and Prop.\,\ref{P:IMPROVEMENTOFAUXILIARYBOOTSTRAP} to deduce the bound
	$|\Lunit \upmu + \upmu \mytr_{\gtorus} \angk| \lesssim 1$ for the integrand factor in the last
	integral on RHS\,\eqref{E:ENERGYNULLFLUXINTEGRALIDENTITIESTRANSPORT},
	and we use \eqref{E:CLOSEDKEYJACOBIANDETERMINANTESTIMATECHOVDOUBLENULLTOGEO} to deduce that
	the integrand factors
	$\MagnitueofinnerproductofLunitanduLunit = \frac{1}{\Lunit \ubar} \approx 1$.
	We also use
	Lemma\,\ref{L:ALLL2CONTROLLINGQUANTITIESINITIALLYSMALL}
	to bound the data-dependent terms 
	$ \ingoingfluxtransport{\blowuprateoftransportWRTnumberofcommutations(N)}[(\tander^N\vortrenormalized,\tander^N\GradEnt)](\leftubar,u) + \outgoingfluxtransport{\blowuprateoftransportWRTnumberofcommutations(N)}[(\tander^N \vortrenormalized,\tander^N\GradEnt)](\ubar,\moreinterestingu_1) $
	on RHS\,\eqref{E:ENERGYNULLFLUXINTEGRALIDENTITIESTRANSPORT} by $\lesssim \initialsmall^2$. Then by Def.\,\ref{D:MAINCOERCIVE},
	the estimate \eqref{E:CLOSEDKEYJACOBIANDETERMINANTESTIMATECHOVDOUBLENULLTOGEO},
	and
	Lemma\,\ref{L:COERCIVENESSOFL2CONTROLLINGQUANITIESFORVORTICITYANDENTROPY},
	Young's inequality, and the pointwise estimate \eqref{E:COMMUTEDTRANSPORTPOINTWISEESTIMATESFORSPECIFICVORTICITYANDENTROPYGRADIENT}
	to bound the integrand factors of $\left|\upmu \Transport \tander^{N} (\vortrenormalized,\GradEnt) \right|$
	on RHS\,\eqref{E:ENERGYNULLFLUXINTEGRALIDENTITIESTRANSPORT},
	we deduce:
		\begin{align}  \label{E:FIRSTSTEPPROOFBELOWTOPORDERENERGYINTEGRALINEQUALITIESSPECIFICVORTICITYANDENTROPYGRADIENT}
		\begin{split}
		& \fluxcontrolVort_{N}(\ubar,u)
		+
		\fluxcontrolGradEnt_{N}(\ubar,u) + \blowuprateoftransportWRTnumberofcommutations(N) \left( \bulkcontrolVort_N(\ubar,u) + \bulkcontrolGradEnt_N(\ubar,u) \right)
		 \\
		&  \ \  \lesssim
		\initialsmall^2
		+
		\int_{u' = \moreinterestingu_1}^u
			\left\lbrace
				\fluxcontrolVort_{\le N}(\ubar,u')
				+
				\fluxcontrolGradEnt_{\le N}(\ubar,u')
			\right\rbrace
		\, \rmd u'
			\\
	& \ \ \ + \max_{\tander^N \in \mathfrak{P}^{(N)}}
		\int_{\characteristicdiamondtwoarg{[\leftubar,\ubar)}{[\moreinterestingu_1,u]}} \left| \tander^{N}(\vortrenormalized,\GradEnt)\right|
			\left\{ \fundbootsmall |\newuL \tander^{[1,N]} \velocityarray|
			+
			\fundbootsmall \upmu |\tander^{[1,N + 1]} \velocityarray|
			+
			\fundbootsmall |\tandersmall^{[1,N]} \badcontrolvars|\right\}  \weight^{\blowuprateoftransportWRTnumberofcommutations(N)}
		\, \voldiamond.
	\end{split}
	\end{align}
	We note that the error terms generated by $\tander^{\le N}(\vortrenormalized,\GradEnt)$ on RHS\,\eqref{E:COMMUTEDTRANSPORTPOINTWISEESTIMATESFORSPECIFICVORTICITYANDENTROPYGRADIENT} are bounded by the first integral $\int_{u' = \moreinterestingu_1}^u\{\fluxcontrolVort_{\le N}(\ubar,u')
				+
				\fluxcontrolGradEnt_{\le N}(\ubar,u') \}
		\, \rmd u'$ on RHS\,\eqref{E:FIRSTSTEPPROOFBELOWTOPORDERENERGYINTEGRALINEQUALITIESSPECIFICVORTICITYANDENTROPYGRADIENT}.
		
	Suppose that $\Ntop - 2 \le N' \le \Ntop$, which implies that $\blowuprateoftransportWRTnumberofcommutations(N') > 0$, $\blowuprateofwaveWRTnumberofcommutations(N')-1 = \blowuprateoftransportWRTnumberofcommutations(N')$, and $\blowuprateoftransportWRTnumberofcommutations(N') > \blowuprateofwaveWRTnumberofcommutations(N'')$ for all $1 \le N'' < N'$ by Def.\,\ref{D:NOTATIONFORPOWERSOFWEIGHT}. Then, it follows that:
	\begin{align}
		\begin{split} \label{E:SECONDSTEPPROOFBELOWTOPORDERENERGYINTEGRALINEQUALITIESSPECIFICVORTICITYANDENTROPYGRADIENT}
			& \int_{\characteristicdiamondtwoarg{[\leftubar,\ubar)}{[\moreinterestingu_1,u]}} \left| \tander^{N'}(\vortrenormalized,\GradEnt)\right| \left\{ \fundbootsmall  |\newuL \tander^{[1,N']} \velocityarray|\right\} \weight^{\blowuprateoftransportWRTnumberofcommutations(N)}
			\, \voldiamond \\
			& \ \  \lesssim \int_{u' = \moreinterestingu_1}^u
			\left\lbrace
				\fluxcontrolVort_{N'}(\ubar,u')
				+
				\fluxcontrolGradEnt_{N'}(\ubar,u')
			\right\rbrace
			\, \rmd u'  + \fundbootsmall^2 \int_{\characteristicdiamondtwoarg{[\leftubar,\ubar)}{[\moreinterestingu_1,u]}} \left|\newuL \tander^{N'} \velocityarray\right|^2 \weight^{\blowuprateofwaveWRTnumberofcommutations(N')-1} \, \voldiamond \\
			& \ \ \ \ 	+
			\fundbootsmall^2 
			\int_{\ubar' = \leftubar}^{\ubar} 
			\totalfluxcontrolvelocity_{[1,N'-1]}(\ubar',u) 
			\, \rmd \ubar' \\
			& \ \  \lesssim \int_{u' = \moreinterestingu_1}^u\left\lbrace
				\fluxcontrolVort_{N'}(\ubar,u')
				+
				\fluxcontrolGradEnt_{N'}(\ubar,u')
			\right\rbrace
			\, \rmd u' 
			+ \fundbootsmall^2 \weakspacetimetransvelocitycontrol_{N'}(\ubar,u) \\
			& \ \ \ \ + \fundbootsmall^2 
			\int_{\ubar' = \leftubar}^{\ubar} 
			\totalfluxcontrolvelocity_{[1,N'-1]}(\ubar',u) 
			\, \rmd \ubar' ,			 
		\end{split}
	\end{align}
where we used \eqref{E:WAVENEWSPACETIMEL2TRANSVERSALCONTROLLINGQUANTITY} in the last inequality in \eqref{E:SECONDSTEPPROOFBELOWTOPORDERENERGYINTEGRALINEQUALITIESSPECIFICVORTICITYANDENTROPYGRADIENT}. Under the same assumption that $\blowuprateoftransportWRTnumberofcommutations(N') > 0$, we now focus on the integrals generated by the terms $\fundbootsmall\upmu|\tander^{[1,N'+1]} \velocityarray|$ in RHS\,\eqref{E:FIRSTSTEPPROOFBELOWTOPORDERENERGYINTEGRALINEQUALITIESSPECIFICVORTICITYANDENTROPYGRADIENT}. By \eqref{E:POINTWISESMOOTHTORIDIFFERENTIALNORMINTERMSOFDOULBENULLTORIDIFFERENTIALINORMANDL} and Cor.\,\ref{C:IMPROVEAUX}, we bound: 
	\begin{align}
		\fundbootsmall\upmu|\tander^{[1,N'+1]} \velocityarray| & \lesssim \fundbootsmall \upmu |\nullangrmD \tander^{[1,N']} \velocityarray|_{\gnulltori} + \fundbootsmall \upmu |\Lunit \tander^{[1,N']} \velocityarray|. \label{E:THIRDSTEPPROOFBELOWTOPORDERENERGYINTEGRALINEQUALITIESSPECIFICVORTICITYANDENTROPYGRADIENT}
	\end{align}
Using now \eqref{E:WAVENEWSPACETIMENULLANGULARL2CONTROLLINGQUANTITY}, the error integral generated by $ \fundbootsmall \upmu |\nullangrmD \tander^{[1,N']} \velocityarray|_{\gnulltori}$ in RHS\,\eqref{E:THIRDSTEPPROOFBELOWTOPORDERENERGYINTEGRALINEQUALITIESSPECIFICVORTICITYANDENTROPYGRADIENT} can be bounded exactly as in \eqref{E:SECONDSTEPPROOFBELOWTOPORDERENERGYINTEGRALINEQUALITIESSPECIFICVORTICITYANDENTROPYGRADIENT} with $ \fundbootsmall^2\weakspacetimeangularvelocitycontrol_{N'}(\ubar,u)$ in place of $\fundbootsmall^2 \weakspacetimetransvelocitycontrol_{N'}(\ubar,u)$, where we emphasize that we must utilize (at least) half of the $\upmu$ present on RHS\,\eqref{E:THIRDSTEPPROOFBELOWTOPORDERENERGYINTEGRALINEQUALITIESSPECIFICVORTICITYANDENTROPYGRADIENT} to apply the coercive bound \eqref{E:WAVENEWSPACETIMENULLANGULARL2CONTROLLINGQUANTITY}. 

Suppose now $N' = \Ntop -3$, which implies that $\blowuprateoftransportWRTnumberofcommutations(N') = 0$ but $0 < \blowuprateofwaveWRTnumberofcommutations(N') = 0.13 < 1$. Then, we bound: 
	\begin{align}
		\begin{split} \label{E:FIFTHSTEPPROOFBELOWTOPORDERENERGYINTEGRALINEQUALITIESSPECIFICVORTICITYANDENTROPYGRADIENT}
			& \int_{\characteristicdiamondtwoarg{[\leftubar,\ubar)}{[\moreinterestingu_1,u]}} \left| \tander^{N'}(\vortrenormalized,\GradEnt)\right| \left\{ \fundbootsmall  |\newuL \tander^{N'} \velocityarray|\right\} 
			\, \voldiamond \\
			& \ \  \lesssim \int_{\characteristicdiamondtwoarg{[\leftubar,\ubar)}{[\moreinterestingu_1,u]}} \left| \tander^{N'}(\vortrenormalized,\GradEnt)\right|^2 \, \voldiamond + \fundbootsmall^2 \int_{\characteristicdiamondtwoarg{[\leftubar,\ubar)}{[\moreinterestingu_1,u]}}  \frac{1}{\weight^{\blowuprateofwaveWRTnumberofcommutations(N')}} \left| \newuL \tander^{N'}\velocityarray\right|^2 \weight^{\blowuprateofwaveWRTnumberofcommutations(N')} \, \voldiamond \\
			& \ \ \lesssim \int_{u' = \moreinterestingu_1}^u \left\lbrace
				\fluxcontrolVort_{N'}(\ubar,u')
				+
				\fluxcontrolGradEnt_{N'}(\ubar,u')
			\right\rbrace
			\, \rmd u' + \fundbootsmall^2 \int_{\ubar' = \leftubar}^{\ubar} \frac{1}{|\ubar'|^{\blowuprateofwaveWRTnumberofcommutations(N')}} \totalfluxcontrolvelocity_{N'}(\ubar',u) \, \rmd \ubar'.
		\end{split}
	\end{align}
The integral generated by $\fundbootsmall\upmu |\nullangrmD \tander^{N'} \velocityarray|_{\gnulltori}$ in \eqref{E:THIRDSTEPPROOFBELOWTOPORDERENERGYINTEGRALINEQUALITIESSPECIFICVORTICITYANDENTROPYGRADIENT} is handled similarly and we omit the details.

We now bound the error integrals from $\fundbootsmall\upmu |\Lunit \tander^{N'} \velocityarray|$ in RHS\,\eqref{E:THIRDSTEPPROOFBELOWTOPORDERENERGYINTEGRALINEQUALITIESSPECIFICVORTICITYANDENTROPYGRADIENT}. If $\blowuprateoftransportWRTnumberofcommutations(N') > 0$, we capitalized on the additional factor of $\upmu \approx \ReciprocaluLunitAppliedtoTimeFunction$ (see \eqref{E:RECIPROCALULUNITAPPLIEDTOTIMEFUNCTIONAPPROXIMATELYMU}) that is present, multiply and divide by $\weight^{1/2}$, 
and apply Young's inequality to bound: 
\begin{align}
		\begin{split} \label{E:FOURTHSTEPPROOFBELOWTOPORDERENERGYINTEGRALINEQUALITIESSPECIFICVORTICITYANDENTROPYGRADIENT}
			&  \int_{\characteristicdiamondtwoarg{[\leftubar,\ubar)}{[\moreinterestingu_1,u]}} \fundbootsmall  \upmu \left| \tander^{N'} (\vortrenormalized,\GradEnt)\right| \left|\Lunit \tander^{[1,N']} \velocityarray \right| \weight^{\blowuprateoftransportWRTnumberofcommutations(N')}\,  \voldiamond \\
			& \ \ \lesssim \varsigma_0  \int_{\characteristicdiamondtwoarg{[\leftubar,\ubar)}{[\moreinterestingu_1,u]}}   \ReciprocaluLunitAppliedtoTimeFunction \left| \tander^{N'} (\vortrenormalized,\GradEnt)\right|^2 \weight^{\blowuprateoftransportWRTnumberofcommutations(N')-1}\,  \voldiamond + (1+\varsigma_0^{-1})\fundbootsmall^2  \int_{\characteristicdiamondtwoarg{[\leftubar,\ubar)}{[\moreinterestingu_1,u]}}  \left| \Lunit \tander^{N'} \velocityarray\right|^2 \weight^{\blowuprateoftransportWRTnumberofcommutations(N')+1}\,  \voldiamond \\
			& \ \ \ \  + \fundbootsmall^2 
			\int_{u' = \moreinterestingu_1}^{u} 
			\totalfluxcontrolvelocity_{[1,N'-1]}(\ubar,u') 
			\, \rmd u' \\
			& \ \  \lesssim \varsigma_0 \left( \bulkcontrolVort_{N'}(\ubar,u) + \bulkcontrolGradEnt_{N'}(\ubar,u) \right) + (1+\varsigma_0^{-1})\fundbootsmall^2 
			\int_{u' = \moreinterestingu_1}^{u} 
			\totalfluxcontrolvelocity_{[1,N']}(\ubar,u') 
			\, \rmd u'.
		\end{split}
	\end{align}
By making $\varsigma_0$ potentially even smaller, the two terms $\varsigma_0 \left( \bulkcontrolVort_{N'}(\ubar,u) + \bulkcontrolGradEnt_{N'}(\ubar,u) \right)$ can be absorbed by the terms $ \blowuprateoftransportWRTnumberofcommutations(N) \left( \bulkcontrolVort_N(\ubar,u) + \bulkcontrolGradEnt_N(\ubar,u) \right)$ on LHS \eqref{E:FIRSTSTEPPROOFBELOWTOPORDERENERGYINTEGRALINEQUALITIESSPECIFICVORTICITYANDENTROPYGRADIENT}. In the case where $ \blowuprateoftransportWRTnumberofcommutations(N') =0$ but $ 0 < \blowuprateofwaveWRTnumberofcommutations(N') < 1$, we bound this by: 
\begin{align}
		\begin{split} \label{E:SIXTHSTEPPROOFBELOWTOPORDERENERGYINTEGRALINEQUALITIESSPECIFICVORTICITYANDENTROPYGRADIENT}
			&  \int_{\characteristicdiamondtwoarg{[\leftubar,\ubar)}{[\moreinterestingu_1,u]}} \fundbootsmall  \upmu \left| \tander^{N'} (\vortrenormalized,\GradEnt)\right| \left|\Lunit \tander^{[1,N']} \velocityarray \right| \,  \voldiamond \\
			& \ \ \lesssim  \int_{\characteristicdiamondtwoarg{[\leftubar,\ubar)}{[\moreinterestingu_1,u]}}   \ReciprocaluLunitAppliedtoTimeFunction  \left| \tander^{N'} (\vortrenormalized,\GradEnt)\right|^2 \weight^{-\blowuprateofwaveWRTnumberofcommutations(N')}\,  \voldiamond + \fundbootsmall^2  \int_{\characteristicdiamondtwoarg{[\leftubar,\ubar)}{[\moreinterestingu_1,u]}}  \left| \Lunit \tander^{N'} \velocityarray\right|^2 \weight^{\blowuprateofwaveWRTnumberofcommutations(N')}\,  \voldiamond \\
			& \ \ \ \  + \fundbootsmall^2 
			\int_{u' = \moreinterestingu_1}^{u} 
			\totalfluxcontrolvelocity_{[1,N'-1]}(\ubar,u') 
			\, \rmd u' \\			
			& \ \ \lesssim \int_{\ubar' = \leftubar}^{\ubar} \frac{1}{|\ubar|^{\blowuprateofwaveWRTnumberofcommutations(N')}}	\left\lbrace
				\fluxcontrolVort_{N'}(\ubar',u)
				+
				\fluxcontrolGradEnt_{N'}(\ubar',u)
			\right\rbrace
			\, \rmd \ubar' +   \fundbootsmall^2 
			\int_{u' = \moreinterestingu_1}^{u} 
			\totalfluxcontrolvelocity_{[1,N']}(\ubar,u') 
			\, \rmd u'.
		\end{split}
	\end{align}We omit the proofs for the analogous bounds in \eqref{E:SECONDSTEPPROOFBELOWTOPORDERENERGYINTEGRALINEQUALITIESSPECIFICVORTICITYANDENTROPYGRADIENT}--\eqref{E:SIXTHSTEPPROOFBELOWTOPORDERENERGYINTEGRALINEQUALITIESSPECIFICVORTICITYANDENTROPYGRADIENT} in the remaining cases where $\blowuprateoftransportWRTnumberofcommutations(N') = \blowuprateofwaveWRTnumberofcommutations(N') = 0$ because they are strictly easier and do not require careful attention to the weights.

Finally, to handle the error integrals generated by the $\tandersmall^{[1,N]} \badcontrolvars$ on RHS\,\eqref{E:FIRSTSTEPPROOFBELOWTOPORDERENERGYINTEGRALINEQUALITIESSPECIFICVORTICITYANDENTROPYGRADIENT}, we use that $\blowuprateoftransportWRTnumberofcommutations(N) = \blowuprateofacousticgeoWRTnumberofcommutations(N-1)$ for all $1 \le N' \le \Ntop$ and so we bound: 
\begin{align}
		\begin{split} \label{E:SEVENTHSTEPPROOFBELOWTOPORDERENERGYINTEGRALINEQUALITIESSPECIFICVORTICITYANDENTROPYGRADIENT}
			& \int_{\characteristicdiamondtwoarg{[\leftubar,\ubar)}{[\moreinterestingu_1,u]}} \left| \tander^{N'}(\vortrenormalized,\GradEnt)\right| \left\{ \fundbootsmall  \left|\tandersmall^{N'} (\upmu,\Lunit^i) \right|\right\} \weight^{\blowuprateoftransportWRTnumberofcommutations(N)} 
			\, \voldiamond \\
			& \ \ \lesssim \int_{u' = \moreinterestingu_1}^u \left\lbrace
				\fluxcontrolVort_{N'}(\ubar,u')
				+
				\fluxcontrolGradEnt_{N'}(\ubar,u')
			\right\rbrace
			\, \rmd u' + \fundbootsmall^2 \int_{\ubar =\leftubar}^{\ubar}  \left\| \tander^{N'}(\upmu,\Lunit^i)\right\|_{L^2_{\blowuprateofacousticgeoWRTnumberofcommutations(N-1)}\left(\ingoingcharacteristicsurfacetwoarg{\ubar'}{[\moreinterestingu_1,u]}\right)}^2 \rmd \ubar'. 
		\end{split}
	\end{align}In total, inserting \eqref{E:PRELIMINARYBELOWTOPORDERL2ESTIMATESFOREIKONALFUINCTIONQUANTITIES}--\eqref{E:PRELIMINARYBELOWTOPORDERL2ESTIMATESFOREIKONALFUINCTIONQUANTITIESWITHZ} into RHS\,\eqref{E:SEVENTHSTEPPROOFBELOWTOPORDERENERGYINTEGRALINEQUALITIESSPECIFICVORTICITYANDENTROPYGRADIENT},  summing up \eqref{E:SECONDSTEPPROOFBELOWTOPORDERENERGYINTEGRALINEQUALITIESSPECIFICVORTICITYANDENTROPYGRADIENT}--\eqref{E:SEVENTHSTEPPROOFBELOWTOPORDERENERGYINTEGRALINEQUALITIESSPECIFICVORTICITYANDENTROPYGRADIENT}, applying Gr\"onwall's inequality on the first integral on RHS\,\eqref{E:SIXTHSTEPPROOFBELOWTOPORDERENERGYINTEGRALINEQUALITIESSPECIFICVORTICITYANDENTROPYGRADIENT} (which relies on the integrability of $\frac{1}{|\ubar|^{\blowuprateofwaveWRTnumberofcommutations(N')}}$ whenever $0 < \blowuprateofwaveWRTnumberofcommutations(N) < 1$, we conclude the proof of \eqref{E:BELOWTOPORDERENERGYINTEGRALINEQUALITIESSPECIFICVORTICITYANDENTROPYGRADIENT}.

	The estimate \eqref{E:BELOWTOPORDERENERGYINTEGRALINEQUALITIESFORMODIFIEDFLUID}
	can be proved using similar arguments
	based on the pointwise estimates 
	\eqref{E:COMMUTEDTRANSPORTPOINTWISEESTIMATESFORMODIFIEDCURLOFVORT}--\eqref{E:COMMUTEDTRANSPORTPOINTWISEESTIMATESFORMODIFIEDDIVERGENCEOFENTROPYGRADIENT}
	for $0 \leq N \leq \Ntop-1$ and $\blowuprateoftransportWRTnumberofcommutations(N) = \blowuprateofacousticgeoWRTnumberofcommutations(N-1)$. 
	Although we do not provide full details, we point out the additional ingredients that play a role in the proof. Firstly, since $N \le \Ntop-1$, it follows that the highest number of derivatives on the velocity variables are $\newuL \tander^{\Ntop-1}\velocityarray$ and $\tander^{\Ntop} \velocityarray$. For these terms, the required power of the weight to be able to apply the coercivity estimates of Lemma\,\ref{L:COERCIVENESSOFL2CONTROLLINGQUANITIESUNIFIEDCOMMUTATOR} is $\blowuprateofwaveWRTnumberofcommutations(\Ntop-1) = \blowuprateoftransportWRTnumberofcommutations(\Ntop) -1$. In particular, there is never a need to multiply and divide the integrants by powers of $\weight$ (nor a need to apply Gronwall's inequality) to conclude \eqref{E:BELOWTOPORDERENERGYINTEGRALINEQUALITIESFORMODIFIEDFLUID}. In addition, since there is no power of $\upmu$ in the terms $\tander^{\Ntop} \velocityarray$ on RHS\,\eqref{E:COMMUTEDTRANSPORTPOINTWISEESTIMATESFORMODIFIEDCURLOFVORT}--\eqref{E:COMMUTEDTRANSPORTPOINTWISEESTIMATESFORMODIFIEDDIVERGENCEOFENTROPYGRADIENT}, we do not appeal to \eqref{E:POINTWISESMOOTHTORIDIFFERENTIALNORMINTERMSOFDOULBENULLTORIDIFFERENTIALINORMANDL} as we did in \eqref{E:THIRDSTEPPROOFBELOWTOPORDERENERGYINTEGRALINEQUALITIESSPECIFICVORTICITYANDENTROPYGRADIENT}. Instead, we use the aforementioned gain in weights when $\tander^{\Ntop} \velocityarray =\Lunit \tander^{\Ntop-1} \velocityarray$ and bound the error integrals generated by $\tander^{\Ntop} \velocityarray = \Yvf{A} \tander^{\Ntop-1} \velocityarray$ using the bulk energy (after applying Young's inequality) by $\fundbootsmall^2 \strongangularcontrolvelocity_{\Ntop-1}$.\end{proof}

\subsection{Proof of \eqref{E:APRIORIL2ESTIMATESFORVORTICITYFLUX}--\eqref{E:APRIORIL2BELOWTOPESTIMATESFORDIVGRADENTFLUX}}
\label{SS:PROOFOFBELOWTOPORDERTRANSPORTENERGYESTIMATES}
We first prove \eqref{E:APRIORIL2ESTIMATESFORVORTICITYFLUX}--\eqref{E:APRIORIL2ESTIMATESFORGRADENTFLUX}. We set
$\mathbb{T}(\ubar,u)  
\eqdef 	
	\fluxcontrolVort_{\leq \Ntop}(\ubar,u) 
	+ 
	\fluxcontrolGradEnt_{\leq \Ntop}(\ubar,u) +  \bulkcontrolVort_{\le \Ntop} (\ubar,u) + \bulkcontrolGradEnt_{\Ntop} (\ubar,u)$.
From
\eqref{E:BELOWTOPORDERENERGYINTEGRALINEQUALITIESSPECIFICVORTICITYANDENTROPYGRADIENT}, the integrability $\int_{\ubar' = \leftubar}^{\ubar} \frac{1}{|\ubar'|^{\blowuprateofwaveWRTnumberofcommutations(N')}} \, \mathrm{d} \ubar' \lesssim 1$,
the wave energy bootstrap assumptions \eqref{E:MAINWAVEENERGYBOOTSTRAP},
and \eqref{E:NONLINEARINEQUALITYRELATINGDATAEPSILONANDBOOTSTRAPEPSILON},
we find that
$
\mathbb{T}(\ubar,u)
\leq 
C \initialsmall^2
+
C
\int_{u' = \moreinterestingu_1}^u 
	\mathbb{T}(\ubar,u')
\, \rmd u' + \int_{\ubar' = \leftubar}^{\ubar} \frac{1}{|\ubar'|^{\blowuprateofwaveWRTnumberofcommutations(N')}} 
	\mathbb{T}(\ubar',u) 
\, \rmd \ubar''
$.
Applying Gr\"{o}nwall's inequality, we conclude that $\mathbb{T}(\ubar,u)
\leq 
C \initialsmall^2$, which yields \eqref{E:APRIORIL2ESTIMATESFORVORTICITYFLUX}--\eqref{E:APRIORIL2ESTIMATESFORGRADENTFLUX}.

Similarly, to prove \eqref{E:APRIORIL2BELOWTOPESTIMATESFORVORTVORTFLUX}--\eqref{E:APRIORIL2BELOWTOPESTIMATESFORDIVGRADENTFLUX}, 
 set
$\mathbb{T}(\ubar,u)  
\eqdef 	
	\fluxcontrolVortVort_{\leq \Ntop-1}(\ubar,u) 
+ 
\fluxcontrolDivGradEnt_{\leq \Ntop-1}(\ubar,u) +  \bulkcontrolVortVort_{\le \Ntop-1}(\ubar,u) + \bulkcontrolDivGradEnt_{\Ntop-1}(\ubar,u)$.
From
\eqref{E:BELOWTOPORDERENERGYINTEGRALINEQUALITIESFORMODIFIEDFLUID},
the wave energy bootstrap assumptions \eqref{E:MAINWAVEENERGYBOOTSTRAP},
\eqref{E:NONLINEARINEQUALITYRELATINGDATAEPSILONANDBOOTSTRAPEPSILON},
and the already proved estimates \eqref{E:APRIORIL2ESTIMATESFORVORTICITYFLUX}--\eqref{E:APRIORIL2ESTIMATESFORGRADENTFLUX},
we find that
$
\mathbb{T}(\ubar,u)
\leq 
C \initialsmall^2
+
C
\int_{u' = \moreinterestingu_1}^u 
	\mathbb{T}(\ubar,u')
\, \rmd u' + \int_{\ubar' = \leftubar}^{\ubar} 
	\mathbb{T}(\ubar',u)
\, \rmd u'
$.
Applying Gr\"{o}nwall's inequality, we conclude that $\mathbb{T}(\ubar,u)
\leq 
C \initialsmall^2$, which yields \eqref{E:APRIORIL2BELOWTOPESTIMATESFORVORTVORTFLUX}--\eqref{E:APRIORIL2BELOWTOPESTIMATESFORDIVGRADENTFLUX}.

\hfill $\qed$

\subsection{Proof of
\eqref{E:MAINTORIL2VORTGRADENTBELOWTOPORDERBLOWUP}--\eqref{E:MAINTORIL2MODIFIEDFLUIDVARIABLESBELOWTOPORDERBLOWUP}}
\label{SS:PROOFOFBELOWTOPORDERROUGHTORIENERGYESTIMATES}
We fix any integer $N$ with $0 \leq N \leq \Ntop-1$.
Using \eqref{E:IDENTITYUBARDERIVATIVEOFDOUBLENULLORUSINTEGRAL}, with $f = \left| \tander^N (\vortrenormalized,\GradEnt) \right|^2  \weight^{\blowuprateoftransportWRTnumberofcommutations(N)}$, we have that:
\begin{align} 
	\begin{split} \label{E:FIRSTSTEPPROOFOFBELOWTOPORDERESTIMATES}
	\int_{\doublenulltoritwoarg{\ubar}{u}} \left| \tander^N (\vortrenormalized,\GradEnt) \right|^2 \weight^{\blowuprateoftransportWRTnumberofcommutations(N+1)} \, \voldoublenulltori 
	& = 	\int_{\doublenulltoritwoarg{\leftubar }{u}} \left| \tander^N (\vortrenormalized,\GradEnt) \right|^2 \weight^{\blowuprateoftransportWRTnumberofcommutations(N+1)} \, \voldoublenulltori + \blowuprateoftransportWRTnumberofcommutations(N) 	\int_{\outgoingcharacteristicsurfacetwoarg{u}{[\leftubar,\ubar]}}  \left| \tander^N (\vortrenormalized,\GradEnt) \right|^2  (\newL \weight) \weight^{\blowuprateoftransportWRTnumberofcommutations(N)} \, \volingoingnullhypersurface \\
	& \ \ + \int_{\outgoingcharacteristicsurfacetwoarg{u}{[\leftubar,\ubar]}} \newL \left( \left| \tander^N (\vortrenormalized,\GradEnt)\right|^2\right) \weight^{\blowuprateoftransportWRTnumberofcommutations(N+1)} \, \volingoingnullhypersurface +  \int_{\outgoingcharacteristicsurfacetwoarg{u}{[\leftubar,\ubar]}}  \left| \tander^N (\vortrenormalized,\GradEnt) \right|^2   \frac{1}{2} f \mytr_{\gnulltori} \deform{\newL} \, \weight^{\blowuprateoftransportWRTnumberofcommutations(N+1)} \, \volingoingnullhypersurface,
	\end{split}
\end{align}
where we note that the second integral on RHS\,\eqref{E:FIRSTSTEPPROOFOFBELOWTOPORDERESTIMATES} is not present whenever its coefficient $\blowuprateoftransportWRTnumberofcommutations(N) = 0$. Note that $\newL \weight = -1$ and so, even when $\blowuprateoftransportWRTnumberofcommutations(N) >0$, this integral is $\le 0$. From this, using the pointwise bounds \eqref{E:POINTWISEBOUNDFORDOUBLENULLOROIDALTRACEOFDEFORMATIONTENSOROFNEWLVECTORFIELD}, the definitions  \eqref{E:NULLTORIVORTICITYL2CONTROLLINGQUANTITY} and \eqref{E:NULLTORIGRADENTL2CONTROLLINGQUANTITY}, and the already proven estimates
	\eqref{E:APRIORIL2ESTIMATESFORVORTICITYFLUX}--\eqref{E:APRIORIL2ESTIMATESFORGRADENTFLUX}
	for
	$
	\fluxcontrolVort_{\le N+1}(\ubar,u)
	$
	and
	$
	\fluxcontrolGradEnt_{\le N+1}(\ubar,u)
	$,
	the desired estimate \eqref{E:MAINTORIL2VORTGRADENTBELOWTOPORDERBLOWUP} follows from:
\begin{align}
	\begin{split} \label{E:SECONDSTEPPROOFOFBELOWTOPORDERESTIMATES}
	\int_{\doublenulltoritwoarg{\ubar}{u}} \left| \tander^N (\vortrenormalized,\GradEnt) \right|^2 \weight^{\blowuprateoftransportWRTnumberofcommutations(N+1)} \, \voldoublenulltori 
	\lesssim 
	\initialsmall^2 +	\fluxcontrolVort_{\le N+1}(\ubar,u)
	+
	\fluxcontrolGradEnt_{\le N+1}(\ubar,u) \lesssim \initialsmall^2
	\end{split}
	\end{align}
	From \eqref{E:FIRSTSTEPPROOFOFBELOWTOPORDERESTIMATES}
	and the already proven estimates
	\eqref{E:APRIORIL2ESTIMATESFORVORTICITYFLUX}--\eqref{E:APRIORIL2ESTIMATESFORGRADENTFLUX}
	for
	$
	\fluxcontrolVort_{\le N+1}(\ubar,u)
	$
	and
	$
	\fluxcontrolGradEnt_{\le N+1}(\ubar,u)
	$,
	we conclude,
	in view of definitions~\eqref{E:NULLTORIVORTICITYL2CONTROLLINGQUANTITY} and \eqref{E:NULLTORIGRADENTL2CONTROLLINGQUANTITY},
	the desired estimate \eqref{E:MAINTORIL2VORTGRADENTBELOWTOPORDERBLOWUP}.
	
	The estimates
	\eqref{E:MAINTORIL2MODIFIEDFLUIDVARIABLESBELOWTOPORDERBLOWUP}
	can be proved via similar arguments based on the data-estimates
	\eqref{E:SMALLDATAOFVORTVORTANDDIVGRADENTROPYONDOUBLENULLTORIONPU0}, the already proven estimates
	\eqref{E:APRIORIL2BELOWTOPESTIMATESFORVORTVORTFLUX}--\eqref{E:APRIORIL2BELOWTOPESTIMATESFORDIVGRADENTFLUX},
	and definitions~\eqref{E:NULLTORIVORTVORTL2CONTROLLINGQUANTITY},\eqref{E:NULLTORIDIVGRADENTL2CONTROLLINGQUANTITY}.
	
\hfill $\qed$

\section{Top-order elliptic-hyperbolic $L^2$ estimates for the specific vorticity and entropy gradient}
\label{S:TOPORDERELIPTICHYPERBOLICL2ESTIMATESFORSPECIFICVORTICITYANDENTROPYGRADIENT}
We continue to work under the assumptions of Sect.\,\ref{SS:SILENTFACTS}.
Our main goal in this section is to derive the top-order $L^2$ estimates of Prop.\,\ref{P:TOPORDERMODIFIEDAPRIORIL2ESTIMATESANDDOUBLENULLTORIESTIMATES} for the modified fluid variables $\VortVort$ and $\DivGradEnt$. The proposition also provides the top order energy estimates for $\vortrenormalized$ and $\GradEnt$ along the double-null tori, which are actually coupled to the energy estimates of $\VortVort$ and $\DivGradEnt$. 

In Sect.\,\ref{S:WAVEANDACOUSTICGEOMETRYAPRIORIESTIMATES}, we will use the estimates for $\VortVort$ and $\DivGradEnt$
that we derive in this section in our proof of the wave a priori estimates, 
which we stated as Prop.\,\ref{P:APRIORIL2ESTIMATESWAVEVARIABLES}.
Hence, we highlight that for the logic of the paper, it is important that 
\textbf{the estimates we derive in this section do not rely on the wave estimates of Prop.\,\ref{P:APRIORIL2ESTIMATESWAVEVARIABLES}}. The proofs of this section will instead rely on the bootstrap assumptions
\eqref{E:MAINWAVEENERGYBOOTSTRAP}
for the wave energies, which are \emph{weaker} than the estimates that we derive in Prop.\,\ref{P:APRIORIL2ESTIMATESWAVEVARIABLES}.

To explain the main challenges in the analysis,
we recall that the transport equations \eqref{E:EVOLUTIONEQUATIONFLATCURLRENORMALIZEDVORTICITY} 
and
\eqref{E:TRANSPORTFLATDIVGRADENT}
satisfied by
$\VortVort$ and $\DivGradEnt$
feature some difficult source terms,
denoted by $\mainnullform_{(\VortVort)}^i$ and $\mainnullform_{(\DivGradEnt)}$,
that depend on the general first-order derivatives
of $\vortrenormalized$ and $\GradEnt$. 
These source terms have the potential to cause the loss of a derivative
at the top-order because they cannot be bounded using pure transport estimates. 
In the below-top-order estimates of 
Sect.\,\ref{S:BELOWTOPORDERHYPERBOLICL2ESTIMATESFORSPECIFICVORTICITYANDENTROPYGRADIENT} where $N \le \Ntop -1$, we allowed
the loss of a derivative, as is signified by the terms
$
\int_{u' = \moreinterestingu_1}^u 
			\left\lbrace
				\fluxcontrolVort_{\leq N+1}(\ubar,u')
				+
				\fluxcontrolGradEnt_{\leq N+1}(\ubar,u')
			\right\rbrace
		\, \rmd u'	
$
on RHS\,\eqref{E:BELOWTOPORDERENERGYINTEGRALINEQUALITIESFORMODIFIEDFLUID}.
To avoid the loss at the top-order, we handle the difficult source terms in a different way,
one that is based on combining the elliptic-hyperbolic integral identity provided by 
Prop.\,\ref{P:INTEGRALIDENTITYFORELLIPTICHYPERBOLICCURRENT}
with pointwise estimates that take into account the special structure of the equations 
of Theorem~\ref{T:GEOMETRICWAVETRANSPORTSYSTEM}, sharp estimates for the double-null geometry, and the below-top-order estimates that we already derived in
Sect.\,\ref{S:BELOWTOPORDERHYPERBOLICL2ESTIMATESFORSPECIFICVORTICITYANDENTROPYGRADIENT}.

We organize this section this as follows: 
\begin{itemize}
\item
	In Sect.\,\ref{SS:PRELIMINARYENERGYINTEGRALIDENTITIESFORTOPDERIVATIVESOFMODIFIEDFLUID},
	we derive some preliminary energy integral inequalities for $\VortVort$ and $\DivGradEnt$, 
	which feature the difficult source terms; this part of the proof is not more difficult than the
	proof of the below-top-order energy integral inequalities we derived in 
	Lemma\,\ref{L:BELOWTOPORDERENERGYINTEGRALINEQUALITIESSPECIFICVORTICITYANDENTROPYGRADIENT}.
\item In Sect.\,\ref{SS:MAINENERGYINTEGRALINEQUALITIESFORTOPORDERMODIFIEDFLUIDCONDITIONALONSOURCETERMPROPOSITION},
	we use the preliminary energy integral inequalities to derive the main energy integral inequalities for the top-order
	derivatives of $\VortVort$ and $\DivGradEnt$. These main energy integral inequalities
	are conditional on having $L^2$ estimates for the difficult source terms, which we prove independently as
	Prop.\,\ref{P:ELLIPTICHYPERBOLICINTEGRALINEQUALITIES}
	in Sect.\,\ref{SS:MAINELLIPTICHYPERBOLICINTEGRALINEQUALITIES}.
\item In Sect.\,\ref{SS:PROOFOFMAINVORTVORTDIVGRADENTTOPORDERBLOWUP},
	we use the main energy integral inequalities
	to prove the top-order $L^2$ estimate \eqref{E:TOPORDERVORTVORTL2APRIORIESTIMATES}--\eqref{E:TOPORDERDIVGRADENTL2APRIORIESTIMATES}.
 \item In Sect.\,\ref{SS:ELLIPTICHYPERBOLICIDENTITYL2ESTIMATESFORERRORINTEGRALS}, 
	in service of the proof of Prop.\,\ref{P:ELLIPTICHYPERBOLICINTEGRALINEQUALITIES},
	we derive estimates for the error integrals appearing in
	the elliptic-hyperbolic integral identity
	\eqref{E:INTEGRALIDENTITYFORELLIPTICHYPERBOLICCURRENT}.
\item In Sect.\,\ref{SS:MAINELLIPTICHYPERBOLICINTEGRALINEQUALITIES},
	we prove Prop.\,\ref{P:ELLIPTICHYPERBOLICINTEGRALINEQUALITIES}.
\item Finally, in Sect.\,\ref{SS:PROOFOFMAINTORIVORTGRADENTTOPORDERBLOWUP}, 
	we prove the top-order estimate \eqref{E:MAINTORIL2VORTGRADENTTOPORDER}
	for $\vortrenormalized$ and $\GradEnt$ along the double-null tori
	as well as the top-order estimate \eqref{E:MAINTORIL2MODIFIEDFLUIDVARIABLESTOPORDER}
	for $\VortVort$ and $\DivGradEnt$ along the double-null tori.
\end{itemize}

\subsection{Preliminary energy integral inequalities for the top-order derivatives of $\VortVort$ and $\DivGradEnt$}
\label{SS:PRELIMINARYENERGYINTEGRALIDENTITIESFORTOPDERIVATIVESOFMODIFIEDFLUID}

\begin{lemma}[Preliminary energy integral inequalities for the top-order derivatives of $\VortVort$ and $\DivGradEnt$]
\label{L:PRELIMINARYINTEGRALINEQUALITIESFORTOPDERIVATIVESOFMODIFIEDFLUID}
	Let $\varsigma_0$ be the constant from Lemmas\,\ref{L:PRELIMINARYBELOWTOPORDERL2ESTIMATESFOREIKONALFUINCTIONQUANTITIES} and\,\ref{L:BELOWTOPORDERENERGYINTEGRALINEQUALITIESSPECIFICVORTICITYANDENTROPYGRADIENT} (see also Remark\,\ref{R:SMALLNESSCONSTANTFORBELOWTOPORDERACOUSTICVARIABLES}). Then for any $\varsigma \in (0,1]$, 
	the following integral inequalities hold for
	$(\ubar,u) \in [\leftubar,\ubarboot) \times [\moreinterestingu_1,\moreinterestingu_2]$,
	where the implicit constants are \textbf{independent} of $\varsigma$ (whereas the dependence on $\varsigma_0$ in the inequality is made explicit):
	\begin{align} 
	\begin{split} \label{E:PRELIMINARYINTEGRALINEQUALITIESFORTOPDERIVATIVESOFMODIFIEDFLUIDSTATEMENT}
		& \fluxcontrolVortVort_{\Ntop}(\ubar,u)
		 +
		\fluxcontrolDivGradEnt_{\Ntop}(\ubar,u) 
		+ \blowupratetoporderacoustic \left(\bulkcontrolVortVort_{\Ntop}(\ubar,u)
		+
		\bulkcontrolDivGradEnt_{\Ntop}(\ubar,u)\right)
		\\
		& \ \ \lesssim
		\initialsmall^2
		+
		\varsigma \max_{\tander^{\Ntop}\in \mathfrak{P}^{(\Ntop)}} 
		\int_{\characteristicdiamondtwoarg{[\leftubar,\ubar)}{[\moreinterestingu_1,u]}}
			\left\lbrace \ellipticCoerciveQuadratic[\pmb{\partial} \tander^{\Ntop} \vortrenormalized,\pmb{\partial}\tander^{\Ntop}  \vortrenormalized] 
			+  \ellipticCoerciveQuadratic[\pmb{\partial} \tander^{\Ntop} \GradEnt,\pmb{\partial} \tander^{\Ntop} \GradEnt] 
			\right\rbrace \weight^{\blowupratetoporderacoustic}\ReciprocalLunitAppliedtoTimeFunction
		\, \voldiamond
			\\
	& \ \ \ \ 
		+
		\left(1 + \frac{1}{\varsigma} \right)
		\int_{u' = \moreinterestingu_1}^u
			\left\lbrace
				\fluxcontrolVortVort_{\Ntop}(\ubar,u')
				+
				\fluxcontrolDivGradEnt_{\Ntop}(\ubar,u')
			\right\rbrace
		\, \rmd u'
			\\
	& \ \ \ \ 
		+
		\int_{u' = \moreinterestingu_1}^u
			\left\lbrace \fluxcontrolVortVort_{\leq \Ntop-1}(\ubar,u') + \fluxcontrolDivGradEnt_{\leq \Ntop-1}(\ubar,u')\right\rbrace
		\, \rmd u'
			\\
	& \ \ \ \ 
		+
		\int_{u' = \moreinterestingu_1}^u 
			\left\lbrace
				\fluxcontrolVort_{\leq \Ntop}(\ubar,u')
				+
				\fluxcontrolGradEnt_{\leq \Ntop}(\ubar,u')
			\right\rbrace
		\, \rmd u'	
			\\
	& \ \ \ \ 
		+
		(1+\varsigma_0^{-1})\fundbootsmall^2 
		\int_{\ubar' = \leftubar}^{\ubar} 
			\totalfluxcontrolvelocity_{[1,\Ntop]}(\ubar',u) 
		\, \rmd \ubar'
		+
		(1+\varsigma_0^{-1})\fundbootsmall^2 
		\int_{u' = \moreinterestingu_1}^u 
			\totalfluxcontrolvelocity_{[1,\Ntop]}(\ubar,u')
		\, \rmd u'
		\\
	& \ \ \ \ +
		(1+\varsigma_0^{-1})\fundbootsmall^2 \strongangularcontrolvelocity_{[1,\Ntop]}(\ubar,u).
	\end{split}
	\end{align}
	
\end{lemma}

	\eqref{E:COMMUTEDTRANSPORTPOINTWISEESTIMATESFORMODIFIEDCURLOFVORT}--\eqref{E:COMMUTEDTRANSPORTPOINTWISEESTIMATESFORMODIFIEDDIVERGENCEOFENTROPYGRADIENT}

\begin{proof}
	We consider the transport identity \eqref{E:ENERGYNULLFLUXINTEGRALIDENTITIESTRANSPORT} with $\tander^{\Ntop}(\VortVort,\DivGradEnt)$ in the place of $f$. Using the data estimates of Lemma\,\ref{L:ALLL2CONTROLLINGQUANTITIESINITIALLYSMALL}, we  bound the data-dependent terms 
	\begin{align}
		 \ingoingfluxtransport{\blowupratetoporderacoustic}[\tander^{\Ntop}(\VortVort,\DivGradEnt)](\leftubar,u) + \outgoingfluxtransport{\blowupratetoporderacoustic}[\tander^{\Ntop}(\VortVort,\DivGradEnt)](\ubar,\moreinterestingu_1) \label{E:PRELIMINARYINTEGRALINEQUALITIESFORTOPDERIVATIVESOFMODIFIEDFLUIDSTEPFORDATA}
	\end{align}
	on RHS\,\eqref{E:ENERGYNULLFLUXINTEGRALIDENTITIESTRANSPORT} by $\lesssim \initialsmall^2$. By the bootstrap assumptions and Prop.\,\ref{P:IMPROVEMENTOFAUXILIARYBOOTSTRAP}, we have  the bound
	$|\Lunit \upmu + \upmu \mytr_{\gtorus} \angk| \lesssim 1$ for the integrand factor in the last
	integral on RHS\,\eqref{E:ENERGYNULLFLUXINTEGRALIDENTITIESTRANSPORT}. From this, inserting the pointwise estimates 	\eqref{E:COMMUTEDTRANSPORTPOINTWISEESTIMATESFORMODIFIEDCURLOFVORT}--\eqref{E:COMMUTEDTRANSPORTPOINTWISEESTIMATESFORMODIFIEDDIVERGENCEOFENTROPYGRADIENT} for $\upmu \Transport(\tander^{\Ntop}(\VortVort,\DivGradEnt))$, we have:
	\begin{align} 
		\begin{split} \label{E:PRELIMINARYINTEGRALINEQUALITIESFORTOPDERIVATIVESOFMODIFIEDFLUID}
			& \fluxcontrolVortVort_{\Ntop}(\ubar,u)
			 +
			\fluxcontrolDivGradEnt_{\Ntop}(\ubar,u) + \blowupratetoporderacoustic \left(\bulkcontrolVortVort_{\Ntop}(\ubar,u)
			+
			\bulkcontrolDivGradEnt_{\Ntop}(\ubar,u)\right)
			\\
			& \ \ \lesssim
			\initialsmall^2 + \int_{u' = \moreinterestingu_1}^u \left\{  \fluxcontrolVortVort_{\Ntop}(\ubar,u')
			 +
			\fluxcontrolDivGradEnt_{\Ntop}(\ubar,u') \right\} \rmd u' \\
			& \ \ \ \ + \sup_{\substack{\tander^{\Ntop} \in \mathfrak{P}^{(\Ntop)} \\ \tander^{\Ntop+1} \in \mathfrak{P}^{(\Ntop+1)}}} \int_{\characteristicdiamondtwoarg{[\leftubar,\ubar)}{[\moreinterestingu_1,u]}} \left|\tander^{\Ntop}(\VortVort,\DivGradEnt)\right| \Big\{ \left| \tander^{\le \Ntop}(\VortVort,\DivGradEnt)\right| + \left| \tander^{\Ntop+1} (\vortrenormalized,\GradEnt)\right| +  \left| \tander^{\le \Ntop} (\vortrenormalized,\GradEnt)\right|   \\ 
			& \ \ \ \  + \fundbootsmall \left|\newuL \tander^{[1,\Ntop]} \velocityarray \right|
			+
			\fundbootsmall \left| \tander^{[1,\Ntop + 1]} \velocityarray \right|
			+
			\fundbootsmall \left|\tandersmall^{[1,\Ntop]} \badcontrolvars \right|\Big\} \weight^{\blowupratetoporderacoustic} \ReciprocalLunitAppliedtoTimeFunction\, \voldiamond		
		\end{split}
	\end{align}Using Lemma\,\ref{L:COERCIVENESSOFL2CONTROLLINGQUANITIESFORVORTICITYANDENTROPY} and $\weight \lesssim 1$, we bound:
 \begin{align} 
		\begin{split} \label{E:PRELIMINARYINTEGRALINEQUALITIESFORTOPDERIVATIVESOFMODIFIEDFLUIDSTEP2}
			& \int_{\characteristicdiamondtwoarg{[\leftubar,\ubar)}{[\moreinterestingu_1,u]}} \left|\tander^{\Ntop}(\VortVort,\DivGradEnt)\right| \Big\{ \left| \tander^{\le \Ntop}(\VortVort,\DivGradEnt)\right| + \left| \tander^{\Ntop+1} (\vortrenormalized,\GradEnt)\right| +  \left| \tander^{\le \Ntop} (\vortrenormalized,\GradEnt)\right| \Big\} \weight^{\blowupratetoporderacoustic}\ReciprocalLunitAppliedtoTimeFunction\,  \voldiamond \\
			 & \ \ \lesssim \initialsmall^2 +  \varsigma  \int_{\characteristicdiamondtwoarg{[\leftubar,\ubar)}{[\moreinterestingu_1,u]}} \left|\tander^{\Ntop+1}(\vortrenormalized,\GradEnt) \right|^2 \weight^{\blowupratetoporderacoustic} \ReciprocalLunitAppliedtoTimeFunction \, \voldiamond 
			 +  (1+\varsigma^{-1})  \int_{u' = \moreinterestingu_1}^u \left\{  \fluxcontrolVortVort_{\le \Ntop}(\ubar,u')
			 +
			\fluxcontrolDivGradEnt_{\le \Ntop}(\ubar,u') \right\} \rmd u' \\
			& \ \ \ \  + \int_{u' = \moreinterestingu_1}^u \left\{  \fluxcontrolVort_{\le \Ntop}(\ubar,u')
			 +
			\fluxcontrolGradEnt_{\le \Ntop}(\ubar,u') \right\} \rmd u'.
		\end{split}
	\end{align} 
	
 Using Lemmas\,\ref{L:COERCIVENESSOFL2CONTROLLINGQUANITIESUNIFIEDCOMMUTATOR} and\,\ref{L:PRELIMINARYBELOWTOPORDERL2ESTIMATESFOREIKONALFUINCTIONQUANTITIES} as well as $\weight \lesssim 1$, we similarly bound: 
 	\begin{align}
		\begin{split} \label{E:PRELIMINARYINTEGRALINEQUALITIESFORTOPDERIVATIVESOFMODIFIEDFLUIDSTEP3}
			& \int_{\characteristicdiamondtwoarg{[\leftubar,\ubar)}{[\moreinterestingu_1,u]}} \left|\tander^{\Ntop}(\VortVort,\DivGradEnt)\right| \Big\{
			 \fundbootsmall \left|\newuL \tander^{[1,\Ntop]} \velocityarray \right|
			+
			\fundbootsmall \left| \tander^{[1,\Ntop + 1]} \velocityarray \right|
			+
			\fundbootsmall \left|\tandersmall^{[1,\Ntop]} \badcontrolvars \right|\Big\} \weight^{\blowupratetoporderacoustic} \, \voldiamond \\
			& \ \ \lesssim    \int_{u' = \moreinterestingu_1}^u \left\{  \fluxcontrolVortVort_{\le \Ntop}(\ubar,u')
			 +
			\fluxcontrolDivGradEnt_{\le \Ntop}(\ubar,u') \right\} \rmd u' + \   (1 + \varsigma_0^{-1})\fundbootsmall^2 \int_{\ubar' = \leftubar}^{\ubar} \totalfluxcontrolvelocity_{[1,\Ntop]}(\ubar',u) \, \mathrm{d} \ubar'   \\
			& \ \  \ \ +  (1 + \varsigma_0^{-1}) \fundbootsmall^2 \int_{u' = \moreinterestingu_1}^u \totalfluxcontrolvelocity_{[1,\Ntop]}(\ubar,u') \, \mathrm{d} u' + (1 + \varsigma_0^{-1}) \fundbootsmall^2 \strongangularcontrolvelocity_{[1,N+1]}(\ubar,u) \\
			& \ \  \ \ + \int_{u' = \moreinterestingu_1}^u \left\{\fluxcontrolVort_{\le N+1}(\ubar,u')+ \fluxcontrolGradEnt_{\le N+1}(\ubar,u')\right\} \, \mathrm{d} u'.
		\end{split}
	\end{align}
	
All error integrals in \eqref{E:PRELIMINARYINTEGRALINEQUALITIESFORTOPDERIVATIVESOFMODIFIEDFLUID}--\eqref{E:PRELIMINARYINTEGRALINEQUALITIESFORTOPDERIVATIVESOFMODIFIEDFLUIDSTEP3} except the one featuring $\tander^{\Ntop+1}(\vortrenormalized,\DivGradEnt)$ are present on RHS\,\eqref{E:PRELIMINARYINTEGRALINEQUALITIESFORTOPDERIVATIVESOFMODIFIEDFLUID}

 To conclude the prove, note that Lemma\,\ref{L:RELATIONSHIPBETWEENCARTESIANPARTIALDERIVATIVESANDSMOOTHGEOMETRICCOMMUTATORS}, Prop.\,\ref{P:SCHEMATICSTRUCTUREOFVARIOUSTENSORSINTERMSOFCONTROLVARS}, the bootstrap assumptions, and coercivity estimate \eqref{E:COERCIVENESSOFELLIPTICHYPERBOLICQUADRATICFORM} imply: 
	\begin{align} 
		\begin{split} \label{E:PRELIMINARYINTEGRALINEQUALITIESFORTOPDERIVATIVESOFMODIFIEDFLUIDSTEP4}
			\left|\tander^{\Ntop + 1} (\vortrenormalized,\GradEnt)\right| & \lesssim \max_{\tander^{\Ntop} \in \mathfrak{P}^{(\Ntop)}} \sum_{\alpha=0}^3 \left| \p_\alpha \tander^{\Ntop} \vortrenormalized\right| + \left| \p_\alpha \tander^{\Ntop} \GradEnt \right| \\
			& \lesssim  \max_{\tander^{\Ntop} \in \mathfrak{P}^{(\Ntop)}} 	\ellipticCoerciveQuadratic[\pmb{\partial} \tander^{\Ntop}\vortrenormalized,\pmb{\partial}\tander^{\Ntop} \vortrenormalized] 
			 + 	\ellipticCoerciveQuadratic[\pmb{\partial} \tander^{\Ntop} \GradEnt,\pmb{\partial}\tander^{\Ntop} \GradEnt]
		\end{split}
	\end{align}\end{proof}

\subsection{The main energy integral inequalities for the top-order derivatives of $\VortVort$ and $\DivGradEnt$, conditional on
Prop.\,\ref{P:ELLIPTICHYPERBOLICINTEGRALINEQUALITIES}}
\label{SS:MAINENERGYINTEGRALINEQUALITIESFORTOPORDERMODIFIEDFLUIDCONDITIONALONSOURCETERMPROPOSITION}
Most of our effort in Sect.\,\ref{S:TOPORDERELIPTICHYPERBOLICL2ESTIMATESFORSPECIFICVORTICITYANDENTROPYGRADIENT} 
is dedicated towards bounding the
spacetime integrals:
$$
\varsigma 
\int_{\characteristicdiamondtwoarg{[\leftubar,\ubarboot)}{[\moreinterestingu_1,\moreinterestingu_2]}}
			\left\lbrace
				\ellipticCoerciveQuadratic[\pmb{\partial} \tander^{\Ntop} \vortrenormalized,\pmb{\partial} \tander^{\Ntop} \vortrenormalized] 
				 + 	\ellipticCoerciveQuadratic[\pmb{\partial} \tander^{\Ntop} \GradEnt,\pmb{\partial} \tander^{\Ntop} \GradEnt]
			\right\rbrace \weight^{\blowupratetoporderacoustic}\ReciprocalLunitAppliedtoTimeFunction
		\, \voldiamond,
	$$which appear on RHS\,\eqref{E:PRELIMINARYINTEGRALINEQUALITIESFORTOPDERIVATIVESOFMODIFIEDFLUIDSTATEMENT}.
	We derive the needed estimates in Prop.\,\ref{P:ELLIPTICHYPERBOLICINTEGRALINEQUALITIES}.
	Given Prop.\,\ref{P:ELLIPTICHYPERBOLICINTEGRALINEQUALITIES} 
	and Lemma\,\ref{L:PRELIMINARYINTEGRALINEQUALITIESFORTOPDERIVATIVESOFMODIFIEDFLUID}, 
	it is easy to derive energy integral inequalities
	that can be used to obtain the desired top-order energy estimates
	for $\VortVort$ and $\DivGradEnt$; we derive these integral inequalities in the next lemma.
	
\begin{lemma}[The main energy integral inequalities for the top-order derivatives of $\VortVort$ and $\DivGradEnt$]
\label{L:MAININTEGRALINEQUALITIESFORTOPDERIVATIVESOFMODIFIEDFLUID}
	Assuming the results of Prop.\,\ref{P:ELLIPTICHYPERBOLICINTEGRALINEQUALITIES}, for any $\varsigma \in (0,1]$, 
	the following integral inequalities hold for
	$(\ubar,u) \in [\leftubar,\ubarboot) \times [\moreinterestingu_1,\moreinterestingu_2]$,
	where the implicit constants are \textbf{independent} of $\varsigma$:
	\begin{align}
	\begin{split}  \label{E:MAININTEGRALINEQUALITIESFORTOPDERIVATIVESOFMODIFIEDFLUID}
		& \fluxcontrolVortVort_{\Ntop}(\ubar,u)
		+
		\fluxcontrolDivGradEnt_{\Ntop}(\ubar,u)  + \blowupratetoporderacoustic \left(\bulkcontrolVortVort_{\Ntop}(\ubar,u)
			+
			\bulkcontrolDivGradEnt_{\Ntop}(\ubar,u)\right) \\
		& \ \ \lesssim
		+  \initialsmall^2 +
		\varsigma 
		\left\lbrace
			\fluxcontrolVortVort_{\Ntop}(\ubar,u)
			+
			\fluxcontrolDivGradEnt_{\Ntop}(\ubar,u)
		\right\rbrace \\
		& \ \ \ \ +  \fluxcontrolVortVort_{\le \Ntop-1}(\ubar,u) + \fluxcontrolDivGradEnt_{\le \Ntop-1}(\ubar,u) +	\fluxcontrolVort_{\le \Ntop}(\ubar,u) + \fluxcontrolGradEnt_{\le \Ntop}(\ubar,u) \\
		& \ \ \ \ + \fundbootsmall^2
			\totalfluxcontrolvelocity_{[1,\Ntop]}(\ubar,u)  + (1 + \varsigma_0^{-1}) \int_{u' = \moreinterestingu_1}^u \left\{ \fluxcontrolVortVort_{\le \Ntop} (\ubar,u') +   \fluxcontrolDivGradEnt_{\le \Ntop }(\ubar,u')\right\}\,  \rmd u'   \\
		& \ \ \ \ 
			+ (1 + \varsigma_0^{-1}) \int_{u' = \moreinterestingu_1}^u \left\{ \fluxcontrolVort_{\le \Ntop}(\ubar,u') + \fluxcontrolGradEnt_{\le N} (\ubar,u') \right\}\rmd u' \\
		& \ \ \ \  +
		\left(
			1 
			+ \varsigma_0^{-1}
		\right)\fundbootsmall^2 
		\int_{\ubar' = \leftubar}^{\ubar} 
			\frac{1}{|\ubar'|^{0.5}}
			\totalfluxcontrolvelocity_{[1,\Ntop]}(\ubar',u) 
		\, \rmd \ubar' \\
		& \ \ \ \   +
		\left(
			1 
			+ \varsigma_0^{-1}
		\right)\fundbootsmall^2 
		\int_{u' = \moreinterestingu_1}^{u} 
			\totalfluxcontrolvelocity_{[1,\Ntop]}(\ubar,u') 
		\, \rmd u' \\
		& \ \ \ \ + (1 + \varsigma_0^{-1}) \fundbootsmall^2\strongangularcontrolvelocity_{[1,\Ntop]}(\ubar,u),
		\end{split}	
	\end{align}	
\end{lemma}

\begin{proof}
	We start by considering inequality \eqref{E:PRELIMINARYINTEGRALINEQUALITIESFORTOPDERIVATIVESOFMODIFIEDFLUIDSTATEMENT}. Thanks to the estimates
	\eqref{E:ELLIPTICHYPERBOLICINTEGRALINEQUALITYFORSPECIFICVORTICITY}--\eqref{E:ELLIPTICHYPERBOLICINTEGRALINEQUALITYFORENTROPYGRADIENT},
	which we prove independently in Sect.\,\ref{SS:MAINELLIPTICHYPERBOLICINTEGRALINEQUALITIES},
	the spacetime integral
	\[
\varsigma
\int_{\characteristicdiamondtwoarg{[\leftubar,\ubarboot)}{[\moreinterestingu_1,\moreinterestingu_2]}}
\left\lbrace \ellipticCoerciveQuadratic[\pmb{\partial} \tander^{\Ntop} \vortrenormalized,\pmb{\partial}\tander^{\Ntop}  \vortrenormalized] 
			+  \ellipticCoerciveQuadratic[\pmb{\partial} \tander^{\Ntop} \GradEnt,\pmb{\partial} \tander^{\Ntop} \GradEnt] 
			\right\rbrace \weight^{\blowupratetoporderacoustic}\ReciprocalLunitAppliedtoTimeFunction 		\, \voldiamond
	\]
	on the first line of RHS\,\eqref{E:PRELIMINARYINTEGRALINEQUALITIESFORTOPDERIVATIVESOFMODIFIEDFLUIDSTATEMENT}
	is bounded by
$
\varsigma
	\left\lbrace
		\mbox{RHS\,\eqref{E:ELLIPTICHYPERBOLICINTEGRALINEQUALITYFORSPECIFICVORTICITY}}
		+
		\mbox{RHS\,\eqref{E:ELLIPTICHYPERBOLICINTEGRALINEQUALITYFORENTROPYGRADIENT}}
	\right\rbrace
$.	
\end{proof}

\subsection{Proof of the main top-order energy estimates \eqref{E:TOPORDERVORTVORTL2APRIORIESTIMATES}--\eqref{E:TOPORDERDIVGRADENTL2APRIORIESTIMATES}
 of Prop.\,\ref{P:TOPORDERMODIFIEDAPRIORIL2ESTIMATESANDDOUBLENULLTORIESTIMATES}}
\label{SS:PROOFOFMAINVORTVORTDIVGRADENTTOPORDERBLOWUP}
Given Lemma\,\ref{L:MAININTEGRALINEQUALITIESFORTOPDERIVATIVESOFMODIFIEDFLUID},
we are now ready to prove our main top-order a priori energy estimates \eqref{E:TOPORDERVORTVORTL2APRIORIESTIMATES}--\eqref{E:TOPORDERDIVGRADENTL2APRIORIESTIMATES}
for the modified fluid variables.
We again emphasize that Lemma\,\ref{L:MAININTEGRALINEQUALITIESFORTOPDERIVATIVESOFMODIFIEDFLUID}
is conditional on the estimates of Prop.\,\ref{P:ELLIPTICHYPERBOLICINTEGRALINEQUALITIES},
which we prove independently below.

To proceed, we choose and fix $\varsigma > 0$ sufficiently small such that the term
$
\varsigma 
		\left\lbrace
			\fluxcontrolVortVort_{\Ntop}(\ubar,u)
			+
			\fluxcontrolDivGradEnt_{\Ntop}(\ubar,u)
		\right\rbrace
$
on RHS\,\eqref{E:MAININTEGRALINEQUALITIESFORTOPDERIVATIVESOFMODIFIEDFLUID}
can be absorbed back into the LHS at the expense of increasing the implicit
constants. Next, we 
use the already proven estimates
\eqref{E:APRIORIL2ESTIMATESFORVORTICITYFLUX}--\eqref{E:APRIORIL2BELOWTOPESTIMATESFORDIVGRADENTFLUX}, the bootstrap assumptions \eqref{E:MAINWAVEENERGYBOOTSTRAP},
and \eqref{E:NONLINEARINEQUALITYRELATINGDATAEPSILONANDBOOTSTRAPEPSILON}, and the integrability of $\int_{\ubar' = \leftubar}^u \frac{1}{|\ubar'|^{0.5}} \, \rmd \ubar ' < \infty$ 
to bound all terms on RHS\,\eqref{E:MAININTEGRALINEQUALITIESFORTOPDERIVATIVESOFMODIFIEDFLUID} by $C \initialsmall^2$
except for the integral
$
\int_{u' = \moreinterestingu_1}^u
			\left\lbrace
				\fluxcontrolVortVort_{\Ntop}(\ubar,u')
				+
				\fluxcontrolDivGradEnt_{\Ntop}(\ubar,u')
			\right\rbrace
		\, \rmd u'
$.
In total, this leads to the following inequality, where $C$ depends on the fixed value of $\varsigma$ and $\varsigma_0$:
\begin{align} \label{E:GRONWALLREADYMAINTOPORDERENERGYESTIMATESMODIFIEDFLUIDVARIABLESBLOWUP}
\fluxcontrolVortVort_{\Ntop}(\ubar,u)
+
\fluxcontrolDivGradEnt_{\Ntop}(\ubar,u) +   \blowupratetoporderacoustic \left(\bulkcontrolVortVort_{\Ntop}(\ubar,u)
			+
			\bulkcontrolDivGradEnt_{\Ntop}(\ubar,u)\right) 
& 
	\leq 
	C
	\initialsmall^2
		+
	C 
		\int_{u' = \moreinterestingu_1}^u
			\left\lbrace
				\fluxcontrolVortVort_{\Ntop}(\ubar,u')
				+
				\fluxcontrolDivGradEnt_{\Ntop}(\ubar,u')
			\right\rbrace
		\, \rmd u'.
\end{align}
From \eqref{E:GRONWALLREADYMAINTOPORDERENERGYESTIMATESMODIFIEDFLUIDVARIABLESBLOWUP} and
Gr\"{o}nwall's inequality,
we conclude
that 
$
\fluxcontrolVortVort_{\Ntop}(\ubar,u)
+
\fluxcontrolDivGradEnt_{\Ntop}(\ubar,u) +  \blowupratetoporderacoustic \left(\bulkcontrolVortVort_{\Ntop}(\ubar,u)
			+
			\bulkcontrolDivGradEnt_{\Ntop}(\ubar,u)\right) 
\leq 
	C
	\initialsmall^2 
$,
which yields the desired bounds \eqref{E:TOPORDERVORTVORTL2APRIORIESTIMATES}--\eqref{E:TOPORDERDIVGRADENTL2APRIORIESTIMATES}.

\hfill $\qed$

\subsection{$L^2$ estimates for the error terms}
\label{SS:ELLIPTICHYPERBOLICIDENTITYL2ESTIMATESFORERRORINTEGRALS}
In the next lemma, we derive $L^2$ estimates for the error integrals
in the elliptic-hyperbolic identities provided by Prop.\,\ref{P:INTEGRALIDENTITYFORELLIPTICHYPERBOLICCURRENT}.

\begin{lemma}[$L^2$ estimates for the error terms in the elliptic-hyperbolic identities]
	\label{L:L2ESTIMTAESFORELLIPTICHYPERBOLICIDENTITYERORTERMS}
	Let $\varsigma \in [0,1)$,
	let $(\ubar,u) \in [\leftubar,\ubarboot) \times [\moreinterestingu_1,\moreinterestingu_2]$,
	let
	$
	\EllipticHyperbolicCurrentIntegralIdentityTotalSpacetimeErrorTerm[\cdot,\cdot]
	$
	be the error term defined by \eqref{E:ELLIPTICHYPERBOLICINTEGRALIDENTITYBULKERRORTERM},
	and let $\ellipticCoerciveQuadratic[\cdot,\cdot]$ be the coercive quadratic form defined in 
	\eqref{E:NULLHYPERSURFACEADAPTEDCOERCIVEQUADRATICFORM}.
	Then the following spacetime integral estimates hold,
	where the implicit constants are \textbf{independent} of $\varsigma$:
	\begin{subequations}
	\begin{align} 
	\begin{split} \label{E:SPACETIMEERRORTERMESTIMATEFORSPECIFICVORTICITYELLIPTICHYPERBOLICIDENTITY}
		\left|
		\int_{\characteristicdiamondtwoarg{[\leftubar,\ubar)}{[\moreinterestingu_1,u]}}
			\EllipticHyperbolicCurrentIntegralIdentityTotalSpacetimeErrorTerm[ \tander^{\Ntop} \vortrenormalized, \pmb{\partial} 			
			\tander^{\Ntop} \vortrenormalized] \weight^{\blowupratetoporderacoustic}\ReciprocalLunitAppliedtoTimeFunction 
		\, \voldiamond
		\right|
		& \lesssim 
			\varsigma 
			\int_{\characteristicdiamondtwoarg{[\leftubar,\ubar)}{[\moreinterestingu_1,\moreinterestingu_2]}}
				\ellipticCoerciveQuadratic[\pmb{\partial}\tander^{\Ntop} \vortrenormalized,\pmb{\partial} \tander^{\Ntop} \vortrenormalized] 
				 \weight^{\blowupratetoporderacoustic}\ReciprocalLunitAppliedtoTimeFunction 
			\, \voldiamond
				\\
	& \ \ + (1 + \varsigma^{-1}) \initialsmall^2 + (1 + \varsigma^{-1}) \int_{u' = \moreinterestingu_1}^u \fluxcontrolVortVort_{\le \Ntop }(\ubar,u') \rmd u' \\
	& \ \  +(1 + \varsigma^{-1}) \int_{u' = \moreinterestingu_1}^u \left\{ \fluxcontrolVort_{\le \Ntop}(\ubar,u') + \fluxcontrolGradEnt_{\le N} (\ubar,u') \right\}\rmd u' \\
	& \ \ +
		\left(
			1 
			+ \varsigma^{-1} + \varsigma_0^{-1}
		\right)\fundbootsmall^2 
		\int_{\ubar' = \leftubar}^{\ubar} 
			\frac{1}{|\ubar'|^{0.5}}
			\totalfluxcontrolvelocity_{[1,\Ntop]}(\ubar',u) 
		\, \rmd \ubar' \\
		& \ \  +
		\left(
			1 
			+ \varsigma^{-1} + \varsigma_0^{-1} +  \varsigma^{-1}\varsigma_0^{-1}
		\right)\fundbootsmall^2 
		\int_{u' = \moreinterestingu_1}^{u} 
			\totalfluxcontrolvelocity_{[1,\Ntop]}(\ubar,u') 
		\, \rmd u' \\
		& \ \ + (1 + \varsigma^{-1} + \varsigma_0^{-1} + \varsigma^{-1} \varsigma_0^{-1}) \fundbootsmall^2\strongangularcontrolvelocity_{[1,\Ntop]}(\ubar,u),
	\end{split}
			\\
	\begin{split}  \label{E:SPACETIMEERRORTERMESTIMATEFORENTROPYGRADIENTELLIPTICHYPERBOLICIDENTITY}
		\left|
		\int_{\characteristicdiamondtwoarg{[\leftubar,\ubar)}{[\moreinterestingu_1,u]}}
			\EllipticHyperbolicCurrentIntegralIdentityTotalSpacetimeErrorTerm[ \tander^{\Ntop} \GradEnt, \pmb{\partial} \tander^{\Ntop} \GradEnt] \weight^{ \blowupratetoporderacoustic}\ReciprocalLunitAppliedtoTimeFunction 
		\, \voldiamond
		\right|
		& \lesssim  
			\varsigma
			\int_{\characteristicdiamondtwoarg{[\leftubar,\ubar)}{[\moreinterestingu_1,\moreinterestingu_2]}}
				\ellipticCoerciveQuadratic[\pmb{\partial}\tander^{\Ntop} \GradEnt,\pmb{\partial} \tander^{\Ntop} \GradEnt]
				 \weight^{\blowupratetoporderacoustic}\ReciprocalLunitAppliedtoTimeFunction 
			\, \voldiamond
				\\
	& \ \ +  (1 + \varsigma^{-1}) \initialsmall^2 + (1 + \varsigma^{-1}) \int_{u' = \moreinterestingu_1}^u \fluxcontrolDivGradEnt_{\le \Ntop }(\ubar,u') \rmd u' \\
	& \ \  +(1 + \varsigma^{-1}) \int_{u' = \moreinterestingu_1}^u \left\{ \fluxcontrolVort_{\le \Ntop}(\ubar,u') + \fluxcontrolGradEnt_{\le N} (\ubar,u') \right\}\rmd u' \\
	& \ \ +
		\left(
			1 
			+ \varsigma^{-1} + \varsigma_0^{-1}
		\right)\fundbootsmall^2 
		\int_{\ubar' = \leftubar}^{\ubar} 
			\frac{1}{|\ubar'|^{0.5}}
			\totalfluxcontrolvelocity_{[1,\Ntop]}(\ubar',u) 
		\, \rmd \ubar' \\
		& \ \  +
		\left(
			1 
			+ \varsigma^{-1} + \varsigma_0^{-1} +  \varsigma^{-1}\varsigma_0^{-1}
		\right)\fundbootsmall^2 
		\int_{u' = \moreinterestingu_1}^{u} 
			\totalfluxcontrolvelocity_{[1,\Ntop]}(\ubar,u') 
		\, \rmd u' \\
		& \ \ + (1 + \varsigma^{-1} + \varsigma_0^{-1}  +  \varsigma^{-1} \varsigma_0^{-1}) \fundbootsmall^2\strongangularcontrolvelocity_{[1,\Ntop]}(\ubar,u),	\end{split}
	\end{align}
	\end{subequations}
	
	Moreover, the error terms
	$
	^{(\SigmatTan)}\Currentboundaryerrorhavetocontrolprincipalweighted{\blowupratetoporderacoustic}, \,
	^{(\SigmatTan)}\Currentboundaryerrorhavetocontrollowerorderweighted{\blowupratetoporderacoustic}
	$
	defined by \eqref{E:PRINCIPALERRORTERMHAVETOCONTROLKEYIDPUTANGENTCURRENTCONTRACTEDAGAINSTVECTORFIELDFORVORTICITYWEIGHTED}--\eqref{E:LOWERORDERERRORTERMHAVETOCONTROLKEYIDPUTANGENTCURRENTCONTRACTEDAGAINSTVECTORFIELDFORVORTICITYWEIGHTED} with $\SigmatTan = \tander^{\Ntop} \vortrenormalized$ or $\tander^{\Ntop}\GradEnt$ 
	verify the following ingoing characteristic hypersurface integral estimates:
	\begin{subequations}
	\begin{align} 
	\begin{split} \label{E:PRINCIPALNULLHYPERRORTERMESTIMATEFORSPECIFICVORTICITYELLIPTICHYPERBOLICIDENTITY}
		 \left|
			\int_{\ingoingcharacteristicsurfacetwoarg{\ubar}{[\moreinterestingu_1,\moreinterestingu_2]} } \,
				^{(\tander^{\Ntop} \vortrenormalized)}\Currentboundaryerrorhavetocontrolprincipalweighted{\blowupratetoporderacoustic} \, \volingoingnullhypersurface
		\right| 
		&  
		 \lesssim \initialsmall^2 + \fluxcontrolVortVort_{\Ntop}(\ubar,u) +	\fluxcontrolVort_{\le \Ntop}(\ubar,u) \\
		& \ \ + \fundbootsmall^2
			\totalfluxcontrolvelocity_{[1,\Ntop]}(\ubar,u)
			+ (1 + \varsigma_0^{-1}) \fundbootsmall^2 \int_{\ubar' = \leftubar}^{\ubar} \totalfluxcontrolvelocity_{[1,\Ntop]}(\ubar',u) \, \mathrm{d} \ubar'   \\
			& \ \ +  (1 + \varsigma_0^{-1}) \fundbootsmall^2 \int_{u' = \moreinterestingu_1}^u \totalfluxcontrolvelocity_{[1,\Ntop]}(\ubar,u') \, \mathrm{d} u' + (1 + \varsigma_0^{-1})\fundbootsmall^2 \strongangularcontrolvelocity_{[1,\Ntop]}(\ubar,u) \\
			& \ \ +  \int_{u' = \moreinterestingu_1}^u \left\{\fluxcontrolVort_{\le \Ntop}(\ubar,u')+ \fluxcontrolGradEnt_{\le \Ntop }(\ubar,u')\right\} \, \mathrm{d} u',
		\end{split}
					 \\
		\begin{split} \label{E:PRINCIPALNULLHYPERRORTERMESTIMATEFORENTROPYGRADIENTELLIPTICHYPERBOLICIDENTITY} 
		 \left|
			\int_{\ingoingcharacteristicsurfacetwoarg{\ubar}{[\moreinterestingu_1,\moreinterestingu_2]} }
				 \,
				^{(\tander^{\Ntop} \GradEnt)}\Currentboundaryerrorhavetocontrolprincipalweighted{\blowupratetoporderacoustic} 
			\, \volingoingnullhypersurface
		\right| 
		& 
		  \lesssim 	
			\fluxcontrolDivGradEnt_{\Ntop}(\ubar,u) + \fluxcontrolGradEnt_{\Ntop}(\ubar,u) 
			+ 
			\fundbootsmall^2 
			\totalfluxcontrolvelocity_{\Ntop}(\ubar,u),
	\end{split}
	\end{align}
	\end{subequations}	\begin{subequations}
	\begin{align} \label{E:LOWERORDERNULLHYPERRORTERMESTIMATEFORSPECIFICVORTICITYELLIPTICHYPERBOLICIDENTITY}
		\begin{split}
			\left|
				\int_{\ingoingcharacteristicsurfacetwoarg{\ubar}{[\moreinterestingu_1,\moreinterestingu_2]} }
					\, ^{(\tander^N \vortrenormalized)}\Currentboundaryerrorhavetocontrollowerorderweighted{\blowupratetoporderacoustic} \, \volingoingnullhypersurface
			\right|
			& \lesssim \initialsmall^2 + \fluxcontrolVortVort_{\le \Ntop-1} (\ubar,u) +  \fluxcontrolDivGradEnt_{\le \Ntop-1} (\ubar,u) + \fluxcontrolVort_{\le \Ntop} (\ubar,u) +  \fluxcontrolGradEnt_{\le \Ntop} (\ubar,u) \\
			& \ \ + \fundbootsmall^2
				\totalfluxcontrolvelocity_{[1,\Ntop]}(\ubar,u)
			+ (1 + \varsigma_0^{-1}) \fundbootsmall^2 \int_{\ubar' = \leftubar}^{\ubar} \totalfluxcontrolvelocity_{[1,\Ntop-1]}(\ubar',u) \, \mathrm{d} \ubar'   \\
			& \ \ +  (1 + \varsigma_0^{-1}) \fundbootsmall^2 \int_{u' = \moreinterestingu_1}^u \totalfluxcontrolvelocity_{[1,\Ntop-1]}(\ubar,u') \, \mathrm{d} u' + (1 + \varsigma_0^{-1})\fundbootsmall^2 \strongangularcontrolvelocity_{[1,\Ntop-1]}(\ubar,u) \\
			& \ \ +  \int_{u' = \moreinterestingu_1}^u \left\{\fluxcontrolVort_{\le \Ntop-1}(\ubar,u')+ \fluxcontrolGradEnt_{\le \Ntop -1}(\ubar,u')\right\} \, \mathrm{d} u',
		\end{split} \\
		\begin{split}
			\left|
			\int_{\ingoingcharacteristicsurfacetwoarg{\ubar}{[\moreinterestingu_1,\moreinterestingu_2]} }
				\, ^{(\tander^N \GradEnt)}\Currentboundaryerrorhavetocontrollowerorderweighted{\blowupratetoporderacoustic} 
			\, \volingoingnullhypersurface
			\right|
			& \lesssim \initialsmall^2 + \fluxcontrolVortVort_{\le \Ntop-1} (\ubar,u) +  \fluxcontrolDivGradEnt_{\le \Ntop-1} (\ubar,u) + \fluxcontrolVort_{\le \Ntop} (\ubar,u) +  \fluxcontrolGradEnt_{\le \Ntop} (\ubar,u) \\
			& \ \ + \fundbootsmall^2
				\totalfluxcontrolvelocity_{[1,\Ntop]}(\ubar,u)
			+ (1 + \varsigma_0^{-1}) \fundbootsmall^2 \int_{\ubar' = \leftubar}^{\ubar} \totalfluxcontrolvelocity_{[1,\Ntop-1]}(\ubar',u) \, \mathrm{d} \ubar'   \\
			& \ \ +  (1 + \varsigma_0^{-1}) \fundbootsmall^2 \int_{u' = \moreinterestingu_1}^u \totalfluxcontrolvelocity_{[1,\Ntop-1]}(\ubar,u') \, \mathrm{d} u' + (1 + \varsigma_0^{-1})\fundbootsmall^2 \strongangularcontrolvelocity_{[1,\Ntop-1]}(\ubar,u) \\
			& \ \ +  \int_{u' = \moreinterestingu_1}^u \left\{\fluxcontrolVort_{\le \Ntop-1}(\ubar,u')+ \fluxcontrolGradEnt_{\le \Ntop -1}(\ubar,u')\right\} \, \mathrm{d} u'.
		\end{split}			\label{E:LOWERORDERNULLHYPERRORTERMESTIMATEFORENTROPYGRADIENTELLIPTICHYPERBOLICIDENTITY}
	\end{align}
	\end{subequations}
	
	Finally, the error terms
	$\coercivequadraticformintoriboundaryterms[\tander^{\Ntop} \vortrenormalized, \tander^{\Ntop} \vortrenormalized],
	\cdots,
	\coercivequadraticformintoriboundaryterms[\tander^{\Ntop} \GradEnt, \tander^{\Ntop} \GradEnt]
	$
	defined by \eqref{E:TORIBOUNDARYTERMSFORELLIPTICHYPERBOLICIDENTITIES}
	verify the following double-null tori integral estimates:
	\begin{subequations}
	\begin{align} 
	 \label{E:TOPORDERSMALLNESSOFVORTICITYANDENTROPYGRADIENTDATAALONGINITIALUDOUBLENULLTORI}
		\int_{\doublenulltoritwoarg{\ubar}{\moreinterestingu_1}} \frac{1}{\upmu}
			\coercivequadraticformintoriboundaryterms[\tander^{\Ntop} \vortrenormalized, \tander^{\Ntop} \vortrenormalized] \weight^{\blowupratetoporderacoustic}
		\, \voldoublenulltori,
			\,
		\int_{\doublenulltoritwoarg{\ubar}{\moreinterestingu_1}}
			 \frac{1}{\upmu} \coercivequadraticformintoriboundaryterms[\tander^{\Ntop} \GradEnt, \tander^{\Ntop} \GradEnt]
		\, \voldoublenulltori
		& \lesssim
			\initialsmall^2,
				\\
		\int_{\doublenulltoritwoarg{\leftubar}{u}}
			 \frac{1}{\upmu} \coercivequadraticformintoriboundaryterms[\tander^{\Ntop} \vortrenormalized, \tander^{\Ntop} \vortrenormalized] \weight^{\blowupratetoporderacoustic}
		\, \voldoublenulltori,
			\,
		\int_{\doublenulltoritwoarg{\leftubar}{u}}
			 \frac{1}{\upmu} \coercivequadraticformintoriboundaryterms[\tander^{\Ntop} \GradEnt, \tander^{\Ntop} \GradEnt] \weight^{\blowupratetoporderacoustic}
		\, \voldoublenulltori
		& \lesssim \initialsmall^2.
			\label{E:TOPORDERSMALLNESSOFVORTICITYANDENTROPYGRADIENTDATAALONGINITIALUBARDOUBLENULLTORI}
\end{align}
\end{subequations}	
	
\end{lemma}

\begin{proof}

We first prove \eqref{E:SPACETIMEERRORTERMESTIMATEFORSPECIFICVORTICITYELLIPTICHYPERBOLICIDENTITY}. From \eqref{E:ANTISYMMETRICVORTICITYELLIPTICHYPERBOLICIDENTITYERORTERMPOINTWISE}, \eqref{E:DIVVORTICITYELLIPTICHYPERBOLICIDENTITYERORTERMPOINTWISE}, \eqref{E:DERIVATIVEOFONEOVERMUWEIGHTVORTICITYELLIPTICHYPERBOLICIDENTITYERORTERMPOINTWISE}, \eqref{E:FIRSTABSORBTERMVORTICITYELLIPTICHYPERBOLICIDENTITYERORTERMPOINTWISE}, \eqref{E:SECONDABSORBTERMVORTICITYELLIPTICHYPERBOLICIDENTITYERORTERMPOINTWISE}, \eqref{E:QUADRATICINMATERIALDERIVATIVETERMVORTICITYELLIPTICHYPERBOLICIDENTITYERORTERMPOINTWISE}, \eqref{E:NULLGEOMETRYVORTICITYELLIPTICHYPERBOLICIDENTITYERORTERMPOINTWISE}, it suffices to bound: 
\begin{align} 
		\begin{split} \label{E:SPACETIMEERRORTERMESTIMATEFORSPECIFICVORTICITYELLIPTICHYPERBOLICIDENTITYSTEP1}
			 \int_{\characteristicdiamondtwoarg{[\leftubar,\ubar)}{[\moreinterestingu_1,u]}}	& 
			\left\{	\varsigma 
			\ellipticCoerciveQuadratic[\pmb{\partial}\tander^{\Ntop} \vortrenormalized,\pmb{\partial} \tander^{\Ntop} \vortrenormalized]
				+
			\left(1 + \frac{1}{\varsigma} \right)
			\left|
				\tander^{\Ntop} \VortVort 
			\right|^2
			+
			\left(1 + \frac{1}{\varsigma} \right)
			\frac{1}{\upmu^2}
			\left|
				\tander^{\le \Ntop -1} \VortVort 
			\right|^2 \right. \\
			& \ \ 
			+ 
			\left(1 + \frac{1}{\varsigma} \right)
			\frac{1}{\upmu^2}
			|\tander^{\le \Ntop} (\vortrenormalized,\GradEnt)|^2  
			+
			\left(1 + \frac{1}{\varsigma} \right)
			\frac{\fundbootsmall^2}{\upmu^2}
			|\newuL \tander^{[1,\Ntop]} \velocityarray|^2 \\
			& \ \ \left. 
			+
			\left(1 + \frac{1}{\varsigma} \right)
			\fundbootsmall^2  | \tander^{[1,\Ntop+1]} \velocityarray|^2
			+
			\left(1 + \frac{1}{\varsigma} \right)
			\frac{\fundbootsmall^2}{\upmu^2}
			|\tandersmall^{[1,\Ntop]} \badcontrolvars|^2 \right\} \weight^{\blowupratetoporderacoustic}\ReciprocalLunitAppliedtoTimeFunction 
			\, \voldiamond
		\end{split}
	\end{align}
The first integral in \eqref{E:SPACETIMEERRORTERMESTIMATEFORSPECIFICVORTICITYELLIPTICHYPERBOLICIDENTITYSTEP1} is present on RHS\,\eqref{E:SPACETIMEERRORTERMESTIMATEFORSPECIFICVORTICITYELLIPTICHYPERBOLICIDENTITY}. The second error integral is  bounded by $(1 + \varsigma^{-1}) \int_{u' = \moreinterestingu_1}^u \fluxcontrolVortVort_{N}(\ubar,u') \rmd u'$. Since $\weight^{\blowupratetoporderacoustic-2} \lesssim  \weight^{\blowuprateofacousticgeoWRTnumberofcommutations(\Ntop-1)} = \weight^{\blowuprateoftransportWRTnumberofcommutations(\Ntop)}$, estimate $\frac{1}{\upmu} \lesssim \frac{1}{\weight}$ (see \eqref{E:KEYESTIMATECONTROLLINGINVERSEMUBYINVERSEWEIGHT}) implies that the third and fourth error integrals in \eqref{E:SPACETIMEERRORTERMESTIMATEFORSPECIFICVORTICITYELLIPTICHYPERBOLICIDENTITYSTEP1} are $\lesssim$ RHS\,\eqref{E:SPACETIMEERRORTERMESTIMATEFORSPECIFICVORTICITYELLIPTICHYPERBOLICIDENTITY}. Similarly, since $\weight^{\blowupratetoporderacoustic-2} = \weight^{\blowuprateofwaveWRTnumberofcommutations(\Ntop)-0.5}$, we bound the fifth integral above by $(1+\varsigma^{-1})\fundbootsmall^2 \int_{\ubar' = \leftubar}^{\ubar} \frac{1}{|\ubar'|^{0.5}}\totalfluxcontrolvelocity_{[1,\Ntop]} (\ubar',u) \, \rmd \ubar'$.  The sixth integral above is bounded by $(1+\varsigma^{-1})\fundbootsmall^2 \int_{u' = \moreinterestingu_1}^u \totalfluxcontrolvelocity_{[1,\Ntop]} (\ubar,u') \, \rmd u'$ when $\tander^{\Ntop+1}\velocityarray = \Lunit\tander^{\Ntop} \velocityarray$ or by $ (1 + \varsigma^{-1})\fundbootsmall^2 \strongangularcontrolvelocity_{[1,\Ntop]}(\ubar,u)$ when $\tander^{\Ntop+1} \velocityarray = \Yvf{A} \tander^{\Ntop} \velocityarray$. To bound the final integral, we use Lemma\,\ref{L:PRELIMINARYBELOWTOPORDERL2ESTIMATESFOREIKONALFUINCTIONQUANTITIES}. The proof of \eqref{E:SPACETIMEERRORTERMESTIMATEFORENTROPYGRADIENTELLIPTICHYPERBOLICIDENTITY} follows from identical arguments, using instead \eqref{E:ANTISYMMETRICENTROPYGRADIENTELLIPTICHYPERBOLICIDENTITYERORTERMPOINTWISE}, \eqref{E:DIVENTROPYGRADIENTELLIPTICHYPERBOLICIDENTITYERORTERMPOINTWISE}, \eqref{E:DERIVATIVEOFONEOVERMUWEIGHTENTROPYGRADIENTELLIPTICHYPERBOLICIDENTITYERORTERMPOINTWISE}, \eqref{E:FIRSTABSORBTERMENTROPYGRADIENTELLIPTICHYPERBOLICIDENTITYERORTERMPOINTWISE}, \eqref{E:SECONDABSORBTERMENTROPYGRADIENTELLIPTICHYPERBOLICIDENTITYERORTERMPOINTWISE}, \eqref{E:QUADRATICINMATERIALDERIVATIVETERMVORTICITYELLIPTICHYPERBOLICIDENTITYERORTERMPOINTWISE}, \eqref{E:NULLGEOMETRYENTROPYGRADIENTELLIPTICHYPERBOLICIDENTITYERORTERMPOINTWISE}.

 Estimates \eqref{E:PRINCIPALNULLHYPERRORTERMESTIMATEFORSPECIFICVORTICITYELLIPTICHYPERBOLICIDENTITY} follows from nearly identica arguments as in the previous paragraph, using instead \eqref{E:PRINCIPALSIGAMTILDEELLIPTICHYPERBOLICVORTICITYPOINTWISE} and multiplying both sides by $\weight^{\blowupratetoporderacoustic}$. We clarify, however, that to estimate $\frac{1}{\upmu} |\tander^{\le \Ntop} (\vortrenormalized,\GradEnt)| \weight^{\blowupratetoporderacoustic}$, we multiply and divide by $\ReciprocaluLunitAppliedtoTimeFunction$, apply the bound \eqref{E:RECIPROCALULUNITAPPLIEDTOTIMEFUNCTIONAPPROXIMATELYMU}, then use \eqref{E:KEYESTIMATECONTROLLINGINVERSEMUBYINVERSEWEIGHT} to bound the integral as: 
 \begin{align}
		 \begin{split} \label{E:PRINCIPALNULLHYPERRORTERMESTIMATEFORSPECIFICVORTICITYELLIPTICHYPERBOLICIDENTITYSTEP1}
 			\int_{\ingoingcharacteristicsurfacetwoarg{\ubar}{[\moreinterestingu_1,\moreinterestingu_2]}} \frac{1}{\upmu} \left| \tander^{\le \Ntop} (\vortrenormalized, \GradEnt)\right|^2  \weight^{\blowupratetoporderacoustic} \, \volingoingnullhypersurface & \lesssim \int_{\ingoingcharacteristicsurfacetwoarg{\ubar}{[\moreinterestingu_1,\moreinterestingu_2]}}\ReciprocaluLunitAppliedtoTimeFunction \left| \tander^{\le \Ntop}  (\vortrenormalized, \GradEnt)\right|^2   \weight^{\blowupratetoporderacoustic-2} \, \volingoingnullhypersurface \\
			& \lesssim \int_{u' = \moreinterestingu_1}^u \left\{\fluxcontrolVort_{\le \Ntop}(\ubar,u')+ \fluxcontrolGradEnt_{\le \Ntop }(\ubar,u')\right\} \, \mathrm{d} u',
		\end{split}
	\end{align}
where we used $\weight^{\blowupratetoporderacoustic-2} \lesssim  \weight^{\blowuprateoftransportWRTnumberofcommutations(\Ntop)}$ and \eqref{E:COERCIVENESSOFCONTROLVORT}. We clarify that we used \eqref{E:PRELIMINARYBELOWTOPORDERL2ESTIMATESFOREIKONALFUINCTIONQUANTITIESWITHZ} to bound the error integrals generated by the last terms $\frac{\fundbootsmall^2}{\upmu} \left| \comdersmall^{[1,N];\le 1} \controlvars \right|^2$. The bound \eqref{E:PRINCIPALNULLHYPERRORTERMESTIMATEFORENTROPYGRADIENTELLIPTICHYPERBOLICIDENTITY} follows from similar arguments based on \eqref{E:PRINCIPALSIGAMTILDEELLIPTICHYPERBOLICENTROPYGRADIENTPOINTWISE}, but they are strictly easier as there are no terms $\frac{\fundbootsmall^2}{\upmu} \left| \comdersmall^{[1,N];\le 1} \controlvars \right|^2$ on RHS\,\eqref{E:PRELIMINARYBELOWTOPORDERL2ESTIMATESFOREIKONALFUINCTIONQUANTITIESWITHZ}.

The bounds \eqref{E:LOWERORDERNULLHYPERRORTERMESTIMATEFORSPECIFICVORTICITYELLIPTICHYPERBOLICIDENTITY}--\eqref{E:LOWERORDERNULLHYPERRORTERMESTIMATEFORENTROPYGRADIENTELLIPTICHYPERBOLICIDENTITY} follow from nearly identical arguments used to prove \eqref{E:PRINCIPALNULLHYPERRORTERMESTIMATEFORSPECIFICVORTICITYELLIPTICHYPERBOLICIDENTITY}--\eqref{E:PRINCIPALNULLHYPERRORTERMESTIMATEFORENTROPYGRADIENTELLIPTICHYPERBOLICIDENTITY} using instead \eqref{E:LOWERORDERSIGAMTILDEELLIPTICHYPERBOLICVORTICITYPOINTWISE}--\eqref{E:LOWERORDERSIGAMTILDEELLIPTICHYPERBOLIENTROPYGRADIENTPOINTWISE}. We omit the details except for the main changes. Firstly, instead of having an extra $\frac{1}{2}$-power of $\weight$ as in \eqref{E:PRINCIPALNULLHYPERRORTERMESTIMATEFORSPECIFICVORTICITYELLIPTICHYPERBOLICIDENTITYSTEP1} (note $\weight^{\blowupratetoporderacoustic - 2} = \weight^{\blowuprateoftransportWRTnumberofcommutations(\Ntop) + \frac{1}{2}} \lesssim \weight^{\blowuprateoftransportWRTnumberofcommutations(\Ntop)}$), due to the singular factor $\upmu^{-3/2}$, we have a sharp bound: $\weight^{\blowupratetoporderacoustic - 2.5} = \weight^{\blowuprateoftransportWRTnumberofcommutations(\Ntop)}$. 
Second, in order to control the $|\tander^{[1,N]} \velocityarray|$ term on RHS\,\eqref{E:LOWERORDERSIGAMTILDEELLIPTICHYPERBOLICVORTICITYPOINTWISE}--\eqref{E:LOWERORDERSIGAMTILDEELLIPTICHYPERBOLIENTROPYGRADIENTPOINTWISE}, which does not feature a $\ingoingcharacteristicsurfacetwoarg{\ubar}{[\moreinterestingu_1,\moreinterestingu_2]}$-tangent differentiation as its last derivative operator, we use \eqref{E:DERIVATIVELOSINGL2ESTIMATEFORTANGENTIALDERIVATIVES}.

Estimates \eqref{E:TOPORDERSMALLNESSOFVORTICITYANDENTROPYGRADIENTDATAALONGINITIALUDOUBLENULLTORI}--\eqref{E:TOPORDERSMALLNESSOFVORTICITYANDENTROPYGRADIENTDATAALONGINITIALUBARDOUBLENULLTORI} follow immediately from the coercive bounds \eqref{E:KEYCOERCIVITYONTORIBOUNDARYTERMSFORELLIPTICHYPERBOLICIDENTITIES} and the initial data estimates \eqref{E:SMALLDATAOFTOPORDERVORTICITYANDGRADENTROPYONDOUBLENULLTORIONUBAR0FORELLIPTICHYPERBOLICIDENTITIES}--\eqref{E:SMALLDATAOFTOPORDERVORTICITYANDGRADENTROPYONDOUBLENULLTORIONU1FORELLIPTICHYPERBOLICIDENTITIES}.

\end{proof}

\subsection{The main elliptic-hyperbolic integral inequalities}
\label{SS:MAINELLIPTICHYPERBOLICINTEGRALINEQUALITIES}
Thanks to the availability of 
Lemma\,\ref{L:L2ESTIMTAESFORELLIPTICHYPERBOLICIDENTITYERORTERMS},
we are now ready to prove Prop.\,\ref{P:ELLIPTICHYPERBOLICINTEGRALINEQUALITIES}.

\begin{proposition}[The main elliptic-hyperbolic integral inequalities]
		\label{P:ELLIPTICHYPERBOLICINTEGRALINEQUALITIES}
	Let $\ellipticCoerciveQuadratic[\pmb{\partial} \SigmatTan,\pmb{\partial} \SigmatTan]$
	be the quadratic form from Def.\,\ref{D:NULLHYPERSURFACEADAPTEDCOERCIVEQUADRATICFORM},
	and let
	$\coercivequadraticformintoriboundaryterms[\SigmatTan,\SigmatTan]$
	be the quadratic form defined by \eqref{E:TORIBOUNDARYTERMSFORELLIPTICHYPERBOLICIDENTITIES}.
	Then the following spacetime integral estimates hold for 
$(\ubar,u) \in [\leftubar,\ubarboot) \times [\moreinterestingu_1,\moreinterestingu_2]$:
	\begin{subequations}
	\begin{align}	
	\begin{split} \label{E:ELLIPTICHYPERBOLICINTEGRALINEQUALITYFORSPECIFICVORTICITY}
		&
		\int_{\characteristicdiamondtwoarg{[\leftubar,\ubar)}{[\moreinterestingu_1,u]}}
			\ellipticCoerciveQuadratic[\pmb{\partial} \tander^{\Ntop} \vortrenormalized, \pmb{\partial} \tander^{\Ntop} \vortrenormalized]\weight^{\blowupratetoporderacoustic}\ReciprocalLunitAppliedtoTimeFunction 
		\, \voldiamond
		+
		\int_{\doublenulltoritwoarg{\ubar}{u}}
		 \frac{1}{\upmu}
			\coercivequadraticformintoriboundaryterms[\tander^{\Ntop} \vortrenormalized, \tander^{\Ntop} \vortrenormalized] \weight^{\blowupratetoporderacoustic}
		\, \voldoublenulltori
			\\
		& \lesssim
			(1 + \varsigma^{-1}) \initialsmall^2 + \fluxcontrolVortVort_{\le \Ntop}(\ubar,u) +  \fluxcontrolDivGradEnt_{\le \Ntop-1}(\ubar,u) +	\fluxcontrolVort_{\le \Ntop}(\ubar,u) + \fluxcontrolGradEnt_{\le \Ntop}(\ubar,u) + \fundbootsmall^2
			\totalfluxcontrolvelocity_{[1,\Ntop]}(\ubar,u) \\
			& \ \ + (1 + \varsigma^{-1}) \int_{u' = \moreinterestingu_1}^u \fluxcontrolVortVort_{\le \Ntop }(\ubar,u') \rmd u'  +(1 + \varsigma^{-1}) \int_{u' = \moreinterestingu_1}^u \left\{ \fluxcontrolVort_{\le \Ntop}(\ubar,u') + \fluxcontrolGradEnt_{\le N} (\ubar,u') \right\}\rmd u' \\
	& \ \ +
		\left(
			1 
			+ \varsigma^{-1} +  \varsigma_0^{-1}
		\right)\fundbootsmall^2 
		\int_{\ubar' = \leftubar}^{\ubar} 
			\frac{1}{|\ubar'|^{0.5}}
			\totalfluxcontrolvelocity_{[1,\Ntop]}(\ubar',u) 
		\, \rmd \ubar' \\
		& \ \  +
		\left(
			1 
			+ \varsigma^{-1} + \varsigma_0^{-1} + \varsigma^{-1}\varsigma_0^{-1}
		\right)\fundbootsmall^2 
		\int_{u' = \moreinterestingu_1}^{u} 
			\totalfluxcontrolvelocity_{[1,\Ntop]}(\ubar,u') 
		\, \rmd u' \\
		& \ \ + (1 + \varsigma^{-1} + \varsigma_0^{-1} + \varsigma^{-1} \varsigma_0^{-1}) \fundbootsmall^2\strongangularcontrolvelocity_{[1,\Ntop]}(\ubar,u),	\end{split}	
			\\
	\begin{split} \label{E:ELLIPTICHYPERBOLICINTEGRALINEQUALITYFORENTROPYGRADIENT} 
	&
		\int_{\characteristicdiamondtwoarg{[\leftubar,\ubar)}{[\moreinterestingu_1,u]}}
			\ellipticCoerciveQuadratic[\pmb{\partial} \tander^{\Ntop} \GradEnt, \pmb{\partial} \tander^{\Ntop} \GradEnt] 
			\weight^{\blowupratetoporderacoustic} \ReciprocalLunitAppliedtoTimeFunction 
		\, \voldiamond
		+
		\int_{\doublenulltoritwoarg{\ubar}{u}}
			 \frac{1}{\upmu}
			\coercivequadraticformintoriboundaryterms[\tander^{\Ntop} \GradEnt, \tander^{\Ntop} \GradEnt] \weight^{\blowupratetoporderacoustic}
		\, \voldoublenulltori
				\\
		& \lesssim
			(1 + \varsigma^{-1}) \initialsmall^2 + \fluxcontrolDivGradEnt_{\le \Ntop}(\ubar,u) +  \fluxcontrolVortVort_{\le \Ntop-1}(\ubar,u) +	\fluxcontrolVort_{\le \Ntop}(\ubar,u) + \fluxcontrolGradEnt_{\le \Ntop}(\ubar,u) + \fundbootsmall^2
			\totalfluxcontrolvelocity_{[1,\Ntop]}(\ubar,u) \\
			& \ \ + (1 + \varsigma^{-1}) \int_{u' = \moreinterestingu_1}^u \fluxcontrolDivGradEnt_{\le \Ntop }(\ubar,u') \rmd u'  
			+ (1 + \varsigma^{-1}) \int_{u' = \moreinterestingu_1}^u \left\{ \fluxcontrolVort_{\le \Ntop}(\ubar,u') + \fluxcontrolGradEnt_{\le N} (\ubar,u') \right\}\rmd u' \\
	& \ \ +
		\left(
			1 
			+ \varsigma^{-1} + \varsigma_0^{-1}
		\right)\fundbootsmall^2 
		\int_{\ubar' = \leftubar}^{\ubar} 
			\frac{1}{|\ubar'|^{0.5}}
			\totalfluxcontrolvelocity_{[1,\Ntop]}(\ubar',u) 
		\, \rmd \ubar' \\
		& \ \  +
		\left(
			1 
			+ \varsigma^{-1} + \varsigma_0^{-1} + \varsigma^{-1}\varsigma_0^{-1}
		\right)\fundbootsmall^2 
		\int_{u' = \moreinterestingu_1}^{u} 
			\totalfluxcontrolvelocity_{[1,\Ntop]}(\ubar,u') 
		\, \rmd u' \\
		& \ \ + (1 + \varsigma^{-1} + \varsigma_0^{-1} +  \varsigma^{-1} \varsigma_0^{-1}) \fundbootsmall^2\strongangularcontrolvelocity_{[1,\Ntop]}(\ubar,u),	\end{split}	
\end{align}
\end{subequations}
\end{proposition}

\begin{proof}
	We consider the integral identity
	\eqref{E:INTEGRALIDENTITYFORELLIPTICHYPERBOLICCURRENT}
	with $\tander^{\Ntop} \vortrenormalized$ and $\tander^{\Ntop} \GradEnt$ in the role of $\SigmatTan$.
	Using Lemma\,\ref{L:L2ESTIMTAESFORELLIPTICHYPERBOLICIDENTITYERORTERMS} with 
	$\ubar_1 \eqdef \leftubar$,
	$\ubar_2 \eqdef \ubar$,
	$u_1 \eqdef \moreinterestingu_1$,
	and
	$u_2 \eqdef u$,
	as well as the data-estimates 
	\eqref{E:WAVEL2CONTROLLINGINITIALLYSMALL}--\eqref{E:FLUXFORGRADENTINITIALLSMALL},
	we bound the integrals on RHS\,\eqref{E:INTEGRALIDENTITYFORELLIPTICHYPERBOLICCURRENT}
	(the data-estimates are used to control the data-hypersurface integrals
	$\int_{\ingoingcharacteristicsurfacetwoarg{\leftubar}{[\moreinterestingu_1,u]}} \cdots$),
	where we can discard the integrals
	$
		-
		\int_{\doublenulltoritwoarg{\leftubar}{\moreinterestingu_1}}
				 \frac{1}{\upmu} \weight^{\blowupratetoporderacoustic} \coercivequadraticformintoriboundaryterms[\SigmatTan,\SigmatTan] 	\, \voldoublenulltori
	$
	because they are non-positive in view of
	\eqref{E:KEYCOERCIVITYONTORIBOUNDARYTERMSFORELLIPTICHYPERBOLICIDENTITIES}.
	Finally, by choosing and fixing $\varsigma > 0$ to be sufficiently small,
we can absorb the first error integral
$
\varsigma
			\int_{\characteristicdiamondtwoarg{[\leftubar,\ubar)}{[\moreinterestingu_1,\moreinterestingu_2]}}
				\ellipticCoerciveQuadratic[\pmb{\partial} \tander^{\Ntop} \vortrenormalized, \pmb{\partial} \tander^{\Ntop} \vortrenormalized] 	\weight^{\blowupratetoporderacoustic} \ReciprocalLunitAppliedtoTimeFunction 
			\, \voldiamond
$
on RHS\,\eqref{E:SPACETIMEERRORTERMESTIMATEFORSPECIFICVORTICITYELLIPTICHYPERBOLICIDENTITY}. We cal also absorb the first error integral $
\varsigma
			\int_{\characteristicdiamondtwoarg{[\leftubar,\ubar)}{[\moreinterestingu_1,\moreinterestingu_2]}}
				\ellipticCoerciveQuadratic[\pmb{\partial} \tander^{\Ntop} \GradEnt, \pmb{\partial} \tander^{\Ntop} \GradEnt] 	\weight^{\blowupratetoporderacoustic} \ReciprocalLunitAppliedtoTimeFunction 
			\, \voldiamond
$
on RHS\,\eqref{E:SPACETIMEERRORTERMESTIMATEFORENTROPYGRADIENTELLIPTICHYPERBOLICIDENTITY}
into LHS~\eqref{E:ELLIPTICHYPERBOLICINTEGRALINEQUALITYFORSPECIFICVORTICITY}
and
LHS~\eqref{E:ELLIPTICHYPERBOLICINTEGRALINEQUALITYFORENTROPYGRADIENT} respectively.
This yields the desired estimates
\eqref{E:ELLIPTICHYPERBOLICINTEGRALINEQUALITYFORSPECIFICVORTICITY}--\eqref{E:ELLIPTICHYPERBOLICINTEGRALINEQUALITYFORENTROPYGRADIENT}.

		
\end{proof}

\subsection{Proof of the top-order double-null tori energy estimates \eqref{E:MAINTORIL2VORTGRADENTTOPORDER} and
\eqref{E:MAINTORIL2MODIFIEDFLUIDVARIABLESTOPORDER}}
\label{SS:PROOFOFMAINTORIVORTGRADENTTOPORDERBLOWUP}
We first prove the estimate \eqref{E:MAINTORIL2VORTGRADENTTOPORDER}.
We consider the estimates
\eqref{E:ELLIPTICHYPERBOLICINTEGRALINEQUALITYFORSPECIFICVORTICITY}--\eqref{E:ELLIPTICHYPERBOLICINTEGRALINEQUALITYFORENTROPYGRADIENT}.
In view of 
 \eqref{E:NULLTORIVORTICITYL2CONTROLLINGQUANTITY}, \eqref{E:NULLTORIGRADENTL2CONTROLLINGQUANTITY}, and 
the quantitative positive definiteness estimate
\eqref{E:KEYCOERCIVITYONTORIBOUNDARYTERMSFORELLIPTICHYPERBOLICIDENTITIES}, we have 
	\begin{align}
		\toricontrolVort_{\Ntop}(\ubar,u)  + \toricontrolGradEnt_{\Ntop}(\ubar,u) \lesssim (\text{RHS\,\eqref{E:ELLIPTICHYPERBOLICINTEGRALINEQUALITYFORSPECIFICVORTICITY} + \eqref{E:ELLIPTICHYPERBOLICINTEGRALINEQUALITYFORENTROPYGRADIENT}}). \label{E:MAINTORIL2VORTGRADENTANDMODIFIEDFLUIDVARIABLESTOPORDERSTEP1}
	\end{align}
 Next, we 
use the already proven estimates
\eqref{E:APRIORIL2ESTIMATESFORVORTICITYFLUX}--\eqref{E:APRIORIL2BELOWTOPESTIMATESFORDIVGRADENTFLUX} and \eqref{E:TOPORDERVORTVORTL2APRIORIESTIMATES}--\eqref{E:TOPORDERDIVGRADENTL2APRIORIESTIMATES}, the bootstrap assumptions \eqref{E:MAINWAVEENERGYBOOTSTRAP},
the nonlinear relationship between the bootstrap parameters \eqref{E:NONLINEARINEQUALITYRELATINGDATAEPSILONANDBOOTSTRAPEPSILON}, and the integrability of $\int_{\ubar' = \leftubar}^u \frac{1}{|\ubar'|^{0.5}} \, \rmd \ubar ' < \infty$ to find that $ (\text{RHS\,\eqref{E:ELLIPTICHYPERBOLICINTEGRALINEQUALITYFORSPECIFICVORTICITY} + \eqref{E:ELLIPTICHYPERBOLICINTEGRALINEQUALITYFORENTROPYGRADIENT}}) \lesssim \initialsmall^2$, proving the estimates.

The estimate \eqref{E:MAINTORIL2MODIFIEDFLUIDVARIABLESTOPORDER} follows from combining the same
arguments we used to prove \eqref{E:FIRSTSTEPPROOFOFBELOWTOPORDERESTIMATES}
with the data-estimates
\eqref{E:SMALLDATAOFTOPORDERVORTICITYANDGRADENTROPYONDOUBLENULLTORIONU1FORELLIPTICHYPERBOLICIDENTITIES}, the already proven estimates \eqref{E:TOPORDERVORTVORTL2APRIORIESTIMATES}--\eqref{E:TOPORDERDIVGRADENTL2APRIORIESTIMATES}, and definitions \eqref{E:NULLTORIVORTVORTL2CONTROLLINGQUANTITY}, \eqref{E:NULLTORIDIVGRADENTL2CONTROLLINGQUANTITY}.

\hfill $\qed$


\section{Proof of the $L^2$ estimates for the wave variables and the characteristic geometry}  
\label{S:WAVEANDACOUSTICGEOMETRYAPRIORIESTIMATES}
We continue to work under the assumptions of Sect.\,\ref{SS:SILENTFACTS}.
In this section, we prove Props.\,\ref{P:APRIORIL2ESTIMATESWAVEVARIABLES} and\,\ref{P:APRIORIL2ESTIMATESACOUSTICGEOMETRY}, 
which provide the main a priori energy estimates for the wave variables and the characteristic geometry
along the double-null foliations.  We recall that the fundamental $L^2$-controlling
quantities, such as $\mastercontroltop$ and $\totalcontrolwave_{[1,N]}$, 
are defined in Sect.\,\ref{SS:FUNDAMENTALL2CONTROLLINGWAVEQUANTITIES}
(see in particular Def.\,\ref{D:MASTERCONTROLLINGQUANTITY}).

\subsection{Statement of the integral inequalities used in proving a priori $L^2$ estimates for the wave variables} 
\label{SS:INTEGRALINEQUALITIIESFORWAVEENERGIES}
In this section, we state Prop.\,\ref{P:MAINWAVEENERGYINTEGRALINEQUALITIES},
which provides an integral inequality
for the top-order master controlling energies $\mastercontroltop$.  In the same proposition we provide analogous below-top-order estimates for the wave controlling quantities.
As we will see in Sect.\,\ref{SSS:PROOFOFTOPORDERAPRIORIL2ESTIMATEMASTERCONTROLANDALLORDERWAVES}, 
these integral inequalities are the main ingredients in our proof of the $L^2$ 
a priori estimates of Props.\,\ref{P:TOPORDERAPRIORIL2ESTIMATESMASTERCONTROL}--\ref{P:APRIORIL2ESTIMATESACOUSTICGEOMETRY}.
Most of our effort in Sect.\,\ref{S:WAVEANDACOUSTICGEOMETRYAPRIORIESTIMATES} is dedicated towards
proving preliminary estimates that we will use in proving Prop.\,\ref{P:MAINWAVEENERGYINTEGRALINEQUALITIES}.

We now state the proposition. Its proof is located in Sect.\,\ref{SS:PROOFOFMAINWAVEENERGYINTEGRALINEQUALITIES}. Throughout all of Sect.\ref{S:WAVEANDACOUSTICGEOMETRYAPRIORIESTIMATES}, we will silently use the estimate \eqref{E:BOUNDSOFLUBARINTERESTINGREGION}, which implies $\ReciprocalLunitAppliedtoTimeFunction = \frac{1}{\Lunit \ubar} \approx 1$.

\begin{proposition}[The system of integral inequalities satisfied by the energy controlling quantities] 
\label{P:MAINWAVEENERGYINTEGRALINEQUALITIES} 
 Let $\varsigma_0$ be the constant featured in Lemma\,\ref{L:PRELIMINARYBELOWTOPORDERL2ESTIMATESFOREIKONALFUINCTIONQUANTITIES} (see also Remark\,\ref{R:SMALLNESSCONSTANTFORBELOWTOPORDERACOUSTICVARIABLES}). If $\varsigma_1 \in (0,1]$ is sufficiently small, then there exist constants
$C > 0$ and $C_* > 0$ 
that are \textbf{independent} of $\varsigma_0,\varsigma_1$ 
such that the following estimates hold
for $(\ubar,u) \in [\leftubar,\ubarboot) \times [\moreinterestingu_1,\moreinterestingu_2]$: 

\medskip

\noindent \underline{\textbf{Top-order integral inequalities for $\mastercontroltop$}.}
We have the following estimates for the $L^2$-controlling quantity $\mastercontroltop(\ubar,u)$
defined in \eqref{E:MASTERCONTROLLINGQUANTITY}:
\begin{align}
		\begin{split}  \label{E:TOPORDERWAVEL2CONTROLLINGINTEGRALINEQUALITY}
			\mastercontroltop(\ubar,u) &\le \boxed{2.81 + 5.18} \weakspacetimetransvelocitycontrol_N(\ubar,u)  + \boxed{2.33 + 1.01} \totalbulkcontrolprecisefullymodifiedchi(\ubar,u) \\
			& \ \ +  \boxed{1.02 + 2.07} \weakspacetimeangularvelocitycontrol_N(\ubar,u) +  \boxed{4.13}  \totalbulkcontrolpartialmodifiedchi(\ubar,u) \\
			& \ \  + \boxed{\frac{2}{3}} \totalfluxcontrolvelocity_N(\ubar,u) + \boxed{\frac{3}{4}(1.01)} \totalfluxcontrolpartialmodifiedchi(\ubar,u) \\
			& \ \ + \Errortoparg{N}(\ubar,u),
		\end{split}
	\end{align}\noindent where we define $\Errortoparg{N}(\ubar,u)$ to be any term that satisfies the following estimate,
in which the implicit constants are \textbf{independent} of $\varsigma_0,\varsigma_1 \in (0,1]$:
	\begin{align}
		\begin{split} \label{E:ERRORTOPORDERWAVEESTIMATES}
			\left|\Errortoparg{N}(\ubar,u)\right| &  \lesssim  (1 + \varsigma_0^{-1}+\varsigma_1^{-1}) \initialsmalldoublenull^2  + (\mr \upalpha + \fundbootsmall + \varsigma_1) \totalfluxcontrolvelocity_N(\ubar,u) +  (1 + \varsigma_0^{-1} + \varsigma_1^{-1}) \int_{\ubar' = \leftubar}^{\ubar} \totalfluxcontrolvelocity_{[1,N]}(\ubar',u) \, \mathrm{d} \ubar' \\ 
			& \ \  + (1+ \varsigma_0^{-1} + \varsigma_1^{-1}) \int_{u' = \moreinterestingu_1}^u  \totalfluxcontrolvelocity_{[1,N]}(\ubar,u') \, \rmd u' + \varsigma_1 (1 + \varsigma_0^{-1})\strongangularcontrolvelocity_{N}(\ubar,u) \\
			& \ \ + \fundbootsmall (1 + \varsigma_0^{-1}) \strongangularcontrolvelocity_{[1,N]}(\ubar,u) +  (1 + \varsigma_1^{-1}) \int_{\ubar' = \leftubar}^{\ubar} \frac{1}{|\ubar'|^{1/2}} \totalfluxcontrolvelocity_{[1,N]}(\ubar',u) \, \mathrm{d} \ubar' \\
			& \ \ +   \int_{u' = \moreinterestingu_1}^u \left\{\fluxcontrolVort_{\le N}(\ubar,u')+ \fluxcontrolGradEnt_{\le N }(\ubar,u')\right\} \, \mathrm{d} u +   \int_{u' = \moreinterestingu_1}^u \left\{\fluxcontrolVortVort_{\le N}(\ubar,u')+ \fluxcontrolDivGradEnt_{\le N }(\ubar,u')\right\} \, \mathrm{d} u' \\	
			& \ \ + (1+\varsigma_1^{-1}) \fluxcontrolVort(\ubar,u) + (1+\varsigma_1^{-1}) \bulkcontrolVort(\ubar,u) + \bulkcontrolVortVort_N(\ubar,u) +  \bulkcontrolDivGradEnt_N(\ubar,u)  \\
			& \ \  + (1 +\varsigma_1^{-1})  \int_{\ubar' = \leftubar}^{\ubar} \frac{1}{|\ubar'|^{0.5}}\totalfluxcontrolprecisefullymodifiedchi(\ubar',u) \, \rmd \ubar'   +  (1 + \varsigma_0^{-1} +  \varsigma_1^{-1}) \int_{\ubar' = \leftubar}^{\ubar}  \totalfluxcontrolimprecisefullymodifiedchi(\ubar',u) \, \mathrm{d} \ubar' \\
			& \ \ +  \fundbootsmall \totalbulkcontrolimprecisefullymodifiedchi(\ubar,u) +   \fundbootsmall \totalbulkcontrolprecisefullymodifiedchi(\ubar,u) +(\mr \upalpha + \fundbootsmall) \totalbulkcontrolpartialmodifiedchi(\ubar,u) +  (\mr \upalpha + \fundbootsmall  + \varsigma_1)\weakspacetimeangularvelocitycontrol_N(\ubar,u) \\
			& \ \ + (1 + \varsigma_1^{-1}) \int_{\ubar' = \leftubar}^{\ubar} \totalfluxcontrolpartialmodifiedchi(\ubar',u) \, \rmd \ubar'.
		\end{split}
	\end{align}

\medskip

\noindent \underline{\textbf{Below-top-order integral inequalities for $\velocityarray$}.}
Finally, if $2 \leq N \leq \Ntop$, 
then we have the following estimates for the $L^2$-controlling quantity $\totalcontrolwave_{[1,N-1]}$
defined by \eqref{E:VELOCITYFLUXANDSTRONGANGULARL2CONTROLLINGQUANTITY} and Def.\,\ref{D:SUMMEDL2CONTROLLINGQUANTITIES}:

	\begin{align}
		\begin{split} \label{E:MAINWAVEBELOWTOPINTEGRALINEQUALITIES}
			\totalcontrolwave_{[1,N-1]}(\ubar,u) & \le \varsigma_1(1 + \varsigma_0^{-1}) \int_{u' = \moreinterestingu_1}^u \totalfluxcontrolvelocity_{N}(\ubar,u') \, \rmd u'  +  \varsigma_1(1 + \varsigma_0^{-1}) \strongangularcontrolvelocity_N(\ubar,u) \\
			& \ \ + \Errorsubcriticalarg{N-1}(\ubar,u),
		\end{split}
	\end{align}
where, for any integer $M \ge 1$, 
$\Errorsubcriticalarg{N-1}(\ubar,u)$ is defined to be any term that satisfies the following estimate,
in which the implicit constants are \textbf{independent} of $\varsigma_0,\varsigma_1 \in (0,1]$:
\begin{align}
	\begin{split} \label{E:ERRORBELOWTOPORDERWAVEESTIMATES}
		\left|
			\Errorsubcriticalarg{M}
		\right|(\ubar,u) &  \lesssim \initialsmalldoublenull^2  + \varsigma_1^{-1} \int_{\ubar' = \leftubar}^{\ubar} \frac{1}{|\ubar'|^{1/2}}\totalfluxcontrolvelocity_{[1,M]}(\ubar',u) \, \rmd \ubar'  \\
		& \ \ + \varsigma_1^{-1} \int_{u' = \moreinterestingu_1}^u \totalfluxcontrolvelocity_{[1,M]}(\ubar,u') \, \rmd u' + \varsigma_1(1 + \varsigma_0^{-1}) \strongangularcontrolvelocity_{[1,M]}(\ubar,u) \\
		& \ \ + \fundbootsmall  \weakspacetimetransvelocitycontrol_{M}(\ubar,u)  + \fundbootsmall \weakspacetimeangularvelocitycontrol_M (\ubar,u) \\
		& \ \ + \int_{\ubar' = \leftubar}^{\ubar} \fluxcontrolVortVort_{\le M} (\ubar',u) \, \rmd \ubar' + \bulkcontrolVortVort_{M}(\ubar,u).
	\end{split}
\end{align}\end{proposition}

\subsection{Estimates for the easiest error integrals}
\label{SS:ESTIMATESFOREASIESTERRORINTEGRALS}

\subsubsection{Estimates for the error integrals generated by the error terms $\HarmlessWave{N}$}
\label{SSS:ESTIMATESFOREERRORINTEGRALSGENERATEDBYHARMLESSWAVETERMS}
In the following lemma, we derive bounds for all the wave equation error integrals 
that involve the $\HarmlessWave{N}$ terms defined in \eqref{E:HARMLESSWAVE}. 

\begin{lemma}[Bounds for error integrals involving $\HarmlessWave{N}$ terms] 
\label{L:HARMLESSWAVETERMSERRORINTEGRALBOUNDS}
Let $1 \leq N \leq \Ntop$ and let $v \in \velocityarray = \{v^1,v^2,v^3\}$.
Let $\tander^N \in \mathfrak{P}^{(N)}$, where $\mathfrak{P}^{(N)}$
is the set of order $N$ $\nullhyparg{u}$-tangential commutator operators from
Def.\,\ref{D:STRINGSOFCOMMUTATIONVECTORFIELDS}. Let $\blowuprateofwaveWRTnumberofcommutations(N)$ denote the blowup rate of the wave variables defined in \eqref{E:BLOWUPRATEOFWAVEWRTNUMBEROFCOMMUTATORS}. Let $\varsigma_0$ be the constant featured in Lemma\,\ref{L:PRELIMINARYBELOWTOPORDERL2ESTIMATESFOREIKONALFUINCTIONQUANTITIES} (see also Remark\,\ref{R:SMALLNESSCONSTANTFORBELOWTOPORDERACOUSTICVARIABLES}).
Recall that terms of type $\HarmlessWave{[1,N]}$ are defined in Def.\,\ref{D:HARMLESSWAVE}.
Then if $\varsigma_1\in (0,1]$ is sufficiently small, the following estimates hold for 
$(\ubar,u) \in [\leftubar,\ubarboot) \times [\moreinterestingu_1,\moreinterestingu_2]$, 
where the implicit constants are \textbf{independent} of $\varsigma_0, \, \varsigma_1$: 
\begin{align}
	\begin{split}  \label{E:ENERGYESTIMATEHARMLESSWAVEERRORTERMS}
		 & \int_{\characteristicdiamondtwoarg{[\leftubar,\ubar]}{[\moreinterestingu_1,u]}} 	
			\left| \begin{pmatrix} 
				(1 + 2 \upmu) \Lunit \tander^N v 
				\\ 
				2 \muX \tander^N v 
			\end{pmatrix} 
			\right|  
			\left| \HarmlessWave{[1,N]} \right| \weight^{\blowuprateofwaveWRTnumberofcommutations(N)} \ReciprocalLunitAppliedtoTimeFunction \, \voldiamond  \\
		& \ \     \lesssim \initialsmalldoublenull^2 +  (1 + \varsigma_0^{-1} + \varsigma_1^{-1}) \int_{\ubar' = \leftubar}^{\ubar} \totalfluxcontrolvelocity_{[1,N]}(\ubar',u) \, \mathrm{d} \ubar' \\
		& \ \ +  (1 + \varsigma_0^{-1} + \varsigma_1^{-1}) \int_{u' = \moreinterestingu_1}^u \totalfluxcontrolvelocity_{[1,N]}(\ubar,u') \, \mathrm{d} u' + \varsigma_1(1 + \varsigma_0^{-1}) \strongangularcontrolvelocity_{[1,N]}(\ubar,u) \\
		& \ \ +   \int_{u' = \moreinterestingu_1}^u \left\{\fluxcontrolVort_{\le N}(\ubar,u')+ \fluxcontrolGradEnt_{\le N }(\ubar,u')\right\} \, \mathrm{d} u'. 
	\end{split}
\end{align}\end{lemma}

\begin{proof}
We control the integral of the quadratic error terms generated by LHS\,\eqref{E:ENERGYESTIMATEHARMLESSWAVEERRORTERMS} using Def.\,\ref{D:HARMLESSWAVE}. Using \eqref{E:ARBITRARYCOMMUTATORSTRINGESTIMATE}, one easily sees that it suffices to control the spacetime integral of $\Lunit \tander^N v$ and $\muX \tander^N v$ times $ \left|\newuL \tander^{[1,N]} \velocityarray \right| + \upmu \left| \nullangrmD \tander^{[1,N]} \velocityarray \right|_{\gnulltori} + \left| \Lunit \tander^{\le N} \velocityarray \right| + \left| \angrmd \tander^{\le N} \velocityarray  \right|_{\gtorus} + \left|  \tander^{\le N}(\vortrenormalized,\GradEnt)\right| +   \left|\comdersmall^{[1,N];1} \controlvars \right| +  \left| \tandersmall^{[1,N]} \badcontrolvars \right|$.

Consider first the terms generated by $\Lunit \tander^N v$. Since $\ReciprocalLunitAppliedtoTimeFunction \approx 1$ and $\weight = |\ubar| < 1$, Young's inequality, Prop.\,\ref{P:IMPROVEMENTOFAUXILIARYBOOTSTRAP}, and \eqref{E:COERCIVENESSOFNULLFLUXCONTROLWAVE} imply:
\begin{align} \label{E:ENERGYESTIMATEFORHARMLESSERRORTERMSSTEP1} 
	\begin{split}
		& \int_{\characteristicdiamondtwoarg{[\leftubar,\ubar]}{[\moreinterestingu_1,u]}} 	
			(1 + 2\upmu) \left| \Lunit \tander^N v \right| \left(  \left|\newuL \tander^{[1,N]} \velocityarray \right| + \upmu \left| \nullangrmD \tander^{[1,N]} \velocityarray \right|_{\gnulltori} +  \left| \Lunit \tander^{\le N} v \right| \right)  \weight^{\blowuprateofwaveWRTnumberofcommutations(N)} \ReciprocalLunitAppliedtoTimeFunction \, \voldiamond \\
			& \ \   \lesssim \int_{\characteristicdiamondtwoarg{[\leftubar,\ubar]}{[\moreinterestingu_1,u]}} \left|\Lunit \tander^{\le N} \velocityarray \right|^2 \weight^{\blowuprateofwaveWRTnumberofcommutations(N)} \ReciprocalLunitAppliedtoTimeFunction \, \voldiamond  + \int_{\characteristicdiamondtwoarg{[\leftubar,\ubar]}{[\moreinterestingu_1,u]}} \left|\newuL \tander^{[1,N]} \velocityarray \right|^2 \weight^{\blowuprateofwaveWRTnumberofcommutations(N)} \ReciprocalLunitAppliedtoTimeFunction \, \voldiamond \\
			& \ \ +  \int_{\characteristicdiamondtwoarg{[\leftubar,\ubar]}{[\moreinterestingu_1,u]}} \upmu \left| \nullangrmD \tander^{[1,N]} \velocityarray \right|_{\gnulltori}^2 \weight^{\blowuprateofwaveWRTnumberofcommutations(N)} \ReciprocalLunitAppliedtoTimeFunction \, \voldiamond \\
			& \ \ \lesssim  \int_{u' = \moreinterestingu_1}^u \totalfluxcontrolvelocity_{[1,N]}(\ubar,u') \, \mathrm{d} u' + \int_{\ubar' = \leftubar}^{\ubar} \totalfluxcontrolvelocity_{[1,N]}(\ubar',u) \, \mathrm{d} \ubar'.
	\end{split}
\end{align}
We clarify that we used Def.\,\ref{D:VOLFORMS} to express the spacetime integral over the characteristic diamond as an iterated integral over the ingoing and outgoing characteristic surfaces. Using \eqref{E:UNIFIEDCOMMUTATORSPACETIMESTRONGANGULARVELOCITYCONTROL} and \eqref{E:COERCIVITYOLDSPACETIMETERM}, one sees easily from Young's inequality that: 
\begin{align} \label{E:ENERGYESTIMATEFORHARMLESSERRORTERMSSTEP2}
	 \int_{\characteristicdiamondtwoarg{[\leftubar,\ubar]}{[\moreinterestingu_1,u]}} 	
			(1 + 2\upmu) \left| \Lunit \tander^N v \right| \cdot \left| \angrmd \tander^{[1,N]} \velocityarray  \right|_{\gtorus}  \weight^{\blowuprateofwaveWRTnumberofcommutations(N)} \ReciprocalLunitAppliedtoTimeFunction \, \voldiamond  \lesssim  \varsigma_1^{-1}  \int_{u' = \moreinterestingu_1}^u \totalfluxcontrolvelocity_{[1,N]}(\ubar,u') \, \mathrm{d} u' + \varsigma_1\strongangularcontrolvelocity_{[1,N]}(\ubar,u).
\end{align}
We now control the error integral generated by the product of $\Lunit \tander^N v$ and $\tandersmall^{[1,N]} \upmu$ as a representative example of the quadratic terms generated by $\left|\comdersmall^{[1,N];1} \controlvars \right| +  \left| \tandersmall^{[1,N]} \badcontrolvars \right|$. Using again  $\ReciprocalLunitAppliedtoTimeFunction \approx 1$, $\weight = |\ubar| < 1$, Young's inequality, as well as \eqref{E:PRELIMINARYBELOWTOPORDERL2ESTIMATESFOREIKONALFUINCTIONQUANTITIES} we see: 
\begin{align} \label{E:ENERGYESTIMATEFORHARMLESSERRORTERMSSTEP3} 
	\begin{split}
		& \int_{\characteristicdiamondtwoarg{[\leftubar,\ubar]}{[\moreinterestingu_1,u]}} 	
				(1 + 2\upmu) \left| \Lunit \tander^N v \right| \cdot \left| \tandersmall^{[1,N]} \upmu \right| \weight^{\blowuprateofwaveWRTnumberofcommutations(N)} \ReciprocalLunitAppliedtoTimeFunction \, \voldiamond    \\
		& \ \ \lesssim \varsigma_1^{-1}  \int_{u' = \moreinterestingu_1}^u \totalfluxcontrolvelocity_{[1,N]}(\ubar,u') \, \mathrm{d} u' + \varsigma_1 \int_{\ubar' = \leftubar}^{\ubar} \int_{\ingoingcharacteristicsurfacetwoarg{\ubar'}{[\moreinterestingu_1,u]}}  
			\left| \tandersmall^{[1,N]} \upmu \right|^2 \weight^{\blowuprateofacousticgeoWRTnumberofcommutations(N)} \volingoingnullhypersurface \mathrm{d} \ubar' \\
		& \ \ \lesssim \initialsmalldoublenull^2 + (1 + \varsigma_0^{-1} + \varsigma_1^{-1}) \int_{u' = \moreinterestingu_1}^u \totalfluxcontrolvelocity_{[1,N]}(\ubar,u') \, \mathrm{d} u' +\varsigma_1 (1 + \varsigma_0^{-1}) \strongangularcontrolvelocity_{[1,N]}(\ubar,u).
	\end{split}
\end{align}
We clarify that in the last inequality on RHS\,\eqref{E:ENERGYESTIMATEFORHARMLESSERRORTERMSSTEP3} we used that the $L^2$-controlling quantities are monotonically increasing in their arguments and that the $\initialsmalldoublenull^2$ term came from the first term RHS\,\eqref{E:PRELIMINARYBELOWTOPORDERL2ESTIMATESFOREIKONALFUINCTIONQUANTITIES}. 

Similarly to \eqref{E:ENERGYESTIMATEFORHARMLESSERRORTERMSSTEP3}, since $\weight^{\blowuprateofwaveWRTnumberofcommutations(N)} \lesssim \weight^{\blowuprateoftransportWRTnumberofcommutations(N)}$, we easily bound $\int_{\characteristicdiamondtwoarg{[\leftubar,\ubar]}{[\moreinterestingu_1,u]}} 	
				(1 + 2\upmu) \left| \Lunit \tander^N v \right| \left| \tander^{\le N}(\vortrenormalized,\GradEnt)\right| \weight^{\blowuprateofwaveWRTnumberofcommutations(N)} \ReciprocalLunitAppliedtoTimeFunction \, \voldiamond  \lesssim  \int_{u' = \moreinterestingu_1}^u \totalfluxcontrolvelocity_{[1,N]}(\ubar,u') \rmd u' +  \int_{u' = \moreinterestingu_1}^u \left\{\fluxcontrolVort_{\le N}(\ubar,u')+ \fluxcontrolGradEnt_{\le N }(\ubar,u')\right\} \, \mathrm{d} u'$.

Using \eqref{E:L2ESTIMATESFORMUXTRANSVERSALDERIVATIVESOFWAVEVARIABLES}, the analysis for the error integrals generated by $\muX \tander^N v$ is carried out in a similar way. We omit the details. 

\end{proof}

\subsubsection{Estimates for the error integrals generated by the multiplier vectorfield}
\label{SSS:ESTIMATESFORERRORTERMSGENERATEDBYMULTIPLIERVECTORFIELD}
In the next lemma, we bound the error integrals that are generated by the multiplier vectorfield. 

\begin{lemma}[Estimates for the error integrals generated by the multiplier vectorfield]
\label{L:ESTIMATESFORERRORTERMSGENERATEDBYMULTIPLIERVECTORFIELD}
Assume that $1 \leq N \leq \Ntop$,
$v \in \velocityarray = \{v^1,v^2,v^3\}$.
Let $\tander^N \in \mathfrak{P}^{(N)}$, where $\mathfrak{P}^{(N)}$
is the set of order $N$ $\nullhyparg{u}$-tangential commutator operators from
Def.\,\ref{D:STRINGSOFCOMMUTATIONVECTORFIELDS}. Let $\blowuprateofwaveWRTnumberofcommutations(N)$ denote the blowup rate of the wave variables defined in \eqref{E:BLOWUPRATEOFWAVEWRTNUMBEROFCOMMUTATORS}. Let $\varsigma_0$ be the constant featured in Lemma\,\ref{L:PRELIMINARYBELOWTOPORDERL2ESTIMATESFOREIKONALFUINCTIONQUANTITIES} (see also Remark\,\ref{R:SMALLNESSCONSTANTFORBELOWTOPORDERACOUSTICVARIABLES}).
Recall that the multiplier vectorfield $\multipliervectorfield$ is defined in \eqref{E:MULTIPLIERVECTORFIELD}
and that ${^{(\multipliervectorfield)}\mathfrak{B}}[\tander^N v]$ 
is the error term defined by 
\eqref{E:WAVEENERGYIDENTITYBULKERRORTERM} and \eqref{E:WAVEENERGYIDENTITYBULKERRORTERM1}--\eqref{E:WAVEENERGYIDENTITYBULKERRORTERM6}
and appearing on RHS\,\eqref{E:ENERGYNULLFLUXINTEGRALIDENTITIESWAVE}.
Then if $\varsigma_1 \in (0,1]$ is sufficiently small, the following estimates hold for 
$(\ubar,u) \in [\leftubar,\ubarboot) \times [\moreinterestingu_1,\moreinterestingu_2]$, 
where the implicit constants are \textbf{independent} of $\varsigma_0,\varsigma_1$:
\begin{align}
	\begin{split} \label{E:WAVEENERGYESTIMATEMULTIPLIERERRORTERMS}
		\int_{\characteristicdiamondtwoarg{[\leftubar,\ubar]}{[\moreinterestingu_1,u]}} 	
			\left|
				{^{(\multipliervectorfield)}\mathfrak{B}}[\tander^N v] 
			\right| \weight^{\blowuprateofwaveWRTnumberofcommutations(N)} 
		\ReciprocalLunitAppliedtoTimeFunction \, \voldiamond & \lesssim (1 + \varsigma_1^{-1}) \int_{\ubar' = \leftubar}^{\ubar} \frac{1}{|\ubar'|^{1/2}} \totalfluxcontrolvelocity_{[1,N]}(\ubar',u) \, \mathrm{d} \ubar' \\
		& \ \  +  (1  +  \varsigma_1^{-1}) \int_{u' = \moreinterestingu_1}^u \totalfluxcontrolvelocity_{[1,N]}(\ubar,u') \, \mathrm{d} u' \\
		& \ \ + \varsigma_1 \strongangularcontrolvelocity_{[1,N]}(\ubar,u).
	\end{split}
\end{align}\end{lemma}

\begin{proof}

We first bound the error integral generated by the first term on RHS\,\eqref{E:WAVEENERGYIDENTITYBULKERRORTERM3}. Using  \eqref{E:POINTWISESTIMATEFORBREVEXINTERMSOFNEWUL}--\eqref{E:POINTWISESMOOTHTORIDIFFERENTIALNORMINTERMSOFDOULBENULLTORIDIFFERENTIALINORMANDL}, the crucial estimate \eqref{E:NEWULMUUBOUNDEDBYSQRTMU}, and  Prop.\,\ref{P:IMPROVEMENTOFAUXILIARYBOOTSTRAP}, we find that the term $(\muX \upmu) \left| \angrmd \tander^N v \right|_{\gtorus}^2$  from \eqref{E:WAVEENERGYIDENTITYBULKERRORTERM3} is bounded in magnitude by $\lesssim \sqrt{\upmu}  \left| \nullangrmD \tander^N v \right|_{\gtorus}^2 + \left| \Lunit \tander^N v  \right|^2$. Hence, using  \eqref{E:KEYESTIMATECONTROLLINGINVERSEMUBYINVERSEWEIGHT} and the coercivity \eqref{E:COERCIVENESSOFNULLFLUXCONTROLWAVE}, we can bound the integral of this term over the characteristic diamond $\characteristicdiamondtwoarg{[\leftubar,\ubar]}{[\moreinterestingu_1,u]}$ by $\lesssim \int_{\ubar' = \leftubar}^{\ubar} \frac{1}{|\ubar'|^{1/2}} \totalfluxcontrolvelocity_{[1,N]}(\ubar',u) \, \mathrm{d} \ubar' +  \int_{u' = \moreinterestingu_1}^u \totalfluxcontrolvelocity_{[1,N]}(\ubar,u') \, \mathrm{d} u',$ as desired.

We now handle the remaining terms in the definition \eqref{E:WAVEENERGYIDENTITYBULKERRORTERM3}
of ${^{(\multipliervectorfield)}\mathfrak{B}_{(3)}[\tander^N v]}$
as well as the
remaining bulk terms ${^{(\multipliervectorfield)}\mathfrak{B}_{(i)}[\tander^N v]}$ with $i \in\{1,2,4,5,6\}$,
i.e., the terms defined in 
\eqref{E:WAVEENERGYIDENTITYBULKERRORTERM1},
\eqref{E:WAVEENERGYIDENTITYBULKERRORTERM2},
\eqref{E:WAVEENERGYIDENTITYBULKERRORTERM4},
\eqref{E:WAVEENERGYIDENTITYBULKERRORTERM5},
and
\eqref{E:WAVEENERGYIDENTITYBULKERRORTERM6}.
To this end, we first use Prop.\,\ref{P:SCHEMATICSTRUCTUREOFVARIOUSTENSORSINTERMSOFCONTROLVARS}, 
the bootstrap assumptions, 
and Young's inequality to deduce the following pointwise estimates,
valid for any $\varsigma_1 \in (0,1]$ with implicit constants that are \textbf{independent} of $\varsigma_1$:
\begin{align} 
	\begin{split} \label{E:POINTWISEESTIMATESFOREASYTERMSGENERATEDBYMULTIPLIER}
		& \left| 2 \upmu  \Lunit \upmu + 
			\frac{1}{2}
			\upmu \mytr_{\gtorus}\upchi 
			+ 
			\upmu^2 \mytr_{\gtorus} \angktan 
			+ 
			\upmu \mytr_{\gtorus} \angktrans 
			\right|  \left| \angrmd \tander^N v \right|_{\gtorus}^2, \\
		& \ \ \max_{i \in\{1,2,4,5,6\}} \left|{^{(\multipliervectorfield)}\mathfrak{B}_{(i)}[\tander^N v]} \right| \\
		& \ \ \lesssim (1 + \varsigma_1)^{-1} (\Lunit \tander^N v)^2 + (1+ \varsigma_1^{-1}) (\muX \tander^N v )^2 \\
		& \ \ + \upmu   \left| \angrmd \tander^N v \right|_{\gtorus}^2  + \varsigma_1   \left| \angrmd \tander^N v \right|_{\gtorus}^2.
	\end{split}
\end{align}
Using \eqref{E:COERCIVENESSOFNULLFLUXCONTROLWAVE}, \eqref{E:COERCIVITYOLDSPACETIMETERM}, \eqref{E:POINTWISESTIMATEFORBREVEXINTERMSOFNEWUL}--\eqref{E:POINTWISESMOOTHTORIDIFFERENTIALNORMINTERMSOFDOULBENULLTORIDIFFERENTIALINORMANDL}, the integral of the product of RHS\,\eqref{E:POINTWISEESTIMATESFOREASYTERMSGENERATEDBYMULTIPLIER} and $\weight^{\blowuprateofwaveWRTnumberofcommutations(N)}$   over the characteristic diamond $\characteristicdiamondtwoarg{[\leftubar,\ubar]}{[\moreinterestingu_1,u]}$ is $\lesssim$ as desired.

\end{proof}

\subsubsection{Estimates for the error integrals generated by the weight}
\label{SSS:ESTIMATESFORERRORTERMSGENERATEDBYTHEWEIGHT}

In the next lemma, we bound the error integrals that are generated by the weight in our fundamental energy identity.

\begin{lemma}[Estimates for the error integrals generated by the weight]\label{L:ESTIMATESFORERRORTERMSGENERATEDBYTHEWEIGHT}
Assume that $1 \le N \le \Ntop, \, v \in \velocityarray = \{v^1,v^2,v^3\}$. Let $\tander^N \in \mathfrak{P}^{(N)}$, where $\mathfrak{P}^{(N)}$
is the set of order $N$ $\nullhyparg{u}$-tangential commutator operators from
Def.\,\ref{D:STRINGSOFCOMMUTATIONVECTORFIELDS}. Let $\blowuprateofwaveWRTnumberofcommutations(N)$ denote the blowup rate of the wave variables defined in \eqref{E:BLOWUPRATEOFWAVEWRTNUMBEROFCOMMUTATORS}.
Recall that the vectorfield  $\ToriTangentVectorfieldAssociatedToDoubleNullFolliations$ is defined in \eqref{E:TORITANGENTVECTORFIELDASSOCIATEDTODOUBLENULLFRAME}. Then the following estimates hold for 
$(\ubar,u) \in [\leftubar,\ubarboot) \times [\moreinterestingu_1,\moreinterestingu_2]$:
	\begin{align}
		\begin{split} \label{E:ESTIMATESFORERRORTERMSGENERATEDBYTHEWEIGHT}
			\int_{\characteristicdiamondtwoarg{[\leftubar,\ubar]}{[\moreinterestingu_1,u]}} 	
				\left| \MagnitueofinnerproductofnewLandnewuL (\ToriTangentVectorfieldAssociatedToDoubleNullFolliations \tander^N v) (\newuL \tander^N v) (\newL \weight) \right|  \weight^{\blowuprateofwaveWRTnumberofcommutations(N)-1} \voldiamond \lesssim \fundbootsmall  \weakspacetimetransvelocitycontrol_N(\ubar,u)  + \fundbootsmall\weakspacetimeangularvelocitycontrol_N(\ubar,u)
		\end{split}
	\end{align}
\end{lemma}

\begin{proof}
	We first note that the $\doublenulltoritwoarg{\ubar}{u}$-tangency of $\ToriTangentVectorfieldAssociatedToDoubleNullFolliations$ implies that $\nulltorusproject \ToriTangentVectorfieldAssociatedToDoubleNullFolliations = \ToriTangentVectorfieldAssociatedToDoubleNullFolliations$. Hence, using the $\gnulltori$-Schwarz inequality, \eqref{E:MAGNITUDEOFINNERPRODUCTOFNEWLANDNEWULAPPROXIMATELYMU}--\eqref{E:RECIPROCALULUNITAPPLIEDTOTIMEFUNCTIONAPPROXIMATELYMU} and  \eqref{E:POINTWISEESTIMATEFORCARTESIANCOMPONENTSOFDOUBLENULLTORITANGENTVECTORFIELD}, we have $|\MagnitueofinnerproductofnewLandnewuL \ToriTangentVectorfieldAssociatedToDoubleNullFolliations \tander^N \Psi|\lesssim \fundbootsmall( \upmu + 2 \upmu \ReciprocaluLunitAppliedtoTimeFunction)  \left| \nullangrmD \tander^N \Psi \right|_{\gnulltori}$. Using $\newL \weight = 1$, \eqref{E:WAVENEWSPACETIMENULLANGULARL2CONTROLLINGQUANTITY}--\eqref{E:WAVENEWSPACETIMEL2TRANSVERSALCONTROLLINGQUANTITY}, and a straightforward application of Young's inequality, the result \eqref{E:ESTIMATESFORERRORTERMSGENERATEDBYTHEWEIGHT} follows.
	
\end{proof}

\subsection{Imprecise top-order $L^2$ estimates for $\upchi$}
\label{SS:IMPRECISETOPORDERL2ESTIMATSFORCHI}
In this section, we derive $L^2$ \emph{spacetime} estimates for the top-order terms 
$\upmu \tander^{\Ntop} \mytr_{\gtorus} \upchi$ and $\upmu \angLie_{\tander}^{\Ntop} \upchi$ 
in terms of the $L^2$-controlling quantities. The key ingredient will be a preliminary estimate for the fully modified quantity $\fullymodquant{\tander^N}$ from Def.\,\ref{D:FULLYANDPARTIALLYMODIFIEDQUANTITIES}, which will be provided in Sect.\,\ref{SSS:IMPRECISETOPORDERESTIMATESFORFULLYMODIFIEDCHI}. By ``imprecise'' we mean that we do not carefully keep track of the constants.

\subsubsection{Statement of the imprecise top-order $L^2$ estimates for $\upchi$}
\label{SSS:STATEMENTOFIMPRECISETOPORDERL2ESTIMATSFORCHI}
In the next proposition, we state the estimates. Its proof is located in Sect.\,\ref{SSS:PROOFOFIMPRECISETOPORDERL2ESTIMATSFORCHI}.

\begin{proposition}[Top-order imprecise $L^2$ estimates for $\upchi$] \label{P:TOPORDERIMPRECISEL2ESTIMATEMUCHI}
Let $N = \Ntop$, and let
$\mathfrak{P}^{(N)}$, $\angLie_{\mathfrak{P}}^{(N)}$
be the sets of order $N$ $\nullhyparg{u}$-tangential commutator operators from
Def.\,\ref{D:STRINGSOFCOMMUTATIONVECTORFIELDS}. Let $\blowuprateofwaveWRTnumberofcommutations(N) = \blowupratetoporderwave$ denote the blowup rate of the wave variables defined in \eqref{E:BLOWUPRATEOFWAVEWRTNUMBEROFCOMMUTATORS}. Let $\varsigma_0$ be the constant featured in Lemma\,\ref{L:PRELIMINARYBELOWTOPORDERL2ESTIMATESFOREIKONALFUINCTIONQUANTITIES} (see also Remark\,\ref{R:SMALLNESSCONSTANTFORBELOWTOPORDERACOUSTICVARIABLES}).
Then if $\varsigma_1 \in(0,1]$ is sufficiently small, for $(\ubar,u) \in [\leftubar,\ubarboot) \times [\moreinterestingu_1,\moreinterestingu_2]$, the following estimate holds for the imprecise fully modified controlling quantities $\totalfluxcontrolimprecisefullymodifiedchi(\ubar,u), \,  \totalbulkcontrolimprecisefullymodifiedchi(\ubar,u)$ from \eqref{E:IMPRECISEFULLYMODQUANTFLUXOTALCONTROL}--\eqref{E:IMPRECISEFULLYMODQUANTBULKTOTALCONTROL}, where the implicit constants are \textbf{independent} of $\varsigma_0,\varsigma_1$:
\begin{align}
	\begin{split} \label{E:FINALIMPRECISETOPORDERESTIMATESFORFULLYMODIFIEDCHI}
		 \totalfluxcontrolimprecisefullymodifiedchi(\ubar,u)  + \totalbulkcontrolimprecisefullymodifiedchi(\ubar,u) &  \lesssim  (1 + \varsigma_0^{-1}) \initialsmalldoublenull^2   +  (1 + \varsigma_0^{-1} +  \varsigma_1^{-1}) \int_{\ubar' = \leftubar}^{\ubar}  \totalfluxcontrolimprecisefullymodifiedchi(\ubar',u) \, \mathrm{d} \ubar'  \\
		& \ \ + (1 + \varsigma_0^{-1}) \int_{\ubar' = \leftubar}^{\ubar} \totalfluxcontrolvelocity_{[1,N]}(\ubar',u) \, \mathrm{d} \ubar'\   + (1+ \varsigma_0^{-1}) \int_{u' = \moreinterestingu_1}^u  \totalfluxcontrolvelocity_N(\ubar,u') \, \rmd u' \\
		& \ \ + \varsigma_1 (1 + \varsigma_0^{-1})\strongangularcontrolvelocity_{N}(\ubar,u) \\
		& \ \  + \bulkcontrolVortVort_N(\ubar,u) +  \bulkcontrolDivGradEnt_N(\ubar,u) 
			  + \int_{u' = \moreinterestingu_1}^u
			\left\lbrace
				\fluxcontrolVortVort_{\le N-1}(\ubar,u')
				+
				\fluxcontrolDivGradEnt_{\le N-1}(\ubar,u')
			\right\rbrace
		\, \rmd u' \\
		& \ \ +  \int_{u' = \moreinterestingu_1}^u \left\{\fluxcontrolVort_{\le N}(\ubar,u')+ \fluxcontrolGradEnt_{\le N }(\ubar,u')\right\} \, \mathrm{d} u', 
 	\end{split}.
\end{align}
In addition, the following estimates hold for $\upmu \tander^N \mytr_{\gtorus}\upchi$ and $\upmu |\angLie_{\tander}^N \upchi|_{\gtorus}$:
\begin{align}
	\begin{split} \label{E:TOPORDERL2ESTIMATEMUCHI}
		& \max_{\tander^N \in \mathfrak{P}^{(N)}}  \int_{\characteristicdiamondtwoarg{[\leftubar,\ubar)}{[\moreinterestingu_1,u]}}  
			 \left| \upmu \tander^N \mytr_{\gtorus}\upchi \right|^2 \weight^{\blowupratetoporderwave} \voldiamond,
			 \, \max_{\angLie_{\tander}^N \in \angLie_{\mathfrak{P}}^{(N)}} \int_{\characteristicdiamondtwoarg{[\leftubar,\ubar)}{[\moreinterestingu_1,u]}} 
			 \left| \upmu \angLie_{\tander}^N \upchi \right|_{\gtorus}^2 \weight^{\blowupratetoporderwave} \voldiamond \\
		&  \lesssim  \initialsmalldoublenull^2 + \totalbulkcontrolimprecisefullymodifiedchi(\ubar,u) +   (1 + \varsigma_0^{-1}) \int_{\ubar' = \leftubar}^{\ubar} \totalfluxcontrolvelocity_{[1,N]}(\ubar',u) \, \mathrm{d} \ubar'   \\
			& \ \ +  (1 + \varsigma_0^{-1}) \int_{u' = \moreinterestingu_1}^u \totalfluxcontrolvelocity_{[1,N]}(\ubar,u') \, \mathrm{d} u' + (1 + \varsigma_0^{-1}) \strongangularcontrolvelocity_{[1,N]}(\ubar,u)  \\
			& \ \ +  \int_{u' = \moreinterestingu_1}^u \left\{\fluxcontrolVort_{\le N}(\ubar,u')+ \fluxcontrolGradEnt_{\le N }(\ubar,u')\right\} \, \mathrm{d} u'.
	\end{split}
\end{align}\end{proposition}

\subsubsection{Preliminary estimates for $\upchi$}
\label{SSS:PRELIMINARYESTIAMTSFORTOPORDERL2ESTIMATSFORCHI}
In the following lemma, we derive preliminary $L^2$ estimates for $\upmu \tander^N \mytr_{\gtorus}\upchi, \,  \upmu \angLie_{\tander}^N \upchi$, and $\fullymodquant{\tander^N}$ that we will use in the proof of Prop.\,\ref{P:TOPORDERIMPRECISEL2ESTIMATEMUCHI}. 

\begin{lemma}[Preliminary top-order $L^2$ estimates for $\upchi$]  \label{L:PRELIMINARYESTIAMTSFORTOPORDERL2ESTIMATSFORCHI}
Let $N = \Ntop$, and let
$\mathfrak{P}^{(N)}$
be the set of order $N$ $\nullhyparg{u}$-tangential commutator operators from
Def.\,\ref{D:STRINGSOFCOMMUTATIONVECTORFIELDS}. Let $\blowuprateofwaveWRTnumberofcommutations(N) = \blowupratetoporderwave$ denote the blowup rate of the wave variables defined in \eqref{E:BLOWUPRATEOFWAVEWRTNUMBEROFCOMMUTATORS}. Let $\varsigma_0$ be the constant featured in Lemma\,\ref{L:PRELIMINARYBELOWTOPORDERL2ESTIMATESFOREIKONALFUINCTIONQUANTITIES} (see also Remark\,\ref{R:SMALLNESSCONSTANTFORBELOWTOPORDERACOUSTICVARIABLES}).
Then the following estimates hold for
$(\ubar,u) \in [\leftubar,\ubarboot) \times [\moreinterestingu_1,\moreinterestingu_2]$, where the implicit constants are \textbf{independent} of $\varsigma_0$:\begin{align}
	\begin{split} \label{E:PRELIMINARYTOPORDERTRACECHI}
		 \int_{\characteristicdiamondtwoarg{[\leftubar,\ubar)}{[\moreinterestingu_1,u]}}  
		 \left| \upmu \tander^N \mytr_{\gtorus} \upchi \right|^2 \weight^{\blowupratetoporderwave} \voldiamond
		 & \lesssim  \initialsmalldoublenull^2 + \totalbulkcontrolimprecisefullymodifiedchi(\ubar,u) +   (1 + \varsigma_0^{-1}) \int_{\ubar' = \leftubar}^{\ubar} \totalfluxcontrolvelocity_{[1,N]}(\ubar',u) \, \mathrm{d} \ubar'   \\
			& \ \ +  (1 + \varsigma_0^{-1}) \int_{u' = \moreinterestingu_1}^u \totalfluxcontrolvelocity_{[1,N]}(\ubar,u') \, \mathrm{d} u' + (1 + \varsigma_0^{-1}) \strongangularcontrolvelocity_{[1,N]}(\ubar,u) \\
			& \ \ +  \int_{u' = \moreinterestingu_1}^u \left\{\fluxcontrolVort_{\le N}(\ubar,u')+ \fluxcontrolGradEnt_{\le N }(\ubar,u')\right\} \, \mathrm{d} u',
	\end{split} \\
	\begin{split}\label{E:PRELIMINARYTOPORDERLIECHI}
		 \int_{\characteristicdiamondtwoarg{[\leftubar,\ubar)}{[\moreinterestingu_1,u]}}  
				 \left| \upmu \angLie_{\tander}^N \upchi \right|_{\gtorus}^2 \weight^{\blowupratetoporderwave} \voldiamond & \lesssim  \initialsmalldoublenull^2 + (1 + \varsigma_0^{-1}) \int_{u' = \moreinterestingu_1}^u \totalfluxcontrolvelocity_{[1,N]}(\ubar,u') \, \mathrm{d} u' + (1 + \varsigma_0^{-1}) \strongangularcontrolvelocity_{[1,N]}(\ubar,u) \\
				 & \ \ +   \int_{u' = \moreinterestingu_1}^u \left\{\fluxcontrolVort_{\le N}(\ubar,u')+ \fluxcontrolGradEnt_{\le N }(\ubar,u')\right\} \, \mathrm{d} u' \\
				 & \ \ + \max_{\tander^N \in \mathfrak{P}^{(N)}}  \int_{\characteristicdiamondtwoarg{[\leftubar,\ubar)}{[\moreinterestingu_1,u]}}  
				 \left| \upmu \tander^N \mytr_{\gtorus}\upchi \right|^2 \weight^{\blowupratetoporderwave} \voldiamond.  
	\end{split}
\end{align}
\end{lemma}

\begin{proof}

\hfill

\noindent \textbf{Proof of \eqref{E:PRELIMINARYTOPORDERTRACECHI}}:

Using the definition of the fully modified quantity \eqref{E:FULLYMODIFIEDQUANTITY}, we  easily bound $ \int_{\characteristicdiamondtwoarg{[\leftubar,\ubar)}{[\moreinterestingu_1,u]}}  
		 \left| \upmu \tander^N \mytr_{\gtorus} \upchi \right|^2 \weight^{\blowupratetoporderwave} \voldiamond \lesssim  \int_{\characteristicdiamondtwoarg{[\leftubar,\ubar)}{[\moreinterestingu_1,u]}}  \left| \, \fullymodquant{\tander^N} \right|^2 \weight^{\blowupratetoporderwave} \voldiamond +  \int_{\characteristicdiamondtwoarg{[\leftubar,\ubar)}{[\moreinterestingu_1,u]}}  \left|\tander^N \mathfrak{X} \right|^2 \weight^{\blowupratetoporderwave} \voldiamond$. The first integral is $\lesssim$ RHS\,\eqref{E:PRELIMINARYTOPORDERTRACECHI} by \eqref{E:IMPRECISEFULLYMODQUANTBULKTOTALCONTROL} and \eqref{E:COERCIVITYNEWSPACETIMETERM}. To bound the square integral of $\tander^N \mathfrak{X}$, we insert the pointwise bound \eqref{E:POINTWISEESTIMATETANGENTDERIVATIVESOFMODQUANTINHOM}. We use Lemma\,\ref{L:COERCIVENESSOFL2CONTROLLINGQUANITIESUNIFIEDCOMMUTATOR} to bound the resulting integrals from the first two terms on RHS\,\eqref{E:POINTWISEESTIMATETANGENTDERIVATIVESOFMODQUANTINHOM} by the two $\totalfluxcontrolvelocity_{[1,N]}(\ubar,u)$ integrals on RHS\,\eqref{E:PRELIMINARYTOPORDERTRACECHI}. To bound the spacetime square integral of the $\comdersmall^{[1,N]}\velocityarray$ terms present on RHS\,\eqref{E:POINTWISEESTIMATETANGENTDERIVATIVESOFMODQUANTINHOM}, we use  \eqref{E:ARBITRARYCOMMUTATORSTRINGESTIMATE} with $N$ in place of $N+1$ and \eqref{E:COERCIVENESSOFNULLFLUXCONTROLWAVE}, \eqref{E:COERCIVITYOLDSPACETIMETERM} to conclude that the integrals of the first three terms on RHS\,\eqref{E:ARBITRARYCOMMUTATORSTRINGESTIMATE} are $\lesssim$ RHS\,\eqref{E:PRELIMINARYTOPORDERTRACECHI}. To control the integrals generated by the vorticity and gradient entropy, since $\weight^{\blowupratetoporderwave} \lesssim \weight^{\blowuprateoftransportWRTnumberofcommutations(N)}$ for any $N$, we easily bound the spacetime integrals by $\int_{u' = \moreinterestingu_1}^u \left\{\fluxcontrolVort_{\le N}(\ubar,u')+ \fluxcontrolGradEnt_{\le N }(\ubar,u')\right\} \, \mathrm{d} u'$.
It remains to control the spacetime integrals $ \int_{\characteristicdiamondtwoarg{[\leftubar,\ubar)}{[\moreinterestingu_1,u]}} \left\{ |\tander^{[1,N]} \controlvars|^2 + |\tandersmall^{[1,N]}\badcontrolvars|^2\right\} \weight^{\blowupratetoporderwave} \voldiamond $. To this end, we integrate \eqref{E:PRELIMINARYBELOWTOPORDERL2ESTIMATESFOREIKONALFUINCTIONQUANTITIES}--\eqref{E:PRELIMINARYBELOWTOPORDERL2ESTIMATESFOREIKONALFUINCTIONQUANTITIESWITHZ} with respect to $\ubar$, and use the fact that RHS of \eqref{E:PRELIMINARYBELOWTOPORDERL2ESTIMATESFOREIKONALFUINCTIONQUANTITIES}--\eqref{E:PRELIMINARYBELOWTOPORDERL2ESTIMATESFOREIKONALFUINCTIONQUANTITIESWITHZ} are increasing with respect to $\ubar$.
		 
\hfill

\noindent \textbf{Proof of \eqref{E:PRELIMINARYTOPORDERLIECHI}}:

We apply the elliptic estimate \eqref{E:ELLIPTICESTIMATECHI} with 
$\upxi \eqdef \angLie_{\tander}^{N-1}\upchi$ to deduce:
\begin{align}  
\begin{split}	\label{E:PRELIMINARYTOPORDERLIECHIINTERMEDIATESTEP1}  
	 \int_{\characteristicdiamondtwoarg{[\leftubar,\ubar)}{[\moreinterestingu_1,u]}} 
		\upmu^2 |\angLie_{\tander}^N \upchi |_{\gtorus}^2  \weight^{\blowupratetoporderwave}
	\, \voldiamond 
	& 
	\lesssim   \int_{\characteristicdiamondtwoarg{[\leftubar,\ubar)}{[\moreinterestingu_1,u]}} 
  	\upmu^2 |\angLie_{\Lunit} \upxi|_{\gtorus}^2 \weight^{\blowupratetoporderwave}
	\, \voldiamond 
	+ 
 	 \int_{\characteristicdiamondtwoarg{[\leftubar,\ubar)}{[\moreinterestingu_1,u]}} 
		\upmu^2 |\angdiv \upxi|_{\gtorus}^2 \weight^{\blowupratetoporderwave}
	\, \voldiamond 
		\\
	& \ \
	+ 
	\sum_{A=2,3} 
	 \int_{\characteristicdiamondtwoarg{[\leftubar,\ubar)}{[\moreinterestingu_1,u]}} 
 		\upmu^2 (\Yvf{A} \mytr_{\gtorus}\upxi)^2 \weight^{\blowupratetoporderwave}
	\, \voldiamond 
	+ 
	\fundbootsmall 
	 \int_{\characteristicdiamondtwoarg{[\leftubar,\ubar)}{[\moreinterestingu_1,u]}} 
 		|\upxi|_{\gtorus}^2 \weight^{\blowupratetoporderwave}
	\, \voldiamond. 
\end{split}
\end{align}
Using the pointwise estimate \eqref{E:ANGLIELTANGENTIALCHIPOINTWISE} and arguments similar to the ones used in the proof of \eqref{E:PRELIMINARYTOPORDERTRACECHI},
the first integral on 
on RHS\,\eqref{E:PRELIMINARYTOPORDERLIECHIINTERMEDIATESTEP1}
is 
$
\lesssim
$ RHS\,\eqref{E:PRELIMINARYTOPORDERLIECHI}.
Those arguments also imply, based on \eqref{E:PRELIMINARYBELOWTOPORDERL2ESTIMATESFOREIKONALFUINCTIONQUANTITIES},
that the last integral on RHS\,\eqref{E:PRELIMINARYTOPORDERLIECHIINTERMEDIATESTEP1} is $\lesssim$ RHS  \eqref{E:PRELIMINARYTOPORDERLIECHI}.

To handle the third term on RHS\,\eqref{E:PRELIMINARYTOPORDERLIECHIINTERMEDIATESTEP1},
we start with the following pointwise triangle inequality estimate:
\begin{align} \label{E:PRELIMINARYTOPORDERLIECHIINTERMEDIATESTEP2}
\upmu^2 
\left|
	\Yvf{A} \mytr_{\gtorus} \angLie_{\tander}^{N-1} \upchi
\right|^2 \weight^{\blowupratetoporderwave}
& \lesssim
\upmu^2 \left|\Yvf{A} \tander^{N-1} \mytr_{\gtorus} \upchi\right|^2 \weight^{\blowupratetoporderwave}
+
\upmu^2
\left|
	\Yvf{A}
	\left(
		\mytr_{\gtorus} \angLie_{\tander}^{N-1} \upchi 
		- 
		\tander^{N-1} \mytr_{\gtorus} \upchi 
	\right)
\right|^2 \weight^{\blowupratetoporderwave}.
\end{align}
The integral of the first term $\upmu^2 \left|\Yvf{A} \tander^{N-1} \mytr_{\gtorus} \angLie_{\tander}^{N-1}\upchi \right|^2\weight^{\blowupratetoporderwave}$ over $\characteristicdiamondtwoarg{[\leftubar,\ubar)}{[\moreinterestingu_1,u]}$ on RHS\,\eqref{E:PRELIMINARYTOPORDERLIECHIINTERMEDIATESTEP2}
is $\lesssim$ the last term on RHS\,\eqref{E:PRELIMINARYTOPORDERLIECHI}.
To handle the second term on RHS\,\eqref{E:PRELIMINARYTOPORDERLIECHIINTERMEDIATESTEP2},
we first note the following pointwise commutator estimate,
which follows easily from the Leibniz rule, 
the bootstrap assumptions,
and \eqref{E:TANDERGANDCHIESTIMATE}:
\begin{align} \label{E:LIEANDTRACECOMMUTATOR}
	\left|
		\Yvf{A} 
		\left(\mytr_{\gtorus} \angLie_{\tander}^{N-1}\upchi 
		- 
		\tander^{N-1} \mytr_{\gtorus} \upchi\right)
	\right| 
	& \lesssim 	
	\left| \tander^{[1,N]} \controlvars \right|. 
\end{align}
Hence, the integral over $\characteristicdiamondtwoarg{[\leftubar,\ubar)}{[\moreinterestingu_1,u]}$ of the second term on RHS\,\eqref{E:PRELIMINARYTOPORDERLIECHIINTERMEDIATESTEP2} is $\lesssim$ RHS\,\eqref{E:PRELIMINARYTOPORDERLIECHI} using \eqref{E:PRELIMINARYBELOWTOPORDERL2ESTIMATESFOREIKONALFUINCTIONQUANTITIES}.

It remains for us to bound the $\angdiv \angLie_{\tander}^{N-1} \upchi$-involving
integral on RHS\,\eqref{E:PRELIMINARYTOPORDERLIECHIINTERMEDIATESTEP1}.
We start with the following pointwise triangle inequality estimate:
\begin{align} \label{E:PRELIMINARYTOPORDERLIECHIINTERMEDIATESTEP3}
\upmu^2
\left|
	\angdiv \angLie_{\tander}^{N-1} \upchi
\right|_{\gtorus}^2  \weight^{\blowupratetoporderwave}
& \lesssim
\upmu^2
\left|
	\angrmd \tander^{N-1} \mytr_{\gtorus} \upchi
\right|_{\gtorus} \weight^{\blowupratetoporderwave}
+
\upmu^2
\left| 
	\angdiv \angLie_{\tander}^{N-1} \upchi 
	- 
	\angrmd \tander^{N-1} \mytr_{\gtorus} \upchi 
\right|_{\gtorus} \weight^{\blowupratetoporderwave}.
\end{align}
By \eqref{E:ANGDFPOINTWISEBOUNDEDBYCOMMUTATORVECTORFIELDS}, the integral of the first term 
$
\upmu^2
\left|
	\angrmd \tander^{N-1} \mytr_{\gtorus} \upchi
\right|_{\gtorus}^2 \weight^{\blowupratetoporderwave}$
on RHS\,\eqref{E:PRELIMINARYTOPORDERLIECHIINTERMEDIATESTEP3} over $\characteristicdiamondtwoarg{[\leftubar,\ubar)}{[\moreinterestingu_1,u]}$
is $\lesssim$ the last term on RHS\,\eqref{E:PRELIMINARYTOPORDERLIECHI}.
To handle the second term on RHS\,\eqref{E:PRELIMINARYTOPORDERLIECHIINTERMEDIATESTEP3},
we multiply the estimate \eqref{E:CODAZZICOMMUTATORESTIMATES} by $\upmu \weight^{\blowupratetoporderwave/2}$, square the resulting inequality, and integrate over $\characteristicdiamondtwoarg{[\leftubar,\ubar)}{[\moreinterestingu_1,u]}$. The resulting error integral is handled using arguments similar to the above. 
We have therefore proved \eqref{E:PRELIMINARYTOPORDERLIECHI}.
\end{proof}

\subsubsection{Imprecise estimates for the top-order fully modified quantities} \label{SSS:IMPRECISETOPORDERESTIMATESFORFULLYMODIFIEDCHI}

In this section we derive imprecise top-order $L^2$ estimates for the fully modified quantities $\fullymodquant{\tander^N}$. By imprecise we mean that we will not carefully track the constants in the error terms as they do not affect the blowup rates. The reader should compare the result of Lemma\,\ref{L:IMPRECISETOPORDERESTIMATESFORFULLYMODIFIEDCHI} with that of Lemma\,\ref{L:L2ESTIMATESFORMOSTDIFFICULTPRODUCT}. 

\begin{lemma}[Imprecise $L^2$ estimates for the top-order fully modified quantities] \label{L:IMPRECISETOPORDERESTIMATESFORFULLYMODIFIEDCHI}

Let $N = \Ntop$,  let
$\mathfrak{P}^{(N)}$, $\angLie_{\mathfrak{P}}^{(N)}$ 
be the set of order $N$ $\nullhyparg{u}$-tangential commutator operators from
Def.\,\ref{D:STRINGSOFCOMMUTATIONVECTORFIELDS}. Let $\blowuprateofwaveWRTnumberofcommutations(N) = \blowupratetoporderwave$ denote the blowup rate of the wave variables defined in \eqref{E:BLOWUPRATEOFWAVEWRTNUMBEROFCOMMUTATORS}. Let $\varsigma_0$ be the constant featured in Lemma\,\ref{L:PRELIMINARYBELOWTOPORDERL2ESTIMATESFOREIKONALFUINCTIONQUANTITIES} (see also Remark\,\ref{R:SMALLNESSCONSTANTFORBELOWTOPORDERACOUSTICVARIABLES}).
Then if $\varsigma_1 \in (0,1]$ is sufficiently small, the following estimates hold for
$(\ubar,u) \in [\leftubar,\ubarboot) \times [\moreinterestingu_1,\moreinterestingu_2]$, where the implicit constants are \textbf{independent} of $\varsigma_0$ and $\varsigma_1$:

	\begin{align}
		\begin{split} \label{E:IMPRECISETOPORDERESTIMATESFORFULLYMODIFIEDCHI}
			& \totalfluxcontrolimprecisefullymodifiedchi(\ubar,u)  +  (\blowupratetoporderwave + 1) \totalbulkcontrolimprecisefullymodifiedchi(\ubar,u) \lesssim  (1 + \varsigma_0^{-1}) \initialsmalldoublenull^2  \\
			& \ \ +  (1 + \varsigma_0^{-1} +  \varsigma_1^{-1}) \int_{\ubar' = \leftubar}^{\ubar}  \totalfluxcontrolimprecisefullymodifiedchi(\ubar',u) \, \mathrm{d} \ubar'  + \varsigma_1\fundbootsmall \max_{\angLie_{\tander}^N \in \angLie_{\mathfrak{P}}^{(N)}}   \int_{ \characteristicdiamondtwoarg{[\leftubar,\ubar)}{[\moreinterestingu_1,u]}}
			\left|\upmu \angLie_{\tander}^N \upchi \right|^2  \weight^{\blowupratetoporderwave} \voldiamond \\
			& \ \ + (1 + \varsigma_0^{-1}) \int_{\ubar' = \leftubar}^{\ubar} \totalfluxcontrolvelocity_{[1,N]}(\ubar',u) \, \mathrm{d} \ubar'\   + (1+ \varsigma_0^{-1}) \int_{u' = \moreinterestingu_1}^u  \totalfluxcontrolvelocity_{[1,N]}(\ubar,u') \, \rmd u' + \varsigma_1 (1 + \varsigma_0^{-1})\strongangularcontrolvelocity_{N}(\ubar,u) \\
			 & \ \ + \bulkcontrolVortVort_N(\ubar,u) +  \bulkcontrolDivGradEnt_N(\ubar,u) 
			  + \int_{u' = \moreinterestingu_1}^u
			\left\lbrace
				\fluxcontrolVortVort_{\le N-1}(\ubar,u')
				+
				\fluxcontrolDivGradEnt_{\le N-1}(\ubar,u')
			\right\rbrace
		\, \rmd u' \\
		& \ \ 
+   \int_{u' = \moreinterestingu_1}^u \left\{\fluxcontrolVort_{\le N}(\ubar,u')+ \fluxcontrolGradEnt_{\le N }(\ubar,u')\right\} \, \mathrm{d} u' \\. 
 		\end{split}
	\end{align}
\end{lemma}

\begin{proof}
The starting point of the proof is the identity \eqref{E:SPACETIMEENERGYNULLFLUXIDENTITYFORIMPRECISETOPORDERCHI}. Using the data estimates of Sect.\,\ref{SSS:QUANTITATIVEASSUMPTIONSONDATAAWAYFROMSYMMETRY}, the bound $|\mytr_{\gtorus}\upchi| \lesssim 1$, we find: 
\begin{align} 
	\begin{split} \label{E:PRELIMINARYTOPORDERMODIFIEDCHISTEP2}
			 \int_{\ingoingcharacteristicsurfacetwoarg{\ubar}{[\moreinterestingu_1,u]}}  
			 \left| \fullymodquant{\tander^N} \right|^2 \weight^{\blowupratetoporderwave + 1} \volingoingnullhypersurface  +(\blowupratetoporderwave + 1) \newspacetimecoercive{\blowupratetoporderwave}\left[ \fullymodquant{\tander^N}\right](\ubar,u)  & \lesssim  \initialsmalldoublenull^2 + \int_{ \characteristicdiamondtwoarg{[\leftubar,\ubar)}{[\moreinterestingu_1,u]}} \, \fullymodquant{\tander^N}\Lunit \left(\fullymodquant{\tander^N} \right)  \weight^{\blowupratetoporderwave + 1} \voldiamond \\
			 & + \int_{ \characteristicdiamondtwoarg{[\leftubar,\ubar)}{[\moreinterestingu_1,u]}}  \left| \fullymodquant{\tander^N} \right|^2 \weight^{\blowupratetoporderwave + 1} \voldiamond.
	\end{split} 
\end{align}
We claim that RHS\,\eqref{E:PRELIMINARYTOPORDERMODIFIEDCHISTEP2} is $\lesssim$ RHS\,\eqref{E:IMPRECISETOPORDERESTIMATESFORFULLYMODIFIEDCHI}, thereby concluding the proof upon taking the maximum over $\tander^N \in \mathfrak{P}^{(N)}$. 

We bound the first integral on RHS\,\eqref{E:PRELIMINARYTOPORDERMODIFIEDCHISTEP2} using \eqref{E:TRANSPORTEQUATIONFORFULLYMODIFIEDQUANTITY}. Using \eqref{E:RELATIONBETWEENCARTESIANNORMALIZEDNULLVECTORFIELDSANDEIKONALFUNCTIONORMALIZEDNULLVECTORFIELDS}, \eqref{E:KEYESTIMATECONTROLLINGINVERSEMUBYINVERSEWEIGHT}, \eqref{E:BOUNDSONNEWLMUINTERESTINGREGION}, and \eqref{E:COERCIVITYNEWSPACETIMETERM}, the error integral generated by the first term on RHS\,\eqref{E:TRANSPORTEQUATIONFORFULLYMODIFIEDQUANTITY} is: 
\begin{align} \label{E:PRELIMINARYTOPORDERMODIFIEDCHISTEP2.5}
	2 \int_{ \characteristicdiamondtwoarg{[\leftubar,\ubar)}{[\moreinterestingu_1,u]}} \frac{\Lunit \upmu}{\upmu} \left|  \fullymodquant{\tander^N}\right|^2   \weight^{\blowupratetoporderwave + 1} \ReciprocalLunitAppliedtoTimeFunction \voldiamond \le - \newspacetimecoercive{\blowupratetoporderwave}\left[ \fullymodquant{\tander^N}\right](\ubar,u).
\end{align}Using \eqref{E:POINTWISEESTIMATETANGENTDERIVATIVESOFMODQUANTINHOM}, we bound the error integral generated by the second term on RHS\,\eqref{E:TRANSPORTEQUATIONFORFULLYMODIFIEDQUANTITY} by:
	\begin{align}
		\begin{split} \label{E:PRELIMINARYTOPORDERMODIFIEDCHISTEP3}
			  \int_{ \characteristicdiamondtwoarg{[\leftubar,\ubar)}{[\moreinterestingu_1,u]}} \frac{1}{\upmu} \left| \fullymodquant{\tander^N} \right| & \left\{ \left|\newuL \tander^N\velocityarray\right| + \upmu  \left|\tander^{[1,N+1]} \velocityarray\right|
			+ 
			\left|\comdersmall^{[1,N];1}\velocityarray\right| \right. \\ 
			& \ \left. \qquad + \left| \tander^{\le N}(\vortrenormalized,\GradEnt)\right| 
			+
			\left|\tander^{[1,N]}\controlvars\right| 
			+ 
			\left|\tandersmall^{[1,N]}\badcontrolvars \right| \right\} \weight^{\blowupratetoporderwave + 1} \voldiamond.
		\end{split}
	\end{align}
Using Young's inequality, \eqref{E:KEYESTIMATECONTROLLINGINVERSEMUBYINVERSEWEIGHT}, and \eqref{E:COERCIVITYNEWSPACETIMETERM}, the first product on RHS\,\eqref{E:PRELIMINARYTOPORDERMODIFIEDCHISTEP3} is:
\begin{align}  \label{E:PRELIMINARYTOPORDERMODIFIEDCHISTEP3.5}
	\lesssim \varsigma_1 \newspacetimecoercive{\blowupratetoporderwave}\left[ \fullymodquant{\tander^N}\right](\ubar,u)  + \varsigma_1^{-1} \int_{\ubar' = \leftubar}^{\ubar} \totalfluxcontrolvelocity_{[1,N]}(\ubar',u) \, \mathrm{d} \ubar'.
\end{align}
Upon observing the cancellation of the $\upmu$'s, the second product on RHS\,\eqref{E:PRELIMINARYTOPORDERMODIFIEDCHISTEP3} can be bounded in magnitude by
$\lesssim \int_{ \characteristicdiamondtwoarg{[\leftubar,\ubar)}{[\moreinterestingu_1,u]}}  \left| \fullymodquant{\tander^N} \right|^2 \weight^{\blowupratetoporderwave + 1} \voldiamond  + \int_{ \characteristicdiamondtwoarg{[\leftubar,\ubar)}{[\moreinterestingu_1,u]}}   \upmu  \left|\tander^{[1,N+1]} \velocityarray\right|^2 \weight^{\blowupratetoporderwave} \voldiamond \lesssim$ RHS\,\eqref{E:IMPRECISETOPORDERESTIMATESFORFULLYMODIFIEDCHI}. Using the estimate $\weight^{\blowupratetoporderwave + 1} \lesssim  \weight^{\blowuprateofacousticgeoWRTnumberofcommutations(N-1)}$, the remaining error integrals from the third, fourth, and fifth products in \eqref{E:PRELIMINARYTOPORDERMODIFIEDCHISTEP3} are: 
 \begin{align}
	\begin{split}  \label{E:PRELIMINARYTOPORDERMODIFIEDCHISTEP4}
		&  \lesssim  \varsigma_1^{-1}\int_{ \characteristicdiamondtwoarg{[\leftubar,\ubar)}{[\moreinterestingu_1,u]}}  \left| \fullymodquant{\tander^N} \right|^2 \weight^{\blowupratetoporderwave + 1} \voldiamond +  \varsigma_1  \int_{ \characteristicdiamondtwoarg{[\leftubar,\ubar)}{[\moreinterestingu_1,u]}}   \upmu \left|\comdersmall^{[1,N];1}\velocityarray\right|^2  \weight^{\blowupratetoporderwave + 1} \voldiamond  \\
		 & \ \ +  \varsigma_1  \int_{\ubar' = \leftubar}^{\ubar} \left\{ \left\| \tandersmall^{[1,N]}\upmu\right\|_{L^2_{\blowuprateofacousticgeoWRTnumberofcommutations(N-1)} \left(\ingoingcharacteristicsurfacetwoarg{\ubar'}{[\moreinterestingu_1,u]}\right)}^2 + \sum_{i=1}^3\left\| \tander^{[1,N]}\Lunit^i \right\|_{L^2_{\blowuprateofacousticgeoWRTnumberofcommutations(N-1)} \left(\ingoingcharacteristicsurfacetwoarg{\ubar'}{[\moreinterestingu_1,u]}\right)}^2  \right\} (\ubar') \, \rmd \ubar'.
	\end{split}
\end{align}
Using \eqref{E:ARBITRARYCOMMUTATORSTRINGESTIMATE}, Lemma\,\ref{L:PRELIMINARYBELOWTOPORDERL2ESTIMATESFOREIKONALFUINCTIONQUANTITIES}, Cor.\,\ref{C:NONSINGULARL2ESTIMATESFORWAVEVARIABLESTHATLOSEONEDERIVATIVE}, the RHS\,\eqref{E:PRELIMINARYTOPORDERMODIFIEDCHISTEP4} is $\lesssim$ \eqref{E:IMPRECISETOPORDERESTIMATESFORFULLYMODIFIEDCHI}. We clarify that we use the fact that the RHS\,\eqref{E:PRELIMINARYBELOWTOPORDERL2ESTIMATESFOREIKONALFUINCTIONQUANTITIES}--\eqref{E:PRELIMINARYBELOWTOPORDERL2ESTIMATESFOREIKONALFUINCTIONQUANTITIESWITHZ} is increasing with respect to $(\ubar,u)$ as the last two error terms in RHS\,\eqref{E:PRELIMINARYTOPORDERMODIFIEDCHISTEP4} feature an additional $\ubar$-integration.

Next we bound the error integral generated by the $\upmu[\Lunit,\tander^N]\mytr_{\gtorus}\upchi$ term on RHS\,\eqref{E:TRANSPORTEQUATIONFORFULLYMODIFIEDQUANTITY} in \eqref{E:PRELIMINARYTOPORDERMODIFIEDCHISTEP2}. Using 
\eqref{E:COMMUTATOROFTANGENTIALANDTANGENTIALCOMMUTATORS},
Lemma\,\ref{L:CRUDEPOINTWISEESTIMATESFORTENSORFIELDS},
Prop.\,\ref{P:POINTWISETRANSPORTINEQUALITIESFOREIKFUNCTIONQUANTITIES},
Prop.\,\ref{P:IMPROVEMENTOFAUXILIARYBOOTSTRAP}, 
Cor\,\ref{C:IMPROVEAUX}, and Young's inequality,
we deduce that this error integral is: 
\begin{align}
	  \label{E:PRELIMINARYTOPORDERMODIFIEDCHISTEP5}
	& \lesssim (1 + \varsigma_1^{-1}) \int_{ \characteristicdiamondtwoarg{[\leftubar,\ubar)}{[\moreinterestingu_1,u]}}  \left| \fullymodquant{\tander^N} \right|^2 \weight^{\blowupratetoporderwave + 1} \voldiamond  +  \int_{ \characteristicdiamondtwoarg{[\leftubar,\ubar)}{[\moreinterestingu_1,u]}} \upmu \left| \comdersmall^{[1,N+1];1} \velocityarray \right|^2 \weight^{\blowupratetoporderwave}  \voldiamond \\	
		& \ \ +  \int_{\characteristicdiamondtwoarg{[\leftubar,\ubar)}{[\moreinterestingu_1,u]}} \left| \tander^{\le N}(\vortrenormalized,\GradEnt)\right|^2\weight^{\blowupratetoporderwave}  \voldiamond  
		+  \varsigma_1 \int_{ \characteristicdiamondtwoarg{[\leftubar,\ubar)}{[\moreinterestingu_1,u]}} \left\{
			\left|\tander^{[1,N]} \controlvars \right|^2
			+
		\left|\tandersmall^{[1,N]}\badcontrolvars \right|^2 \right\}\weight^{\blowupratetoporderwave}  \voldiamond 	 \label{E:PRELIMINARYTOPORDERMODIFIEDCHISTEP5.5} \\
	& \ \ + \max_{\angLie_{\tander}^N \in \angLie_{\mathfrak{P}}^{(N)}}  \varsigma_1 \fundbootsmall \int_{ \characteristicdiamondtwoarg{[\leftubar,\ubar)}{[\moreinterestingu_1,u]}}
		\left|\upmu \angLie_{\tander}^N \upchi \right|^2  \weight^{\blowupratetoporderwave} \voldiamond.\label{E:PRELIMINARYTOPORDERMODIFIEDCHISTEP6}
\end{align}
Using \eqref{E:ARBITRARYCOMMUTATORSTRINGESTIMATE}, Lemmas\,\ref{L:COERCIVITYOFSPACETIMEL2INTEGRALS}--\ref{L:COERCIVENESSOFL2CONTROLLINGQUANITIESFORVORTICITYANDENTROPY} and\,\ref{L:PRELIMINARYBELOWTOPORDERL2ESTIMATESFOREIKONALFUINCTIONQUANTITIES}, the second error integral in \eqref{E:PRELIMINARYTOPORDERMODIFIEDCHISTEP5} and the ones in \eqref{E:PRELIMINARYTOPORDERMODIFIEDCHISTEP5.5} are $\lesssim$ RHS\,\eqref{E:IMPRECISETOPORDERESTIMATESFORFULLYMODIFIEDCHI}.
The integral on \eqref{E:PRELIMINARYTOPORDERMODIFIEDCHISTEP6} is clearly bounded by the second term on RHS\,\eqref{E:IMPRECISETOPORDERESTIMATESFORFULLYMODIFIEDCHI}.

The same arguments given above imply that the terms on RHS\,\eqref{E:PRELIMINARYTOPORDERMODIFIEDCHISTEP2}
generated by the three commutator terms on the second line of RHS\,\eqref{E:TRANSPORTEQUATIONFORFULLYMODIFIEDQUANTITY} are $\lesssim$ \eqref{E:IMPRECISETOPORDERESTIMATESFORFULLYMODIFIEDCHI}. 

To handle the terms on RHS\,\eqref{E:PRELIMINARYTOPORDERMODIFIEDCHISTEP2}
arising from the term $\tander^N(\upmu |\upchi|_{\gtorus}^2)$ on RHS\,\eqref{E:TRANSPORTEQUATIONFORFULLYMODIFIEDQUANTITY},
we expand this term using the Leibniz rule for the operators $\angLie_{\tander}$.
Arguing as in the proof of \eqref{E:PRELIMINARYTOPORDERMODIFIEDCHISTEP6},
we find that:
\begin{align} \label{E:PRELIMINARYTOPORDERMODIFIEDCHISTEP7}
\left|
\tander^N(\upmu |\upchi|_{\gtorus}^2) 
\right| 
& 
\lesssim    
\fundbootsmall
\left|\upmu \angLie_{\tander}^N \upchi \right|_{\gtorus} 
+ 
\left|\tander^{[1,N]}\controlvars \right| 
+ 
\left|\tandersmall^{[1,N]}\badcontrolvars \right|.
\end{align}
This implies that the error integral generated from these terms is also $\lesssim$ \eqref{E:PRELIMINARYTOPORDERMODIFIEDCHISTEP5}--\eqref{E:PRELIMINARYTOPORDERMODIFIEDCHISTEP6}. 

The same arguments as above show that the all error integrals arising from the term $\tander^N \mathfrak{A}$ on RHS\,\eqref{E:TRANSPORTEQUATIONFORFULLYMODIFIEDQUANTITY}
except the $\upmu(\VortVort,\GradEnt)$ are $\lesssim$ \eqref{E:IMPRECISETOPORDERESTIMATESFORFULLYMODIFIEDCHI}. To control these, we easily bound:
\begin{align}
		\begin{split} \label{E:PRELIMINARYTOPORDERMODIFIEDCHISTEP8}
			&  \int_{ \characteristicdiamondtwoarg{[\leftubar,\ubar)}{[\moreinterestingu_1,u]}}  \left| \fullymodquant{\tander^N}\right| \left| \upmu \tander^{\le N}(\VortVort,\DivGradEnt)\right|   \weight^{\blowupratetoporderwave+1} \voldiamond \\
			& \ \  \lesssim  \int_{ \characteristicdiamondtwoarg{[\leftubar,\ubar)}{[\moreinterestingu_1,u]}}    \left| \fullymodquant{\tander^N}\right|^2  \weight^{\blowupratetoporderwave +1} \voldiamond +   \int_{ \characteristicdiamondtwoarg{[\leftubar,\ubar)}{[\moreinterestingu_1,u]}}\upmu \left| \tander^N (\VortVort,\DivGradEnt)\right|^2\weight^{\blowupratetoporderwave +1} \voldiamond \\
			& \ \ \ \ +  \int_{ \characteristicdiamondtwoarg{[\leftubar,\ubar)}{[\moreinterestingu_1,u]}}\upmu \left| \tander^{\le N-1} (\VortVort,\DivGradEnt)\right|^2\weight^{\blowupratetoporderwave +1} \voldiamond \\
			& \ \ \lesssim \varsigma_1^{-1}  \int_{\ubar' = \leftubar}^{\ubar}  \totalfluxcontrolimprecisefullymodifiedchi(\ubar',u) \, \mathrm{d} \ubar' +  \bulkcontrolVortVort_N(\ubar,u) +  \bulkcontrolDivGradEnt_N(\ubar,u) \\
			& \ \ \ \ + \int_{u' = \moreinterestingu_1}^u
			\left\lbrace
				\fluxcontrolVortVort_{\le N-1}(\ubar,u')
				+
				\fluxcontrolDivGradEnt_{\le N-1}(\ubar,u')
			\right\rbrace
		\, \rmd u',
		\end{split}
	\end{align}
where we used the trivial bound $\weight^{\blowupratetoporderwave+1} \lesssim \weight^{\blowuprateofacousticgeoWRTnumberofcommutations(N)-1}$ to bound $  \int_{ \characteristicdiamondtwoarg{[\leftubar,\ubar)}{[\moreinterestingu_1,u]}}\upmu \left| \tander^N (\VortVort,\DivGradEnt)\right|^2\weight^{\blowupratetoporderwave +1} \voldiamond \lesssim \bulkcontrolVortVort_N(\ubar,u) +  \bulkcontrolDivGradEnt_N(\ubar,u)$, see \eqref{E:RECIPROCALULUNITAPPLIEDTOTIMEFUNCTIONAPPROXIMATELYMU}, \eqref{E:COERCIVENESSOFNEWBULKVORTVORT}, and  \eqref{E:COERCIVENESSOFNEWBULKDIVGRADENT}.

To conclude the claim that LHS\,\eqref{E:PRELIMINARYTOPORDERMODIFIEDCHISTEP2} is $\lesssim$ RHS\,\eqref{E:IMPRECISETOPORDERESTIMATESFORFULLYMODIFIEDCHI}, we make $\varsigma_0$ sufficiently so that the $(-1 + C\varsigma_1) \newspacetimecoercive{\blowupratetoporderwave }\left[ \fullymodquant{\tander^N}\right](\ubar,u)$ integrals from \eqref{E:PRELIMINARYTOPORDERMODIFIEDCHISTEP2.5} and \eqref{E:PRELIMINARYTOPORDERMODIFIEDCHISTEP3.5} are $\le 0$.

\end{proof}

\subsubsection{Proof of Prop.\,\ref{P:TOPORDERIMPRECISEL2ESTIMATEMUCHI}}
\label{SSS:PROOFOFIMPRECISETOPORDERL2ESTIMATSFORCHI}

The estimate \eqref{E:TOPORDERL2ESTIMATEMUCHI}
for 
$
\max_{\tander^N \in \mathfrak{P}^{(N)}}  \int_{\characteristicdiamondtwoarg{[\leftubar,\ubar)}{[\moreinterestingu_1,u]}}  
			 \left| \upmu \tander^N \mytr_{\gtorus}\upchi \right|^2 \weight^{\blowupratetoporderwave} \voldiamond$ is merely a restatement of \eqref{E:PRELIMINARYTOPORDERTRACECHI}. The estimate \eqref{E:TOPORDERL2ESTIMATEMUCHI} for $\max_{\angLie_{\tander}^N \in \angLie_{\mathfrak{P}}^{(N)}} \int_{\characteristicdiamondtwoarg{[\leftubar,\ubar)}{[\moreinterestingu_1,u]}} 
			 \left| \upmu \angLie_{\tander}^N \upchi \right|_{\gtorus}^2 \weight^{\blowupratetoporderwave} \voldiamond$ follows from inserting \eqref{E:PRELIMINARYTOPORDERTRACECHI} into \eqref{E:PRELIMINARYTOPORDERLIECHI}. Finally, \eqref{E:FINALIMPRECISETOPORDERESTIMATESFORFULLYMODIFIEDCHI} follows from inserting the bound \eqref{E:TOPORDERL2ESTIMATEMUCHI} derived for the term  $\max_{\angLie_{\tander}^N \in \angLie_{\mathfrak{P}}^{(N)}} \int_{\characteristicdiamondtwoarg{[\leftubar,\ubar)}{[\moreinterestingu_1,u]}} 
			 \left| \upmu \angLie_{\tander}^N \upchi \right|_{\gtorus}^2 \weight^{\blowupratetoporderwave} \voldiamond$ into \eqref{E:IMPRECISETOPORDERESTIMATESFORFULLYMODIFIEDCHI}, and taking $\fundbootsmall$ sufficiently small to absorb the $\totalbulkcontrolimprecisefullymodifiedchi(\ubar,u)$ term from RHS\,\eqref{E:TOPORDERL2ESTIMATEMUCHI} into the LHS \eqref{E:IMPRECISETOPORDERESTIMATESFORFULLYMODIFIEDCHI}. 
\hfill $\qed$

\subsection{Estimates for the easy top-order eikonal function-involving error integrals}
\label{SS:ESTIMATESFOREASYTOPORDEREIKONALFUNCTIONERRORINTEGRALS}
In Lemma\,\ref{L:ESTIMATESFOREASYTOPORDEREIKONALFUNCTIONERRORINTEGRALS}, 
we use Prop.\,\ref{P:TOPORDERIMPRECISEL2ESTIMATEMUCHI} to control ``easy'' 
top-order wave equation error integrals
that are generated by the first product $\angrmd^{\sharp} v^i \cdot \upmu \angrmd \tanderY^{N-1} \mytr_{\gtorus}\upchi $
on RHS\,\eqref{E:TOPCOMMUTEDWAVELFIRSTTHENALLYS}
and the second product 
$\Yvfsmallcoeff{A} \angrmd^{\sharp} v_i \cdot \upmu \angrmd \tanderY^{N-1} \mytr_{\gtorus}\upchi$
on RHS\,\eqref{E:TOPCOMMUTEDWAVEALLYS}.
The corresponding error integrals are easy because of the helpful factor of $\upmu$.
In Sect.\,\ref{SS:ESTIMATESFORMOSTDIFFICULTEIKONALFUNCTIONERRORINTEGRALS}, we
will control the analogous -- but much more difficult -- error integral generated by the first product
$
(\muX v^i) \tanderY^{N-1} \Yvf{A} \mytr_{\gtorus} \upchi 
$
on RHS\,\eqref{E:TOPCOMMUTEDWAVEALLYS},
which lacks the factor of $\upmu$.

\begin{lemma}[Estimates for the easy top-order eikonal function-involving error integrals]
\label{L:ESTIMATESFOREASYTOPORDEREIKONALFUNCTIONERRORINTEGRALS}
Let $N = \Ntop$, let $v \in \velocityarray = \{v^1,v^2,v^3\}$,
and let $\tander^N \in \mathfrak{P}^{(N)}$, where $\mathfrak{P}^{(N)}$
is defined in Def.\,\ref{D:STRINGSOFCOMMUTATIONVECTORFIELDS}.
Then the following estimates hold for 
$(\ubar,u) \in [\leftubar,\ubarboot) \times [\moreinterestingu_1,\moreinterestingu_2]$:
\begin{align}
\begin{split}  \label{E:ESTIMATESFOREASYTOPORDEREIKONALFUNCTIONERRORINTEGRALS}
	& 
	 \int_{\characteristicdiamondtwoarg{[\leftubar,\ubar]}{[\moreinterestingu_1,u]}} 	
		\left| 
			\begin{pmatrix} 
				(1 + 2 \upmu) \Lunit \tander^N v
				\\ 
				2 \muX \tander^N v 
			\end{pmatrix} 
		\right|  
		\left| 
			\begin{pmatrix} 
				2 \upmu (\muX v) \tander^N \mytr_{\gtorus} \upchi \\
				(\angrmd^{\sharp} v) \cdot \upmu \angrmd \tanderY^{N-1} \mytr_{\gtorus} \upchi
				\\ 
				\Yvfsmallcoeff{A} \angrmd^{\sharp} v \cdot \upmu \angrmd \tanderY^{N-1}\mytr_{\gtorus}\upchi 
			\end{pmatrix} 
		\right|  \weight^{\blowupratetoporderwave} \ReciprocalLunitAppliedtoTimeFunction
	\, \voldiamond 
		\\
	& =
		\Errortoparg{N}(\ubar,u),
\end{split}
\end{align}
where $\Errortoparg{N}(\ubar,u)$ satisfies the estimate \eqref{E:ERRORTOPORDERWAVEESTIMATES}. 
\end{lemma}

\begin{proof}
	The bootstrap assumptions, \eqref{E:ANGDFPOINTWISEBOUNDEDBYCOMMUTATORVECTORFIELDS}, and \eqref{E:POINTWISESTIMATEFORBREVEXINTERMSOFNEWUL}
	imply that the integrand on LHS\,\eqref{E:ESTIMATESFOREASYTOPORDEREIKONALFUNCTIONERRORINTEGRALS} is
	pointwise bounded by $
	\lesssim
	\left(
		|\newuL \tander^N v|
		+
		|\Lunit \tander^N v| + \upmu |\angrmd \tander^N v|_{\gtorus}
	\right)
	\cdot
	\upmu |\tanderY^N \mytr_{\gtorus} \upchi|
	$.
	Hence, integrating over $\characteristicdiamondtwoarg{[\leftubar,\ubar]}{[\moreinterestingu_1,u]}$
	and using \eqref{E:COERCIVENESSOFNULLFLUXCONTROLWAVE}
	and Young's inequality,
	we bound LHS\,\eqref{E:ESTIMATESFOREASYTOPORDEREIKONALFUNCTIONERRORINTEGRALS} by:
	\begin{align}
		\lesssim \varsigma_1^{-1} \int_{\ubar' = \leftubar}^{\ubar}\totalfluxcontrolvelocity_N(\ubar',u) \, \rmd \ubar' + \varsigma_1^{-1} \int_{u' = \moreinterestingu_1}^u \totalfluxcontrolvelocity_N(\ubar',u) \, \rmd u' + \varsigma_1 \int_{\characteristicdiamondtwoarg{[\leftubar,\ubar)}{[\moreinterestingu_1,u]}}  
			 \left| \upmu \tander^N \mytr_{\gtorus}\upchi \right|^2 \weight^{\blowupratetoporderwave} \voldiamond. \label{E:ESTIMATESFOREASYTOPORDEREIKONALFUNCTIONERRORINTEGRALSSTEP1}
	\end{align}
	By \eqref{E:ERRORTOPORDERWAVEESTIMATES} and \eqref{E:TOPORDERL2ESTIMATEMUCHI}, we have finished the proof.
\end{proof}

\subsection{Estimates for the most difficult top-order eikonal function-involving error integrals}
\label{SS:ESTIMATESFORMOSTDIFFICULTEIKONALFUNCTIONERRORINTEGRALS}
In this section, we control the most difficult eikonal function-involving 
terms appearing in the commuted wave equations of Prop.\,\ref{P:MOSTDIFFICULTWAVETERMS}. 
More precisely, the most difficult product is the first one
$(\muX v^i) \tanderY^{N-1} \Yvf{A} \mytr_{\gtorus} \upchi$
on RHS\,\eqref{E:TOPCOMMUTEDWAVEALLYS}.
In our energy estimates, these terms are multiplied by $\multipliervectorfield \tanderY^N  v^i$, where 
$\multipliervectorfield$ is defined in \eqref{E:MULTIPLIERVECTORFIELD}.
This leads to the following difficult error integrals:
\begin{subequations}
\begin{align}
\int_{\characteristicdiamondtwoarg{[\leftubar,\ubar)}{[\moreinterestingu_1,u]}} 
	\left\lbrace 
		2 \muX \tanderY^N v^i
	\right\rbrace 
	\left\lbrace 
		(\muX v^i ) \tanderY^N \mytr_{\gtorus} \upchi 
	\right\rbrace \weight^{\blowupratetoporderwave} \ReciprocalLunitAppliedtoTimeFunction
	\, \voldiamond, 
		\label{E:OVERVIEWMUXPSIDIFFICULTERRORINTEGAL} 
			\\
	\int_{\characteristicdiamondtwoarg{[\leftubar,\ubar)}{[\moreinterestingu_1,u]}} 
		\left\lbrace 
			\Lunit \tanderY^N v^i
		\right\rbrace 
		\left\lbrace 
			(\muX v^i) \tanderY^N \mytr_{\gtorus}\upchi 
		\right\rbrace \weight^{\blowupratetoporderwave} \ReciprocalLunitAppliedtoTimeFunction
		\, \voldiamond.
		\label{E:OVERVIEWLPSIDIFFICULTERRORINTEGAL}
\end{align}
\end{subequations}
We bound the integral \eqref{E:OVERVIEWMUXPSIDIFFICULTERRORINTEGAL} 
in Sect.\,\ref{SSS:EIKONALTOPORDERENERGYESTIMATESWITHOUTIBP} 
and the integral \eqref{E:OVERVIEWLPSIDIFFICULTERRORINTEGAL},
which we control via a further integration by parts with respect to $\Lunit$,
in Sect.\,\ref{SSS:EIKONALTOPORDERENERGYESTIMATEWITHIBP}. 
It turns out that these two integrals are the main ones driving
top-order wave energy blowup-rate, i.e., power of the weight $ \weight^{\blowupratetoporderwave}$.
We provide the main estimates for these two integrals in Lemmas\,\ref{L:BOUNDSFORMOSTDIFFICULTWAVEERRORNTEGRALS}
and\,\ref{L:ESTIMATESFORDIFFICULTSPACETIMEERRORINTEGRALSINVOLVINGIBPWRTL}.

\subsubsection{Estimates that do not involve integration by parts}  
\label{SSS:EIKONALTOPORDERENERGYESTIMATESWITHOUTIBP}
We begin our analysis of the integral \eqref{E:OVERVIEWMUXPSIDIFFICULTERRORINTEGAL} 
with the following lemma, which provides precise $L^2$ a priori estimates for the ``precise'' fully-modified quantity $\ReciprocalLunitAppliedtoTimeFunction \sumofbrevexvisquared \, \fullymodquant{\tander^N}$ in terms of the controlling quantities. From the transport identity \eqref{E:SPACETIMEENERGYNULLFLUXIDENTITYFORFULLYMODIFIEDQUANTITY}, the most difficult terms will be the spacetime integral involving $\sumofbrevexvisquared^2 \, \fullymodquant{\tander^N} \Lunit \, \fullymodquant{\tander^N}$. The proof will show that the most difficult error integrals generated from $\sumofbrevexvisquared^2$ (see definition \eqref{E:FACTORINTOPACOUSTICENERGYNEEDEDFORSHARPCONSTANTS}) are those involving $(\muX v^1)^2  \, \fullymodquant{\tander^N} \Lunit \, \fullymodquant{\tander^N}$. The other integrals featuring $\muX v^2$ and $\muX v^3$ are much easier to handle due to the smallness bounds \eqref{E:LINFINITYIMPROVEMENTAUXTRANSVERSALPDERIVATIVESPARTIALWAVEARRAYSMALL}.

We also point out here that we will not estimate \eqref{E:OVERVIEWMUXPSIDIFFICULTERRORINTEGAL} individually; instead we bound the \emph{sum} of \eqref{E:OVERVIEWMUXPSIDIFFICULTERRORINTEGAL} over $i = 1,2,3$. We do this because the velocity controlling quantities of Def.\,\ref{D:MAINCOERCIVE} are defined with the sum over $i = 1,2,3$ and we will need to account for all three velocity simultaneously in order to absorb certain bulk integrals.

\begin{lemma}[Precise $L^2$ estimates for the top-order fully modified quantities]
\label{L:L2ESTIMATESFORMOSTDIFFICULTPRODUCT}
Let $N = \Ntop$,  let
$\tander^N \in \mathfrak{P}^{(N)}$ be the set of order $N$ $\ell_{t,u}$-tangential commutator operators from
Def.\,\ref{D:STRINGSOFCOMMUTATIONVECTORFIELDS}. Let $\blowuprateofwaveWRTnumberofcommutations(N) = \blowupratetoporderwave$ denote the blowup rate of the wave variables defined in \eqref{E:BLOWUPRATEOFWAVEWRTNUMBEROFCOMMUTATORS}. Let $\varsigma_0$ be the constant featured in Lemma\,\ref{L:PRELIMINARYBELOWTOPORDERL2ESTIMATESFOREIKONALFUINCTIONQUANTITIES} (see also Remark\,\ref{R:SMALLNESSCONSTANTFORBELOWTOPORDERACOUSTICVARIABLES}). Let $\sumofbrevexvisquared$ be the function from \eqref{E:FACTORINTOPACOUSTICENERGYNEEDEDFORSHARPCONSTANTS}. 
Then if $\varsigma_1 \in(0,1]$ is sufficiently small, the following estimates hold for
$(\ubar,u) \in [\leftubar,\ubarboot) \times [\moreinterestingu_1,\moreinterestingu_2]$, where the implicit constants are \textbf{independent} of $\varsigma_0$ and $\varsigma_1$:
\begin{align}
		\begin{split} \label{E:L2ESTIMATESFORMOSTDIFFICULTPRODUCT}
			& \int_{\ingoingcharacteristicsurfacetwoarg{\ubar}{[\moreinterestingu_1,u]}} \ReciprocalLunitAppliedtoTimeFunction^2\sumofbrevexvisquared^2 \left( \, \fullymodquant{\tander^N}\right)^2 \weight^{\blowupratetoporderwave} \, \volingoingnullhypersurface + \blowupratetoporderwave 
			\int_{\characteristicdiamondtwoarg{[\leftubar,\ubar)}{[\moreinterestingu_1,u]}} \ReciprocalLunitAppliedtoTimeFunction^2 \sumofbrevexvisquared^2 \left( \, \fullymodquant{\tander^N}\right)^2 \weight^{\blowupratetoporderwave - 1} \voldiamond \\
			 & \ \ \le \boxed{2.81} \weakspacetimetransvelocitycontrol_N(\ubar,u) + \boxed{2.33} \totalbulkcontrolprecisefullymodifiedchi(\ubar,u) \\
			 & \ \ + \Errortoparg{N}(\ubar,u).
		\end{split}
	\end{align}\end{lemma}
\begin{proof}
We begin with the identity \eqref{E:SPACETIMEENERGYNULLFLUXIDENTITYFORFULLYMODIFIEDQUANTITY}. From \eqref{E:LINFINITYIMPROVEMENTAUXWAVEARRAY}, we have that $| \sumofbrevexvisquared \Lunit \sumofbrevexvisquared |\lesssim \fundbootsmall$. Hence, using the data estimates of Sect.\,\ref{SSS:QUANTITATIVEASSUMPTIONSONDATAAWAYFROMSYMMETRY}, the bounds $|\mytr_{\gtorus}\upchi| \lesssim 1$, and $|\Lunit \ReciprocalLunitAppliedtoTimeFunction| = \left| \Lunit (\frac{1}{\Lunit \ubar})\right|  \lesssim 1$, (see \eqref{E:CLOSEDVERSIONLUBARESTIMATECHOVGEOTODOUBLENULL} and \eqref{E:CLOSEDVERSIONC11BOUNDFORCHOVDOUBLENULLTOGEO}), we have that RHS\,\eqref{E:SPACETIMEENERGYNULLFLUXIDENTITYFORFULLYMODIFIEDQUANTITY} is: 
\begin{align} 
	\begin{split} \label{E:L2ESTIMATESFORMOSTDIFFICULTPRODUCTSTEP1}
	& \le 2 \int_{ \characteristicdiamondtwoarg{[\leftubar,\ubar)}{[\moreinterestingu_1,u]}}  \ReciprocalLunitAppliedtoTimeFunction^2 \sumofbrevexvisquared ^2 \,  \fullymodquant{\tander^N}  \Lunit \left( \fullymodquant{\tander^N}\right) \weight^{\blowupratetoporderwave} \ReciprocalLunitAppliedtoTimeFunction \, \voldiamond  \\
	& \ \ + C \initialsmalldoublenull^2 + C \fundbootsmall \totalbulkcontrolimprecisefullymodifiedchi(\ubar,u) + C \int_{\ubar' = \leftubar}^{\ubar} \totalfluxcontrolprecisefullymodifiedchi(\ubar',u) \, \rmd \ubar'.
	\end{split}
\end{align}
The second line of \eqref{E:L2ESTIMATESFORMOSTDIFFICULTPRODUCTSTEP1} is clearly $ = \Errortoparg{N}(\ubar,u)$.

We now analyze the first integral on RHS\,\eqref{E:L2ESTIMATESFORMOSTDIFFICULTPRODUCTSTEP1} using \eqref{E:TRANSPORTEQUATIONFORFULLYMODIFIEDQUANTITY}. The error integral generated by the first term on RHS\,\eqref{E:TRANSPORTEQUATIONFORFULLYMODIFIEDQUANTITY} is precisely: 
\begin{align} \label{E:L2ESTIMATESFORMOSTDIFFICULTPRODUCTSTEP1.3}
		4 \int_{ \characteristicdiamondtwoarg{[\leftubar,\ubar)}{[\moreinterestingu_1,u]}} \frac{\Lunit \upmu}{\upmu} \ReciprocalLunitAppliedtoTimeFunction^2 \sumofbrevexvisquared^2  \left(\fullymodquant{\tander^N} \right)^2 \weight^{\blowupratetoporderwave} \ReciprocalLunitAppliedtoTimeFunction \, \voldiamond.
	\end{align}
	We decompose the spacetime integral in \eqref{E:L2ESTIMATESFORMOSTDIFFICULTPRODUCTSTEP1.3} using the identity:
\begin{align} \label{E:DIVIDINGTHEDOMAINSFORSHARPCONSTANTS}
		1
		= 
		\mathbf{1}_{\characteristicdiamondtwoarg{[\leftubar,\ubar)}{[\moreinterestingu_1,u]} \cap \smallneighborhoodofcreasearg{[\leftubar,\ubarboot]}}+
		\mathbf{1}_{\characteristicdiamondtwoarg{[\leftubar,\ubar)}{[\moreinterestingu_1,u]} \setminus \smallneighborhoodofcreasearg{[\leftubar,\ubarboot]}},
	\end{align}
where $\smallneighborhoodofcreasearg{[\leftubar,\ubarboot]}$ is the region from \eqref{E:SMALLNEIGHBORHOOD}. Using $\newL = \ReciprocalLunitAppliedtoTimeFunction \Lunit$ and \eqref{E:NEWLMUISALMOSTMINUSONEINSMALLNEIGHBORHOOD}, we have:
\begin{align}
		\begin{split} \label{E:L2ESTIMATESFORMOSTDIFFICULTPRODUCTSTEP1.6}
			& 4 \int_{ \characteristicdiamondtwoarg{[\leftubar,\ubar)}{[\moreinterestingu_1,u]} \cap \smallneighborhoodofcreasearg{[\leftubar,\ubarboot]}} 
		 \frac{\Lunit \upmu}{\upmu} \ReciprocalLunitAppliedtoTimeFunction^2 \sumofbrevexvisquared^2  \left(\fullymodquant{\tander^N} \right)^2 \weight^{\blowupratetoporderwave} \ReciprocalLunitAppliedtoTimeFunction \, \voldiamond  \\
		 	&   \ \ \le \boxed{-4(0.98)} \int_{ \characteristicdiamondtwoarg{[\leftubar,\ubar)}{[\moreinterestingu_1,u]} \cap \smallneighborhoodofcreasearg{[\leftubar,\ubarboot]}} 
		 \frac{1}{\upmu} \ReciprocalLunitAppliedtoTimeFunction^2 \sumofbrevexvisquared^2 \left(\fullymodquant{\tander^N} \right)^2 \weight^{\blowupratetoporderwave} \ReciprocalLunitAppliedtoTimeFunction \, \voldiamond .
		\end{split}
	\end{align}On the other hand, by applying first \eqref{E:OUTSIDEOFSMALLNEIGHBORHOODLOWERBOUNDFORMU} and then overestimating the integral over the entire characteristic diamond, we bound: 
\begin{align}
		\begin{split} \label{E:L2ESTIMATESFORMOSTDIFFICULTPRODUCTSTEP1.9}
			 4 \int_{ \characteristicdiamondtwoarg{[\leftubar,\ubar)}{[\moreinterestingu_1,u]} \setminus \smallneighborhoodofcreasearg{[\leftubar,\ubarboot]}} 
		 \frac{\Lunit \upmu}{\upmu} \ReciprocalLunitAppliedtoTimeFunction^2 \sumofbrevexvisquared^2  \left(\fullymodquant{\tander^N} \right)^2 \weight^{\blowupratetoporderwave} \ReciprocalLunitAppliedtoTimeFunction \, \voldiamond & \le  4 \int_{ \characteristicdiamondtwoarg{[\leftubar,\ubar)}{[\moreinterestingu_1,u]}}
		 |\Lunit \upmu| \ReciprocalLunitAppliedtoTimeFunction^2 \sumofbrevexvisquared^2 \left(\fullymodquant{\tander^N} \right)^2 \weight^{\blowupratetoporderwave} \ReciprocalLunitAppliedtoTimeFunction \, \voldiamond  \\
		 & \le \int_{\ubar' = \leftubar}^{\ubar} \totalfluxcontrolprecisefullymodifiedchi(\ubar',u) \, \rmd \ubar',
		\end{split}
	\end{align}
which concludes the analysis of the integral in \eqref{E:L2ESTIMATESFORMOSTDIFFICULTPRODUCTSTEP1.3}. 

We now focus on the error integral generated by the second term on RHS\,\eqref{E:TRANSPORTEQUATIONFORFULLYMODIFIEDQUANTITY} is: 
\begin{align}
	\begin{split} \label{E:L2ESTIMATESFORMOSTDIFFICULTPRODUCTSTEP2}
		-4 \int_{ \characteristicdiamondtwoarg{[\leftubar,\ubar)}{[\moreinterestingu_1,u]}} \frac{\newL \upmu}{\upmu} \ReciprocalLunitAppliedtoTimeFunction \sumofbrevexvisquared^2 \left(\fullymodquant{\tander^N} \right) (\tanderY^N \mathfrak{X}) \weight^{\blowupratetoporderwave} \ReciprocalLunitAppliedtoTimeFunction \, \voldiamond,
	\end{split}
\end{align}
where we used  $\newL = \ReciprocalLunitAppliedtoTimeFunction \Lunit$ to express $\ReciprocalLunitAppliedtoTimeFunction \Lunit \upmu = \ReciprocalLunitAppliedtoTimeFunction \newL \upmu$. Using $\sumofbrevexvisquared^2 = (\muX v^1)^2 + (\muX v^2)^2 + (\muX v^3)^2$, we focus first on the term generated by $(\muX v^1)^2$, namely:
\begin{align}
		\begin{split} \label{E:L2ESTIMATESFORMOSTDIFFICULTPRODUCTSTEP3}
			-4 \int_{ \characteristicdiamondtwoarg{[\leftubar,\ubar)}{[\moreinterestingu_1,u]}} \frac{\newL \upmu}{\upmu} \ReciprocalLunitAppliedtoTimeFunction  (\muX v^1)^2 \left(\fullymodquant{\tander^N} \right) (\tander^N \mathfrak{X}) \weight^{\blowupratetoporderwave} \ReciprocalLunitAppliedtoTimeFunction \, \voldiamond.
		\end{split}
	\end{align}By adding and subtracting $2(\Lunit \upmu) \muX \tander^N v^1$ in the integrand and applying pointwise estimate \eqref{E:POINTWISESUMOFMODIFIEDQUANTITYINHOMANDGLL}, it suffices to bound: 
\begin{align}
		\begin{split} \label{E:L2ESTIMATESFORMOSTDIFFICULTPRODUCTSTEP4}
			& 8 \int_{ \characteristicdiamondtwoarg{[\leftubar,\ubar)}{[\moreinterestingu_1,u]}} \left| \frac{(\newL \upmu)^2}{\upmu}   (\muX v^1) \left(\fullymodquant{\tander^N} \right) \muX \tander^N v^1 \right| \weight^{\blowupratetoporderwave} \ReciprocalLunitAppliedtoTimeFunction \, \voldiamond 
		\end{split} \\
		\begin{split} \label{E:L2ESTIMATESFORMOSTDIFFICULTPRODUCTSTEP4.5}
			&  \ \ + C  \int_{ \characteristicdiamondtwoarg{[\leftubar,\ubar)}{[\moreinterestingu_1,u]}} \left| \frac{\newuL \upmu}{\upmu} (\muX v^1)^2 \left(\fullymodquant{\tander^N} \right)\right| \cdot \Big\{   \fundbootsmall \left| \newuL \tander^{[1,N]} \velocityarray\right| + \fundbootsmall \upmu \left| \tander^{[1,N+1]}\velocityarray\right|    \\
			& \qquad \qquad  + \upmu \left| \tander^{\le N}(\vortrenormalized,\GradEnt)\right| +  \left| \tander^{[1,N]}\controlvars\right| +  \left|\tandersmall^{[1,N]}\badcontrolvars\right| \Big\} \weight^{\blowupratetoporderwave} \ReciprocalLunitAppliedtoTimeFunction \, \voldiamond. 
		\end{split}
	\end{align}We focus first on \eqref{E:L2ESTIMATESFORMOSTDIFFICULTPRODUCTSTEP4}. Dividing the characteristic diamond with \eqref{E:DIVIDINGTHEDOMAINSFORSHARPCONSTANTS} once again, using \eqref{E:POINTWISESTIMATEFORBREVEXINTERMSOFNEWUL}, Cor\,\ref{C:IMPROVEAUX}, \eqref{E:KEYESTIMATECONTROLLINGINVERSEMUBYINVERSEWEIGHT}, \eqref{E:NEWLMUISALMOSTMINUSONEINSMALLNEIGHBORHOOD}, \eqref{E:COERCIVENESSOFNULLFLUXCONTROLWAVE}, \eqref{E:COERCIVITYOLDSPACETIMETERM}, the trivial identity $(\muX v^1)^2 \le \sumofbrevexvisquared^2$, and Young's inequality $8 ab \le \frac{8}{3} a^2 + 6 b^2$, we have:
\begin{align}
	\begin{split} \label{E:L2ESTIMATESFORMOSTDIFFICULTPRODUCTSTEP5}
		& 8 \int_{ \characteristicdiamondtwoarg{[\leftubar,\ubar)}{[\moreinterestingu_1,u]}\cap \smallneighborhoodofcreasearg{[\leftubar,\ubarboot]}}  \frac{|\newL \upmu|^2 }{\upmu} \muX \tander^N v^1 (\muX v^1) \fullymodquant{\tander^N}  \weight^{\blowupratetoporderwave} \ReciprocalLunitAppliedtoTimeFunction \, \voldiamond \\ 
		& \le \boxed{2.81}  \int_{ \characteristicdiamondtwoarg{[\leftubar,\ubar)}{[\moreinterestingu_1,u]}} |\newuL \tander^N v^1|^2 \weight^{\blowupratetoporderwave - 1} \, \voldiamond + \boxed{6(1.02)^2}   \int_{\characteristicdiamondtwoarg{[\leftubar,\ubar)}{[\moreinterestingu_1,u]} \cap \smallneighborhoodofcreasearg{[\leftubar,\ubarboot]}}   \frac{1}{\upmu} \ReciprocalLunitAppliedtoTimeFunction^2 \sumofbrevexvisquared^2 \left( \, \fullymodquant{\tander^N}\right)^2 \weight^{\blowupratetoporderwave} \voldiamond \\
		& \ \ + C \int_{u' = \moreinterestingu_1}^u\totalfluxcontrolvelocity_N(\ubar,u') \, \rmd u' + C \fundbootsmall \strongangularcontrolvelocity_N(\ubar,u).
	\end{split}
\end{align} 
We clarify that in deriving \eqref{E:L2ESTIMATESFORMOSTDIFFICULTPRODUCTSTEP5}, we used \eqref{E:KEYESTIMATECONTROLLINGINVERSEMUBYINVERSEWEIGHT} to bound $\upmu^{-1} \le (1.01) \weight^{-1}$ and overestimated the integral featuring $\newuL \tander^N v^1$ over $\characteristicdiamondtwoarg{[\leftubar,\ubar)}{[\moreinterestingu_1,u]} \cap \smallneighborhoodofcreasearg{[\leftubar,\ubarboot]}$ to the whole characteristic diamond. From $\frac{8}{3}(1.02)^2(1.01) < 2.81$, where the two factors of $(1.02)$ came from \eqref{E:NEWLMUISALMOSTMINUSONEINSMALLNEIGHBORHOOD}, the first integral on RHS\,\eqref{E:L2ESTIMATESFORMOSTDIFFICULTPRODUCTSTEP5} is bounded by the first integral on RHS\,\eqref{E:L2ESTIMATESFORMOSTDIFFICULTPRODUCT}. Similarly, using the bound $6(1.02)^2 - 4(0.98) < 2.33$, the second integral in \eqref{E:L2ESTIMATESFORMOSTDIFFICULTPRODUCTSTEP5} plus the one in RHS\,\eqref{E:L2ESTIMATESFORMOSTDIFFICULTPRODUCTSTEP1.6} are easily bounded by the second integral on RHS\,\eqref{E:L2ESTIMATESFORMOSTDIFFICULTPRODUCT}. Using similar but strictly easier arguments (since we do not need to keep track of the constants) based on \eqref{E:OUTSIDEOFSMALLNEIGHBORHOODLOWERBOUNDFORMU}, the respective integrals over $\characteristicdiamondtwoarg{[\leftubar,\ubar)}{[\moreinterestingu_1,u]} \setminus \smallneighborhoodofcreasearg{[\leftubar,\ubarboot]}$  as in \eqref{E:L2ESTIMATESFORMOSTDIFFICULTPRODUCTSTEP5} are $ \lesssim \int_{\ubar' = \leftubar}^{\ubar} \totalfluxcontrolvelocity_N(\ubar',u) \, \rmd \ubar ' +  
 \int_{u' = \moreinterestingu_1}^u\totalfluxcontrolvelocity_N(\ubar,u') \, \rmd u' +  \int_{\ubar' = \leftubar}^{\ubar} \totalfluxcontrolprecisefullymodifiedchi(\ubar',u) \, \rmd \ubar' =  \Errortoparg{N}(\ubar,u) $. This concludes the analysis of the integral on \eqref{E:L2ESTIMATESFORMOSTDIFFICULTPRODUCTSTEP4}. 
 
 We now bound the error integrals in \eqref{E:L2ESTIMATESFORMOSTDIFFICULTPRODUCTSTEP4.5}. By \eqref{E:KEYESTIMATECONTROLLINGINVERSEMUBYINVERSEWEIGHT}, the first product is easily bounded by $\lesssim \fundbootsmall \weakspacetimetransvelocitycontrol_N(\ubar,u) + \fundbootsmall \totalbulkcontrolprecisefullymodifiedchi(\ubar,u) =   \Errortoparg{N}(\ubar,u) $. Similarly, by accounting for the cancelation of the factor of $\frac{1}{\upmu}$ by the terms $\fundbootsmall \upmu| \tander^{[1,N+1]}\velocityarray|$ in \eqref{E:L2ESTIMATESFORMOSTDIFFICULTPRODUCTSTEP4.5}, the corresponding error integral is bounded by $ \fundbootsmall\int_{\ubar' = \leftubar}^{\ubar} \totalfluxcontrolprecisefullymodifiedchi(\ubar',u) \, \rmd \ubar' 
  +  
\fundbootsmall  \int_{u' = \moreinterestingu_1}^u\totalfluxcontrolvelocity_N(\ubar,u') \, \rmd u' 
 + \fundbootsmall \strongangularcontrolvelocity_N(\ubar,u) =  \Errortoparg{N}(\ubar,u)$. Since the vorticity and entropy terms in \eqref{E:L2ESTIMATESFORMOSTDIFFICULTPRODUCTSTEP4.5} also feature a factor of $\upmu$, the contributing error integrals are similarly bounded by $\int_{\ubar' = \leftubar}^{\ubar} \totalfluxcontrolprecisefullymodifiedchi(\ubar',u) \, \rmd \ubar' 
  +  
 \int_{u' = \moreinterestingu_1}^u \left\{\fluxcontrolVort_{\le N}(\ubar,u') + \fluxcontrolGradEnt_{\le N}(\ubar,u') \right\} \, \rmd u =  \Errortoparg{N}(\ubar,u)$  using $\weight^{\blowupratetoporderwave} \lesssim \weight^{\blowuprateoftransportWRTnumberofcommutations(N)}$. Next, we control the most difficult error integrals from the last two terms in \eqref{E:L2ESTIMATESFORMOSTDIFFICULTPRODUCTSTEP4.5} using \eqref{E:KEYESTIMATECONTROLLINGINVERSEMUBYINVERSEWEIGHT}, and Young's inequality: 
	\begin{align}
		\begin{split} \label{E:L2ESTIMATESFORMOSTDIFFICULTPRODUCTSTEP6}
			&  \int_{ \characteristicdiamondtwoarg{[\leftubar,\ubar)}{[\moreinterestingu_1,u]}} \left| \frac{(\newL \upmu)^2}{\upmu}   (\muX v^1) \left(\fullymodquant{\tander^N} \right) \muX \tander^N v^1 \right| \left\{ \sum_{i=1}^3 \left| \tander^N L^i\right|+ \left| \tandersmall^N \upmu\right| \right\} \, \weight^{\blowupratetoporderwave} \voldiamond \\
			 & \ \ \lesssim \varsigma_1^{-1} \int_{ \characteristicdiamondtwoarg{[\leftubar,\ubar)}{[\moreinterestingu_1,u]}} \sumofbrevexvisquared^2 \left(\, \fullymodquant{\tander^N} \right)^2 \weight^{\blowupratetoporderwave} \voldiamond + \varsigma_1 \int_{ \characteristicdiamondtwoarg{[\leftubar,\ubar)}{[\moreinterestingu_1,u]}}  \left\{ \sum_{i=1}^3 \left| \tander^N L^i\right|^2+ \left| \tandersmall^N \upmu\right|^2 \right\} \, \weight^{\blowupratetoporderwave - 2} \voldiamond 
		\end{split}
	\end{align}
 
\noindent Using \eqref{E:PRELIMINARYJUSTBELOWTOPORDERL2SPACETIMEESTIMATESFORLSMALL}--\eqref{E:PRELIMINARYJUSTBELOWTOPORDERL2SPACETIMEESTIMATESFORMU}, and the fact that $\blowupratetoporderwave -2  =  \blowuprateofacousticgeoWRTnumberofcommutations(N-1)-1$, due to the factor of $\varsigma_1$ in the spacetime integrals featuring $\tander^N L^i$ and $\tandersmall \upmu$, the RHS\,\eqref{E:L2ESTIMATESFORMOSTDIFFICULTPRODUCTSTEP6} is $\lesssim  \Errortoparg{N}(\ubar,u)$. The error integrals for the rest of the terms $\tander^{[1,N-1]}\controlvars$ and $\tandersmall^{[1,N-1]}\badcontrolvars$ in \eqref{E:L2ESTIMATESFORMOSTDIFFICULTPRODUCTSTEP4.5} are bounded instead with \eqref{E:PRELIMINARYBELOWTOPORDERL2ESTIMATESFOREIKONALFUINCTIONQUANTITIES}--\eqref{E:PRELIMINARYBELOWTOPORDERL2ESTIMATESFOREIKONALFUINCTIONQUANTITIESWITHZ}. We omit the details.

This concludes the analysis of the integral in \eqref{E:L2ESTIMATESFORMOSTDIFFICULTPRODUCTSTEP3}, which accounts for the $(\muX v^1)^2$ terms in $\sumofbrevexvisquared^2 = (\muX v^1)^2 +  (\muX v^2)^2 +  (\muX v^3)^2$ in \eqref{E:L2ESTIMATESFORMOSTDIFFICULTPRODUCTSTEP2}. To control the remaining integrals, we use the fact that $|\muX v^2|, |\muX v^3| \lesssim \fundbootsmall$ (see \eqref{E:IMPROVEAUXWAVEARRAYPARTIALLINFINITY}) to bound:
\begin{align}
		\begin{split} \label{E:L2ESTIMATESFORMOSTDIFFICULTPRODUCTSTEP7}
			& -4 \int_{ \characteristicdiamondtwoarg{[\leftubar,\ubar)}{[\moreinterestingu_1,u]}} \frac{\newL \upmu}{\upmu} \ReciprocalLunitAppliedtoTimeFunction  \left((\muX v^2)^2 + (\muX v^3)^2\right) \left(\fullymodquant{\tander^N} \right) (\tanderY^N \mathfrak{X}) \weight^{\blowupratetoporderwave} \ReciprocalLunitAppliedtoTimeFunction \, \voldiamond  \\
			& \lesssim \fundbootsmall \int_{ \characteristicdiamondtwoarg{[\leftubar,\ubar)}{[\moreinterestingu_1,u]}} \frac{1}{\upmu}\sumofbrevexvisquared^2 \left| \, \fullymodquant{\tander^N}\right| \left|  \tander^N \mathfrak{X} \right| \weight^{\blowupratetoporderwave} \ReciprocalLunitAppliedtoTimeFunction \, \voldiamond.
		\end{split}
	\end{align}
Using the smallness factor of $\fundbootsmall$ in RHS\,\eqref{E:L2ESTIMATESFORMOSTDIFFICULTPRODUCTSTEP7}, strictly easier but very similar arguments used to bound \eqref{E:L2ESTIMATESFORMOSTDIFFICULTPRODUCTSTEP4}--\eqref{E:L2ESTIMATESFORMOSTDIFFICULTPRODUCTSTEP6} by $\lesssim  \Errortoparg{N}(\ubar,u)$ using  \eqref{E:POINTWISEESTIMATETANGENTDERIVATIVESOFMODQUANTINHOM} instead of \eqref{E:POINTWISESUMOFMODIFIEDQUANTITYINHOMANDGLL} prove that RHS\,\eqref{E:L2ESTIMATESFORMOSTDIFFICULTPRODUCTSTEP7} is also $\lesssim  \Errortoparg{N}(\ubar,u)$. We clarify that the arguments are strictly easier because the factor of $\fundbootsmall$ allows us to be use imprecise constants as opposed to, say, \eqref{E:L2ESTIMATESFORMOSTDIFFICULTPRODUCTSTEP5}.

To finish the proof that the first integral on RHS\,\eqref{E:L2ESTIMATESFORMOSTDIFFICULTPRODUCTSTEP1} is $\lesssim$ \eqref{E:L2ESTIMATESFORMOSTDIFFICULTPRODUCT}, it remains to consider the remaining terms on RHS\,\eqref{E:TRANSPORTEQUATIONFORFULLYMODIFIEDQUANTITY}. Using arguments similar to those used in the paragraphs preceeding the proof of \eqref{E:PRELIMINARYTOPORDERMODIFIEDCHISTEP5}--\eqref{E:PRELIMINARYTOPORDERMODIFIEDCHISTEP6}, we bound the $\upmu[\Lunit,\tander^N ]\mytr_{\gtorus}\upchi$ term on RHS\,\eqref{E:TRANSPORTEQUATIONFORFULLYMODIFIEDQUANTITY} in \eqref{E:L2ESTIMATESFORMOSTDIFFICULTPRODUCTSTEP1} by:
\begin{align}
	\begin{split} \label{E:L2ESTIMATESFORMOSTDIFFICULTPRODUCTSTEP13}
		\lesssim \max_{\angLie_{\tander}^N \in \angLie_{\mathfrak{P}}^{(N)}}  \fundbootsmall \int_{ \characteristicdiamondtwoarg{[\leftubar,\ubar)}{[\moreinterestingu_1,u]}}
		\left|\upmu \angLie_{\tander}^N \upchi \right|^2  \weight^{\blowupratetoporderwave} \voldiamond +  \Errortoparg{N}(\ubar,u).
	\end{split}
\end{align}
Estimate \eqref{E:TOPORDERL2ESTIMATEMUCHI} implies that RHS\,\eqref{E:L2ESTIMATESFORMOSTDIFFICULTPRODUCTSTEP13} is $=  \Errortoparg{N}(\ubar,u)$ thanks to the presence of the small $\fundbootsmall$ coefficient in \eqref{E:L2ESTIMATESFORMOSTDIFFICULTPRODUCTSTEP13}. Similarly, all of the error integrals generated by the commutator terms on RHS\,\eqref{E:TRANSPORTEQUATIONFORFULLYMODIFIEDQUANTITY} are $=  \Errortoparg{N}(\ubar,u)$; we omit the details. To handle the integrals arising from the term $\tander^N(\upmu |\upchi|_{\gtorus}^2)$ on RHS\,\eqref{E:TRANSPORTEQUATIONFORFULLYMODIFIEDQUANTITY}, we use \eqref{E:PRELIMINARYTOPORDERMODIFIEDCHISTEP7} and similar arguments used to show that RHS\,\eqref{E:L2ESTIMATESFORMOSTDIFFICULTPRODUCTSTEP13} is $= \Errortoparg{N}(\ubar,u)$.

From \eqref{E:L2ESTIMATESFORMOSTDIFFICULTPRODUCTSTEP5}, all of the contributing error integrals from inserting the $\tander^N \mathfrak{A}$ term on RHS\,\eqref{E:TRANSPORTEQUATIONFORFULLYMODIFIEDQUANTITY} to the $\Lunit \, \fullymodquant{\tander^N}$ factor in \eqref{E:L2ESTIMATESFORMOSTDIFFICULTPRODUCTSTEP1} have already been controlled except for $\upmu \tander^N(\VortVort,\DivGradEnt)$ and $\tander^{\le N-1}(\VortVort,\DivGradEnt)$. To control the top order terms, we use \eqref{E:RECIPROCALULUNITAPPLIEDTOTIMEFUNCTIONAPPROXIMATELYMU}, the definition of the bulk $L^2$ controlling quantities  \eqref{E:VORTVORTNEWSPACETIMEL2CONTROLLINGQUANTITY} and \eqref{E:DIVGRADENTNEWSPACETIMEL2CONTROLLINGQUANTITY}, multiply and divide by $\weight^{\frac{1}{4}}$ in the integrand, apply Young's inequality, and use $\blowupratetoporderwave + \frac{1}{2} = \blowuprateofacousticgeoWRTnumberofcommutations(N) - 1$ to bound:
\begin{align}
		\begin{split} \label{E:L2ESTIMATESFORMOSTDIFFICULTPRODUCTSTEP14}
			 & \int_{ \characteristicdiamondtwoarg{[\leftubar,\ubar)}{[\moreinterestingu_1,u]}} \ReciprocalLunitAppliedtoTimeFunction^2\sumofbrevexvisquared^2 \left| \, \fullymodquant{\tander^N}\right| \left|  \upmu \tander^N (\VortVort,\DivGradEnt)\right|  \weight^{\blowupratetoporderwave} \ReciprocalLunitAppliedtoTimeFunction \, \voldiamond \\
			 & \ \ \lesssim  \int_{ \characteristicdiamondtwoarg{[\leftubar,\ubar)}{[\moreinterestingu_1,u]}} \ReciprocalLunitAppliedtoTimeFunction^2\sumofbrevexvisquared^2 \left| \, \fullymodquant{\tander^N}\right|^2   \weight^{\blowupratetoporderwave- \frac{1}{2}}  \, \voldiamond 
			 +   \int_{ \characteristicdiamondtwoarg{[\leftubar,\ubar)}{[\moreinterestingu_1,u]}} \ReciprocaluLunitAppliedtoTimeFunction
			    \left|   \tander^N (\VortVort,\DivGradEnt)\right|^2   \weight^{\blowupratetoporderwave+ \frac{1}{2}} \ReciprocalLunitAppliedtoTimeFunction \, \voldiamond \\
			    & \ \ \lesssim   \int_{\ubar' = \leftubar}^{\ubar} \frac{1}{|\ubar'|^{0.5}}  \totalfluxcontrolprecisefullymodifiedchi(\ubar',u) \, \rmd \ubar' 
			    +\bulkcontrolVortVort_{N}(\ubar,u)
			    +
   			   \bulkcontrolDivGradEnt_{N}(\ubar,u) \\
			   & \ \ =  \Errortoparg{N}(\ubar,u).
		\end{split}
	\end{align}
\noindent To estimate the integrals featuring the lower order terms $\tander^{\le N-1}(\VortVort,\DivGradEnt)$, we use that $\blowupratetoporderwave > \blowuprateofacousticgeoWRTnumberofcommutations(N')$ for all integers $N' \le N-1$ to bound:
	\begin{align}
		\begin{split} \label{E:L2ESTIMATESFORMOSTDIFFICULTPRODUCTSTEP15}
			& \int_{ \characteristicdiamondtwoarg{[\leftubar,\ubar)}{[\moreinterestingu_1,u]}} \ReciprocalLunitAppliedtoTimeFunction^2\sumofbrevexvisquared^2 \left| \, \fullymodquant{\tander^N}\right| \left|   \tander^{\le N-1} (\VortVort,\DivGradEnt)\right|  \weight^{\blowupratetoporderwave} \ReciprocalLunitAppliedtoTimeFunction \, \voldiamond \\
			& \ \ 
			 \lesssim  \int_{\ubar' = \leftubar}^{\ubar} \totalfluxcontrolprecisefullymodifiedchi(\ubar',u) \, \rmd \ubar' 
			+ \int_{u' =\moreinterestingu_1}^{u} \left\{ \fluxcontrolVortVort_{\le N-1} (\ubar,u') +  \fluxcontrolDivGradEnt_{\le N-1} (\ubar,u') \right\} \, \rmd u' \\
			& \ \ =  \Errortoparg{N}(\ubar,u).
		\end{split}
	\end{align}

\end{proof}

\begin{lemma}[Bounds for the most difficult error integrals in the wave equation energy estimates] 
\label{L:BOUNDSFORMOSTDIFFICULTWAVEERRORNTEGRALS}
Let $N = \Ntop$,  let
$\mathfrak{Y}^{(N)}$ be the set of order $N$ $\ell_{t,u}$-tangential commutator operators from
Def.\,\ref{D:STRINGSOFCOMMUTATIONVECTORFIELDS}. Let $\blowuprateofwaveWRTnumberofcommutations(N) = \blowupratetoporderwave$ denote the blowup rate of the wave variables defined in \eqref{E:BLOWUPRATEOFWAVEWRTNUMBEROFCOMMUTATORS}. Let $\sumofbrevexvisquared$ be the function from \eqref{E:FACTORINTOPACOUSTICENERGYNEEDEDFORSHARPCONSTANTS}. Let $\varsigma_0$ be the constant featured in Lemma\,\ref{L:PRELIMINARYBELOWTOPORDERL2ESTIMATESFOREIKONALFUINCTIONQUANTITIES} (see also Remark\,\ref{R:SMALLNESSCONSTANTFORBELOWTOPORDERACOUSTICVARIABLES}).
Then if $\varsigma_1 \in(0,1]$ is sufficiently small, the following estimates hold for
$(\ubar,u) \in [\leftubar,\ubarboot) \times [\moreinterestingu_1,\moreinterestingu_2]$, where the implicit constants are \textbf{independent} of $\varsigma_0$ and $\varsigma_1$:
	\begin{align} 
		\begin{split} \label{E:SPACETIMEBOUNDSMOSTDIFFICULTWAVEPRODUCT}
			2 \left| \sum_{i=1}^3 \int_{ \characteristicdiamondtwoarg{[\leftubar,\ubar)}{[\moreinterestingu_1,u]}} (\muX \tanderY^N v^i) (\muX v^i) \tanderY^N \mytr_{\gtorus}\upchi \weight^{\blowupratetoporderwave} \ReciprocalLunitAppliedtoTimeFunction \voldiamond \right| & \le \boxed{5.18} \weakspacetimetransvelocitycontrol_N(\ubar,u) +  \boxed{1.01} \totalbulkcontrolprecisefullymodifiedchi(\ubar,u) \\
			& \ \ +   \Errortoparg{N}(\ubar,u) 
		\end{split}
	\end{align}
\end{lemma}

\begin{proof}
By using the definition \eqref{E:FULLYMODIFIEDQUANTITY} of the fully modified quantity, we express the integral in question as: 
	\begin{align}
		\begin{split} \label{E:SPACETIMEBOUNDSMOSTDIFFICULTWAVEPRODUCTSTEP1}
			2 \sum_{i=1}^3 \int_{ \characteristicdiamondtwoarg{[\leftubar,\ubar)}{[\moreinterestingu_1,u]}} \frac{1}{\upmu} (\muX \tanderY^N v^i) (\muX v^i) \fullymodquant{\tanderY^N} \weight^{\blowupratetoporderwave} \ReciprocalLunitAppliedtoTimeFunction \voldiamond  - 2 \sum_{i=1}^3 \int_{ \characteristicdiamondtwoarg{[\leftubar,\ubar)}{[\moreinterestingu_1,u]}} \frac{1}{\upmu} (\muX \tanderY^N v^i) (\muX v^i) \tanderY^N \mathfrak{X} \weight^{\blowupratetoporderwave} \ReciprocalLunitAppliedtoTimeFunction \voldiamond 
		\end{split}
	\end{align}
Using Young's inequality $ab \le \frac{1}{2} a^2 + \frac{1}{2} b^2$, \eqref{E:POINTWISESTIMATEFORBREVEXINTERMSOFNEWUL}, \eqref{E:KEYESTIMATECONTROLLINGINVERSEMUBYINVERSEWEIGHT},  Cor\,\ref{C:IMPROVEAUX}, \eqref{E:UNIFIEDCOMMUTATORVELOCITYSPACETIMEL2TRANSCONTROLLINGQUANTITY}, \eqref{E:PRECISEFULLYMODQUANTBULKTOTALCONTROL}, \eqref{E:COERCIVITYOLDSPACETIMETERM}, \eqref{E:COERCIVENESSOFNULLFLUXCONTROLWAVE}, \eqref{E:PRECISEFULLYMODQUANTBULKTOTALCONTROL}, and $ \sumofbrevexvisquared^2 = \sum_{i=1}^3 (\muX v^i)^2$ (see \eqref{E:FACTORINTOPACOUSTICENERGYNEEDEDFORSHARPCONSTANTS}), the first integral on RHS\,\eqref{E:SPACETIMEBOUNDSMOSTDIFFICULTWAVEPRODUCTSTEP1} may be bounded as: 
\begin{align}
		\begin{split} \label{E:SPACETIMEBOUNDSMOSTDIFFICULTWAVEPRODUCTSTEP2}
			 & \le  \boxed{1.01} \weakspacetimetransvelocitycontrol_N(\ubar,u) +  \boxed{1.01} \totalbulkcontrolprecisefullymodifiedchi(\ubar,u)  + C \int_{u' = \moreinterestingu_1}^u \totalfluxcontrolvelocity_N(\ubar,u') \, \rmd u' + C \fundbootsmall \strongangularcontrolvelocity_N(\ubar,u),
		\end{split}
	\end{align}
where we note that the last two integrals on RHS\,\eqref{E:SPACETIMEBOUNDSMOSTDIFFICULTWAVEPRODUCTSTEP2} are $=  \Errortoparg{N}(\ubar,u)$.
	
Next, using \eqref{E:POINTWISESUMOFMODIFIEDQUANTITYINHOMANDGLL}, we may bound the $i = 1$ summand in the second integral in RHS\,\eqref{E:SPACETIMEBOUNDSMOSTDIFFICULTWAVEPRODUCTSTEP1} as :	
	\begin{align}
		& \leq 4 \int_{ \characteristicdiamondtwoarg{[\leftubar,\ubar)}{[\moreinterestingu_1,u]}}  \frac{|\Lunit \upmu|}{\upmu} \left|\muX \tanderY^N v^1\right|^2
		\weight^{\blowupratetoporderwave} \ReciprocalLunitAppliedtoTimeFunction \, \voldiamond \label{E:SPACETIMEBOUNDSMOSTDIFFICULTWAVEPRODUCTSTEP3} \\
			\begin{split}
				& \ \ + C  \int_{ \characteristicdiamondtwoarg{[\leftubar,\ubar)}{[\moreinterestingu_1,u]}}    \frac{1}{\upmu} \left| \muX \tanderY^N v^1\right| 
		\left\{   \fundbootsmall \left| \newuL \tander^{[1,N]} \velocityarray\right| + \fundbootsmall \upmu \left| \tander^{[1,N+1]}\velocityarray\right| \right. \\
				& \qquad \qquad \qquad \qquad  \left.  + \upmu \left| \tander^{\le N}(\vortrenormalized,\GradEnt)\right| +  \left| \tander^{[1,N]}\controlvars\right| +  \left|\tandersmall^{[1,N]}\badcontrolvars\right| \right\} \weight^{\blowupratetoporderwave} \ReciprocalLunitAppliedtoTimeFunction \, \voldiamond. \label{E:SPACETIMEBOUNDSMOSTDIFFICULTWAVEPRODUCTSTEP4}
			\end{split}
	\end{align}
Due to the small $\fundbootsmall$ coefficients in the first two terms of \eqref{E:SPACETIMEBOUNDSMOSTDIFFICULTWAVEPRODUCTSTEP4},  arguments nearly identical to the ones used to bound the integral in \eqref{E:L2ESTIMATESFORMOSTDIFFICULTPRODUCTSTEP4.5} (e.g. see \eqref{E:L2ESTIMATESFORMOSTDIFFICULTPRODUCTSTEP6} and the surrounding discussion), imply that  \eqref{E:SPACETIMEBOUNDSMOSTDIFFICULTWAVEPRODUCTSTEP4}  $=  \Errortoparg{N}(\ubar,u)$. 

To  bound \eqref{E:SPACETIMEBOUNDSMOSTDIFFICULTWAVEPRODUCTSTEP3}, we consider \eqref{E:DIVIDINGTHEDOMAINSFORSHARPCONSTANTS}. Using \eqref{E:POINTWISESTIMATEFORBREVEXINTERMSOFNEWUL}, Cor\,\ref{C:IMPROVEAUX}, \eqref{E:KEYESTIMATECONTROLLINGINVERSEMUBYINVERSEWEIGHT}, $\newL = \ReciprocalLunitAppliedtoTimeFunction \Lunit$ with \eqref{E:NEWLMUISALMOSTMINUSONEINSMALLNEIGHBORHOOD}, and \eqref{E:COERCIVENESSOFNULLFLUXCONTROLWAVE}, we bound the integral resulting from the first term in \eqref{E:DIVIDINGTHEDOMAINSFORSHARPCONSTANTS} as:
\begin{align}
		\begin{split} \label{E:SPACETIMEBOUNDSMOSTDIFFICULTWAVEPRODUCTSTEP5}
			& 4 \int_{\characteristicdiamondtwoarg{[\leftubar,\ubar)}{[\moreinterestingu_1,u]} \cap \smallneighborhoodofcreasearg{[\leftubar,\ubarboot]}} \frac{|\newL \upmu|}{\upmu}  \left| \muX \tanderY^N v^1\right|^2 \weight^{\blowupratetoporderwave} \ReciprocalLunitAppliedtoTimeFunction \, \voldiamond \\
			& \ \ \le \boxed{4(1.04)} \weakspacetimetransvelocitycontrol_N(\ubar,u) +  C \int_{u' = \moreinterestingu_1}^u\totalfluxcontrolvelocity_N(\ubar,u') \, \rmd u' + C \fundbootsmall \strongangularcontrolvelocity_N(\ubar,u),
		\end{split}
	\end{align}
where we note that the last two integrals on RHS\,\eqref{E:SPACETIMEBOUNDSMOSTDIFFICULTWAVEPRODUCTSTEP2} are $=  \Errortoparg{N}(\ubar,u)$. 
We clarify that in deriving \eqref{E:SPACETIMEBOUNDSMOSTDIFFICULTWAVEPRODUCTSTEP5}, we used \eqref{E:KEYESTIMATECONTROLLINGINVERSEMUBYINVERSEWEIGHT} and \eqref{E:NEWLMUISALMOSTMINUSONEINSMALLNEIGHBORHOOD} to bound $|\newL \upmu| \mathbf{1}_{\characteristicdiamondtwoarg{[\leftubar,\ubar)}{[\moreinterestingu_1,u]} \cap \smallneighborhoodofcreasearg{[\leftubar,\ubarboot]}}  \upmu^{-1} \le (1.02) \upmu^{-1} \le 1.04 \weight^{-1}$. Notice that the first terms in \eqref{E:SPACETIMEBOUNDSMOSTDIFFICULTWAVEPRODUCTSTEP2} and \eqref{E:SPACETIMEBOUNDSMOSTDIFFICULTWAVEPRODUCTSTEP5} correspond to the first term on RHS\,\eqref{E:SPACETIMEBOUNDSMOSTDIFFICULTWAVEPRODUCT} since $1.01 + 4(1.04) < 5.18$.

Next, using \eqref{E:OUTSIDEOFSMALLNEIGHBORHOODLOWERBOUNDFORMU} in place of \eqref{E:NEWLMUISALMOSTMINUSONEINSMALLNEIGHBORHOOD} and estimate \eqref{E:L2ESTIMATESFORMUXTRANSVERSALDERIVATIVESOFWAVEVARIABLES}, the portion of the error integral \eqref{E:SPACETIMEBOUNDSMOSTDIFFICULTWAVEPRODUCTSTEP3} generated from the complement of the region in \eqref{E:SPACETIMEBOUNDSMOSTDIFFICULTWAVEPRODUCTSTEP5} is bounded by $ \int_{\ubar' = \leftubar}^{\ubar} \totalfluxcontrolvelocity_N(\ubar',u) \, \rmd \ubar' + C\int_{\moreinterestingu_1}^u \totalfluxcontrolvelocity_N(\ubar,u') \, \rmd u' =  \Errortoparg{N}(\ubar,u)$. 

This concludes the proof of the $i = 1$ summand in the second integral of \eqref{E:SPACETIMEBOUNDSMOSTDIFFICULTWAVEPRODUCTSTEP1}. For $i = 2,3$, we use the fact that $|\muX v^2|,|\muX v^3| \lesssim \fundbootsmall$ (see \eqref{E:LINFINITYIMPROVEMENTAUXTRANSVERSALPDERIVATIVESPARTIALWAVEARRAYSMALL}) and pointwise bound \eqref{E:POINTWISEESTIMATETANGENTDERIVATIVESOFMODQUANTINHOM} to estimate the integral by: 
\begin{align}
		\begin{split} \label{E:SPACETIMEBOUNDSMOSTDIFFICULTWAVEPRODUCTSTEP6}
			\fundbootsmall \int_{ \characteristicdiamondtwoarg{[\leftubar,\ubar)}{[\moreinterestingu_1,u]}} \frac{1}{\upmu} \left| \muX \tanderY^N v^i \right| & \left\{\left|\newuL \tander^N\velocityarray \right| + \upmu  \left|\tander^{[1,N+1]} \velocityarray\right|
			+ 
			\left|\comdersmall^{[1,N];1}\velocityarray \right| \right. \\
			& \ \  \left. + \left| \tander^{\le N}(\vortrenormalized,\GradEnt)\right| 
			+
			\left|\tander^{[1,N]}\controlvars\right| 
			+ 
			\left|\tandersmall^{[1,N]}\badcontrolvars \right| \right\} \weight^{\blowupratetoporderwave} \ReciprocalLunitAppliedtoTimeFunction \, \voldiamond 
		\end{split}
	\end{align}
Thanks to the small factor of $\fundbootsmall$ in \eqref{E:SPACETIMEBOUNDSMOSTDIFFICULTWAVEPRODUCTSTEP6} and pointwise estimates \eqref{E:POINTWISESTIMATEFORBREVEXINTERMSOFNEWUL}, arguments strictly easier yet nearly identical to the ones used to bound the integral in \eqref{E:L2ESTIMATESFORMOSTDIFFICULTPRODUCTSTEP4.5} (e.g. there is no need to multiply and divide by $\varsigma_1^{-1/2}$ before applying Young's inequality as in \eqref{E:L2ESTIMATESFORMOSTDIFFICULTPRODUCTSTEP6}) and the surrounding discussion imply that  \eqref{E:SPACETIMEBOUNDSMOSTDIFFICULTWAVEPRODUCTSTEP6}  $=  \Errortoparg{N}(\ubar,u)$.

\end{proof}

\subsubsection{Estimates involving integration by parts with respect to $\Lunit$}
\label{SSS:EIKONALTOPORDERENERGYESTIMATEWITHIBP}
In this section, we bound the difficult top-order integrals highlighted in \eqref{E:OVERVIEWLPSIDIFFICULTERRORINTEGAL}.
We first establish a series of preliminary lemmas before establishing the main result in 
Lemma\,\ref{L:MAINESTIMATESFORDIFFICULTSPACETIMEERRORINTEGRALSINVOLVINGIBPWRTL}.

\begin{lemma}[Precise $L^2$ estimates for the partially modified quantities]
\label{L:PRECISEL2ESTIMATEFORPARTIALLYMODIFIEDQUANTITIES}
Let $N = \Ntop$,  let 
$\tanderY^N \in \mathfrak{Y}^{(N)}$ be the set of order $N$ $\ell_{t,u}$-tangential commutator operators from
Def.\,\ref{D:STRINGSOFCOMMUTATIONVECTORFIELDS}. Let $\blowuprateofwaveWRTnumberofcommutations(N) = \blowupratetoporderwave$ denote the blowup rate of the wave variables defined in \eqref{E:BLOWUPRATEOFWAVEWRTNUMBEROFCOMMUTATORS}. Let $\varsigma_0$ be the constant featured in Lemma\,\ref{L:PRELIMINARYBELOWTOPORDERL2ESTIMATESFOREIKONALFUINCTIONQUANTITIES} (see also Remark\,\ref{R:SMALLNESSCONSTANTFORBELOWTOPORDERACOUSTICVARIABLES}). Let $\sumofbrevexvisquared$ be the function from \eqref{E:FACTORINTOPACOUSTICENERGYNEEDEDFORSHARPCONSTANTS}.
Then if $\varsigma_1 \in(0,1]$ is sufficiently small, the following estimates hold for
$(\ubar,u) \in [\leftubar,\ubarboot) \times [\moreinterestingu_1,\moreinterestingu_2]$, where the implicit constants are \textbf{independent} of $\varsigma_0$ and $\varsigma_1$:\begin{align}
	\begin{split} \label{E:PRECISEL2ESTIMATEFORPARTIALLYMODIFIEDQUANTITIES}
		& \left\| \sumofbrevexvisquared \, \partialmodquant{\tanderY^{N-1}} \right\|^2_{L^2_{\blowupratetoporderwave - 1}\left( \ingoingcharacteristicsurfacetwoarg{\ubar}{[\moreinterestingu_1,u]}\right)} +(\blowupratetoporderwave - 1)\newspacetimecoercive{\blowupratetoporderwave - 2}\left[ \sumofbrevexvisquared \, \partialmodquant{\tanderY^{N-1}}\right](\ubar,u) \\
		& \le \boxed{4.13}  \totalbulkcontrolpartialmodifiedchi(\ubar,u) + \boxed{1.02} \weakspacetimeangularvelocitycontrol_N(\ubar,u)  + \Errortoparg{N}(\ubar,u).
	\end{split}
\end{align}\end{lemma}

\begin{proof}
We begin with the identity \eqref{E:SPACETIMEENERGYNULLFLUXIDENTITYFORPARTIALLYMODIFIEDQUANTITY}. Using the data estimates of Sect.\,\ref{SSS:QUANTITATIVEASSUMPTIONSONDATAAWAYFROMSYMMETRY}, the $L^\infty$ estimates of Prop\,\ref{P:IMPROVEMENTOFAUXILIARYBOOTSTRAP}, the bounds $|\mytr_{\gtorus}\upchi| \lesssim 1$, $|\sumofbrevexvisquared (\Lunit \sumofbrevexvisquared)| \lesssim \fundbootsmall$ (see \eqref{E:LINFINITYIMPROVEMENTAUXWAVEARRAY}), and $|\Lunit \ReciprocalLunitAppliedtoTimeFunction| = \left| \Lunit (\frac{1}{\Lunit \ubar})\right|  \lesssim 1$, (see \eqref{E:CLOSEDVERSIONLUBARESTIMATECHOVGEOTODOUBLENULL} and \eqref{E:CLOSEDVERSIONC11BOUNDFORCHOVDOUBLENULLTOGEO}), we have that RHS\,\eqref{E:SPACETIMEENERGYNULLFLUXIDENTITYFORPARTIALLYMODIFIEDQUANTITY} is:
\begin{align} 
	\begin{split} \label{E:PRECISEL2ESTIMATEFORPARTIALLYMODIFIEDQUANTITIESSTEP1}
	& \le 2 \int_{ \characteristicdiamondtwoarg{[\leftubar,\ubar)}{[\moreinterestingu_1,u]}}    \sumofbrevexvisquared^2  \,   \partialmodquant{\tanderY^{N-1}} \Lunit \left( \partialmodquant{\tanderY^{N-1}}\right)\weight^{\blowupratetoporderwave - 1} \ReciprocalLunitAppliedtoTimeFunction \, \voldiamond  \\
	& \ \ + C \fundbootsmall  \int_{ \characteristicdiamondtwoarg{[\leftubar,\ubar)}{[\moreinterestingu_1,u]}} \left( \partialmodquant{\tanderY^{N-1}}\right)^2\weight^{\blowupratetoporderwave - 1}  \, \voldiamond \\
	& \ \ + C \initialsmalldoublenull^2 + C \int_{\ubar' = \leftubar}^{\ubar} \totalfluxcontrolpartialmodifiedchi(\ubar',u) \, \rmd \ubar',
	\end{split}
\end{align}
where we note that the last two terms on RHS\,\eqref{E:PRECISEL2ESTIMATEFORPARTIALLYMODIFIEDQUANTITIESSTEP1} are $= \Errortoparg{N}(\ubar,u)$. Using \eqref{E:PARTIALMODIFIEDQUANTITY}--\eqref{E:PARTIALMODIFIEDQUANTITYINHOMDEFINITION} and \eqref{E:PRELIMINARYBELOWTOPORDERL2ESTIMATESFOREIKONALFUINCTIONQUANTITIES}, we bound the second integral on RHS\,\eqref{E:PRECISEL2ESTIMATEFORPARTIALLYMODIFIEDQUANTITIESSTEP1} as:
\begin{align} 
	\begin{split} \label{E:PRECISEL2ESTIMATEFORPARTIALLYMODIFIEDQUANTITIESSTEP2}
		& \lesssim \fundbootsmall \int_{ \characteristicdiamondtwoarg{[\leftubar,\ubar)}{[\moreinterestingu_1,u]}} \left| \tanderY^{N-1} \mytr_{\gtorus}\upchi\right| \, \weight^{\blowupratetoporderwave - 1}  \, \voldiamond + \int_{ \characteristicdiamondtwoarg{[\leftubar,\ubar)}{[\moreinterestingu_1,u]}} \left( \left| \Lunit \tanderY^{N-1} \wavearray\right|^2 + \upmu \left| \angrmd \tanderY^{N-1} \wavearray\right|_{\gtorus}^2 \right) \weight^{\blowupratetoporderwave - 1}  \, \voldiamond \\
		& \lesssim \Errortoparg{N}(\ubar,u).
	\end{split}
\end{align}It remains to estimate the first term on RHS\,\eqref{E:PRECISEL2ESTIMATEFORPARTIALLYMODIFIEDQUANTITIESSTEP1}. Using \eqref{E:TRANSPORTEQUATIONFORPARTIALMODIFIEDQUANTITY} and \eqref{E:FACTORINTOPACOUSTICENERGYNEEDEDFORSHARPCONSTANTS}, this is integral is:
	\begin{align} 
		& =  \int_{ \characteristicdiamondtwoarg{[\leftubar,\ubar)}{[\moreinterestingu_1,u]}}   \left( (\muX v^1)^2 + (\muX v^2)^2 + (\muX v^3)^2\right)  \,   \partialmodquant{\tanderY^{N-1}} \left( \vec{G}_{\Lunit \Lunit}\diamond \angLap \tanderY^{N-1} \wavearray \right) \weight^{\blowupratetoporderwave - 1} \ReciprocalLunitAppliedtoTimeFunction \, \voldiamond  \label{E:PRECISEL2ESTIMATEFORPARTIALLYMODIFIEDQUANTITIESSTEP3} \\
		 & \ \ +   \int_{ \characteristicdiamondtwoarg{[\leftubar,\ubar)}{[\moreinterestingu_1,u]}} \sumofbrevexvisquared^2  \,   \partialmodquant{\tanderY^{N-1}} \left(^{(\tanderY^{N-1})}\mathfrak{B}\right) \weight^{\blowupratetoporderwave - 1} \ReciprocalLunitAppliedtoTimeFunction \, \voldiamond. \label{E:PRECISEL2ESTIMATEFORPARTIALLYMODIFIEDQUANTITIESSTEP4}
	\end{align}

We focus our attention first on the error integrals featuring only $(\muX v^1)^2$ in \eqref{E:PRECISEL2ESTIMATEFORPARTIALLYMODIFIEDQUANTITIESSTEP3}. Using the trilinear structure \eqref{E:TRILINEARSTRUCTUREFORPUREANGULARLMULTIPLIERTERM}, Cor.\,\ref{C:IMPROVEAUX}, and $\newL = \ReciprocalLunitAppliedtoTimeFunction \Lunit$, this corresponding integral is bounded in magnitude by: 
\begin{align}
			& \le \int_{ \characteristicdiamondtwoarg{[\leftubar,\ubar)}{[\moreinterestingu_1,u]}} \left| (\newL \upmu) (\muX v^1) \left( \, \partialmodquant{\tanderY^{N-1}}\right) \angLap \tanderY^{N-1} v^1 \right|  \weight^{\blowupratetoporderwave - 1} \voldiamond  \label{E:PRECISEL2ESTIMATEFORPARTIALLYMODIFIEDQUANTITIESSTEP5} \\
			& \ \ + C \int_{ \characteristicdiamondtwoarg{[\leftubar,\ubar)}{[\moreinterestingu_1,u]}} \left| (\muX v^1) \, \partialmodquant{\tanderY^{N-1}} \right| \left\{ \left| \Lunit \tander^{[1,N]}\velocityarray\right| +  \fundbootsmall \left| \tander^{[1,N+1]} \velocityarray\right|  + \left| \tander^{[1,N]}(\vortrenormalized,\GradEnt)\right| + \fundbootsmall \left|\tander^{[1,N]}\controlvars\right| \right\} 
			\weight^{\blowupratetoporderwave - 1}  \label{E:PRECISEL2ESTIMATEFORPARTIALLYMODIFIEDQUANTITIESSTEP6}
	\end{align}To control \eqref{E:PRECISEL2ESTIMATEFORPARTIALLYMODIFIEDQUANTITIESSTEP5}, we first decompose the characteristic diamond using \eqref{E:DIVIDINGTHEDOMAINSFORSHARPCONSTANTS}. On the first region featuring $\smallneighborhoodofcreasearg{[\leftubar,\ubarboot]}$, we use \eqref{E:SMOOTHANGULARHESSIANOFFPOINTWISEBOUNDEDBYCOMMUTATORVECTORFIELDS}, \eqref{E:POINTWISESMOOTHTORIDIFFERENTIALNORMINTERMSOFDOULBENULLTORIDIFFERENTIALINORMANDL}, the pointwise estimates of Prop\,\ref{P:IMPROVEMENTOFAUXILIARYBOOTSTRAP}, Cor\,\ref{C:IMPROVEAUX}, multiply and divide by $\sqrt{\upmu}$ alongside  \eqref{E:KEYESTIMATECONTROLLINGINVERSEMUBYINVERSEWEIGHT} in the $\nullangrmD \tanderY^N v^1$ terms, multiply and divide by $\sqrt{\weight}$ in the $\Lunit \tanderY^N \velocityarray$ terms, \eqref{E:NEWLMUISALMOSTMINUSONEINSMALLNEIGHBORHOOD}, and $|\muX v^1| \le \sumofbrevexvisquared$:
\begin{align} 
	\begin{split} \label{E:PRECISEL2ESTIMATEFORPARTIALLYMODIFIEDQUANTITIESSTEP7}
		& 2 \int_{ \characteristicdiamondtwoarg{[\leftubar,\ubar)}{[\moreinterestingu_1,u]}\cap \smallneighborhoodofcreasearg{[\leftubar,\ubarboot]}} \left|  (\newL \upmu)    (\muX v^1)  \left(\,\partialmodquant{\tanderY^{N-1}} \right) \left( \angLap \tanderY^{N-1}  v^1\right)\right| \weight^{\blowupratetoporderwave - 1}\voldiamond \\
		& \le \left(\boxed{2(1.02)} + \mathcal{O}(\mr \upalpha) +  \mathcal{O}(\fundbootsmall)\right)  \int_{ \characteristicdiamondtwoarg{[\leftubar,\ubar)}{[\moreinterestingu_1,u]}}  \sumofbrevexvisquared  \left| \frac{1}{\sqrt{\upmu}}\,\partialmodquant{\tanderY^{N-1}} \right| \cdot \sqrt{2\upmu \sum_{A=2}^3  \left|\nullangrmD \Yvf{A}\tanderY^{N-1} v^1\right|_{\gnulltori}^2} \,\weight^{\blowupratetoporderwave - 1}\voldiamond \\
		& \ \ + C\fundbootsmall \int_{ \characteristicdiamondtwoarg{[\leftubar,\ubar)}{[\moreinterestingu_1,u]}}   \sumofbrevexvisquared \left| \frac{1}{\sqrt{\weight}}\,\partialmodquant{\tanderY^{N-1}} \right|\cdot \left|\sqrt{\weight} \Lunit \tanderY^N \velocityarray\right| \weight^{\blowupratetoporderwave - 1}\voldiamond \\
		& \ \ + C \fundbootsmall \int_{ \characteristicdiamondtwoarg{[\leftubar,\ubar)}{[\moreinterestingu_1,u]}}  \sumofbrevexvisquared  \left| \,\partialmodquant{\tanderY^{N-1}} \right|\cdot \left| \angrmD \tanderY^{N-1} \velocityarray\right|_{\gtorus}  \weight^{\blowupratetoporderwave - 1}\voldiamond.
	\end{split} 
\end{align}

By applying Young's inequality, the last two integrals on RHS\,\eqref{E:PRECISEL2ESTIMATEFORPARTIALLYMODIFIEDQUANTITIESSTEP7} can be bounded as $\lesssim \fundbootsmall\totalbulkcontrolpartialmodifiedchi(\ubar,u) + \int_{\ubar'= \leftubar}^{\ubar} \totalfluxcontrolpartialmodifiedchi(\ubar',u) \, \rmd \ubar'  + \int_{u' = \moreinterestingu_1}^{u} \totalfluxcontrolvelocity_{[N-1,N]}(\ubar,u') \, \rmd u' = \Errortoparg{N}(\ubar,u)$. From \eqref{E:KEYESTIMATECONTROLLINGINVERSEMUBYINVERSEWEIGHT}, \eqref{E:WAVENEWSPACETIMENULLANGULARL2CONTROLLINGQUANTITY}, and Young's inequality $2 a(\sqrt{2}b) \le 4 a^2 + \frac{1}{2} b^2$, we bound the first integral on RHS\,\eqref{E:PRECISEL2ESTIMATEFORPARTIALLYMODIFIEDQUANTITIESSTEP7} by: 
\begin{align} 
	\begin{split} \label{E:PRECISEL2ESTIMATEFORPARTIALLYMODIFIEDQUANTITIESSTEP8}
		& \le  \left(\boxed{4.13} + \mathcal{O}(\mr \upalpha) + \mathcal{O}(\fundbootsmall)\right)  \int_{ \characteristicdiamondtwoarg{[\leftubar,\ubar)}{[\moreinterestingu_1,u]}}  \sumofbrevexvisquared^2  \left(\,   \partialmodquant{\tanderY^{N-1}}\right)^2\weight^{\blowupratetoporderwave-2} \, \voldiamond \\ 
		& \ \ +\left( \boxed{\frac{1.02}{2}} + \mathcal{O}(\mr \upalpha) + \mathcal{O}(\fundbootsmall)\right)  \int_{ \characteristicdiamondtwoarg{[\leftubar,\ubar)}{[\moreinterestingu_1,u]}} \upmu \sum_{A=2}^3 \left|\nullangrmD \Yvf{A}\tanderY^{N-1} v^1\right|_{\gnulltori}^2  \weight^{\blowupratetoporderwave - 1} \, \voldiamond  \\
	 	& \le \boxed{4(1.02)^2}  \totalbulkcontrolpartialmodifiedchi(\ubar,u) + \boxed{1.02} \weakspacetimeangularvelocitycontrol_N(\ubar,u) \\
		& \ \ + C(\mr \upalpha + \fundbootsmall) \totalbulkcontrolpartialmodifiedchi(\ubar,u) + C (\mr \upalpha + \fundbootsmall ) \weakspacetimeangularvelocitycontrol_N(\ubar,u).
	\end{split}
\end{align}
where we note the last two integrals on RHS\,\eqref{E:PRECISEL2ESTIMATEFORPARTIALLYMODIFIEDQUANTITIESSTEP7} are $= \Errortoparg{N}(\ubar,u)$ and the two contributing summands from $A = 2,3$ in the first line of \eqref{E:PRECISEL2ESTIMATEFORPARTIALLYMODIFIEDQUANTITIESSTEP8} contributed to a factor of 2 in the coefficient of $ \weakspacetimeangularvelocitycontrol_N(\ubar,u)$.

This concludes the estimate for \eqref{E:PRECISEL2ESTIMATEFORPARTIALLYMODIFIEDQUANTITIESSTEP5} in the first spacetime region defined by \eqref{E:DIVIDINGTHEDOMAINSFORSHARPCONSTANTS}. To control the second spacetime region, we use \eqref{E:OUTSIDEOFSMALLNEIGHBORHOODLOWERBOUNDFORMU} and arguments similar to those used to prove \eqref{E:PRECISEL2ESTIMATEFORPARTIALLYMODIFIEDQUANTITIESSTEP7} to see:
 \begin{align} 
	\begin{split} \label{E:PRECISEL2ESTIMATEFORPARTIALLYMODIFIEDQUANTITIESSTEP9}
		& 2 \int_{ \characteristicdiamondtwoarg{[\leftubar,\ubar)}{[\moreinterestingu_1,u]} \setminus \smallneighborhoodofcreasearg{[\leftubar,\ubarboot]}} \left| (\newL \upmu)    \sumofbrevexvisquared  \left(\,\partialmodquant{\tanderY^{N-1}} \right) \left( \angLap \tanderY^{N-1}  v^1\right)\right| \weight^{\blowupratetoporderwave - 1}\voldiamond \\
		& \lesssim   \int_{ \characteristicdiamondtwoarg{[\leftubar,\ubar)}{[\moreinterestingu_1,u]}} \sumofbrevexvisquared  \left| \,\partialmodquant{\tanderY^{N-1}} \right| \cdot  \left| \sqrt{\upmu}\nullangrmD \tanderY^{N} \wavearray\right|_{\gnulltori} \weight^{\blowupratetoporderwave - 1}\voldiamond  \\
		& \ \  +  \int_{ \characteristicdiamondtwoarg{[\leftubar,\ubar)}{[\moreinterestingu_1,u]}} \sumofbrevexvisquared   \left| \,\partialmodquant{\tanderY^{N-1}} \right| \cdot  \left| \Lunit \tanderY^{N} \velocityarray \right| \weight^{\blowupratetoporderwave - 1}\voldiamond  + \Errortoparg{N}(\ubar,u) \\
		& \lesssim \varsigma_1^{-1} \int_{\ubar' = \leftubar}^{\ubar} \totalfluxcontrolpartialmodifiedchi(\ubar',u) \, \rmd \ubar' + \varsigma_1 \weakspacetimeangularvelocitycontrol_N(\ubar,u)  + \varsigma_1 \totalbulkcontrolpartialmodifiedchi(\ubar,u) + \int_{u' = \moreinterestingu_1}^u \totalfluxcontrolvelocity_N(\ubar,u') \, \rmd u' + \Errortoparg{N}(\ubar,u) \\
		& =   \Errortoparg{N}(\ubar,u).
	\end{split} 
\end{align}This concludes the analysis of spacetime integral \eqref{E:PRECISEL2ESTIMATEFORPARTIALLYMODIFIEDQUANTITIESSTEP5}. We now bound \eqref{E:PRECISEL2ESTIMATEFORPARTIALLYMODIFIEDQUANTITIESSTEP6}. Similar to \eqref{E:PRECISEL2ESTIMATEFORPARTIALLYMODIFIEDQUANTITIESSTEP7} and \eqref{E:PRECISEL2ESTIMATEFORPARTIALLYMODIFIEDQUANTITIESSTEP9}, the first two error integrals are bounded by $\lesssim  \varsigma_1 \totalbulkcontrolpartialmodifiedchi(\ubar,u)  +   \int_{u' = \moreinterestingu_1}^u \totalfluxcontrolvelocity_N(\ubar,u') \, \rmd u' + \fundbootsmall  \strongangularcontrolvelocity_N(\ubar,u) = \Errortoparg{N}(\ubar,u)$. Since $\blowupratetoporderwave -1 = \blowuprateoftransportWRTnumberofcommutations(N)$, the third spacetime integrals are bounded using Young's inequality by 
$\int_{\ubar' = \leftubar}^{\ubar} \totalfluxcontrolpartialmodifiedchi(\ubar',u) \, \rmd \ubar' 
+ \int_{u' = \moreinterestingu_1}^u	\left\lbrace \fluxcontrolVort_{\le N}(\ubar,u')  + \fluxcontrolGradEnt_{\le N}(\ubar,u') \right\rbrace \, \rmd u' =  \Errortoparg{N}(\ubar,u)$. Finally, thanks to the small coefficient $\fundbootsmall$ in the last term on RHS\,\eqref{E:PRECISEL2ESTIMATEFORPARTIALLYMODIFIEDQUANTITIESSTEP6}, this integral is also $ =  \Errortoparg{N}(\ubar,u)$ using Lemma\,\ref{L:PRELIMINARYBELOWTOPORDERL2ESTIMATESFOREIKONALFUINCTIONQUANTITIES}. 

This concludes the proof of \eqref{E:PRECISEL2ESTIMATEFORPARTIALLYMODIFIEDQUANTITIESSTEP3} for the terms featuring $\muX v^1$. We now handle the other terms. Using \eqref{E:LINFINITYIMPROVEMENTAUXTRANSVERSALPDERIVATIVESPARTIALWAVEARRAYSMALL}, we bound this integral by: 
\begin{align} 
		\begin{split} \label{E:PRECISEL2ESTIMATEFORPARTIALLYMODIFIEDQUANTITIESSTEP10}
			&  \int_{ \characteristicdiamondtwoarg{[\leftubar,\ubar)}{[\moreinterestingu_1,u]}} \left|  \left( (\muX v^2)^2 + (\muX v^3)^2\right)  \left( \,   \partialmodquant{\tanderY^{N-1}} \right) \left( \vec{G}_{\Lunit \Lunit}\diamond \angLap \tanderY^{N-1} \wavearray \right) \right|  \weight^{\blowupratetoporderwave - 1} \ReciprocalLunitAppliedtoTimeFunction \, \voldiamond \\
			 & \ \ \lesssim  	 
			 \fundbootsmall \int_{ \characteristicdiamondtwoarg{[\leftubar,\ubar)}{[\moreinterestingu_1,u]}} \left|  \sumofbrevexvisquared \left( \,   \partialmodquant{\tanderY^{N-1}} \right) \left( \vec{G}_{\Lunit \Lunit}\diamond \angLap \tanderY^{N-1} \wavearray \right) \right|  \weight^{\blowupratetoporderwave - 1} \ReciprocalLunitAppliedtoTimeFunction \, \voldiamond.
		\end{split}
	\end{align}
By \eqref{E:SMOOTHANGULARLAPLACIANPOINTWISEBOUNDEDBYCOMMUTATORVECTORFIELDS} and \eqref{E:HIGHERORDERTANGENTIALDERIVATIVESOFDENSITY}, it follows that $\left| \vec{G}_{\Lunit \Lunit}\diamond \angLap \tanderY^{N-1} \wavearray \right| \lesssim \left| \tander^{[1,N+1]} \velocityarray\right| + \left| \tander^{\le N}(\vortrenormalized,\GradEnt)\right| + \fundbootsmall \left|\tander^{[1,N]} \controlvars\right|$. From this bound and small coefficient of $\fundbootsmall$ in \eqref{E:PRECISEL2ESTIMATEFORPARTIALLYMODIFIEDQUANTITIESSTEP9}, the corresponding integral was already estimated as part of \eqref{E:PRECISEL2ESTIMATEFORPARTIALLYMODIFIEDQUANTITIESSTEP6}. 

Finally, we return to \eqref{E:PRECISEL2ESTIMATEFORPARTIALLYMODIFIEDQUANTITIESSTEP4}. Using \eqref{E:POINTWISETANGENTDERIVATIVEOFBINHOM} and \eqref{E:PRELIMINARYBELOWTOPORDERL2ESTIMATESFOREIKONALFUINCTIONQUANTITIES}, we bound the integral by:
\begin{align} 
	\begin{split} \label{E:PRECISEL2ESTIMATEFORPARTIALLYMODIFIEDQUANTITIESSTEP11}
		\lesssim \varsigma^{-1} \int_{\ubar' = \leftubar}^{\ubar} \totalfluxcontrolpartialmodifiedchi(\ubar',u) \, \rmd \ubar' + \varsigma_1 \int_{ \characteristicdiamondtwoarg{[\leftubar,\ubar)}{[\moreinterestingu_1,u]}}  \left| \tander^{[1,N]} \controlvars\right|^2 \weight^{\blowupratetoporderwave - 1}\voldiamond \lesssim   \Errortoparg{N}(\ubar,u).
	\end{split}
\end{align}\end{proof}

\begin{lemma}[Estimates for difficult top-order error integrals related to integration by parts
 with respect to $\newL$] 
\label{L:ESTIMATESFORDIFFICULTSPACETIMEERRORINTEGRALSINVOLVINGIBPWRTL}
Let $N = \Ntop$, 
and let $\tanderY^N \in \mathfrak{Y}^{(N)}$, where $\mathfrak{Y}^{(N)}$
is defined in Def.\,\ref{D:STRINGSOFCOMMUTATIONVECTORFIELDS}. Let $\tanderY^{N-1} \in \mathfrak{Y}^{(N-1)}$
be such that $\tanderY^N = \Yvf{A} \tanderY^{N-1}$ for some $\Yvf{A} \in \Angularset$,
and let $\partialmodquant{\tanderY^{N-1}}$ be the corresponding partially modified quantity defined by
\eqref{E:PARTIALMODIFIEDQUANTITY}.
Then the following estimates hold for 
$(\ubar,u) \in [\leftubar,\ubarboot) \times [\moreinterestingu_1,\moreinterestingu_2]$:
	\begin{align} 
		\label{E:MAINSPACETIMEWAVEESTIMATEIBP}   \begin{split} 
			 \left| \int_{ \characteristicdiamondtwoarg{[\leftubar,\ubar)}{[\moreinterestingu_1,u]}} 
			(\Yvf{A} \tanderY^N v^1) 
			(\muX v^1) 
			\newL
			\partialmodquant{\tanderY^{N-1}} \, \weight^{\blowupratetoporderwave} \voldiamond \right| 
			& \le \boxed{2.07} \weakspacetimeangularvelocitycontrol_{N}(\ubar,u) +  \Errortoparg{N}(\ubar,u),
		\end{split} \\
		 \label{E:MAINHYPERSURFACEESTIMATEIBP} \begin{split} 
			 \left|  \sum_{i=1}^3 \int_{\ingoingcharacteristicsurfacetwoarg{\ubar}{[\moreinterestingu_1,u]}} 
			( \muX v^i )\left\{ \frac{1}{\Speed^2} (\gtorus^{-1})^{AB} \gtorusdoublenullCOV_B^C \nullgeop{x^C} \tanderY^N  v^i \right\}  \,\partialmodquant{\tanderY^{N-1}} \, \weight^{\blowupratetoporderwave}
			\, \volingoingnullhypersurface \right| 
			&  \le \boxed{\frac{2}{3}} \totalfluxcontrolvelocity_N(\ubar,u) + \boxed{\frac{3}{4}(1.01)}\totalfluxcontrolpartialmodifiedchi(\ubar,u) \\
			& \ \  + \Errortoparg{N}(\ubar,u).
		\end{split}		
	\end{align}
	
\noindent Moreover, for any $v \in \velocityarraypartial = \{ v^2,v^3\}$, we have the following estimates:
	\begin{align}
		 \label{E:MAINSPACETIMEPARTIALWAVEESTIMATEIBP}   \begin{split}
			& \left| \sum_{i=1}^3 \int_{ \characteristicdiamondtwoarg{[\leftubar,\ubar)}{[\moreinterestingu_1,u]}} 
			(\Yvf{A} \tanderY^N v^i) 
			(\muX v^i)  
			\newL
			\partialmodquant{\tanderY^{N-1}} \, \weight^{\blowupratetoporderwave} \voldiamond \right| \lesssim  \Errortoparg{N}(\ubar,u).
		\end{split}		
	\end{align}
	
\noindent Finally, we also have the following data estimate:
	\begin{align}
		\begin{split}
			 \label{E:MAINHYPERSURFACEDATAESTIMATEIBP} 
				& \left|  \int_{\ingoingcharacteristicsurfacetwoarg{\leftubar}{[\moreinterestingu_1,u]}} 
			\sum_{i=1}^3 \int_{\ingoingcharacteristicsurfacetwoarg{\ubar}{[\moreinterestingu_1,u]}} 
			( \muX v^i )\left\{ \frac{1}{\Speed^2} (\gtorus^{-1})^{AB} \gtorusdoublenullCOV_B^C \nullgeop{x^C} \tanderY^N  v^i \right\}  \,\partialmodquant{\tanderY^{N-1}}\, \weight^{\blowupratetoporderwave}
				\, \volingoingnullhypersurface \right| \\
				& \ \ \lesssim  \Errortoparg{N}(\ubar,u).
		\end{split}
	\end{align}\end{lemma}

\begin{proof} \hfill

\noindent \underline{\textbf{Proof of \eqref{E:MAINSPACETIMEWAVEESTIMATEIBP}:}}

Using \eqref{E:TRANSPORTEQUATIONFORPARTIALMODIFIEDQUANTITY} and $\newL = \ReciprocalLunitAppliedtoTimeFunction \Lunit$, 
the integral in question is:  
	\begin{align}
		& = \frac{1}{2} \int_{ \characteristicdiamondtwoarg{[\leftubar,\ubar)}{[\moreinterestingu_1,u]}} 
		(\Yvf{A} \tanderY^N v^1) 
		(\muX v^1)  \left( \vec{G}_{\Lunit \Lunit}\diamond \angLap \tanderY^{N-1} \wavearray \right) \weight^{\blowupratetoporderwave} \ReciprocalLunitAppliedtoTimeFunction \, \voldiamond  \label{E:MAINSPACETIMEWAVEESTIMATEIBPSTEP1} \\
		 & \ \ +  \sum_{i=1}^3  \int_{ \characteristicdiamondtwoarg{[\leftubar,\ubar)}{[\moreinterestingu_1,u]}} 
		(\Yvf{A} \tanderY^N v^i) 
		(\muX v^i)  \left(^{(\tanderY^{N-1})}\mathfrak{B}\right) \weight^{\blowupratetoporderwave} \ReciprocalLunitAppliedtoTimeFunction \, \voldiamond. \label{E:MAINSPACETIMEWAVEESTIMATEIBPSTEP2}			
	\end{align}To estimate \eqref{E:MAINSPACETIMEWAVEESTIMATEIBPSTEP1}, we consider the trilinear structure bound  \eqref{E:TRILINEARSTRUCTUREFORPUREANGULARLMULTIPLIERTERM}  to bound the integral by:
\begin{align}
			& \le \int_{ \characteristicdiamondtwoarg{[\leftubar,\ubar)}{[\moreinterestingu_1,u]} }  \left| (\newL \upmu)  \left(  \Yvf{A} \tanderY  v^1 \right) \left( \angLap \tanderY^{N-1} v^1\right) \right|  \weight^{\blowupratetoporderwave} \voldiamond  \label{E:MAINSPACETIMEWAVEESTIMATEIBPSTEP3} \\
			& \ \ + C \int_{ \characteristicdiamondtwoarg{[\leftubar,\ubar)}{[\moreinterestingu_1,u]}}  \left| (\newL \upmu)  \Yvf{A} \tanderY  v^1  \right|  \left\{ \left| \Lunit \tander^{[1,N]}\velocityarray\right| +  \fundbootsmall \left| \tander^{[1,N+1]} \velocityarray\right|  + \left| \tander^{[1,N]}(\vortrenormalized,\GradEnt)\right| + \fundbootsmall \left|\tander^{[1,N]}\controlvars\right| \right\} 
			\weight^{\blowupratetoporderwave}  \label{E:MAINSPACETIMEWAVEESTIMATEIBPSTEP4}
	\end{align}Dividing the integral over the characteristic diamond in \eqref{E:MAINSPACETIMEWAVEESTIMATEIBPSTEP3} with \eqref{E:DIVIDINGTHEDOMAINSFORSHARPCONSTANTS}, the portion of the integral in \eqref{E:MAINSPACETIMEWAVEESTIMATEIBPSTEP3} featuring $ \smallneighborhoodofcreasearg{[\leftubar,\ubarboot]}$ is bounded using 
\eqref{E:ANGDFPOINTWISEBOUNDEDBYCOMMUTATORVECTORFIELDS}--\eqref{E:SMOOTHANGULARLAPLACIANPOINTWISEBOUNDEDBYCOMMUTATORVECTORFIELDS}, \eqref{E:POINTWISESMOOTHTORIDIFFERENTIALNORMINTERMSOFDOULBENULLTORIDIFFERENTIALINORMANDL}, the pointwise estimates of Prop\,\ref{P:IMPROVEMENTOFAUXILIARYBOOTSTRAP}, Cor\,\ref{C:IMPROVEAUX}, multiplying and dividing by $\upmu$ alongside  \eqref{E:KEYESTIMATECONTROLLINGINVERSEMUBYINVERSEWEIGHT} in the $\nullangrmD \tanderY^N v^1$ terms, \eqref{E:NEWLMUISALMOSTMINUSONEINSMALLNEIGHBORHOOD}, and Young's inequality as follows:
\begin{align} 
	\begin{split} \label{E:MAINSPACETIMEWAVEESTIMATEIBPSTEP5}
		& \int_{ \characteristicdiamondtwoarg{[\leftubar,\ubar)}{[\moreinterestingu_1,u]} \cap \smallneighborhoodofcreasearg{[\leftubar,\ubarboot]}}  \left| (\newL \upmu)  \left(  \Yvf{A} \tanderY  v^1 \right) \left( \angLap \tanderY^{N-1} v^1\right) \right|  \weight^{\blowupratetoporderwave} \voldiamond  \\
		& \le \left(\boxed{1.02} + \mathcal{O}(\mr \upalpha)\right) \int_{ \characteristicdiamondtwoarg{[\leftubar,\ubar)}{[\moreinterestingu_1,u]}} \left| \angrmd \tanderY^N v^1 \right|_{\gtorus} \cdot \sqrt{2 \sum_{B=2}^3 \left| \angrmd \Yvf{B} \tanderY^{N-1} v^1\right|^2_{\gtorus} }\weight^{\blowupratetoporderwave}\voldiamond \\
		& \ \ + C\fundbootsmall  \int_{ \characteristicdiamondtwoarg{[\leftubar,\ubar)}{[\moreinterestingu_1,u]}} \left| \angrmd \tanderY^{[N-1,N]} v^1 \right|_{\gtorus}^2 \weight^{\blowupratetoporderwave}\voldiamond  \\
		& \le  \left(\boxed{2.07} + \mathcal{O}(\mr \upalpha) + \mathcal{O}(\fundbootsmall) \right) \max_{\tanderY^N \in \mathfrak{Y}^{(N)}}\int_{ \characteristicdiamondtwoarg{[\leftubar,\ubar)}{[\moreinterestingu_1,u]}}  \upmu  \left| \nullangrmD \tanderY^N v^1 \right|_{\gnulltori}^2 \weight^{\blowupratetoporderwave - 1}\voldiamond \\
		& \ \ + C \fundbootsmall \int_{ \characteristicdiamondtwoarg{[\leftubar,\ubar)}{[\moreinterestingu_1,u]}}  \left| \Lunit \tanderY^N v^1 \right|^2 \, \weight^{\blowupratetoporderwave}\voldiamond +  \Errortoparg{N}(\ubar,u) \\
		& \le \boxed{2.07} \weakspacetimeangularvelocitycontrol_{N}(\ubar,u) +   \Errortoparg{N}(\ubar,u).
	\end{split}
\end{align} Using 
\eqref{E:ANGDFPOINTWISEBOUNDEDBYCOMMUTATORVECTORFIELDS}--\eqref{E:SMOOTHANGULARHESSIANOFFPOINTWISEBOUNDEDBYCOMMUTATORVECTORFIELDS}, Cor\,\ref{C:IMPROVEAUX},   \eqref{E:OUTSIDEOFSMALLNEIGHBORHOODLOWERBOUNDFORMU}, and \eqref{E:COERCIVENESSOFNULLFLUXCONTROLWAVE}, we bound the integral from the remaining spacetime region of \eqref{E:DIVIDINGTHEDOMAINSFORSHARPCONSTANTS} by:
\begin{align} 
	\begin{split} \label{E:MAINSPACETIMEWAVEESTIMATEIBPSTEP6}
		&\int_{ \characteristicdiamondtwoarg{[\leftubar,\ubar)}{[\moreinterestingu_1,u]} \setminus \smallneighborhoodofcreasearg{[\leftubar,\ubarboot]}}   \left| (\newL \upmu)      \left((\Yvf{A} \tanderY^N v^1) \right) \left( \angLap \tanderY^{N-1}  v^1\right) \right|  \weight^{\blowupratetoporderwave}\voldiamond \\
		& \lesssim  \int_{ \characteristicdiamondtwoarg{[\leftubar,\ubar)}{[\moreinterestingu_1,u]}} \upmu \left|\angrmd \tanderY^{[N-1,N]} v^1\right|_{\gtorus} \,  \weight^{\blowupratetoporderwave}\voldiamond  = \Errortoparg{N}(\ubar,u).
	\end{split}
\end{align}Using strictly easier arguments used to prove \eqref{E:MAINSPACETIMEWAVEESTIMATEIBPSTEP5}--\eqref{E:MAINSPACETIMEWAVEESTIMATEIBPSTEP6}, the integrals in \eqref{E:MAINSPACETIMEWAVEESTIMATEIBPSTEP4} are all $= \Errortoparg{N}(\ubar,u)$. The same applies to the integral in \eqref{E:MAINSPACETIMEWAVEESTIMATEIBPSTEP2} after using  \eqref{E:POINTWISETANGENTDERIVATIVEOFBINHOM} and \eqref{E:PRELIMINARYBELOWTOPORDERL2ESTIMATESFOREIKONALFUINCTIONQUANTITIES}. We omit the details.

\noindent \underline{\textbf{Proof of \eqref{E:MAINHYPERSURFACEESTIMATEIBP} and \eqref{E:MAINHYPERSURFACEDATAESTIMATEIBP}:}}

We use \eqref{E:CHOVCOEFFICIENTSSMOOTHANGULARDERIVATIVESINTERMSOFDOUBLENULLONESANDL}, the bootstrap assumptions, estimate \eqref{E:SMALLC01ESTIMATESFORUBAR}, \eqref{E:NULLCOORDINATEPARTIALINTERMSOFMETRICNORM}, and Cor.\,\ref{C:IMPROVEAUX} to bound: 
	\begin{align} \label{E:MAINHYPERSURFACEESTIMATEIBPSTEP1}
		\left| \frac{1}{\Speed^2} (\gtorus^{-1})^{AB} \gtorusdoublenullCOV_B^C \nullgeop{x^C} \tanderY^N  v^i \right| \le \left| \frac{1}{\Speed^2} (\gtorus^{-1})^{AB}\nullgeop{x^C} \tanderY^N  v^i\right| + \left(\mathcal{O}(\mr\upalpha) + \mathcal{O}(\fundbootsmall)\right)\left| \nullangrmD \tanderY^N \velocityarray \right|_{\gnulltori}.
	\end{align}
Next, using \eqref{E:SMOOTHGINVERSEABEXPRESSION}, $\Xsmall^A = X^A$, \eqref{E:SCHEMATICSTRUCTUREOFXSMALL}, \eqref{E:NULLCOORDINATEPARTIALINTERMSOFMETRICNORM}, and the bootstrap assumptions, we have:
	\begin{align} \label{E:MAINHYPERSURFACEESTIMATEIBPSTEP2}
		\left| \frac{1}{\Speed^2} (\gtorus^{-1})^{AB}\nullgeop{x^C} \tanderY^N  v^i \right|  \le \left| \nullangrmD \tanderY^N v^i\right|_{\gnulltori} + \left(\mathcal{O}(\mr\upalpha) + \mathcal{O}(\fundbootsmall)\right)\left| \nullangrmD \velocityarray \right|_{\gnulltori}.
	\end{align}
Inserting \eqref{E:MAINHYPERSURFACEESTIMATEIBPSTEP1}--\eqref{E:MAINHYPERSURFACEESTIMATEIBPSTEP2} into the sum  in \eqref{E:MAINHYPERSURFACEESTIMATEIBP}, multiplying and dividing by $\sqrt{\upmu}$, \eqref{E:KEYESTIMATECONTROLLINGINVERSEMUBYINVERSEWEIGHT}, Young's inequality $ab \le \frac{1}{3} a^2 + \frac{3}{4} b^2$, and  we have the bound:  
	\begin{align}
		\begin{split} \label{E:MAINHYPERSURFACEESTIMATEIBPSTEP3}
			&  \left| \sum_{i=1}^3  \int_{\ingoingcharacteristicsurfacetwoarg{\ubar}{[\moreinterestingu_1,u]}} 
			 (\muX v^i)\left\{ \frac{1}{\Speed^2} (\gtorus^{-1})^{AB} \gtorusdoublenullCOV_B^C \nullgeop{x^C} \tanderY^N  v^i \right\} \,\partialmodquant{\belowtopordercomforpartialmodifiedchi^{N-1}} \, \weight^{\blowupratetoporderwave}
			\, \volingoingnullhypersurface \right| \\
			& \ \  \lesssim \boxed{\frac{1}{3}} \sum_{i=1}^3 \int_{\ingoingcharacteristicsurfacetwoarg{\ubar}{[\moreinterestingu_1,u]}}  \upmu \left| \nullangrmD \tanderY^N v^i\right|_{\gnulltori}^2 \weight^{\blowupratetoporderwave}
			\, \volingoingnullhypersurface + \boxed{\frac{3}{4}(1.01)} \sum_{i=1}^3 \int_{\ingoingcharacteristicsurfacetwoarg{\ubar}{[\moreinterestingu_1,u]}} (\muX v^i)^2 \left( \partialmodquant{\belowtopordercomforpartialmodifiedchi^{N-1}}\right)^2  \weight^{\blowupratetoporderwave-1}.
		\end{split}
	\end{align}
The desired estimate \eqref{E:MAINHYPERSURFACEESTIMATEIBP} then follows from 
the first inequality in \eqref{E:COERCIVENESSOFNULLFLUXCONTROLWAVE} (which we emphasize has a $\frac{1}{2}$, which explains the $\boxed{\frac{2}{3}}$ in the first term in \eqref{E:MAINHYPERSURFACEESTIMATEIBP}), \eqref{E:FACTORINTOPACOUSTICENERGYNEEDEDFORSHARPCONSTANTS}, and \eqref{E:PARTIALMODQUANTFLUXTOTALCONTROL}. Finally, \eqref{E:MAINHYPERSURFACEDATAESTIMATEIBP} follows from the same details but instead the data assumptions of Sect.\,\ref{SSS:QUANTITATIVEASSUMPTIONSONDATAAWAYFROMSYMMETRY}.

\noindent \underline{\textbf{Proof of \eqref{E:MAINSPACETIMEPARTIALWAVEESTIMATEIBP}:}}

The proof of estimate \eqref{E:MAINSPACETIMEPARTIALWAVEESTIMATEIBP} uses similar ideas as in the proofs of \eqref{E:MAINSPACETIMEWAVEESTIMATEIBP}--\eqref{E:MAINHYPERSURFACEESTIMATEIBP}, but the details are much easier due to the smallness $|\muX v^A| \lesssim \fundbootsmall$, see \eqref{E:LINFINITYIMPROVEMENTAUXTRANSVERSALPDERIVATIVESPARTIALWAVEARRAYSMALL}. We omit the details.

\end{proof}

\begin{lemma}[Estimates for easy error integrals that arise during integration by parts with respect to $\newL$]
	\label{L:IBPEASYSPACETIMEERRRORINTEGRALSESITMATES}
	Let $N = \Ntop$, 
	and let $\tanderY^N \in \mathfrak{Y}^{(N)}$, where $\mathfrak{Y}^{(N)}$
is defined in Def.\,\ref{D:STRINGSOFCOMMUTATIONVECTORFIELDS}. Let $\tanderY^{N-1} \in \mathfrak{Y}^{(N-1)}$
be such that $\tanderY^N = \Yvf{A} \tanderY^{N-1}$ for some $\Yvf{A} \in \Angularset$.
Let $\partialmodquant{\tanderY^{N-1}}$ be the partially modified quantity defined by
\eqref{E:PARTIALMODIFIEDQUANTITY},
and let
$\ErrorIBP^{AG;L}_1[\tander^N v; \, \partialmodquant{\tanderY^{N-1}}; \Yvf{A}]$ and $\ErrorIBP^{AG;L}_2[\tander^N v; \, \partialmodquant{\tanderY^{N-1}}; \Yvf{A}]$ 
	be the error terms defined in
	\eqref{E:ERRORTERM1KEYIBPIDENTIFYFORWAVEEQUATIONENERGYESTIMATES}--\eqref{E:ERRORTERM2KEYIBPIDENTIFYFORWAVEEQUATIONENERGYESTIMATES}
	(with $v$ in the role of $\varphi$ and $\partialmodquant{\tanderY^{N-1}}$ in the role of $\upeta$), respectively.
	Then the following estimate holds for 
	$(\ubar,u) \in [\leftubar,\ubarboot) \times [\moreinterestingu_1,\moreinterestingu_2]$:
	\begin{align} 
		\begin{split}
			& \left| \sum_{i=1}^3 \int_{\characteristicdiamondtwoarg{[\leftubar,\ubar)}{[\moreinterestingu_1,u]}}
			\ErrorIBP^{AG;L}_1[\tander^N v^i; \, \partialmodquant{\tanderY^{N-1}}; \Yvf{A}]  \,\voldiamond \right|, \, 
			\left| \sum_{i=1}^3 \int_{\characteristicdiamondtwoarg{[\leftubar,\ubar)}{[\moreinterestingu_1,u]}}
			\ErrorIBP^{AG;L}_2[\tander^N v^i; \, \partialmodquant{\tanderY^{N-1}}; \Yvf{A}] \,\voldiamond \right| \\
			 & \lesssim \Errortoparg{N}(\ubar,u).\label{E:IBPFIRSTEASYSPACETIMELERRRORINTEGRALESITMATE}
		\end{split}
	\end{align}
	
\end{lemma}

\begin{proof}
We first prove \eqref{E:IBPFIRSTEASYSPACETIMELERRRORINTEGRALESITMATE} for $\ErrorIBP_1^{AG;L}$. Using the identity 
$
\mytr_{\gtorus} \, \angdeform{\Yvf{A}}
= 
(\gtorus^{-1})^{\alpha \beta}
\Yvf{A} \gfour_{\alpha \beta}
+
2
\smoothtorusproject_{\kappa}^{\ \lambda}
\partial_{\lambda} \Yvf{A}^{\kappa}
$,
Cor\,\ref{C:ELLTUPROJECTEDVERSIONOFCARTESIANPARTIALDERIVATIVES},
Prop.\,\ref{P:SCHEMATICSTRUCTUREOFVARIOUSTENSORSINTERMSOFCONTROLVARS},
\eqref{E:POINTWISEBOUNDCOMMUTATORSTANGENTIALANDTANGENTIALDERIVATIVESONSCALARFUNCTION},\eqref{E:PARTIALMODIFIEDQUANTITY}, 
\eqref{E:POINTWISEBELOWTOPORDERPARTIALMODQUANTINHOM}, 
and the bootstrap assumptions,
we deduce the following pointwise bound:
	\begin{align} \label{E:IBPFIRSTEASYSPACETIMELERRRORINTEGRALESITMATESTEP1}
		\begin{split}
			\left| \ErrorIBP_1^{AG;L}[\tander^N v^i; \, \partialmodquant{\tanderY^{N-1}}; \Yvf{A}] \right| & \lesssim \left| \tander^{[1,N+1]} \velocityarray \right| \left| \partialmodquant{\tanderY^{N-1}}\right| \weight^{\blowupratetoporderwave}\\
			&  \lesssim \left| \tander^{[1,N+1]} \velocityarray \right| \left\{\left| \tanderY^{N-1}\mytr_{\gtorus}\upchi\right| + \left| \tander^{[1,N]}\velocityarray \right| \right\}\weight^{\blowupratetoporderwave} 
		\end{split}
	\end{align}
Multiplying and dividing by $\sqrt{\upmu}$ in RHS\,\eqref{E:IBPFIRSTEASYSPACETIMELERRRORINTEGRALESITMATESTEP1}, using \eqref{E:KEYESTIMATECONTROLLINGINVERSEMUBYINVERSEWEIGHT}, and applying Young's inequality, we have that RHS\,\eqref{E:IBPFIRSTEASYSPACETIMELERRRORINTEGRALESITMATESTEP1} is: 
	\begin{align} \label{E:IBPFIRSTEASYSPACETIMELERRRORINTEGRALESITMATESTEP2}
		\lesssim \varsigma_1^{-1} \upmu \left| \tander^{[1,N+1]}\velocityarray\right|^2 \weight^{\blowupratetoporderwave}+  \varsigma_1 \left| \tanderY^{N-1} \mytr_{\gtorus}\upchi\right| \weight^{\blowupratetoporderwave - 1}+ \upmu \left| \tander^{[1,N]}\velocityarray\right|^2 \weight^{\blowupratetoporderwave - 2}.
	\end{align}
	With the coercive bound \eqref{E:COERCIVENESSOFNULLFLUXCONTROLWAVE} and \eqref{E:PRELIMINARYBELOWTOPORDERL2ESTIMATESFOREIKONALFUINCTIONQUANTITIES}, the integral of RHS\,\eqref{E:IBPFIRSTEASYSPACETIMELERRRORINTEGRALESITMATESTEP2} over $\characteristicdiamondtwoarg{[\leftubar,\ubar)}{[\moreinterestingu_1,u]}$ is $\lesssim \Errortoparg{N}(\ubar,u)$.
	
We now prove \eqref{E:IBPFIRSTEASYSPACETIMELERRRORINTEGRALESITMATE} for $\ErrorIBP_2^{AG;L}$. Multiplying and dividing the first term in \eqref{E:ERRORTERM2KEYIBPIDENTIFYFORWAVEEQUATIONENERGYESTIMATES} by $\sqrt{\weight}$, using \eqref{E:POINTWISEBELOWTOPORDERPARTIALMODQUANTINHOM}, the bound $|\Yvf{A}\weight|\lesssim \fundbootsmall$ (which is a consequence of \eqref{E:Y2INTERMSOFGEOMETRICCOORDINATEVECTORFIELDS}--\eqref{E:Y3INTERMSOFGEOMETRICCOORDINATEVECTORFIELDS}, \eqref{E:SMALLC01ESTIMATESFORUBAR}, and Cor\,\ref{C:IMPROVEAUX}), the definition of $\sumofbrevexvisquared$ in \eqref{E:FACTORINTOPACOUSTICENERGYNEEDEDFORSHARPCONSTANTS}, and Young's inequality, it follows that: 
	\begin{align} \label{E:IBPFIRSTEASYSPACETIMELERRRORINTEGRALESITMATESTEP3}
		\left| \sum_{i=1}^3 \blowupratetoporderwave \, \weight^{\blowupratetoporderwave - 1} (\Yvf{A} \weight) (\muX v^i) (\newL \tander^N v^i) \partialmodquant{\tanderY^N} \right| \lesssim \left| \Lunit \tander^N \velocityarray \right|^2  \weight^{\blowupratetoporderwave} + \fundbootsmall \sumofbrevexvisquared^2 \left( \partialmodquant{\tanderY^N}\right)^2 \weight^{\blowupratetoporderwave - 2}.
	\end{align}
By \eqref{E:COERCIVENESSOFNULLFLUXCONTROLWAVE} and \eqref{E:PARTIALMODQUANTBULKTOTALCONTROL}, the integral of RHS\,\eqref{E:IBPFIRSTEASYSPACETIMELERRRORINTEGRALESITMATESTEP3} over $\characteristicdiamondtwoarg{[\leftubar,\ubar)}{[\moreinterestingu_1,u]}$ is $\lesssim \int_{u' = \moreinterestingu_1}^u \totalfluxcontrolvelocity_{N}(\ubar,u') \rmd u' + \fundbootsmall \totalbulkcontrolpartialmodifiedchi_{N}(\ubar,u) \lesssim \Errortoparg{N}(\ubar,u)$.

To handle the second term in \eqref{E:ERRORTERM2KEYIBPIDENTIFYFORWAVEEQUATIONENERGYESTIMATES}, we use \eqref{E:SMOOTHTORUSNORMCOMPARBLETOTANGENTIALCONTRACTIONS}, Cor.\,\ref{C:IMPROVEAUX}, multiply and divide by $\varsigma_1^{1/2}$, \eqref{E:FACTORINTOPACOUSTICENERGYNEEDEDFORSHARPCONSTANTS}, \eqref{E:DIVIDINGTHEDOMAINSFORSHARPCONSTANTS}, and Young's inequality to bound:
\begin{align} \label{E:IBPFIRSTEASYSPACETIMELERRRORINTEGRALESITMATESTEP4}
		\begin{split}
			& \left| \sum_{i=1}^3 \blowupratetoporderwave \, \weight^{\blowupratetoporderwave - 1} (\newL  \weight)(\muX v^i) (\Yvf{A} \tander^N v^i) \partialmodquant{\tanderY^N} \right| 
			\lesssim \varsigma_1 \left|\angrmD \tander^N \velocityarray \right|^2_{\gtorus} \weight^{\blowupratetoporderwave - 1} \mathbf{1}_{\characteristicdiamondtwoarg{[\leftubar,\ubar)}{[\moreinterestingu_1,u]} \cap \smallneighborhoodofcreasearg{[\leftubar,\ubarboot]}} \\
			& \ \ \ \ + \left|\angrmD \tander^N \velocityarray \right|^2_{\gtorus} \weight^{\blowupratetoporderwave - 1} \mathbf{1}_{\characteristicdiamondtwoarg{[\leftubar,\ubar)}{[\moreinterestingu_1,u]} \setminus \smallneighborhoodofcreasearg{[\leftubar,\ubarboot]}} + (1 +  \varsigma_1^{-1}) \sumofbrevexvisquared^2 \left( \partialmodquant{\tanderY^N}\right)^2\weight^{\blowupratetoporderwave - 1} 	\end{split}
	\end{align} 
	
\noindent Applying \eqref{E:OUTSIDEOFSMALLNEIGHBORHOODLOWERBOUNDFORMU} to the term featuring
$ \mathbf{1}_{\characteristicdiamondtwoarg{[\leftubar,\ubar)}{[\moreinterestingu_1,u]} \setminus \smallneighborhoodofcreasearg{[\leftubar,\ubarboot]}}$, the 
 integral of RHS\,\eqref{E:IBPFIRSTEASYSPACETIMELERRRORINTEGRALESITMATESTEP3} over the characteristic diamond is then $\lesssim \varsigma_1\strongangularcontrolvelocity_N(\ubar,u) +  \int_{u' = \moreinterestingu_1}^u \totalfluxcontrolvelocity_{N}(\ubar,u') \rmd u' + \varsigma_1^{-1}\int_{\ubar' = \leftubar}^{\ubar} \totalfluxcontrolpartialmodifiedchi(\ubar',u) \, \rmd \ubar' = \Errortoparg{N}(\ubar,u)$.

\end{proof}

We are now ready to derive the main estimates for the top-order error integrals highlighted 
in \eqref{E:OVERVIEWLPSIDIFFICULTERRORINTEGAL}.

\begin{lemma}[The main estimates for difficult top-order spacetime error integrals requiring integration by parts in $\newL$] 
\label{L:MAINESTIMATESFORDIFFICULTSPACETIMEERRORINTEGRALSINVOLVINGIBPWRTL}
Let $N = \Ntop$, 
and let $\tanderY^N \in \mathfrak{Y}^{(N)}$ be of the form $\tanderY^N = \Yvf{A} \tanderY^{N-1}$ for some $A = 2,3$, where $\tanderY^{N-1} \in \mathfrak{Y}^{(N-1)}$. 
Then the following estimate holds for 
$(\ubar,u) \in [\leftubar,\ubarboot) \times [\moreinterestingu_1,\moreinterestingu_2]$:
	\begin{align}
		\begin{split} \label{E:FINALMAINSPACETIMEWAVEESTIMATEIBP}
			& \left| \sum_{i=1}^3 \int_{\characteristicdiamondtwoarg{[\leftubar,\ubar)}{[\moreinterestingu_1,u]}} 
			 (\Lunit \tanderY^N v^i)
			(\muX v^i) \tanderY^N \mytr_{\gtorus}\upchi
			\weight^{\blowupratetoporderwave} \ReciprocalLunitAppliedtoTimeFunction
			\, \voldiamond\right| \\
			&  \le  \boxed{\frac{2}{3}} \totalfluxcontrolvelocity_N(\ubar,u) + \boxed{\frac{3}{4}(1.01)}\totalfluxcontrolpartialmodifiedchi(\ubar,u)  +  \boxed{2.07} \weakspacetimeangularvelocitycontrol_N(\ubar,u) \\
			& \ \ + \Errortoparg{N}(\ubar,u).
		\end{split}
	\end{align}
where $\Errortoparg{N}(\timefunction,u)$ satisfies \eqref{E:ERRORTOPORDERWAVEESTIMATES}.

\end{lemma}

\begin{proof}
We begin by using \eqref{E:PARTIALMODIFIEDQUANTITY} 
to decompose the factor $\tanderY^N \mytr_{\gtorus}\upchi$ on LHS~\eqref{E:FINALMAINSPACETIMEWAVEESTIMATEIBP}
as follows:
$\tanderY^N \mytr_{\gtorus}\upchi 
= 
\Yvf{A} 
\tanderY^{N-1} \mytr_{\gtorus} \upchi 
= 
\Yvf{A} \, \partialmodquant{\tanderY^{N-1}} 
- 
\Yvf{A} \, \partialmodquantinhom{\tanderY^{N-1}}$.
We insert this decomposition into LHS~\eqref{E:FINALMAINSPACETIMEWAVEESTIMATEIBP}
and will handle each of the two integrals separately, 
starting with the one generated by $\Yvf{A}\, \partialmodquantinhom{\tanderY^{N-1}}$,
which is easier. Specifically, 
the estimate \eqref{E:POINTWISELANDYDERIVATIVESOFTOPORDERPARTIALMODQUANTINHOM}, and the $L^\infty$ bounds of Sect.\,\ref{S:MAINLINFINITYESTIMATES} imply that the products
$
(\muX v^i)
\Yvf{A} \, \partialmodquantinhom{\tanderY^{N-1}}
$ is of type
$\HarmlessWave{N}$.
Hence, \eqref{E:ENERGYESTIMATEHARMLESSWAVEERRORTERMS}
implies that the corresponding integral
$- \int_{\characteristicdiamondtwoarg{[\leftubar,\ubar)}{[\moreinterestingu_1,u]}} 
	\left\lbrace 
		\Lunit \tanderY^N v^i
	\right\rbrace 
	\left\lbrace 
		(\muX v^i)  
		\Yvf{A} \partialmodquantinhom{\tanderY^{N-1}} 
	\right\rbrace  \weight^{\blowupratetoporderwave} \ReciprocalLunitAppliedtoTimeFunction
			\, \voldiamond
$
is of type $\Errortoparg{N}(\ubar,u)$ as desired.

To complete the proof of \eqref{E:FINALMAINSPACETIMEWAVEESTIMATEIBP},
it remains for us to bound the following spacetime integral:
\begin{align} \label{E:PROOFOFMAINL2ESTIMATESANNOYINGSPACETIMEINTEGRALINVOLVINGIBP}
	\sum_{i=1}^3
	\int_{\characteristicdiamondtwoarg{[\leftubar,\ubar)}{[\moreinterestingu_1,u]}} 
	\left\lbrace 
		\Lunit \tanderY^N v^i
	\right\rbrace 
	\left\lbrace 
		(\muX v^i) 
		\Yvf{A} \, \partialmodquant{\tanderY^{N-1}} 
	\right\rbrace \weight^{\blowupratetoporderwave} \ReciprocalLunitAppliedtoTimeFunction
			\, \voldiamond
\end{align}
To bound \eqref{E:PROOFOFMAINL2ESTIMATESANNOYINGSPACETIMEINTEGRALINVOLVINGIBP}, we use $\newL = \ReciprocalLunitAppliedtoTimeFunction\Lunit$, integrate by parts using the identity \eqref{E:KEYIBPIDENTIFYFORWAVEEQUATIONENERGYESTIMATES}
with $v^i$ in the role of $\varphi$, $\partialmodquant{\tanderY^{N-1}}$ in the role of $\upeta$,
 $\tanderY^N$ in the role of $\tander^N$, and sum over $i = 1,2,3$.
The first two integrals on RHS\,\eqref{E:KEYIBPIDENTIFYFORWAVEEQUATIONENERGYESTIMATES},
namely:
\begin{align*}
&
 \sum_{i=1}^3  \int_{\characteristicdiamondtwoarg{[\leftubar,\ubar)}{[\moreinterestingu_1,u]}} 
		(\muX v^i) 
		(\Yvf{A} \tanderY^N v^i) 
		\newL
		\partialmodquant{\tanderY^{N-1}}  \weight^{\blowupratetoporderwave} \voldiamond
		\\
&
- 
\int_{\ingoingcharacteristicsurfacetwoarg{\leftubar}{[\moreinterestingu_1,u]}} 
				 ( \muX v^i) \left\{ \frac{1}{\Speed^2} (\gtorus^{-1})^{AB} \gtorusdoublenullCOV_B^C \nullgeop{x^C}\tander^N v^i \right\}  \partialmodquant{\tanderY^{N-1}} \, \weight^{\blowupratetoporderwave}
				\, \volingoingnullhypersurface
\end{align*}
are the main ones,
and in \eqref{E:MAINSPACETIMEWAVEESTIMATEIBP}--\eqref{E:MAINHYPERSURFACEESTIMATEIBP},
we showed that they are bounded by $\leq \mbox{RHS\,\eqref{E:FINALMAINSPACETIMEWAVEESTIMATEIBP}}$.
The third integral on RHS\,\eqref{E:KEYIBPIDENTIFYFORWAVEEQUATIONENERGYESTIMATES}, namely
$  \sum_{i=1}^3 \int_{\ingoingcharacteristicsurfacetwoarg{\leftubar}{[\moreinterestingu_1,u]}} 
				 ( \muX v^i) \left\{ \frac{1}{\Speed^2} (\gtorus^{-1})^{AB} \gtorusdoublenullCOV_B^C \nullgeop{x^C}\tander^N v^i \right\}  \partialmodquant{\tanderY^{N-1}} \, \weight^{\blowupratetoporderwave}
				\, \volingoingnullhypersurface
$,
was shown to be of type $\Errortop(\ubar,u)$
in Lemma\,\ref{L:ESTIMATESFORDIFFICULTSPACETIMEERRORINTEGRALSINVOLVINGIBPWRTL}.
The last integral on RHS\,\eqref{E:KEYIBPIDENTIFYFORWAVEEQUATIONENERGYESTIMATES},
whose integrand is of type 
$ \ErrorIBP_1^{AG;L}[\tander^Nv^i ;\partialmodquant{\tanderY^{N-1}}] + \ErrorIBP_2^{AG;L}[\tander^N v^i;\, \partialmodquant{\tanderY^{N-1}}] +   \ErrorIBP_2[\tander^N\varphi;\partialmodquant{\tanderY^{N-1}}]$
was shown to be of type $\Errortop(\ubar,u)$
in Lemma\,\ref{L:IBPEASYSPACETIMEERRRORINTEGRALSESITMATES}.

\end{proof}

\subsection{Estimates for wave equation error integrals involving top order derivatives of vorticity} \label{SS:L2ESTIMATESFORTOPORDERVORTICITYINRHSOFWAVEEQUATIONS}

In this section, we derive a priori estimates for the integrals of Sect.\,\ref{S:IBPIDENTITIESNEEDEDFORTOPORDERVORTICITY} which requires a perfect divergence structure in the nonlinear interaction between vorticity and the velocity variables. These appear in the commuted wave equations satisfied by the velocity from Prop.\,\ref{P:MOSTDIFFICULTWAVETERMS} in the terms $\tander^N \mathfrak{G}_i$. Much like in the estimates of Sect.\,\ref{SS:ESTIMATESFORMOSTDIFFICULTEIKONALFUNCTIONERRORINTEGRALS},  this perfect divergence structure is only visible when summing over all the energy identities \eqref{E:ENERGYNULLFLUXINTEGRALIDENTITIESWAVE} for $v^1,v^2$, and $v^3$. Here we will not analyze the full contribution of $\tander^N \mathfrak{G}_i$, but only the error integrals requiring a perfect divergence. The main result is Prop.\,\ref{P:L2CONTROLOFTOPORDERINTERACTIONOFVORTICITYANDVELOCITY}.

\begin{lemma}[$L^2$ estimates for the terms featuring a perfect divergence structure]
\label{L:L2ESTIMATESFORTHETERMSFEATURINGADIVERGENCESTRUCTURE}

Let $N = \Ntop$,  let 
$\tander^N \in \mathfrak{P}^{(N)}$ be the set of top order $\nullhyparg{u}$-tangential commutator operators from
Def.\,\ref{D:STRINGSOFCOMMUTATIONVECTORFIELDS}. Let $\blowuprateofwaveWRTnumberofcommutations(N) = \blowupratetoporderwave$ denote the blowup rate of the wave variables defined in \eqref{E:BLOWUPRATEOFWAVEWRTNUMBEROFCOMMUTATORS}. For $\Sigma_t$-tangent vectorfields $\SigmatTan_1, \, \SigmatTan_2$, denote their Cartesian inner product by $\SigmatTan_1\cdot\SigmatTan_2 = \sum_{i=1}^3\SigmatTan_1^i\cdot\SigmatTan_2^i$. Then the following estimate holds for
$(\ubar,u) \in [\leftubar,\ubarboot) \times [\moreinterestingu_1,\moreinterestingu_2]$:
\begin{align}
		\begin{split} \label{E:L2ESTIMATESFORTHETERMSFEATURINGADIVERGENCESTRUCTURE}
			 - \int_{\characteristicdiamondtwoarg{[\leftubar,\ubar)}{[\moreinterestingu_1,u]}} \Speed^2 \exp(2\LogDensity) (\multipliervectorfield \tander^N \vortrenormalized)\cdot (\upmu \tander^N \vortrenormalized) \weight^{\blowupratetoporderwave} \ReciprocalLunitAppliedtoTimeFunction \, \voldiamond \le  C \, \Errortop(\ubar,u),
		\end{split}
	\end{align}
where we note that \eqref{E:L2ESTIMATESFORTHETERMSFEATURINGADIVERGENCESTRUCTURE} does not feature an absolute value on the left hand side.

\end{lemma}

\begin{proof}
	The starting point of the proof is identity \eqref{E:PERFECTDIVERGENCESTRUCTUREFOROVERDIFFERENTIATEDVORTICITY}. Using \eqref{E:COMMUTEDTRANSPORTPOINTWISEESTIMATESFORSPECIFICVORTICITYANDENTROPYGRADIENT}, Lemmas\,\ref{L:COERCIVENESSOFL2CONTROLLINGQUANITIESUNIFIEDCOMMUTATOR} and\,\ref{L:PRELIMINARYBELOWTOPORDERL2ESTIMATESFOREIKONALFUINCTIONQUANTITIES}, we see that the spacetime integral on RHS\,\eqref{E:PERFECTDIVERGENCESTRUCTUREFOROVERDIFFERENTIATEDVORTICITY} featuring $\upmu \Transport \tander^N \vortrenormalized$ is bounded in magnitude by: 
\begin{align} 
		\begin{split} \label{E:L2ESTIMATESFORTHETERMSFEATURINGADIVERGENCESTRUCTURESTEP1}
			& \lesssim \initialsmalldoublenull^2 + \int_{\ubar' = \leftubar}^{\ubar} \totalfluxcontrolvelocity_{[1,N]}(\ubar',u) \, \mathrm{d} \ubar'  +  \int_{u' = \moreinterestingu_1}^u \totalfluxcontrolvelocity_{[1,N]}(\ubar,u') \, \mathrm{d} u' \\
			& \ \ \ \ + \fundbootsmall(1 + \varsigma_0^{-1}) \strongangularcontrolvelocity_{[1,N]}(\ubar,u) 
			+ \int_{u' = \moreinterestingu_1}^u \left\{\fluxcontrolVort_{\le N}(\ubar,u')+ \fluxcontrolGradEnt_{\le N}(\ubar,u')\right\} \, \mathrm{d} u' \\
			& = \Errortop(\ubar,u).
		\end{split}
	\end{align}
	
By the bootstrap assumptions, the $L^\infty$-estimates of Prop.\,\ref{P:IMPROVEMENTOFAUXILIARYBOOTSTRAP}, and the identity $\blowuprateoftransportWRTnumberofcommutations(N) = \blowupratetoporderwave-1$, all the remaining spacetime integrals on RHS\,\eqref{E:PERFECTDIVERGENCESTRUCTUREFOROVERDIFFERENTIATEDVORTICITY} are bounded by $ \int_{\characteristicdiamondtwoarg{[\leftubar,\ubar)}{[\moreinterestingu_1,u]}} \{ \upmu |\tander^N \vortrenormalized|^2 + \weight |\tander^N \vortrenormalized|^2\} \weight^{\blowupratetoporderwave-1} \voldiamond \lesssim \int_{\ubar' = \leftubar}^{\ubar} \fluxcontrolVort_{ N}(\ubar',u) \, \rmd \ubar' + \int_{u' = \moreinterestingu_1}^u \fluxcontrolVort_{ N}(\ubar,u') \, \rmd u' = \Errortop(\ubar,u)$. 

Next, we note that the integral on RHS\,\eqref{E:PERFECTDIVERGENCESTRUCTUREFOROVERDIFFERENTIATEDVORTICITY} over the ingoing characteristic $\ingoingcharacteristicsurfacetwoarg{\ubar}{[\moreinterestingu_1,u]}$ features an overall negative sign, and hence it is trivially bounded above by 0. The remaining integral over $\ingoingcharacteristicsurfacetwoarg{\leftubar}{[\moreinterestingu_1,u]}$ is bounded by $\initialsmall$ by the initial data assumption \eqref{E:L2NORMSOFVORTANDGRADENTSMALLALONGINITIALINGOINGNULLSURFACE} since $\upmu \approx \ReciprocaluLunitAppliedtoTimeFunction$ (see \eqref{E:RECIPROCALULUNITAPPLIEDTOTIMEFUNCTIONAPPROXIMATELYMU}).

\end{proof}

\begin{lemma}[$L^2$ estimates for the easy over-differentiated vorticity integrals] \label{L:L2ESTIMATESFORTHEEASYOVERLYDIFFERENTIATEDINTEGRALS}

Let $N = \Ntop$,  let 
$\tander^N \in \mathfrak{P}^{(N)}$ be the set of top order $\nullhyparg{u}$-tangential commutator operators from
Def.\,\ref{D:STRINGSOFCOMMUTATIONVECTORFIELDS}. Let $\blowuprateofwaveWRTnumberofcommutations(N) = \blowupratetoporderwave$ denote the blowup rate of the wave variables defined in \eqref{E:BLOWUPRATEOFWAVEWRTNUMBEROFCOMMUTATORS}. 
Then the following estimates hold for
$(\ubar,u) \in [\leftubar,\ubarboot) \times [\moreinterestingu_1,\moreinterestingu_2]$:
\begin{align}
		\begin{split} \label{E:L2ESTIMATESFORTHEEASYOVERLYDIFFERENTIATEDINTEGRALS}
			& \left|  \int_{\characteristicdiamondtwoarg{[\leftubar,\ubar)}{[\moreinterestingu_1,u]}} \left\{\Speed^2 \exp(\LogDensity) (\multipliervectorfield \tander^N \vortrenormalized)\cdot [\upmu,\tander^N] (\Flatcurl v) + \Speed^2\exp(\LogDensity) (\multipliervectorfield \tander^N \vortrenormalized)\cdot [\tander^N,\upmu \Flatcurl] v \right\} \weight^{\blowupratetoporderwave} \ReciprocalLunitAppliedtoTimeFunction \, \voldiamond \right|, \\
			& \left|   \int_{\characteristicdiamondtwoarg{[\leftubar,\ubar)}{[\moreinterestingu_1,u]}} \Speed^2 \exp(\LogDensity) (\multipliervectorfield \tander^N \vortrenormalized^i ) \left\{\sum_{\substack{ \tander^{N_1}\tander^{N_2} = \tander^N \\ N_1 \le N-1}} \upmu (\tander^{N_1} \vortrenormalized^i) \tander^{N_2}(\exp(\LogDensity)) \right\} \weight^{\blowupratetoporderwave} \ReciprocalLunitAppliedtoTimeFunction \, \voldiamond \right| \\
			& \qquad \qquad \lesssim \Errortop(\ubar,u).
		\end{split}
	\end{align}
	
\end{lemma} 

\begin{proof}

We being with the integration by parts identity of \eqref{E:HARMLESSOVERDIFFERENTIATEDIBPFORVORTICITY}. Using the commutator estimates of Prop.\,\ref{P:COMMUTATORESTIMATES}, pointwise estimates \eqref{E:POINTWISEESTIMATESFORMUWEIGHTEDCARTESIANCOMMUTATOR}, writing $\Flatcurl v = \exp(\LogDensity) \vortrenormalized)$,  \eqref{E:COMMUTEDTRANSPORTPOINTWISEESTIMATESFORSPECIFICVORTICITYANDENTROPYGRADIENT}, the bootstrap assumptions, the $L^\infty$ estimates of Prop.\,\ref{P:IMPROVEMENTOFAUXILIARYBOOTSTRAP}, and Lemmas\,\ref{L:COERCIVENESSOFL2CONTROLLINGQUANITIESUNIFIEDCOMMUTATOR} and\,\ref{L:PRELIMINARYBELOWTOPORDERL2ESTIMATESFOREIKONALFUINCTIONQUANTITIES}, the integrals on RHS\,\eqref{E:HARMLESSOVERDIFFERENTIATEDIBPFORVORTICITY} featuring $\upmu \Transport \tander^N \vortrenormalized$ are: 
\begin{align} 
		\begin{split} \label{E:L2ESTIMATESFORTHEEASYOVERLYDIFFERENTIATEDINTEGRALSSTEP1}
			& \lesssim \initialsmalldoublenull^2 + \int_{\ubar' = \leftubar}^{\ubar} \totalfluxcontrolvelocity_{[1,N]}(\ubar',u) \, \mathrm{d} \ubar'  +  \int_{u' = \moreinterestingu_1}^u \totalfluxcontrolvelocity_{[1,N]}(\ubar,u') \, \mathrm{d} u' \\
			& \ \ \ \ + \fundbootsmall(1 + \varsigma_0^{-1}) \strongangularcontrolvelocity_{[1,N]}(\ubar,u) 
			+ \int_{u' = \moreinterestingu_1}^u \left\{\fluxcontrolVort_{\le N}(\ubar,u')+ \fluxcontrolGradEnt_{\le N}(\ubar,u')\right\} \, \mathrm{d} u' \\
			& = \Errortop(\ubar,u).
		\end{split}
	\end{align} 

Similarly, by the Leibniz rule, the bootstrap assumptions, commutator estimates of Prop.\,\ref{P:COMMUTATORESTIMATES}, the $L^\infty$ estimates of Prop.\,\ref{P:IMPROVEMENTOFAUXILIARYBOOTSTRAP}, pointwise estimate \eqref{E:POINTWISEESTIMATESFORMUWEIGHTEDCARTESIANCOMMUTATOR}, the spacetime integrals on RHS\,\eqref{E:HARMLESSOVERDIFFERENTIATEDIBPFORVORTICITY} featuring the terms $\ErrorIBP_1^{\vort;\text{below-top}}[\tander^N]$, as well as those featuring an outer-most $\Lunit$-differentiation (e.g. $\Lunit\left([\upmu,\tander^N](\Flatcurl v)\right)$), are also $= \Errortop(\ubar,u)$. We omit the details.

Similarly, since $\blowuprateoftransportWRTnumberofcommutations(N) = \blowupratetoporderwave-1$, the spacetime integral of $\ErrorIBP_1^{\vort;\text{below-top}}[\tander^N]$ is also $= \Errortop(\ubar,u)$. 

It remains to control the boundary integrals on the outgoing characteristics. Using similar arguments as above and Young's inequality, these integrals are bounded by: 
	\begin{align} \label{E:L2ESTIMATESFORTHEEASYOVERLYDIFFERENTIATEDINTEGRALSSTEP2}
		& \int_{\ingoingcharacteristicsurfacetwoarg{\ubar}{[\moreinterestingu_1,u]}} \left\{ \varsigma_1^{-1} \left|\tander^{\le N} \vortrenormalized\right|^2 + \varsigma_1  \left|\comdersmall^{[1,N];\le 1} \velocityarray\right|^2 + \varsigma_1 \left| \tandersmall^{[1,N]} \badcontrolvars\right|^2 \right\} \weight^{\blowupratetoporderwave}\,  \volingoingnullhypersurface.
	\end{align}
	
\noindent Thanks to the small coefficient $\varsigma_1$, Lemma \,\ref{L:PRELIMINARYBELOWTOPORDERL2ESTIMATESFOREIKONALFUINCTIONQUANTITIES}, Cor.\,\ref{C:NONSINGULARL2ESTIMATESFORWAVEVARIABLESTHATLOSEONEDERIVATIVE} imply that \eqref{E:L2ESTIMATESFORTHEEASYOVERLYDIFFERENTIATEDINTEGRALSSTEP2} is also $= \Errortop(\ubar,u)$. Nearly identical arguments hold for the integrals over $\ingoingcharacteristicsurfacetwoarg{\leftubar}{[\moreinterestingu_1,u]}$, after which the data assumptions of Sect.\,\ref{S:ASSUMPTIONSONTHEDATA} imply they are $\lesssim \initialsmall^2$ (which is $=  \Errortop(\ubar,u)$).

\end{proof}

\begin{proposition}[$L^2$ control of the nonlinear interaction between vorticity and velocity at top order] \label{P:L2CONTROLOFTOPORDERINTERACTIONOFVORTICITYANDVELOCITY}
Let $N = \Ntop$,  let 
$\tander^N \in \mathfrak{P}^{(N)}$ be the set of top order $\nullhyparg{u}$-tangential commutator operators from
Def.\,\ref{D:STRINGSOFCOMMUTATIONVECTORFIELDS}. Let $\blowuprateofwaveWRTnumberofcommutations(N) = \blowupratetoporderwave$ denote the blowup rate of the wave variables defined in \eqref{E:BLOWUPRATEOFWAVEWRTNUMBEROFCOMMUTATORS}. For $\Sigma_t$-tangent vectorfields $\SigmatTan_1, \, \SigmatTan_2$, denote their Cartesian inner product by $\SigmatTan_1\cdot\SigmatTan_2 = \sum_{i=1}^3\SigmatTan_1^i\cdot\SigmatTan_2^i$. For $\Sigma_t$-tangent vectorfields $\SigmatTan_1, \, \SigmatTan_2$, denote their Cartesian inner product by $\SigmatTan_1\cdot\SigmatTan_2 = \sum_{i=1}^3\SigmatTan_1^i\cdot\SigmatTan_2^i$. Then there exists an implicit constant such that the following estimate holds for
$(\ubar,u) \in [\leftubar,\ubarboot) \times [\moreinterestingu_1,\moreinterestingu_2]$:
\begin{align} \label{E:L2CONTROLOFTOPORDERINTERACTIONOFVORTICITYANDVELOCITY}
		\left| \int_{\characteristicdiamondtwoarg{[\leftubar,\ubar)}{[\moreinterestingu_1,u]}} \Speed^2 \exp(\LogDensity) \upmu \left\{ \Flatcurl(\tander^N \vortrenormalized)\right\}  \cdot (\multipliervectorfield \tander^N v) \weight^{\blowupratetoporderwave} \ReciprocalLunitAppliedtoTimeFunction \, \voldiamond\right| \lesssim   \,  \Errortop(\ubar,u),
	\end{align}\end{proposition} 

\begin{proof}
The starting point of the proof is the integration by parts identity \eqref{E:IBPIDENTITYNEEDEDFORTOPORDERVORTICITY}. By the already proved \eqref{E:L2ESTIMATESFORTHETERMSFEATURINGADIVERGENCESTRUCTURE} and \eqref{E:L2ESTIMATESFORTHEEASYOVERLYDIFFERENTIATEDINTEGRALS}, except for the spacetime integrals featuring $ \ErrorIBP_1^{\vortrenormalized;\text{top}}[\tander^N]$, $\ErrorIBP_2^{\vortrenormalized;\text{top}}[\tander^N]$ and the boundary integrals, the RHS\,\eqref{E:IBPIDENTITYNEEDEDFORTOPORDERVORTICITY} has already been shown to be $\le C\,  \Errortop(\ubar,u)$. 

We now bound the spacetime integrals on RHS\,\eqref{E:IBPIDENTITYNEEDEDFORTOPORDERVORTICITY} featuring $ \ErrorIBP_1^{\vortrenormalized;\text{top}}[\tander^N]$. Multiplying \eqref{E:CARTESIANPARTIAL1TOCOMMUTATORS}--\eqref{E:CARTESIANPARTIAL23TOCOMMUTATORS} by $\upmu$, using Prop.\,\ref{P:SCHEMATICSTRUCTUREOFVARIOUSTENSORSINTERMSOFCONTROLVARS}, \eqref{E:TANGENTIALSIMPLECOMMUTATORIDENTITY}--\eqref{E:TRANSVERSALTANGENTIALSIMPLECOMMUTATORIDENTITY}, the bootstrap assumptions,  \eqref{E:POINTWISESTIMATEFORBREVEXINTERMSOFNEWUL}, and Young's inequality, we have the following pointwise bound:
\begin{align} \label{E:L2CONTROLOFTOPORDERINTERACTIONOFVORTICITYANDVELOCITYSTEP1}
		\left| \ErrorIBP_1^{\vortrenormalized;\text{top}}[\tander^N]\right| \lesssim  \weight^{\blowupratetoporderwave} \left\{\left| \newuL \tander^N \velocityarray \right|^2 +  \left| \Lunit \tander^N \velocityarray\right|^2 + \varsigma_1 \left|\angrmD \tander^N \velocityarray\right|_{\gtorus}^2 + (1+\varsigma_1^{-1})\left|\tander^N \vortrenormalized\right|^2 \right\},
	\end{align}\noindent where we clarify that the term $ \left|\angrmD \tander^N \velocityarray\right|_{\gtorus}^2$ in RHS\,\eqref{E:L2CONTROLOFTOPORDERINTERACTIONOFVORTICITYANDVELOCITYSTEP1} came from $[\multipliervectorfield, \upepsilon_{iab} \upmu \p_a] \tander^N v^i$ present in $\ErrorIBP_1^{\vortrenormalized;\text{top}}[\tander^N]$, applying \eqref{E:TANGENTIALSIMPLECOMMUTATORIDENTITY}--\eqref{E:TRANSVERSALTANGENTIALSIMPLECOMMUTATORIDENTITY}, and \eqref{E:SMOOTHTORUSNORMCOMPARBLETOTANGENTIALCONTRACTIONS}.
	
By \eqref{E:COERCIVENESSOFNULLFLUXCONTROLWAVE} and \eqref{E:COERCIVENESSOFCONTROLVORT}, except for the $\varsigma_1 \left|\angrmD \tander^N \velocityarray\right|_{\gtorus}^2$ term, the spacetime integral of RHS\,\eqref{E:L2CONTROLOFTOPORDERINTERACTIONOFVORTICITYANDVELOCITYSTEP1} is bounded by
$ \int_{\ubar' = \leftubar}^{\ubar} \totalfluxcontrolvelocity_{N} (\ubar',u) \, \mathrm{d} \ubar'  
+  \int_{u' = \moreinterestingu_1}^u \totalfluxcontrolvelocity_{N} (\ubar,u') \, \mathrm{d} u' 
+ \int_{u' = \moreinterestingu_1}^u \fluxcontrolVort_{ N}(\ubar,u') \, \rmd u' =  \Errortop(\ubar,u)$. For this remaining term, we decompose the integral over the characteristic diamond using \eqref{E:DIVIDINGTHEDOMAINSFORSHARPCONSTANTS}, apply \eqref{E:OUTSIDEOFSMALLNEIGHBORHOODLOWERBOUNDFORMU} on $\characteristicdiamondtwoarg{[\leftubar,\ubarboot]}{[\moreinterestingu_1,u]} \backslash \smallneighborhoodofcreasearg{[\leftubar,\ubar]}$, and use the coercivity bounds \eqref{E:COERCIVENESSOFNULLFLUXCONTROLWAVE}, \eqref{E:UNIFIEDCOMMUTATORMIXANGULARSPACEITMECOERCIVITYWITHANGULARLAPLACIANINTHEINTEGRAND} to bound the spacetime integral by $\varsigma_1\strongangularcontrolvelocity_N(\ubar,u) +  \int_{u' = \moreinterestingu_1}^u \totalfluxcontrolvelocity_{N} (\ubar,u') \, \mathrm{d} u'  =  \Errortop(\ubar,u)$. 

We now bound the spacetime integrals of $\ErrorIBP_2^{\vortrenormalized;\text{top}}[\tander^N]$. Multiplying \eqref{E:CARTESIANPARTIAL1TOCOMMUTATORS}--\eqref{E:CARTESIANPARTIAL23TOCOMMUTATORS} by $\upmu$, we have that $|\upmu \p_a \weight| \lesssim \upmu \approx \ReciprocaluLunitAppliedtoTimeFunction$. Hence, using similar arguments as in the previous two paragraphs, we bound: 
\begin{align}
		\begin{split} \label{E:L2CONTROLOFTOPORDERINTERACTIONOFVORTICITYANDVELOCITYSTEP2}
			& \left| \blowupratetoporderwave\weight^{\blowupratetoporderwave-1} \upepsilon_{iab}(\upmu \p_a \weight)\Speed^2 \exp(\LogDensity) \left(\tander^N \vortrenormalized^b\right) \multipliervectorfield  \tander^N v^i  \right| \\
			& \ \  \lesssim  \weight^{\blowupratetoporderwave} \left\{\left| \newuL \tander^N \velocityarray \right|^2 +  \left| \Lunit \tander^N \velocityarray\right|^2 + \upmu \left|\nullangrmD \tander^N \velocityarray\right|_{\gnulltori}^2 + \ReciprocaluLunitAppliedtoTimeFunction \weight^{-2}\left|\tander^N \vortrenormalized\right|^2 \right\},
		\end{split}
	\end{align}
\notag Since $\blowuprateoftransportWRTnumberofcommutations(N) -1 = \blowupratetoporderwave-2$, the spacetime integral of \eqref{E:L2CONTROLOFTOPORDERINTERACTIONOFVORTICITYANDVELOCITYSTEP2} is bounded by $ \int_{\ubar' = \leftubar}^{\ubar} \totalfluxcontrolvelocity_{N} (\ubar',u) \, \mathrm{d} \ubar'  + \bulkcontrolVort_N(\ubar,u)  =  \Errortop(\ubar,u)$.

To bound the integral of the second term on RHS\,\eqref{E:ERRORTERM2IBPIDENTITYNEEDEDFORTOPORDERVORTICITY}, notice that $|\multipliervectorfield \weight| = |(1+2\upmu)\Lunit \weight + \muX \weight| \lesssim 1$. By Prop.\,\ref{P:SCHEMATICSTRUCTUREOFVARIOUSTENSORSINTERMSOFCONTROLVARS}, \eqref{E:TORITANGENTVECTORFIELDASSOCIATEDTODOUBLENULLFRAME}, \eqref{E:MUWEIGHTEDCARTESIANINTERMSOFDOUBLENULLFRAME}, and the bootstrap assumptions, we have:
	\begin{align}
		\begin{split} \label{E:L2CONTROLOFTOPORDERINTERACTIONOFVORTICITYANDVELOCITYSTEP3}
			& \left| \blowupratetoporderwave\weight^{\blowupratetoporderwave-1}( \multipliervectorfield \weight)\Speed^2 \exp(\LogDensity) \left(\tander^N \vortrenormalized^b\right) \upepsilon_{iab}\upmu \p_a (\tander^N v^i) \right| \\
			& \ \ \lesssim \weight^{\blowupratetoporderwave-1} \left| \tander^N \vortrenormalized\right|  \left\{\left| \newuL \tander^N \velocityarray\right| + \upmu \left| \Lunit \tander^N \velocityarray\right| + \upmu \left| \nullangrmD \tander^N \velocityarray\right|_{\gnulltori}\right\} \\
			& \ \ \lesssim (1+\varsigma_1^{-1}) \ReciprocaluLunitAppliedtoTimeFunction \weight^{\blowupratetoporderwave-2} \left| \tander^N \vortrenormalized\right|^2 + \varsigma_1 \weight^{\blowupratetoporderwave-1} \left\{ \left| \newuL \tander^N \velocityarray\right|^2 + \upmu \left| \nullangrmD \tander^N \velocityarray\right|_{\gnulltori}^2 \right\} + \weight^{\blowupratetoporderwave} \left| \Lunit \tander^N\velocityarray\right|^2.
		\end{split}
	\end{align}\noindent The spacetime integral of RHS\,\eqref{E:L2CONTROLOFTOPORDERINTERACTIONOFVORTICITYANDVELOCITYSTEP3} is then bounded by $(1+\varsigma_1^{-1}) \bulkcontrolVort_N(\ubar,u) + \varsigma_1\weakspacetimetransvelocitycontrol_N(\ubar,u) + \varsigma_1\weakspacetimeangularvelocitycontrol_N(\ubar,u) +  \int_{u' = \moreinterestingu_1}^u \totalfluxcontrolvelocity_{N} (\ubar,u') \, \mathrm{d} u'  =  \Errortop(\ubar,u)$.

It remains to control the boundary integrals of RHS\,\eqref{E:IBPIDENTITYNEEDEDFORTOPORDERVORTICITY}. From \eqref{E:TORITANGENTVECTORFIELDASSOCIATEDTODOUBLENULLFRAME}, the bootstrap assumptions, and the estimates for the ingoing eikonal function of Lemma\,\ref{L:LINFTYESTIMATESFORINGOINGEIKONALFUNCTIONANDDERIVATIVES} and Cor.\,\ref{C:PRELIMINARYESTIMATESFORTHEDOUBLENULLACOUSTICSCALARS},  \eqref{E:NULLCOORDINATEPARTIALINTERMSOFMETRICNORM}, \eqref{E:KEYESTIMATECONTROLLINGINVERSEMUBYINVERSEWEIGHT}, \eqref{E:COERCIVENESSOFNULLFLUXCONTROLWAVE}, \eqref{E:COERCIVENESSOFCONTROLVORT} and Young's inequality it follows that: 
\begin{align} 
			& \left|   \int_{\ingoingcharacteristicsurfacetwoarg{\ubar}{[\moreinterestingu_1,u]}} \left\{ \Speed^2 \exp(\LogDensity)(\tander^N\vortrenormalized^b) \left( 2 \newuL \tander^N v^i + 2\MagnitueofinnerproductofnewLandnewuL \ToriTangentVectorfieldAssociatedToDoubleNullFolliations \tander^N v^i\right) \upepsilon_{iab}\left\{\left(\frac{\MagnitueofinnerproductofnewLandnewuL}{\ingoingmu} - \upmu\right)\frac{X^a}{\Speed^2} + \upmu \angularcoeffinmuweightedspatialcartesian_a^A \frac{\Yvf{A} \ubar}{\Lunit \ubar}\right\}  \right.  \right.\\
			& \qquad  \left. \left. - \Speed^2  \exp(\LogDensity)(\tander^N\vortrenormalized^b) \upepsilon_{iab}  \left(\frac{X^a}{\Speed^2} \newuL v^i  + \frac{X^a}{\Speed^2} \MagnitueofinnerproductofnewLandnewuL \ToriTangentVectorfieldAssociatedToDoubleNullFolliations \tander^N v^i + \frac{\upmu}{\Speed^2}\angularcoeffinmuweightedspatialcartesian_a^A (\gtorus^{-1})^{AB}  \gtorusdoublenullCOV_B^C \nullgeop{x^C} \tander^N v^i\right) \left(1 + \frac{2 \MagnitueofinnerproductofnewLandnewuL}{\ingoingmu}\right) \right\} \weight^{\blowupratetoporderwave} \volingoingnullhypersurface \right| \\
			& \lesssim \varsigma_1 \int_{\ingoingcharacteristicsurfacetwoarg{\ubar}{[\moreinterestingu_1,u]}} \left\{ \left|\newuL \tander^N \velocityarray \right|^2 + \upmu \left| \nullangrmD \tander^N \velocityarray\right|_{\gnulltori}^2\right\} \weight^{\blowupratetoporderwave} \volingoingnullhypersurface + \varsigma_1^{-1}  \int_{\ingoingcharacteristicsurfacetwoarg{\ubar}{[\moreinterestingu_1,u]}}  \ReciprocaluLunitAppliedtoTimeFunction \left|\tander^N \vortrenormalized\right|^2  \weight^{\blowupratetoporderwave-1} \volingoingnullhypersurface \\
			& \lesssim  \varsigma_1\totalfluxcontrolvelocity_{N}(\ubar,u) + \varsigma_1^{-1} \fluxcontrolVort_N(\ubar,u)  \\
			& =   \Errortop(\ubar,u).
	\end{align}
By identical arguments and the initial data assumptions, the boundary integral over $\ingoingcharacteristicsurfacetwoarg{\leftubar}{[\moreinterestingu_1,u]}$ is $\lesssim \initialsmall^2 =    \Errortop(\ubar,u)$. 
							
Similarly, by \eqref{E:SMOOTHTORUSNORMCOMPARBLETOTANGENTIALCONTRACTIONS}, the boundary integrals over $\outgoingcharacteristicsurfacetwoarg{u}{[\leftubar,\ubar)}$  are also bounded in magnitude by $\lesssim  \varsigma_1\totalfluxcontrolvelocity_{N}(\ubar,u) + \varsigma_1^{-1} \fluxcontrolVort_N(\ubar,u)  =   \Errortop(\ubar,u)$.

\end{proof}

\subsection{Estimates for wave equation error integrals involving a loss of one derivative}
\label{SS:ESTIMATESFORWAVEERRORINTEGRALSTHATLOSEADERIVATIVE}
In Lemma\,\ref{L:ESTIMATESFORWAVEERRORINTEGRALSTHATLOSEADERIVATIVE},
we control below-top-order error integrals involving the $\tanderY^N \mytr_{\gtorus}\upchi$-involving products on
RHSs~\eqref{E:TOPCOMMUTEDWAVELFIRSTTHENALLYS}--\eqref{E:TOPCOMMUTEDWAVEALLYS}
using an approach that loses one derivative 
(which is permissible below top-order)
but leads to estimates that are less singular 
with respect to powers of $|\timefunction|$
compared to
the top-order error integrals we bounded in
Lemmas~\ref{L:BOUNDSFORMOSTDIFFICULTWAVEERRORNTEGRALS} and
\ref{L:MAINESTIMATESFORDIFFICULTSPACETIMEERRORINTEGRALSINVOLVINGIBPWRTL}.
The estimates are crucial for the ``descent scheme'' of
Prop.\,\ref{P:APRIORIL2ESTIMATESWAVEVARIABLES}, which states that the wave energies are weighted by one less factor of $\weight^2$
at every level of descent below top-order
until one reaches a level
where the energies are bounded without a weight.

\begin{lemma}[Estimates for error integrals involving a loss of one derivative]
\label{L:ESTIMATESFORWAVEERRORINTEGRALSTHATLOSEADERIVATIVE}
Assume that $2 \leq N \leq \Ntop$ and
$v \in \velocityarray = \{v^1,v^2,v^3\}$,
and recall that $\multipliervectorfield$ is defined in \eqref{E:MULTIPLIERVECTORFIELD}. Let $\varsigma_0$ be the constant featured in Lemma\,\ref{L:PRELIMINARYBELOWTOPORDERL2ESTIMATESFOREIKONALFUINCTIONQUANTITIES} (see also Remark\,\ref{R:SMALLNESSCONSTANTFORBELOWTOPORDERACOUSTICVARIABLES}).
Then if $\varsigma_1 \in (0,1]$ is sufficiently small, the following estimate holds
 for
$(\ubar,u) \in [\leftubar,\ubarboot) \times [\moreinterestingu_1,\moreinterestingu_2]$, where the implicit constant is \textbf{independent} if $\varsigma_0$ and $\varsigma_1$:
\begin{align}
\begin{split} \label{E:ESTIMATESFORWAVEERRORINTEGRALSTHATLOSEADERIVATIVE}
	&  \int_{\characteristicdiamondtwoarg{[\leftubar,\ubar)}{[\moreinterestingu_1,u]}} 
			\left| 
				\multipliervectorfield \tander^{N-1} v 
			\right| 
			\left|
				\begin{pmatrix}
					(\muX v) \tander^{N-1} \mytr_{\gtorus}\upchi 
						\\
					(\angrmd^{\#} v) \cdot \upmu \angrmd \tander^{N-2} \mytr_{\gtorus} \upchi
						\\
					\Yvfsmallcoeff{A} (\angrmd^{\#} v) \cdot \upmu \angrmd \tander^{N-2} \mytr_{\gtorus}\upchi
				\end{pmatrix}
		  \right|_{\gtorus} \weight^{\blowuprateofwaveWRTnumberofcommutations(N-1)}\ReciprocalLunitAppliedtoTimeFunction
		\, \voldiamond  
			\\
	&  \lesssim \varsigma_1(1 + \varsigma_0^{-1}) \int_{\moreinterestingu_1}^u \totalfluxcontrolvelocity_{N}(\ubar,u') \, \rmd u  +  \varsigma_1(1 + \varsigma_0^{-1}) \strongangularcontrolvelocity_N(\ubar,u) 
			\\
	& \ \ 
		+ 
		\Errorsubcriticalarg{N-1}(\timefunction,u),
\end{split}
\end{align}
where $\Errorsubcriticalarg{N-1}(\timefunction,u)$ 
satisfies \eqref{E:ERRORBELOWTOPORDERWAVEESTIMATES} with $N-1$ in the role of $M$.

\end{lemma}

\begin{proof}
We prove \eqref{E:ESTIMATESFORWAVEERRORINTEGRALSTHATLOSEADERIVATIVE} only
for the term generated by the first entry $(\muX v) \tander^{N-1} \mytr_{\gtorus} \upchi$ on 
LHS~\eqref{E:ESTIMATESFORWAVEERRORINTEGRALSTHATLOSEADERIVATIVE}; 
the remaining terms can be handled similarly. 
Since $\multipliervectorfield \tander^{N-1} v = (1 + 2 \upmu) \Lunit \tander^{N-1} v + 2\muX \tander^{N-1} v$, 
we can use 
\eqref{E:COERCIVENESSOFNULLFLUXCONTROLWAVE}, \eqref{E:L2ESTIMATESFORMUXTRANSVERSALDERIVATIVESOFWAVEVARIABLES},
\eqref{E:PRELIMINARYBELOWTOPORDERL2ESTIMATESFOREIKONALFUINCTIONQUANTITIES},
the bootstrap assumptions,
Young's inequality, the trivial inequality $\weight^{\blowuprateofwaveWRTnumberofcommutations(N-1)} \lesssim \weight^{\blowuprateofacousticgeoWRTnumberofcommutations(N-1)}$, 
and the fact that the error terms featured in RHS\,\eqref{E:PRELIMINARYBELOWTOPORDERL2ESTIMATESFOREIKONALFUINCTIONQUANTITIES} are increasing in their arguments
to deduce:
\begin{align}
		\begin{split} \label{E:PROOFENERGYWAVEESTIMATESBELOWTOPORDERTERMS}
			& \int_{\characteristicdiamondtwoarg{[\leftubar,\ubar)}{[\moreinterestingu_1,u]}} 
			\left| \multipliervectorfield \tander^{N-1} v \right| 
			\left| (\muX v) \tander^{N-1} \mytr_{\gtorus}\upchi \right| 
			\weight^{\blowuprateofwaveWRTnumberofcommutations(N-1)}\ReciprocalLunitAppliedtoTimeFunction
			\, \voldiamond  \\
			& \lesssim \varsigma_1^{-1} \int_{\characteristicdiamondtwoarg{[\leftubar,\ubar)}{[\moreinterestingu_1,u]}} \left\{ \left| \Lunit \tander^{N-1}\velocityarray\right|^2 + \left| \muX \tander^N \velocityarray\right|^2\right\}\weight^{\blowuprateofwaveWRTnumberofcommutations(N-1)}\ReciprocalLunitAppliedtoTimeFunction
			\, \voldiamond + \varsigma_1\int_{\characteristicdiamondtwoarg{[\leftubar,\ubar)}{[\moreinterestingu_1,u]}} \left| \tander^{N-1}\mytr_{\gtorus}\upchi\right|^2 \weight^{\blowuprateofacousticgeoWRTnumberofcommutations(N-1)} \voldiamond \\
			& \lesssim \varsigma_1(1 + \varsigma_0^{-1}) \int_{\moreinterestingu_1}^u \totalfluxcontrolvelocity_{N}(\ubar,u') \, \rmd u  +  \varsigma_1(1 + \varsigma_0^{-1}) \strongangularcontrolvelocity_N(\ubar,u) + \initialsmalldoublenull^2 \\
			& \ \ + \varsigma_1^{-1} \int_{\ubar' = \leftubar}^{\ubar} \totalfluxcontrolvelocity_{[1,N-1]}(\ubar',u) \, \rmd \ubar' + \varsigma_1^{-1} \int_{\moreinterestingu_1}^u \totalfluxcontrolvelocity_{[1,N-1]}(\ubar,u') \, \rmd u' + \varsigma_1(1 + \varsigma_0^{-1}) \strongangularcontrolvelocity_{[1,N-1]}(\ubar,u) \\
			& \ \ + \int_{u' = \moreinterestingu_1}^u \left\{\fluxcontrolVort_{\le N}(\ubar,u')+ \fluxcontrolGradEnt_{\le N}(\ubar,u')\right\} \, \mathrm{d} u'.
		\end{split}
	\end{align}
 
\end{proof}

\subsection{Proof of Prop.\,\ref{P:MAINWAVEENERGYINTEGRALINEQUALITIES}}
\label{SS:PROOFOFMAINWAVEENERGYINTEGRALINEQUALITIES}
\ \\

\noindent \textbf{Proof of \eqref{E:TOPORDERWAVEL2CONTROLLINGINTEGRALINEQUALITY}:}
We will prove estimate \eqref{E:TOPORDERWAVEL2CONTROLLINGINTEGRALINEQUALITY} independently for each term featured in the definition of  the master controlling quantity $\mastercontroltop(\ubar,u)$ defined in \eqref{E:MASTERCONTROLLINGQUANTITY}. Some of these terms, e.g. $\totalfluxcontrolpartialmodifiedchi(\ubar,u)$, generate boxed-constant-multiplied integrals on RHS\,\eqref{E:TOPORDERWAVEL2CONTROLLINGINTEGRALINEQUALITY}. We stress that \emph{the boxed-constant-multiplied integrals featured on RHS\,\eqref{E:TOPORDERWAVEL2CONTROLLINGINTEGRALINEQUALITY} are equal to the sum of all of the boxed-constant-multiplied integrals obtained when estimating each term in $\mastercontroltop(\ubar,u)$.}  In the rest of the proof given below, $\Errortoparg{N}(\ubar,u)$ denotes  
any term satisfying \eqref{E:ERRORTOPORDERWAVEESTIMATES}.

First, it is easy to see that RHS\,\eqref{E:FINALIMPRECISETOPORDERESTIMATESFORFULLYMODIFIEDCHI} is of type $\Errortoparg{N}(\ubar,u) \le$ RHS\,\eqref{E:TOPORDERWAVEL2CONTROLLINGINTEGRALINEQUALITY}.

Next, consider the $L^2$ estimates of Lemmas\,\ref{L:L2ESTIMATESFORMOSTDIFFICULTPRODUCT} and\,\ref{L:PRECISEL2ESTIMATEFORPARTIALLYMODIFIEDQUANTITIES}. Taking the supremum of \eqref{E:L2ESTIMATESFORMOSTDIFFICULTPRODUCT} and \eqref{E:PRECISEL2ESTIMATEFORPARTIALLYMODIFIEDQUANTITIES} over the relevant values of $\ubar$ and $u$, taking the maximum of \eqref{E:L2ESTIMATESFORMOSTDIFFICULTPRODUCT} over $\tander^N \in \mathfrak{P}^{(N)}$ and of \eqref{E:PRECISEL2ESTIMATEFORPARTIALLYMODIFIEDQUANTITIES}  over $\tanderY^{N-1} \in \mathfrak{Y}^{(N-1)}$, it is easy to see that $\totalfluxcontrolprecisefullymodifiedchi(\ubar,u),\, \blowupratetoporderwave \totalbulkcontrolprecisefullymodifiedchi(\ubar,u), \, \totalfluxcontrolpartialmodifiedchi(\ubar,u)$, and $(\blowupratetoporderwave - 1)\totalbulkcontrolpartialmodifiedchi(\ubar,u)$ are all $\le$ RHS\,\eqref{E:TOPORDERWAVEL2CONTROLLINGINTEGRALINEQUALITY}.

Now let $N = \Ntop$ and fix any $\tander^N \in \mathfrak{P}^{(N)}$, where  
$\mathfrak{P}^{(N)}$ is defined in Def.\,\ref{D:STRINGSOFCOMMUTATIONVECTORFIELDS}, and let $(\ubar,u) \in [\leftubar,\ubarboot)\times[\moreinterestingu_1,\moreinterestingu_2]$. Using \eqref{E:ENERGYNULLFLUXINTEGRALIDENTITIESWAVE}, we deduce:
	\begin{align} 
		\begin{split} \label{E:TOPORDERWAVEL2CONTROLLINGINTEGRALINEQUALITYPROOFSTEP1}
			& \sum_{i=1}^3 \ingoingfluxwave{\blowupratetoporderwave}[\tander^N  v^i](\ubar,u) + \outgoingfluxwave{\blowupratetoporderwave}[\tander^N  v^i](\ubar,u) + \frac{1}{2} \spacetimecoercive{\blowupratetoporderwave}[ \angrmD \tander^N  v^i](\ubar,u)  \\
			& + \frac{\blowupratetoporderwave}{2} \newspacetimecoercive{\blowupratetoporderwave - 1}[ \sqrt{\upmu + 2 \upmu\ReciprocaluLunitAppliedtoTimeFunction} \nullangrmD \tander^N  v^i](\ubar,u) + 2\blowupratetoporderwave \newspacetimecoercive{\blowupratetoporderwave - 1}[\newuL \tander^N  v^i](\ubar,u)  \\
			& = \sum_{i=1}^3 \ingoingfluxwave{\blowupratetoporderwave}[\tander^N  v^i](\leftubar,u) + \outgoingfluxwave{\blowupratetoporderwave}[\tander^N  v^i](\ubar,\moreinterestingu_1) - \int_{\characteristicdiamondtwoarg{[\leftubar,\ubar)}{[\moreinterestingu_1,u]}} \left\{ (1+2\upmu) L\tander^N  v^i + 2 \muX \tander^N  v^i \right\} \upmu \Box_{\gfour} \tander^N  v^i \, \weight^{\blowupratetoporderwave} \ReciprocalLunitAppliedtoTimeFunction \voldiamond  \\
			& \ \
			+ \int_{\characteristicdiamondtwoarg{[\leftubar,\ubar)}{[\moreinterestingu_1,u]}}  
			{^{(\multipliervectorfield)}\mathfrak{B}}[\tander^N  v^i]\, \weight^{\blowupratetoporderwave} \ReciprocalLunitAppliedtoTimeFunction \voldiamond  + 12.6 \int_{\characteristicdiamondtwoarg{[\leftubar,\ubar)}{[\moreinterestingu_1,u]} }\MagnitueofinnerproductofnewLandnewuL \ToriTangentVectorfieldAssociatedToDoubleNullFolliations \tander^N  v^i (\newuL \tander^N  v^i) (\newL \weight) \, \weight^{\blowupratetoporderwave - 1} \voldiamond.
		\end{split}
	\end{align}We will show that RHS\,\eqref{E:TOPORDERWAVEL2CONTROLLINGINTEGRALINEQUALITYPROOFSTEP1} is $\le$ RHS\,\eqref{E:TOPORDERWAVEL2CONTROLLINGINTEGRALINEQUALITY}. After that, we will take the supremum of the resulting estimate over the relevant values of $\ubar$ and $u$, and finally the maximum over all $\tander^N \in \mathfrak{P}^{(N)}$.

First, we note that \eqref{E:WAVEL2CONTROLLINGINITIALLYSMALL}  and implies $\sum_{i=1}^3\ingoingfluxwave{\blowupratetoporderwave}[\tander^N  v^i](\leftubar,u) + \outgoingfluxwave{\blowupratetoporderwave}[\tander^N  v^i](\ubar,\moreinterestingu_1)  \lesssim \initialsmalldoublenull^2 = \Errortoparg{N}(\ubar,u)$. 

Second, using Lemma\,\ref{L:ESTIMATESFORERRORTERMSGENERATEDBYMULTIPLIERVECTORFIELD}, we see that the ${^{(\multipliervectorfield)}\mathfrak{B}}[\tander^N  v^i]$-involving integrals on RHS\,\eqref{E:TOPORDERWAVEL2CONTROLLINGINTEGRALINEQUALITYPROOFSTEP1} is of type $\Errortoparg{N}(\ubar,u)$. Similarly, Lemma\,\ref{L:ESTIMATESFORERRORTERMSGENERATEDBYTHEWEIGHT} implies that the $\ToriTangentVectorfieldAssociatedToDoubleNullFolliations \tander^N v^i$-involving integrals on RHS\,\eqref{E:TOPORDERWAVEL2CONTROLLINGINTEGRALINEQUALITYPROOFSTEP1} is of type $\Errortoparg{N}(\ubar,u)$.

It remains to bound the last integral $- \sum_{i=1}^3 \int_{\characteristicdiamondtwoarg{[\leftubar,\ubar)}{[\moreinterestingu_1,u]}} \left\{ (1+2\upmu) L\tander^N  v^i + 2 \muX \tander^N  v^i \right\}( \upmu \Box_{\gfour} \tander^N  v^i )\, \weight^{\blowupratetoporderwave} \ReciprocalLunitAppliedtoTimeFunction \voldiamond$. Let $\mathfrak{G}_i$ denote the RHS of the wave equation satisfied by $v^i$, i.e. $\upmu \Box_{\gfour} v^i = \mathfrak{G}_i$. By Prop.\,\ref{P:MOSTDIFFICULTWAVETERMS}, we have the following three types of terms on the RHS $\upmu\Box_{\gfour} \tander^N v^i$: terms involving top-order derivatives of $\tander^N \mathfrak{G}_i$, the harmless terms $\HarmlessWave{N}$, and top-order derivatives of  $\mytr_{\gtorus}\upchi$.

We split the argument into Steps 0, 1, 2, 3.

\medskip

\noindent \emph{Step~0: $\tander^N \mathfrak{G}_i$ and $\HarmlessWave{N}$ terms}. By Lemma\,\ref{L:HARMLESSWAVETERMSERRORINTEGRALBOUNDS}, the spacetime integrals generated by the  $\HarmlessWave{N}$ terms are of  type $\Errortoparg{N}(\ubar,u)$.
 
Expressing $(1 +2\upmu)\Lunit + 2\muX  = \multipliervectorfield$ and using the covariant wave equation \eqref{E:VELOCITYWAVEEQUATIONWITHCURLCURL}, we have that the spacetime integral generated by the $\mathfrak{G}_i$ terms are:
 \begin{align}
	 	\begin{split} \label{E:TOPORDERWAVEL2CONTROLLINGINTEGRALINEQUALITYPROOFSTEP2}
			 & - \sum_{i=1}^3 \int_{\characteristicdiamondtwoarg{[\leftubar,\ubar)}{[\moreinterestingu_1,u]}} \left\{ \multipliervectorfield \tander^N  v^i \right\} \left\{ \upmu \tander^N  \mathfrak{G}_i \right\} \, \weight^{\blowupratetoporderwave} \ReciprocalLunitAppliedtoTimeFunction \voldiamond  = \sum_{i=1}^3 \int_{\characteristicdiamondtwoarg{[\leftubar,\ubar)}{[\moreinterestingu_1,u]}} \left\{ \multipliervectorfield \tander^N v^i\right\} \left\{ \Speed^2 \exp(\LogDensity) \upmu \left(\Flatcurl \, \tander^N \vortrenormalized \right)^i \right\} \, \weight^{\blowupratetoporderwave} \ReciprocalLunitAppliedtoTimeFunction \voldiamond \\
			 &  \ \ + \sum_{i=1}^3  \int_{\characteristicdiamondtwoarg{[\leftubar,\ubar)}{[\moreinterestingu_1,u]}} \left\{ \multipliervectorfield \tander^N v^i\right\} \left\{ \Speed^2 \exp(\LogDensity) \left([\tander^N, \upmu \Flatcurl ]\, \vortrenormalized\right)^i \right\} \, \weight^{\blowupratetoporderwave} \ReciprocalLunitAppliedtoTimeFunction \voldiamond \\
			 & \ \ + \sum_{i=1}^3   \int_{\characteristicdiamondtwoarg{[\leftubar,\ubar)}{[\moreinterestingu_1,u]}} \left\{ \multipliervectorfield \tander^N v^i\right\} \left\{ \sum_{\substack{ \tander^{N_1}\tander^{N_2} = \tander^{N} \\ N_2 \le N-1}} \tander^{N_1}(\Speed^2 \exp(\LogDensity)) \tander^{N_2}\left( \upmu (\Flatcurl \vortrenormalized)^i\right) \right\} \, \weight^{\blowupratetoporderwave} \ReciprocalLunitAppliedtoTimeFunction \voldiamond  \\
			 & \ \ + \sum_{i=1}^3   \int_{\characteristicdiamondtwoarg{[\leftubar,\ubar)}{[\moreinterestingu_1,u]}} \left\{ \multipliervectorfield \tander^N v^i\right\} \left\{ 	\tander^N \left(  \upmu \exp(-\LogDensity) \frac{p_{;\Ent}}{\overline{\varrho}} \GradEnt^a \p_a v^i\right) -\tander^N\left( \upmu \exp(-\LogDensity) \frac{p_{;\Ent}}{\overline{\varrho}} (\Flatdiv v) \GradEnt^i\right) \right\}  \, \weight^{\blowupratetoporderwave} \ReciprocalLunitAppliedtoTimeFunction \voldiamond  \\ 
			 & \ \ - \sum_{i=1}^3   \int_{\characteristicdiamondtwoarg{[\leftubar,\ubar)}{[\moreinterestingu_1,u]}} \left\{ \multipliervectorfield \tander^N v^i\right\} \left\{ \tander^N (\upmu \nullform_{(v)}^i) + \tander^N (\upmu 
			\mathfrak{L}_{(v)}^i) \right\}  \, \weight^{\blowupratetoporderwave} \ReciprocalLunitAppliedtoTimeFunction \voldiamond.
		\end{split}
	\end{align}
By Prop.\,\ref{P:SCHEMATICSTRUCTUREOFVARIOUSTENSORSINTERMSOFCONTROLVARS}, the bootstrap assumptions, pointwise estimates \eqref{E:POINTWISEESIMATESFORDERIVATIVESOFNULLFORMSTRUCTUREWAVEVARIABLES}, \eqref{E:POINTWISEESTIMATESFORALLLINEARTERMS}, \eqref{E:POINTWISEESTIMATESFORMUWEIGHTEDCARTESIANCOMMUTATOR}, and \eqref{E:COMMUTEDTRANSVERSALDERVIATVIESOFVORTICITYINTERMSOFTANGENTIAL}, all integrals on the RHS\,\eqref{E:TOPORDERWAVEL2CONTROLLINGINTEGRALINEQUALITYPROOFSTEP2} \emph{except for} the first one featuring $ \upmu \left(\Flatcurl \, \tander^N \vortrenormalized \right)^i $  are $= \Errortoparg{N}(\ubar,u)$. 

To control this remaining integral, we use the \eqref{E:L2CONTROLOFTOPORDERINTERACTIONOFVORTICITYANDVELOCITY} to conclude that the remaining error integral satisfies: 
\begin{align}
		\sum_{i=1}^3 \int_{\characteristicdiamondtwoarg{[\leftubar,\ubar)}{[\moreinterestingu_1,u]}} \left\{ \multipliervectorfield \tander^N v^i\right\} \left\{ \Speed^2 \exp(\LogDensity) \upmu \left(\Flatcurl \, \tander^N \vortrenormalized \right)^i \right\} \, \weight^{\blowupratetoporderwave} \ReciprocalLunitAppliedtoTimeFunction \voldiamond = \Errortoparg{N}(\ubar,u). \label{E:TOPORDERWAVEL2CONTROLLINGINTEGRALINEQUALITYPROOFSTEP3}
	\end{align}\medskip

\noindent \emph{Step~1: the case $\tander^N \notin \left\lbrace \tanderY^{N-1} \Lunit, \, \tanderY^N \right\rbrace$}.
By substituting \eqref{E:TOPCOMMUTEDWAVENOTDIFFICULT} 
for $\upmu \Box_{\gfour} \tander^N  v$ on 
RHS\,\eqref{E:TOPORDERWAVEL2CONTROLLINGINTEGRALINEQUALITYPROOFSTEP1}, using Lemma\,\ref{L:HARMLESSWAVETERMSERRORINTEGRALBOUNDS}, and Step~0, we find that the corresponding error integrals
are of type $\Errortoparg{N}(\ubar,u)$. 

\medskip

\noindent \emph{Step~2: the case $\tander^N = \tanderY^{N-1} \Lunit$}.
We use \eqref{E:TOPCOMMUTEDWAVELFIRSTTHENALLYS} to substitute for $\upmu \Box_{\gfour} \tander^N  v$ on 
RHS\,\eqref{E:TOPORDERWAVEL2CONTROLLINGINTEGRALINEQUALITYPROOFSTEP1}. All error integrals
except for the one generated by the first product on RHS\,\eqref{E:TOPCOMMUTEDWAVELFIRSTTHENALLYS} 
can be handled using the same arguments given in
Step~1.
The error integral generated by the first product on RHS\,\eqref{E:TOPCOMMUTEDWAVELFIRSTTHENALLYS}
is: 
\begin{align} \label{E:EASIERTOPORDERERRORINTEGRALLFIRSTPROOFOFWAVEL2ESTIMATES}
	&
	\int_{\characteristicdiamondtwoarg{[\leftubar,\ubar)}{[\moreinterestingu_1,u]}}
		\left\lbrace
			(1 + 2 \upmu) \Lunit \tanderY^{N-1}\Lunit  v 
			+ 
			2 \muX \tanderY^{N-1}\Lunit  v 
		\right\rbrace 
		(\angrmd^{\sharp}  v) 
		\cdot \upmu \angrmd \tanderY^{N-1} \mytr_{\gtorus}\upchi  
	\, \weight^{\blowupratetoporderwave} \ReciprocalLunitAppliedtoTimeFunction  \voldiamond,
\end{align}
and the estimate \eqref{E:ESTIMATESFOREASYTOPORDEREIKONALFUNCTIONERRORINTEGRALS}
implies that the integral is also of type $\Errortoparg{N}(\ubar,u)$.

\medskip

\noindent \emph{Step 3A: the case $\tander^N = \tanderY^N$}.

We substitute RHS\,\eqref{E:TOPCOMMUTEDWAVEALLYS} 
for $\upmu \Box_{\gfour} \tanderY^N  v$ on RHS\,\eqref{E:TOPORDERWAVEL2CONTROLLINGINTEGRALINEQUALITYPROOFSTEP1}.
All error integrals
except for the ones generated by the first two products on
RHS\,\eqref{E:TOPCOMMUTEDWAVEALLYS} can be handled using the same arguments given in
Step~1.
Writing $\tanderY^N = \tanderY^{N-1}\Yvf{A}$, the error integral generated by the second product on
RHS\,\eqref{E:TOPCOMMUTEDWAVEALLYS} is: 
\begin{align} \label{E:EASIERTOPORDERERRORINTEGRALALLYSPROOFOFWAVEL2ESTIMATES}
	&
	- \sum_{i=1}^3 \int_{\characteristicdiamondtwoarg{[\leftubar,\ubar)}{[\moreinterestingu_1,u]}}
		\left\lbrace
			(1 + 2 \upmu) \Lunit \tanderY^N v^i 
			+ 
			2 \muX \tanderY^N v^i
		\right\rbrace 
		\Yvfsmallcoeff{A}
		(\angrmd^{\sharp} v^i)
		\cdot 
		\upmu 
		\angrmd \tanderY^{N-1} \mytr_{\gtorus} \upchi
	\, \weight^{\blowupratetoporderwave} \ReciprocalLunitAppliedtoTimeFunction  \voldiamond,
\end{align}
and the estimate \eqref{E:ESTIMATESFOREASYTOPORDEREIKONALFUNCTIONERRORINTEGRALS}
implies that the integral is also of type $\Errortoparg{N}(\ubar,u)$.

The first product on
RHS\,\eqref{E:TOPCOMMUTEDWAVEALLYS} generates the following two difficult error integrals:
\begin{align} \label{E:DIFFICULTTOPORDERERRORINTEGRAL1ALLYSPROOFOFWAVEL2ESTIMATES}
	&
	-\sum_{i=1}^3\int_{\characteristicdiamondtwoarg{[\leftubar,\ubar)}{[\moreinterestingu_1,u]}}
		\left\lbrace
			(1 + 2 \upmu) \Lunit \tanderY^N v^i 
		\right\rbrace 
		(\muX v^i) 
		\tanderY^{N-1} \Yvf{A} \mytr_{\gtorus} \upchi
	\, \weight^{\blowupratetoporderwave} \ReciprocalLunitAppliedtoTimeFunction  \voldiamond,
			\\
		& -2\sum_{i=1}^3
		\int_{\characteristicdiamondtwoarg{[\leftubar,\ubar)}{[\moreinterestingu_1,u]}}
		\left\lbrace
			 \muX \tanderY^N v^i
		\right\rbrace 
		(\muX v^i) 
		\tanderY^{N-1} \Yvf{A} \mytr_{\gtorus} \upchi
	\, \weight^{\blowupratetoporderwave} \ReciprocalLunitAppliedtoTimeFunction  \voldiamond,
	\label{E:DIFFICULTTOPORDERERRORINTEGRAL2ALLYSPROOFOFWAVEL2ESTIMATES}
\end{align} 
and in \eqref{E:ESTIMATESFOREASYTOPORDEREIKONALFUNCTIONERRORINTEGRALS}, \eqref{E:FINALMAINSPACETIMEWAVEESTIMATEIBP},
and
\eqref{E:SPACETIMEBOUNDSMOSTDIFFICULTWAVEPRODUCT}
we proved that
\eqref{E:DIFFICULTTOPORDERERRORINTEGRAL1ALLYSPROOFOFWAVEL2ESTIMATES}--\eqref{E:DIFFICULTTOPORDERERRORINTEGRAL2ALLYSPROOFOFWAVEL2ESTIMATES}
are bounded in magnitude by $\leq$ 
RHS\,\eqref{E:TOPORDERWAVEL2CONTROLLINGINTEGRALINEQUALITY} as desired.
\emph{It is precisely this step that
generates all the boxed-constant-multiplied integrals on RHS\,\eqref{E:TOPORDERWAVEL2CONTROLLINGINTEGRALINEQUALITY}.}

We have therefore proved \eqref{E:TOPORDERWAVEL2CONTROLLINGINTEGRALINEQUALITY}.

\medskip

\noindent \textbf{Proof of \eqref{E:MAINWAVEBELOWTOPINTEGRALINEQUALITIES}:}

We fix $v \in \velocityarray = \{v^1,v^2,v^3\}$ and consider \eqref{E:TOPORDERWAVEL2CONTROLLINGINTEGRALINEQUALITYPROOFSTEP1} with $N'$ in the role of $N$ and $\blowuprateofwaveWRTnumberofcommutations(N')$ as in \eqref{E:BLOWUPRATEOFWAVEWRTNUMBEROFCOMMUTATORS}:
\begin{align} 
		\begin{split} \label{E:MAINWAVEINTEGRALIDENTITIESBELOWTOPORDERINTERMEDIATESTEP1}
			& \ingoingfluxwave{\blowuprateofwaveWRTnumberofcommutations(N')}[\tander^{N'}  v](\ubar,u) + \outgoingfluxwave{\blowuprateofwaveWRTnumberofcommutations(N')}[\tander^{N'}  v](\ubar,u) + \frac{1}{2} \spacetimecoercive{\blowuprateofwaveWRTnumberofcommutations(N')}[ \angrmD\tander^{N'}  v](\ubar,u)  \\
			& + \frac{\blowuprateofwaveWRTnumberofcommutations(N')}{2} \newspacetimecoercive{\blowuprateofwaveWRTnumberofcommutations(N')-1}[ \sqrt{\upmu + 2 \upmu\ReciprocaluLunitAppliedtoTimeFunction} \nullangrmD \tander^{N'}  v](\ubar,u) + 2 \blowuprateofwaveWRTnumberofcommutations(N') \newspacetimecoercive{\blowuprateofwaveWRTnumberofcommutations(N')-1}[\newuL \tander^{N'} v](\ubar,u)  \\
			& = \ingoingfluxwave{\blowuprateofwaveWRTnumberofcommutations(N')}[\tander^{N'}  v](\leftubar,u) + \outgoingfluxwave{\blowuprateofwaveWRTnumberofcommutations(N')}[\tander^{N'} v](\ubar,\moreinterestingu_1) - \int_{\characteristicdiamondtwoarg{[\leftubar,\ubar)}{[\moreinterestingu_1,u]}} \left\{ (1+2\upmu) L\tander^{N'}  v + 2 \muX \tander^{N'}  v \right\} \upmu \Box_{\gfour} \tander^{N'}  v \, \weight^{\blowuprateofwaveWRTnumberofcommutations(N')} \ReciprocalLunitAppliedtoTimeFunction \voldiamond  \\
			& \ \
			+ \int_{\characteristicdiamondtwoarg{[\leftubar,\ubar)}{[\moreinterestingu_1,u]}}  
			{^{(\multipliervectorfield)}\mathfrak{B}}[\tander^{N'}  v]\, \weight^{\blowuprateofwaveWRTnumberofcommutations(N')} \ReciprocalLunitAppliedtoTimeFunction \voldiamond  + 2\blowuprateofwaveWRTnumberofcommutations(N') \int_{\characteristicdiamondtwoarg{[\leftubar,\ubar)}{[\moreinterestingu_1,u]} }\MagnitueofinnerproductofnewLandnewuL \ToriTangentVectorfieldAssociatedToDoubleNullFolliations \tander^{N'}  v (\newuL \tander^{N'}  v) (\newL \weight) \, \weight^{\blowuprateofwaveWRTnumberofcommutations(N')-1} \voldiamond.
		\end{split}
	\end{align}
We will show that if $2 \leq N \leq \Ntop$ and
$1 \leq N' \leq N-1$, 
then 
$
\mbox{RHS\,\eqref{E:MAINWAVEINTEGRALIDENTITIESBELOWTOPORDERINTERMEDIATESTEP1}} 
\leq 
\mbox{RHS\,\eqref{E:MAINWAVEBELOWTOPINTEGRALINEQUALITIES}}$. 
Then for the same reasons given just below \eqref{E:TOPORDERWAVEL2CONTROLLINGINTEGRALINEQUALITYPROOFSTEP1},
we see, taking into account the Def.\,\ref{D:SUMMEDL2CONTROLLINGQUANTITIES}
of $\totalcontrolwave_{[1,N-1]} (\ubar,u)$,
that would imply \eqref{E:MAINWAVEBELOWTOPINTEGRALINEQUALITIES}.

 First, we note that \eqref{E:WAVEL2CONTROLLINGINITIALLYSMALL}  and implies $\ingoingfluxwave{\blowuprateofwaveWRTnumberofcommutations(N')}[\tander^{N'}  v](\leftubar,u) + \outgoingfluxwave{\blowuprateofwaveWRTnumberofcommutations(N')}[\tander^{N'}  v](\ubar,\moreinterestingu_1)  \lesssim \initialsmalldoublenull^2 = 			\Errorsubcriticalarg{N'}$

Next, whenever $\blowuprateofwaveWRTnumberofcommutations(N') > 0$, we use \eqref{E:ESTIMATESFORERRORTERMSGENERATEDBYTHEWEIGHT} to show that the last integral on RHS\,\eqref{E:MAINWAVEINTEGRALIDENTITIESBELOWTOPORDERINTERMEDIATESTEP1} featuring $\MagnitueofinnerproductofnewLandnewuL \ToriTangentVectorfieldAssociatedToDoubleNullFolliations \tander^{N'}  v $ is $= \Errorsubcriticalarg{N'}$; when $\blowuprateofwaveWRTnumberofcommutations(N') = 0$ there is nothing to bound.

 The spacetime integral featuring ${^{(\multipliervectorfield)}\mathfrak{B}}[\tander^{N'}  v]$ is also $= \Errorsubcriticalarg{N'}$ by \eqref{E:WAVEENERGYESTIMATEMULTIPLIERERRORTERMS}.
 
 It remains to control the integrals featuring $\upmu \Box_{\gfour}\tander^{N'} v$ on RHS\,\eqref{E:MAINWAVEINTEGRALIDENTITIESBELOWTOPORDERINTERMEDIATESTEP1}, which we do using Prop.\,\ref{P:MOSTDIFFICULTWAVETERMS}. Strictly easier arguments used to prove Step 0 above proves that the spacetime integrals from the  inhomogeneity satisfy: 
\begin{align}
		\begin{split} \label{E:MAINWAVEINTEGRALIDENTITIESBELOWTOPORDERINTERMEDIATESTEP2}
			\int_{\characteristicdiamondtwoarg{[\leftubar,\ubar)}{[\moreinterestingu_1,u]}} \left\{ (1+2\upmu) L\tander^{N'}  v + 2 \muX \tander^{N'}  v \right\} \left\{ \tander^{N'} (\upmu \mathfrak{G})\right\}\, \weight^{\blowuprateofwaveWRTnumberofcommutations(N')} \ReciprocalLunitAppliedtoTimeFunction \voldiamond =  \Errorsubcriticalarg{N'}(\ubar,u).
		\end{split}
	\end{align}
The reason that the proof of \eqref{E:MAINWAVEINTEGRALIDENTITIESBELOWTOPORDERINTERMEDIATESTEP2} is strictly easier than Step 0 above for $\Ntop$ is that $\blowuprateofacousticgeoWRTnumberofcommutations(N') = \blowuprateofwaveWRTnumberofcommutations(N') + 1$ whenever $N' \le \Ntop -1$. Specifically, we do not have to rely on decomposing $\VortVort^i$ as in  \eqref{E:VELOCITYWAVEEQUATIONWITHCURLCURL} and apply the integration by parts arguments of Sect.\,\ref{SS:L2ESTIMATESFORTOPORDERVORTICITYINRHSOFWAVEEQUATIONS}. That is, we simply apply Young's inequality as follows:
\begin{align}
		\begin{split} \label{E:MAINWAVEINTEGRALIDENTITIESBELOWTOPORDERINTERMEDIATESTEP3}
			& \left| \int_{\characteristicdiamondtwoarg{[\leftubar,\ubar)}{[\moreinterestingu_1,u]}}  \left\{ (1+2\upmu) L\tander^{N'}  v + 2 \muX \tander^{N'}  v \right\} \left\{ \upmu \tander^{N'} \VortVort \right\} \, \weight^{\blowuprateofwaveWRTnumberofcommutations(N')} \ReciprocalLunitAppliedtoTimeFunction \voldiamond \right| \\
			& \ \ \lesssim \int_{\characteristicdiamondtwoarg{[\leftubar,\ubar)}{[\moreinterestingu_1,u]}} \left\{ \left| \Lunit \tander^{N'} \velocityarray\right|^2  +\left| \muX \tander^{N'} \velocityarray\right|^2 \right\} \weight^{\blowuprateofwaveWRTnumberofcommutations(N')} \ReciprocalLunitAppliedtoTimeFunction \voldiamond + \int_{\characteristicdiamondtwoarg{[\leftubar,\ubar)}{[\moreinterestingu_1,u]}}  \upmu \left| \tander^{N'} \VortVort\right|^2   \weight^{\blowuprateofwaveWRTnumberofcommutations(N')} \ReciprocalLunitAppliedtoTimeFunction \voldiamond
		\end{split}.
	\end{align}Using \eqref{E:COERCIVENESSOFNULLFLUXCONTROLWAVE} and \eqref{E:L2ESTIMATESFORMUXTRANSVERSALDERIVATIVESOFWAVEVARIABLES}, the first integral on RHS\,\eqref{E:MAINWAVEINTEGRALIDENTITIESBELOWTOPORDERINTERMEDIATESTEP3} is always $\lesssim \int_{\ubar' = \leftubar}^{\ubar} \totalfluxcontrolvelocity_{N'} (\ubar',u) \, \rmd \ubar' + \int_{\moreinterestingu_1}^u \totalfluxcontrolvelocity_{N'}(\ubar,u') \, \rmd u' =  \Errorsubcriticalarg{N'}(\ubar,u) $.  Whenever $\blowuprateofacousticgeoWRTnumberofcommutations(N') -1 > 0$, the second integral is bounded by $\bulkcontrolVortVort_{N'}(\ubar,u) = \Errorsubcriticalarg{N'}(\ubar,u)$. Whenever $\blowuprateofacousticgeoWRTnumberofcommutations(N') - 1 = 0$, it is bounded by $\int_{\ubar' = \leftubar}^{\ubar} \fluxcontrolVortVort_{N'} (\ubar',u) \, \rmd \ubar'  = \Errorsubcriticalarg{N'}(\ubar,u)$. Combining \eqref{E:MAINWAVEINTEGRALIDENTITIESBELOWTOPORDERINTERMEDIATESTEP2}--\eqref{E:MAINWAVEINTEGRALIDENTITIESBELOWTOPORDERINTERMEDIATESTEP3}, the analogous of Step 0  above for $N' \le \Ntop -1$ is $\Errorsubcriticalarg{N'}(\ubar,u)$.

We now clarify the reason as to why controlling the nonlinear effects from the derivatives of $\VortVort$ where $N' \le \Ntop -1$ was so much easier than $N = \Ntop$. This was because $\blowuprateofacousticgeoWRTnumberofcommutations(\Ntop) = \blowuprateofwaveWRTnumberofcommutations(\Ntop) + 1.5$, and hence we would not be able to estimate the RHS\,\eqref{E:MAINWAVEINTEGRALIDENTITIESBELOWTOPORDERINTERMEDIATESTEP3} for $N = \Ntop$ as we did in the previous paragraph. This is why we had to find the perfect divergence structure in Sect.\,\ref{S:IBPIDENTITIESNEEDEDFORTOPORDERVORTICITY} (see also the $L^2$ estimates of Sect.\,\ref{SS:L2ESTIMATESFORTOPORDERVORTICITYINRHSOFWAVEEQUATIONS}).

The analogs of all the estimates through Step~1 above can be carried out
just as before; our arguments show
that since $1 \leq N' \leq N-1$,
all the corresponding error integrals are of type $\Errorsubcriticalarg{N-1}$,
i.e., that they are $\lesssim \mbox{RHS\,\eqref{E:ERRORBELOWTOPORDERWAVEESTIMATES}}$
with $M = N-1$.
The big difference occurs in Steps~2 and 3, where we now control the error integrals generated by the $\mytr_{\gtorus}\upchi$ terms on 
RHSs~\eqref{E:TOPCOMMUTEDWAVELFIRSTTHENALLYS}--\eqref{E:TOPCOMMUTEDWAVEALLYS}
using the derivative-losing estimate
\eqref{E:ESTIMATESFORWAVEERRORINTEGRALSTHATLOSEADERIVATIVE},
which shows that they are bounded by:
\begin{align} \label{E:MAINWAVEBELOWTOPINTEGRALINEQUALITIESSTEP1}
	 \lesssim \varsigma_1(1 + \varsigma_0^{-1}) \int_{\moreinterestingu_1}^u \totalfluxcontrolvelocity_{N'+1}(\ubar,u') \, \rmd u  +  \varsigma_1(1 + \varsigma_0^{-1}) \strongangularcontrolvelocity_{N'+1}(\ubar,u) 
		+ 
		\Errorsubcriticalarg{N'}(\ubar,u).
\end{align}
Since $1 \leq N' \leq N-1$, we have $\Errorsubcriticalarg{N'} \lesssim \Errorsubcriticalarg{N-1}$,
which is $\leq \mbox{RHS\,\eqref{E:MAINWAVEBELOWTOPINTEGRALINEQUALITIES}}$ as desired. If $N' = N-1$, then the two integrals featured in RHS\,\eqref{E:MAINWAVEBELOWTOPINTEGRALINEQUALITIESSTEP1} are directly featured in \eqref{E:MAINWAVEBELOWTOPINTEGRALINEQUALITIES}. If $N' \le N -2$, then the integrals are of type $\Errorsubcriticalarg{N-1}(\ubar,u)$, which concludes the proof.

\subsubsection{The proof of Props.\,\ref{P:TOPORDERAPRIORIL2ESTIMATESMASTERCONTROL}--\ref{P:APRIORIL2ESTIMATESWAVEVARIABLES}}
\label{SSS:PROOFOFTOPORDERAPRIORIL2ESTIMATEMASTERCONTROLANDALLORDERWAVES} \hfill

\medskip

\noindent \textbf{Estimate for $ \mastercontroltop$:}

Throughout this proof, we silently use the trivial inequalities $\totalfluxcontrolvelocity_{[1,N]}  \le \totalcontrolwave_{[1,N]}(\ubar,u)$ and $  \strongangularcontrolvelocity_{[1,N]}(\ubar,u) \le  \totalcontrolwave_{[1,N]}(\ubar,u)$, $\weakspacetimeangularvelocitycontrol_N{[1,N]}(\ubar,u) \le  \totalcontrolwave_{[1,N]}(\ubar,u)$, and $ \weakspacetimetransvelocitycontrol_{[1,N]}(\ubar,u) \le  \totalcontrolwave_{[1,N]}(\ubar,u)$, see \eqref{E:VELOCITYFLUXANDSTRONGANGULARL2CONTROLLINGQUANTITY}. Similar remarks apply to the individual terms present in $\mastercontroltop(\ubar,u)$.

The starting point of the proof is the a priori $L^2$ estimate \eqref{E:TOPORDERWAVEL2CONTROLLINGINTEGRALINEQUALITY}. Expanding $\mastercontroltop(\ubar,u)$ using \eqref{E:MASTERCONTROLLINGQUANTITY}, we have: 
\begin{align}
		\begin{split} \label{E:PROOFOFTOPORDERAPRIORIL2ESTIMATEMASTERCONTROLANDALLORDERWAVESSTEP1}
			&   \totalfluxcontrolvelocity_N(\ubar,u) +
				\frac{1}{2} \strongangularcontrolvelocity_N(\ubar,u) + 2\blowupratetoporderwave \weakspacetimetransvelocitycontrol_{\Ntop}(\ubar,u)  +\frac{\blowupratetoporderwave}{2}  \weakspacetimeangularvelocitycontrol_{\Ntop}(\ubar,u) \\
				& \ \  + \totalfluxcontrolprecisefullymodifiedchi(\ubar,u)
				 + \blowupratetoporderwave \totalbulkcontrolprecisefullymodifiedchi(\ubar,u)  \\
			&  \ \ + \totalfluxcontrolimprecisefullymodifiedchi(\ubar,u) + (\blowupratetoporderwave + 1) \totalbulkcontrolimprecisefullymodifiedchi(\ubar,u) \\
			& \ \ + \totalfluxcontrolpartialmodifiedchi(\ubar,u) +(\blowupratetoporderwave - 1)\totalbulkcontrolpartialmodifiedchi(\ubar,u) \\
			& \le   \boxed{2.81 + 5.18} \weakspacetimetransvelocitycontrol_N(\ubar,u)  + \boxed{2.33 + 1.01} \totalbulkcontrolprecisefullymodifiedchi(\ubar,u) \\
			& \ \ +  \boxed{1.02 + 2.07} \weakspacetimeangularvelocitycontrol_N(\ubar,u) +  \boxed{4.13}  \totalbulkcontrolpartialmodifiedchi(\ubar,u) \\
			& \ \  + \boxed{\frac{2}{3}} \totalfluxcontrolvelocity_N(\ubar,u) + \boxed{\frac{3}{4}(1.01)} \totalfluxcontrolpartialmodifiedchi(\ubar,u) \\
			& \ \ + \Errortoparg{N}(\ubar,u),
		\end{split} 
	\end{align}
We now highlight the following crucial fact: \emph{all of the boxed-constants on RHS\,\eqref{E:PROOFOFTOPORDERAPRIORIL2ESTIMATEMASTERCONTROLANDALLORDERWAVESSTEP1} are \underline{\textbf{strictly less}} than the respective term on LHS\,\eqref{E:PROOFOFTOPORDERAPRIORIL2ESTIMATEMASTERCONTROLANDALLORDERWAVESSTEP1}.}  Hence, all of the integrals on RHS\,\eqref{E:PROOFOFTOPORDERAPRIORIL2ESTIMATEMASTERCONTROLANDALLORDERWAVESSTEP1} featuring boxed-constants terms may be absorbed onto the LHS\,\eqref{E:PROOFOFTOPORDERAPRIORIL2ESTIMATEMASTERCONTROLANDALLORDERWAVESSTEP1}. Taking into account the definition of $\Errortoparg{N}(\ubar,u)$ in \eqref{E:ERRORTOPORDERWAVEESTIMATES}, we thus proven that there exists a universal constant $C$ such that for all $\varsigma_0,\varsigma_1 \in (0,1]$ sufficiently small, we have:
\begin{align}
		\begin{split} \label{E:PROOFOFTOPORDERAPRIORIL2ESTIMATEMASTERCONTROLANDALLORDERWAVESSTEP2}
			\mastercontroltop(\ubar,u) & \le C  (1 + \varsigma_0^{-1}+\varsigma_1^{-1}) \initialsmalldoublenull^2 \\
			& \ \ + C (1 + \varsigma_0^{-1}+\varsigma_1^{-1}) \int_{\ubar'=\leftubar}^{\ubar} \frac{1}{|\ubar'|^{1/2}} \mastercontroltop(\ubar',u) \, \rmd \ubar' +  C (1 + \varsigma_0^{-1}+\varsigma_1^{-1}) \int_{u'=\moreinterestingu_1}^{u} \mastercontroltop(\ubar,u') \, \rmd u' \\
			& \ \ + C\left\{\varsigma_1\left(1+\varsigma_0^{-1}\right) + \fundbootsmall \left(1+\varsigma_0^{-1}\right) + \mr{\upalpha}\right\} \mastercontroltop(\ubar,u) \\
			&  \ \ +   \int_{u' = \moreinterestingu_1}^u \left\{\fluxcontrolVort_{\le N}(\ubar,u')+ \fluxcontrolGradEnt_{\le N }(\ubar,u')\right\} \, \mathrm{d} u +   \int_{u' = \moreinterestingu_1}^u \left\{\fluxcontrolVortVort_{\le N}(\ubar,u')+ \fluxcontrolDivGradEnt_{\le N }(\ubar,u')\right\} \, \mathrm{d} u' \\	
			& \ \ + (1+\varsigma_1^{-1}) \fluxcontrolVort(\ubar,u) + (1+\varsigma_1^{-1}) \bulkcontrolVort(\ubar,u) + \bulkcontrolVortVort_N(\ubar,u) +  \bulkcontrolDivGradEnt_N(\ubar,u). 
		\end{split}
	\end{align}Recall that $\varsigma_0 \in (0,1]$ is the fixed constant featured in Lemma\,\ref{L:PRELIMINARYBELOWTOPORDERL2ESTIMATESFOREIKONALFUINCTIONQUANTITIES} (see also Remark\,\ref{R:SMALLNESSCONSTANTFORBELOWTOPORDERACOUSTICVARIABLES}). Taking $\varsigma_1, \, \fundbootsmall, \mr\upalpha$ to be sufficiently small, we may absorb the $ C\left\{\varsigma_1\left(1+\varsigma_0^{-1}\right) + \fundbootsmall \left(1+\varsigma_0^{-1}\right) + \mr{\upalpha}\right\} \mastercontroltop(\ubar,u)$ term into the LHS of \eqref{E:PROOFOFTOPORDERAPRIORIL2ESTIMATEMASTERCONTROLANDALLORDERWAVESSTEP2}. Next, using the already proved $L^2$ estimates for the transport variables of Props.\,\ref{P:APRIORIL2ESTIMATESFORTRANSPORTVARIABLES}--\ref{P:TOPORDERMODIFIEDAPRIORIL2ESTIMATESANDDOUBLENULLTORIESTIMATES}, \emph{which only relied on the bootstrap assumptions for the wave variables}, the last two lines of RHS\,\eqref{E:PROOFOFTOPORDERAPRIORIL2ESTIMATEMASTERCONTROLANDALLORDERWAVESSTEP2} are $\lesssim \initialsmall^2$. 

 Applying Gr\"onwall's inequality on the remaining terms concludes the proof of Prop.\,\ref{P:TOPORDERAPRIORIL2ESTIMATESMASTERCONTROL}. We clarify that the constant $C> 0$ present on RHS\,\eqref{E:MAINTOPORDERL2ESTIMATESMASTERCONTROLQUANTITY} is allowed to depend on $\varsigma_0,\, \varsigma_1$.

\medskip

\noindent \textbf{Estimate for $ \totalcontrolwave_{[1,\Ntop-1]}$:}

Consider the $L^2$ a priori estimate \eqref{E:MAINWAVEBELOWTOPINTEGRALINEQUALITIES}. Using the precise form \eqref{E:ERRORBELOWTOPORDERWAVEESTIMATES}, this implies: 
	\begin{align}
		\begin{split} \label{E:PROOFOFBELOWTOPORDERAPRIORIL2ESTIMATEWAVESTEP1}
			 \totalcontrolwave_{[1,\Ntop-1]}(\ubar,u) & \le C  \initialsmalldoublenull^2 + C\varsigma_1(1+\varsigma_0^{-1})\int_{u' = \moreinterestingu_1}^u \mastercontroltop(\ubar,u') \, \mathrm{d} u' + C\varsigma_1(1+\varsigma_0^{-1})\mastercontroltop(\ubar,u) \\
			 & \ \ + C\varsigma_1^{-1} \int_{u' = \moreinterestingu_1}^u  \totalcontrolwave_{[1,\Ntop-1]}(\ubar,u') \, \rmd u' + C\varsigma_1(1 + \varsigma_0^{-1}) \totalcontrolwave_{[1,\Ntop-1]}(\ubar,u) \\
			& \ \ + \int_{\ubar' = \leftubar}^{\ubar} \fluxcontrolVortVort_{\le \Ntop-1} (\ubar',u) \, \rmd \ubar' + \bulkcontrolVortVort_{\le \Ntop-1}(\ubar,u).
		\end{split}
	\end{align}
where $C$ is independent of $\varsigma_0$ and $\varsigma_1$, provided that $\varsigma_1 \in (0,1]$ is sufficiently small. Since $\varsigma_0 \in (0,1]$ is the fixed constant featured in Lemma\,\ref{L:PRELIMINARYBELOWTOPORDERL2ESTIMATESFOREIKONALFUINCTIONQUANTITIES}, by making $\varsigma_1$ even smaller if necessary, we may absorb the term $ C\varsigma_1(1 + \varsigma_0^{-1}) \totalcontrolwave_{[1,\Ntop-1]}(\ubar,u) $ on RHS\,\eqref{E:PROOFOFBELOWTOPORDERAPRIORIL2ESTIMATEWAVESTEP1} over to the LHS\,\eqref{E:PROOFOFBELOWTOPORDERAPRIORIL2ESTIMATEWAVESTEP1}: 
\begin{align}
		\begin{split} \label{E:PROOFOFBELOWTOPORDERAPRIORIL2ESTIMATEWAVESTEP2}
			 \totalcontrolwave_{[1,\Ntop-1]}(\ubar,u) & \le C  \initialsmalldoublenull^2 + C\varsigma_1(1+\varsigma_0^{-1})\int_{u' = \moreinterestingu_1}^u \mastercontroltop(\ubar,u') \, \mathrm{d} u' + C\varsigma_1(1+\varsigma_0^{-1})\mastercontroltop(\ubar,u) \\
			 & \ \ + C\varsigma_1^{-1} \int_{u' = \moreinterestingu_1}^u  \totalcontrolwave_{[1,\Ntop-1]}(\ubar,u') \, \rmd u' \\
			& \ \ + \int_{\ubar' = \leftubar}^{\ubar} \fluxcontrolVortVort_{\le \Ntop-1} (\ubar',u) \, \rmd \ubar' + \bulkcontrolVortVort_{\le \Ntop-1}(\ubar,u).			 
		\end{split}
	\end{align}
Since $\mastercontroltop(\ubar,u)$ is increasing in its arguments, inserting the already proved estimate \eqref{E:MAINTOPORDERL2ESTIMATESMASTERCONTROLQUANTITY} into \eqref{E:PROOFOFBELOWTOPORDERAPRIORIL2ESTIMATEWAVESTEP2}, using the already proved $L^2$ estimates for the transport variables of Props\,\ref{P:APRIORIL2ESTIMATESFORTRANSPORTVARIABLES}, and applying Gr\"onwall's inequality proves the desired result. 

\hfill $\qed$

\subsection{Proof of Prop.\,\ref{P:APRIORIL2ESTIMATESACOUSTICGEOMETRY}}
\label{SS:PROOFOFENERGYESTIMATESFORACOUSTICGEOMETRY} Throughout this proof, we silently use the trivial inequalities $\totalfluxcontrolvelocity_{[1,N]}  \le \totalcontrolwave_{[1,N]}(\ubar,u)$ and $\strongangularcontrolvelocity_{[1,N]}(\ubar,u) \le  \totalcontrolwave_{[1,N]}(\ubar,u)$, see \eqref{E:VELOCITYFLUXANDSTRONGANGULARL2CONTROLLINGQUANTITY}. Similar remarks apply to the individual terms present in $\mastercontroltop(\ubar,u)$. 

To prove \eqref{E:MAINACOUSTICBELOWTOPESTIMATES}, we insert the already proved estimates for $ \totalcontrolwave_{[1,\Ntop]}(\ubar,u)$ into RHS\,\eqref{E:PRELIMINARYBELOWTOPORDERL2ESTIMATESFOREIKONALFUINCTIONQUANTITIES}. To prove \eqref{E:MAINACOUSTICTOPORDERESTIMATES}, we first use \eqref{E:KEYESTIMATECONTROLLINGINVERSEMUBYINVERSEWEIGHT} to bound $\left\| \tander^{\Ntop} \mytr_{\gtorus}\upchi \right\|_{L^2_{\blowupratetoporderwave+2}\left(\characteristicdiamondtwoarg{[\leftubar,\ubar)}{[\moreinterestingu_1,u]}\right)} \lesssim \left\| \upmu \tander^{\Ntop} \mytr_{\gtorus}\upchi \right\|_{L^2_{\blowupratetoporderwave}\left(\characteristicdiamondtwoarg{[\leftubar,\ubar)}{[\moreinterestingu_1,u]}\right)}$. The desired result follows from inserting the already proved bounds for $\mastercontroltop(\ubar,u)$ into RHS\,\eqref{E:TOPORDERL2ESTIMATEMUCHI}. An identical argument holds for $\angLie_{\tander}^{\Ntop}\upchi$. 

\hfill $\qed$

\section{Improvements of the fundamental quantitative \texorpdfstring{$L^{\infty}$}{sup-norm} bootstrap assumptions} 
\label{S:IMPROVEMENTSOFFUNDAMENTALQUANTITATIVEBOOTSTRAPASSUMPTIONS}
We continue to work under the assumptions of Sect.\,\ref{SS:SILENTFACTS}.
In this short section, 
we derive $L^{\infty}$ estimates that yield an improvement of the fundamental quantitative bootstrap assumptions
stated in Sect.\,\ref{SSS:FUNDAMENTALQUANTITATIVE}.  

\begin{proposition}[Improvement of the fundamental quantitative bootstrap assumptions]
\label{P:IMPROVEMENTOFFUNDAMENTALQUANTITATIVEBOOTSTRAPASSUMPTIONS}
The following estimates hold for 
		$(\ubar,u) \in [\leftubar,\ubarboot) \times [\moreinterestingu_1,\moreinterestingu_2]$:
	\begin{align} \label{E:RECOVEREDLINFINITYESTIMATE}
		\left\| 
	\tander^{[1,\Ntop-6]} \wavearray
\right\|_{L^{\infty}\left(\doublenulltoritwoarg{\ubar}{u}\right)}, \,   \left\| 
	\tander^{\le \Ntop-6} (\vortrenormalized,\GradEnt)
\right\|_{L^{\infty}\left(\doublenulltoritwoarg{\ubar}{u}\right)},  \,   \left\| 
	\tander^{\le \Ntop-7} (\VortVort,\DivGradEnt)
\right\|_{L^{\infty}\left(\doublenulltoritwoarg{\ubar}{u}\right)}, \,  & \leq C \initialsmall
	\end{align}
In particular, if $C \initialsmall < \fundbootsmall$,
then \eqref{E:RECOVEREDLINFINITYESTIMATE} yields a strict improvement of 
\eqref{E:FUNDAMENTALQUANTITATIVEBOOT}.

\end{proposition}

\begin{proof}

We begin with the proof for $\tander^{[1,\Ntop-6]} (\vortrenormalized,\GradEnt)$. Using \eqref{E:H2LINFINITYSOBOLEVEMBEDDINGDOUBLENULLTORUS} and \eqref{E:MAINTORIL2VORTGRADENTBELOWTOPORDERBLOWUP}, we have the estimates $  \left\| 
	\tander^{[1,\Ntop-6]} (\vortrenormalized,\GradEnt)
\right\|_{L^{\infty}\left(\doublenulltoritwoarg{\ubar}{u}\right)} \lesssim \left\| \tander^{\le \Ntop-4} (\vortrenormalized,\GradEnt)\right\|_{L^{2}_0\left(\doublenulltoritwoarg{\ubar}{u}\right)}  \lesssim \initialsmall$. The proof for $\tander^{\le 7}(\VortVort,\DivGradEnt)$ is nearly identical but based instead on \eqref{E:MAINTORIL2MODIFIEDFLUIDVARIABLESBELOWTOPORDERBLOWUP}. 

We now improve the bootstrap assumptions for $\velocityarray = \{v^1,v^2,v^3\} \subset \wavearray = \{v^2,v^2,v^3,\LogDensity,\Ent\}$.

We use the estimate \eqref{E:H2LINFINITYFUNDAMENTALTHEOREMOFCALCULUSPLUSSOBOLEVEMBEDDINGONDOUBLENULLTORUS}, 
the data-assumption \eqref{E:SMALLDATAOFPSIONINITIALDOUBLENULLTORI}, 
Lemma\,\ref{L:COERCIVENESSOFL2CONTROLLINGQUANITIESUNIFIEDCOMMUTATOR},
and the estimate \eqref{E:NONWEIGHTEDL2ESTIMATESFORWAVE}
to conclude that
$\left \|\tander^{[1,\Ntop-6]} \velocityarray \right\|_{L^{\infty}\left(\doublenulltoritwoarg{\ubar}{u}\right)}^2 
\lesssim 
\initialsmalldoublenull^2 
+ 
 \left\| \Lunit \tander^{[1,\Ntop-4]} \velocityarray \right\|_{L^2_0\left(\outgoingcharacteristicsurfacetwoarg{u}{[\leftubar,\ubar)}\right)}^2
\lesssim 
\initialsmalldoublenull^2$
as desired. The estimate \eqref{E:RECOVEREDLINFINITYESTIMATE} for $\LogDensity \in \wavearray$ follows now from \eqref{E:HIGHERORDERTANGENTIALDERIVATIVESOFDENSITY} and the already proved pointwise bounds for $\velocityarray$ and $(\vortrenormalized,\GradEnt)$. We note that the last term on RHS\,\eqref{E:HIGHERORDERTANGENTIALDERIVATIVESOFDENSITY} is handled using the bootstrap assumptions; namely \eqref{E:NONLINEARINEQUALITYRELATINGDATAEPSILONANDBOOTSTRAPEPSILON} implies $\fundbootsmall |\tander^{[1,\Ntop-7]}\controlvars| \lesssim \fundbootsmall^2 \lesssim \initialsmall^{4/3} \lesssim \initialsmall$.
\end{proof}

\section{Existence up to the Cauchy horizon via continuation criteria} \label{S:EXISTENCEUPTOCAUCHYHORIZONBYCONTINUATIONCRITERIA}

In main goal of this section is to state and prove our main theorem, Theorem\,\ref{T:EXISTENCEUPTOCAUCHYHORIZONBYCONTINUATIONCRITERIA}. We show that the solution exists on the characteristic diamond $\characteristicdiamondtwoarg{[\leftubar,0]}{[\moreinterestingu_1,\moreinterestingu_2]}$, which in particular contains a portion of the Cauchy horizon. The proof of the theorem depends on a local well-posedness result for the compressible Euler equations where the data is given on a pair of transversely intersecting $\gfour$-null hypersurfaces. This goes by the name of the \emph{characteristic initial value problem}, and is adapted from the work of \cite{speckyucharacteristicIVP}. We also depend on a continuation criteria for the characteristic IVP, also adapted from \cite{speckyucharacteristicIVP} and are provided in Sect.\,\ref{SS:CONTINUATIONCRITERIACHARACTERISTICIVP}. 

We now define the portion of the Cauchy horizon that we obtain.

\begin{definition}[The Cauchy horizon] \label{D:DEFCH}We define the \textbf{Cauchy horizon} to be the set:
\begin{align} \label{E:DEFCH}
	\CH \eqdef  \ingoingcharacteristicsurfacetwoarg{0}{[\moreinterestingu_1,\moreinterestingu_2]} \cap \{ \newuL \upmu \ge 0\}.
\end{align}
\end{definition}

\begin{theorem}[Existence up to the Cauchy horizon via continuation criteria] \label{T:EXISTENCEUPTOCAUCHYHORIZONBYCONTINUATIONCRITERIA}
Under the assumptions on $\Sigma_0$ of Theorem\,\ref{T:MAINRESULTSFROMSINGULARBOUNDARYPAPER}, which in particular imply that the data-assumptions of Sect.\,\ref{S:ASSUMPTIONSONTHEDATA} hold, we have the following: 

	\medskip
	
\noindent \underline{\textbf{The ingoing eikonal function, ingoing characteristic surfaces, and classical existence in geometric coordinates}}.

\begin{itemize}
	\item There exists a solution to the eikonal equation initial value problem \eqref{E:INGOINGEIKONALEQUATIONINITIALVALUEPROBLEM} with range $[\leftubar,0]$ (see \eqref{E:MOSTNEGATIVEVALUEOFUBAR}),  and we denote its level-set portions, viewed as substes of geometric coordinate space $\R_t \times \R_u \times \T^2$, as follows: $\ingoingcharacteristicsurfacetwoarg{\ubar}{[u_1,u_2]} = \{ (t,u,x^2,x^3) \ | \ \ubar(t,u,x^2,x^3) = \ubar, \, u_1 \le u \le u_2, \, (x^2,x^3) \in \T^2\}$.  More precisely, with $\moreinterestingu_1,\moreinterestingu_2$ as in \eqref{E:MOSTNEGATIVEMOREINTERESTINGVALUEOFU}--\eqref{E:MOSTPOSITIVEMOREINTERESTINGVALUEOFU}, $\ubar$ is defined on the characteristic diamond $\characteristicdiamondtwoarg{[\leftubar,0]}{[\moreinterestingu_1,\moreinterestingu_2]} = \cup_{\ubar \in [\leftubar,0]} \ingoingcharacteristicsurfacetwoarg{\ubar}{[\moreinterestingu_1,\moreinterestingu_2]}$.
	\item The change of variables map $\CHOVgeotodoublenull(t,u,x^2,x^3) \eqdef (\ubar,u,x^2,x^3)$ is a $C^{1,1}_{\textnormal{geo}}$-diffeomorphism from $\characteristicdiamondtwoarg{[\leftubar,0]}{[\moreinterestingu_1,\moreinterestingu_2]}$ onto its image $[\leftubar,0]\times[\moreinterestingu_1,\moreinterestingu_2] \times \T^2$ satisfying $\left\| \CHOVgeotodoublenull\right\|_{C_{\textnormal{geo}}^{1,1}\left(\characteristicdiamondtwoarg{[\leftubar,0]}{[\moreinterestingu_1,\moreinterestingu_2]}\right)} \le C$. Moreover, $\geop{t} \ubar \approx 1$ on $\characteristicdiamondtwoarg{[\leftubar,0]}{[\moreinterestingu_1,\moreinterestingu_2]}$.
	\item For any $\ubar \in [\leftubar,0]$, the ingoing characteristic surfaces $ \ingoingcharacteristicsurfacetwoarg{\ubar}{[\moreinterestingu_1,\moreinterestingu_2]}$ are embedded $C^{1,1}_{\textnormal{geo}}$ submanifolds of spacetime relative to the differential structure of the geometric coordinates. In particular, there exists a $C^{1,1}$ scalar function $\Cartesiantisafunctiononlevelsetsofingoingcharacteristicsarg{\ubar}: 
					[\moreinterestingu_1,\moreinterestingu_2] \times \mathbb{T}^2 \rightarrow \mathbb{R}$,
					depending on $\ubar$, such that $\left\| \Cartesiantisafunctiononlevelsetsofingoingcharacteristicsarg{\ubar} \right\|_{C^{1,1}([\moreinterestingu_1,\moreinterestingu_2] \times \mathbb{T}^2)}  
						 \leq C$ and $\ingoingcharacteristicsurfacetwoarg{\ubar}{[\moreinterestingu_1,\moreinterestingu_2]}
							 =
						\left\lbrace
							(t,u,x^2,x^3)
							\ | \
							t = \Cartesiantisafunctiononlevelsetsofingoingcharacteristicsarg{\ubar}(u,x^2,x^3),
								\,
							(u,x^2,x^3) \in [\moreinterestingu_1,\moreinterestingu_2] \times \mathbb{T}^2
						\right\rbrace$.
	\item The following quantities extend as solutions to the compact set $\characteristicdiamondtwoarg{[\leftubar,0]}{[\moreinterestingu_1,\moreinterestingu_2]}$ as elements of the following spacetime H\"older spaces\footnote{Actually, thanks to the estimates offered by Prop.\,\ref{P:IMPROVEMENTOFAUXILIARYBOOTSTRAP}, the solution enjoys additional regularity in the directions tangent to the outgoing characteristics compared to what we have stated here; we have stated simpler, sub-optimal regularity conclusions only to avoid cluttering the presentation.} with respect to the geometric coordinates, and their corresponding spacetime H\"older norms on $\characteristicdiamondtwoarg{[\leftubar,0]}{[\moreinterestingu_1,\moreinterestingu_2]}$ are $\le C$:
	\begin{align} \label{E:MAINTHEOREMHOLDERREGULARITY}
		\wavearray \in C^{3,1}_{\textnormal{geo}}\left( \characteristicdiamondtwoarg{[\leftubar,0]}{[\moreinterestingu_1,\moreinterestingu_2]}\right), & & \Upsilon \in C^{3,1}_{\textnormal{geo}}\left(\characteristicdiamondtwoarg{[\leftubar,0]}{[\moreinterestingu_1,\moreinterestingu_2]}\right), & & \Lunit^i, \, \upmu \in C^{2,1}_{\textnormal{geo}}\left(\characteristicdiamondtwoarg{[\leftubar,0]}{[\moreinterestingu_1,\moreinterestingu_2]}\right).
	\end{align}
	
	\item The $L^\infty$ estimates of Prop.\,\ref{P:IMPROVEMENTOFAUXILIARYBOOTSTRAP} with $\fundbootsmall$ replaced by $C\initialsmalldoublenull$, and the $L^2$ energy--null-flux estimates of Props.\,\ref{P:TOPORDERAPRIORIL2ESTIMATESMASTERCONTROL},\,\ref{P:APRIORIL2ESTIMATESWAVEVARIABLES}, and\,\ref{P:APRIORIL2ESTIMATESACOUSTICGEOMETRY} hold with $\ubarboot = 0$, i.e., they hold on $\characteristicdiamondtwoarg{[\leftubar,0]}{[\moreinterestingu_1,\moreinterestingu_2]}$.
	\item For $\ubar \in [\leftubar,0]$,
		 we have: 
		\begin{align} \label{E:MAINRESULTSMINVALUEOFMUONFOLIATION}
			\min_{ \ingoingcharacteristicsurfacetwoarg{\ubar}{[\moreinterestingu_1,\moreinterestingu_2]}} \upmu
			& \ge \frac{1}{1.01} |\ubar|.
		\end{align}
		Moreover, within $\ingoingcharacteristicsurfacetwoarg{\ubar}{[\moreinterestingu_1,\moreinterestingu_2]}$,
		the minimum value
		in \eqref{E:MAINRESULTSMINVALUEOFMUONFOLIATION} is achieved by $\upmu$ precisely on the
		set $\argubarnewultorus{\ubar}$ from definition~\ref{E:DEFOFNEWULADAPTEDTORI},
		which is a $C^{0,1}$-embedded torus for $\ubar \in [\leftubar,0)$ and the $C^{1,1}$-embedded torus $\argubarnewultorus{0} = \crease$ when $\ubar = 0$.
		In particular, in $\characteristicdiamondtwoarg{[\leftubar,0]}{[\moreinterestingu_1,\moreinterestingu_2]}$,
		$\upmu$ vanishes precisely along the crease $\argubarnewultorus{0} = \crease$,
		which is a subset of $\ingoingcharacteristicsurfacetwoarg{0}{[\moreinterestingu_1,\moreinterestingu_2]}$.

\end{itemize}

\medskip
	
\noindent \underline{\textbf{A description of the Cauchy horizon and its causal character in geometric coordinates.}}
	\begin{itemize}
		\item The vectorfield $\newuL$ defined by \eqref{E:DOUBLENULLEIKONALFUNCTIONNORMALIZEDNULLVECTORFIELDS} enjoys the following properties:
			\begin{itemize}
				\item $\newuL$ is well-defined on and tangent everywhere on $\CH$, viewed as a subset of geometric coordinate space.
				\item $\newuL$ satisfies $\gfour(\newuL,\newuL) = 0$ and $\newuL u = 1$ on $\CH$.
			\end{itemize}
		\item $\CH$ is a three-dimensional $C^{1,1}_{\textnormal{geo}}$ embedded sub-manifold-with-boundary of $\ingoingcharacteristicsurfacetwoarg{0}{[\moreinterestingu_1,\moreinterestingu_2]}$, and by the previous bullet point, is a $\gfour$-null hypersurface. 
		\item $\upmu$ is increasing along integral curves of $\newuL$, and in particular, is positive everywhere on $\CH \setminus \crease$. In particular, the past boundary of $\newuL$ is the crease $\crease$.
	\end{itemize}

\medskip

\noindent \underline{\textbf{Properties of the change of variables map $\Upsilon$}}.

\begin{itemize}
\item  On $\characteristicdiamondtwoarg{[\leftubar,0]}{[\moreinterestingu_1,\moreinterestingu_2]}$,
		the change of variables map 
		$\Upsilon(t,u,x^2,x^3) = (t,x^1,x^2,x^3)$ is an injection onto its image in Cartesian coordinate space
		satisfying $\| \Upsilon \|_{C_{\textnormal{geo}}^{3,1}(\characteristicdiamondtwoarg{[\leftubar,0]}{[\moreinterestingu_1,\moreinterestingu_2]})} \leq C$.
		In particular, $\Upsilon$ is a homeomorphism from the compact set
		$\characteristicdiamondtwoarg{[\leftubar,0]}{[\moreinterestingu_1,\moreinterestingu_2]}$ onto its image.
	\item With $\rmd_{\textnormal{geo}} \Upsilon$ denoting the Jacobian matrix of $\Upsilon$,
		we have:
		\begin{align} \label{E:MAINTHEOREMGEOTOCARTESIANJACOBIANDETERMINANTESTIMATE}
			\mydet \rmd_{\textnormal{geo}} \Upsilon
			& \approx
				- 
				\upmu.
		\end{align}
		Hence, on $\characteristicdiamondtwoarg{[\leftubar,0]}{[\moreinterestingu_1,\moreinterestingu_2]} \backslash \argubarnewultorus{0}$,
		$\Upsilon$ is a diffeomorphism.
\end{itemize}

\medskip 
\noindent \underline{\textbf{A description of the Cauchy horizon and its causal character in Cartesian coordinates.}} The vectorfield $\uLunit$ defined by \eqref{E:DOUBLENULLCARTESIANTIMENORMALIZEDNULLVECTORFIELDS} enjoys the following properties:
		\begin{itemize}			
			\item In the language of differential geometry, since $\Upsilon$ is a homeomorphism even along the crease, the vectorfield $\uLunit$ is $\Upsilon$-related to the vectorfield $\Lunit + \ingoingmu\Lunit \ubar X - \ingoingmu \angD^A \ubar (\p_A - \Yvfsmallcoeff{A}X^1\p_1)$. That is, for any $(t,x^1,x^2,x^3) \in \Upsilon(\CH)$, we have that the pushforward $\Upsilon_*\uLunit|_{(t,x^1,x^2,x^3)}$ is well defined and satisfies:
			\begin{align} \label{E:PUSHFORWARDOFTNORMALIZEDLBARINCARTESIANCOORDINATES}
				\begin{split}
					\Upsilon_* \uLunit \big|_{(t,x^1,x^2,x^3)} \eqdef \left. \rmd_{\textnormal{geo}}\Upsilon \cdot \uLunit \right|_{\Upsilon^{-1}(t,x^1,x^2,x^3)} = \left. \Lunit + \ingoingmu\Lunit \ubar X - \ingoingmu \angD^A \ubar (\p_A - \Yvfsmallcoeff{A}X^1\p_1)\right|_{(t,x^1,x^2,x^3)}.
				\end{split}
			\end{align}
			Consequently, $\Upsilon_* \uLunit$ is a well-defined $C^{0,1}$-vectorfield and tangent everywhere on $\Upsilon(\CH)$, viewed as a subset of Cartesian coordinate space.
			\item The pushforward $\Upsilon_* \uLunit$  satisfies $\gfour(\Upsilon_* \uLunit,\Upsilon_* \uLunit) = 0$ and $(\Upsilon_* \uLunit) t = 1$ along $\Upsilon(\CH)$. 
			\item $\Upsilon(\CH\setminus \crease)$ is a three-dimensional sub-manifold-with-boundary of class $C^{1,1}$, whereas $\Upsilon(\CH)$ is of class $C^{0,1}$. In particular, it is a $\gfour$-null hypersurface.
			\item $\upmu$ is increasing along integral curves of $\Upsilon_* \uLunit$, and in particular, is positive everywhere on $\Upsilon (\CH \setminus \crease)$.
		\end{itemize}

	\medskip

\noindent \underline{\textbf{A description of the solution's 
singular and regular behavior with respect to Cartesian coordinates}}.
	\begin{itemize}
		\item (Region without singularities).
		On the subset
		$\Upsilon\left(\characteristicdiamondtwoarg{[\leftubar,0]}{[\moreinterestingu_1,\moreinterestingu_2]}
 \backslash \crease \right)$
		of Cartesian coordinate space,
		the solution exists classically with respect to the Cartesian coordinates. In particular, the solution is \textbf{regular with respect to the Cartesian coordinates on the portion of the Cauchy horizon $\CH\setminus \crease$}. 
	\item (The wave singularity).
		The following lower bound holds in 
		$\Upsilon\left(\characteristicdiamondtwoarg{[\leftubar,0]}{[\moreinterestingu_1,\moreinterestingu_2]} \right)$: 
		\begin{align} \label{E:MAINTHEOREMBLOWUPLOWERBOUND}
			|X \RRiemann|
			&
			\geq
			\frac{\blowupdeltadoublenull}{|\bar{\Speed}_{;\LogDensity} + 1|}\frac{1}{\upmu}
		\end{align}
		where $\blowupdeltadoublenull > 0$ is the data-parameter from \eqref{E:DELTASTARDOUBLENULLDEF}.		In particular, if $q \in \Upsilon\left( \argubarnewultorus{0} \right)$,
		then since 
		$\crease = \argubarnewultorus{0} \subset \characteristicdiamondtwoarg{[\leftubar,0]}{[\moreinterestingu_1,\moreinterestingu_2]}$ 
		by \eqref{E:IMPROVEDLEVELSETSTRUCTUREANDLOCATIONOFMIN}, 
		and since $\upmu = 0$ along 
		$\Upsilon\left(\argubarnewultorus{0}\right)$,
		it follows that
		$|X \RRiemann|(q') \to \infty$ as $q' \rightarrow q$ in
		$\Upsilon\left(\characteristicdiamondtwoarg{[\leftubar,0]}{[\moreinterestingu_1,\moreinterestingu_2]}
		\backslash \argubarnewultorus{0} \right)$.
		\item (Regular behavior\footnote{Here we have only highlighted some of the quantities that remain $L^\infty$-bounded on $\Upsilon\left( \characteristicdiamondtwoarg{[\leftubar,0]}{[\moreinterestingu_1,\moreinterestingu_2]}\right)$. We refer to Prop.\,\ref{P:IMPROVEMENTOFAUXILIARYBOOTSTRAP} for more comprehensive results.} along the outgoing characteristics) The derivatives of $\wavearray$ up to order $\Ntop -6$ with respect to the vectorfields in the $\nullhyparg{u}$-tangent commutation set $\Tanset$ defined in \eqref{E:COMMUTATIONVECTORFIELDS} are $L^\infty$-bounded on $\Upsilon\left(\characteristicdiamondtwoarg{[\leftubar,0]}{[\moreinterestingu_1,\moreinterestingu_2]}
\right)$. 
	\end{itemize}

\end{theorem}

\begin{proof} \hfill

\noindent \underline{\textbf{Classical existence up to the Cauchy horizon in geometric coordinates:}}
Local well-posedness for the characteristic initial value problem provided in Theorem\,\ref{T:CONTINUATIONCRITERIACHARACTERISTICIVP} implies that there exists a $\ubar_{\textnormal{Local}} \in (\frac{3}{4} \leftubar,0)$
	such that the solution variables 
$\wavearray$,
$u$, 
and
$\ubar$
are classical solutions on
$\characteristicdiamondtwoarg{[\leftubar,\ubar_{\textnormal{Local}})}{[\moreinterestingu_1,\moreinterestingu_2]}$
and such that all of the bootstrap assumptions from Sect.\,\ref{S:BOOTSTRAPEVERYTHINGEXCEPTENERGIES}
hold on $\characteristicdiamondtwoarg{[\leftubar,\ubar_{\textnormal{Local}})}{[\moreinterestingu_1,\moreinterestingu_2]}$.
Let $\ubar_{\textnormal{Max}}$ be the supremum over all such $\ubar_{\textnormal{Local}}$.
Then the solution exists classically 
and satisfies all the bootstrap assumptions from Sect.\,\ref{S:BOOTSTRAPEVERYTHINGEXCEPTENERGIES}
on $\characteristicdiamondtwoarg{[\leftubar,\ubar_{\textnormal{Max}})}{[\moreinterestingu_1,\moreinterestingu_2]}$.
If it were true that $\ubar_{\textnormal{Max}} < 0$, then \eqref{E:KEYESTIMATECONTROLLINGINVERSEMUBYINVERSEWEIGHT}
would imply that $\upmu$ is uniformly positive on 
$\characteristicdiamondtwoarg{[\leftubar,\ubar_{\textnormal{Max}})}{[\moreinterestingu_1,\moreinterestingu_2]}$,
and then Lemma\,\ref{L:RELATIONSHIPBETWEENCARTESIANPARTIALDERIVATIVESANDSMOOTHGEOMETRICCOMMUTATORS} 
would imply that the solution's Cartesian partial derivatives are uniformly bounded 
on $\characteristicdiamondtwoarg{[\leftubar,\ubar_{\textnormal{Max}})}{[\moreinterestingu_1,\moreinterestingu_2]}$.
Moreover,
in the previous sections, 
we derived strict 
(assuming that $\initialsmalldoublenull$ and $\mathring{\upalpha}$ are sufficiently small)
improvements of all the bootstrap assumptions.
By the continuation criteria of Theorem\,\ref{T:CONTINUATIONCRITERIACHARACTERISTICIVP}, there is a $\Delta > 0$ 
with $\ubar_{\textnormal{Max}} + \Delta < 0$
such that
the solution exists classically 
and satisfies the bootstrap assumptions from Sect.\,\ref{S:BOOTSTRAPEVERYTHINGEXCEPTENERGIES}
on $\characteristicdiamondtwoarg{[\leftubar,\ubar_{\textnormal{Max}} + \Delta)}{[\moreinterestingu_1,\moreinterestingu_2]}$.
That is impossible in view of the definition of $\ubar_{\textnormal{Max}}$. Hence, 
$\ubar_{\textnormal{Max}} = 0$.

We now prove that $\argubarnewultorus{0} = \crease$. Using \eqref{E:CONVENIENTIDENTITYFORULUNIT} and \eqref{E:MAGNITUDEOFINNERPRODUCTOFNEWLANDNEWULAPPROXIMATELYMU}--\eqref{E:RECIPROCALULUNITAPPLIEDTOTIMEFUNCTIONAPPROXIMATELYMU}, we see that $\newuL \upmu \big|_{\upmu = 0} = \muX \upmu \big|_{\upmu = 0}$. By \eqref{E:DATAFORINGOINGEIKONALINITIALVALUEPROBLEM}, it follows that $\crease \subset \argubarnewultorus{0}$. To prove the other direction, suppose $p \in \argubarnewultorus{0}$ and so $\newuL \upmu (p)= 0$. By \eqref{E:CONVENIENTIDENTITYFORULUNIT}, we have $0 = \left\{\ReciprocaluLunitAppliedtoTimeFunction \Lunit \upmu + \muX \upmu - \MagnitueofinnerproductofnewLandnewuL \angD^A \ubar \geop{x^A}\upmu \right\}(p)$. By \eqref{E:MUTRANSVERSALCONVEXITY}, it follows that $\newuL$ is transverse to level sets of $\newuL \upmu$ within $\ingoingcharacteristicsurfacetwoarg{0}{[\moreinterestingu_1,\moreinterestingu_2]}$. By the fundamental theorem of calculus, it then follows that $\upmu > 0$ on $\ingoingcharacteristicsurfacetwoarg{0}{[\moreinterestingu_1,\moreinterestingu_2]} \cap \{\newuL \upmu > 0\}$ and $\ingoingcharacteristicsurfacetwoarg{0}{[\moreinterestingu_1,\moreinterestingu_2]} \cap \{\newuL \upmu < 0\}$. Since $\upmu$ achieves its minima on $\argubarnewultorus{0}$, the desired inclusion $\argubarnewultorus{0} \subset \crease$ once gain follows from \eqref{E:DATAFORINGOINGEIKONALINITIALVALUEPROBLEM} and \eqref{E:MAGNITUDEOFINNERPRODUCTOFNEWLANDNEWULAPPROXIMATELYMU}--\eqref{E:RECIPROCALULUNITAPPLIEDTOTIMEFUNCTIONAPPROXIMATELYMU}.

The remaining conclusions of the theorem on the characteristic diamond $\characteristicdiamondtwoarg{[\leftubar,0]}{[\moreinterestingu_1,\moreinterestingu_2]}$ now follow from
\eqref{E:KEYESTIMATECONTROLLINGINVERSEMUBYINVERSEWEIGHT} and the statements just below it,
\eqref{E:LOWERBOUNDONMAGNITUDEOFXRPLUS},
Lemmas\,\ref{L:PROPERTIESANDDIFFEOMORPHICEXTENSIONOFDOUBLENULLCOORDINATES} and\,\ref{L:CONTINUOUSEXTNESION},
Props.\,\ref{P:IMPROVEMENTOFAUXILIARYBOOTSTRAP},
\ref{P:HOMEOMORPHICANDDIFFEOMORPHICEXTENSIONOFCARTESIANCOORDINATES},\,\ref{P:TOPORDERAPRIORIL2ESTIMATESMASTERCONTROL}, 
\ref{P:APRIORIL2ESTIMATESWAVEVARIABLES},
\ref{P:APRIORIL2ESTIMATESACOUSTICGEOMETRY},
and\,\ref{P:IMPROVEMENTOFFUNDAMENTALQUANTITATIVEBOOTSTRAPASSUMPTIONS},
with $\ubarboot= 0$ in all these results.

\medskip

\noindent \underline{\textbf{Description of the Cauchy horizon in geometric coordinates:}}
The properties of $\newuL$ follow from \eqref{E:DOUBLENULLALLTHENEWVECTORFIELDSARENULL} and \eqref{E:GEOMETRICCOORDINATENULLVECTORFIELDSAPPLIEDTOGEOMETRICCOORDINATES}. The $C^{1,1}_{\textnormal{geo}}$ regularity of $\CH$ follows from $\left\| \CHOVgeotodoublenull\right\|_{C_{\textnormal{geo}}^{1,1}\left(\characteristicdiamondtwoarg{[\leftubar,0]}{[\moreinterestingu_1,\moreinterestingu_2]}\right)} \le C$ since $\CH \subset \ingoingcharacteristicsurfacetwoarg{0}{[\moreinterestingu_1,\moreinterestingu_2]}$ and the fact that $0$ is a regular value of $\ubar$, which is a consequence of the diffeomorphism property of $\CHOVgeotodoublenull$. By \eqref{E:BEHAVIOROFMUATTHECREASEINDIRECTIONSOFNEWUL} and \eqref{E:MUTRANSVERSALCONVEXITY}, it follows that $\newuL \upmu > 0$ on $\CH \setminus \crease$. The monotonicity properties of $\upmu$ easily follow.

\medskip

\noindent \underline{\textbf{Properties of $\Upsilon$:}} These results follow from Prop.\,\ref{P:HOMEOMORPHICANDDIFFEOMORPHICEXTENSIONOFCARTESIANCOORDINATES} with $\ubarboot = 0$.

\medskip

\noindent \underline{\textbf{Description of the Cauchy horizon in Cartesian coordinates:}} Consider the identity \eqref{E:CONVENIENTIDENTITYFORULUNIT}. Since $\frac{1}{\ReciprocaluLunitAppliedtoTimeFunction} \muX = \frac{\ingoingmu}{\upmu \ReciprocalLunitAppliedtoTimeFunction} \muX =  \frac{\ingoingmu \Lunit \ubar}{\upmu}\muX $ by \eqref{E:RATIOOFNULLGEOSICINNERPRODUCTANDFOLIATIONDENSITY}, \eqref{E:PUSHFORWARDOFTNORMALIZEDLBARINCARTESIANCOORDINATES} follows from \eqref{E:GEOMETRICVECTORFIELDSINTERMSOFCARTESIANONES} and \eqref{E:CONVENIENTIDENTITYFORULUNIT}. The regularity $C^{0,1}$ regularity of $\Upsilon_* \uLunit$ follows from Lemma\,\ref{L:PROPERTIESANDDIFFEOMORPHICEXTENSIONOFDOUBLENULLCOORDINATES} and Prop.\,\ref{P:HOMEOMORPHICANDDIFFEOMORPHICEXTENSIONOFCARTESIANCOORDINATES}. The tangency of $\Upsilon_* \uLunit$ to $\Upsilon(\CH)$, as well as $\gfour(\Upsilon_* \uLunit,\Upsilon_* \uLunit) = 0$, follow from \eqref{E:DOUBLENULLALLTHENEWVECTORFIELDSARENULL} and  \eqref{E:LUNITANDULUNITAPPLIEDTOEIKONALANDCARTESIANTIME}. The $C^{1,1}$ regularity of $\Upsilon(\CH \setminus \crease)$ follows from Lemma\,\ref{L:PROPERTIESANDDIFFEOMORPHICEXTENSIONOFDOUBLENULLCOORDINATES} and Prop.\,\ref{P:HOMEOMORPHICANDDIFFEOMORPHICEXTENSIONOFCARTESIANCOORDINATES}, specifically the fact that $\Upsilon$ is a diffeomorphism away from the crease. The $C^{0,1}$ regularity of $\Upsilon(\CH)$ follows from the fact that $\CH$ is the flowout of the crease $\crease$ by the Lipschitz vectorfield $\Upsilon_* \uLunit$. We remark that flows of $C^{0,1}$ vectorfields are $C^{1,1}$ in the flow-parameter, which in this case may be identified with $t$, but are only $C^{0,1}$ maps of the data \cite{rampazzo2007frobenius}.  

\medskip

\noindent \underline{\textbf{Description of the singular and regular behavior of the solution in Cartesian coordinates:}} These results follow from the conclusions proved above in geometric coordinates, specifically \eqref{E:LOWERBOUNDONMAGNITUDEOFXRPLUS} and \eqref{E:RECOVEREDLINFINITYESTIMATE}, as well as invoking Lemma\,\ref{L:PROPERTIESANDDIFFEOMORPHICEXTENSIONOFDOUBLENULLCOORDINATES}. 

\end{proof}

\subsection{Continuation criteria for the characteristic initial value problem} \label{SS:CONTINUATIONCRITERIACHARACTERISTICIVP}

We now state the results needed from \cite{speckyucharacteristicIVP}, adapted to the current setting.\footnote{The results of 
\cite{speckyucharacteristicIVP} require the characteristic data to satisfy constraint equations on the transversal characteristic surfaces as
well as corner compatibility conditions on their intersection. These are automatically satisfied in the context of the present paper because
our characteristic data is induced by an existing $C^{\infty}$ solution that originally arose from spacelike data, and
the constraints and compatibility conditions are consequences of the Euler equations being satisfied to all orders along the
characteristic hypersurfaces.}

\begin{theorem}[Existence and continuation criteria for the characteristic initial value problem; \cite{speckyucharacteristicIVP}] \label{T:CONTINUATIONCRITERIACHARACTERISTICIVP} \hfill

\medskip
	
\noindent \underline{\textbf{Local existence for the characteristic IVP}}. 
Recall that the classical solution is $C^{\infty}$ on $\twoargMrough{[\timefunction_0,0],[- \rightu,\leftu]}{0}$
and therefore induces $C^{\infty}$ data on the characteristic surfaces
$\ingoingcharacteristicsurfacetwoarg{\leftubar}{[\moreinterestingu_1,\moreinterestingu_2]}, \, \outgoingcharacteristicsurfacetwoarg{\moreinterestingu_1}{[\leftubar,0]}$; see Appendix~\ref{A:DATASSUMPTIONS}.
Then there exists a  
$\ubarlocal \in (\leftubar,0)$, depending on the data, such the characteristic data launch a unique $C^{\infty}$ solution
on $\characteristicdiamondtwoarg{[\leftubar,\ubarlocal)}{[\moreinterestingu_1,\moreinterestingu_2]}$,
and the change of variables map $(u,\ubar,x^2,x^3) \rightarrow (t,x^1,x^2,x^3)$
is a diffeomorphism from $\characteristicdiamondtwoarg{[\leftubar,\ubarlocal)}{[\moreinterestingu_1,\moreinterestingu_2]}$
onto its image in Cartesian coordinate space.

\medskip
	
\noindent \underline{\textbf{Continuation criteria for the characteristic IVP}}.

Assume the following: 
\begin{itemize}
	\item The assumptions of Theorem\,\ref{T:EXISTENCEUPTOCAUCHYHORIZONBYCONTINUATIONCRITERIA} hold.
	\item $\ubarboot < 0$.
	\item The ingoing eikonal function $\ubar$, the wave variables $\wavearray$, and the acoustic geometry variables $u, \, \chi$, etc. are classical solutions on $\characteristicdiamondtwoarg{[\leftubar,\ubarboot)}{[\moreinterestingu_1,\moreinterestingu_2]}$. 
	\item The bootstrap assumptions of Sects.\,\ref{S:BOOTSTRAPEVERYTHINGEXCEPTENERGIES} and\,\ref{SS:BOOTSTRAPASSUMPTIONSFORTHEWAVEENERGIES} hold for $(\ubar,u) \in [\leftubar,\ubarboot)\times[\moreinterestingu_1,\moreinterestingu_2]$. 
\end{itemize}

Then there exists a $\Delta \in (0,|\ubarboot|)$ such that the ingoing eikonal function $\ubar$, the solution variables $\wavearray, \, u$, and all of the other geometric quantities defined throughout the article can be uniquely extended to a strictly larger region of the form $\characteristicdiamondtwoarg{[\leftubar,\ubarboot + \Delta)}{[\moreinterestingu_1,\moreinterestingu_2]}$ on which all of the bootstrap
assumptions of Sects.\,\ref{S:BOOTSTRAPEVERYTHINGEXCEPTENERGIES} and\,\ref{SS:BOOTSTRAPASSUMPTIONSFORTHEWAVEENERGIES} hold.

\end{theorem}

\begin{proof}[Discussion of the proof]
The existence aspects of the theorem are proved in \cite{speckyucharacteristicIVP}.
In \cite{speckyucharacteristicIVP}, it was also shown that the continuation criteria allow one to extend the solution
as a $C^{\infty}$ solution to a larger characteristic diamond of the form 
$\characteristicdiamondtwoarg{[\leftubar,\ubarboot + \Delta)}{[\moreinterestingu_1,\moreinterestingu_2]}$.
That result requires as a hypothesis that the change of variables map $(u,\ubar,x^2,x^3) \rightarrow (t,x^1,x^2,x^3)$ extends to be a diffeomorphism  
on the closure $\characteristicdiamondtwoarg{[\leftubar,\ubarlocal]}{[\moreinterestingu_1,\moreinterestingu_2]}$.
In the present context, this hypothesis follows from combining 
Prop.\,\ref{P:SHARPCONTROLOFMUANDDERIVATIVES}, which shows that $\upmu$ is uniformly positive on
$\characteristicdiamondtwoarg{[\leftubar,\ubarboot)}{[\moreinterestingu_1,\moreinterestingu_2]}$
(the positivity depends on $\ubarboot$, which is by assumption not equal to $0$)
with the results for $\Upsilon$ proved in Prop.\,\ref{P:HOMEOMORPHICANDDIFFEOMORPHICEXTENSIONOFCARTESIANCOORDINATES} 
and the results for $\CHOVgeotodoublenull$ proved in Lemma~\ref{L:PROPERTIESANDDIFFEOMORPHICEXTENSIONOFDOUBLENULLCOORDINATES}.
Moreover,
in the previous sections, 
we derived strict 
(assuming that $\initialsmalldoublenull$ and $\mathring{\upalpha}$ are sufficiently small)
improvements of all the bootstrap assumptions
on $\characteristicdiamondtwoarg{[\leftubar,\ubarlocal)}{[\moreinterestingu_1,\moreinterestingu_2]}$.
Hence, the fact that the bootstrap assumptions of Sects.\,\ref{S:BOOTSTRAPEVERYTHINGEXCEPTENERGIES} and\,\ref{SS:BOOTSTRAPASSUMPTIONSFORTHEWAVEENERGIES} hold for $(\ubar,u) \in [\leftubar,\ubarboot+ \Delta)\times[\moreinterestingu_1,\moreinterestingu_2]$
follows from continuity and the fact that the solution is $C^{\infty}$ on 
$\characteristicdiamondtwoarg{[\leftubar,\ubarboot + \Delta)}{[\moreinterestingu_1,\moreinterestingu_2]}$.

\end{proof}

\appendix 

\section{The characteristic data assumptions hold} \label{A:DATASSUMPTIONS}

In this appendix, we prove that the initial data assumptions of Sect.\,\ref{S:ASSUMPTIONSONTHEDATA} hold provided the statements and conclusions of Theorems\,\ref{T:MAINRESULTSFROMSINGULARBOUNDARYPAPER} and\,\ref{T:CONSTRUCTIONOFTHEINGOINGEIKONALFUNCTIONCONTRACTIONMAPPING} hold. 

\begin{proposition}[The classical solution on $\twoargMrough{[\timefunction_0,0],[- \rightu,\leftu]}{0}$ induces characteristic data]
Assume the assumptions and conclusions of Theorems\,\ref{T:MAINRESULTSFROMSINGULARBOUNDARYPAPER} and\,\ref{T:CONSTRUCTIONOFTHEINGOINGEIKONALFUNCTIONCONTRACTIONMAPPING} hold. Then there exist constants $\leftubar, \, \moreinterestingu_1, \, \moreinterestingu_2$ satisfying \eqref{E:MOSTNEGATIVEVALUEOFUBAR}--\eqref{E:MOSTPOSITIVEMOREINTERESTINGVALUEOFU} such that \eqref{E:INITIALINGOINGCHARACTERISTICSURFACE}--\eqref{E:INITIALOUTGOINGCHARACTERISTICSURFACE} hold and such that the classical solution on $\twoargMrough{[\timefunction_0,0],[- \rightu,\leftu]}{0}$ induces data on the initial characteristic surfaces $\ingoingcharacteristicsurfacetwoarg{\leftubar}{[\moreinterestingu_1,\moreinterestingu_2]}, \, \outgoingcharacteristicsurfacetwoarg{\moreinterestingu_1}{[\leftubar,0]}$ satisfying all of the qualitative and quantitative assumptions of Sect.\,\ref{S:ASSUMPTIONSONTHEDATA} with $\initialsmalldoublenull = C \initialsmall$.

\end{proposition}

\begin{proof} We will prove there exist constants $\leftubar,\, \moreinterestingu_1,\, \moreinterestingu_2$ such that \eqref{E:MOSTNEGATIVEVALUEOFUBAR}--\eqref{E:CONDITIONTHATGUARANTEESNONTRIVIALPARTOFCAUCHYHORIZON} hold. We will then prove that the classical solution on $\twoargMrough{[\timefunction_0,0],[- \rightu,\leftu]}{0}$ induces the $L^2$ assumptions for $\wavearray$ \eqref{E:TRANSVERSALDERIVATIVEOFTANGENTIALL2NORMSOFWAVEVARIABLESSMALLALONGINITIALINGOINGNULLSURFACE}--\eqref{E:SMALLDATAOFPSIONINITIALDOUBLENULLTORI} and the $L^2$ assumptions for the top order eikonal function quantities \eqref{E:L2DATASSUMPTIONFORFULLYMODIFIEDCHITOPORDER}--\eqref{E:L2DATASSUMPTIONFORPARTIALLYMODIFIEDCHITOPORDER}. The bounds for the below-top-order eikonal function quantities can be derived from \eqref{E:TRANSVERSALDERIVATIVEOFTANGENTIALL2NORMSOFWAVEVARIABLESSMALLALONGINITIALINGOINGNULLSURFACE}--\eqref{E:SMALLDATAOFPSIONINITIALDOUBLENULLTORI} and \eqref{E:L2DATASSUMPTIONFORFULLYMODIFIEDCHITOPORDER}--\eqref{E:L2DATASSUMPTIONFORPARTIALLYMODIFIEDCHITOPORDER}, see Footnote\,\ref{FN:INTUITIONONMUDATAASSUMPTIONS}. Finally, the $L^\infty$ estimates on $\doublenulltoritwoarg{\leftubar}{u}$ of the wave variables, as well as the localized assumptions from Sect.\,\ref{SSS:LOCALIZEDDATAASSUMPTIONSFORMUANDDERIVATIVES}, follow from the $C^{k,1}_{\textnormal{geo}}$ estimates of Theorems\,\ref{T:MAINRESULTSFROMSINGULARBOUNDARYPAPER} and\,\ref{T:CONSTRUCTIONOFTHEINGOINGEIKONALFUNCTIONCONTRACTIONMAPPING}. 

We now begin to prove the $L^2$ data assumptions on through a series of steps.

 \hfill

\noindent \underline{\textbf{A weighted energy scheme on $\twoargMrough{[\timefunction_0,0],[- \rightu,\leftu]}{0}$:}}

A technical issue in proving that the classical solution on $\twoargMrough{[\timefunction_0,0],[- \rightu,\leftu]}{0}$ induces data on initial characteristic surfaces is that the $L^2$-analysis of \cite{abbrescia2022emergence} features \emph{singular unweighted} energy estimates, whereas the $L^2$-analysis of the current paper features \emph{regular weighted} energy estimates. We begin by sketching a proof of how one would adapt the weighted energy scheme from the current paper to $\twoargMrough{[\timefunction_0,0],[- \rightu,\leftu]}{0}$. First, we recall that the $L^2$-analysis of \cite{abbrescia2022emergence} used energies on the \emph{spacelike} hypersurfaces $ \hypthreearg{\timefunction}{[- \rightu,\leftu]}{0}$, fluxes on the truncated outgoing characteristics $\nullhypthreearg{0}{u}{[\timefunction_0,0]}$ (see \eqref{E:NULLHYPERSURFACEROUGHTRUNCATED}), and spacetime bulk energies on $\twoargMrough{[\timefunction_0,0],[- \rightu,\leftu]}{0}$: \begin{subequations}
	\begin{align}
		\mathbb{E}[\tander^N\Psi](\timefunction,u) & \approx \int_{\hypthreearg{\timefunction}{[- \rightu,u]}{0}} 2 (\muX \tander^N \Psi)^2 + \upmu (\Lunit \tander^N \Psi)^2 + \upmu \left| \angrmd \tander^N \Psi \right|_{\gtorus}^2, \label{E:ENERGYONROUGHFOLIATIONS} \\
		\mathbb{F}[\tander^N \Psi](\timefunction,u) & \approx \int_{\nullhypthreearg{0}{u}{[\timefunction_0,\timefunction]}}  (\Lunit \tander^N \Psi)^2 + \upmu \left| \angrmd \tander^N \Psi \right|_{\gtorus}^2, \label{E:NULLFLUXONROUGHFOLIATIONTRUNCATEDOUTGOINGCHARACTERISTICS} \\
		\mathbb{K}[\tander^N \Psi] (\timefunction,u) & \approx \int_{\twoargMrough{[\timefunction_0,\timefunction],[- \rightu,u]}{0}} \left|\angrmd \tander^N \Psi\right|_{\gtorus}^2, \label{E:BLUKENERGYONROUGHCLASSICALDEVELOPMENT}
	\end{align} 
\end{subequations}
where we have suppressed the volume forms in \eqref{E:ENERGYONROUGHFOLIATIONS}--\eqref{E:BLUKENERGYONROUGHCLASSICALDEVELOPMENT} for convenience. In particular, the $L^2$-controlling quantities \emph{do not feature a singular weight that vanishes at the crease} like the ones we introduced in Defs.\,\ref{D:NULLFLUXESFORTHEWAVEANDTRANSPORTVARIABLES} and\,\ref{D:SPACETIMECOERCIVEINTEGRALS}. In fact, at the top and mid-to-top order, the estimates we derived in \cite{abbrescia2022emergence} allowed the energies to potentially blow up. That is, for some inform positive constant $M_*$ and $N_* = \lceil \tfrac{1}{2}M_*\rceil$, we have:
\begin{subequations}
	\begin{align}
		\mathbb{E}[\tander^{\Ntop}\Psi](\timefunction,u) + \mathbb{F}[\tander^{\Ntop} \Psi](\timefunction,u) + \mathbb{K}[\tander^{\Ntop} \Psi] (\timefunction,u) & \lesssim \initialsmall^2 |\timefunction|^{-M_*}, \label{E:TOPORDERENERGYBLOWUPROUGHFOLIATIONS} \\
		\mathbb{E}[\tander^{\Ntop-1}\Psi](\timefunction,u) + \mathbb{F}[\tander^{\Ntop-1} \Psi](\timefunction,u) + \mathbb{K}[\tander^{\Ntop-1} \Psi] (\timefunction,u) & \lesssim \initialsmall^2 |\timefunction|^{-M_* + 2}, \label{E:BELOWTOPORDERENERGYBLOWUPROUGHFOLIATIONS} \\
		& \vdots \notag \\
		\mathbb{E}[\tander^{\Ntop-N_*}\Psi](\timefunction,u) + \mathbb{F}[\tander^{\Ntop-N_*} \Psi](\timefunction,u) + \mathbb{K}[\tander^{\Ntop-N_*} \Psi] (\timefunction,u) & \lesssim \initialsmall^2, \label{E:REGULARENERGYESTIMATESROUGHFOLIATIONS}
	\end{align}
\end{subequations}
see \cite{abbrescia2022emergence}*{Prop.~25.1}. Based on \eqref{E:MAINRESULTSMINVALUEOFMUONROUGHFOLIATION} and \eqref{E:TOPORDERENERGYBLOWUPROUGHFOLIATIONS}--\eqref{E:REGULARENERGYESTIMATESROUGHFOLIATIONS}, one sees that if one wanted to use the \emph{weighted} $L^2$-analysis used in the current paper on the region $\twoargMrough{[\timefunction_0,0],[- \rightu,\leftu]}{0}$, the appropriate weight would be $|\timefunction|$. That is, for $\blowuprateofwaveWRTnumberofcommutations(N),\blowuprateoftransportWRTnumberofcommutations(N)$ as in \eqref{E:BLOWUPRATEOFWAVEWRTNUMBEROFCOMMUTATORS}--\eqref{E:BLOWUPRATEOFTRANSPORTWRTNUMBEROFCOMMUTATORS}, we would have\footnote{In deriving weighted estimates on $\twoargMrough{[\timefunction_0,\timefunction],[- \rightu,u]}{0}$, we technically would not need \eqref{E:WEIGHTEDIMPRECISEFULLYMODQUANTFLUXROUGHCLASSICALDEVELOPMENT}--\eqref{E:WEIGHTEDIMPRECISEFULLYMODQUANTBULKROUGHCLASSICALDEVELOPMENT}. The reason for this is that deriving $L^2$-estimates for $\upmu \tander^N \mytr_{\gtorus}\upchi$ is much easier on $\hypthreearg{\timefunction}{[- \rightu,u]}{0}$ than it is on $\ingoingcharacteristicsurfacetwoarg{\ubar}{[\moreinterestingu_1,u]}$ due to the coercive control of $\upmu( \Lunit \tander^N \Psi)^2$ in \eqref{E:WEIGHTEDENERGYONROUGHFOLIATIONS}. However, we still consider it in \eqref{E:WEIGHTEDIMPRECISEFULLYMODQUANTFLUXROUGHCLASSICALDEVELOPMENT}--\eqref{E:WEIGHTEDIMPRECISEFULLYMODQUANTBULKROUGHCLASSICALDEVELOPMENT} because it makes it slightly easier to prove the first inequality in \eqref{E:L2DATASSUMPTIONFORFULLYMODIFIEDCHITOPORDER}.} \footnote{We use \eqref{E:WEIGHTEDENERGYONROUGHFOLIATIONSWITHLOSS} as an analog for \eqref{E:DERIVATIVELOSINGL2ESTIMATEFORTANGENTIALDERIVATIVES}.}\begin{subequations}
	\begin{align}
		\mathbb{E}[\tander^N\velocityarray](\timefunction,u) & \approx  \sum_{i=1}^3 \int_{\hypthreearg{\timefunction}{[- \rightu,u]}{0}} \left\{ 2 (\muX \tander^N v^i)^2 + \upmu (\Lunit \tander^N v^i)^2 + \upmu \left| \angrmd \tander^N v^i \right|_{\gtorus}^2 \right\} |\timefunction|^{\blowuprateofwaveWRTnumberofcommutations(N)}, \label{E:WEIGHTEDENERGYONROUGHFOLIATIONS} \\
		\mathbb{F}[\tander^N \velocityarray](\timefunction,u) & \approx   \sum_{i=1}^3 \int_{\nullhypthreearg{0}{u}{[\timefunction_0,\timefunction]}} \left\{  (\Lunit \tander^N v^i)^2 + \upmu \left| \angrmd \tander^N v^i \right|_{\gtorus}^2\right\} |\timefunction|^{\blowuprateofwaveWRTnumberofcommutations(N)}, \label{E:WEIGHTEDNULLFLUXONROUGHFOLIATIONTRUNCATEDOUTGOINGCHARACTERISTICS} \\
		\mathbb{K}[\tander^N \velocityarray] (\timefunction,u) & \approx    \sum_{i=1}^3 \int_{\twoargMrough{[\timefunction_0,\timefunction],[- \rightu,u]}{0}}  \left|\angrmd \tander^N v^i \right|_{\gtorus}^2  |\timefunction|^{\blowuprateofwaveWRTnumberofcommutations(N)},\label{E:WEIGHTEDBLUKENERGYONROUGHCLASSICALDEVELOPMENT} \\
		\mathbb{E}^{(\textnormal{Loss})}[\tander^N\velocityarray] (\timefunction,u) & =  \sum_{i=1}^3 \int_{\hypthreearg{\timefunction}{[- \rightu,u]}{0}} |\tander^N v^i|^2 |\timefunction|^{\blowuprateofwaveWRTnumberofcommutations(N)}, \label{E:WEIGHTEDENERGYONROUGHFOLIATIONSWITHLOSS} \\
		\mathbb{E}[\tander^N(\vortrenormalized,\GradEnt)](\timefunction,u) & \approx \int_{\hypthreearg{\timefunction}{[- \rightu,u]}{0}}  \upmu \left| \tander^N (\vortrenormalized,\GradEnt)\right|^2 |\timefunction|^{\blowuprateoftransportWRTnumberofcommutations(N)}, \label{E:WEIGHTEDENERGYONROUGHFOLIATIONSVORTGRADENT} \\
		\mathbb{F}[\tander^N (\vortrenormalized,\GradEnt)](\timefunction,u) & \approx  \int_{\nullhypthreearg{0}{u}{[\timefunction_0,\timefunction]}} \left| \tander^N (\vortrenormalized,\GradEnt)\right|^2 |\timefunction|^{\blowuprateoftransportWRTnumberofcommutations(N)}, \label{E:WEIGHTEDNULLFLUXONROUGHFOLIATIONTRUNCATEDOUTGOINGCHARACTERISTICSVORTGRADENT} \\
		\mathbb{E}[\tander^N(\VortVort,\DivGradEnt)](\timefunction,u) & \approx 	\int_{\hypthreearg{\timefunction}{[- \rightu,u]}{0}}  \upmu \left| \tander^N (\VortVort,\DivGradEnt)\right|^2 |\timefunction|^{\blowuprateofacousticgeoWRTnumberofcommutations(N)}, \label{E:WEIGHTEDENERGYONROUGHFOLIATIONSVORTVORTDIVGRADENT} \\
		\mathbb{F}[\tander^N (\VortVort,\DivGradEnt)](\timefunction,u) & \approx  \int_{\nullhypthreearg{0}{u}{[\timefunction_0,\timefunction]}} \left| \tander^N (\VortVort,\DivGradEnt)\right|^2 |\timefunction|^{\blowuprateofacousticgeoWRTnumberofcommutations(N)}, \label{E:WEIGHTEDNULLFLUXONROUGHFOLIATIONTRUNCATEDOUTGOINGCHARACTERISTICSVORTVORTDIVGRADENT} \\				
		\mathbb{X}^{(\textnormal{Precise})}(\timefunction,u) & \approx \max_{\tander^{\Ntop} \in \mathfrak{P}^{(\Ntop)}}   \int_{\hypthreearg{\timefunction}{[- \rightu,u]}{0}} \sumofbrevexvisquared^2 \left(\fullymodquant{\tander^{\Ntop}}\right)^2 |\timefunction|^{\blowupratetoporderwave}, \label{E:WEIGHTEDPRECISEFULLYMODQUANTFLUXROUGHCLASSICALDEVELOPMENT} \\
		\mathbb{KX}^{(\textnormal{Precise})}(\timefunction,u) & \approx  \max_{\tander^{\Ntop} \in \mathfrak{P}^{(\Ntop)}}  \int_{\twoargMrough{[\timefunction_0,\timefunction],[- \rightu,u]}{0}}   \sumofbrevexvisquared^2 \left(\fullymodquant{\tander^{\Ntop}}\right)^2 |\timefunction|^{\blowupratetoporderwave - 1}, \label{E:WEIGHTEDPRECISEFULLYMODQUANTBULKROUGHCLASSICALDEVELOPMENT} \\
		\mathbb{X}^{(\textnormal{Imprecise})}(\timefunction,u) & \approx \max_{\tander^N \in \mathfrak{P}^{(N)}}  \int_{\hypthreearg{\timefunction}{[- \rightu,u]}{0}}  \left(\fullymodquant{\tander^{\Ntop}}\right)^2 |\timefunction|^{\blowupratetoporderwave + 1}, \label{E:WEIGHTEDIMPRECISEFULLYMODQUANTFLUXROUGHCLASSICALDEVELOPMENT} \\
		\mathbb{KX}^{(\textnormal{Imprecise})}(\timefunction,u) & \approx \max_{\tander^N \in \mathfrak{P}^{(N)}}   \int_{\twoargMrough{[\timefunction_0,\timefunction],[- \rightu,u]}{0}}     \left(\fullymodquant{\tanderY^{\Ntop}}\right)^2 |\timefunction|^{\blowupratetoporderwave}, \label{E:WEIGHTEDIMPRECISEFULLYMODQUANTBULKROUGHCLASSICALDEVELOPMENT} \\
		\widetilde{\mathbb{X}}(\timefunction,u) & \approx \max_{\tanderY^{\Ntop-1} \in \mathfrak{Y}^{(\Ntop-1)}}  \int_{\hypthreearg{\timefunction}{[- \rightu,u]}{0}} \sumofbrevexvisquared^2 \left(\partialmodquant{\tanderY^{\Ntop-1}}\right)^2 |\timefunction|^{\blowupratetoporderwave - 1}, \label{E:WEIGHTEDPARTIALMODQUANTFLUXROUGHCLASSICALDEVELOPMENT} \\
		\mathbb{K}\widetilde{\mathbb{X}}(\timefunction,u) & \approx  \max_{\tanderY^{\Ntop-1} \in \mathfrak{Y}^{(\Ntop-1)}}  \int_{\twoargMrough{[\timefunction_0,\timefunction],[- \rightu,u]}{0}}   \sumofbrevexvisquared^2 \left(\partialmodquant{\tanderY^{\Ntop-1}}\right)^2 |\timefunction|^{\blowupratetoporderwave - 2}, \label{E:WEIGHTEDPARTIALMODQUANTBULKROUGHCLASSICALDEVELOPMENT}.
	\end{align} 
\end{subequations}Then the analogs of Props.\,\ref{P:TOPORDERAPRIORIL2ESTIMATESMASTERCONTROL},\,\ref{P:APRIORIL2ESTIMATESWAVEVARIABLES}, and\,\ref{P:APRIORIL2ESTIMATESACOUSTICGEOMETRY} with $\initialsmalldoublenull$ replaced with $\initialsmall$, $\ingoingcharacteristicsurfacetwoarg{\ubar}{[\moreinterestingu_1,u]}$ replaced with $\hypthreearg{\timefunction}{[- \rightu,u]}{0}$, and $\characteristicdiamondtwoarg{[\leftubar,\ubar)}{[\moreinterestingu_1,u]}$ replaced with $\twoargMrough{[\timefunction_0,\timefunction],[- \rightu,u]}{0}$ hold. That is, for any $(\timefunction,u) \in [\timefunction_0,0]\times [-\rightu,\leftu]$, we have:
\begin{subequations}
	\begin{align} 
		\max_{\substack{N \in \{1,\dots,\Ntop\} \\ \tander^N \in \mathfrak{P}^{(N)}}}  \left\{ \mathbb{E}[\tander^N \velocityarray](\timefunction,u) + \mathbb{F}[\tander^N \velocityarray](\timefunction,u) + \mathbb{K}[\tander^N \velocityarray] (\timefunction,u)\right\} & \le C \initialsmall^2, \label{E:BOUNDEDWAVEWEIGHTEDENERGYESTIMATESONROUGHCLASSICALDEV} \\
		\max_{N \in \{1,\dots,\Ntop\}}  \left\{ \mathbb{E}[\tander^N (\vortrenormalized,\GradEnt)](\timefunction,u) + \mathbb{F}[\tander^N (\vortrenormalized,\GradEnt)](\timefunction,u) \right\} & \le C \initialsmall^2, \label{E:BOUNDEDVORTGRADENTWEIGHTEDENERGYESTIMATESONROUGHCLASSICALDEV} \\
		\max_{N \in \{1,\dots,\Ntop\} } \left\{ \mathbb{E}[\tander^N (\VortVort,\DivGradEnt)](\timefunction,u) + \mathbb{F}[\tander^N (\VortVort,\DivGradEnt)](\timefunction,u) \right\} & \le C \initialsmall^2, \label{E:BOUNDEDVORTVORTDIVGRADENTWEIGHTEDENERGYESTIMATESONROUGHCLASSICALDEV} \\			
		\mathbb{X}^{(\textnormal{Precise})}(\timefunction,u) + \mathbb{KX}^{(\textnormal{Precise})}(\timefunction,u) & \le C \initialsmall^2 \label{E:BOUNDEDPRECISEFULLYMODCHIWEIGHTEDENERGYESTIMATESONROUGHCLASSICALDEV}, \\
		\mathbb{X}^{(\textnormal{Imprecise})}(\timefunction,u) + \mathbb{KX}^{(\textnormal{Imprecise})}(\timefunction,u) & \le C \initialsmall^2, \label{E:BOUNDEDIMPRECISEFULLYMODCHIWEIGHTEDENERGYESTIMATESONROUGHCLASSICALDEV} \\
		\widetilde{\mathbb{X}}(\timefunction,u) + \mathbb{K}\widetilde{\mathbb{X}}(\timefunction,u) & \le C \initialsmall^2. \label{E:BOUNDEDPARTIALMODCHIWEIGHTEDENERGYESTIMATESONROUGHCLASSICALDEV}
	\end{align}
\end{subequations}
Next, define the following:
	\begin{align}
		\begin{split} \label{E:INTEGRANDFORCHEBYCHEV}
			\mathbf{M} & \eqdef \sum_{ \substack{ N \in \{1,\dots,\Ntop\} \\ \tander^N \in \mathfrak{P}^{(N)} \\ i  \in \{1,2,3\} }} \left\{ (\muX \tander^N v^i)^2 + \upmu (\Lunit \tander^N v^i)^2 + (1 + \upmu) |\angrmd \tander^N v^i|_{\gtorus}^2 + |\tander^N v^i|^2 \right\} |\timefunction|^{\blowuprateofwaveWRTnumberofcommutations(N)} \\
			& \ \ + \sum_{\substack{ N \in \{1,\dots,\Ntop\} \\ \tander^N \in \mathfrak{P}^{(N)} }} \upmu \left| \tander^N (\vortrenormalized,\GradEnt)\right|^2  |\timefunction|^{\blowuprateoftransportWRTnumberofcommutations(N)} +  \upmu \left| \tander^N (\VortVort,\DivGradEnt)\right|^2  |\timefunction|^{\blowuprateofacousticgeoWRTnumberofcommutations(N)}  \\	
			& \ \ +  \sum_{\tander^{\Ntop} \in \mathfrak{P}^{(\Ntop)}}  \left| \pmb{\p} \tander^{\Ntop} (\vortrenormalized,\GradEnt)\right|^2 |\timefunction|^{\blowupratetoporderwave} \\
			& \ \ + \sum_{\tander^{\Ntop} \in \mathfrak{P}^{(\Ntop)}}   \left(\fullymodquant{\tanderY^{\Ntop}}\right)^2 |\timefunction|^{\blowupratetoporderwave + 1} + \sum_{\tanderY^{\Ntop} \in \mathfrak{Y}^{(\Ntop)}} \sumofbrevexvisquared^2 \left(\fullymodquant{\tanderY^{\Ntop}}\right)^2 |\timefunction|^{\blowupratetoporderwave}  \\
			& \ \  +  \sum_{\tanderY^{\Ntop-1} \in \mathfrak{Y}^{(\Ntop-1)}} \sumofbrevexvisquared^2 \left(\partialmodquant{\tanderY^{\Ntop-1}}\right)^2 |\timefunction|^{\blowupratetoporderwave - 1}.
		\end{split}
	\end{align}
Since $\int_{\twoargMrough{[\timefunction_0,0],[- \rightu,u]}{0}} f = \int_{\timefunction' = \timefunction_0}^0 \int_{\hypthreearg{\timefunction'}{[- \rightu,u]}{0}} f \, \rmd \timefunction'$ by \cite{abbrescia2022emergence}*{(9.12d)}, integrating \eqref{E:BOUNDEDWAVEWEIGHTEDENERGYESTIMATESONROUGHCLASSICALDEV}--\eqref{E:BOUNDEDPARTIALMODCHIWEIGHTEDENERGYESTIMATESONROUGHCLASSICALDEV}, as well as the $L^2$-estimates for $|\pmb{\p}\tander^{\Ntop} (\vortrenormalized,\GradEnt)|^2$ (i.e. inserting the bounds from \cite{abbrescia2022emergence}*{Prop.\,24.1-24.4} into \cite{abbrescia2022emergence}*{Prop.\,27.5}) provides the following $L^2$-spacetime bound:
\begin{align} \label{E:SPACETIMEINTEGRALFORCHEBYCHEV}
	\int_{\twoargMrough{[\timefunction_0,0],[- \rightu,u]}{0}}\mathbf{M} \lesssim \initialsmall^2.
\end{align}
\medskip

\noindent \underline{\textbf{$L^2$ data estimates on $\ingoingcharacteristicsurfacetwoarg{\leftubar}{[\moreinterestingu_1,\moreinterestingu_2]}$ by a Chebyshev-type argument:}}  

Recall the embedding $\embeddingdatahypfortimefunctionarg{0}\colon \domainforembeddingdatahypfortimefunctiontwoarg{0}{[0,\mupositive]}  \to \datahypfortimefunctiontwoarg{0}{[\timefunction_0,0]}$ from Prop.\,\ref{P:PROPERTIESOFMUXMUZEROLEVELSET}, the flowout domain $\underline{\mathscr{D}}$ on which $\ubar$ was constructed in Theorem\,\ref{T:CONSTRUCTIONOFTHEINGOINGEIKONALFUNCTIONCONTRACTIONMAPPING}, and the flow map $\flowmapofnewularg{\Delta s}$ of $\newuL$ from \eqref{E:FLOWMAPOFNEWUL}. From the definitions \eqref{E:FLOWOUTOFDATASURFACE}--\eqref{E:DEFINITIONOFFLOWOUTDOMAINFORUBAR}, it follows that: 
	\begin{align} \label{E:FLOWOUTDOMAINFOLIATED}
		\underline{\mathscr{D}}  = \bigcup_{\timefunction\in [\timefunction_0,0]} \left\{ \flowmapofnewularg{\Delta s} \circ \embeddingdatahypfortimefunctionarg{0}(\Cartesiantisafunctiononmumxtoriarg{-\timefunction}{0}(x^2,x^3),x^2,x^3) \ \bigg| \ (x^2,x^3) \in \T^2, \text{ and } [\Delta \underline{a}(\Cartesiantisafunctiononmumxtoriarg{-\timefunction}{0}(x^2,x^3),x^2,x^3),\Delta \underline{b}(\Cartesiantisafunctiononmumxtoriarg{-\timefunction}{0}(x^2,x^3),x^2,x^3)] \right\}.
	\end{align}
	That is, $\underline{\mathscr{D}}$ is foliated by the sets on RHS\,\eqref{E:FLOWOUTDOMAINFOLIATED}, which are subsets of $\{ \ubar  \equiv - \mulevelsetvalue\}$ for $\mulevelsetvalue \in [0,\mupositive]$. By abusing notation, we denote those sets by $\ingoingcharacteristicsurfacetwoarg{\timefunction}{[\Delta \underline{a},\Delta \underline{b}]}$. By Fubini's theorem and \eqref{E:SPACETIMEINTEGRALFORCHEBYCHEV}, we bound the integral of $\mathbf{M}$ over $\underline{\mathscr{D}}$ as: 
	\begin{align} \label{E:FUBINITHEOREMAPPLIEDTOFLOWOUTDOMAIN}
		 \int_{\timefunction' = \timefunction_0}^0 \int_{\ingoingcharacteristicsurfacetwoarg{\timefunction}{[\Delta \underline{a},\Delta \underline{b}]}}  \mathbf{M} = \int_{\mathscr{D}} \mathbf{M} \le \int_{\twoargMrough{[\timefunction_0,0],[- \rightu,u]}{0}}\mathbf{M}  \lesssim \initialsmall^2.
	\end{align}
From \eqref{E:FUBINITHEOREMAPPLIEDTOFLOWOUTDOMAIN}, and Chebychev's inequality, it follows that there exists a $\leftubar \in [\frac{3}{4}\timefunction_0,\frac{1}{4}\timefunction_0]$ such that $ \int_{\ingoingcharacteristicsurfacetwoarg{\leftubar}{[\Delta \underline{a},\Delta \underline{b}]}} \mathbf{M}  \lesssim \initialsmall^2$. 

Next, using arguments similar to those used in the proof of \eqref{E:REQUIREMENTOFFLOWINTERVALFORNEWULTOMAKELEFTFLOWOUTENDINROUGHHYPZERO}, it follows that there exists a non-negative function $\Delta \underline{c}$ such that:
\begin{align}
	\timefunction\left(  \flowmapofnewularg{\Delta \underline{c}(\Cartesiantisafunctiononmumxtoriarg{-\leftubar}{0}(x^2,x^3),x^2,x^3)}\circ  \embeddingdatahypfortimefunctionarg{0}(\Cartesiantisafunctiononmumxtoriarg{-\leftubar}{0}(x^2,x^3),x^2,x^3)\right) = \frac{1}{100} \timefunction_0. \label{E:ENDLEFTFLOWOUTBEFOREFINALROUGHHYP}
\end{align}
Then by setting:
\begin{subequations}
	\begin{align}
		\moreinterestingu_2 & \eqdef \inf_{(x^2,x^3) \in \T^2} u \circ \flowmapofnewularg{\Delta \underline{c}(\Cartesiantisafunctiononmumxtoriarg{-\leftubar}{0}(x^2,x^3),x^2,x^3)
} \circ \embeddingdatahypfortimefunctionarg{0}(\Cartesiantisafunctiononmumxtoriarg{-\leftubar}{0}(x^2,x^3),x^2,x^3), \label{E:DEFINITIONOFLEFTMOREINTERESTINGUVALUE} \\
		\moreinterestingu_1 & \eqdef \sup_{(x^2,x^3) \in \T^2} u \circ \flowmapofnewularg{\Delta a(\Cartesiantisafunctiononmumxtoriarg{-\leftubar}{0}(x^2,x^3),x^2,x^3)
} \circ \embeddingdatahypfortimefunctionarg{0}(\Cartesiantisafunctiononmumxtoriarg{-\leftubar}{0}(x^2,x^3),x^2,x^3),\label{E:DEFINITIONOFRIGHTMOREINTERESTINGUVALUE}
	\end{align}
\end{subequations}
it follows that the desired bounds \eqref{E:MOSTNEGATIVEMOREINTERESTINGVALUEOFU}--\eqref{E:CONDITIONTHATGUARANTEESNONTRIVIALPARTOFCAUCHYHORIZON} hold using  \eqref{E:BOUNDSFORITERATESOFNEWULFLOWOUTDOMAININTERVALDEFININGFUNCTIONS} and choosing $\initialsmall$ to be sufficiently smaller than $|\timefunction_0|$. In total, we have proven that: 
\begin{align}
	\outgoingcharacteristicsurfacetwoarg{\moreinterestingu_1}{[\leftubar,0]}, \, \ingoingcharacteristicsurfacetwoarg{\leftubar}{[\moreinterestingu_1,\moreinterestingu_2]} \subset \underline{\mathscr{D}} \cap \twoargMrough{[\timefunction_0,\frac{1}{100}\timefunction_0],[- \rightu,u]}{0}, \label{E:INITIALINGOINGCHARACTERISTICSURFACEDOESNOTTOUCHFINALROUGHHYP} \\ 
	\int_{\ingoingcharacteristicsurfacetwoarg{\leftubar}{[\moreinterestingu_1,\moreinterestingu_2]}} \mathbf{M} \lesssim  \int_{\ingoingcharacteristicsurfacetwoarg{\leftubar}{[\Delta \underline{a},\Delta \underline{b}]}} \mathbf{M}  \lesssim \initialsmall^2. \label{E:COERCIVECONTROLOFMOSTTHINGSONINITIALINGOINGCHARACTERISTIC} 
\end{align}
We are now ready to prove that the $L^2$ data estimates $\ingoingcharacteristicsurfacetwoarg{\leftubar}{[\moreinterestingu_1,\moreinterestingu_2]}$  from Sect.\,\ref{SSS:QUANTITATIVEASSUMPTIONSONDATAAWAYFROMSYMMETRY} hold. From \eqref{E:MOSTNEGATIVEVALUEOFUBAR} and \eqref{E:INITIALINGOINGCHARACTERISTICSURFACEDOESNOTTOUCHFINALROUGHHYP}, for any $p \in \ingoingcharacteristicsurfacetwoarg{\leftubar}{[\moreinterestingu_1,\moreinterestingu_2]}$, the following estimates hold:
\begin{align}
	\weight (p) = |\leftubar| \approx |\timefunction_0| \approx | \timefunction|(p). \label{E:UBARAPPROXTIMEFUNCTIONONINITIALINGOINGCHARACTERISTIC}
\end{align}
From \eqref{E:UBARAPPROXTIMEFUNCTIONONINITIALINGOINGCHARACTERISTIC}, the last three sums in \eqref{E:INTEGRANDFORCHEBYCHEV}, and \eqref{E:COERCIVECONTROLOFMOSTTHINGSONINITIALINGOINGCHARACTERISTIC}, the desired bounds \eqref{E:L2DATASSUMPTIONFORFULLYMODIFIEDCHITOPORDER}--\eqref{E:L2DATASSUMPTIONFORPARTIALLYMODIFIEDCHITOPORDER} hold with $\initialsmalldoublenull = C \initialsmall$. Similarly, from \eqref{E:GEOP2TOCOMMUTATORS}--\eqref{E:GEOP3TOCOMMUTATORS}, \eqref{E:CONVENIENTIDENTITYFORULUNIT}, Prop\,\ref{P:SCHEMATICSTRUCTUREOFVARIOUSTENSORSINTERMSOFCONTROLVARS}, \eqref{E:MAGNITUDEOFINNERPRODUCTOFNEWLANDNEWULAPPROXIMATELYMU}--\eqref{E:RECIPROCALULUNITAPPLIEDTOTIMEFUNCTIONAPPROXIMATELYMU}, \eqref{E:ANGDFPOINTWISEBOUNDEDBYCOMMUTATORVECTORFIELDS}, and \eqref{E:POINTWISEDOULBENULLTORIDIFFERENTIALINORMINTERMSOFSMOOTHTORIDIFFERENTIALNORMANDL}, and the first three sums in \eqref{E:INTEGRANDFORCHEBYCHEV}, we have the following bounds for any $N \in \{0,\dots,\Ntop\}$ and $\tander^N \in \mathfrak{P}^{(N)}$:
	\begin{subequations}
		\begin{align}
			\sum_{i=1}^3 \int_{\ingoingcharacteristicsurfacetwoarg{\leftubar}{[\moreinterestingu_1,\moreinterestingu_2]}} \left\{ (\newuL \tander^N v^i)^2 + \upmu |\nullangrmD \tander^N v^i |_{\gnulltori}^2 + |\tander^N v^i |^2 \right\} |\weight|^{\blowuprateofwaveWRTnumberofcommutations(N)} & \lesssim \initialsmall^2  \label{E:INITIALINGOINGCHARACTERISTICSURFACEL2WAVEBOUNDS} \\
			 \int_{\ingoingcharacteristicsurfacetwoarg{\leftubar}{[\moreinterestingu_1,\moreinterestingu_2]}} \upmu \left| \tander^{N} (\vortrenormalized,\GradEnt)\right|^2  |\weight|^{\blowuprateoftransportWRTnumberofcommutations(N)} & \lesssim \initialsmall^2,\label{E:INITIALINGOINGCHARACTERISTICSURFACEL2VORTGRADENTBOUNDS} \\
			 \int_{\ingoingcharacteristicsurfacetwoarg{\leftubar}{[\moreinterestingu_1,\moreinterestingu_2]}} \upmu \left| \tander^{N} (\VortVort,\DivGradEnt)\right|^2  |\weight|^{\blowuprateofacousticgeoWRTnumberofcommutations(N)} & \lesssim \initialsmall^2.\label{E:INITIALINGOINGCHARACTERISTICSURFACEL2VORTVORTDIVGRADENTBOUNDS}		 
		\end{align}
	\end{subequations}
The first two terms in \eqref{E:INITIALINGOINGCHARACTERISTICSURFACEL2WAVEBOUNDS} conclude the proof of \eqref{E:TRANSVERSALDERIVATIVEOFTANGENTIALL2NORMSOFWAVEVARIABLESSMALLALONGINITIALINGOINGNULLSURFACE}--\eqref{E:NULLANGULARDERIVATIVEOFTANGENTIALL2NORMSOFWAVEVARIABLESSMALLALONGINITIALROUGHHYPERSURFACE}, while \eqref{E:INITIALINGOINGCHARACTERISTICSURFACEL2VORTGRADENTBOUNDS}--\eqref{E:INITIALINGOINGCHARACTERISTICSURFACEL2VORTVORTDIVGRADENTBOUNDS} respectively prove \eqref{E:L2NORMSOFVORTANDGRADENTSMALLALONGINITIALINGOINGNULLSURFACE}--\eqref{E:L2NORMSOFVORTVORTANDDIVGRADENTSMALLALONGINITIALINGOINGNULLSURFACE}  with $\initialsmalldoublenull = C \initialsmall$.

\medskip

\noindent \underline{\textbf{$L^2$ below-top-order data estimates on $\doublenulltoritwoarg{\leftubar}{u}$ by a Chebyshev-type argument:}}  

We now prove \eqref{E:SMALLDATAOFPSIONINITIALDOUBLENULLTORI}--\eqref{E:SMALLDATAOFVORTVORTANDDIVGRADENTROPYONDOUBLENULLTORIONPU0}. Let $N \in \{1,\dots,\Ntop\}$. Since $ \int_{\ingoingcharacteristicsurfacetwoarg{\leftubar}{[\moreinterestingu_1,\moreinterestingu_2]}} \smoothfunction = \int_{u' = \moreinterestingu_1}^{\moreinterestingu_2} \int_{\doublenulltoritwoarg{\leftubar}{u'}} \smoothfunction \, \rmd u'$, \eqref{E:INITIALINGOINGCHARACTERISTICSURFACEL2WAVEBOUNDS}--\eqref{E:INITIALINGOINGCHARACTERISTICSURFACEL2VORTVORTDIVGRADENTBOUNDS} and Chebyshev's inequality imply the existence of a $u_* \in [\moreinterestingu_1,\moreinterestingu_2]$ such that: 
\begin{align} \label{E:SPECIALNULLTORIL2WAVEBOUNDS} 
	\int_{\doublenulltoritwoarg{\leftubar}{u_*}} \left\{|\tander^N \velocityarray|^2 |\weight|^{\blowuprateofwaveWRTnumberofcommutations(N)} + \upmu |\tander^{N-1}(\vortrenormalized,\GradEnt)|^2 |\weight|^{\blowuprateoftransportWRTnumberofcommutations(N)} + \upmu |\tander^{N-1}(\VortVort,\DivGradEnt)|^2 |\weight|^{\blowuprateofacousticgeoWRTnumberofcommutations(N)}\right\} \lesssim \initialsmall^2. 
\end{align}
Since $\newuL \weight = 0$, estimate \eqref{E:POINTWISEBOUNDFORDOUBLENULLTOROIDALTRACEOFDEFORMATIONTENSOROFNEWULVECTORFIELD}, identity  \eqref{E:NEWULINTEGRALONDOUBLENULLTORI}, Gr\"onwall's inequality, and \eqref{E:INITIALINGOINGCHARACTERISTICSURFACEL2WAVEBOUNDS} conclude the proof of \eqref{E:SMALLDATAOFPSIONINITIALDOUBLENULLTORI} with $\initialsmalldoublenull = C \initialsmall$. 

To prove \eqref{E:SMALLDATAOFVORTICITYANDGRADENTROPYONDOUBLENULLTORIONPU0}--\eqref{E:SMALLDATAOFVORTVORTANDDIVGRADENTROPYONDOUBLENULLTORIONPU0}, by Gr\"onwall's inequality, the same arguments as in the previous paragraph imply that for any $u \in [\moreinterestingu_1,\moreinterestingu_2]$, we have: 
\begin{align} 
		\begin{split} \label{E:GRONWALLTYPEESTIMATEFORPROVINGBELOWTOPORDERDATASSUMPTIONSONTORIFORVORTICITY}
			& \sum_{N \le \Ntop} \int_{\doublenulltoritwoarg{\leftubar}{u}} \left\{ \upmu |\tander^{N-1}(\vortrenormalized,\GradEnt)|^2 |\weight|^{\blowuprateoftransportWRTnumberofcommutations(N)} + \upmu |\tander^{N-1}(\VortVort,\DivGradEnt)|^2 |\weight|^{\blowuprateofacousticgeoWRTnumberofcommutations(N)}\right\} \\
			& \ \  \lesssim \sum_{N \le \Ntop} \initialsmall^2 + \int_{\ingoingcharacteristicsurfacetwoarg{\leftubar}{[u_*,u]}} \left\{ \upmu \left| \newuL \tander^{N-1}(\vortrenormalized,\GradEnt)\right|^2 |\timefunction|^{\blowuprateoftransportWRTnumberofcommutations(N)} + \upmu \left| \newuL \tander^{N-1}(\VortVort,\DivGradEnt)\right|^2 |\timefunction|^{\blowuprateofacousticgeoWRTnumberofcommutations(N)} \right\}
		\end{split} 
	\end{align}
Using the second identity in \eqref{E:CONVENIENTIDENTITYFORULUNIT}, the pointwise bounds  \eqref{E:COMMUTEDTRANSPORTPOINTWISEESTIMATESFORMODIFIEDCURLOFVORT}--\eqref{E:COMMUTEDTRANSVERSALDERVIATVIESOFENTROPYGRADIENTINTERMSOFTANGENTIAL}, the identity $\muX = \upmu \Transport - \upmu \Lunit$, and over-estimating the integral to the larger region $\ingoingcharacteristicsurfacetwoarg{\leftubar}{[\moreinterestingu_1,\moreinterestingu_2]}$, the desired bounds  \eqref{E:SMALLDATAOFVORTICITYANDGRADENTROPYONDOUBLENULLTORIONPU0}--\eqref{E:SMALLDATAOFVORTVORTANDDIVGRADENTROPYONDOUBLENULLTORIONPU0} follow from \eqref{E:COERCIVECONTROLOFMOSTTHINGSONINITIALINGOINGCHARACTERISTIC}.

\medskip

\noindent \underline{\textbf{$L^2$ data estimates on $\outgoingcharacteristicsurfacetwoarg{\moreinterestingu_1}{[\leftubar,0]}$:}}

From \eqref{E:WEIGHTEDNULLFLUXONROUGHFOLIATIONTRUNCATEDOUTGOINGCHARACTERISTICS}, \eqref{E:BOUNDEDWAVEWEIGHTEDENERGYESTIMATESONROUGHCLASSICALDEV}, and \eqref{E:INITIALINGOINGCHARACTERISTICSURFACEDOESNOTTOUCHFINALROUGHHYP}, we have that: 
\begin{align} 
	\begin{split} \label{E:DATAESTIMATEL2BOUNDONOUTGOINGCHARACTERISTIC}
		\sum_{i=1}^3 \int_{\outgoingcharacteristicsurfacetwoarg{\moreinterestingu_1}{[\leftubar,0]}} \left\{ (\Lunit \tander^N v^i )^2 + \upmu \left| \angrmd  \tander^N v^i \right|_{\gtorus}^2\right\} \weight^{\blowuprateofwaveWRTnumberofcommutations(N)} & \approx \int_{\outgoingcharacteristicsurfacetwoarg{\moreinterestingu_1}{[\leftubar,0]}}  \left\{ (\Lunit \tander^N v^i )^2 + \upmu \left| \angrmd \tander^N v^i \right|_{\gtorus}^2\right\}  |\timefunction| ^{\blowuprateofwaveWRTnumberofcommutations(N)}  \\
		& \le \int_{\nullhypthreearg{0}{\moreinterestingu_1}{[\timefunction_0,0]} }  \left\{ (\Lunit \tander^N v^i )^2 + \upmu \left| \angrmd 
		\tander^N v^i \right|_{\gtorus}^2\right\}  |\timefunction|^{\blowuprateofwaveWRTnumberofcommutations(N)} \\
		& \lesssim \initialsmall^2.
	\end{split}
\end{align}
The desired the data estimates \eqref{E:LDERIVATIVEOFTANGENTIALL2NORMSOFWAVEVARIABLESSMALLALONGINITIALOUTGOINGCHARACTERISTICHYPERSURFACE}--\eqref{E:SMOOTHANGULARDERIVATIVEOFTANGENTIALL2NORMSOFWAVEVARIABLESSMALLALONGINITIALOUTGOINGCHARACTERISTICHYPERSURFACE} then follow from \eqref{E:DATAESTIMATEL2BOUNDONOUTGOINGCHARACTERISTIC} with $\initialsmalldoublenull = C \initialsmall$. 

Identical arguments hold for \eqref{E:SMALLDATAOFVORTICITYANDGRADENTROPYONDOUBLENULLTORIONPU0}--\eqref{E:L2NORMSOFVORTVORTANDDIVGRADENTSMALLALONGINITIALOUTGOINGNULLSURFACE} using instead \eqref{E:BOUNDEDVORTGRADENTWEIGHTEDENERGYESTIMATESONROUGHCLASSICALDEV}--\eqref{E:BOUNDEDVORTVORTDIVGRADENTWEIGHTEDENERGYESTIMATESONROUGHCLASSICALDEV}.

\noindent \underline{\textbf{Top order $L^2$ data estimates on double-null tori for the elliptic hyperbolic identities:}}

We now prove \eqref{E:SMALLDATAOFTOPORDERVORTICITYANDGRADENTROPYONDOUBLENULLTORIONUBAR0FORELLIPTICHYPERBOLICIDENTITIES}--\eqref{E:SMALLDATAOFTOPORDERVORTICITYANDGRADENTROPYONDOUBLENULLTORIONU1FORELLIPTICHYPERBOLICIDENTITIES}. Throughout this proof, we note that $\upmu \ge C > 0$ on both the outgoing $\outgoingcharacteristicsurfacetwoarg{\moreinterestingu_1}{[\leftubar,0]}$ and ingoing $\ingoingcharacteristicsurfacetwoarg{\leftubar}{[\moreinterestingu_1,\moreinterestingu_2]}$ characteristics. 

Since $|\timefunction|^{\blowupratetoporderacoustic} \lesssim |\timefunction|^{\blowuprateoftransportWRTnumberofcommutations(\Ntop)}$ on $\ingoingcharacteristicsurfacetwoarg{\leftubar}{[\moreinterestingu_1,\moreinterestingu_2]}$, the bound \eqref{E:COERCIVECONTROLOFMOSTTHINGSONINITIALINGOINGCHARACTERISTIC} easily implies the following estimates:
	\begin{align}
		\int_{\ingoingcharacteristicsurfacetwoarg{\leftubar}{[\moreinterestingu_1,\moreinterestingu_2]}} \left| \pmb{\p}\tander^{\Ntop} (\vortrenormalized,\GradEnt)\right|^2 |\timefunction|^{\blowupratetoporderacoustic},  \int_{\ingoingcharacteristicsurfacetwoarg{\leftubar}{[\moreinterestingu_1,\moreinterestingu_2]}} \left| \tander^{\Ntop} (\vortrenormalized,\GradEnt)\right|^2 |\timefunction|^{\blowupratetoporderacoustic}  \lesssim  \initialsmall^2. \label{E:BOUNDONINITIALINGOINGCHARACTERISTICFORTOPORDERVORTANDGRADENT} 
	\end{align}
From the second bound in \eqref{E:BOUNDONINITIALINGOINGCHARACTERISTICFORTOPORDERVORTANDGRADENT} and Chebychev's inequality, there exists a $U_*$ such that $\int_{\doublenulltoritwoarg{\leftubar}{U_*}} \left| \tander^{\Ntop} (\vortrenormalized,\GradEnt)\right|^2 |\timefunction|^{\blowupratetoporderacoustic}  \lesssim  \initialsmall^2$. From this, the first bound in \eqref{E:BOUNDONINITIALINGOINGCHARACTERISTICFORTOPORDERVORTANDGRADENT}, and a similar Gr\"onwall inequality argument as in \eqref{E:GRONWALLTYPEESTIMATEFORPROVINGBELOWTOPORDERDATASSUMPTIONSONTORIFORVORTICITY}, we have that the following estimate holds:
	\begin{align}
		\int_{\doublenulltoritwoarg{\leftubar}{u}} \frac{1}{\upmu} \left| \tander^{\Ntop} (\vortrenormalized,\GradEnt)\right|^2 |\timefunction|^{\blowupratetoporderacoustic}  \lesssim  \initialsmall^2, & & \forall u \in [\moreinterestingu_1,\moreinterestingu_2]. \label{E:ANYNULLTORIONINITIALINGOINGCHARACTERISTICBOUNDWITHTIMEFUNCTION}
	\end{align}
Next, since  $|\timefunction|^{\blowupratetoporderacoustic} \lesssim |\timefunction|^{\blowuprateoftransportWRTnumberofcommutations(\Ntop)}$ on $\outgoingcharacteristicsurfacetwoarg{\moreinterestingu_1}{[\leftubar,0]}$, \eqref{E:BOUNDEDVORTGRADENTWEIGHTEDENERGYESTIMATESONROUGHCLASSICALDEV} implies $ \int_{\outgoingcharacteristicsurfacetwoarg{\moreinterestingu_1}{[\leftubar,0]}} \frac{1}{\upmu} \left| \tander^{\Ntop} (\vortrenormalized,\GradEnt)\right|^2 |\timefunction|^{\blowupratetoporderacoustic}  \lesssim  \initialsmall^2$. By Chebychev, there exists a $\underline{\mathfrak{U}}_* \in [\leftubar,0]$ such that: 
	\begin{align}
		 \int_{\doublenulltoritwoarg{\underline{\mathfrak{U}}_*}{\moreinterestingu_1}} \frac{1}{\upmu} \left| \tander^{\Ntop} (\vortrenormalized,\GradEnt)\right|^2 |\timefunction|^{\blowupratetoporderacoustic}  \lesssim  \initialsmall^2. \label{E:ASPECIALNULLTORIONINITIALOUTGOINGCHARACTERISTICBOUNDWITHTIMEFUNCTION}
	\end{align}By setting $\mathfrak{U}_* := \inf_{(x^2,x^3) \in \T^2} \Eikonalisafunctiononubarnewultori_0(x^2,x^3) - \uplambda$ where $\Eikonalisafunctiononubarnewultori_0(x^2,x^3)$ denotes the $u$-coordinate of the crease as a graph over $\T^2$ (see \eqref{E:BAUBARNEWULMUTORISTRUCTURE} and \eqref{E:C01BOUNDFORUDOUBLENULLCOORDINATEGRAPHTORI}) and $\uplambda \ll |\timefunction_0|$, it follows that $\characteristicdiamondtwoarg{[\leftubar,0]}{[\moreinterestingu_1,\mathfrak{U}_*]}$ is a characteristic diamond completely contained in the region of classical existence $\twoargMrough{[\timefunction_0,0],[- \rightu,\leftu]}{0}$ and uniformly bounded away from the crease. Since: 
	\begin{align} \label{E:CHEBYCHEVAGAINTOFINDAGOODOUTGOINGCHARACTERISTICFORDATA}
		\int_{u' = \moreinterestingu_1}^{\mathfrak{U}_*} \int_{\outgoingcharacteristicsurfacetwoarg{u'}{[\leftubar,0]}} \left| \pmb{\p} \tander^{\Ntop}(\vortrenormalized,\GradEnt)\right|^2 |\timefunction|^{\blowupratetoporderacoustic} = \int_{\characteristicdiamondtwoarg{[\leftubar,0]}{[\moreinterestingu_1,\mathfrak{U}_*]}} \left| \pmb{\p} \tander^{\Ntop}(\vortrenormalized,\GradEnt)\right|^2 |\timefunction|^{\blowupratetoporderacoustic} \lesssim \initialsmall^2,
	\end{align}
by Chebychev again there exists a $\moreinterestingu_* \in [\moreinterestingu_1,\mathfrak{U}_*]$ such that $\int_{\outgoingcharacteristicsurfacetwoarg{\moreinterestingu_*}{[\leftubar,0]}} \left| \pmb{\p} \tander^{\Ntop}(\vortrenormalized,\GradEnt)\right|^2 |\timefunction|^{\blowupratetoporderacoustic} \lesssim \initialsmall^2$. Similarly, by the same arguments used to prove \eqref{E:ASPECIALNULLTORIONINITIALOUTGOINGCHARACTERISTICBOUNDWITHTIMEFUNCTION}, there exists a $\underline{\mathfrak{U}}_*' \in [\leftubar,0]$ such that $ \int_{\doublenulltoritwoarg{\underline{\mathfrak{U}}_*'}{\moreinterestingu_*}} \frac{1}{\upmu} \left| \tander^{\Ntop} (\vortrenormalized,\GradEnt)\right|^2 |\timefunction|^{\blowupratetoporderacoustic}  \lesssim  \initialsmall^2$. From this estimate and using a similar Gr\"onwall inequality argument as in \eqref{E:GRONWALLTYPEESTIMATEFORPROVINGBELOWTOPORDERDATASSUMPTIONSONTORIFORVORTICITY} (with bound \eqref{E:DOUBLENULLTORUSINNULLHYPERSURFACEL2FUNDAMENTALTHEOREMOFCALCULUSESTIMATE} in place of \eqref{E:NEWULINTEGRALONDOUBLENULLTORI}), we have the following:
\begin{align} \label{E:ANENTIREFAMILYOFGOODDOUBLENULLTORIDATAESTIMATES}
		 \int_{\doublenulltoritwoarg{\ubar}{\moreinterestingu_*}} \frac{1}{\upmu} \left| \tander^{\Ntop} (\vortrenormalized,\GradEnt)\right|^2 |\timefunction|^{\blowupratetoporderacoustic} \lesssim \initialsmall^2, & & \forall \ubar \in [\leftubar,0].
	\end{align}
For $\ubar \in [\leftubar, \underline{\mathfrak{U}}_*)$, set  $\ubar_1 = \ubar$ and $\ubar_2 = \underline{\mathfrak{U}}_*$. Whenever $\ubar \in ( \underline{\mathfrak{U}}_*,0]$, set $\ubar_1 = \underline{\mathfrak{U}}_*$ and $\ubar_2 = \ubar$. Set $u_1 := \moreinterestingu_1$ and $u_2 = \moreinterestingu_*$. Applying then  \eqref{E:ANENTIREFAMILYOFGOODDOUBLENULLTORIDATAESTIMATES} and \eqref{E:ASPECIALNULLTORIONINITIALOUTGOINGCHARACTERISTICBOUNDWITHTIMEFUNCTION} as ``data estimates'', we can apply the elliptic-hyperbolic integral identity of Prop.\,\ref{P:INTEGRALIDENTITYFORELLIPTICHYPERBOLICCURRENT}\footnote{Note that Prop.\,\ref{P:INTEGRALIDENTITYFORELLIPTICHYPERBOLICCURRENT} also holds true for any arbitrary weight, up to error terms that fall on the weight. These error terms are trivially bounded using the bounds on the rough time function $\timefunction$ and the fact that the identity is applied on $\characteristicdiamondtwoarg{[\ubar_1,\ubar_2]}{[u_1,u_2]} \subset \twoargMrough{[\timefunction_0,0],[- \rightu,\leftu]}{0}$.} and algebraically solve for the double-null-tori integral $ \int_{\doublenulltoritwoarg{\ubar}{\moreinterestingu_1}} \frac{1}{\upmu} \left| \tander^{\Ntop} (\vortrenormalized,\GradEnt)\right|^2 |\timefunction|^{\blowupratetoporderacoustic}$. This  integral in question is then given in terms of error terms that are all $\lesssim \initialsmall^2$ from nearly identical proofs of Prop.\,\ref{P:TOPORDERAPRIORIL2ESTIMATESMASTERCONTROL}--\ref{P:TOPORDERMODIFIEDAPRIORIL2ESTIMATESANDDOUBLENULLTORIESTIMATES}. In particular, we have proven that: 
\begin{align} \label{E:DESIREDNULLTORIFAMILYESTIMATES}
		 \int_{\doublenulltoritwoarg{\ubar}{\moreinterestingu_1}} \frac{1}{\upmu} \left| \tander^{\Ntop} (\vortrenormalized,\GradEnt)\right|^2 |\timefunction|^{\blowupratetoporderacoustic} \lesssim \initialsmall^2, & & \forall \ubar \in [\leftubar,0].
	\end{align}
Since $|\weight| \le |\timefunction_0| \approx |\timefunction|$ on $\outgoingcharacteristicsurfacetwoarg{\moreinterestingu_1}{[\leftubar,0]}$, the desired data bounds \eqref{E:SMALLDATAOFTOPORDERVORTICITYANDGRADENTROPYONDOUBLENULLTORIONUBAR0FORELLIPTICHYPERBOLICIDENTITIES}--\eqref{E:SMALLDATAOFTOPORDERVORTICITYANDGRADENTROPYONDOUBLENULLTORIONU1FORELLIPTICHYPERBOLICIDENTITIES} follow from \eqref{E:ANYNULLTORIONINITIALINGOINGCHARACTERISTICBOUNDWITHTIMEFUNCTION} and \eqref{E:DESIREDNULLTORIFAMILYESTIMATES}, as desired.

\end{proof}

\section{A translational dictionary and comparison of the approaches of \cite{dC2007,abbrescia2022emergence,shkoller2024geometry}} 
\label{S:DICTIONARY}
 Our first goal in this appendix is to provide a translational ``dictionary'' between the analytic and geometric methods of \cite{dC2007,abbrescia2022emergence} and \cite{shkoller2024geometry}. We hope that the traditional fluids community benefits from our making the  ``geometric'' language of \cite{dC2007,abbrescia2022emergence} more accessible, while the geometric waves community benefits from our making the more ``fluids'' language of \cite{shkoller2024geometry} more accessible. We clarify that we are not literally translating the full scope of 
\cite{dC2007,abbrescia2022emergence,shkoller2024geometry}. Rather, our aim is for readers to be able to use this appendix to understand the 
key commonalities of the approaches.

Our second goal is to compare the approaches. The works share many equivalent underlying methods. 
In particular, the two ``languages'' are connected by two coordinate systems that we prove to be identical; see Prop.\,\ref{P:ALETOGEOCOORD}.
There are also some pronounced differences between the works. Below we will discuss these similarities and differences.
We focus our attention on the following fundamental issues:
 
 \begin{itemize}
 	\item \textbf{Coordinates; Sect.\,\ref{SSS:ACOUSTICALCOORDINATESANDALECOORDINATES}:} The ``geometric'' and ``fluids'' languages respectively involve 	
	coordinate systems that we prove to be identical; see Prop.\,\ref{P:ALETOGEOCOORD}. The key point is that in all works, these coordinate systems	
	regularize the singularity and form the backbone of the PDE analysis.
	\item \textbf{Fundamental unknowns; Sec.\,\ref{SS:SOLUTIONVARIABLESUSEDINENERGYMETHOD}:} In \cite{dC2007,abbrescia2022emergence,shkoller2024geometry}
	and in this paper, energy estimates are derived for certain nonlinear combinations of the first derivatives of the physical fluid variables. 
	While the nonlinear combinations used in each of the works are different,
	they are closely related, and we will refer to all of them as the ``$L^2$-differentiated solution variables'' to distinguish them from the (un-differentiated) physical fluid variables. 
		A unifying theme in all the works is that the energies are \emph{anisotropic} with respect to the
	way that the $\upmu$ weight appears in their expression. More precisely, in all the energy definitions, 
	relative to the coordinates from the previous point,
	differentiations or contractions in directions
	that are \emph{transversal} to the characteristics require an extra  
	power of $\upmu$ compared to tangential differentiations or contractions, i.e., the transversal directions are ``worse.''
	\item \textbf{Consequences of multiplying the evolution equations first by $\upmu$ before commuting with good derivatives; Sect.\,\ref{SS:MULTIPLYINGFIRSTBYMU}:} The first step in deriving higher order energy estimates for the $L^2$ solution variables 
	is to derive evolution equations that they satisfy. In \cite{dC2007}, Christodoulou observed that there is a specific order of operations one needs to follow to accomplish this: one \emph{first} multiplies the base evolution equations -- without any commutations -- of the $L^2$-solution variables by $\upmu$. One \emph{then} commutes these $\upmu$-weighted equations with smooth (with respect to the good coordinates of the first bullet point) vectorfield operators. This fundamental idea is carried out in \cite{abbrescia2022emergence, shkoller2024geometry} and this paper, which we explain in both the ``geometric'' and ``fluids'' approaches in Sect.\,\ref{SS:MULTIPLYINGFIRSTBYMU}. 
	This strategy for commuting the equations generates many common difficulties in \cite{dC2007,abbrescia2022emergence,shkoller2024geometry}, 
	two of which we highlight in this appendix:
	\begin{itemize}
		\item Some error terms generated by the commutations are ``dangerous'' in that 
		they feature insufficient $\upmu$-weights to be directly controlled by the anisotropic energies described in the second bullet point, 
		even though they feature a small coefficient. Fortunately, there exist ``compensating'' 
		error terms generated by the commutations that feature a good sign (i.e., they are \emph{coercive}) and compatible with the missing $\upmu$-weights.  
		These ``friction integrals'' can be used to absorb (thanks to small coefficients) the dangerous integrals. 
		We discuss these signed spacetime integrals in Sect.\,\ref{SS:COERCIVEBULKANDAMPINGTERMS}.
		\item Some error error terms that arise in the commuted equations appear to lose derivatives;  
		see Sect.\,\ref{SSS:AVOIDINGDERIVATIVELOSSINTHEGEOMETRY}.
		In the $3D$ setting, this loss was overcome in \cite{dC2007,abbrescia2022emergence} 
		and this paper by combining Raychaudhuri's equation with elliptic estimates
		to show that the harmful terms satisfy a better-than-expected estimate. 
		In contrast, in the $2D$ isentropic case,
		the paper \cite{shkoller2024geometry} features a novel divergence identity
		that allows one to handle these terms in a different way, based on integrating by parts.
	\end{itemize}
\end{itemize}
	
Despite the similar themes mentioned above, the detailed implementation of the ideas is different in the various papers. 
Moreover, we would like to again highlight the fundamental difference between the solutions studied in \cite{abbrescia2022emergence,dC2007,shkoller2024geometry} by recalling the table from Sect.\,\ref{SSS:NONLINEARGEOMETRICOPTICS}:

\begin{table}[h!]
\centering
\begin{tabular}{|c|c|c|c|c|c|c|}
\hline
 &   \makecell{$2D$ isentropic \\ Euler up to $\singhyp$} &   \makecell{$2D$ isentropic \\ Euler up to $\Cauchyhor$} &  \makecell{$3D$ quasilinear scalar \\  equations up to $\singhyp$} & \makecell{$3D$ full Euler \\ system up to $\singhyp$ } & \makecell{$3D$ full Euler\\ system up to $\Cauchyhor$} \\ \hline
 Strict convexity & \cite{shkoller2024geometry} & \cite{shkoller2024geometry} & \cite{dC2007} & &\\ \hline
Transversal convexity  &    &  &  & \cite{abbrescia2022emergence}  & This paper \\ \hline
\end{tabular}
\end{table}

Above all else, the distinctions between transversal convexity and strict convexity are the most fundamental. Strict convexity, if it holds at all, usually does so only in a small neighborhood of the point of first gradient-blowup.  Transversal convexity is much more stable because it ``causes'' the gradient-blowup to happen in the first place in a non-degenerate way. Without any kind of convexity, the structure of the boundary of the MGHD can dramatically change, and the shock development problem can end up ill-posed, even in $1D$.

\subsection{Characteristics and co-characteristics of the compressible Euler system} \label{SS:CHARSANDCOCHARSFOR2DEULER} To facilitate the discussion below, we first provide a standard definition of characteristics for the compressible Euler system.
Because \cite{dC2007,abbrescia2022emergence} treats $3D$ fluids while \cite{shkoller2024geometry} treats $2D$ isentropic fluids, for the purposes of this appendix and this appendix only, we sometimes restrict our discussion to the $2D$ isentropic case to facilitate our comparison of the works.

For any $p \in \R \times \Sigma$,\footnote{Here, $\Sigma$ corresponds to the spatial slices of spacetime at a fixed Cartesian time value. In \cite{dC2007}, $\Sigma = \R^3$, in \cite{shkoller2024geometry}, $\Sigma = \T^2$, and in \cite{abbrescia2022emergence}, $\Sigma = \R \times \T^2$.}  given a spacetime co-vector $\bm{\upxi} = (\upxi_0,\upxi_1,\upxi_2) \in T^*_p(\R \times \Sigma)$, the symbol of the $2D$ isentropic Euler equations is: 
	\begin{align} \label{E:SYMBOLOF2DISENTROPICEULER}
		\upsigma(p,\upxi) = \begin{pmatrix}
			\upxi_0 + v^a \upxi_a & \varrho \upxi_1 & \varrho \upxi_2 & 0 \\
			\frac{\Speed^2}{\varrho} \upxi_1 & \upxi_0 + v^a \upxi_a & 0 & 0 \\
			\frac{\Speed^2}{\varrho} \upxi_2 & 0 & \upxi_0 + v^a \upxi_a & 0 \\
			0 & 0 & 0 & \upxi_0 + v^a \upxi_a
			\end{pmatrix}.
	\end{align}
Using \eqref{E:SYMBOLOF2DISENTROPICEULER}, we compute the \emph{characteristic polynomial}:\footnote{In the $3D$ non-isentropic setting, $\det \upsigma(p,\upxi) = (\upxi_0 + v^a \upxi_a)^3 ((\upxi_0 + v^a \upxi_a)^2 - \Speed^2 |\upxi|^2)$, and thus the material derivative characteristics have multiplicity $3$.} 
	\begin{align}
		\det\upsigma(p,\bm{\upxi}) = (\upxi_0 + v^a \upxi_a)((\xi_0 + v^a \upxi_a)^2 - \Speed^2 |\upxi|^2), \label{E:CHARACTERISTICPOLYNOMIAL}
	\end{align}
where $|\upxi|^2 := \upxi_1^2 + \upxi_2^2$ is the Euclidean length of the spatial components of $\bm{\upxi}$. The \emph{co-characteristics of the $2D$-isentropic Euler system}\footnote{The subset of the co-tangent bundle $\{(p,\upxi) \in T^*(\R\times\Sigma) , | \, \det \sigma(p,\upxi) = 0\}$ is often referred to as the \emph{characteristic variety}.} at $p \in \R\times\Sigma$ are defined to be the set of $\bm{\upxi} \in T_p^*(\R\times\Sigma)$ such that $\det \upsigma(p,\bm{\upxi}) = 0$. From \eqref{E:CHARACTERISTICPOLYNOMIAL}, we see that the co-characteristics are precisely the plane:
\begin{align} \label{E:COTANGENTSPACEMATERIALDERIVATIVEPLANE}
	P^*_p
	&:=\{\bm{\upxi} \in T^*_p (\R \times \Sigma) \, | \upxi_0 = - v^i\upxi_i\}, 
\end{align}
representing the co-characteristics of the material derivative vectorfield, 
together with the cone:
\begin{align} \label{E:COTANGENTSPACESOUNDCONES}
	C^*_p & :=\{\bm{\upxi} \in T_p^*(\R\times\Sigma)\, | \, \upxi_0 = - v^i\upxi_i \pm \Speed |\upxi|\}, 
\end{align}
representing the (co-)ingoing and (co-)outgoing sound cones. 

Next, we note the following crucial fact, which follows from a straightforward computation based 
on \eqref{E:ACOUSTICALMETRICANDMATERIALDER}: 
any spacetime one-form of the form 
$\bm{\upxi} := \upxi_\alpha \mathrm{d}x^\alpha$ satisfying:
\begin{align} \label{E:XI0CONDITIONFORNULL}
	\upxi_0 & = - v^a \upxi_a \pm \Speed |\upxi|, 
\end{align}
is necessarily $\gfour$-\emph{null}:
	\begin{align} 
		(\gfour^{-1})^{\alpha\beta}\upxi_\alpha \upxi_\beta = 0. \label{E:COCHARACTERISTICSARENULL}
	\end{align}

We are now ready to define what we mean by the ``characteristics'' of the compressible Euler system.

\begin{definition}[The characteristics for the compressible Euler system] \label{D:CHARACTERISTICSURFACESCANONICAL} 
 A surface $S \subset \R\times\Sigma$ (or curve $\ell \subset \R\times\Sigma$) is a \emph{characteristic} if at every point $p \in S$ (or $p \in \ell$), the co-normal space to $S$ (or $\ell$) at $p$, defined to be $\{\bm{\upxi} \in T_p^*(\R\times\Sigma) \, | \, \bm{\upxi}(V) = 0 \text{ for all } V \in  T_pS \}$ (or  $\{\bm{\upxi} \in T_p^*(\R\times\Sigma) \, | \, \bm{\upxi}(V) = 0 \text{ for all } V \in  T_p\ell \}$), is a subset of $C_p^*$ (or $P_p^*$).


\end{definition}


\begin{lemma}[Level sets of eikonal functions are characteristics as in Def.\,\ref{D:CHARACTERISTICSURFACESCANONICAL}] \label{L:LEVELSETSOFEIKONALEQUATIONSARECANONICALCHARS}
Let $U$ be any $C^1$ eikonal function, that is, any $C^1$ solution to $(\mathbf{\gfour}^{-1})^{\alpha\beta}\partial_\alpha U \partial_\beta U = 0$. 
Then the level sets of $U$ are characteristics in the sense of Def.\,\ref{D:CHARACTERISTICSURFACESCANONICAL}. In particular, 
the outgoing and ingoing characteristic surfaces $\outgoingcharacteristicsurfacetwoarg{u}{I}$ and $\ingoingcharacteristicsurfacetwoarg{\ubar}{J}$ defined by  \eqref{E:OUTGOINGCHRACTERISTICSURFACE}--\eqref{E:INGOINGCHARACTERISTICSURFACE} are both characteristics.
\end{lemma}
\begin{proof}
It is a standard fact from differential geometry that the co-normal spaces to the level sets of $U$ are spanned by the one-form $\bm{\upxi} = \upxi_\alpha \mathrm{d} x^\alpha$ with components $\upxi_\alpha := \partial_\alpha U$. We compute that:
	\begin{align}
			0 = (\gfour^{-1})^{\alpha\beta}\p_\alpha U \p_\beta U  = - (\partial_t U + v^a \partial_a U)^2 
			+ 
			\sum_{a=1}^2 \Speed^2 (\partial_a U)^2 & & \Longrightarrow & & \partial_t U + v^a\p_a U = \pm \Speed |\p U|,
	\end{align}
where $|\p U| = (\p_{x^1} U)^2 + (\p_{x^2}U)^2$. In particular, $\bm{\upxi} = \mathrm{d} U$ is a sound cone co-characteristic, as desired.  
\end{proof}

\subsection{Trivial notational differences in the physical fluid variables} In this short section, we record some notational differences between \cite{abbrescia2022emergence,shkoller2024geometry}. Since \cite{shkoller2024geometry} treats only adiabatic
equations of state $p(\varrho) = \tfrac{1}{\gamma} \varrho^\gamma$ for $\gamma > 1$, it follows that the speed of sound in \cite{shkoller2024geometry} is 
$\Speed = \sqrt{p'(\varrho)} = \varrho^{\frac{\gamma-1}{2}}$. In \cite{shkoller2024geometry}*{Sect.1.5.1}, the authors 
often use the rescaled sound speed $\sigma:= \frac{1}{\alpha} \Speed$ where $\alpha = \tfrac{\gamma-1}{2}$ instead of $\Speed$. We thus have the following table,
where the ``$\sim$'' in the table is explained in the subsequent paragraph and in the discussion surrounding \eqref{E:ANALOGYBETWEENRIEMANNINVARIANTS}:
\begin{table}[h!]
\centering
\begin{tabular}{|c|c|c|c|}
\hline
 &   \cite{abbrescia2022emergence} and this paper & & \cite{shkoller2024geometry} \\ \hline
 Fluid velocity & $v$  & = & $u$ \\ \hline
Sound speed  &  $\Speed$ & =  & $\alpha\sigma$  \\ \hline
Riemann invariants & $\RRiemann,\LRiemann$  & $\sim$ & $w,z$ \\
\hline \end{tabular}
\end{table}

For the $2D$-isentropic Euler system, there are three physical fluid unknowns. The choice of these unknowns corresponding to our results are $\RRiemann, \LRiemann, \LogDensity$, where $\RRiemann, \LRiemann$ are the ``almost'' Riemann invariants of \eqref{E:ALMOSTRIEMANNINVARIANTS}, and $\LogDensity$ is the logarithmic density as in \eqref{E:LOGDENS}. Instead of $\LogDensity$, \cite{shkoller2024geometry} adopt the \emph{tangential velocity} $u \cdot \tau = u^i \tau^i$ as their third unknown, 
where $\tau^i$ is a basis for the tangent space of the outgoing characteristics restricted to constant time-values; see \eqref{E:DEFINITIONOFTANGENTVECTORSFORSTEVEOUTGOINGCHAR}. We also note that $w,z$ defined in \cite{shkoller2024geometry}*{Equation (1.5)} are not exactly equal to $\RRiemann,\LRiemann$. Rather, in \cite{shkoller2024geometry},
the following definitions are adopted:
	\begin{align}
		w & := u\cdot n + \sigma, \qquad z := u \cdot n - \sigma, \label{E:STEVERIEMANN}
	\end{align}
where $n$ is the normal vector of the outgoing characteristics restricted to constant time-values; see \eqref{E:DEFINITIONOFNORMALVECTORSFORSTEVEOUTGOINGCHAR}. 
If we were to restrict our results to adiabatic equations of state,\footnote{Recall that we allow for any equation of state other than the Chaplygin gas.} then the analysis of Sect.\,\ref{SSS:ACOUSTICALCOORDINATESANDALECOORDINATES} (see \eqref{E:OUTGOINGCOCHARCHOICEDATA}) implies that at time $0$, 
the different definitions of Riemann invariants agree up to a constant
(e.g., $\RRiemann\restriction_{t=0} =( v^1 + \sigma)\restriction_{t=0} +C  = w\restriction_{t=0} + C$). However, 
we remind the reader that the data in \cite{shkoller2024geometry} are fundamentally different than the data we treat here, due to our assuming only transversal convexity instead of strict convexity; see Sect.\,\ref{SS:RESULTSFROMSINGULARBOUNDARYPAPER}. In this appendix, we will mostly use our notation for the fluid velocity, sound speed, and Riemann invariants. We will refer to $w,z$ only when making direct comparisons between \cite{shkoller2024geometry} and our work.

\subsection{Coordinates that regularize the singularity} \label{SSS:ACOUSTICALCOORDINATESANDALECOORDINATES}
Any method designed to study the solution until the formation of a shock must grapple with infinite gradients that occur in the standard Cartesian coordinate system $(t,x^1,x^2,x^3)$. In \cite{dC2007,abbrescia2022emergence}, the solution is studied in geometric coordinates $(t,u,x^2,x^3)$, relative to which it remains rather smooth, and the singularity manifests in a degeneration between the change of variables map $\Upsilon: (t,u,x^2,x^3) \to (t,x^1,x^2,x^3)$. Similarly, in \cite{shkoller2024geometry}, the solution is studied relative to \emph{Arbitrary Eulerian-Lagrangian} (ALE) coordinates, which is also a coordinate system in which the solution remains rather smooth, and the singularity again manifests in the change of variables back to Cartesian coordinates. Our main result in this section is a proof that the geometric coordinates and the ALE coordinates are identical; see Prop.\,\ref{P:ALETOGEOCOORD}.\footnote{Strictly speaking, only the geometric coordinates from \cite{abbrescia2022emergence} are equivalent to the ALE coordinates used in \cite{shkoller2024geometry}. Christodoulou defined geometric spherical variables $\vartheta^A$ to better capture the geometry of the null cones in \cite{dC2007}, which are foliated by spheres. Since the characteristic surfaces in \cite{abbrescia2022emergence,shkoller2024geometry} are foliated by tori, to facilitate our comparison of these two works, we restrict our attention to the tori case, which corresponds to $(x^2,x^3) \in \T^2$.}

We begin by outlining the construction of the characteristic surfaces 
(see Def.\,\ref{D:CHARACTERISTICSURFACESCANONICAL})
given in \cite{shkoller2024geometry}. \cite{shkoller2024geometry} constructs a vector-valued function 
$\eta(t,x) = (\eta^1(t,x),\eta^2(t,x))$ that solves the following nonlinear initial value problem for $i=1,2$
(see \cite{shkoller2024geometry}*{Equation\,(1.12)}):
	\begin{align}
		\geop{t} \eta^i(t,x^1,x^2) 
		& = v^i(t,\eta(t,x^1,x^2)) + \Speed(t,\eta(t,x^1,x^2)) n_i(t,\eta(t,x^1,x^2)), & & \eta^i(t,x^1,x^2) = x^i,\label{E:DEFINITIONOFSTEVEFLOWMAP}
	\end{align}
where $n = (n_1,n_2)$ is the co-vector defined as follows, with $|\geop{x^2} \eta|$ denoting the Euclidean length of $\geop{x^2} \eta$:
	\begin{align}
		n(t,x^1,x^2) := \frac{1}{|\geop{x^2} \eta|} \left(\geop{x^2} \eta^2(t,x^1,x^2),- \geop{x^2} \eta^1(t,x^1,x^2)\right). \label{E:DEFINITIONOFNORMALVECTORSFORSTEVEOUTGOINGCHAR}
	\end{align}
 \cite{shkoller2024geometry} also defines the vector $\tau = (\tau^1,\tau^2)$ as follows ((see \cite{shkoller2024geometry}*{Equation\,(1.11a)})):
	\begin{align}
		\tau(t,x^1,x^2) := \frac{1}{|\geop{x^2} \eta|} \left(\geop{x^2} \eta^1(t,x^1,x^2), \geop{x^2} \eta^2(t,x^1,x^2)\right). 
		\label{E:DEFINITIONOFTANGENTVECTORSFORSTEVEOUTGOINGCHAR}
	\end{align}
In Prop.\,\ref{P:ALETOGEOCOORD}, we will prove that $\eta$ represents the spatial components of a null-geodesic flow map. It will then become clear that $n$ and $\tau$ are respectively spatial co-normal and tangent vectors of the outgoing characteristics at a fixed Cartesian time.

To define a characteristic surface as in Def.\,\ref{D:CHARACTERISTICSURFACESCANONICAL}, 
we consider the spacetime one-form $\mathbf{n}$ with the following Cartesian components:
\begin{align} \label{E:CARTESIANCOMPONENTSOFSPACETIMEONEFORM}
	(n_0,n_1,n_2),
\end{align}
where:
		\begin{align}
			n_0 &: = - v^a n_a - \Speed = - v \cdot n - \Speed. \label{E:OUTGOINGCOCHARCHOICE} 
		\end{align}
The discussion just above \eqref{E:COCHARACTERISTICSARENULL} implies that
$\mathbf{n}$ is co-characteristic, i.e., $\mathbf{n} \restriction_p \in C_p^*$. 
Next, we note that \eqref{E:DEFINITIONOFSTEVEFLOWMAP}--\eqref{E:OUTGOINGCOCHARCHOICE}  
imply that $\mathbf{n}$ has the following initial condition:
		\begin{align}
			n_0\restriction_{t=0}= -( v^1 - \Speed)\restriction_{t=0}, \qquad (n_1,n_2)\restriction_{t=0} & := (1,0). \label{E:OUTGOINGCOCHARCHOICEDATA}
		\end{align}
	Next, we use \eqref{E:ACOUSTICALMETRICANDMATERIALDER} to compute that the $\gfour$-dual vectorfield to $\mathbf{n}$ is,
	up to a harmless conformal factor of $- \Speed$, the following vectorfield:
	\begin{align}
		\Lunit^{ALE} 
		& := \frac{1}{- \Speed} \gfour^{-1}\cdot \mathbf{n} = \p_t + (v^a + \Speed n_a) \partial_a. \label{E:ALELUNIT}
	\end{align}
From \eqref{E:COCHARACTERISTICSARENULL}, it follows immediately that $\Lunit^{ALE}$ is $\gfour$-null:
	\begin{align}
		\gfour(\Lunit^{ALE},\Lunit^{ALE}) = 0.  \label{E:ALELUNITISNULL}
	\end{align}
The authors of \cite{shkoller2024geometry} denote the Cartesian spatial components of the vectorfield $\Lunit^{ALE}$ by $\mathcal{V}_3$, 
i.e., $\Lunit^{ALE} = \partial_t + \mathcal{V}_3^a \partial_a$, where
$\mathcal{V}_3^i := v^i + \Speed n_i$; see \cite{shkoller2024geometry}*{Equation\,(1.7)}.

Let $\Lambda_{ALE} = (\Lambda_{ALE}^0,\Lambda_{ALE}^1,\Lambda_{ALE}^2)$ denote the flow map of $\Lunit^{ALE}$, i.e.,
	\begin{align}
		\geop{t} \Lambda_{ALE}^\alpha(t,x^1,x^2) 
		&= (\Lunit^{ALE})^\alpha \circ \Lambda_{ALE}(t,x^1,x^2), & & \Lambda_{ALE}(0,x^1,x^2) = (0,x^1,x^2) \in \{0\}\times\T^2. \label{E:FLOWMAPOFALELUNIT}
	\end{align}
Note that \eqref{E:DEFINITIONOFSTEVEFLOWMAP} and \eqref{E:FLOWMAPOFALELUNIT} imply that $\eta^i$ and $\Lambda_{ALE}^i$ solve the same ODE with the same initial conditions. From ODE uniqueness, we conclude that $\Lambda_{ALE}^i = \eta^i$ for $i = 1,2$. 
Then, given $x^1 \in \T$, \cite{shkoller2024geometry} defines the truncated(-in time) characteristic surfaces $\Gamma_{x^1}(T)$ to be the following image set: 
	\begin{align}
		\Gamma_{x^1}(T) := \{ \Lambda_{ALE}(t,x^1,x^2) \, : \, (t,x^2) \in [0,T]\times \T \}. \label{E:CHARACTERISTICS}
	\end{align}
In Prop.\,\ref{P:ALETOGEOCOORD}, we prove that $\Gamma_{x^1}(T)$ are characteristic surfaces 
in the sense of Def.\,\ref{D:CHARACTERISTICSURFACESCANONICAL}. 

We highlight now the following crucial point: 
	\begin{quote}
		By defining the characteristic surfaces by \eqref{E:CHARACTERISTICS}, the authors of \cite{shkoller2024geometry} \emph{are no longer using $x^1$ to denote the first spatial Cartesian coordinate of the ambient spacetime}.\footnote{The spatial manifold used in \cite{shkoller2024geometry} is $\Sigma = \T^2$, whereas a $2D$-analog of our work \cite{abbrescia2022emergence} would have $\Sigma = \R\times\T$. This minor point is not important for the present discussion.} Instead, the symbol $x^1$ in \cite{shkoller2024geometry} denotes a label that identifies the characteristic surface emanating from the torus $\{x^1\}\times\T \subset \lbrace t = 0 \rbrace$ under the flow of $\Lunit^{ALE}$. The notation given in \cite{shkoller2024geometry} for the ambient Cartesian spatial coordinate directions is $e_1 := (1,0)$ and $e_2 := (0,1)$. In particular, the partial differentiation in $x^a$ in \eqref{E:ALELUNIT} is partial differentiation in the Cartesian direction $e_a$, \emph{not} with respect to 
		$x^a$ in the ALE-coordinate differential structure. To prevent further confusion with the Cartesian coordinate $x^1$ that we use in \cite{abbrescia2022emergence} and this paper, we denote the ``non-Cartesian'' spatial coordinate introduced in \cite{shkoller2024geometry} by $x^1_{ALE}$. We note that $x^1_{ALE} \restriction_{t=0}  = x^1 \restriction_{t=0}$, i.e., initially, the ALE coordinate system is the same as the Cartesian one.
	\end{quote}
	
The \emph{ALE coordinates} are defined to be $(t,x^1_{ALE},x^2)$. For each $t \in [0,T]$, the ALE diffeomorphism $\psi_t : \T^2_{(x^1_{ALE},x^2)} \to \T^2_{(e_1,e_2)}$ is defined by $\psi_t(x^1_{ALE},x^2) = (h(t,x^1_{ALE},x^2),x^2)$; see \cite{shkoller2024geometry}*{Section~1.6}. 
Here, for each fixed $(t,x^1_{ALE}) \in [0,T]\times\T$, the function $h(t,x^1_{ALE},\cdot)$ is defined by $x^2 \rightarrow h(t,x^1_{ALE},x^2) := \Lambda^1(t,x^1_{ALE},x^2)$, 
i.e., at each fixed pair $(t,x^1_{ALE})$,
$h(t,x^1_{ALE},\cdot)$ is the ``$x^1$ Cartesian-height'' (in the $e_1$-direction) of $\Lambda_{ALE}(t,x^1_{ALE},\cdot)$, viewed as a graph over $\T$. 

\begin{remark}[Notation for partial differentiation] In light of the forth-coming equivalence of the geometric coordinates and ALE coordinates, we adopt the same notation for partial differentiation, that is $\left\{\geop{t},\geop{u},\geop{x^2} \right\} = \left\{\geop{t},\geop{x^1_{ALE}},\geop{x^2} \right\}$. We continue to denote partial differentiation with respect to the Cartesian differential structure by $\{\p_{t} = \p_{0},\p_{1},\p_{2}\}$.
\end{remark}

We will now prove that ALE coordinates are identical to the acoustical geometric coordinates afforded by Christodoulou's framework \cite{dC2007}. 
We will also prove that the surfaces $\Gamma_{x^1}(T)$ defined by \eqref{E:CHARACTERISTICS} are identical to the truncated characteristics
$\nullhyparg{u}^{[0,T]}$ that we use here and in \cite{abbrescia2022emergence},
which are portions of level sets of a solution $u$ to the eikonal equation \eqref{E:REGULAREIKONALEQUATION}.

\begin{proposition}[ALE coordinates are identical to acoustical geometric coordinates] \label{P:ALETOGEOCOORD}
	 Let $\eta, n, \tau$ be as in \eqref{E:DEFINITIONOFSTEVEFLOWMAP}--\eqref{E:DEFINITIONOFTANGENTVECTORSFORSTEVEOUTGOINGCHAR}, where we recall that equation \eqref{E:DEFINITIONOFSTEVEFLOWMAP} is solved in conjunction with the compressible Euler equations. Let $\mathbf{n}= (n_0,n_1,n_2)$ be the $\mathbf{g}$-null one-form from 
	\eqref{E:OUTGOINGCOCHARCHOICE}, let $\Lunit^{ALE}$ be as in \eqref{E:ALELUNIT}, 
	and let $\Lambda_{ALE}$ be the flow-map of $\Lunit^{ALE}$ as in \eqref{E:FLOWMAPOFALELUNIT}. 
	Then for $C^1$ solutions to the compressible Euler equations, the following hold:
			\begin{enumerate}
				\item \underline{\textbf{The flow of $\Lunit^{ALE}$}}. The flow map $\Lambda_{ALE}$ satisfies $\Lambda_{ALE}^0 = t$ and $\Lambda_{ALE}^i = \eta^i$ for $i = 1,2$.
				\item \underline{\textbf{$\Gamma_{x^1}(T)$ are characteristics in the sense of Def.\,\ref{D:CHARACTERISTICSURFACESCANONICAL}}}. 
				Recall that $\Gamma_{x^1}(T)$ is the surface defined in \eqref{E:CHARACTERISTICS}.
				For any $\Gamma_{x^1}(T)$-tangent vectorfield $W$, we have $\mathbf{n}_\alpha W^\alpha = 0$. Hence, $\Gamma_{x^1}(T)$ is a characteristic surface as in Def.\,\ref{D:CHARACTERISTICSURFACESCANONICAL}.
				\item \underline{\textbf{$\Gamma_{x^1}(T)$ are $\gfour$-null hypersurfaces}}. There exists a $\Gamma_{x^1}(T)$-tangent vectorfield $V$ 
					that is $\gfour$-orthogonal to the tangent space of $\Gamma_{x^1}(T)$. In particular, $\gfour(V,V) = 0$.
				\item \underline{\textbf{The geometric coordinates are identical to ALE coordinates}}. Suppose that $u$ is the solution to the eikonal equation $(\gfour^{-1})^{\alpha\beta}\p_\alpha u \p_\beta u = 0$ with $u \restriction_{t=0} = x^1$ and\footnote{For the proposition to be correct, $\partial_t u < 0$ is the right sign; $\partial_t u > 0$ would correspond to constructing the ingoing characteristics. We point out that  we have $\partial_t u > 0$ in \eqref{E:EIKONALBRANCHINTRO}. The reason is that in the bulk of the paper, we chose the initial data \eqref{E:EIKONALDATAINTRO} to have the opposite sign compared to this proposition.} 
				$\partial_t u \restriction_{t=0} < 0$. 
				Then $x^1_{ALE} = u$. In addition, the vectorfield $\Lunit$ defined by \eqref{E:LUNIT} satisfies $\Lunit = \Lunit^{ALE}$.				
				\item \underline{\textbf{$\Lunit^{ALE}$ is a non-affinely  parametrized geodesic vectorfield}}. There is a non-vanishing function $\smoothfunction$ such that $\Dfour_{\Lunit^{ALE}}\Lunit^{ALE} = \smoothfunction \Lunit^{ALE}$, where $\Dfour$ is the Levi-Civita covariant derivative defined by $\gfour$. In particular, the evolution equation \eqref{E:DEFINITIONOFSTEVEFLOWMAP} satisfied by $\eta$ is a rewriting of $\frac{1}{\smoothfunction} \Dfour_{\Lunit^{ALE}}(\Lunit^{ALE})^i = (\Lunit^{ALE})^i$. 
			\end{enumerate}
		
	\end{proposition}

	\begin{proof} \hfill \\
	\noindent \underline{\emph{Proof of 1.}} $\Lambda_{ALE}^0 = t$ follows $\Lunit^{ALE}t = 1$ and $\Lambda_{ALE}(0,x^1_{ALE},x^2) =(0, x^1_{ALE},x^2)$. The statement that $\Lambda_{ALE}^i = \eta^i$ was already proved below \eqref{E:FLOWMAPOFALELUNIT}. 
	
	\medskip

	\noindent \underline{\emph{Proof of 2.}} 
	Recall that  $\mathbf{n}$ is, by construction, a co-characteristic. Since $\Gamma_{x^1}(T)$ is a co-dimension 1 submanifold of $\R \times \Sigma$, it follows that its co-normal space at any point in  $\Gamma_{x^1}(T)$ is one dimensional. Therefore, it suffices to prove that
	for any $\Gamma_{x^1}(T)$-tangent vectorfield $W$, we have
	$\mathbf{n} \cdot W := n_\alpha  W^\alpha = 0$. We will prove this relation for $W$ belonging to a basis for the tangent space of the surface.
	From the parametrization of \eqref{E:CHARACTERISTICS} by $(t,x^2) \in [0,T]\times \T$, it follows that $\frac{\p}{\p x^2} \Lambda_{ALE}^\alpha(t,x^1_{ALE},x^2)$ are the components of the tangent vector which is the push-forward of $\geop{x^2}\restriction_{t=0} = \p_{x^2}\restriction_{t=0}$ by $\Lambda_{ALE}$, which is tangent to $\{(0,x^1_{ALE} \restriction_{t=0})\}\times\T$. We note further that since $\Lambda_{ALE}^0(t,x^1_{ALE},x^2) = t$, the vectorfield $\frac{\p}{\p x^2} \Lambda_{ALE}(t,x^1_{ALE},x^2)$ is $\Sigma_t$-tangent, i.e., $\frac{\p}{\p x^2} \Lambda_{ALE}^0(t,x^1_{ALE},x^2) = 0$. Using that $\Lambda_{ALE}^i = \eta^i$ for $i = 1,2$, we compute,
	with the help of \eqref{E:DEFINITIONOFNORMALVECTORSFORSTEVEOUTGOINGCHAR}--\eqref{E:DEFINITIONOFTANGENTVECTORSFORSTEVEOUTGOINGCHAR}, that:
	\begin{align}
		\mathbf{n}\cdot \geop{x^2} \Lambda_{ALE} = n_a \geop{x^2} \Lambda_{ALE}^a = \left|\geop{x^2} \eta\right| n_a \tau^a =0.
	\end{align}

Moreover, from the definition \eqref{E:FLOWMAPOFALELUNIT} of the flow map, it follows 
that $\Lunit^{ALE}$ is $\Gamma_{x^1}(T)$-tangent. From \eqref{E:ALELUNIT}--\eqref{E:ALELUNITISNULL}, we have $\mathbf{n} \cdot \Lunit^{ALE} = 0$.

Since $\{\Lunit^{ALE},\geop{x^2} \Lambda_{ALE}\}$ span the tangent space to $\Gamma_{x^1_{ALE}}(T)$ at any point, we have proved that $\Gamma_{x^1_{ALE}}(T)$ is a characteristic surface in the sense of Def.\,\ref{D:CHARACTERISTICSURFACESCANONICAL}.
	
	\medskip

	\noindent \underline{\emph{Proof of 3.}} Since $\Lunit^{ALE} = -\frac{1}{\Speed} \gfour^{-1} \cdot \mathbf{n}$ (see \eqref{E:ALELUNIT}), 
	the proof of Point 2 implies that $\gfour(\Lunit^{ALE},\geop{x^2} \Lambda_{ALE}) = 0$. From \eqref{E:ALELUNITISNULL} 
	and the aforementioned fact that $\{\Lunit^{ALE},\geop{x^2} \Lambda_{ALE}\}$ spans the tangent spaces to $\Gamma_{x^1_{ALE}}(T)$,
	it follows that $\Lunit^{ALE}$ is $\mathbf{g}$-null and $\gfour$-orthogonal to the tangent spaces of $\Gamma_{x^1_{ALE}}(T)$, 
	i.e., $\Gamma_{x^1_{ALE}}(T)$ is a $\gfour$-null hypersurface.
	
	\medskip
	
	\noindent \underline{\emph{Proof of 4.}} Since $x^1_{ALE}$ is constant along each characteristic $\Gamma_{x^1_{ALE}}(T)$, it follows that the one-form $\mathrm{d} x^1_{ALE}$ is co-normal to it. By Point 3, it follows that $\mathrm{d} x^1_{ALE}$ must be parallel to the $\gfour$-dual of 
$\Lunit^{ALE}$, namely it must be parallel to the co-characteristic $\mathbf{n}$. By \eqref{E:COCHARACTERISTICSARENULL}, we have that $0 = (\gfour^{-1})(\mathrm{d} x^1_{ALE},\mathrm{d} x^1_{ALE})$, 
and hence $x^1_{ALE}$ is an eikonal function. Since $\Lunit^{ALE} x^1_{ALE} = 0$ and $\partial_{x^1} x^1_{ALE}\restriction_{t=0} = 1$, it follows from \eqref{E:ALELUNIT} that 
for subsonic solutions\footnote{The solutions studied here and in \cite{shkoller2024geometry} are subsonic. \label{FN:SUBSONIC}} 
(in which $|v|$ from \eqref{E:ALELUNIT} is smaller than $\Speed$), we have
$\partial_t x^1_{ALE}\restriction_{t=0} < 0$.
In particular, $u$ and $x^1_{ALE}$ are $C^1$ solutions to the eikonal equation with the same initial data and with negative time derivatives at time $0$.
Standard arguments, based on differentiating the eikonal equation and then deriving 
estimates for the differences $\partial_{\alpha}(u - x^1_{ALE})$, can be used to show that $u - x^1_{ALE} = 0$, i.e., that $u = x^1_{ALE}$, as is desired.
Moreover, it follows that the one-forms $\mathrm{d} u$ and $\mathbf{n}$ are both co-normal to the characteristics and hence must be parallel.
Therefore, the vectorfields $\Lunit^{ALE}$ (defined in \eqref{E:ALELUNIT}) and $\Lunit$ (defined in \eqref{E:LUNIT})
must be parallel. Since $\Lunit t = 1 = \Lunit^{ALE} t$, we in fact must have $\Lunit = \Lunit^{ALE}$, as is desired.

\medskip
	
\noindent \underline{\emph{Proof of 5.}}
	Let $u$ be as in Point 4, and recall that we showed that $u = x^1_{ALE}$. A straightforward calculation shows that its $\gfour$-gradient defined by\footnote{Since $\partial_t u > 0$ in \eqref{E:REGULAREIKONALEQUATION}, here we have changed the minus sign compared to \eqref{E:LGEO} in order to be consistent with the sign conventions of this appendix.}  $\Lgeo^\alpha := (\gfour^{-1})^{\alpha\beta} \p_\beta u$ is an affinely parameterized null geodesic. That is, $\Dfour_{\Lgeo}\Lgeo = 0$. 
	By $\Lunit := \upmu \Lgeo$ and Point 4, we compute $\Dfour_{\Lunit^{ALE}} \Lunit^{ALE} = \Dfour_{\Lunit} \Lunit = \upmu  \Dfour_{\Lgeo} (\upmu \Lgeo) = \upmu \Lgeo(\upmu) \Lgeo  = (\Lunit \upmu) \Lgeo = \frac{\Lunit \upmu}{\upmu} \Lunit = \frac{\Lunit \upmu}{\upmu} \Lunit^{ALE} $, as is desired.
	\end{proof}

	Next, we note that the relationship between our change of variables map $\Upsilon$ from geometric coordinates to Cartesian coordinates (see \eqref{E:CHOVGEOTOCARTESIAN} for the $3D$ version) and the map $\psi_t$ from \cite{shkoller2024geometry} is 
$\Upsilon(t,u,x^2) = (t,\psi_t(u,x^2))$. 
The truncated characteristic surfaces $\Gamma_{x^1_{ALE}}(T)$ are precisely the same as our $\nullhyparg{u}^{[0,T]}$; 
see \eqref{E:TRUNCATEDNULLHYPERSURFACES}. That is the reason that  
the solution stays relatively smooth in the differential structure of the ALE coordinates $(t,x^1_{ALE},x^2)$ in \cite{shkoller2024geometry}, just like in the case of the geometric acoustical eikonal coordinates $(t,u,x^2)$ in \cite{dC2007,abbrescia2022emergence}. 
The gradient singularity is also captured by the same mechanism: the vanishing of $\upmu \sim \det D\Upsilon(t,u,x^2) = \det D \psi_t(x^1_{ALE},x^2)$, which in \cite{shkoller2024geometry} is denoted by $J_g$; 
see \cite{shkoller2024geometry}*{Equations\,(2.12) and (2.16a)}.

\subsection{Anisotropic $L^2$-differentiated solution variables used in the energy method}   \label{SS:SOLUTIONVARIABLESUSEDINENERGYMETHOD}
For convenience, we now denote the $2D$ isentropic fluid variables as $\wavearray := \{\Psi^0,\Psi^1,\Psi^2\} = \{\LogDensity, v^1,v^2\}$. All of \cite{dC2007,abbrescia2022emergence,shkoller2024geometry}, as well as this paper, derive $L^2$-based energy estimates for quantities that are at the level of regularity level $\pmb{\partial}\wavearray$. For convenience, we henceforth refer to these quantities as the ``$L^2$-differentiated solution variables''. The energy estimates in \emph{all} of these works are directionally anisotropic with respect to powers of $\upmu = J_g$. By this, we mean that tensorial contractions of the expressions $\p_{\alpha}\Psi^\beta$ (which make up the $L^2$-differentiated solution variables) with directions that are \emph{transverse} to the characteristics $\nullhyparg{u}$ feature an additional weight of $\upmu = J_g$ compared to \emph{tangential} contractions. We note, however, that there are substantial differences in the energy methods used to control the vorticity scalar in $2D$ \cite{shkoller2024geometry} compared to the $3D$ vorticity vectorfield of \cite{abbrescia2022emergence} and this paper. 
Here, we also developed new energy methods to incorporate the dynamics of the entropy $\Ent$.

To aid, simplify, and unify the notion for the anisotropic energies featured in \cite{dC2007,abbrescia2022emergence,shkoller2024geometry}, in this section, we restrict our attention to two types of regions:
	\begin{subequations}
		\begin{align}
			M_{t,u} & := \bigcup_{t' \in [0,t]} \Sigma_{t'}^{  [u_1,u_2]} = \bigcup_{u \in [u_1,u_2]} \nullhyparg{u}^{[0,t]} , \label{E:SPACELIKEORCHARACTERISTICFOLIATIONS} \\
			\mathcal{M}_t & :=  \bigcup_{t' \in [0,t]} \Sigma_{t'}. \label{E:SPACELIKEFOLIATIONSOFSPACETIME}
		\end{align}
	\end{subequations}
We view $M_{t,u}$ as regions of spacetime that can be foliated by either spacelike hypersurfaces $\Sigma_t^{[u_1,u_2]}$, whose boundaries are $\partial \Sigma_t^{[u_1,u_2]} =  \ell_{t,u_1} \cup \ell_{t,u}$,  or by outgoing characteristics $P_u^{[0,t]}$, whose boundaries are $\partial \nullhyparg{u}^{[0,t]} = \ell_{0,u}\cup\ell_{t,u}$. We view $\mathcal{M}_t$ as regions of spacetime that are foliated by spacelike hypersurfaces \emph{without boundary}. In particular, this is the case considered in \cite{shkoller2024geometry}, as their spatial topology is homeomorphic to $\T^2$.

\subsubsection{The $L^2$-differentiated solution variables in \cite{dC2007}}  \label{SSS:L2SOLUTIONVARIABLESFORCHRISTODOULOU}
We now explain how the anisotropicity highlighted above manifests in \cite{dC2007}. 
 Christodoulou's monograph focused on the $3D$ irrotational and isentropic relativistic Euler equations.  
In that setting, the equations of motion reduce to a quasilinear wave equation
$(\gfour^{-1})^{\alpha\beta}(\partial \Phi) \p_\alpha\p_\beta\Phi = 0$ 
for a potential function $\Phi$ satisfying, schematically,
$\pmb{\partial} \Phi \sim \wavearray$. 
Christodoulou then took a rectangular coordinate of this equation and derived the following \cite{dC2007}*{Equation (1.86)}:
	\begin{align}
		0 & = \Box_{\tilde{\gfour}} (\partial_\kappa \Phi) = |\mbox{det}(\tilde{\gfour})|^{-1/2} 
		\partial_{\alpha}\left(|\mbox{det}(\tilde{\gfour})|^{1/2}\tilde{\gfour}^{\alpha\beta}\p_\beta (\partial_\kappa\Phi)\right), \label{E:CHRISTODOULOUGEOMETRICEQUATION} 
	\end{align}
where $\tilde{\gfour} = F(\partial \Phi) \gfour$ is a conformal rescaling of the acoustical metric for a well-chosen $F > 0$, 
and $\kappa \in \{0,1,2,3\}$. For solutions to equation \eqref{E:CHRISTODOULOUGEOMETRICEQUATION}, the $L^2$-energy method yields the following identity 
(see \cite{dC2007}*{Chapter 5} or \eqref{E:ENERGYNULLFLUXINTEGRALIDENTITIESWAVE}) on $M_{t,u}$:
	\begin{align}
		\begin{split} \label{E:CHRISTODOULOUENERGYAPPENDIX}
			 & \int_{\Sigma_t^{[u_1,u_2]}}\left\{\ ( \newuL  \p_\kappa \Phi)^2 + \upmu (Y \p_\kappa \Phi)^2 + \upmu (\Lunit \p_\kappa \Phi)^2\right\} + \int_{\nullhyparg{u}^{[0,t]}} \left\{ (\Lunit   \p_\kappa \phi)^2 + \upmu (Y \p_\kappa \phi)^2 \right\}  \\
			 & \ \ + \int_{M_{t,u}} (-\Lunit \upmu) (Y  \p_\kappa \Phi)^2  = \textnormal{data} + \cdots,
		\end{split}
	\end{align}
where 
we have suppressed the volume forms for simplicity. Notice that \eqref{E:CHRISTODOULOUENERGYAPPENDIX} provides $L^2$-estimates for $\Lunit \partial_\kappa \Phi$ on the outgoing characteristic $\nullhyparg{u}^{[0,t]}$ \emph{without a $\upmu$-weight}, and for $\newuL \partial_\kappa \Phi = \upmu \uLunit \partial_\kappa \Phi$, which features a $\upmu$ weight.\footnote{Christodoulou's notation for $\newuL$ is in fact ``$\uLunit$''; compare \eqref{E:MUXINTERMSOFGEOMETRICCOORDINATEVECTORFIELDS},\eqref{E:CONVENIENTIDENTITYFORULUNIT}, and \cite{dC2007}*{Equation (5.39)}. Hence, the $L^2$-hypersurface integral \cite{dC2007}*{Equation (5.47)} does in fact measure  $\newuL \partial_\kappa \Phi$ in $L^2$.} In addition, the spacetime integral on LHS~\eqref{E:CHRISTODOULOUGEOMETRICEQUATION} provides $L^2_{t,u}$-control over
$Y \p_\kappa \Phi$ \emph{without a $\upmu$-weight}. Note that $\partial_\kappa\Phi$ is at the same regularity level of the fluid unknowns $\wavearray$ and that differentiation with respect to\footnote{In two spatial dimensions, there is only one $\ell_{t,u}$-tangent vectorfield compared to the two $\Yvf{2},\Yvf{3}$ present in this paper.} $\Lunit$ or $Y$ is equivalent to differentiating with respect to Cartesian derivatives $\p_\alpha \p_\kappa \Phi$ and then contracting with \emph{tangential} directions $L^\alpha$ or $Y^\alpha$. Similarly, differentiating with resepct to $\newuL$ is equivalent to differentiating with respect to Cartesian derivatives $\p_\alpha \p_\kappa \Phi$, contracting with \emph{transverse} directions $\uLunit^\alpha$, and multiplying by $\upmu= J_g$.

\subsubsection{The $L^2$-differentiated solution variables in \cite{abbrescia2022emergence} and this paper} \label{SSS:L2SOLUTIONVARIABLESFORLEOJARED}
We now briefly explain how the $L^2$-differentiated solution variables are constructed in \cite{abbrescia2022emergence} and this paper. For the rest of this appendix, when referencing the present paper and our work \cite{abbrescia2022emergence}, it should be understood that we are discussing the full $3D$ compressible Euler system. 
We work directly with the geometric-wave-transport formulation of the compressible Euler system as in \eqref{E:COVWAVEEQv}--\eqref{E:MODIFIEDFLUIDVARIABLES}. 
The $L^2$ energy method for the wave equations $\Box_{\gfour} \wavearray = \cdots$ satisfied by the wave variables $\wavearray = \{\LogDensity,v^1,v^2,v^3,\Ent\}$ on $M_{t,u}$ yields the following identity (as in \eqref{E:ENERGYNULLFLUXINTEGRALIDENTITIESWAVE}):
	\begin{align}
		\begin{split} \label{E:EULERENERGYAPPENDIX}
			 & \int_{\Sigma_t^{[u_1,u_2]}}\left\{\ (\upmu \newuL \wavearray)^2 + \upmu (Y  \wavearray)^2 + \upmu (\Lunit  \wavearray)^2\right\}  + \int_{\nullhyparg{u}^{[0,t]}} \left\{ (\Lunit  \wavearray)^2 + \upmu (Y  \wavearray)^2 \right\} \\
			 & \ \ + \int_{M_{t,u}} (-\Lunit \upmu) (Y    \wavearray )^2   = \textnormal{data} + \cdots,
		\end{split}
	\end{align}

As we explained in detail in Sects.\,\ref{SS:RESULTSFROMSINGULARBOUNDARYPAPER}--\ref{SS:NEWIDEASANDMETHODS}, the $L^2$-solution variables in 
\eqref{E:ENERGYNULLFLUXINTEGRALIDENTITIESWAVE} are $\Lunit \wavearray$,
$\newuL \wavearray = \upmu \uLunit \wavearray$, 
and $\sqrt{\upmu} Y \wavearray$,
similar to Sect.\,\ref{SSS:L2SOLUTIONVARIABLESFORCHRISTODOULOU}. Our complete energy method 
relies on numerous technical innovations to complement \eqref{E:EULERENERGYAPPENDIX},
such as the elliptic-hyperbolic integral identities of Sect.\,\ref{S:ELLIPTICHYPERBOLICIDENTITIES}; see Sect.\,\ref{SS:NEWIDEASANDMETHODS} for a detailed overview.

\subsubsection{The $L^2$-differentiated solution variables in \cite{shkoller2024geometry}} \label{SSS:DIFFERENTIATEDSOLUTIONVARSANDMUWEIGHTEDEQUATIONS}
In this section, we explain how the analysis of \cite{shkoller2024geometry} also uses the same fundamental properties of ``anisotropic $L^2$-differentiated solution variables'' 
introduced at the start of Sect.\,\ref{SS:SOLUTIONVARIABLESUSEDINENERGYMETHOD}. Recall the notation introduced above in Sect.\,\ref{SSS:ACOUSTICALCOORDINATESANDALECOORDINATES}.

A technical difference between \cite{dC2007,abbrescia2022emergence} and \cite{shkoller2024geometry} is that the latter does \emph{not} treat the compressible Euler equations as a second order geometric wave-transport-Hodge system. 
Instead, the fundamental $L^2$-differentiated solution variables in \cite{shkoller2024geometry} solve a \emph{first}-order system. To derive these equations, the authors of \cite{shkoller2024geometry} differentiate \eqref{E:INTROTRANSPORTVI}--\eqref{E:INTROTRANSPORTDENSITY} with respect to the Cartesian derivative $\partial_k$. They then re-write the resulting differentiated equations 
in terms of the characteristic vectorfield $\Lunit^{ALE}$ by adding and subtracting the contribution of $c n_j \partial_j$, where $n_j$ is as in \eqref{E:DEFINITIONOFNORMALVECTORSFORSTEVEOUTGOINGCHAR},
thereby obtaining:
	\begin{subequations} \label{E:DERIVATIVELOSINGEQUATIONS}
		\begin{align}
			\Lunit^{ALE} (\partial_k v^i)  - \Speed  n_j\p_j (\partial_k v^i) & = - \frac{1}{\varrho} \Speed^2 \partial_k  \p_i \varrho  - 2 \Speed_{;\varrho}\Speed \partial_k\varrho  \partial_i \varrho + \frac{1}{\varrho^2} \Speed^2 \partial_k\varrho \partial_i\varrho - \partial_k v^j \partial_j v^i, \label{E:VELOCITYEQREWRITTENINTERMSOFLUNITALE} \\
			\Lunit^{ALE} (\partial_k \varrho) -  \Speed n_j\p_j (\partial_k \varrho) & = - \varrho \p_k\partial_i v^i - (\partial_k \varrho) \partial_i v^i - (\partial_k v^i) \partial_i \varrho. \label{E:DENSITYEQREWRITTENINTERMSOFLUNITALE}
		\end{align}
	\end{subequations}
Now, instead of solving \eqref{E:VELOCITYEQREWRITTENINTERMSOFLUNITALE}--\eqref{E:DENSITYEQREWRITTENINTERMSOFLUNITALE} for $\{\partial_k v^1, \partial_k v^2, \partial_k \varrho\}_{k \in \{1,2\}}$, Shkoller-Vicol study the corresponding equations satisfied by the following ``differentiated Riemann variables,'' for $k = 1,2$: 
	\begin{align}
		W_k :=  n_i \partial_k v^i + \partial_k \sigma, && Z_k := n_i \partial_k v^i - \partial_k\sigma, && A_k:=  \tau_i \partial_k v^i,  \label{E:SHKOLLERDIFFERENTIATEDRIEMANNVAIRABLES}
	\end{align}
where $\tau$ is as in \eqref{E:DEFINITIONOFTANGENTVECTORSFORSTEVEOUTGOINGCHAR} and $\sigma = \frac{1}{\alpha} \Speed$ (recall \cite{shkoller2024geometry} considered only equations of states of the form $p = \frac{1}{2\alpha + 1} \varrho^{2\alpha + 1})$. From \eqref{E:DERIVATIVELOSINGEQUATIONS} and straightforward calculations, the following kinds of evolution equations are derived in \cite{shkoller2024geometry}:
	\begin{align}
		\Lunit^{ALE} W_k = \cdots,  && \Lunit^{ALE} Z_k = \cdots,  && \Lunit^{ALE} A_k = \cdots. \label{E:SHKOLLERDIFFERENTIATEDRIEMANNVAIRABLESEQUATIONS}
	\end{align}
The choice of $W_k$, $Z_k$, and $A_k$ as 
$L^2$-differentiated solution variables is convenient in \cite{shkoller2024geometry} because there, the initial data is a perturbation of a fully non-degenerate (strictly convex) 
shock-forming profile for the right Riemann invariant $n_a v^a + \sigma\restriction_{t=0}$. Formally, readers could draw the following analogy with our definition of almost Riemann invariants 
(see \eqref{E:ALMOSTRIEMANNINVARIANTS}) and the differentiated Riemann variables of \cite{shkoller2024geometry}:
\begin{align} \label{E:ANALOGYBETWEENRIEMANNINVARIANTS}
	 \pmb{\p} \RRiemann \sim W_k , & & \pmb{\p} \LRiemann \sim Z_k, & & \pmb{\partial}v^2 \sim A_k.
\end{align}
The schematic identities of \eqref{E:CHRISTODOULOUENERGYAPPENDIX} and \eqref{E:EULERENERGYAPPENDIX} (see the rigorous statement in a characteristic diamond in Prop.\,\ref{P:ENERGYNULLFLUXINTEGRALIDENTITIES}, as well as Sects.\,\ref{SS:RESULTSFROMSINGULARBOUNDARYPAPER}--\ref{SS:NEWIDEASANDMETHODS}), show that the fundamental energy identity for geometric wave equations of the form $\Box_{\gfour} F = G$ automatically provide energy estimates for \emph{all of} the $L^2$-differentiated solution variables $\Lunit F$, $\muX F$, and $Y F$. 
In contrast, being a first order system along $\Lunit = \Lunit^{ALE}$, the $L^2$-estimates for the transport equations \eqref{E:SHKOLLERDIFFERENTIATEDRIEMANNVAIRABLESEQUATIONS} only provide estimates for $W_k, Z_k$, and $A_k$, and do not directly reveal any further tensorial structure; see \eqref{E:SPACETIMEENERGYNULLFLUXIDENTITYFORBELOWTOPORDERACOUSTIC}. 
Hence, to properly identify and study the correct tensorial directions of the solution's gradient that blow-up, and to account for the singular powers of $\upmu = J_g$, 
in \cite{shkoller2024geometry}, the evolution equations are manually contracted with $n_k$ and $\tau^k$ and multiplied by $\upmu$, which yields the following system:
\begin{subequations}
	\begin{align}
		\mathbf{W}_N & :=  n_a W_a, & \mathbf{Z}_N & := n_a Z_a, & \mathbf{A}_N & := n_a A_a,  \label{E:TRANSVERSALDIFFERENTIATEDSHKOLLERVARS} \\
		\mathbf{W}_T & := \tau^a W_a, & \mathbf{Z}_T & := \tau^a Z_a, &  \mathbf{A}_T & := \tau^a A_a,  \label{E:TANGENTIALDIFFERENTIATEDSHKOLLERVARS} \\
		\Lunit^{ALE} (J_g \mathbf{W}_N) & = \cdots, &  \Lunit^{ALE} (J_g\mathbf{Z}_N) & =\cdots, & \Lunit^{ALE}(J_g \mathbf{A}_N) & = \cdots, \label{E:TRANSVERSALDIFFERENTIATEDSHKOLLERVARSEQUATIONS} \\
		\Lunit^{ALE}  \mathbf{W}_T & = \cdots, &  \Lunit^{ALE} \mathbf{Z}_T & =\cdots, & \Lunit^{ALE} \mathbf{A}_T & = \cdots.  \label{E:TANGENTIALDIFFERENTIATEDSHKOLLERVARSEQUATIONS}
	\end{align}
\end{subequations}
	
Since we have shown that $\Lunit = \Lunit^{ALE}$, the standard transport equation energy identity for solutions to \eqref{E:TRANSVERSALDIFFERENTIATEDSHKOLLERVARSEQUATIONS}--\eqref{E:TANGENTIALDIFFERENTIATEDSHKOLLERVARSEQUATIONS} takes the following form:\footnote{The power of $J_g = \upmu$ comes from the following identity for the volume element on $\Sigma_t$: 
$\mathrm{d} x^1 \mathrm{d} x^2 = J_g \mathrm{d}x^1_{ALE}  \mathrm{d} x^2$.} 
	\begin{subequations}
		\begin{align}
			\int_{\Sigma_t}  J_g |(J_g \mathbf{W}_N, J_g\mathbf{Z}_N, J_g \mathbf{A}_N)|^2  & = \textnormal{data} + \cdots, \label{E:ENERGYIDENTITYTRANSVERSALSTEVE} \\
			\int_{\Sigma_t}  J_g | (\mathbf{W}_T , \mathbf{Z}_T ,  \mathbf{A}_T)|^2  & = \textnormal{data} + \cdots, \label{E:ENERGYIDENTITYTANGENTIALSTEVE} 
		\end{align}
	\end{subequations}
In light of \eqref{E:ENERGYIDENTITYTRANSVERSALSTEVE}--\eqref{E:ENERGYIDENTITYTRANSVERSALSTEVE}, we see that the fundamental $L^2$-differentiated solution variables used in \cite{shkoller2024geometry} are precisely $\{J_g \mathbf{W}_N, J_g\mathbf{Z}_N, J_g \mathbf{A}_N,  \mathbf{W}_T , \mathbf{Z}_T ,  \mathbf{A}_T$\}. One could draw the following analogies between the variables used in 
\cite{abbrescia2022emergence} and the present paper:
	\begin{subequations}
		\begin{align}
			J_g \mathbf{W}_N  &\sim \muX \RRiemann, & J_g \mathbf{Z}_N & \sim \muX \LRiemann, & J_g \mathbf{A}_N & \sim \muX v^2,   
			\\
			\mathbf{W}_T & \sim Y \RRiemann, &  \mathbf{Z}_T & \sim Y \LRiemann,  &  \mathbf{A}_T &\sim Y v^2.
		\end{align}
	\end{subequations}

\subsection{Multiplying by $\upmu=J_g$ first and then applying smooth differential operators and regular commutators} \label{SS:MULTIPLYINGFIRSTBYMU}
For equations whose solutions form shocks, including \eqref{E:CHRISTODOULOUGEOMETRICEQUATION} from \cite{dC2007}, \eqref{E:COVARIANTWAVEEQUATIONSWAVEVARIABLES}--\eqref{E:TRANSPORTFLATDIVGRADENT} from \cite{abbrescia2022emergence} or this paper, or \eqref{E:TRANSVERSALDIFFERENTIATEDSHKOLLERVARSEQUATIONS}--\eqref{E:TANGENTIALDIFFERENTIATEDSHKOLLERVARSEQUATIONS} from \cite{shkoller2024geometry}, commuting with generic derivatives leads to uncontrollable error terms featuring powers of $\upmu^{-1}$. For example, given $P  \in\{\Lunit,\Yvf{2},\Yvf{3}\}$, \eqref{E:SCHEMATICGRADIENTBLOWUP} implies $[P,\partial_{x^\alpha}]\sim \frac{1}{\upmu}\partial_{x^\alpha}$. Hence, commuting $P$ with $\Box_{\gfour} \sim \gfour^{\alpha\beta} \partial_{\alpha}\partial_{\beta}$ would lead to powers of $\frac{1}{\upmu}$. For these reasons, \emph{all of} \cite{dC2007,abbrescia2022emergence,shkoller2024geometry} and this paper multiply the evolution equations \emph{first} by $\upmu = J_g$ before commuting. For example, after multiplying by $\upmu$, we can express $\upmu \Box_{\gfour} \sim Z^{\le 2}$ where $Z \in \{\Lunit,\muX,\Yvf{2},\Yvf{3}\}$. In Lemma\,\ref{L:SIMPLECOMMUTATORIDENTITY}, we prove that $[P,Z]\sim \Yvf{A}$, 
which is much better than $[P,\p_{x^\alpha}] \sim \frac{1}{\upmu} \p_{x^\alpha}$. 

Similarly, RHSs~\eqref{E:TRANSVERSALDIFFERENTIATEDSHKOLLERVARSEQUATIONS}--\eqref{E:TANGENTIALDIFFERENTIATEDSHKOLLERVARSEQUATIONS} feature error terms with coefficients given by $J_g^{-1}$; see \cite{shkoller2024geometry}*{Equations\,(1.20)--(1.21)}. Differentiating these equations would lead to worse powers of $J_g^{-1}$, but multiplying \emph{first} by $\upmu = J_g$ regularizes the equations, see the additional power of $J_g = \upmu$ on LHS of \cite{shkoller2024geometry}*{Equations\,(3.35b)--(3.35c)}. We also note that all of the commutators from \cite{dC2007,abbrescia2022emergence,shkoller2024geometry} and this paper are differential operators which are \emph{smooth} relative to the geometric coordinates (equivalently, to the ALE coordinates); see Sect.\,\ref{SSS:ENERGYDESCENTSCHEMEINROUGHCOORDINATES}.

\subsubsection{Coercive bulk energies and damping norms} \label{SS:COERCIVEBULKANDAMPINGTERMS}
For the sake of simplicity and clarity, consider the problem of only controlling the solution up until the first Cartesian time of gradient blowup. In this case, an appropriate weight in the energy method is $\weight = \upmu = J_g$. Now, there are many error integrals that arise upon commuting the energy identities 
\eqref{E:CHRISTODOULOUENERGYAPPENDIX}--\eqref{E:EULERENERGYAPPENDIX}, \eqref{E:ENERGYIDENTITYTRANSVERSALSTEVE}--\eqref{E:ENERGYIDENTITYTANGENTIALSTEVE} (where the commutations must follow the order of operations developed by Christodoulou \cite{dC2007} described in Sect.\,\ref{SS:MULTIPLYINGFIRSTBYMU}). 
Here, we depict one such error integral, first in the style of the present paper, and then in the style of \cite{shkoller2024geometry}.
One such error integral is depicted in RHSs~\eqref{E:SCHEMATICWAVECOMMUTATIONS}--\eqref{E:SCHEMATICFIRSTORDERECOMMUTATIONS} below, see \cite{shkoller2024geometry}*{Equations (1.25)--(1.26)}:\footnote{The papers \cite{dC2007,abbrescia2022emergence,shkoller2024geometry} relied on different values of $M_*$ to close the estimates.}
	\begin{subequations}
		\begin{align}
			\begin{split}
				 & \int_{\Sigma_t^{[u_1,u_2]}} \upmu^{M_*} \left\{ (\upmu \newuL \tander^{\Ntop} \wavearray)^2 + \upmu (Y  \tander^{\Ntop} \wavearray)^2 + \upmu (\Lunit  \tander^{\Ntop} \wavearray)^2\right\}  + \int_{\nullhyparg{u}^{[0,t]}} \upmu^{M_*} \left\{ (\Lunit \tander^{\Ntop} \wavearray)^2 + \upmu (Y\tander^{\Ntop}  \wavearray )^2 \right\} \\
				 & \ \ + (M_* + 1)\int_{M_{t,u}} (-\Lunit \upmu) \upmu^{M_*} (Y  \tander^{\Ntop}  \wavearray )^2   \lesssim \mathcal{O}(\varepsilon)  \int_{M_{t,u}}  \upmu^{M_*} (Y  \tander^N  \wavearray )^2,
			\end{split}  \label{E:SCHEMATICWAVECOMMUTATIONS} \\
			\begin{split} \label{E:SCHEMATICFIRSTORDERECOMMUTATIONS}
				& \int_{\Sigma_t}  J_g^{M_* + 1} |\tander^{\Ntop}(J_g \mathbf{W}_N, J_g\mathbf{Z}_N, J_g \mathbf{A}_N)|^2  + (M_* +1) \int_{M_t} \left(- \geop{t} J_g\right)  J_g^{M_*} |\tander^{\Ntop}(J_g \mathbf{W}_N, J_g\mathbf{Z}_N, J_g \mathbf{A}_N)|^2\\
				& \ \ 
				\lesssim 	\mathcal{O}(\varepsilon) \int_{M_t}  J_g^{M_*} |\tander^{\Ntop}(J_g \mathbf{W}_N, J_g\mathbf{Z}_N, J_g \mathbf{A}_N)|^2. 
			\end{split}
		\end{align}
	\end{subequations}
Observe that the coercive control afforded by the $L^2$-energies and null-fluxes on $\Sigma_t^{[u_1,u_]}, \nullhyparg{u}^{[0,t]}$ or $\Sigma_t$ on the LHS are insufficient with respect to the singular powers of $\upmu$ yo control the integrals on the RHS via a standard method based on Gronwall's inequality. 
However, since $\geop{t} J_g = \geop{t} \upmu \approx \Lunit \upmu \approx -1$ (see \eqref{E:NEWLMUISALMOSTMINUSONEINSMALLNEIGHBORHOOD}, \eqref{E:SPACETIMECOERCIVEINTEGRALS}, and \cite{shkoller2024geometry}*{Equation (1.23a)}), the spacetime integrals on LHS \eqref{E:SCHEMATICWAVECOMMUTATIONS}--\eqref{E:SCHEMATICFIRSTORDERECOMMUTATIONS} provide  stronger coercive control  by an entire power of $\upmu = J_g$. We call these \emph{bulk integrals} in \cite{abbrescia2022emergence} and this paper, while \cite{shkoller2024geometry}*{Sect.\,1.10.1} calls them \emph{damping norms}. Observe that the underlying mechanism for the coercivity is the same in both spacetime integrals, namely the compressive sign of $\Lunit \upmu < 0$. If $\epsilon$ is sufficiently small relative to $M_* + 1$, it follows that the error integrals on 
the RHS can be absorbed by the bulk integrals on the LHS.


\subsubsection{Avoiding derivative-loss in the acoustic geometry} \label{SSS:AVOIDINGDERIVATIVELOSSINTHEGEOMETRY}
Since the coordinate systems used in \cite{dC2007,abbrescia2022emergence} and \cite{shkoller2024geometry} are the same
(see Prop.\,\ref{P:ALETOGEOCOORD}), 
all these papers must grapple with a potential loss of derivatives in certain error integrals in the energy estimates
that feature top-order commutator terms involving the acoustic geometry. 
In the former papers, as we explained in Sect.\,\ref{SSS:GAUSSCURVATUREESTIMATESNOTREQUIREDIN2D}, 
the difficult terms involve the top-order derivatives of the tensor $\upchi$. The potential loss of regularity was
handled by using Raychaudhuri's equation \cite{aR1955} to derive a transport equation for 
a modified version of $\mytr_{\gtorus}\upchi$, which allows one to bound it terms of error terms that have a consistent level of regularity, 
and also employing elliptic estimates to handle the non-trace part of $\upchi$. We point out that 
the elliptic estimates are not needed in $2D$ because $\upchi$ is determined by its trace in $2D$.
In the $2D$ isentropic setting, the authors of \cite{shkoller2024geometry} used a different approach. They 
found a novel perfect divergence structure in the derivative-losing acoustic geometry terms, 
and in the energy estimates, the unwanted derivative can integrated away; 
see \cite{shkoller2024geometry}*{Sect\,10}. Identifying this structure was an important contribution of \cite{shkoller2024geometry}.
It remains an interesting question whether a similar divergence structure is present in the general $3D$ case with vorticity and entropy.

\bibliographystyle{alpha} 

\bibliography{JBib}

\end{document}